\documentclass[numbers=enddot,12pt,final,onecolumn,notitlepage]{scrartcl}%
\usepackage[headsepline,footsepline,manualmark]{scrlayer-scrpage}
\usepackage[all,cmtip]{xy}
\usepackage{amssymb}
\usepackage{amsmath}
\usepackage{amsthm}
\usepackage{framed}
\usepackage{comment}
\usepackage{color}
\usepackage[breaklinks=True]{hyperref}
\usepackage[sc]{mathpazo}
\usepackage[T1]{fontenc}
\usepackage{needspace}
\usepackage{tabls}
\usepackage{ytableau}
\usepackage{tikz}
\usepackage{pgfplots}
\usepackage[type={CC}, modifier={zero}, version={1.0},]{doclicense}
\providecommand{\U}[1]{\protect\rule{.1in}{.1in}}
\pgfplotsset{compat=1.13}
\theoremstyle{definition}
\newtheorem{theo}{Theorem}[subsection]
\newtheorem{exer}{Exercise}[subsubsection]
\newenvironment{theorem}[1][]
{\begin{theo}[#1]\begin{leftbar}}
{\end{leftbar}\end{theo}}
\newtheorem{lem}[theo]{Lemma}
\newenvironment{lemma}[1][]
{\begin{lem}[#1]\begin{leftbar}}
{\end{leftbar}\end{lem}}
\newtheorem{prop}[theo]{Proposition}
\newenvironment{proposition}[1][]
{\begin{prop}[#1]\begin{leftbar}}
{\end{leftbar}\end{prop}}
\newtheorem{defi}[theo]{Definition}
\newenvironment{definition}[1][]
{\begin{defi}[#1]\begin{leftbar}}
{\end{leftbar}\end{defi}}
\newtheorem{remk}[theo]{Remark}
\newenvironment{remark}[1][]
{\begin{remk}[#1]\begin{leftbar}}
{\end{leftbar}\end{remk}}
\newtheorem{coro}[theo]{Corollary}
\newenvironment{corollary}[1][]
{\begin{coro}[#1]\begin{leftbar}}
{\end{leftbar}\end{coro}}
\newtheorem{conv}[theo]{Convention}
\newenvironment{convention}[1][]
{\begin{conv}[#1]\begin{leftbar}}
{\end{leftbar}\end{conv}}
\newtheorem{quest}[theo]{Question}
\newenvironment{question}[1][]
{\begin{quest}[#1]\begin{leftbar}}
{\end{leftbar}\end{quest}}
\newtheorem{warn}[theo]{Warning}

\newtheorem{conj}[theo]{Conjecture}
\newenvironment{conjecture}[1][]
{\begin{conj}[#1]\begin{leftbar}}
{\end{leftbar}\end{conj}}
\newtheorem{exam}[theo]{Example}
\newenvironment{example}[1][]
{\begin{exam}[#1]\begin{leftbar}}
{\end{leftbar}\end{exam}}
\newtheorem{exmp}[exer]{Exercise}
\newenvironment{exercise}[1][]
{\begin{exmp}[#1]\begin{leftbar}}
{\end{leftbar}\end{exmp}}
\newenvironment{statement}{\begin{quote}}{\end{quote}}
\newenvironment{fineprint}{\begin{small}}{\end{small}}
\let\sumnonlimits\sum
\let\prodnonlimits\prod
\let\cupnonlimits\bigcup
\let\capnonlimits\bigcap
\renewcommand{\sum}{\sumnonlimits\limits}
\renewcommand{\prod}{\prodnonlimits\limits}
\renewcommand{\bigcup}{\cupnonlimits\limits}
\renewcommand{\bigcap}{\capnonlimits\limits}
\usetikzlibrary{arrows,arrows.meta,decorations.markings}
\newenvironment{noncompile}{}{}
\excludecomment{verlong}
\includecomment{vershort}
\excludecomment{noncompile}

\newcommand{\arinj}{\ar@{_{(}->}}
\newcommand{\arinjrev}{\ar@{^{(}->}}
\newcommand{\arsurj}{\ar@{->>}}
\newcommand{\arelem}{\ar@{|->}}
\newcommand{\arback}{\ar@{<-}}

\definecolor{dbluecolor}{rgb}{0.01,0.02,0.7}
\definecolor{dgreencolor}{rgb}{0.2,0.4,0.0}
\definecolor{darkred}{rgb}{0.7,0,0}
\newtheoremstyle{plainsl}
{8pt plus 2pt minus 4pt}
{8pt plus 2pt minus 4pt}
{\slshape}
{0pt}
{\bfseries}
{.}
{5pt plus 1pt minus 1pt}
{}
\theoremstyle{plainsl}
\ihead{Introduction to Algebraic Combinatorics, version 2026-11-21}
\ohead{page \thepage}
\begin{document}

\title{An Introduction to Algebraic Combinatorics\\{\normalsize [Math 701, Spring 2021 lecture notes]}\\{\normalsize [Math 531, Winter 2024 lecture notes]}}
\author{Darij Grinberg}
\date{November 21, 2026}
\maketitle

\begin{abstract}
\textbf{Abstract.} This is an introduction to algebraic combinatorics, written
for a quarter-long graduate course. It starts with a rigorous introduction to
formal power series with some combinatorial applications, then discusses
integer partitions (proving Jacobi's triple product identity), permutations
(Lehmer codes, cycles) and subtractive methods (alternating sums,
cancellations and inclusion-exclusion principles, with a particular focus on
sign-reversing involutions and determinants). The last chapter introduces
symmetric polynomials and proves the Littlewood--Richardson rule using
Bender--Knuth involutions (a la Stembridge).

The appendix contains over 200 exercises (without solutions).

\end{abstract}
\tableofcontents

\doclicenseThis

\newpage

\section{What is this?}

These are the notes for an introductory course on algebraic combinatorics held
in the Spring Quarter 2021 at Drexel University\footnote{The website of this
course is \url{https://darijgrinberg.gitlab.io/t/21s/}}. The topics covered are

\begin{itemize}
\item formal power series and their use as generating functions (Chapter
\ref{chap.gf});

\item integer partitions and $q$-binomial coefficients (Chapter
\ref{chap.pars});

\item permutations and their lengths, inversions and cycles (Chapter
\ref{chap.perm});

\item alternating sums, the use of sign-reversing involutions and the
combinatorial view on determinants (Chapter \ref{chap.sign});

\item the basics of symmetric polynomials, particularly Schur polynomials
(Chapter \ref{chap.sf}).
\end{itemize}

Almost all of this text is in a finished state (the only exception being
Section \ref{sec.gf.multivar}, which is currently an overview). However, it is
improvable both in detail and in coverage. They might grow in later iterations
of this course (in particular, various further topics could get included).
Errors and confusions will be fixed whenever I become aware of them (any
assistance is greatly appreciated!\footnote{Please send comments to
\href{mailto:darijgrinberg@gmail.com}{\texttt{darijgrinberg@gmail.com}}}).

Exercises of varying difficulty appear at the end of the text (Chapter
\ref{chap.hw}).

\subsection*{Acknowledgments}

Thanks to the students in my Math 701 course for what was in essence an
alpha-test of these notes. Some exercises have been adapted from collections
by Richard P. Stanley, Martin Aigner, Donald Knuth, Mikl\'{o}s B\'{o}na, Mark
Wildon and Igor Pak. A
\href{https://math.stackexchange.com/users/556825/mindlack}{math.stackexchange
user named Mindlack} has contributed the proof of Proposition
\ref{prop.fps.prodrule-inf-inf}. Andrew Solomon and Mikey Becht have reported
typos. Your name could most likely stand here.

\subsection{AI use}

Various versions of GPT have found errors and suggested a few exercises as
well as a correct proof of Proposition \ref{prop.fps.prods-mulable-rules.SW1}
(a technical point whose subtlety I had initially underestimated). All the
writing is my own.

\newpage

\section{\label{chap.intro}Before we start...}

\subsection{What is this?}

This is a course on algebraic combinatorics. This subject can be viewed either
as a continuation of enumerative combinatorics by other means (specifically,
algebraic ones), or as the part of algebra where one studies concrete
polynomials (more precisely, families of polynomials). For example, the Schur
polynomials can be viewed on the one hand as a tool for enumerating certain
kinds of tableaux (essentially, tabular arrangements of numbers that increase
along rows and columns), while on the other hand they form a family of
polynomials with a myriad surprising properties, generalizing (e.g.) the
Vandermonde determinant. I hope to cover both aspects of the subject to a
reasonable extent in this course.

\subsection{Prerequisites}

To understand this course, you are assumed to speak the language of rings and
fields (we will mostly need the basic properties of polynomials and linear
maps; we will define what we need about power series), and to have some basic
knowledge of enumerative combinatorics (see below). My notes \cite{23wa}, and
the references I gave therein, can help refresh your knowledge of the former.
As for the latter, there are dozens of sources available (I made a list at
\url{https://math.stackexchange.com/a/1454420/} , focussing mostly on texts
available online), including my own notes \cite{22fco}.

\subsection{Notations and elementary facts}

We will use the following notations and conventions:

\begin{itemize}
\item The symbol $\mathbb{N}$ will denote the set $\left\{  0,1,2,3,\ldots
\right\}  $ of nonnegative integers.

\item The size (i.e., cardinality) of a set $A$ will be denoted by $\left\vert
A\right\vert $.

\item The symbol \textquotedblleft\#\textquotedblright\ means
\textquotedblleft number\textquotedblright. For example, the size $\left\vert
A\right\vert $ of a set $A$ is the \# of elements of $A$.

\item As usual in combinatorics, $a^{0}=1$ for any $a$, including $0$.

\item We often (but not always) use the notation $\left[  n\right]  $ for the
set $\left\{  1,2,\ldots,n\right\}  $, where $n$ is an arbitrary integer. (If
$n\leq0$, then this set is empty.)
\end{itemize}

We will need some basics from enumerative combinatorics (see, e.g.,
\cite[\S 8.1]{Newste19} for details, and \cite[Chapters 1 and 2]{19fco} for
more details):

\begin{itemize}
\item \textbf{addition principle = sum rule:} If $A$ and $B$ are two disjoint
sets, then $\left\vert A\cup B\right\vert =\left\vert A\right\vert +\left\vert
B\right\vert $.

\item \textbf{multiplication principle = product rule:} If $A$ and $B$ are any
two sets, then $\left\vert A\times B\right\vert =\left\vert A\right\vert
\cdot\left\vert B\right\vert $.

\item \textbf{bijection principle:} There is a bijection (= bijective map =
invertible map = one-to-one correspondence) between two sets $X$ and $Y$ if
and only if $\left\vert X\right\vert =\left\vert Y\right\vert $.

\item A set with $n$ elements has $2^{n}$ subsets, and has $\dbinom{n}{k}$
size-$k$ subsets for any $k\in\mathbb{R}$.

\item A set with $n$ elements has $n!$ permutations (= bijective maps from
this set to itself).

\item \textbf{dependent product rule:} Consider a situation in which you have
to make $n$ choices (sequentially). Assume that you have $a_{1}$ options
available in choice $1$, and then (after making choice $1$) you have $a_{2}$
options available in choice $2$ (no matter which option you chose in choice
$1$), and then (after both choices $1$ and $2$) you have $a_{3}$ options
available in choice $3$ (no matter which options you chose in choices $1$ and
$2$), and so on. Then, the total \# of ways to make all $n$ choices is
$a_{1}a_{2}\cdots a_{n}$. (This is formalized in \cite[Theorem 8.1.19]%
{Newste19}.)
\end{itemize}

A few words about binomial coefficients are in order:

\begin{definition}
\label{def.binom.binom}For any numbers $n$ and $k$, we set%
\begin{equation}
\dbinom{n}{k}=%
\begin{cases}
\dfrac{n\left(  n-1\right)  \left(  n-2\right)  \cdots\left(  n-k+1\right)
}{k!}, & \text{if }k\in\mathbb{N};\\
0, & \text{else.}%
\end{cases}
\label{eq.def.binom.binom.eq}%
\end{equation}
Note that \textquotedblleft numbers\textquotedblright\ is to be understood
fairly liberally here. In particular, $n$ can be any integer, rational, real
or complex number (or, more generally, any element in a $\mathbb{Q}$-algebra),
whereas $k$ can be anything (although the only nonzero values of $\dbinom
{n}{k}$ will be achieved for $k\in\mathbb{N}$, by the above definition).
\end{definition}

\begin{example}
\label{exa.binom.-1choosek}For any $k\in\mathbb{N}$, we have
\begin{align*}
\dbinom{-1}{k}  &  =\dfrac{\left(  -1\right)  \left(  -1-1\right)  \left(
-1-2\right)  \cdots\left(  -1-k+1\right)  }{k!}\\
&  =\dfrac{\left(  -1\right)  \left(  -2\right)  \left(  -3\right)
\cdots\left(  -k\right)  }{k!}=\dfrac{\left(  -1\right)  ^{k}k!}{k!}=\left(
-1\right)  ^{k}.
\end{align*}

\end{example}

If $n,k\in\mathbb{N}$ and $n\geq k$, then
\begin{equation}
\dbinom{n}{k}=\dfrac{n!}{k!\left(  n-k\right)  !}. \label{eq.binom.fac-form}%
\end{equation}
But this formula only applies to the case when $n,k\in\mathbb{N}$ and $n\geq
k$. Our above definition is more general than it. The combinatorial meaning of
the binomial coefficient $\dbinom{n}{k}$ (as the \# of $k$-element subsets of
a given $n$-element set) also cannot be used for negative or non-integer
values of $n$.

\begin{example}
\label{exa.binom.2n-choose-n}Let $n\in\mathbb{N}$. Then, $\dbinom{2n}%
{n}=\dfrac{1\cdot3\cdot5\cdot\cdots\cdot\left(  2n-1\right)  }{n!}\cdot2^{n}$.
\end{example}

\begin{proof}
[Proof of Example \ref{exa.binom.2n-choose-n}.]We have%
\begin{align*}
\left(  2n\right)  !  &  =1\cdot2\cdot\cdots\cdot\left(  2n\right) \\
&  =\left(  1\cdot3\cdot5\cdot\cdots\cdot\left(  2n-1\right)  \right)
\cdot\underbrace{\left(  2\cdot4\cdot6\cdot\cdots\cdot\left(  2n\right)
\right)  }_{\substack{=\left(  2\cdot1\right)  \cdot\left(  2\cdot2\right)
\cdot\cdots\cdot\left(  2\cdot n\right)  \\=2^{n}\left(  1\cdot2\cdot
\cdots\cdot n\right)  }}\\
&  =\left(  1\cdot3\cdot5\cdot\cdots\cdot\left(  2n-1\right)  \right)
\cdot2^{n}\underbrace{\left(  1\cdot2\cdot\cdots\cdot n\right)  }_{=n!}\\
&  =\left(  1\cdot3\cdot5\cdot\cdots\cdot\left(  2n-1\right)  \right)
\cdot2^{n}n!.
\end{align*}
Now, (\ref{eq.binom.fac-form}) yields%
\begin{align*}
\dbinom{2n}{n}  &  =\dfrac{\left(  2n\right)  !}{n!\left(  2n-n\right)
!}=\dfrac{\left(  2n\right)  !}{n!n!}=\dfrac{\left(  1\cdot3\cdot5\cdot
\cdots\cdot\left(  2n-1\right)  \right)  \cdot2^{n}n!}{n!\cdot n!}\\
&  \ \ \ \ \ \ \ \ \ \ \ \ \ \ \ \ \ \ \ \ \left(  \text{since }\left(
2n\right)  !=\left(  1\cdot3\cdot5\cdot\cdots\cdot\left(  2n-1\right)
\right)  \cdot2^{n}n!\right) \\
&  =\dfrac{\left(  1\cdot3\cdot5\cdot\cdots\cdot\left(  2n-1\right)  \right)
\cdot2^{n}}{n!}=\dfrac{1\cdot3\cdot5\cdot\cdots\cdot\left(  2n-1\right)  }%
{n!}\cdot2^{n}.
\end{align*}
This proves Example \ref{exa.binom.2n-choose-n}.
\end{proof}

Entire books have been written about binomial coefficients and their
properties. See \cite{Spivey19} for a recent text; elementary introductions
can be found in \cite[Chapter 5]{GKP} and \cite[Chapter 3]{detnotes} and
\cite[\S 1.2.6]{Knuth-TAoCP1} and \cite[Chapter 2]{Wildon19}. Here are two
more basic facts that we will need (\cite[Theorem 1.3.8]{19fco} and
\cite[Proposition 1.3.6]{19fco}, respectively):

\begin{proposition}
[\emph{Pascal's identity}, aka \emph{recurrence of the binomial coefficients}%
]\label{prop.binom.rec}We have%
\begin{equation}
\dbinom{m}{n}=\dbinom{m-1}{n-1}+\dbinom{m-1}{n} \label{eq.binom.rec.m}%
\end{equation}
for any numbers $m$ and $n$.
\end{proposition}

\begin{proposition}
\label{prop.binom.0}Let $m,n\in\mathbb{N}$ satisfy $m<n$. Then, $\dbinom{m}%
{n}=0$.
\end{proposition}

Note that Proposition \ref{prop.binom.0} really requires $m\in\mathbb{N}$. For
example, $1.5<2$ but $\dbinom{1.5}{2}=0.375\neq0$.

Yet another useful property of the binomial coefficients is the following
(\cite[Theorem 1.3.11]{19fco}):

\begin{theorem}
[\emph{Symmetry of the binomial coefficients}]\label{thm.binom.sym}Let
$n\in\mathbb{N}$ and $k\in\mathbb{R}$. Then,%
\[
\dbinom{n}{k}=\dbinom{n}{n-k}.
\]

\end{theorem}

Note the $n\in\mathbb{N}$ requirement. Convince yourself that Theorem
\ref{thm.binom.sym} would fail for $n=-1$ and $k=0$.

\newpage

\section{\label{chap.gf}Generating functions}

In this first chapter, we will discuss generating functions: first informally,
then on a rigorous footing. You may have seen generating functions already, as
their usefulness extends far beyond combinatorics; but they are so important
to this course that they are worth covering twice in case of doubt.

Rigorous introductions to generating functions (and formal power series in
general) can also be found in \cite[Chapter 7 (in the 1st edition)]{Loehr-BC},
in \cite[Chapter 1]{Henric74}, in \cite{Sambal22}, and (to some extent) in
\cite[Chapter 7]{19s}.\footnote{Bourbaki's \cite[\S IV.4]{Bourba03} contains
what might be the most rigorous and honest treatment of formal power series
available in the literature; however, it is not the most readable source, as
the notation is dense and heavily relies on other volumes by the same author.}
A quick overview is given in \cite{Niven69}, and many applications are found
in \cite{Wilf09} and \cite[Chapters 4--5]{Riorda79}. There are furthermore
numerous books that explore enumerative combinatorics through the lens of
generating functions (\cite{GouJac83}, \cite{Wagner08}, \cite{Lando03} and others).

\subsection{\label{sec.gf.exas}Examples}

Let me first show what generating functions are good for. Then, starting in
the next section, I will explain how to rigorously define them. For now, I
will work informally; please suspend your disbelief until the next section.

The \textbf{idea} behind generating functions is easy: Any sequence $\left(
a_{0},a_{1},a_{2},\ldots\right)  $ of numbers gives rise to a
\textquotedblleft power series\textquotedblright\ $a_{0}+a_{1}x+a_{2}%
x^{2}+\cdots$, which is called the \emph{generating function} of this
sequence. This \textquotedblleft power series\textquotedblright\ is an
infinite sum (an \textquotedblleft infinite polynomial\textquotedblright\ in
an indeterminate $x$), so it is not immediately clear what it means and what
we are allowed to do with it; but before we answer such questions, let us
first play around with these power series and hope for the best. The following
four examples show how they can be useful.

\subsubsection{Example 1: The Fibonacci sequence}

\textbf{Example 1.} The \emph{Fibonacci sequence} is the sequence $\left(
f_{0},f_{1},f_{2},\ldots\right)  $ of integers defined recursively by%
\[
f_{0}=0,\ \ \ \ \ \ \ \ \ \ f_{1}=1,\ \ \ \ \ \ \ \ \ \ f_{n}=f_{n-1}%
+f_{n-2}\text{ for each }n\geq2.
\]
Its entries are known as the \emph{Fibonacci numbers}. Here are the first few
of them:%
\[%
\begin{tabular}
[c]{|c||c|c|c|c|c|c|c|c|c|c|c|c|}\hline
$n$ & $0$ & $1$ & $2$ & $3$ & $4$ & $5$ & $6$ & $7$ & $8$ & $9$ & $10$ &
$11$\\\hline
$f_{n}$ & $0$ & $1$ & $1$ & $2$ & $3$ & $5$ & $8$ & $13$ & $21$ & $34$ & $55$
& $89$\\\hline
\end{tabular}
\ \ \
\]

Let us see what we can learn about this sequence by considering its generating
function%
\begin{align*}
F\left(  x\right)   &  :=f_{0}+f_{1}x+f_{2}x^{2}+f_{3}x^{3}+\cdots\\
&  =0+1x+1x^{2}+2x^{3}+3x^{4}+5x^{5}+\cdots.
\end{align*}

We have%
\begin{align*}
F\left(  x\right)   &  =f_{0}+f_{1}x+f_{2}x^{2}+f_{3}x^{3}+f_{4}x^{4}+\cdots\\
&  =\underbrace{0+1x}_{=x}+\left(  f_{1}+f_{0}\right)  x^{2}+\left(
f_{2}+f_{1}\right)  x^{3}+\left(  f_{3}+f_{2}\right)  x^{4}+\cdots\\
&  \ \ \ \ \ \ \ \ \ \ \ \ \ \ \ \ \ \ \ \ \left(  \text{since }f_{0}=0\text{
and }f_{1}=1\text{ and }f_{n}=f_{n-1}+f_{n-2}\text{ for each }n\geq2\right) \\
&  =x+\underbrace{\left(  f_{1}+f_{0}\right)  x^{2}+\left(  f_{2}%
+f_{1}\right)  x^{3}+\left(  f_{3}+f_{2}\right)  x^{4}+\cdots}%
_{\substack{=\left(  f_{1}x^{2}+f_{2}x^{3}+f_{3}x^{4}+\cdots\right)  +\left(
f_{0}x^{2}+f_{1}x^{3}+f_{2}x^{4}+\cdots\right)  \\\text{(here we are hoping
that this manipulation of}\\\text{infinite sums is indeed legitimate)}}}\\
&  =x+\underbrace{\left(  f_{1}x^{2}+f_{2}x^{3}+f_{3}x^{4}+\cdots\right)
}_{\substack{=x\left(  f_{1}x+f_{2}x^{2}+f_{3}x^{3}+\cdots\right)  \\=x\left(
F\left(  x\right)  -f_{0}\right)  =xF\left(  x\right)  \\\text{(since }%
f_{0}=0\text{)}}}+\underbrace{\left(  f_{0}x^{2}+f_{1}x^{3}+f_{2}x^{4}%
+\cdots\right)  }_{\substack{=x^{2}\left(  f_{0}+f_{1}x+f_{2}x^{2}+f_{3}%
x^{3}+\cdots\right)  \\=x^{2}F\left(  x\right)  }}\\
&  =x+xF\left(  x\right)  +x^{2}F\left(  x\right)  =x+\left(  x+x^{2}\right)
F\left(  x\right)  .
\end{align*}
Solving this equation for $F\left(  x\right)  $ (assuming that we are allowed
to divide by $1-x-x^{2}$), we get%
\begin{align}
F\left(  x\right)   &  =\dfrac{x}{1-x-x^{2}}\label{eq.sec.gf.exas.1.Fx=1}\\
&  =\dfrac{x}{\left(  1-\phi_{+}x\right)  \left(  1-\phi_{-}x\right)
},\nonumber
\end{align}
where $\phi_{+}=\dfrac{1+\sqrt{5}}{2}$ and $\phi_{-}=\dfrac{1-\sqrt{5}}{2}$
are the two roots of the quadratic polynomial $x^{2}-x-1$ (note that $\phi
_{+}$ and $\phi_{-}$ are sometimes known as the \textquotedblleft golden
ratios\textquotedblright; we have $\phi_{+}\approx1.618$ and $\phi_{-}%
\approx-0.618$). Hence,%
\begin{align}
F\left(  x\right)   &  =\dfrac{x}{\left(  1-\phi_{+}x\right)  \left(
1-\phi_{-}x\right)  }\nonumber\\
&  =\dfrac{1}{\sqrt{5}}\cdot\dfrac{1}{1-\phi_{+}x}-\dfrac{1}{\sqrt{5}}%
\cdot\dfrac{1}{1-\phi_{-}x} \label{eq.sec.gf.exas.1.Fx=2}%
\end{align}
(by partial fraction decomposition).

Now, what are the coefficients of the power series $\dfrac{1}{1-\alpha x}$ for
an $\alpha\in\mathbb{C}$ ? Let me first answer this question for $\alpha=1$.
Namely, I claim that%
\begin{equation}
\dfrac{1}{1-x}=1+x+x^{2}+x^{3}+\cdots. \label{eq.sec.gf.exas.1.1/1-x}%
\end{equation}
Indeed, this follows by observing that
\begin{align*}
&  \left(  1-x\right)  \left(  1+x+x^{2}+x^{3}+\cdots\right) \\
&  =\left(  1+x+x^{2}+x^{3}+\cdots\right)  -x\left(  1+x+x^{2}+x^{3}%
+\cdots\right) \\
&  =\left(  1+x+x^{2}+x^{3}+\cdots\right)  -\left(  x+x^{2}+x^{3}+x^{4}%
+\cdots\right) \\
&  =1
\end{align*}
(again, we are hoping that these manipulations of infinite sums are allowed).
Note that the equality (\ref{eq.sec.gf.exas.1.1/1-x}) is a version of the
\emph{geometric series formula} familiar from real analysis. Now, for any
$\alpha\in\mathbb{C}$, we can substitute $\alpha x$ for $x$ in the equality
(\ref{eq.sec.gf.exas.1.1/1-x}), and thus obtain%
\begin{align}
\dfrac{1}{1-\alpha x}  &  =1+\alpha x+\left(  \alpha x\right)  ^{2}+\left(
\alpha x\right)  ^{3}+\cdots\nonumber\\
&  =1+\alpha x+\alpha^{2}x^{2}+\alpha^{3}x^{3}+\cdots.
\label{eq.sec.gf.exas.1.1/1-ax}%
\end{align}
Hence, our above formula (\ref{eq.sec.gf.exas.1.Fx=2}) becomes%
\begin{align*}
F\left(  x\right)   &  =\dfrac{1}{\sqrt{5}}\cdot\dfrac{1}{1-\phi_{+}x}%
-\dfrac{1}{\sqrt{5}}\cdot\dfrac{1}{1-\phi_{-}x}\\
&  =\dfrac{1}{\sqrt{5}}\cdot\left(  1+\phi_{+}x+\phi_{+}^{2}x^{2}+\phi_{+}%
^{3}x^{3}+\cdots\right)  -\dfrac{1}{\sqrt{5}}\cdot\left(  1+\phi_{-}x+\phi
_{-}^{2}x^{2}+\phi_{-}^{3}x^{3}+\cdots\right) \\
&  \ \ \ \ \ \ \ \ \ \ \ \ \ \ \ \ \ \ \ \ \left(  \text{by
(\ref{eq.sec.gf.exas.1.1/1-ax}), applied to }\alpha=\phi_{+}\text{ and again
to }\alpha=\phi_{-}\right) \\
&  =\dfrac{1}{\sqrt{5}}\sum_{k\geq0}\phi_{+}^{k}x^{k}-\dfrac{1}{\sqrt{5}}%
\sum_{k\geq0}\phi_{-}^{k}x^{k}\\
&  =\sum_{k\geq0}\left(  \dfrac{1}{\sqrt{5}}\cdot\phi_{+}^{k}-\dfrac{1}%
{\sqrt{5}}\cdot\phi_{-}^{k}\right)  x^{k}.
\end{align*}
Now, for any given $n\in\mathbb{N}$, the coefficient of $x^{n}$ in the power
series on the left hand side of this equality is $f_{n}$ (since $F\left(
x\right)  =f_{0}+f_{1}x+f_{2}x^{2}+f_{3}x^{3}+\cdots$), whereas the
coefficient of $x^{n}$ on the right hand side is clearly $\dfrac{1}{\sqrt{5}%
}\cdot\phi_{+}^{n}-\dfrac{1}{\sqrt{5}}\cdot\phi_{-}^{n}$. Thus, comparing
coefficients before $x^{n}$, we obtain%
\begin{align*}
f_{n}  &  =\dfrac{1}{\sqrt{5}}\cdot\phi_{+}^{n}-\dfrac{1}{\sqrt{5}}\cdot
\phi_{-}^{n}\\
&  \ \ \ \ \ \ \ \ \ \ \ \ \ \ \ \ \ \ \ \ \left(
\begin{array}
[c]{c}%
\text{assuming that \textquotedblleft comparing coefficients\textquotedblright%
\ is allowed, i.e.,}\\
\text{that equal power series really have equal coefficients}%
\end{array}
\right) \\
&  =\dfrac{1}{\sqrt{5}}\cdot\left(  \dfrac{1+\sqrt{5}}{2}\right)  ^{n}%
-\dfrac{1}{\sqrt{5}}\cdot\left(  \dfrac{1-\sqrt{5}}{2}\right)  ^{n}%
\end{align*}
for any $n\in\mathbb{N}$. This formula is known as \emph{Binet's formula}. It
has many consequences; for example, it implies easily that $\lim
\limits_{n\rightarrow\infty}\dfrac{f_{n+1}}{f_{n}}=\phi_{+}=\dfrac{1+\sqrt{5}%
}{2}\approx1.618\ldots$. Thus, $f_{n}\sim\dfrac{1}{\sqrt{5}}\phi_{+}^{n}$ in
the asymptotical sense.

\subsubsection{Example 2: Dyck words and Catalan numbers}

Before the next example, let us address a warmup question: What is the number
of $2n$-tuples that contain $n$ entries equal to $0$ and $n$ entries equal to
$1$ (where $n$ is a given nonnegative integer)?

(For example, for $n=2$, these $2n$-tuples are $\left(  1,1,0,0\right)  $,
$\left(  1,0,1,0\right)  $, $\left(  1,0,0,1\right)  $, $\left(
0,1,1,0\right)  $, $\left(  0,1,0,1\right)  $, $\left(  0,0,1,1\right)  $, so
there are $6$ of them.)

\textbf{Answer:} The number is $\dbinom{2n}{n}$, since choosing a $2n$-tuple
that contains $n$ entries equal to $0$ and $n$ entries equal to $1$ is
tantamount to choosing an $n$-element subset of $\left\{  1,2,\ldots
,2n\right\}  $ (and we know that there are $\dbinom{2n}{n}$ ways to choose the latter).

\bigskip

\textbf{Example 2.} A \emph{Dyck word} of length $2n$ (where $n\in\mathbb{N}$)
is a $2n$-tuple that contains $n$ entries equal to $0$ and $n$ entries equal
to $1$, and has the additional property that for each $k$, we have
\begin{align}
&  \left(  \text{\# of }0\text{'s among its first }k\text{ entries}\right)
\nonumber\\
&  \leq\left(  \text{\# of }1\text{'s among its first }k\text{ entries}%
\right)  . \label{eq.sec.gf.exas.2.dyckword.leq}%
\end{align}
(The symbol \textquotedblleft\#\textquotedblright\ means \textquotedblleft
number\textquotedblright.)

Some examples: The tuples
\[
\left(  1,0,1,0\right)  ,\ \ \left(  1,1,0,0\right)  ,\ \ \left(
1,1,0,1,0,0\right)  ,\ \ \left(  {}\right)  ,\ \ \left(  1,0\right)
\]
are Dyck words. The tuples
\[
\left(  0,1,1,0\right)  ,\ \ \left(  1,0,0,1\right)  ,\ \ \left(
1,1,0\right)  ,\ \ \left(  1\right)  ,\ \ \left(  1,1,1,0\right)
\]
are not Dyck words.

A \emph{Dyck path} of length $2n$ is a path from the point $\left(
0,0\right)  $ to the point $\left(  2n,0\right)  $ in the Cartesian plane that
moves only using \textquotedblleft NE-steps\textquotedblright\ (i.e., steps of
the form $\left(  x,y\right)  \rightarrow\left(  x+1,y+1\right)  $) and
\textquotedblleft SE-steps\textquotedblright\ (i.e., steps of the form
$\left(  x,y\right)  \rightarrow\left(  x+1,y-1\right)  $) and never falls
below the x-axis (i.e., does not contain any point $\left(  x,y\right)  $ with
$y<0$).

Examples: For $n=2$, the Dyck paths from $\left(  0,0\right)  $ to $\left(
2n,0\right)  $ are%
\[%
\begin{tabular}
[c]{ccc}%
\begin{tikzpicture}
\draw(0, 0) -- (1, 1) -- (2, 0) -- (3, 1) -- (4, 0);
\filldraw(0, 0) circle [fill=red, radius=0.1];
\filldraw(1, 1) circle [fill=red, radius=0.1];
\filldraw(2, 0) circle [fill=red, radius=0.1];
\filldraw(3, 1) circle [fill=red, radius=0.1];
\filldraw(4, 0) circle [fill=red, radius=0.1];
\end{tikzpicture}%
& $\ \ \ \ \ \ \ \ \ \ \ \ \ \ \ \ \ \ \ \ $ &
\begin{tikzpicture}
\draw(0, 0) -- (1, 1) -- (2, 2) -- (3, 1) -- (4, 0);
\filldraw(0, 0) circle [fill=red, radius=0.1];
\filldraw(1, 1) circle [fill=red, radius=0.1];
\filldraw(2, 2) circle [fill=red, radius=0.1];
\filldraw(3, 1) circle [fill=red, radius=0.1];
\filldraw(4, 0) circle [fill=red, radius=0.1];
\end{tikzpicture}%
\end{tabular}
\]

A Dyck path can be viewed as the \textquotedblleft skyline\textquotedblright%
\ of a \textquotedblleft mountain range\textquotedblright. For example:%
\[%
\begin{tabular}
[c]{|c|c|}\hline
Dyck path & \textquotedblleft mountain range\textquotedblright\\\hline
\ \
\begin{tikzpicture}
\draw(0, 0) -- (1, 1) -- (2, 2) -- (3, 1) -- (4, 0) -- (5, 1) -- (6, 0);
\filldraw(0, 0) circle [fill=red, radius=0.1];
\filldraw(1, 1) circle [fill=red, radius=0.1];
\filldraw(2, 2) circle [fill=red, radius=0.1];
\filldraw(3, 1) circle [fill=red, radius=0.1];
\filldraw(4, 0) circle [fill=red, radius=0.1];
\filldraw(5, 1) circle [fill=red, radius=0.1];
\filldraw(6, 0) circle [fill=red, radius=0.1];
\end{tikzpicture}%
\ \  & \ \
\begin{tikzpicture}
\draw(0, 0) -- (1, 1) -- (2, 2) -- (3, 1) -- (4, 0) -- (5, 1) -- (6, 0);
\fill
[green!60!black] (0, 0) -- (1, 1) -- (2, 2) -- (3, 1) -- (4, 0) -- (5, 1) -- (6, 0);
\filldraw(0, 0) circle [fill=red, radius=0.1];
\filldraw(1, 1) circle [fill=red, radius=0.1];
\filldraw(2, 2) circle [fill=red, radius=0.1];
\filldraw(3, 1) circle [fill=red, radius=0.1];
\filldraw(4, 0) circle [fill=red, radius=0.1];
\filldraw(5, 1) circle [fill=red, radius=0.1];
\filldraw(6, 0) circle [fill=red, radius=0.1];
\end{tikzpicture}%
\ \ \\\hline
\end{tabular}
\]

The names \textquotedblleft NE-steps\textquotedblright\ and \textquotedblleft
SE-steps\textquotedblright\ in the definition of a Dyck path refer to compass
directions: If we treat the Cartesian plane as a map with the x-axis directed
eastwards and the y-axis directed northwards, then an NE-step moves to the
northeast, and an SE-step moves to the southeast.

Note that any NE-step and any SE-step increases the x-coordinate by $1$ (that
is, the step goes from a point with x-coordinate $k$ to a point with
x-coordinate $k+1$). Thus, any Dyck path from $\left(  0,0\right)  $ to
$\left(  2n,0\right)  $ has precisely $2n$ steps. Of these $2n$ steps, exactly
$n$ are NE-steps while the remaining $n$ are SE-steps (because any NE-step
increases the y-coordinate by $1$, while any SE-step decreases the
y-coordinate by $1$). Since a Dyck path must never fall below the x-axis, we
see that the number of SE-steps up to any given point can never be larger than
the number of NE-steps up to this point. But this is exactly the condition
(\ref{eq.sec.gf.exas.2.dyckword.leq}) from the definition of a Dyck word,
except that we are talking about NE-steps and SE-steps instead of $1$'s and
$0$'s. Thus, there is a simple bijection between Dyck words of length $2n$ and
Dyck paths from $\left(  0,0\right)  $ to $\left(  2n,0\right)  $:

\begin{itemize}
\item send each $1$ in the Dyck word to a NE-step in the Dyck path;

\item send each $0$ in the Dyck word to a SE-step in the Dyck path.
\end{itemize}

So the \# of Dyck words (of length $2n$) equals the \# of Dyck paths (from
$\left(  0,0\right)  $ to $\left(  2n,0\right)  $). But what is this number?

Example: For $n=3$, this number is $5$. Indeed, here are all Dyck paths from
$\left(  0,0\right)  $ to $\left(  6,0\right)  $, and their corresponding Dyck
words:%
\[%
\begin{tabular}
[c]{|c|c|}\hline
Dyck path & Dyck word\\\hline
\ \
\begin{tikzpicture}
\draw(0, 0) -- (1, 1) -- (2, 2) -- (3, 1) -- (4, 0) -- (5, 1) -- (6, 0);
\filldraw(0, 0) circle [fill=red, radius=0.1];
\filldraw(1, 1) circle [fill=red, radius=0.1];
\filldraw(2, 2) circle [fill=red, radius=0.1];
\filldraw(3, 1) circle [fill=red, radius=0.1];
\filldraw(4, 0) circle [fill=red, radius=0.1];
\filldraw(5, 1) circle [fill=red, radius=0.1];
\filldraw(6, 0) circle [fill=red, radius=0.1];
\end{tikzpicture}%
\ \  & $\left(  1,1,0,0,1,0\right)  $\\\hline
\ \
\begin{tikzpicture}
\draw(0, 0) -- (1, 1) -- (2, 2) -- (3, 3) -- (4, 2) -- (5, 1) -- (6, 0);
\filldraw(0, 0) circle [fill=red, radius=0.1];
\filldraw(1, 1) circle [fill=red, radius=0.1];
\filldraw(2, 2) circle [fill=red, radius=0.1];
\filldraw(3, 3) circle [fill=red, radius=0.1];
\filldraw(4, 2) circle [fill=red, radius=0.1];
\filldraw(5, 1) circle [fill=red, radius=0.1];
\filldraw(6, 0) circle [fill=red, radius=0.1];
\end{tikzpicture}%
\ \  & $\left(  1,1,1,0,0,0\right)  $\\\hline
\ \
\begin{tikzpicture}
\draw(0, 0) -- (1, 1) -- (2, 0) -- (3, 1) -- (4, 0) -- (5, 1) -- (6, 0);
\filldraw(0, 0) circle [fill=red, radius=0.1];
\filldraw(1, 1) circle [fill=red, radius=0.1];
\filldraw(2, 0) circle [fill=red, radius=0.1];
\filldraw(3, 1) circle [fill=red, radius=0.1];
\filldraw(4, 0) circle [fill=red, radius=0.1];
\filldraw(5, 1) circle [fill=red, radius=0.1];
\filldraw(6, 0) circle [fill=red, radius=0.1];
\end{tikzpicture}%
\ \  & $\left(  1,0,1,0,1,0\right)  $\\\hline
\ \
\begin{tikzpicture}
\draw(0, 0) -- (1, 1) -- (2, 0) -- (3, 1) -- (4, 2) -- (5, 1) -- (6, 0);
\filldraw(0, 0) circle [fill=red, radius=0.1];
\filldraw(1, 1) circle [fill=red, radius=0.1];
\filldraw(2, 0) circle [fill=red, radius=0.1];
\filldraw(3, 1) circle [fill=red, radius=0.1];
\filldraw(4, 2) circle [fill=red, radius=0.1];
\filldraw(5, 1) circle [fill=red, radius=0.1];
\filldraw(6, 0) circle [fill=red, radius=0.1];
\end{tikzpicture}%
\ \  & $\left(  1,0,1,1,0,0\right)  $\\\hline
\ \
\begin{tikzpicture}
\draw(0, 0) -- (1, 1) -- (2, 2) -- (3, 1) -- (4, 2) -- (5, 1) -- (6, 0);
\filldraw(0, 0) circle [fill=red, radius=0.1];
\filldraw(1, 1) circle [fill=red, radius=0.1];
\filldraw(2, 2) circle [fill=red, radius=0.1];
\filldraw(3, 1) circle [fill=red, radius=0.1];
\filldraw(4, 2) circle [fill=red, radius=0.1];
\filldraw(5, 1) circle [fill=red, radius=0.1];
\filldraw(6, 0) circle [fill=red, radius=0.1];
\end{tikzpicture}%
\ \  & $\left(  1,1,0,1,0,0\right)  $\\\hline
\end{tabular}
\]
(We will soon stop writing the commas and parentheses when writing down words.
For example, the word $\left(  1,1,0,0,1,0\right)  $ will just become $110010$.)

Back to the general question.

For each $n\in\mathbb{N}$, we define
\begin{align*}
c_{n}:=  &  \ \left(  \text{\# of Dyck paths from }\left(  0,0\right)  \text{
to }\left(  2n,0\right)  \right) \\
=  &  \ \left(  \text{\# of Dyck words of length }2n\right)
\ \ \ \ \ \ \ \ \ \ \left(  \text{as we have seen above}\right)  .
\end{align*}
Then, $c_{0}=1$ (since the only Dyck path from $\left(  0,0\right)  $ to
$\left(  0,0\right)  $ is the trivial path) and $c_{1}=1$ and $c_{2}=2$ and
$c_{3}=5$ and $c_{4}=14$ and so on. These numbers $c_{n}$ are known as the
\emph{Catalan numbers}. Entire books have been written about them, such as
\cite{Stanle15}.

Let us first find a recurrence relation for $c_{n}$. The argument below is
best understood by following an example; namely, consider the following Dyck
path from $\left(  0,0\right)  $ to $\left(  16,0\right)  $ (so the
corresponding $n$ is $8$):%
\[%
\begin{tikzpicture}
\draw
(0, 0) -- (1, 1) -- (2, 2) -- (3, 1) -- (4, 2) -- (5, 1) -- (6, 0) -- (7, 1) -- (8, 2) -- (9, 1) -- (10, 0) -- (11, 1) -- (12, 0) -- (13, 1) -- (14, 2) -- (15, 1) -- (16, 0);
\filldraw(0, 0) circle [fill=red, radius=0.1];
\filldraw(1, 1) circle [fill=red, radius=0.1];
\filldraw(2, 2) circle [fill=red, radius=0.1];
\filldraw(3, 1) circle [fill=red, radius=0.1];
\filldraw(4, 2) circle [fill=red, radius=0.1];
\filldraw(5, 1) circle [fill=red, radius=0.1];
\filldraw(6, 0) circle [fill=red, radius=0.1];
\filldraw(7, 1) circle [fill=red, radius=0.1];
\filldraw(8, 2) circle [fill=red, radius=0.1];
\filldraw(9, 1) circle [fill=red, radius=0.1];
\filldraw(10, 0) circle [fill=red, radius=0.1];
\filldraw(11, 1) circle [fill=red, radius=0.1];
\filldraw(12, 0) circle [fill=red, radius=0.1];
\filldraw(13, 1) circle [fill=red, radius=0.1];
\filldraw(14, 2) circle [fill=red, radius=0.1];
\filldraw(15, 1) circle [fill=red, radius=0.1];
\filldraw(16, 0) circle [fill=red, radius=0.1];
\end{tikzpicture}%
\]

Fix a positive integer $n$. If $d$ is a Dyck path from $\left(  0,0\right)  $
to $\left(  2n,0\right)  $, then the \emph{first return} of $d$ (this is short
for \textquotedblleft first return of $d$ to the x-axis\textquotedblright)
shall mean the first point on $d$ that lies on the x-axis but is not the
origin (i.e., that has the form $\left(  i,0\right)  $ for some integer
$i>0$). For instance, in the example that we just gave, the first return is
the point $\left(  6,0\right)  $. If $d$ is a Dyck path from $\left(
0,0\right)  $ to $\left(  2n,0\right)  $, and if $\left(  i,0\right)  $ is its
first return, then $i$ is even\footnote{\textit{Proof.} The number of NE-steps
before the first return must equal the number of SE-steps before the first
return (because these steps have altogether taken us from the origin to a
point on the x-axis, and thus must have increased and decreased the
y-coordinate an equal number of times). This shows that the total number of
steps before the first return is even. In other words, $i$ is even (because
the total number of steps before the first return is $i$).}, and therefore we
have $i=2k$ for some $k\in\left\{  1,2,\ldots,n\right\}  $. Hence, for any
Dyck path from $\left(  0,0\right)  $ to $\left(  2n,0\right)  $, the first
return is a point of the form $\left(  2k,0\right)  $ for some $k\in\left\{
1,2,\ldots,n\right\}  $. Thus,%
\begin{align*}
&  \left(  \text{\# of Dyck paths from }\left(  0,0\right)  \text{ to }\left(
2n,0\right)  \right) \\
&  =\sum_{k=1}^{n}\left(  \text{\# of Dyck paths from }\left(  0,0\right)
\text{ to }\left(  2n,0\right)  \text{ whose first return is }\left(
2k,0\right)  \right)  .
\end{align*}

Now, let us fix some $k\in\left\{  1,2,\ldots,n\right\}  $. We shall compute
the \# of Dyck paths from $\left(  0,0\right)  $ to $\left(  2n,0\right)  $
whose first return is $\left(  2k,0\right)  $. Any such Dyck path has a
natural \textquotedblleft two-part\textquotedblright\ structure: Its first
$2k$ steps form a path from $\left(  0,0\right)  $ to $\left(  2k,0\right)  $,
while its last (i.e., remaining) $2\left(  n-k\right)  $ steps form a path
from $\left(  2k,0\right)  $ to $\left(  2n,0\right)  $. Thus, in order to
construct such a path, we

\begin{itemize}
\item first choose its first $2k$ steps: They have to form a Dyck path from
$\left(  0,0\right)  $ to $\left(  2k,0\right)  $ that never returns to the
x-axis until $\left(  2k,0\right)  $. Hence, they begin with a NE-step and end
with a SE-step (since any other steps here would cause the path to fall below
the x-axis). Between these two steps, the remaining $2k-2=2\left(  k-1\right)
$ steps form a path that not only never falls below the x-axis, but also never
touches it (since $\left(  2k,0\right)  $ is the first return of our Dyck
path, so that our Dyck path does not touch the x-axis between $\left(
0,0\right)  $ and $\left(  2k,0\right)  $). In other words, these $2\left(
k-1\right)  $ steps form a path from $\left(  1,1\right)  $ to $\left(
2k-1,1\right)  $ that never falls below the $y=1$ line (i.e., below the x-axis
shifted by $1$ upwards). This means that it is a Dyck path from $\left(
0,0\right)  $ to $\left(  2\left(  k-1\right)  ,0\right)  $ (shifted by
$\left(  1,1\right)  $). Thus, there are $c_{k-1}$ possibilities for this
path. Hence, there are $c_{k-1}$ choices for the first $2k$ steps of our Dyck path.

\item then choose its last $2\left(  n-k\right)  $ steps: They have to form a
path from $\left(  2k,0\right)  $ to $\left(  2n,0\right)  $ that never falls
below the x-axis (but is allowed to touch it any number of times). Thus, they
form a Dyck path from $\left(  0,0\right)  $ to $\left(  2\left(  n-k\right)
,0\right)  $ (shifted by $\left(  2k,0\right)  $). So there are $c_{n-k}$
choices for these last $2\left(  n-k\right)  $ steps.
\end{itemize}

Thus, there are $c_{k-1}c_{n-k}$ many options for such a Dyck path from
$\left(  0,0\right)  $ to $\left(  2n,0\right)  $ (since choosing the first
$2k$ steps and choosing the last $2\left(  n-k\right)  $ steps are independent).

Let me illustrate this reasoning on the Dyck path from $\left(  0,0\right)  $
to $\left(  16,0\right)  $ shown above. This Dyck path has first return
$\left(  6,0\right)  $; thus, the corresponding $k$ is $3$. Since this Dyck
path does not return to the x-axis before $\left(  2k,0\right)  =\left(
6,0\right)  $, its first $2k$ steps stay above (or on) the yellow trapezoid
shown here:%
\[%
\begin{tikzpicture}
\draw
(0, 0) -- (1, 1) -- (2, 2) -- (3, 1) -- (4, 2) -- (5, 1) -- (6, 0) -- (7, 1) -- (8, 2) -- (9, 1) -- (10, 0) -- (11, 1) -- (12, 0) -- (13, 1) -- (14, 2) -- (15, 1) -- (16, 0);
\fill[yellow] (0, 0) -- (1, 1) -- (5, 1) -- (6, 0) -- cycle;
\filldraw(0, 0) circle [fill=red, radius=0.1];
\filldraw(1, 1) circle [fill=red, radius=0.1];
\filldraw(2, 2) circle [fill=red, radius=0.1];
\filldraw(3, 1) circle [fill=red, radius=0.1];
\filldraw(4, 2) circle [fill=red, radius=0.1];
\filldraw(5, 1) circle [fill=red, radius=0.1];
\filldraw(6, 0) circle [fill=red, radius=0.1];
\filldraw(7, 1) circle [fill=red, radius=0.1];
\filldraw(8, 2) circle [fill=red, radius=0.1];
\filldraw(9, 1) circle [fill=red, radius=0.1];
\filldraw(10, 0) circle [fill=red, radius=0.1];
\filldraw(11, 1) circle [fill=red, radius=0.1];
\filldraw(12, 0) circle [fill=red, radius=0.1];
\filldraw(13, 1) circle [fill=red, radius=0.1];
\filldraw(14, 2) circle [fill=red, radius=0.1];
\filldraw(15, 1) circle [fill=red, radius=0.1];
\filldraw(16, 0) circle [fill=red, radius=0.1];
\end{tikzpicture}%
\]
In particular, the first and the last of these $2k$ steps are uniquely
determined, while the steps between them form a diagonally shifted Dyck path
that is filled in green here:%
\[%
\begin{tikzpicture}
\draw
(0, 0) -- (1, 1) -- (2, 2) -- (3, 1) -- (4, 2) -- (5, 1) -- (6, 0) -- (7, 1) -- (8, 2) -- (9, 1) -- (10, 0) -- (11, 1) -- (12, 0) -- (13, 1) -- (14, 2) -- (15, 1) -- (16, 0);
\fill[yellow] (0, 0) -- (1, 1) -- (5, 1) -- (6, 0) -- cycle;
\fill
[green!60!black] (1, 1) -- (2, 2) -- (3, 1) -- (4, 2) -- (5, 1) -- (1, 1);
\filldraw(0, 0) circle [fill=red, radius=0.1];
\filldraw(1, 1) circle [fill=red, radius=0.1];
\filldraw(2, 2) circle [fill=red, radius=0.1];
\filldraw(3, 1) circle [fill=red, radius=0.1];
\filldraw(4, 2) circle [fill=red, radius=0.1];
\filldraw(5, 1) circle [fill=red, radius=0.1];
\filldraw(6, 0) circle [fill=red, radius=0.1];
\filldraw(7, 1) circle [fill=red, radius=0.1];
\filldraw(8, 2) circle [fill=red, radius=0.1];
\filldraw(9, 1) circle [fill=red, radius=0.1];
\filldraw(10, 0) circle [fill=red, radius=0.1];
\filldraw(11, 1) circle [fill=red, radius=0.1];
\filldraw(12, 0) circle [fill=red, radius=0.1];
\filldraw(13, 1) circle [fill=red, radius=0.1];
\filldraw(14, 2) circle [fill=red, radius=0.1];
\filldraw(15, 1) circle [fill=red, radius=0.1];
\filldraw(16, 0) circle [fill=red, radius=0.1];
\end{tikzpicture}%
\]
Finally, the last $2\left(  n-k\right)  $ steps form a horizontally shifted
Dyck path that is filled in purple here:%
\[%
\begin{tikzpicture}
\draw
(0, 0) -- (1, 1) -- (2, 2) -- (3, 1) -- (4, 2) -- (5, 1) -- (6, 0) -- (7, 1) -- (8, 2) -- (9, 1) -- (10, 0) -- (11, 1) -- (12, 0) -- (13, 1) -- (14, 2) -- (15, 1) -- (16, 0);
\fill[yellow] (0, 0) -- (1, 1) -- (5, 1) -- (6, 0) -- cycle;
\fill
[green!60!black] (1, 1) -- (2, 2) -- (3, 1) -- (4, 2) -- (5, 1) -- (1, 1);
\fill
[blue!60!white] (6, 0) -- (7, 1) -- (8, 2) -- (9, 1) -- (10, 0) -- (11, 1) -- (12, 0) -- (13, 1) -- (14, 2) -- (15, 1) -- (16, 0);
\filldraw(0, 0) circle [fill=red, radius=0.1];
\filldraw(1, 1) circle [fill=red, radius=0.1];
\filldraw(2, 2) circle [fill=red, radius=0.1];
\filldraw(3, 1) circle [fill=red, radius=0.1];
\filldraw(4, 2) circle [fill=red, radius=0.1];
\filldraw(5, 1) circle [fill=red, radius=0.1];
\filldraw(6, 0) circle [fill=red, radius=0.1];
\filldraw(7, 1) circle [fill=red, radius=0.1];
\filldraw(8, 2) circle [fill=red, radius=0.1];
\filldraw(9, 1) circle [fill=red, radius=0.1];
\filldraw(10, 0) circle [fill=red, radius=0.1];
\filldraw(11, 1) circle [fill=red, radius=0.1];
\filldraw(12, 0) circle [fill=red, radius=0.1];
\filldraw(13, 1) circle [fill=red, radius=0.1];
\filldraw(14, 2) circle [fill=red, radius=0.1];
\filldraw(15, 1) circle [fill=red, radius=0.1];
\filldraw(16, 0) circle [fill=red, radius=0.1];
\end{tikzpicture}%
\]
Our above argument shows that there are $c_{k-1}$ choices for the green Dyck
path and $c_{n-k}$ choices for the purple Dyck path, therefore $c_{k-1}%
c_{n-k}$ options in total.

Forget that we fixed $k$. Our counting argument above shows that%
\begin{align}
&  \left(  \text{\# of Dyck paths from }\left(  0,0\right)  \text{ to }\left(
2n,0\right)  \text{ whose first return is }\left(  2k,0\right)  \right)
\nonumber\\
&  =c_{k-1}c_{n-k} \label{eq.sec.gf.exas.2.dyckword.knum}%
\end{align}
for each $k\in\left\{  1,2,\ldots,n\right\}  $. Now,
\begin{align*}
c_{n}  &  =\left(  \text{\# of Dyck paths from }\left(  0,0\right)  \text{ to
}\left(  2n,0\right)  \right) \\
&  =\sum_{k=1}^{n}\underbrace{\left(  \text{\# of Dyck paths from }\left(
0,0\right)  \text{ to }\left(  2n,0\right)  \text{ whose first return is
}\left(  2k,0\right)  \right)  }_{\substack{=c_{k-1}c_{n-k}\\\text{(by
(\ref{eq.sec.gf.exas.2.dyckword.knum}))}}}\\
&  =\sum_{k=1}^{n}c_{k-1}c_{n-k}=c_{0}c_{n-1}+c_{1}c_{n-2}+c_{2}c_{n-3}%
+\cdots+c_{n-1}c_{0}.
\end{align*}
This is a recurrence equation for $c_{n}$. Combining it with $c_{0}=1$, we can
use it to compute any value of $c_{n}$ recursively. Let us, however, try to
digest it using generating functions!

Let%
\[
C\left(  x\right)  :=\sum_{n\geq0}c_{n}x^{n}=c_{0}+c_{1}x+c_{2}x^{2}%
+c_{3}x^{3}+\cdots.
\]
Thus,%
\begin{align*}
C\left(  x\right)   &  =c_{0}+c_{1}x+c_{2}x^{2}+c_{3}x^{3}+\cdots\\
&  =1+\left(  c_{0}c_{0}\right)  x+\left(  c_{0}c_{1}+c_{1}c_{0}\right)
x^{2}+\left(  c_{0}c_{2}+c_{1}c_{1}+c_{2}c_{0}\right)  x^{3}+\cdots\\
&  \ \ \ \ \ \ \ \ \ \ \ \ \ \ \ \ \ \ \ \ \left(
\begin{array}
[c]{c}%
\text{since }c_{0}=1\\
\text{and }c_{n}=c_{0}c_{n-1}+c_{1}c_{n-2}+c_{2}c_{n-3}+\cdots+c_{n-1}c_{0}\\
\text{for each }n>0
\end{array}
\right) \\
&  =1+x\underbrace{\left(  \left(  c_{0}c_{0}\right)  +\left(  c_{0}%
c_{1}+c_{1}c_{0}\right)  x+\left(  c_{0}c_{2}+c_{1}c_{1}+c_{2}c_{0}\right)
x^{2}+\cdots\right)  }_{\substack{=\left(  c_{0}+c_{1}x+c_{2}x^{2}%
+\cdots\right)  ^{2}\\\text{(because if we multiply out }\left(  c_{0}%
+c_{1}x+c_{2}x^{2}+\cdots\right)  ^{2}\\\text{and collect like powers of
}x\text{, we obtain}\\\text{exactly }\left(  c_{0}c_{0}\right)  +\left(
c_{0}c_{1}+c_{1}c_{0}\right)  x+\left(  c_{0}c_{2}+c_{1}c_{1}+c_{2}%
c_{0}\right)  x^{2}+\cdots\text{)}}}\\
&  =1+x\left(  \underbrace{c_{0}+c_{1}x+c_{2}x^{2}+\cdots}_{=C\left(
x\right)  }\right)  ^{2}=1+x\left(  C\left(  x\right)  \right)  ^{2}.
\end{align*}
This is a quadratic equation in $C\left(  x\right)  $. Let us solve it by the
quadratic formula (assuming for now that this is allowed -- i.e., that the
quadratic formula really does apply to our \textquotedblleft power
series\textquotedblright, whatever they are). Thus, we get%
\begin{equation}
C\left(  x\right)  =\dfrac{1\pm\sqrt{1-4x}}{2x}. \label{eq.sec.gf.exas.2.-+}%
\end{equation}
The $\pm$ sign here cannot be a $+$ sign, because if it were a $+$, then the
power series on top of the fraction would not be divisible by $2x$ (as its
constant term would be $2$ and thus nonzero\footnote{If you find this
unconvincing, here is a cleaner way to argue this: Multiplying the equality
(\ref{eq.sec.gf.exas.2.-+}) by $2x$, we obtain $2xC\left(  x\right)
=1\pm\sqrt{1-4x}$. The left hand side of this equality has constant term $0$,
but the right hand side has constant term $1\pm1$ (here, we are making the
assumption that $\sqrt{1-4x}$ is a power series with constant term $1$; this
is plausible because $\sqrt{1-4\cdot0}=1$ and will also be justified further
below). Thus, $0=1\pm1$; this shows that the $\pm$ sign is a $-$ sign.}).
Thus, (\ref{eq.sec.gf.exas.2.-+}) becomes%
\begin{equation}
C\left(  x\right)  =\dfrac{1-\sqrt{1-4x}}{2x}=\dfrac{1}{2x}\left(  1-\left(
1-4x\right)  ^{1/2}\right)  . \label{eq.sec.gf.exas.2.-}%
\end{equation}

How do we find the coefficients of the power series $\left(  1-4x\right)
^{1/2}$ ?

For each $n\in\mathbb{N}$, the binomial formula yields%
\begin{equation}
\left(  1+x\right)  ^{n}=\sum_{k=0}^{n}\dbinom{n}{k}x^{k}=\sum_{k\geq0}%
\dbinom{n}{k}x^{k}. \label{eq.sec.gf.exas.2.(1+x)n}%
\end{equation}
(Here, we have replaced the $\sum_{k=0}^{n}$ sign by a $\sum_{k\geq0}$ sign,
thus extending the summation from all $k\in\left\{  0,1,\ldots,n\right\}  $ to
all $k\in\mathbb{N}$. This does not change the value of the sum, since all the
newly appearing addends are $0$, as you can easily check.)

Let us pretend that the formula (\ref{eq.sec.gf.exas.2.(1+x)n}) holds not only
for $n\in\mathbb{N}$, but also for $n=1/2$. That is, we have%
\begin{equation}
\left(  1+x\right)  ^{1/2}=\sum_{k\geq0}\dbinom{1/2}{k}x^{k}.
\label{eq.sec.gf.exas.2.(1+x)1/2}%
\end{equation}
Now, substitute $-4x$ for $x$ in this equality (here we are making the rather
plausible assumption that we can substitute $-4x$ for $x$ in a power series);
then, we get%
\begin{align*}
\left(  1-4x\right)  ^{1/2}  &  =\sum_{k\geq0}\dbinom{1/2}{k}\left(
-4x\right)  ^{k}=\sum_{k\geq0}\dbinom{1/2}{k}\left(  -4\right)  ^{k}x^{k}\\
&  =\underbrace{\dbinom{1/2}{0}}_{=1}\underbrace{\left(  -4\right)  ^{0}}%
_{=1}\underbrace{x^{0}}_{=1}+\sum\limits_{k\geq1}\dbinom{1/2}{k}\left(
-4\right)  ^{k}x^{k}\\
&  =1+\sum\limits_{k\geq1}\dbinom{1/2}{k}\left(  -4\right)  ^{k}x^{k}.
\end{align*}
Hence,%
\[
1-\left(  1-4x\right)  ^{1/2}=1-\left(  1+\sum\limits_{k\geq1}\dbinom{1/2}%
{k}\left(  -4\right)  ^{k}x^{k}\right)  =-\sum\limits_{k\geq1}\dbinom{1/2}%
{k}\left(  -4\right)  ^{k}x^{k}.
\]
Thus, (\ref{eq.sec.gf.exas.2.-}) becomes%
\begin{align*}
C\left(  x\right)   &  =\dfrac{1}{2x}\underbrace{\left(  1-\left(
1-4x\right)  ^{1/2}\right)  }_{=-\sum\limits_{k\geq1}\dbinom{1/2}{k}\left(
-4\right)  ^{k}x^{k}}=\dfrac{1}{2x}\left(  -\sum\limits_{k\geq1}\dbinom
{1/2}{k}\left(  -4\right)  ^{k}x^{k}\right) \\
&  =\sum\limits_{k\geq1}\dbinom{1/2}{k}\underbrace{\dfrac{-\left(  -4\right)
^{k}x^{k}}{2x}}_{=2\left(  -4\right)  ^{k-1}x^{k-1}}=\sum\limits_{k\geq
1}\dbinom{1/2}{k}2\left(  -4\right)  ^{k-1}x^{k-1}\\
&  =\sum\limits_{k\geq0}\dbinom{1/2}{k+1}2\left(  -4\right)  ^{k}x^{k}\\
&  \ \ \ \ \ \ \ \ \ \ \ \ \ \ \ \ \ \ \ \ \left(  \text{here, we have
substituted }k+1\text{ for }k\text{ in the sum}\right)  .
\end{align*}
Comparing coefficients before $x^{n}$ in this equality gives%
\begin{equation}
c_{n}=\dbinom{1/2}{n+1}2\left(  -4\right)  ^{n}.
\label{eq.sec.gf.exas.cn=frac}%
\end{equation}

This is an explicit formula for $c_{n}$ (and makes computation of $c_{n}$
pretty easy!), but it turns out that it can be simplified further. Indeed, the
definition of $\dbinom{1/2}{n+1}$ yields%
\begin{align*}
\dbinom{1/2}{n+1}  &  =\dfrac{\left(  1/2\right)  \left(  1/2-1\right)
\left(  1/2-2\right)  \cdots\left(  1/2-n\right)  }{\left(  n+1\right)  !}\\
&  =\dfrac{\vphantom{\dfrac{f}{f}}\dfrac{1}{2}\cdot\dfrac{-1}{2}\cdot
\dfrac{-3}{2}\cdot\dfrac{-5}{2}\cdot\cdots\cdot\dfrac{-\left(  2n-1\right)
}{2}}{\left(  n+1\right)  !}\\
&  =\dfrac{\left(  1\cdot\left(  -1\right)  \cdot\left(  -3\right)
\cdot\left(  -5\right)  \cdot\cdots\cdot\left(  -\left(  2n-1\right)  \right)
\right)  /2^{n+1}}{\left(  n+1\right)  !}\\
&  =\dfrac{\left(  \left(  -1\right)  \cdot\left(  -3\right)  \cdot\left(
-5\right)  \cdot\cdots\cdot\left(  -\left(  2n-1\right)  \right)  \right)
/2^{n+1}}{\left(  n+1\right)  !}\\
&  =\dfrac{\left(  -1\right)  ^{n}\left(  1\cdot3\cdot5\cdot\cdots\cdot\left(
2n-1\right)  \right)  /2^{n+1}}{\left(  n+1\right)  !}.
\end{align*}
Thus, (\ref{eq.sec.gf.exas.cn=frac}) rewrites as
\begin{align*}
c_{n}  &  =\dfrac{\left(  -1\right)  ^{n}\left(  1\cdot3\cdot5\cdot\cdots
\cdot\left(  2n-1\right)  \right)  /2^{n+1}}{\left(  n+1\right)  !}%
\cdot2\left(  -4\right)  ^{n}\\
&  =\underbrace{\dfrac{1\cdot3\cdot5\cdot\cdots\cdot\left(  2n-1\right)
}{\left(  n+1\right)  !}}_{\substack{=\dfrac{1}{n+1}\cdot\dfrac{1\cdot
3\cdot5\cdot\cdots\cdot\left(  2n-1\right)  }{n!}\\\text{(since }\left(
n+1\right)  !=\left(  n+1\right)  \cdot n!\text{)}}}\cdot\,\underbrace{\dfrac
{\left(  -1\right)  ^{n}\cdot2\left(  -4\right)  ^{n}}{2^{n+1}}}_{=2^{n}}\\
&  =\dfrac{1}{n+1}\cdot\underbrace{\dfrac{1\cdot3\cdot5\cdot\cdots\cdot\left(
2n-1\right)  }{n!}\cdot2^{n}}_{\substack{=\dbinom{2n}{n}\\\text{(by Example
\ref{exa.binom.2n-choose-n})}}}=\dfrac{1}{n+1}\dbinom{2n}{n}.
\end{align*}
Hence, we have shown that
\begin{equation}
c_{n}=\dfrac{1}{n+1}\dbinom{2n}{n}. \label{eq.sec.gf.exas.cn=1/(n+1)}%
\end{equation}
Moreover, we can rewrite this further as%
\begin{equation}
c_{n}=\dbinom{2n}{n}-\dbinom{2n}{n-1} \label{eq.sec.gf.exas.cn=diff}%
\end{equation}
(since another binomial coefficient manipulation\footnote{See Exercise
\ref{exe.gf.cn-easy-manip} for this.} yields $\dfrac{1}{n+1}\dbinom{2n}%
{n}=\dbinom{2n}{n}-\dbinom{2n}{n-1}$).

Here is the upshot: The \# of Dyck words of length $2n$ is $c_{n}=\dfrac
{1}{n+1}\dbinom{2n}{n}$. In other words, a $2n$-tuple that consists of $n$
entries equal to $0$ and $n$ entries equal to $1$ (chosen uniformly at random)
is a Dyck word with probability $\dfrac{1}{n+1}$.

(There are also combinatorial ways to prove this; see, e.g., \cite[\S 7.5,
discussion at the end of Example 4]{GKP} or \cite[\S 1.6]{Stanle15} or
\cite{Martin-Dyck} or \cite[Lecture 29, Theorem 5.6.7]{22fco} or \cite[Theorem
1.56]{Loehr-BC}\footnote{Note that the \textquotedblleft Dyck
paths\textquotedblright\ in \cite{Loehr-BC} differ from ours in that they use
N-steps (i.e., steps $\left(  i,j\right)  \mapsto\left(  i,j+1\right)  $) and
E-steps (i.e., steps $\left(  i,j\right)  \mapsto\left(  i+1,j\right)  $)
instead of NE-steps and SE-steps, and stay above the $x=y$ line instead of
above the x-axis. But this notion of Dyck paths is equivalent to ours, since a
clockwise rotation by $45^{\circ}$ followed by a $\sqrt{2}$-homothety
transforms it into ours.} or \cite[\S 8.5, proofs of Identity 244]%
{Spivey19}\footnote{Again, \cite{Spivey19} works not directly with Dyck paths,
but rather with paths that use E-steps (i.e., steps $\left(  i,j\right)
\mapsto\left(  i+1,j\right)  $) and N-steps (i.e., steps $\left(  i,j\right)
\mapsto\left(  i,j+1\right)  $) instead of NE-steps and SE-steps, and stay
below the $x=y$ line instead of above the x-axis. But this kind of Dyck paths
is equivalent to our Dyck paths, since a reflection across the $x=y$ line,
followed by a clockwise rotation by $45^{\circ}$ followed by a $\sqrt{2}%
$-homothety transforms it into ours.}.)

Here is a table of the first $12$ Catalan numbers $c_{n}$:%
\[%
\begin{tabular}
[c]{|c||c|c|c|c|c|c|c|c|c|c|c|c|}\hline
$n$ & $0$ & $1$ & $2$ & $3$ & $4$ & $5$ & $6$ & $7$ & $8$ & $9$ & $10$ &
$11$\\\hline
$c_{n}$ & $1$ & $1$ & $2$ & $5$ & $14$ & $42$ & $132$ & $429$ & $1430$ &
$4862$ & $16\,796$ & $58\,786$\\\hline
\end{tabular}
\ \ \ \ \ .
\]

\subsubsection{Example 3: The Vandermonde convolution}

\textbf{Example 3:} The \emph{Vandermonde convolution identity} (also known as
the \emph{Chu--Vandermonde identity}) says that%
\[
\dbinom{a+b}{n}=\sum_{k=0}^{n}\dbinom{a}{k}\dbinom{b}{n-k}%
\ \ \ \ \ \ \ \ \ \ \text{for any numbers }a,b\text{ and any }n\in\mathbb{N}%
\]
(where \textquotedblleft numbers\textquotedblright\ can mean, e.g.,
\textquotedblleft complex numbers\textquotedblright).

Let us prove this using generating functions. For now, we shall only prove
this for $a,b\in\mathbb{N}$; later I will explain why it also holds for
arbitrary (rational, real or complex) numbers $a,b$ as well.

Indeed, fix $a,b\in\mathbb{N}$. Recall (from (\ref{eq.sec.gf.exas.2.(1+x)n}))
that
\[
\left(  1+x\right)  ^{n}=\sum_{k\geq0}\dbinom{n}{k}x^{k}%
\]
for each $n\in\mathbb{N}$. Hence,%
\begin{align*}
\left(  1+x\right)  ^{a}  &  =\sum_{k\geq0}\dbinom{a}{k}x^{k}%
\ \ \ \ \ \ \ \ \ \ \text{and}\\
\left(  1+x\right)  ^{b}  &  =\sum_{k\geq0}\dbinom{b}{k}x^{k}%
\ \ \ \ \ \ \ \ \ \ \text{and}\\
\left(  1+x\right)  ^{a+b}  &  =\sum_{k\geq0}\dbinom{a+b}{k}x^{k}.
\end{align*}
Thus,
\begin{align*}
\underbrace{\left(  1+x\right)  ^{a}}_{=\sum_{k\geq0}\dbinom{a}{k}x^{k}%
}\ \ \underbrace{\left(  1+x\right)  ^{b}}_{=\sum_{k\geq0}\dbinom{b}{k}x^{k}}
&  =\left(  \sum_{k\geq0}\dbinom{a}{k}x^{k}\right)  \left(  \sum_{k\geq
0}\dbinom{b}{k}x^{k}\right) \\
&  =\left(  \sum_{k\geq0}\dbinom{a}{k}x^{k}\right)  \left(  \sum_{\ell\geq
0}\dbinom{b}{\ell}x^{\ell}\right) \\
&  =\sum_{k\geq0}\ \ \sum_{\ell\geq0}\dbinom{a}{k}x^{k}\dbinom{b}{\ell}%
x^{\ell}\\
&  =\sum_{k\geq0}\ \ \sum_{\ell\geq0}\dbinom{a}{k}\dbinom{b}{\ell}x^{k+\ell}\\
&  =\sum_{n\geq0}\left(  \sum_{k=0}^{n}\dbinom{a}{k}\dbinom{b}{n-k}\right)
x^{n}%
\end{align*}
(here, we have merged addends in which $x$ appears in the same power). Hence,%
\begin{align*}
\sum_{n\geq0}\left(  \sum_{k=0}^{n}\dbinom{a}{k}\dbinom{b}{n-k}\right)  x^{n}
&  =\left(  1+x\right)  ^{a}\left(  1+x\right)  ^{b}=\left(  1+x\right)
^{a+b}=\sum_{k\geq0}\dbinom{a+b}{k}x^{k}\\
&  =\sum_{n\geq0}\dbinom{a+b}{n}x^{n}.
\end{align*}
Comparing coefficients in this equality, we obtain%
\[
\sum_{k=0}^{n}\dbinom{a}{k}\dbinom{b}{n-k}=\dbinom{a+b}{n}%
\ \ \ \ \ \ \ \ \ \ \text{for each }n\in\mathbb{N}.
\]
This completes the proof of the Vandermonde convolution identity for
$a,b\in\mathbb{N}$.

\subsubsection{Example 4: Solving a recurrence}

\textbf{Example 4.} The following example is from \cite[\S 1.2]{Wilf04}.
Define a sequence $\left(  a_{0},a_{1},a_{2},\ldots\right)  $ of numbers
recursively by%
\[
a_{0}=1,\ \ \ \ \ \ \ \ \ \ a_{n+1}=2a_{n}+n\text{ for all }n\geq0.
\]
Thus, the first entries of this sequence are $1,2,5,12,27,58,121,\ldots$. This
sequence appears in \href{https://oeis.org/}{the OEIS (= Online Encyclopedia
of Integer Sequences)} as \href{https://oeis.org/A000325}{A000325}, up to a
shift of the index.

Can we find an explicit formula for $a_{n}$ (without looking it up in the OEIS)?

Again, generating functions are helpful. Set%
\[
A\left(  x\right)  =a_{0}+a_{1}x+a_{2}x^{2}+a_{3}x^{3}+\cdots.
\]
Then,%
\begin{align}
A\left(  x\right)   &  =a_{0}+a_{1}x+a_{2}x^{2}+a_{3}x^{3}+\cdots\nonumber\\
&  =1+\left(  2a_{0}+0\right)  x+\left(  2a_{1}+1\right)  x^{2}+\left(
2a_{2}+2\right)  x^{3}+\cdots\nonumber\\
&  \ \ \ \ \ \ \ \ \ \ \ \ \ \ \ \ \ \ \ \ \left(  \text{since }a_{0}=1\text{
and }a_{n+1}=2a_{n}+n\text{ for all }n\geq0\right) \nonumber\\
&  =1+2\underbrace{\left(  a_{0}x+a_{1}x^{2}+a_{2}x^{3}+\cdots\right)
}_{=xA\left(  x\right)  }+\underbrace{\left(  0x+1x^{2}+2x^{3}+\cdots\right)
}_{=x\left(  0+1x+2x^{2}+3x^{3}+\cdots\right)  }\nonumber\\
&  =1+2xA\left(  x\right)  +x\left(  0+1x+2x^{2}+3x^{3}+\cdots\right)  .
\label{eq.sec.gf.exas.4.Ax=1}%
\end{align}

Thus, it would clearly be helpful to find a simple expression for
$0+1x+2x^{2}+3x^{3}+\cdots$. Here are two ways to do so: \medskip

\textit{First way:} We assume that our power series (whatever they actually
are) can be differentiated (as if they were functions). We furthermore assume
that these derivatives satisfy the same basic rules (sum rule, product rule,
quotient rule, chain rule) as the derivatives in real analysis. (Again, these
assumptions shall be justified later on.)

Denoting the derivative of a power series $f$ by $f^{\prime}$, we then have%
\[
\left(  1+x+x^{2}+x^{3}+\cdots\right)  ^{\prime}=1+2x+3x^{2}+4x^{3}+\cdots.
\]
Hence,%
\[
1+2x+3x^{2}+4x^{3}+\cdots=\left(  1+x+x^{2}+x^{3}+\cdots\right)  ^{\prime
}=\left(  \dfrac{1}{1-x}\right)  ^{\prime}%
\]
(since (\ref{eq.sec.gf.exas.1.1/1-x}) yields $1+x+x^{2}+x^{3}+\cdots=\dfrac
{1}{1-x}$). Using the quotient rule, we can easily find that $\left(
\dfrac{1}{1-x}\right)  ^{\prime}=\dfrac{1}{\left(  1-x\right)  ^{2}}$, so that%
\begin{equation}
1+2x+3x^{2}+4x^{3}+\cdots=\left(  \dfrac{1}{1-x}\right)  ^{\prime}=\dfrac
{1}{\left(  1-x\right)  ^{2}}. \label{eq.sec.gf.exas.4.1/(1-x)2}%
\end{equation}

The left hand side of this looks very similar to the power series
$0+1x+2x^{2}+3x^{3}+\cdots$ that we want to simplify. And indeed, we have the
following:%
\begin{align}
0+1x+2x^{2}+3x^{3}+\cdots &  =x\underbrace{\left(  1+2x+3x^{2}+4x^{3}%
+\cdots\right)  }_{=\dfrac{1}{\left(  1-x\right)  ^{2}}}=x\cdot\dfrac
{1}{\left(  1-x\right)  ^{2}}\nonumber\\
&  =\dfrac{x}{\left(  1-x\right)  ^{2}}. \label{eq.sec.gf.exas.4.x/(1-x)2}%
\end{align}
\medskip

\textit{Second way:} We rewrite $0+1x+2x^{2}+3x^{3}+\cdots$ as an infinite sum
of infinite sums:%
\begin{align}
&  \ 0+1x+2x^{2}+3x^{3}+\cdots\nonumber\\
&
\begin{array}
[c]{cccccccccc}%
= & x^{1} & + & x^{2} & + & x^{3} & + & x^{4} & + & \cdots\\
&  & + & x^{2} & + & x^{3} & + & x^{4} & + & \cdots\\
&  &  &  & + & x^{3} & + & x^{4} & + & \cdots\\
&  &  &  &  &  & + & x^{4} & + & \cdots\\
&  &  &  &  &  &  &  & \ddots & \ddots
\end{array}
\nonumber\\
&  =\sum_{k\geq1}\underbrace{\left(  x^{k}+x^{k+1}+x^{k+2}+\cdots\right)
}_{\substack{=x^{k}\cdot\left(  1+x+x^{2}+x^{3}+\cdots\right)  \\=x^{k}%
\cdot\dfrac{1}{1-x}\\\text{(by (\ref{eq.sec.gf.exas.1.1/1-x}))}}}=\sum
_{k\geq1}x^{k}\cdot\dfrac{1}{1-x}=\dfrac{1}{1-x}\cdot\underbrace{\sum_{k\geq
1}x^{k}}_{\substack{=x^{1}+x^{2}+x^{3}+\cdots\\=x\cdot\left(  1+x+x^{2}%
+x^{3}+\cdots\right)  \\=x\cdot\dfrac{1}{1-x}\\\text{(by
(\ref{eq.sec.gf.exas.1.1/1-x}))}}}\nonumber\\
&  =\dfrac{1}{1-x}\cdot x\cdot\dfrac{1}{1-x}=\dfrac{x}{\left(  1-x\right)
^{2}}. \label{eq.sec.gf.exas.4.x/(1-x)2alt}%
\end{align}
(We used some unstated assumptions here about infinite sums -- specifically,
we assumed that we can rearrange them without worrying about absolute
convergence or similar issues -- but we will later see that these assumptions
are well justified. Besides, we obtained the same result as by our first way,
which is reassuring.) \medskip

Having computed $0+1x+2x^{2}+3x^{3}+\cdots$, we can now simplify
(\ref{eq.sec.gf.exas.4.Ax=1}), obtaining%
\begin{align*}
A\left(  x\right)   &  =1+2xA\left(  x\right)  +x\underbrace{\left(
0+1x+2x^{2}+3x^{3}+\cdots\right)  }_{=\dfrac{x}{\left(  1-x\right)  ^{2}}}\\
&  =1+2xA\left(  x\right)  +x\cdot\dfrac{x}{\left(  1-x\right)  ^{2}}.
\end{align*}
This is a linear equation in $A\left(  x\right)  $. Solving it yields%
\begin{align*}
A\left(  x\right)   &  =\dfrac{1-2x+2x^{2}}{\left(  1-x\right)  ^{2}\left(
1-2x\right)  }\\
&  =\underbrace{\dfrac{2}{1-2x}}_{\substack{=2\sum_{k\geq0}2^{k}%
x^{k}\\\text{(by (\ref{eq.sec.gf.exas.1.1/1-ax}))}}}-\underbrace{\dfrac
{1}{\left(  1-x\right)  ^{2}}}_{\substack{=1+2x+3x^{2}+4x^{3}+\cdots
\\\text{(by (\ref{eq.sec.gf.exas.4.1/(1-x)2}), or alternatively}\\\text{by
dividing the}\\\text{equality (\ref{eq.sec.gf.exas.4.x/(1-x)2alt}) by
}x\text{)}}}\ \ \ \ \ \ \ \ \ \ \left(  \text{by partial fraction
decomposition}\right) \\
&  =\underbrace{2\sum_{k\geq0}2^{k}x^{k}}_{=\sum_{k\geq0}2^{k+1}x^{k}%
}-\underbrace{\left(  1+2x+3x^{2}+4x^{3}+\cdots\right)  }_{=\sum_{k\geq
0}\left(  k+1\right)  x^{k}}\\
&  =\sum_{k\geq0}2^{k+1}x^{k}-\sum_{k\geq0}\left(  k+1\right)  x^{k}%
=\sum_{k\geq0}\left(  2^{k+1}-\left(  k+1\right)  \right)  x^{k}.
\end{align*}
Comparing coefficients, we obtain%
\[
a_{n}=2^{n+1}-\left(  n+1\right)  \ \ \ \ \ \ \ \ \ \ \text{for each }%
n\in\mathbb{N}.
\]
This is also easy to prove directly (by induction on $n$).

\subsection{\label{sec.gf.defs}Definitions}

The four examples above should have convinced you that generating functions
can be useful. Thus, it is worthwhile to put them on a rigorous footing by
first \textbf{defining} generating functions and then \textbf{justifying} the
manipulations we have been doing to them in the previous section (e.g.,
dividing them, solving quadratic equations, taking infinite sums, taking
derivatives, ...). We are next going to sketch how this can be done (see
\cite[Chapter 7 (in the 1st edition)]{Loehr-BC} and \cite[Chapter 7]{19s} for
some details).

First things first: Generating functions are not actually functions. They are
so-called \emph{formal power series} (short \emph{FPSs}). Roughly speaking, a
formal power series is a \textquotedblleft formal\textquotedblright\ infinite
sum of the form $a_{0}+a_{1}x+a_{2}x^{2}+\cdots$, where $x$ is an
\textquotedblleft indeterminate\textquotedblright\ (we shall soon see what
this all means). You cannot substitute $x=2$ into such a power series. (For
example, substituting $x=2$ into $\dfrac{1}{1-x}=1+x+x^{2}+x^{3}+\cdots$ would
lead to the absurd equality $\dfrac{1}{-1}=1+2+4+8+16+\cdots$.) The word
\textquotedblleft function\textquotedblright\ in \textquotedblleft generating
function\textquotedblright\ is somewhat of a historical artifact.

\subsubsection{\label{subsec.gf.defs.commrings}Reminder: Commutative rings}

In order to obtain a precise understanding of what FPSs are, we go back to
abstract algebra. We begin by recalling the concept of a commutative ring.
This is defined in any textbook on abstract algebra, but we recall the
definition for the sake of completeness.

Informally, a \emph{commutative ring} is a set $K$ equipped with binary
operations $\oplus$, $\ominus$ and $\odot$ and elements $\mathbf{0}$ and
$\mathbf{1}$ that \textquotedblleft behave\textquotedblright\ like addition,
subtraction and multiplication (of numbers) and the numbers $0$ and $1$,
respectively. For example, they should satisfy rules like $\left(  a\oplus
b\right)  \odot c=\left(  a\odot c\right)  \oplus\left(  b\odot c\right)  $.

Formally, commutative rings are defined as follows:

\begin{definition}
\label{def.alg.commring}A \emph{commutative ring} means a set $K$ equipped
with three maps%
\begin{align*}
\oplus &  :K\times K\rightarrow K,\\
\ominus &  :K\times K\rightarrow K,\\
\odot &  :K\times K\rightarrow K
\end{align*}
and two elements $\mathbf{0}\in K$ and $\mathbf{1}\in K$ satisfying the
following axioms:

\begin{enumerate}
\item \emph{Commutativity of addition:} We have $a\oplus b=b\oplus a$ for all
$a,b\in K$.

(Here and in the following, we write the three maps $\oplus$, $\ominus$ and
$\odot$ infix -- i.e., we denote the image of a pair $\left(  a,b\right)  \in
K\times K$ under the map $\oplus$ by $a\oplus b$ rather than by $\oplus\left(
a,b\right)  $.)

\item \emph{Associativity of addition:} We have $a\oplus\left(  b\oplus
c\right)  =\left(  a\oplus b\right)  \oplus c$ for all $a,b,c\in K$.

\item \emph{Neutrality of zero:} We have $a\oplus\mathbf{0}=\mathbf{0}\oplus
a=a$ for all $a\in K$.

\item \emph{Subtraction undoes addition:} Let $a,b,c\in K$. We have $a\oplus
b=c$ if and only if $a=c\ominus b$.

\item \emph{Commutativity of multiplication:} We have $a\odot b=b\odot a$ for
all $a,b\in K$.

\item \emph{Associativity of multiplication:} We have $a\odot\left(  b\odot
c\right)  =\left(  a\odot b\right)  \odot c$ for all $a,b,c\in K$.

\item \emph{Distributivity:} We have%
\[
a\odot\left(  b\oplus c\right)  =\left(  a\odot b\right)  \oplus\left(  a\odot
c\right)  \ \ \ \ \ \ \ \ \ \ \text{and}\ \ \ \ \ \ \ \ \ \ \left(  a\oplus
b\right)  \odot c=\left(  a\odot c\right)  \oplus\left(  b\odot c\right)
\]
for all $a,b,c\in K$.

\item \emph{Neutrality of one:} We have $a\odot\mathbf{1}=\mathbf{1}\odot a=a$
for all $a\in K$.

\item \emph{Annihilation:} We have $a\odot\mathbf{0}=\mathbf{0}\odot
a=\mathbf{0}$ for all $a\in K$.
\end{enumerate}

[\textbf{Note:} Most authors do not include $\ominus$ in the definition of a
commutative ring, but instead require the existence of additive inverses for
all $a\in K$. This is equivalent to our definition, because if additive
inverses exist, then we can define $a\ominus b$ to be $a\oplus\overline{b}$
where $\overline{b}$ is the additive inverse of $b$.]

The operations $\oplus$, $\ominus$ and $\odot$ are called the \emph{addition},
the \emph{subtraction} and the \emph{multiplication} of the ring $K$. This
does not imply that they have any connection with the usual addition,
subtraction and multiplication of numbers; it merely means that they play
similar roles to the latter and behave similarly. When confusion is unlikely,
we will denote these three operations $\oplus$, $\ominus$ and $\odot$ by $+$,
$-$ and $\cdot$, respectively, and we will abbreviate $a\odot b=a\cdot b$ by
$ab$.

The elements $\mathbf{0}$ and $\mathbf{1}$ are called the \emph{zero} and the
\emph{unity} (or the \emph{one}) of the ring $K$. Again, this does not imply
that they equal the numbers $0$ and $1$, but merely that they play analogous
roles. We will simply call these elements $0$ and $1$ when confusion with the
corresponding numbers is unlikely.

We will use \emph{PEMDAS conventions} for the three operations $\oplus$,
$\ominus$ and $\odot$. These imply that the operation $\odot$ has higher
precedence than $\oplus$ and $\ominus$, while the operations $\oplus$ and
$\ominus$ are left-associative. Thus, for example, \textquotedblleft%
$ab+ac$\textquotedblright\ means $\left(  ab\right)  +\left(  ac\right)  $
(that is, $\left(  a\odot b\right)  \oplus\left(  a\odot c\right)  $).
Likewise, \textquotedblleft$a-b+c$\textquotedblright\ means $\left(
a-b\right)  +c=\left(  a\ominus b\right)  \oplus c$.
\end{definition}

Here are some examples of commutative rings:

\begin{itemize}
\item The sets $\mathbb{Z}$, $\mathbb{Q}$, $\mathbb{R}$ and $\mathbb{C}$ are
commutative rings. (Of course, the operations $\oplus$, $\ominus$ and $\odot$
of these rings are just the usual operations $+$, $-$ and $\cdot$ known from
high school.)

\item The set $\mathbb{N}$ is not a commutative ring, since it has no
subtraction. (It is, however, something called a \emph{commutative semiring}.)

\item The matrix ring $\mathbb{Q}^{m\times m}$ (this is the ring of all
$m\times m$-matrices with rational entries) is not a commutative ring for
$m>1$ (because it fails the \textquotedblleft commutativity of
multiplication\textquotedblright\ axiom). However, it satisfies all axioms
other than \textquotedblleft commutativity of multiplication\textquotedblright%
. This makes it a \emph{noncommutative ring}.

\item The set
\[
\mathbb{Z}\left[  \sqrt{5}\right]  =\left\{  a+b\sqrt{5}\ \mid\ a,b\in
\mathbb{Z}\right\}
\]
is a commutative ring with operations $+$, $-$ and $\cdot$ inherited from
$\mathbb{R}$. This is because any $a,b,c,d\in\mathbb{Z}$ satisfy%
\begin{align*}
\left(  a+b\sqrt{5}\right)  +\left(  c+d\sqrt{5}\right)   &  =\left(
a+c\right)  +\left(  b+d\right)  \sqrt{5}\in\mathbb{Z}\left[  \sqrt{5}\right]
;\\
\left(  a+b\sqrt{5}\right)  -\left(  c+d\sqrt{5}\right)   &  =\left(
a-c\right)  +\left(  b-d\right)  \sqrt{5}\in\mathbb{Z}\left[  \sqrt{5}\right]
;\\
\left(  a+b\sqrt{5}\right)  \left(  c+d\sqrt{5}\right)   &  =\left(
ac+5bd\right)  +\left(  ad+bc\right)  \sqrt{5}\in\mathbb{Z}\left[  \sqrt
{5}\right]  .
\end{align*}
This is called a \emph{subring} of $\mathbb{R}$ (i.e., a subset of
$\mathbb{R}$ that contains $0$ and $1$ and is closed under the operations $+$,
$-$ and $\cdot$, and therefore constitutes a commutative ring with these
operations inherited from $\mathbb{R}$).

\item For each $m\in\mathbb{Z}$, the set
\begin{align*}
\mathbb{Z}/m  &  =\left\{  \text{all residue classes modulo }m\right\} \\
&  =\left\{
\begin{array}
[c]{c}%
\text{equivalence classes of integers with respect to}\\
\text{the equivalence \textquotedblleft congruent modulo }%
m\text{\textquotedblright}\\
\text{(that is, \textquotedblleft differ by a multiple of }%
m\text{\textquotedblright)}%
\end{array}
\right\}
\end{align*}
is a commutative ring, with its operations defined by%
\[
\overline{a}+\overline{b}=\overline{a+b},\ \ \ \ \ \ \ \ \ \ \overline
{a}-\overline{b}=\overline{a-b},\ \ \ \ \ \ \ \ \ \ \overline{a}\cdot
\overline{b}=\overline{ab}.
\]
If $m>0$, then this ring $\mathbb{Z}/m$ is finite and has size $m$. It is also
known as $\mathbb{Z}/m\mathbb{Z}$ or $\mathbb{Z}_{m}$ (careful with the latter
notation; it can mean different things to different people). When $m$ is
prime, the ring $\mathbb{Z}/m$ is actually a finite field and is called
$\mathbb{F}_{m}$. (But, e.g., the ring $\mathbb{Z}/4$ is not a field, and not
the same as $\mathbb{F}_{4}$.)

\item In the examples we have seen so far, the elements of the commutative
ring either are numbers or (as in the case of matrices or residue classes)
consist of numbers. For a contrast, here is an example where they are sets:

For any two sets $X$ and $Y$, we define the \emph{symmetric difference}
$X\bigtriangleup Y$ of $X$ and $Y$ to be the set%
\begin{align*}
\left(  X\cup Y\right)  \setminus\left(  X\cap Y\right)   &  =\left(
X\setminus Y\right)  \cup\left(  Y\setminus X\right) \\
&  =\left\{  \text{all elements that belong to exactly one of }X\text{ and
}Y\right\}  .
\end{align*}

Fix a set $S$. Consider the power set $\mathcal{P}\left(  S\right)  $ of $S$
(that is, the set of all subsets of $S$). This power set $\mathcal{P}\left(
S\right)  $ is a commutative ring if we equip it with the operation
$\bigtriangleup$ as addition (that is, $X\oplus Y=X\bigtriangleup Y$ for any
subsets $X$ and $Y$ of $S$), with the same operation $\bigtriangleup$ as
subtraction (that is, $X\ominus Y=X\bigtriangleup Y$), and with the operation
$\cap$ as multiplication (that is, $X\odot Y=X\cap Y$), and with the elements
$\varnothing$ and $S$ as zero and unity (that is, with $\mathbf{0}%
=\varnothing$ and $\mathbf{1}=S$). Indeed, it is straightforward to see that
all the axioms in Definition \ref{def.alg.commring} hold for this ring. (For
example, distributivity holds because any three sets $A,B,C$ satisfy
$A\cap\left(  B\bigtriangleup C\right)  =\left(  A\cap B\right)
\bigtriangleup\left(  A\cap C\right)  $ and $\left(  A\bigtriangleup B\right)
\cap C=\left(  A\cap C\right)  \bigtriangleup\left(  B\cap C\right)  $.) This
is an example of a \emph{Boolean ring} (i.e., a ring in which $aa=a$ for each
element $a$ of the ring).

\item Here is another example of a semiring, which is rather useful in
combinatorics. Let $\mathbb{T}$ be the set $\mathbb{Z}\cup\left\{
-\infty\right\}  $, where $-\infty$ is just some extra symbol. Define two
operations $\oplus$ and $\odot$ on this set $\mathbb{T}$ by setting%
\[
a\oplus b=\max\left\{  a,b\right\}  \ \ \ \ \ \ \ \ \ \ \left(  \text{where
}\max\left\{  n,-\infty\right\}  :=n\text{ for each }n\in\mathbb{T}\right)
\]
and%
\[
a\odot b=a+b\ \ \ \ \ \ \ \ \ \ \left(  \text{where }n+\left(  -\infty\right)
:=\left(  -\infty\right)  +n:=-\infty\text{ for each }n\in\mathbb{T}\right)
.
\]
Then, $\mathbb{T}$ is a commutative semiring (i.e., it would satisfy
Definition \ref{def.alg.commring} if not for the lack of subtraction). It is
called \href{https://en.wikipedia.org/wiki/Tropical_semiring}{the
\emph{tropical semiring}} of $\mathbb{Z}$.
\end{itemize}

More examples can be found in algebra textbooks (or in \cite[\S 6.1]{detnotes}
or \cite[\S 5.2]{19s} or \cite[\S 2.1.2, \S 2.3.2, \S 2.3.3]{23wa}). \medskip

\textbf{Good news:} In any commutative ring $K$, the standard rules of
computation apply:

\begin{itemize}
\item You can compute finite sums (of elements of $K$) without specifying the
order of summation or the placement of parentheses. For example, for any
$a,b,c,d,e\in K$, we have%
\[
\left(  \left(  a+b\right)  +\left(  c+d\right)  \right)  +e=\left(  a+\left(
b+c\right)  \right)  +\left(  d+e\right)  ,
\]
so you can write the sum $a+b+c+d+e$ without putting parentheses around
anything. (This is called \textquotedblleft general
associativity\textquotedblright.)

Also, finite sums do not depend on the order of addends. For example, for any
$a,b,c,d,e\in K$, we have%
\[
a+b+c+d+e=d+b+a+e+c.
\]
(This is called \textquotedblleft general commutativity\textquotedblright.)

More formally: If $\left(  a_{s}\right)  _{s\in S}$ is any finite family of
elements of a commutative ring $K$ (this means that $S$ is a finite set, and
$a_{s}$ is an element of $K$ for each $s\in S$), then the finite sum%
\[
\sum\limits_{s\in S}a_{s}%
\]
is a well-defined element of $K$. Furthermore, such sums satisfy the usual
rules of sums (see \cite[\S 1.4.2]{detnotes}); for instance:

\begin{itemize}
\item If $S=X\cup Y$ and $X\cap Y=\varnothing$, then $\sum\limits_{s\in
S}a_{s}=\sum\limits_{s\in X}a_{s}+\sum\limits_{s\in Y}a_{s}$.

\item We have $\sum\limits_{s\in S}\left(  a_{s}+b_{s}\right)  =\sum
\limits_{s\in S}a_{s}+\sum\limits_{s\in S}b_{s}$.

\item If $S$ and $W$ are two finite sets, and if $f:S\rightarrow W$ is a map,
then $\sum_{s\in S}a_{s}=\sum_{w\in W}\sum_{\substack{s\in S;\\f\left(
s\right)  =w}}a_{s}$. (That is, you can subdivide a finite sum into a finite
sum of finite sums by bunching its addends arbitrarily.)
\end{itemize}

(See \cite[\S 2.14]{detnotes} for very detailed proofs\footnote{These proofs
are stated for numbers rather than elements of an arbitrary commutative ring
$K$, but the exact same reasoning works in an arbitrary ring $K$.}.
Alternatively, you can treat them as exercises on induction.)

If $S=\varnothing$, then $\sum\limits_{s\in S}a_{s}=0$ by definition. Such a
sum $\sum\limits_{s\in S}a_{s}$ with $S=\varnothing$ is called an empty sum.

\item The same holds for finite products. If $S=\varnothing$, then
$\prod\limits_{s\in S}a_{s}=1$ by definition.

\item If $a\in K$, then $-a$ denotes $0-a=\mathbf{0}-a\in K$.

\item If $n\in\mathbb{Z}$ and $a\in K$, then we can define an element $na\in
K$ by%
\[
na=%
\begin{cases}
\underbrace{a+a+\cdots+a}_{n\text{ times}}, & \text{if }n\geq0;\\
-\left(  \underbrace{a+a+\cdots+a}_{-n\text{ times}}\right)  , & \text{if
}n<0.
\end{cases}
\]
This generalizes the classical definition of multiplication (for integers) as
repeated addition.

\item If $n\in\mathbb{N}$ and $a\in K$, then we can define an element%
\[
a^{n}=\underbrace{aa\cdots a}_{n\text{ times}}\in K.
\]
In particular, $a^{0}=\underbrace{aa\cdots a}_{0\text{ times}}=1$ (since an
empty product is always $1$).

\item Standard rules hold:%
\begin{align*}
-\left(  a+b\right)   &  =\left(  -a\right)  +\left(  -b\right)
\ \ \ \ \ \ \ \ \ \ \text{for any }a,b\in K;\\
-\left(  -a\right)   &  =a\ \ \ \ \ \ \ \ \ \ \text{for any }a\in K;\\
\left(  n+m\right)  a  &  =na+ma\ \ \ \ \ \ \ \ \ \ \text{for any }a\in
K\text{ and }n,m\in\mathbb{Z};\\
\left(  nm\right)  a  &  =n\left(  ma\right)  \ \ \ \ \ \ \ \ \ \ \text{for
any }a\in K\text{ and }n,m\in\mathbb{Z};\\
a\left(  b-c\right)   &  =\left(  ab\right)  -\left(  ac\right)
\ \ \ \ \ \ \ \ \ \ \text{for any }a,b,c\in K;\\
\left(  ab\right)  ^{n}  &  =a^{n}b^{n}\ \ \ \ \ \ \ \ \ \ \text{for any
}a,b\in K\text{ and }n\in\mathbb{N};\\
a^{n+m}  &  =a^{n}a^{m}\ \ \ \ \ \ \ \ \ \ \text{for any }a\in K\text{ and
}n,m\in\mathbb{N};\\
a^{nm}  &  =\left(  a^{n}\right)  ^{m}\ \ \ \ \ \ \ \ \ \ \text{for any }a\in
K\text{ and }n,m\in\mathbb{N};\\
&  \ldots;\\
\left(  a+b\right)  ^{n}  &  =\sum_{k=0}^{n}\dbinom{n}{k}a^{k}b^{n-k}%
\ \ \ \ \ \ \ \ \ \ \text{for any }a,b\in K\text{ and }n\in\mathbb{N}.
\end{align*}
(The latter equality is known as the \emph{binomial theorem} or \emph{binomial
formula}.)
\end{itemize}

A further concept will be useful. Namely, if $K$ is a commutative ring, then
the notion of a $K$\emph{-module} is the straightforward generalization of the
concept of a $K$-vector space to cases where $K$ is not a field (but just a
commutative ring). Here is the definition of a $K$-module in detail:

\begin{definition}
\label{def.alg.module}Let $K$ be a commutative ring.

A $K$\emph{-module} means a set $M$ equipped with three maps
\begin{align*}
\oplus &  :M\times M\rightarrow M,\\
\ominus &  :M\times M\rightarrow M,\\
\rightharpoonup &  :K\times M\rightarrow M
\end{align*}
(notice that the third map has domain $K\times M$, not $M\times M$) and an
element $\overrightarrow{0}\in M$ satisfying the following axioms:

\begin{enumerate}
\item \emph{Commutativity of addition:} We have $a\oplus b=b\oplus a$ for all
$a,b\in M$.

(Here and in the following, we write the three maps $\oplus$, $\ominus$ and
$\rightharpoonup$ infix, just as for a commutative ring.)

\item \emph{Associativity of addition:} We have $a\oplus\left(  b\oplus
c\right)  =\left(  a\oplus b\right)  \oplus c$ for all $a,b,c\in M$.

\item \emph{Neutrality of zero:} We have $a\oplus\overrightarrow{0}%
=\overrightarrow{0}\oplus a=a$ for all $a\in M$.

\item \emph{Subtraction undoes addition:} Let $a,b,c\in M$. We have $a\oplus
b=c$ if and only if $a=c\ominus b$.

\item \emph{Associativity of scaling:} We have $u\rightharpoonup\left(
v\rightharpoonup a\right)  =\left(  uv\right)  \rightharpoonup a$ for all
$u,v\in K$ and $a\in M$.

\item \emph{Left distributivity:} We have $u\rightharpoonup\left(  a\oplus
b\right)  =\left(  u\rightharpoonup a\right)  \oplus\left(  u\rightharpoonup
b\right)  $ for all $u\in K$ and $a,b\in M$.

\item \emph{Right distributivity:} We have $\left(  u+v\right)
\rightharpoonup a=\left(  u\rightharpoonup a\right)  \oplus\left(
v\rightharpoonup a\right)  $ for all $u,v\in K$ and $a\in M$.

\item \emph{Neutrality of one:} We have $1\rightharpoonup a=a$ for all $a\in
M$.

\item \emph{Left annihilation:} We have $0\rightharpoonup a=\overrightarrow{0}%
$ for all $a\in M$.

\item \emph{Right annihilation:} We have $u\rightharpoonup\overrightarrow{0}%
=\overrightarrow{0}$ for all $u\in K$.
\end{enumerate}

[\textbf{Note:} Most authors do not include $\ominus$ in the definition of a
$K$-module. Instead, they may require the existence of additive inverses for
all $a\in M$. Just like for commutative rings, this is equivalent. Actually,
even the existence of additive inverses is not necessary, since the additive
inverse of an $a\in M$ can be constructed using the operation $\rightharpoonup
$ as $\left(  -1\right)  \rightharpoonup a$. Thus, you will sometimes see the
notion of a $K$-module defined entirely without any reference to subtraction
or additive inverses; it is still equivalent to ours.]

The operations $\oplus$, $\ominus$ and $\rightharpoonup$ are called the
\emph{addition}, the \emph{subtraction} and the \emph{scaling} (or the
$K$\emph{-action}) of the $K$-module $M$. When confusion is unlikely, we will
denote these three operations $\oplus$, $\ominus$ and $\rightharpoonup$ by
$+$, $-$ and $\cdot$, respectively, and we will abbreviate $a\rightharpoonup
b=a\cdot b$ by $ab$.

The element $\overrightarrow{0}$ is called the \emph{zero} (or the \emph{zero
vector}) of the $K$-module $M$. We will usually just call it $0$.

When $M$ is a $K$-module, the elements of $M$ are called \emph{vectors}, while
the elements of $K$ are called \emph{scalars}.

We will use \emph{PEMDAS conventions} for the three operations $\oplus$,
$\ominus$ and $\rightharpoonup$, with the operation $\rightharpoonup$ having
higher precedence than $\oplus$ and $\ominus$.
\end{definition}

This all having been said, we can now define formal power series.

\subsubsection{\label{subsec.gf.defs.fps}The definition of formal power
series}

Until the end of Chapter \ref{chap.gf}, the following convention will be in place:

\begin{convention}
Fix a commutative ring $K$. (For example, $K$ can be $\mathbb{Z}$ or
$\mathbb{Q}$ or $\mathbb{C}$.)
\end{convention}

\begin{definition}
\label{def.fps.fps}A \emph{formal power series} (or, short, \emph{FPS}) in $1$
indeterminate over $K$ means a sequence $\left(  a_{0},a_{1},a_{2}%
,\ldots\right)  =\left(  a_{n}\right)  _{n\in\mathbb{N}}\in K^{\mathbb{N}}$ of
elements of $K$.
\end{definition}

Examples of FPSs over $\mathbb{Z}$ are $\left(  0,0,0,\ldots\right)  $ and
$\left(  1,0,0,0,\ldots\right)  $ and $\left(  1,1,1,1,\ldots\right)  $ and
$\left(  1,2,3,4,\ldots\right)  $.

Definition \ref{def.fps.fps} technically answers the question
\textquotedblleft what is an FPS\textquotedblright; however, the questions
\textquotedblleft what can we do with an FPS\textquotedblright\ and
\textquotedblleft why do the examples in Section \ref{sec.gf.exas}
work\textquotedblright\ or \textquotedblleft what is $x$\textquotedblright%
\ remain open. These questions will take us a while.

First, let us define some operations on FPSs.

\begin{definition}
\label{def.fps.ops}\textbf{(a)} The \emph{sum} of two FPSs $\mathbf{a}=\left(
a_{0},a_{1},a_{2},\ldots\right)  $ and $\mathbf{b}=\left(  b_{0},b_{1}%
,b_{2},\ldots\right)  $ is defined to be the FPS%
\[
\left(  a_{0}+b_{0},\ \ a_{1}+b_{1},\ \ a_{2}+b_{2},\ \ \ldots\right)  .
\]
It is denoted by $\mathbf{a}+\mathbf{b}$. \medskip

\textbf{(b)} The \emph{difference} of two FPSs $\mathbf{a}=\left(  a_{0}%
,a_{1},a_{2},\ldots\right)  $ and $\mathbf{b}=\left(  b_{0},b_{1},b_{2}%
,\ldots\right)  $ is defined to be the FPS%
\[
\left(  a_{0}-b_{0},\ \ a_{1}-b_{1},\ \ a_{2}-b_{2},\ \ \ldots\right)  .
\]
It is denoted by $\mathbf{a}-\mathbf{b}$. \medskip

\textbf{(c)} If $\lambda\in K$ and if $\mathbf{a}=\left(  a_{0},a_{1}%
,a_{2},\ldots\right)  $ is an FPS, then we define an FPS
\[
\lambda\mathbf{a}:=\left(  \lambda a_{0},\lambda a_{1},\lambda a_{2}%
,\ldots\right)  .
\]

\textbf{(d)} The \emph{product} of two FPSs $\mathbf{a}=\left(  a_{0}%
,a_{1},a_{2},\ldots\right)  $ and $\mathbf{b}=\left(  b_{0},b_{1},b_{2}%
,\ldots\right)  $ is defined to be the FPS $\left(  c_{0},c_{1},c_{2}%
,\ldots\right)  $, where
\begin{align*}
c_{n}  &  =\sum_{i=0}^{n}a_{i}b_{n-i}=\sum_{\substack{\left(  i,j\right)
\in\mathbb{N}^{2};\\i+j=n}}a_{i}b_{j}\\
&  =a_{0}b_{n}+a_{1}b_{n-1}+a_{2}b_{n-2}+\cdots+a_{n}b_{0}%
\ \ \ \ \ \ \ \ \ \ \text{for each }n\in\mathbb{N}.
\end{align*}
This product is denoted by $\mathbf{a}\cdot\mathbf{b}$ or just by
$\mathbf{ab}$. \medskip

\textbf{(e)} For each $a\in K$, we define $\underline{a}$ to be the FPS
$\left(  a,0,0,0,\ldots\right)  $. An FPS of the form $\underline{a}$ for some
$a\in K$ (that is, an FPS $\left(  a_{0},a_{1},a_{2},\ldots\right)  $
satisfying $a_{1}=a_{2}=a_{3}=\cdots=0$) is said to be \emph{constant}.
\medskip

\textbf{(f)} The set of all FPSs (in $1$ indeterminate over $K$) is denoted
$K\left[  \left[  x\right]  \right]  $.
\end{definition}

The following theorem is crucial: essentially it says that the operations on
FPSs that we have just defined behave as such operations should:

\begin{theorem}
\label{thm.fps.ring}\textbf{(a)} The set $K\left[  \left[  x\right]  \right]
$ is a commutative ring (with its operations $+$, $-$ and $\cdot$ defined in
Definition \ref{def.fps.ops}). Its zero and its unity are the FPSs
$\underline{0}=\left(  0,0,0,\ldots\right)  $ and $\underline{1}=\left(
1,0,0,0,\ldots\right)  $. This means, concretely, that the following facts hold:

\begin{enumerate}
\item \emph{Commutativity of addition:} We have $\mathbf{a}+\mathbf{b}%
=\mathbf{b}+\mathbf{a}$ for all $\mathbf{a},\mathbf{b}\in K\left[  \left[
x\right]  \right]  $.

\item \emph{Associativity of addition:} We have $\mathbf{a}+\left(
\mathbf{b}+\mathbf{c}\right)  =\left(  \mathbf{a}+\mathbf{b}\right)
+\mathbf{c}$ for all $\mathbf{a},\mathbf{b},\mathbf{c}\in K\left[  \left[
x\right]  \right]  $.

\item \emph{Neutrality of zero:} We have $\mathbf{a}+\underline{0}%
=\underline{0}+\mathbf{a}=\mathbf{a}$ for all $\mathbf{a}\in K\left[  \left[
x\right]  \right]  $.

\item \emph{Subtraction undoes addition:} Let $\mathbf{a},\mathbf{b}%
,\mathbf{c}\in K\left[  \left[  x\right]  \right]  $. We have $\mathbf{a}%
+\mathbf{b}=\mathbf{c}$ if and only if $\mathbf{a}=\mathbf{c}-\mathbf{b}$.

\item \emph{Commutativity of multiplication:} We have $\mathbf{ab}%
=\mathbf{ba}$ for all $\mathbf{a},\mathbf{b}\in K\left[  \left[  x\right]
\right]  $.

\item \emph{Associativity of multiplication:} We have $\mathbf{a}\left(
\mathbf{bc}\right)  =\left(  \mathbf{ab}\right)  \mathbf{c}$ for all
$\mathbf{a},\mathbf{b},\mathbf{c}\in K\left[  \left[  x\right]  \right]  $.

\item \emph{Distributivity:} We have%
\[
\mathbf{a}\left(  \mathbf{b}+\mathbf{c}\right)  =\mathbf{ab}+\mathbf{ac}%
\ \ \ \ \ \ \ \ \ \ \text{and}\ \ \ \ \ \ \ \ \ \ \left(  \mathbf{a}%
+\mathbf{b}\right)  \mathbf{c}=\mathbf{ac}+\mathbf{bc}%
\]
for all $\mathbf{a},\mathbf{b},\mathbf{c}\in K\left[  \left[  x\right]
\right]  $.

\item \emph{Neutrality of one:} We have $\mathbf{a}\cdot\underline{1}%
=\underline{1}\cdot\mathbf{a}=\mathbf{a}$ for all $\mathbf{a}\in K\left[
\left[  x\right]  \right]  $.

\item \emph{Annihilation:} We have $\mathbf{a}\cdot\underline{0}%
=\underline{0}\cdot\mathbf{a}=\underline{0}$ for all $\mathbf{a}\in K\left[
\left[  x\right]  \right]  $.
\end{enumerate}

\textbf{(b)} Furthermore, $K\left[  \left[  x\right]  \right]  $ is a
$K$-module (with its scaling being the map that sends each $\left(
\lambda,\mathbf{a}\right)  \in K\times K\left[  \left[  x\right]  \right]  $
to the FPS $\lambda\mathbf{a}$ defined in Definition \ref{def.fps.ops}
\textbf{(c)}). Its zero vector is $\underline{0}$. Concretely, this means that:

\begin{enumerate}
\item[10.] \emph{Associativity of scaling:} We have $\lambda\left(
\mu\mathbf{a}\right)  =\left(  \lambda\mu\right)  \mathbf{a}$ for all
$\lambda,\mu\in K$ and $\mathbf{a}\in K\left[  \left[  x\right]  \right]  $.

\item[11.] \emph{Left distributivity:} We have $\lambda\left(  \mathbf{a}%
+\mathbf{b}\right)  =\lambda\mathbf{a}+\lambda\mathbf{b}$ for all $\lambda\in
K$ and $\mathbf{a},\mathbf{b}\in K\left[  \left[  x\right]  \right]  $.

\item[12.] \emph{Right distributivity:} We have $\left(  \lambda+\mu\right)
\mathbf{a}=\lambda\mathbf{a}+\mu\mathbf{a}$ for all $\lambda,\mu\in K$ and
$\mathbf{a}\in K\left[  \left[  x\right]  \right]  $.

\item[13.] \emph{Neutrality of one:} We have $1\mathbf{a}=\mathbf{a}$ for all
$\mathbf{a}\in K\left[  \left[  x\right]  \right]  $.

\item[14.] \emph{Left annihilation:} We have $0\mathbf{a}=\underline{0}$ for
all $\mathbf{a}\in K\left[  \left[  x\right]  \right]  $.

\item[15.] \emph{Right annihilation:} We have $\lambda\underline{0}%
=\underline{0}$ for all $\lambda\in K$.
\end{enumerate}

(Also, some of the facts from part \textbf{(a)} are included in this
statement.) \medskip

\textbf{(c)} We have $\lambda\left(  \mathbf{a}\cdot\mathbf{b}\right)
=\left(  \lambda\mathbf{a}\right)  \cdot\mathbf{b}=\mathbf{a}\cdot\left(
\lambda\mathbf{b}\right)  $ for all $\lambda\in K$ and $\mathbf{a}%
,\mathbf{b}\in K\left[  \left[  x\right]  \right]  $. \medskip

\textbf{(d)} Finally, we have $\lambda\mathbf{a}=\underline{\lambda}%
\cdot\mathbf{a}$ for all $\lambda\in K$ and $\mathbf{a}\in K\left[  \left[
x\right]  \right]  $.
\end{theorem}

Theorem \ref{thm.fps.ring} allows us to calculate with FPSs as we do with
numbers, at least as far as the operations $+$, $-$ and $\cdot$ are concerned.
Hence, e.g., we know that:

\begin{itemize}
\item Sums and products in $K\left[  \left[  x\right]  \right]  $ need no
parentheses and do not depend on the order of addends/factors. For example,
for any $\mathbf{a},\mathbf{b},\mathbf{c},\mathbf{d}\in K\left[  \left[
x\right]  \right]  $, we have $\left(  \left(  \mathbf{ab}\right)
\mathbf{c}\right)  \mathbf{d}=\mathbf{a}\left(  \left(  \mathbf{bc}\right)
\mathbf{d}\right)  =\left(  \mathbf{ab}\right)  \left(  \mathbf{cd}\right)  $,
so that we can write $\mathbf{abcd}$ for each of these products; and moreover,
we have $\mathbf{abcd}=\mathbf{bdac}=\mathbf{dacb}$.

\item Finite sums and products (such as $\sum_{i=1}^{k}\mathbf{a}_{i}$ or
$\sum_{i\in I}\mathbf{a}_{i}$ or $\prod_{i=1}^{k}\mathbf{a}_{i}$ or
$\prod_{i\in I}\mathbf{a}_{i}$, where $I$ is a finite set) make sense and
behave as one would expect.

\item Powers exist: that is, you can take $\mathbf{a}^{n}$ for each FPS
$\mathbf{a}$ and each $n\in\mathbb{N}$.

\item Standard rules hold: e.g., we have $\mathbf{a}^{n+m}=\mathbf{a}%
^{n}\mathbf{a}^{m}$ and $\left(  \mathbf{ab}\right)  ^{n}=\mathbf{a}%
^{n}\mathbf{b}^{n}$ for any $\mathbf{a},\mathbf{b}\in K\left[  \left[
x\right]  \right]  $ and any $n,m\in\mathbb{N}$.

\item The binomial theorem holds: For any $\mathbf{a},\mathbf{b}\in K\left[
\left[  x\right]  \right]  $ and any $n\in\mathbb{N}$, we have%
\[
\left(  \mathbf{a}+\mathbf{b}\right)  ^{n}=\sum_{k=0}^{n}\dbinom{n}%
{k}\mathbf{a}^{k}\mathbf{b}^{n-k}.
\]

\end{itemize}

Before we prove Theorem \ref{thm.fps.ring}, we introduce one more notation:

\begin{definition}
\label{def.fps.coeff}If $n\in\mathbb{N}$, and if $\mathbf{a}=\left(
a_{0},a_{1},a_{2},\ldots\right)  \in K\left[  \left[  x\right]  \right]  $ is
an FPS, then we define an element $\left[  x^{n}\right]  \mathbf{a}\in K$ by%
\[
\left[  x^{n}\right]  \mathbf{a}:=a_{n}.
\]
This is called the \emph{coefficient of }$x^{n}$\emph{ in }$\mathbf{a}$, or
the $n$\emph{-th coefficient} of $\mathbf{a}$, or the $x^{n}$%
\emph{-coefficient} of $\mathbf{a}$.
\end{definition}

Thus, the definition of the sum of two FPSs (Definition \ref{def.fps.ops}
\textbf{(a)}) rewrites as follows: For any $\mathbf{a},\mathbf{b}\in K\left[
\left[  x\right]  \right]  $ and any $n\in\mathbb{N}$, we have
\begin{equation}
\left[  x^{n}\right]  \left(  \mathbf{a}+\mathbf{b}\right)  =\left[
x^{n}\right]  \mathbf{a}+\left[  x^{n}\right]  \mathbf{b}.
\label{pf.thm.fps.ring.xn(a+b)=}%
\end{equation}
Similarly, for any $\mathbf{a},\mathbf{b}\in K\left[  \left[  x\right]
\right]  $ and any $n\in\mathbb{N}$, we have
\begin{equation}
\left[  x^{n}\right]  \left(  \mathbf{a}-\mathbf{b}\right)  =\left[
x^{n}\right]  \mathbf{a}-\left[  x^{n}\right]  \mathbf{b}.
\label{pf.thm.fps.ring.xn(a-b)=}%
\end{equation}
Meanwhile, the definition of the product of two FPSs (Definition
\ref{def.fps.ops} \textbf{(d)}) rewrites as follows: For any $\mathbf{a}%
,\mathbf{b}\in K\left[  \left[  x\right]  \right]  $ and any $n\in\mathbb{N}$,
we have
\begin{align}
&  \left[  x^{n}\right]  \left(  \mathbf{ab}\right) \nonumber\\
&  =\left[  x^{0}\right]  \mathbf{a}\cdot\left[  x^{n}\right]  \mathbf{b}%
+\left[  x^{1}\right]  \mathbf{a}\cdot\left[  x^{n-1}\right]  \mathbf{b}%
+\left[  x^{2}\right]  \mathbf{a}\cdot\left[  x^{n-2}\right]  \mathbf{b}%
+\cdots+\left[  x^{n}\right]  \mathbf{a}\cdot\left[  x^{0}\right]
\mathbf{b}\nonumber\\
&  =\sum_{i=0}^{n}\left[  x^{i}\right]  \mathbf{a}\cdot\left[  x^{n-i}\right]
\mathbf{b}\label{pf.thm.fps.ring.xn(ab)=2}\\
&  =\sum_{j=0}^{n}\left[  x^{n-j}\right]  \mathbf{a}\cdot\left[  x^{j}\right]
\mathbf{b} \label{pf.thm.fps.ring.xn(ab)=3}%
\end{align}
(here, we have substituted $n-j$ for $i$ in the sum). For $n=0$, this equality
simplifies to%
\begin{equation}
\left[  x^{0}\right]  \left(  \mathbf{ab}\right)  =\left[  x^{0}\right]
\mathbf{a}\cdot\left[  x^{0}\right]  \mathbf{b}.
\label{pf.thm.fps.ring.x0(ab)=}%
\end{equation}
In other words, when we multiply two FPSs, their constant terms get
multiplied. Here and in the following, the \emph{constant term} of an FPS
$\mathbf{a}\in K\left[  \left[  x\right]  \right]  $ is defined to be its
$0$-th coefficient $\left[  x^{0}\right]  \mathbf{a}$.

Finally, Definition \ref{def.fps.ops} \textbf{(c)} rewrites as follows: For
any $\lambda\in K$ and $\mathbf{a}\in K\left[  \left[  x\right]  \right]  $
and any $n\in\mathbb{N}$, we have%
\begin{equation}
\left[  x^{n}\right]  \left(  \lambda\mathbf{a}\right)  =\lambda\cdot\left[
x^{n}\right]  \mathbf{a}. \label{pf.thm.fps.ring.xn(la)=}%
\end{equation}

\begin{proof}
[Proof of Theorem \ref{thm.fps.ring}.]Most parts of Theorem \ref{thm.fps.ring}
are straightforward to verify. Let us check the associativity of multiplication.

Let $\mathbf{a},\mathbf{b},\mathbf{c}\in K\left[  \left[  x\right]  \right]
$. We must prove that $\mathbf{a}\left(  \mathbf{bc}\right)  =\left(
\mathbf{ab}\right)  \mathbf{c}$. Let $n\in\mathbb{N}$. Consider the two
equalities%
\begin{align*}
\left[  x^{n}\right]  \left(  \left(  \mathbf{ab}\right)  \mathbf{c}\right)
&  =\sum_{j=0}^{n}\underbrace{\left[  x^{n-j}\right]  \left(  \mathbf{ab}%
\right)  }_{\substack{=\sum_{i=0}^{n-j}\left[  x^{i}\right]  \mathbf{a}%
\cdot\left[  x^{n-j-i}\right]  \mathbf{b}\\\text{(by
(\ref{pf.thm.fps.ring.xn(ab)=2}),}\\\text{applied to }n-j\text{ instead of
}n\text{)}}}\cdot\left[  x^{j}\right]  \mathbf{c}\\
&  \ \ \ \ \ \ \ \ \ \ \ \ \ \ \ \ \ \ \ \ \left(
\begin{array}
[c]{c}%
\text{by (\ref{pf.thm.fps.ring.xn(ab)=3}), applied to }\mathbf{ab}\text{ and
}\mathbf{c}\\
\text{instead of }\mathbf{a}\text{ and }\mathbf{b}%
\end{array}
\right) \\
&  =\sum_{j=0}^{n}\left(  \sum_{i=0}^{n-j}\left[  x^{i}\right]  \mathbf{a}%
\cdot\left[  x^{n-j-i}\right]  \mathbf{b}\right)  \cdot\left[  x^{j}\right]
\mathbf{c}\\
&  =\sum_{j=0}^{n}\ \ \sum_{i=0}^{n-j}\left[  x^{i}\right]  \mathbf{a}%
\cdot\left[  x^{n-j-i}\right]  \mathbf{b}\cdot\left[  x^{j}\right]  \mathbf{c}%
\end{align*}
and%
\begin{align*}
\left[  x^{n}\right]  \left(  \mathbf{a}\left(  \mathbf{bc}\right)  \right)
&  =\sum_{i=0}^{n}\left[  x^{i}\right]  \mathbf{a}\cdot\underbrace{\left[
x^{n-i}\right]  \left(  \mathbf{bc}\right)  }_{\substack{=\sum_{j=0}%
^{n-i}\left[  x^{n-i-j}\right]  \mathbf{b}\cdot\left[  x^{j}\right]
\mathbf{c}\\\text{(by (\ref{pf.thm.fps.ring.xn(ab)=3}), applied to }n-i\text{,
}\mathbf{b}\text{ and }\mathbf{c}\\\text{instead of }n\text{, }\mathbf{a}%
\text{ and }\mathbf{b}\text{)}}}\\
&  \ \ \ \ \ \ \ \ \ \ \ \ \ \ \ \ \ \ \ \ \left(  \text{by
(\ref{pf.thm.fps.ring.xn(ab)=2}), applied to }\mathbf{bc}\text{ instead of
}\mathbf{b}\right) \\
&  =\sum_{i=0}^{n}\left[  x^{i}\right]  \mathbf{a}\cdot\left(  \sum
_{j=0}^{n-i}\left[  x^{n-i-j}\right]  \mathbf{b}\cdot\left[  x^{j}\right]
\mathbf{c}\right) \\
&  =\sum_{i=0}^{n}\ \ \sum_{j=0}^{n-i}\left[  x^{i}\right]  \mathbf{a}%
\cdot\left[  x^{n-i-j}\right]  \mathbf{b}\cdot\left[  x^{j}\right]
\mathbf{c}.
\end{align*}
The right hand sides of these two equalities are equal, since\footnote{The
first equality we are about to state is an equality of summation signs. Such
an equality means that whatever you put inside the summation signs, they
produce equal results. For example, $\sum_{i\in\left\{  1,2,3\right\}  }%
=\sum_{i=1}^{3}$ and $\sum_{i\in\mathbb{N}}=\sum_{i\geq0}$.}%
\[
\sum_{j=0}^{n}\ \ \sum_{i=0}^{n-j}=\sum_{\substack{\left(  i,j\right)
\in\mathbb{N}^{2};\\i+j\leq n}}=\sum_{i=0}^{n}\ \ \sum_{j=0}^{n-i}%
\]
and $n-j-i=n-i-j$. Thus, their left hand sides are equal as well. In other
words,%
\[
\left[  x^{n}\right]  \left(  \left(  \mathbf{ab}\right)  \mathbf{c}\right)
=\left[  x^{n}\right]  \left(  \mathbf{a}\left(  \mathbf{bc}\right)  \right)
.
\]

Now, forget that we fixed $n$. We thus have shown that $\left[  x^{n}\right]
\left(  \left(  \mathbf{ab}\right)  \mathbf{c}\right)  =\left[  x^{n}\right]
\left(  \mathbf{a}\left(  \mathbf{bc}\right)  \right)  $ for each
$n\in\mathbb{N}$. In other words, each entry of $\left(  \mathbf{ab}\right)
\mathbf{c}$ equals the corresponding entry of $\mathbf{a}\left(
\mathbf{bc}\right)  $. This entails $\left(  \mathbf{ab}\right)
\mathbf{c}=\mathbf{a}\left(  \mathbf{bc}\right)  $ (since a FPS is just the
sequence of its entries). In other words, $\mathbf{a}\left(  \mathbf{bc}%
\right)  =\left(  \mathbf{ab}\right)  \mathbf{c}$. This concludes the proof of
associativity of multiplication.

The remaining claims of Theorem \ref{thm.fps.ring} are
LTTR\footnote{\textquotedblleft LTTR\textquotedblright\ means
\textquotedblleft left to the reader\textquotedblright.} (their proofs follow
the same pattern, but are easier to execute).
\end{proof}

Since $K\left[  \left[  x\right]  \right]  $ is a commutative ring, any finite
sum of FPSs is well-defined. Sometimes, however, infinite sums of FPSs make
sense as well: for example, it stands to reason that%
\begin{align}
&  \ \ \ \ \left(  1,1,1,1,\ldots\right) \nonumber\\
&  +\left(  0,1,1,1,\ldots\right) \nonumber\\
&  +\left(  0,0,1,1,\ldots\right) \nonumber\\
&  +\left(  0,0,0,1,\ldots\right) \nonumber\\
&  +\cdots\nonumber\\
&  =\left(  1,2,3,4,\ldots\right)  , \label{eq.fps.summable.exa1}%
\end{align}
because FPSs are added entrywise. Let us rigorously define such sums. First,
we define \textquotedblleft essentially finite\textquotedblright\ sums of
elements of $K$:

\begin{definition}
\label{def.infsum.essfin}\textbf{(a)} A family $\left(  a_{i}\right)  _{i\in
I}\in K^{I}$ of elements of $K$ is said to be \emph{essentially finite} if all
but finitely many $i\in I$ satisfy $a_{i}=0$ (in other words, if the set
$\left\{  i\in I\ \mid\ a_{i}\neq0\right\}  $ is finite). \medskip

\textbf{(b)} Let $\left(  a_{i}\right)  _{i\in I}\in K^{I}$ be an essentially
finite family of elements of $K$. Then, the infinite sum $\sum_{i\in I}a_{i}$
is defined to equal the finite sum $\sum_{\substack{i\in I;\\a_{i}\neq0}%
}a_{i}$. Such an infinite sum is said to be \emph{essentially finite}.
\end{definition}

For example, the family $\left(  \left\lfloor \dfrac{5}{2^{n}}\right\rfloor
\right)  _{n\in\mathbb{N}}$ of integers is essentially finite\footnote{Here
and in the following, the notation $\left\lfloor x\right\rfloor $ (where $x$
is a real number) stands for the largest integer that is $\leq x$. For
instance, $\left\lfloor \pi\right\rfloor =3$ and $\left\lfloor 3\right\rfloor
=3$.}, and its sum is%
\begin{align*}
\sum_{n\in\mathbb{N}}\left\lfloor \dfrac{5}{2^{n}}\right\rfloor  &
=\left\lfloor \dfrac{5}{2^{0}}\right\rfloor +\left\lfloor \dfrac{5}{2^{1}%
}\right\rfloor +\left\lfloor \dfrac{5}{2^{2}}\right\rfloor +\left\lfloor
\dfrac{5}{2^{3}}\right\rfloor +\left\lfloor \dfrac{5}{2^{4}}\right\rfloor
+\cdots\\
&  =5+2+1+\underbrace{0+0+0+\cdots}_{\text{throw these away}}=5+2+1=8.
\end{align*}

Note the following:

\begin{itemize}
\item A family $\left(  a_{i}\right)  _{i\in I}\in K^{I}$ is always
essentially finite if $I$ is finite; thus, essentially finite families are a
wider class than finite families.

\item Any essentially finite sum of real or complex numbers is convergent in
the sense of analysis, but the converse is not true; for instance, the
infinite sum $\sum_{n\in\mathbb{N}}\dfrac{1}{2^{n}}=\dfrac{1}{1}+\dfrac{1}%
{2}+\dfrac{1}{4}+\dfrac{1}{8}+\cdots$ is convergent in the sense of analysis,
but not essentially finite. Essential finiteness is a crude algebraic
imitation of convergence. One of its advantages is that it works for sums of
elements of any commutative ring, not just of $\mathbb{R}$ or $\mathbb{C}$.
\end{itemize}

So the idea behind Definition \ref{def.infsum.essfin} \textbf{(b)} is that
addends that equal $0$ can be discarded in a sum, even when there are
infinitely many of them. \medskip

Sums of essentially finite families satisfy the usual rules for sums (such as
the breaking-apart rule $\sum_{s\in S}a_{s}=\sum_{s\in X}a_{s}+\sum_{s\in
Y}a_{s}$ when a set $S$ is the union of two disjoint sets $X$ and $Y$). See
\cite[\S 2.14.15]{detnotes} for details\footnote{Note that \cite[\S 2.14.15]%
{detnotes} uses the words \textquotedblleft\emph{finitely supported}%
\textquotedblright\ instead of \textquotedblleft essentially
finite\textquotedblright.}. There is only one \textbf{caveat}: Interchange of
summation signs (e.g., replacing $\sum_{i\in I}\ \ \sum_{j\in J}a_{i,j}$ by
$\sum_{j\in J}\ \ \sum_{i\in I}a_{i,j}$) works only if the family $\left(
a_{i,j}\right)  _{\left(  i,j\right)  \in I\times J}$ is essentially finite
(i.e., all but finitely many \textbf{pairs} $\left(  i,j\right)  \in I\times
J$ satisfy $a_{i,j}=0$); it does not suffice that the sums $\sum_{i\in
I}\ \ \sum_{j\in J}a_{i,j}$ and $\sum_{j\in J}\ \ \sum_{i\in I}a_{i,j}$
themselves are essentially finite (i.e., that the families $\left(
a_{i,j}\right)  _{j\in J}$ for all $i\in I$, the families $\left(
a_{i,j}\right)  _{i\in I}$ for all $j\in J$, and the families $\left(
\sum_{j\in J}a_{i,j}\right)  _{i\in I}$ and $\left(  \sum_{i\in I}%
a_{i,j}\right)  _{j\in J}$ are essentially finite).

\begin{fineprint}
For a counterexample, consider the family $\left(  a_{i,j}\right)  _{\left(
i,j\right)  \in I\times J}$ of integers with $I=\left\{  1,2,3,\ldots\right\}
$ and $J=\left\{  1,2,3,\ldots\right\}  $, where $a_{i,j}$ is given by the
following table:%
\[%
\begin{tabular}
[c]{|c||c|c|c|c|c|c|}\hline
$a_{i,j}$ & $j=1$ & $j=2$ & $j=3$ & $j=4$ & $j=5$ & $\cdots$\\\hline\hline
$i=1$ & $1$ & $-1$ &  &  &  & $\cdots$\\\hline
$i=2$ & $\phantom{-1}$ & $1$ & $-1$ &  &  & $\cdots$\\\hline
$i=3$ &  &  & $1$ & $-1$ &  & $\cdots$\\\hline
$i=4$ &  &  &  & $1$ & $-1$ & $\cdots$\\\hline
$i=5$ &  &  &  &  & $1$ & $\cdots$\\\hline
$\vdots$ & $\vdots$ & $\vdots$ & $\vdots$ & $\vdots$ & $\vdots$ & $\ddots
$\\\hline
\end{tabular}
\ \ \ \ \ \
\]
(where all the entries in the empty cells are $0$). For this family $\left(
a_{i,j}\right)  _{\left(  i,j\right)  \in I\times J}$, both sums $\sum_{i\in
I}\ \ \sum_{j\in J}a_{i,j}$ and $\sum_{j\in J}\ \ \sum_{i\in I}a_{i,j}$ are
essentially finite, but they are not equal (indeed, the former sum is
$\sum_{i\in I}\ \ \underbrace{\sum_{j\in J}a_{i,j}}_{=0}=\sum_{i\in I}0=0$,
whereas the latter sum is $\sum_{j\in J}\ \ \sum_{i\in I}a_{i,j}%
=\underbrace{\sum_{i\in I}a_{i,1}}_{=1}+\sum_{j>1}\ \ \underbrace{\sum_{i\in
I}a_{i,j}}_{=0}=1+\sum_{j>1}0=1$). And indeed, this family $\left(
a_{i,j}\right)  _{\left(  i,j\right)  \in I\times J}$ is not essentially finite.
\end{fineprint}

\bigskip

We have now made sense of infinite sums of elements of $K$ when all but
finitely many addends are $0$. Of course, we can do the same for $K\left[
\left[  x\right]  \right]  $ instead of $K$ (since $K\left[  \left[  x\right]
\right]  $, too, is a commutative ring). However, this does not help make
sense of the sum on the left hand side of (\ref{eq.fps.summable.exa1}),
because this sum is not essentially finite (it is a sum of infinitely many
nonzero FPSs). Thus, for sums of FPSs, we need a weaker version of essential
finiteness. Here is its definition:

\Needspace{19pc}

\begin{definition}
\label{def.fps.summable}A (possibly infinite) family $\left(  \mathbf{a}%
_{i}\right)  _{i\in I}$ of FPSs is said to be \emph{summable} (or
\emph{entrywise essentially finite}) if%
\[
\text{for each }n\in\mathbb{N}\text{, all but finitely many }i\in I\text{
satisfy }\left[  x^{n}\right]  \mathbf{a}_{i}=0.
\]
In this case, the sum $\sum_{i\in I}\mathbf{a}_{i}$ is defined to be the FPS
with%
\begin{equation}
\left[  x^{n}\right]  \left(  \sum_{i\in I}\mathbf{a}_{i}\right)
=\underbrace{\sum_{i\in I}\left[  x^{n}\right]  \mathbf{a}_{i}}%
_{\substack{\text{an essentially}\\\text{finite sum}}%
}\ \ \ \ \ \ \ \ \ \ \text{for all }n\in\mathbb{N}\text{.}
\label{eq.def.fps.summable.sum}%
\end{equation}
(By abuse of notation, we will sometimes call the sum $\sum_{i\in I}%
\mathbf{a}_{i}$ summable instead of the family $\left(  \mathbf{a}_{i}\right)
_{i\in I}$; but this will always mean that the latter family is summable.)
\end{definition}

\begin{remark}
The condition \textquotedblleft all but finitely many $i\in I$ satisfy
$\left[  x^{n}\right]  \mathbf{a}_{i}=0$\textquotedblright\ in Definition
\ref{def.fps.summable} is \textbf{not} equivalent to \textquotedblleft
infinitely many $i\in I$ satisfy $\left[  x^{n}\right]  \mathbf{a}_{i}%
=0$\textquotedblright.
\end{remark}

Any essentially finite family of FPSs is summable, but the converse is not
generally the case.

\begin{example}
Let us see how Definition \ref{def.fps.summable} justifies
(\ref{eq.fps.summable.exa1}). Consider the family $\left(  \mathbf{a}%
_{i}\right)  _{i\in\mathbb{N}}\in K\left[  \left[  x\right]  \right]
^{\mathbb{N}}$ of FPSs, where%
\[
\mathbf{a}_{i}:=\left(  \underbrace{0,0,\ldots,0}_{i\text{ zeroes}%
},1,1,1,\ldots\right)  \ \ \ \ \ \ \ \ \ \ \text{for each }i\in\mathbb{N}.
\]
The left hand side of (\ref{eq.fps.summable.exa1}) is precisely $\mathbf{a}%
_{0}+\mathbf{a}_{1}+\mathbf{a}_{2}+\cdots=\sum_{i\in\mathbb{N}}\mathbf{a}_{i}%
$, so let us check that the family $\left(  \mathbf{a}_{i}\right)
_{i\in\mathbb{N}}$ is summable and that its sum $\sum_{i\in\mathbb{N}%
}\mathbf{a}_{i}$ really equals the right hand side of
(\ref{eq.fps.summable.exa1}).

For each $n\in\mathbb{N}$, all but finitely many $i\in\mathbb{N}$ satisfy
$\left[  x^{n}\right]  \mathbf{a}_{i}=0$ (indeed, all $i>n$ satisfy this
equality). Thus, the family $\left(  \mathbf{a}_{i}\right)  _{i\in\mathbb{N}}$
is summable. For each $n\in\mathbb{N}$, we have%
\begin{align*}
\left[  x^{n}\right]  \left(  \sum_{i\in\mathbb{N}}\mathbf{a}_{i}\right)   &
=\sum_{i\in\mathbb{N}}\left[  x^{n}\right]  \mathbf{a}_{i}%
\ \ \ \ \ \ \ \ \ \ \left(  \text{by (\ref{eq.def.fps.summable.sum})}\right)
\\
&  =\sum_{i=0}^{n}\underbrace{\left[  x^{n}\right]  \mathbf{a}_{i}%
}_{\substack{=1\\\text{(since }i\leq n\text{)}}}+\sum_{i>n}\underbrace{\left[
x^{n}\right]  \mathbf{a}_{i}}_{\substack{=0\\\text{(since }i>n\text{)}}%
}=\sum_{i=0}^{n}1+\underbrace{\sum_{i>n}0}_{=0}\\
&  =\sum_{i=0}^{n}1=n+1.
\end{align*}
Thus, $\sum_{i\in\mathbb{N}}\mathbf{a}_{i}=\left(  1,2,3,4,\ldots\right)  $.
This is precisely the right hand side of (\ref{eq.fps.summable.exa1}). Thus,
(\ref{eq.fps.summable.exa1}) has been justified rigorously.
\end{example}

You can think of summable infinite sums of FPSs as a crude algebraic imitation
of uniformly convergent infinite sums of holomorphic functions in complex
analysis. (However, the former are a much simpler concept than the latter. In
particular, complex analysis is completely unnecessary in their study.)

The following fact is nearly obvious:\footnote{A \emph{subfamily} of a family
$\left(  f_{i}\right)  _{i\in I}$ means a family of the form $\left(
f_{i}\right)  _{i\in J}$, where $J$ is a subset of $I$.}

\begin{proposition}
\label{prop.fps.summable.sub}Let $\left(  \mathbf{a}_{i}\right)  _{i\in I}$ be
a summable family of FPSs. Then, any subfamily of $\left(  \mathbf{a}%
_{i}\right)  _{i\in I}$ is summable as well.
\end{proposition}

\begin{proof}
[Proof of Proposition \ref{prop.fps.summable.sub}.]Let $J$ be a subset of $I$.
We must prove that the subfamily $\left(  \mathbf{a}_{i}\right)  _{i\in J}$ is summable.

Let $n\in\mathbb{N}$. Then, all but finitely many $i\in I$ satisfy $\left[
x^{n}\right]  \mathbf{a}_{i}=0$ (since the family $\left(  \mathbf{a}%
_{i}\right)  _{i\in I}$ is summable). Hence, all but finitely many $i\in J$
satisfy $\left[  x^{n}\right]  \mathbf{a}_{i}=0$ (since $J$ is a subset of
$I$). Since we have proved this for each $n\in\mathbb{N}$, we thus conclude
that the family $\left(  \mathbf{a}_{i}\right)  _{i\in J}$ is summable. This
proves Proposition \ref{prop.fps.summable.sub}.
\end{proof}

Just as with essentially finite families, we can work with summable sums of
FPSs \textquotedblleft as if they were finite\textquotedblright\ most of the time:

\begin{proposition}
\label{prop.fps.summable-sums-rule}Sums of summable families of FPSs satisfy
the usual rules for sums (such as the breaking-apart rule $\sum_{s\in S}%
a_{s}=\sum_{s\in X}a_{s}+\sum_{s\in Y}a_{s}$ when a set $S$ is the union of
two disjoint sets $X$ and $Y$). See \cite[Proposition 7.2.11]{19s} for
details. Again, the only \textbf{caveat} is about interchange of summation
signs: The equality
\[
\sum_{i\in I}\ \ \sum_{j\in J}\mathbf{a}_{i,j}=\sum_{j\in J}\ \ \sum_{i\in
I}\mathbf{a}_{i,j}%
\]
holds when the family $\left(  \mathbf{a}_{i,j}\right)  _{\left(  i,j\right)
\in I\times J}$ is summable (i.e., when for each $n\in\mathbb{N}$, all but
finitely many \textbf{pairs} $\left(  i,j\right)  \in I\times J$ satisfy
$\left[  x^{n}\right]  \mathbf{a}_{i,j}=0$); it does not generally hold if we
merely assume that the sums $\sum_{i\in I}\ \ \sum_{j\in J}\mathbf{a}_{i,j}$
and $\sum_{j\in J}\ \ \sum_{i\in I}\mathbf{a}_{i,j}$ are summable.
\end{proposition}

\begin{proof}
[Proof of Proposition \ref{prop.fps.summable-sums-rule}.]The proof is tedious
(as there are many rules to check), but fairly straightforward (the idea is
always to focus on a single coefficient, and then to reduce the infinite sums
to finite sums). For example, consider the \textquotedblleft discrete Fubini
rule\textquotedblright, which says that%
\begin{equation}
\sum_{i\in I}\ \ \sum_{j\in J}\mathbf{a}_{i,j}=\sum_{\left(  i,j\right)  \in
I\times J}\mathbf{a}_{i,j}=\sum_{j\in J}\ \ \sum_{i\in I}\mathbf{a}_{i,j}
\label{pf.prop.fps.summable-sums-rule.fub}%
\end{equation}
whenever $\left(  \mathbf{a}_{i,j}\right)  _{\left(  i,j\right)  \in I\times
J}$ is a summable family of FPSs. In order to prove this rule, we fix a
summable family $\left(  \mathbf{a}_{i,j}\right)  _{\left(  i,j\right)  \in
I\times J}$ of FPSs. It is easy to see that the families $\left(
\mathbf{a}_{i,j}\right)  _{j\in J}$ for all $i\in I$ are summable as well, as
are the families $\left(  \mathbf{a}_{i,j}\right)  _{i\in I}$ for all $j\in
J$, and the families $\left(  \sum_{j\in J}\mathbf{a}_{i,j}\right)  _{i\in I}$
and $\left(  \sum_{i\in I}\mathbf{a}_{i,j}\right)  _{j\in J}$. Hence, all sums
in (\ref{pf.prop.fps.summable-sums-rule.fub}) are well-defined. Now, in order
to prove (\ref{pf.prop.fps.summable-sums-rule.fub}), it suffices to check that%
\[
\left[  x^{n}\right]  \left(  \sum_{i\in I}\ \ \sum_{j\in J}\mathbf{a}%
_{i,j}\right)  =\left[  x^{n}\right]  \left(  \sum_{\left(  i,j\right)  \in
I\times J}\mathbf{a}_{i,j}\right)  =\left[  x^{n}\right]  \left(  \sum_{j\in
J}\ \ \sum_{i\in I}\mathbf{a}_{i,j}\right)
\]
for each $n\in\mathbb{N}$. Fix $n\in\mathbb{N}$; then, we have $\left[
x^{n}\right]  \left(  \mathbf{a}_{i,j}\right)  =0$ for all but finitely many
$\left(  i,j\right)  \in I\times J$ (since the family $\left(  \mathbf{a}%
_{i,j}\right)  _{\left(  i,j\right)  \in I\times J}$ is summable). That is,
the set of all pairs $\left(  i,j\right)  \in I\times J$ satisfying $\left[
x^{n}\right]  \left(  \mathbf{a}_{i,j}\right)  \neq0$ is finite. Hence, the
set $I^{\prime}:=\left\{  i\ \mid\ \left(  i,j\right)  \in I\times J\text{
with }\left[  x^{n}\right]  \left(  \mathbf{a}_{i,j}\right)  \neq0\right\}  $
of the first entries of all these pairs is also finite, and so is the set
$J^{\prime}:=\left\{  j\ \mid\ \left(  i,j\right)  \in I\times J\text{ with
}\left[  x^{n}\right]  \left(  \mathbf{a}_{i,j}\right)  \neq0\right\}  $ of
the second entries of all these pairs. Now, the definitions of $I^{\prime}$
and $J^{\prime}$ ensure that any pair $\left(  i,j\right)  \in I\times J$
satisfies $\left[  x^{n}\right]  \mathbf{a}_{i,j}=0$ unless $i\in I^{\prime}$
and $j\in J^{\prime}$. Hence, we easily obtain the three equalities%
\[
\left[  x^{n}\right]  \left(  \sum_{i\in I}\ \ \sum_{j\in J}\mathbf{a}%
_{i,j}\right)  =\sum_{i\in I}\ \ \sum_{j\in J}\left[  x^{n}\right]
\mathbf{a}_{i,j}=\sum_{i\in I^{\prime}}\ \ \sum_{j\in J^{\prime}}\left[
x^{n}\right]  \mathbf{a}_{i,j}%
\]
and%
\[
\left[  x^{n}\right]  \left(  \sum_{\left(  i,j\right)  \in I\times
J}\mathbf{a}_{i,j}\right)  =\sum_{\left(  i,j\right)  \in I\times J}\left[
x^{n}\right]  \mathbf{a}_{i,j}=\sum_{\left(  i,j\right)  \in I^{\prime}\times
J^{\prime}}\left[  x^{n}\right]  \mathbf{a}_{i,j}%
\]
and%
\[
\left[  x^{n}\right]  \left(  \sum_{j\in J}\ \ \sum_{i\in I}\mathbf{a}%
_{i,j}\right)  =\sum_{j\in J}\ \ \sum_{i\in I}\left[  x^{n}\right]
\mathbf{a}_{i,j}=\sum_{j\in J^{\prime}}\ \ \sum_{i\in I^{\prime}}\left[
x^{n}\right]  \mathbf{a}_{i,j}.
\]
However, the right hand sides of these three equalities are equal (since the
sums appearing in them are finite sums, and thus satisfy the usual rules for
sums). Thus, the left hand sides are equal, exactly as we needed to show. See
\cite[proof of Proposition 7.2.11]{19s} for more details of this proof.
Proving the other properties of sums is easier.
\end{proof}

A few conventions about infinite sums will be used rather often:

\begin{convention}
\label{conv.fps.infsum}\textbf{(a)} For any given integer $m\in\mathbb{Z}$,
the summation sign $\sum_{k\geq m}$ is to be understood as $\sum_{k\in\left\{
m,m+1,m+2,\ldots\right\}  }$. We also write $\sum_{k=m}^{\infty}$ for this
summation sign. \medskip

\textbf{(b)} For any given integer $m\in\mathbb{Z}$, the summation sign
$\sum_{k>m}$ is to be understood as $\sum_{k\in\left\{  m+1,m+2,m+3,\ldots
\right\}  }$. \medskip

\textbf{(c)} Let $I$ be a set, and let $\mathcal{A}\left(  i\right)  $ be a
logical statement for each $i\in I$. (For example, $I$ can be $\mathbb{N}$,
and $\mathcal{A}\left(  i\right)  $ can be the statement \textquotedblleft$i$
is odd\textquotedblright.) Then, the summation sign $\sum_{\substack{i\in
I;\\\mathcal{A}\left(  i\right)  }}$ is to be understood as $\sum
_{i\in\left\{  j\in I\ \mid\ \mathcal{A}\left(  j\right)  \right\}  }$. (For
example, the summation sign $\sum_{\substack{i\in\mathbb{N};\\i\text{ is odd}%
}}$ means $\sum_{i\in\left\{  j\in\mathbb{N}\ \mid\ j\text{ is odd}\right\}
}$, that is, a sum over all odd elements of $\mathbb{N}$.)
\end{convention}

We can now define the $x$ that figured so prominently in our informal
exploration of formal power series back in Section \ref{sec.gf.exas}:

\begin{definition}
\label{def.fps.x}Let $x$ denote the FPS $\left(  0,1,0,0,0,\ldots\right)  $.
In other words, let $x$ denote the FPS with $\left[  x^{1}\right]  x=1$ and
$\left[  x^{i}\right]  x=0$ for all $i\neq1$.
\end{definition}

The following simple lemma follows almost immediately from the definition of
multiplication of FPSs:

\begin{lemma}
\label{lem.fps.xa}Let $\mathbf{a}=\left(  a_{0},a_{1},a_{2},\ldots\right)  $
be an FPS. Then, $x\cdot\mathbf{a}=\left(  0,a_{0},a_{1},a_{2},\ldots\right)
$.
\end{lemma}

In other words, multiplying an FPS $\mathbf{a}$ by $x$ is tantamount to
inserting a $0$ at the front of $\mathbf{a}$ (and shifting all the previously
existing entries of $\mathbf{a}$ to the right by one position).

\begin{proof}
[Proof of Lemma \ref{lem.fps.xa}.]If $n$ is a positive integer, then%
\begin{align*}
\left[  x^{n}\right]  \left(  x\cdot\mathbf{a}\right)   &  =\sum_{i=0}%
^{n}\left[  x^{i}\right]  x\cdot\underbrace{\left[  x^{n-i}\right]
\mathbf{a}}_{\substack{=a_{n-i}\\\text{(since }\mathbf{a}=\left(  a_{0}%
,a_{1},a_{2},\ldots\right)  \text{)}}}\\
&  \ \ \ \ \ \ \ \ \ \ \ \ \ \ \ \ \ \ \ \ \left(  \text{by
(\ref{pf.thm.fps.ring.xn(ab)=2}), applied to }x\text{ and }\mathbf{a}\text{
instead of }\mathbf{a}\text{ and }\mathbf{b}\right) \\
&  =\sum_{i=0}^{n}\left[  x^{i}\right]  x\cdot a_{n-i}=\underbrace{\left[
x^{0}\right]  x}_{=0}\cdot\,a_{n-0}+\underbrace{\left[  x^{1}\right]  x}%
_{=1}\cdot\,a_{n-1}+\sum_{i=2}^{n}\underbrace{\left[  x^{i}\right]
x}_{\substack{=0\\\text{(since }i\geq2>1\text{)}}}\cdot\,a_{n-i}\\
&  \ \ \ \ \ \ \ \ \ \ \ \ \ \ \ \ \ \ \ \ \left(
\begin{array}
[c]{c}%
\text{here, we have split off the addends}\\
\text{for }i=0\text{ and }i=1\text{ from the sum}%
\end{array}
\right) \\
&  =\underbrace{0\cdot a_{n-0}}_{=0}+\,1\cdot a_{n-1}+\underbrace{\sum
_{i=2}^{n}0\cdot a_{n-i}}_{=0}=1\cdot a_{n-1}=a_{n-1}.
\end{align*}
A similar argument can be used for $n=0$ (except that now, the sum $\sum
_{i=0}^{n}\left[  x^{i}\right]  x\cdot a_{n-i}$ has no $\left[  x^{1}\right]
x\cdot a_{n-1}$ addend), and results in the conclusion that $\left[
x^{n}\right]  \left(  x\cdot\mathbf{a}\right)  =0$ in this case. Thus, for
each $n\in\mathbb{N}$, we have%
\[
\left[  x^{n}\right]  \left(  x\cdot\mathbf{a}\right)  =%
\begin{cases}
a_{n-1}, & \text{if }n>0;\\
0, & \text{if }n=0.
\end{cases}
\]
In other words, $x\cdot\mathbf{a}=\left(  0,a_{0},a_{1},a_{2},\ldots\right)
$. This proves Lemma \ref{lem.fps.xa}.
\end{proof}

Recall that $x^{k}=\underbrace{xx\cdots x}_{k\text{ times}}$ for each
$k\in\mathbb{N}$ (by the definition of powers in any commutative ring). In
particular, $x^{0}=\underline{1}$ (since $\underline{1}$ is the unity of the
ring $K\left[  \left[  x\right]  \right]  $). The following proposition
describes $x^{k}$ explicitly for each $k\in\mathbb{N}$:

\begin{proposition}
\label{prop.fps.xk}We have%
\[
x^{k}=\left(  \underbrace{0,0,\ldots,0}_{k\text{ times}},1,0,0,0,\ldots
\right)  \ \ \ \ \ \ \ \ \ \ \text{for each }k\in\mathbb{N}.
\]

\end{proposition}

\begin{proof}
[Proof of Proposition \ref{prop.fps.xk}.]We induct on $k$.

\textit{Induction base:} We have $x^{0}=\underline{1}=\left(  1,0,0,0,0,\ldots
\right)  =\left(  \underbrace{0,0,\ldots,0}_{0\text{ times}},1,0,0,0,\ldots
\right)  $. In other words, Proposition \ref{prop.fps.xk} holds for $k=0$.

\textit{Induction step:} Let $m\in\mathbb{N}$. Assume that Proposition
\ref{prop.fps.xk} holds for $k=m$. We must prove that Proposition
\ref{prop.fps.xk} holds for $k=m+1$.

We have $x^{m}=\left(  \underbrace{0,0,\ldots,0}_{m\text{ times}%
},1,0,0,0,\ldots\right)  $ (since Proposition \ref{prop.fps.xk} holds for
$k=m$). Thus, Lemma \ref{lem.fps.xa} (applied to $\mathbf{a}=x^{m}$ and
\newline$\left(  a_{0},a_{1},a_{2},\ldots\right)  =\left(
\underbrace{0,0,\ldots,0}_{m\text{ times}},1,0,0,0,\ldots\right)  $) yields%
\[
x\cdot x^{m}=\left(  0,\underbrace{0,0,\ldots,0}_{m\text{ times}%
},1,0,0,0,\ldots\right)  =\left(  \underbrace{0,0,\ldots,0}_{m+1\text{ times}%
},1,0,0,0,\ldots\right)  .
\]
In other words, $x^{m+1}=\left(  \underbrace{0,0,\ldots,0}_{m+1\text{ times}%
},1,0,0,0,\ldots\right)  $ (since $x\cdot x^{m}=x^{m+1}$). In other words,
Proposition \ref{prop.fps.xk} holds for $k=m+1$. This completes the induction
step, thus proving Proposition \ref{prop.fps.xk}.
\end{proof}

Finally, the following corollary allows us to rewrite any FPS $\left(
a_{0},a_{1},a_{2},\ldots\right)  $ in the familiar form $a_{0}+a_{1}%
x+a_{2}x^{2}+a_{3}x^{3}+\cdots$:

\begin{corollary}
\label{cor.fps.sumakxk}Any FPS $\left(  a_{0},a_{1},a_{2},\ldots\right)  \in
K\left[  \left[  x\right]  \right]  $ satisfies%
\[
\left(  a_{0},a_{1},a_{2},\ldots\right)  =a_{0}+a_{1}x+a_{2}x^{2}+a_{3}%
x^{3}+\cdots=\sum_{n\in\mathbb{N}}a_{n}x^{n}.
\]
In particular, the right hand side here is well-defined, i.e., the family
$\left(  a_{n}x^{n}\right)  _{n\in\mathbb{N}}$ is summable.
\end{corollary}

\begin{proof}
[Proof of Corollary \ref{cor.fps.sumakxk} (sketched).](See \cite[Corollary
7.2.16]{19s} for details.) By Proposition \ref{prop.fps.xk}, we have%
\begin{align*}
&  a_{0}+a_{1}x+a_{2}x^{2}+a_{3}x^{3}+\cdots\\
&  =\ \ \ \ a_{0}\left(  1,0,0,0,\ldots\right) \\
&  \ \ \ \ +a_{1}\left(  0,1,0,0,\ldots\right) \\
&  \ \ \ \ +a_{2}\left(  0,0,1,0,\ldots\right) \\
&  \ \ \ \ +a_{3}\left(  0,0,0,1,\ldots\right) \\
&  \ \ \ \ +\cdots\\
&  =\ \ \ \ \left(  a_{0},0,0,0,\ldots\right) \\
&  \ \ \ \ +\left(  0,a_{1},0,0,\ldots\right) \\
&  \ \ \ \ +\left(  0,0,a_{2},0,\ldots\right) \\
&  \ \ \ \ +\left(  0,0,0,a_{3},\ldots\right) \\
&  \ \ \ \ +\cdots\\
&  =\left(  a_{0},a_{1},a_{2},a_{3},\ldots\right)  .
\end{align*}
This proves Corollary \ref{cor.fps.sumakxk} (since $a_{0}+a_{1}x+a_{2}%
x^{2}+a_{3}x^{3}+\cdots=\sum_{n\in\mathbb{N}}a_{n}x^{n}$ holds for obvious reasons).
\end{proof}

So we have \textquotedblleft found\textquotedblright\ our $x$ and given a
rigorous justification for writing $\left(  a_{0},a_{1},a_{2},\ldots\right)  $
as $a_{0}+a_{1}x+a_{2}x^{2}+a_{3}x^{3}+\cdots$. Note that we did not use any
analysis (real or complex) in the process; in particular, we did not have to
worry about convergence (we did have to worry about summability, but this is
much simpler than convergence). It is easy to come up with FPSs that don't
converge when any nonzero real number is substituted for $x$ (for example,
$\sum_{n\in\mathbb{N}}n!x^{n}=1+x+2x^{2}+6x^{3}+24x^{4}+\cdots$ is such an
FPS); they are nevertheless completely legitimate FPSs.

We can now also answer the question \textquotedblleft what is a generating
function\textquotedblright, albeit the answer is somewhat anticlimactic:

\begin{definition}
Let $\left(  a_{0},a_{1},a_{2},\ldots\right)  $ be a sequence of elements of
$K$. Then, the \emph{(ordinary) generating function} of $\left(  a_{0}%
,a_{1},a_{2},\ldots\right)  $ will mean the FPS $\left(  a_{0},a_{1}%
,a_{2},\ldots\right)  =a_{0}+a_{1}x+a_{2}x^{2}+a_{3}x^{3}+\cdots$.
\end{definition}

For future use, let us observe the following fact (in a way, part of
Proposition \ref{prop.fps.summable-sums-rule}), which extends the law of
distributivity in $K\left[  \left[  x\right]  \right]  $ to infinite sums (as
long as they are summable):

\begin{proposition}
\label{prop.fps.sum-dist-inf}Let $\left(  \mathbf{a}_{i}\right)  _{i\in I}$ be
a summable family of FPSs. Let $\mathbf{b}$ be a further FPS. Then, the family
$\left(  \mathbf{ba}_{i}\right)  _{i\in I}$ is again summable, and satisfies%
\[
\sum_{i\in I}\mathbf{ba}_{i}=\mathbf{b}\sum_{i\in I}\mathbf{a}_{i}.
\]

\end{proposition}

\begin{proof}
[Proof sketch.]The \textquotedblleft hard\textquotedblright\ part is showing
that the family $\left(  \mathbf{ba}_{i}\right)  _{i\in I}$ is summable. For
this purpose, we must show that for each $n\in\mathbb{N}$, all but finitely
many $i\in I$ satisfy $\left[  x^{n}\right]  \left(  \mathbf{ba}_{i}\right)
=0$. So let us fix $n\in\mathbb{N}$. Then:

\begin{itemize}
\item All but finitely many $i\in I$ satisfy $\left[  x^{0}\right]
\mathbf{a}_{i}=0$ (since the family $\left(  \mathbf{a}_{i}\right)  _{i\in I}$
is summable).

\item All but finitely many $i\in I$ satisfy $\left[  x^{1}\right]
\mathbf{a}_{i}=0$ (since the family $\left(  \mathbf{a}_{i}\right)  _{i\in I}$
is summable).

\item All but finitely many $i\in I$ satisfy $\left[  x^{2}\right]
\mathbf{a}_{i}=0$ (since the family $\left(  \mathbf{a}_{i}\right)  _{i\in I}$
is summable).

\item $\ldots$

\item All but finitely many $i\in I$ satisfy $\left[  x^{n}\right]
\mathbf{a}_{i}=0$ (since the family $\left(  \mathbf{a}_{i}\right)  _{i\in I}$
is summable).
\end{itemize}

Combining all these bullet points, we conclude that all but finitely many
$i\in I$ satisfy all the $n+1$ equalities%
\[
\left[  x^{m}\right]  \mathbf{a}_{i}=0\ \ \ \ \ \ \ \ \ \ \text{for all }%
m\in\left\{  0,1,\ldots,n\right\}
\]
(since a union of $n+1$ finite sets is still a finite set). All these $i$ must
then satisfy%
\begin{align*}
\left[  x^{n}\right]  \left(  \mathbf{ba}_{i}\right)   &  =\sum_{j=0}%
^{n}\left[  x^{n-j}\right]  \mathbf{b}\cdot\underbrace{\left[  x^{j}\right]
\mathbf{a}_{i}}_{\substack{=0\\\text{(since }\left[  x^{m}\right]
\mathbf{a}_{i}=0\\\text{for all }m\in\left\{  0,1,\ldots,n\right\}  \text{)}%
}}\ \ \ \ \ \ \ \ \ \ \left(
\begin{array}
[c]{c}%
\text{by (\ref{pf.thm.fps.ring.xn(ab)=3}), applied to }\mathbf{b}\text{ and
}\mathbf{a}_{i}\\
\text{instead of }\mathbf{a}\text{ and }\mathbf{b}%
\end{array}
\right) \\
&  =\sum_{j=0}^{n}\left[  x^{n-j}\right]  \mathbf{b}\cdot0=0.
\end{align*}

Forget that we fixed $n$. Thus we have shown that for each $n\in\mathbb{N}$,
all but finitely many $i\in I$ satisfy $\left[  x^{n}\right]  \left(
\mathbf{ba}_{i}\right)  =0$. In other words, the family $\left(
\mathbf{ba}_{i}\right)  _{i\in I}$ is summable.

It remains to prove that $\sum_{i\in I}\mathbf{ba}_{i}=\mathbf{b}\sum_{i\in
I}\mathbf{a}_{i}$. This can be checked coefficient by coefficient: For each
$n\in\mathbb{N}$, we have%
\begin{align*}
\left[  x^{n}\right]  \left(  \sum_{i\in I}\mathbf{ba}_{i}\right)   &
=\sum_{i\in I}\underbrace{\left[  x^{n}\right]  \left(  \mathbf{ba}%
_{i}\right)  }_{\substack{=\sum_{j=0}^{n}\left[  x^{n-j}\right]
\mathbf{b}\cdot\left[  x^{j}\right]  \mathbf{a}_{i}\\\text{(by
(\ref{pf.thm.fps.ring.xn(ab)=3}) again)}}}\\
&  =\sum_{i\in I}\ \ \sum_{j=0}^{n}\left[  x^{n-j}\right]  \mathbf{b}%
\cdot\left[  x^{j}\right]  \mathbf{a}_{i}\\
&  =\sum_{j=0}^{n}\left[  x^{n-j}\right]  \mathbf{b}\cdot\underbrace{\sum
_{i\in I}\left[  x^{j}\right]  \mathbf{a}_{i}}_{=\left[  x^{j}\right]  \left(
\sum_{i\in I}\mathbf{a}_{i}\right)  }\\
&  \ \ \ \ \ \ \ \ \ \ \ \ \ \ \ \ \ \ \ \ \left(
\begin{array}
[c]{c}%
\text{this interchange of summation signs was}\\
\text{allowed, since the summability of }\left(  \mathbf{a}_{i}\right)  _{i\in
I}\\
\text{ensures that for each }j\in\left\{  0,1,\ldots,n\right\}  \text{,}\\
\text{all but finitely many }i\in I\text{ satisfy }\left[  x^{j}\right]
\mathbf{a}_{i}=0\text{,}\\
\text{and thus all but finitely many}\\
\text{pairs }\left(  i,j\right)  \in I\times\left\{  0,1,\ldots,n\right\} \\
\text{satisfy }\left[  x^{j}\right]  \mathbf{a}_{i}=0
\end{array}
\right) \\
&  =\sum_{j=0}^{n}\left[  x^{n-j}\right]  \mathbf{b}\cdot\left[  x^{j}\right]
\left(  \sum_{i\in I}\mathbf{a}_{i}\right) \\
&  =\left[  x^{n}\right]  \left(  \mathbf{b}\sum_{i\in I}\mathbf{a}%
_{i}\right)
\end{align*}
(since another application of (\ref{pf.thm.fps.ring.xn(ab)=3}) shows that
$\left[  x^{n}\right]  \left(  \mathbf{b}\sum_{i\in I}\mathbf{a}_{i}\right)
=\sum_{j=0}^{n}\left[  x^{n-j}\right]  \mathbf{b}\cdot\left[  x^{j}\right]
\left(  \sum_{i\in I}\mathbf{a}_{i}\right)  $). This completes the proof of
Proposition \ref{prop.fps.sum-dist-inf}.
\end{proof}

\subsubsection{\label{subsec.gf.defs.cvi}The Chu--Vandermonde identity}

What we have done so far is sufficient to justify Example 3 from Section
\ref{sec.gf.exas}. Thus, let us record the result of Example 3 as a proposition:

\begin{proposition}
\label{prop.binom.vandermonde.NN}Let $a,b\in\mathbb{N}$, and let
$n\in\mathbb{N}$. Then,%
\begin{equation}
\dbinom{a+b}{n}=\sum_{k=0}^{n}\dbinom{a}{k}\dbinom{b}{n-k}.
\label{eq.prop.binom.vandermonde.NN.eq}%
\end{equation}

\end{proposition}

We have yet to justify Examples 1, 2 and 4; we shall do so later. First,
however, let us generalize Proposition \ref{prop.binom.vandermonde.NN} to
arbitrary numbers $a,b$ (as opposed to merely $a,b\in\mathbb{N}$). That is, we
shall prove the following:

\begin{theorem}
[\emph{Vandermonde convolution identity}, aka \emph{Chu--Vandermonde
identity}]\label{thm.binom.vandermonde.CC}Let $a,b\in\mathbb{C}$, and let
$n\in\mathbb{N}$. Then,%
\begin{equation}
\dbinom{a+b}{n}=\sum_{k=0}^{n}\dbinom{a}{k}\dbinom{b}{n-k}.
\label{eq.prop.binom.vandermonde.CC.eq}%
\end{equation}

\end{theorem}

(Actually, the \textquotedblleft Let $a,b\in\mathbb{C}$\textquotedblright\ in
Theorem \ref{thm.binom.vandermonde.CC} can be replaced by \textquotedblleft
Let $a,b$ be any numbers\textquotedblright, where \textquotedblleft
numbers\textquotedblright\ is understood appropriately\footnote{Namely, a
\textquotedblleft number\textquotedblright\ should here be viewed as an
element of a commutative $\mathbb{Q}$-algebra. This includes complex numbers,
polynomials over complex numbers, power series over complex numbers and even
commuting matrices with complex entries.}. We just wrote \textquotedblleft Let
$a,b\in\mathbb{C}$\textquotedblright\ for simplicity.)

To recall, back in Example 3, we proved Proposition
\ref{prop.binom.vandermonde.NN} by multiplying out the identity $\left(
1+x\right)  ^{a+b}=\left(  1+x\right)  ^{a}\left(  1+x\right)  ^{b}$ using the
binomial formula. If we wanted to extend this argument to arbitrary
$a,b\in\mathbb{C}$, then we would need to make sense of powers like $\left(
1+x\right)  ^{-3}$ or $\left(  1+x\right)  ^{1/2}$ or $\left(  1+x\right)
^{\pi}$. This is indeed possible (in fact, we will briefly outline this later
on), but there is a much shorter way.

In fact, there is a slick trick to automatically extend a claim like
(\ref{eq.prop.binom.vandermonde.NN.eq}) from nonnegative integers to complex
numbers. It is sometimes known as the \emph{polynomial identity trick}, and is
used a lot in algebra. The proof of Theorem \ref{thm.binom.vandermonde.CC}
that we shall sketch below should illustrate this trick; you can find more
details and further examples in \cite[\S 7.5.3]{20f}.

\begin{proof}
[Proof of Theorem \ref{thm.binom.vandermonde.CC} (sketched).]Fix
$n\in\mathbb{N}$ and $b\in\mathbb{N}$, but forget that $a$ was fixed. Then,
Proposition \ref{prop.binom.vandermonde.NN} (which we have already proved)
says that the equality%
\begin{equation}
\dbinom{a+b}{n}=\sum_{k=0}^{n}\dbinom{a}{k}\dbinom{b}{n-k}
\label{pf.thm.binom.vandermonde.CC.1}%
\end{equation}
holds for each $a\in\mathbb{N}$. However, both sides of this equality are
polynomials (more precisely, polynomial functions) in $a$: indeed,%
\begin{align*}
\dbinom{a+b}{n}  &  =\dfrac{\left(  a+b\right)  \left(  a+b-1\right)  \left(
a+b-2\right)  \cdots\left(  a+b-n+1\right)  }{n!}%
\ \ \ \ \ \ \ \ \ \ \text{and}\\
\sum_{k=0}^{n}\dbinom{a}{k}\dbinom{b}{n-k}  &  =\sum_{k=0}^{n}\dfrac{a\left(
a-1\right)  \cdots\left(  a-k+1\right)  }{k!}\dbinom{b}{n-k}.
\end{align*}
If two univariate polynomials $p$ and $q$ (with rational, real or complex
coefficients) are equal for each $a\in\mathbb{N}$ (that is: if we have
$p\left(  a\right)  =q\left(  a\right)  $ for each $a\in\mathbb{N}$), then
they must be identical (because two univariate polynomials that are equal at
infinitely many points must necessarily be identical\footnote{Quick reminder
on why this is true: If $p$ and $q$ are two univariate polynomials (with
rational, real or complex coefficients) that are equal at infinitely many
points (i.e., if there exist infinitely many numbers $z$ satisfying $p\left(
z\right)  =q\left(  z\right)  $), then $p=q$ (because the assumption entails
that the difference $p-q$ has infinitely many roots, but this entails $p-q=0$
and thus $p=q$). See \cite[Corollary 7.5.7]{20f} for this argument in more
detail.}). Hence, the two sides of (\ref{pf.thm.binom.vandermonde.CC.1}) must
be identical as polynomials in $a$. Thus, the equality
(\ref{pf.thm.binom.vandermonde.CC.1}) holds not only for each $a\in\mathbb{N}%
$, but also for each $a\in\mathbb{C}$.

Now, forget that $b$ was fixed. Instead, let us fix $a\in\mathbb{C}$. As we
just have proved, the equality (\ref{pf.thm.binom.vandermonde.CC.1}) holds for
each $b\in\mathbb{N}$. We want to show that it holds for each $b\in\mathbb{C}%
$. But this can be achieved by the same argument that we just used to extend
it from $a\in\mathbb{N}$ to $a\in\mathbb{C}$: We view both sides of the
equality as polynomials (but this time in $b$, not in $a$), and argue that
these polynomials must be identical because they are equal at infinitely many
points. The upshot is that the equality (\ref{pf.thm.binom.vandermonde.CC.1})
holds for all $a,b\in\mathbb{C}$; thus, Theorem \ref{thm.binom.vandermonde.CC}
is proven. (See \cite[proofs of Lemma 7.5.8 and Theorem 7.5.3]{20f} or
\cite[\S 2.17.3]{19s} for this proof in more detail. Alternatively, see
\cite[\S 3.3.2 and \S 3.3.3]{detnotes} for two other proofs.)
\end{proof}

\subsubsection{What next?}

Let us now return to our quest of justifying Examples 1, 2 and 4 from Section
\ref{sec.gf.exas}. In order to do so, we need to know

\begin{itemize}
\item what we can substitute into a FPS;

\item when and why we can divide FPSs by FPSs;

\item when and why we can take the square root of an FPS and solve a quadratic
equation using the quadratic formula.
\end{itemize}

So we need to do more. The following sections are devoted to this.

\subsection{Dividing FPSs}

\subsubsection{Conventions}

We shall make ourselves at home in the ring $K\left[  \left[  x\right]
\right]  $ a bit more. (Recall that $K$ is a fixed commutative ring.)

\begin{convention}
\label{conv.fps.const}From now on, we identify each $a\in K$ with the constant
FPS $\underline{a}=\left(  a,0,0,0,0,\ldots\right)  \in K\left[  \left[
x\right]  \right]  $.
\end{convention}

This is motivated by the fact that $\underline{a}=a+0x+0x^{2}+0x^{3}+\cdots$
for any $a\in K$. Convention \ref{conv.fps.const} does not cause any dangerous
ambiguities, because we have%
\begin{align*}
\underline{a+b}  &  =\underline{a}+\underline{b}\ \ \ \ \ \ \ \ \ \ \text{and}%
\\
\underline{a-b}  &  =\underline{a}-\underline{b}\ \ \ \ \ \ \ \ \ \ \text{and}%
\\
\underline{a\cdot b}  &  =\underline{a}\cdot\underline{b}%
\ \ \ \ \ \ \ \ \ \ \text{for any }a,b\in K
\end{align*}
(check this!), and because the zero and the unity of the ring $K\left[
\left[  x\right]  \right]  $ are $\underline{0}$ and $\underline{1}$, respectively.

Furthermore, I will stop using boldfaced letters (like $\mathbf{a}%
,\mathbf{b},\mathbf{c}$) for FPSs. (I did this above for the sake of
convenience, but this is rarely done in the literature.)

\subsubsection{Inverses in commutative rings}

We recall the notion of an \emph{inverse} in a commutative ring:

\begin{definition}
\label{def.commring.inverse}Let $L$ be a commutative ring. Let $a\in L$. Then:
\medskip

\textbf{(a)} An \emph{inverse} (or \emph{multiplicative inverse}) of $a$ means
an element $b\in L$ such that $ab=ba=1$ (where the \textquotedblleft%
$1$\textquotedblright\ means the unity of $L$). \medskip

\textbf{(b)} We say that $a$ is \emph{invertible} in $L$ (or a \emph{unit} of
$L$) if $a$ has an inverse.
\end{definition}

Note that the condition \textquotedblleft$ab=ba=1$\textquotedblright\ in
Definition \ref{def.commring.inverse} \textbf{(a)} can be restated as
\textquotedblleft$ab=1$\textquotedblright, because we automatically have
$ab=ba$ (since $L$ is a commutative ring). I have chosen to write
\textquotedblleft$ab=ba=1$\textquotedblright\ in order to state the definition
in a form that applies verbatim to noncommutative rings as well.

\begin{example}
\textbf{(a)} In the ring $\mathbb{Z}$, the only two invertible elements are
$1$ and $-1$. Each of these two elements is its own inverse. \medskip

\textbf{(b)} In the ring $\mathbb{Q}$, every nonzero element is invertible.
The same holds for the rings $\mathbb{R}$ and $\mathbb{C}$ (and, more
generally, for any field).
\end{example}

Our next goal is to study inverses of FPSs in $K\left[  \left[  x\right]
\right]  $, answering in particular the natural question \textquotedblleft
which elements of $K\left[  \left[  x\right]  \right]  $ have
inverses\textquotedblright. But let us first prove their uniqueness in the
generality of an arbitrary commutative ring:

\begin{theorem}
\label{thm.commring.inverse-uni}Let $L$ be a commutative ring. Let $a\in L$.
Then, there is \textbf{at most one} inverse of $a$.
\end{theorem}

\begin{proof}
Let $b$ and $c$ be two inverses of $a$. We must prove that $b=c$.

Since $b$ is an inverse of $a$, we have $ab=ba=1$. Since $c$ is an inverse of
$a$, we have $ac=ca=1$. Now, we have $b\underbrace{\left(  ac\right)  }%
_{=1}=b\cdot1=b$ and $\underbrace{\left(  ba\right)  }_{=1}c=1\cdot c=c$.
However, because of the \textquotedblleft associativity of
multiplication\textquotedblright\ axiom in Definition \ref{def.alg.commring},
we have $b\left(  ac\right)  =\left(  ba\right)  c$. Hence, $b=b\left(
ac\right)  =\left(  ba\right)  c=c$. This proves Theorem
\ref{thm.commring.inverse-uni}.
\end{proof}

Theorem \ref{thm.commring.inverse-uni} allows us to make the following definition:

\begin{definition}
\label{def.commring.fracs}Let $L$ be a commutative ring. Let $a\in L$. Assume
that $a$ is invertible. Then: \medskip

\textbf{(a)} The inverse of $a$ is called $a^{-1}$. (This notation is
unambiguous, since Theorem \ref{thm.commring.inverse-uni} shows that the
inverse of $a$ is unique.) \medskip

\textbf{(b)} For any $b\in L$, the product $b\cdot a^{-1}$ is called
$\dfrac{b}{a}$ (or $b/a$). \medskip

\textbf{(c)} For any negative integer $n$, we define $a^{n}$ to be $\left(
a^{-1}\right)  ^{-n}$. Thus, the $n$-th power $a^{n}$ is defined for each
$n\in\mathbb{Z}$.
\end{definition}

The following facts are easy to check:

\begin{proposition}
\label{prop.commring.fracs.1}Let $L$ be a commutative ring. Then: \medskip

\textbf{(a)} Any invertible element $a\in L$ satisfies $a^{-1}=1/a$. \medskip

\textbf{(b)} For any invertible elements $a,b\in L$, the element $ab$ is
invertible as well, and satisfies $\left(  ab\right)  ^{-1}=b^{-1}%
a^{-1}=a^{-1}b^{-1}$. \medskip

\textbf{(c)} If $a\in L$ is invertible, and if $n\in\mathbb{Z}$ is arbitrary,
then $a^{-n}=\left(  a^{-1}\right)  ^{n}=\left(  a^{n}\right)  ^{-1}$.
\medskip

\textbf{(d)} Laws of exponents hold for negative exponents as well: If $a$ and
$b$ are invertible elements of $L$, then%
\begin{align*}
a^{n+m}  &  =a^{n}a^{m}\ \ \ \ \ \ \ \ \ \ \text{for all }n,m\in\mathbb{Z};\\
\left(  ab\right)  ^{n}  &  =a^{n}b^{n}\ \ \ \ \ \ \ \ \ \ \text{for all }%
n\in\mathbb{Z};\\
\left(  a^{n}\right)  ^{m}  &  =a^{nm}\ \ \ \ \ \ \ \ \ \ \text{for all
}n,m\in\mathbb{Z}.
\end{align*}

\textbf{(e)} Laws of fractions hold: If $a$ and $c$ are two invertible
elements of $L$, and if $b$ and $d$ are any two elements of $L$, then
\[
\dfrac{b}{a}+\dfrac{d}{c}=\dfrac{bc+ad}{ac}\ \ \ \ \ \ \ \ \ \ \text{and}%
\ \ \ \ \ \ \ \ \ \ \dfrac{b}{a}\cdot\dfrac{d}{c}=\dfrac{bd}{ac}.
\]

\textbf{(f)} Division undoes multiplication: If $a,b,c$ are three elements of
$L$ with $a$ being invertible, then the equality $\dfrac{c}{a}=b$ is
equivalent to $c=ab$.
\end{proposition}

\begin{proof}
Exercise. (See, e.g., \cite[solution to Exercise 4.1.1]{19s} for a proof of
parts \textbf{(c)} and \textbf{(d)} in the special case where $L=\mathbb{C}$;
essentially the same argument works in the general case. The remaining parts
of Proposition \ref{prop.commring.fracs.1} are even easier to check. Note that
parts \textbf{(a)} and \textbf{(c)} as well as the $\left(  ab\right)
^{-1}=b^{-1}a^{-1}$ part of part \textbf{(b)} would hold even if $L$ were a
noncommutative ring.)
\end{proof}

\subsubsection{Inverses in $K\left[  \left[  x\right]  \right]  $}

Now, which FPSs are invertible in the ring $K\left[  \left[  x\right]
\right]  $ ? For example, we know from (\ref{eq.sec.gf.exas.1.1/1-x}) that the
FPS $1-x$ is invertible, with inverse $1+x+x^{2}+x^{3}+\cdots$. On the other
hand, the FPS $x$ is not invertible, since Lemma \ref{lem.fps.xa} shows that
any product of $x$ with an FPS must begin with a $0$ (but the unity of
$K\left[  \left[  x\right]  \right]  $ does not begin with a $0$). (Strictly
speaking, this is only true if the ring $K$ is nontrivial -- i.e., if not all
elements of $K$ are equal. If $K$ is trivial, then $K\left[  \left[  x\right]
\right]  $ is trivial, and thus any FPS in $K\left[  \left[  x\right]
\right]  $ is invertible, but this does not make an interesting statement.)

It turns out that we can characterize invertible FPSs in $K\left[  \left[
x\right]  \right]  $ in a rather simple way:

\begin{proposition}
\label{prop.fps.invertible}Let $a\in K\left[  \left[  x\right]  \right]  $.
Then, the FPS $a$ is invertible in $K\left[  \left[  x\right]  \right]  $ if
and only if its constant term $\left[  x^{0}\right]  a$ is invertible in $K$.
\end{proposition}

\begin{proof}
The statement we are proving is an \textquotedblleft if and only
if\textquotedblright\ statement. We shall prove its \textquotedblleft only
if\textquotedblright\ (i.e., \textquotedblleft$\Longrightarrow$%
\textquotedblright) and its \textquotedblleft if\textquotedblright\ (i.e.,
\textquotedblleft$\Longleftarrow$\textquotedblright) directions separately:

\begin{enumerate}
\item[$\Longrightarrow:$] Assume that $a$ is invertible in $K\left[  \left[
x\right]  \right]  $. That is, $a$ has an inverse $b\in K\left[  \left[
x\right]  \right]  $. Consider this $b$. Since $b$ is an inverse of $a$, we
have $ab=ba=1$ (where \textquotedblleft$1$\textquotedblright\ means
$\underline{1}$ by Convention \ref{conv.fps.const}). However,
(\ref{pf.thm.fps.ring.x0(ab)=}) (applied to $\mathbf{a}=a$ and $\mathbf{b}=b$)
yields $\left[  x^{0}\right]  \left(  ab\right)  =\left[  x^{0}\right]
a\cdot\left[  x^{0}\right]  b$. Comparing this with $\left[  x^{0}\right]
\underbrace{\left(  ab\right)  }_{=1}=\left[  x^{0}\right]  1=1$, we find
$\left[  x^{0}\right]  a\cdot\left[  x^{0}\right]  b=1$. Thus, $\left[
x^{0}\right]  b$ is an inverse of $\left[  x^{0}\right]  a$ in $K$ (since
$\left[  x^{0}\right]  b\cdot\left[  x^{0}\right]  a=\left[  x^{0}\right]
a\cdot\left[  x^{0}\right]  b=1$, so that $\left[  x^{0}\right]  a\cdot\left[
x^{0}\right]  b=\left[  x^{0}\right]  b\cdot\left[  x^{0}\right]  a=1$).
Therefore, $\left[  x^{0}\right]  a$ is invertible in $K$ (with inverse
$\left[  x^{0}\right]  b$). This proves the \textquotedblleft$\Longrightarrow
$\textquotedblright\ direction of Proposition \ref{prop.fps.invertible}.

\item[$\Longleftarrow:$] Assume that $\left[  x^{0}\right]  a$ is invertible
in $K$. Write the FPS $a$ in the form $a=\left(  a_{0},a_{1},a_{2}%
,\ldots\right)  $. Thus, $\left[  x^{0}\right]  a=a_{0}$, so that $a_{0}$ is
invertible in $K$ (since $\left[  x^{0}\right]  a$ is invertible in $K$).
Thus, its inverse $a_{0}^{-1}$ is well-defined.

Now, we want to prove that $a$ is invertible in $K\left[  \left[  x\right]
\right]  $. We thus try to find an inverse of $a$.

We work backwards at first: We assume that $b=\left(  b_{0},b_{1},b_{2}%
,\ldots\right)  \in K\left[  \left[  x\right]  \right]  $ is an inverse for
$a$, and we try to figure out how this inverse looks like.

Since $b$ is an inverse of $a$, we have $ab=\underline{1}=\left(
1,0,0,0,\ldots\right)  $. However, from $a=\left(  a_{0},a_{1},a_{2}%
,\ldots\right)  $ and $b=\left(  b_{0},b_{1},b_{2},\ldots\right)  $, we have%
\begin{align*}
ab  &  =\left(  a_{0},a_{1},a_{2},\ldots\right)  \left(  b_{0},b_{1}%
,b_{2},\ldots\right) \\
&  =\left(  a_{0}b_{0},\ \ a_{0}b_{1}+a_{1}b_{0},\ \ a_{0}b_{2}+a_{1}%
b_{1}+a_{2}b_{0},\ \ a_{0}b_{3}+a_{1}b_{2}+a_{2}b_{1}+a_{3}b_{0}%
,\ \ \ldots\right)
\end{align*}
(by the definition of the product of FPSs). Comparing this with $ab=\left(
1,0,0,0,\ldots\right)  $, we obtain%
\begin{align*}
&  \left(  1,0,0,0,\ldots\right) \\
&  =\left(  a_{0}b_{0},\ \ a_{0}b_{1}+a_{1}b_{0},\ \ a_{0}b_{2}+a_{1}%
b_{1}+a_{2}b_{0},\ \ a_{0}b_{3}+a_{1}b_{2}+a_{2}b_{1}+a_{3}b_{0}%
,\ \ \ldots\right)  .
\end{align*}
This can be rewritten as the following system of equations:%
\begin{equation}
\left\{
\begin{array}
[c]{l}%
1=a_{0}b_{0},\\
0=a_{0}b_{1}+a_{1}b_{0},\\
0=a_{0}b_{2}+a_{1}b_{1}+a_{2}b_{0},\\
0=a_{0}b_{3}+a_{1}b_{2}+a_{2}b_{1}+a_{3}b_{0},\\
\ldots.
\end{array}
\right.  \label{pf.prop.fps.invertible.sys}%
\end{equation}
I claim that this system of equations uniquely determines $\left(  b_{0}%
,b_{1},b_{2},\ldots\right)  $. Indeed, we can solve the first equation
($1=a_{0}b_{0}$) for $b_{0}$, thus obtaining $b_{0}=a_{0}^{-1}$ (since $a_{0}$
is invertible). Having thus found $b_{0}$, we can solve the second equation
($0=a_{0}b_{1}+a_{1}b_{0}$) for $b_{1}$, thus obtaining $b_{1}=-a_{0}%
^{-1}\left(  a_{1}b_{0}\right)  $ (again because $a_{0}$ is invertible).
Having thus found both $b_{0}$ and $b_{1}$, we can solve the third equation
($0=a_{0}b_{2}+a_{1}b_{1}+a_{2}b_{0}$) for $b_{2}$, thus obtaining
$b_{2}=-a_{0}^{-1}\left(  a_{1}b_{1}+a_{2}b_{0}\right)  $. Proceeding like
this, we obtain recursive expressions for all coefficients $b_{0},b_{1}%
,b_{2},\ldots$ of $b$, namely%
\begin{equation}
\left\{
\begin{array}
[c]{l}%
b_{0}=a_{0}^{-1},\\
b_{1}=-a_{0}^{-1}\left(  a_{1}b_{0}\right)  ,\\
b_{2}=-a_{0}^{-1}\left(  a_{1}b_{1}+a_{2}b_{0}\right)  ,\\
b_{3}=-a_{0}^{-1}\left(  a_{1}b_{2}+a_{2}b_{1}+a_{3}b_{0}\right)  ,\\
\ldots.
\end{array}
\right.  \label{pf.prop.fps.invertible.sol}%
\end{equation}
(This procedure for solving systems of linear equations is well-known from
linear algebra -- it is a form of Gaussian elimination, but a particularly
simple one because our system is triangular with invertible coefficients on
the diagonal. The only complication is that it has infinitely many variables
and infinitely many equations.)

So we have shown that if $b$ is an inverse of $a$, then the entries $b_{i}$ of
the FPS $b$ are given recursively by (\ref{pf.prop.fps.invertible.sol}). This
yields that $b$ is unique; alas, this is not what we want to prove. Instead,
we want to prove that $b$ exists.

Fortunately, we can achieve this by simply turning our above argument around:
Forget that we fixed $b$. Instead, we define a sequence $\left(  b_{0}%
,b_{1},b_{2},\ldots\right)  $ of elements of $K$ recursively by
(\ref{pf.prop.fps.invertible.sol}), and we define the FPS $b=\left(
b_{0},b_{1},b_{2},\ldots\right)  \in K\left[  \left[  x\right]  \right]
$.\ Then, the equalities (\ref{pf.prop.fps.invertible.sys}) hold (because they
are just equivalent restatements of the equalities
(\ref{pf.prop.fps.invertible.sol})). In other words, we have%
\begin{align*}
&  \left(  1,0,0,0,\ldots\right) \\
&  =\left(  a_{0}b_{0},\ \ a_{0}b_{1}+a_{1}b_{0},\ \ a_{0}b_{2}+a_{1}%
b_{1}+a_{2}b_{0},\ \ a_{0}b_{3}+a_{1}b_{2}+a_{2}b_{1}+a_{3}b_{0}%
,\ \ \ldots\right)  .
\end{align*}
However, as before, we can show that%
\[
ab=\left(  a_{0}b_{0},\ \ a_{0}b_{1}+a_{1}b_{0},\ \ a_{0}b_{2}+a_{1}%
b_{1}+a_{2}b_{0},\ \ a_{0}b_{3}+a_{1}b_{2}+a_{2}b_{1}+a_{3}b_{0}%
,\ \ \ldots\right)  .
\]
Comparing these two equalities, we find $ab=\left(  1,0,0,0,\ldots\right)
=\underline{1}$. Thus, $ba=ab=\underline{1}$, so that $ab=ba=\underline{1}$.
This shows that $b$ is an inverse of $a$, so that $a$ is invertible. This
proves the \textquotedblleft$\Longleftarrow$\textquotedblright\ direction of
Proposition \ref{prop.fps.invertible}.
\end{enumerate}
\end{proof}

We note a particularly simple corollary of Proposition
\ref{prop.fps.invertible} when $K$ is a field:

\begin{corollary}
\label{cor.fps.invertible.field}Assume that $K$ is a field. Let $a\in K\left[
\left[  x\right]  \right]  $. Then, the FPS $a$ is invertible in $K\left[
\left[  x\right]  \right]  $ if and only if $\left[  x^{0}\right]  a\neq0$.
\end{corollary}

\begin{proof}
An element of $K$ is invertible in $K$ if and only if it is nonzero (since $K$
is a field). Hence, Corollary \ref{cor.fps.invertible.field} follows from
Proposition \ref{prop.fps.invertible}.
\end{proof}

\subsubsection{Newton's binomial formula}

Let us now return to considering specific FPSs. We have already seen that the
FPS $1-x$ is invertible, with inverse $1+x+x^{2}+x^{3}+\cdots$. We shall now
show an analogous result for the FPS $1+x$. Its invertibility follows from
Proposition \ref{prop.fps.invertible}, but it is better to derive it by hand,
as this also gives a formula for the inverse:

\begin{proposition}
\label{prop.fps.invertible.1+x}The FPS $1+x\in K\left[  \left[  x\right]
\right]  $ is invertible, and its inverse is%
\[
\left(  1+x\right)  ^{-1}=1-x+x^{2}-x^{3}+x^{4}-x^{5}\pm\cdots=\sum
_{n\in\mathbb{N}}\left(  -1\right)  ^{n}x^{n}.
\]

\end{proposition}

\begin{proof}
[First proof of Proposition \ref{prop.fps.invertible.1+x}.]We have%
\begin{align*}
&  \left(  1+x\right)  \cdot\left(  1-x+x^{2}-x^{3}+x^{4}-x^{5}\pm
\cdots\right) \\
&  =\underbrace{1\cdot\left(  1-x+x^{2}-x^{3}+x^{4}-x^{5}\pm\cdots\right)
}_{=1-x+x^{2}-x^{3}+x^{4}-x^{5}\pm\cdots}+\underbrace{x\cdot\left(
1-x+x^{2}-x^{3}+x^{4}-x^{5}\pm\cdots\right)  }_{=x-x^{2}+x^{3}-x^{4}%
+x^{5}-x^{6}\pm\cdots}\\
&  =\left(  1-x+x^{2}-x^{3}+x^{4}-x^{5}\pm\cdots\right)  +\left(
x-x^{2}+x^{3}-x^{4}+x^{5}-x^{6}\pm\cdots\right) \\
&  =1
\end{align*}
(since all powers of $x$ other than $x^{0}=1$ cancel out). This shows that
$1-x+x^{2}-x^{3}+x^{4}-x^{5}\pm\cdots$ is an inverse of $1+x$ (since $K\left[
\left[  x\right]  \right]  $ is a commutative ring). Thus, $1+x$ is
invertible, and its inverse is $\left(  1+x\right)  ^{-1}=1-x+x^{2}%
-x^{3}+x^{4}-x^{5}\pm\cdots=\sum_{n\in\mathbb{N}}\left(  -1\right)  ^{n}x^{n}%
$. This proves Proposition \ref{prop.fps.invertible.1+x}.
\end{proof}

\begin{proof}
[Second proof of Proposition \ref{prop.fps.invertible.1+x}.]Using Proposition
\ref{prop.fps.sum-dist-inf}, we have%
\begin{align*}
&  \left(  1+x\right)  \cdot\left(  1-x+x^{2}-x^{3}+x^{4}-x^{5}\pm
\cdots\right) \\
&  =\underbrace{\left(  1+x\right)  \cdot1}_{=1+x}-\underbrace{\left(
1+x\right)  \cdot x}_{=x+x^{2}}+\underbrace{\left(  1+x\right)  \cdot x^{2}%
}_{=x^{2}+x^{3}}-\underbrace{\left(  1+x\right)  \cdot x^{3}}_{=x^{3}+x^{4}%
}+\underbrace{\left(  1+x\right)  \cdot x^{4}}_{=x^{4}+x^{5}}%
-\underbrace{\left(  1+x\right)  \cdot x^{5}}_{=x^{5}+x^{6}}\pm\cdots\\
&  =\left(  1+x\right)  -\left(  x+x^{2}\right)  +\left(  x^{2}+x^{3}\right)
-\left(  x^{3}+x^{4}\right)  +\left(  x^{4}+x^{5}\right)  -\left(  x^{5}%
+x^{6}\right)  \pm\cdots\\
&  =1
\end{align*}
(since we have a telescoping sum in front of us, in which all powers of $x$
other than $x^{0}=1$ cancel out). This shows that $1-x+x^{2}-x^{3}+x^{4}%
-x^{5}\pm\cdots$ is an inverse of $1+x$ (since $K\left[  \left[  x\right]
\right]  $ is a commutative ring). Thus, $1+x$ is invertible, and its inverse
is $\left(  1+x\right)  ^{-1}=1-x+x^{2}-x^{3}+x^{4}-x^{5}\pm\cdots=\sum
_{n\in\mathbb{N}}\left(  -1\right)  ^{n}x^{n}$. This proves Proposition
\ref{prop.fps.invertible.1+x}.
\end{proof}

Proposition \ref{prop.fps.invertible.1+x} shows that the FPS $1+x$ is
invertible; thus, its powers $\left(  1+x\right)  ^{n}$ are defined for all
$n\in\mathbb{Z}$ (by Definition \ref{def.commring.fracs} \textbf{(c)}). The
following formula -- known as \emph{Newton's binomial theorem}\footnote{or
\emph{Newton's binomial formula}} -- describes these powers explicitly:

\begin{theorem}
\label{thm.fps.newton-binom}For each $n\in\mathbb{Z}$, we have%
\[
\left(  1+x\right)  ^{n}=\sum_{k\in\mathbb{N}}\dbinom{n}{k}x^{k}.
\]

\end{theorem}

Note that the sum $\sum_{k\in\mathbb{N}}\dbinom{n}{k}x^{k}$ is summable for
each $n\in\mathbb{Z}$ (indeed, it equals the FPS $\left(  \dbinom{n}%
{0},\dbinom{n}{1},\dbinom{n}{2},\ldots\right)  $). If $n\in\mathbb{N}$, then
it is essentially finite.

The reader may want to check that the particular case $n=-1$ of Theorem
\ref{thm.fps.newton-binom} agrees with Proposition
\ref{prop.fps.invertible.1+x}. (Recall Example \ref{exa.binom.-1choosek}!)

Of course, Theorem \ref{thm.fps.newton-binom} should look familiar -- an
identical-looking formula appears in real analysis under the same name.
However, the result in real analysis is concerned with infinite sums of real
numbers, while our Theorem \ref{thm.fps.newton-binom} is an identity between
FPSs over an arbitrary commutative ring. Thus, the two facts are not the same.

We will prove Theorem \ref{thm.fps.newton-binom} in a somewhat roundabout way,
since this gives us an opportunity to establish some auxiliary results that
are of separate interest (and usefulness). The first of these auxiliary
results is a fundamental property of binomial coefficients, known as the
\emph{upper negation formula} (see, e.g., \cite[Proposition 1.3.7]{19fco}):

\begin{theorem}
\label{thm.binom.upneg-n}Let $n\in\mathbb{C}$ and $k\in\mathbb{Z}$. Then,%
\[
\dbinom{-n}{k}=\left(  -1\right)  ^{k}\dbinom{k+n-1}{k}.
\]

\end{theorem}

\begin{proof}
[Proof of Theorem \ref{thm.binom.upneg-n} (sketched).]If $k<0$, then this is
trivial because both $\dbinom{-n}{k}$ and $\dbinom{k+n-1}{k}$ are $0$ (by
(\ref{eq.def.binom.binom.eq})). Thus, we WLOG assume that $k\geq0$. Hence,
$k\in\mathbb{N}$. Thus, (\ref{eq.def.binom.binom.eq}) yields%
\begin{align*}
\dbinom{-n}{k}  &  =\dfrac{\left(  -n\right)  \left(  -n-1\right)  \left(
-n-2\right)  \cdots\left(  -n-k+1\right)  }{k!}\ \ \ \ \ \ \ \ \ \ \text{and}%
\\
\dbinom{k+n-1}{k}  &  =\dfrac{\left(  k+n-1\right)  \left(  k+n-2\right)
\left(  k+n-3\right)  \cdots\left(  k+n-k\right)  }{k!}.
\end{align*}
A moment of thought reveals that the right hand sides of these two equalities
are equal up to a factor of $\left(  -1\right)  ^{k}$. Thus, so are the left
hand sides. In other words, $\dbinom{-n}{k}=\left(  -1\right)  ^{k}%
\dbinom{k+n-1}{k}$. This proves Theorem \ref{thm.binom.upneg-n}.
\end{proof}

(Quick exercise: Rederive Example \ref{exa.binom.-1choosek} from Theorem
\ref{thm.binom.upneg-n}.)

Next, we show a formula for negative powers of $1+x$:

\begin{proposition}
\label{prop.fps.anti-newton-binom}For each $n\in\mathbb{N}$, we have%
\[
\left(  1+x\right)  ^{-n}=\sum_{k\in\mathbb{N}}\left(  -1\right)  ^{k}%
\dbinom{n+k-1}{k}x^{k}.
\]

\end{proposition}

\begin{proof}
[Proof of Proposition \ref{prop.fps.anti-newton-binom}.]We proceed by
induction on $n$:

\textit{Induction base:} Comparing
\[
\left(  1+x\right)  ^{-0}=\left(  1+x\right)  ^{0}=\underline{1}=\left(
1,0,0,0,\ldots\right)  =1
\]
with%
\begin{align*}
&  \sum_{k\in\mathbb{N}}\left(  -1\right)  ^{k}\dbinom{0+k-1}{k}x^{k}\\
&  =\underbrace{\left(  -1\right)  ^{0}}_{=1}\underbrace{\dbinom{0+0-1}{0}%
}_{=1}\underbrace{x^{0}}_{=1}+\sum_{\substack{k\in\mathbb{N};\\k>0}}\left(
-1\right)  ^{k}\underbrace{\dbinom{0+k-1}{k}}_{\substack{=\dbinom{k-1}%
{k}=0\\\text{(by Proposition \ref{prop.binom.0}}\\\text{(applied to
}m=k-1\text{ and }n=k\text{),}\\\text{since }k-1<k\text{ and }k-1\in
\mathbb{N}\text{)}}}x^{k}\\
&  \ \ \ \ \ \ \ \ \ \ \ \ \ \ \ \ \ \ \ \ \left(  \text{here, we have split
off the addend for }k=0\text{ from the sum}\right) \\
&  =1+\underbrace{\sum_{\substack{k\in\mathbb{N};\\k>0}}\left(  -1\right)
^{k}0x^{k}}_{=0}=1,
\end{align*}
we obtain $\left(  1+x\right)  ^{-0}=\sum_{k\in\mathbb{N}}\left(  -1\right)
^{k}\dbinom{0+k-1}{k}x^{k}$. In other words, Proposition
\ref{prop.fps.anti-newton-binom} holds for $n=0$.

\textit{Induction step:} Let $j\in\mathbb{N}$. Assume that Proposition
\ref{prop.fps.anti-newton-binom} holds for $n=j$. We must prove that
Proposition \ref{prop.fps.anti-newton-binom} holds for $n=j+1$.

We have assumed that Proposition \ref{prop.fps.anti-newton-binom} holds for
$n=j$. In other words, we have%
\begin{equation}
\left(  1+x\right)  ^{-j}=\sum_{k\in\mathbb{N}}\left(  -1\right)  ^{k}%
\dbinom{j+k-1}{k}x^{k}. \label{pf.prop.fps.anti-newton-binom.IH}%
\end{equation}

Now, we want to prove that Proposition \ref{prop.fps.anti-newton-binom} holds
for $n=j+1$. In other words, we want to prove that%
\[
\left(  1+x\right)  ^{-\left(  j+1\right)  }=\sum_{k\in\mathbb{N}}\left(
-1\right)  ^{k}\dbinom{\left(  j+1\right)  +k-1}{k}x^{k}.
\]
In view of $\left(  1+x\right)  ^{-\left(  j+1\right)  }=\left(  1+x\right)
^{-j}\cdot\left(  1+x\right)  ^{-1}$ and $\left(  j+1\right)  +k-1=j+k$, this
equality can be rewritten as%
\[
\left(  1+x\right)  ^{-j}\cdot\left(  1+x\right)  ^{-1}=\sum_{k\in\mathbb{N}%
}\left(  -1\right)  ^{k}\dbinom{j+k}{k}x^{k}.
\]
Since $1+x$ is invertible, we can equivalently transform this equality by
multiplying both sides with $1+x$; thus, it becomes%
\[
\left(  1+x\right)  ^{-j}=\left(  \sum_{k\in\mathbb{N}}\left(  -1\right)
^{k}\dbinom{j+k}{k}x^{k}\right)  \cdot\left(  1+x\right)  .
\]
So this is the equality we must prove.

We do this by simplifying its right hand side: Using Proposition
\ref{prop.fps.sum-dist-inf}, we have%
\begin{align*}
&  \left(  \sum_{k\in\mathbb{N}}\left(  -1\right)  ^{k}\dbinom{j+k}{k}%
x^{k}\right)  \cdot\left(  1+x\right) \\
&  =\sum_{k\in\mathbb{N}}\underbrace{\left(  -1\right)  ^{k}\dbinom{j+k}%
{k}x^{k}\left(  1+x\right)  }_{\substack{=\left(  -1\right)  ^{k}\dbinom
{j+k}{k}\left(  x^{k}+x^{k+1}\right)  \\=\left(  -1\right)  ^{k}\dbinom
{j+k}{k}x^{k}+\left(  -1\right)  ^{k}\dbinom{j+k}{k}x^{k+1}}}\\
&  =\sum_{k\in\mathbb{N}}\left(  \left(  -1\right)  ^{k}\dbinom{j+k}{k}%
x^{k}+\left(  -1\right)  ^{k}\dbinom{j+k}{k}x^{k+1}\right) \\
&  =\sum_{k\in\mathbb{N}}\left(  -1\right)  ^{k}\underbrace{\dbinom{j+k}{k}%
}_{\substack{=\dbinom{j+k-1}{k-1}+\dbinom{j+k-1}{k}\\\text{(by Proposition
\ref{prop.binom.rec},}\\\text{applied to }m=j+k\text{ and }n=k\text{)}}%
}x^{k}+\underbrace{\sum_{k\in\mathbb{N}}\left(  -1\right)  ^{k}\dbinom{j+k}%
{k}x^{k+1}}_{\substack{=\sum_{k\geq1}\left(  -1\right)  ^{k-1}\dbinom
{j+\left(  k-1\right)  }{k-1}x^{k}\\\text{(here, we have substituted
}k-1\text{ for }k\\\text{in the sum)}}}\\
&  =\underbrace{\sum_{k\in\mathbb{N}}\left(  -1\right)  ^{k}\left(
\dbinom{j+k-1}{k-1}+\dbinom{j+k-1}{k}\right)  x^{k}}_{=\sum_{k\in\mathbb{N}%
}\left(  -1\right)  ^{k}\dbinom{j+k-1}{k-1}x^{k}+\sum_{k\in\mathbb{N}}\left(
-1\right)  ^{k}\dbinom{j+k-1}{k}x^{k}}+\sum_{k\geq1}\underbrace{\left(
-1\right)  ^{k-1}}_{=-\left(  -1\right)  ^{k}}\underbrace{\dbinom{j+\left(
k-1\right)  }{k-1}}_{=\dbinom{j+k-1}{k-1}}x^{k}\\
&  =\sum_{k\in\mathbb{N}}\left(  -1\right)  ^{k}\dbinom{j+k-1}{k-1}x^{k}%
+\sum_{k\in\mathbb{N}}\left(  -1\right)  ^{k}\dbinom{j+k-1}{k}x^{k}%
+\sum_{k\geq1}\left(  -\left(  -1\right)  ^{k}\right)  \dbinom{j+k-1}%
{k-1}x^{k}\\
&  =\sum_{k\in\mathbb{N}}\left(  -1\right)  ^{k}\dbinom{j+k-1}{k-1}x^{k}%
+\sum_{k\in\mathbb{N}}\left(  -1\right)  ^{k}\dbinom{j+k-1}{k}x^{k}%
-\sum_{k\geq1}\left(  -1\right)  ^{k}\dbinom{j+k-1}{k-1}x^{k}\\
&  =\underbrace{\left(  \sum_{k\in\mathbb{N}}\left(  -1\right)  ^{k}%
\dbinom{j+k-1}{k-1}x^{k}-\sum_{k\geq1}\left(  -1\right)  ^{k}\dbinom
{j+k-1}{k-1}x^{k}\right)  }_{\substack{=\left(  -1\right)  ^{0}\dbinom
{j+0-1}{0-1}x^{0}\\\text{(since the two sums differ only in their }k=0\text{
addend)}}}+\underbrace{\sum_{k\in\mathbb{N}}\left(  -1\right)  ^{k}%
\dbinom{j+k-1}{k}x^{k}}_{\substack{=\left(  1+x\right)  ^{-j}\\\text{(by
(\ref{pf.prop.fps.anti-newton-binom.IH}))}}}\\
&  =\left(  -1\right)  ^{0}\underbrace{\dbinom{j+0-1}{0-1}}%
_{\substack{=0\\\text{(by (\ref{eq.def.binom.binom.eq}), since }%
0-1\notin\mathbb{N}\text{)}}}x^{0}+\left(  1+x\right)  ^{-j}=\left(
1+x\right)  ^{-j}.
\end{align*}
Multiplying both sides of this equality by $\left(  1+x\right)  ^{-1}$, we
obtain%
\[
\sum_{k\in\mathbb{N}}\left(  -1\right)  ^{k}\dbinom{j+k}{k}x^{k}=\left(
1+x\right)  ^{-j}\cdot\left(  1+x\right)  ^{-1}=\left(  1+x\right)  ^{\left(
-j\right)  +\left(  -1\right)  }=\left(  1+x\right)  ^{-\left(  j+1\right)
},
\]
so that
\[
\left(  1+x\right)  ^{-\left(  j+1\right)  }=\sum_{k\in\mathbb{N}}\left(
-1\right)  ^{k}\underbrace{\dbinom{j+k}{k}}_{=\dbinom{\left(  j+1\right)
+k-1}{k}}x^{k}=\sum_{k\in\mathbb{N}}\left(  -1\right)  ^{k}\dbinom{\left(
j+1\right)  +k-1}{k}x^{k}.
\]
In other words, Proposition \ref{prop.fps.anti-newton-binom} holds for
$n=j+1$. This completes the induction step. Thus, Proposition
\ref{prop.fps.anti-newton-binom} is proved.
\end{proof}

We can rewrite Proposition \ref{prop.fps.anti-newton-binom} using negative
binomial coefficients:

\begin{corollary}
\label{cor.fps.anti-newton-binom-2}For each $n\in\mathbb{N}$, we have%
\[
\left(  1+x\right)  ^{-n}=\sum_{k\in\mathbb{N}}\dbinom{-n}{k}x^{k}.
\]

\end{corollary}

\begin{proof}
[Proof of Corollary \ref{cor.fps.anti-newton-binom-2}.]Proposition
\ref{prop.fps.anti-newton-binom} yields%
\[
\left(  1+x\right)  ^{-n}=\sum_{k\in\mathbb{N}}\underbrace{\left(  -1\right)
^{k}\dbinom{n+k-1}{k}}_{\substack{=\left(  -1\right)  ^{k}\dbinom{k+n-1}%
{k}\\=\dbinom{-n}{k}\\\text{(by Theorem \ref{thm.binom.upneg-n})}}}x^{k}%
=\sum_{k\in\mathbb{N}}\dbinom{-n}{k}x^{k}.
\]
This proves Corollary \ref{cor.fps.anti-newton-binom-2}.
\end{proof}

We can now easily prove Newton's binomial formula:

\begin{proof}
[Proof of Theorem \ref{thm.fps.newton-binom}.]Let $n\in\mathbb{Z}$. We must
prove that $\left(  1+x\right)  ^{n}=\sum_{k\in\mathbb{N}}\dbinom{n}{k}x^{k}$.

If $n\in\mathbb{N}$, then this follows by comparing%
\begin{align*}
\left(  1+x\right)  ^{n}  &  =\left(  x+1\right)  ^{n}=\sum_{k=0}^{n}%
\dbinom{n}{k}x^{k}\underbrace{1^{n-k}}_{=1}\ \ \ \ \ \ \ \ \ \ \left(
\text{by the binomial theorem}\right) \\
&  =\sum_{k=0}^{n}\dbinom{n}{k}x^{k}%
\end{align*}
with%
\begin{align*}
\sum_{k\in\mathbb{N}}\dbinom{n}{k}x^{k}  &  =\sum_{k=0}^{n}\dbinom{n}{k}%
x^{k}+\sum_{k>n}\underbrace{\dbinom{n}{k}}_{\substack{=0\\\text{(by
Proposition \ref{prop.binom.0},}\\\text{since }n<k\text{)}}}x^{k}=\sum
_{k=0}^{n}\dbinom{n}{k}x^{k}+\underbrace{\sum_{k>n}0x^{k}}_{=0}\\
&  =\sum_{k=0}^{n}\dbinom{n}{k}x^{k}.
\end{align*}
Hence, for the rest of this proof, we WLOG assume that $n\notin\mathbb{N}$.

Hence, $n$ is a negative integer, so that $-n\in\mathbb{N}$. Thus, Corollary
\ref{cor.fps.anti-newton-binom-2} (applied to $-n$ instead of $n$) yields%
\[
\left(  1+x\right)  ^{-\left(  -n\right)  }=\sum_{k\in\mathbb{N}}%
\dbinom{-\left(  -n\right)  }{k}x^{k}.
\]
Since $-\left(  -n\right)  =n$, this rewrites as $\left(  1+x\right)
^{n}=\sum_{k\in\mathbb{N}}\dbinom{n}{k}x^{k}$. Thus, Theorem
\ref{thm.fps.newton-binom} is proven.
\end{proof}

We thus have a formula for $\left(  1+x\right)  ^{n}$ for each
\textbf{integer} $n$. We don't yet have such a formula for $\left(
1+x\right)  ^{1/2}$ (nor do we have a proper definition of $\left(
1+x\right)  ^{1/2}$), but this was clearly a step forward.

\subsubsection{Dividing by $x$}

Let us see how this all helps us justify our arguments in Section
\ref{sec.gf.exas}. Proposition \ref{prop.fps.invertible} justifies the
fractions that appear in (\ref{eq.sec.gf.exas.1.Fx=2}), but it does not
justify dividing by the FPS $2x$ in (\ref{eq.sec.gf.exas.2.-}), since the
constant term $\left[  x^{0}\right]  \left(  2x\right)  $ is surely not
invertible. And indeed, the FPS $2x$ is not invertible; the fraction
$\dfrac{1}{2x}$ is not a well-defined FPS.

However, it is easy to see directly which FPSs can be divided by $x$ (and thus
by $2x$, if $K=\mathbb{Q}$), and what it means to divide them by $x$. In fact,
Lemma \ref{lem.fps.xa} shows that multiplying an FPS by $x$ means moving all
its entries by one position to the right, and putting a $0$ into the newly
vacated starting position. Thus, it is rather clear what dividing by $x$
should be:

\begin{definition}
\label{def.fps.div-by-x}Let $a=\left(  a_{0},a_{1},a_{2},\ldots\right)  $ be
an FPS whose constant term $a_{0}$ is $0$. Then, $\dfrac{a}{x}$ is defined to
be the FPS $\left(  a_{1},a_{2},a_{3},\ldots\right)  $.
\end{definition}

The following is almost trivial:

\begin{proposition}
\label{prop.fps.div-by-x-inverts}Let $a\in K\left[  \left[  x\right]  \right]
$ and $b\in K\left[  \left[  x\right]  \right]  $ be two FPSs. Then, $a=xb$ if
and only if $\left(  \left[  x^{0}\right]  a=0\text{ and }b=\dfrac{a}%
{x}\right)  $.
\end{proposition}

\begin{proof}
Exercise.
\end{proof}

Having defined $\dfrac{a}{x}$ in Definition \ref{def.fps.div-by-x} (when $a$
has constant term $0$), we can also define $\dfrac{a}{2x}$ when $2$ is
invertible in $K$ (just set $\dfrac{a}{2x}=\dfrac{1}{2}\cdot\dfrac{a}{x}$).
Thus, the fraction $\dfrac{1\pm\sqrt{1-4x}}{2x}$ in (\ref{eq.sec.gf.exas.2.-+}%
) makes sense when the $\pm$ sign is a $-$ sign (but not when it is a $+$
sign), at least if we interpret the square root $\sqrt{1-4x}$ as $\sum
_{k\geq0}\dbinom{1/2}{k}\left(  -4x\right)  ^{k}$.

\subsubsection{A few lemmas}

Let us use this occasion to state a few simple lemmas (vaguely related to
Definition \ref{def.fps.div-by-x}) that will be used later on:

\begin{lemma}
\label{lem.fps.g=xh}Let $a\in K\left[  \left[  x\right]  \right]  $ be an FPS
with $\left[  x^{0}\right]  a=0$. Then, there exists an $h\in K\left[  \left[
x\right]  \right]  $ such that $a=xh$.
\end{lemma}

\begin{proof}
[Proof of Lemma \ref{lem.fps.g=xh}.]Write the FPS $a$ in the form $a=\left(
a_{0},a_{1},a_{2},\ldots\right)  $. Thus, $a_{0}=\left[  x^{0}\right]  a=0$.
Hence, the FPS $\dfrac{a}{x}$ is well-defined. Moreover, it is easy to see
that $a=x\cdot\dfrac{a}{x}$\ \ \ \ \footnote{\textit{Proof.} We have
$\dfrac{a}{x}=\left(  a_{1},a_{2},a_{3},\ldots\right)  $ (by Definition
\ref{def.fps.div-by-x}). Thus, Lemma \ref{lem.fps.xa} (applied to $\dfrac
{a}{x}$ and $a_{i+1}$ instead of $\mathbf{a}$ and $a_{i}$) yields%
\begin{align*}
x\cdot\dfrac{a}{x}  &  =\left(  0,a_{1},a_{2},a_{3},\ldots\right)  =\left(
a_{0},a_{1},a_{2},a_{3},\ldots\right)  \ \ \ \ \ \ \ \ \ \ \left(  \text{since
}0=a_{0}\right) \\
&  =\left(  a_{0},a_{1},a_{2},\ldots\right)  =a.
\end{align*}
Thus, $a=x\cdot\dfrac{a}{x}$.}. Hence, there exists an $h\in K\left[  \left[
x\right]  \right]  $ such that $a=xh$ (namely, $h=\dfrac{a}{x}$). This proves
Lemma \ref{lem.fps.g=xh}.
\end{proof}

\begin{lemma}
\label{lem.fps.first-n-coeffs-of-xna}Let $k\in\mathbb{N}$. Let $a\in K\left[
\left[  x\right]  \right]  $ be any FPS. Then, the first $k$ coefficients of
the FPS $x^{k}a$ are $0$.
\end{lemma}

\begin{proof}
[Proof of Lemma \ref{lem.fps.first-n-coeffs-of-xna}.]We must show that
$\left[  x^{m}\right]  \left(  x^{k}a\right)  =0$ for any nonnegative integer
$m<k$. But we can do this directly: If $m$ is a nonnegative integer such that
$m<k$, then (\ref{pf.thm.fps.ring.xn(ab)=2}) (applied to $x^{k}$, $a$ and $m$
instead of $\mathbf{a}$, $\mathbf{b}$ and $n$) yields%
\[
\left[  x^{m}\right]  \left(  x^{k}a\right)  =\sum_{i=0}^{m}%
\underbrace{\left[  x^{i}\right]  \left(  x^{k}\right)  }%
_{\substack{=0\\\text{(since }i\leq m<k\\\text{and thus }i\neq k\text{)}%
}}\cdot\left[  x^{m-i}\right]  a=\sum_{i=0}^{m}0\cdot\left[  x^{m-i}\right]
a=0,
\]
which is exactly what we wanted to show. Thus, Lemma
\ref{lem.fps.first-n-coeffs-of-xna} is proved.

(Alternatively, we could prove Lemma \ref{lem.fps.first-n-coeffs-of-xna} by
writing $a$ in the form $a=\left(  a_{0},a_{1},a_{2},\ldots\right)  $ and
observing that $x^{k}a=\left(  \underbrace{0,0,\ldots,0}_{k\text{ times}%
},a_{0},a_{1},a_{2},\ldots\right)  $. This follows by applying Lemma
\ref{lem.fps.xa} a total of $k$ times, or more formally by induction on $k$.)
\end{proof}

Lemma \ref{lem.fps.first-n-coeffs-of-xna} has a converse; here is a statement
that combines it with this converse:

\begin{lemma}
\label{lem.fps.muls-of-xn}Let $k\in\mathbb{N}$. Let $f\in K\left[  \left[
x\right]  \right]  $ be any FPS. Then, the first $k$ coefficients of the FPS
$f$ are $0$ if and only if $f$ is a multiple of $x^{k}$.
\end{lemma}

Here, we use the following notation:

\begin{definition}
Let $g\in K\left[  \left[  x\right]  \right]  $ be an FPS. Then, a
\emph{multiple} of $g$ means an FPS of the form $ga$ with $a\in K\left[
\left[  x\right]  \right]  $.
\end{definition}

(This is just a particular case of the usual concept of multiples in a
commutative ring.)

\begin{proof}
[Proof of Lemma \ref{lem.fps.muls-of-xn}.]The statement we are proving is an
\textquotedblleft if and only if\textquotedblright\ statement. We shall prove
its \textquotedblleft only if\textquotedblright\ (i.e., \textquotedblleft%
$\Longrightarrow$\textquotedblright) and its \textquotedblleft
if\textquotedblright\ (i.e., \textquotedblleft$\Longleftarrow$%
\textquotedblright) directions separately:

\begin{enumerate}
\item[$\Longrightarrow:$] Assume that the first $k$ coefficients of the FPS
$f$ are $0$. We must show that $f$ is a multiple of $x^{k}$.

Write $f$ as $f=\left(  f_{0},f_{1},f_{2},\ldots\right)  $. Then, the first
$k$ coefficients of the FPS $f$ are $f_{0},f_{1},\ldots,f_{k-1}$. Hence, these
$k$ coefficients $f_{0},f_{1},\ldots,f_{k-1}$ are $0$ (since we have assumed
that the first $k$ coefficients of the FPS $f$ are $0$). In other words,
$f_{n}=0$ for each $n\in\left\{  0,1,\ldots,k-1\right\}  $. Hence, $\sum
_{n=0}^{k-1}\underbrace{f_{n}}_{=0}x^{n}=\sum_{n=0}^{k-1}0x^{n}=0$.

Now,
\begin{align*}
f  &  =\left(  f_{0},f_{1},f_{2},\ldots\right)  =\sum_{n\in\mathbb{N}}%
f_{n}x^{n}=\underbrace{\sum_{n=0}^{k-1}f_{n}x^{n}}_{=0}+\sum_{n=k}^{\infty
}f_{n}\underbrace{x^{n}}_{\substack{=x^{k}x^{n-k}\\\text{(since }n\geq
k\text{)}}}=\sum_{n=k}^{\infty}f_{n}x^{k}x^{n-k}\\
&  =x^{k}\sum_{n=k}^{\infty}f_{n}x^{n-k}\ \ \ \ \ \ \ \ \ \ \left(  \text{by
Proposition \ref{prop.fps.sum-dist-inf}}\right)  .
\end{align*}
In other words, $f=x^{k}a$ for $a=\sum_{n=k}^{\infty}f_{n}x^{n-k}$. This shows
that $f$ is a multiple of $x^{k}$. Thus, the \textquotedblleft$\Longrightarrow
$\textquotedblright\ direction of Lemma \ref{lem.fps.muls-of-xn} is proved.

\item[$\Longleftarrow:$] Assume that $f$ is a multiple of $x^{k}$. In other
words, $f=x^{k}a$ for some $a\in K\left[  \left[  x\right]  \right]  $.
Consider this $a$. Now, Lemma \ref{lem.fps.first-n-coeffs-of-xna} yields that
the first $k$ coefficients of the FPS $x^{k}a$ are $0$. In other words, the
first $k$ coefficients of the FPS $f$ are $0$ (since $f=x^{k}a$). This proves
the \textquotedblleft$\Longleftarrow$\textquotedblright\ direction of Lemma
\ref{lem.fps.muls-of-xn}.
\end{enumerate}

The proof of Lemma \ref{lem.fps.muls-of-xn} is now complete, as both
directions have been proved.
\end{proof}

Another lemma that will prove its usefulness much later concerns FPSs that are
equal up until a certain coefficient. It says that if $f$ and $g$ are two FPSs
whose first $n+1$ coefficients agree (for a certain $n\in\mathbb{N}$), then
the same is true of the FPSs $af$ and $ag$ whenever $a$ is any further FPS. In
more detail:

\begin{lemma}
\label{lem.fps.prod.irlv.fg}Let $a,f,g\in K\left[  \left[  x\right]  \right]
$ be three FPSs. Let $n\in\mathbb{N}$. Assume that
\begin{equation}
\left[  x^{m}\right]  f=\left[  x^{m}\right]  g\ \ \ \ \ \ \ \ \ \ \text{for
each }m\in\left\{  0,1,\ldots,n\right\}  . \label{eq.lem.fps.prod.irlv.fg.ass}%
\end{equation}
Then,
\[
\left[  x^{m}\right]  \left(  af\right)  =\left[  x^{m}\right]  \left(
ag\right)  \ \ \ \ \ \ \ \ \ \ \text{for each }m\in\left\{  0,1,\ldots
,n\right\}  .
\]

\end{lemma}

\begin{proof}
[Proof of Lemma \ref{lem.fps.prod.irlv.fg}.]Let $m\in\left\{  0,1,\ldots
,n\right\}  $. Then, $m\leq n$. Hence, each $j\in\left\{  0,1,\ldots
,m\right\}  $ satisfies $j\leq m\leq n$ and thus $j\in\left\{  0,1,\ldots
,n\right\}  $ and therefore%
\begin{equation}
\left[  x^{j}\right]  f=\left[  x^{j}\right]  g
\label{pf.lem.fps.prod.irlv.fg.1}%
\end{equation}
(by (\ref{eq.lem.fps.prod.irlv.fg.ass}), applied to $j$ instead of $m$).
However, (\ref{pf.thm.fps.ring.xn(ab)=3}) (applied to $m$, $a$ and $f$ instead
of $n$, $\mathbf{a}$ and $\mathbf{b}$) yields%
\[
\left[  x^{m}\right]  \left(  af\right)  =\sum_{j=0}^{m}\left[  x^{m-j}%
\right]  a\cdot\underbrace{\left[  x^{j}\right]  f}_{\substack{=\left[
x^{j}\right]  g\\\text{(by (\ref{pf.lem.fps.prod.irlv.fg.1}))}}}=\sum
_{j=0}^{m}\left[  x^{m-j}\right]  a\cdot\left[  x^{j}\right]  g.
\]
On the other hand, (\ref{pf.thm.fps.ring.xn(ab)=3}) (applied to $m$, $a$ and
$g$ instead of $n$, $\mathbf{a}$ and $\mathbf{b}$) yields%
\[
\left[  x^{m}\right]  \left(  ag\right)  =\sum_{j=0}^{m}\left[  x^{m-j}%
\right]  a\cdot\left[  x^{j}\right]  g.
\]
Comparing these two equalities, we obtain $\left[  x^{m}\right]  \left(
af\right)  =\left[  x^{m}\right]  \left(  ag\right)  $. This proves Lemma
\ref{lem.fps.prod.irlv.fg}.
\end{proof}

A consequence of Lemma \ref{lem.fps.prod.irlv.fg} is the following fact:

\begin{lemma}
\label{lem.fps.prod.irlv.mul}Let $u,v\in K\left[  \left[  x\right]  \right]  $
be two FPSs such that $v$ is a multiple of $u$. Let $n\in\mathbb{N}$. Assume
that
\begin{equation}
\left[  x^{m}\right]  u=0\ \ \ \ \ \ \ \ \ \ \text{for each }m\in\left\{
0,1,\ldots,n\right\}  . \label{eq.lem.fps.prod.irlv.mul.ass}%
\end{equation}
Then,
\[
\left[  x^{m}\right]  v=0\ \ \ \ \ \ \ \ \ \ \text{for each }m\in\left\{
0,1,\ldots,n\right\}  .
\]

\end{lemma}

\begin{proof}
[Proof of Lemma \ref{lem.fps.prod.irlv.mul}.]We have assumed that $v$ is a
multiple of $u$. In other words, $v=ua$ for some $a\in K\left[  \left[
x\right]  \right]  $. Consider this $a$.

For each $m\in\left\{  0,1,\ldots,n\right\}  $, we have%
\begin{align*}
\left[  x^{m}\right]  u  &  =0\ \ \ \ \ \ \ \ \ \ \left(  \text{by
(\ref{eq.lem.fps.prod.irlv.mul.ass})}\right) \\
&  =\left[  x^{m}\right]  0\ \ \ \ \ \ \ \ \ \ \left(  \text{since the FPS
}0\text{ satisfies }\left[  x^{m}\right]  0=0\right)  .
\end{align*}
Hence, Lemma \ref{lem.fps.prod.irlv.fg} (applied to $f=u$ and $g=0$) yields
that%
\begin{equation}
\left[  x^{m}\right]  \left(  au\right)  =\left[  x^{m}\right]  \left(
a\cdot0\right)  \ \ \ \ \ \ \ \ \ \ \text{for each }m\in\left\{
0,1,\ldots,n\right\}  . \label{pf.lem.fps.prod.irlv.mul.2}%
\end{equation}
Now, for each $m\in\left\{  0,1,\ldots,n\right\}  $, we have%
\begin{align*}
\left[  x^{m}\right]  v  &  =\left[  x^{m}\right]  \left(  au\right)
\ \ \ \ \ \ \ \ \ \ \left(  \text{since }v=ua=au\right) \\
&  =\left[  x^{m}\right]  \underbrace{\left(  a\cdot0\right)  }_{=0}%
\ \ \ \ \ \ \ \ \ \ \left(  \text{by (\ref{pf.lem.fps.prod.irlv.mul.2}%
)}\right) \\
&  =\left[  x^{m}\right]  0=0.
\end{align*}
This proves Lemma \ref{lem.fps.prod.irlv.mul}.
\end{proof}

We can derive a further useful consequence from Lemma
\ref{lem.fps.prod.irlv.mul}:

\begin{lemma}
\label{lem.fps.prod.irlv.cong-mul}Let $a,b,c,d\in K\left[  \left[  x\right]
\right]  $ be four FPSs. Let $n\in\mathbb{N}$. Assume that
\begin{equation}
\left[  x^{m}\right]  a=\left[  x^{m}\right]  b\ \ \ \ \ \ \ \ \ \ \text{for
each }m\in\left\{  0,1,\ldots,n\right\}  .
\label{pf.lem.fps.prod.irlv.cong-mul.assab}%
\end{equation}
Assume further that%
\begin{equation}
\left[  x^{m}\right]  c=\left[  x^{m}\right]  d\ \ \ \ \ \ \ \ \ \ \text{for
each }m\in\left\{  0,1,\ldots,n\right\}  .
\label{pf.lem.fps.prod.irlv.cong-mul.asscd}%
\end{equation}
Then,
\[
\left[  x^{m}\right]  \left(  ac\right)  =\left[  x^{m}\right]  \left(
bd\right)  \ \ \ \ \ \ \ \ \ \ \text{for each }m\in\left\{  0,1,\ldots
,n\right\}  .
\]

\end{lemma}

\begin{proof}
[Proof of Lemma \ref{lem.fps.prod.irlv.cong-mul}.]For each $m\in\left\{
0,1,\ldots,n\right\}  $, we have%
\begin{align*}
\left[  x^{m}\right]  \left(  a-b\right)   &  =\left[  x^{m}\right]  a-\left[
x^{m}\right]  b\ \ \ \ \ \ \ \ \ \ \left(  \text{by
(\ref{pf.thm.fps.ring.xn(a-b)=})}\right) \\
&  =0\ \ \ \ \ \ \ \ \ \ \left(  \text{by
(\ref{pf.lem.fps.prod.irlv.cong-mul.assab})}\right)  .
\end{align*}
Moreover, the FPS $ac-bc$ is a multiple of $a-b$ (since $ac-bc=\left(
a-b\right)  c$). Hence, Lemma \ref{lem.fps.prod.irlv.mul} (applied to $u=a-b$
and $v=ac-bc$) shows that
\begin{equation}
\left[  x^{m}\right]  \left(  ac-bc\right)  =0\ \ \ \ \ \ \ \ \ \ \text{for
each }m\in\left\{  0,1,\ldots,n\right\}
\label{pf.lem.fps.prod.irlv.cong-mul.2}%
\end{equation}
(since we have $\left[  x^{m}\right]  \left(  a-b\right)  =0$ for each
$m\in\left\{  0,1,\ldots,n\right\}  $).

For each $m\in\left\{  0,1,\ldots,n\right\}  $, we have%
\begin{align*}
\left[  x^{m}\right]  \left(  c-d\right)   &  =\left[  x^{m}\right]  c-\left[
x^{m}\right]  d\ \ \ \ \ \ \ \ \ \ \left(  \text{by
(\ref{pf.thm.fps.ring.xn(a-b)=})}\right) \\
&  =0\ \ \ \ \ \ \ \ \ \ \left(  \text{by
(\ref{pf.lem.fps.prod.irlv.cong-mul.asscd})}\right)  .
\end{align*}
Moreover, the FPS $bc-bd$ is a multiple of $c-d$ (since $bc-bd=b\left(
c-d\right)  =\left(  c-d\right)  b$). Hence, Lemma \ref{lem.fps.prod.irlv.mul}
(applied to $u=c-d$ and $v=bc-bd$) shows that
\begin{equation}
\left[  x^{m}\right]  \left(  bc-bd\right)  =0\ \ \ \ \ \ \ \ \ \ \text{for
each }m\in\left\{  0,1,\ldots,n\right\}
\label{pf.lem.fps.prod.irlv.cong-mul.4}%
\end{equation}
(since we have $\left[  x^{m}\right]  \left(  c-d\right)  =0$ for each
$m\in\left\{  0,1,\ldots,n\right\}  $).

Now, let $m\in\left\{  0,1,\ldots,n\right\}  $. Then,
(\ref{pf.thm.fps.ring.xn(a-b)=}) yields $\left[  x^{m}\right]  \left(
ac-bc\right)  =\left[  x^{m}\right]  \left(  ac\right)  -\left[  x^{m}\right]
\left(  bc\right)  $. Comparing this with
(\ref{pf.lem.fps.prod.irlv.cong-mul.2}), we obtain $\left[  x^{m}\right]
\left(  ac\right)  -\left[  x^{m}\right]  \left(  bc\right)  =0$. In other
words, $\left[  x^{m}\right]  \left(  ac\right)  =\left[  x^{m}\right]
\left(  bc\right)  $. On the other hand, (\ref{pf.thm.fps.ring.xn(a-b)=})
yields $\left[  x^{m}\right]  \left(  bc-bd\right)  =\left[  x^{m}\right]
\left(  bc\right)  -\left[  x^{m}\right]  \left(  bd\right)  $. Comparing this
with (\ref{pf.lem.fps.prod.irlv.cong-mul.4}), we obtain $\left[  x^{m}\right]
\left(  bc\right)  -\left[  x^{m}\right]  \left(  bd\right)  =0$. In other
words, $\left[  x^{m}\right]  \left(  bc\right)  =\left[  x^{m}\right]
\left(  bd\right)  $. Hence, $\left[  x^{m}\right]  \left(  ac\right)
=\left[  x^{m}\right]  \left(  bc\right)  =\left[  x^{m}\right]  \left(
bd\right)  $. This proves Lemma \ref{lem.fps.prod.irlv.cong-mul}.
\end{proof}

\subsection{Polynomials}

\subsubsection{Definition}

Let us take a little side trip to relate FPSs to polynomials. As should be
clear enough from the definitions, we can think of an FPS as a
\textquotedblleft polynomial with (potentially) infinitely many nonzero
coefficients\textquotedblright. This can be easily made precise. Indeed, we
can \textbf{define} polynomials as FPSs that have only finitely many nonzero coefficients:

\begin{definition}
\label{def.fps.pol}\textbf{(a)} An FPS $a\in K\left[  \left[  x\right]
\right]  $ is said to be a \emph{polynomial} if all but finitely many
$n\in\mathbb{N}$ satisfy $\left[  x^{n}\right]  a=0$ (that is, if all but
finitely many coefficients of $a$ are $0$). \medskip

\textbf{(b)} We let $K\left[  x\right]  $ be the set of all polynomials $a\in
K\left[  \left[  x\right]  \right]  $. This set $K\left[  x\right]  $ is a
subring of $K\left[  \left[  x\right]  \right]  $ (according to Theorem
\ref{thm.fps.pol.ring} below), and is called the \emph{univariate polynomial
ring} over $K$.
\end{definition}

For example, $2+3x+7x^{5}$ is a polynomial, whereas $1+x+x^{2}+x^{3}+\cdots$
is not (unless $K$ is a trivial ring).

The definition of a \textquotedblleft polynomial\textquotedblright\ that you
have seen in your abstract algebra course might be superficially different
from that in Definition \ref{def.fps.pol}; but it necessarily is equivalent.
(In fact, Definition \ref{def.fps.pol} \textbf{(a)} can be restated as
\textquotedblleft a polynomial means a $K$-linear combination of the monomials
$x^{0},x^{1},x^{2},\ldots$\textquotedblright, and it is clear that the
monomials $x^{0},x^{1},x^{2},\ldots$ in $K\left[  \left[  x\right]  \right]  $
are $K$-linearly independent; thus, the polynomial ring $K\left[  x\right]  $
as we have defined it in Definition \ref{def.fps.pol} \textbf{(b)} is a free
$K$-module with basis $\left(  x^{0},x^{1},x^{2},\ldots\right)  $. The same is
true for the polynomial ring $K\left[  x\right]  $ that you know from abstract
algebra. Moreover, the rules for adding, subtracting and multiplying
polynomials known from abstract algebra agree with the formulas for
$\mathbf{a}+\mathbf{b}$, $\mathbf{a}-\mathbf{b}$ and $\mathbf{a}%
\cdot\mathbf{b}$ that we gave in Definition \ref{def.fps.ops}.)

We owe a theorem:

\begin{theorem}
\label{thm.fps.pol.ring}The set $K\left[  x\right]  $ is a subring of
$K\left[  \left[  x\right]  \right]  $ (that is, it is closed under addition,
subtraction and multiplication, and contains the zero $\underline{0}$ and the
unity $\underline{1}$) and is a $K$-submodule of $K\left[  \left[  x\right]
\right]  $ (that is, it is closed under addition and scaling by elements of
$K$).
\end{theorem}

\begin{proof}
[Proof of Theorem \ref{thm.fps.pol.ring} (sketched).]This is a rather easy
exercise. The hardest part is to show that $K\left[  x\right]  $ is closed
under multiplication. But this, too, is easy: Let $a,b\in K\left[  x\right]
$. Then, all but finitely many $n\in\mathbb{N}$ satisfy $\left[  x^{n}\right]
a=0$ (since $a\in K\left[  x\right]  $). In other words, there exists a finite
subset $I$ of $\mathbb{N}$ such that%
\begin{equation}
\left[  x^{i}\right]  a=0\text{ for all }i\in\mathbb{N}\setminus I.
\label{pf.thm.fps.pol.ring.xia=0}%
\end{equation}
Similarly, there exists a finite subset $J$ of $\mathbb{N}$ such that
\begin{equation}
\left[  x^{j}\right]  b=0\text{ for all }j\in\mathbb{N}\setminus J.
\label{pf.thm.fps.pol.ring.xjb=0}%
\end{equation}
Consider these $I$ and $J$. Now, let $S$ be the subset $\left\{
i+j\ \mid\ i\in I\text{ and }j\in J\right\}  $ of $\mathbb{N}$. This set $S$
is again finite (since $I$ and $J$ are finite), and we can easily see (using
(\ref{pf.thm.fps.ring.xn(ab)=2})) that%
\[
\left[  x^{n}\right]  \left(  ab\right)  =0\text{ for all }n\in\mathbb{N}%
\setminus S
\]
\footnote{\textit{Proof.} Let $n\in\mathbb{N}\setminus S$. Thus,
$n\in\mathbb{N}$ and $n\notin S$. Now, (\ref{pf.thm.fps.ring.xn(ab)=2}) yields
$\left[  x^{n}\right]  \left(  ab\right)  =\sum_{i=0}^{n}\left[  x^{i}\right]
a\cdot\left[  x^{n-i}\right]  b$. We shall next show that each addend in this
sum is $0$.
\par
Indeed, let $i\in\left\{  0,1,\ldots,n\right\}  $ be arbitrary. Let $j=n-i$.
Thus, $n=i+j$. Hence, we cannot have $i\in I$ and $j\in J$ simultaneously
(because if we did, then we would have $n=\underbrace{i}_{\in I}%
+\underbrace{j}_{\in J}\in S$ (by the definition of $S$), which would
contradict $n\notin S$). Hence, we must have $i\notin I$ or $j\notin J$ (or
both). In the former case, we have $\left[  x^{i}\right]  a=0$ (by
(\ref{pf.thm.fps.pol.ring.xia=0}), since $i\notin I$ entails $i\in
\mathbb{N}\setminus I$). In the latter case, we have $\left[  x^{j}\right]
b=0$ (by (\ref{pf.thm.fps.pol.ring.xjb=0}), since $j\notin J$ entails
$j\in\mathbb{N}\setminus J$). Thus, in either case, at least one of the two
coefficients $\left[  x^{i}\right]  a$ and $\left[  x^{j}\right]  b$ is $0$,
so that their product $\left[  x^{i}\right]  a\cdot\left[  x^{j}\right]  b$ is
$0$. In other words, $\left[  x^{i}\right]  a\cdot\left[  x^{n-i}\right]  b$
is $0$ (since $j=n-i$).
\par
Forget that we fixed $i$. We thus have shown that $\left[  x^{i}\right]
a\cdot\left[  x^{n-i}\right]  b$ is $0$ for each $i\in\left\{  0,1,\ldots
,n\right\}  $. In other words, all addends of the sum $\sum_{i=0}^{n}\left[
x^{i}\right]  a\cdot\left[  x^{n-i}\right]  b$ are $0$. Hence, the whole sum
is $0$. In other words, $\left[  x^{n}\right]  \left(  ab\right)  =0$ (since
$\left[  x^{n}\right]  \left(  ab\right)  =\sum_{i=0}^{n}\left[  x^{i}\right]
a\cdot\left[  x^{n-i}\right]  b$), qed.}. Thus, all but finitely many
$n\in\mathbb{N}$ satisfy $\left[  x^{n}\right]  \left(  ab\right)  =0$ (since
$S$ is finite). This shows that $ab\in K\left[  x\right]  $. Hence, we have
shown that $K\left[  x\right]  $ is closed under multiplication. The remaining
claims of Theorem \ref{thm.fps.pol.ring} are similar but easier.
\end{proof}

\subsubsection{Reminders on rings and $K$-algebras}

As we now know, polynomials are just a special case of FPSs. However, they
have some features that FPSs don't have in general. The most important of
these features is \emph{substitution}. To wit, we can substitute an element of
$K$, or more generally an element of any $K$-algebra, into a polynomial (but
generally not into an FPS). Before we explain how, let us recall the notions
of rings and $K$-algebras (see \cite[Chapter 2 and \S 3.11]{23wa} for details
and more about them):

\begin{definition}
\label{def.alg.ring}The notion of a \emph{ring} (also known as a
\emph{noncommutative ring}) is defined in the exact same way as we defined the
notion of a commutative ring in Definition \ref{def.alg.commring}, except that
the \textquotedblleft Commutativity of multiplication\textquotedblright\ axiom
is removed.
\end{definition}

Examples of noncommutative rings\footnote{Note that the word \textquotedblleft
noncommutative ring\textquotedblright\ does not imply that the ring is not
commutative; it merely means that commutativity is not required. Thus, any
commutative ring is a noncommutative ring.} abound in linear algebra:

\begin{itemize}
\item For any $n\in\mathbb{N}$, the matrix ring $\mathbb{R}^{n\times n}$ (that
is, the ring of all $n\times n$-matrices with real entries) is a ring. This
ring is commutative if $n\leq1$, but not if $n>1$.

More generally, if $K$ is any ring (commutative or not), then the matrix ring
$K^{n\times n}$ is a ring for every $n\in\mathbb{N}$.

\item The ring $\mathbb{H}$ of
\href{https://en.wikipedia.org/wiki/Quaternion}{quaternions} is a ring that is
not commutative.

\item If $M$ is an abelian group, then the ring of all endomorphisms of $M$
(that is, the ring of all $\mathbb{Z}$-linear maps from $M$ to $M$) is a
noncommutative ring. (Its multiplication is composition of endomorphisms.)
\end{itemize}

Next, let us recall the notion of a $K$-algebra (\cite[\S 3.11]{23wa}). Recall
that $K$ is a fixed commutative ring.

\begin{definition}
\label{def.alg.Kalg}A $K$\emph{-algebra} is a set $A$ equipped with four maps%
\begin{align*}
\oplus &  :A\times A\rightarrow A,\\
\ominus &  :A\times A\rightarrow A,\\
\odot &  :A\times A\rightarrow A,\\
\rightharpoonup &  :K\times A\rightarrow A
\end{align*}
and two elements $\overrightarrow{0}\in A$ and $\overrightarrow{1}\in A$
satisfying the following properties:

\begin{enumerate}
\item The set $A$, equipped with the maps $\oplus$, $\ominus$ and $\odot$ and
the two elements $\overrightarrow{0}$ and $\overrightarrow{1}$, is a
(noncommutative) ring.

\item The set $A$, equipped with the maps $\oplus$, $\ominus$ and
$\rightharpoonup$ and the element $\overrightarrow{0}$, is a $K$-module.

\item We have%
\begin{equation}
\lambda\rightharpoonup\left(  a\odot b\right)  =\left(  \lambda\rightharpoonup
a\right)  \odot b=a\odot\left(  \lambda\rightharpoonup b\right)
\label{eq.def.alg.Kalg.scaleinv}%
\end{equation}
for all $\lambda\in K$ and $a,b\in A$.
\end{enumerate}

(Thus, in a nutshell, a $K$-algebra is a set $A$ that is simultaneously a ring
and a $K$-module, with the property that the ring $A$ and the $K$-module $A$
have the same addition, the same subtraction and the same zero, and satisfy
the additional compatibility property (\ref{eq.def.alg.Kalg.scaleinv}).)

Consequently, a $K$-algebra is automatically a ring and a $K$-module. Thus,
all the notations and shorthands that we have introduced for rings and for
$K$-modules will also be used for $K$-algebras. For example, if $A$ is a
$K$-algebra, then both maps $\odot:A\times A\rightarrow A$ and $\left.
\rightharpoonup\right.  :K\times A\rightarrow A$ will be denoted by $\cdot$
unless there is a risk of confusion. (There is rarely a risk of confusion,
since the two maps act on different inputs: $a\cdot b$ means $a\odot b$ if $a$
belongs to $A$, and means $a\rightharpoonup b$ if $a$ belongs to $K$. Often,
even when an element $a$ belongs to both $A$ and $K$, the elements $a\odot b$
and $a\rightharpoonup b$ are equal, so confusion cannot arise.)
\end{definition}

Examples of $K$-algebras include:

\begin{itemize}
\item the ring $K$ itself;

\item the ring $K\left[  \left[  x\right]  \right]  $ of FPSs (we have defined
the relevant maps in Definition \ref{def.fps.ops}, and claimed the relevant
properties in Theorem \ref{thm.fps.ring});

\item its subring $K\left[  x\right]  $ (all its maps are inherited from
$K\left[  \left[  x\right]  \right]  $);

\item the matrix ring $K^{n\times n}$ for each $n\in\mathbb{N}$;

\item any quotient ring of $K$ (that is, any ring of the form $K/I$ where $I$
is an ideal of $K$);

\item any commutative ring that contains $K$ as a subring.
\end{itemize}

Note that the axiom (\ref{eq.def.alg.Kalg.scaleinv}) in the definition of
$K$-algebra can be rewritten as%
\[
\lambda\left(  ab\right)  =\left(  \lambda a\right)  b=a\left(  \lambda
b\right)  \ \ \ \ \ \ \ \ \ \ \text{for all }\lambda\in K\text{ and }a,b\in A
\]
using our conventions (to write $ab$ for $a\odot b$ and to write $\lambda c$
for $\lambda\rightharpoonup c$). It says that scaling a product in $A$ by a
scalar $\lambda\in K$ is equivalent to scaling either of its two factors by
$\lambda$.

\subsubsection{Evaluation aka substitution into polynomials}

We can now define what it means to substitute an element of a $K$-algebra into
a polynomial:

\begin{definition}
\label{def.pol.subs}Let $f\in K\left[  x\right]  $ be a polynomial. Let $A$ be
any $K$-algebra. Let $a\in A$ be any element. We then define an element
$f\left[  a\right]  \in A$ as follows:

Write $f$ in the form $f=\sum_{n\in\mathbb{N}}f_{n}x^{n}$ with $f_{0}%
,f_{1},f_{2},\ldots\in K$. (That is, $f_{n}=\left[  x^{n}\right]  f$ for each
$n\in\mathbb{N}$.) Then, set%
\[
f\left[  a\right]  :=\sum_{n\in\mathbb{N}}f_{n}a^{n}.
\]
(This sum is essentially finite, since $f$ is a polynomial.)

The element $f\left[  a\right]  $ is also denoted by $f\circ a$ or by
$f\left(  a\right)  $, and is called the \emph{value} of $f$ at $a$ (or the
\emph{evaluation} of $f$ at $a$, or the \emph{result of substituting }$a$ for
$x$ in $f$).
\end{definition}

The most popular notation for the value $f\left[  a\right]  $ that we have
just defined is $f\left(  a\right)  $. Unfortunately, this notation can lead
to ambiguities (for example, $f\left(  x+1\right)  $ could mean either the
value of $f$ at $x+1$, or the product of $f$ with $x+1$). By writing $f\left[
a\right]  $ or $f\circ a$ instead, I will avoid these ambiguities, although I
will not be fully consistent and will occasionally revert to writing $f\left(
a\right)  $ when confusion is unlikely.

For example, if $f=4x^{3}+2x+7$, then $f\left[  a\right]  =4a^{3}+2a+7$. For
another example, if $f=\left(  x+5\right)  ^{3}$, then $f\left[  a\right]
=\left(  a+5\right)  ^{3}$ (although this is not obvious; it follows from
Theorem \ref{thm.pol.eval.a+b} below). \medskip

If $f$ and $g$ are two polynomials in $K\left[  x\right]  $, then the value
$f\left[  g\right]  =f\circ g$ (this is the value of $f$ at $g$; it is
well-defined because $K\left[  x\right]  $ is a $K$-algebra) is also known as
the \emph{composition} of $f$ with $g$. We note that any polynomial $f\in
K\left[  x\right]  $ satisfies
\begin{align*}
f\left[  x\right]   &  =f\ \ \ \ \ \ \ \ \ \ \text{and}\\
f\left[  0\right]   &  =\left[  x^{0}\right]  f=\left(  \text{the constant
term of }f\right)  \ \ \ \ \ \ \ \ \ \ \text{and}\\
f\left[  1\right]   &  =\left[  x^{0}\right]  f+\left[  x^{1}\right]
f+\left[  x^{2}\right]  f+\cdots=\left(  \text{the sum of all coefficients of
}f\right)  .
\end{align*}
It is fairly common to write $f\left[  x\right]  $ instead of $f$ for a
polynomial, just to stress the fact that it is a polynomial in an
indeterminate called $x$. The equality $f\left[  x\right]  =f$ justifies this.
\medskip

Definition \ref{def.pol.subs} is rather versatile. For example, if
$f\in\mathbb{Z}\left[  x\right]  $ is a polynomial with integer coefficients,
then it allows evaluating $f$ at integers, at complex numbers, at residue
classes in $\mathbb{Z}/n$, at square matrices, at other polynomials and at
FPSs. More generally, a polynomial $f\in\mathbb{Z}\left[  x\right]  $ can be
evaluated at any element of any ring, since any ring is automatically a
$\mathbb{Z}$-algebra. Evaluating polynomials at square matrices is an
important idea in linear algebra (e.g.,
\href{https://en.wikipedia.org/wiki/Cayley-Hamilton_theorem}{the
Cayley--Hamilton theorem} is concerned with the characteristic polynomial of a
square matrix, evaluated at this matrix itself).

The following theorem surveys the most basic properties of values of polynomials:

\begin{theorem}
\label{thm.pol.eval.a+b}Let $A$ be a $K$-algebra. Let $a\in A$. Then: \medskip

\textbf{(a)} Any $f,g\in K\left[  x\right]  $ satisfy%
\[
\left(  f+g\right)  \left[  a\right]  =f\left[  a\right]  +g\left[  a\right]
\ \ \ \ \ \ \ \ \ \ \text{and}\ \ \ \ \ \ \ \ \ \ \left(  fg\right)  \left[
a\right]  =f\left[  a\right]  \cdot g\left[  a\right]  .
\]

\textbf{(b)} Any $\lambda\in K$ and $f\in K\left[  x\right]  $ satisfy%
\[
\left(  \lambda f\right)  \left[  a\right]  =\lambda\cdot f\left[  a\right]
.
\]

\textbf{(c)} Any $\lambda\in K$ satisfies $\underline{\lambda}\left[
a\right]  =\lambda\cdot1_{A}$, where $1_{A}$ is the unity of the ring $A$.
(This is often written as \textquotedblleft$\underline{\lambda}\left[
a\right]  =\lambda$\textquotedblright, but keep in mind that the
\textquotedblleft$\lambda$\textquotedblright\ on the right hand side of this
equality is understood to be \textquotedblleft coerced into $A$%
\textquotedblright, so it actually means \textquotedblleft the element of $A$
corresponding to $\lambda$\textquotedblright, which is $\lambda\cdot1_{A}$.)
\medskip

\textbf{(d)} We have $x\left[  a\right]  =a$. \medskip

\textbf{(e)} We have $x^{i}\left[  a\right]  =a^{i}$ for each $i\in\mathbb{N}%
$. \medskip

\textbf{(f)} Any $f,g\in K\left[  x\right]  $ satisfy $f\left[  g\left[
a\right]  \right]  =\left(  f\left[  g\right]  \right)  \left[  a\right]  $.
\end{theorem}

\begin{proof}
See \cite[Theorem 7.6.3]{19s} for parts \textbf{(a)}, \textbf{(b)},
\textbf{(c)}, \textbf{(d)} and \textbf{(e)}. Part \textbf{(f)} is
\cite[Proposition 7.6.14]{19s}.
\end{proof}

\subsection{Substitution and evaluation of power series}

\subsubsection{Defining substitution}

Definition \ref{def.pol.subs} shows that if $f\in K\left[  x\right]  $ is a
polynomial, then almost anything (to be more precise: any element of a
$K$-algebra) can be substituted into $f$.

In contrast, if $f\in K\left[  \left[  x\right]  \right]  $ is a FPS, then
there are far fewer things that can be substituted into $f$. Even elements of
$K$ itself cannot always be substituted into $f$. For example, if we try to
substitute $1$ for $x$ in the FPS $1+x+x^{2}+x^{3}+\cdots$, then we get%
\[
1+1+1^{2}+1^{3}+\cdots=1+1+1+1+\cdots,
\]
which is undefined. Real analysis can help make sense of certain values of
FPSs (for example, substituting $\dfrac{1}{2}$ for $x$ into the FPS
$1+x+x^{2}+x^{3}+\cdots$ yields the convergent series $1+\dfrac{1}{2}%
+\dfrac{1}{2^{2}}+\dfrac{1}{2^{3}}+\cdots=2$), but this is subtle and specific
to certain numbers and certain FPSs.\footnote{For instance, it is not hard to
see that there is no nonzero complex number that can be substituted into the
FPS $\sum_{n\in\mathbb{N}}n!x^{n}$ to obtain a convergent result. Thus, even
though some complex numbers can be substituted into some FPSs, there is no
complex number other than $0$ that can be substituted into \textbf{every}
FPS.}

Thus, polynomials have a clear advantage over FPSs.

However, let us not give up on FPSs yet. \textbf{Some} things can be
substituted into an FPS. For example:

\begin{itemize}
\item We can always substitute $0$ for $x$ in an FPS $a_{0}+a_{1}x+a_{2}%
x^{2}+a_{3}x^{3}+\cdots$. The result is%
\[
a_{0}+a_{1}0+a_{2}0^{2}+a_{3}0^{3}+\cdots=a_{0}+0+0+0+\cdots=a_{0}.
\]

\item We can always substitute $x$ for $x$ in an FPS $a_{0}+a_{1}x+a_{2}%
x^{2}+a_{3}x^{3}+\cdots$. The result is the same FPS $a_{0}+a_{1}x+a_{2}%
x^{2}+a_{3}x^{3}+\cdots$ that we started with (obviously).

\item We can always substitute $2x$ for $x$ in an FPS $a_{0}+a_{1}x+a_{2}%
x^{2}+a_{3}x^{3}+\cdots$. The result is%
\[
a_{0}+a_{1}\left(  2x\right)  +a_{2}\left(  2x\right)  ^{2}+a_{3}\left(
2x\right)  ^{3}+\cdots=a_{0}+2a_{1}x+4a_{2}x^{2}+8a_{3}x^{3}+\cdots.
\]

\item We can always substitute $x^{2}+x$ for $x$ in an FPS $a_{0}+a_{1}%
x+a_{2}x^{2}+a_{3}x^{3}+\cdots$. This is less obvious, so let me explain why.
If we try to substitute $x^{2}+x$ for $x$ in an FPS $a_{0}+a_{1}x+a_{2}%
x^{2}+a_{3}x^{3}+\cdots$, then we obtain%
\begin{align*}
&  \left(  a_{0}+a_{1}x+a_{2}x^{2}+a_{3}x^{3}+\cdots\right)  \left[
x+x^{2}\right] \\
&  =a_{0}+a_{1}\left(  x+x^{2}\right)  +a_{2}\left(  x+x^{2}\right)
^{2}+a_{3}\left(  x+x^{2}\right)  ^{3}+\cdots\\
&  =a_{0}+a_{1}\left(  x+x^{2}\right)  +a_{2}\left(  x^{2}+2x^{3}%
+x^{4}\right)  +a_{3}\left(  x^{3}+3x^{4}+3x^{5}+x^{6}\right)  +\cdots\\
&  =a_{0}+a_{1}x+\left(  a_{1}+a_{2}\right)  x^{2}+\left(  2a_{2}%
+a_{3}\right)  x^{3}+\left(  a_{2}+3a_{3}+a_{4}\right)  x^{4}+\cdots.
\end{align*}
I claim that the right hand side here is well-defined. To prove this, I need
to show that for each $n\in\mathbb{N}$, the coefficient of $x^{n}$ on this
right hand side is a \textbf{finite} sum of $a_{i}$'s. Indeed, fix
$n\in\mathbb{N}$. Recall that the right hand side is obtained by expanding the
infinite sum%
\[
a_{0}+a_{1}\left(  x+x^{2}\right)  +a_{2}\left(  x+x^{2}\right)  ^{2}%
+a_{3}\left(  x+x^{2}\right)  ^{3}+\cdots.
\]
Only the first $n+1$ addends of this infinite sum (i.e., only the addends
$a_{k}\left(  x+x^{2}\right)  ^{k}$ with $k\leq n$) can contribute to the
coefficient of $x^{n}$, since any of the remaining addends is a multiple of
$x^{n+1}$ (because it has the form $a_{k}\left(  x+x^{2}\right)  ^{k}%
=a_{k}\left(  x\left(  1+x\right)  \right)  ^{k}=a_{k}x^{k}\left(  1+x\right)
^{k}$ with $k\geq n+1$) and thus has a zero coefficient of $x^{n}$. Hence, the
coefficient of $x^{n}$ in this infinite sum equals the coefficient of $x^{n}$
in the \textbf{finite} sum
\[
a_{0}+a_{1}\left(  x+x^{2}\right)  +a_{2}\left(  x+x^{2}\right)  ^{2}%
+a_{3}\left(  x+x^{2}\right)  ^{3}+\cdots+a_{n}\left(  x+x^{2}\right)  ^{n}.
\]
But the latter coefficient is clearly a \textbf{finite} sum of $a_{i}$'s.
Thus, my claim is proved, and it follows that the result of substituting
$x^{2}+x$ for $x$ in an FPS $a_{0}+a_{1}x+a_{2}x^{2}+a_{3}x^{3}+\cdots$ is well-defined.
\end{itemize}

The idea of the last example can be generalized; there was nothing special
about $x+x^{2}$ that we used other than the fact that $x+x^{2}$ is a multiple
of $x$ (that is, an FPS whose constant term is $0$). Thus, generalizing our
reasoning from this example, we can convince ourselves that any FPS $g$ that
is a multiple of $x$ (that is, whose constant term is $0$) can be substituted
into any FPS. Let us introduce a notation for this, exactly like we did for
substituting things into polynomials:

\begin{definition}
\label{def.fps.subs}Let $f$ and $g$ be two FPSs in $K\left[  \left[  x\right]
\right]  $. Assume that $\left[  x^{0}\right]  g=0$ (that is, $g=g_{1}%
x^{1}+g_{2}x^{2}+g_{3}x^{3}+\cdots$ for some $g_{1},g_{2},g_{3},\ldots\in K$).

We then define an FPS $f\left[  g\right]  \in K\left[  \left[  x\right]
\right]  $ as follows:

Write $f$ in the form $f=\sum_{n\in\mathbb{N}}f_{n}x^{n}$ with $f_{0}%
,f_{1},f_{2},\ldots\in K$. (That is, $f_{n}=\left[  x^{n}\right]  f$ for each
$n\in\mathbb{N}$.) Then, set%
\begin{equation}
f\left[  g\right]  :=\sum_{n\in\mathbb{N}}f_{n}g^{n}.
\label{eq.def.fps.subs.eq}%
\end{equation}
(This sum is well-defined, as we will see in Proposition
\ref{prop.fps.subs.wd} \textbf{(b)} below.)

This FPS $f\left[  g\right]  $ is also denoted by $f\circ g$, and is called
the \emph{composition} of $f$ with $g$, or the result of \emph{substituting}
$g$ for $x$ in $f$.
\end{definition}

Once again, it is not uncommon to see this FPS $f\left[  g\right]  $ denoted
by $f\left(  g\right)  $, but I will eschew the latter notation (since it can
be confused with a product).

In order to prove that Definition \ref{def.fps.subs} makes sense, we need to
ensure that the infinite sum $\sum_{n\in\mathbb{N}}f_{n}g^{n}$ in
(\ref{eq.def.fps.subs.eq}) is well-defined. The proof of this fact is
analogous to the reasoning I used in the last example; let me present it again
in the general case:

\begin{proposition}
\label{prop.fps.subs.wd}Let $f$ and $g$ be two FPSs in $K\left[  \left[
x\right]  \right]  $. Assume that $\left[  x^{0}\right]  g=0$. Write $f$ in
the form $f=\sum_{n\in\mathbb{N}}f_{n}x^{n}$ with $f_{0},f_{1},f_{2},\ldots\in
K$. Then: \medskip

\textbf{(a)} For each $n\in\mathbb{N}$, the first $n$ coefficients of the FPS
$g^{n}$ are $0$. \medskip

\textbf{(b)} The sum $\sum_{n\in\mathbb{N}}f_{n}g^{n}$ is well-defined, i.e.,
the family $\left(  f_{n}g^{n}\right)  _{n\in\mathbb{N}}$ is summable.
\medskip

\textbf{(c)} We have $\left[  x^{0}\right]  \left(  \sum_{n\in\mathbb{N}}%
f_{n}g^{n}\right)  =f_{0}$.
\end{proposition}

\begin{proof}
[Proof of Proposition \ref{prop.fps.subs.wd}.]\textbf{(a)} This is easily
proved by induction on $n$. Here is a shorter alternative argument:

The FPS $g$ has constant term $\left[  x^{0}\right]  g=0$. Hence, Lemma
\ref{lem.fps.g=xh} (applied to $a=g$) yields that there exists an $h\in
K\left[  \left[  x\right]  \right]  $ such that $g=xh$. Consider this $h$.

Now, let $n\in\mathbb{N}$. From $g=xh$, we obtain $g^{n}=\left(  xh\right)
^{n}=x^{n}h^{n}$. However, Lemma \ref{lem.fps.first-n-coeffs-of-xna} (applied
to $k=n$ and $a=h^{n}$) yields that the first $n$ coefficients of the FPS
$x^{n}h^{n}$ are $0$. In other words, the first $n$ coefficients of the FPS
$g^{n}$ are $0$ (since $g^{n}=x^{n}h^{n}$). Thus, Proposition
\ref{prop.fps.subs.wd} \textbf{(a)} is proved. \medskip

\textbf{(b)} This follows from part \textbf{(a)}. Here are the details.

We must prove that the family $\left(  f_{n}g^{n}\right)  _{n\in\mathbb{N}}$
is summable. In other words, we must prove that the family $\left(  f_{i}%
g^{i}\right)  _{i\in\mathbb{N}}$ is summable (since $\left(  f_{i}%
g^{i}\right)  _{i\in\mathbb{N}}=\left(  f_{n}g^{n}\right)  _{n\in\mathbb{N}}%
$). In other words, we must prove that for each $n\in\mathbb{N}$, all but
finitely many $i\in\mathbb{N}$ satisfy $\left[  x^{n}\right]  \left(
f_{i}g^{i}\right)  =0$ (by the definition of \textquotedblleft
summable\textquotedblright). So let us prove this.

Fix $n\in\mathbb{N}$. We must prove that all but finitely many $i\in
\mathbb{N}$ satisfy $\left[  x^{n}\right]  \left(  f_{i}g^{i}\right)  =0$.

Indeed, let $i\in\mathbb{N}$ satisfy $i>n$. Then, $n<i$. Now, the first $i$
coefficients of the FPS $g^{i}$ are $0$ (by Proposition \ref{prop.fps.subs.wd}
\textbf{(a)}, applied to $i$ instead of $n$). However, the coefficient
$\left[  x^{n}\right]  \left(  g^{i}\right)  $ of $g^{i}$ is one of these
first $i$ coefficients (because $n<i$). Thus, this coefficient $\left[
x^{n}\right]  \left(  g^{i}\right)  $ must be $0$. Now, $f_{i}\in K$; thus,
(\ref{pf.thm.fps.ring.xn(la)=}) (applied to $\lambda=f_{i}$ and $\mathbf{a}%
=g^{i}$) yields $\left[  x^{n}\right]  \left(  f_{i}g^{i}\right)  =f_{i}%
\cdot\underbrace{\left[  x^{n}\right]  \left(  g^{i}\right)  }_{=0}=0$.

Forget that we fixed $i$. We thus have shown that all $i\in\mathbb{N}$
satisfying $i>n$ satisfy $\left[  x^{n}\right]  \left(  f_{i}g^{i}\right)
=0$. Hence, all but finitely many $i\in\mathbb{N}$ satisfy $\left[
x^{n}\right]  \left(  f_{i}g^{i}\right)  =0$ (because all but finitely many
$i\in\mathbb{N}$ satisfy $i>n$). This is precisely what we wanted to prove.
Thus, Proposition \ref{prop.fps.subs.wd} \textbf{(b)} is proved. \medskip

\textbf{(c)} Let $n$ be a positive integer. We shall first show that $\left[
x^{0}\right]  \left(  f_{n}g^{n}\right)  =0$.

Indeed, Proposition \ref{prop.fps.subs.wd} \textbf{(a)} shows that the first
$n$ coefficients of the FPS $g^{n}$ are $0$. However, the coefficient $\left[
x^{0}\right]  \left(  g^{n}\right)  $ is one of these first $n$ coefficients
(since $n$ is positive). Thus, this coefficient $\left[  x^{0}\right]  \left(
g^{n}\right)  $ must be $0$. Now, $f_{n}\in K$; thus,
(\ref{pf.thm.fps.ring.xn(la)=}) (applied to $f_{n}$, $g^{n}$ and $0$ instead
of $\lambda$, $\mathbf{a}$ and $n$) yields $\left[  x^{0}\right]  \left(
f_{n}g^{n}\right)  =f_{n}\cdot\underbrace{\left[  x^{0}\right]  \left(
g^{n}\right)  }_{=0}=0$.

Forget that we fixed $n$. We thus have shown that
\begin{equation}
\left[  x^{0}\right]  \left(  f_{n}g^{n}\right)
=0\ \ \ \ \ \ \ \ \ \ \text{for each positive integer }n.
\label{pf.prop.fps.subs.wd.c.1}%
\end{equation}

Now,%
\begin{align*}
\left[  x^{0}\right]  \left(  \sum_{n\in\mathbb{N}}f_{n}g^{n}\right)   &
=\sum_{n\in\mathbb{N}}\left[  x^{0}\right]  \left(  f_{n}g^{n}\right)
\ \ \ \ \ \ \ \ \ \ \left(  \text{by (\ref{eq.def.fps.summable.sum})}\right)
\\
&  =\left[  x^{0}\right]  \left(  f_{0}\underbrace{g^{0}}_{=\underline{1}%
}\right)  +\sum_{n>0}\underbrace{\left[  x^{0}\right]  \left(  f_{n}%
g^{n}\right)  }_{\substack{=0\\\text{(by (\ref{pf.prop.fps.subs.wd.c.1}))}}}\\
&  \ \ \ \ \ \ \ \ \ \ \ \ \ \ \ \ \ \ \ \ \left(
\begin{array}
[c]{c}%
\text{here, we have split off the}\\
\text{addend for }n=0\text{ from the sum}%
\end{array}
\right) \\
&  =\left[  x^{0}\right]  \underbrace{\left(  f_{0}\underline{1}\right)
}_{\substack{=f_{0}\left(  1,0,0,0,\ldots\right)  \\=\left(  f_{0}%
,0,0,0,\ldots\right)  }}+\underbrace{\sum_{n>0}0}_{=0}=\left[  x^{0}\right]
\left(  f_{0},0,0,0,\ldots\right)  =f_{0}.
\end{align*}
This proves Proposition \ref{prop.fps.subs.wd} \textbf{(c)}.
\end{proof}

\begin{example}
\label{exa.fps.subs.fibonacci}The FPS $x+x^{2}$ has constant term $\left[
x^{0}\right]  \left(  x+x^{2}\right)  =0$. Hence, according to Definition
\ref{def.fps.subs}, we can substitute it for $x$ into $1+x+x^{2}+x^{3}+\cdots
$. The result is%
\begin{align*}
&  \left(  1+x+x^{2}+x^{3}+\cdots\right)  \left[  x+x^{2}\right] \\
&  =1+\left(  x+x^{2}\right)  +\left(  x+x^{2}\right)  ^{2}+\left(
x+x^{2}\right)  ^{3}+\left(  x+x^{2}\right)  ^{4}+\left(  x+x^{2}\right)
^{5}+\cdots\\
&  =1+x+2x^{2}+3x^{3}+5x^{4}+8x^{5}+\cdots.
\end{align*}
The right hand side appears to be $f_{1}+f_{2}x+f_{3}x^{2}+f_{4}x^{3}+\cdots$,
where $\left(  f_{0},f_{1},f_{2},\ldots\right)  $ is the Fibonacci sequence
(as defined in Section \ref{sec.gf.exas}). Let me show that this is indeed the case.

In Example 1 in Section \ref{sec.gf.exas}, we had shown that%
\[
f_{0}+f_{1}x+f_{2}x^{2}+f_{3}x^{3}+\cdots=\dfrac{x}{1-x-x^{2}}.
\]
Thus,%
\begin{align*}
\dfrac{x}{1-x-x^{2}}  &  =\underbrace{f_{0}}_{=0}+\,f_{1}x+f_{2}x^{2}%
+f_{3}x^{3}+\cdots\\
&  =f_{1}x+f_{2}x^{2}+f_{3}x^{3}+f_{4}x^{4}+\cdots\\
&  =x\left(  f_{1}+f_{2}x+f_{3}x^{2}+f_{4}x^{3}+\cdots\right)  .
\end{align*}
Cancelling $x$ from this equality (this is indeed allowed -- make sure you
understand why!), we obtain%
\[
\dfrac{1}{1-x-x^{2}}=f_{1}+f_{2}x+f_{3}x^{2}+f_{4}x^{3}+\cdots.
\]
However, it appears reasonable to expect that
\begin{equation}
\dfrac{1}{1-x-x^{2}}=\dfrac{1}{1-x}\left[  x+x^{2}\right]  ,
\label{eq.exa.fps.subs.fibonacci.subs-plausible}%
\end{equation}
because substituting $x+x^{2}$ for $x$ in the expression $\dfrac{1}{1-x}$
results in $\dfrac{1}{1-x-x^{2}}$. This is plausible but not obvious -- after
all, we defined $\dfrac{1}{1-x}\left[  x+x^{2}\right]  $ to be the result of
substituting $x+x^{2}$ for $x$ into the \textbf{expanded} version of
$\dfrac{1}{1-x}$ (which is $1+x+x^{2}+x^{3}+\cdots$), not into the fractional
expression $\dfrac{1}{1-x}$. Nevertheless,
(\ref{eq.exa.fps.subs.fibonacci.subs-plausible}) is true (and will soon be
proved). If we take this fact for granted, then our claim easily follows:%
\begin{align*}
f_{1}+f_{2}x+f_{3}x^{2}+f_{4}x^{3}+\cdots &  =\dfrac{1}{1-x-x^{2}}=\dfrac
{1}{1-x}\left[  x+x^{2}\right] \\
&  =\left(  1+x+x^{2}+x^{3}+\cdots\right)  \left[  x+x^{2}\right]
\end{align*}
(since $\dfrac{1}{1-x}=1+x+x^{2}+x^{3}+\cdots$).
\end{example}

\subsubsection{Laws of substitution}

The plausible but nontrivial statement
(\ref{eq.exa.fps.subs.fibonacci.subs-plausible}) that we have just used
follows from part \textbf{(c)} of the following proposition:\footnote{We are
treating the symbol \textquotedblleft$\circ$\textquotedblright\ similarly to
the multiplication sign $\cdot$ in our PEMDAS convention. Thus, an expression
like \textquotedblleft$f_{1}\circ g+f_{2}\circ g$\textquotedblright\ is
understood to mean $\left(  f_{1}\circ g\right)  +\left(  f_{2}\circ g\right)
$.}

\begin{proposition}
\label{prop.fps.subs.rules}Composition of FPSs satisfies the rules you would
expect it to satisfy: \medskip

\textbf{(a)} If $f_{1},f_{2},g\in K\left[  \left[  x\right]  \right]  $
satisfy $\left[  x^{0}\right]  g=0$, then $\left(  f_{1}+f_{2}\right)  \circ
g=f_{1}\circ g+f_{2}\circ g$. \medskip

\textbf{(b)} If $f_{1},f_{2},g\in K\left[  \left[  x\right]  \right]  $
satisfy $\left[  x^{0}\right]  g=0$, then $\left(  f_{1}\cdot f_{2}\right)
\circ g=\left(  f_{1}\circ g\right)  \cdot\left(  f_{2}\circ g\right)  $.
\medskip

\textbf{(c)} If $f_{1},f_{2},g\in K\left[  \left[  x\right]  \right]  $
satisfy $\left[  x^{0}\right]  g=0$, then $\dfrac{f_{1}}{f_{2}}\circ
g=\dfrac{f_{1}\circ g}{f_{2}\circ g}$, as long as $f_{2}$ is invertible. (In
particular, $f_{2}\circ g$ is automatically invertible under these
assumptions.) \medskip

\textbf{(d)} If $f,g\in K\left[  \left[  x\right]  \right]  $ satisfy $\left[
x^{0}\right]  g=0$, then $f^{k}\circ g=\left(  f\circ g\right)  ^{k}$ for each
$k\in\mathbb{N}$. \medskip

\textbf{(e)} If $f,g,h\in K\left[  \left[  x\right]  \right]  $ satisfy
$\left[  x^{0}\right]  g=0$ and $\left[  x^{0}\right]  h=0$, then $\left[
x^{0}\right]  \left(  g\circ h\right)  =0$ and $\left(  f\circ g\right)  \circ
h=f\circ\left(  g\circ h\right)  $. \medskip

\textbf{(f)} We have $\underline{a}\circ g=\underline{a}$ for each $a\in K$
and $g\in K\left[  \left[  x\right]  \right]  $. \medskip

\textbf{(g)} We have $x\circ g=g\circ x=g$ for each $g\in K\left[  \left[
x\right]  \right]  $. \medskip

\textbf{(h)} If $\left(  f_{i}\right)  _{i\in I}\in K\left[  \left[  x\right]
\right]  ^{I}$ is a summable family of FPSs, and if $g\in K\left[  \left[
x\right]  \right]  $ is an FPS satisfying $\left[  x^{0}\right]  g=0$, then
the family $\left(  f_{i}\circ g\right)  _{i\in I}\in K\left[  \left[
x\right]  \right]  ^{I}$ is summable as well and we have $\left(  \sum_{i\in
I}f_{i}\right)  \circ g=\sum_{i\in I}f_{i}\circ g$.
\end{proposition}

For our proof of Proposition \ref{prop.fps.subs.rules}, we will need the
following lemma:

\begin{lemma}
\label{lem.fps.fg-coeffs-0}Let $f,g\in K\left[  \left[  x\right]  \right]  $
satisfy $\left[  x^{0}\right]  g=0$. Let $k\in\mathbb{N}$ be such that the
first $k$ coefficients of $f$ are $0$. Then, the first $k$ coefficients of
$f\circ g$ are $0$.
\end{lemma}

\begin{proof}
[Proof of Lemma \ref{lem.fps.fg-coeffs-0}.]This is very similar to the proof
of Proposition \ref{prop.fps.subs.wd} \textbf{(a)}.

We have $\left[  x^{0}\right]  g=0$. Hence, Lemma \ref{lem.fps.g=xh} (applied
to $a=g$) yields that there exists an $h\in K\left[  \left[  x\right]
\right]  $ such that $g=xh$. Consider this $h$.

Write the FPS $f$ in the form $f=\left(  f_{0},f_{1},f_{2},\ldots\right)  $.
Then, the first $k$ coefficients of $f$ are $f_{0},f_{1},\ldots,f_{k-1}$.
Hence, these coefficients $f_{0},f_{1},\ldots,f_{k-1}$ are $0$ (since the
first $k$ coefficients of $f$ are $0$). In other words,%
\begin{equation}
f_{n}=0\ \ \ \ \ \ \ \ \ \ \text{for each }n<k.
\label{pf.lem.fps.fg-coeffs-0.gn=0}%
\end{equation}

Now, $f=\left(  f_{0},f_{1},f_{2},\ldots\right)  =\sum_{n\in\mathbb{N}}%
f_{n}x^{n}$ with $f_{0},f_{1},f_{2},\ldots\in K$. Hence, Definition
\ref{def.fps.subs} yields%
\begin{align*}
f\left[  g\right]   &  =\sum_{n\in\mathbb{N}}f_{n}g^{n}=\sum_{\substack{n\in
\mathbb{N};\\n<k}}\underbrace{f_{n}}_{\substack{=0\\\text{(by
(\ref{pf.lem.fps.fg-coeffs-0.gn=0}))}}}g^{n}+\underbrace{\sum_{\substack{n\in
\mathbb{N};\\n\geq k}}}_{=\sum_{\substack{n\in\mathbb{N};\\k\leq n}}}%
f_{n}\underbrace{g^{n}}_{\substack{=\left(  xh\right)  ^{n}\\\text{(since
}g=xh\text{)}}}=\underbrace{\sum_{\substack{n\in\mathbb{N};\\n<k}}0g^{n}}%
_{=0}+\sum_{\substack{n\in\mathbb{N};\\k\leq n}}f_{n}\left(  xh\right)  ^{n}\\
&  =\sum_{\substack{n\in\mathbb{N};\\k\leq n}}f_{n}\underbrace{\left(
xh\right)  ^{n}}_{=x^{n}h^{n}}=\sum_{\substack{n\in\mathbb{N};\\k\leq n}%
}f_{n}x^{n}h^{n}.
\end{align*}
Thus,
\[
f\circ g=f\left[  g\right]  =\sum_{\substack{n\in\mathbb{N};\\k\leq
n}}\underbrace{f_{n}x^{n}}_{=x^{n}f_{n}}h^{n}=\sum_{\substack{n\in
\mathbb{N};\\k\leq n}}x^{n}f_{n}h^{n}.
\]
But this ensures that the first $k$ coefficients of $f\circ g$ are
$0$\ \ \ \ \footnote{\textit{Proof.} We must show that $\left[  x^{m}\right]
\left(  f\circ g\right)  =0$ for any nonnegative integer $m<k$. But we can do
this directly: If $m$ is a nonnegative integer such that $m<k$, then%
\begin{align*}
\left[  x^{m}\right]  \left(  f\circ g\right)   &  =\left[  x^{m}\right]
\left(  \sum_{\substack{n\in\mathbb{N};\\k\leq n}}x^{n}f_{n}h^{n}\right)
\ \ \ \ \ \ \ \ \ \ \left(  \text{since }f\circ g=\sum_{\substack{n\in
\mathbb{N};\\k\leq n}}x^{n}f_{n}h^{n}\right) \\
&  =\sum_{\substack{n\in\mathbb{N};\\k\leq n}}\underbrace{\left[
x^{m}\right]  \left(  x^{n}f_{n}h^{n}\right)  }_{\substack{=\sum_{i=0}%
^{m}\left[  x^{i}\right]  \left(  x^{n}\right)  \cdot\left[  x^{m-i}\right]
\left(  f_{n}h^{n}\right)  \\\text{(by (\ref{pf.thm.fps.ring.xn(ab)=2}%
),}\\\text{applied to }x^{n}\text{, }f_{n}h^{n}\text{ and }m\\\text{instead of
}\mathbf{a}\text{, }\mathbf{b}\text{ and }n\text{)}}}=\sum_{\substack{n\in
\mathbb{N};\\k\leq n}}\ \ \sum_{i=0}^{m}\underbrace{\left[  x^{i}\right]
\left(  x^{n}\right)  }_{\substack{=0\\\text{(since }i\leq m<k\leq
n\\\text{and thus }i\neq n\text{)}}}\cdot\left[  x^{m-i}\right]  \left(
f_{n}h^{n}\right) \\
&  =\sum_{\substack{n\in\mathbb{N};\\k\leq n}}\ \ \sum_{i=0}^{m}0\cdot\left[
x^{m-i}\right]  \left(  f_{n}h^{n}\right)  =0,
\end{align*}
exactly as we wanted to show.}. Thus, Lemma \ref{lem.fps.fg-coeffs-0} follows.
\end{proof}

Our proof of Proposition \ref{prop.fps.subs.rules} will furthermore use the
\emph{Kronecker delta notation}:

\begin{definition}
\label{def.kron-delta}If $i$ and $j$ are any objects, then $\delta_{i,j}$
means the element $%
\begin{cases}
1, & \text{if }i=j;\\
0, & \text{if }i\neq j
\end{cases}
$\ \ \ of $K$.
\end{definition}

For example, $\delta_{2,2}=1$ and $\delta_{3,8}=0$.

\begin{proof}
[Proof of Proposition \ref{prop.fps.subs.rules}.]The proof is long and not
particularly combinatorial. I am merely writing it down because it is so
rarely explained in the literature. \medskip

\textbf{(a)} This is an easy consequence of the definitions, and also appears
in \cite[Theorem 7.62]{Loehr-BC} and \cite[Proposition 2.2.2]{Brewer14}.

Here are the details: Let $f_{1},f_{2},g\in K\left[  \left[  x\right]
\right]  $ satisfy $\left[  x^{0}\right]  g=0$. Write the FPSs $f_{1}$ and
$f_{2}$ as
\begin{equation}
f_{1}=\sum_{n\in\mathbb{N}}f_{1,n}x^{n}\ \ \ \ \ \ \ \ \ \ \text{and}%
\ \ \ \ \ \ \ \ \ \ f_{2}=\sum_{n\in\mathbb{N}}f_{2,n}x^{n}
\label{pf.prop.fps.subs.rules.a.1}%
\end{equation}
with $f_{1,0},f_{1,1},f_{1,2},\ldots\in K$ and $f_{2,0},f_{2,1},f_{2,2}%
,\ldots\in K$. Then, adding the two equalities in
(\ref{pf.prop.fps.subs.rules.a.1}) together, we find%
\[
f_{1}+f_{2}=\sum_{n\in\mathbb{N}}f_{1,n}x^{n}+\sum_{n\in\mathbb{N}}%
f_{2,n}x^{n}=\sum_{n\in\mathbb{N}}\underbrace{\left(  f_{1,n}x^{n}%
+f_{2,n}x^{n}\right)  }_{=\left(  f_{1,n}+f_{2,n}\right)  x^{n}}=\sum
_{n\in\mathbb{N}}\left(  f_{1,n}+f_{2,n}\right)  x^{n}.
\]
Thus, Definition \ref{def.fps.subs} (applied to $f=f_{1}+f_{2}$) yields%
\begin{equation}
\left(  f_{1}+f_{2}\right)  \left[  g\right]  =\sum_{n\in\mathbb{N}}\left(
f_{1,n}+f_{2,n}\right)  g^{n} \label{pf.prop.fps.subs.rules.a.3}%
\end{equation}
(since $f_{1,n}+f_{2,n}\in K$ for each $n\in\mathbb{N}$).

On the other hand, we have $f_{1}=\sum_{n\in\mathbb{N}}f_{1,n}x^{n}$. Thus,
Definition \ref{def.fps.subs} yields $f_{1}\left[  g\right]  =\sum
_{n\in\mathbb{N}}f_{1,n}g^{n}$ (since $f_{1,n}\in K$ for each $n\in\mathbb{N}%
$). Similarly, $f_{2}\left[  g\right]  =\sum_{n\in\mathbb{N}}f_{2,n}g^{n}$.
Adding these two equalities together, we obtain%
\[
f_{1}\left[  g\right]  +f_{2}\left[  g\right]  =\sum_{n\in\mathbb{N}}%
f_{1,n}g^{n}+\sum_{n\in\mathbb{N}}f_{2,n}g^{n}=\sum_{n\in\mathbb{N}%
}\underbrace{\left(  f_{1,n}g^{n}+f_{2,n}g^{n}\right)  }_{=\left(
f_{1,n}+f_{2,n}\right)  g^{n}}=\sum_{n\in\mathbb{N}}\left(  f_{1,n}%
+f_{2,n}\right)  g^{n}.
\]
Comparing this with (\ref{pf.prop.fps.subs.rules.a.3}), we obtain $\left(
f_{1}+f_{2}\right)  \left[  g\right]  =f_{1}\left[  g\right]  +f_{2}\left[
g\right]  $. In other words, $\left(  f_{1}+f_{2}\right)  \circ g=f_{1}\circ
g+f_{2}\circ g$ (since the notation $f\circ g$ is synonymous to $f\left[
g\right]  $). This proves Proposition \ref{prop.fps.subs.rules} \textbf{(a)}.

\bigskip

\textbf{(f)} This is a near-trivial consequence of the definitions. To wit:
Let $a\in K$. Then,
\begin{align}
\sum_{n\in\mathbb{N}}a\delta_{n,0}x^{n}  &  =a\underbrace{\delta_{0,0}%
}_{\substack{=1\\\text{(since }0=0\text{)}}}\underbrace{x^{0}}_{=\underline{1}%
}+\sum_{\substack{n\in\mathbb{N};\\n\neq0}}a\underbrace{\delta_{n,0}%
}_{\substack{=0\\\text{(since }n\neq0\text{)}}}x^{n}=a\cdot\underline{1}%
+\underbrace{\sum_{\substack{n\in\mathbb{N};\\n\neq0}}a\cdot0x^{n}}%
_{=0}\nonumber\\
&  =a\cdot\underline{1}=a\cdot\left(  1,0,0,0,\ldots\right)  =\left(
a\cdot1,a\cdot0,a\cdot0,a\cdot0,\ldots\right) \nonumber\\
&  =\left(  a,0,0,0,\ldots\right)  =\underline{a}.
\label{pf.prop.fps.subs.rules.f.1}%
\end{align}

Let $g\in K\left[  \left[  x\right]  \right]  $. We must prove that
$\underline{a}\circ g=\underline{a}$. From (\ref{pf.prop.fps.subs.rules.f.1}),
we obtain $\underline{a}=\sum_{n\in\mathbb{N}}a\delta_{n,0}x^{n}$. Hence,
Definition \ref{def.fps.subs} (or Definition \ref{def.pol.subs}) yields
$\underline{a}\left[  g\right]  =\sum_{n\in\mathbb{N}}a\delta_{n,0}g^{n}$
(since $a\delta_{n,0}\in K$ for each $n\in\mathbb{N}$). Thus,%
\begin{align*}
\underline{a}\left[  g\right]   &  =\sum_{n\in\mathbb{N}}a\delta_{n,0}%
g^{n}=a\underbrace{\delta_{0,0}}_{\substack{=1\\\text{(since }0=0\text{)}%
}}\underbrace{g^{0}}_{=\underline{1}}+\sum_{\substack{n\in\mathbb{N};\\n\neq
0}}a\underbrace{\delta_{n,0}}_{\substack{=0\\\text{(since }n\neq0\text{)}%
}}g^{n}=a\cdot\underline{1}+\underbrace{\sum_{\substack{n\in\mathbb{N}%
;\\n\neq0}}a\cdot0g^{n}}_{=0}\\
&  =a\cdot\underline{1}=\underline{a}.
\end{align*}
In other words, $\underline{a}\circ g=\underline{a}$ (since $\underline{a}%
\circ g$ is a synonym for $\underline{a}\left[  g\right]  $). This proves
Proposition \ref{prop.fps.subs.rules} \textbf{(f)}.

\bigskip

\textbf{(g)} This is easy, too. Indeed, we have%
\begin{align}
\sum_{n\in\mathbb{N}}\delta_{n,1}x^{n}  &  =\underbrace{\delta_{1,1}%
}_{\substack{=1\\\text{(since }1=1\text{)}}}\underbrace{x^{1}}_{=x}%
+\sum_{\substack{n\in\mathbb{N};\\n\neq1}}\underbrace{\delta_{n,1}%
}_{\substack{=0\\\text{(since }n\neq1\text{)}}}x^{n}=x+\underbrace{\sum
_{\substack{n\in\mathbb{N};\\n\neq1}}0x^{n}}_{=0}\nonumber\\
&  =x. \label{pf.prop.fps.subs.rules.g.1}%
\end{align}

Let $g\in K\left[  \left[  x\right]  \right]  $. We must prove that $x\circ
g=g\circ x=g$. From (\ref{pf.prop.fps.subs.rules.g.1}), we obtain
$x=\sum_{n\in\mathbb{N}}\delta_{n,1}x^{n}$. Hence, Definition
\ref{def.fps.subs} (or Definition \ref{def.pol.subs}) yields $x\left[
g\right]  =\sum_{n\in\mathbb{N}}\delta_{n,1}g^{n}$ (since $\delta_{n,1}\in K$
for each $n\in\mathbb{N}$). Thus,%
\[
x\left[  g\right]  =\sum_{n\in\mathbb{N}}\delta_{n,1}g^{n}=\underbrace{\delta
_{1,1}}_{\substack{=1\\\text{(since }1=1\text{)}}}\underbrace{g^{1}}_{=g}%
+\sum_{\substack{n\in\mathbb{N};\\n\neq1}}\underbrace{\delta_{n,1}%
}_{\substack{=0\\\text{(since }n\neq1\text{)}}}g^{n}=g+\underbrace{\sum
_{\substack{n\in\mathbb{N};\\n\neq1}}0g^{n}}_{=0}=g.
\]
In other words, $x\circ g=g$ (since $x\circ g$ is a synonym for $x\left[
g\right]  $).

Next, let us write $g$ in the form $g=\sum_{n\in\mathbb{N}}g_{n}x^{n}$ for
some $g_{0},g_{1},g_{2},\ldots\in K$. Then, Definition \ref{def.fps.subs}
yields $g\left[  x\right]  =\sum_{n\in\mathbb{N}}g_{n}x^{n}=g$. Thus, $g\circ
x=g\left[  x\right]  =g$. Combining this with $x\circ g=g$, we obtain $x\circ
g=g\circ x=g$. This proves Proposition \ref{prop.fps.subs.rules} \textbf{(g)}.

\bigskip

\textbf{(b)} This appears in \cite[Theorem 7.62]{Loehr-BC} and
\cite[Proposition 2.2.2]{Brewer14}. Here is the proof:

Let $f_{1},f_{2},g\in K\left[  \left[  x\right]  \right]  $ satisfy $\left[
x^{0}\right]  g=0$. Write the FPSs $f_{1}$ and $f_{2}$ as
\begin{equation}
f_{1}=\sum_{n\in\mathbb{N}}f_{1,n}x^{n}\ \ \ \ \ \ \ \ \ \ \text{and}%
\ \ \ \ \ \ \ \ \ \ f_{2}=\sum_{n\in\mathbb{N}}f_{2,n}x^{n}
\label{pf.prop.fps.subs.rules.b.1}%
\end{equation}
with $f_{1,0},f_{1,1},f_{1,2},\ldots\in K$ and $f_{2,0},f_{2,1},f_{2,2}%
,\ldots\in K$. Thus, Definition \ref{def.fps.subs} yields%
\[
f_{1}\left[  x\right]  =\sum_{n\in\mathbb{N}}f_{1,n}x^{n}=f_{1}%
\ \ \ \ \ \ \ \ \ \ \text{and}\ \ \ \ \ \ \ \ \ \ f_{2}\left[  x\right]
=\sum_{n\in\mathbb{N}}f_{2,n}x^{n}=f_{2}%
\]
and
\[
f_{1}\left[  g\right]  =\sum_{n\in\mathbb{N}}f_{1,n}g^{n}%
\ \ \ \ \ \ \ \ \ \ \text{and}\ \ \ \ \ \ \ \ \ \ f_{2}\left[  g\right]
=\sum_{n\in\mathbb{N}}f_{2,n}g^{n}.
\]
Hence,%
\begin{align*}
f_{1}\left[  g\right]   &  =\sum_{n\in\mathbb{N}}f_{1,n}g^{n}=\sum
_{i\in\mathbb{N}}f_{1,i}g^{i}\ \ \ \ \ \ \ \ \ \ \text{and}\\
f_{2}\left[  g\right]   &  =\sum_{n\in\mathbb{N}}f_{2,n}g^{n}=\sum
_{j\in\mathbb{N}}f_{2,j}g^{j}.
\end{align*}
Multiplying these two equalities together, we find%
\begin{align}
f_{1}\left[  g\right]  \cdot f_{2}\left[  g\right]   &  =\left(  \sum
_{i\in\mathbb{N}}f_{1,i}g^{i}\right)  \left(  \sum_{j\in\mathbb{N}}%
f_{2,j}g^{j}\right)  =\sum_{i\in\mathbb{N}}\ \ \sum_{j\in\mathbb{N}%
}\underbrace{f_{1,i}g^{i}f_{2,j}g^{j}}_{\substack{=f_{1,i}f_{2,j}g^{i}%
g^{j}\\=f_{1,i}f_{2,j}g^{i+j}}}\nonumber\\
&  =\underbrace{\sum_{i\in\mathbb{N}}\ \ \sum_{j\in\mathbb{N}}}_{=\sum
_{\left(  i,j\right)  \in\mathbb{N}^{2}}}f_{1,i}f_{2,j}g^{i+j}%
=\underbrace{\sum_{\left(  i,j\right)  \in\mathbb{N}^{2}}}_{=\sum
_{n\in\mathbb{N}}\ \ \sum_{\substack{\left(  i,j\right)  \in\mathbb{N}%
^{2};\\i+j=n}}}f_{1,i}f_{2,j}g^{i+j}\nonumber\\
&  =\sum_{n\in\mathbb{N}}\ \ \sum_{\substack{\left(  i,j\right)  \in
\mathbb{N}^{2};\\i+j=n}}f_{1,i}f_{2,j}\underbrace{g^{i+j}}_{\substack{=g^{n}%
\\\text{(since }i+j=n\text{)}}}=\sum_{n\in\mathbb{N}}\ \ \sum
_{\substack{\left(  i,j\right)  \in\mathbb{N}^{2};\\i+j=n}}f_{1,i}f_{2,j}%
g^{n}\nonumber\\
&  =\sum_{n\in\mathbb{N}}\left(  \sum_{\substack{\left(  i,j\right)
\in\mathbb{N}^{2};\\i+j=n}}f_{1,i}f_{2,j}\right)  g^{n}.
\label{pf.prop.fps.subs.rules.b.3}%
\end{align}

However, we can apply the same computations to $x$ instead of $g$ (since $x$
is also an FPS with $\left[  x^{0}\right]  x=0$). Thus, we obtain%
\[
f_{1}\left[  x\right]  \cdot f_{2}\left[  x\right]  =\sum_{n\in\mathbb{N}%
}\left(  \sum_{\substack{\left(  i,j\right)  \in\mathbb{N}^{2};\\i+j=n}%
}f_{1,i}f_{2,j}\right)  x^{n}.
\]
In view of $f_{1}\left[  x\right]  =f_{1}$ and $f_{2}\left[  x\right]  =f_{2}%
$, this rewrites as%
\[
f_{1}\cdot f_{2}=\sum_{n\in\mathbb{N}}\left(  \sum_{\substack{\left(
i,j\right)  \in\mathbb{N}^{2};\\i+j=n}}f_{1,i}f_{2,j}\right)  x^{n}.
\]
Hence, Definition \ref{def.fps.subs} yields%
\[
\left(  f_{1}\cdot f_{2}\right)  \left[  g\right]  =\sum_{n\in\mathbb{N}%
}\left(  \sum_{\substack{\left(  i,j\right)  \in\mathbb{N}^{2};\\i+j=n}%
}f_{1,i}f_{2,j}\right)  g^{n}%
\]
(since $\sum_{\substack{\left(  i,j\right)  \in\mathbb{N}^{2};\\i+j=n}%
}f_{1,i}f_{2,j}\in K$ for each $n\in\mathbb{N}$). Comparing this with
(\ref{pf.prop.fps.subs.rules.b.3}), we obtain%
\[
\left(  f_{1}\cdot f_{2}\right)  \left[  g\right]  =f_{1}\left[  g\right]
\cdot f_{2}\left[  g\right]  .
\]
In other words, $\left(  f_{1}\cdot f_{2}\right)  \circ g=\left(  f_{1}\circ
g\right)  \cdot\left(  f_{2}\circ g\right)  $ (since the notation $f\circ g$
is synonymous to $f\left[  g\right]  $). This proves Proposition
\ref{prop.fps.subs.rules} \textbf{(b)}. \smallskip

\begin{fineprint}
Did you notice it? I have cheated. The above proof of Proposition
\ref{prop.fps.subs.rules} \textbf{(b)} relied on some manipulations of
infinite sums that need to be justified. Namely, we replaced \textquotedblleft%
$\sum_{i\in\mathbb{N}}\ \ \sum_{j\in\mathbb{N}}$\textquotedblright\ by
\textquotedblleft$\sum_{\left(  i,j\right)  \in\mathbb{N}^{2}}$%
\textquotedblright. This is an application of the \textquotedblleft discrete
Fubini rule\textquotedblright, and as we said above, this rule can only be
used if we know that the family $\left(  f_{1,i}f_{2,j}g^{i+j}\right)
_{\left(  i,j\right)  \in\mathbb{N}\times\mathbb{N}}$ is summable. In other
words, we need to show the following statement:

\begin{statement}
\textit{Statement 1:} For each $m\in\mathbb{N}$, all but finitely many pairs
$\left(  i,j\right)  \in\mathbb{N}\times\mathbb{N}$ satisfy $\left[
x^{m}\right]  \left(  f_{1,i}f_{2,j}g^{i+j}\right)  =0$.
\end{statement}

We shall achieve this by proving the following statement:

\begin{statement}
\textit{Statement 2:} For any three nonnegative integers $m,i,j$ with $m<i+j$,
we have $\left[  x^{m}\right]  \left(  g^{i+j}\right)  =0$.
\end{statement}

[\textit{Proof of Statement 2:} Let $m,i,j$ be three nonnegative integers with
$m<i+j$. We must show that $\left[  x^{m}\right]  \left(  g^{i+j}\right)  =0$.

We have $\left[  x^{0}\right]  g=0$. Hence, Lemma \ref{lem.fps.g=xh} (applied
to $a=g$) shows that there exists an $h\in K\left[  \left[  x\right]  \right]
$ such that $g=xh$. Consider this $h$.

Now, from $g=xh$, we obtain $g^{i+j}=\left(  xh\right)  ^{i+j}=x^{i+j}h^{i+j}%
$. However, Lemma \ref{lem.fps.first-n-coeffs-of-xna} (applied to $i+j$ and
$h^{i+j}$ instead of $k$ and $a$) shows that the first $i+j$ coefficients of
the FPS $x^{i+j}h^{i+j}$ are $0$. In other words, the first $i+j$ coefficients
of the FPS $g^{i+j}$ are $0$ (since $g^{i+j}=x^{i+j}h^{i+j}$). But $\left[
x^{m}\right]  \left(  g^{i+j}\right)  $ is one of these first $i+j$
coefficients (since $m<i+j$). Thus, we conclude that $\left[  x^{m}\right]
\left(  g^{i+j}\right)  =0$. This proves Statement 2.]

[\textit{Proof of Statement 1:} Let $m\in\mathbb{N}$. If $\left(  i,j\right)
\in\mathbb{N}\times\mathbb{N}$ is a pair satisfying $m<i+j$, then%
\begin{align*}
\left[  x^{m}\right]  \left(  f_{1,i}f_{2,j}g^{i+j}\right)   &  =f_{1,i}%
f_{2,j}\underbrace{\left[  x^{m}\right]  \left(  g^{i+j}\right)
}_{\substack{=0\\\text{(by Statement 2)}}}\ \ \ \ \ \ \ \ \ \ \left(  \text{by
(\ref{pf.thm.fps.ring.xn(la)=})}\right) \\
&  =0.
\end{align*}
Thus, all but finitely many pairs $\left(  i,j\right)  \in\mathbb{N}%
\times\mathbb{N}$ satisfy $\left[  x^{m}\right]  \left(  f_{1,i}f_{2,j}%
g^{i+j}\right)  =0$ (because all but finitely many such pairs satisfy
$m<i+j$). This proves Statement 1.]

As explained above, Statement 1 shows that the family \newline$\left(
f_{1,i}f_{2,j}g^{i+j}\right)  _{\left(  i,j\right)  \in\mathbb{N}%
\times\mathbb{N}}$ is summable, and thus our interchange of summation signs
made above is justified. This completes our proof of Proposition
\ref{prop.fps.subs.rules} \textbf{(b)}.
\end{fineprint}

\bigskip

\textbf{(c)} This follows easily from parts \textbf{(b)} and \textbf{(f)}. In
detail: Let $f_{1},f_{2},g\in K\left[  \left[  x\right]  \right]  $ be such
that $\left[  x^{0}\right]  g=0$. Assume that $f_{2}$ is invertible. Let us
first show that $f_{2}\circ g$ is invertible.

Indeed, consider the inverse $f_{2}^{-1}$ of $f_{2}$. This inverse exists
(since $f_{2}$ is invertible) and satisfies $f_{2}^{-1}\cdot f_{2}%
=\underline{1}$. Now, Proposition \ref{prop.fps.subs.rules} \textbf{(b)}
(applied to $f_{2}^{-1}$ instead of $f_{1}$) yields $\left(  f_{2}^{-1}\cdot
f_{2}\right)  \circ g=\left(  f_{2}^{-1}\circ g\right)  \cdot\left(
f_{2}\circ g\right)  $. Hence,%
\[
\left(  f_{2}^{-1}\circ g\right)  \cdot\left(  f_{2}\circ g\right)
=\underbrace{\left(  f_{2}^{-1}\cdot f_{2}\right)  }_{=\underline{1}}%
\circ\,g=\underline{1}\circ g=\underline{1}%
\]
(by Proposition \ref{prop.fps.subs.rules} \textbf{(f)}, applied to $a=1$).
Thus, the FPS $f_{2}^{-1}\circ g$ is an inverse of $f_{2}\circ g$. Hence,
$f_{2}\circ g$ is invertible. The expression $\dfrac{f_{1}\circ g}{f_{2}\circ
g}$ is therefore well-defined.

It now remains to prove that $\dfrac{f_{1}}{f_{2}}\circ g=\dfrac{f_{1}\circ
g}{f_{2}\circ g}$. To this purpose, we argue as follows: The expression
$\dfrac{f_{1}}{f_{2}}$ is well-defined, since $f_{2}$ is invertible.
Proposition \ref{prop.fps.subs.rules} \textbf{(b)} (applied to $\dfrac{f_{1}%
}{f_{2}}$ instead of $f_{1}$) yields $\left(  \dfrac{f_{1}}{f_{2}}\cdot
f_{2}\right)  \circ g=\left(  \dfrac{f_{1}}{f_{2}}\circ g\right)  \cdot\left(
f_{2}\circ g\right)  $. In view of $\dfrac{f_{1}}{f_{2}}\cdot f_{2}=f_{1}$,
this rewrites as $f_{1}\circ g=\left(  \dfrac{f_{1}}{f_{2}}\circ g\right)
\cdot\left(  f_{2}\circ g\right)  $. We can divide both sides of this equality
by $f_{2}\circ g$ (since $f_{2}\circ g$ is invertible), and thus obtain
$\dfrac{f_{1}\circ g}{f_{2}\circ g}=\dfrac{f_{1}}{f_{2}}\circ g$. In other
words, $\dfrac{f_{1}}{f_{2}}\circ g=\dfrac{f_{1}\circ g}{f_{2}\circ g}$. Thus,
Proposition \ref{prop.fps.subs.rules} \textbf{(c)} is proven.

\bigskip

\textbf{(d)} Let $f,g\in K\left[  \left[  x\right]  \right]  $ satisfy
$\left[  x^{0}\right]  g=0$. We must prove that $f^{k}\circ g=\left(  f\circ
g\right)  ^{k}$ for each $k\in\mathbb{N}$.

We prove this by induction on $k$:

\textit{Induction base:} We have $\underbrace{f^{0}}_{=\underline{1}}%
\circ\,g=\underline{1}\circ g=\underline{1}$ (by Proposition
\ref{prop.fps.subs.rules} \textbf{(f)}, applied to $a=1$). Comparing this with
$\left(  f\circ g\right)  ^{0}=\underline{1}$, we find $f^{0}\circ g=\left(
f\circ g\right)  ^{0}$. In other words, $f^{k}\circ g=\left(  f\circ g\right)
^{k}$ holds for $k=0$.

\textit{Induction step:} Let $m\in\mathbb{N}$. Assume that $f^{k}\circ
g=\left(  f\circ g\right)  ^{k}$ holds for $k=m$. We must prove that
$f^{k}\circ g=\left(  f\circ g\right)  ^{k}$ holds for $k=m+1$.

We have assumed that $f^{k}\circ g=\left(  f\circ g\right)  ^{k}$ holds for
$k=m$. In other words, we have $f^{m}\circ g=\left(  f\circ g\right)  ^{m}$.
Now,%
\begin{align*}
\underbrace{f^{m+1}}_{=f\cdot f^{m}}\circ\,g  &  =\left(  f\cdot f^{m}\right)
\circ g=\left(  f\circ g\right)  \cdot\underbrace{\left(  f^{m}\circ g\right)
}_{=\left(  f\circ g\right)  ^{m}}\\
&  \ \ \ \ \ \ \ \ \ \ \left(  \text{by Proposition \ref{prop.fps.subs.rules}
\textbf{(b)}, applied to }f_{1}=f\text{ and }f_{2}=f^{m}\right) \\
&  =\left(  f\circ g\right)  \cdot\left(  f\circ g\right)  ^{m}=\left(  f\circ
g\right)  ^{m+1}.
\end{align*}
In other words, $f^{k}\circ g=\left(  f\circ g\right)  ^{k}$ holds for
$k=m+1$. This completes the induction step. Thus, we have proven that
$f^{k}\circ g=\left(  f\circ g\right)  ^{k}$ for each $k\in\mathbb{N}$.
Proposition \ref{prop.fps.subs.rules} \textbf{(d)} is now proven.

\bigskip

\textbf{(h)} This is just a generalization of Proposition
\ref{prop.fps.subs.rules} \textbf{(a)} to (potentially) infinite sums. The
proof follows the same method, but unfortunately requires some technical
reasoning about summability. I will give the full proof for the sake of
completeness, but be warned that it contains nothing of interest.

Let $\left(  f_{i}\right)  _{i\in I}\in K\left[  \left[  x\right]  \right]
^{I}$ be a summable family of FPSs. Let $g\in K\left[  \left[  x\right]
\right]  $ be an FPS satisfying $\left[  x^{0}\right]  g=0$.

First, we shall prove that the family $\left(  f_{i}\circ g\right)  _{i\in
I}\in K\left[  \left[  x\right]  \right]  ^{I}$ is summable. Indeed, we recall
that the family $\left(  f_{i}\right)  _{i\in I}\in K\left[  \left[  x\right]
\right]  ^{I}$ is summable. In other words,%
\[
\text{for each }n\in\mathbb{N}\text{, all but finitely many }i\in I\text{
satisfy }\left[  x^{n}\right]  f_{i}=0
\]
(by the definition of \textquotedblleft summable\textquotedblright). In other
words, for each $n\in\mathbb{N}$, there exists a finite subset $I_{n}$ of $I$
such that%
\begin{equation}
\text{all }i\in I\setminus I_{n}\text{ satisfy }\left[  x^{n}\right]  f_{i}=0.
\label{pf.prop.fps.subs.rules.h.1}%
\end{equation}
Consider this subset $I_{n}$. Thus, all the sets $I_{0},I_{1},I_{2},\ldots$
are finite subsets of $I$.

Now, let $n\in\mathbb{N}$ be arbitrary. The set $I_{0}\cup I_{1}\cup\cdots\cup
I_{n}$ is a union of $n+1$ finite subsets of $I$ (because all the sets
$I_{0},I_{1},I_{2},\ldots$ are finite subsets of $I$), and thus itself is a
finite subset of $I$. Moreover,
\begin{equation}
\text{all }i\in I\setminus\left(  I_{0}\cup I_{1}\cup\cdots\cup I_{n}\right)
\text{ satisfy }\left[  x^{n}\right]  \left(  f_{i}\circ g\right)  =0.
\label{pf.prop.fps.subs.rules.h.2}%
\end{equation}

\begin{fineprint}
[\textit{Proof of (\ref{pf.prop.fps.subs.rules.h.2}):} Let $i\in
I\setminus\left(  I_{0}\cup I_{1}\cup\cdots\cup I_{n}\right)  $. We must show
that $\left[  x^{n}\right]  \left(  f_{i}\circ g\right)  =0$.

Let $m\in\left\{  0,1,\ldots,n\right\}  $. Then, $I_{m}\subseteq I_{0}\cup
I_{1}\cup\cdots\cup I_{n}$, so that $I_{0}\cup I_{1}\cup\cdots\cup
I_{n}\supseteq I_{m}$ and thus $I\setminus\underbrace{\left(  I_{0}\cup
I_{1}\cup\cdots\cup I_{n}\right)  }_{\supseteq I_{m}}\subseteq I\setminus
I_{m}$. Hence, $i\in I\setminus\left(  I_{0}\cup I_{1}\cup\cdots\cup
I_{n}\right)  \subseteq I\setminus I_{m}$. Therefore,
(\ref{pf.prop.fps.subs.rules.h.1}) (applied to $m$ instead of $n$) yields
$\left[  x^{m}\right]  f_{i}=0$.

Forget that we fixed $m$. We thus have shown that $\left[  x^{m}\right]
f_{i}=0$ for each $m\in\left\{  0,1,\ldots,n\right\}  $. In other words, the
first $n+1$ coefficients of $f_{i}$ are $0$. Hence, Lemma
\ref{lem.fps.fg-coeffs-0} (applied to $f_{i}$ and $n+1$ instead of $f$ and
$k$) shows that the first $n+1$ coefficients of $f_{i}\circ g$ are $0$.
However, $\left[  x^{n}\right]  \left(  f_{i}\circ g\right)  $ is one of these
first $n+1$ coefficients (indeed, it is the last of them); thus, this
coefficient $\left[  x^{n}\right]  \left(  f_{i}\circ g\right)  $ must be $0$.
This proves (\ref{pf.prop.fps.subs.rules.h.2}).] \smallskip
\end{fineprint}

Now, recall that $I_{0}\cup I_{1}\cup\cdots\cup I_{n}$ is a finite subset of
$I$. Hence, thanks to (\ref{pf.prop.fps.subs.rules.h.2}), we know that there
exists a finite subset $J$ of $I$ such that all $i\in I\setminus J$ satisfy
$\left[  x^{n}\right]  \left(  f_{i}\circ g\right)  =0$ (namely, $J=I_{0}\cup
I_{1}\cup\cdots\cup I_{n}$). In other words, all but finitely many $i\in I$
satisfy $\left[  x^{n}\right]  \left(  f_{i}\circ g\right)  =0$.

Forget that we fixed $n$. We thus have shown that%
\[
\text{for each }n\in\mathbb{N}\text{, all but finitely many }i\in I\text{
satisfy }\left[  x^{n}\right]  \left(  f_{i}\circ g\right)  =0.
\]
In other words, the family $\left(  f_{i}\circ g\right)  _{i\in I}\in K\left[
\left[  x\right]  \right]  ^{I}$ is summable (by the definition of
\textquotedblleft summable\textquotedblright).

It now remains to prove that $\left(  \sum_{i\in I}f_{i}\right)  \circ
g=\sum_{i\in I}f_{i}\circ g$.

For each $i\in I$, we write the FPS $f_{i}$ in the form $f_{i}=\sum
_{n\in\mathbb{N}}f_{i,n}x^{n}$ with \newline$f_{i,0},f_{i,1},f_{i,2},\ldots\in
K$. First, we shall show that
\begin{equation}
\text{the family }\left(  f_{i,m}g^{m}\right)  _{\left(  i,m\right)  \in
I\times\mathbb{N}}\text{ is summable.} \label{pf.prop.fps.subs.rules.h.3}%
\end{equation}

\begin{fineprint}
[\textit{Proof of (\ref{pf.prop.fps.subs.rules.h.3}):} Fix an $n\in\mathbb{N}%
$. Let $J$ denote the set $I_{0}\cup I_{1}\cup\cdots\cup I_{n}$. Hence, $J$ is
a union of $n+1$ finite subsets of $I$ (because all the sets $I_{0}%
,I_{1},I_{2},\ldots$ are finite subsets of $I$), and thus itself is a finite
subset of $I$. The set $J\times\left\{  0,1,\ldots,n\right\}  $ must be finite
(since it is the product of the two finite sets $J$ and $\left\{
0,1,\ldots,n\right\}  $).

Now, let $\left(  i,m\right)  \in\left(  I\times\mathbb{N}\right)
\setminus\left(  J\times\left\{  0,1,\ldots,n\right\}  \right)  $. We shall
prove that $\left[  x^{n}\right]  \left(  f_{i,m}g^{m}\right)  =0$.

We have $\left(  i,m\right)  \notin J\times\left\{  0,1,\ldots,n\right\}  $
(since $\left(  i,m\right)  \in\left(  I\times\mathbb{N}\right)
\setminus\left(  J\times\left\{  0,1,\ldots,n\right\}  \right)  $).

We note that $f_{i,m}\in K$ and thus $\left[  x^{n}\right]  \left(
f_{i,m}g^{m}\right)  =f_{i,m}\cdot\left[  x^{n}\right]  \left(  g^{m}\right)
$ (by (\ref{pf.thm.fps.ring.xn(la)=})). However, we have $\left[
x^{0}\right]  g=0$. Thus, Proposition \ref{prop.fps.subs.wd} \textbf{(a)}
(applied to $f_{i}$ and $f_{i,j}$ and $m$ instead of $f$ and $f_{j}$ and $n$)
yields that the first $m$ coefficients of the FPS $g^{m}$ are $0$. In other
words, we have
\begin{equation}
\left[  x^{k}\right]  \left(  g^{m}\right)  =0\ \ \ \ \ \ \ \ \ \ \text{for
each }k\in\left\{  0,1,\ldots,m-1\right\}  .
\label{pf.prop.fps.subs.rules.h.3.pf.5}%
\end{equation}

Now, if we have $n\in\left\{  0,1,\ldots,m-1\right\}  $, then we have $\left[
x^{n}\right]  \left(  g^{m}\right)  =0$ (by
(\ref{pf.prop.fps.subs.rules.h.3.pf.5}), applied to $k=n$) and therefore
$\left[  x^{n}\right]  \left(  f_{i,m}g^{m}\right)  =f_{i,m}\cdot
\underbrace{\left[  x^{n}\right]  \left(  g^{m}\right)  }_{=0}=0$. Hence,
$\left[  x^{n}\right]  \left(  f_{i,m}g^{m}\right)  =0$ is proved in the case
when $n\in\left\{  0,1,\ldots,m-1\right\}  $. Thus, for the rest of this proof
of $\left[  x^{n}\right]  \left(  f_{i,m}g^{m}\right)  =0$, we WLOG assume
that $n\notin\left\{  0,1,\ldots,m-1\right\}  $. Hence, $n>m-1$ (since
$n\in\mathbb{N}$), so that $n\geq m$ (since $n$ and $m$ are integers).
Therefore, $m\leq n$, so that $m\in\left\{  0,1,\ldots,n\right\}  $ and
therefore $I_{m}\subseteq I_{0}\cup I_{1}\cup\cdots\cup I_{n}=J$ (since we
defined $J$ to be $I_{0}\cup I_{1}\cup\cdots\cup I_{n}$).

If we had $i\in I_{m}$, then we would have $\left(  i,m\right)  \in
J\times\left\{  0,1,\ldots,n\right\}  $ (since $i\in I_{m}\subseteq J$ and
$m\in\left\{  0,1,\ldots,n\right\}  $), which would contradict the fact that
$\left(  i,m\right)  \notin J\times\left\{  0,1,\ldots,n\right\}  $. Thus, we
cannot have $i\in I_{m}$. Hence, $i\notin I_{m}$, so that $i\in I\setminus
I_{m}$ (since $i\in I$). Thus, (\ref{pf.prop.fps.subs.rules.h.1}) (applied to
$m$ instead of $n$) yields $\left[  x^{m}\right]  f_{i}=0$. However, from
$f_{i}=\sum_{n\in\mathbb{N}}f_{i,n}x^{n}$, we see that $\left[  x^{m}\right]
f_{i}=f_{i,m}$. Thus, $f_{i,m}=\left[  x^{m}\right]  f_{i}=0$. Consequently,
$\left[  x^{n}\right]  \left(  f_{i,m}g^{m}\right)  =\underbrace{f_{i,m}}%
_{=0}\cdot\left[  x^{n}\right]  \left(  g^{m}\right)  =0$.

Forget that we fixed $\left(  i,m\right)  $. We thus have shown that
\[
\text{all }\left(  i,m\right)  \in\left(  I\times\mathbb{N}\right)
\setminus\left(  J\times\left\{  0,1,\ldots,n\right\}  \right)  \text{ satisfy
}\left[  x^{n}\right]  \left(  f_{i,m}g^{m}\right)  =0.
\]
Therefore, all but finitely many $\left(  i,m\right)  \in I\times\mathbb{N}$
satisfy $\left[  x^{n}\right]  \left(  f_{i,m}g^{m}\right)  =0$ (since
$J\times\left\{  0,1,\ldots,n\right\}  $ is a finite subset of $I\times
\mathbb{N}$).

Forget that we fixed $n$. We thus have shown that%
\[
\text{for each }n\in\mathbb{N}\text{, all but finitely many }\left(
i,m\right)  \in I\times\mathbb{N}\text{ satisfy }\left[  x^{n}\right]  \left(
f_{i,m}g^{m}\right)  =0.
\]
In other words, the family $\left(  f_{i,m}g^{m}\right)  _{\left(  i,m\right)
\in I\times\mathbb{N}}$ is summable. This proves
(\ref{pf.prop.fps.subs.rules.h.3}).]
\end{fineprint}

Now, we have shown that the family $\left(  f_{i,m}g^{m}\right)  _{\left(
i,m\right)  \in I\times\mathbb{N}}$ is summable. Renaming the index $\left(
i,m\right)  $ as $\left(  i,n\right)  $, we thus conclude that the family
$\left(  f_{i,n}g^{n}\right)  _{\left(  i,n\right)  \in I\times\mathbb{N}}$ is
summable. The same argument (but with $g$ replaced by $x$) shows that the
family $\left(  f_{i,n}x^{n}\right)  _{\left(  i,n\right)  \in I\times
\mathbb{N}}$ is summable (since the FPS $x$ satisfies $\left[  x^{0}\right]
x=0$).

The proof of $\left(  \sum_{i\in I}f_{i}\right)  \circ g=\sum_{i\in I}%
f_{i}\circ g$ is now just a matter of computation: Summing the equalities
$f_{i}=\sum_{n\in\mathbb{N}}f_{i,n}x^{n}$ over all $i\in I$, we obtain%
\[
\sum_{i\in I}f_{i}=\sum_{i\in I}\ \ \sum_{n\in\mathbb{N}}f_{i,n}x^{n}%
=\sum_{n\in\mathbb{N}}\ \ \sum_{i\in I}f_{i,n}x^{n}%
\]
(here, we have been able to interchange the summation signs, since the family
$\left(  f_{i,n}x^{n}\right)  _{\left(  i,n\right)  \in I\times\mathbb{N}}$ is
summable). Thus,%
\[
\sum_{i\in I}f_{i}=\sum_{n\in\mathbb{N}}\ \ \sum_{i\in I}f_{i,n}x^{n}%
=\sum_{n\in\mathbb{N}}\left(  \sum_{i\in I}f_{i,n}\right)  x^{n}.
\]
Hence, Definition \ref{def.fps.subs} (applied to $f=\sum_{i\in I}f_{i}$)
yields%
\begin{equation}
\left(  \sum_{i\in I}f_{i}\right)  \left[  g\right]  =\sum_{n\in\mathbb{N}%
}\left(  \sum_{i\in I}f_{i,n}\right)  g^{n} \label{pf.prop.fps.subs.rules.h.5}%
\end{equation}
(since $\sum_{i\in I}f_{i,n}\in K$ for each $n\in\mathbb{N}$).

On the other hand, for each $i\in I$, we have $f_{i}\left[  g\right]
=\sum_{n\in\mathbb{N}}f_{i,n}g^{n}$ (by Definition \ref{def.fps.subs}, since
$f_{i}=\sum_{n\in\mathbb{N}}f_{i,n}x^{n}$ with $f_{i,0},f_{i,1},f_{i,2}%
,\ldots\in K$). Summing these equalities over all $i\in I$, we find%
\[
\sum_{i\in I}f_{i}\left[  g\right]  =\sum_{i\in I}\ \ \sum_{n\in\mathbb{N}%
}f_{i,n}g^{n}=\sum_{n\in\mathbb{N}}\ \ \sum_{i\in I}f_{i,n}g^{n}%
\]
(again, we have been able to interchange the summation signs, since the family
$\left(  f_{i,n}g^{n}\right)  _{\left(  i,n\right)  \in I\times\mathbb{N}}$ is
summable). Thus,%
\[
\sum_{i\in I}f_{i}\left[  g\right]  =\sum_{n\in\mathbb{N}}\ \ \sum_{i\in
I}f_{i,n}g^{n}=\sum_{n\in\mathbb{N}}\left(  \sum_{i\in I}f_{i,n}\right)
g^{n}.
\]
Comparing this with (\ref{pf.prop.fps.subs.rules.h.5}), we find $\left(
\sum_{i\in I}f_{i}\right)  \left[  g\right]  =\sum_{i\in I}f_{i}\left[
g\right]  $. In other words, $\left(  \sum_{i\in I}f_{i}\right)  \circ
g=\sum_{i\in I}f_{i}\circ g$ (since the notation $f\circ g$ is synonymous to
$f\left[  g\right]  $). Thus, Proposition \ref{prop.fps.subs.rules}
\textbf{(h)} is proven.

\bigskip

\textbf{(e)} This is \cite[Theorem 7.63]{Loehr-BC} and \cite[Proposition
2.2.5]{Brewer14}.\footnote{See also \cite[Proposition 7.6.14]{19s} for a
similar property for polynomials.} Again, let us give the proof:

Write the FPS $g$ in the form $g=\sum_{n\in\mathbb{N}}g_{n}x^{n}$ for some
$g_{0},g_{1},g_{2},\ldots\in K$. Then, $g_{0}=\left[  x^{0}\right]  g=0$.
Moreover, $g\circ h=g\left[  h\right]  =\sum_{n\in\mathbb{N}}g_{n}h^{n}$ (by
Definition \ref{def.fps.subs}, because $g=\sum_{n\in\mathbb{N}}g_{n}x^{n}$
with $g_{0},g_{1},g_{2},\ldots\in K$). But Proposition \ref{prop.fps.subs.wd}
\textbf{(c)} (applied to $g$, $h$ and $g_{n}$ instead of $f$, $g$ and $f_{n}$)
yields $\left[  x^{0}\right]  \left(  \sum_{n\in\mathbb{N}}g_{n}h^{n}\right)
=g_{0}=0$. In view of $g\circ h=\sum_{n\in\mathbb{N}}g_{n}h^{n}$, this
rewrites as $\left[  x^{0}\right]  \left(  g\circ h\right)  =0$. Hence, the
composition $f\circ\left(  g\circ h\right)  $ is well-defined.

It remains to show that $\left(  f\circ g\right)  \circ h=f\circ\left(  g\circ
h\right)  $.

Write the FPS $f$ in the form $f=\sum_{n\in\mathbb{N}}f_{n}x^{n}$ for some
$f_{0},f_{1},f_{2},\ldots\in K$. Thus, Definition \ref{def.fps.subs} yields%
\[
f\left[  g\right]  =\sum_{n\in\mathbb{N}}f_{n}g^{n}%
\ \ \ \ \ \ \ \ \ \ \text{and}\ \ \ \ \ \ \ \ \ \ f\left[  g\circ h\right]
=\sum_{n\in\mathbb{N}}f_{n}\cdot\left(  g\circ h\right)  ^{n}.
\]
Moreover, the family $\left(  f_{n}g^{n}\right)  _{n\in\mathbb{N}}$ is
summable (by Proposition \ref{prop.fps.subs.wd} \textbf{(b)}). Hence,
Proposition \ref{prop.fps.subs.rules} \textbf{(h)} (applied to $\left(
f_{n}g^{n}\right)  _{n\in\mathbb{N}}$ and $h$ instead of $\left(
f_{i}\right)  _{i\in I}$ and $g$) yields that the family $\left(  \left(
f_{n}g^{n}\right)  \circ h\right)  _{n\in\mathbb{N}}\in K\left[  \left[
x\right]  \right]  ^{\mathbb{N}}$ is summable as well and that we have
\begin{equation}
\left(  \sum_{n\in\mathbb{N}}f_{n}g^{n}\right)  \circ h=\sum_{n\in\mathbb{N}%
}\left(  f_{n}g^{n}\right)  \circ h. \label{pf.prop.fps.subs.rules.e.3}%
\end{equation}
In view of $f\circ g=f\left[  g\right]  =\sum_{n\in\mathbb{N}}f_{n}g^{n}$, we
can rewrite (\ref{pf.prop.fps.subs.rules.e.3}) as%
\begin{equation}
\left(  f\circ g\right)  \circ h=\sum_{n\in\mathbb{N}}\left(  f_{n}%
g^{n}\right)  \circ h. \label{pf.prop.fps.subs.rules.e.4}%
\end{equation}
However, for each $n\in\mathbb{N}$, we have $f_{n}g^{n}=\underline{f_{n}}\cdot
g^{n}$ (by Theorem \ref{thm.fps.ring} \textbf{(d)}, applied to $\lambda=f_{n}$
and $\mathbf{a}=g^{n}$) and thus%
\begin{align}
\left(  f_{n}g^{n}\right)  \circ h  &  =\left(  \underline{f_{n}}\cdot
g^{n}\right)  \circ h=\underbrace{\left(  \underline{f_{n}}\circ h\right)
}_{\substack{=\underline{f_{n}}\\\text{(by Proposition
\ref{prop.fps.subs.rules} \textbf{(f)},}\\\text{applied to }f_{n}\text{ and
}h\\\text{instead of }a\text{ and }g\text{)}}}\cdot\underbrace{\left(
g^{n}\circ h\right)  }_{\substack{=\left(  g\circ h\right)  ^{n}\\\text{(by
Proposition \ref{prop.fps.subs.rules} \textbf{(d)},}\\\text{applied to
}g\text{ and }h\\\text{instead of }f\text{ and }g\text{)}}}\nonumber\\
&  \ \ \ \ \ \ \ \ \ \ \ \ \ \ \ \ \ \ \ \ \left(
\begin{array}
[c]{c}%
\text{by Proposition \ref{prop.fps.subs.rules} \textbf{(b)},}\\
\text{applied to }\underline{f_{n}}\text{, }g^{n}\text{ and }h\text{ instead
of }f_{1}\text{, }f_{2}\text{ and }g
\end{array}
\right) \nonumber\\
&  =\underline{f_{n}}\cdot\left(  g\circ h\right)  ^{n}=f_{n}\cdot\left(
g\circ h\right)  ^{n} \label{pf.prop.fps.subs.rules.e.5}%
\end{align}
(since Theorem \ref{thm.fps.ring} \textbf{(d)} (applied to $\lambda=f_{n}$ and
$\mathbf{a}=\left(  g\circ h\right)  ^{n}$) yields $f_{n}\cdot\left(  g\circ
h\right)  ^{n}=\underline{f_{n}}\cdot\left(  g\circ h\right)  ^{n}$). Thus,
(\ref{pf.prop.fps.subs.rules.e.4}) becomes%
\begin{align*}
\left(  f\circ g\right)  \circ h  &  =\sum_{n\in\mathbb{N}}\underbrace{\left(
f_{n}g^{n}\right)  \circ h}_{\substack{=f_{n}\cdot\left(  g\circ h\right)
^{n}\\\text{(by (\ref{pf.prop.fps.subs.rules.e.5}))}}}=\sum_{n\in\mathbb{N}%
}f_{n}\cdot\left(  g\circ h\right)  ^{n}\\
&  =f\left[  g\circ h\right]  \ \ \ \ \ \ \ \ \ \ \left(  \text{since
}f\left[  g\circ h\right]  =\sum_{n\in\mathbb{N}}f_{n}\cdot\left(  g\circ
h\right)  ^{n}\right) \\
&  =f\circ\left(  g\circ h\right)  .
\end{align*}
This completes the proof of Proposition \ref{prop.fps.subs.rules} \textbf{(e)}.
\end{proof}

\begin{example}
Let us use Proposition \ref{prop.fps.subs.rules} \textbf{(c)} to justify the
equality (\ref{eq.exa.fps.subs.fibonacci.subs-plausible}) that we used in
Example \ref{exa.fps.subs.fibonacci}. Indeed, we know that the FPS $1-x$ is
invertible. Thus, applying Proposition \ref{prop.fps.subs.rules} \textbf{(c)}
to $f_{1}=1$ and $f_{2}=1-x$ and $g=x+x^{2}$, we obtain%
\[
\dfrac{1}{1-x}\circ\left(  x+x^{2}\right)  =\dfrac{1\circ\left(
x+x^{2}\right)  }{\left(  1-x\right)  \circ\left(  x+x^{2}\right)  }.
\]
Using the notation $f\left[  g\right]  $ instead of $f\circ g$, we can rewrite
this as%
\[
\dfrac{1}{1-x}\left[  x+x^{2}\right]  =\dfrac{1\left[  x+x^{2}\right]
}{\left(  1-x\right)  \left[  x+x^{2}\right]  }.
\]
In view of $1\left[  x+x^{2}\right]  =1$ and $\left(  1-x\right)  \left[
x+x^{2}\right]  =1-\left(  x+x^{2}\right)  =1-x-x^{2}$, this rewrites as
$\dfrac{1}{1-x}\left[  x+x^{2}\right]  =\dfrac{1}{1-x-x^{2}}$. Thus,
(\ref{eq.exa.fps.subs.fibonacci.subs-plausible}) is proved.
\end{example}

Let us summarize: If $f,g\in K\left[  \left[  x\right]  \right]  $ are two
FPSs, then the composition $f\left[  g\right]  $ is not always defined.
However, it is defined at least in the following two cases:

\begin{itemize}
\item in the case when $f$ is a polynomial (that is, $f\in K\left[  x\right]
$), and

\item in the case when $g$ has constant term $0$ (that is, $\left[
x^{0}\right]  g=0$).
\end{itemize}

This justifies some more of the things we did back in Section
\ref{sec.gf.exas}; in particular, Example 1 from that section is now fully
justified. But we still have not defined (e.g.) the square root of an FPS,
which we used in Example 2.

Before I explain square roots, let me quickly survey differentiation of FPSs.

\subsection{Derivatives of FPSs}

Our definition of the derivative of a FPS copycats the well-known formula for
the derivative of a power series in analysis:

\begin{definition}
\label{def.fps.deriv}Let $f\in K\left[  \left[  x\right]  \right]  $ be an
FPS. Then, the \emph{derivative} $f^{\prime}$ of $f$ is an FPS defined as
follows: Write $f$ as $f=\sum_{n\in\mathbb{N}}f_{n}x^{n}$ (with $f_{0}%
,f_{1},f_{2},\ldots\in K$), and set%
\[
f^{\prime}:=\sum_{n>0}nf_{n}x^{n-1}.
\]

\end{definition}

To make sure that this derivative behaves nicely, we need to check that it
satisfies the familiar properties of derivatives. And indeed, it does:

\begin{theorem}
\label{thm.fps.deriv.rules}\textbf{(a)} We have $\left(  f+g\right)  ^{\prime
}=f^{\prime}+g^{\prime}$ for any $f,g\in K\left[  \left[  x\right]  \right]
$. \medskip

\textbf{(b)} If $\left(  f_{i}\right)  _{i\in I}$ is a summable family of
FPSs, then the family $\left(  f_{i}^{\prime}\right)  _{i\in I}$ is summable
as well, and we have%
\[
\left(  \sum_{i\in I}f_{i}\right)  ^{\prime}=\sum_{i\in I}f_{i}^{\prime}.
\]

\textbf{(c)} We have $\left(  cf\right)  ^{\prime}=cf^{\prime}$ for any $c\in
K$ and $f\in K\left[  \left[  x\right]  \right]  $. \medskip

\textbf{(d)} We have $\left(  fg\right)  ^{\prime}=f^{\prime}g+fg^{\prime}$
for any $f,g\in K\left[  \left[  x\right]  \right]  $. (This is known as the
\emph{Leibniz rule}.) \medskip

\textbf{(e)} If $f,g\in K\left[  \left[  x\right]  \right]  $ are two FPSs
such that $g$ is invertible, then%
\[
\left(  \dfrac{f}{g}\right)  ^{\prime}=\dfrac{f^{\prime}g-fg^{\prime}}{g^{2}%
}.
\]
(This is known as the \emph{quotient rule}.) \medskip

\textbf{(f)} If $g\in K\left[  \left[  x\right]  \right]  $ is an FPS, then
$\left(  g^{n}\right)  ^{\prime}=ng^{\prime}g^{n-1}$ for any $n\in\mathbb{N}$
(where the expression $ng^{\prime}g^{n-1}$ is to be understood as $0$ if
$n=0$). \medskip

\textbf{(g)} Given two FPSs $f,g\in K\left[  \left[  x\right]  \right]  $, we
have
\[
\left(  f\circ g\right)  ^{\prime}=\left(  f^{\prime}\circ g\right)  \cdot
g^{\prime}%
\]
if $f$ is a polynomial or if $\left[  x^{0}\right]  g=0$. (This is known as
the \emph{chain rule}.) \medskip

\textbf{(h)} If $K$ is a $\mathbb{Q}$-algebra, and if two FPSs $f,g\in
K\left[  \left[  x\right]  \right]  $ satisfy $f^{\prime}=g^{\prime}$, then
$f-g$ is constant.
\end{theorem}

Theorem \ref{thm.fps.deriv.rules} justifies Example 4 in Section
\ref{sec.gf.exas} (specifically, Theorem \ref{thm.fps.deriv.rules}
\textbf{(e)} is the quotient rule that we used to compute $\left(  \dfrac
{1}{1-x}\right)  ^{\prime}$).

\begin{proof}
[Proof of Theorem \ref{thm.fps.deriv.rules} (sketched).]\textbf{(a)} This is
part of \cite[Exercise 5 \textbf{(b)}]{19s-mt3s} (specifically, it is
Statement 1 in \cite[solution to Exercise 5 \textbf{(b)}]{19s-mt3s}). Anyway,
the proof is very easy. \medskip

\textbf{(b)} This is just the natural generalization of Theorem
\ref{thm.fps.deriv.rules} \textbf{(a)} to (potentially) infinite sums. The
proof follows the same idea, but requires some straightforward technical
verifications (mainly to check that the summation signs can be interchanged).
\medskip

\textbf{(c)} This is part of \cite[Exercise 5 \textbf{(b)}]{19s-mt3s}
(specifically, it is Statement 3 in \cite[solution to Exercise 5 \textbf{(b)}%
]{19s-mt3s}). \medskip

\textbf{(d)} This is \cite[Exercise 5 \textbf{(c)}]{19s-mt3s} and
\cite[Proposition 0.2 \textbf{(c)}]{logexp}. \medskip

\textbf{(e)} Let $f,g\in K\left[  \left[  x\right]  \right]  $ be two FPSs
such that $g$ is invertible. Then, Theorem \ref{thm.fps.deriv.rules}
\textbf{(d)} (applied to $\dfrac{f}{g}$ instead of $f$) yields $\left(
\dfrac{f}{g}\cdot g\right)  ^{\prime}=\left(  \dfrac{f}{g}\right)  ^{\prime
}\cdot g+\dfrac{f}{g}\cdot g^{\prime}$. In view of $\dfrac{f}{g}\cdot g=f$,
this rewrites as $f^{\prime}=\left(  \dfrac{f}{g}\right)  ^{\prime}\cdot
g+\dfrac{f}{g}\cdot g^{\prime}$. Solving this for $\left(  \dfrac{f}%
{g}\right)  ^{\prime}$, we find $\left(  \dfrac{f}{g}\right)  ^{\prime}%
=\dfrac{f^{\prime}g-fg^{\prime}}{g^{2}}$. This proves Theorem
\ref{thm.fps.deriv.rules} \textbf{(e)}. \medskip

\textbf{(f)} This follows by induction on $n$, using part \textbf{(d)} (in the
induction step) and $1^{\prime}=0$ (in the induction base). \medskip

\textbf{(g)} Let $f,g\in K\left[  \left[  x\right]  \right]  $ be two FPSs
such that $f$ is a polynomial or $\left[  x^{0}\right]  g=0$. Write the FPS
$f$ in the form $f=\sum_{n\in\mathbb{N}}f_{n}x^{n}$ with $f_{0},f_{1}%
,f_{2},\ldots\in K$. Then, either Definition \ref{def.fps.subs} or Definition
\ref{def.pol.subs} (depending on whether we have $\left[  x^{0}\right]  g=0$
or $f$ is a polynomial) yields $f\left[  g\right]  =\sum_{n\in\mathbb{N}}%
f_{n}g^{n}$. Hence,%
\begin{align}
\left(  f\left[  g\right]  \right)  ^{\prime}  &  =\left(  \sum_{n\in
\mathbb{N}}f_{n}g^{n}\right)  ^{\prime}=\sum_{n\in\mathbb{N}}%
\underbrace{\left(  f_{n}g^{n}\right)  ^{\prime}}_{\substack{=f_{n}\left(
g^{n}\right)  ^{\prime}\\\text{(by Theorem \ref{thm.fps.deriv.rules}
\textbf{(c)})}}}\ \ \ \ \ \ \ \ \ \ \left(  \text{by Theorem
\ref{thm.fps.deriv.rules} \textbf{(b)}}\right) \nonumber\\
&  =\sum_{n\in\mathbb{N}}f_{n}\left(  g^{n}\right)  ^{\prime}=f_{0}%
\underbrace{\left(  g^{0}\right)  ^{\prime}}_{=1^{\prime}=0}+\sum_{n>0}%
f_{n}\underbrace{\left(  g^{n}\right)  ^{\prime}}_{\substack{=ng^{\prime
}g^{n-1}\\\text{(by Theorem \ref{thm.fps.deriv.rules} \textbf{(f)})}%
}}=\underbrace{f_{0}\cdot0}_{=0}+\sum_{n>0}f_{n}ng^{\prime}g^{n-1}\nonumber\\
&  =\sum_{n>0}f_{n}n\underbrace{g^{\prime}g^{n-1}}_{=g^{n-1}g^{\prime}}%
=\sum_{n>0}nf_{n}g^{n-1}g^{\prime}. \label{pf.thm.fps.deriv.rules.g.3}%
\end{align}

On the other hand, from $f=\sum_{n\in\mathbb{N}}f_{n}x^{n}$, we obtain
$f^{\prime}=\sum_{n>0}nf_{n}x^{n-1}=\sum_{m\in\mathbb{N}}\left(  m+1\right)
f_{m+1}x^{m}$ (here, we have substituted $m+1$ for $n$ in the sum). Hence,
Definition \ref{def.fps.subs} or Definition \ref{def.pol.subs} (depending on
whether we have $\left[  x^{0}\right]  g=0$ or $f$ is a polynomial) yields
\[
f^{\prime}\left[  g\right]  =\sum_{m\in\mathbb{N}}\left(  m+1\right)
f_{m+1}g^{m}=\sum_{n>0}nf_{n}g^{n-1}%
\]
(here, we have substituted $n-1$ for $m$ in the sum). Multiplying both sides
of this equality by $g^{\prime}$, we find%
\[
f^{\prime}\left[  g\right]  \cdot g^{\prime}=\left(  \sum_{n>0}nf_{n}%
g^{n-1}\right)  \cdot g^{\prime}=\sum_{n>0}nf_{n}g^{n-1}g^{\prime}.
\]
Comparing this with (\ref{pf.thm.fps.deriv.rules.g.3}), we find $\left(
f\left[  g\right]  \right)  ^{\prime}=f^{\prime}\left[  g\right]  \cdot
g^{\prime}$. In other words, $\left(  f\circ g\right)  ^{\prime}=\left(
f^{\prime}\circ g\right)  \cdot g^{\prime}$ (since $f\circ g$ is a synonym for
$f\left[  g\right]  $). This proves Theorem \ref{thm.fps.deriv.rules}
\textbf{(g)}. (Note that this proof is done in \cite[proof of Theorem 7.57
(d)]{Loehr-BC} in the case when $f$ is a polynomial.) \medskip

\textbf{(h)} The proof is easy and can be found in \cite[Lemma 0.3]{logexp}.
Note that the condition that $K$ be a $\mathbb{Q}$-algebra is crucial, since
it allows dividing by positive integers. (For example, if $K$ could be
$\mathbb{Z}/2$, then we would easily get a counterexample, e.g., by taking
$f=x^{2}$ and $g=0$.)
\end{proof}

\subsection{\label{subsec.fps.exp}Exponentials and logarithms}

\begin{convention}
\label{conv.fps.exp.K-Q-alg}Throughout Section \ref{subsec.fps.exp}, we assume
that $K$ is not just a commutative ring, but actually a commutative
$\mathbb{Q}$-algebra.
\end{convention}

\subsubsection{Definitions}

Convention \ref{conv.fps.exp.K-Q-alg} entails that elements of $K$ can be
divided by the positive integers $1,2,3,\ldots$. We can use this to define
three specific (and particularly important) FPSs over $K$:

\begin{definition}
\label{def.fps.exp-log}Define three FPSs $\exp$, $\overline{\log}$ and
$\overline{\exp}$ in $K\left[  \left[  x\right]  \right]  $ by%
\begin{align*}
\exp &  :=\sum_{n\in\mathbb{N}}\dfrac{1}{n!}x^{n},\\
\overline{\log}  &  :=\sum_{n\geq1}\dfrac{\left(  -1\right)  ^{n-1}}{n}%
x^{n},\\
\overline{\exp}  &  :=\exp-1=\sum_{n\geq1}\dfrac{1}{n!}x^{n}.
\end{align*}

(The last equality sign here follows from $\exp=\sum_{n\in\mathbb{N}}\dfrac
{1}{n!}x^{n}=\underbrace{\dfrac{1}{0!}}_{=1}\underbrace{x^{0}}_{=1}%
+\sum_{n\geq1}\dfrac{1}{n!}x^{n}=1+\sum_{n\geq1}\dfrac{1}{n!}x^{n}$.)
\end{definition}

Note that the FPS $\exp$ is the usual exponential series from analysis, but
now manifesting itself as a FPS. Likewise, $\overline{\log}$ is the Mercator
series for $\log\left(  1+x\right)  $, where $\log$ stands for the natural
logarithm function. The natural logarithm function itself cannot be
interpreted as an FPS, since $\log0$ is undefined.

\subsubsection{The exponential and the logarithm are inverse}

I will prove that
\begin{equation}
\overline{\exp}\circ\overline{\log}=\overline{\log}\circ\overline{\exp}=x.
\label{eq.fps.exp-log-inverse.1}%
\end{equation}
This is an algebraic analogue of the well-known fact from analysis which
states that the exponential and logarithm functions are mutually inverse.

There is a short way of proving (\ref{eq.fps.exp-log-inverse.1}), which I will
not take: Namely, one can show that any equality between holomorphic functions
on an open disk around the origin leads to an equality between their Taylor
series (viewed as FPSs). Thus, if you have proved in complex analysis that
$\log\circ\exp=\operatorname*{id}$ on an open disk around $0$ and $\exp
\circ\log=\operatorname*{id}$ on an open disk around $1$, then you
automatically get (\ref{eq.fps.exp-log-inverse.1}) (indeed, $\overline{\exp}$
and $\overline{\log}$ are the Taylor series of the functions $\exp$ and $\log$
around $0$ and $1$, with the reservation that the point $1$ has been moved to
the origin by a shift). This approach uses nontrivial results from complex
analysis, so I will not follow it and instead start from scratch.

The main tool in the proof of (\ref{eq.fps.exp-log-inverse.1}) will be the
following useful proposition (\cite[Lemma 0.4]{logexp}):

\begin{proposition}
\label{prop.fps.exp-log-der}Let $g\in K\left[  \left[  x\right]  \right]  $
with $\left[  x^{0}\right]  g=0$. Then: \medskip

\textbf{(a)} We have%
\[
\left(  \overline{\exp}\circ g\right)  ^{\prime}=\left(  \exp\circ g\right)
^{\prime}=\left(  \exp\circ g\right)  \cdot g^{\prime}.
\]

\textbf{(b)} We have%
\[
\left(  \overline{\log}\circ g\right)  ^{\prime}=\left(  1+g\right)
^{-1}\cdot g^{\prime}.
\]

\end{proposition}

\begin{proof}
[Proof of Proposition \ref{prop.fps.exp-log-der}.]\textbf{(a)} Let us first
show that $\overline{\exp}^{\prime}=\exp^{\prime}=\exp$. Indeed,
$\overline{\exp}=\exp-1$, so that $\exp=\overline{\exp}+1$ and therefore%
\begin{align}
\exp^{\prime}  &  =\left(  \overline{\exp}+1\right)  ^{\prime}=\overline{\exp
}^{\prime}+\underbrace{1^{\prime}}_{=0}\ \ \ \ \ \ \ \ \ \ \left(  \text{by
Theorem \ref{thm.fps.deriv.rules} \textbf{(a)}}\right) \nonumber\\
&  =\overline{\exp}^{\prime}. \label{pf.prop.fps.exp-log-der.a.1}%
\end{align}
Next, we recall that $\exp=\sum_{n\in\mathbb{N}}\dfrac{1}{n!}x^{n}$. Hence,
the definition of a derivative yields%
\begin{align}
\exp^{\prime}  &  =\sum_{n\geq1}\underbrace{n\cdot\dfrac{1}{n!}}%
_{\substack{=\dfrac{1}{\left(  n-1\right)  !}\\\text{(since }n!=n\cdot\left(
n-1\right)  !\text{)}}}x^{n-1}=\sum_{n\geq1}\dfrac{1}{\left(  n-1\right)
!}x^{n-1}=\sum_{n\in\mathbb{N}}\dfrac{1}{n!}x^{n}\nonumber\\
&  \ \ \ \ \ \ \ \ \ \ \ \ \ \ \ \ \ \ \ \ \left(  \text{here, we have
substituted }n\text{ for }n-1\text{ in the sum}\right) \nonumber\\
&  =\exp. \label{pf.prop.fps.exp-log-der.a.2}%
\end{align}
Comparing this with (\ref{pf.prop.fps.exp-log-der.a.1}), we find
\begin{equation}
\overline{\exp}^{\prime}=\exp. \label{pf.prop.fps.exp-log-der.a.3}%
\end{equation}

Now, we can apply the chain rule (Theorem \ref{thm.fps.deriv.rules}
\textbf{(g)}) to $f=\overline{\exp}$ (since $\left[  x^{0}\right]  g=0$), and
thus obtain%
\[
\left(  \overline{\exp}\circ g\right)  ^{\prime}=\left(  \underbrace{\overline
{\exp}^{\prime}}_{=\exp}\circ\,g\right)  \cdot g^{\prime}=\left(  \exp\circ
g\right)  \cdot g^{\prime}.
\]
The same computation (but with $\overline{\exp}$ replaced by $\exp$) yields
$\left(  \exp\circ g\right)  ^{\prime}=\left(  \exp\circ g\right)  \cdot
g^{\prime}$. Combining these two formulas, we obtain $\left(  \overline{\exp
}\circ g\right)  ^{\prime}=\left(  \exp\circ g\right)  ^{\prime}=\left(
\exp\circ g\right)  \cdot g^{\prime}$. Thus, we have proved Proposition
\ref{prop.fps.exp-log-der} \textbf{(a)}.

\textbf{(b)} We have $\overline{\log}=\sum_{n\geq1}\dfrac{\left(  -1\right)
^{n-1}}{n}x^{n}$. Thus,
\begin{align*}
\overline{\log}^{\prime}  &  =\left(  \sum_{n\geq1}\dfrac{\left(  -1\right)
^{n-1}}{n}x^{n}\right)  ^{\prime}\\
&  =\sum_{n\geq1}\dfrac{\left(  -1\right)  ^{n-1}}{n}\underbrace{\left(
x^{n}\right)  ^{\prime}}_{\substack{=nx^{\prime}x^{n-1}\\\text{(by Theorem
\ref{thm.fps.deriv.rules} \textbf{(f)},}\\\text{applied to }x\text{ instead of
}g\text{)}}}\ \ \ \ \ \ \ \ \ \ \left(  \text{by Theorem
\ref{thm.fps.deriv.rules} \textbf{(b)}}\right) \\
&  =\sum_{n\geq1}\underbrace{\dfrac{\left(  -1\right)  ^{n-1}}{n}n}_{=\left(
-1\right)  ^{n-1}}\underbrace{x^{\prime}}_{=1}x^{n-1}=\sum_{n\geq1}\left(
-1\right)  ^{n-1}x^{n-1}=\sum_{n\in\mathbb{N}}\left(  -1\right)  ^{n}x^{n}%
\end{align*}
(here, we have substituted $n$ for $n-1$ in the sum). On the other hand,
Proposition \ref{prop.fps.invertible.1+x} yields $\left(  1+x\right)
^{-1}=\sum_{n\in\mathbb{N}}\left(  -1\right)  ^{n}x^{n}$. Comparing these two
equalities, we find%
\begin{equation}
\overline{\log}^{\prime}=\left(  1+x\right)  ^{-1}.
\label{pf.prop.fps.exp-log-der.b.2}%
\end{equation}

Now, we can apply the chain rule (Theorem \ref{thm.fps.deriv.rules}
\textbf{(g)}) to $f=\overline{\log}$ (since $\left[  x^{0}\right]  g=0$), and
thus obtain%
\begin{equation}
\left(  \overline{\log}\circ g\right)  ^{\prime}=\left(  \underbrace{\overline
{\log}^{\prime}}_{=\left(  1+x\right)  ^{-1}}\circ\,g\right)  \cdot g^{\prime
}=\left(  \left(  1+x\right)  ^{-1}\circ g\right)  \cdot g^{\prime}.
\label{pf.prop.fps.exp-log-der.b.3}%
\end{equation}

However, we claim that $\left(  1+x\right)  ^{-1}\circ g=\left(  1+g\right)
^{-1}$. Indeed, Proposition \ref{prop.fps.subs.rules} \textbf{(c)} (applied to
$f_{1}=1$ and $f_{2}=1+x$) yields $\dfrac{1}{1+x}\circ g=\dfrac{1\circ
g}{\left(  1+x\right)  \circ g}$ (since the FPS $1+x$ is invertible). In view
of $\dfrac{1}{1+x}=\left(  1+x\right)  ^{-1}$ and
\[
\underbrace{1}_{=\underline{1}}\circ\,g=\underline{1}\circ g=\underline{1}%
\ \ \ \ \ \ \ \ \ \ \left(  \text{by Proposition \ref{prop.fps.subs.rules}
\textbf{(f)}, applied to }a=1\right)
\]
and%
\[
\left(  1+x\right)  \circ g=1+g\ \ \ \ \ \ \ \ \ \ \left(  \text{this follows
easily from Definition \ref{def.fps.subs}}\right)  ,
\]
this rewrites as $\left(  1+x\right)  ^{-1}\circ g=\dfrac{\underline{1}}%
{1+g}=\left(  1+g\right)  ^{-1}$. Hence, (\ref{pf.prop.fps.exp-log-der.b.3})
becomes%
\[
\left(  \overline{\log}\circ g\right)  ^{\prime}=\underbrace{\left(  \left(
1+x\right)  ^{-1}\circ g\right)  }_{=\left(  1+g\right)  ^{-1}}\cdot
\,g^{\prime}=\left(  1+g\right)  ^{-1}\cdot g^{\prime}.
\]
Thus, we have proved Proposition \ref{prop.fps.exp-log-der} \textbf{(b)}.
\end{proof}

We will need a very simple lemma, which says (in particular) that if two FPSs
have constant terms $0$, then so does their composition:

\begin{lemma}
\label{lem.fps.compos-cst-term-0}Let $f,g\in K\left[  \left[  x\right]
\right]  $ be two FPSs with $\left[  x^{0}\right]  g=0$. Then, $\left[
x^{0}\right]  \left(  f\circ g\right)  =\left[  x^{0}\right]  f$.
\end{lemma}

\begin{proof}
[Proof of Lemma \ref{lem.fps.compos-cst-term-0}.]Write $f$ in the form
$f=\sum_{n\in\mathbb{N}}f_{n}x^{n}$ with $f_{0},f_{1},f_{2},\ldots\in K$.
Thus, $f_{0}=\left[  x^{0}\right]  f$. Now, Definition \ref{def.fps.subs}
yields $f\left[  g\right]  =\sum_{n\in\mathbb{N}}f_{n}g^{n}$. However,
Proposition \ref{prop.fps.subs.wd} \textbf{(c)} yields $\left[  x^{0}\right]
\left(  \sum_{n\in\mathbb{N}}f_{n}g^{n}\right)  =f_{0}=\left[  x^{0}\right]
f$. In view of $f\circ g=f\left[  g\right]  =\sum_{n\in\mathbb{N}}f_{n}g^{n}$,
we can rewrite this as $\left[  x^{0}\right]  \left(  f\circ g\right)
=\left[  x^{0}\right]  f$. This proves Lemma \ref{lem.fps.compos-cst-term-0}.
\end{proof}

Now, we can prove the equalities we promised (\cite[Theorem 0.1]{logexp}):

\begin{theorem}
\label{thm.fps.exp-log-inv}We have%
\[
\overline{\exp}\circ\overline{\log}=x\ \ \ \ \ \ \ \ \ \ \text{and}%
\ \ \ \ \ \ \ \ \ \ \overline{\log}\circ\overline{\exp}=x.
\]

\end{theorem}

\begin{proof}
[Proof of Theorem \ref{thm.fps.exp-log-inv}.]Let us first prove that
$\overline{\log}\circ\overline{\exp}=x$.

Indeed, the idea of this proof is to show that $\overline{\log}\circ
\overline{\exp}$ and $x$ are two FPSs with the same constant term (namely,
$0$) and with the same derivative. Once this is proved, Theorem
\ref{thm.fps.deriv.rules} \textbf{(h)} will easily yield that they are equal.

Let us fill in the details. We have $\left[  x^{0}\right]  \overline{\exp}=0$
(since $\overline{\exp}=\sum_{n\geq1}\dfrac{1}{n!}x^{n}$). Hence, Lemma
\ref{lem.fps.compos-cst-term-0} (applied to $f=\overline{\log}$ and
$g=\overline{\exp}$) yields $\left[  x^{0}\right]  \left(  \overline{\log
}\circ\overline{\exp}\right)  =\left[  x^{0}\right]  \overline{\log}=0$ (since
$\overline{\log}=\sum_{n\geq1}\dfrac{\left(  -1\right)  ^{n-1}}{n}x^{n}$).
Now, (\ref{pf.thm.fps.ring.xn(a-b)=}) yields%
\begin{equation}
\left[  x^{0}\right]  \left(  \overline{\log}\circ\overline{\exp}-x\right)
=\underbrace{\left[  x^{0}\right]  \left(  \overline{\log}\circ\overline{\exp
}\right)  }_{=0}-\underbrace{\left[  x^{0}\right]  x}_{=0}=0.
\label{pf.thm.fps.exp-log-inv.a.1}%
\end{equation}
However, $\overline{\exp}=\exp-1$ and thus $1+\overline{\exp}=\exp$. Now,
Proposition \ref{prop.fps.exp-log-der} \textbf{(b)} (applied to $g=\overline
{\exp}$) yields
\[
\left(  \overline{\log}\circ\overline{\exp}\right)  ^{\prime}=\left(
\underbrace{1+\overline{\exp}}_{=\exp}\right)  ^{-1}\cdot\underbrace{\overline
{\exp}^{\prime}}_{\substack{=\exp\\\text{(by
(\ref{pf.prop.fps.exp-log-der.a.3}))}}}=\exp^{-1}\cdot\exp=1=x^{\prime}%
\]
(since $x^{\prime}=1$). Hence, Theorem \ref{thm.fps.deriv.rules} \textbf{(h)}
(applied to $f=\overline{\log}\circ\overline{\exp}$ and $g=x$) yields that
$\overline{\log}\circ\overline{\exp}-x$ is constant. In other words,
$\overline{\log}\circ\overline{\exp}-x=\underline{a}$ for some $a\in K$.
Consider this $a$. From $\overline{\log}\circ\overline{\exp}-x=\underline{a}$,
we obtain $\left[  x^{0}\right]  \left(  \overline{\log}\circ\overline{\exp
}-x\right)  =\left[  x^{0}\right]  \underline{a}=a$. Comparing this with
(\ref{pf.thm.fps.exp-log-inv.a.1}), we find $a=0$. Hence, $\overline{\log
}\circ\overline{\exp}-x=\underline{a}$ rewrites as $\overline{\log}%
\circ\overline{\exp}-x=\underline{0}$. In other words, $\overline{\log}%
\circ\overline{\exp}=x$.

Now it remains to prove that $\overline{\exp}\circ\overline{\log}=x$. There
are (at least) two ways to do this:

\begin{itemize}
\item \textit{1st way:} A homework exercise (Exercise \ref{exe.fps.comp-inv})
says that any FPS $f$ with $\left[  x^{0}\right]  f=0$ and with $\left[
x^{1}\right]  f$ invertible has a unique compositional inverse (i.e., there is
a unique FPS $g$ with $\left[  x^{0}\right]  g=0$ and $f\circ g=g\circ f=x$).
We can apply this to $f=\overline{\log}$ (since $\left[  x^{0}\right]
\overline{\log}=0$ and since $\left[  x^{1}\right]  \overline{\log}=1$ is
invertible), and thus see that $\overline{\log}$ has a unique compositional
inverse $g$. This compositional inverse $g$ must be $\overline{\exp}$, since
$\overline{\log}\circ\overline{\exp}=x$ (indeed, comparing
$\underbrace{\left(  g\circ\overline{\log}\right)  }_{=x}\circ\,\overline
{\exp}=x\circ\overline{\exp}=\overline{\exp}$ with
\begin{align*}
\left(  g\circ\overline{\log}\right)  \circ\overline{\exp}  &  =g\circ
\underbrace{\left(  \overline{\log}\circ\overline{\exp}\right)  }%
_{=x}\ \ \ \ \ \ \ \ \ \ \left(  \text{by Proposition
\ref{prop.fps.subs.rules} \textbf{(e)}}\right) \\
&  =g\circ x=g
\end{align*}
yields $g=\overline{\exp}$). But this entails that $\overline{\exp}%
\circ\overline{\log}=x$ as well.

\item \textit{2nd way:} Here is a more direct argument. We shall first show
that $\exp\circ\overline{\log}=1+x$.

To wit: The FPS $1+x$ is invertible (by Proposition
\ref{prop.fps.invertible.1+x}). Thus, applying the quotient rule (Theorem
\ref{thm.fps.deriv.rules} \textbf{(e)}) to $f=\exp\circ\,\overline{\log}$ and
$g=1+x$, we obtain%
\[
\left(  \dfrac{\exp\circ\,\overline{\log}}{1+x}\right)  ^{\prime}%
=\dfrac{\left(  \exp\circ\,\overline{\log}\right)  ^{\prime}\cdot\left(
1+x\right)  -\left(  \exp\circ\,\overline{\log}\right)  \cdot\left(
1+x\right)  ^{\prime}}{\left(  1+x\right)  ^{2}}.
\]
In view of%
\begin{align*}
\left(  \exp\circ\,\overline{\log}\right)  ^{\prime}  &  =\left(  \exp
\circ\,\overline{\log}\right)  \cdot\underbrace{\overline{\log}^{\prime}%
}_{\substack{=\left(  1+x\right)  ^{-1}\\\text{(by
(\ref{pf.prop.fps.exp-log-der.b.2}))}}}\\
&  \ \ \ \ \ \ \ \ \ \ \ \ \ \ \ \ \ \ \ \ \left(  \text{by Proposition
\ref{prop.fps.exp-log-der} \textbf{(a)}, applied to }g=\overline{\log}\right)
\\
&  =\left(  \exp\circ\,\overline{\log}\right)  \cdot\left(  1+x\right)  ^{-1}%
\end{align*}
and $\left(  1+x\right)  ^{\prime}=1$, we can rewrite this as%
\begin{align*}
\left(  \dfrac{\exp\circ\,\overline{\log}}{1+x}\right)  ^{\prime}  &
=\dfrac{\left(  \exp\circ\,\overline{\log}\right)  \cdot\left(  1+x\right)
^{-1}\cdot\left(  1+x\right)  -\left(  \exp\circ\,\overline{\log}\right)
\cdot1}{\left(  1+x\right)  ^{2}}\\
&  =\dfrac{\left(  \exp\circ\,\overline{\log}\right)  -\left(  \exp
\circ\,\overline{\log}\right)  }{\left(  1+x\right)  ^{2}}=0=0^{\prime}.
\end{align*}
Thus, Theorem \ref{thm.fps.deriv.rules} \textbf{(h)} (applied to
$f=\dfrac{\exp\circ\,\overline{\log}}{1+x}$ and $g=0$) yields that
$\dfrac{\exp\circ\,\overline{\log}}{1+x}-0$ is constant. In other words,
$\dfrac{\exp\circ\,\overline{\log}}{1+x}$ is constant. In other words,
$\dfrac{\exp\circ\,\overline{\log}}{1+x}=\underline{a}$ for some $a\in K$.
Consider this $a$. From $\dfrac{\exp\circ\,\overline{\log}}{1+x}%
=\underline{a}$, we obtain $\exp\circ\,\overline{\log}=\underline{a}\left(
1+x\right)  =a\left(  1+x\right)  $. Thus,%
\[
\left[  x^{0}\right]  \left(  \exp\circ\,\overline{\log}\right)  =\left[
x^{0}\right]  \left(  a\left(  1+x\right)  \right)  =a.
\]
However, it is easy to see that $\left[  x^{0}\right]  \left(  \exp
\circ\,\overline{\log}\right)  =1$\ \ \ \ \footnote{\textit{Proof.} Recall
that $\left[  x^{0}\right]  \overline{\log}=0$. Hence, Lemma
\ref{lem.fps.compos-cst-term-0} (applied to $f=\exp$ and $g=\overline{\log}$)
yields
\begin{align*}
\left[  x^{0}\right]  \left(  \exp\circ\,\overline{\log}\right)   &  =\left[
x^{0}\right]  \exp=\dfrac{1}{0!}\ \ \ \ \ \ \ \ \ \ \left(  \text{since }%
\exp=\sum_{n\in\mathbb{N}}\dfrac{1}{n!}x^{n}\right) \\
&  =\dfrac{1}{1}=1.
\end{align*}
}. Comparing these two equalities, we find $a=1$. Thus, $\exp\circ
\,\overline{\log}=\underbrace{a}_{=1}\left(  1+x\right)  =1+x$.

Now, $\overline{\exp}=\exp-1=\exp+\,\underline{-1}$. Hence,
\begin{align*}
&  \overline{\exp}\circ\overline{\log}\\
&  =\left(  \exp+\,\underline{-1}\right)  \circ\overline{\log}%
=\underbrace{\exp\circ\,\overline{\log}}_{=1+x}+\underbrace{\underline{-1}%
\circ\overline{\log}}_{\substack{=\underline{-1}\\\text{(by Proposition
\ref{prop.fps.subs.rules} \textbf{(f)},}\\\text{applied to }-1\text{ and
}\overline{\log}\text{ instead of }a\text{ and }g\text{)}}}\\
&  \ \ \ \ \ \ \ \ \ \ \ \ \ \ \ \ \ \ \ \ \left(
\begin{array}
[c]{c}%
\text{by Proposition \ref{prop.fps.subs.rules} \textbf{(a)},}\\
\text{applied to }f_{1}=\exp\text{ and }f_{2}=\underline{-1}\text{ and
}g=\overline{\log}%
\end{array}
\right) \\
&  =1+x+\underline{-1}=\underbrace{\left(  1+\underline{-1}\right)  }%
_{=0}+\,x=x.
\end{align*}

\end{itemize}

Either way, we have shown that $\overline{\exp}\circ\overline{\log}=x$. Thus,
the proof of Theorem \ref{thm.fps.exp-log-inv} is complete.
\end{proof}

\subsubsection{The exponential and the logarithm of an FPS}

In Definition \ref{def.fps.exp-log}, we have found algebraic versions of the
exponential and logarithm functions as FPSs. Next, we shall define analogues
of these functions as operators acting on FPSs (i.e., analogues not of the
functions $\exp$ and $\log$ themselves, but rather of composition with these functions):

\begin{definition}
\label{def.fps.Exp-Log-maps}\textbf{(a)} We let $K\left[  \left[  x\right]
\right]  _{0}$ denote the set of all FPSs $f\in K\left[  \left[  x\right]
\right]  $ with $\left[  x^{0}\right]  f=0$. \medskip

\textbf{(b)} We let $K\left[  \left[  x\right]  \right]  _{1}$ denote the set
of all FPSs $f\in K\left[  \left[  x\right]  \right]  $ with $\left[
x^{0}\right]  f=1$. \medskip

\textbf{(c)} We define two maps%
\begin{align*}
\operatorname*{Exp}:K\left[  \left[  x\right]  \right]  _{0}  &  \rightarrow
K\left[  \left[  x\right]  \right]  _{1},\\
g  &  \mapsto\exp\circ g
\end{align*}
and%
\begin{align*}
\operatorname*{Log}:K\left[  \left[  x\right]  \right]  _{1}  &  \rightarrow
K\left[  \left[  x\right]  \right]  _{0},\\
f  &  \mapsto\overline{\log}\circ\left(  f-1\right)  .
\end{align*}
(These two maps are well-defined according to parts \textbf{(c)} and
\textbf{(d)} of Lemma \ref{lem.fps.Exp-Log-maps-wd} below.)
\end{definition}

The maps $\operatorname*{Exp}$ and $\operatorname*{Log}$ are algebraic
analogues of the maps in complex analysis that take any holomorphic function
$f$ to its exponential and logarithm, respectively (at least within certain
regions in which these things are well-defined). As one would hope, and as we
will soon see, they are mutually inverse. Let us first check that their
definition is justified:

\begin{lemma}
\label{lem.fps.Exp-Log-maps-wd}\textbf{(a)} For any $f,g\in K\left[  \left[
x\right]  \right]  _{0}$, we have $f\circ g\in K\left[  \left[  x\right]
\right]  _{0}$. \medskip

\textbf{(b)} For any $f\in K\left[  \left[  x\right]  \right]  _{1}$ and $g\in
K\left[  \left[  x\right]  \right]  _{0}$, we have $f\circ g\in K\left[
\left[  x\right]  \right]  _{1}$. \medskip

\textbf{(c)} For any $g\in K\left[  \left[  x\right]  \right]  _{0}$, we have
$\exp\circ g\in K\left[  \left[  x\right]  \right]  _{1}$. \medskip

\textbf{(d)} For any $f\in K\left[  \left[  x\right]  \right]  _{1}$, we have
$f-1\in K\left[  \left[  x\right]  \right]  _{0}$ and $\overline{\log}%
\circ\left(  f-1\right)  \in K\left[  \left[  x\right]  \right]  _{0}$.
\end{lemma}

\begin{proof}
[Proof of Lemma \ref{lem.fps.Exp-Log-maps-wd}.]\textbf{(a)} Let $f,g\in
K\left[  \left[  x\right]  \right]  _{0}$. In view of the definition of
$K\left[  \left[  x\right]  \right]  _{0}$, this entails that $\left[
x^{0}\right]  f=0$ and $\left[  x^{0}\right]  g=0$. Hence, Lemma
\ref{lem.fps.compos-cst-term-0} yields $\left[  x^{0}\right]  \left(  f\circ
g\right)  =\left[  x^{0}\right]  f=0$. In other words, $f\circ g\in K\left[
\left[  x\right]  \right]  _{0}$ (by the definition of $K\left[  \left[
x\right]  \right]  _{0}$). This proves Lemma \ref{lem.fps.Exp-Log-maps-wd}
\textbf{(a)}.

\textbf{(b)} This is analogous to the proof of Lemma
\ref{lem.fps.Exp-Log-maps-wd} \textbf{(a)}.

\textbf{(c)} Let $g\in K\left[  \left[  x\right]  \right]  _{0}$. From
$\exp=\sum_{n\in\mathbb{N}}\dfrac{1}{n!}x^{n}$, we obtain $\left[
x^{0}\right]  \exp=\dfrac{1}{0!}=1$, so that $\exp\in K\left[  \left[
x\right]  \right]  _{1}$. Hence, Lemma \ref{lem.fps.Exp-Log-maps-wd}
\textbf{(b)} (applied to $f=\exp$) yields $\exp\circ g\in K\left[  \left[
x\right]  \right]  _{1}$. This proves Lemma \ref{lem.fps.Exp-Log-maps-wd}
\textbf{(c)}.

\textbf{(d)} Let $f\in K\left[  \left[  x\right]  \right]  _{1}$. Thus,
$\left[  x^{0}\right]  f=1$. Now, (\ref{pf.thm.fps.ring.xn(a-b)=}) yields
$\left[  x^{0}\right]  \left(  f-1\right)  =\underbrace{\left[  x^{0}\right]
f}_{=1}-\underbrace{\left[  x^{0}\right]  1}_{=1}=1-1=0$, so that $f-1\in
K\left[  \left[  x\right]  \right]  _{0}$. Furthermore, $\left[  x^{0}\right]
\overline{\log}=0$ (since $\overline{\log}=\sum_{n\geq1}\dfrac{\left(
-1\right)  ^{n-1}}{n}x^{n}$) and thus $\overline{\log}\in K\left[  \left[
x\right]  \right]  _{0}$. Hence, Lemma \ref{lem.fps.Exp-Log-maps-wd}
\textbf{(a)} (applied to $\overline{\log}$ and $f-1$ instead of $f$ and $g$)
yields $\overline{\log}\circ\left(  f-1\right)  \in K\left[  \left[  x\right]
\right]  _{0}$. Thus, Lemma \ref{lem.fps.Exp-Log-maps-wd} \textbf{(d)} is proven.
\end{proof}

\begin{lemma}
\label{lem.fps.Exp-Log-maps-inv}The maps $\operatorname*{Exp}$ and
$\operatorname*{Log}$ are mutually inverse bijections between $K\left[
\left[  x\right]  \right]  _{0}$ and $K\left[  \left[  x\right]  \right]
_{1}$.
\end{lemma}

\begin{proof}
[Proof of Lemma \ref{lem.fps.Exp-Log-maps-inv}.]First, we make a simple
auxiliary observation: Each $g\in K\left[  \left[  x\right]  \right]  _{0}$
satisfies\footnote{As before, the \textquotedblleft$\circ$\textquotedblright%
\ operation behaves like multiplication in the sense of PEMDAS conventions.
Thus, the expression \textquotedblleft$\overline{\exp}\circ g+1$%
\textquotedblright\ means $\left(  \overline{\exp}\circ g\right)  +1$.}%
\begin{equation}
\exp\circ g=\overline{\exp}\circ g+1. \label{pf.lem.fps.Exp-Log-maps-inv.1}%
\end{equation}

[\textit{Proof of (\ref{pf.lem.fps.Exp-Log-maps-inv.1}):} Let $g\in K\left[
\left[  x\right]  \right]  _{0}$. Recall that $\overline{\exp}=\exp-1$, so
that $\exp=\overline{\exp}+1=\overline{\exp}+\underline{1}$. Hence,%
\[
\exp\circ g=\left(  \overline{\exp}+\underline{1}\right)  \circ g=\overline
{\exp}\circ g+\underline{1}\circ g
\]
(by Proposition \ref{prop.fps.subs.rules} \textbf{(a)}, applied to
$f_{1}=\overline{\exp}$ and $f_{2}=\underline{1}$). However, Proposition
\ref{prop.fps.subs.rules} \textbf{(f)} (applied to $a=1$) yields
$\underline{1}\circ g=\underline{1}=1$. Hence, $\exp\circ g=\overline{\exp
}\circ g+\underbrace{\underline{1}\circ g}_{=1}=\overline{\exp}\circ g+1$.
This proves (\ref{pf.lem.fps.Exp-Log-maps-inv.1}).]

Now, let us show that $\operatorname*{Exp}\circ\operatorname*{Log}%
=\operatorname*{id}$. Indeed, we fix some $f\in K\left[  \left[  x\right]
\right]  _{1}$. Then, $f-1\in K\left[  \left[  x\right]  \right]  _{0}$ (by
Lemma \ref{lem.fps.Exp-Log-maps-wd} \textbf{(d)}). Hence, Proposition
\ref{prop.fps.subs.rules} \textbf{(e)} (applied to $\overline{\exp}$,
$\overline{\log}$ and $f-1$ instead of $f$, $g$ and $h$) yields $\left(
\overline{\exp}\circ\overline{\log}\right)  \circ\left(  f-1\right)
=\overline{\exp}\circ\left(  \overline{\log}\circ\left(  f-1\right)  \right)
$. Thus,%
\begin{align}
\overline{\exp}\circ\left(  \overline{\log}\circ\left(  f-1\right)  \right)
&  =\underbrace{\left(  \overline{\exp}\circ\overline{\log}\right)
}_{\substack{=x\\\text{(by Theorem \ref{thm.fps.exp-log-inv})}}}\circ\left(
f-1\right)  =x\circ\left(  f-1\right) \nonumber\\
&  =f-1 \label{pf.lem.fps.Exp-Log-maps-inv.a.1}%
\end{align}
(by Proposition \ref{prop.fps.subs.rules} \textbf{(g)}, applied to $g=f-1$).
However,%
\begin{align*}
\left(  \operatorname*{Exp}\circ\operatorname*{Log}\right)  \left(  f\right)
&  =\operatorname*{Exp}\left(  \operatorname*{Log}f\right) \\
&  =\exp\circ\left(  \operatorname*{Log}f\right)  \ \ \ \ \ \ \ \ \ \ \left(
\text{by the definition of }\operatorname*{Exp}\right) \\
&  =\overline{\exp}\circ\underbrace{\left(  \operatorname*{Log}f\right)
}_{\substack{=\overline{\log}\circ\left(  f-1\right)  \\\text{(by the
definition of }\operatorname*{Log}\text{)}}}+\,1\\
&  \ \ \ \ \ \ \ \ \ \ \ \ \ \ \ \ \ \ \ \ \left(  \text{by
(\ref{pf.lem.fps.Exp-Log-maps-inv.1}), applied to }g=\operatorname*{Log}%
f\right) \\
&  =\underbrace{\overline{\exp}\circ\left(  \overline{\log}\circ\left(
f-1\right)  \right)  }_{\substack{=f-1\\\text{(by
(\ref{pf.lem.fps.Exp-Log-maps-inv.a.1}))}}}+\,1\\
&  =\left(  f-1\right)  +1=f=\operatorname*{id}\left(  f\right)  .
\end{align*}

Forget that we fixed $f$. We thus have shown that $\left(  \operatorname*{Exp}%
\circ\operatorname*{Log}\right)  \left(  f\right)  =\operatorname*{id}\left(
f\right)  $ for each $f\in K\left[  \left[  x\right]  \right]  _{1}$. In other
words, $\operatorname*{Exp}\circ\operatorname*{Log}=\operatorname*{id}$.

Using a similar argument, we can show that $\operatorname*{Log}\circ
\operatorname*{Exp}=\operatorname*{id}$. Indeed, let us fix some $g\in
K\left[  \left[  x\right]  \right]  _{0}$. Hence, Proposition
\ref{prop.fps.subs.rules} \textbf{(e)} (applied to $\overline{\log}$,
$\overline{\exp}$ and $g$ instead of $f$, $g$ and $h$) yields $\left(
\overline{\log}\circ\overline{\exp}\right)  \circ g=\overline{\log}%
\circ\left(  \overline{\exp}\circ g\right)  $. Thus,%
\begin{align}
\overline{\log}\circ\left(  \overline{\exp}\circ g\right)   &
=\underbrace{\left(  \overline{\log}\circ\overline{\exp}\right)
}_{\substack{=x\\\text{(by Theorem \ref{thm.fps.exp-log-inv})}}}\circ
\,g=x\circ g\nonumber\\
&  =g \label{pf.lem.fps.Exp-Log-maps-inv.b.1}%
\end{align}
(by Proposition \ref{prop.fps.subs.rules} \textbf{(g)}). But the definition of
$\operatorname*{Exp}$ yields $\operatorname*{Exp}g=\exp\circ g=\overline{\exp
}\circ g+1$ (by (\ref{pf.lem.fps.Exp-Log-maps-inv.1})). Hence,
$\operatorname*{Exp}g-1=\overline{\exp}\circ g$. Now,%
\begin{align*}
\left(  \operatorname*{Log}\circ\operatorname*{Exp}\right)  \left(  g\right)
&  =\operatorname*{Log}\left(  \operatorname*{Exp}g\right) \\
&  =\overline{\log}\circ\left(  \operatorname*{Exp}g-1\right)
\ \ \ \ \ \ \ \ \ \ \left(  \text{by the definition of }\operatorname*{Log}%
\right) \\
&  =\overline{\log}\circ\left(  \overline{\exp}\circ g\right)
\ \ \ \ \ \ \ \ \ \ \left(  \text{since }\operatorname*{Exp}g-1=\overline
{\exp}\circ g\right) \\
&  =g\ \ \ \ \ \ \ \ \ \ \left(  \text{by
(\ref{pf.lem.fps.Exp-Log-maps-inv.b.1})}\right) \\
&  =\operatorname*{id}\left(  g\right)  .
\end{align*}

Forget that we fixed $g$. We thus have shown that $\left(  \operatorname*{Log}%
\circ\operatorname*{Exp}\right)  \left(  g\right)  =\operatorname*{id}\left(
g\right)  $ for each $g\in K\left[  \left[  x\right]  \right]  _{0}$. In other
words, $\operatorname*{Log}\circ\operatorname*{Exp}=\operatorname*{id}$.
Combining this with $\operatorname*{Exp}\circ\operatorname*{Log}%
=\operatorname*{id}$, we see that the maps $\operatorname*{Exp}$ and
$\operatorname*{Log}$ are mutually inverse bijections between $K\left[
\left[  x\right]  \right]  _{0}$ and $K\left[  \left[  x\right]  \right]
_{1}$. This proves Lemma \ref{lem.fps.Exp-Log-maps-inv}.
\end{proof}

\subsubsection{Addition to multiplication}

We will now prove another lemma, which says that the $\operatorname*{Exp}$ and
$\operatorname*{Log}$ maps deserve their names:

\begin{lemma}
\label{lem.fps.Exp-Log-additive}\textbf{(a)} For any $f,g\in K\left[  \left[
x\right]  \right]  _{0}$, we have%
\[
\operatorname*{Exp}\left(  f+g\right)  =\operatorname*{Exp}f\cdot
\operatorname*{Exp}g.
\]

\textbf{(b)} For any $f,g\in K\left[  \left[  x\right]  \right]  _{1}$, we
have%
\[
\operatorname*{Log}\left(  fg\right)  =\operatorname*{Log}%
f+\operatorname*{Log}g.
\]

\end{lemma}

\begin{proof}
[Proof of Lemma \ref{lem.fps.Exp-Log-additive} (sketched).]\textbf{(a)} Like
many of our arguments involving FPSs, this will be a short computation
followed by lengthy technical arguments justifying the interchanges of
summation signs. (In this aspect, our algebraic replica of the analysis of
infinite sums doesn't differ that much from the original.) We begin with the
computation; the justifying arguments will be sketched afterwards.

Let $f,g\in K\left[  \left[  x\right]  \right]  _{0}$. Thus, $\left[
x^{0}\right]  f=0$ and $\left[  x^{0}\right]  g=0$. Hence, $f+g\in K\left[
\left[  x\right]  \right]  _{0}$ (since (\ref{pf.thm.fps.ring.xn(a+b)=})
yields $\left[  x^{0}\right]  \left(  f+g\right)  =\underbrace{\left[
x^{0}\right]  f}_{=0}+\underbrace{\left[  x^{0}\right]  g}_{=0}=0$).

By the definition of $\operatorname*{Exp}$, we have%
\[
\operatorname*{Exp}f=\exp\circ f=\exp\left[  f\right]  =\sum_{n\in\mathbb{N}%
}\dfrac{1}{n!}f^{n}%
\]
(by Definition \ref{def.fps.subs}, since $\exp=\sum_{n\in\mathbb{N}}\dfrac
{1}{n!}x^{n}$). Similarly,%
\[
\operatorname*{Exp}g=\sum_{n\in\mathbb{N}}\dfrac{1}{n!}g^{n}%
\]
and%
\[
\operatorname*{Exp}\left(  f+g\right)  =\sum_{n\in\mathbb{N}}\dfrac{1}%
{n!}\left(  f+g\right)  ^{n}.
\]

Now, the latter equality becomes%
\begin{align*}
\operatorname*{Exp}\left(  f+g\right)   &  =\sum_{n\in\mathbb{N}}\dfrac{1}%
{n!}\underbrace{\left(  f+g\right)  ^{n}}_{\substack{=\sum\limits_{k=0}%
^{n}\dbinom{n}{k}f^{k}g^{n-k}\\\text{(by the binomial theorem)}}}=\sum
_{n\in\mathbb{N}}\dfrac{1}{n!}\sum\limits_{k=0}^{n}\dbinom{n}{k}f^{k}g^{n-k}\\
&  =\sum_{n\in\mathbb{N}}\ \ \sum\limits_{k=0}^{n}\underbrace{\dfrac{1}%
{n!}\dbinom{n}{k}}_{\substack{=\dfrac{1}{k!\left(  n-k\right)  !}\\\text{(by
(\ref{eq.binom.fac-form}))}}}f^{k}g^{n-k}=\underbrace{\sum_{n\in\mathbb{N}%
}\ \ \sum\limits_{k=0}^{n}}_{\substack{=\sum_{\substack{\left(  n,k\right)
\in\mathbb{N}\times\mathbb{N};\\k\leq n}}\\=\sum_{k\in\mathbb{N}}%
\ \ \sum\limits_{n\geq k}}}\dfrac{1}{k!\left(  n-k\right)  !}f^{k}g^{n-k}\\
&  =\sum_{k\in\mathbb{N}}\ \ \sum\limits_{n\geq k}\dfrac{1}{k!\left(
n-k\right)  !}f^{k}g^{n-k}=\sum_{k\in\mathbb{N}}\ \ \sum\limits_{\ell
\in\mathbb{N}}\dfrac{1}{k!\ell!}f^{k}g^{\ell}\\
&  \ \ \ \ \ \ \ \ \ \ \ \ \ \ \ \ \ \ \ \ \left(  \text{here, we have
substituted }\ell\text{ for }n-k\text{ in the second sum}\right)  .
\end{align*}
Comparing this with%
\begin{align*}
\operatorname*{Exp}f\cdot\operatorname*{Exp}g  &  =\left(  \sum_{k\in
\mathbb{N}}\dfrac{1}{k!}f^{k}\right)  \cdot\left(  \sum_{\ell\in\mathbb{N}%
}\dfrac{1}{\ell!}g^{\ell}\right) \\
&  \ \ \ \ \ \ \ \ \ \ \ \ \ \ \ \ \ \ \ \ \left(
\begin{array}
[c]{c}%
\text{since }\operatorname*{Exp}f=\sum_{n\in\mathbb{N}}\dfrac{1}{n!}f^{n}%
=\sum_{k\in\mathbb{N}}\dfrac{1}{k!}f^{k}\\
\text{and }\operatorname*{Exp}g=\sum_{n\in\mathbb{N}}\dfrac{1}{n!}g^{n}%
=\sum_{\ell\in\mathbb{N}}\dfrac{1}{\ell!}g^{\ell}%
\end{array}
\right) \\
&  =\sum_{k\in\mathbb{N}}\ \ \sum_{\ell\in\mathbb{N}}\dfrac{1}{k!}f^{k}%
\cdot\dfrac{1}{\ell!}g^{\ell}=\sum_{k\in\mathbb{N}}\ \ \sum_{\ell\in
\mathbb{N}}\dfrac{1}{k!\ell!}f^{k}g^{\ell},
\end{align*}
we obtain $\operatorname*{Exp}\left(  f+g\right)  =\operatorname*{Exp}%
f\cdot\operatorname*{Exp}g$.

\begin{fineprint}
This is sufficient to prove Lemma \ref{lem.fps.Exp-Log-additive} \textbf{(a)}
if we can justify the above manipulations of infinite sums. Actually, there is
just one manipulation that we need to justify, and that is our replacement of
\textquotedblleft$\sum_{n\in\mathbb{N}}\ \ \sum\limits_{k=0}^{n}%
$\textquotedblright\ by \textquotedblleft$\sum_{k\in\mathbb{N}}\ \ \sum
\limits_{n\geq k}$\textquotedblright. This is an application of the
\textquotedblleft discrete Fubini rule\textquotedblright\ (specifically, of a
version thereof in which the summation is over all pairs $\left(  n,k\right)
\in\mathbb{N}\times\mathbb{N}$ satisfying $k\leq n$). In order to justify this
manipulation, we need to show that the family $\left(  \dfrac{1}{k!\left(
n-k\right)  !}f^{k}g^{n-k}\right)  _{\left(  n,k\right)  \in\mathbb{N}%
\times\mathbb{N}\text{ satisfying }k\leq n}$ is summable. In other words, we
need to show the following statement:

\begin{statement}
\textit{Statement 1:} For each $m\in\mathbb{N}$, all but finitely many pairs
$\left(  n,k\right)  \in\mathbb{N}\times\mathbb{N}$ satisfying $k\leq n$
satisfy $\left[  x^{m}\right]  \left(  \dfrac{1}{k!\left(  n-k\right)  !}%
f^{k}g^{n-k}\right)  =0$.
\end{statement}

We shall achieve this by proving the following statement:

\begin{statement}
\textit{Statement 2:} For any three nonnegative integers $m,k,\ell$ with
$m<k+\ell$, we have $\left[  x^{m}\right]  \left(  f^{k}g^{\ell}\right)  =0$.
\end{statement}

[\textit{Proof of Statement 2:} Let $m,k,\ell$ be three nonnegative integers
with $m<k+\ell$. We must show that $\left[  x^{m}\right]  \left(  f^{k}%
g^{\ell}\right)  =0$.

We have $\left[  x^{0}\right]  f=0$. Hence, Lemma \ref{lem.fps.g=xh} (applied
to $a=f$) shows that there exists an $h\in K\left[  \left[  x\right]  \right]
$ such that $f=xh$. Consider this $h$ and denote it by $u$. Thus, $u\in
K\left[  \left[  x\right]  \right]  $ and $f=xu$.

We have $\left[  x^{0}\right]  g=0$. Hence, Lemma \ref{lem.fps.g=xh} (applied
to $a=g$) shows that there exists an $h\in K\left[  \left[  x\right]  \right]
$ such that $g=xh$. Consider this $h$ and denote it by $v$. Thus, $v\in
K\left[  \left[  x\right]  \right]  $ and $g=xv$.

Now, from $f=xu$ and $g=xv$, we obtain $f^{k}g^{\ell}=\left(  xu\right)
^{k}\left(  xv\right)  ^{\ell}=x^{k}u^{k}x^{\ell}v^{\ell}=x^{k+\ell}%
u^{k}v^{\ell}$. However, Lemma \ref{lem.fps.first-n-coeffs-of-xna} (applied to
$k+\ell$ and $u^{k}v^{\ell}$ instead of $k$ and $a$) shows that the first
$k+\ell$ coefficients of the FPS $x^{k+\ell}u^{k}v^{\ell}$ are $0$. In other
words, the first $k+\ell$ coefficients of the FPS $f^{k}g^{\ell}$ are $0$
(since $f^{k}g^{\ell}=x^{k+\ell}u^{k}v^{\ell}$). But $\left[  x^{m}\right]
\left(  f^{k}g^{\ell}\right)  $ is one of these first $k+\ell$ coefficients
(since $m<k+\ell$). Thus, we conclude that $\left[  x^{m}\right]  \left(
f^{k}g^{\ell}\right)  =0$. This proves Statement 2.]

Note that Statement 2 entails that the family $\left(  f^{k}g^{\ell}\right)
_{\left(  k,\ell\right)  \in\mathbb{N}\times\mathbb{N}}$ is summable (because
when $m\in\mathbb{N}$ is given, all but finitely many pairs $\left(
k,\ell\right)  \in\mathbb{N}\times\mathbb{N}$ satisfy $m<k+\ell$). However, we
need to prove Statement 1, so let us do this:

[\textit{Proof of Statement 1:} Let $m\in\mathbb{N}$. If $\left(  n,k\right)
\in\mathbb{N}\times\mathbb{N}$ is a pair satisfying $k\leq n$ and $m<n$, then
\begin{align*}
\left[  x^{m}\right]  \left(  \dfrac{1}{k!\left(  n-k\right)  !}f^{k}%
g^{n-k}\right)   &  =\dfrac{1}{k!\left(  n-k\right)  !}\underbrace{\left[
x^{m}\right]  \left(  f^{k}g^{n-k}\right)  }_{\substack{=0\\\text{(by
Statement 2}\\\text{(applied to }\ell=n-k\text{),}\\\text{since }m<n=k+\left(
n-k\right)  \text{)}}}\ \ \ \ \ \ \ \ \ \ \left(  \text{by
(\ref{pf.thm.fps.ring.xn(la)=})}\right) \\
&  =0.
\end{align*}
Thus, all but finitely many pairs $\left(  n,k\right)  \in\mathbb{N}%
\times\mathbb{N}$ satisfying $k\leq n$ satisfy \newline$\left[  x^{m}\right]
\left(  \dfrac{1}{k!\left(  n-k\right)  !}f^{k}g^{n-k}\right)  =0$ (because
all but finitely many such pairs satisfy $m<n$). This proves Statement 1.]

As explained above, Statement 1 shows that the family \newline$\left(
\dfrac{1}{k!\left(  n-k\right)  !}f^{k}g^{n-k}\right)  _{\left(  n,k\right)
\in\mathbb{N}\times\mathbb{N}\text{ satisfying }k\leq n}$ is summable, and
thus our interchange of summation signs made above is justified. This
completes our proof of Lemma \ref{lem.fps.Exp-Log-additive} \textbf{(a)}.
\end{fineprint}

\bigskip

\textbf{(b)} This easily follows from part \textbf{(a)}, since we know that
$\operatorname*{Log}$ is inverse to $\operatorname*{Exp}$. Here are the details:

Let $f,g\in K\left[  \left[  x\right]  \right]  _{1}$. Set
$u=\operatorname*{Log}f$ and $v=\operatorname*{Log}g$; then, $u,v\in K\left[
\left[  x\right]  \right]  _{0}$ (since $\operatorname*{Log}$ is a map from
$K\left[  \left[  x\right]  \right]  _{1}$ to $K\left[  \left[  x\right]
\right]  _{0}$). Hence, Lemma \ref{lem.fps.Exp-Log-additive} \textbf{(a)}
(applied to $u$ and $v$ instead of $f$ and $g$) yields $\operatorname*{Exp}%
\left(  u+v\right)  =\operatorname*{Exp}u\cdot\operatorname*{Exp}v$.

However, Lemma \ref{lem.fps.Exp-Log-maps-inv} says that the maps
$\operatorname*{Exp}$ and $\operatorname*{Log}$ are mutually inverse
bijections between $K\left[  \left[  x\right]  \right]  _{0}$ and $K\left[
\left[  x\right]  \right]  _{1}$. Hence, $\operatorname*{Exp}\circ
\operatorname*{Log}=\operatorname*{id}$ and $\operatorname*{Log}%
\circ\operatorname*{Exp}=\operatorname*{id}$.

Now, from $u=\operatorname*{Log}f$, we obtain $\operatorname*{Exp}%
u=\operatorname*{Exp}\left(  \operatorname*{Log}f\right)  =\underbrace{\left(
\operatorname*{Exp}\circ\operatorname*{Log}\right)  }_{=\operatorname*{id}%
}\left(  f\right)  =\operatorname*{id}\left(  f\right)  =f$. Similarly,
$\operatorname*{Exp}v=g$. Multiplying these two equalities, we find
$\operatorname*{Exp}u\cdot\operatorname*{Exp}v=fg$. Now, we have%
\[
\operatorname*{Log}\left(  \operatorname*{Exp}\left(  u+v\right)  \right)
=\underbrace{\left(  \operatorname*{Log}\circ\operatorname*{Exp}\right)
}_{=\operatorname*{id}}\left(  u+v\right)  =\operatorname*{id}\left(
u+v\right)  =u+v=\operatorname*{Log}f+\operatorname*{Log}g
\]
(since $u=\operatorname*{Log}f$ and $v=\operatorname*{Log}g$). In view of
$\operatorname*{Exp}\left(  u+v\right)  =\operatorname*{Exp}u\cdot
\operatorname*{Exp}v=fg$, this rewrites as $\operatorname*{Log}\left(
fg\right)  =\operatorname*{Log}f+\operatorname*{Log}g$. This proves Lemma
\ref{lem.fps.Exp-Log-additive} \textbf{(b)}.
\end{proof}

We can neatly pack the last few lemmas into a single theorem through the use
of group isomorphisms. To this purpose, we need to observe that $K\left[
\left[  x\right]  \right]  _{0}$ is a group under addition and $K\left[
\left[  x\right]  \right]  _{1}$ is a group under multiplication:

\begin{proposition}
\label{prop.fps.Exp-Log-groups}\textbf{(a)} The subset $K\left[  \left[
x\right]  \right]  _{0}$ of $K\left[  \left[  x\right]  \right]  $ is closed
under addition and subtraction and contains $0$, and thus forms a group
$\left(  K\left[  \left[  x\right]  \right]  _{0},+,0\right)  $. \medskip

\textbf{(b)} The subset $K\left[  \left[  x\right]  \right]  _{1}$ of
$K\left[  \left[  x\right]  \right]  $ is closed under multiplication and
division and contains $1$, and thus forms a group $\left(  K\left[  \left[
x\right]  \right]  _{1},\cdot,1\right)  $.
\end{proposition}

\begin{proof}
[Proof of Proposition \ref{prop.fps.Exp-Log-groups}.]\textbf{(a)} It is clear
that the set $K\left[  \left[  x\right]  \right]  _{0}$ contains the FPS $0$
(since $\left[  x^{0}\right]  0=0$). Thus, it remains to show that $K\left[
\left[  x\right]  \right]  _{0}$ is closed under addition and subtraction. But
this is easy: If $f,g\in K\left[  \left[  x\right]  \right]  _{0}$, then
$\left[  x^{0}\right]  f=0$ and $\left[  x^{0}\right]  g=0$, and therefore
$f+g\in K\left[  \left[  x\right]  \right]  _{0}$ (since
(\ref{pf.thm.fps.ring.xn(a+b)=}) yields $\left[  x^{0}\right]  \left(
f+g\right)  =\underbrace{\left[  x^{0}\right]  f}_{=0}+\underbrace{\left[
x^{0}\right]  g}_{=0}=0$) and $f-g\in K\left[  \left[  x\right]  \right]
_{0}$ (by a similar argument using (\ref{pf.thm.fps.ring.xn(a-b)=})). Thus,
$K\left[  \left[  x\right]  \right]  _{0}$ is closed under addition and
subtraction. This proves Proposition \ref{prop.fps.Exp-Log-groups}
\textbf{(a)}.

\textbf{(b)} Any $a\in K\left[  \left[  x\right]  \right]  _{1}$ is invertible
in $K\left[  \left[  x\right]  \right]  $ (indeed, $a\in K\left[  \left[
x\right]  \right]  _{1}$ shows that $\left[  x^{0}\right]  a=1$; thus,
$\left[  x^{0}\right]  a$ is invertible in $K$; therefore, Proposition
\ref{prop.fps.invertible} entails that $a$ is invertible in $K\left[  \left[
x\right]  \right]  $). Hence, $\dfrac{f}{g}$ is well-defined for any $f,g\in
K\left[  \left[  x\right]  \right]  _{1}$.

Next, we claim that $K\left[  \left[  x\right]  \right]  _{1}$ is closed under
multiplication. Indeed, if $f,g\in K\left[  \left[  x\right]  \right]  _{1}$,
then $\left[  x^{0}\right]  f=1$ and $\left[  x^{0}\right]  g=1$, and
therefore $fg\in K\left[  \left[  x\right]  \right]  _{1}$ (since
(\ref{pf.thm.fps.ring.x0(ab)=}) yields $\left[  x^{0}\right]  \left(
fg\right)  =\underbrace{\left[  x^{0}\right]  f}_{=1}\cdot\underbrace{\left[
x^{0}\right]  g}_{=1}=1$). This shows that $K\left[  \left[  x\right]
\right]  _{1}$ is closed under multiplication.

It remains to prove that $K\left[  \left[  x\right]  \right]  _{1}$ is closed
under division. Indeed, if $f,g\in K\left[  \left[  x\right]  \right]  _{1}$,
then $\left[  x^{0}\right]  f=1$ and $\left[  x^{0}\right]  g=1$, and
therefore $\dfrac{f}{g}\in K\left[  \left[  x\right]  \right]  _{1}$ (because
we have $f=\dfrac{f}{g}\cdot g$ and thus
\begin{align*}
\left[  x^{0}\right]  f  &  =\left[  x^{0}\right]  \left(  \dfrac{f}{g}\cdot
g\right)  =\left[  x^{0}\right]  \dfrac{f}{g}\cdot\underbrace{\left[
x^{0}\right]  g}_{=1}\ \ \ \ \ \ \ \ \ \ \left(  \text{by
(\ref{pf.thm.fps.ring.x0(ab)=})}\right) \\
&  =\left[  x^{0}\right]  \dfrac{f}{g}%
\end{align*}
and thus $\left[  x^{0}\right]  \dfrac{f}{g}=\left[  x^{0}\right]  f=1$, so
that $\dfrac{f}{g}\in K\left[  \left[  x\right]  \right]  _{1}$). This shows
that $K\left[  \left[  x\right]  \right]  _{1}$ is closed under division.
Thus, Proposition \ref{prop.fps.Exp-Log-groups} \textbf{(b)} is proven.
\end{proof}

The two groups in Proposition \ref{prop.fps.Exp-Log-groups} can now be
connected through $\operatorname*{Exp}$ and $\operatorname*{Log}$:

\begin{theorem}
\label{thm.fps.Exp-Log-group-iso}The maps
\[
\operatorname*{Exp}:\left(  K\left[  \left[  x\right]  \right]  _{0}%
,+,0\right)  \rightarrow\left(  K\left[  \left[  x\right]  \right]  _{1}%
,\cdot,1\right)
\]
and%
\[
\operatorname*{Log}:\left(  K\left[  \left[  x\right]  \right]  _{1}%
,\cdot,1\right)  \rightarrow\left(  K\left[  \left[  x\right]  \right]
_{0},+,0\right)
\]
are mutually inverse group isomorphisms.
\end{theorem}

\begin{proof}
[Proof of Theorem \ref{thm.fps.Exp-Log-group-iso} (sketched).]Lemma
\ref{lem.fps.Exp-Log-additive} yields that these two maps are group
homomorphisms\footnote{Here, we are using the following fact: If $\left(
G,\ast,e_{G}\right)  $ and $\left(  H,\ast,e_{H}\right)  $ are any two groups,
and if $\Phi:G\rightarrow H$ is a map such that every $f,g\in G$ satisfy
$\Phi\left(  f\ast g\right)  =\Phi\left(  f\right)  \ast\Phi\left(  g\right)
$, then $\Phi$ is a group homomorphism.}. Lemma \ref{lem.fps.Exp-Log-maps-inv}
shows that they are mutually inverse. Combining these results, we conclude
that these two maps are mutually inverse group isomorphisms. This proves
Theorem \ref{thm.fps.Exp-Log-group-iso}.
\end{proof}

Theorem \ref{thm.fps.Exp-Log-group-iso} helps us turn addition into
multiplication and vice versa when it comes to FPSs, at least if the constant
terms are the right ones. This will come useful rather soon.

\subsubsection{The logarithmic derivative}

For future use, we shall define and briefly study one more concept related to logarithms:

\begin{definition}
\label{def.fps.loder.1}In this definition, we do \textbf{not} use Convention
\ref{conv.fps.exp.K-Q-alg}; thus, $K$ can be an arbitrary commutative ring.
However, we set $K\left[  \left[  x\right]  \right]  _{1}=\left\{  f\in
K\left[  \left[  x\right]  \right]  \ \mid\ \left[  x^{0}\right]  f=1\right\}
$.

For any FPS $f\in K\left[  \left[  x\right]  \right]  _{1}$, we define the
\emph{logarithmic derivative} $\operatorname*{loder}f\in K\left[  \left[
x\right]  \right]  $ to be the FPS $\dfrac{f^{\prime}}{f}$. (This is
well-defined, since $f$ is easily seen to be invertible.\footnotemark)
\end{definition}

\footnotetext{Indeed: Let $f\in K\left[  \left[  x\right]  \right]  _{1}$.
Then, $\left[  x^{0}\right]  f=1$. Thus, $\left[  x^{0}\right]  f$ is
invertible. Hence, Proposition \ref{prop.fps.invertible} (applied to $a=f$)
yields that $f$ is invertible.}The reason why this FPS $\operatorname*{loder}%
f$ is called the logarithmic derivative of $f$ is made clear by the following
simple fact:

\begin{proposition}
\label{prop.fps.loder.log}Let $K$ be a commutative $\mathbb{Q}$-algebra. Let
$f\in K\left[  \left[  x\right]  \right]  _{1}$ be an FPS. Then,
$\operatorname*{loder}f=\left(  \operatorname*{Log}f\right)  ^{\prime}$.
\end{proposition}

\begin{proof}
[Proof of Proposition \ref{prop.fps.loder.log}.]The definition of
$\operatorname*{Log}$ yields $\operatorname*{Log}f=\overline{\log}\circ\left(
f-1\right)  $.

From $f\in K\left[  \left[  x\right]  \right]  _{1}$, we obtain $\left[
x^{0}\right]  f=1=\left[  x^{0}\right]  1$. Now, $\left[  x^{0}\right]
\left(  f-1\right)  =\left[  x^{0}\right]  f-\left[  x^{0}\right]  1=0$ (since
$\left[  x^{0}\right]  f=\left[  x^{0}\right]  1$). Hence, Proposition
\ref{prop.fps.exp-log-der} \textbf{(b)} (applied to $g=f-1$) yields
\[
\left(  \overline{\log}\circ\left(  f-1\right)  \right)  ^{\prime}=\left(
\underbrace{1+\left(  f-1\right)  }_{=f}\right)  ^{-1}\cdot\left(  f-1\right)
^{\prime}=f^{-1}\cdot\left(  f-1\right)  ^{\prime}=\dfrac{\left(  f-1\right)
^{\prime}}{f}.
\]
In view of $\operatorname*{Log}f=\overline{\log}\circ\left(  f-1\right)  $, we
can rewrite this as%
\[
\left(  \operatorname*{Log}f\right)  ^{\prime}=\dfrac{\left(  f-1\right)
^{\prime}}{f}.
\]

Since $\left(  f-1\right)  ^{\prime}=f^{\prime}$ (which is easy to
see\footnote{\textit{Proof.} Theorem \ref{thm.fps.deriv.rules} \textbf{(a)}
(applied to $g=-1$) yields $\left(  f+\left(  -1\right)  \right)  ^{\prime
}=f^{\prime}+\underbrace{\left(  -1\right)  ^{\prime}}_{=0}=f^{\prime}$. In
other words, $\left(  f-1\right)  ^{\prime}=f^{\prime}$.}), we can rewrite
this further as $\left(  \operatorname*{Log}f\right)  ^{\prime}=\dfrac
{f^{\prime}}{f}=\operatorname*{loder}f$ (since $\operatorname*{loder}f$ is
defined as $\dfrac{f^{\prime}}{f}$). This proves Proposition
\ref{prop.fps.loder.log}.
\end{proof}

Note that Proposition \ref{prop.fps.loder.log} only makes sense when $K$ is a
$\mathbb{Q}$-algebra (otherwise, $\operatorname*{Log}f$ would make no sense),
which is why we are not using it as a definition of $\operatorname*{loder}f$.
The logarithmic derivative is defined even when the logarithm is not!

If you have seen the logarithmic derivative in analysis, you will likely
expect the following property:

\begin{proposition}
\label{prop.fps.loder.prod}Let $f,g\in K\left[  \left[  x\right]  \right]
_{1}$ be two FPSs. Then, $\operatorname*{loder}\left(  fg\right)
=\operatorname*{loder}f+\operatorname*{loder}g$.

(Here, we do \textbf{not} use Convention \ref{conv.fps.exp.K-Q-alg}; thus, $K$
can be an arbitrary commutative ring.)
\end{proposition}

\begin{proof}
[Proof of Proposition \ref{prop.fps.loder.prod}.]The definition of a
logarithmic derivative yields $\operatorname*{loder}f=\dfrac{f^{\prime}}{f}$
and $\operatorname*{loder}g=\dfrac{g^{\prime}}{g}$, but it also yields%
\begin{align*}
\operatorname*{loder}\left(  fg\right)   &  =\dfrac{\left(  fg\right)
^{\prime}}{fg}=\dfrac{f^{\prime}g+fg^{\prime}}{fg}\ \ \ \ \ \ \ \ \ \ \left(
\begin{array}
[c]{c}%
\text{since Theorem \ref{thm.fps.deriv.rules} \textbf{(d)}}\\
\text{says that }\left(  fg\right)  ^{\prime}=f^{\prime}g+fg^{\prime}%
\end{array}
\right) \\
&  =\underbrace{\dfrac{f^{\prime}}{f}}_{=\operatorname*{loder}f}%
+\underbrace{\dfrac{g^{\prime}}{g}}_{=\operatorname*{loder}g}%
=\operatorname*{loder}f+\operatorname*{loder}g.
\end{align*}
This proves Proposition \ref{prop.fps.loder.prod}.
\end{proof}

\begin{corollary}
\label{cor.fps.loder.prodk}Let $f_{1},f_{2},\ldots,f_{k}$ be any $k$ FPSs in
$K\left[  \left[  x\right]  \right]  _{1}$. Then,
\[
\operatorname*{loder}\left(  f_{1}f_{2}\cdots f_{k}\right)
=\operatorname*{loder}\left(  f_{1}\right)  +\operatorname*{loder}\left(
f_{2}\right)  +\cdots+\operatorname*{loder}\left(  f_{k}\right)  .
\]

(Here, we do \textbf{not} use Convention \ref{conv.fps.exp.K-Q-alg}; thus, $K$
can be an arbitrary commutative ring.)
\end{corollary}

\begin{proof}
[Proof of Corollary \ref{cor.fps.loder.prodk}.]Induct on $k$. The \textit{base
case} ($k=0$) requires showing that $\operatorname*{loder}1=0$, which is easy
to see from the definition of a logarithmic derivative (since $1^{\prime}=0$).
The \textit{induction step} follows easily from Proposition
\ref{prop.fps.loder.prod}.
\end{proof}

\begin{corollary}
\label{cor.fps.loder.inv}Let $f$ be any FPS in $K\left[  \left[  x\right]
\right]  _{1}$. Then, $\operatorname*{loder}\left(  f^{-1}\right)
=-\operatorname*{loder}f$.

(Here, we do \textbf{not} use Convention \ref{conv.fps.exp.K-Q-alg}; thus, $K$
can be an arbitrary commutative ring.)
\end{corollary}

\begin{proof}
[Proof of Corollary \ref{cor.fps.loder.inv}.]First of all, $f$ is invertible
(as we showed in Definition \ref{def.fps.loder.1}). This shows that $f^{-1}$
is well-defined.

Proposition \ref{prop.fps.loder.prod} (applied to $g=f^{-1}$) yields
$\operatorname*{loder}\left(  ff^{-1}\right)  =\operatorname*{loder}%
f+\operatorname*{loder}\left(  f^{-1}\right)  $. However, we also have
$\operatorname*{loder}\left(  \underbrace{ff^{-1}}_{=1}\right)
=\operatorname*{loder}1=0$ (since $1^{\prime}=0$). Comparing these two
equalities, we obtain $\operatorname*{loder}f+\operatorname*{loder}\left(
f^{-1}\right)  =0$. In other words, $\operatorname*{loder}\left(
f^{-1}\right)  =-\operatorname*{loder}f$. This proves Corollary
\ref{cor.fps.loder.inv}.
\end{proof}

\subsection{\label{sec.gf.nips}Non-integer powers}

\subsubsection{Definition}

Now, let us again recall Example 2 from Section \ref{sec.gf.exas}. In order to
fully justify that example, we still need to explain what $\sqrt{1-4x}$ is.

More generally, let us try to define non-integer powers of FPSs (since square
roots are just $1/2$-th powers). Thus, we are trying to solve the following problem:

\begin{statement}
\textit{Problem:} Devise a reasonable definition of the $c$-th power $f^{c}$
for any FPS $f\in K\left[  \left[  x\right]  \right]  $ and any $c\in K$.
\end{statement}

Here, \textquotedblleft reasonable\textquotedblright\ means that it should
have some of the properties we would expect:

\begin{itemize}
\item It should not conflict with the existing notion of $f^{c}$ for
$c\in\mathbb{N}$. That is, if $c\in\mathbb{N}$, then our new definition of
$f^{c}$ should yield the same result as the existing meaning that $f^{c}$ has
in this case (namely, $\underbrace{ff\cdots f}_{c\text{ times}}$). The same
should hold for $c\in\mathbb{Z}$ when $f$ is invertible.

\item Rules of exponents should hold: i.e., we should have%
\begin{equation}
f^{a+b}=f^{a}f^{b},\ \ \ \ \ \ \ \ \ \ \left(  fg\right)  ^{a}=f^{a}%
g^{a},\ \ \ \ \ \ \ \ \ \ \left(  f^{a}\right)  ^{b}=f^{ab}
\label{eq.sec.gf.nips.rules-of-exps}%
\end{equation}
for all $a,b\in K$ and $f,g\in K\left[  \left[  x\right]  \right]  $.

\item For any positive integer $n$ and any FPS $f\in K\left[  \left[
x\right]  \right]  $, the $1/n$-th power $f^{1/n}$ should be an $n$-th root of
$f$ (that is, an FPS whose $n$-th power is $f$). (If $K$ is a field of
characteristic $0$, then this actually follows from the previous two
properties, since we can apply the rule $\left(  f^{a}\right)  ^{b}=f^{ab}$ to
$a=1/n$ and $b=n$.)
\end{itemize}

Clearly, we cannot solve the above problem in full generality:

\begin{itemize}
\item The power $0^{-1}$ cannot be reasonably defined (unless $K$ is trivial).
Indeed, $0^{-1}\cdot0^{1}$ would have to equal $0^{-1+1}=0^{0}=1$, but this
would contradict $0^{-1}\cdot0^{1}=0^{-1}\cdot0=0$.

\item The power $x^{1/2}$ cannot be reasonably defined either (unless $K$ is
trivial). Indeed, there is no FPS whose square is $x$. This will be proved in
Exercise \ref{exe.fps.not-squares} \textbf{(a)}.

\item Even the power $\left(  -1\right)  ^{1/2}$ cannot always be defined:
There is no guarantee that $K$ contains a square root of $-1$ (and if $K$ does
not, then it is easy to see that $K\left[  \left[  x\right]  \right]  $ does
not either).
\end{itemize}

However, all we want is to make sense of $\sqrt{1-4x}$, so let us restrict
ourselves to FPSs whose constant term is $1$. Using the notation from
Definition \ref{def.fps.Exp-Log-maps} \textbf{(b)}, we are thus moving on to
the following problem:

\begin{statement}
\textit{More realistic problem:} Devise a reasonable definition of the $c$-th
power $f^{c}$ for any FPS $f\in K\left[  \left[  x\right]  \right]  _{1}$ and
any $c\in K$.
\end{statement}

Besides imposing the above wishlist of properties, we want this $c$-th power
$f^{c}$ itself to belong to $K\left[  \left[  x\right]  \right]  _{1}$, since
otherwise the iterated power $\left(  f^{a}\right)  ^{b}$ in our rules of
exponents might be undefined.

It turns out that this is still too much to ask. Indeed, if $K=\mathbb{Z}/2$,
then the FPS $1+x\in K\left[  \left[  x\right]  \right]  _{1}$ has no square
root (you get to prove this in Exercise \ref{exe.fps.not-squares}
\textbf{(c)}), so its $1/2$-th power $\left(  1+x\right)  ^{1/2}$ cannot be
reasonably defined.

However, if we assume (as in Convention \ref{conv.fps.exp.K-Q-alg}) that $K$
is a commutative $\mathbb{Q}$-algebra, then we get lucky: Our
\textquotedblleft more realistic problem\textquotedblright\ can be solved in
(at least) two ways: \medskip

\textit{1st solution:} We define
\[
\left(  1+x\right)  ^{c}:=\sum_{k\in\mathbb{N}}\dbinom{c}{k}x^{k}%
\ \ \ \ \ \ \ \ \ \ \text{for each }c\in K,
\]
in order to make Newton's binomial formula (Theorem \ref{thm.fps.newton-binom}%
) hold for arbitrary exponents\footnote{Note that $\dbinom{c}{k}%
=\dfrac{c\left(  c-1\right)  \left(  c-2\right)  \cdots\left(  c-k+1\right)
}{k!}$ is well-defined since $K$ is a commutative $\mathbb{Q}$-algebra.}.
Subsequently, we define%
\begin{equation}
f^{c}:=\left(  1+x\right)  ^{c}\left[  f-1\right]
\ \ \ \ \ \ \ \ \ \ \text{for any }f\in K\left[  \left[  x\right]  \right]
_{1}\text{ and }c\in K \label{eq.sec.gf.nips.newton-def2}%
\end{equation}
(in order to have $\left(  1+g\right)  ^{c}=\sum_{k\in\mathbb{N}}\dbinom{c}%
{k}g^{k}$ hold not only for $g=x$, but also for all $g\in K\left[  \left[
x\right]  \right]  _{0}$).

It is clear that the FPS $f^{c}$ is well-defined in this way. However, proving
that this definition satisfies all our wishlist (particularly the rules of
exponents (\ref{eq.sec.gf.nips.rules-of-exps})) is highly nontrivial. Some of
this is done in \cite[\S 7.12]{Loehr-BC}, but it is still a lot of work.

Thus, we shall discard this definition of $f^{c}$, and instead take a
different way: \medskip

\textit{2nd solution:} Recall the mutually inverse group isomorphisms%
\begin{align*}
\operatorname*{Exp}  &  :\left(  K\left[  \left[  x\right]  \right]
_{0},+,0\right)  \rightarrow\left(  K\left[  \left[  x\right]  \right]
_{1},\cdot,1\right)  \ \ \ \ \ \ \ \ \ \ \text{and}\\
\operatorname*{Log}  &  :\left(  K\left[  \left[  x\right]  \right]
_{1},\cdot,1\right)  \rightarrow\left(  K\left[  \left[  x\right]  \right]
_{0},+,0\right)
\end{align*}
from Theorem \ref{thm.fps.Exp-Log-group-iso}. Thus, for any $f\in K\left[
\left[  x\right]  \right]  _{1}$ and any $c\in\mathbb{Z}$, the equation%
\[
\operatorname*{Log}\left(  f^{c}\right)  =c\operatorname*{Log}f
\]
holds (since $\operatorname*{Log}$ is a group homomorphism). This suggests
that we define $f^{c}$ for all $c\in K$ by the same equation. In other words,
we define $f^{c}$ for all $c\in K$ by setting $f^{c}=\operatorname*{Exp}%
\left(  c\operatorname*{Log}f\right)  $ (since the map $\operatorname*{Exp}$
is inverse to $\operatorname*{Log}$). And this is what we shall do now:

\begin{definition}
\label{def.fps.power-c}Assume that $K$ is a commutative $\mathbb{Q}$-algebra.
Let $f\in K\left[  \left[  x\right]  \right]  _{1}$ and $c\in K$. Then, we
define an FPS%
\[
f^{c}:=\operatorname*{Exp}\left(  c\operatorname*{Log}f\right)  \in K\left[
\left[  x\right]  \right]  _{1}.
\]

\end{definition}

This definition of $f^{c}$ does not conflict with our original definition of
$f^{c}$ when $c\in\mathbb{Z}$ because (as we said) the original definition of
$f^{c}$ already satisfies $\operatorname*{Log}\left(  f^{c}\right)
=c\operatorname*{Log}f$ and therefore $f^{c}=\operatorname*{Exp}\left(
c\operatorname*{Log}f\right)  $.

Moreover, Definition \ref{def.fps.power-c} makes the rules of exponents hold:

\begin{theorem}
\label{thm.fps.power-c.rules}Assume that $K$ is a commutative $\mathbb{Q}%
$-algebra. For any $a,b\in K$ and $f,g\in K\left[  \left[  x\right]  \right]
_{1}$, we have
\[
f^{a+b}=f^{a}f^{b},\ \ \ \ \ \ \ \ \ \ \left(  fg\right)  ^{a}=f^{a}%
g^{a},\ \ \ \ \ \ \ \ \ \ \left(  f^{a}\right)  ^{b}=f^{ab}.
\]

\end{theorem}

\begin{proof}
Easy exercise (Exercise \ref{exe.fps.power-c.rules}).
\end{proof}

Now, let us return to Example 2 from Section \ref{sec.gf.exas}. In that
example, we had to solve the quadratic equation
\[
C\left(  x\right)  =1+x\left(  C\left(  x\right)  \right)  ^{2}%
\ \ \ \ \ \ \ \ \ \ \text{for an FPS }C\left(  x\right)  \in\mathbb{Q}\left[
\left[  x\right]  \right]  .
\]
Let us write $C$ for $C\left(  x\right)  $; thus, this quadratic equation
becomes%
\[
C=1+xC^{2}.
\]
By completing the square, we can rewrite this equation in the equivalent form%
\[
\left(  1-2xC\right)  ^{2}=1-4x.
\]
Taking both sides of this equation to the $1/2$-th power, we obtain%
\[
\left(  \left(  1-2xC\right)  ^{2}\right)  ^{1/2}=\left(  1-4x\right)  ^{1/2}%
\]
(since both sides are FPSs with constant term $1$). However, the FPS $1-2xC$
has constant term $1$; thus, the rules of exponents yield $\left(  \left(
1-2xC\right)  ^{2}\right)  ^{1/2}=\left(  1-2xC\right)  ^{2\cdot1/2}=1-2xC$.
Hence,%
\[
1-2xC=\left(  \left(  1-2xC\right)  ^{2}\right)  ^{1/2}=\left(  1-4x\right)
^{1/2}.
\]
This is a linear equation in $C$; solving it for $C$ (using Proposition
\ref{prop.fps.div-by-x-inverts}) yields%
\[
C=\dfrac{1}{2x}\left(  1-\left(  1-4x\right)  ^{1/2}\right)  .
\]
This is precisely the \textquotedblleft square-root\textquotedblright%
\ expression for $C=C\left(  x\right)  $ that we have obtained back in Section
\ref{sec.gf.exas}, but now we have proved it rigorously.

\subsubsection{The Newton binomial formula for arbitrary exponents}

Is Example 2 from Section \ref{sec.gf.exas} fully justified now? No, because
we still need to prove the identity (\ref{eq.sec.gf.exas.2.(1+x)1/2}) that we
used back there. Since we are defining powers in the 2nd way (i.e., using
Definition \ref{def.fps.power-c} rather than using
(\ref{eq.sec.gf.nips.newton-def2})), it is not immediately obvious.
Nevertheless, it can be proved. More generally, we can prove the following:

\begin{theorem}
[Generalized Newton binomial formula]\label{thm.fps.gen-newton}Assume that $K$
is a commutative $\mathbb{Q}$-algebra. Let $c\in K$. Then,%
\[
\left(  1+x\right)  ^{c}=\sum_{k\in\mathbb{N}}\dbinom{c}{k}x^{k}.
\]

\end{theorem}

The following proof illustrates a technique that will probably appear
preposterous if you are seeing it for the first time, but is in fact both
legitimate and rather useful.

\begin{proof}
[Proof of Theorem \ref{thm.fps.gen-newton} (sketched).]The definition of
$\operatorname*{Log}$ yields
\[
\operatorname*{Log}\left(  1+x\right)  =\overline{\log}\circ\left(
\underbrace{\left(  1+x\right)  -1}_{=x}\right)  =\overline{\log}\circ
x=\overline{\log}%
\]
(by Proposition \ref{prop.fps.subs.rules} \textbf{(g)}, applied to
$g=\overline{\log}$).

Now, let us obstinately compute $\left(  1+x\right)  ^{c}$ using Definition
\ref{def.fps.power-c} and the definitions of $\operatorname*{Exp}$ and
$\operatorname*{Log}$. To wit: Let $\mathbb{P}$ denote the set $\left\{
1,2,3,\ldots\right\}  $. By Definition \ref{def.fps.power-c}, we have%
\begin{align}
&  \left(  1+x\right)  ^{c}\nonumber\\
&  =\operatorname*{Exp}\left(  c\operatorname*{Log}\left(  1+x\right)
\right)  =\operatorname*{Exp}\left(  c\,\overline{\log}\right)
\ \ \ \ \ \ \ \ \ \ \left(  \text{since }\operatorname*{Log}\left(
1+x\right)  =\overline{\log}\right) \nonumber\\
&  =\exp\circ\left(  c\,\overline{\log}\right)  \ \ \ \ \ \ \ \ \ \ \left(
\text{by the definition of }\operatorname*{Exp}\right) \nonumber\\
&  =\exp\circ\left(  c\sum_{n\geq1}\dfrac{\left(  -1\right)  ^{n-1}}{n}%
x^{n}\right)  \ \ \ \ \ \ \ \ \ \ \left(  \text{since }\overline{\log}%
=\sum_{n\geq1}\dfrac{\left(  -1\right)  ^{n-1}}{n}x^{n}\right) \nonumber\\
&  =\exp\circ\left(  \sum_{n\geq1}\dfrac{\left(  -1\right)  ^{n-1}}{n}%
cx^{n}\right) \nonumber\\
&  =\sum_{m\in\mathbb{N}}\dfrac{1}{m!}\left(  \sum_{n\geq1}\dfrac{\left(
-1\right)  ^{n-1}}{n}cx^{n}\right)  ^{m} \label{pf.thm.fps.gen-newton.1}%
\end{align}
(by Definition \ref{def.fps.subs}, since $\exp=\sum_{n\in\mathbb{N}}\dfrac
{1}{n!}x^{n}=\sum_{m\in\mathbb{N}}\dfrac{1}{m!}x^{m}$).

Now, fix $m\in\mathbb{N}$. We shall expand $\left(  \sum_{n\geq1}%
\dfrac{\left(  -1\right)  ^{n-1}}{n}cx^{n}\right)  ^{m}$. Indeed, we can
replace the \textquotedblleft$\sum_{n\geq1}$\textquotedblright\ sign by an
\textquotedblleft$\sum_{n\in\mathbb{P}}$\textquotedblright\ sign, since
$\mathbb{P}=\left\{  1,2,3,\ldots\right\}  $. Thus,%
\begin{align*}
&  \left(  \sum_{n\geq1}\dfrac{\left(  -1\right)  ^{n-1}}{n}cx^{n}\right)
^{m}\\
&  =\left(  \sum_{n\in\mathbb{P}}\dfrac{\left(  -1\right)  ^{n-1}}{n}%
cx^{n}\right)  ^{m}\\
&  =\underbrace{\left(  \sum_{n\in\mathbb{P}}\dfrac{\left(  -1\right)  ^{n-1}%
}{n}cx^{n}\right)  \left(  \sum_{n\in\mathbb{P}}\dfrac{\left(  -1\right)
^{n-1}}{n}cx^{n}\right)  \cdots\left(  \sum_{n\in\mathbb{P}}\dfrac{\left(
-1\right)  ^{n-1}}{n}cx^{n}\right)  }_{m\text{ times}}\\
&  =\left(  \sum_{n_{1}\in\mathbb{P}}\dfrac{\left(  -1\right)  ^{n_{1}-1}%
}{n_{1}}cx^{n_{1}}\right)  \left(  \sum_{n_{2}\in\mathbb{P}}\dfrac{\left(
-1\right)  ^{n_{2}-1}}{n_{2}}cx^{n_{2}}\right)  \cdots\left(  \sum_{n_{m}%
\in\mathbb{P}}\dfrac{\left(  -1\right)  ^{n_{m}-1}}{n_{m}}cx^{n_{m}}\right) \\
&  \ \ \ \ \ \ \ \ \ \ \ \ \ \ \ \ \ \ \ \ \left(  \text{here, we have renamed
the summation indices}\right) \\
&  =\sum_{\left(  n_{1},n_{2},\ldots,n_{m}\right)  \in\mathbb{P}^{m}}\left(
\dfrac{\left(  -1\right)  ^{n_{1}-1}}{n_{1}}cx^{n_{1}}\right)  \left(
\dfrac{\left(  -1\right)  ^{n_{2}-1}}{n_{2}}cx^{n_{2}}\right)  \cdots\left(
\dfrac{\left(  -1\right)  ^{n_{m}-1}}{n_{m}}cx^{n_{m}}\right)
\end{align*}
(by a product rule for the product of $m$ sums\footnote{This product rule says
that
\[
\left(  \sum_{n_{1}\in A_{1}}a_{1,n_{1}}\right)  \left(  \sum_{n_{2}\in A_{2}%
}a_{2,n_{2}}\right)  \cdots\left(  \sum_{n_{m}\in A_{m}}a_{m,n_{m}}\right)
=\sum_{\left(  n_{1},n_{2},\ldots,n_{m}\right)  \in A_{1}\times A_{2}%
\times\cdots\times A_{m}}a_{1,n_{1}}a_{2,n_{2}}\cdots a_{m,n_{m}}%
\]
for any $m$ sets $A_{1},A_{2},\ldots,A_{m}$ and any elements $a_{i,j}\in
K\left[  \left[  x\right]  \right]  $, provided that all the sums on the left
hand side of this equality are summable. We leave it to the reader to convince
himself of this rule (intuitively, it just says that we can expand a product
of sums in the usual way, even when the sums are infinite) and to check that
the sums we are applying it to are indeed summable.}). Hence,%
\begin{align}
&  \left(  \sum_{n\geq1}\dfrac{\left(  -1\right)  ^{n-1}}{n}cx^{n}\right)
^{m}\nonumber\\
&  =\underbrace{\sum_{\left(  n_{1},n_{2},\ldots,n_{m}\right)  \in
\mathbb{P}^{m}}}_{=\sum_{k\in\mathbb{N}}\ \ \sum_{\substack{\left(
n_{1},n_{2},\ldots,n_{m}\right)  \in\mathbb{P}^{m};\\n_{1}+n_{2}+\cdots
+n_{m}=k}}}\underbrace{\left(  \dfrac{\left(  -1\right)  ^{n_{1}-1}}{n_{1}%
}cx^{n_{1}}\right)  \left(  \dfrac{\left(  -1\right)  ^{n_{2}-1}}{n_{2}%
}cx^{n_{2}}\right)  \cdots\left(  \dfrac{\left(  -1\right)  ^{n_{m}-1}}{n_{m}%
}cx^{n_{m}}\right)  }_{=\dfrac{\left(  -1\right)  ^{n_{1}+n_{2}+\cdots
+n_{m}-m}}{n_{1}n_{2}\cdots n_{m}}c^{m}x^{n_{1}+n_{2}+\cdots+n_{m}}%
}\nonumber\\
&  =\sum_{k\in\mathbb{N}}\ \ \sum_{\substack{\left(  n_{1},n_{2},\ldots
,n_{m}\right)  \in\mathbb{P}^{m};\\n_{1}+n_{2}+\cdots+n_{m}=k}}\dfrac{\left(
-1\right)  ^{n_{1}+n_{2}+\cdots+n_{m}-m}}{n_{1}n_{2}\cdots n_{m}}%
c^{m}\underbrace{x^{n_{1}+n_{2}+\cdots+n_{m}}}_{\substack{=x^{k}\\\text{(since
}n_{1}+n_{2}+\cdots+n_{m}=k\text{)}}}\nonumber\\
&  =\sum_{k\in\mathbb{N}}\ \ \sum_{\substack{\left(  n_{1},n_{2},\ldots
,n_{m}\right)  \in\mathbb{P}^{m};\\n_{1}+n_{2}+\cdots+n_{m}=k}}\dfrac{\left(
-1\right)  ^{n_{1}+n_{2}+\cdots+n_{m}-m}}{n_{1}n_{2}\cdots n_{m}}c^{m}x^{k}.
\label{pf.thm.fps.gen-newton.2}%
\end{align}

Now, forget that we fixed $m$. We thus have proved
(\ref{pf.thm.fps.gen-newton.2}) for each $m\in\mathbb{N}$.

Now, (\ref{pf.thm.fps.gen-newton.1}) becomes%
\begin{align}
&  \left(  1+x\right)  ^{c}\nonumber\\
&  =\sum_{m\in\mathbb{N}}\dfrac{1}{m!}\left(  \sum_{n\geq1}\dfrac{\left(
-1\right)  ^{n-1}}{n}cx^{n}\right)  ^{m}\nonumber\\
&  =\sum_{m\in\mathbb{N}}\dfrac{1}{m!}\sum_{k\in\mathbb{N}}\ \ \sum
_{\substack{\left(  n_{1},n_{2},\ldots,n_{m}\right)  \in\mathbb{P}^{m}%
;\\n_{1}+n_{2}+\cdots+n_{m}=k}}\dfrac{\left(  -1\right)  ^{n_{1}+n_{2}%
+\cdots+n_{m}-m}}{n_{1}n_{2}\cdots n_{m}}c^{m}x^{k}\ \ \ \ \ \ \ \ \ \ \left(
\text{by (\ref{pf.thm.fps.gen-newton.2})}\right) \nonumber\\
&  =\sum_{k\in\mathbb{N}}\ \ \sum_{m\in\mathbb{N}}\ \ \sum_{\substack{\left(
n_{1},n_{2},\ldots,n_{m}\right)  \in\mathbb{P}^{m};\\n_{1}+n_{2}+\cdots
+n_{m}=k}}\dfrac{1}{m!}\cdot\dfrac{\left(  -1\right)  ^{n_{1}+n_{2}%
+\cdots+n_{m}-m}}{n_{1}n_{2}\cdots n_{m}}c^{m}x^{k}\nonumber\\
&  =\sum_{k\in\mathbb{N}}\left(  \sum_{m\in\mathbb{N}}\ \ \sum
_{\substack{\left(  n_{1},n_{2},\ldots,n_{m}\right)  \in\mathbb{P}^{m}%
;\\n_{1}+n_{2}+\cdots+n_{m}=k}}\dfrac{1}{m!}\cdot\dfrac{\left(  -1\right)
^{n_{1}+n_{2}+\cdots+n_{m}-m}}{n_{1}n_{2}\cdots n_{m}}c^{m}\right)  x^{k}.
\label{pf.thm.fps.gen-newton.3}%
\end{align}

Now, let $k\in\mathbb{N}$. Let us rewrite the \textquotedblleft middle
sum\textquotedblright\ $\sum_{m\in\mathbb{N}}\ \ \sum_{\substack{\left(
n_{1},n_{2},\ldots,n_{m}\right)  \in\mathbb{P}^{m};\\n_{1}+n_{2}+\cdots
+n_{m}=k}}\dfrac{1}{m!}\cdot\dfrac{\left(  -1\right)  ^{n_{1}+n_{2}%
+\cdots+n_{m}-m}}{n_{1}n_{2}\cdots n_{m}}c^{m}$ on the right hand side as a
finite sum. Indeed, a \emph{composition} of $k$ shall mean a tuple $\left(
n_{1},n_{2},\ldots,n_{m}\right)  $ of positive integers satisfying
$n_{1}+n_{2}+\cdots+n_{m}=k$. (For example, $\left(  1,3,1\right)  $ is a
composition of $5$. We will study compositions in more detail in Section
\ref{sec.fps.intcomps}.) Let $\operatorname*{Comp}\left(  k\right)  $ denote
the set of all compositions of $k$. It is easy to see that this set
$\operatorname*{Comp}\left(  k\right)  $ is finite\footnote{\textit{Proof.}
Let $\left(  n_{1},n_{2},\ldots,n_{m}\right)  \in\operatorname*{Comp}\left(
k\right)  $. Thus, $\left(  n_{1},n_{2},\ldots,n_{m}\right)  $ is a
composition of $k$. In other words, $\left(  n_{1},n_{2},\ldots,n_{m}\right)
$ is a finite tuple of positive integers satisfying $n_{1}+n_{2}+\cdots
+n_{m}=k$. Hence, all its $m$ entries $n_{1},n_{2},\ldots,n_{m}$ are positive
integers and thus are $\geq1$; therefore, $n_{1}+n_{2}+\cdots+n_{m}%
\geq\underbrace{1+1+\cdots+1}_{m\text{ times}}=m$, so that $m\leq n_{1}%
+n_{2}+\cdots+n_{m}=k$. Thus, $m\in\left\{  0,1,\ldots,k\right\}  $.
\par
Furthermore, the sum $n_{1}+n_{2}+\cdots+n_{m}$ is $\geq$ to each of its $m$
addends (since its $m$ addends $n_{1},n_{2},\ldots,n_{m}$ are positive). In
other words, we have $n_{1}+n_{2}+\cdots+n_{m}\geq n_{i}$ for each
$i\in\left\{  1,2,\ldots,m\right\}  $. Thus, for each $i\in\left\{
1,2,\ldots,m\right\}  $, we have $n_{i}\leq n_{1}+n_{2}+\cdots+n_{m}=k$ and
therefore $n_{i}\in\left\{  1,2,\ldots,k\right\}  $ (since $n_{i}$ is a
positive integer). Hence,
\[
\left(  n_{1},n_{2},\ldots,n_{m}\right)  \in\left\{  1,2,\ldots,k\right\}
^{m}\subseteq\bigcup_{\ell\in\left\{  0,1,\ldots,k\right\}  }\left\{
1,2,\ldots,k\right\}  ^{\ell}%
\]
(since $m\in\left\{  0,1,\ldots,k\right\}  $).
\par
Now, forget that we fixed $\left(  n_{1},n_{2},\ldots,n_{m}\right)  $. We thus
have shown that $\left(  n_{1},n_{2},\ldots,n_{m}\right)  \in\bigcup_{\ell
\in\left\{  0,1,\ldots,k\right\}  }\left\{  1,2,\ldots,k\right\}  ^{\ell}$ for
each $\left(  n_{1},n_{2},\ldots,n_{m}\right)  \in\operatorname*{Comp}\left(
k\right)  $. In other words, $\operatorname*{Comp}\left(  k\right)
\subseteq\bigcup_{\ell\in\left\{  0,1,\ldots,k\right\}  }\left\{
1,2,\ldots,k\right\}  ^{\ell}$. Since the set $\bigcup_{\ell\in\left\{
0,1,\ldots,k\right\}  }\left\{  1,2,\ldots,k\right\}  ^{\ell}$ is clearly
finite (having size $\sum_{\ell\in\left\{  0,1,\ldots,k\right\}  }k^{\ell}$),
this entails that the set $\operatorname*{Comp}\left(  k\right)  $ is finite
as well, qed.
\par
(Incidentally, we will see in Section \ref{sec.fps.intcomps} that this set
$\operatorname*{Comp}\left(  k\right)  $ has size $2^{k-1}$ for $k\geq1$, and
size $1$ for $k=0$.)}. Now, we can rewrite the double summation sign
\textquotedblleft$\sum_{m\in\mathbb{N}}\ \ \sum_{\substack{\left(  n_{1}%
,n_{2},\ldots,n_{m}\right)  \in\mathbb{P}^{m};\\n_{1}+n_{2}+\cdots+n_{m}=k}%
}$\textquotedblright\ as a single summation sign \textquotedblleft%
$\sum_{\left(  n_{1},n_{2},\ldots,n_{m}\right)  \in\operatorname*{Comp}\left(
k\right)  }$\textquotedblright\ (since $\operatorname*{Comp}\left(  k\right)
$ is precisely the set of all tuples $\left(  n_{1},n_{2},\ldots,n_{m}\right)
\in\mathbb{P}^{m}$ satisfying $n_{1}+n_{2}+\cdots+n_{m}=k$). Hence, we obtain
\begin{align}
&  \sum_{m\in\mathbb{N}}\ \ \sum_{\substack{\left(  n_{1},n_{2},\ldots
,n_{m}\right)  \in\mathbb{P}^{m};\\n_{1}+n_{2}+\cdots+n_{m}=k}}\dfrac{1}%
{m!}\cdot\dfrac{\left(  -1\right)  ^{n_{1}+n_{2}+\cdots+n_{m}-m}}{n_{1}%
n_{2}\cdots n_{m}}c^{m}\nonumber\\
&  =\sum_{\left(  n_{1},n_{2},\ldots,n_{m}\right)  \in\operatorname*{Comp}%
\left(  k\right)  }\dfrac{1}{m!}\cdot\dfrac{\left(  -1\right)  ^{n_{1}%
+n_{2}+\cdots+n_{m}-m}}{n_{1}n_{2}\cdots n_{m}}c^{m}.
\label{pf.thm.fps.gen-newton.4}%
\end{align}

Forget that we fixed $k$. Thus, for each $k\in\mathbb{N}$, we have defined a
finite set $\operatorname*{Comp}\left(  k\right)  $ and shown that
(\ref{pf.thm.fps.gen-newton.4}) holds.

Using (\ref{pf.thm.fps.gen-newton.4}), we can rewrite
(\ref{pf.thm.fps.gen-newton.3}) as%
\begin{align}
&  \left(  1+x\right)  ^{c}\nonumber\\
&  =\sum_{k\in\mathbb{N}}\left(  \sum_{\left(  n_{1},n_{2},\ldots
,n_{m}\right)  \in\operatorname*{Comp}\left(  k\right)  }\dfrac{1}{m!}%
\cdot\dfrac{\left(  -1\right)  ^{n_{1}+n_{2}+\cdots+n_{m}-m}}{n_{1}n_{2}\cdots
n_{m}}c^{m}\right)  x^{k}. \label{pf.thm.fps.gen-newton.5}%
\end{align}

Now, recall that our goal is to prove that this equals
\[
\sum_{k\in\mathbb{N}}\dbinom{c}{k}x^{k}.
\]
This is equivalent to proving that the equality%
\begin{equation}
\sum_{\left(  n_{1},n_{2},\ldots,n_{m}\right)  \in\operatorname*{Comp}\left(
k\right)  }\dfrac{1}{m!}\cdot\dfrac{\left(  -1\right)  ^{n_{1}+n_{2}%
+\cdots+n_{m}-m}}{n_{1}n_{2}\cdots n_{m}}c^{m}=\dbinom{c}{k}
\label{pf.thm.fps.gen-newton.polyeq}%
\end{equation}
holds for each $k\in\mathbb{N}$ (because two FPSs $u=\sum_{k\in\mathbb{N}%
}u_{k}x^{k}$ and $v=\sum_{k\in\mathbb{N}}v_{k}x^{k}$ (with $u_{k}\in K$ and
$v_{k}\in K$) are equal if and only if the equality $u_{k}=v_{k}$ holds for
each $k\in\mathbb{N}$).

Thus, we have reduced our original goal (which was to prove $\left(
1+x\right)  ^{c}=\sum_{k\in\mathbb{N}}\dbinom{c}{k}x^{k}$) to the auxiliary
goal of proving the equality (\ref{pf.thm.fps.gen-newton.polyeq}) for each
$k\in\mathbb{N}$. However, this doesn't look very useful, since
(\ref{pf.thm.fps.gen-newton.polyeq}) is too messy an equality to have a simple
proof. We are seemingly stuck.

However, it turns out that we are almost there -- we just need to take a
bird's eye view. Here is the plan: We fix $k\in\mathbb{N}$. Instead of trying
to prove the equality (\ref{pf.thm.fps.gen-newton.polyeq}) directly, we
observe that both sides of this equality are polynomials (with rational
coefficients) in $c$. (Indeed, the left hand side is clearly a polynomial in
$c$, since it is a finite sum of \textquotedblleft rational number times a
power of $c$\textquotedblright\ expressions. The right hand side is a
polynomial in $c$ because $\dbinom{c}{k}=\dfrac{c\left(  c-1\right)  \left(
c-2\right)  \cdots\left(  c-k+1\right)  }{k!}$.) Thus, the polynomial identity
trick (which we learnt in Subsection \ref{subsec.gf.defs.cvi}) tells us that
if we can prove this equality (\ref{pf.thm.fps.gen-newton.polyeq}) for each
$c\in\mathbb{N}$, then it will automatically hold for each $c\in K$ (since the
two polynomials that yield its left and right hand sides will have to be
equal, having infinitely many equal values). Hence, in order to prove
(\ref{pf.thm.fps.gen-newton.polyeq}) for each $c\in K$, it suffices to prove
it for each $c\in\mathbb{N}$. Now, how can we prove it for each $c\in
\mathbb{N}$ ? We forget that we fixed $k$, and we remember that the equality
(\ref{pf.thm.fps.gen-newton.polyeq}) (for all $k\in\mathbb{N}$) is just an
equivalent restatement of the FPS equality $\left(  1+x\right)  ^{c}%
=\sum_{k\in\mathbb{N}}\dbinom{c}{k}x^{k}$ (that is, the equality we have
originally set out to prove). However, we know for sure that this equality
holds for each $c\in\mathbb{N}$ (by Theorem \ref{thm.fps.newton-binom},
applied to $n=c$). Hence, the equality (\ref{pf.thm.fps.gen-newton.polyeq})
also holds for each $c\in\mathbb{N}$ (and each $k\in\mathbb{N}$). And this is
precisely what we needed to show!

Let me explain this argument in detail now, as it is somewhat
vertigo-inducing. We forget that we fixed $K$ and $c$. Now, fix $c\in
\mathbb{N}$. Thus, $c\in\mathbb{N}\subseteq\mathbb{Z}$. Hence, in the ring
$\mathbb{Q}\left[  \left[  x\right]  \right]  $, we have%
\[
\left(  1+x\right)  ^{c}=\sum_{k\in\mathbb{N}}\dbinom{c}{k}x^{k}%
\ \ \ \ \ \ \ \ \ \ \left(  \text{by Theorem \ref{thm.fps.newton-binom},
applied to }n=c\right)  .
\]
However, (\ref{pf.thm.fps.gen-newton.5}) (applied to $K=\mathbb{Q}$) shows
that%
\[
\left(  1+x\right)  ^{c}=\sum_{k\in\mathbb{N}}\left(  \sum_{\left(
n_{1},n_{2},\ldots,n_{m}\right)  \in\operatorname*{Comp}\left(  k\right)
}\dfrac{1}{m!}\cdot\dfrac{\left(  -1\right)  ^{n_{1}+n_{2}+\cdots+n_{m}-m}%
}{n_{1}n_{2}\cdots n_{m}}c^{m}\right)  x^{k}.
\]
Comparing these two equalities, we obtain%
\[
\sum_{k\in\mathbb{N}}\left(  \sum_{\left(  n_{1},n_{2},\ldots,n_{m}\right)
\in\operatorname*{Comp}\left(  k\right)  }\dfrac{1}{m!}\cdot\dfrac{\left(
-1\right)  ^{n_{1}+n_{2}+\cdots+n_{m}-m}}{n_{1}n_{2}\cdots n_{m}}c^{m}\right)
x^{k}=\sum_{k\in\mathbb{N}}\dbinom{c}{k}x^{k}.
\]
Comparing coefficients in this equality, we see that%
\begin{equation}
\sum_{\left(  n_{1},n_{2},\ldots,n_{m}\right)  \in\operatorname*{Comp}\left(
k\right)  }\dfrac{1}{m!}\cdot\dfrac{\left(  -1\right)  ^{n_{1}+n_{2}%
+\cdots+n_{m}-m}}{n_{1}n_{2}\cdots n_{m}}c^{m}=\dbinom{c}{k}
\label{pf.thm.fps.gen-newton.finale.1}%
\end{equation}
for each $k\in\mathbb{N}$. This is an equality between two rational numbers.

Now, forget that we fixed $c$. We thus have shown that
(\ref{pf.thm.fps.gen-newton.finale.1}) holds for each $k\in\mathbb{N}$ and
each $c\in\mathbb{N}$.

Let us now fix $k\in\mathbb{N}$. We have just shown that the equality
(\ref{pf.thm.fps.gen-newton.finale.1}) holds for each $c\in\mathbb{N}$. In
other words, the two polynomials%
\[
f:=\sum_{\left(  n_{1},n_{2},\ldots,n_{m}\right)  \in\operatorname*{Comp}%
\left(  k\right)  }\dfrac{1}{m!}\cdot\dfrac{\left(  -1\right)  ^{n_{1}%
+n_{2}+\cdots+n_{m}-m}}{n_{1}n_{2}\cdots n_{m}}x^{m}\in\mathbb{Q}\left[
x\right]
\]
and%
\[
g:=\dbinom{x}{k}=\dfrac{x\left(  x-1\right)  \left(  x-2\right)  \cdots\left(
x-k+1\right)  }{k!}\in\mathbb{Q}\left[  x\right]
\]
satisfy $f\left[  c\right]  =g\left[  c\right]  $ for each $c\in\mathbb{N}$
(because $f\left[  c\right]  $ is the left hand side of
(\ref{pf.thm.fps.gen-newton.finale.1}), while $g\left[  c\right]  $ is the
right hand side of (\ref{pf.thm.fps.gen-newton.finale.1})). Thus, each
$c\in\mathbb{N}$ satisfies $\left(  f-g\right)  \left[  c\right]
=\underbrace{f\left[  c\right]  }_{=g\left[  c\right]  }-\,g\left[  c\right]
=g\left[  c\right]  -g\left[  c\right]  =0$. In other words, each
$c\in\mathbb{N}$ is a root of $f-g$. Hence, the polynomial $f-g$ has
infinitely many roots in $\mathbb{Q}$ (since there are infinitely many
$c\in\mathbb{N}$). Since $f-g$ is a polynomial with rational coefficients,
this is impossible unless $f-g=0$. We thus must have $f-g=0$, so that $f=g$.
In other words,
\begin{equation}
\sum_{\left(  n_{1},n_{2},\ldots,n_{m}\right)  \in\operatorname*{Comp}\left(
k\right)  }\dfrac{1}{m!}\cdot\dfrac{\left(  -1\right)  ^{n_{1}+n_{2}%
+\cdots+n_{m}-m}}{n_{1}n_{2}\cdots n_{m}}x^{m}=\dbinom{x}{k}
\label{pf.thm.fps.gen-newton.finale.4}%
\end{equation}
holds in the polynomial ring $\mathbb{Q}\left[  x\right]  $.

Now, forget that we fixed $k$. We thus have shown that the equality
(\ref{pf.thm.fps.gen-newton.finale.4}) holds for each $k\in\mathbb{N}$.

Now, fix a commutative $\mathbb{Q}$-algebra $K$ and an arbitrary element $c\in
K$. For each $k\in\mathbb{N}$, we then have%
\begin{equation}
\sum_{\left(  n_{1},n_{2},\ldots,n_{m}\right)  \in\operatorname*{Comp}\left(
k\right)  }\dfrac{1}{m!}\cdot\dfrac{\left(  -1\right)  ^{n_{1}+n_{2}%
+\cdots+n_{m}-m}}{n_{1}n_{2}\cdots n_{m}}c^{m}=\dbinom{c}{k}
\label{pf.thm.fps.gen-newton.finale.5}%
\end{equation}
(by substituting $c$ for $x$ on both sides of the equality
(\ref{pf.thm.fps.gen-newton.finale.4})). Consequently, we can rewrite
(\ref{pf.thm.fps.gen-newton.5}) as%
\[
\left(  1+x\right)  ^{c}=\sum_{k\in\mathbb{N}}\dbinom{c}{k}x^{k}.
\]
This proves Theorem \ref{thm.fps.gen-newton}.
\end{proof}

The method we used in the above proof is worth recapitulating in broad strokes:

\begin{itemize}
\item We had to prove a fairly abstract statement (namely, the identity
$\left(  1+x\right)  ^{c}=\sum_{k\in\mathbb{N}}\dbinom{c}{k}x^{k}$).

\item We translated this statement into an awkward but more concrete statement
(namely, the equality (\ref{pf.thm.fps.gen-newton.polyeq})).

\item We then argued that this concrete statement needs only to be proven in a
special case (viz., for all $c\in\mathbb{N}$ rather than for all $c\in K$),
because it is an equality between two polynomials with rational coefficients.

\item To prove this concrete statement in this special case, we translated it
back into the abstract language of FPSs, and realized that in this special
case it is already known (as a consequence of Theorem
\ref{thm.fps.newton-binom}).
\end{itemize}

Thus, by strategically switching between the abstract and the concrete, we
have managed to use the advantages of both sides.

Now that Theorem \ref{thm.fps.gen-newton} is proved, Example 2 from Section
\ref{sec.gf.exas} is fully justified (since we can obtain
(\ref{eq.sec.gf.exas.2.(1+x)1/2}) by applying Theorem \ref{thm.fps.gen-newton}
to $K=\mathbb{Q}$ and $c=1/2$).

\subsubsection{Another application}

Let us show yet another application of powers with non-integer exponents and
the generalized Newton formula. We shall show the following binomial identity:

\begin{proposition}
\label{prop.binom.nCk-2i-qedmo.CN}Let $n\in\mathbb{C}$ and $k\in\mathbb{N}$.
Then,%
\[
\sum_{i=0}^{k}\dbinom{n+i-1}{i}\dbinom{n}{k-2i}=\dbinom{n+k-1}{k}.
\]

\end{proposition}

Proposition \ref{prop.binom.nCk-2i-qedmo.CN} can be proved in various ways.
For example, a mostly combinatorial proof is found in \cite[Exercise 2.10.7
and Exercise 2.10.8]{19fco}\footnote{Specifically, \cite[Exercise
2.10.7]{19fco} proves Proposition \ref{prop.binom.nCk-2i-qedmo.CN} in the
particular case when $n\in\mathbb{N}$; then, \cite[Exercise 2.10.8]{19fco}
extends it to the case when $n\in\mathbb{R}$. However, the latter argument can
just as well be used to extend it to arbitrary $n\in\mathbb{C}$.}. We shall
give a proof using generating functions instead.

\begin{proof}
[Proof of Proposition \ref{prop.binom.nCk-2i-qedmo.CN}.]Define two FPSs
$f,g\in\mathbb{C}\left[  \left[  x\right]  \right]  $ by
\begin{align}
f  &  =\sum_{i\in\mathbb{N}}\dbinom{n+i-1}{i}x^{2i}%
\label{pf.prop.binom.nCk-2i-qedmo.CN.f=}\\
\text{and}\ \ \ \ \ \ \ \ \ \ g  &  =\sum_{j\in\mathbb{N}}\dbinom{n}{j}x^{j}.
\label{pf.prop.binom.nCk-2i-qedmo.CN.g=}%
\end{align}
(We will soon see why we chose to define them this way.) Multiplying the two
equalities (\ref{pf.prop.binom.nCk-2i-qedmo.CN.f=}) and
(\ref{pf.prop.binom.nCk-2i-qedmo.CN.g=}), we find%
\begin{align*}
fg  &  =\left(  \sum_{i\in\mathbb{N}}\dbinom{n+i-1}{i}x^{2i}\right)  \left(
\sum_{j\in\mathbb{N}}\dbinom{n}{j}x^{j}\right) \\
&  =\underbrace{\sum_{i\in\mathbb{N}}\ \ \sum_{j\in\mathbb{N}}}_{=\sum
_{\left(  i,j\right)  \in\mathbb{N}\times\mathbb{N}}}\ \ \underbrace{\dbinom
{n+i-1}{i}x^{2i}\dbinom{n}{j}x^{j}}_{=\dbinom{n+i-1}{i}\dbinom{n}{j}x^{2i+j}%
}\\
&  =\underbrace{\sum_{\left(  i,j\right)  \in\mathbb{N}\times\mathbb{N}}%
}_{=\sum_{m\in\mathbb{N}}\ \ \sum_{\substack{\left(  i,j\right)  \in
\mathbb{N}\times\mathbb{N};\\2i+j=m}}}\dbinom{n+i-1}{i}\dbinom{n}{j}x^{2i+j}\\
&  =\sum_{m\in\mathbb{N}}\ \ \sum_{\substack{\left(  i,j\right)  \in
\mathbb{N}\times\mathbb{N};\\2i+j=m}}\dbinom{n+i-1}{i}\dbinom{n}%
{j}\underbrace{x^{2i+j}}_{\substack{=x^{m}\\\text{(since }2i+j=m\text{)}}}\\
&  =\sum_{m\in\mathbb{N}}\ \ \sum_{\substack{\left(  i,j\right)  \in
\mathbb{N}\times\mathbb{N};\\2i+j=m}}\dbinom{n+i-1}{i}\dbinom{n}{j}x^{m}\\
&  =\sum_{m\in\mathbb{N}}\left(  \sum_{\substack{\left(  i,j\right)
\in\mathbb{N}\times\mathbb{N};\\2i+j=m}}\dbinom{n+i-1}{i}\dbinom{n}{j}\right)
x^{m}.
\end{align*}
Hence, the $x^{k}$-coefficient of this FPS $fg$ is%
\begin{equation}
\left[  x^{k}\right]  \left(  fg\right)  =\sum_{\substack{\left(  i,j\right)
\in\mathbb{N}\times\mathbb{N};\\2i+j=k}}\dbinom{n+i-1}{i}\dbinom{n}{j}.
\label{pf.prop.binom.nCk-2i-qedmo.CN.3}%
\end{equation}

Now, a pair $\left(  i,j\right)  \in\mathbb{N}\times\mathbb{N}$ satisfying
$2i+j=k$ is uniquely determined by its first entry $i$, since its second entry
$j$ is given by $j=k-2i$. Hence, we can substitute $\left(  i,k-2i\right)  $
for $\left(  i,j\right)  $ in the sum $\sum_{\substack{\left(  i,j\right)
\in\mathbb{N}\times\mathbb{N};\\2i+j=k}}\dbinom{n+i-1}{i}\dbinom{n}{j}$, thus
rewriting this sum as $\sum_{\substack{i\in\mathbb{N};\\k-2i\in\mathbb{N}%
}}\dbinom{n+i-1}{i}\dbinom{n}{k-2i}$. Hence,%
\begin{align*}
\sum_{\substack{\left(  i,j\right)  \in\mathbb{N}\times\mathbb{N}%
;\\2i+j=k}}\dbinom{n+i-1}{i}\dbinom{n}{j}  &  =\sum_{\substack{i\in
\mathbb{N};\\k-2i\in\mathbb{N}}}\dbinom{n+i-1}{i}\dbinom{n}{k-2i}\\
&  =\sum_{\substack{i\in\mathbb{N};\\2i\leq k}}\dbinom{n+i-1}{i}\dbinom
{n}{k-2i}\\
&  \ \ \ \ \ \ \ \ \ \ \ \ \ \ \ \ \ \ \ \ \left(
\begin{array}
[c]{c}%
\text{since an }i\in\mathbb{N}\text{ satisfies }k-2i\in\mathbb{N}\\
\text{if and only if it satisfies }2i\leq k
\end{array}
\right) \\
&  =\sum_{\substack{i\in\mathbb{N};\\i\leq k}}\dbinom{n+i-1}{i}\dbinom
{n}{k-2i}%
\end{align*}
(here, we have replaced the condition \textquotedblleft$2i\leq k$%
\textquotedblright\ under the summation sign by the weaker condition
\textquotedblleft$i\leq k$\textquotedblright, thus extending the range of the
sum; but this did not change the sum, since all the newly introduced addends
are $0$ because of the vanishing $\dbinom{n}{k-2i}$ factor). Thus,
(\ref{pf.prop.binom.nCk-2i-qedmo.CN.3}) becomes%
\begin{align}
\left[  x^{k}\right]  \left(  fg\right)   &  =\sum_{\substack{\left(
i,j\right)  \in\mathbb{N}\times\mathbb{N};\\2i+j=k}}\dbinom{n+i-1}{i}%
\dbinom{n}{j}\nonumber\\
&  =\sum_{\substack{i\in\mathbb{N};\\i\leq k}}\dbinom{n+i-1}{i}\dbinom
{n}{k-2i}\nonumber\\
&  =\sum_{i=0}^{k}\dbinom{n+i-1}{i}\dbinom{n}{k-2i}.
\label{pf.prop.binom.nCk-2i-qedmo.CN.5}%
\end{align}
Note that the right hand side here is precisely the left hand side of the
identity we are trying to prove. This is why we defined $f$ and $g$ as we did.
With a bit of experience, the computation above can easily be
reverse-engineered, and the definitions of $f$ and $g$ are essentially forced
by the goal of making (\ref{pf.prop.binom.nCk-2i-qedmo.CN.5}) hold.

Anyway, it is now clear that a simple expression for $fg$ would move us
forward. So let us try to simplify $f$ and $g$. For $g$, the answer is
easiest: We have%
\[
g=\sum_{j\in\mathbb{N}}\dbinom{n}{j}x^{j}=\left(  1+x\right)  ^{n},
\]
because Theorem \ref{thm.fps.gen-newton} (applied to $c=n$) yields $\left(
1+x\right)  ^{n}=\sum_{j\in\mathbb{N}}\dbinom{n}{j}x^{j}$. For $f$, we need a
few more steps. Proposition \ref{prop.fps.anti-newton-binom} yields%
\begin{equation}
\left(  1+x\right)  ^{-n}=\sum_{i\in\mathbb{N}}\left(  -1\right)  ^{i}%
\dbinom{n+i-1}{i}x^{i}. \label{pf.prop.binom.nCk-2i-qedmo.CN.9}%
\end{equation}
Substituting $-x^{2}$ for $x$ on both sides of this equality, we obtain%
\begin{align*}
\left(  1-x^{2}\right)  ^{-n}  &  =\sum_{i\in\mathbb{N}}\left(  -1\right)
^{i}\dbinom{n+i-1}{i}\underbrace{\left(  -x^{2}\right)  ^{i}}_{=\left(
-1\right)  ^{i}x^{2i}}=\sum_{i\in\mathbb{N}}\left(  -1\right)  ^{i}%
\dbinom{n+i-1}{i}\left(  -1\right)  ^{i}x^{2i}\\
&  =\sum_{i\in\mathbb{N}}\underbrace{\left(  -1\right)  ^{i}\left(  -1\right)
^{i}}_{=1}\dbinom{n+i-1}{i}x^{2i}=\sum_{i\in\mathbb{N}}\dbinom{n+i-1}{i}%
x^{2i}=f
\end{align*}
and thus $f=\left(  1-x^{2}\right)  ^{-n}$. Multiplying this equality by
$g=\left(  1+x\right)  ^{n}$, we obtain%
\begin{align*}
fg  &  =\left(  1-x^{2}\right)  ^{-n}\left(  1+x\right)  ^{n}=\left(
\dfrac{1+x}{1-x^{2}}\right)  ^{n}=\left(  \dfrac{1-x^{2}}{1+x}\right)  ^{-n}\\
&  =\left(  1-x\right)  ^{-n}\ \ \ \ \ \ \ \ \ \ \left(  \text{since }%
\dfrac{1-x^{2}}{1+x}=\dfrac{\left(  1-x\right)  \left(  1+x\right)  }%
{1+x}=1-x\right) \\
&  =\sum_{i\in\mathbb{N}}\left(  -1\right)  ^{i}\dbinom{n+i-1}{i}%
\underbrace{\left(  -x\right)  ^{i}}_{=\left(  -1\right)  ^{i}x^{i}}\\
&  \ \ \ \ \ \ \ \ \ \ \ \ \ \ \ \ \ \ \ \ \left(  \text{this follows by
substituting }-x\text{ for }x\text{ on both sides of
(\ref{pf.prop.binom.nCk-2i-qedmo.CN.9})}\right) \\
&  =\sum_{i\in\mathbb{N}}\left(  -1\right)  ^{i}\dbinom{n+i-1}{i}\left(
-1\right)  ^{i}x^{i}=\sum_{i\in\mathbb{N}}\underbrace{\left(  -1\right)
^{i}\left(  -1\right)  ^{i}}_{=1}\dbinom{n+i-1}{i}x^{i}\\
&  =\sum_{i\in\mathbb{N}}\dbinom{n+i-1}{i}x^{i}.
\end{align*}
Hence, the $x^{k}$-coefficient in $fg$ is $\left[  x^{k}\right]  \left(
fg\right)  =\dbinom{n+k-1}{k}$. Comparing this with
(\ref{pf.prop.binom.nCk-2i-qedmo.CN.5}), we obtain%
\[
\sum_{i=0}^{k}\dbinom{n+i-1}{i}\dbinom{n}{k-2i}=\dbinom{n+k-1}{k}.
\]
Thus, Proposition \ref{prop.binom.nCk-2i-qedmo.CN} is proved.
\end{proof}

\subsection{\label{sec.fps.intcomps}Integer compositions}

\subsubsection{Compositions}

Next, let us count certain simple combinatorial objects known as \emph{integer
compositions}. There are easy combinatorial ways to do this (see, e.g.,
\cite[\S 2.10.1]{19fco}), but we shall employ generating functions, in order
to see one more example of how these can be used.

\begin{definition}
\label{def.fps.comps}\textbf{(a)} An \emph{(integer) composition} means a
(finite) tuple of positive integers. \medskip

\textbf{(b)} The \emph{size} of an integer composition $\alpha=\left(
\alpha_{1},\alpha_{2},\ldots,\alpha_{m}\right)  $ is defined to be the integer
$\alpha_{1}+\alpha_{2}+\cdots+\alpha_{m}$. It is written $\left\vert
\alpha\right\vert $. \medskip

\textbf{(c)} The \emph{length} of an integer composition $\alpha=\left(
\alpha_{1},\alpha_{2},\ldots,\alpha_{m}\right)  $ is defined to be the integer
$m$. It is written $\ell\left(  \alpha\right)  $. \medskip

\textbf{(d)} Let $n\in\mathbb{N}$. A \emph{composition of }$n$ means a
composition whose size is $n$. \medskip

\textbf{(e)} Let $n\in\mathbb{N}$ and $k\in\mathbb{N}$. A \emph{composition of
}$n$\emph{ into }$k$\emph{ parts} is a composition whose size is $n$ and whose
length is $k$.
\end{definition}

\begin{example}
\label{exa.fps.comps.1}The tuple $\left(  3,8,6\right)  $ is a composition
with size $3+8+6=17$ and length $3$. In other words, it is a composition of
$17$ into $3$ parts.

The empty tuple $\left(  {}\right)  $ is a composition of $0$ into $0$ parts.
It is the only composition of $0$, and also is the only composition with
length $0$.
\end{example}

The following \textbf{questions} arise quite naturally:

\begin{enumerate}
\item How many compositions of $n$ exist for a given $n\in\mathbb{N}$ ?

\item How many compositions of $n$ into $k$ parts exist for given
$n,k\in\mathbb{N}$ ?
\end{enumerate}

Let us use generating functions to answer question 2.

\begin{proof}
[Approach to question 2.]Fix $k\in\mathbb{N}$, but don't fix $n$. Let
\begin{equation}
a_{n,k}=\left(  \text{\# of compositions of }n\text{ into }k\text{
parts}\right)  . \label{pf.thm.fps.comps.n-into-k-parts.1st.ank=}%
\end{equation}
We want to find $a_{n,k}$. We define the generating function%
\begin{equation}
A_{k}:=\sum_{n\in\mathbb{N}}a_{n,k}x^{n}=\left(  a_{0,k},a_{1,k}%
,a_{2,k},\ldots\right)  \in\mathbb{Q}\left[  \left[  x\right]  \right]  .
\label{pf.thm.fps.comps.n-into-k-parts.1st.Ak=}%
\end{equation}

Let us write $\mathbb{P}$ for the set $\left\{  1,2,3,\ldots\right\}  $. Then,
a composition of $n$ into $k$ parts is nothing but a $k$-tuple $\left(
\alpha_{1},\alpha_{2},\ldots,\alpha_{k}\right)  \in\mathbb{P}^{k}$ satisfying
$\alpha_{1}+\alpha_{2}+\cdots+\alpha_{k}=n$. Hence,
(\ref{pf.thm.fps.comps.n-into-k-parts.1st.ank=}) can be rewritten as%
\begin{align}
a_{n,k}  &  =\left(  \text{\# of all }k\text{-tuples }\left(  \alpha
_{1},\alpha_{2},\ldots,\alpha_{k}\right)  \in\mathbb{P}^{k}\text{ satisfying
}\alpha_{1}+\alpha_{2}+\cdots+\alpha_{k}=n\right) \nonumber\\
&  =\sum\limits_{\substack{\left(  \alpha_{1},\alpha_{2},\ldots,\alpha
_{k}\right)  \in\mathbb{P}^{k};\\\alpha_{1}+\alpha_{2}+\cdots+\alpha_{k}=n}}1.
\label{pf.thm.fps.comps.n-into-k-parts.1st.ank=2}%
\end{align}
Thus, we can rewrite the equality $A_{k}=\sum_{n\in\mathbb{N}}a_{n,k}x^{n}$
as
\begin{align*}
A_{k}  &  =\sum_{n\in\mathbb{N}}\left(  \sum\limits_{\substack{\left(
\alpha_{1},\alpha_{2},\ldots,\alpha_{k}\right)  \in\mathbb{P}^{k};\\\alpha
_{1}+\alpha_{2}+\cdots+\alpha_{k}=n}}1\right)  x^{n}=\sum_{n\in\mathbb{N}%
}\ \ \sum\limits_{\substack{\left(  \alpha_{1},\alpha_{2},\ldots,\alpha
_{k}\right)  \in\mathbb{P}^{k};\\\alpha_{1}+\alpha_{2}+\cdots+\alpha_{k}%
=n}}\ \ \underbrace{x^{n}}_{\substack{=x^{\alpha_{1}+\alpha_{2}+\cdots
+\alpha_{k}}\\\text{(since }\alpha_{1}+\alpha_{2}+\cdots+\alpha_{k}=n\text{)}%
}}\\
&  =\underbrace{\sum_{n\in\mathbb{N}}\ \ \sum\limits_{\substack{\left(
\alpha_{1},\alpha_{2},\ldots,\alpha_{k}\right)  \in\mathbb{P}^{k};\\\alpha
_{1}+\alpha_{2}+\cdots+\alpha_{k}=n}}}_{=\sum\limits_{\left(  \alpha
_{1},\alpha_{2},\ldots,\alpha_{k}\right)  \in\mathbb{P}^{k}}}%
\ \ \underbrace{x^{\alpha_{1}+\alpha_{2}+\cdots+\alpha_{k}}}_{=x^{\alpha_{1}%
}x^{\alpha_{2}}\cdots x^{\alpha_{k}}}=\sum\limits_{\left(  \alpha_{1}%
,\alpha_{2},\ldots,\alpha_{k}\right)  \in\mathbb{P}^{k}}x^{\alpha_{1}%
}x^{\alpha_{2}}\cdots x^{\alpha_{k}}\\
&  =\left(  \sum_{\alpha_{1}\in\mathbb{P}}x^{\alpha_{1}}\right)  \left(
\sum_{\alpha_{2}\in\mathbb{P}}x^{\alpha_{2}}\right)  \cdots\left(
\sum_{\alpha_{k}\in\mathbb{P}}x^{\alpha_{k}}\right) \\
&  \ \ \ \ \ \ \ \ \ \ \ \ \ \ \ \ \ \ \ \ \left(
\begin{array}
[c]{c}%
\text{by the same product rule that we used back}\\
\text{in the proof of Theorem \ref{thm.fps.gen-newton}}%
\end{array}
\right) \\
&  =\left(  \sum_{n\in\mathbb{P}}x^{n}\right)  ^{k}%
\end{align*}
(here, we have renamed all the $k$ summation indices as $n$, and realized that
all $k$ sums are identical). Since%
\[
\sum_{n\in\mathbb{P}}x^{n}=x^{1}+x^{2}+x^{3}+\cdots=x\underbrace{\left(
1+x+x^{2}+\cdots\right)  }_{=\dfrac{1}{1-x}}=x\cdot\dfrac{1}{1-x}=\dfrac
{x}{1-x},
\]
this can be rewritten further as%
\begin{equation}
A_{k}=\left(  \dfrac{x}{1-x}\right)  ^{k}=x^{k}\left(  1-x\right)  ^{-k}.
\label{pf.thm.fps.comps.num-comps-n-k.Ak=}%
\end{equation}

In order to simplify this, we need to expand $\left(  1-x\right)  ^{-k}$. This
is routine by now: Theorem \ref{thm.fps.newton-binom} (applied to $-k$ instead
of $n$) yields%
\[
\left(  1+x\right)  ^{-k}=\sum_{j\in\mathbb{N}}\dbinom{-k}{j}x^{j}.
\]
Substituting $-x$ for $x$ on both sides of this equality (i.e., applying the
map $f\mapsto f\circ\left(  -x\right)  $), we obtain%
\begin{align}
\left(  1-x\right)  ^{-k}  &  =\sum_{j\in\mathbb{N}}\underbrace{\dbinom{-k}%
{j}}_{\substack{=\left(  -1\right)  ^{j}\dbinom{j+k-1}{j}\\\text{(by Theorem
\ref{thm.binom.upneg-n},}\\\text{applied to }k\text{ and }j\\\text{instead of
}n\text{ and }k\text{)}}}\underbrace{\left(  -x\right)  ^{j}}_{=\left(
-1\right)  ^{j}x^{j}}=\sum_{j\in\mathbb{N}}\left(  -1\right)  ^{j}%
\dbinom{j+k-1}{j}\cdot\left(  -1\right)  ^{j}x^{j}\nonumber\\
&  =\sum_{j\in\mathbb{N}}\underbrace{\left(  \left(  -1\right)  ^{j}\right)
^{2}}_{=1}\dbinom{j+k-1}{j}x^{j}=\sum_{j\in\mathbb{N}}\dbinom{j+k-1}{j}%
x^{j}\label{pf.thm.fps.comps.num-comps-n-k.1-x-k}\\
&  =\sum\limits_{n\geq k}\dbinom{n-1}{n-k}x^{n-k}\nonumber
\end{align}
(here, we have substituted $n-k$ for $j$ in the sum). Hence, our above
computation of $A_{k}$ can be completed as follows:
\begin{align*}
A_{k}  &  =x^{k}\underbrace{\left(  1-x\right)  ^{-k}}_{=\sum\limits_{n\geq
k}\dbinom{n-1}{n-k}x^{n-k}}=x^{k}\sum\limits_{n\geq k}\dbinom{n-1}{n-k}%
x^{n-k}\\
&  =\sum\limits_{n\geq k}\dbinom{n-1}{n-k}\underbrace{x^{k}x^{n-k}}_{=x^{n}%
}=\sum\limits_{n\geq k}\dbinom{n-1}{n-k}x^{n}=\sum\limits_{n\in\mathbb{N}%
}\dbinom{n-1}{n-k}x^{n}%
\end{align*}
(here, we have extended the range of the summation from all $n\geq k$ to all
$n\in\mathbb{N}$; this did not change the sum, since all the newly introduced
addends with $n<k$ are $0$). Comparing coefficients, we thus obtain%
\begin{equation}
a_{n,k}=\dbinom{n-1}{n-k}\ \ \ \ \ \ \ \ \ \ \text{for each }n\in\mathbb{N}.
\label{pf.thm.fps.comps.n-into-k-parts.1st.rhs1}%
\end{equation}

If $n>0$, then we can rewrite the right hand side of this equality as
$\dbinom{n-1}{k-1}$ (using Theorem \ref{thm.binom.sym}). However, if $n=0$,
then this right hand side equals $\delta_{k,0}$ instead (where we are using
Definition \ref{def.kron-delta}). Thus, we can rewrite
(\ref{pf.thm.fps.comps.n-into-k-parts.1st.rhs1}) as%
\begin{equation}
a_{n,k}=%
\begin{cases}
\dbinom{n-1}{k-1}, & \text{if }n>0;\\
\delta_{k,0}, & \text{if }n=0
\end{cases}
\ \ \ \ \ \ \ \ \ \ \text{for each }n\in\mathbb{N}.
\label{pf.thm.fps.comps.n-into-k-parts.1st.rhs2}%
\end{equation}

\end{proof}

We have thus answered our Question 2. Let us summarize the two answers we have
found ((\ref{pf.thm.fps.comps.n-into-k-parts.1st.rhs1}) and
(\ref{pf.thm.fps.comps.n-into-k-parts.1st.rhs2})) in the following theorem
(\cite[Theorem 2.10.1]{19fco}):

\begin{theorem}
\label{thm.fps.comps.num-comps-n-k}Let $n,k\in\mathbb{N}$. Then, the \# of
compositions of $n$ into $k$ parts is%
\[
\dbinom{n-1}{n-k}=%
\begin{cases}
\dbinom{n-1}{k-1}, & \text{if }n>0;\\
\delta_{k,0}, & \text{if }n=0.
\end{cases}
\]

\end{theorem}

This theorem has other proofs as well. See \cite[Proof of Theorem
2.10.1]{19fco} for a proof by bijection and \cite[solution to Exercise
2.10.2]{19fco} for a proof by induction.

As an easy consequence of Theorem \ref{thm.fps.comps.num-comps-n-k}, we can
now answer Question 1 as well:

\begin{theorem}
\label{thm.fps.comps.num-comps-n}Let $n\in\mathbb{N}$. Then, the \# of
compositions of $n$ is%
\[%
\begin{cases}
2^{n-1}, & \text{if }n>0;\\
1, & \text{if }n=0.
\end{cases}
\]

\end{theorem}

\begin{proof}
[Proof of Theorem \ref{thm.fps.comps.num-comps-n} (sketched).]If $n=0$, then
the \# of compositions of $n$ is $1$ (since the empty tuple $\left(
{}\right)  $ is the only composition of $0$). Thus, for the rest of this
proof, we WLOG assume that $n\neq0$. Hence, we must prove that the \# of
compositions of $n$ is $2^{n-1}$.

If $\left(  n_{1},n_{2},\ldots,n_{m}\right)  $ is a composition of $n$, then
$m\in\left\{  1,2,\ldots,n\right\}  $ (why?). In other words, any composition
of $n$ is a composition of $n$ into $k$ parts for some $k\in\left\{
1,2,\ldots,n\right\}  $. Hence,%
\begin{align*}
&  \left(  \text{\# of compositions of }n\right) \\
&  =\sum_{k\in\left\{  1,2,\ldots,n\right\}  }\underbrace{\left(  \text{\# of
compositions of }n\text{ into }k\text{ parts}\right)  }_{\substack{=\dbinom
{n-1}{n-k}\\\text{(by Theorem \ref{thm.fps.comps.num-comps-n-k})}}}\\
&  =\sum_{k\in\left\{  1,2,\ldots,n\right\}  }\dbinom{n-1}{n-k}=\sum_{k=1}%
^{n}\dbinom{n-1}{n-k}=\sum_{k=0}^{n-1}\dbinom{n-1}{k}%
\end{align*}
(here, we have substituted $k$ for $n-k$ in the sum). Comparing this with%
\begin{align*}
2^{n-1}  &  =\left(  1+1\right)  ^{n-1}=\sum_{k=0}^{n-1}\dbinom{n-1}%
{k}\underbrace{1^{k}1^{n-1-k}}_{=1}\ \ \ \ \ \ \ \ \ \ \left(  \text{by the
binomial theorem}\right) \\
&  =\sum_{k=0}^{n-1}\dbinom{n-1}{k},
\end{align*}
we obtain $\left(  \text{\# of compositions of }n\right)  =2^{n-1}$. This
proves Theorem \ref{thm.fps.comps.num-comps-n}.
\end{proof}

\subsubsection{Weak compositions}

A variant of compositions are the \emph{weak compositions}. These are like
compositions, but their entries have to only be nonnegative rather than
positive. For the sake of completeness, let us give their definition in
full:\footnote{Beware that the word \textquotedblleft weak
composition\textquotedblright\ does not have a unique meaning in the
literature.}

\begin{definition}
\label{def.fps.wcomps}\textbf{(a)} An \emph{(integer) weak composition} means
a (finite) tuple of nonnegative integers. \medskip

\textbf{(b)} The \emph{size} of a weak composition $\alpha=\left(  \alpha
_{1},\alpha_{2},\ldots,\alpha_{m}\right)  $ is defined to be the integer
$\alpha_{1}+\alpha_{2}+\cdots+\alpha_{m}$. It is written $\left\vert
\alpha\right\vert $. \medskip

\textbf{(c)} The \emph{length} of a weak composition $\alpha=\left(
\alpha_{1},\alpha_{2},\ldots,\alpha_{m}\right)  $ is defined to be the integer
$m$. It is written $\ell\left(  \alpha\right)  $. \medskip

\textbf{(d)} Let $n\in\mathbb{N}$. A \emph{weak composition of }$n$ means a
weak composition whose size is $n$. \medskip

\textbf{(e)} Let $n\in\mathbb{N}$ and $k\in\mathbb{N}$. A \emph{weak
composition of }$n$\emph{ into }$k$\emph{ parts} is a weak composition whose
size is $n$ and whose length is $k$.
\end{definition}

\begin{example}
\label{exa.fps.wcomps.1}The tuple $\left(  3,0,1,2\right)  $ is a weak
composition with size $3+0+1+2=6$ and length $4$. In other words, it is a weak
composition of $6$ into $4$ parts. It is not a composition, since one of its
entries is a $0$.
\end{example}

Weak compositions are a rather natural analogue of compositions, but behave
dissimilarly in one important way: Any $n\in\mathbb{N}$ has infinitely many
weak compositions. Indeed, all the tuples $\left(  n\right)  ,\ \left(
n,0\right)  ,\ \left(  n,0,0\right)  ,\ \left(  n,0,0,0\right)  ,\ \ldots$ are
weak compositions of $n$ (and of course, there are many more weak compositions
of $n$, unless $n=0$). Thus, it makes no sense to look for an analogue of
Theorem \ref{thm.fps.comps.num-comps-n} for weak compositions. However, an
analogue of Theorem \ref{thm.fps.comps.num-comps-n-k} does exist:

\begin{theorem}
\label{thm.fps.comps.num-wcomps-n-k}Let $n,k\in\mathbb{N}$. Then, the \# of
weak compositions of $n$ into $k$ parts is
\[
\dbinom{n+k-1}{n}=%
\begin{cases}
\dbinom{n+k-1}{k-1}, & \text{if }k>0;\\
\delta_{n,0}, & \text{if }k=0.
\end{cases}
\]

\end{theorem}

\begin{proof}
[Proof of Theorem \ref{thm.fps.comps.num-wcomps-n-k} (sketched).](See
\cite[Theorem 2.10.5]{19fco} for details and alternative proofs.) Adding $1$
to a nonnegative integer yields a positive integer. Furthermore, if we add $1$
to each entry of a $k$-tuple, then the sum of all entries of the $k$-tuple
increases by $k$.

Thus, if $\left(  \alpha_{1},\alpha_{2},\ldots,\alpha_{k}\right)  $ is a weak
composition of $n$ into $k$ parts, then \newline$\left(  \alpha_{1}%
+1,\alpha_{2}+1,\ldots,\alpha_{k}+1\right)  $ is a composition of $n+k$ into
$k$ parts. Hence, we can define a map%
\begin{align*}
\left\{  \text{weak compositions of }n\text{ into }k\text{ parts}\right\}   &
\rightarrow\left\{  \text{compositions of }n+k\text{ into }k\text{
parts}\right\}  ,\\
\left(  \alpha_{1},\alpha_{2},\ldots,\alpha_{k}\right)   &  \mapsto\left(
\alpha_{1}+1,\alpha_{2}+1,\ldots,\alpha_{k}+1\right)  .
\end{align*}
Furthermore, it is easy to see that this map is a bijection (indeed, its
inverse is rather easy to construct). Thus, by the bijection principle, we
have%
\begin{align*}
&  \left\vert \left\{  \text{weak compositions of }n\text{ into }k\text{
parts}\right\}  \right\vert \\
&  =\left\vert \left\{  \text{compositions of }n+k\text{ into }k\text{
parts}\right\}  \right\vert \\
&  =\left(  \text{\# of compositions of }n+k\text{ into }k\text{ parts}\right)
\\
&  =\dbinom{n+k-1}{n+k-k}\ \ \ \ \ \ \ \ \ \ \left(
\begin{array}
[c]{c}%
\text{by the first equality sign in Theorem \ref{thm.fps.comps.num-comps-n-k}%
,}\\
\text{applied to }n+k\text{ instead of }n
\end{array}
\right) \\
&  =\dbinom{n+k-1}{n}.
\end{align*}
Thus, we have shown that the \# of weak compositions of $n$ into $k$ parts is
$\dbinom{n+k-1}{n}$. It remains to prove that this equals $%
\begin{cases}
\dbinom{n+k-1}{k-1}, & \text{if }k>0;\\
\delta_{n,0}, & \text{if }k=0
\end{cases}
$ as well. This is similar to how we obtained
(\ref{pf.thm.fps.comps.n-into-k-parts.1st.rhs2}): If $k=0$, then it is clear
by inspection; otherwise it follows from Theorem \ref{thm.binom.sym}. Theorem
\ref{thm.fps.comps.num-wcomps-n-k} is proven.
\end{proof}

\subsubsection{Weak compositions with entries from $\left\{  0,1,\ldots
,p-1\right\}  $}

Theorems \ref{thm.fps.comps.num-comps-n-k} and
\ref{thm.fps.comps.num-wcomps-n-k} may stir up hopes that other tuple-counting
problems also have simple answers. Let us see if this hope holds up.

\begin{proof}
[An attempt.]We fix three nonnegative integers $n$, $k$ and $p$. We are
looking for the \# of $k$-tuples $\left(  \alpha_{1},\alpha_{2},\ldots
,\alpha_{k}\right)  \in\left\{  0,1,\ldots,p-1\right\}  ^{k}$ satisfying
$\alpha_{1}+\alpha_{2}+\cdots+\alpha_{k}=n$. In other words, we are looking
for the \# of all weak compositions of $n$ into $k$ parts with the property
that each entry is $<p$. Let us denote this \# by $w_{n,k,p}$.

Just as when counting compositions, we invoke a generating function. Forget
that we fixed $n$, and define the FPS%
\[
W_{k,p}:=\sum_{n\in\mathbb{N}}w_{n,k,p}x^{n}\in\mathbb{Q}\left[  \left[
x\right]  \right]  .
\]
For each $n\in\mathbb{N}$, we have%
\begin{align*}
w_{n,k,p}  &  =\left(  \text{\# of all }k\text{-tuples }\left(  \alpha
_{1},\alpha_{2},\ldots,\alpha_{k}\right)  \in\left\{  0,1,\ldots,p-1\right\}
^{k}\text{ }\right. \\
&  \ \ \ \ \ \ \ \ \ \ \ \ \ \ \ \ \ \ \ \ \ \left.  \text{satisfying }%
\alpha_{1}+\alpha_{2}+\cdots+\alpha_{k}=n\right) \\
&  =\sum\limits_{\substack{\left(  \alpha_{1},\alpha_{2},\ldots,\alpha
_{k}\right)  \in\left\{  0,1,\ldots,p-1\right\}  ^{k};\\\alpha_{1}+\alpha
_{2}+\cdots+\alpha_{k}=n}}1.
\end{align*}
Thus, we can rewrite the equality $W_{k,p}=\sum_{n\in\mathbb{N}}w_{n,k,p}%
x^{n}$ as%
\begin{align*}
W_{k,p}  &  =\sum_{n\in\mathbb{N}}\left(  \sum\limits_{\substack{\left(
\alpha_{1},\alpha_{2},\ldots,\alpha_{k}\right)  \in\left\{  0,1,\ldots
,p-1\right\}  ^{k};\\\alpha_{1}+\alpha_{2}+\cdots+\alpha_{k}=n}}1\right)
x^{n}\\
&  =\sum_{n\in\mathbb{N}}\ \ \sum\limits_{\substack{\left(  \alpha_{1}%
,\alpha_{2},\ldots,\alpha_{k}\right)  \in\left\{  0,1,\ldots,p-1\right\}
^{k};\\\alpha_{1}+\alpha_{2}+\cdots+\alpha_{k}=n}}\ \ \underbrace{x^{n}%
}_{\substack{=x^{\alpha_{1}+\alpha_{2}+\cdots+\alpha_{k}}\\\text{(since
}\alpha_{1}+\alpha_{2}+\cdots+\alpha_{k}=n\text{)}}}\\
&  =\underbrace{\sum_{n\in\mathbb{N}}\ \ \sum\limits_{\substack{\left(
\alpha_{1},\alpha_{2},\ldots,\alpha_{k}\right)  \in\left\{  0,1,\ldots
,p-1\right\}  ^{k};\\\alpha_{1}+\alpha_{2}+\cdots+\alpha_{k}=n}}}%
_{=\sum\limits_{\left(  \alpha_{1},\alpha_{2},\ldots,\alpha_{k}\right)
\in\left\{  0,1,\ldots,p-1\right\}  ^{k}}}\ \ \underbrace{x^{\alpha_{1}%
+\alpha_{2}+\cdots+\alpha_{k}}}_{=x^{\alpha_{1}}x^{\alpha_{2}}\cdots
x^{\alpha_{k}}}\\
&  =\sum\limits_{\left(  \alpha_{1},\alpha_{2},\ldots,\alpha_{k}\right)
\in\left\{  0,1,\ldots,p-1\right\}  ^{k}}x^{\alpha_{1}}x^{\alpha_{2}}\cdots
x^{\alpha_{k}}\\
&  =\left(  \sum_{\alpha_{1}\in\left\{  0,1,\ldots,p-1\right\}  }x^{\alpha
_{1}}\right)  \left(  \sum_{\alpha_{2}\in\left\{  0,1,\ldots,p-1\right\}
}x^{\alpha_{2}}\right)  \cdots\left(  \sum_{\alpha_{k}\in\left\{
0,1,\ldots,p-1\right\}  }x^{\alpha_{k}}\right) \\
&  \ \ \ \ \ \ \ \ \ \ \ \ \ \ \ \ \ \ \ \ \left(
\begin{array}
[c]{c}%
\text{by the same product rule that we used back}\\
\text{in the proof of Theorem \ref{thm.fps.gen-newton}}%
\end{array}
\right) \\
&  =\left(  \sum_{n\in\left\{  0,1,\ldots,p-1\right\}  }x^{n}\right)  ^{k}%
\end{align*}
(here, we have renamed all the $k$ summation indices as $n$, and realized that
all $k$ sums are identical). Since%
\[
\sum_{n\in\left\{  0,1,\ldots,p-1\right\}  }x^{n}=x^{0}+x^{1}+\cdots
+x^{p-1}=\dfrac{1-x^{p}}{1-x}%
\]
(the last equality sign here is easy to check\footnote{Just show that $\left(
1-x\right)  \left(  x^{0}+x^{1}+\cdots+x^{p-1}\right)  =1-x^{p}$.}), this can
be rewritten further as%
\begin{equation}
W_{k,p}=\left(  \dfrac{1-x^{p}}{1-x}\right)  ^{k}=\left(  1-x^{p}\right)
^{k}\left(  1-x\right)  ^{-k}. \label{pf.thm.fps.comps.num-wpcomps-n-k.3}%
\end{equation}
In order to expand the right hand side, let us expand $\left(  1-x^{p}\right)
^{k}$ and $\left(  1-x\right)  ^{-k}$ separately.

The binomial theorem yields%
\[
\left(  1-x^{p}\right)  ^{k}=\sum_{j=0}^{k}\left(  -1\right)  ^{j}\dbinom
{k}{j}\underbrace{\left(  x^{p}\right)  ^{j}}_{=x^{pj}}=\sum_{j=0}^{k}\left(
-1\right)  ^{j}\dbinom{k}{j}x^{pj}=\sum_{j\in\mathbb{N}}\left(  -1\right)
^{j}\dbinom{k}{j}x^{pj}%
\]
(here, we have extended the range of the summation from $j\in\left\{
0,1,\ldots,k\right\}  $ to all $j\in\mathbb{N}$; this did not change the sum,
since all the newly introduced addends are $0$). Multiplying this by
(\ref{pf.thm.fps.comps.num-comps-n-k.1-x-k}), we obtain%
\begin{align*}
&  \left(  1-x^{p}\right)  ^{k}\left(  1-x\right)  ^{-k}\\
&  =\left(  \sum_{j\in\mathbb{N}}\left(  -1\right)  ^{j}\dbinom{k}{j}%
x^{pj}\right)  \left(  \sum_{j\in\mathbb{N}}\dbinom{j+k-1}{j}x^{j}\right) \\
&  =\left(  \sum_{j\in\mathbb{N}}\left(  -1\right)  ^{j}\dbinom{k}{j}%
x^{pj}\right)  \left(  \sum_{i\in\mathbb{N}}\dbinom{i+k-1}{i}x^{i}\right) \\
&  \ \ \ \ \ \ \ \ \ \ \ \ \ \ \ \ \ \ \ \ \left(
\begin{array}
[c]{c}%
\text{here, we have renamed the summation index }j\\
\text{as }i\text{ in the second sum}%
\end{array}
\right) \\
&  =\underbrace{\sum_{\left(  i,j\right)  \in\mathbb{N}\times\mathbb{N}}%
}_{=\sum_{n\in\mathbb{N}}\ \ \sum_{\substack{\left(  i,j\right)  \in
\mathbb{N}\times\mathbb{N};\\pj+i=n}}}\left(  -1\right)  ^{j}\dbinom{k}%
{j}\dbinom{i+k-1}{i}x^{pj+i}\\
&  =\sum_{n\in\mathbb{N}}\ \ \sum_{\substack{\left(  i,j\right)  \in
\mathbb{N}\times\mathbb{N};\\pj+i=n}}\left(  -1\right)  ^{j}\dbinom{k}%
{j}\dbinom{i+k-1}{i}\underbrace{x^{pj+i}}_{\substack{=x^{n}\\\text{(since
}pj+i=n\text{)}}}\\
&  =\sum_{n\in\mathbb{N}}\ \ \sum_{\substack{\left(  i,j\right)  \in
\mathbb{N}\times\mathbb{N};\\pj+i=n}}\left(  -1\right)  ^{j}\dbinom{k}%
{j}\dbinom{i+k-1}{i}x^{n}\\
&  =\sum_{n\in\mathbb{N}}\left(  \sum_{\substack{\left(  i,j\right)
\in\mathbb{N}\times\mathbb{N};\\pj+i=n}}\left(  -1\right)  ^{j}\dbinom{k}%
{j}\dbinom{i+k-1}{i}\right)  x^{n}.
\end{align*}
Thus, (\ref{pf.thm.fps.comps.num-wpcomps-n-k.3}) becomes
\[
W_{k,p}=\left(  1-x^{p}\right)  ^{k}\left(  1-x\right)  ^{-k}=\sum
_{n\in\mathbb{N}}\left(  \sum_{\substack{\left(  i,j\right)  \in
\mathbb{N}\times\mathbb{N};\\pj+i=n}}\left(  -1\right)  ^{j}\dbinom{k}%
{j}\dbinom{i+k-1}{i}\right)  x^{n}.
\]
Comparing coefficients in this equality, we find that each $n\in\mathbb{N}$
satisfies%
\begin{align*}
w_{n,k,p}  &  =\sum_{\substack{\left(  i,j\right)  \in\mathbb{N}%
\times\mathbb{N};\\pj+i=n}}\left(  -1\right)  ^{j}\dbinom{k}{j}\dbinom
{i+k-1}{i}\\
&  =\sum_{\substack{j\in\mathbb{N};\\pj\leq n}}\left(  -1\right)  ^{j}%
\dbinom{k}{j}\dbinom{n-pj+k-1}{n-pj}\\
&  \ \ \ \ \ \ \ \ \ \ \ \ \ \ \ \ \ \ \ \ \left(
\begin{array}
[c]{c}%
\text{here, we have substituted }\left(  n-pj,j\right)  \text{ for }\left(
i,j\right)  \text{ in the sum,}\\
\text{since any pair }\left(  i,j\right)  \in\mathbb{N}\times\mathbb{N}\text{
satisfying }pj+i=n\\
\text{is uniquely determined by its second entry }j
\end{array}
\right) \\
&  =\sum_{j\in\mathbb{N}}\left(  -1\right)  ^{j}\dbinom{k}{j}\dbinom
{n-pj+k-1}{n-pj}\\
&  \ \ \ \ \ \ \ \ \ \ \ \ \ \ \ \ \ \ \ \ \left(
\begin{array}
[c]{c}%
\text{here, we have extended the range of summation by}\\
\text{dropping the \textquotedblleft}pj\leq n\text{\textquotedblright%
\ requirement; this does not change}\\
\text{the sum, since all newly introduced addends are }0
\end{array}
\right) \\
&  =\sum_{j=0}^{k}\left(  -1\right)  ^{j}\dbinom{k}{j}\dbinom{n-pj+k-1}{n-pj}%
\end{align*}
(here, we have removed all addends with $j>k$ from the sum; this does not
change the sum, since all these addends are $0$).
\end{proof}

Thus, we have proved the following fact:

\begin{theorem}
\label{thm.fps.comps.num-wpcomps-n-k}Let $n,k,p\in\mathbb{N}$. Then, the \# of
$k$-tuples $\left(  \alpha_{1},\alpha_{2},\ldots,\alpha_{k}\right)
\in\left\{  0,1,\ldots,p-1\right\}  ^{k}$ satisfying $\alpha_{1}+\alpha
_{2}+\cdots+\alpha_{k}=n$ is%
\[
\sum_{j=0}^{k}\left(  -1\right)  ^{j}\dbinom{k}{j}\dbinom{n-pj+k-1}{n-pj}.
\]

\end{theorem}

In general, this expression is the simplest we can get. A combinatorial proof
of Theorem \ref{thm.fps.comps.num-wpcomps-n-k} can be found in \cite[Exercise
2.10.6]{19fco}.

However, the particular case when $p=2$ is worth exploring, as it allows for a
much simpler expression. Indeed, the $k$-tuples $\left(  \alpha_{1},\alpha
_{2},\ldots,\alpha_{k}\right)  \in\left\{  0,1\right\}  ^{k}$ are just the
\textquotedblleft\emph{binary }$k$\emph{-strings}\textquotedblright, i.e., the
$k$-tuples formed of $0$s and $1$s. Imposing the condition $\alpha_{1}%
+\alpha_{2}+\cdots+\alpha_{k}=n$ on such a $k$-tuple is tantamount to
requiring that it contain exactly $n$ many $1$s. Therefore, the \# of all
$k$-tuples $\left(  \alpha_{1},\alpha_{2},\ldots,\alpha_{k}\right)
\in\left\{  0,1\right\}  ^{k}$ satisfying $\alpha_{1}+\alpha_{2}+\cdots
+\alpha_{k}=n$ is $\dbinom{k}{n}$, since we have to choose which $n$ of its
$k$ positions will be occupied by $1$s (and then all remaining $k-n$ positions
will be occupied by $0$s). However, Theorem
\ref{thm.fps.comps.num-wpcomps-n-k} (applied to $p=2$) yields that this \#
equals $\sum_{j=0}^{k}\left(  -1\right)  ^{j}\dbinom{k}{j}\dbinom
{n-2j+k-1}{n-2j}$. Comparing these two results, we obtain the following identity:

\begin{proposition}
\label{prop.fps.comps.num-w2comps-n-k-id}Let $n,k\in\mathbb{N}$. Then,%
\[
\dbinom{k}{n}=\sum_{j=0}^{k}\left(  -1\right)  ^{j}\dbinom{k}{j}%
\dbinom{n-2j+k-1}{n-2j}.
\]

\end{proposition}

\subsection{\label{sec.gf.xneq}$x^{n}$-equivalence}

We now return to general properties of FPSs.

\begin{definition}
\label{def.fps.xneq}Let $n\in\mathbb{N}$. Let $f,g\in K\left[  \left[
x\right]  \right]  $ be two FPSs. We write $f\overset{x^{n}}{\equiv}g$ if and
only if%
\[
\text{each }m\in\left\{  0,1,\ldots,n\right\}  \text{ satisfies }\left[
x^{m}\right]  f=\left[  x^{m}\right]  g.
\]

Thus, we have defined a binary relation $\overset{x^{n}}{\equiv}$ on the set
$K\left[  \left[  x\right]  \right]  $. We say that an FPS $f$ is $x^{n}%
$\emph{-equivalent} to an FPS $g$ if and only if $f\overset{x^{n}}{\equiv}g$.
\end{definition}

Thus, an FPS $f$ is $x^{n}$-equivalent to an FPS $g$ if and only if the first
$n+1$ coefficients of $f$ agree with the first $n+1$ coefficients of $g$. Here
are some examples:

\begin{example}
\textbf{(a)} Consider the two FPSs
\[
\left(  1+x\right)  ^{3}=1+3x+3x^{2}+x^{3}=x^{0}+3x^{1}+3x^{2}+1x^{3}%
+0x^{4}+0x^{5}+\cdots
\]
and%
\[
\dfrac{1}{1-3x}=1+3x+\left(  3x\right)  ^{2}+\left(  3x\right)  ^{3}%
+\cdots=x^{0}+3x^{1}+9x^{2}+27x^{3}+\cdots.
\]
Thus, $\left(  1+x\right)  ^{3}\overset{x^{1}}{\equiv}\dfrac{1}{1-3x}$, since
each $m\in\left\{  0,1\right\}  $ satisfies $\left[  x^{m}\right]  \left(
\left(  1+x\right)  ^{3}\right)  =\left[  x^{m}\right]  \dfrac{1}{1-3x}$
(indeed, we have $\left[  x^{0}\right]  \left(  \left(  1+x\right)
^{3}\right)  =1=\left[  x^{0}\right]  \dfrac{1}{1-3x}$ and $\left[
x^{1}\right]  \left(  \left(  1+x\right)  ^{3}\right)  =3=\left[
x^{1}\right]  \dfrac{1}{1-3x}$). Of course, this also shows that $\left(
1+x\right)  ^{3}\overset{x^{0}}{\equiv}\dfrac{1}{1-3x}$. However, we don't
have $\left(  1+x\right)  ^{3}\overset{x^{2}}{\equiv}\dfrac{1}{1-3x}$ (at
least not for $K=\mathbb{Z}$), since $\left[  x^{2}\right]  \left(  \left(
1+x\right)  ^{3}\right)  =3$ does not equal $\left[  x^{2}\right]  \dfrac
{1}{1-3x}=9$. \medskip

\textbf{(b)} More generally, $\left(  1+x\right)  ^{n}\overset{x^{1}}{\equiv
}\dfrac{1}{1-nx}$ and $\left(  1+x\right)  ^{n}\overset{x^{1}}{\equiv}1+nx$
for each $n\in\mathbb{Z}$. \medskip

\textbf{(c)} The generating function $F=F\left(  x\right)  =0+1x+1x^{2}%
+2x^{3}+3x^{4}+5x^{5}+\cdots$ from Section \ref{sec.gf.exas} satisfies
$F\overset{x^{2}}{\equiv}x+x^{2}$ and $F\overset{x^{3}}{\equiv}x+x^{2}+2x^{3}%
$. \medskip

\textbf{(d)} If $f\in K\left[  \left[  x\right]  \right]  $ is any FPS, and if
$n\in\mathbb{N}$, then there exists a polynomial $p\in K\left[  x\right]  $
such that $f\overset{x^{n}}{\equiv}p$. Indeed, we can take $p=\sum_{k=0}%
^{n}\left(  \left[  x^{k}\right]  f\right)  \cdot x^{k}$.
\end{example}

One way to get an intuition for the relation $\overset{x^{n}}{\equiv}$ is to
think of it as a kind of \textquotedblleft approximate
equality\textquotedblright\ up to degree $n$. (This makes the most sense if
one thinks of $x$ as an infinitesimal quantity, in which case a term $\lambda
x^{k}$ (with $\lambda\in K$) is the more \textquotedblleft
important\textquotedblright\ the lower $k$ is. From this viewpoint,
$f\overset{x^{n}}{\equiv}g$ means that the FPSs $f$ and $g$ agree in their
$n+1$ most \textquotedblleft important\textquotedblright\ terms and differ at
most in their \textquotedblleft error terms\textquotedblright.) For this
reason, the statement \textquotedblleft$f\overset{x^{n}}{\equiv}%
g$\textquotedblright\ is sometimes written as \textquotedblleft$f=g+o\left(
x^{n}\right)  $\textquotedblright\ (an algebraic imitation of Landau's
little-o notation from asymptotic analysis). Another intuition comes from
elementary number theory: The relation $\overset{x^{n}}{\equiv}$ is similar to
congruence of integers modulo a given integer. This is more than a similarity;
the relation $\overset{x^{n}}{\equiv}$ can in fact be restated as a
divisibility in the same fashion as for congruences of integers (see
Proposition \ref{prop.fps.xneq-multiple} below). For this reason, the
statement \textquotedblleft$f\overset{x^{n}}{\equiv}g$\textquotedblright\ is
sometimes written as \textquotedblleft$f\equiv g\operatorname{mod}x^{n+1}%
$\textquotedblright. We shall, however, eschew both of these alternative
notations, and use the original notation \textquotedblleft$f\overset{x^{n}%
}{\equiv}g$\textquotedblright\ from Definition \ref{def.fps.xneq}, as both
intuitions (while useful) would distract from the simplicity of Definition
\ref{def.fps.xneq}\footnote{Case in point: Definition \ref{def.fps.xneq} can
be generalized to multivariate FPSs, but the two intuitions are no longer
available (or, worse, give the \textquotedblleft wrong\textquotedblright%
\ concepts) when extended to this generality.}.

Here are some basic properties of the relation $\overset{x^{n}}{\equiv}$ (some
of which will be used without explicit reference):

\begin{theorem}
\label{thm.fps.xneq.props}Let $n\in\mathbb{N}$. \medskip

\textbf{(a)} The relation $\overset{x^{n}}{\equiv}$ on $K\left[  \left[
x\right]  \right]  $ is an equivalence relation. In other words:

\begin{itemize}
\item This relation is reflexive (i.e., we have $f\overset{x^{n}}{\equiv}f$
for each $f\in K\left[  \left[  x\right]  \right]  $).

\item This relation is transitive (i.e., if three FPSs $f,g,h\in K\left[
\left[  x\right]  \right]  $ satisfy $f\overset{x^{n}}{\equiv}g$ and
$g\overset{x^{n}}{\equiv}h$, then $f\overset{x^{n}}{\equiv}h$).

\item This relation is symmetric (i.e., if two FPSs $f,g\in K\left[  \left[
x\right]  \right]  $ satisfy $f\overset{x^{n}}{\equiv}g$, then
$g\overset{x^{n}}{\equiv}f$).
\end{itemize}

\textbf{(b)} If $a,b,c,d\in K\left[  \left[  x\right]  \right]  $ are four
FPSs satisfying $a\overset{x^{n}}{\equiv}b$ and $c\overset{x^{n}}{\equiv}d$,
then we also have%
\begin{align}
&  a+c\overset{x^{n}}{\equiv}b+d;\label{eq.thm.fps.xneq.props.b.+}\\
&  a-c\overset{x^{n}}{\equiv}b-d;\label{eq.thm.fps.xneq.props.b.-}\\
&  ac\overset{x^{n}}{\equiv}bd. \label{eq.thm.fps.xneq.props.b.*}%
\end{align}

\textbf{(c)} If $a,b\in K\left[  \left[  x\right]  \right]  $ are two FPSs
satisfying $a\overset{x^{n}}{\equiv}b$, then $\lambda a\overset{x^{n}}{\equiv
}\lambda b$ for each $\lambda\in K$. \medskip

\textbf{(d)} If $a,b\in K\left[  \left[  x\right]  \right]  $ are two
invertible FPSs satisfying $a\overset{x^{n}}{\equiv}b$, then $a^{-1}%
\overset{x^{n}}{\equiv}b^{-1}$. \medskip

\textbf{(e)} If $a,b,c,d\in K\left[  \left[  x\right]  \right]  $ are four
FPSs satisfying $a\overset{x^{n}}{\equiv}b$ and $c\overset{x^{n}}{\equiv}d$,
and if the FPSs $c$ and $d$ are invertible, then we also have%
\begin{equation}
\dfrac{a}{c}\overset{x^{n}}{\equiv}\dfrac{b}{d}.
\label{eq.thm.fps.xneq.props.e./}%
\end{equation}

\textbf{(f)} Let $S$ be a finite set. Let $\left(  a_{s}\right)  _{s\in S}\in
K\left[  \left[  x\right]  \right]  ^{S}$ and $\left(  b_{s}\right)  _{s\in
S}\in K\left[  \left[  x\right]  \right]  ^{S}$ be two families of FPSs such
that
\begin{equation}
\text{each }s\in S\text{ satisfies }a_{s}\overset{x^{n}}{\equiv}b_{s}.
\label{eq.thm.fps.xneq.props.e.ass}%
\end{equation}
Then, we have%
\begin{align}
&  \sum_{s\in S}a_{s}\overset{x^{n}}{\equiv}\sum_{s\in S}b_{s}%
;\label{eq.thm.fps.xneq.props.e.+}\\
&  \prod_{s\in S}a_{s}\overset{x^{n}}{\equiv}\prod_{s\in S}b_{s}.
\label{eq.thm.fps.xneq.props.e.*}%
\end{align}

\end{theorem}

\begin{proof}
[Proof of Theorem \ref{thm.fps.xneq.props} (sketched).]All of these properties
are analogous to familiar properties of integer congruences, except for
Theorem \ref{thm.fps.xneq.props} \textbf{(d)}, which is moot for integers
(since there are not many integers that are invertible in $\mathbb{Z}$). The
proofs are similarly simple (using (\ref{pf.thm.fps.ring.xn(a+b)=}),
(\ref{pf.thm.fps.ring.xn(a-b)=}), (\ref{pf.thm.fps.ring.xn(ab)=2}) and
(\ref{pf.thm.fps.ring.xn(la)=})). Thus, we shall only give some hints for the
proof of Theorem \ref{thm.fps.xneq.props} \textbf{(d)} here; detailed proofs
of all parts of Theorem \ref{thm.fps.xneq.props} can be found in Section
\ref{sec.details.gf.xneq}.

\textbf{(d)} Let $a,b\in K\left[  \left[  x\right]  \right]  $ be two
invertible FPSs satisfying $a\overset{x^{n}}{\equiv}b$. We want to show that
$a^{-1}\overset{x^{n}}{\equiv}b^{-1}$.

The FPS $a$ is invertible; thus, its constant term $\left[  x^{0}\right]  a$
is invertible in $K$ (by Proposition \ref{prop.fps.invertible}).

Recall that $a\overset{x^{n}}{\equiv}b$. In other words,
\begin{equation}
\text{each }m\in\left\{  0,1,\ldots,n\right\}  \text{ satisfies }\left[
x^{m}\right]  a=\left[  x^{m}\right]  b. \label{pf.thm.fps.xneq.props.d.ass}%
\end{equation}

Now, we want to prove that $a^{-1}\overset{x^{n}}{\equiv}b^{-1}$. In other
words, we want to prove that each $m\in\left\{  0,1,\ldots,n\right\}  $
satisfies $\left[  x^{m}\right]  \left(  a^{-1}\right)  =\left[  x^{m}\right]
\left(  b^{-1}\right)  $. We shall prove this by strong induction on $m$: We
fix some $k\in\left\{  0,1,\ldots,n\right\}  $, and we assume (as an induction
hypothesis) that
\begin{equation}
\left[  x^{m}\right]  \left(  a^{-1}\right)  =\left[  x^{m}\right]  \left(
b^{-1}\right)  \ \ \ \ \ \ \ \ \ \ \text{for each }m\in\left\{  0,1,\ldots
,k-1\right\}  . \label{pf.thm.fps.xneq.props.d.IH}%
\end{equation}
We must now prove that $\left[  x^{k}\right]  \left(  a^{-1}\right)  =\left[
x^{k}\right]  \left(  b^{-1}\right)  $. We know that%
\begin{align*}
\left[  x^{k}\right]  \left(  aa^{-1}\right)   &  =\sum_{i=0}^{k}\left[
x^{i}\right]  a\cdot\left[  x^{k-i}\right]  \left(  a^{-1}\right)
\ \ \ \ \ \ \ \ \ \ \left(  \text{by (\ref{pf.thm.fps.ring.xn(ab)=2})}\right)
\\
&  =\left[  x^{0}\right]  a\cdot\left[  x^{k}\right]  \left(  a^{-1}\right)
+\sum_{i=1}^{k}\left[  x^{i}\right]  a\cdot\left[  x^{k-i}\right]  \left(
a^{-1}\right)
\end{align*}
(here, we have split off the addend for $i=0$ from the sum). Thus,%
\[
\left[  x^{0}\right]  a\cdot\left[  x^{k}\right]  \left(  a^{-1}\right)
+\sum_{i=1}^{k}\left[  x^{i}\right]  a\cdot\left[  x^{k-i}\right]  \left(
a^{-1}\right)  =\left[  x^{k}\right]  \underbrace{\left(  aa^{-1}\right)
}_{=1}=\left[  x^{k}\right]  1.
\]
We can solve this equation for $\left[  x^{k}\right]  \left(  a^{-1}\right)  $
(since $\left[  x^{0}\right]  a$ is invertible), and thus obtain%
\[
\left[  x^{k}\right]  \left(  a^{-1}\right)  =\dfrac{1}{\left[  x^{0}\right]
a}\cdot\left(  \left[  x^{k}\right]  1-\sum_{i=1}^{k}\left[  x^{i}\right]
a\cdot\left[  x^{k-i}\right]  \left(  a^{-1}\right)  \right)  .
\]
The same argument (applied to $b$ instead of $a$) yields%
\[
\left[  x^{k}\right]  \left(  b^{-1}\right)  =\dfrac{1}{\left[  x^{0}\right]
b}\cdot\left(  \left[  x^{k}\right]  1-\sum_{i=1}^{k}\left[  x^{i}\right]
b\cdot\left[  x^{k-i}\right]  \left(  b^{-1}\right)  \right)  .
\]
The right hand sides of the latter two equalities are equal (since each
$i\in\left\{  1,2,\ldots,k\right\}  $ satisfies $\left[  x^{i}\right]
a=\left[  x^{i}\right]  b$ as a consequence of
(\ref{pf.thm.fps.xneq.props.d.ass}), and satisfies $\left[  x^{k-i}\right]
\left(  a^{-1}\right)  =\left[  x^{k-i}\right]  \left(  b^{-1}\right)  $ as a
consequence of (\ref{pf.thm.fps.xneq.props.d.IH}), and since we have $\left[
x^{0}\right]  a=\left[  x^{0}\right]  b$ as a consequence of
(\ref{pf.thm.fps.xneq.props.d.ass})). Hence, the left hand sides must also be
equal. In other words, $\left[  x^{k}\right]  \left(  a^{-1}\right)  =\left[
x^{k}\right]  \left(  b^{-1}\right)  $, which is precisely what we wanted to
prove. Thus, the induction step is complete, so that $a^{-1}\overset{x^{n}%
}{\equiv}b^{-1}$ is proved. Thus, Theorem \ref{thm.fps.xneq.props}
\textbf{(d)} is proved. (See Section \ref{sec.details.gf.xneq} for more details.)
\end{proof}

Let us next characterize $x^{n}$-equivalence of FPSs in terms of divisibility:

\begin{proposition}
\label{prop.fps.xneq-multiple}Let $n\in\mathbb{N}$. Let $f,g\in K\left[
\left[  x\right]  \right]  $ be two FPSs. Then, we have $f\overset{x^{n}%
}{\equiv}g$ if and only if the FPS $f-g$ is a multiple of $x^{n+1}$.
\end{proposition}

\begin{proof}
[Proof of Proposition \ref{prop.fps.xneq-multiple}.]See Section
\ref{sec.details.gf.xneq} for this proof (a simple consequence of Lemma
\ref{lem.fps.muls-of-xn}).
\end{proof}

Finally, here is a subtler property of $x^{n}$-equivalence similar to the ones
in Theorem \ref{thm.fps.xneq.props} \textbf{(b)}:

\begin{proposition}
\label{prop.fps.xneq.comp}Let $n\in\mathbb{N}$. Let $a,b,c,d\in K\left[
\left[  x\right]  \right]  $ be four FPSs satisfying $a\overset{x^{n}}{\equiv
}b$ and $c\overset{x^{n}}{\equiv}d$ and $\left[  x^{0}\right]  c=0$ and
$\left[  x^{0}\right]  d=0$. Then,%
\[
a\circ c\overset{x^{n}}{\equiv}b\circ d.
\]

\end{proposition}

\begin{proof}
[Proof of Proposition \ref{prop.fps.xneq.comp} (sketched).]Write $a$ and $b$
as $a=\sum_{i\in\mathbb{N}}a_{i}x^{i}$ and $b=\sum_{i\in\mathbb{N}}b_{i}x^{i}$
(with $a_{i},b_{i}\in K$). Then, $a\overset{x^{n}}{\equiv}b$ means that
$a_{i}=b_{i}$ for all $i\leq n$. Combine this with $c^{i}\overset{x^{n}%
}{\equiv}d^{i}$ (which holds for all $i\in\mathbb{N}$ as a consequence of
$c\overset{x^{n}}{\equiv}d$ and of (\ref{eq.thm.fps.xneq.props.e.*})) to
obtain the relation $a_{i}c^{i}\overset{x^{n}}{\equiv}b_{i}d^{i}$ for all
$i\leq n$. But this relation also holds for all $i>n$, since all such $i$
satisfy $\left[  x^{m}\right]  \left(  c^{i}\right)  =\left[  x^{m}\right]
\left(  d^{i}\right)  =0$ for all $m\in\left\{  0,1,\ldots,n\right\}  $ (a
consequence of Lemma \ref{lem.fps.muls-of-xn} using $\left[  x^{0}\right]
c=0$ and $\left[  x^{0}\right]  d=0$). Thus, the relation $a_{i}%
c^{i}\overset{x^{n}}{\equiv}b_{i}d^{i}$ holds for all $i\in\mathbb{N}$.
Summing over all $i$, we find $a\circ c\overset{x^{n}}{\equiv}b\circ d$.

See Section \ref{sec.details.gf.xneq} for the details of this argument.
\end{proof}

\subsection{\label{sec.gf.prod}Infinite products}

Let us now extend our FPS playground somewhat. We have made sense of infinite
sums. What about infinite products?

\subsubsection{\label{subsec.gf.prod.exa}An example}

We start with a \textbf{motivating example} (due to Euler, in
\cite[\S 328--329]{Euler48}), which we shall first discuss informally.

Assume (for the time being) that the infinite product%
\begin{equation}
\prod_{i\in\mathbb{N}}\left(  1+x^{2^{i}}\right)  =\left(  1+x^{1}\right)
\left(  1+x^{2}\right)  \left(  1+x^{4}\right)  \left(  1+x^{8}\right)
\cdots\label{eq.fps.prod.binary.1}%
\end{equation}
in the ring $K\left[  \left[  x\right]  \right]  $ is meaningful, and that
such products behave as nicely as finite products. Can we simplify this product?

We can observe that each $i\in\mathbb{N}$ satisfies $1+x^{2^{i}}%
=\dfrac{1-x^{2^{i+1}}}{1-x^{2^{i}}}$ (since $1-x^{2^{i+1}}=1-\left(  x^{2^{i}%
}\right)  ^{2}=\left(  1-x^{2^{i}}\right)  \left(  1+x^{2^{i}}\right)  $).
Multiplying these equalities over all $i\in\mathbb{N}$, we obtain%
\[
\prod_{i\in\mathbb{N}}\left(  1+x^{2^{i}}\right)  =\prod_{i\in\mathbb{N}%
}\dfrac{1-x^{2^{i+1}}}{1-x^{2^{i}}}=\dfrac{1-x^{2}}{1-x^{1}}\cdot
\dfrac{1-x^{4}}{1-x^{2}}\cdot\dfrac{1-x^{8}}{1-x^{4}}\cdot\dfrac{1-x^{16}%
}{1-x^{8}}\cdot\cdots.
\]
The product on the right hand side here is a \emph{telescoping product} --
meaning that each numerator is cancelled by the denominator of the following
fraction. Assuming (somewhat plausibly, but far from rigorously) that we are
allowed to cancel infinitely many factors from an infinite product, we thus
end up with a single $1-x^{1}$ factor in the denominator. That is, our product
simplifies to $\dfrac{1}{1-x^{1}}$. Thus, we obtain%
\begin{equation}
\prod_{i\in\mathbb{N}}\left(  1+x^{2^{i}}\right)  =\dfrac{1}{1-x^{1}}%
=\dfrac{1}{1-x}=1+x+x^{2}+x^{3}+\cdots. \label{eq.fps.prod.binary.2}%
\end{equation}

This was not very rigorous, so let us try to compute the product $\prod
_{i\in\mathbb{N}}\left(  1+x^{2^{i}}\right)  $ in a different way. Namely, we
recall a simple fact about finite products: If $a_{0},a_{1},\ldots,a_{m}$ are
finitely many elements of a commutative ring, then the product
\begin{equation}
\prod_{i=0}^{m}\left(  1+a_{i}\right)  =\left(  1+a_{0}\right)  \left(
1+a_{1}\right)  \cdots\left(  1+a_{m}\right)
\label{eq.fps.prod.binary.prod-fin}%
\end{equation}
equals the sum\footnote{The indices $i_{1},i_{2},\ldots,i_{k}$ of the sum are
supposed to be nonnegative integers. Their number $k$ is allowed to vary; in
particular it can be $0$, in which case the chain of inequalities $i_{1}%
<i_{2}<\cdots<i_{k}\leq m$ is regarded as a tautology (i.e., automatically
true).}
\[
\sum_{i_{1}<i_{2}<\cdots<i_{k}\leq m}a_{i_{1}}a_{i_{2}}\cdots a_{i_{k}}%
\]
of all the $2^{m+1}$ \textquotedblleft sub-products\textquotedblright\ of the
product $a_{0}a_{1}\cdots a_{m}$ (because this sum is what we obtain if we
expand $\prod_{i=0}^{m}\left(  1+a_{i}\right)  $ by repeatedly applying
distributivity). For example, for $m=2$, this is saying that%
\begin{align*}
&  \left(  1+a_{0}\right)  \left(  1+a_{1}\right)  \left(  1+a_{2}\right) \\
&  =1+a_{0}+a_{1}+a_{2}+a_{0}a_{1}+a_{0}a_{2}+a_{1}a_{2}+a_{0}a_{1}a_{2}.
\end{align*}
Now, it is plausible to expect the same formula $\prod_{i=0}^{m}\left(
1+a_{i}\right)  =\sum_{I\subseteq\left\{  0,1,\ldots,m\right\}  }%
\ \ \prod_{i\in I}a_{i}$ to hold if \textquotedblleft$m$ is $\infty
$\textquotedblright\ (that is, if the product ranges over all $i\in\mathbb{N}%
$), provided that the product is meaningful. In other words, it is plausible
to expect that%
\begin{equation}
\prod_{i\in\mathbb{N}}\left(  1+a_{i}\right)  =\sum_{i_{1}<i_{2}<\cdots<i_{k}%
}a_{i_{1}}a_{i_{2}}\cdots a_{i_{k}} \label{eq.fps.prod.binary.prod-inf}%
\end{equation}
for any infinite sequence $a_{0},a_{1},a_{2},\ldots$ as long as $\prod
_{i\in\mathbb{N}}\left(  1+a_{i}\right)  $ makes sense. (It's a little bit
more complicated than that, but we aren't trying to be fully rigorous yet. The
correct condition is that the sequence $\left(  a_{0},a_{1},a_{2}%
,\ldots\right)  $ is summable.) If we now apply
(\ref{eq.fps.prod.binary.prod-inf}) to $a_{i}=x^{2^{i}}$, then we obtain%
\begin{align}
\prod_{i\in\mathbb{N}}\left(  1+x^{2^{i}}\right)   &  =\sum_{i_{1}%
<i_{2}<\cdots<i_{k}}x^{2^{i_{1}}}x^{2^{i_{2}}}\cdots x^{2^{i_{k}}}=\sum
_{i_{1}<i_{2}<\cdots<i_{k}}x^{2^{i_{1}}+2^{i_{2}}+\cdots+2^{i_{k}}}\nonumber\\
&  =\sum_{n\in\mathbb{N}}q_{n}x^{n}, \label{eq.fps.prod.binary.5}%
\end{align}
where $q_{n}$ is the \# of ways to write the integer $n$ as a sum $2^{i_{1}%
}+2^{i_{2}}+\cdots+2^{i_{k}}$ with nonnegative integers $i_{1},i_{2}%
,\ldots,i_{k}$ satisfying $i_{1}<i_{2}<\cdots<i_{k}$. Comparing this with
(\ref{eq.fps.prod.binary.2}), we obtain%
\[
\sum_{n\in\mathbb{N}}q_{n}x^{n}=1+x+x^{2}+x^{3}+\cdots=\sum_{n\in\mathbb{N}%
}x^{n},
\]
at least if our assumptions were valid. Comparing coefficients, this would
mean that $q_{n}=1$ for each $n\in\mathbb{N}$. In other words, each
$n\in\mathbb{N}$ can be written in \textbf{exactly one} way as a sum
$2^{i_{1}}+2^{i_{2}}+\cdots+2^{i_{k}}$ with nonnegative integers $i_{1}%
,i_{2},\ldots,i_{k}$ satisfying $i_{1}<i_{2}<\cdots<i_{k}$. In other words,
each $n\in\mathbb{N}$ can be written uniquely as a finite sum of distinct
powers of $2$.

Is this true? Yes, because this is just saying that each $n\in\mathbb{N}$ has
a unique binary representation. For example, $21=2^{4}+2^{2}+2^{0}$
corresponds to the binary representation $21=\left(  10101\right)  _{2}$.

Thus, the two results we have obtained in (\ref{eq.fps.prod.binary.2}) and
(\ref{eq.fps.prod.binary.5}) are actually equal, which is reassuring. Yet,
this does not replace a formal definition of infinite products that rigorously
justifies the above arguments.

\subsubsection{\label{subsec.gf.prod.def}A rigorous definition}

One way of rigorously defining infinite products of FPSs can be found in
\cite[\S 7.5]{Loehr-BC}. However, this definition only defines infinite
products of the form $\prod_{i\in\mathbb{N}}$ or $\prod_{i=1}^{\infty}$, but
not (for example) of the form $\prod_{I\subseteq\mathbb{N}}$ or $\prod
_{\left(  i,j\right)  \in\mathbb{N}\times\mathbb{N}}$. Another definition of
infinite products uses the $\operatorname*{Log}$ and $\operatorname*{Exp}$
bijections from Definition \ref{def.fps.Exp-Log-maps} to turn products into
sums; but this requires $K$ to be a $\mathbb{Q}$-algebra (since
$\operatorname*{Log}$ and $\operatorname*{Exp}$ aren't defined otherwise).
Thus, we shall give a different definition here.

We recall our definition of infinite sums of FPSs (Definition
\ref{def.fps.summable}):

\begin{statement}
\textbf{Definition \ref{def.fps.summable} (repeated).} A (possibly infinite)
family $\left(  \mathbf{a}_{i}\right)  _{i\in I}$ of FPSs is said to be
\emph{summable} if%
\[
\text{for each }n\in\mathbb{N}\text{, all but finitely many }i\in I\text{
satisfy }\left[  x^{n}\right]  \mathbf{a}_{i}=0.
\]
In this case, the sum $\sum_{i\in I}\mathbf{a}_{i}$ is defined to be the FPS
with%
\[
\left[  x^{n}\right]  \left(  \sum_{i\in I}\mathbf{a}_{i}\right)
=\underbrace{\sum_{i\in I}\left[  x^{n}\right]  \mathbf{a}_{i}}%
_{\substack{\text{an essentially}\\\text{finite sum}}%
}\ \ \ \ \ \ \ \ \ \ \text{for all }n\in\mathbb{N}\text{.}%
\]

\end{statement}

This is how we defined infinite sums of FPSs. We cannot use the same
definition for infinite products, because usually%
\[
\text{we \textbf{don't} expect to have }\left[  x^{n}\right]  \left(
\prod_{i\in I}\mathbf{a}_{i}\right)  =\prod_{i\in I}\left[  x^{n}\right]
\mathbf{a}_{i}%
\]
(after all, multiplication of FPSs is not defined coefficientwise). The
condition \textquotedblleft all but finitely many $i\in I$ satisfy $\left[
x^{n}\right]  \mathbf{a}_{i}=0$\textquotedblright\ is therefore not what we
are looking for.

Let us instead go back to the idea behind Definition \ref{def.fps.summable}.
Let us fix some $n\in\mathbb{N}$. What was the actual purpose of the
\textquotedblleft all but finitely many $i\in I$ satisfy $\left[
x^{n}\right]  \mathbf{a}_{i}=0$\textquotedblright\ condition? The purpose was
to ensure that the coefficient $\left[  x^{n}\right]  \left(  \sum_{i\in
I}\mathbf{a}_{i}\right)  $ is determined by \textbf{finitely many} of the
$\mathbf{a}_{i}$'s. In other words, the purpose was to ensure that there is a
\textbf{finite} partial sum of $\sum_{i\in I}\mathbf{a}_{i}$ such that if we
add any further $\mathbf{a}_{i}$'s to this partial sum, then the coefficient
of $x^{n}$ does not change any more. Here is a way to restate this condition
more rigorously: There is a \textbf{finite} subset $M$ of $I$ such that every
finite subset $J$ of $I$ satisfying $M\subseteq J\subseteq I$ satisfies%
\[
\left[  x^{n}\right]  \left(  \sum_{i\in M}\mathbf{a}_{i}\right)  =\left[
x^{n}\right]  \left(  \sum_{i\in J}\mathbf{a}_{i}\right)  .
\]
(The subset $M$ here is the indexing set of our finite partial sum, and the
set $J$ is what it becomes if we add some further $\mathbf{a}_{i}$'s to this
partial sum.)

This condition is a mouthful; this is why we found it easier to boil it down
to the simple \textquotedblleft all but finitely many $i\in I$ satisfy
$\left[  x^{n}\right]  \mathbf{a}_{i}=0$\textquotedblright\ condition in the
case of infinite sums. However, in the case of infinite products, we cannot
boil it down to something this simple; thus, we have to live with it.

Fortunately, we can simplify our life by giving this condition a name. To
highlight the analogy, let us define it both for sums and for products:

\begin{definition}
\label{def.fps.determines-xn-coeff}Let $\left(  \mathbf{a}_{i}\right)  _{i\in
I}\in K\left[  \left[  x\right]  \right]  ^{I}$ be a (possibly infinite)
family of FPSs. Let $n\in\mathbb{N}$. Let $M$ be a finite subset of $I$.
\medskip

\textbf{(a)} We say that $M$ \emph{determines the }$x^{n}$\emph{-coefficient
in the sum of }$\left(  \mathbf{a}_{i}\right)  _{i\in I}$ if every finite
subset $J$ of $I$ satisfying $M\subseteq J\subseteq I$ satisfies%
\[
\left[  x^{n}\right]  \left(  \sum_{i\in J}\mathbf{a}_{i}\right)  =\left[
x^{n}\right]  \left(  \sum_{i\in M}\mathbf{a}_{i}\right)  .
\]
(You can think of this condition as saying \textquotedblleft If you add any
further $\mathbf{a}_{i}$s to the sum $\sum_{i\in M}\mathbf{a}_{i}$, then the
$x^{n}$-coefficient stays unchanged\textquotedblright, or, more informally:
\textquotedblleft If you want to know the $x^{n}$-coefficient of $\sum_{i\in
I}\mathbf{a}_{i}$, it suffices to take the partial sum over all $i\in
M$\textquotedblright.) \medskip

\textbf{(b)} We say that $M$ \emph{determines the }$x^{n}$\emph{-coefficient
in the product of }$\left(  \mathbf{a}_{i}\right)  _{i\in I}$ if every finite
subset $J$ of $I$ satisfying $M\subseteq J\subseteq I$ satisfies%
\[
\left[  x^{n}\right]  \left(  \prod_{i\in J}\mathbf{a}_{i}\right)  =\left[
x^{n}\right]  \left(  \prod_{i\in M}\mathbf{a}_{i}\right)  .
\]
(You can think of this condition as saying \textquotedblleft If you multiply
any further $\mathbf{a}_{i}$s to the product $\prod_{i\in M}\mathbf{a}_{i}$,
then the $x^{n}$-coefficient stays unchanged\textquotedblright, or, more
informally: \textquotedblleft If you want to know the $x^{n}$-coefficient of
$\prod_{i\in I}\mathbf{a}_{i}$, it suffices to take the partial product over
all $i\in M$\textquotedblright.)
\end{definition}

\begin{example}
\textbf{(a)} Consider the family
\begin{align*}
&  \left(  \left(  x+x^{2}\right)  ^{i}\right)  _{i\in\mathbb{N}}\\
&  =\left(  \left(  x+x^{2}\right)  ^{0},\ \ \left(  x+x^{2}\right)
^{1},\ \ \left(  x+x^{2}\right)  ^{2},\ \ \left(  x+x^{2}\right)
^{3},\ \ \left(  x+x^{2}\right)  ^{4},\ \ \ldots\right) \\
&  =\left(  1,\ \ x+x^{2},\ \ x^{2}+2x^{3}+x^{4},\ \ x^{3}+3x^{4}+3x^{5}%
+x^{6},\ \ x^{4}+4x^{5}+6x^{6}+4x^{7}+x^{8},\ \ \ldots\right)
\end{align*}
of FPSs. The subset $\left\{  2,3\right\}  $ of $\mathbb{N}$ determines the
$x^{3}$-coefficient in the sum of this family $\left(  \left(  x+x^{2}\right)
^{i}\right)  _{i\in\mathbb{N}}$, because every finite subset $J$ of
$\mathbb{N}$ satisfying $\left\{  2,3\right\}  \subseteq J\subseteq\mathbb{N}$
satisfies%
\[
\left[  x^{3}\right]  \left(  \sum_{i\in J}\left(  x+x^{2}\right)
^{i}\right)  =\left[  x^{3}\right]  \left(  \sum_{i\in\left\{  2,3\right\}
}\left(  x+x^{2}\right)  ^{i}\right)
\]
(this is simply a consequence of the fact that the only two entries of our
family that have a nonzero $x^{3}$-coefficient are the entries $\left(
x+x^{2}\right)  ^{i}$ for $i\in\left\{  2,3\right\}  $). Thus, any finite
subset of $\mathbb{N}$ that contains $\left\{  2,3\right\}  $ as a subset
determines the $x^{3}$-coefficient in the sum of this family $\left(  \left(
x+x^{2}\right)  ^{i}\right)  _{i\in\mathbb{N}}$. \medskip

\textbf{(b)} Consider the family%
\[
\left(  1+x^{i}\right)  _{i\in\mathbb{N}}=\left(  1+1,\ \ 1+x,\ \ 1+x^{2}%
,\ \ 1+x^{3},\ \ 1+x^{4},\ \ \ldots\right)
\]
of FPSs over $K=\mathbb{Z}$. The subset $\left\{  0,1,2,3\right\}  $ of
$\mathbb{N}$ determines the $x^{3}$-coefficient in the product of this family
$\left(  1+x^{i}\right)  _{i\in\mathbb{N}}$, because every finite subset $J$
of $\mathbb{N}$ satisfying $\left\{  0,1,2,3\right\}  \subseteq J\subseteq
\mathbb{N}$ satisfies%
\[
\left[  x^{3}\right]  \left(  \prod_{i\in J}\left(  1+x^{i}\right)  \right)
=\left[  x^{3}\right]  \left(  \prod_{i\in\left\{  0,1,2,3\right\}  }\left(
1+x^{i}\right)  \right)  .
\]
(This is because multiplying an FPS by any of the polynomials $1+x^{4}%
,\ \ 1+x^{5},\ \ 1+x^{6},\ \ \ldots$ leaves its $x^{3}$-coefficient
unchanged.) Thus, any finite subset of $\mathbb{N}$ that contains $\left\{
0,1,2,3\right\}  $ as a subset determines the $x^{3}$-coefficient in the
product of this family $\left(  1+x^{i}\right)  _{i\in\mathbb{N}}$.

On the other hand, the subset $\left\{  0,3\right\}  $ of $\mathbb{N}$ does
\textbf{not} determine the $x^{3}$-coefficient in the product of $\left(
1+x^{i}\right)  _{i\in\mathbb{N}}$. To see this, it suffices to notice that%
\[
\underbrace{\left[  x^{3}\right]  \left(  \prod_{i\in\left\{  0,1,2,3\right\}
}\left(  1+x^{i}\right)  \right)  }_{\substack{=\left[  x^{3}\right]  \left(
\left(  1+x^{0}\right)  \left(  1+x^{1}\right)  \left(  1+x^{2}\right)
\left(  1+x^{3}\right)  \right)  \\=4}}\neq\underbrace{\left[  x^{3}\right]
\left(  \prod_{i\in\left\{  0,3\right\}  }\left(  1+x^{i}\right)  \right)
}_{\substack{=\left[  x^{3}\right]  \left(  \left(  1+x^{0}\right)  \left(
1+x^{3}\right)  \right)  \\=2}}.
\]
(The philosophical reason is that, even though the monomial $x^{3}$ itself
does not appear in any of the entries $1+x^{1}$ and $1+x^{2}$, it does emerge
in the product of these two entries with the constant term of $\prod
_{i\in\left\{  0,3\right\}  }\left(  1+x^{i}\right)  =\left(  1+1\right)
\left(  1+x^{3}\right)  $.) \medskip

\textbf{(c)} Here is a simple but somewhat slippery example: Let $I=\left\{
1,2,3\right\}  $, and define the FPS $\mathbf{a}_{i}=1+\left(  -1\right)
^{i}x$ over $K=\mathbb{Z}$ for each $i\in I$ (so that $\mathbf{a}%
_{1}=\mathbf{a}_{3}=1-x$ and $\mathbf{a}_{2}=1+x$). Consider the finite family
$\left(  \mathbf{a}_{i}\right)  _{i\in I}$. Its product is already
well-defined by dint of its finiteness:%
\[
\prod_{i\in I}\mathbf{a}_{i}=\left(  1-x\right)  \left(  1+x\right)  \left(
1-x\right)  =1-x-x^{2}+x^{3}.
\]
But let us nevertheless see which subsets $M$ of $I$ determine the $x^{1}%
$-coefficient in the product of $\left(  \mathbf{a}_{i}\right)  _{i\in I}$.
The subset $I$ itself clearly does (since the only finite subset $J$ of $I$
satisfying $I\subseteq J\subseteq I$ is $I$ itself). The subset $\left\{
1\right\}  $ of $I$, however, does not, even though the product $\prod_{i\in
I}\mathbf{a}_{i}$ has the same $x^{1}$-coefficient as its subproduct
$\prod_{i\in\left\{  1\right\}  }\mathbf{a}_{i}$. (To see why $\left\{
1\right\}  $ does not determine the $x^{1}$-coefficient in the product of
$\left(  \mathbf{a}_{i}\right)  _{i\in I}$, it suffices to check that not
every finite subset $J$ of $I$ satisfying $\left\{  1\right\}  \subseteq
J\subseteq I$ satisfies the equality%
\[
\left[  x^{1}\right]  \left(  \prod_{i\in J}\mathbf{a}_{i}\right)  =\left[
x^{1}\right]  \left(  \prod_{i\in\left\{  1\right\}  }\mathbf{a}_{i}\right)
.
\]
For example, $J=\left\{  1,2\right\}  $ does not, since $\left[  x^{1}\right]
\left(  \prod_{i\in\left\{  1\right\}  }\mathbf{a}_{i}\right)  =-1$ but
$\left[  x^{1}\right]  \left(  \prod_{i\in\left\{  1,2\right\}  }%
\mathbf{a}_{i}\right)  =0$.)

This shows that, in order for a subset $M$ of $I$ to determine the $x^{n}%
$-coefficient in the product of a family $\left(  \mathbf{a}_{i}\right)
_{i\in I}$, it does \textbf{not} suffice to check that $\left[  x^{n}\right]
\left(  \prod_{i\in I}\mathbf{a}_{i}\right)  =\left[  x^{n}\right]  \left(
\prod_{i\in M}\mathbf{a}_{i}\right)  $; it rather needs to be shown that
$\left[  x^{n}\right]  \left(  \prod_{i\in J}\mathbf{a}_{i}\right)  =\left[
x^{n}\right]  \left(  \prod_{i\in M}\mathbf{a}_{i}\right)  $ for each finite
subset $J$ of $I$ satisfying $M\subseteq J\subseteq I$.
\end{example}

\begin{definition}
\label{def.fps.xn-coeff-fin-determined}Let $\left(  \mathbf{a}_{i}\right)
_{i\in I}\in K\left[  \left[  x\right]  \right]  ^{I}$ be a (possibly
infinite) family of FPSs. Let $n\in\mathbb{N}$. \medskip

\textbf{(a)} We say that \emph{the }$x^{n}$\emph{-coefficient in the sum of
}$\left(  \mathbf{a}_{i}\right)  _{i\in I}$\emph{ is finitely determined} if
there is a finite subset $M$ of $I$ that determines the $x^{n}$-coefficient in
the sum of $\left(  \mathbf{a}_{i}\right)  _{i\in I}$. \medskip

\textbf{(b)} We say that \emph{the }$x^{n}$\emph{-coefficient in the product
of }$\left(  \mathbf{a}_{i}\right)  _{i\in I}$\emph{ is finitely determined}
if there is a finite subset $M$ of $I$ that determines the $x^{n}$-coefficient
in the product of $\left(  \mathbf{a}_{i}\right)  _{i\in I}$.
\end{definition}

Using these concepts, we can now reword our definition of infinite sums as follows:

\begin{proposition}
\label{prop.fps.summable=fin-det}Let $\left(  \mathbf{a}_{i}\right)  _{i\in
I}\in K\left[  \left[  x\right]  \right]  ^{I}$ be a (possibly infinite)
family of FPSs. Then: \medskip

\textbf{(a)} The family $\left(  \mathbf{a}_{i}\right)  _{i\in I}$ is summable
if and only if each coefficient in its sum is finitely determined (i.e., for
each $n\in\mathbb{N}$, the $x^{n}$-coefficient in the sum of $\left(
\mathbf{a}_{i}\right)  _{i\in I}$ is finitely determined). \medskip

\textbf{(b)} If the family $\left(  \mathbf{a}_{i}\right)  _{i\in I}$ is
summable, then its sum $\sum_{i\in I}\mathbf{a}_{i}$ is the FPS whose $x^{n}%
$-coefficient (for any $n\in\mathbb{N}$) can be computed as follows: If
$n\in\mathbb{N}$, and if $M$ is a finite subset of $I$ that determines the
$x^{n}$-coefficient in the sum of $\left(  \mathbf{a}_{i}\right)  _{i\in I}$,
then
\[
\left[  x^{n}\right]  \left(  \sum_{i\in I}\mathbf{a}_{i}\right)  =\left[
x^{n}\right]  \left(  \sum_{i\in M}\mathbf{a}_{i}\right)  .
\]

\end{proposition}

\begin{proof}
Easy and LTTR.
\end{proof}

Inspired by Proposition \ref{prop.fps.summable=fin-det}, we can now define
infinite products of FPSs at last:

\begin{definition}
\label{def.fps.multipliable}Let $\left(  \mathbf{a}_{i}\right)  _{i\in I}$ be
a (possibly infinite) family of FPSs. Then: \medskip

\textbf{(a)} The family $\left(  \mathbf{a}_{i}\right)  _{i\in I}$ is said to
be \emph{multipliable} if and only if each coefficient in its product is
finitely determined. \medskip

\textbf{(b)} If the family $\left(  \mathbf{a}_{i}\right)  _{i\in I}$ is
multipliable, then its \emph{product} $\prod_{i\in I}\mathbf{a}_{i}$ is
defined to be the FPS whose $x^{n}$-coefficient (for any $n\in\mathbb{N}$) can
be computed as follows: If $n\in\mathbb{N}$, and if $M$ is a finite subset of
$I$ that determines the $x^{n}$-coefficient in the product of $\left(
\mathbf{a}_{i}\right)  _{i\in I}$, then
\[
\left[  x^{n}\right]  \left(  \prod_{i\in I}\mathbf{a}_{i}\right)  =\left[
x^{n}\right]  \left(  \prod_{i\in M}\mathbf{a}_{i}\right)  .
\]

\end{definition}

\begin{proposition}
\label{prop.fps.multipliable.prod-wd}This definition of $\prod_{i\in
I}\mathbf{a}_{i}$ is well-defined -- i.e., the coefficient $\left[
x^{n}\right]  \left(  \prod_{i\in M}\mathbf{a}_{i}\right)  $ does not depend
on $M$ (as long as $M$ is a finite subset of $I$ that determines the $x^{n}%
$-coefficient in the product of $\left(  \mathbf{a}_{i}\right)  _{i\in I}$).
\end{proposition}

\begin{proof}
Let $n\in\mathbb{N}$. We need to check that the coefficient $\left[
x^{n}\right]  \left(  \prod_{i\in M}\mathbf{a}_{i}\right)  $ does not depend
on $M$ (as long as $M$ is a finite subset of $I$ that determines the $x^{n}%
$-coefficient in the product of $\left(  \mathbf{a}_{i}\right)  _{i\in I}$).
In other words, we need to check that if $M_{1}$ and $M_{2}$ are two finite
subsets of $I$ that each determine the $x^{n}$-coefficient in the product of
$\left(  \mathbf{a}_{i}\right)  _{i\in I}$, then%
\begin{equation}
\left[  x^{n}\right]  \left(  \prod_{i\in M_{1}}\mathbf{a}_{i}\right)
=\left[  x^{n}\right]  \left(  \prod_{i\in M_{2}}\mathbf{a}_{i}\right)  .
\label{pf.prop.fps.multipliable.prod-wd.goal}%
\end{equation}

So let us prove this. Let $M_{1}$ and $M_{2}$ be two finite subsets of $I$
that each determine the $x^{n}$-coefficient in the product of $\left(
\mathbf{a}_{i}\right)  _{i\in I}$. Thus, in particular, $M_{1}$ determines the
$x^{n}$-coefficient in the product of $\left(  \mathbf{a}_{i}\right)  _{i\in
I}$. In other words, every finite subset $J$ of $I$ satisfying $M_{1}\subseteq
J\subseteq I$ satisfies%
\[
\left[  x^{n}\right]  \left(  \prod_{i\in J}\mathbf{a}_{i}\right)  =\left[
x^{n}\right]  \left(  \prod_{i\in M_{1}}\mathbf{a}_{i}\right)  .
\]
Applying this to $J=M_{1}\cup M_{2}$, we obtain%
\begin{equation}
\left[  x^{n}\right]  \left(  \prod_{i\in M_{1}\cup M_{2}}\mathbf{a}%
_{i}\right)  =\left[  x^{n}\right]  \left(  \prod_{i\in M_{1}}\mathbf{a}%
_{i}\right)  \label{pf.prop.fps.multipliable.prod-wd.1}%
\end{equation}
(since $M_{1}\cup M_{2}$ is a subset of $I$ satisfying $M_{1}\subseteq
M_{1}\cup M_{2}\subseteq I$). The same argument (with the roles of $M_{1}$ and
$M_{2}$ swapped) yields%
\begin{equation}
\left[  x^{n}\right]  \left(  \prod_{i\in M_{2}\cup M_{1}}\mathbf{a}%
_{i}\right)  =\left[  x^{n}\right]  \left(  \prod_{i\in M_{2}}\mathbf{a}%
_{i}\right)  . \label{pf.prop.fps.multipliable.prod-wd.2}%
\end{equation}
The left hand sides of the equalities
(\ref{pf.prop.fps.multipliable.prod-wd.1}) and
(\ref{pf.prop.fps.multipliable.prod-wd.2}) are equal (since $M_{1}\cup
M_{2}=M_{2}\cup M_{1}$). Thus, the right hand sides are equal as well. In
other words, $\left[  x^{n}\right]  \left(  \prod_{i\in M_{1}}\mathbf{a}%
_{i}\right)  =\left[  x^{n}\right]  \left(  \prod_{i\in M_{2}}\mathbf{a}%
_{i}\right)  $. Thus, we have proved
(\ref{pf.prop.fps.multipliable.prod-wd.goal}), and with it Proposition
\ref{prop.fps.multipliable.prod-wd}.
\end{proof}

The attentive (and pedantic) reader might notice that there is one more thing
that needs to be checked in order to make sure that Definition
\ref{def.fps.multipliable} \textbf{(b)} is legitimate. In fact, this
definition does not merely define (some) infinite products $\prod_{i\in
I}\mathbf{a}_{i}$ of FPSs, but also \textquotedblleft
accidentally\textquotedblright\ gives a new meaning to \textbf{finite}
products $\prod_{i\in I}\mathbf{a}_{i}$ (since a finite family $\left(
\mathbf{a}_{i}\right)  _{i\in I}$ of FPSs is always multipliable). We
therefore need to check that this new meaning does not conflict with the
original definition of a finite product of elements of a commutative ring. In
other words, we need to prove the following:

\begin{proposition}
\label{prop.fps.multipliable.prod-wd2}Let $\left(  \mathbf{a}_{i}\right)
_{i\in I}$ be a finite family of FPSs. Then, the product $\prod_{i\in
I}\mathbf{a}_{i}$ defined according to Definition \ref{def.fps.multipliable}
\textbf{(b)} equals the finite product $\prod_{i\in I}\mathbf{a}_{i}$ defined
in the usual way (i.e., defined as in any commutative ring).
\end{proposition}

\begin{proof}
Argue that $I$ itself is a subset of $I$ that determines all coefficients in
the product of $\left(  \mathbf{a}_{i}\right)  _{i\in I}$. See Section
\ref{sec.details.gf.prod} for a detailed proof.
\end{proof}

We shall apply Convention \ref{conv.fps.infsum} to infinite products just like
we have been applying it to infinite sums. For instance, the product sign
$\prod\limits_{k=m}^{\infty}$ (for a fixed $m\in\mathbb{Z}$) means
$\prod\limits_{k\in\left\{  m,m+1,m+2,\ldots\right\}  }$.

\subsubsection{Why $\prod_{i\in\mathbb{N}}\left(  1+x^{2^{i}}\right)  $ works
and $\prod_{i\in\mathbb{N}}\left(  1+ix\right)  $ doesn't}

Let us now see how Definition \ref{def.fps.multipliable} legitimizes our
product $\prod_{i\in\mathbb{N}}\left(  1+x^{2^{i}}\right)  $ from Subsection
\ref{subsec.gf.prod.exa}. Indeed,%
\[
\prod_{i\in\mathbb{N}}\left(  1+x^{2^{i}}\right)  =\left(  1+x^{1}\right)
\left(  1+x^{2}\right)  \left(  1+x^{4}\right)  \left(  1+x^{8}\right)
\cdots.
\]
If you want to compute the $x^{6}$-coefficient in this product, you only need
to multiply the first 3 factors $\left(  1+x^{1}\right)  \left(
1+x^{2}\right)  \left(  1+x^{4}\right)  $; none of the other factors will
change this coefficient in any way, because multiplying an FPS by $1+x^{m}$
(for some $m>0$) does not change its first $m$ coefficients\footnote{For
example, let us check this for $m=3$: If we multiply an FPS $a_{0}x^{0}%
+a_{1}x^{1}+a_{2}x^{2}+\cdots$ by $1+x^{3}$, then we obtain%
\begin{align*}
&  \left(  a_{0}x^{0}+a_{1}x^{1}+a_{2}x^{2}+\cdots\right)  \left(
1+x^{3}\right) \\
&  =a_{0}x^{0}+a_{1}x^{1}+a_{2}x^{2}+\left(  a_{3}+a_{0}\right)  x^{3}+\left(
a_{4}+a_{1}\right)  x^{4}+\left(  a_{5}+a_{2}\right)  x^{5}+\cdots,
\end{align*}
and so the first $3$ coefficients are left unchanged.}. Likewise, if you want
to compute the $x^{13}$-coefficient of the above product, then you only need
to multiply the first $4$ factors; none of the others will have any effect on
this coefficient. The same logic applies to the $x^{n}$-coefficient for any
$n\in\mathbb{N}$; it is determined by the first $\left\lfloor \log
_{2}n\right\rfloor +1$ factors of the product (if $n>0$; otherwise it equals
$1$). Thus, each coefficient in the product is finitely determined. This means
that the family is multipliable; thus, its product makes sense.

In contrast, the product%
\[
\left(  1+0x\right)  \left(  1+1x\right)  \left(  1+2x\right)  \left(
1+3x\right)  \left(  1+4x\right)  \cdots=\prod_{i\in\mathbb{N}}\left(
1+ix\right)
\]
does not make sense. Indeed, its $x^{1}$-coefficient is not finitely
determined (any of the factors other than $1+0x$ affects it), so the family
$\left(  1+ix\right)  _{i\in\mathbb{N}}$ is not multipliable.

Likewise, the product%
\[
\left(  1+\dfrac{x}{1^{2}}\right)  \left(  1+\dfrac{x}{2^{2}}\right)  \left(
1+\dfrac{x}{3^{2}}\right)  \cdots=\prod_{i\in\left\{  1,2,3,\ldots\right\}
}\left(  1+\dfrac{x}{i^{2}}\right)
\]
does not make sense in $K\left[  \left[  x\right]  \right]  $, even for
$K=\mathbb{C}$ (although the analogous product in complex analysis defines a
holomorphic function).

\subsubsection{A general criterion for multipliability}

Recall our reasoning that we used above to prove that the family $\left(
1+x^{2^{i}}\right)  _{i\in\mathbb{N}}$ is multipliable. The core of this
reasoning was the observation that multiplying an FPS by $1+x^{m}$ (for some
$m>0$) does not change its first $m$ coefficients. This can be generalized: If
$f\in K\left[  \left[  x\right]  \right]  $ is an FPS whose first $m$
coefficients are $0$ (for example, $f$ can be $x^{m}$, in which case we
recover the statement in our preceding sentence), then multiplying an FPS $a$
by $1+f$ does not change its first $m$ coefficients (that is, the first $m$
coefficients of $a\left(  1+f\right)  $ are the first $m$ coefficients of
$a$). This is a useful fact, so let us state it as a lemma (renaming $m$ as
$n+1$):

\begin{lemma}
\label{lem.fps.prod.irlv.1}Let $a,f\in K\left[  \left[  x\right]  \right]  $
be two FPSs. Let $n\in\mathbb{N}$. Assume that
\begin{equation}
\left[  x^{m}\right]  f=0\ \ \ \ \ \ \ \ \ \ \text{for each }m\in\left\{
0,1,\ldots,n\right\}  . \label{eq.lem.fps.prod.irlv.1.ass}%
\end{equation}
Then,
\[
\left[  x^{m}\right]  \left(  a\left(  1+f\right)  \right)  =\left[
x^{m}\right]  a\ \ \ \ \ \ \ \ \ \ \text{for each }m\in\left\{  0,1,\ldots
,n\right\}  .
\]

\end{lemma}

\begin{proof}
[Proof of Lemma \ref{lem.fps.prod.irlv.1}.]The FPS $af$ is a multiple of $f$
(since $af=fa$). Hence, Lemma \ref{lem.fps.prod.irlv.mul} (applied to $u=f$
and $v=af$) yields that%
\begin{equation}
\left[  x^{m}\right]  \left(  af\right)  =0\ \ \ \ \ \ \ \ \ \ \text{for each
}m\in\left\{  0,1,\ldots,n\right\}  \label{pf.lem.fps.prod.irlv.1.4}%
\end{equation}
(since we have assumed that $\left[  x^{m}\right]  f=0$ for each $m\in\left\{
0,1,\ldots,n\right\}  $).

Now, for each $m\in\left\{  0,1,\ldots,n\right\}  $, we have
\begin{align*}
\left[  x^{m}\right]  \left(  \underbrace{a\left(  1+f\right)  }%
_{=a+af}\right)   &  =\left[  x^{m}\right]  \left(  a+af\right)  =\left[
x^{m}\right]  a+\underbrace{\left[  x^{m}\right]  \left(  af\right)
}_{\substack{=0\\\text{(by (\ref{pf.lem.fps.prod.irlv.1.4}))}}%
}\ \ \ \ \ \ \ \ \ \ \left(  \text{by (\ref{pf.thm.fps.ring.xn(a+b)=})}\right)
\\
&  =\left[  x^{m}\right]  a.
\end{align*}
This proves Lemma \ref{lem.fps.prod.irlv.1}.
\end{proof}

For convenience, let us extend Lemma \ref{lem.fps.prod.irlv.1} to products of
several factors:

\begin{lemma}
\label{lem.fps.prod.irlv.fin}Let $a\in K\left[  \left[  x\right]  \right]  $
be an FPS. Let $\left(  f_{i}\right)  _{i\in J}\in K\left[  \left[  x\right]
\right]  ^{J}$ be a finite family of FPSs. Let $n\in\mathbb{N}$. Assume that
each $i\in J$ satisfies
\begin{equation}
\left[  x^{m}\right]  \left(  f_{i}\right)  =0\ \ \ \ \ \ \ \ \ \ \text{for
each }m\in\left\{  0,1,\ldots,n\right\}  .
\label{eq.lem.fps.prod.irlv.fin.ass}%
\end{equation}
Then,
\[
\left[  x^{m}\right]  \left(  a\prod_{i\in J}\left(  1+f_{i}\right)  \right)
=\left[  x^{m}\right]  a\ \ \ \ \ \ \ \ \ \ \text{for each }m\in\left\{
0,1,\ldots,n\right\}  .
\]

\end{lemma}

\begin{proof}
[Proof of Lemma \ref{lem.fps.prod.irlv.fin}.]This is just Lemma
\ref{lem.fps.prod.irlv.1}, applied several times (specifically, $\left\vert
J\right\vert $ many times). See Section \ref{sec.details.gf.prod} for a
detailed proof.
\end{proof}

Now, using Lemma \ref{lem.fps.prod.irlv.fin}, we can obtain the following
convenient criterion for multipliability:

\begin{theorem}
\label{thm.fps.1+f-mulable}Let $\left(  f_{i}\right)  _{i\in I}\in K\left[
\left[  x\right]  \right]  ^{I}$ be a (possibly infinite) summable family of
FPSs. Then, the family $\left(  1+f_{i}\right)  _{i\in I}$ is multipliable.
\end{theorem}

\begin{proof}
[Proof of Theorem \ref{thm.fps.1+f-mulable}.]This is an easy consequence of
Lemma \ref{lem.fps.prod.irlv.fin}. See Section \ref{sec.details.gf.prod} for a
detailed proof.
\end{proof}

We notice two simple sufficient (if rarely satisfied) criteria for multipliability:

\begin{proposition}
\label{prop.fps.1-mulable}If all but finitely many entries of a family
$\left(  \mathbf{a}_{i}\right)  _{i\in I}\in K\left[  \left[  x\right]
\right]  ^{I}$ equal $1$ (that is, if all but finitely many $i\in I$ satisfy
$\mathbf{a}_{i}=1$), then this family is multipliable.
\end{proposition}

\begin{proof}
LTTR. (See Section \ref{sec.details.gf.prod} for a detailed proof.)
\end{proof}

\begin{remark}
\label{rmk.fps.0-mulable}If a family $\left(  \mathbf{a}_{i}\right)  _{i\in
I}\in K\left[  \left[  x\right]  \right]  ^{I}$ contains $0$ as an entry
(i.e., if there exists an $i\in I$ such that $\mathbf{a}_{i}=0$), then this
family is automatically multipliable, and its product is $0$.
\end{remark}

\begin{proof}
Assume that the family $\left(  \mathbf{a}_{i}\right)  _{i\in I}$ contains $0$
as an entry. That is, there exists some $j\in I$ such that $\mathbf{a}_{j}=0$.
Consider this $j$. Now, it is easy to see that the subset $\left\{  j\right\}
$ of $I$ determines all coefficients in the product of $\left(  \mathbf{a}%
_{i}\right)  _{i\in I}$. The details are LTTR.
\end{proof}

\subsubsection{$x^{n}$-approximators}

Working with multipliable families gets slightly easier using the following notion:

\begin{definition}
\label{def.fps.infprod-approx}Let $\left(  \mathbf{a}_{i}\right)  _{i\in I}\in
K\left[  \left[  x\right]  \right]  ^{I}$ be a family of FPSs. Let
$n\in\mathbb{N}$. An $x^{n}$\emph{-approximator} for $\left(  \mathbf{a}%
_{i}\right)  _{i\in I}$ means a finite subset $M$ of $I$ that determines the
first $n+1$ coefficients in the product of $\left(  \mathbf{a}_{i}\right)
_{i\in I}$. (In other words, $M$ has to determine the $x^{m}$-coefficient in
the product of $\left(  \mathbf{a}_{i}\right)  _{i\in I}$ for each
$m\in\left\{  0,1,\ldots,n\right\}  $.)
\end{definition}

The name \textquotedblleft$x^{n}$-approximator\textquotedblright\ is supposed
to hint at the fact that if $M$ is an $x^{n}$-approximator for a multipliable
family $\left(  \mathbf{a}_{i}\right)  _{i\in I}$, then the (finite)
subproduct $\prod_{i\in M}\mathbf{a}_{i}$ \textquotedblleft
approximates\textquotedblright\ the full product $\prod_{i\in I}\mathbf{a}%
_{i}$ up until the $x^{n}$-coefficient (i.e., the first $n+1$ coefficients of
$\prod_{i\in M}\mathbf{a}_{i}$ equal the respective coefficients of
$\prod_{i\in I}\mathbf{a}_{i}$). See Proposition
\ref{prop.fps.infprod-approx-xneq} \textbf{(b)} below for the precise
statement of this fact.

Clearly, an $x^{n}$-approximator for a family $\left(  \mathbf{a}_{i}\right)
_{i\in I}$ always determines the $x^{n}$-coefficient in the product of
$\left(  \mathbf{a}_{i}\right)  _{i\in I}$. But the converse is not true, as
the following example shows:

\begin{example}
Consider the family
\[
\left(  1+x^{2^{i}}\right)  _{i\in\mathbb{N}}=\left(  1+x^{1},\ \ 1+x^{2}%
,\ \ 1+x^{4},\ \ 1+x^{8},\ \ \ldots\right)
\]
of FPSs. The finite subset $\left\{  1,2\right\}  $ of $\mathbb{N}$ determines
the $x^{6}$-coefficient in the product of this family (indeed, the $x^{6}%
$-coefficient of the product $\left(  1+x^{2}\right)  \left(  1+x^{4}\right)
$ is $1$, and this does not change if we multiply any further factors onto
this product), but is not an $x^{6}$-approximator for this family (since,
e.g., it does not determine the $x^{5}$-coefficient in its product).
\end{example}

\begin{lemma}
\label{lem.fps.mulable.approx}Let $\left(  \mathbf{a}_{i}\right)  _{i\in I}\in
K\left[  \left[  x\right]  \right]  ^{I}$ be a multipliable family of FPSs.
Let $n\in\mathbb{N}$. Then, there exists an $x^{n}$-approximator for $\left(
\mathbf{a}_{i}\right)  _{i\in I}$.
\end{lemma}

\begin{proof}
[Proof of Lemma \ref{lem.fps.mulable.approx} (sketched).]This is an easy
consequence of the fact that a union of finitely many finite sets is finite. A
detailed proof can be found in Section \ref{sec.details.gf.prod}.
\end{proof}

As promised above, we can use $x^{n}$-approximators to \textquotedblleft
approximate\textquotedblright\ infinite products of FPSs (in the sense of:
compute the first $n+1$ coefficients of these products). Here is why this
works:\footnote{See Definition \ref{def.fps.xneq} for the meaning of the
symbol \textquotedblleft$\overset{x^{n}}{\equiv}$\textquotedblright\ appearing
in this proposition.}

\begin{proposition}
\label{prop.fps.infprod-approx-xneq}Let $\left(  \mathbf{a}_{i}\right)  _{i\in
I}\in K\left[  \left[  x\right]  \right]  ^{I}$ be a family of FPSs. Let
$n\in\mathbb{N}$. Let $M$ be an $x^{n}$-approximator for $\left(
\mathbf{a}_{i}\right)  _{i\in I}$. Then: \medskip

\textbf{(a)} Every finite subset $J$ of $I$ satisfying $M\subseteq J\subseteq
I$ satisfies%
\[
\prod_{i\in J}\mathbf{a}_{i}\overset{x^{n}}{\equiv}\prod_{i\in M}%
\mathbf{a}_{i}.
\]

\textbf{(b)} If the family $\left(  \mathbf{a}_{i}\right)  _{i\in I}$ is
multipliable, then
\[
\prod_{i\in I}\mathbf{a}_{i}\overset{x^{n}}{\equiv}\prod_{i\in M}%
\mathbf{a}_{i}.
\]

\end{proposition}

\begin{proof}
This follows easily from Definition \ref{def.fps.infprod-approx} and
Definition \ref{def.fps.multipliable} \textbf{(b)}. See Section
\ref{sec.details.gf.prod} for a detailed proof.
\end{proof}

\subsubsection{Properties of infinite products}

We shall now establish some general properties of infinite products. By and
large, these properties are analogous to corresponding properties of finite
products, although they need a few technical requirements (see Remark
\ref{rmk.fps.subfamily-not-mulable} for why).

We begin with the first property, which says in essence that a product can be
broken into two parts:

\begin{proposition}
\label{prop.fps.union-mulable}Let $\left(  \mathbf{a}_{i}\right)  _{i\in I}\in
K\left[  \left[  x\right]  \right]  ^{I}$ be a family of FPSs. Let $J$ be a
subset of $I$. Assume that the subfamilies $\left(  \mathbf{a}_{i}\right)
_{i\in J}$ and $\left(  \mathbf{a}_{i}\right)  _{i\in I\setminus J}$ are
multipliable. Then: \medskip

\textbf{(a)} The entire family $\left(  \mathbf{a}_{i}\right)  _{i\in I}$ is
multipliable. \medskip

\textbf{(b)} We have
\[
\prod_{i\in I}\mathbf{a}_{i}=\left(  \prod_{i\in J}\mathbf{a}_{i}\right)
\cdot\left(  \prod_{i\in I\setminus J}\mathbf{a}_{i}\right)  .
\]

\end{proposition}

\begin{proof}
[Proof of Proposition \ref{prop.fps.union-mulable} (sketched).]Here is the
idea: Fix $n\in\mathbb{N}$. Lemma \ref{lem.fps.mulable.approx} (applied to $J$
instead of $I$) shows that there exists an $x^{n}$-approximator $U$ for
$\left(  \mathbf{a}_{i}\right)  _{i\in J}$. Consider this $U$. Lemma
\ref{lem.fps.mulable.approx} (applied to $I\setminus J$ instead of $I$) shows
that there exists an $x^{n}$-approximator $V$ for $\left(  \mathbf{a}%
_{i}\right)  _{i\in I\setminus J}$. Consider this $V$. Note that $U\cup V$ is
finite (since $U$ and $V$ are finite). Now, it is not hard to see that $U\cup
V$ determines the $x^{n}$-coefficient in the product of $\left(
\mathbf{a}_{i}\right)  _{i\in I}$ (indeed, it is not much harder to see that
$U\cup V$ is an $x^{n}$-approximator for $\left(  \mathbf{a}_{i}\right)
_{i\in I}$). Hence, the $x^{n}$-coefficient in the product of $\left(
\mathbf{a}_{i}\right)  _{i\in I}$ is finitely determined (since $U\cup V$ is
finite). Now, forget that we fixed $n$, and conclude that the family $\left(
\mathbf{a}_{i}\right)  _{i\in I}$ is multipliable. This proves part
\textbf{(a)}. Part \textbf{(b)} easily follows using Proposition
\ref{prop.fps.infprod-approx-xneq} \textbf{(b)}.

The details of this proof can be found in Section \ref{sec.details.gf.prod}.
\end{proof}

A useful particular case of Proposition \ref{prop.fps.union-mulable} is
obtained when the subset $J$ is a single-element set $\left\{  j\right\}  $.
In this case, the proposition says that%
\[
\prod_{i\in I}\mathbf{a}_{i}=\mathbf{a}_{j}\cdot\prod_{i\in I\setminus\left\{
j\right\}  }\mathbf{a}_{i}%
\]
(assuming that the family $\left(  \mathbf{a}_{i}\right)  _{i\in
I\setminus\left\{  j\right\}  }$ is multipliable\footnote{The one-element
family $\left(  \mathbf{a}_{i}\right)  _{i\in\left\{  j\right\}  }$ is, of
course, always multipliable.}). This rule allows us to split off any factor
from a multipliable product, as long as the rest of the product is still
multipliable. \medskip

Our next property generalizes the classical rule $\prod_{i\in I}\left(
a_{i}b_{i}\right)  =\left(  \prod_{i\in I}a_{i}\right)  \cdot\left(
\prod_{i\in I}b_{i}\right)  $ of finite products to infinite ones:

\begin{proposition}
\label{prop.fps.prod-mulable}Let $\left(  \mathbf{a}_{i}\right)  _{i\in I}\in
K\left[  \left[  x\right]  \right]  ^{I}$ and $\left(  \mathbf{b}_{i}\right)
_{i\in I}\in K\left[  \left[  x\right]  \right]  ^{I}$ be two multipliable
families of FPSs. Then: \medskip

\textbf{(a)} The family $\left(  \mathbf{a}_{i}\mathbf{b}_{i}\right)  _{i\in
I}$ is multipliable. \medskip

\textbf{(b)} We have%
\[
\prod_{i\in I}\left(  \mathbf{a}_{i}\mathbf{b}_{i}\right)  =\left(
\prod_{i\in I}\mathbf{a}_{i}\right)  \cdot\left(  \prod_{i\in I}\mathbf{b}%
_{i}\right)  .
\]

\end{proposition}

\begin{proof}
[Proof of Proposition \ref{prop.fps.prod-mulable} (sketched).]Here is the
idea: Fix $n\in\mathbb{N}$. Lemma \ref{lem.fps.mulable.approx} shows that
there exists an $x^{n}$-approximator $U$ for $\left(  \mathbf{a}_{i}\right)
_{i\in I}$. Consider this $U$. Lemma \ref{lem.fps.mulable.approx} (applied to
$\mathbf{b}_{i}$ instead of $\mathbf{a}_{i}$) shows that there exists an
$x^{n}$-approximator $V$ for $\left(  \mathbf{b}_{i}\right)  _{i\in I}$.
Consider this $V$. Note that $U\cup V$ is finite (since $U$ and $V$ are
finite). Now, it is not hard to see that $U\cup V$ is an $x^{n}$-approximator
for $\left(  \mathbf{a}_{i}\mathbf{b}_{i}\right)  _{i\in I}$. From here,
proceed as in the proof of Proposition \ref{prop.fps.union-mulable}.

The details of this proof can be found in Section \ref{sec.details.gf.prod}.
\end{proof}

Next comes an analogous property for quotients instead of products:

\begin{proposition}
\label{prop.fps.div-mulable}Let $\left(  \mathbf{a}_{i}\right)  _{i\in I}\in
K\left[  \left[  x\right]  \right]  ^{I}$ and $\left(  \mathbf{b}_{i}\right)
_{i\in I}\in K\left[  \left[  x\right]  \right]  ^{I}$ be two multipliable
families of FPSs. Assume that the FPS $\mathbf{b}_{i}$ is invertible for each
$i\in I$. Then: \medskip

\textbf{(a)} The family $\left(  \dfrac{\mathbf{a}_{i}}{\mathbf{b}_{i}%
}\right)  _{i\in I}$ is multipliable. \medskip

\textbf{(b)} We have%
\[
\prod_{i\in I}\dfrac{\mathbf{a}_{i}}{\mathbf{b}_{i}}=\dfrac{\prod\limits_{i\in
I}\mathbf{a}_{i}}{\prod\limits_{i\in I}\mathbf{b}_{i}}.
\]

\end{proposition}

\begin{proof}
[Proof of Proposition \ref{prop.fps.div-mulable} (sketch).]This is similar to
the proof of Proposition \ref{prop.fps.prod-mulable}, but using
(\ref{eq.thm.fps.xneq.props.e./}) instead of Lemma
\ref{lem.fps.prod.irlv.cong-mul}. The details of this proof can be found in
Section \ref{sec.details.gf.prod}.
\end{proof}

Now we come to an annoying technicality. In Proposition
\ref{prop.fps.summable.sub}, we have learnt that any subfamily of a summable
family is again summable. Alas, the analogous fact for multipliability is false:

\begin{remark}
\label{rmk.fps.subfamily-not-mulable}Not every subfamily of a multipliable
family is itself multipliable. For example, for $K=\mathbb{Z}$, the family
$\left(  0,1,2,3,\ldots\right)  $ is multipliable, but its subfamily $\left(
1,2,3,\ldots\right)  $ is not.
\end{remark}

Nevertheless, if a family of \textbf{invertible} FPSs is multipliable, then so
is any subfamily of it. In other words:

\begin{proposition}
\label{prop.fps.prods-mulable-subfams}Let $\left(  \mathbf{a}_{i}\right)
_{i\in I}\in K\left[  \left[  x\right]  \right]  ^{I}$ be a multipliable
family of invertible FPSs. Then, any subfamily of $\left(  \mathbf{a}%
_{i}\right)  _{i\in I}$ is multipliable.
\end{proposition}

\begin{proof}
[Proof of Proposition \ref{prop.fps.prods-mulable-subfams} (sketched).]This is
another proof in the tradition of the proofs of Proposition
\ref{prop.fps.union-mulable} and Proposition \ref{prop.fps.prod-mulable}. We
must show that the family $\left(  \mathbf{a}_{i}\right)  _{i\in J}$ is
multipliable whenever $J$ is a subset of $I$. The idea is to show that if $U$
is an $x^{n}$-approximator for $\left(  \mathbf{a}_{i}\right)  _{i\in I}$,
then $U\cap J$ determines the $x^{n}$-coefficient in the product of $\left(
\mathbf{a}_{i}\right)  _{i\in J}$ (and, in fact, is an $x^{n}$-approximator
for $\left(  \mathbf{a}_{i}\right)  _{i\in J}$). This relies on the
invertibility of $\prod_{i\in U\setminus J}\mathbf{a}_{i}$; this is why the
FPS $\mathbf{a}_{i}$ are required to be invertible in the proposition.

The details of this proof can be found in Section \ref{sec.details.gf.prod}.
\end{proof}

As Remark \ref{rmk.fps.subfamily-not-mulable} shows, Proposition
\ref{prop.fps.prods-mulable-subfams} would not hold without the word
\textquotedblleft invertible\textquotedblright. \medskip

Now let us state a trivial rule that is nevertheless worth stating. Recall
that finite products can be reindexed using a bijection -- i.e., if
$f:S\rightarrow T$ is a bijection between two finite sets $S$ and $T$, then
any product $\prod_{t\in T}\mathbf{a}_{t}$ can be rewritten as $\prod_{s\in
S}\mathbf{a}_{f\left(  s\right)  }$. The same holds for infinite products:

\begin{proposition}
\label{prop.fps.prods-mulable-rules.reindex}Let $S$ and $T$ be two sets. Let
$f:S\rightarrow T$ be a bijection. Let $\left(  \mathbf{a}_{t}\right)  _{t\in
T}\in K\left[  \left[  x\right]  \right]  ^{T}$ be a multipliable family of
FPSs. Then,
\[
\prod_{t\in T}\mathbf{a}_{t}=\prod_{s\in S}\mathbf{a}_{f\left(  s\right)  }%
\]
(and, in particular, the product on the right hand side is well-defined, i.e.,
the family $\left(  \mathbf{a}_{f\left(  s\right)  }\right)  _{s\in S}$ is multipliable).
\end{proposition}

In other words, if we reindex a multipliable family of FPSs (using a
bijection), then the resulting family will still be multipliable and have the
same product as the original family.

\begin{proof}
[Proof of Proposition \ref{prop.fps.prods-mulable-rules.reindex}
(sketched).]This is a trivial consequence of the definitions of
multipliability and infinite products.
\end{proof}

Next comes a rule that allows us to break an infinite product into any amount
(possibly infinite!) of subproducts:

\begin{proposition}
\label{prop.fps.prods-mulable-rules.SW1}Let $\left(  \mathbf{a}_{s}\right)
_{s\in S}\in K\left[  \left[  x\right]  \right]  ^{S}$ be a multipliable
family of FPSs. Let $W$ be a set. Let $f:S\rightarrow W$ be a map. Assume that
for each $w\in W$, the family $\left(  \mathbf{a}_{s}\right)  _{s\in
S;\ f\left(  s\right)  =w}$ is multipliable.\footnotemark\ Then,%
\begin{equation}
\prod_{s\in S}\mathbf{a}_{s}=\prod_{w\in W}\ \ \prod_{\substack{s\in
S;\\f\left(  s\right)  =w}}\mathbf{a}_{s}.
\label{eq.prop.fps.prods-mulable-rules.SW1.eq}%
\end{equation}
(In particular, the right hand side is well-defined -- i.e., the family
$\left(  \prod_{\substack{s\in S;\\f\left(  s\right)  =w}}\mathbf{a}%
_{s}\right)  _{w\in W}$ is multipliable.)
\end{proposition}

\footnotetext{Note that this assumption automatically holds if we assume that
all FPSs $\mathbf{a}_{s}$ are invertible. Indeed, the family $\left(
\mathbf{a}_{s}\right)  _{s\in S;\ f\left(  s\right)  =w}$ is a subfamily of
the multipliable family $\left(  \mathbf{a}_{s}\right)  _{s\in S}$ and thus
must itself be multipliable if all the $\mathbf{a}_{s}$ are invertible (by
Proposition \ref{prop.fps.prods-mulable-subfams}).}

\begin{proof}
[Proof of Proposition \ref{prop.fps.prods-mulable-rules.SW1} (sketched).]This
can be derived from the analogous property of finite products, since all
coefficients in a multipliable product are finitely determined (conveniently
using $x^{n}$-approximators, which determine several coefficients at the same
time). See Section \ref{sec.details.gf.prod} for the details of this proof.
\end{proof}

The next rule allows (under certain technical conditions) to interchange
product signs, and to rewrite a nested product as a product over pairs:

\begin{proposition}
[Fubini rule for infinite products of FPSs]%
\label{prop.fps.prods-mulable-rules.fubini1}Let $I$ and $J$ be two sets. Let
$\left(  \mathbf{a}_{\left(  i,j\right)  }\right)  _{\left(  i,j\right)  \in
I\times J}\in K\left[  \left[  x\right]  \right]  ^{I\times J}$ be a
multipliable family of FPSs. Assume that for each $i\in I$, the family
$\left(  \mathbf{a}_{\left(  i,j\right)  }\right)  _{j\in J}$ is multipliable.
Assume that for each $j\in J$, the family $\left(  \mathbf{a}_{\left(
i,j\right)  }\right)  _{i\in I}$ is multipliable.\ Then,%
\[
\prod_{i\in I}\ \ \prod_{j\in J}\mathbf{a}_{\left(  i,j\right)  }%
=\prod_{\left(  i,j\right)  \in I\times J}\mathbf{a}_{\left(  i,j\right)
}=\prod_{j\in J}\ \ \prod_{i\in I}\mathbf{a}_{\left(  i,j\right)  }.
\]
(In particular, all the products appearing in this equality are well-defined.)
\end{proposition}

\begin{proof}
[Proof of Proposition \ref{prop.fps.prods-mulable-rules.fubini1}
(sketched).]The first equality follows by applying Proposition
\ref{prop.fps.prods-mulable-rules.SW1} to $S=I\times J$ and $W=I$ and
$f\left(  i,j\right)  =i$ (and appropriately reindexing the products). The
second equality is analogous. See Section \ref{sec.details.gf.prod} for the
details of this proof.
\end{proof}

The above rules (Proposition \ref{prop.fps.union-mulable}, Proposition
\ref{prop.fps.prod-mulable}, Proposition \ref{prop.fps.div-mulable},
Proposition \ref{prop.fps.prods-mulable-subfams}, Proposition
\ref{prop.fps.prods-mulable-rules.reindex}, Proposition
\ref{prop.fps.prods-mulable-rules.SW1}, Proposition
\ref{prop.fps.prods-mulable-rules.fubini1}) show that infinite products (as we
have defined them) are well-behaved -- i.e., satisfy the usual rules that
finite products satisfy, with only one minor caveat: Subfamilies of
multipliable families might fail to be multipliable (as we saw in Remark
\ref{rmk.fps.subfamily-not-mulable}). If we restrict ourselves to multipliable
families of \textbf{invertible} FPSs, then even this caveat is avoided:

\begin{proposition}
[Fubini rule for infinite products of FPSs, invertible case]%
\label{prop.fps.prods-mulable-rules.fubini}Let $I$ and $J$ be two sets. Let
$\left(  \mathbf{a}_{\left(  i,j\right)  }\right)  _{\left(  i,j\right)  \in
I\times J}\in K\left[  \left[  x\right]  \right]  ^{I\times J}$ be a
multipliable family of invertible FPSs. Then,%
\[
\prod_{i\in I}\ \ \prod_{j\in J}\mathbf{a}_{\left(  i,j\right)  }%
=\prod_{\left(  i,j\right)  \in I\times J}\mathbf{a}_{\left(  i,j\right)
}=\prod_{j\in J}\ \ \prod_{i\in I}\mathbf{a}_{\left(  i,j\right)  }.
\]
(In particular, all the products appearing in this equality are well-defined.)
\end{proposition}

\begin{proof}
[Proof of Proposition \ref{prop.fps.prods-mulable-rules.fubini} (sketched).]%
See Section \ref{sec.details.gf.prod}.
\end{proof}

These rules (and some similar ones, which the reader can easily invent and
prove\footnote{For example, if $\left(  \mathbf{a}_{i}\right)  _{i\in I}$ is a
multipliable family of FPSs, and if $k\in\mathbb{N}$, then $\prod_{i\in
I}\mathbf{a}_{i}^{k}=\left(  \prod_{i\in I}\mathbf{a}_{i}\right)  ^{k}$. If
the $\mathbf{a}_{i}$ are moreover invertible, then this equality also holds
for all $k\in\mathbb{Z}$.}) allow us to work with infinite products almost as
comfortably as with finite ones. Of course, we need to check that our families
are multipliable, and occasionally verify that a few other technical
requirements are met (such as multipliability of subfamilies), but usually
such verifications are straightforward and easy and can be done in one's head.

Does this justify our manipulations in Subsection \ref{subsec.gf.prod.exa}? To
some extent. We need to be careful with the telescope principle, whose
infinite analogue is rather subtle and needs some qualifications. Here is an
example of how not to use the telescope principle:
\[
\dfrac{1}{2}\cdot\dfrac{2}{2}\cdot\dfrac{2}{2}\cdot\dfrac{2}{2}\cdot\cdots
\neq1.
\]
It is tempting to argue that the infinitely many $2$'s in these fractions
cancel each other out, and yet the $1$ that remains is not the right result.
See Exercise \ref{exe.fps.prods.telescope1} \textbf{(b)} for how an infinite
telescope principle should actually look like. Anyway, our computations in
Subsection \ref{subsec.gf.prod.exa} did not truly need the telescope
principle; they can just as well be made using more fundamental
rules\footnote{To wit, we can argue as follows: We have%
\[
\prod_{i\in\mathbb{N}}\left(  1+x^{2^{i}}\right)  =\prod_{i\in\mathbb{N}%
}\dfrac{1-x^{2^{i+1}}}{1-x^{2^{i}}}=\dfrac{\prod_{i\in\mathbb{N}}\left(
1-x^{2^{i+1}}\right)  }{\prod_{i\in\mathbb{N}}\left(  1-x^{2^{i}}\right)  },
\]
where the last step used Proposition \ref{prop.fps.div-mulable} and relied on
the fact that both families $\left(  1-x^{2^{i+1}}\right)  _{i\in\mathbb{N}}$
and $\left(  1-x^{2^{i}}\right)  _{i\in\mathbb{N}}$ are multipliable (this is
important, but very easy to check in this case) and that each FPS $1-x^{2^{i}%
}$ is invertible (because its constant term is $1$). However, splitting off
the factor for $i=0$ from the product $\prod_{i\in\mathbb{N}}\left(
1-x^{2^{i}}\right)  $, we obtain
\[
\prod_{i\in\mathbb{N}}\left(  1-x^{2^{i}}\right)  =\underbrace{\left(
1-x^{2^{0}}\right)  }_{=1-x^{1}=1-x}\cdot\underbrace{\prod_{i>0}\left(
1-x^{2^{i}}\right)  }_{\substack{=\prod_{i\in\mathbb{N}}\left(  1-x^{2^{i+1}%
}\right)  \\\text{(here, we substituted }i+1\\\text{for }i\text{ in the
product)}}}=\left(  1-x\right)  \cdot\prod_{i\in\mathbb{N}}\left(
1-x^{2^{i+1}}\right)  ,
\]
so that%
\[
\dfrac{\prod_{i\in\mathbb{N}}\left(  1-x^{2^{i+1}}\right)  }{\prod
_{i\in\mathbb{N}}\left(  1-x^{2^{i}}\right)  }=\dfrac{1}{1-x}.
\]
Hence, $\prod_{i\in\mathbb{N}}\left(  1+x^{2^{i}}\right)  =\dfrac{\prod
_{i\in\mathbb{N}}\left(  1-x^{2^{i+1}}\right)  }{\prod_{i\in\mathbb{N}}\left(
1-x^{2^{i}}\right)  }=\dfrac{1}{1-x}$. So we don't need the telescope
principle to justify this equality.}.

There is one more rule that we have used in Subsection
\ref{subsec.gf.prod.exa} and have not justified yet: the equality
(\ref{eq.fps.prod.binary.prod-inf}). We will justify it next.

\subsubsection{Product rules (generalized distributive laws)}

The equality (\ref{eq.fps.prod.binary.prod-inf}) is an instance of a
\emph{product rule} -- a statement of the form \textquotedblleft a product of
sums can be expanded into one big sum\textquotedblright. The simplest product
rules are the distributive laws $a\left(  b+c\right)  =ab+ac$ and $\left(
a+b\right)  c=ac+bc$ (here, one of the sums being multiplied is a one-addend
sum); one of the next-simplest is $\left(  a+b\right)  \left(  c+d\right)
=ac+ad+bc+bd$. As far as finite sums and finite products are concerned, the
following product rule is one of the most general:\footnote{Keep in mind that
an empty Cartesian product (i.e., a Cartesian product of $0$ sets) is always a
$1$-element set; its only element is the $0$-tuple $\left(  {}\right)  $.
Thus, a sum ranging over an empty Cartesian product has exactly $1$ addend.}

\begin{proposition}
\label{prop.fps.prodrule-fin-fin}Let $L$ be a commutative ring. For every
$n\in\mathbb{N}$, let $\left[  n\right]  $ denote the set $\left\{
1,2,\ldots,n\right\}  $.

Let $n\in\mathbb{N}$. For every $i\in\left[  n\right]  $, let $p_{i,1}%
,p_{i,2},\ldots,p_{i,m_{i}}$ be finitely many elements of $L$. Then,%
\begin{equation}
\prod_{i=1}^{n}\ \ \sum_{k=1}^{m_{i}}p_{i,k}=\sum_{\left(  k_{1},k_{2}%
,\ldots,k_{n}\right)  \in\left[  m_{1}\right]  \times\left[  m_{2}\right]
\times\cdots\times\left[  m_{n}\right]  }\ \ \prod_{i=1}^{n}p_{i,k_{i}}.
\label{eq.lem.fps.prodrule-fin-fin.eq}%
\end{equation}

\end{proposition}

We can rewrite (\ref{eq.lem.fps.prodrule-fin-fin.eq}) in a less abstract way
as follows:%
\begin{align*}
&  \left(  p_{1,1}+p_{1,2}+\cdots+p_{1,m_{1}}\right)  \left(  p_{2,1}%
+p_{2,2}+\cdots+p_{2,m_{2}}\right)  \cdots\left(  p_{n,1}+p_{n,2}%
+\cdots+p_{n,m_{n}}\right) \\
&  =p_{1,1}p_{2,1}\cdots p_{n,1}+p_{1,1}p_{2,1}\cdots p_{n-1,1}p_{n,2}%
+\cdots+p_{1,m_{1}}p_{2,m_{2}}\cdots p_{n,m_{n}},
\end{align*}
where the right hand side is the sum of all $m_{1}m_{2}\cdots m_{n}$ many ways
to multiply one addend from each of the factors on the left hand side.

See \cite[solution to Exercise 6.9]{detnotes} for a formal proof of
Proposition \ref{prop.fps.prodrule-fin-fin}. (The idea is to reduce it to the
case $n=2$ by induction, then to use the discrete Fubini rule.)

Let us now move on to product rules for infinite sums and products. First, let
us extend Proposition \ref{prop.fps.prodrule-fin-fin} to a finite product of
infinite sums (which are now required to be in $K\left[  \left[  x\right]
\right]  $ in order to have a notion of summability):

\begin{proposition}
\label{prop.fps.prodrule-fin-inf}For every $n\in\mathbb{N}$, let $\left[
n\right]  $ denote the set $\left\{  1,2,\ldots,n\right\}  $.

Let $n\in\mathbb{N}$. For every $i\in\left[  n\right]  $, let $\left(
p_{i,k}\right)  _{k\in S_{i}}$ be a summable family of elements of $K\left[
\left[  x\right]  \right]  $. Then,%
\begin{equation}
\prod_{i=1}^{n}\ \ \sum_{k\in S_{i}}p_{i,k}=\sum_{\left(  k_{1},k_{2}%
,\ldots,k_{n}\right)  \in S_{1}\times S_{2}\times\cdots\times S_{n}}%
\ \ \prod_{i=1}^{n}p_{i,k_{i}}. \label{eq.lem.fps.prodrule-fin-inf.eq}%
\end{equation}
In particular, the family $\left(  \prod_{i=1}^{n}p_{i,k_{i}}\right)
_{\left(  k_{1},k_{2},\ldots,k_{n}\right)  \in S_{1}\times S_{2}\times
\cdots\times S_{n}}$ is summable.
\end{proposition}

\begin{proof}
Same method as for Proposition \ref{prop.fps.prodrule-fin-fin}, but now using
the discrete Fubini rule for infinite sums.
\end{proof}

Proposition \ref{prop.fps.prodrule-fin-inf} is rather general, but is only
concerned with \textbf{finite} products. Thus, it cannot be directly used to
justify (\ref{eq.fps.prod.binary.prod-inf}), since the product in
(\ref{eq.fps.prod.binary.prod-inf}) is infinite. Thus, we need a product rule
for infinite products of sums. Such rules are subtle and require particular
care: Not only do our sums have to be summable and our product multipliable,
but we also must avoid cases like $\left(  1-1\right)  \left(  1-1\right)
\left(  1-1\right)  \cdots$, which would produce non-summable infinite sums
when expanded (despite being multipliable). We also need to come clear about
what addends we get when we expand our products: For example, when expanding
the product
\[
\left(  1+a_{0}\right)  \left(  1+a_{1}\right)  \left(  1+a_{2}\right)
\left(  1+a_{3}\right)  \left(  1+a_{4}\right)  \cdots,
\]
we should get addends like $a_{0}\cdot1\cdot a_{2}\cdot\underbrace{1\cdot
1\cdot1\cdot\cdots}_{\text{infinitely many }1\text{s}}$ or $1\cdot a_{1}\cdot
a_{2}\cdot1\cdot a_{4}\cdot\underbrace{1\cdot1\cdot1\cdot\cdots}%
_{\text{infinitely many }1\text{s}}$, but not addends like $a_{0}\cdot
a_{1}\cdot a_{2}\cdot a_{3}\cdot\cdots$; otherwise, the right hand side of
(\ref{eq.fps.prod.binary.prod-inf}) would have to include infinite products of
$a_{i}$'s. To filter out the latter kind of addends, let us define the notion
of \textquotedblleft essentially finite\textquotedblright\ sequences or families:

\begin{definition}
\label{def.fps.prodrule.ess-fin}\textbf{(a)} A sequence $\left(  k_{1}%
,k_{2},k_{3},\ldots\right)  $ is said to be \emph{essentially finite} if all
but finitely many $i\in\left\{  1,2,3,\ldots\right\}  $ satisfy $k_{i}=0$.
\medskip

\textbf{(b)} A family $\left(  k_{i}\right)  _{i\in I}$ is said to be
\emph{essentially finite} if all but finitely many $i\in I$ satisfy $k_{i}=0$.
\end{definition}

For example, the sequence $\left(  2,4,1,0,0,0,0,\ldots\right)  $ is
essentially finite, whereas the sequence $\left(  0,1,0,1,0,1,\ldots\right)  $
(which alternates between $0$s and $1$s) is not.

Of course, Definition \ref{def.fps.prodrule.ess-fin} \textbf{(a)} is a
particular case of Definition \ref{def.fps.prodrule.ess-fin} \textbf{(b)},
since a sequence $\left(  k_{1},k_{2},k_{3},\ldots\right)  $ is the same as a
family $\left(  k_{i}\right)  _{i\in\left\{  1,2,3,\ldots\right\}  }$ indexed
by the positive integers. We also remark that Definition
\ref{def.infsum.essfin} \textbf{(a)} is a particular case of Definition
\ref{def.fps.prodrule.ess-fin} \textbf{(b)}, since a family of elements of $K$
is one particular type of family.

Now, we can finally state a version of the product rule for infinite products
of potentially infinite sums. This will help us derive
(\ref{eq.fps.prod.binary.prod-inf}) (even though the sums being multiplied in
(\ref{eq.fps.prod.binary.prod-inf}) are finite).

\begin{proposition}
\label{prop.fps.prodrule-inf-infN}Let $S_{1},S_{2},S_{3},\ldots$ be infinitely
many sets that all contain the number $0$. Set%
\[
\overline{S}=\left\{  \left(  i,k\right)  \ \mid\ i\in\left\{  1,2,3,\ldots
\right\}  \text{ and }k\in S_{i}\text{ and }k\neq0\right\}  .
\]

For any $i\in\left\{  1,2,3,\ldots\right\}  $ and any $k\in S_{i}$, let
$p_{i,k}$ be an element of $K\left[  \left[  x\right]  \right]  $. Assume
that
\begin{equation}
p_{i,0}=1\ \ \ \ \ \ \ \ \ \ \text{for any }i\in\left\{  1,2,3,\ldots\right\}
. \label{eq.prop.fps.prodrule-inf-infN.pi0=1}%
\end{equation}
Assume further that the family $\left(  p_{i,k}\right)  _{\left(  i,k\right)
\in\overline{S}}$ is summable. Then, the product $\prod_{i=1}^{\infty}%
\ \ \sum_{k\in S_{i}}p_{i,k}$ is well-defined (i.e., the family $\left(
p_{i,k}\right)  _{k\in S_{i}}$ is summable for each $i\in\left\{
1,2,3,\ldots\right\}  $, and the family $\left(  \sum_{k\in S_{i}}%
p_{i,k}\right)  _{i\in\left\{  1,2,3,\ldots\right\}  }$ is multipliable), and
we have%
\begin{equation}
\prod_{i=1}^{\infty}\ \ \sum_{k\in S_{i}}p_{i,k}=\sum_{\substack{\left(
k_{1},k_{2},k_{3},\ldots\right)  \in S_{1}\times S_{2}\times S_{3}\times
\cdots\\\text{is essentially finite}}}\ \ \prod_{i=1}^{\infty}p_{i,k_{i}}.
\label{eq.prop.fps.prodrule-inf-infN.eq}%
\end{equation}
In particular, the family $\left(  \prod_{i=1}^{\infty}p_{i,k_{i}}\right)
_{\left(  k_{1},k_{2},k_{3},\ldots\right)  \in S_{1}\times S_{2}\times
S_{3}\times\cdots\text{ is essentially finite}}$ is summable.
\end{proposition}

Note that the assumption (\ref{eq.prop.fps.prodrule-inf-infN.pi0=1}) in
Proposition \ref{prop.fps.prodrule-inf-infN} ensures that each of the sums
being multiplied contains an addend that equals $1$ (and that this addend is
named $p_{i,0}$; but this clearly does not restrict the generality of the
proposition). The equality (\ref{eq.prop.fps.prodrule-inf-infN.eq}) shows how
we can expand an infinite product of such sums. The result is a huge sum of
products (the right hand side of (\ref{eq.prop.fps.prodrule-inf-infN.eq}));
each of these products is formed by picking an addend from each sum
$\sum_{k\in S_{i}}p_{i,k}$ (just as in the finite case). The picking has to be
done in such a way that the addend $1$ gets picked from all but finitely many
of our sums (for instance, when expanding $\left(  1+x^{1}\right)  \left(
1+x^{2}\right)  \left(  1+x^{3}\right)  \cdots$, we don't want to pick the
$x^{i}$'s from all sums, as this would lead to $x^{1}x^{2}x^{3}\cdots
=$\textquotedblleft$x^{\infty}$\textquotedblright); this is why the sum on the
right hand side of (\ref{eq.prop.fps.prodrule-inf-infN.eq}) is ranging only
over the \textbf{essentially finite} sequences $\left(  k_{1},k_{2}%
,k_{3},\ldots\right)  \in S_{1}\times S_{2}\times S_{3}\times\cdots$. Finally,
the condition that the family $\left(  p_{i,k}\right)  _{\left(  i,k\right)
\in\overline{S}}$ be summable is our way of ruling out cases like $\left(
1-1\right)  \left(  1-1\right)  \left(  1-1\right)  \cdots$, which cannot be
expanded. Proposition \ref{prop.fps.prodrule-inf-infN} is not the most general
version of the product rule, but it is sufficient for many of our needs (and
we will see some more general versions below).

The rest of this subsection should probably be skipped at first reading -- it
is largely a succession of technical arguments about finite and infinite sets
in service of making rigorous what is already intuitively clear.

Before we sketch a proof of Proposition \ref{prop.fps.prodrule-inf-infN}, let
us show how (\ref{eq.fps.prod.binary.prod-inf}) can be derived from it:

\begin{proof}
[Proof of (\ref{eq.fps.prod.binary.prod-inf}) using Proposition
\ref{prop.fps.prodrule-inf-infN}.]Let $\left(  a_{0},a_{1},a_{2}%
,\ldots\right)  =\left(  a_{n}\right)  _{n\in\mathbb{N}}$ be a summable
sequence of FPSs in $K\left[  \left[  x\right]  \right]  $. We must prove the
equality (\ref{eq.fps.prod.binary.prod-inf}).

Set $S_{i}=\left\{  0,1\right\}  $ for each $i\in\left\{  1,2,3,\ldots
\right\}  $. Thus,%
\[
S_{1}\times S_{2}\times S_{3}\times\cdots=\left\{  0,1\right\}  \times\left\{
0,1\right\}  \times\left\{  0,1\right\}  \times\cdots=\left\{  0,1\right\}
^{\infty}.
\]

Define the set $\overline{S}$ as in Proposition
\ref{prop.fps.prodrule-inf-infN}. Then, $\overline{S}=\left\{  \left(
1,1\right)  ,\ \left(  2,1\right)  ,\ \left(  3,1\right)  ,\ \ldots\right\}
$. Now, set $p_{i,k}=a_{i-1}^{k}$ for each $i\in\left\{  1,2,3,\ldots\right\}
$ and each $k\in\left\{  0,1\right\}  $. Thus, the family $\left(
p_{i,k}\right)  _{\left(  i,k\right)  \in\overline{S}}$ is
summable\footnote{Indeed, this family can be rewritten as $\left(
p_{1,1},p_{2,1},p_{3,1},\ldots\right)  =\left(  a_{0}^{1},a_{1}^{1},a_{2}%
^{1},\ldots\right)  =\left(  a_{0},a_{1},a_{2},\ldots\right)  $, but we have
assumed that the latter family is summable.}, and the statement
(\ref{eq.prop.fps.prodrule-inf-infN.pi0=1}) holds as well (indeed, we have
$p_{i,0}=a_{i-1}^{0}=1$ for any $i\in\left\{  1,2,3,\ldots\right\}  $). Hence,
Proposition \ref{prop.fps.prodrule-inf-infN} can be applied.

Moreover, for each $i\in\left\{  1,2,3,\ldots\right\}  $, we have
$S_{i}=\left\{  0,1\right\}  $. Thus, for each $i\in\left\{  1,2,3,\ldots
\right\}  $, we have%
\[
\sum_{k\in S_{i}}p_{i,k}=\sum_{k\in\left\{  0,1\right\}  }p_{i,k}%
=\underbrace{p_{i,0}}_{=a_{i-1}^{0}=1}+\underbrace{p_{i,1}}_{=a_{i-1}%
^{1}=a_{i-1}}=1+a_{i-1}.
\]
Hence,%
\[
\prod_{i=1}^{\infty}\ \ \underbrace{\sum_{k\in S_{i}}p_{i,k}}_{=1+a_{i-1}%
}=\prod_{i=1}^{\infty}\left(  1+a_{i-1}\right)  =\prod_{i\in\mathbb{N}}\left(
1+a_{i}\right)
\]
(here, we have substituted $i$ for $i-1$ in the product). Thus,%
\begin{align}
\prod_{i\in\mathbb{N}}\left(  1+a_{i}\right)   &  =\prod_{i=1}^{\infty
}\ \ \sum_{k\in S_{i}}p_{i,k}\nonumber\\
&  =\underbrace{\sum_{\substack{\left(  k_{1},k_{2},k_{3},\ldots\right)  \in
S_{1}\times S_{2}\times S_{3}\times\cdots\\\text{is essentially finite}}%
}}_{\substack{=\sum_{\substack{\left(  k_{1},k_{2},k_{3},\ldots\right)
\in\left\{  0,1\right\}  ^{\infty}\\\text{is essentially finite}%
}}\\\text{(since }S_{1}\times S_{2}\times S_{3}\times\cdots=\left\{
0,1\right\}  ^{\infty}\text{)}}}\ \ \prod_{i=1}^{\infty}\underbrace{p_{i,k_{i}%
}}_{=a_{i-1}^{k_{i}}}\ \ \ \ \ \ \ \ \ \ \left(  \text{by
(\ref{eq.prop.fps.prodrule-inf-infN.eq})}\right) \nonumber\\
&  =\sum_{\substack{\left(  k_{1},k_{2},k_{3},\ldots\right)  \in\left\{
0,1\right\}  ^{\infty}\\\text{is essentially finite}}}\ \ \prod_{i=1}^{\infty
}a_{i-1}^{k_{i}}=\sum_{\substack{\left(  k_{0},k_{1},k_{2},\ldots\right)
\in\left\{  0,1\right\}  ^{\infty}\\\text{is essentially finite}%
}}\ \ \underbrace{\prod_{i=1}^{\infty}a_{i-1}^{k_{i-1}}}_{=\prod
_{i\in\mathbb{N}}a_{i}^{k_{i}}}\nonumber\\
&  \ \ \ \ \ \ \ \ \ \ \ \ \ \ \ \ \ \ \ \ \left(
\begin{array}
[c]{c}%
\text{here, we have renamed the summation}\\
\text{index }\left(  k_{1},k_{2},k_{3},\ldots\right)  \text{ as }\left(
k_{0},k_{1},k_{2},\ldots\right)
\end{array}
\right) \nonumber\\
&  =\sum_{\substack{\left(  k_{0},k_{1},k_{2},\ldots\right)  \in\left\{
0,1\right\}  ^{\infty}\\\text{is essentially finite}}}\ \ \prod_{i\in
\mathbb{N}}a_{i}^{k_{i}}. \label{pf.eq.fps.prod.binary.prod-inf.5}%
\end{align}

Now, an essentially finite sequence $\left(  k_{0},k_{1},k_{2},\ldots\right)
\in\left\{  0,1\right\}  ^{\infty}$ is just an infinite sequence of $0$s and
$1$s that contains only finitely many $1$s. Thus, such a sequence is uniquely
determined by the positions of the $1$s in it, and these positions form a
finite set, so they can be uniquely labeled as $i_{1},i_{2},\ldots,i_{k}$ with
$i_{1}<i_{2}<\cdots<i_{k}$. To put this more formally: There is a bijection%
\begin{align*}
&  \text{from }\left\{  \text{essentially finite sequences }\left(
k_{0},k_{1},k_{2},\ldots\right)  \in\left\{  0,1\right\}  ^{\infty}\right\} \\
&  \text{to }\left\{  \text{finite lists }\left(  i_{1},i_{2},\ldots
,i_{k}\right)  \text{ of nonnegative integers such that }i_{1}<i_{2}%
<\cdots<i_{k}\right\}
\end{align*}
that sends each sequence $\left(  k_{0},k_{1},k_{2},\ldots\right)  \in\left\{
0,1\right\}  ^{\infty}$ to the list of all $i\in\mathbb{N}$ satisfying
$k_{i}=1$ (written in increasing order). Furthermore, if this bijection sends
an essentially finite sequence $\left(  k_{0},k_{1},k_{2},\ldots\right)
\in\left\{  0,1\right\}  ^{\infty}$ to a finite list $\left(  i_{1}%
,i_{2},\ldots,i_{k}\right)  $, then
\begin{align*}
\prod_{i\in\mathbb{N}}a_{i}^{k_{i}}  &  =\left(  \prod_{\substack{i\in
\mathbb{N};\\k_{i}=0}}\ \ \underbrace{a_{i}^{k_{i}}}_{=a_{i}^{0}=1}\right)
\left(  \prod_{\substack{i\in\mathbb{N};\\k_{i}=1}}\ \ \underbrace{a_{i}%
^{k_{i}}}_{=a_{i}^{1}=a_{i}}\right)  \ \ \ \ \ \ \ \ \ \ \left(  \text{since
each }k_{i}\text{ is either }0\text{ or }1\right) \\
&  =\underbrace{\left(  \prod_{\substack{i\in\mathbb{N};\\k_{i}=0}}1\right)
}_{=1}\left(  \prod_{\substack{i\in\mathbb{N};\\k_{i}=1}}a_{i}\right)
=\prod_{\substack{i\in\mathbb{N};\\k_{i}=1}}a_{i}=a_{i_{1}}a_{i_{2}}\cdots
a_{i_{k}}%
\end{align*}
(since $\left(  i_{1},i_{2},\ldots,i_{k}\right)  $ is the list of all
$i\in\mathbb{N}$ satisfying $k_{i}=1$). Thus, we can use this bijection to
re-index the sum on the right hand side of
(\ref{pf.eq.fps.prod.binary.prod-inf.5}); we obtain%
\[
\sum_{\substack{\left(  k_{0},k_{1},k_{2},\ldots\right)  \in\left\{
0,1\right\}  ^{\infty}\\\text{is essentially finite}}}\ \ \prod_{i\in
\mathbb{N}}a_{i}^{k_{i}}=\sum_{i_{1}<i_{2}<\cdots<i_{k}}a_{i_{1}}a_{i_{2}%
}\cdots a_{i_{k}}.
\]
Thus, (\ref{pf.eq.fps.prod.binary.prod-inf.5}) rewrites as%
\[
\prod_{i\in\mathbb{N}}\left(  1+a_{i}\right)  =\sum_{i_{1}<i_{2}<\cdots<i_{k}%
}a_{i_{1}}a_{i_{2}}\cdots a_{i_{k}}.
\]
This proves (\ref{eq.fps.prod.binary.prod-inf}).
\end{proof}

We note that (\ref{eq.fps.prod.binary.prod-inf}) can be rewritten as%
\begin{equation}
\prod_{i\in\mathbb{N}}\left(  1+a_{i}\right)  =\sum_{\substack{J\text{ is a
finite}\\\text{subset of }\mathbb{N}}}\ \ \prod_{i\in J}a_{i}.
\label{eq.fps.prod.binary.prod-inf2}%
\end{equation}
Indeed, the right hand side of this equality is precisely the right hand side
of (\ref{eq.fps.prod.binary.prod-inf}).

We will prove Proposition \ref{prop.fps.prodrule-inf-infN} soon. First,
however, let us generalize Proposition \ref{prop.fps.prodrule-inf-infN} to
products indexed by arbitrary sets instead of $\left\{  1,2,3,\ldots\right\}
$:

\begin{proposition}
\label{prop.fps.prodrule-inf-inf}Let $I$ be a set. For any $i\in I$, let
$S_{i}$ be a set that contains the number $0$. Set%
\[
\overline{S}=\left\{  \left(  i,k\right)  \ \mid\ i\in I\text{ and }k\in
S_{i}\text{ and }k\neq0\right\}  .
\]

For any $i\in I$ and any $k\in S_{i}$, let $p_{i,k}$ be an element of
$K\left[  \left[  x\right]  \right]  $. Assume that
\begin{equation}
p_{i,0}=1\ \ \ \ \ \ \ \ \ \ \text{for any }i\in I.
\label{eq.prop.fps.prodrule-inf-inf.pi0=1}%
\end{equation}
Assume further that the family $\left(  p_{i,k}\right)  _{\left(  i,k\right)
\in\overline{S}}$ is summable. Then, the product $\prod_{i\in I}\ \ \sum_{k\in
S_{i}}p_{i,k}$ is well-defined (i.e., the family $\left(  p_{i,k}\right)
_{k\in S_{i}}$ is summable for each $i\in I$, and the family $\left(
\sum_{k\in S_{i}}p_{i,k}\right)  _{i\in I}$ is multipliable), and we have%
\begin{equation}
\prod_{i\in I}\ \ \sum_{k\in S_{i}}p_{i,k}=\sum_{\substack{\left(
k_{i}\right)  _{i\in I}\in\prod_{i\in I}S_{i}\\\text{is essentially finite}%
}}\ \ \prod_{i\in I}p_{i,k_{i}}. \label{eq.prop.fps.prodrule-inf-inf.eq}%
\end{equation}
In particular, the family $\left(  \prod_{i\in I}p_{i,k_{i}}\right)  _{\left(
k_{i}\right)  _{i\in I}\in\prod_{i\in I}S_{i}\text{ is essentially finite}}$
is summable.
\end{proposition}

Clearly, Proposition \ref{prop.fps.prodrule-inf-infN} is the particular case
of Proposition \ref{prop.fps.prodrule-inf-inf} for $I=\left\{  1,2,3,\ldots
\right\}  $. The proof of Proposition \ref{prop.fps.prodrule-inf-inf} is a
long-winded reduction to the finite case; nothing substantial is going on in
it. Thus, I recommend skipping it unless specifically interested. For the sake
of convenience, before proving Proposition \ref{prop.fps.prodrule-inf-inf},
let me restate Proposition \ref{prop.fps.prodrule-fin-inf} in a form that
makes it particularly easy to use:

\begin{proposition}
\label{prop.fps.prodrule-fin-infJ}Let $N$ be a finite set. For any $i\in N$,
let $\left(  p_{i,k}\right)  _{k\in S_{i}}$ be a summable family of elements
of $K\left[  \left[  x\right]  \right]  $. Then,%
\begin{equation}
\prod_{i\in N}\ \ \sum_{k\in S_{i}}p_{i,k}=\sum_{\left(  k_{i}\right)  _{i\in
N}\in\prod_{i\in N}S_{i}}\ \ \prod_{i\in N}p_{i,k_{i}}.
\label{eq.lem.fps.prodrule-fin-infJ.eq}%
\end{equation}
In particular, the family $\left(  \prod_{i\in N}p_{i,k_{i}}\right)  _{\left(
k_{i}\right)  _{i\in N}\in\prod_{i\in N}S_{i}}$ is summable.
\end{proposition}

\begin{proof}
[Proof of Proposition \ref{prop.fps.prodrule-fin-infJ}.]This is just
Proposition \ref{prop.fps.prodrule-fin-inf}, with the indexing set $\left\{
1,2,\ldots,n\right\}  $ replaced by $N$. It can be proved by reindexing the
products (or directly by induction on $\left\vert N\right\vert $).
\end{proof}

Let me furthermore state an infinite version of Lemma
\ref{lem.fps.prod.irlv.fin}, which will be used in the proof of Proposition
\ref{prop.fps.prodrule-inf-inf} given below:

\begin{lemma}
\label{lem.fps.prod.irlv.inf}Let $a\in K\left[  \left[  x\right]  \right]  $
be an FPS. Let $\left(  f_{i}\right)  _{i\in J}\in K\left[  \left[  x\right]
\right]  ^{J}$ be a summable family of FPSs. Let $n\in\mathbb{N}$. Assume that
each $i\in J$ satisfies
\begin{equation}
\left[  x^{m}\right]  \left(  f_{i}\right)  =0\ \ \ \ \ \ \ \ \ \ \text{for
each }m\in\left\{  0,1,\ldots,n\right\}  .
\label{eq.lem.fps.prod.irlv.inf.ass}%
\end{equation}
Then,
\[
\left[  x^{m}\right]  \left(  a\prod_{i\in J}\left(  1+f_{i}\right)  \right)
=\left[  x^{m}\right]  a\ \ \ \ \ \ \ \ \ \ \text{for each }m\in\left\{
0,1,\ldots,n\right\}  .
\]

\end{lemma}

\begin{proof}
[Proof of Lemma \ref{lem.fps.prod.irlv.inf} (sketched).]The idea here is to
argue that the first $n+1$ coefficients of $a\prod_{i\in J}\left(
1+f_{i}\right)  $ agree with those of $a\prod_{i\in M}\left(  1+f_{i}\right)
$ for some finite subset $M$ of $J$. Then, apply Lemma
\ref{lem.fps.prod.irlv.fin} to this subset $M$. The details of this argument
can be found in Section \ref{sec.details.gf.prod}.
\end{proof}

With this lemma proved, we have all necessary prerequisites for the proof of
Proposition \ref{prop.fps.prodrule-inf-inf}. The proof, however, is rather
long due to the bookkeeping required, and therefore has been banished to the
appendix (Section \ref{sec.details.gf.prod}, to be specific).

We notice that Proposition \ref{prop.fps.prodrule-inf-inf} (and thus also
Proposition \ref{prop.fps.prodrule-inf-infN}) can be generalized slightly by
lifting the requirement that all sets $S_{i}$ contain $0$ (this means that the
sums being multiplied no longer need to contain $1$ as an addend). See
Exercise \ref{exe.fps.prodrule-inf-inf-0} for this generalization. This lets
us expand products such as $x\left(  1+x^{1}\right)  \left(  1+x^{2}\right)
\left(  1+x^{3}\right)  \cdots$ (but of course, we could just as well expand
this particular product by splitting off the $x$ factor and applying
Proposition \ref{prop.fps.prodrule-inf-infN}).

\subsubsection{Another example}

As another bit of practice with infinite products of FPSs, let us prove the
following identity by Euler (\cite[\S 326]{Euler48}):

\begin{proposition}
[Euler]\label{prop.gf.prod.euler-odd}We have%
\[
\prod_{i>0}\left(  1-x^{2i-1}\right)  ^{-1}=\prod_{k>0}\left(  1+x^{k}\right)
.
\]

\end{proposition}

\begin{proof}
First, it is easy to see that the families $\left(  1+x^{k}\right)  _{k>0}$
and $\left(  1-x^{k}\right)  _{k>0}$ and $\left(  1-x^{2k}\right)  _{k>0}$ and
$\left(  1-x^{2i-1}\right)  _{i>0}$ are multipliable\footnote{Indeed, this
follows from Theorem \ref{thm.fps.1+f-mulable}, since the families $\left(
x^{k}\right)  _{k>0}$ and $\left(  -x^{k}\right)  _{k>0}$ and $\left(
-x^{2k}\right)  _{k>0}$ and $\left(  -x^{2i-1}\right)  _{i>0}$ are summable.}.
This shows that the product on the right hand side of Proposition
\ref{prop.gf.prod.euler-odd} is well-defined.

Proposition \ref{prop.fps.div-mulable} \textbf{(a)} shows that the family
$\left(  \dfrac{1}{1-x^{2i-1}}\right)  _{i>0}$ is multipliable (since the
families $\left(  1\right)  _{i>0}$ and $\left(  1-x^{2i-1}\right)  _{i>0}$
are multipliable, and since the FPS $1-x^{2i-1}$ is invertible for each
$i>0$). In other words, the family $\left(  \left(  1-x^{2i-1}\right)
^{-1}\right)  _{i>0}$ is multipliable. Thus, the product on the left hand side
of Proposition \ref{prop.gf.prod.euler-odd} is well-defined.

For each $k>0$, we have $1+x^{k}=\dfrac{1-x^{2k}}{1-x^{k}}$ (since
$1-x^{2k}=1-\left(  x^{k}\right)  ^{2}=\left(  1-x^{k}\right)  \left(
1+x^{k}\right)  $). Thus,%
\begin{align*}
\prod_{k>0}\underbrace{\left(  1+x^{k}\right)  }_{=\dfrac{1-x^{2k}}{1-x^{k}}}
&  =\prod_{k>0}\dfrac{1-x^{2k}}{1-x^{k}}\\
&  =\dfrac{\prod_{k>0}\left(  1-x^{2k}\right)  }{\prod_{k>0}\left(
1-x^{k}\right)  }\ \ \ \ \ \ \ \ \ \ \left(  \text{by Proposition
\ref{prop.fps.div-mulable} \textbf{(b)}}\right) \\
&  =\dfrac{\left(  1-x^{2}\right)  \left(  1-x^{4}\right)  \left(
1-x^{6}\right)  \left(  1-x^{8}\right)  \cdots}{\left(  1-x^{1}\right)
\left(  1-x^{2}\right)  \left(  1-x^{3}\right)  \left(  1-x^{4}\right)
\cdots}\\
&  =\dfrac{1}{\left(  1-x^{1}\right)  \left(  1-x^{3}\right)  \left(
1-x^{5}\right)  \left(  1-x^{7}\right)  \cdots}\\
&  \ \ \ \ \ \ \ \ \ \ \ \ \ \ \ \ \ \ \ \ \left(
\begin{array}
[c]{c}%
\text{since all }1-x^{i}\text{ factors with }i\text{ even}\\
\text{cancel out (see below for details)}%
\end{array}
\right) \\
&  =\dfrac{1}{1-x^{1}}\cdot\dfrac{1}{1-x^{3}}\cdot\dfrac{1}{1-x^{5}}%
\cdot\dfrac{1}{1-x^{7}}\cdot\cdots\\
&  \ \ \ \ \ \ \ \ \ \ \ \ \ \ \ \ \ \ \ \ \left(  \text{by Proposition
\ref{prop.fps.div-mulable} \textbf{(b)}}\right) \\
&  =\prod_{i>0}\dfrac{1}{1-x^{2i-1}}=\prod_{i>0}\left(  1-x^{2i-1}\right)
^{-1},
\end{align*}
which is precisely the claim of Proposition \ref{prop.gf.prod.euler-odd}.

The step where we cancelled out the $1-x^{i}$ factors with $i$ even can be
justified as follows: We know that the family $\left(  1-x^{k}\right)  _{k>0}$
as well as its subfamilies $\left(  1-x^{k}\right)  _{k>0\text{ is even}}$ and
$\left(  1-x^{k}\right)  _{k>0\text{ are odd}}$ are multipliable. Thus,
Proposition \ref{prop.fps.union-mulable} \textbf{(b)} yields%
\[
\prod_{k>0}\left(  1-x^{k}\right)  =\left(  \prod_{\substack{k>0\\\text{is
even}}}\left(  1-x^{k}\right)  \right)  \cdot\left(  \prod
_{\substack{k>0\\\text{is odd}}}\left(  1-x^{k}\right)  \right)  .
\]
In other words,%
\begin{align*}
&  \left(  1-x^{1}\right)  \left(  1-x^{2}\right)  \left(  1-x^{3}\right)
\left(  1-x^{4}\right)  \cdots\\
&  =\left(  \left(  1-x^{2}\right)  \left(  1-x^{4}\right)  \left(
1-x^{6}\right)  \left(  1-x^{8}\right)  \cdots\right) \\
&  \ \ \ \ \ \ \ \ \ \ \cdot\left(  \left(  1-x^{1}\right)  \left(
1-x^{3}\right)  \left(  1-x^{5}\right)  \left(  1-x^{7}\right)  \cdots\right)
.
\end{align*}
Hence,%
\[
\dfrac{\left(  1-x^{1}\right)  \left(  1-x^{2}\right)  \left(  1-x^{3}\right)
\left(  1-x^{4}\right)  \cdots}{\left(  1-x^{2}\right)  \left(  1-x^{4}%
\right)  \left(  1-x^{6}\right)  \left(  1-x^{8}\right)  \cdots}=\left(
1-x^{1}\right)  \left(  1-x^{3}\right)  \left(  1-x^{5}\right)  \left(
1-x^{7}\right)  \cdots.
\]
Taking reciprocals, we thus obtain%
\[
\dfrac{\left(  1-x^{2}\right)  \left(  1-x^{4}\right)  \left(  1-x^{6}\right)
\left(  1-x^{8}\right)  \cdots}{\left(  1-x^{1}\right)  \left(  1-x^{2}%
\right)  \left(  1-x^{3}\right)  \left(  1-x^{4}\right)  \cdots}=\dfrac
{1}{\left(  1-x^{1}\right)  \left(  1-x^{3}\right)  \left(  1-x^{5}\right)
\left(  1-x^{7}\right)  \cdots},
\]
just as we claimed. Thus, the proof of Proposition
\ref{prop.gf.prod.euler-odd} is complete.
\end{proof}

Let us try to interpret Proposition \ref{prop.gf.prod.euler-odd}
combinatorially, by expanding both sides.

First, we use (\ref{eq.fps.prod.binary.prod-inf}) to expand the right hand
side:%
\begin{align}
\prod_{k>0}\left(  1+x^{k}\right)   &  =\left(  1+x^{1}\right)  \left(
1+x^{2}\right)  \left(  1+x^{3}\right)  \left(  1+x^{4}\right)  \cdots
\nonumber\\
&  =\sum_{\substack{i_{1},i_{2},\ldots,i_{k}\in\left\{  1,2,3,\ldots\right\}
;\\i_{1}<i_{2}<\cdots<i_{k}}}\underbrace{x^{i_{1}}x^{i_{2}}\cdots x^{i_{k}}%
}_{=x^{i_{1}+i_{2}+\cdots+i_{k}}}\ \ \ \ \ \ \ \ \ \ \left(  \text{by
(\ref{eq.fps.prod.binary.prod-inf})}\right) \nonumber\\
&  =\sum_{\substack{i_{1},i_{2},\ldots,i_{k}\in\left\{  1,2,3,\ldots\right\}
;\\i_{1}<i_{2}<\cdots<i_{k}}}x^{i_{1}+i_{2}+\cdots+i_{k}}\nonumber\\
&  =\sum_{n\in\mathbb{N}}d_{n}x^{n}, \label{eq.prop.gf.prod.euler-odd.d}%
\end{align}
where $d_{n}$ is the \# of all strictly increasing tuples $\left(  i_{1}%
<i_{2}<\cdots<i_{k}\right)  $ of positive integers such that $n=i_{1}%
+i_{2}+\cdots+i_{k}$. We can rewrite this definition as follows: $d_{n}$ is
the \# of ways to write $n$ as a sum of distinct positive integers, with no
regard for the order\footnote{\textquotedblleft No regard for the
order\textquotedblright\ means that, for example, $3+4+1$ and $1+3+4$ count as
the same way of writing $8$ as a sum of distinct integers.}.

Next, let us expand the left hand side: For each positive integer $i$, we have%
\begin{align*}
\left(  1-x^{2i-1}\right)  ^{-1}  &  =1+x^{2i-1}+\left(  x^{2i-1}\right)
^{2}+\left(  x^{2i-1}\right)  ^{3}+\cdots\\
&  \ \ \ \ \ \ \ \ \ \ \ \ \ \ \ \ \ \ \ \ \left(
\begin{array}
[c]{c}%
\text{by substituting }x^{2i-1}\text{ for }x\text{ in the}\\
\text{equality }\left(  1-x\right)  ^{-1}=1+x+x^{2}+x^{3}+\cdots
\end{array}
\right) \\
&  =1+x^{2i-1}+x^{2\left(  2i-1\right)  }+x^{3\left(  2i-1\right)  }%
+\cdots=\sum_{k\in\mathbb{N}}x^{k\left(  2i-1\right)  }.
\end{align*}
Multiplying these equalities over all $i>0$, we find%
\begin{align}
&  \prod_{i>0}\left(  1-x^{2i-1}\right)  ^{-1}\nonumber\\
&  =\prod_{i>0}\ \ \sum_{k\in\mathbb{N}}x^{k\left(  2i-1\right)  }\nonumber\\
&  =\sum_{\substack{\left(  k_{1},k_{2},k_{3},\ldots\right)  \in
\mathbb{N}^{\infty}\\\text{is essentially finite}}}\underbrace{\prod
_{i>0}x^{k_{i}\left(  2i-1\right)  }}_{\substack{=x^{k_{1}\cdot1}x^{k_{2}%
\cdot3}x^{k_{3}\cdot5}\cdots\\=x^{k_{1}\cdot1+k_{2}\cdot3+k_{3}\cdot5+\cdots}%
}}\ \ \ \ \ \ \ \ \ \ \left(
\begin{array}
[c]{c}%
\text{by Proposition \ref{prop.fps.prodrule-inf-infN},}\\
\text{applied to }S_{i}=\mathbb{N}\text{ and }p_{i,k}=x^{k\left(  2i-1\right)
}%
\end{array}
\right) \nonumber\\
&  =\sum_{\substack{\left(  k_{1},k_{2},k_{3},\ldots\right)  \in
\mathbb{N}^{\infty}\\\text{is essentially finite}}}x^{k_{1}\cdot1+k_{2}%
\cdot3+k_{3}\cdot5+\cdots}\nonumber\\
&  =\sum_{n\in\mathbb{N}}o_{n}x^{n}, \label{eq.prop.gf.prod.euler-odd.o}%
\end{align}
where $o_{n}$ is the \# of all essentially finite sequences $\left(
k_{1},k_{2},k_{3},\ldots\right)  \in\mathbb{N}^{\infty}$ such that $k_{1}%
\cdot1+k_{2}\cdot3+k_{3}\cdot5+\cdots=n$. I claim that $o_{n}$ is the \# of
ways to write $n$ as a sum of (not necessarily distinct) odd positive
integers, without regard to the order. (Why? Because if we can write $n$ as a
sum of $k_{1}$ many $1$s, $k_{2}$ many $3$s, $k_{3}$ many $5$s and so on, then
$\left(  k_{1},k_{2},k_{3},\ldots\right)  \in\mathbb{N}^{\infty}$ is an
essentially finite sequence such that $k_{1}\cdot1+k_{2}\cdot3+k_{3}%
\cdot5+\cdots=n$; and conversely, if $\left(  k_{1},k_{2},k_{3},\ldots\right)
\in\mathbb{N}^{\infty}$ is such a sequence, then $n$ can be written as a sum
of $k_{1}$ many $1$s, $k_{2}$ many $3$s, $k_{3}$ many $5$s and so on.)

Proposition \ref{prop.gf.prod.euler-odd} tells us that $\prod_{k>0}\left(
1+x^{k}\right)  =\prod_{i>0}\left(  1-x^{2i-1}\right)  ^{-1}$. In view of
(\ref{eq.prop.gf.prod.euler-odd.d}) and (\ref{eq.prop.gf.prod.euler-odd.o}),
we can rewrite this as
\[
\sum_{n\in\mathbb{N}}d_{n}x^{n}=\sum_{n\in\mathbb{N}}o_{n}x^{n}.
\]
Therefore, $d_{n}=o_{n}$ for each $n\in\mathbb{N}$. Thus, we have proven the
following purely combinatorial statement:

\begin{theorem}
[Euler]\label{thm.gf.prod.euler-comb}Let $n\in\mathbb{N}$. Then, $d_{n}=o_{n}%
$, where

\begin{itemize}
\item $d_{n}$ is the \# of ways to write $n$ as a sum of distinct positive
integers, without regard to the order;

\item $o_{n}$ is the \# of ways to write $n$ as a sum of (not necessarily
distinct) odd positive integers, without regard to the order.
\end{itemize}
\end{theorem}

\begin{example}
Let $n=6$. Then, $d_{n}=4$, because the ways to write $n=6$ as a sum of
distinct positive integers are%
\[
6,\ \ \ \ \ \ \ \ \ \ 1+5,\ \ \ \ \ \ \ \ \ \ 2+4,\ \ \ \ \ \ \ \ \ \ 1+2+3
\]
(don't forget the first of these ways, trivial as it looks!). On the other
hand, $o_{n}=4$, because the ways to write $n=6$ as a sum of odd positive
integers are%
\[
1+5,\ \ \ \ \ \ \ \ \ \ 3+3,\ \ \ \ \ \ \ \ \ \ 3+1+1+1,\ \ \ \ \ \ \ \ \ \ 1+1+1+1+1+1.
\]

\end{example}

We will soon learn a bijective proof of Theorem \ref{thm.gf.prod.euler-comb}
as well (see the Second proof of Theorem \ref{thm.pars.odd-dist-equal} below).

\subsubsection{Infinite products and substitution}

Proposition \ref{prop.fps.subs.rules} \textbf{(h)} has an analogue for
products instead of sums:

\begin{proposition}
\label{prop.fps.subs.rule-infprod}If $\left(  f_{i}\right)  _{i\in I}\in
K\left[  \left[  x\right]  \right]  ^{I}$ is a multipliable family of FPSs,
and if $g\in K\left[  \left[  x\right]  \right]  $ is an FPS satisfying
$\left[  x^{0}\right]  g=0$, then the family $\left(  f_{i}\circ g\right)
_{i\in I}\in K\left[  \left[  x\right]  \right]  ^{I}$ is multipliable as well
and we have $\left(  \prod_{i\in I}f_{i}\right)  \circ g=\prod_{i\in I}\left(
f_{i}\circ g\right)  $.
\end{proposition}

\begin{proof}
[Proof of Proposition \ref{prop.fps.subs.rule-infprod} (sketched).]See Section
\ref{sec.details.gf.prod}.
\end{proof}

\subsubsection{Exponentials, logarithms and infinite products}

Recall Definition \ref{def.fps.Exp-Log-maps}. As we saw in Lemma
\ref{lem.fps.Exp-Log-additive}, the exponential map $\operatorname*{Exp}$ and
the logarithm map $\operatorname*{Log}$ can be used to convert sums into
products and vice versa (under certain conditions). This also extends to
infinite sums and infinite products:

\begin{proposition}
\label{prop.fps.Exp-Log-infsum}Assume that $K$ is a commutative $\mathbb{Q}%
$-algebra. Let $\left(  f_{i}\right)  _{i\in I}\in K\left[  \left[  x\right]
\right]  ^{I}$ be a summable family of FPSs in $K\left[  \left[  x\right]
\right]  _{0}$. Then, $\left(  \operatorname*{Exp}f_{i}\right)  _{i\in I}$ is
a multipliable family of FPSs in $K\left[  \left[  x\right]  \right]  _{1}$,
and we have $\sum_{i\in I}f_{i}\in K\left[  \left[  x\right]  \right]  _{0}$
and%
\[
\operatorname*{Exp}\left(  \sum_{i\in I}f_{i}\right)  =\prod_{i\in I}\left(
\operatorname*{Exp}f_{i}\right)  .
\]

\end{proposition}

\begin{proposition}
\label{prop.fps.Exp-Log-infprod}Assume that $K$ is a commutative $\mathbb{Q}%
$-algebra. Let $\left(  f_{i}\right)  _{i\in I}\in K\left[  \left[  x\right]
\right]  ^{I}$ be a multipliable family of FPSs in $K\left[  \left[  x\right]
\right]  _{1}$. Then, $\left(  \operatorname*{Log}f_{i}\right)  _{i\in I}$ is
a summable family of FPSs in $K\left[  \left[  x\right]  \right]  _{0}$, and
we have $\prod\limits_{i\in I}f_{i}\in K\left[  \left[  x\right]  \right]
_{1}$ and%
\[
\operatorname*{Log}\left(  \prod\limits_{i\in I}f_{i}\right)  =\sum_{i\in
I}\left(  \operatorname*{Log}f_{i}\right)  .
\]

\end{proposition}

We leave the proofs of these two propositions as an exercise (Exercise
\ref{exe.fps.Exp-Log-infsum}).

\subsection{\label{sec.gf.weighted-set}The generating function of a weighted
set}

So far, we have built a theory of FPSs, and we have occasionally applied them
to counting problems. However, so far, the latter applications were either
limited to counting tuples of integers (because such tuples naturally appear
as indexes when we multiply several sums) or somewhat indirect, mediated by
integer sequences (i.e., we could not immediately transform our counting
problem into an FPS but rather had to first define an integer sequence and
then take its generating function). For instance, our above proof of Theorem
\ref{thm.gf.prod.euler-comb} was an application of the former type, while
Example 2 in Section \ref{sec.gf.exas} was one of the latter type.

I shall now explain a more direct approach to transforming counting problems
into FPSs. This approach is the theory of so-called \textquotedblleft%
\emph{combinatorial classes}\textquotedblright\ or \textquotedblleft%
\emph{finite-type weighted sets}\textquotedblright, and can be used
(particularly in some of its more advanced variants) to prove deep results.
However, I will not use it much in these notes, so I will only give a brief
introduction. Much more can be learned from \cite[Part A]{FlaSed09};
introductions can also be found in \cite[\S 3.3--\S 3.4]{Fink17}, \cite[Parts
1 and 2]{Melcze24} and \cite[Chapter 6]{Loehr-BC}. The theory is not magic,
and in principle can always be eschewed (any proof using it can be rewritten
without its use), but it offers the major advantages of systematizing the
combinatorial objects being counted and making the proofs more insightful. In
a sense, it allows us to perform certain computations not just with the
generating functions of the objects, but with the objects themselves -- a
middle ground between algebra and combinatorics.

I will now outline the beginnings of this theory, following Fink's
introduction \cite[\S 3.3--\S 3.4]{Fink17} (but speaking of \emph{finite-type
weighted sets} where \cite{Fink17}, \cite{Melcze24} and \cite{FlaSed09} speak
of \textquotedblleft combinatorial classes\textquotedblright).

\subsubsection{The theory}

\begin{definition}
\label{def.gf-ws.weighted-sets}\textbf{(a)} A \emph{weighted set} is a set $A$
equipped with a function $w:A\rightarrow\mathbb{N}$, which is called the
\emph{weight function} of this weighted set. For each $a\in A$, the value
$w\left(  a\right)  $ is denoted $\left\vert a\right\vert $ and is called the
\emph{weight} of $a$ (in our weighted set). \medskip

Usually, instead of explicitly specifying the weight function $w$ as a
function, we will simply specify the weight $\left\vert a\right\vert $ for
each $a\in A$. The weighted set consisting of the set $A$ and the weight
function $w$ will be called $\left(  A,w\right)  $ or simply $A$ when the
weight function $w$ is clear from the context. \medskip

\textbf{(b)} A weighted set $A$ is said to be \emph{finite-type} if for each
$n\in\mathbb{N}$, there are only finitely many $a\in A$ having weight
$\left\vert a\right\vert =n$. \medskip

\textbf{(c)} If $A$ is a finite-type weighted set, then the \emph{weight
generating function} of $A$ is defined to be the FPS%
\[
\sum_{a\in A}x^{\left\vert a\right\vert }=\sum_{n\in\mathbb{N}}\left(
\text{\# of }a\in A\text{ having weight }n\right)  \cdot x^{n}\in
\mathbb{Z}\left[  \left[  x\right]  \right]  .
\]
This FPS is denoted by $\overline{A}$. \medskip

\textbf{(d)} An \emph{isomorphism} between two weighted sets $A$ and $B$ means
a bijection $\rho:A\rightarrow B$ that preserves the weight (i.e., each $a\in
A$ satisfies $\left\vert \rho\left(  a\right)  \right\vert =\left\vert
a\right\vert $). \medskip

\textbf{(e)} We say that two weighted sets $A$ and $B$ are \emph{isomorphic}
(this is written $A\cong B$) if there exists an isomorphism between $A$ and
$B$.
\end{definition}

\begin{example}
\label{exa.ws.bin-string1}Let $B$ be the weighted set of all \emph{binary
strings}, i.e., finite tuples consisting of $0$s and $1$s. Thus,%
\[
B=\left\{  \left(  {}\right)  ,\ \left(  0\right)  ,\ \left(  1\right)
,\ \left(  0,0\right)  ,\ \left(  0,1\right)  ,\ \left(  1,0\right)
,\ \left(  1,1\right)  ,\ \left(  0,0,0\right)  ,\ \left(  0,0,1\right)
,\ \ldots\right\}  .
\]
The weight of an $n$-tuple is defined to be $n$. In other words, the weight of
a binary string is the length of this string. This weighted set $B$ is
finite-type, since for each $n\in\mathbb{N}$, there are only finitely many
binary strings of length $n$ (namely, there are $2^{n}$ such strings). The
weight generating function of $B$ is%
\[
\overline{B}=\sum_{n\in\mathbb{N}}\underbrace{\left(  \text{\# of }a\in
B\text{ having weight }n\right)  }_{\substack{=\left(  \text{\# of binary
strings of length }n\right)  \\=2^{n}}}\cdot\,x^{n}=\sum_{n\in\mathbb{N}}%
2^{n}x^{n}=\dfrac{1}{1-2x}.
\]

\end{example}

\begin{example}
Let $B^{\prime}$ be the weighted set of all binary strings, i.e., finite
tuples consisting of $0$s and $1$s. The weight of an $n$-tuple $\left(
a_{1},a_{2},\ldots,a_{n}\right)  \in B^{\prime}$ is defined to be $a_{1}%
+a_{2}+\cdots+a_{n}$. (Thus, $B^{\prime}$ is the same set as the $B$ in
Example \ref{exa.ws.bin-string1}, but equipped with a different weight
function.) This weighted set $B^{\prime}$ is \textbf{not} finite-type, since
there are infinitely many $a\in B^{\prime}$ having weight $0$ (namely, all
tuples of the form $\left(  0,0,\ldots,0\right)  $, with any number of zeroes.)
\end{example}

Those familiar with graded vector spaces (or graded modules) will recognize
weighted sets as their combinatorial analogues: Weighted sets are to sets what
graded vector spaces are to vector spaces. (In particular, the weight
generating function $\overline{A}$ of a finite-type weighted set $A$ is an
analogue of the Hilbert series of a finite-type graded vector space.) \medskip

Let us now prove a few basic properties of weighted sets and their weight
generating functions.

\begin{proposition}
\label{prop.gf-ws.iso}Let $A$ and $B$ be two isomorphic finite-type weighted
sets. Then, $\overline{A}=\overline{B}$.
\end{proposition}

\begin{proof}
This is almost trivial: The weighted sets $A$ and $B$ are isomorphic; thus,
there exists an isomorphism $\rho:A\rightarrow B$. Consider this $\rho$. Then,
$\rho$ is a bijection and preserves the weight (since $\rho$ is an isomorphism
of weighted sets). The latter property says that we have
\begin{equation}
\left\vert \rho\left(  a\right)  \right\vert =\left\vert a\right\vert
\ \ \ \ \ \ \ \ \ \ \text{for each }a\in A. \label{pf.prop.gf-ws.iso.wtpres}%
\end{equation}
Now, the definition of $\overline{B}$ yields%
\begin{align*}
\overline{B}  &  =\sum_{b\in B}x^{\left\vert b\right\vert }=\sum_{a\in
A}\underbrace{x^{\left\vert \rho\left(  a\right)  \right\vert }}%
_{\substack{=x^{\left\vert a\right\vert }\\\text{(by
(\ref{pf.prop.gf-ws.iso.wtpres}))}}}\ \ \ \ \ \ \ \ \ \ \left(
\begin{array}
[c]{c}%
\text{here, we have substituted }\rho\left(  a\right)  \text{ for }b\\
\text{in the sum, since }\rho\text{ is a bijection}%
\end{array}
\right) \\
&  =\sum_{a\in A}x^{\left\vert a\right\vert }.
\end{align*}
Comparing this with%
\[
\overline{A}=\sum_{a\in A}x^{\left\vert a\right\vert }%
\ \ \ \ \ \ \ \ \ \ \left(  \text{by the definition of }\overline{A}\right)
,
\]
we obtain $\overline{A}=\overline{B}$. This proves Proposition
\ref{prop.gf-ws.iso}.
\end{proof}

Note that Proposition \ref{prop.gf-ws.iso} has a converse: If $\overline
{A}=\overline{B}$, then $A\cong B$. \medskip

Recall that the \emph{disjoint union} of two sets $A$ and $B$ is informally
understood to be \textquotedblleft the union of $A$ and $B$ where we pretend
that the sets $A$ and $B$ are disjoint (even if they aren't)\textquotedblright%
. Formally, it is defined as the set $\left(  \left\{  0\right\}  \times
A\right)  \cup\left(  \left\{  1\right\}  \times B\right)  $, but we think of
the elements $\left(  0,a\right)  $ of this set as copies of the respective
elements $a\in A$, and we think of the elements $\left(  1,b\right)  $ of this
set as copies of the respective elements $b\in B$. (Thus, even if the original
sets $A$ and $B$ have some elements in common, say $c\in A\cap B$, the two
copies $\left(  0,c\right)  $ and $\left(  1,c\right)  $ of such common
elements will not coincide.)

If $A$ and $B$ are two finite sets, then their disjoint union always has size
$\left\vert A\right\vert +\left\vert B\right\vert $.

Let us now extend the concept of a disjoint union to weighted sets:

\begin{definition}
\label{def.gf-ws.djun}Let $A$ and $B$ be two weighted sets. Then, the weighted
set $A+B$ is defined to be the disjoint union of $A$ and $B$, with the weight
function inherited from $A$ and $B$ (meaning that each element of $A$ has the
same weight that it had in $A$, and each element of $B$ has the same weight
that it had in $B$). Formally speaking, this means that $A+B$ is the set
$\left(  \left\{  0\right\}  \times A\right)  \cup\left(  \left\{  1\right\}
\times B\right)  $, with the weight function given by
\begin{equation}
\left\vert \left(  0,a\right)  \right\vert =\left\vert a\right\vert
\ \ \ \ \ \ \ \ \ \ \text{for each }a\in A \label{eq.def.gf-ws.djun.wt0}%
\end{equation}
and%
\begin{equation}
\left\vert \left(  1,b\right)  \right\vert =\left\vert b\right\vert
\ \ \ \ \ \ \ \ \ \ \text{for each }b\in B. \label{eq.def.gf-ws.djun.wt1}%
\end{equation}

\end{definition}

\begin{proposition}
\label{prop.gf-ws.djun}Let $A$ and $B$ be two finite-type weighted sets. Then,
$A+B$ is finite-type, too, and satisfies $\overline{A+B}=\overline
{A}+\overline{B}$.
\end{proposition}

\begin{proof}
The formal definition of $A+B$ says that $A+B=\left(  \left\{  0\right\}
\times A\right)  \cup\left(  \left\{  1\right\}  \times B\right)  $. Thus, the
set $A+B$ is the union of the two disjoint sets $\left\{  0\right\}  \times A$
and $\left\{  1\right\}  \times B$ (indeed, these two sets are clearly
disjoint, since $0\neq1$).

Let us first check that $A+B$ is finite-type.

For each $n\in\mathbb{N}$, there are only finitely many $a\in A$ having weight
$\left\vert a\right\vert =n$ (since $A$ is finite-type), and there are only
finitely many $b\in B$ having weight $\left\vert b\right\vert =n$ (since $B$
is finite-type). Hence, there are only finitely many $c\in A+B$ having weight
$\left\vert c\right\vert =n$ (because any such $c$ either has the form
$\left(  0,a\right)  $ for some $a\in A$ having weight $\left\vert
a\right\vert =n$, or has the form $\left(  1,b\right)  $ for some $b\in B$
having weight $\left\vert b\right\vert =n$). In other words, the weighted set
$A+B$ is finite-type.

Let us now check that $\overline{A+B}=\overline{A}+\overline{B}$. Indeed, the
definition of $\overline{A+B}$ yields%
\begin{align*}
\overline{A+B}  &  =\sum_{c\in A+B}x^{\left\vert c\right\vert }\\
&  =\underbrace{\sum_{c\in\left\{  0\right\}  \times A}x^{\left\vert
c\right\vert }}_{\substack{=\sum_{a\in A}x^{\left\vert \left(  0,a\right)
\right\vert }\\\text{(here, we have substituted }\left(  0,a\right)
\\\text{for }c\text{ in the sum, since the}\\\text{map }A\rightarrow\left\{
0\right\}  \times A,\ a\mapsto\left(  0,a\right)  \text{ is a bijection)}%
}}+\underbrace{\sum_{c\in\left\{  1\right\}  \times B}x^{\left\vert
c\right\vert }}_{\substack{=\sum_{b\in B}x^{\left\vert \left(  1,b\right)
\right\vert }\\\text{(here, we have substituted }\left(  1,b\right)
\\\text{for }c\text{ in the sum, since the}\\\text{map }B\rightarrow\left\{
1\right\}  \times B,\ b\mapsto\left(  1,b\right)  \text{ is a bijection)}}}\\
&  \ \ \ \ \ \ \ \ \ \ \ \ \ \ \ \ \ \ \ \ \left(
\begin{array}
[c]{c}%
\text{here, we have split the sum, since the}\\
\text{set }A+B\text{ is the union of the two disjoint}\\
\text{sets }\left\{  0\right\}  \times A\text{ and }\left\{  1\right\}  \times
B
\end{array}
\right) \\
&  =\sum_{a\in A}\underbrace{x^{\left\vert \left(  0,a\right)  \right\vert }%
}_{\substack{=x^{\left\vert a\right\vert }\\\text{(by
(\ref{eq.def.gf-ws.djun.wt0}))}}}+\sum_{b\in B}\underbrace{x^{\left\vert
\left(  1,b\right)  \right\vert }}_{\substack{=x^{\left\vert b\right\vert
}\\\text{(by (\ref{eq.def.gf-ws.djun.wt1}))}}}=\underbrace{\sum_{a\in
A}x^{\left\vert a\right\vert }}_{=\overline{A}}+\underbrace{\sum_{b\in
B}x^{\left\vert b\right\vert }}_{=\overline{B}}=\overline{A}+\overline{B}.
\end{align*}
Thus, the proof of Proposition \ref{prop.gf-ws.djun} is complete.
\end{proof}

We can easily extend Definition \ref{def.gf-ws.djun} and Proposition
\ref{prop.gf-ws.djun} to disjoint unions of any number (even infinite) of
weighted sets. (However, in the case of infinite disjoint unions, we need to
require the disjoint union to be finite-type in order for Proposition
\ref{prop.gf-ws.djun} to make sense.)

We note that disjoint unions of weighted sets do not satisfy associativity
\textquotedblleft on the nose\textquotedblright: If $A$, $B$ and $C$ are three
weighted sets, then the weighted sets $A+B+C$ and $\left(  A+B\right)  +C$ and
$A+\left(  B+C\right)  $ are not literally equal but rather are isomorphic via
canonical isomorphisms. But Proposition \ref{prop.gf-ws.iso} shows that their
weight generating functions are the same, so that we need not distinguish
between them if we are only interested in these functions. \medskip

Another operation on weighted sets is the Cartesian product:

\begin{definition}
\label{def.gf-ws.prod}Let $A$ and $B$ be two weighted sets. Then, the weighted
set $A\times B$ is defined to be the Cartesian product of the sets $A$ and $B$
(that is, the set $\left\{  \left(  a,b\right)  \ \mid\ a\in A\text{ and }b\in
B\right\}  $), with the weight function defined as follows: For any $\left(
a,b\right)  \in A\times B$, we set%
\begin{equation}
\left\vert \left(  a,b\right)  \right\vert =\left\vert a\right\vert
+\left\vert b\right\vert . \label{eq.def.gf-ws.prod.wt}%
\end{equation}

\end{definition}

\begin{proposition}
\label{prop.gf-ws.prod}Let $A$ and $B$ be two finite-type weighted sets. Then,
$A\times B$ is finite-type, too, and satisfies $\overline{A\times B}%
=\overline{A}\cdot\overline{B}$.
\end{proposition}

\begin{proof}
[Proof of Proposition \ref{prop.gf-ws.prod} (sketched).]The proof that
$A\times B$ is finite-type is easy\footnote{Here is the idea: Let
$n\in\mathbb{N}$. We must prove that there are only finitely many pairs
$\left(  a,b\right)  \in A\times B$ having weight $\left\vert \left(
a,b\right)  \right\vert =n$. Since $\left\vert \left(  a,b\right)  \right\vert
=\left\vert a\right\vert +\left\vert b\right\vert $, these pairs have the
property that $a\in A$ has weight $k$ for some $k\in\left\{  0,1,\ldots
,n\right\}  $, and that $b\in B$ has weight $n-k$ for the same $k$. This
leaves only finitely many options for $k$, only finitely many options for $a$
(since $A$ is finite-type and thus has only finitely many elements of weight
$k$), and only finitely many options for $b$ (since $B$ is finite-type and
thus has only finitely many elements of weight $n-k$). Altogether, we thus
obtain only finitely many options for $\left(  a,b\right)  $.}.

Now, the definition of $\overline{A\times B}$ yields%
\begin{align*}
\overline{A\times B}  &  =\sum_{\left(  a,b\right)  \in A\times B}%
\underbrace{x^{\left\vert \left(  a,b\right)  \right\vert }}%
_{\substack{=x^{\left\vert a\right\vert +\left\vert b\right\vert }\\\text{(by
(\ref{eq.def.gf-ws.prod.wt}))}}}\ \ \ \ \ \ \ \ \ \ \left(
\begin{array}
[c]{c}%
\text{since all elements of }A\times B\text{ have}\\
\text{the form }\left(  a,b\right)  \text{ for some }a\text{ and }b
\end{array}
\right) \\
&  =\sum_{\left(  a,b\right)  \in A\times B}\underbrace{x^{\left\vert
a\right\vert +\left\vert b\right\vert }}_{=x^{\left\vert a\right\vert }\cdot
x^{\left\vert b\right\vert }}\\
&  =\sum_{\left(  a,b\right)  \in A\times B}x^{\left\vert a\right\vert }\cdot
x^{\left\vert b\right\vert }=\underbrace{\left(  \sum_{a\in A}x^{\left\vert
a\right\vert }\right)  }_{=\overline{A}}\underbrace{\left(  \sum_{b\in
B}x^{\left\vert b\right\vert }\right)  }_{=\overline{B}}=\overline{A}%
\cdot\overline{B}.
\end{align*}
The proof of Proposition \ref{prop.gf-ws.prod} is thus complete.
\end{proof}

We can easily extend Definition \ref{def.gf-ws.prod} and Proposition
\ref{prop.gf-ws.prod} to Cartesian products of $k$ weighted sets. The weight
function on such a Cartesian product $A_{1}\times A_{2}\times\cdots\times
A_{k}$ is defined by%
\begin{equation}
\left\vert \left(  a_{1},a_{2},\ldots,a_{k}\right)  \right\vert =\left\vert
a_{1}\right\vert +\left\vert a_{2}\right\vert +\cdots+\left\vert
a_{k}\right\vert . \label{eq.gf-ws.prodk.weight}%
\end{equation}
As a particular case of such Cartesian products, we obtain Cartesian powers
when we multiply $k$ copies of the same weighted set:

\begin{definition}
Let $A$ be a weighted set. Then, $A^{k}$ (for $k\in\mathbb{N}$) means the
weighted set $\underbrace{A\times A\times\cdots\times A}_{k\text{ times}}$.
\end{definition}

The analogue of Proposition \ref{prop.gf-ws.prod} for Cartesian products of
$k$ weighted sets then yields the following:

\begin{proposition}
\label{prop.gf-ws.pow}Let $A$ be a finite-type weighted set. Let
$k\in\mathbb{N}$. Then, $A^{k}$ is finite-type, too, and satisfies
$\overline{A^{k}}=\overline{A}^{k}$.
\end{proposition}

Note that the $0$-th Cartesian power $A^{0}$ of a weighted set $A$ always
consists of a single element -- namely, the empty $0$-tuple $\left(
{}\right)  $, which has weight $0$.

If $A$ is a weighted set, then the infinite disjoint union $A^{0}+A^{1}%
+A^{2}+\cdots$ consists of all (finite) tuples of elements of $A$ (including
the $0$-tuple, the $1$-tuples, and so on). The weight of a tuple is the sum of
the weights of its entries (indeed, this is just what
(\ref{eq.gf-ws.prodk.weight}) says).

\subsubsection{Examples}

Now, let us use this theory to revisit some of the things we have already counted:

\begin{itemize}
\item Fix $k\in\mathbb{N}$, and let
\begin{align*}
C_{k}:=  &  \ \left\{  \text{compositions of length }k\right\} \\
=  &  \ \left\{  \left(  a_{1},a_{2},\ldots,a_{k}\right)  \ \mid\ a_{1}%
,a_{2},\ldots,a_{k}\text{ are positive integers}\right\} \\
=  &  \ \mathbb{P}^{k},\ \ \ \ \ \ \ \ \ \ \text{where }\mathbb{P}=\left\{
1,2,3,\ldots\right\}  .
\end{align*}
This set $C_{k}$ becomes a finite-type weighted set if we set $\left\vert
\left(  a_{1},a_{2},\ldots,a_{k}\right)  \right\vert =a_{1}+a_{2}+\cdots
+a_{k}$ for every $\left(  a_{1},a_{2},\ldots,a_{k}\right)  \in C_{k}$. What
is its weight generating function $\overline{C_{k}}$ ? We can turn
$\mathbb{P}$ itself into a weighted set, by defining the weight of a positive
integer $n$ by $\left\vert n\right\vert =n$. Then, $C_{k}=\mathbb{P}^{k}$ not
just as sets, but as weighted sets\footnote{\textit{Proof.} The weight of a
composition $\left(  a_{1},a_{2},\ldots,a_{k}\right)  \in C_{k}$ in $C_{k}$ is%
\[
\left\vert \left(  a_{1},a_{2},\ldots,a_{k}\right)  \right\vert =a_{1}%
+a_{2}+\cdots+a_{k}\ \ \ \ \ \ \ \ \ \ \left(  \text{by definition}\right)  .
\]
The weight of a composition $\left(  a_{1},a_{2},\ldots,a_{k}\right)
\in\mathbb{P}^{k}$ in $\mathbb{P}^{k}$ is%
\begin{align*}
\left\vert \left(  a_{1},a_{2},\ldots,a_{k}\right)  \right\vert  &
=\left\vert a_{1}\right\vert +\left\vert a_{2}\right\vert +\cdots+\left\vert
a_{k}\right\vert \ \ \ \ \ \ \ \ \ \ \left(  \text{by
(\ref{eq.gf-ws.prodk.weight})}\right) \\
&  =a_{1}+a_{2}+\cdots+a_{k}\ \ \ \ \ \ \ \ \ \ \left(  \text{since
}\left\vert n\right\vert =n\text{ for each }n\in\mathbb{P}\right)  .
\end{align*}
These two weights are clearly equal. Thus, the weight functions of $C_{k}$ and
$\mathbb{P}^{k}$ agree. Hence, $C_{k}=\mathbb{P}^{k}$ as weighted sets.}.
Hence,%
\begin{align*}
\overline{C_{k}}  &  =\overline{\mathbb{P}^{k}}=\overline{\mathbb{P}}%
^{k}\ \ \ \ \ \ \ \ \ \ \left(  \text{by Proposition \ref{prop.gf-ws.pow}%
}\right) \\
&  =\left(  x^{1}+x^{2}+x^{3}+\cdots\right)  ^{k}\ \ \ \ \ \ \ \ \ \ \left(
\text{since }\overline{\mathbb{P}}=x^{1}+x^{2}+x^{3}+\cdots\right) \\
&  =\left(  \dfrac{x}{1-x}\right)  ^{k}\ \ \ \ \ \ \ \ \ \ \left(  \text{since
}x^{1}+x^{2}+x^{3}+\cdots=\dfrac{x}{1-x}\right)  .
\end{align*}
This recovers the equality (\ref{pf.thm.fps.comps.num-comps-n-k.Ak=}).

\item Recall the notion of Dyck paths (as defined in Example 2 in Section
\ref{sec.gf.exas}), as well as the Catalan numbers $c_{0},c_{1},c_{2},\ldots$
(defined in the same place). Let
\[
D:=\left\{  \text{Dyck paths from }\left(  0,0\right)  \text{ to }\left(
2n,0\right)  \text{ for some }n\in\mathbb{N}\right\}  .
\]
This set $D$ becomes a finite-type weighted set if we set
\[
\left\vert p\right\vert =n\ \ \ \ \ \ \ \ \ \ \text{whenever }p\text{ is a
Dyck path from }\left(  0,0\right)  \text{ to }\left(  2n,0\right)  .
\]
The weight generating function $\overline{D}$ of this weighted set $D$ is%
\[
\overline{D}=\sum_{n\in\mathbb{N}}\underbrace{\left(  \text{\# of Dyck paths
from }\left(  0,0\right)  \text{ to }\left(  2n,0\right)  \right)
}_{\substack{=c_{n}\\\text{(by the definition of }c_{n}\text{)}}}x^{n}%
=\sum_{n\in\mathbb{N}}c_{n}x^{n}.
\]
This is the generating function we called $C\left(  x\right)  $ in Example 2
of Section \ref{sec.gf.exas}.

Recall that there is only one Dyck path from $\left(  0,0\right)  $ to
$\left(  0,0\right)  $, namely the trivial path. All the other Dyck paths in
$D$ are nontrivial. We let%
\begin{align*}
D_{\operatorname*{triv}}  &  :=\left\{  \text{trivial Dyck paths in
}D\right\}  \ \ \ \ \ \ \ \ \ \ \text{and}\\
D_{\operatorname*{non}}  &  :=\left\{  \text{nontrivial Dyck paths in
}D\right\}  .
\end{align*}
These two sets $D_{\operatorname*{triv}}$ and $D_{\operatorname*{non}}$ are
subsets of $D$, and thus are weighted sets themselves (we define their weight
functions by restricting the one of $D$). The set $D_{\operatorname*{triv}}$
consists of a single Dyck path, which has weight $0$; thus, its weight
generating function is%
\[
\overline{D_{\operatorname*{triv}}}=x^{0}=1.
\]

In Example 2 of Section \ref{sec.gf.exas}, we have seen that any nontrivial
Dyck path $\pi$ has the following structure:\footnote{The colors are referring
to the following picture:%
\[%
\begin{tikzpicture}
\draw
(0, 0) -- (1, 1) -- (2, 2) -- (3, 1) -- (4, 2) -- (5, 1) -- (6, 0) -- (7, 1) -- (8, 2) -- (9, 1) -- (10, 0) -- (11, 1) -- (12, 0) -- (13, 1) -- (14, 2) -- (15, 1) -- (16, 0);
\fill[yellow] (0, 0) -- (1, 1) -- (5, 1) -- (6, 0) -- cycle;
\fill
[green!60!black] (1, 1) -- (2, 2) -- (3, 1) -- (4, 2) -- (5, 1) -- (1, 1);
\fill
[blue!60!white] (6, 0) -- (7, 1) -- (8, 2) -- (9, 1) -- (10, 0) -- (11, 1) -- (12, 0) -- (13, 1) -- (14, 2) -- (15, 1) -- (16, 0);
\filldraw(0, 0) circle [fill=red, radius=0.1];
\filldraw(1, 1) circle [fill=red, radius=0.1];
\filldraw(2, 2) circle [fill=red, radius=0.1];
\filldraw(3, 1) circle [fill=red, radius=0.1];
\filldraw(4, 2) circle [fill=red, radius=0.1];
\filldraw(5, 1) circle [fill=red, radius=0.1];
\filldraw(6, 0) circle [fill=red, radius=0.1];
\filldraw(7, 1) circle [fill=red, radius=0.1];
\filldraw(8, 2) circle [fill=red, radius=0.1];
\filldraw(9, 1) circle [fill=red, radius=0.1];
\filldraw(10, 0) circle [fill=red, radius=0.1];
\filldraw(11, 1) circle [fill=red, radius=0.1];
\filldraw(12, 0) circle [fill=red, radius=0.1];
\filldraw(13, 1) circle [fill=red, radius=0.1];
\filldraw(14, 2) circle [fill=red, radius=0.1];
\filldraw(15, 1) circle [fill=red, radius=0.1];
\filldraw(16, 0) circle [fill=red, radius=0.1];
\end{tikzpicture}%
\]
}

\begin{itemize}
\item a NE-step,

\item followed by a (diagonally shifted) Dyck path (drawn in green),

\item followed by a SE-step,

\item followed by another (horizontally shifted) Dyck path (drawn in purple).
\end{itemize}

If we denote the green Dyck path by $\alpha$ and the purple Dyck path by
$\beta$, then we obtain a pair $\left(  \alpha,\beta\right)  \in D\times D$ of
two Dyck paths. Thus, each nontrivial Dyck path $\pi\in D_{\operatorname*{non}%
}$ gives rise to a pair $\left(  \alpha,\beta\right)  \in D\times D$ of two
Dyck paths. This yields a map%
\begin{align*}
D_{\operatorname*{non}}  &  \rightarrow D\times D,\\
\pi &  \mapsto\left(  \alpha,\beta\right)  ,
\end{align*}
which is easily seen to be a bijection (since any pair $\left(  \alpha
,\beta\right)  \in D\times D$ can be assembled into a single nontrivial Dyck
path $\pi$ that starts with an NE-step, is followed by a shifted copy of
$\alpha$, then by a SE-step, then by a shifted copy of $\beta$).

Alas, this bijection is not an isomorphism of weighted sets, since it fails to
preserve the weight. Indeed, $\left\vert \left(  \alpha,\beta\right)
\right\vert =\left\vert \alpha\right\vert +\left\vert \beta\right\vert
=\left\vert \pi\right\vert -1\neq\left\vert \pi\right\vert $.

Fortunately, we can fix this rather easily. Define a weighted set%
\[
X:=\left\{  1\right\}  \ \ \ \ \ \ \ \ \ \ \text{with }\left\vert 1\right\vert
=1.
\]
This is a one-element set, so the only real difference between the weighted
sets $X\times D\times D$ and $D\times D$ is in the weights. Indeed, the sets
$D\times D$ and $X\times D\times D$ are in bijection (any pair $\left(
\alpha,\beta\right)  \in D\times D$ corresponds to the triple $\left(
1,\alpha,\beta\right)  \in X\times D\times D$), but the weights of
corresponding elements differ by $1$ (namely, $\left\vert \left(
1,\alpha,\beta\right)  \right\vert =\underbrace{\left\vert 1\right\vert }%
_{=1}+\underbrace{\left\vert \alpha\right\vert +\left\vert \beta\right\vert
}_{=\left\vert \left(  \alpha,\beta\right)  \right\vert }=1+\left\vert \left(
\alpha,\beta\right)  \right\vert $).

Thus, by replacing $D\times D$ by $X\times D\times D$, we can fix the degrees
in our above bijection. We thus obtain a bijection%
\begin{align*}
D_{\operatorname*{non}}  &  \rightarrow X\times D\times D,\\
\pi &  \mapsto\left(  1,\alpha,\beta\right)  ,
\end{align*}
which does preserve the weight. This bijection is thus an isomorphism of
weighted sets. Hence, $D_{\operatorname*{non}}\cong X\times D\times D$.

Each Dyck path is either trivial or nontrivial. Hence,%
\[
D\cong D_{\operatorname*{triv}}+\underbrace{D_{\operatorname*{non}}}_{\cong
X\times D\times D}\cong D_{\operatorname*{triv}}+X\times D\times D,
\]
so that%
\begin{align*}
\overline{D}  &  =\overline{D_{\operatorname*{triv}}+X\times D\times
D}\ \ \ \ \ \ \ \ \ \ \left(  \text{by Proposition \ref{prop.gf-ws.iso}%
}\right) \\
&  =\underbrace{\overline{D_{\operatorname*{triv}}}}_{=1}%
+\underbrace{\overline{X}}_{=x}\cdot\underbrace{\overline{D}\cdot\overline{D}%
}_{=\overline{D}^{2}}\ \ \ \ \ \ \ \ \ \ \left(  \text{by Proposition
\ref{prop.gf-ws.djun} and Proposition \ref{prop.gf-ws.prod}}\right) \\
&  =1+x\cdot\overline{D}^{2}.
\end{align*}
This is precisely the quadratic equation $C\left(  x\right)  =1+x\left(
C\left(  x\right)  \right)  ^{2}$ that we obtained in Section
\ref{sec.gf.exas}. But this time, we obtained it in a more conceptual way,
through a combinatorially defined isomorphism $D\cong D_{\operatorname*{triv}%
}+X\times D\times D$ rather than by ad-hoc manipulation of FPSs.
\end{itemize}

\subsubsection{\label{subsec.gf.weighted-set.domino}Domino tilings}

Let us now apply the theory of weight generating functions to something we
haven't already counted. Namely, we shall count the domino tilings of a
rectangle. Informally, these are defined as follows:

\begin{itemize}
\item For any $n,m\in\mathbb{N}$, we let $R_{n,m}$ be a rectangle with width
$n$ and height $m$.

\item A \emph{domino} means a rectangle that is an $R_{1,2}$ or an $R_{2,1}$.

\item A \emph{domino tiling} of a shape $S$ means a tiling of $S$ by dominos
(i.e., a set of dominos that cover $S$ and whose interiors don't intersect).
\end{itemize}

For example, here is a domino tiling of $R_{8,8}$:%
\[%
\begin{tikzpicture}
\filldraw[draw=black, fill=green!20!white] (0, 2) rectangle (1, 4);
\filldraw[draw=black, fill=green!20!white] (0, 4) rectangle (1, 6);
\filldraw[draw=black, fill=green!20!white] (1, 3) rectangle (2, 5);
\filldraw[draw=black, fill=green!20!white] (2, 0) rectangle (3, 2);
\filldraw[draw=black, fill=green!20!white] (2, 6) rectangle (3, 8);
\filldraw[draw=black, fill=green!20!white] (3, 1) rectangle (4, 3);
\filldraw[draw=black, fill=green!20!white] (3, 5) rectangle (4, 7);
\filldraw[draw=black, fill=green!20!white] (4, 2) rectangle (5, 4);
\filldraw[draw=black, fill=green!20!white] (4, 4) rectangle (5, 6);
\filldraw[draw=black, fill=green!20!white] (5, 2) rectangle (6, 4);
\filldraw[draw=black, fill=green!20!white] (5, 4) rectangle (6, 6);
\filldraw[draw=black, fill=green!20!white] (6, 1) rectangle (7, 3);
\filldraw[draw=black, fill=green!20!white] (6, 3) rectangle (7, 5);
\filldraw[draw=black, fill=green!20!white] (6, 5) rectangle (7, 7);
\filldraw[draw=black, fill=green!20!white] (7, 0) rectangle (8, 2);
\filldraw[draw=black, fill=green!20!white] (7, 2) rectangle (8, 4);
\filldraw[draw=black, fill=green!20!white] (7, 4) rectangle (8, 6);
\filldraw[draw=black, fill=green!20!white] (7, 6) rectangle (8, 8);
\filldraw[draw=black, fill=red!20!white] (0, 0) rectangle (2, 1);
\filldraw[draw=black, fill=red!20!white] (0, 1) rectangle (2, 2);
\filldraw[draw=black, fill=red!20!white] (0, 6) rectangle (2, 7);
\filldraw[draw=black, fill=red!20!white] (0, 7) rectangle (2, 8);
\filldraw[draw=black, fill=red!20!white] (1, 2) rectangle (3, 3);
\filldraw[draw=black, fill=red!20!white] (1, 5) rectangle (3, 6);
\filldraw[draw=black, fill=red!20!white] (2, 3) rectangle (4, 4);
\filldraw[draw=black, fill=red!20!white] (2, 4) rectangle (4, 5);
\filldraw[draw=black, fill=red!20!white] (3, 0) rectangle (5, 1);
\filldraw[draw=black, fill=red!20!white] (3, 7) rectangle (5, 8);
\filldraw[draw=black, fill=red!20!white] (4, 1) rectangle (6, 2);
\filldraw[draw=black, fill=red!20!white] (4, 6) rectangle (6, 7);
\filldraw[draw=black, fill=red!20!white] (5, 0) rectangle (7, 1);
\filldraw[draw=black, fill=red!20!white] (5, 7) rectangle (7, 8);
\end{tikzpicture}%
\]
(note that the colors are purely ornamental here: we are coloring a domino
pink if it lies horizontally and green if it stands vertically, for the sake
of convenience).

This sounds geometric, but actually is a combinatorial object hiding behind
geometric language. Our rectangles and dominos all align to a square grid.
Thus, rectangles can be modeled simply as finite sets of grid squares, and
dominos are unordered pairs of adjacent grid squares. Grid squares, in turn,
can be modeled as pairs $\left(  i,j\right)  $ of integers (corresponding to
the Cartesian coordinates of their centers). Thus, we redefine domino tilings
combinatorially as follows:

\begin{definition}
\label{def.domino.shapes-and-tilings}\textbf{(a)} A \emph{shape} means a
subset of $\mathbb{Z}^{2}$.

We draw each $\left(  i,j\right)  \in\mathbb{Z}^{2}$ as a unit square with
center at the point $\left(  i,j\right)  $ (in Cartesian coordinates); thus, a
shape can be drawn as a cluster of squares. \medskip

\textbf{(b)} For any $n,m\in\mathbb{N}$, the shape $R_{n,m}$ (called the
$n\times m$\emph{-rectangle}) is defined to be%
\[
\left\{  1,2,\ldots,n\right\}  \times\left\{  1,2,\ldots,m\right\}  =\left\{
\left(  i,j\right)  \in\mathbb{Z}^{2}\ \mid\ 1\leq i\leq n\text{ and }1\leq
j\leq m\right\}  .
\]

\textbf{(c)} A \emph{domino} means a size-$2$ shape of the form%
\begin{align*}
&  \left\{  \left(  i,j\right)  ,\ \left(  i+1,j\right)  \right\}  \text{ (a
\textquotedblleft\emph{horizontal domino}\textquotedblright)}%
\ \ \ \ \ \ \ \ \ \ \text{or}\\
&  \left\{  \left(  i,j\right)  ,\ \left(  i,j+1\right)  \right\}  \text{ (a
\textquotedblleft\emph{vertical domino}\textquotedblright)}%
\end{align*}
for some $\left(  i,j\right)  \in\mathbb{Z}^{2}$. \medskip

\textbf{(d)} A \emph{domino tiling} of a shape $S$ is a set partition of $S$
into dominos (i.e., a set of disjoint dominos whose union is $S$). \medskip

\textbf{(e)} For any $n,m\in\mathbb{N}$, let $d_{n,m}$ be the \# of domino
tilings of $R_{n,m}$.
\end{definition}

For example, the domino tiling $\raisebox{-3.5ex}{
\begin{tikzpicture}[scale=0.8] \draw (1, 0) rectangle (3, 1); \draw (1, 2) rectangle (3, 1); \draw (0, 0) rectangle (1, 2); \end{tikzpicture}
}$ of the rectangle $R_{3,2}$ is the set partition%
\[
\left\{  \left\{  \left(  1,1\right)  ,\ \left(  1,2\right)  \right\}
,\ \left\{  \left(  2,1\right)  ,\ \left(  3,1\right)  \right\}  ,\ \left\{
\left(  2,2\right)  ,\ \left(  3,2\right)  \right\}  \right\}
\]
of $R_{3,2}$ (here, $\left\{  \left(  1,1\right)  ,\ \left(  1,2\right)
\right\}  $ is the vertical domino, whereas $\left\{  \left(  2,1\right)
,\ \left(  3,1\right)  \right\}  $ is the bottom horizontal domino, and
$\left\{  \left(  2,2\right)  ,\ \left(  3,2\right)  \right\}  $ is the top
horizontal domino). \medskip

Can we compute $d_{n,m}$ ?

The case $m=1$ is a bit too simple (do it!), so let us start with the case
$m=2$. Here are the $d_{n,2}$ for $n\in\left\{  0,1,\ldots,4\right\}  $:%
\[%
\begin{tabular}
[c]{|c|c|l|}\hline
$n$ & $d_{n,2}$ & domino tilings of $R_{n,2}$\\\hline\hline
$0$ & $d_{0,2}=1$ &
$\vphantom{\dfrac{\int}{\int}}\begin{tikzpicture}[scale=0.8] \draw (0, 0) rectangle (0, 2); \end{tikzpicture}$%
\\\hline
$1$ & $d_{1,2}=1$ &
$\vphantom{\dfrac{\int}{\int}}\begin{tikzpicture}[scale=0.8] \draw (0, 0) rectangle (1, 2); \end{tikzpicture}$%
\\\hline
$2$ & $d_{2,2}=2$ &
$\vphantom{\dfrac{\int}{\int}}\begin{tikzpicture}[scale=0.8] \draw (0, 0) rectangle (1, 2); \draw (2, 0) rectangle (1, 2); \end{tikzpicture}\ \ \ ,\ \ \ \begin{tikzpicture}[scale=0.8] \draw (0, 0) rectangle (2, 1); \draw (0, 2) rectangle (2, 1); \end{tikzpicture}$%
\\\hline
$3$ & $d_{3,2}=3$ &
$\vphantom{\dfrac{\int}{\int}}\begin{tikzpicture}[scale=0.8] \draw (0, 0) rectangle (1, 2); \draw (2, 0) rectangle (1, 2); \draw (2, 0) rectangle (3, 2); \end{tikzpicture}\ \ \ ,\ \ \ \begin{tikzpicture}[scale=0.8] \draw (0, 0) rectangle (2, 1); \draw (0, 2) rectangle (2, 1); \draw (2, 0) rectangle (3, 2); \end{tikzpicture}\ \ \ ,\ \ \ \begin{tikzpicture}[scale=0.8] \draw (1, 0) rectangle (3, 1); \draw (1, 2) rectangle (3, 1); \draw (0, 0) rectangle (1, 2); \end{tikzpicture}$%
\\\hline
$4$ & $d_{4,2}=5$ &
$\vphantom{\dfrac{\int}{\int}}\begin{tikzpicture}[scale=0.8] \draw (0, 0) rectangle (1, 2); \draw (2, 0) rectangle (1, 2); \draw (2, 0) rectangle (3, 2); \draw (4, 0) rectangle (3, 2); \end{tikzpicture}\ \ \ ,\ \ \ \begin{tikzpicture}[scale=0.8] \draw (0, 0) rectangle (2, 1); \draw (0, 2) rectangle (2, 1); \draw (2, 0) rectangle (3, 2); \draw (4, 0) rectangle (3, 2); \end{tikzpicture}\ \ \ ,\ \ \ \begin{tikzpicture}[scale=0.8] \draw (1, 0) rectangle (3, 1); \draw (1, 2) rectangle (3, 1); \draw (0, 0) rectangle (1, 2); \draw (4, 0) rectangle (3, 2); \end{tikzpicture}\ \ \ ,$%
\\
&  &
$\vphantom{\dfrac{\int}{\int}}\begin{tikzpicture}[scale=0.8] \draw (2, 0) rectangle (4, 1); \draw (2, 2) rectangle (4, 1); \draw (0, 0) rectangle (1, 2); \draw (2, 0) rectangle (1, 2); \end{tikzpicture}\ \ \ ,\ \ \ \begin{tikzpicture}[scale=0.8] \draw (2, 0) rectangle (4, 1); \draw (2, 2) rectangle (4, 1); \draw (0, 0) rectangle (2, 1); \draw (0, 2) rectangle (2, 1); \end{tikzpicture}$%
\\\hline
\end{tabular}
\ \ \ \ \ \ \ \ \ .
\]

Let us try to compute $d_{n,2}$ in general.

A \emph{height-}$2$ \emph{rectangle} shall mean a rectangle of the form
$R_{n,2}$ with $n\in\mathbb{N}$. Let us define the weighted set%
\begin{align*}
D:=  &  \left\{  \text{domino tilings of height-}2\text{ rectangles}\right\}
\\
=  &  \left\{  \text{domino tilings of }R_{n,2}\text{ with }n\in
\mathbb{N}\right\}  .
\end{align*}
We define the weight of a tiling $T$ of $R_{n,2}$ to be $\left\vert
T\right\vert :=n$ (that is, the width of the rectangle tiled by this tiling).

Thus, $D$ is a finite-type weighted set, with weight generating function%
\[
\overline{D}=\sum_{n\in\mathbb{N}}d_{n,2}x^{n}.
\]

So we want to compute $\overline{D}$. Let us define a new weighted set that
will help us at that.

Namely, we say that a \emph{fault} of a domino tiling $T$ is a vertical line
$\ell$ such that

\begin{itemize}
\item each domino of $T$ lies either left of $\ell$ or right of $\ell$ (but
does not straddle $\ell$), and

\item there is at least one domino of $T$ that lies left of $\ell$, and at
least one domino of $T$ that lies right of $\ell$.
\end{itemize}

For example, here are two domino tilings:%
\[%
\begin{tikzpicture}[scale=0.8]
\draw(3, 0) rectangle (2, 2);
\draw(2, 2) rectangle (3, 4);
\draw(2, 0) rectangle (0, 1);
\draw(2, 1) rectangle (1, 3);
\draw(0, 1) rectangle (1, 3);
\draw(2, 3) rectangle (0, 4);
\end{tikzpicture}%
\ \ \ \ \ \ \ \ \ \ \
\begin{tikzpicture}[scale=0.8]
\draw(0, 3) rectangle (2, 2);
\draw(2, 2) rectangle (4, 3);
\draw(0, 2) rectangle (1, 0);
\draw(1, 2) rectangle (3, 1);
\draw(1, 0) rectangle (3, 1);
\draw(3, 2) rectangle (4, 0);
\end{tikzpicture}%
\ \ \ \ \ \ .
\]
The tiling on the left has a fault (namely, the vertical line separating the
$2$nd from the $3$rd column), but the tiling on the right has none (a fault
must be a vertical line by definition; a horizontal line doesn't count).

A domino tiling will be called \emph{faultfree} if it is nonempty and has no
fault. Thus, the tiling on the right (in the above example) is faultfree.

We now observe a crucial (but trivial) lemma:

\begin{lemma}
\label{lem.gf.weighted-set.domino.fd}Any domino tiling of a height-$2$
rectangle can be decomposed uniquely into a tuple of faultfree tilings of
(usually smaller) height-$2$ rectangles, by cutting it along its faults. For
example:%
\begin{align*}
&  \raisebox{-3.5ex}{
\begin{tikzpicture}[scale=0.8] \draw (1, 0) rectangle (3, 1); \draw (1, 2) rectangle (3, 1); \draw (0, 0) rectangle (1, 2); \draw (4, 0) rectangle (3, 2); \draw (5, 0) rectangle (4, 2); \end{tikzpicture}
}\ \text{ decomposes as}\\
&  \left(  \ \ \raisebox{-3.5ex}{
\begin{tikzpicture}[scale=0.8] \draw (0, 0) rectangle (1, 2); \end{tikzpicture}
}\ \ \ \ \ ,\ \ \raisebox{-3.5ex}{
\begin{tikzpicture}[scale=0.8] \draw (0, 0) rectangle (2, 1); \draw (0, 2) rectangle (2, 1); \end{tikzpicture}
}\ \ \ \ \ ,\ \ \raisebox{-3.5ex}{
\begin{tikzpicture}[scale=0.8] \draw (0, 0) rectangle (1, 2); \end{tikzpicture}
}\ \ \ \ \ ,\ \ \raisebox{-3.5ex}{
\begin{tikzpicture}[scale=0.8] \draw (0, 0) rectangle (1, 2); \end{tikzpicture}
}\ \ \right)  .
\end{align*}
(Note that if the original tiling was faultfree, then it will decompose into a
$1$-tuple. If the original tiling was empty, then it will decompose into a $0$-tuple.)

Moreover, the sum of the weights of the faultfree tilings in the tuple is the
weight of the original tiling. (In other words, if a tiling $T$ decomposes
into the tuple $\left(  T_{1},T_{2},\ldots,T_{k}\right)  $, then $\left\vert
T\right\vert =\left\vert T_{1}\right\vert +\left\vert T_{2}\right\vert
+\cdots+\left\vert T_{k}\right\vert $.)
\end{lemma}

Thus, if we define a new weighted set%
\[
F:=\left\{  \text{\textbf{faultfree} domino tilings of }R_{n,2}\text{ with
}n\in\mathbb{N}\right\}
\]
(with the same weights as in $D$), then we obtain an isomorphism (i.e.,
weight-preserving bijection)%
\[
D\rightarrow\underbrace{F^{0}+F^{1}+F^{2}+F^{3}+\cdots}_{\substack{\text{an
infinite disjoint union}\\\text{of weighted sets}}}
\]
that sends each tiling to the tuple it decomposes into. Hence,%
\[
D\cong F^{0}+F^{1}+F^{2}+F^{3}+\cdots,
\]
and therefore%
\begin{align*}
\overline{D}  &  =\overline{F^{0}+F^{1}+F^{2}+F^{3}+\cdots}=\overline{F^{0}%
}+\overline{F^{1}}+\overline{F^{2}}+\overline{F^{3}}+\cdots\\
&  \ \ \ \ \ \ \ \ \ \ \ \ \ \ \ \ \ \ \ \ \left(  \text{by the infinite
analogue of Proposition \ref{prop.gf-ws.djun}}\right) \\
&  =\overline{F}^{0}+\overline{F}^{1}+\overline{F}^{2}+\overline{F}^{3}%
+\cdots\ \ \ \ \ \ \ \ \ \ \left(  \text{by Proposition \ref{prop.gf-ws.pow}%
}\right) \\
&  =\dfrac{1}{1-\overline{F}}\ \ \ \ \ \ \ \ \ \ \left(  \text{by
(\ref{eq.sec.gf.exas.1.1/1-x}), with }\overline{F}\text{ substituted for
}x\right)  .
\end{align*}
Thus, if we can compute $\overline{F}$, then we can compute $\overline{D}$.

In order to compute $\overline{F}$, let us see how a faultfree domino tiling
of a height-$2$ rectangle looks like. Here are two such tilings:%
\[
\raisebox{-3.5ex}{
\begin{tikzpicture}[scale=0.8] \draw (0, 0) rectangle (1, 2); \end{tikzpicture}
}\ \ \ \ \ \ \ \ \ \ \text{and}\ \ \ \ \ \ \ \ \ \ \raisebox{-3.5ex}{
\begin{tikzpicture}[scale=0.8] \draw (0, 0) rectangle (2, 1); \draw (0, 2) rectangle (2, 1); \end{tikzpicture}
}\ \ \ \ \ .
\]

I claim that these two tilings are the \textbf{only} faultfree tilings of
height-$2$ rectangles. Indeed, consider any faultfree tiling of a height-$2$
rectangle. In this tiling, look at the domino that covers the box $\left(
1,1\right)  $. If it is a vertical domino, then this vertical domino must
constitute the entire tiling, since otherwise there would be a fault to its
right. If it is a horizontal domino, then there must be a second horizontal
domino stacked atop it, and these two dominos must then constitute the entire
tiling, since otherwise there would be a fault to their right. This leads to
the two options we just named.

Thus, the weighted set $F$ consists of just the two tilings shown above: one
tiling of weight $1$ and one tiling of weight $2$. Hence, its weight
generating function is $\overline{F}=x+x^{2}$. So%
\[
\overline{D}=\dfrac{1}{1-\overline{F}}=\dfrac{1}{1-\left(  x+x^{2}\right)
}=\dfrac{1}{1-x-x^{2}}=f_{1}+f_{2}x+f_{3}x^{2}+f_{4}x^{3}+\cdots,
\]
where $\left(  f_{0},f_{1},f_{2},\ldots\right)  $ is the Fibonacci sequence.
Thus, comparing coefficients, we find%
\[
d_{n,2}=f_{n+1}\ \ \ \ \ \ \ \ \ \ \text{for each }n\in\mathbb{N}.
\]

There are, of course, more elementary proofs of this (see \cite[Proposition
1.1.11]{19fco}).

\begin{remark}
Here is an outline of an alternative proof of $\overline{D}=\dfrac
{1}{1-\overline{F}}$:

Any tiling in $D$ is either empty, or can be uniquely split into a pair of a
faultfree tiling and an arbitrary tiling (just split it along its leftmost
fault, or along the right end if there is no fault). Thus,%
\[
D\cong\left\{  0\right\}  +F\times D,
\]
where the set $\left\{  0\right\}  $ here is viewed as a weighted set with
$\left\vert 0\right\vert =0$. Hence,%
\[
\overline{D}=\overline{\left\{  0\right\}  +F\times D}=\overline{\left\{
0\right\}  }+\overline{F}\cdot\overline{D}=1+\overline{F}\cdot\overline{D}.
\]
Solving this for $\overline{D}$, we find $\overline{D}=\dfrac{1}%
{1-\overline{F}}.$
\end{remark}

Now, let us try to solve the analogous problem for height-$3$ rectangles.
Forget about the $D$ and $F$ we defined above. Instead, define a new weighted
set%
\begin{align*}
D:=  &  \left\{  \text{domino tilings of height-}3\text{ rectangles}\right\}
\\
=  &  \left\{  \text{domino tilings of }R_{n,3}\text{ with }n\in
\mathbb{N}\right\}  .
\end{align*}
The weight of a tiling $T$ of $R_{n,3}$ is defined as before (i.e., it is
$\left\vert T\right\vert =n$). Thus, $D$ is a finite-type weighted set, with
generating function%
\[
\overline{D}=\sum_{n\in\mathbb{N}}d_{n,3}x^{n}.
\]
We want to compute this $\overline{D}$.

Set%
\[
F:=\left\{  \text{\textbf{faultfree} domino tilings of }R_{n,3}\text{ with
}n\in\mathbb{N}\right\}  .
\]
Again, we can show $\overline{D}=\dfrac{1}{1-\overline{F}}$ (as we did above
in the case of $R_{n,2}$). We thus need to compute $\overline{F}$.

How does a faultfree domino tiling of a height-$3$ rectangle look like? Let us
classify such tilings according to the kind of dominos that occupy the first
column of the tiling. The proof of the following classification can be found
in the Appendix (Proposition \ref{prop.gf.weighted-set.domino.Rn3.ABC} in
Section \ref{sec.details.domino}):

\begin{itemize}
\item The faultfree domino tilings of a height-$3$ rectangle that contain a
vertical domino in the \textbf{top} two squares of the first column are%
\[%
\begin{tikzpicture}[scale=0.8]
\draw(0, 0) rectangle (2, 1);
\draw(0, 1) rectangle (1, 3);
\draw(1, 3) rectangle (2, 1);
\end{tikzpicture}
\qquad\begin{tikzpicture}[scale=0.8]
\draw(0, 0) rectangle (2, 1);
\draw(2, 1) rectangle (4, 0);
\draw(0, 1) rectangle (1, 3);
\draw(1, 1) rectangle (3, 2);
\draw(1, 3) rectangle (3, 2);
\draw(3, 1) rectangle (4, 3);
\end{tikzpicture}
\qquad\begin{tikzpicture}[scale=0.8]
\draw(0, 0) rectangle (2, 1);
\draw(2, 1) rectangle (4, 0);
\draw(6, 1) rectangle (4, 0);
\draw(0, 1) rectangle (1, 3);
\draw(1, 1) rectangle (3, 2);
\draw(1, 3) rectangle (3, 2);
\draw(3, 1) rectangle (5, 2);
\draw(3, 3) rectangle (5, 2);
\draw(5, 1) rectangle (6, 3);
\node(X) at (7.5, 1.5) {$\Large\cdots$};
\end{tikzpicture}%
\]
(this is an infinite sequence of tilings, each obtained from the previous by
inserting two columns in the middle by a fairly self-explanatory procedure).
The weights of these tilings are $2,4,6,\ldots$, so their total contribution
to the weight generating function $\overline{F}$ of $F$ is $x^{2}+x^{4}%
+x^{6}+\cdots$.

\item The faultfree domino tilings of a height-$3$ rectangle that contain a
vertical domino in the \textbf{bottom} two squares of the first column are%
\[%
\begin{tikzpicture}[scale=0.8]
\draw(0, 3) rectangle (2, 2);
\draw(0, 2) rectangle (1, 0);
\draw(1, 0) rectangle (2, 2);
\end{tikzpicture}
\qquad\begin{tikzpicture}[scale=0.8]
\draw(0, 3) rectangle (2, 2);
\draw(2, 2) rectangle (4, 3);
\draw(0, 2) rectangle (1, 0);
\draw(1, 2) rectangle (3, 1);
\draw(1, 0) rectangle (3, 1);
\draw(3, 2) rectangle (4, 0);
\end{tikzpicture}
\qquad\begin{tikzpicture}[scale=0.8]
\draw(0, 3) rectangle (2, 2);
\draw(2, 2) rectangle (4, 3);
\draw(6, 2) rectangle (4, 3);
\draw(0, 2) rectangle (1, 0);
\draw(1, 2) rectangle (3, 1);
\draw(1, 0) rectangle (3, 1);
\draw(3, 2) rectangle (5, 1);
\draw(3, 0) rectangle (5, 1);
\draw(5, 2) rectangle (6, 0);
\node(X) at (7.5, 1.5) {$\Large\cdots$};
\end{tikzpicture}%
\]
(these are the top-down mirror images of the tilings from the previous bullet
point). The weights of these tilings are $2,4,6,\ldots$, so their total
contribution to the weight generating function $\overline{F}$ of $F$ is
$x^{2}+x^{4}+x^{6}+\cdots$.

\item The faultfree domino tilings of a height-$3$ rectangle that contain
\textbf{no} vertical domino in the first column are%
\[%
\begin{tikzpicture}[scale=0.8]
\draw(0, 0) rectangle (2, 1);
\draw(0, 1) rectangle (2, 2);
\draw(0, 2) rectangle (2, 3);
\end{tikzpicture}%
\]
(yes, there is only one such tiling). The weight of this tiling is $2$, so its
total contribution to the weight generating function $\overline{F}$ of $F$ is
$x^{2}$.
\end{itemize}

This classification of faultfree domino tilings entails%
\begin{align*}
\overline{F}  &  =\left(  x^{2}+x^{4}+x^{6}+\cdots\right)  +\left(
x^{2}+x^{4}+x^{6}+\cdots\right)  +x^{2}\\
&  =x^{2}\cdot\dfrac{1}{1-x^{2}}+x^{2}\cdot\dfrac{1}{1-x^{2}}+x^{2}%
\ \ \ \ \ \ \ \ \ \ \left(  \text{since }x^{2}+x^{4}+x^{6}+\cdots=x^{2}%
\cdot\dfrac{1}{1-x^{2}}\right) \\
&  =\dfrac{3x^{2}-x^{4}}{1-x^{2}}.
\end{align*}
Thus,%
\begin{align*}
\overline{D}  &  =\dfrac{1}{1-\overline{F}}=\dfrac{1}{1-\dfrac{3x^{2}-x^{4}%
}{1-x^{2}}}=\dfrac{1-x^{2}}{1-4x^{2}+x^{4}}\\
&  =1+3x^{2}+11x^{4}+41x^{6}+153x^{8}+\cdots.
\end{align*}

You will notice that only even powers of $x$ appear in this FPS. In other
words,%
\[
d_{n,3}=0\text{ when }n\text{ is odd.}%
\]
This is not surprising, because if $n$ is odd, then the rectangle $R_{n,3}$
has an odd \# of squares, and thus cannot be tiled by dominos.

But we can also compute $d_{n,3}$ for even $n$. Indeed, using the same method
(partial fractions) that we used for the Fibonacci sequence in Section
\ref{sec.gf.exas}, we can expand $\dfrac{1-x^{2}}{1-4x^{2}+x^{4}}$ as a sum of
geometric series:%
\[
\dfrac{1-x^{2}}{1-4x^{2}+x^{4}}=\dfrac{3+\sqrt{3}}{6}\cdot\dfrac{1}{1-\left(
2+\sqrt{3}\right)  x^{2}}+\dfrac{3-\sqrt{3}}{6}\cdot\dfrac{1}{1-\left(
2-\sqrt{3}\right)  x^{2}}.
\]
Thus, we find%
\[
d_{n,3}=\dfrac{3+\sqrt{3}}{6}\left(  2+\sqrt{3}\right)  ^{n/2}+\dfrac
{3-\sqrt{3}}{6}\left(  2-\sqrt{3}\right)  ^{n/2}\ \ \ \ \ \ \ \ \ \ \text{for
any even }n.
\]

Now, what about computing $d_{n,m}$ in general? The above reasoning leading up
to $\overline{D}=\dfrac{1}{1-\overline{F}}$ can be applied for any
$m\in\mathbb{N}$, but describing $F$ becomes harder and harder as $m$ grows
larger. The generating function $\overline{D}$ is still a quotient of two
polynomials for any $m$ (see, e.g., \cite{KlaPol79}), but this requires more
insight to prove. For $m\geq6$, it appears that there is no formula for
$d_{n,m}$ that requires only quadratic irrationalities.

It is worth mentioning a different formula for $d_{n,m}$, found by Kasteleyn
in 1961 (motivated by a theoretical physics model):

\begin{theorem}
[Kasteleyn's formula]Let $n,m\in\mathbb{N}$ be such that $m$ is even and
$n\geq1$. Then,%
\[
d_{n,m}=2^{mn/2}\ \prod_{j=1}^{m/2}\ \ \prod_{k=1}^{n}\sqrt{\left(  \cos
\dfrac{j\pi}{m+1}\right)  ^{2}+\left(  \cos\dfrac{k\pi}{n+1}\right)  ^{2}}.
\]

\end{theorem}

See \cite[Theorem 12.85]{Loehr-BC} or \cite{Stucky15} for proofs of this
rather surprising formula (which, alas, require some more advanced methods
than those introduced in this text). Note that it can indeed be used for exact
computation of $d_{n,m}$ (as there are algorithms for exact manipulation of
\textquotedblleft trigonometric irrationals\textquotedblright\footnote{The
technical term is \textquotedblleft elements of cyclotomic
fields\textquotedblright.}\ such as $\cos\dfrac{j\pi}{m+1}$); for example, it
yields $d_{8,8}=12\ 988\ 816$.

\subsection{\label{sec.gf.lim}Limits of FPSs}

Much of what we have been doing with FPSs has been an analogue (a mockery, as
some might say) of analysis, particularly of complex analysis. We have defined
FPSs, derivatives, infinite sums and products, exponentials and various
related concepts, all by mimicking the respective analytic notions but
avoiding any reference to real numbers or any other peculiarities that our
base ring $K$ may or may not have.

Let me introduce yet another such formal analogue of a classical concept from
analysis: that of a limit. Specifically, I will define (coefficientwise)
limits of sequences of FPSs. The exposition will follow \cite[\S 7.5]%
{Loehr-BC}. These limits are very easy to define, yet quite useful. (In
particular, they can be used to give a simpler definition of infinite
products, although a less natural one.)

\subsubsection{Stabilization of scalars}

We start with the most trivial kind of limit whatsoever: the limit of an
eventually constant sequence. Its main advantage is its generality and the
simplicity of its definition:

\begin{definition}
\label{def.fps.lim.stab}Let $\left(  a_{i}\right)  _{i\in\mathbb{N}}=\left(
a_{0},a_{1},a_{2},\ldots\right)  \in K^{\mathbb{N}}$ be a sequence of elements
of $K$. Let $a\in K$.

We say that the sequence $\left(  a_{i}\right)  _{i\in\mathbb{N}}$
\emph{stabilizes to }$a$ if there exists some $N\in\mathbb{N}$ such that%
\[
\text{all integers }i\geq N\text{ satisfy }a_{i}=a.
\]

Instead of saying \textquotedblleft$\left(  a_{i}\right)  _{i\in\mathbb{N}}$
stabilizes to $a$\textquotedblright, we can also say \textquotedblleft$a_{i}$
stabilizes to $a$ as $i\rightarrow\infty$\textquotedblright. (The words
\textquotedblleft as $i\rightarrow\infty$\textquotedblright\ here signify that
we are talking not about a single term $a_{i}$ but about the entire sequence
$\left(  a_{i}\right)  _{i\in\mathbb{N}}$.)

If $a_{i}$ stabilizes to $a$ as $i\rightarrow\infty$, then we write
$\lim\limits_{i\rightarrow\infty}a_{i}=a$ and say that $a$ is the \emph{limit}
(or \emph{eventual value}) of $\left(  a_{i}\right)  _{i\in\mathbb{N}}$ (or,
less precisely, that $a$ is the limit of the $a_{i}$). This is legitimate,
since $a$ is uniquely determined by the sequence $\left(  a_{i}\right)
_{i\in\mathbb{N}}$.

We can replace $\mathbb{N}$ by $\mathbb{Z}_{\geq q}=\left\{  q,q+1,q+2,\ldots
\right\}  $ in this definition, where $q$ is any integer. That is, everything
we said applies not just to sequences of the form $\left(  a_{i}\right)
_{i\in\mathbb{N}}=\left(  a_{0},a_{1},a_{2},\ldots\right)  $, but also to
sequences of the form $\left(  a_{i}\right)  _{i\geq q}=\left(  a_{q}%
,a_{q+1},a_{q+2},\ldots\right)  $.
\end{definition}

Here are some examples and non-examples of stabilization:

\begin{itemize}
\item The sequence%
\[
\left(  \left\lfloor \dfrac{5}{i}\right\rfloor \right)  _{i\geq1}=\left(
\left\lfloor \dfrac{5}{1}\right\rfloor ,\left\lfloor \dfrac{5}{2}\right\rfloor
,\left\lfloor \dfrac{5}{3}\right\rfloor ,\ldots\right)  =\left(
5,2,1,1,1,0,0,0,0,\ldots\right)  \in\mathbb{Z}^{\left\{  1,2,3,\ldots\right\}
}%
\]
stabilizes to $0$, since all integers $i\geq6$ satisfy $\left\lfloor \dfrac
{5}{i}\right\rfloor =0$.

\item On the other hand, the sequence $\left(  \dfrac{5}{i}\right)  _{i\geq
1}\in\mathbb{Q}^{\left\{  1,2,3,\ldots\right\}  }$ does not stabilize to
anything. (It does converge to $0$ in the sense of real analysis, but this is
a weaker statement than stabilization.)

\item Any constant sequence $\left(  a,a,a,\ldots\right)  \in K^{\mathbb{N}}$
stabilizes to $a$.

\item A well-known (if non-constructive) fact says that every weakly
decreasing sequence of nonnegative integers stabilizes to some nonnegative
integer (since it cannot decrease infinitely often).

\item If $s$ is a nilpotent element of a ring $K$ (that is, an element
satisfying $s^{k}=0$ for some $k\geq0$), then the sequence $\left(
s^{i}\right)  _{i\in\mathbb{N}}=\left(  s^{0},s^{1},s^{2},\ldots\right)  $
stabilizes to $0\in K$. For example, the element $\overline{2}$ of the ring
$\mathbb{Z}/8\mathbb{Z}$ is nilpotent (with $\overline{2}^{3}=\overline{0}$),
so that the sequence $\left(  \overline{2}^{i}\right)  _{i\in\mathbb{N}%
}=\left(  \overline{1},\overline{2},\overline{4},\overline{0},\overline
{0},\overline{0},\ldots\right)  $ stabilizes to $\overline{0}$.

\item If $s$ is an idempotent element of a ring $K$ (that is, an element
satisfying $s^{2}=s$), then the sequence $\left(  s^{i}\right)  _{i\in
\mathbb{N}}=\left(  s^{0},s^{1},s^{2},\ldots\right)  $ stabilizes to $s$ (in
fact, $s^{1}=s^{2}=s^{3}=\cdots=s$).
\end{itemize}

\begin{remark}
If you are familiar with point-set topology, then you will recognize
Definition \ref{def.fps.lim.stab} as an instance of topological convergence:
Namely, a sequence $\left(  a_{i}\right)  _{i\in\mathbb{N}}$ stabilizes to
some $a\in K$ if and only if $\left(  a_{i}\right)  _{i\in\mathbb{N}}$
converges to $a$ in the discrete topology.
\end{remark}

\subsubsection{Coefficientwise stabilization of FPSs}

Now we define a somewhat subtler notion of limits, one specific to FPSs (as
opposed to elements of an arbitrary commutative ring). This is the notion that
we will use in what follows:

\begin{definition}
\label{def.fps.lim.coeff-stab}Let $\left(  f_{i}\right)  _{i\in\mathbb{N}}\in
K\left[  \left[  x\right]  \right]  ^{\mathbb{N}}$ be a sequence of FPSs over
$K$. Let $f\in K\left[  \left[  x\right]  \right]  $ be an FPS.

We say that $\left(  f_{i}\right)  _{i\in\mathbb{N}}$ \emph{coefficientwise
stabilizes to }$f$ if for each $n\in\mathbb{N}$,%
\[
\text{the sequence }\left(  \left[  x^{n}\right]  f_{i}\right)  _{i\in
\mathbb{N}}\text{ stabilizes to }\left[  x^{n}\right]  f.
\]

Instead of saying \textquotedblleft$\left(  f_{i}\right)  _{i\in\mathbb{N}}$
coefficientwise stabilizes to $f$\textquotedblright, we can also say
\textquotedblleft$f_{i}$ coefficientwise stabilizes to $f$ as $i\rightarrow
\infty$\textquotedblright\ or write \textquotedblleft$f_{i}\rightarrow f$ as
$i\rightarrow\infty$\textquotedblright.

If $f_{i}$ coefficientwise stabilizes to $f$ as $i\rightarrow\infty$, then we
write $\lim\limits_{i\rightarrow\infty}f_{i}=f$ and say that $f$ is the
\emph{limit} of $\left(  f_{i}\right)  _{i\in\mathbb{N}}$ (or, less precisely,
that $f$ is the limit of the $f_{i}$). This is legitimate, because $f$ is
uniquely determined by the sequence $\left(  f_{i}\right)  _{i\in\mathbb{N}}$.

Again, we can replace $\mathbb{N}$ by $\mathbb{Z}_{\geq q}=\left\{
q,q+1,q+2,\ldots\right\}  $ in this definition, where $q$ is any integer.
\end{definition}

\begin{example}
\textbf{(a)} The sequence $\left(  x^{i}\right)  _{i\in\mathbb{N}}$ of FPSs
coefficientwise stabilizes to the FPS $0$. (Indeed, for each $n\in\mathbb{N}$,
the sequence $\left(  \left[  x^{n}\right]  \left(  x^{i}\right)  \right)
_{i\in\mathbb{N}}$ consists of a single $1$ and infinitely many $0$s; thus,
this sequence stabilizes to $0$.) In other words, we have $\lim
\limits_{i\rightarrow\infty}x^{i}=0$. \medskip

\textbf{(b)} The sequence $\left(  \dfrac{1}{i}x\right)  _{i\geq1}$ of FPSs
over the ring $\mathbb{Q}$ does not stabilize (since the $x^{1}$-coefficients
of these FPSs $\dfrac{1}{i}x$ never stabilize). Thus, $\lim
\limits_{i\rightarrow\infty}\left(  \dfrac{1}{i}x\right)  $ does not exist.
\medskip

\textbf{(c)} We have%
\[
\lim\limits_{i\rightarrow\infty}\left(  \left(  1+x^{1}\right)  \left(
1+x^{2}\right)  \cdots\left(  1+x^{i}\right)  \right)  =\prod_{k=1}^{\infty
}\left(  1+x^{k}\right)  .
\]
Indeed, the $x^{n}$-coefficient of the product $\left(  1+x^{1}\right)
\left(  1+x^{2}\right)  \cdots\left(  1+x^{i}\right)  $ depends only on those
factors in which the exponent is $\leq n$, and thus stops changing after $i$
surpasses $n$; therefore, it stabilizes to the $x^{n}$-coefficient of the
infinite product $\prod_{k=1}^{\infty}\left(  1+x^{k}\right)  $. \medskip

\textbf{(d)} It would be nice to have $\lim\limits_{n\rightarrow\infty}\left(
1+\dfrac{x}{n}\right)  ^{n}=\exp$ (over the ring $\mathbb{Q}$), as in real
analysis. Unfortunately, this is not the case. In fact, the binomial formula
yields
\[
\left(  1+\dfrac{x}{n}\right)  ^{n}=1+\underbrace{n\cdot\dfrac{x}{n}}%
_{=x}+\,\dbinom{n}{2}\cdot\left(  \dfrac{x}{n}\right)  ^{2}+\cdots
=1+x+\dfrac{\dbinom{n}{2}}{n^{2}}x^{2}+\cdots.
\]
This shows that the $x^{0}$-coefficient and the $x^{1}$-coefficient of
$\left(  1+\dfrac{x}{n}\right)  ^{n}$ stabilize as $n\rightarrow\infty$, but
the $x^{2}$-coefficient does not. Thus, $\lim\limits_{n\rightarrow\infty
}\left(  1+\dfrac{x}{n}\right)  ^{n}$ does not exist (according to our
definition of limit). \medskip

\textbf{(e)} Coefficientwise stabilization is weaker than stabilization: If a
sequence $\left(  f_{0},f_{1},f_{2},\ldots\right)  $ of FPSs stabilizes to an
FPS $f$ (meaning that $f_{i}=f$ for all sufficiently high $i\in\mathbb{N}$),
then it coefficientwise stabilizes to $f$ as well. In particular, any constant
sequence $\left(  f,f,f,\ldots\right)  $ of FPSs stabilizes to $f$.
\end{example}

\begin{remark}
If you are familiar with point-set topology, then you will again recognize
Definition \ref{def.fps.lim.coeff-stab} as an instance of topological
convergence: Namely, we recall that each FPS is an infinite sequence of
elements of $K$ (its coefficients). Thus, the set $K\left[  \left[  x\right]
\right]  $ is the infinite Cartesian product $K\times K\times K\times\cdots$.
If we equip each factor $K$ in this product with the discrete topology, then
the entire product $K\left[  \left[  x\right]  \right]  $ becomes a product
space equipped with the product topology. Now, a sequence $\left(
f_{i}\right)  _{i\in\mathbb{N}}$ of FPSs in $K\left[  \left[  x\right]
\right]  $ coefficientwise stabilizes to some FPS $f\in K\left[  \left[
x\right]  \right]  $ if and only if $\left(  f_{i}\right)  _{i\in\mathbb{N}}$
converges to $f$ in this product topology. (This is not surprising: After all,
the product topology is also known as the \textquotedblleft topology of
pointwise convergence\textquotedblright, and in our case \textquotedblleft
pointwise\textquotedblright\ means \textquotedblleft
coefficientwise\textquotedblright.)
\end{remark}

\subsubsection{Some properties of limits}

We shall now state some basic properties of limits of FPSs. All the proofs are
simple and therefore omitted; some can be found in the Appendix (Section
\ref{sec.details.gf.lim}).

\begin{theorem}
\label{thm.fps.lim.lim-crit}Let $\left(  f_{i}\right)  _{i\in\mathbb{N}}\in
K\left[  \left[  x\right]  \right]  ^{\mathbb{N}}$ be a sequence of FPSs.
Assume that for each $n\in\mathbb{N}$, there exists some $g_{n}\in K$ such
that the sequence $\left(  \left[  x^{n}\right]  f_{i}\right)  _{i\in
\mathbb{N}}$ stabilizes to $g_{n}$. Then, the sequence $\left(  f_{i}\right)
_{i\in\mathbb{N}}$ coefficientwise stabilizes to $\sum_{n\in\mathbb{N}}%
g_{n}x^{n}$. (That is, $\lim\limits_{i\rightarrow\infty}f_{i}=\sum
_{n\in\mathbb{N}}g_{n}x^{n}$.)
\end{theorem}

\begin{proof}
Obvious.
\end{proof}

The following easy lemma connects limits with the notion of $x^{n}%
$-equivalence (as introduced in Definition \ref{def.fps.xneq}):

\begin{lemma}
\label{lem.fps.lim.xn-equiv}Assume that $\left(  f_{i}\right)  _{i\in
\mathbb{N}}$ is a sequence of FPSs, and that $f$ is an FPS such that
$\lim\limits_{i\rightarrow\infty}f_{i}=f$. Then, for each $n\in\mathbb{N}$,
there exists some integer $N\in\mathbb{N}$ such that
\[
\text{all integers }i\geq N\text{ satisfy }f_{i}\overset{x^{n}}{\equiv}f.
\]

\end{lemma}

\begin{proof}
[Proof of Lemma \ref{lem.fps.lim.xn-equiv} (sketched).]Let $n\in\mathbb{N}$.
For each $k\in\left\{  0,1,\ldots,n\right\}  $, we pick an $N_{k}\in
\mathbb{N}$ such that all $i\geq N_{k}$ satisfy $\left[  x^{k}\right]
f_{i}=\left[  x^{k}\right]  f$ (such an $N_{k}$ exists, since $\lim
\limits_{i\rightarrow\infty}f_{i}=f$). Then, all integers $i\geq\max\left\{
N_{0},N_{1},\ldots,N_{n}\right\}  $ satisfy $f_{i}\overset{x^{n}}{\equiv}f$.
See Section \ref{sec.details.gf.lim} for the details of this proof.
\end{proof}

The following proposition is an analogue of the classical \textquotedblleft
limits respect sums and products\textquotedblright\ theorem from real analysis:

\begin{proposition}
\label{prop.fps.lim.sum-prod}Assume that $\left(  f_{i}\right)  _{i\in
\mathbb{N}}$ and $\left(  g_{i}\right)  _{i\in\mathbb{N}}$ are two sequences
of FPSs, and that $f$ and $g$ are two FPSs such that%
\[
\lim\limits_{i\rightarrow\infty}f_{i}=f\ \ \ \ \ \ \ \ \ \ \text{and}%
\ \ \ \ \ \ \ \ \ \ \lim\limits_{i\rightarrow\infty}g_{i}=g.
\]

Then,
\[
\lim\limits_{i\rightarrow\infty}\left(  f_{i}+g_{i}\right)
=f+g\ \ \ \ \ \ \ \ \ \ \text{and}\ \ \ \ \ \ \ \ \ \ \lim
\limits_{i\rightarrow\infty}\left(  f_{i}g_{i}\right)  =fg.
\]
(In other words, the sequences $\left(  f_{i}+g_{i}\right)  _{i\in\mathbb{N}}$
and $\left(  f_{i}g_{i}\right)  _{i\in\mathbb{N}}$ coefficientwise stabilize
to $f+g$ and $fg$, respectively.)
\end{proposition}

\begin{proof}
[Proof of Proposition \ref{prop.fps.lim.sum-prod} (sketched).]Let
$n\in\mathbb{N}$. Then, Lemma \ref{lem.fps.lim.xn-equiv} shows that there
exists some integer $N\in\mathbb{N}$ such that
\[
\text{all integers }i\geq N\text{ satisfy }f_{i}\overset{x^{n}}{\equiv}f.
\]
Similarly, there exists some integer $M\in\mathbb{N}$ such that%
\[
\text{all integers }i\geq M\text{ satisfy }g_{i}\overset{x^{n}}{\equiv}g.
\]
Pick such integers $N$ and $M$, and set $P:=\max\left\{  N,M\right\}  $. Then,
all integers $i\geq P$ satisfy both $f_{i}\overset{x^{n}}{\equiv}f$ and
$g_{i}\overset{x^{n}}{\equiv}g$, and therefore $f_{i}g_{i}\overset{x^{n}%
}{\equiv}fg$ (by (\ref{eq.thm.fps.xneq.props.b.*})), and thus $\left[
x^{n}\right]  \left(  f_{i}g_{i}\right)  =\left[  x^{n}\right]  \left(
fg\right)  $. This shows that the sequence $\left(  \left[  x^{n}\right]
\left(  f_{i}g_{i}\right)  \right)  _{i\in\mathbb{N}}$ stabilizes to $\left[
x^{n}\right]  \left(  fg\right)  $. Since this holds for each $n\in\mathbb{N}%
$, we conclude that $\lim\limits_{i\rightarrow\infty}\left(  f_{i}%
g_{i}\right)  =fg$. The proof of $\lim\limits_{i\rightarrow\infty}\left(
f_{i}+g_{i}\right)  =f+g$ is analogous.

See Section \ref{sec.details.gf.lim} for the details of this proof.
\end{proof}

\begin{corollary}
\label{cor.fps.lim.sum-prod-k}Let $k\in\mathbb{N}$. For each $i\in\left\{
1,2,\ldots,k\right\}  $, let $f_{i}$ be an FPS, and let $\left(
f_{i,n}\right)  _{n\in\mathbb{N}}$ be a sequence of FPSs such that%
\[
\lim\limits_{n\rightarrow\infty}\left(  f_{i,n}\right)  =f_{i}%
\]
(note that it is $n$, not $i$, that goes to $\infty$ here!). Then,%
\[
\lim\limits_{n\rightarrow\infty}\sum_{i=1}^{k}f_{i,n}=\sum_{i=1}^{k}%
f_{i}\ \ \ \ \ \ \ \ \ \ \text{and}\ \ \ \ \ \ \ \ \ \ \lim
\limits_{n\rightarrow\infty}\prod_{i=1}^{k}f_{i,n}=\prod_{i=1}^{k}f_{i}.
\]

\end{corollary}

\begin{proof}
Follows by induction on $k$, using Proposition \ref{prop.fps.lim.sum-prod}.
\end{proof}

\begin{remark}
The analogue of Corollary \ref{cor.fps.lim.sum-prod-k} for infinite sums (or
products) is not true without further requirements. For instance, for each
$i\in\mathbb{N}$, we can define a constant FPS $f_{i,n}=\delta_{i,n}$ (using
Definition \ref{def.kron-delta}). Then, we have $\lim\limits_{n\rightarrow
\infty}f_{i,n}=0$ for any fixed $i\in\mathbb{N}$, but $\lim
\limits_{n\rightarrow\infty}\underbrace{\sum_{i=0}^{\infty}f_{i,n}}_{=1}%
=\lim\limits_{n\rightarrow\infty}1=1$.
\end{remark}

The following proposition says that limits respect quotients:

\begin{proposition}
\label{prop.fps.lim.sum-quot}Assume that $\left(  f_{i}\right)  _{i\in
\mathbb{N}}$ and $\left(  g_{i}\right)  _{i\in\mathbb{N}}$ are two sequences
of FPSs, and that $f$ and $g$ are two FPSs such that%
\[
\lim\limits_{i\rightarrow\infty}f_{i}=f\ \ \ \ \ \ \ \ \ \ \text{and}%
\ \ \ \ \ \ \ \ \ \ \lim\limits_{i\rightarrow\infty}g_{i}=g.
\]
Assume that each FPS $g_{i}$ is invertible. Then, $g$ is also invertible, and
we have%
\[
\lim\limits_{i\rightarrow\infty}\dfrac{f_{i}}{g_{i}}=\dfrac{f}{g}.
\]

\end{proposition}

\begin{proof}
[Proof of Proposition \ref{prop.fps.lim.sum-quot} (sketched).]First use
Proposition \ref{prop.fps.invertible} to show that $g$ is invertible; then use
Theorem \ref{thm.fps.xneq.props} \textbf{(e)} to prove $\lim
\limits_{i\rightarrow\infty}\dfrac{f_{i}}{g_{i}}=\dfrac{f}{g}$. A detailed
proof can be found in Section \ref{sec.details.gf.lim}.
\end{proof}

Limits of FPSs furthermore respect composition:

\begin{proposition}
\label{prop.fps.lim.comp}Assume that $\left(  f_{i}\right)  _{i\in\mathbb{N}}$
and $\left(  g_{i}\right)  _{i\in\mathbb{N}}$ are two sequences of FPSs, and
that $f$ and $g$ are two FPSs such that%
\[
\lim\limits_{i\rightarrow\infty}f_{i}=f\ \ \ \ \ \ \ \ \ \ \text{and}%
\ \ \ \ \ \ \ \ \ \ \lim\limits_{i\rightarrow\infty}g_{i}=g.
\]
Assume that $\left[  x^{0}\right]  g_{i}=0$ for each $i\in\mathbb{N}$. Then,
$\left[  x^{0}\right]  g=0$ and%
\[
\lim\limits_{i\rightarrow\infty}\left(  f_{i}\circ g_{i}\right)  =f\circ g.
\]

\end{proposition}

\begin{proof}
[Proof of Proposition \ref{prop.fps.lim.comp} (sketched).]Similar to
Proposition \ref{prop.fps.lim.sum-quot}, but now using Proposition
\ref{prop.fps.xneq.comp}. See Exercise \ref{exe.fps.lim.sum-prod-quot}
\textbf{(b)} for details.
\end{proof}

Limits of FPSs also respect derivatives (unlike in real analysis, where this
\href{https://en.wikipedia.org/wiki/Uniform_convergence#To_differentiability}{holds
only under subtle additional conditions}):

\begin{proposition}
\label{prop.fps.lim.deriv-lim}Let $\left(  f_{n}\right)  _{n\in\mathbb{N}}$ be
a sequence of FPSs, and let $f$ be an FPS such that%
\[
\lim\limits_{i\rightarrow\infty}f_{i}=f.
\]
Then,%
\[
\lim\limits_{i\rightarrow\infty}f_{i}^{\prime}=f^{\prime}.
\]

\end{proposition}

\begin{proof}
This is nearly trivial (see Exercise \ref{exe.fps.lim.sum-prod-quot}
\textbf{(a)}).
\end{proof}

Next, let us see how infinite sums and infinite products can be written as
limits of finite sums and products, at least as long as they are indexed by
$\mathbb{N}$. The following two theorems and one corollary are easy to prove;
detailed proofs can be found in Section \ref{sec.details.gf.lim}.

\begin{theorem}
\label{thm.fps.lim.sum-lim}Let $\left(  f_{n}\right)  _{n\in\mathbb{N}}$ be a
summable sequence of FPSs. Then,%
\[
\lim\limits_{i\rightarrow\infty}\sum_{n=0}^{i}f_{n}=\sum_{n\in\mathbb{N}}%
f_{n}.
\]
In other words, the infinite sum $\sum_{n\in\mathbb{N}}f_{n}$ is the limit of
the finite partial sums $\sum_{n=0}^{i}f_{n}$.
\end{theorem}

\begin{theorem}
\label{thm.fps.lim.prod-lim}Let $\left(  f_{n}\right)  _{n\in\mathbb{N}}$ be a
multipliable sequence of FPSs. Then,%
\[
\lim\limits_{i\rightarrow\infty}\prod_{n=0}^{i}f_{n}=\prod_{n\in\mathbb{N}%
}f_{n}.
\]
In other words, the infinite product $\prod_{n\in\mathbb{N}}f_{n}$ is the
limit of the finite partial products $\prod_{n=0}^{i}f_{n}$.
\end{theorem}

\begin{corollary}
\label{cor.fps.lim.fps-as-pol}Each FPS is a limit of a sequence of
polynomials. Indeed, if $a=\sum_{n\in\mathbb{N}}a_{n}x^{n}$ (with $a_{n}\in
K$), then%
\[
a=\lim\limits_{i\rightarrow\infty}\sum_{n=0}^{i}a_{n}x^{n}.
\]

\end{corollary}

This corollary can be restated as \textquotedblleft the polynomials are dense
in the FPSs\textquotedblright\ (more formally: $K\left[  x\right]  $ is dense
in $K\left[  \left[  x\right]  \right]  $). This fact is useful, as it allows
you to restrict yourself to polynomials when proving some properties of FPSs.
For example, if you want to prove that some identity (such as the Leibniz rule
$\left(  fg\right)  ^{\prime}=f^{\prime}g+fg^{\prime}$) holds for all FPSs, it
suffices to prove it for polynomials, and then conclude the general case using
a limiting argument.

Both Theorem \ref{thm.fps.lim.sum-lim} and Theorem \ref{thm.fps.lim.prod-lim}
have converses (which are somewhat harder to prove):

\begin{theorem}
\label{thm.fps.lim.sum-lim-conv}Let $\left(  f_{0},f_{1},f_{2},\ldots\right)
$ be a sequence of FPSs such that $\lim\limits_{i\rightarrow\infty}\sum
_{n=0}^{i}f_{n}$ exists. Then, the family $\left(  f_{n}\right)
_{n\in\mathbb{N}}$ is summable, and satisfies%
\[
\sum_{n\in\mathbb{N}}f_{n}=\lim\limits_{i\rightarrow\infty}\sum_{n=0}^{i}%
f_{n}.
\]

\end{theorem}

\begin{theorem}
\label{thm.fps.lim.prod-lim-conv}Let $\left(  f_{0},f_{1},f_{2},\ldots\right)
$ be a sequence of FPSs such that $\lim\limits_{i\rightarrow\infty}\prod
_{n=0}^{i}f_{n}$ exists. Then, the family $\left(  f_{n}\right)
_{n\in\mathbb{N}}$ is multipliable, and satisfies%
\[
\prod_{n\in\mathbb{N}}f_{n}=\lim\limits_{i\rightarrow\infty}\prod_{n=0}%
^{i}f_{n}.
\]

\end{theorem}

We refer to Section \ref{sec.details.gf.lim} for the proofs of these two theorems.

Theorem \ref{thm.fps.lim.sum-lim} and Theorem \ref{thm.fps.lim.sum-lim-conv}
combined show that we could have just as well defined infinite sums of FPSs as
limits of finite partial sums, as long as we were happy with sums of the form
$\sum_{n\in\mathbb{N}}f_{n}$ (as opposed to sums of the form $\sum_{\left(
n,m\right)  \in\mathbb{N}\times\mathbb{N}}$ or other kinds). Likewise, Theorem
\ref{thm.fps.lim.prod-lim} and Theorem \ref{thm.fps.lim.prod-lim-conv} show
the same for products. However, I consider my original definitions to be more
natural. \medskip

It should be clear that all the above-stated properties of limits remain true
if the sequences involved are no longer indexed by $\mathbb{N}$ but rather
indexed by $\mathbb{Z}_{\geq q}=\left\{  q,q+1,q+2,\ldots\right\}  $ for some
integer $q$. Even better: The limit of a sequence $\left(  f_{n}\right)
_{n\geq q}$ does not depend on any chosen finite piece of this sequence, so it
is simultaneously the limit of the sequence $\left(  f_{n}\right)  _{n\geq r}$
for any integer $r$ (as long as $f_{n}$ is well-defined for all $n\geq r$). In
this regard, limits behave exactly as in real analysis.

\subsection{\label{sec.gf.laure}Laurent power series}

\subsubsection{Motivation}

Next, we shall try to extend the concept of FPSs to allow negative powers of
$x$. First, however, let me motivate this by showing an example where the need
for such negative powers arises.

First, we recall the binary positional system for nonnegative integers:

\begin{definition}
A \emph{binary representation} of an integer $n$ means an essentially finite
sequence $\left(  b_{i}\right)  _{i\in\mathbb{N}}=\left(  b_{0},b_{1}%
,b_{2},\ldots\right)  \in\left\{  0,1\right\}  ^{\mathbb{N}}$ such that%
\[
n=\sum_{i\in\mathbb{N}}b_{i}2^{i}.
\]
(Recall that \textquotedblleft essentially finite\textquotedblright\ means
\textquotedblleft all but finitely many $i\in\mathbb{N}$ satisfy $b_{i}%
=0$\textquotedblright.)
\end{definition}

The following theorem is well-known (and has been proved in Subsection
\ref{subsec.gf.prod.exa}):

\begin{theorem}
\label{thm.fps.laure.binary-rep-uniq}Each $n\in\mathbb{N}$ has a unique binary representation.
\end{theorem}

Note that we are encoding the digits (actually, bits) of a binary
representation as essentially finite sequences instead of finite tuples. This
way, we don't have to worry about leading zeros breaking the uniqueness in
Theorem \ref{thm.fps.laure.binary-rep-uniq}. \medskip

Let us now define a variation of binary representation:

\begin{definition}
A \emph{balanced ternary representation} of an integer $n$ means an
essentially finite sequence $\left(  b_{i}\right)  _{i\in\mathbb{N}}=\left(
b_{0},b_{1},b_{2},\ldots\right)  \in\left\{  0,1,-1\right\}  ^{\mathbb{N}}$
such that%
\[
n=\sum_{i\in\mathbb{N}}b_{i}3^{i}.
\]

\end{definition}

Here are some examples:

\begin{itemize}
\item The integer $19$ has a balanced ternary representation $\left(
1,0,-1,1,0,0,0,\ldots\right)  $, because%
\[
19=1-9+27=3^{0}-3^{2}+3^{3}=1\cdot3^{0}+\left(  -1\right)  \cdot3^{2}%
+1\cdot3^{3}.
\]

\item The integer $42$ has a balanced ternary representation $\left(
0,-1,-1,-1,1,0,0,0,\ldots\right)  $, because%
\[
42=81-27-9-3=3^{4}-3^{3}-3^{2}-3^{1}.
\]

\item Note that (unlike with binary representations) even negative integers
can have balanced ternary representations. For example, the integer $-11$ has
a balanced ternary representation $\left(  1,-1,-1,0,0,0,\ldots\right)  $.
\end{itemize}

In the Soviet Union of the 1960s/70s, balanced ternary representations have
been used as a foundation for computers (see
\href{https://en.wikipedia.org/wiki/Setun}{the Setun computer}). The idea was
shelved in the 1970s, but a noticeable number of algorithms have been invented
for working with balanced ternary representations. (See \cite[\S 4.1]%
{Knuth-TAoCP2} for some discussion of these.) The following theorem (which
goes back to Fibonacci) is crucial for making balanced ternary representations useful:

\begin{theorem}
\label{thm.fps.laure.balanced-tern-rep-uniq}Each integer $n$ has a unique
balanced ternary representation.
\end{theorem}

There are various ways to prove this (see, e.g., \cite[solution to Exercise
3.7.8]{20f} for an elementary one). Let us here try to prove it using FPSs
(imitating our proof of Theorem \ref{thm.fps.laure.binary-rep-uniq} in
Subsection \ref{subsec.gf.prod.exa}). Since the $b_{i}$ can be $-1$s, we must
allow for negative powers of $x$.

Let us first argue informally; we will later see whether we can make sense of
what we have done. The following informal argument was proposed by Euler
(\cite[\S 331]{Euler48}). We shall compute the product%
\[
\left(  1+x+x^{-1}\right)  \left(  1+x^{3}+x^{-3}\right)  \left(
1+x^{9}+x^{-9}\right)  \cdots=\prod_{i\geq0}\left(  1+x^{3^{i}}+x^{-3^{i}%
}\right)
\]
in two ways:

\begin{itemize}
\item On the one hand, this product equals%
\begin{align}
&  \prod_{i\geq0}\underbrace{\left(  1+x^{3^{i}}+x^{-3^{i}}\right)  }%
_{=\sum_{b\in\left\{  0,1,-1\right\}  }x^{b\cdot3^{i}}}\nonumber\\
&  =\prod_{i\geq0}\ \ \sum_{b\in\left\{  0,1,-1\right\}  }x^{b\cdot3^{i}%
}\nonumber\\
&  =\sum_{\substack{\left(  b_{0},b_{1},b_{2},\ldots\right)  \in\left\{
0,1,-1\right\}  ^{\mathbb{N}}\\\text{is essentially finite}}}x^{b_{0}3^{0}%
}x^{b_{1}3^{1}}x^{b_{2}3^{2}}\cdots\nonumber\\
&  \ \ \ \ \ \ \ \ \ \ \ \ \ \ \ \ \ \ \ \ \left(
\begin{array}
[c]{c}%
\text{here, we have just expanded the product using
(\ref{eq.prop.fps.prodrule-inf-inf.eq})}\\
\text{(hoping that (\ref{eq.prop.fps.prodrule-inf-inf.eq}) still works in our
setting)}%
\end{array}
\right) \nonumber\\
&  =\sum_{\substack{\left(  b_{0},b_{1},b_{2},\ldots\right)  \in\left\{
0,1,-1\right\}  ^{\mathbb{N}}\\\text{is essentially finite}}}x^{b_{0}%
3^{0}+b_{1}3^{1}+b_{2}3^{2}+\cdots}\nonumber\\
&  =\sum_{n\in\mathbb{Z}}\ \ \sum_{\substack{\left(  b_{0},b_{1},b_{2}%
,\ldots\right)  \in\left\{  0,1,-1\right\}  ^{\mathbb{N}}\\\text{is
essentially finite;}\\b_{0}3^{0}+b_{1}3^{1}+b_{2}3^{2}+\cdots=n}%
}x^{n}\nonumber\\
&  =\sum_{n\in\mathbb{Z}}\left(  \text{\# of balanced ternary representations
of }n\right)  \cdot x^{n}, \label{eq.thm.fps.laure.balanced-tern-rep-uniq.1}%
\end{align}
since a balanced ternary representation of $n$ is precisely an essentially
finite sequence $\left(  b_{0},b_{1},b_{2},\ldots\right)  \in\left\{
0,1,-1\right\}  ^{\mathbb{N}}$ satisfying $b_{0}3^{0}+b_{1}3^{1}+b_{2}%
3^{2}+\cdots=n$.

\item On the other hand, we have $1+x+x^{-1}=\dfrac{1-x^{3}}{x\left(
1-x\right)  }$. Substituting $x^{3^{i}}$ for $x$ in this equality, we obtain%
\[
1+x^{3^{i}}+x^{-3^{i}}=\dfrac{1-\left(  x^{3^{i}}\right)  ^{3}}{x^{3^{i}%
}\left(  1-x^{3^{i}}\right)  }=\dfrac{1-x^{3^{i+1}}}{x^{3^{i}}\left(
1-x^{3^{i}}\right)  }\ \ \ \ \ \ \ \ \ \ \text{for each }i\geq0.
\]
Hence,
\begin{align}
&  \prod_{i\geq0}\left(  1+x^{3^{i}}+x^{-3^{i}}\right) \nonumber\\
&  =\prod_{i\geq0}\dfrac{1-x^{3^{i+1}}}{x^{3^{i}}\left(  1-x^{3^{i}}\right)
}\nonumber\\
&  =\dfrac{1-x^{3}}{x\left(  1-x\right)  }\cdot\dfrac{1-x^{9}}{x^{3}\left(
1-x^{3}\right)  }\cdot\dfrac{1-x^{27}}{x^{9}\left(  1-x^{9}\right)  }%
\cdot\dfrac{1-x^{81}}{x^{27}\left(  1-x^{27}\right)  }\cdot\cdots\nonumber\\
&  =\underbrace{\dfrac{1}{xx^{3}x^{9}x^{27}\cdots}}_{\substack{=x^{-\infty
}\\\text{(whatever this means)}}}\cdot\underbrace{\dfrac{1}{1-x}}%
_{=1+x+x^{2}+x^{3}+\cdots}\nonumber\\
&  \ \ \ \ \ \ \ \ \ \ \ \ \ \ \ \ \ \ \ \ \left(
\begin{array}
[c]{c}%
\text{here we cancelled the factors }1-x^{3},\ \ 1-x^{9},\ \ 1-x^{27}%
,\ \ \ldots\\
\text{by a somewhat daring use of the telescope principle}%
\end{array}
\right) \nonumber\\
&  =x^{-\infty}\left(  1+x+x^{2}+x^{3}+\cdots\right) \nonumber\\
&  =\cdots+x^{-2}+x^{-1}+x^{0}+x^{1}+x^{2}+\cdots\nonumber\\
&  \ \ \ \ \ \ \ \ \ \ \ \ \ \ \ \ \ \ \ \ \left(
\begin{array}
[c]{c}%
\text{with some artistic license,}\\
\text{since }x^{i}\left(  1+x+x^{2}+x^{3}+\cdots\right)  =x^{i}+x^{i+1}%
+x^{i+2}+\cdots\\
\text{for each }i\in\mathbb{Z}%
\end{array}
\right) \nonumber\\
&  =\sum_{n\in\mathbb{Z}}x^{n}.
\label{eq.thm.fps.laure.balanced-tern-rep-uniq.2}%
\end{align}

\end{itemize}

Comparing (\ref{eq.thm.fps.laure.balanced-tern-rep-uniq.1}) with
(\ref{eq.thm.fps.laure.balanced-tern-rep-uniq.2}), we find%
\[
\sum_{n\in\mathbb{Z}}\left(  \text{\# of balanced ternary representations of
}n\right)  \cdot x^{n}=\sum_{n\in\mathbb{Z}}x^{n}.
\]
Comparing coefficients, we thus conclude that%
\[
\left(  \text{\# of balanced ternary representations of }n\right)  =1
\]
for each $n\in\mathbb{Z}$. This \textquotedblleft proves\textquotedblright%
\ Theorem \ref{thm.fps.laure.balanced-tern-rep-uniq}; we just need to make our
computations rigorous -- i.e., define the ring in which we have been
computing, explain what $x$ is, and justify the well-definedness of our
infinite products and sums.

Let us first play around a bit further. We have%
\begin{align*}
&  \left(  1-x\right)  \left(  \cdots+x^{-2}+x^{-1}+x^{0}+x^{1}+x^{2}%
+\cdots\right) \\
&  =\left(  \cdots+x^{-2}+x^{-1}+x^{0}+x^{1}+x^{2}+\cdots\right)
-\underbrace{x\left(  \cdots+x^{-2}+x^{-1}+x^{0}+x^{1}+x^{2}+\cdots\right)
}_{\substack{=\cdots+x^{-1}+x^{0}+x^{1}+x^{2}+x^{3}+\cdots\\=\cdots
+x^{-2}+x^{-1}+x^{0}+x^{1}+x^{2}+\cdots}}\\
&  =\left(  \cdots+x^{-2}+x^{-1}+x^{0}+x^{1}+x^{2}+\cdots\right)  -\left(
\cdots+x^{-2}+x^{-1}+x^{0}+x^{1}+x^{2}+\cdots\right) \\
&  =0.
\end{align*}
Thus, dividing by $1-x$, we obtain%
\[
\cdots+x^{-2}+x^{-1}+x^{0}+x^{1}+x^{2}+\cdots=0.
\]
Comparing coefficients, we conclude that%
\[
1=0\ \ \ \ \ \ \ \ \ \ \text{for each }n\in\mathbb{Z}.
\]
Oops! Looks like we have overtaxed our artistic license.

So we need to be careful with negative powers of $x$. Not everything that
looks like a valid computation actually is one. Thus, we need to be rigorous
and delimit what can and what cannot be done with negative powers of $x$.

\subsubsection{The space $K\left[  \left[  x^{\pm}\right]  \right]  $}

Let us try to define \textquotedblleft FPSs with negative powers of
$x$\textquotedblright\ formally. First, we define the largest possible space
of such FPSs:

\begin{definition}
\label{def.fps.laure.double}Let $K\left[  \left[  x^{\pm}\right]  \right]  $
be the $K$-module $K^{\mathbb{Z}}$ of all families $\left(  a_{n}\right)
_{n\in\mathbb{Z}}=\left(  \ldots,a_{-2},a_{-1},a_{0},a_{1},a_{2}%
,\ldots\right)  $ of elements of $K$. Its addition and its scaling are defined
entrywise:%
\begin{align*}
\left(  a_{n}\right)  _{n\in\mathbb{Z}}+\left(  b_{n}\right)  _{n\in
\mathbb{Z}}  &  =\left(  a_{n}+b_{n}\right)  _{n\in\mathbb{Z}};\\
\lambda\left(  a_{n}\right)  _{n\in\mathbb{Z}}  &  =\left(  \lambda
a_{n}\right)  _{n\in\mathbb{Z}}\ \ \ \ \ \ \ \ \ \ \text{for each }\lambda\in
K.
\end{align*}
An element of $K\left[  \left[  x^{\pm}\right]  \right]  $ will be called a
\emph{doubly infinite power series}. This name is justified by the fact that
we will later use the notation $\sum_{n\in\mathbb{Z}}a_{n}x^{n}$ for a family
$\left(  a_{n}\right)  _{n\in\mathbb{Z}}\in K\left[  \left[  x^{\pm}\right]
\right]  $.
\end{definition}

Now, let us try to define a multiplication on this $K$-module $K\left[
\left[  x^{\pm}\right]  \right]  $, in order to turn it into a $K$-algebra
(like $K\left[  \left[  x\right]  \right]  $). This multiplication should
satisfy%
\[
\left(  a_{n}\right)  _{n\in\mathbb{Z}}\cdot\left(  b_{n}\right)
_{n\in\mathbb{Z}}=\left(  c_{n}\right)  _{n\in\mathbb{Z}}%
,\ \ \ \ \ \ \ \ \ \ \text{where}\ \ \ \ \ \ \ \ \ \ c_{n}=\sum_{i\in
\mathbb{Z}}a_{i}b_{n-i}%
\]
(since this is what we would get if we expanded $\left(  \sum_{n\in\mathbb{Z}%
}a_{n}x^{n}\right)  \left(  \sum_{n\in\mathbb{Z}}b_{n}x^{n}\right)  $ and
combined like powers of $x$). Unfortunately, the sum $\sum_{i\in\mathbb{Z}%
}a_{i}b_{n-i}$ is now infinite (unlike for $K\left[  \left[  x\right]
\right]  $), and is not always well-defined. Thus, the product $\left(
a_{n}\right)  _{n\in\mathbb{Z}}\cdot\left(  b_{n}\right)  _{n\in\mathbb{Z}}$
of two elements of $K\left[  \left[  x^{\pm}\right]  \right]  $ does not
always exist\footnote{The simplest example for this is $\left(  1\right)
_{n\in\mathbb{Z}}\cdot\left(  1\right)  _{n\in\mathbb{Z}}$. If this product
existed, then it should be $\left(  c_{n}\right)  _{n\in\mathbb{Z}}$, where
$c_{n}=\sum_{i\in\mathbb{Z}}1\cdot1$, but the latter sum clearly does not
exist.}. Therefore, the $K$-module $K\left[  \left[  x^{\pm}\right]  \right]
$ is not a $K$-algebra. This explains why our above computations have led us astray.

\subsubsection{Laurent polynomials}

If we cannot multiply two arbitrary elements of $K\left[  \left[  x^{\pm
}\right]  \right]  $, can we perhaps restrict ourselves to a smaller
$K$-submodule of $K\left[  \left[  x^{\pm}\right]  \right]  $ whose elements
can be multiplied?

One such submodule is $K\left[  \left[  x\right]  \right]  $, which is
embedded in $K\left[  \left[  x^{\pm}\right]  \right]  $ in the
\textquotedblleft obvious\textquotedblright\ way (by identifying each FPS
$\left(  a_{n}\right)  _{n\in\mathbb{N}}$ with the doubly infinite power
series $\left(  a_{n}\right)  _{n\in\mathbb{Z}}$, where we set $a_{n}:=0$ for
all $n<0$). But there are also some others. Here is one:\footnote{For proofs
of the statements made (implicitly and explicitly) in the following (or at
least for the less obvious among these proofs), we refer to Section
\ref{sec.details.gf.laure}.}

\begin{definition}
\label{def.fps.laure.laupol}Let $K\left[  x^{\pm}\right]  $ be the
$K$-submodule of $K\left[  \left[  x^{\pm}\right]  \right]  $ consisting of
all \textbf{essentially finite} families $\left(  a_{n}\right)  _{n\in
\mathbb{Z}}$. This is indeed a $K$-submodule (check it!). It should be thought
of as an analogue of the ring of polynomials $K\left[  x\right]  $, but now
allowing for negative powers of $x$.

The elements of $K\left[  x^{\pm}\right]  $ are called \emph{Laurent
polynomials} in the indeterminate $x$ over $K$.

We define a multiplication on $K\left[  x^{\pm}\right]  $ by setting%
\[
\left(  a_{n}\right)  _{n\in\mathbb{Z}}\cdot\left(  b_{n}\right)
_{n\in\mathbb{Z}}=\left(  c_{n}\right)  _{n\in\mathbb{Z}}%
,\ \ \ \ \ \ \ \ \ \ \text{where}\ \ \ \ \ \ \ \ \ \ c_{n}=\sum_{i\in
\mathbb{Z}}a_{i}b_{n-i}.
\]
Note that the sum $\sum_{i\in\mathbb{Z}}a_{i}b_{n-i}$ is now well-defined,
because it is essentially finite.

We define an element $x\in K\left[  x^{\pm}\right]  $ by%
\[
x=\left(  \delta_{i,1}\right)  _{i\in\mathbb{Z}}%
\]
(where we are using the Kronecker delta notation again, as in Definition
\ref{def.kron-delta}).
\end{definition}

\begin{theorem}
\label{thm.fps.laure.laupol-ring}The $K$-module $K\left[  x^{\pm}\right]  $,
equipped with the multiplication we just defined, is a commutative
$K$-algebra. Its unity is $\left(  \delta_{i,0}\right)  _{i\in\mathbb{Z}}$.
The element $x$ is invertible in this $K$-algebra.
\end{theorem}

This $K$-algebra $K\left[  x^{\pm}\right]  $ is called the \emph{Laurent
polynomial ring} in one indeterminate $x$ over $K$. It is often denoted by
$K\left[  x^{\pm1}\right]  $ or $K\left[  x,x^{-1}\right]  $ as well.

\begin{proposition}
\label{prop.fps.laure.a=sumaixi}Any doubly infinite power series $a=\left(
a_{i}\right)  _{i\in\mathbb{Z}}\in K\left[  \left[  x^{\pm}\right]  \right]  $
satisfies%
\[
a=\sum_{i\in\mathbb{Z}}a_{i}x^{i}.
\]
Here, the powers $x^{i}$ are taken in the Laurent polynomial ring $K\left[
x^{\pm}\right]  $, but the infinite sum $\sum_{i\in\mathbb{Z}}a_{i}x^{i}$ is
taken in the $K$-module $K\left[  \left[  x^{\pm}\right]  \right]  $. (The
notions of summable families and infinite sums are defined in $K\left[
\left[  x^{\pm}\right]  \right]  $ in the same way as they are defined in
$K\left[  \left[  x\right]  \right]  $.)
\end{proposition}

Examples of Laurent polynomials are

\begin{itemize}
\item any polynomial in $K\left[  x\right]  $;

\item $x^{-15}$;

\item $x^{2}+3+7x^{-3}$.
\end{itemize}

There are other, equivalent ways to define the Laurent polynomial ring
$K\left[  x^{\pm}\right]  $:

\begin{itemize}
\item as the group algebra of the cyclic group $\mathbb{Z}$ over $K$;

\item as the localization of the polynomial ring $K\left[  x\right]  $ at the
powers of $x$.
\end{itemize}

(These are done in some textbooks on abstract algebra -- e.g., see
\cite[Exercise 3.6.32]{Ford21} for a quick overview.)

\subsubsection{Laurent polynomials are not enough}

Now let us see if we can make our above proof of Theorem
\ref{thm.fps.laure.balanced-tern-rep-uniq} rigorous using Laurent polynomials.
Unfortunately, $K\left[  x^{\pm}\right]  $ is \textquotedblleft too
small\textquotedblright\ to contain the infinite product $\prod_{i\geq
0}\left(  1+x^{3^{i}}+x^{-3^{i}}\right)  $. However, we can try using its
partial products, which are finite. For each $k\in\mathbb{N}$, we have%
\begin{align*}
&  \prod_{i=0}^{k}\left(  1+x^{3^{i}}+x^{-3^{i}}\right) \\
&  =\left(  1+x+x^{-1}\right)  \left(  1+x^{3}+x^{-3}\right)  \cdots\left(
1+x^{3^{k}}+x^{-3^{k}}\right) \\
&  =\dfrac{1-x^{3}}{x\left(  1-x\right)  }\cdot\dfrac{1-x^{9}}{x^{3}\left(
1-x^{3}\right)  }\cdot\cdots\cdot\dfrac{1-x^{3^{k+1}}}{x^{3^{k}}\left(
1-x^{3^{k}}\right)  }\\
&  \ \ \ \ \ \ \ \ \ \ \ \ \ \ \ \ \ \ \ \ \left(
\begin{array}
[c]{c}%
\text{this is somewhat unrigorous, since }1-x^{3^{i}}\text{ are not
invertible}\\
\text{in }K\left[  x^{\pm}\right]  \text{, but this will soon be made
rigorous}%
\end{array}
\right) \\
&  =\underbrace{\dfrac{1}{xx^{3}x^{9}\cdots x^{3^{k}}}}_{=x^{-\left(
3^{0}+3^{1}+\cdots+3^{k}\right)  }}\cdot\underbrace{\dfrac{1-x^{3^{k+1}}}%
{1-x}}_{=1+x+x^{2}+\cdots+x^{3^{k+1}-1}}\ \ \ \ \ \ \ \ \ \ \left(  \text{by
cancelling factors}\right) \\
&  =x^{-\left(  3^{0}+3^{1}+\cdots+3^{k}\right)  }\cdot\left(  1+x+x^{2}%
+\cdots+x^{3^{k+1}-1}\right) \\
&  =x^{-\left(  3^{0}+3^{1}+\cdots+3^{k}\right)  }\cdot\left(  1+x+x^{2}%
+\cdots+x^{2\left(  3^{0}+3^{1}+\cdots+3^{k}\right)  }\right) \\
&  \ \ \ \ \ \ \ \ \ \ \ \ \ \ \ \ \ \ \ \ \left(  \text{since }%
3^{k+1}-1=2\cdot\left(  3^{0}+3^{1}+\cdots+3^{k}\right)  \text{ (check
this!)}\right) \\
&  =x^{-\left(  3^{0}+3^{1}+\cdots+3^{k}\right)  }+x^{-\left(  3^{0}%
+3^{1}+\cdots+3^{k}\right)  +1}+x^{-\left(  3^{0}+3^{1}+\cdots+3^{k}\right)
+2}+\cdots+x^{3^{0}+3^{1}+\cdots+3^{k}}\\
&  =\sum_{\substack{n\in\mathbb{Z};\\\left\vert n\right\vert \leq3^{0}%
+3^{1}+\cdots+3^{k}}}x^{n}.
\end{align*}
On the other hand,%
\begin{align*}
&  \prod_{i=0}^{k}\underbrace{\left(  1+x^{3^{i}}+x^{-3^{i}}\right)  }%
_{=\sum_{b\in\left\{  0,1,-1\right\}  }x^{b\cdot3^{i}}}\\
&  =\prod_{i=0}^{k}\ \ \sum_{b\in\left\{  0,1,-1\right\}  }x^{b\cdot3^{i}}\\
&  =\sum_{\left(  b_{0},b_{1},\ldots,b_{k}\right)  \in\left\{  0,1,-1\right\}
^{k+1}}x^{b_{0}3^{0}}x^{b_{1}3^{1}}\cdots x^{b_{k}3^{k}}%
\ \ \ \ \ \ \ \ \ \ \left(
\begin{array}
[c]{c}%
\text{by expanding the product}\\
\text{using Proposition \ref{prop.fps.prodrule-fin-fin}}%
\end{array}
\right) \\
&  =\sum_{\left(  b_{0},b_{1},\ldots,b_{k}\right)  \in\left\{  0,1,-1\right\}
^{k+1}}x^{b_{0}3^{0}+b_{1}3^{1}+\cdots+b_{k}3^{k}}\\
&  =\sum_{n\in\mathbb{Z}}\ \ \sum_{\substack{\left(  b_{0},b_{1},\ldots
,b_{k}\right)  \in\left\{  0,1,-1\right\}  ^{k+1};\\b_{0}3^{0}+b_{1}%
3^{1}+\cdots+b_{k}3^{k}=n}}x^{n}\\
&  =\sum_{n\in\mathbb{Z}}\left(  \text{\# of }k\text{-bounded balanced ternary
representations of }n\right)  \cdot x^{n},
\end{align*}
where a balanced ternary representation $\left(  b_{0},b_{1},b_{2}%
,\ldots\right)  $ is said to be $k$\emph{-bounded} if $b_{k+1}=b_{k+2}%
=b_{k+3}=\cdots=0$. Comparing the two results, we find%
\begin{align*}
&  \sum_{n\in\mathbb{Z}}\left(  \text{\# of }k\text{-bounded balanced ternary
representations of }n\right)  \cdot x^{n}\\
&  =\sum_{\substack{n\in\mathbb{Z};\\\left\vert n\right\vert \leq3^{0}%
+3^{1}+\cdots+3^{k}}}x^{n}.
\end{align*}
Comparing coefficients, we thus see that each $n\in\mathbb{Z}$ satisfying
$\left\vert n\right\vert \leq3^{0}+3^{1}+\cdots+3^{k}$ has a unique
$k$-bounded balanced ternary representation. Letting $k\rightarrow\infty$ now
quickly yields Theorem \ref{thm.fps.laure.balanced-tern-rep-uniq} (because any
balanced ternary representation is $k$-bounded for any sufficiently large $k$).

Thus we have proved Theorem \ref{thm.fps.laure.balanced-tern-rep-uniq} up to
the fact that we have divided by the polynomials $1-x$, $1-x^{3}$, $1-x^{9}$,
$\ldots$, which are not invertible in the Laurent polynomial ring $K\left[
x^{\pm}\right]  $. In order to fill this gap, we need a new $K$-algebra: The
Laurent polynomial ring $K\left[  x^{\pm}\right]  $ is too small, whereas the
original $K$-module $K\left[  \left[  x^{\pm}\right]  \right]  $ is not a
ring. We need some kind of middle ground: some $K$-module lying between
$K\left[  x^{\pm}\right]  $ and $K\left[  \left[  x^{\pm}\right]  \right]  $
that is a ring but allows division by $1-x$, $1-x^{3}$, $1-x^{9}$, $\ldots$
(and, more generally, by $1-x^{i}$ for each positive integer $i$).

\subsubsection{Laurent series}

This middle ground is called $K\left(  \left(  x\right)  \right)  $ and is
defined as follows:

\begin{definition}
\label{def.fps.laure.lauser}We let $K\left(  \left(  x\right)  \right)  $ be
the subset of $K\left[  \left[  x^{\pm}\right]  \right]  $ consisting of all
families $\left(  a_{i}\right)  _{i\in\mathbb{Z}}\in K\left[  \left[  x^{\pm
}\right]  \right]  $ such that the sequence $\left(  a_{-1},a_{-2}%
,a_{-3},\ldots\right)  $ is essentially finite -- i.e., such that all
sufficiently low $i\in\mathbb{Z}$ satisfy $a_{i}=0$.

The elements of $K\left(  \left(  x\right)  \right)  $ are called
\emph{Laurent series} in one indeterminate $x$ over $K$.
\end{definition}

For example (for $K=\mathbb{Z}$):

\begin{itemize}
\item the \textquotedblleft series\textquotedblright\ $x^{-3}+x^{-2}%
+x^{-1}+x^{0}+x^{1}+\cdots$ belongs to $K\left(  \left(  x\right)  \right)  $;

\item the \textquotedblleft series\textquotedblright\ $1+x^{-1}+x^{-2}%
+x^{-3}+\cdots$ does not belong to $K\left(  \left(  x\right)  \right)  $;

\item the \textquotedblleft series\textquotedblright\ $\sum_{n\in\mathbb{Z}%
}x^{n}=\cdots+x^{-2}+x^{-1}+x^{0}+x^{1}+x^{2}+\cdots$ does not belong to
$K\left(  \left(  x\right)  \right)  $.
\end{itemize}

The subset $K\left(  \left(  x\right)  \right)  $ turns out to behave as
nicely as we might hope for:\footnote{See Section \ref{sec.details.gf.laure}
for a proof of Theorem \ref{thm.fps.laure.lauser-ring}.}

\begin{theorem}
\label{thm.fps.laure.lauser-ring}The subset $K\left(  \left(  x\right)
\right)  $ is a $K$-submodule of $K\left[  \left[  x^{\pm}\right]  \right]  $.
But it has a multiplication (unlike $K\left[  \left[  x^{\pm}\right]  \right]
$). This multiplication is given by the same rule as the multiplication of
$K\left[  x^{\pm}\right]  $: namely,%
\[
\left(  a_{n}\right)  _{n\in\mathbb{Z}}\cdot\left(  b_{n}\right)
_{n\in\mathbb{Z}}=\left(  c_{n}\right)  _{n\in\mathbb{Z}}%
,\ \ \ \ \ \ \ \ \ \ \text{where}\ \ \ \ \ \ \ \ \ \ c_{n}=\sum_{i\in
\mathbb{Z}}a_{i}b_{n-i}.
\]
The sum $\sum_{i\in\mathbb{Z}}a_{i}b_{n-i}$ here is well-defined, because it
is essentially finite (indeed, for all sufficiently low $i\in\mathbb{Z}$, we
have $a_{i}=0$ and thus $a_{i}b_{n-i}=0$; on the other hand, for all
sufficiently high $i\in\mathbb{Z}$, we have $b_{n-i}=0$ and thus $a_{i}%
b_{n-i}=0$).

Equipped with the multiplication we just defined, the $K$-module $K\left(
\left(  x\right)  \right)  $ becomes a commutative $K$-algebra with unity
$\left(  \delta_{i,0}\right)  _{i\in\mathbb{Z}}$.
\end{theorem}

Now, the ring $K\left(  \left(  x\right)  \right)  $ contains both the FPS
ring $K\left[  \left[  x\right]  \right]  $ and the Laurent polynomial ring
$K\left[  x^{\pm}\right]  $ as subrings (and actually as $K$-subalgebras).
This makes it one of the most convenient places for formal manipulation of
FPSs. (However, it has a disadvantage compared to the FPS ring $K\left[
\left[  x\right]  \right]  $: Namely, you cannot easily substitute something
for $x$ in a Laurent series $f\in K\left(  \left(  x\right)  \right)  $.)

Now, our above computation of $\prod_{i=0}^{k}\left(  1+x^{3^{i}}+x^{-3^{i}%
}\right)  $ makes perfect sense in the Laurent series ring $K\left(  \left(
x\right)  \right)  $: Indeed, for each positive integer $i$, the power series
$1-x^{i}$ is invertible in $K\left[  \left[  x\right]  \right]  $ and thus
also invertible in $K\left(  \left(  x\right)  \right)  $. Hence, at last, we
have a rigorous proof of Theorem \ref{thm.fps.laure.balanced-tern-rep-uniq}.

Actually, you \textbf{could} also make sense of our original argument for
proving Theorem \ref{thm.fps.laure.balanced-tern-rep-uniq}, with the infinite
product $\prod_{i\geq0}\left(  1+x^{3^{i}}+x^{-3^{i}}\right)  $, as long as
you made sure to interpret it correctly: First, compute the finite products
$\prod_{i=0}^{k}\left(  1+x^{3^{i}}+x^{-3^{i}}\right)  $ in the ring $K\left(
\left(  x\right)  \right)  $. Then, take their limit $\lim
\limits_{k\rightarrow\infty}\prod_{i=0}^{k}\left(  1+x^{3^{i}}+x^{-3^{i}%
}\right)  $ in $K\left[  \left[  x^{\pm}\right]  \right]  $ (this is not a
ring, but the notion of a limit in $K\left[  \left[  x^{\pm}\right]  \right]
$ is defined just as it was in $K\left[  \left[  x\right]  \right]  $).

More can be said about the $K$-algebra $K\left(  \left(  x\right)  \right)  $
when $K$ is a field: Indeed, in this case, it is itself a field! This fact
(whose proof is Exercise \ref{exe.fps.laure.field}) is not very useful in
combinatorics, but quite so in abstract algebra.

\subsubsection{A $K\left[  x^{\pm}\right]  $-module structure on $K\left[
\left[  x^{\pm}\right]  \right]  $}

One more remark is in order. As we have explained, $K\left[  \left[  x^{\pm
}\right]  \right]  $ is not a ring. However, some semblance of multiplication
can be defined in $K\left[  \left[  x^{\pm}\right]  \right]  $. Namely, a
\textquotedblleft doubly infinite power series\textquotedblright\ $\sum
_{n\in\mathbb{Z}}a_{n}x^{n}=\left(  a_{n}\right)  _{n\in\mathbb{Z}}\in
K\left[  \left[  x^{\pm}\right]  \right]  $ can be multiplied by a Laurent
polynomial $\sum_{n\in\mathbb{Z}}b_{n}x^{n}=\left(  b_{n}\right)
_{n\in\mathbb{Z}}\in K\left[  x^{\pm}\right]  $, because in this case the sums
$\sum_{i\in\mathbb{Z}}a_{i}b_{n-i}$ will be essentially finite. Thus, while we
cannot always multiply two elements of $K\left[  \left[  x^{\pm}\right]
\right]  $, we can always multiply an element of $K\left[  \left[  x^{\pm
}\right]  \right]  $ with an element of $K\left[  x^{\pm}\right]  $. This
makes $K\left[  \left[  x^{\pm}\right]  \right]  $ into a $K\left[  x^{\pm
}\right]  $-module. The equality%
\[
\left(  1-x\right)  \left(  \cdots+x^{-2}+x^{-1}+x^{0}+x^{1}+x^{2}%
+\cdots\right)  =0
\]
reveals that this module has torsion (i.e., a product of two nonzero factors
can be zero). Thus, while we can multiply a doubly infinite power series by
$1-x$, we cannot (in general) divide it by $1-x$.

\subsection{\label{sec.gf.multivar}Multivariate FPSs}%

\[%
\begin{tabular}
[c]{||l||}\hline\hline
\textbf{TODO: This section will get more details in a future version.}%
\\\hline\hline
\end{tabular}
\]

Multivariate FPSs are just FPSs in several variables. Their theory is mostly
analogous to the theory of univariate FPSs (which is what we have been
studying above), but requires more subscripts. I will just discuss the main differences.

For instance, FPSs in two variables $x$ and $y$ have the form $\sum
_{i,j\in\mathbb{N}}a_{i,j}x^{i}y^{j}$, where the summation sign $\sum
_{i,j\in\mathbb{N}}$ means $\sum_{\left(  i,j\right)  \in\mathbb{N}^{2}}$ of
course. Formally, such an FPS is a family $\left(  a_{i,j}\right)  _{\left(
i,j\right)  \in\mathbb{N}^{2}}$ of elements of $K$. The indeterminates $x$ and
$y$ are defined by%
\[
x=\left(  \delta_{\left(  i,j\right)  ,\left(  1,0\right)  }\right)  _{\left(
i,j\right)  \in\mathbb{N}^{2}}\ \ \ \ \ \ \ \ \ \ \text{and}%
\ \ \ \ \ \ \ \ \ \ y=\left(  \delta_{\left(  i,j\right)  ,\left(  0,1\right)
}\right)  _{\left(  i,j\right)  \in\mathbb{N}^{2}}%
\]
(so that each indeterminate has exactly one coefficient equal to $1$, while
all other coefficients are $0$). The rules for addition, subtraction, scaling
and multiplication are essentially as they are for univariate FPSs, except
that now the indexing set is $\mathbb{N}^{2}$ instead of $\mathbb{N}$. For
example, multiplication of FPSs in $x$ and $y$ is defined by the formula%
\[
\left[  x^{n}y^{m}\right]  \left(  ab\right)  =\sum_{\substack{\left(
i,j\right)  ,\ \left(  k,\ell\right)  \in\mathbb{N}^{2};\\i+k=n;\\j+\ell
=m}}\left[  x^{i}y^{j}\right]  a\cdot\left[  x^{k}y^{\ell}\right]  b
\]
(for any two FPSs $a$ and $b$ and any $n,m\in\mathbb{N}$).

More generally, for any $k\in\mathbb{N}$, the FPSs in $k$ variables
$x_{1},x_{2},\ldots,x_{k}$ are defined to be the families $\left(
a_{\mathbf{i}}\right)  _{\mathbf{i}\in\mathbb{N}^{k}}$ of elements of $K$
indexed by $k$-tuples $\mathbf{i}=\left(  i_{1},i_{2},\ldots,i_{k}\right)
\in\mathbb{N}^{k}$. Addition, subtraction and scaling of such families is
defined entrywise. Multiplication is defined by the formula%
\[
\left[  \mathbf{x}^{\mathbf{n}}\right]  \left(  ab\right)  =\sum
_{\substack{\mathbf{i},\mathbf{j}\in\mathbb{N}^{k};\\\mathbf{i}+\mathbf{j}%
=\mathbf{n}}}\left[  \mathbf{x}^{\mathbf{i}}\right]  a\cdot\left[
\mathbf{x}^{\mathbf{j}}\right]  b
\]
(for any two FPSs $a$ and $b$ and any $\mathbf{n}\in\mathbb{N}^{k}$), where

\begin{itemize}
\item the sum $\mathbf{i}+\mathbf{j}$ means the entrywise sum of the
$k$-tuples $\mathbf{i}$ and $\mathbf{j}$ (that is, $\mathbf{i}+\mathbf{j}%
=\left(  i_{1}+j_{1},i_{2}+j_{2},\ldots,i_{k}+j_{k}\right)  $, where
$\mathbf{i}=\left(  i_{1},i_{2},\ldots,i_{k}\right)  $ and $\mathbf{j}=\left(
j_{1},j_{2},\ldots,j_{k}\right)  $);

\item if $\mathbf{m}\in\mathbb{N}^{k}$ and $h$ is an FPS in $x_{1}%
,x_{2},\ldots,x_{k}$, then $\left[  \mathbf{x}^{\mathbf{m}}\right]  h$ is the
$\mathbf{m}$-th entry of the family $h$. Just as in the univariate case, this
entry $\left[  \mathbf{x}^{\mathbf{m}}\right]  h$ is called the
\emph{coefficient of the monomial} $\mathbf{x}^{\mathbf{m}}:=x_{1}^{m_{1}%
}x_{2}^{m_{2}}\cdots x_{k}^{m_{k}}$ in $h$ (where $\mathbf{m}=\left(
m_{1},m_{2},\ldots,m_{k}\right)  $).
\end{itemize}

With this notation, the formula for $ab$ doesn't look any more complicated
than the analogous formula in the univariate case. The only difference is that
the monomials and the coefficients are now indexed not by the nonnegative
integers, but by the $k$-tuples of nonnegative integers instead.

The indeterminates $x_{1},x_{2},\ldots,x_{k}$ are defined by%
\[
x_{i}=\left(  \delta_{\mathbf{n},\left(  0,0,\ldots,0,1,0,0,\ldots,0\right)
}\right)  _{\mathbf{n}\in\mathbb{N}^{k}},
\]
where the tuple $\left(  0,0,\ldots,0,1,0,0,\ldots,0\right)  $ is a $k$-tuple
with a lone $1$ in its $i$-th position. Thus, if $\mathbf{m}=\left(
m_{1},m_{2},\ldots,m_{k}\right)  \in\mathbb{N}^{k}$, then%
\[
x_{1}^{m_{1}}x_{2}^{m_{2}}\cdots x_{k}^{m_{k}}=\left(  \delta_{\mathbf{n}%
,\mathbf{m}}\right)  _{\mathbf{n}\in\mathbb{N}^{k}},
\]
so that each FPS $f=\left(  f_{\mathbf{m}}\right)  _{\mathbf{m}\in
\mathbb{N}^{k}}$ satisfies%
\[
f=\sum_{\mathbf{m}=\left(  m_{1},m_{2},\ldots,m_{k}\right)  \in\mathbb{N}^{k}%
}f_{\mathbf{m}}x_{1}^{m_{1}}x_{2}^{m_{2}}\cdots x_{k}^{m_{k}}.
\]

Most of what we have said about FPSs in one variable applies similarly to FPSs
in multiple variables. The proofs are similar but more laborious due to the
need for subscripts. Instead of the derivative of an FPS, there are now $k$
derivatives (one for each variable); they are called \emph{partial
derivatives}. One needs to be somewhat careful with substitution -- e.g., one
cannot substitute non-commuting elements into a multivariate polynomial. For
example, you cannot substitute two non-commuting matrices $A$ and $B$ for $x$
and $y$ into the polynomial $xy$, at least not without sacrificing the rule
that a value of the product of two polynomials should be the product of their
values (since $\underbrace{x\left[  A,B\right]  }_{=A}\cdot
\underbrace{y\left[  A,B\right]  }_{=B}=AB$ would have to equal
$\underbrace{y\left[  A,B\right]  }_{=B}\cdot\underbrace{x\left[  A,B\right]
}_{=A}=BA$). But you can still substitute $k$ commuting elements for the $k$
indeterminates in a $k$-variable polynomial. You can also compose multivariate
FPSs as long as appropriate summability conditions are satisfied.

\begin{definition}
Let $k\in\mathbb{N}$. The $K$-algebra of all FPSs in $k$ variables
$x_{1},x_{2},\ldots,x_{k}$ over $K$ will be denoted by $K\left[  \left[
x_{1},x_{2},\ldots,x_{k}\right]  \right]  $.
\end{definition}

Sometimes we will use different names for our variables. For example, if we
work with $2$ variables, we will commonly call them $x$ and $y$ instead of
$x_{1}$ and $x_{2}$. Correspondingly, we will use the notation $K\left[
\left[  x,y\right]  \right]  $ (instead of $K\left[  \left[  x_{1}%
,x_{2}\right]  \right]  $) for the $K$-algebra of FPSs in these two variables.

Let me give an example of working with multivariate FPSs.

Let us work in $K\left[  \left[  x,y\right]  \right]  $. On the one hand, we
have%
\begin{align*}
\sum_{n,k\in\mathbb{N}}\dbinom{n}{k}x^{n}y^{k}  &  =\sum_{n\in\mathbb{N}%
}\ \ \sum_{k\in\mathbb{N}}\dbinom{n}{k}x^{n}y^{k}=\sum_{n\in\mathbb{N}}%
x^{n}\underbrace{\sum_{k\in\mathbb{N}}\dbinom{n}{k}y^{k}}_{\substack{=\left(
1+y\right)  ^{n}\\\text{(by the binomial formula)}}}\\
&  =\sum_{n\in\mathbb{N}}x^{n}\left(  1+y\right)  ^{n}=\sum_{n\in\mathbb{N}%
}\left(  x\left(  1+y\right)  \right)  ^{n}=\dfrac{1}{1-x\left(  1+y\right)
}\\
&  \ \ \ \ \ \ \ \ \ \ \left(
\begin{array}
[c]{c}%
\text{here, we have substituted }x\left(  1+y\right)  \text{ for }x\text{
in}\\
\text{the formula }\sum_{n\in\mathbb{N}}x^{n}=\dfrac{1}{1-x}\text{;}\\
\text{this is allowed since }x\left(  1+y\right)  \text{ has constant term }0
\end{array}
\right) \\
&  =\dfrac{1}{1-x}\cdot\dfrac{1}{1-\dfrac{x}{1-x}y}\ \ \ \ \ \ \ \ \ \ \left(
\text{easy to check by computation}\right) \\
&  =\dfrac{1}{1-x}\cdot\sum_{k\in\mathbb{N}}\left(  \dfrac{x}{1-x}y\right)
^{k}\\
&  \ \ \ \ \ \ \ \ \ \ \left(
\begin{array}
[c]{c}%
\text{here, we have substituted }\dfrac{x}{1-x}y\text{ for }x\\
\text{in the formula }\dfrac{1}{1-x}=\sum_{k\in\mathbb{N}}x^{k}%
\end{array}
\right) \\
&  =\dfrac{1}{1-x}\cdot\sum_{k\in\mathbb{N}}\dfrac{x^{k}}{\left(  1-x\right)
^{k}}y^{k}=\sum_{k\in\mathbb{N}}\dfrac{x^{k}}{\left(  1-x\right)  ^{k+1}}%
y^{k}.
\end{align*}
On the other hand, we have%
\[
\sum_{n,k\in\mathbb{N}}\dbinom{n}{k}x^{n}y^{k}=\sum_{k\in\mathbb{N}}%
\ \ \sum_{n\in\mathbb{N}}\dbinom{n}{k}x^{n}y^{k}=\sum_{k\in\mathbb{N}}\left(
\sum_{n\in\mathbb{N}}\dbinom{n}{k}x^{n}\right)  y^{k}.
\]
Comparing these two equalities, we find%
\begin{equation}
\sum_{k\in\mathbb{N}}\dfrac{x^{k}}{\left(  1-x\right)  ^{k+1}}y^{k}=\sum
_{k\in\mathbb{N}}\left(  \sum_{n\in\mathbb{N}}\dbinom{n}{k}x^{n}\right)
y^{k}. \label{eq.fps.mulvar.exa1.xyk}%
\end{equation}

Now, comparing coefficients in front of $x^{i}y^{k}$ in
(\ref{eq.fps.mulvar.exa1.xyk}), we conclude that%
\begin{equation}
\dfrac{x^{k}}{\left(  1-x\right)  ^{k+1}}=\sum_{n\in\mathbb{N}}\dbinom{n}%
{k}x^{n}\ \ \ \ \ \ \ \ \ \ \text{for each }k\in\mathbb{N}.
\label{eq.fps.mulvar.exa1.res1}%
\end{equation}
Let me explain in a bit more detail what I mean by \textquotedblleft comparing
coefficients in front of $x^{i}y^{k}$\textquotedblright. What we have used is
the following simple fact:

\begin{proposition}
\label{prop.fps.mulvar.comp-y-coeff}Let $f_{0},f_{1},f_{2},\ldots$ and
$g_{0},g_{1},g_{2},\ldots$ be FPSs in a single variable $x$ such that%
\begin{equation}
\sum_{k\in\mathbb{N}}f_{k}y^{k}=\sum_{k\in\mathbb{N}}g_{k}y^{k}%
\ \ \ \ \ \ \ \ \ \ \text{in }K\left[  \left[  x,y\right]  \right]  .
\label{eq.prop.fps.mulvar.comp-y-coeff.ass}%
\end{equation}
Then, $f_{k}=g_{k}$ for each $k\in\mathbb{N}$.
\end{proposition}

\begin{proof}
For each $k\in\mathbb{N}$, let us write the two FPSs $f_{k}$ and $g_{k}$ as
$f_{k}=\sum_{n\in\mathbb{N}}f_{k,n}x^{n}$ and $g_{k}=\sum_{n\in\mathbb{N}%
}g_{k,n}x^{n}$ with $f_{k,n},g_{k,n}\in K$. Then, the equality
(\ref{eq.prop.fps.mulvar.comp-y-coeff.ass}) can be rewritten as
\[
\sum_{k\in\mathbb{N}}\ \ \sum_{n\in\mathbb{N}}f_{k,n}x^{n}y^{k}=\sum
_{k\in\mathbb{N}}\ \ \sum_{n\in\mathbb{N}}g_{k,n}x^{n}y^{k}.
\]
Now, comparing coefficients in front of $x^{n}y^{k}$ in this equality, we
obtain $f_{k,n}=g_{k,n}$ for each $k,n\in\mathbb{N}$. Therefore, $f_{k}=g_{k}$
for each $k\in\mathbb{N}$. This proves Proposition
\ref{prop.fps.mulvar.comp-y-coeff}.
\end{proof}

Now, because of (\ref{eq.fps.mulvar.exa1.xyk}), we can apply Proposition
\ref{prop.fps.mulvar.comp-y-coeff} to $f_{k}=\dfrac{x^{k}}{\left(  1-x\right)
^{k+1}}$ and $g_{k}=\sum_{n\in\mathbb{N}}\dbinom{n}{k}x^{n}$. Thus, we obtain
the equality%
\[
\dfrac{x^{k}}{\left(  1-x\right)  ^{k+1}}=\sum_{n\in\mathbb{N}}\dbinom{n}%
{k}x^{n}\ \ \ \ \ \ \ \ \ \ \text{for each }k\in\mathbb{N}.
\]
This is an equality between \textbf{univariate} FPSs, even though we have
obtained it by manipulating \textbf{bivariate} FPSs (i.e., FPSs in two
variables $x$ and $y$). As a homework exercise (Exercise
\ref{exe.fps.sum-n-choose-k-xn-elementary}), you can prove this equality in a
more elementary, univariate way. But the idea to introduce extra variables (in
our case, in order to discover an equality) is highly useful in many
situations (some of which we will see later in this course).

\newpage

\section{\label{chap.pars}Integer partitions and $q$-binomial coefficients}

We have previously counted compositions of an $n\in\mathbb{N}$. These are
(roughly speaking) ways to write $n$ as a sum of finitely many positive
integers, where the order matters. For example, $3=3=2+1=1+2=1+1+1$, so $3$
has $4$ compositions. Formally, compositions are tuples of positive integers.

Now, let us disregard the order. There are two ways to make this rigorous:
either we replace tuples by multisets, or we require the tuples to be weakly
decreasing. These result in the same count, but we will use the 2nd way, just
because tuples are easier to work with than multisets. This will lead to the
notion of \emph{integer partitions}.

\subsection{Partition basics}

\subsubsection{Definitions}

The following definition is built in analogy to Definition \ref{def.fps.comps}:

\begin{definition}
\label{def.pars.parts}\textbf{(a)} An \emph{(integer) partition} means a
(finite) weakly decreasing tuple of positive integers -- i.e., a finite tuple
$\left(  \lambda_{1},\lambda_{2},\ldots,\lambda_{m}\right)  $ of positive
integers such that $\lambda_{1}\geq\lambda_{2}\geq\cdots\geq\lambda_{m}$.

Thus, partitions are the same as weakly decreasing compositions. Hence, the
notions of \emph{size} and \emph{length} of a partition are automatically
defined, since we have defined them for compositions (in Definition
\ref{def.fps.comps}). \medskip

\textbf{(b)} The \emph{parts} of a partition $\left(  \lambda_{1},\lambda
_{2},\ldots,\lambda_{m}\right)  $ are simply its entries $\lambda_{1}%
,\lambda_{2},\ldots,\lambda_{m}$. \medskip

\textbf{(c)} Let $n\in\mathbb{Z}$. A \emph{partition of }$n$ means a partition
whose size is $n$. \medskip

\textbf{(d)} Let $n\in\mathbb{Z}$ and $k\in\mathbb{N}$. A \emph{partition of
}$n$\emph{ into }$k$\emph{ parts} is a partition whose size is $n$ and whose
length is $k$.
\end{definition}

\begin{example}
\label{exa.pars.pars5}The partitions of $5$ are%
\[
\left(  5\right)  ,\ \ \left(  4,1\right)  ,\ \ \left(  3,2\right)
,\ \ \left(  3,1,1\right)  ,\ \ \left(  2,2,1\right)  ,\ \ \left(
2,1,1,1\right)  ,\ \ \left(  1,1,1,1,1\right)  .
\]

\end{example}

\begin{definition}
\label{def.pars.pn-pkn}\textbf{(a)} Let $n\in\mathbb{Z}$ and $k\in\mathbb{N}$.
Then, we set%
\[
p_{k}\left(  n\right)  :=\left(  \text{\# of partitions of }n\text{ into
}k\text{ parts}\right)  .
\]

\textbf{(b)} Let $n\in\mathbb{Z}$. Then, we set%
\[
p\left(  n\right)  :=\left(  \text{\# of partitions of }n\right)  .
\]
This is called the $n$\emph{-th partition number}.
\end{definition}

\begin{example}
Our above list of partitions of $5$ reveals that%
\begin{align*}
p_{0}\left(  5\right)   &  =0;\\
p_{1}\left(  5\right)   &  =1;\\
p_{2}\left(  5\right)   &  =2;\\
p_{3}\left(  5\right)   &  =2;\\
p_{4}\left(  5\right)   &  =1;\\
p_{5}\left(  5\right)   &  =1;\\
p_{k}\left(  5\right)   &  =0\ \ \ \ \ \ \ \ \ \ \text{for any }k>5;
\end{align*}
and finally $p\left(  5\right)  =7$.
\end{example}

Here are the values of $p\left(  n\right)  $ for the first $15$ nonnegative
integers $n$:%
\[%
\begin{tabular}
[c]{|c||c|c|c|c|c|c|c|c|c|c|c|c|c|c|c|}\hline
$n$ & $0$ & $1$ & $2$ & $3$ & $4$ & $5$ & $6$ & $7$ & $8$ & $9$ & $10$ & $11$
& $12$ & $13$ & $14$\\\hline
$p\left(  n\right)  $ & $1$ & $1$ & $2$ & $3$ & $5$ & $7$ & $11$ & $15$ & $22$
& $30$ & $42$ & $56$ & $77$ & $101$ & $135$\\\hline
\end{tabular}
\ \ \ \ .
\]
The sequence $\left(  p\left(  0\right)  ,p\left(  1\right)  ,p\left(
2\right)  ,\ldots\right)  $ is remarkable for being an integer sequence that
grows faster than polynomially, but still considerably slower than
exponentially. (See (\ref{eq.pars.asymptotic-pn}) for an asymptotic
expansion.) This not-too-fast growth (for instance, $p\left(  100\right)
=190\ 569\ 292$ is far smaller than $2^{100}$) makes integer partitions rather
convenient for computer experiments.

\subsubsection{Simple properties of partition numbers}

We will next state some elementary properties of $p_{k}\left(  n\right)  $ and
$p\left(  n\right)  $, but first we introduce a few very basic notations:

\begin{definition}
\label{def.pars.iverson}We will use the \emph{Iverson bracket notation}: If
$\mathcal{A}$ is a logical statement, then $\left[  \mathcal{A}\right]  $
means the \emph{truth value} of $\mathcal{A}$; this is the integer $%
\begin{cases}
1, & \text{if }\mathcal{A}\text{ is true};\\
0, & \text{if }\mathcal{A}\text{ is false}.
\end{cases}
$
\end{definition}

For example, $\left[  2+2=4\right]  =1$ and $\left[  2+2=5\right]  =0$.

Note that the Kronecker delta notation is a particular case of the Iverson
bracket: We have
\[
\delta_{i,j}=\left[  i=j\right]  \ \ \ \ \ \ \ \ \ \ \text{for any objects
}i\text{ and }j.
\]

\begin{definition}
\label{def.pars.floor-ceil}Let $a$ be a real number.

Then, $\left\lfloor a\right\rfloor $ (called the \emph{floor} of $a$) means
the largest integer that is $\leq a$.

Likewise, $\left\lceil a\right\rceil $ (called the \emph{ceiling} of $a$)
means the smallest integer that is $\geq a$.
\end{definition}

For example, the number $\pi\approx3.14$ satisfies $\left\lfloor
\pi\right\rfloor =3$ and $\left\lceil \pi\right\rceil =4$ and $\left\lfloor
-\pi\right\rfloor =-4$ and $\left\lceil -\pi\right\rceil =-3$. For another
example, $\left\lfloor n\right\rfloor =\left\lceil n\right\rceil =n$ for each
$n\in\mathbb{Z}$.

The following proposition collects various basic properties of the numbers
introduced in Definition \ref{def.pars.pn-pkn}:

\begin{proposition}
\label{prop.pars.basics}Let $n\in\mathbb{Z}$ and $k\in\mathbb{N}$. \medskip

\textbf{(a)} We have $p_{k}\left(  n\right)  =0$ whenever $n<0$ and
$k\in\mathbb{N}$. \medskip

\textbf{(b)} We have $p_{k}\left(  n\right)  =0$ whenever $k>n$. \medskip

\textbf{(c)} We have $p_{0}\left(  n\right)  =\left[  n=0\right]  $. \medskip

\textbf{(d)} We have $p_{1}\left(  n\right)  =\left[  n>0\right]  $. \medskip

\textbf{(e)} We have $p_{k}\left(  n\right)  =p_{k}\left(  n-k\right)
+p_{k-1}\left(  n-1\right)  $ whenever $k>0$. \medskip

\textbf{(f)} We have $p_{2}\left(  n\right)  =\left\lfloor n/2\right\rfloor $
whenever $n\in\mathbb{N}$. \medskip

\textbf{(g)} We have $p\left(  n\right)  =p_{0}\left(  n\right)  +p_{1}\left(
n\right)  +\cdots+p_{n}\left(  n\right)  $ whenever $n\in\mathbb{N}$. \medskip

\textbf{(h)} We have $p\left(  n\right)  =0$ whenever $n<0$.
\end{proposition}

\begin{proof}
[Proof of Proposition \ref{prop.pars.basics} (sketched).]\textbf{(a)} The size
of a partition is always nonnegative (being a sum of positive integers). Thus,
a negative number $n$ has no partitions whatsoever. Thus, $p_{k}\left(
n\right)  =0$ whenever $n<0$ and $k\in\mathbb{N}$. \medskip

\textbf{(b)} If $\left(  \lambda_{1},\lambda_{2},\ldots,\lambda_{k}\right)  $
is a partition, then $\lambda_{i}\geq1$ for each $i\in\left\{  1,2,\ldots
,k\right\}  $ (because a partition is a tuple of positive integers, i.e., of
integers $\geq1$). Hence, if $\left(  \lambda_{1},\lambda_{2},\ldots
,\lambda_{k}\right)  $ is a partition of $n$ into $k$ parts, then%
\begin{align*}
n  &  =\lambda_{1}+\lambda_{2}+\cdots+\lambda_{k}\ \ \ \ \ \ \ \ \ \ \left(
\text{since }\left(  \lambda_{1},\lambda_{2},\ldots,\lambda_{k}\right)  \text{
is a partition of }n\right) \\
&  \geq\underbrace{1+1+\cdots+1}_{k\text{ times}}\ \ \ \ \ \ \ \ \ \ \left(
\text{since }\lambda_{i}\geq1\text{ for each }i\in\left\{  1,2,\ldots
,k\right\}  \right) \\
&  =k.
\end{align*}
Thus, a partition of $n$ into $k$ parts cannot satisfy $k>n$. Thus, no such
partitions exist if $k>n$. In other words, $p_{k}\left(  n\right)  =0$ if
$k>n$. \medskip

\textbf{(c)} The integer $0$ has a unique partition into $0$ parts, namely the
empty tuple $\left(  {}\right)  $. A nonzero integer $n$ cannot have any
partitions into $0$ parts, since the empty tuple has size $0\neq n$. Thus,
$p_{0}\left(  n\right)  $ equals $1$ for $n=0$ and equals $0$ for $n\neq0$. In
other words, $p_{0}\left(  n\right)  =\left[  n=0\right]  $. \medskip

\textbf{(d)} Any positive integer $n$ has a unique partition into $1$ part --
namely, the $1$-tuple $\left(  n\right)  $. On the other hand, if $n$ is not
positive, then this $1$-tuple is not a partition, so in this case $n$ has no
partition into $1$ part. Thus, $p_{1}\left(  n\right)  $ equals $1$ if $n$ is
positive and equals $0$ otherwise. In other words, $p_{1}\left(  n\right)
=\left[  n>0\right]  $. \medskip

\textbf{(e)} Assume that $k>0$. We must prove that $p_{k}\left(  n\right)
=p_{k}\left(  n-k\right)  +p_{k-1}\left(  n-1\right)  $.

We consider all partitions of $n$ into $k$ parts. We classify these partitions
into two types:

\begin{itemize}
\item \emph{Type 1} consists of all partitions that have $1$ as a part.

\item \emph{Type 2} consists of all partitions that don't.
\end{itemize}

For example, here are the partitions of $5$ along with their types:%
\[
\underbrace{\left(  4,1\right)  ,\ \ \left(  3,1,1\right)  ,\ \ \left(
2,2,1\right)  ,\ \ \left(  2,1,1,1\right)  ,\ \ \left(  1,1,1,1,1\right)
}_{\text{Type 1}},\ \ \underbrace{\left(  5\right)  ,\ \ \left(  3,2\right)
}_{\text{Type 2}}.
\]

Let us count the type-1 partitions and the type-2 partitions separately.

Any type-1 partition has $1$ as a part, and therefore has $1$ as its last part
(because it is weakly decreasing). Hence, any type-1 partition has the form
$\left(  \lambda_{1},\lambda_{2},\ldots,\lambda_{k-1},1\right)  $. If $\left(
\lambda_{1},\lambda_{2},\ldots,\lambda_{k-1},1\right)  $ is a type-1 partition
(of $n$ into $k$ parts), then $\left(  \lambda_{1},\lambda_{2},\ldots
,\lambda_{k-1}\right)  $ is a partition of $n-1$ into $k-1$ parts. Thus, we
have a map%
\begin{align*}
\left\{  \text{type-1 partitions of }n\text{ into }k\text{ parts}\right\}   &
\rightarrow\left\{  \text{partitions of }n-1\text{ into }k-1\text{
parts}\right\}  ,\\
\left(  \lambda_{1},\lambda_{2},\ldots,\lambda_{k-1},1\right)   &
\mapsto\left(  \lambda_{1},\lambda_{2},\ldots,\lambda_{k-1}\right)  .
\end{align*}
This map is a bijection (since it has an inverse map, which simply inserts a
$1$ at the end of a partition). Thus, the bijection principle shows that%
\begin{align*}
&  \left(  \text{\# of type-1 partitions of }n\text{ into }k\text{
parts}\right) \\
&  =\left(  \text{\# of partitions of }n-1\text{ into }k-1\text{
parts}\right)  =p_{k-1}\left(  n-1\right)
\end{align*}
(by the definition of $p_{k-1}\left(  n-1\right)  $).

Now let us count type-2 partitions. A type-2 partition does not have $1$ as a
part; hence, all its parts are larger than $1$ (because all its parts are
positive integers), and therefore we can subtract $1$ from each part and still
have a partition in front of us. To be more specific: If $\left(  \lambda
_{1},\lambda_{2},\ldots,\lambda_{k}\right)  $ is a type-2 partition of $n$
into $k$ parts, then subtracting $1$ from each of its parts produces the
$k$-tuple $\left(  \lambda_{1}-1,\lambda_{2}-1,\ldots,\lambda_{k}-1\right)  $,
which is a partition of $n-k$ into $k$ parts. Hence, we have a map%
\begin{align*}
\left\{  \text{type-2 partitions of }n\text{ into }k\text{ parts}\right\}   &
\rightarrow\left\{  \text{partitions of }n-k\text{ into }k\text{
parts}\right\}  ,\\
\left(  \lambda_{1},\lambda_{2},\ldots,\lambda_{k}\right)   &  \mapsto\left(
\lambda_{1}-1,\lambda_{2}-1,\ldots,\lambda_{k}-1\right)  .
\end{align*}
This map is a bijection (since it has an inverse map, which simply adds $1$ to
each entry of a partition). Thus, the bijection principle shows that%
\begin{align*}
&  \left(  \text{\# of type-2 partitions of }n\text{ into }k\text{
parts}\right) \\
&  =\left(  \text{\# of partitions of }n-k\text{ into }k\text{ parts}\right)
=p_{k}\left(  n-k\right)
\end{align*}
(by the definition of $p_{k}\left(  n-k\right)  $).

Since any partition of $n$ into $k$ parts is either type-1 or type-2 (but not
both at the same time), we now have%
\begin{align*}
&  \left(  \text{\# of partitions of }n\text{ into }k\text{ parts}\right) \\
&  =\underbrace{\left(  \text{\# of type-1 partitions of }n\text{ into
}k\text{ parts}\right)  }_{=p_{k-1}\left(  n-1\right)  }+\underbrace{\left(
\text{\# of type-2 partitions of }n\text{ into }k\text{ parts}\right)
}_{=p_{k}\left(  n-k\right)  }\\
&  =p_{k-1}\left(  n-1\right)  +p_{k}\left(  n-k\right)  =p_{k}\left(
n-k\right)  +p_{k-1}\left(  n-1\right)  .
\end{align*}
Since the left hand side of this equality is $p_{k}\left(  n\right)  $, we
thus have proved that $p_{k}\left(  n\right)  =p_{k}\left(  n-k\right)
+p_{k-1}\left(  n-1\right)  $. \medskip

\textbf{(f)} Let $n\in\mathbb{N}$. The partitions of $n$ into $2$ parts are%
\[
\left(  n-1,1\right)  ,\ \ \left(  n-2,2\right)  ,\ \ \left(  n-3,3\right)
,\ \ \ldots,\ \ \left(  \underbrace{n-\left\lfloor n/2\right\rfloor
}_{=\left\lceil n/2\right\rceil },\left\lfloor n/2\right\rfloor \right)  .
\]
Thus there are $\left\lfloor n/2\right\rfloor $ of them. In other words,
$p_{2}\left(  n\right)  =\left\lfloor n/2\right\rfloor $. \medskip

\textbf{(g)} Let $n\in\mathbb{N}$. Any partition of $n$ must have $k$ parts
for some $k\in\mathbb{N}$. Thus,%
\begin{align*}
p\left(  n\right)   &  =\sum_{k\in\mathbb{N}}p_{k}\left(  n\right)
=\sum_{k=0}^{n}p_{k}\left(  n\right)  +\sum_{k=n+1}^{\infty}\underbrace{p_{k}%
\left(  n\right)  }_{\substack{=0\\\text{(by Proposition
\ref{prop.pars.basics} \textbf{(b)})}}}\\
&  =\sum_{k=0}^{n}p_{k}\left(  n\right)  =p_{0}\left(  n\right)  +p_{1}\left(
n\right)  +\cdots+p_{n}\left(  n\right)  .
\end{align*}

\textbf{(h)} Same argument as for \textbf{(a)}.
\end{proof}

\subsubsection{The generating function}

Proposition \ref{prop.pars.basics} \textbf{(e)} is a recursive formula that
makes it not too hard to compute $p_{k}\left(  n\right)  $ for reasonably
small values of $n$ and $k$. Then, using Proposition \ref{prop.pars.basics}
\textbf{(g)}, we can compute $p\left(  n\right)  $ from these $p_{k}\left(
n\right)  $'s. However, one might want a better, faster method.

To get there, let me first express the generating function of the numbers
$p\left(  n\right)  $:

\begin{theorem}
\label{thm.pars.main-gf}In the FPS ring $\mathbb{Z}\left[  \left[  x\right]
\right]  $, we have%
\[
\sum_{n\in\mathbb{N}}p\left(  n\right)  x^{n}=\prod_{k=1}^{\infty}\dfrac
{1}{1-x^{k}}.
\]

(The product on the right hand side is well-defined, since multiplying a FPS
by $\dfrac{1}{1-x^{k}}$ does not affect its first $k$ coefficients.)
\end{theorem}

\begin{example}
Let us check the above equality \textquotedblleft up to $x^{5}$%
\textquotedblright, i.e., let us compare the coefficients of $x^{i}$ for
$i<5$. (In doing so, we can ignore all powers of $x$ higher than $x^{4}$.) We
have%
\begin{align*}
\prod_{k=1}^{\infty}\dfrac{1}{1-x^{k}}  &  =\dfrac{1}{1-x^{1}}\cdot\dfrac
{1}{1-x^{2}}\cdot\dfrac{1}{1-x^{3}}\cdot\dfrac{1}{1-x^{4}}\cdot\cdots\\
&  =\ \ \ \left(  1+x+x^{2}+x^{3}+x^{4}+\cdots\right) \\
&  \ \ \ \ \ \cdot\left(  1+x^{2}+x^{4}+\cdots\right) \\
&  \ \ \ \ \ \cdot\left(  1+x^{3}+\cdots\right) \\
&  \ \ \ \ \ \cdot\left(  1+x^{4}+\cdots\right) \\
&  \ \ \ \ \ \cdot\left(  1+\cdots\right) \\
&  \ \ \ \ \ \cdot\left(  1+\cdots\right) \\
&  \ \ \ \ \ \cdot\cdots\\
&  =1+x+2x^{2}+3x^{3}+5x^{4}+\cdots\\
&  =p\left(  0\right)  +p\left(  1\right)  x+p\left(  2\right)  x^{2}+p\left(
3\right)  x^{3}+p\left(  4\right)  x^{4}+\cdots.
\end{align*}

\end{example}

\begin{proof}
[Proof of Theorem \ref{thm.pars.main-gf}.]We have%
\begin{align*}
&  \prod_{k=1}^{\infty}\underbrace{\dfrac{1}{1-x^{k}}}_{=1+x^{k}+x^{2k}%
+x^{3k}+\cdots}\\
&  =\prod_{k=1}^{\infty}\underbrace{\left(  1+x^{k}+x^{2k}+x^{3k}%
+\cdots\right)  }_{=\sum\limits_{u\in\mathbb{N}}x^{ku}}=\prod_{k=1}^{\infty
}\ \ \sum_{u\in\mathbb{N}}x^{ku}\\
&  =\sum_{\substack{\left(  u_{1},u_{2},u_{3},\ldots\right)  \in
\mathbb{N}^{\infty}\text{ is}\\\text{essentially finite}}}x^{1u_{1}}x^{2u_{2}%
}x^{3u_{3}}\cdots\ \ \ \ \ \ \ \ \ \ \left(
\begin{array}
[c]{c}%
\text{here, we expanded the product}\\
\text{using Proposition \ref{prop.fps.prodrule-inf-infN}}%
\end{array}
\right) \\
&  =\sum_{\substack{\left(  u_{1},u_{2},u_{3},\ldots\right)  \in
\mathbb{N}^{\infty}\text{ is}\\\text{essentially finite}}}x^{1u_{1}%
+2u_{2}+3u_{3}+\cdots}\\
&  =\sum_{n\in\mathbb{N}}\left\vert Q_{n}\right\vert x^{n},
\end{align*}
where
\[
Q_{n}=\left\{  \left(  u_{1},u_{2},u_{3},\ldots\right)  \in\mathbb{N}^{\infty
}\text{ essentially finite }\mid\ 1u_{1}+2u_{2}+3u_{3}+\cdots=n\right\}  .
\]

Thus, it will suffice to show that%
\[
\left\vert Q_{n}\right\vert =p\left(  n\right)  \ \ \ \ \ \ \ \ \ \ \text{for
each }n\in\mathbb{N}.
\]

Let us fix $n\in\mathbb{N}$. We want to construct a bijection from $Q_{n}$ to
$\left\{  \text{partitions of }n\right\}  $.

Here is how to do this: For any $\left(  u_{1},u_{2},u_{3},\ldots\right)  \in
Q_{n}$, define a partition%
\begin{align*}
\pi\left(  u_{1},u_{2},u_{3},\ldots\right)  :=  &  \left(  \text{the partition
that contains each }i\text{ exactly }u_{i}\text{ times}\right) \\
=  &  \left(  \ldots,\underbrace{3,3,\ldots,3}_{u_{3}\text{ times}%
},\underbrace{2,2,\ldots,2}_{u_{2}\text{ times}},\underbrace{1,1,\ldots
,1}_{u_{1}\text{ times}}\right)  .
\end{align*}
This is a partition of $n$, since its size is $1u_{1}+2u_{2}+3u_{3}+\cdots=n$
(because $\left(  u_{1},u_{2},u_{3},\ldots\right)  \in Q_{n}$). Thus, we have
defined a partition $\pi\left(  u_{1},u_{2},u_{3},\ldots\right)  $ of $n$ for
each $\left(  u_{1},u_{2},u_{3},\ldots\right)  \in Q_{n}$. In other words, we
have defined a map
\[
\pi:Q_{n}\rightarrow\left\{  \text{partitions of }n\right\}  .
\]
It remains to show that this map $\pi$ is a bijection. We define a map%
\[
\rho:\left\{  \text{partitions of }n\right\}  \rightarrow Q_{n}%
\]
that sends each partition $\lambda$ of $n$ to the sequence%
\[
\left(  \text{\# of }1\text{'s in }\lambda\text{,\ \ \# of }2\text{'s in
}\lambda\text{,\ \ \# of }3\text{'s in }\lambda\text{, \ }\ldots\right)  \in
Q_{n}%
\]
(this is indeed a sequence in $Q_{n}$, since
\[
\sum_{i=1}^{\infty}i\left(  \text{\# of }i\text{'s in }\lambda\right)
=\left(  \text{the sum of all entries of }\lambda\right)  =\left\vert
\lambda\right\vert =n
\]
because $\lambda$ is a partition of $n$). It is now easy to check that the
maps $\pi$ and $\rho$ are mutually inverse, so that $\pi$ is a bijection. The
bijection principle therefore yields $\left\vert Q_{n}\right\vert =\left(
\text{\# of partitions of }n\right)  =p\left(  n\right)  $; but this is
precisely what we wanted to show. The proof of Theorem \ref{thm.pars.main-gf}
is thus complete.
\end{proof}

Theorem \ref{thm.pars.main-gf} has a \textquotedblleft
finite\textquotedblright\ analogue (finite in the sense that the product
$\prod_{k=1}^{\infty}\dfrac{1}{1-x^{k}}$ is replaced by a finite product; the
FPSs are still infinite):

\begin{theorem}
\label{thm.pars.main-gf-parts-n}Let $m\in\mathbb{N}$. For each $n\in
\mathbb{N}$, let $p_{\operatorname*{parts}\leq m}\left(  n\right)  $ be the \#
of partitions $\lambda$ of $n$ such that all parts of $\lambda$ are $\leq m$.
Then,
\[
\sum_{n\in\mathbb{N}}p_{\operatorname*{parts}\leq m}\left(  n\right)
x^{n}=\prod_{k=1}^{m}\dfrac{1}{1-x^{k}}.
\]

\end{theorem}

\begin{proof}
[Proof of Theorem \ref{thm.pars.main-gf-parts-n}.]This proof is mostly
analogous to the above proof of Theorem \ref{thm.pars.main-gf}, and to some
extent even simpler because it uses $m$-tuples instead of infinite sequences.

We have%
\begin{align*}
&  \prod_{k=1}^{m}\underbrace{\dfrac{1}{1-x^{k}}}_{=1+x^{k}+x^{2k}%
+x^{3k}+\cdots}\\
&  =\prod_{k=1}^{m}\underbrace{\left(  1+x^{k}+x^{2k}+x^{3k}+\cdots\right)
}_{=\sum\limits_{u\in\mathbb{N}}x^{ku}}=\prod_{k=1}^{m}\ \ \sum_{u\in
\mathbb{N}}x^{ku}\\
&  =\sum_{\left(  u_{1},u_{2},\ldots,u_{m}\right)  \in\mathbb{N}^{m}}%
x^{1u_{1}}x^{2u_{2}}\cdots x^{mu_{m}}\ \ \ \ \ \ \ \ \ \ \left(
\begin{array}
[c]{c}%
\text{here, we expanded the product}\\
\text{using Proposition \ref{prop.fps.prodrule-fin-inf}}%
\end{array}
\right) \\
&  =\sum_{\left(  u_{1},u_{2},\ldots,u_{m}\right)  \in\mathbb{N}^{m}}%
x^{1u_{1}+2u_{2}+\cdots+mu_{m}}=\sum_{n\in\mathbb{N}}\left\vert Q_{n}%
\right\vert x^{n},
\end{align*}
where
\[
Q_{n}=\left\{  \left(  u_{1},u_{2},\ldots,u_{m}\right)  \in\mathbb{N}%
^{m}\text{ }\mid\ 1u_{1}+2u_{2}+\cdots+mu_{m}=n\right\}  .
\]

Thus, it will suffice to show that%
\[
\left\vert Q_{n}\right\vert =p_{\operatorname*{parts}\leq m}\left(  n\right)
\ \ \ \ \ \ \ \ \ \ \text{for each }n\in\mathbb{N}.
\]

Let us fix $n\in\mathbb{N}$. We want to construct a bijection from $Q_{n}$ to
the set $\left\{  \text{partitions }\lambda\text{ of }n\text{ such that all
parts of }\lambda\text{ are }\leq m\right\}  $.

Here is how to do this: For any $\left(  u_{1},u_{2},\ldots,u_{m}\right)  \in
Q_{n}$, define a partition%
\begin{align*}
\pi\left(  u_{1},u_{2},\ldots,u_{m}\right)  :=  &  \left(  \text{the partition
that contains each }i\text{ exactly }u_{i}\text{ times}\right) \\
=  &  \left(  \underbrace{m,m,\ldots,m}_{u_{m}\text{ times}},\ldots
,\underbrace{2,2,\ldots,2}_{u_{2}\text{ times}},\underbrace{1,1,\ldots
,1}_{u_{1}\text{ times}}\right)  .
\end{align*}
Thus, we have defined a map
\[
\pi:Q_{n}\rightarrow\left\{  \text{partitions }\lambda\text{ of }n\text{ such
that all parts of }\lambda\text{ are }\leq m\right\}  .
\]
It is easy to see that this map $\pi$ is a bijection\footnote{The argument is
analogous to the one used in the proof of Theorem \ref{thm.pars.main-gf}.}.
The bijection principle therefore yields
\[
\left\vert Q_{n}\right\vert =\left(  \text{\# of partitions }\lambda\text{ of
}n\text{ such that all parts of }\lambda\text{ are }\leq m\right)
=p_{\operatorname*{parts}\leq m}\left(  n\right)  ;
\]
but this is precisely what we wanted to show. The proof of Theorem
\ref{thm.pars.main-gf-parts-n} is thus complete.
\end{proof}

Theorem \ref{thm.pars.main-gf-parts-n} can be generalized further: Instead of
counting the partitions of $n$ whose all parts are $\leq m$, we can count the
partitions of $n$ whose parts all belong to a given set $I$ of positive
integers. This leads to the following formula (which specializes back to
Theorem \ref{thm.pars.main-gf-parts-n} when we take $I=\left\{  1,2,\ldots
,m\right\}  $):

\begin{theorem}
\label{thm.pars.main-gf-parts-I}Let $I$ be a subset of $\left\{
1,2,3,\ldots\right\}  $. For each $n\in\mathbb{N}$, let $p_{I}\left(
n\right)  $ be the \# of partitions $\lambda$ of $n$ such that all parts of
$\lambda$ belong to $I$. Then,
\[
\sum_{n\in\mathbb{N}}p_{I}\left(  n\right)  x^{n}=\prod_{k\in I}\dfrac
{1}{1-x^{k}}.
\]

\end{theorem}

\begin{proof}
[Proof of Theorem \ref{thm.pars.main-gf-parts-I} (sketched).]This is analogous
to the proof of Theorem \ref{thm.pars.main-gf-parts-n}, with some minor
changes: The $\prod_{k=1}^{m}$ sign has to be replaced by $\prod_{k\in I}$;
the $m$-tuples $\left(  u_{1},u_{2},\ldots,u_{m}\right)  \in\mathbb{N}^{m}$
must be replaced by the essentially finite families $\left(  u_{i}\right)
_{i\in I}\in\mathbb{N}^{I}$; the bijection $\pi$ has to be replaced by the new
bijection
\[
\pi:Q_{n}\rightarrow\left\{  \text{partitions }\lambda\text{ of }n\text{ such
that all parts of }\lambda\text{ are }\in I\right\}
\]
(where $Q_{n}$ is now the set of all essentially finite families $\left(
u_{i}\right)  _{i\in I}\in\mathbb{N}^{I}$ satisfying $\sum_{i\in I} iu_{i} =
n$) defined by%
\begin{align*}
\pi\left(  \left(  u_{i}\right)  _{i\in I}\right)  :=  &  \left(  \text{the
partition that contains each }i\in I\text{ exactly }u_{i}\text{ times}\right)
\\
=  &  \left(  \ldots,\underbrace{i_{3},i_{3},\ldots,i_{3}}_{u_{i_{3}}\text{
times}},\underbrace{i_{2},i_{2},\ldots,i_{2}}_{u_{i_{2}}\text{ times}%
},\underbrace{i_{1},i_{1},\ldots,i_{1}}_{u_{i_{1}}\text{ times}}\right)  ,
\end{align*}
where $i_{1},i_{2},i_{3},\ldots$ are the elements of $I$ listed in increasing
order. The details are left to the reader.
\end{proof}

\subsubsection{Odd parts and distinct parts}

Next, we shall state a result of Euler that we have already discovered in a
different language.

\begin{definition}
\label{def.pars.odd-dist-parts}Let $n\in\mathbb{Z}$. \medskip

\textbf{(a)} A \emph{partition of }$n$\emph{ into odd parts} means a partition
of $n$ whose all parts are odd. \medskip

\textbf{(b)} A \emph{partition of }$n$\emph{ into distinct parts} means a
partition of $n$ whose parts are distinct. \medskip

\textbf{(c)} Let%
\begin{align*}
p_{\operatorname*{odd}}\left(  n\right)   &  :=\left(  \text{\# of partitions
of }n\text{ into odd parts}\right)  \ \ \ \ \ \ \ \ \ \ \text{and}\\
p_{\operatorname*{dist}}\left(  n\right)   &  :=\left(  \text{\# of partitions
of }n\text{ into distinct parts}\right)  .
\end{align*}

\end{definition}

\begin{example}
We have%
\begin{align*}
p_{\operatorname*{odd}}\left(  7\right)   &  =\left\vert \left\{  \left(
7\right)  ,\ \ \left(  5,1,1\right)  ,\ \ \left(  3,3,1\right)  ,\ \ \left(
3,1,1,1,1\right)  ,\ \ \left(  1,1,1,1,1,1,1\right)  \right\}  \right\vert
=5;\\
p_{\operatorname*{dist}}\left(  7\right)   &  =\left\vert \left\{  \left(
7\right)  ,\ \ \left(  6,1\right)  ,\ \ \left(  5,2\right)  ,\ \ \left(
4,3\right)  ,\ \ \left(  4,2,1\right)  \right\}  \right\vert =5.
\end{align*}

\end{example}

\begin{theorem}
[Euler's odd-distinct identity]\label{thm.pars.odd-dist-equal}We have
$p_{\operatorname*{odd}}\left(  n\right)  =p_{\operatorname*{dist}}\left(
n\right)  $ for each $n\in\mathbb{N}$.
\end{theorem}

We have already encountered this theorem before (as Theorem
\ref{thm.gf.prod.euler-comb}, albeit in less precise language), and we have
proved it using the generating function identity%
\[
\prod_{i>0}\left(  1-x^{2i-1}\right)  ^{-1}=\prod_{k>0}\left(  1+x^{k}\right)
.
\]

Let me outline a different, bijective proof.

\begin{proof}
[Second proof of Theorem \ref{thm.pars.odd-dist-equal} (sketched).]Let
$n\in\mathbb{N}$. We want to construct a bijection%
\[
A:\left\{  \text{partitions of }n\text{ into odd parts}\right\}
\rightarrow\left\{  \text{partitions of }n\text{ into distinct parts}\right\}
.
\]
We shall do this as follows: Given a partition $\lambda$ of $n$ into odd
parts, we repeatedly merge pairs of equal parts in $\lambda$ until no more
equal parts appear (and then we sort the entries of the resulting tuple into
weakly decreasing order). The final result will be $A\left(  \lambda\right)
$. Here are two examples:

\begin{itemize}
\item To compute $A\left(  5,5,3,1,1,1\right)  $, we compute\footnote{The two
entries underlined are the two equal entries that are going to get merged in
the next step. Note that there are usually several candidates, and we just
pick one pair at will.}%
\[
\left(  \underline{5,5},3,1,1,1\right)  \rightarrow\left(
10,3,1,\underline{1,1}\right)  \rightarrow\left(  10,3,2,1\right)  .
\]
Thus, $A\left(  5,5,3,1,1,1\right)  =\left(  10,3,2,1\right)  $.

\item To compute $A\left(  5,3,1,1,1,1\right)  $, we compute%
\[
\left(  5,3,1,1,\underline{1,1}\right)  \rightarrow\left(
5,3,2,\underline{1,1}\right)  \rightarrow\left(  5,3,\underline{2,2}\right)
\rightarrow\left(  5,4,3\right)  .
\]
Thus, $A\left(  5,3,1,1,1,1\right)  =\left(  5,4,3\right)  $.
\end{itemize}

Why is this map $A$ well-defined? We only specified the sort of steps we are
allowed to take when computing $A\left(  \lambda\right)  $; however, there is
often a choice involved in taking these steps (since there are often several
pairs of equal parts).\footnote{For example, we could have also computed
$A\left(  5,5,3,1,1,1\right)  $ as follows:%
\[
\left(  5,5,3,\underline{1},1,\underline{1}\right)  \rightarrow\left(
\underline{5,5},3,2,1\right)  \rightarrow\left(  10,3,2,1\right)  .
\]
} So we have specified a non-deterministic algorithm. Why is the resulting
partition independent of the choices we make?

One way to prove this is using the \emph{diamond lemma}, which is a general
tool for proving that certain non-deterministic algorithms have unique final
outcomes (independent of the choices taken). See, e.g.,
\url{https://mathoverflow.net/questions/289300/} for a list of references on
this lemma.

For the map $A$, we can also proceed differently, by analyzing the algorithm
that we used to define $A$. Namely, we observe what is really going on when we
are merging equal parts. Let us say our original partition $\lambda$ has $p$
many $1$s. Let us first merge them in pairs, so that we get $\left\lfloor
p/2\right\rfloor $ many $2$s and maybe one single $1$. Then, let us merge the
$2$s in pairs, so that we get $\left\lfloor \left\lfloor p/2\right\rfloor
/2\right\rfloor $ many $4$s, maybe a single $2$, and maybe a single $1$.
Proceed until no more than one $1$, no more than one $2$, no more than one
$4$, no more than one $8$, and so on remain. This clears out any duplicate
parts of the form $2^{k}$. Next do the same with parts of the form
$3\cdot2^{k}$ (that is, with parts equal to $3$, $6$, $12$, $24$, and so on),
then with parts of the form $5\cdot2^{k}$, and so on.

The nice thing about this way of proceeding is that we can explicitly describe
the final outcome. Indeed, if the original partition $\lambda$ (a partition of
$n$ into odd parts) contains an odd part $k$ precisely $m$ many times, and if
the binary representation of $m$ is $m=\left(  m_{i}m_{i-1}\cdots m_{1}%
m_{0}\right)  _{2}$ (that is, if $m_{0},m_{1},\ldots,m_{i}\in\left\{
0,1\right\}  $ satisfy $m=\sum_{j=0}^{i}m_{j}2^{j}$), then the partition
$A\left(  \lambda\right)  $ will contain the number $2^{0}k$ exactly $m_{0}$
times, the number $2^{1}k$ exactly $m_{1}$ times, the number $2^{2}k$ exactly
$m_{2}$ times, and so on. Since the binary digits $m_{0},m_{1},\ldots,m_{i}$
are all $\leq1$, this partition $A\left(  \lambda\right)  $ will therefore not
contain any number more than once, i.e., it will be a partition into distinct parts.

It is not hard to check that this map $A$ is indeed a bijection. Indeed, in
order to see this, we construct a map $B$ that will turn out to be its
inverse. Here, we start with a partition $\lambda$ of $n$ into distinct parts.
Let us represent each part of this partition in the form $k\cdot2^{i}$ for
some odd $k\geq1$ and some integer $i\geq0$. (Recall that any positive integer
can be represented uniquely in this form.) Now, replace this part $k\cdot
2^{i}$ by $2^{i}$ many $k$'s. The resulting partition (once all parts have
been replaced) will usually have many equal parts, but all its parts are odd.
We define $B\left(  \lambda\right)  $ to be this resulting partition.
Alternatively, $B\left(  \lambda\right)  $ can also be constructed
step-by-step by a non-deterministic algorithm: Starting with $\lambda$, keep
\textquotedblleft breaking even parts into halves\textquotedblright\ (i.e.,
whenever you see an even part $m$, replace it by two parts $\dfrac{m}{2}$ and
$\dfrac{m}{2}$), until no even parts remain any more. The result is $B\left(
\lambda\right)  $. It is not hard to see that both descriptions of $B\left(
\lambda\right)  $ describe the same partition. It is furthermore easy to see
that this map $B$ is indeed an inverse of $A$, so that $A$ is indeed a
bijection. Thus, the bijection principle yields%
\[
\left\vert \left\{  \text{partitions of }n\text{ into odd parts}\right\}
\right\vert =\left\vert \left\{  \text{partitions of }n\text{ into distinct
parts}\right\}  \right\vert .
\]
In other words, $p_{\operatorname*{odd}}\left(  n\right)
=p_{\operatorname*{dist}}\left(  n\right)  $. This proves Theorem
\ref{thm.pars.odd-dist-equal}.
\end{proof}

Yet another proof of Theorem \ref{thm.pars.odd-dist-equal} will be given in
Subsection \ref{subsec.sign.pie.exas}.

\subsubsection{Partitions with a given largest part}

Here is another situation in which two kinds of partitions are equinumerous:

\begin{proposition}
\label{prop.pars.pkn=dual}Let $n\in\mathbb{N}$ and $k\in\mathbb{N}$. Then,%
\[
p_{k}\left(  n\right)  =\left(  \text{\# of partitions of }n\text{ whose
largest part is }k\right)  .
\]

\end{proposition}

Here and in the following, we use the following convention:

\begin{convention}
\label{conv.pars.largest-part-0}We agree to say that the largest part of the
empty partition $\left(  {}\right)  $ is $0$ (even though this partition has
no parts).
\end{convention}

\begin{example}
For $n=4$ and $k=3$, we have%
\[
p_{k}\left(  n\right)  =p_{3}\left(  4\right)  =1\ \ \ \ \ \ \ \ \ \ \left(
\text{due to the partition }\left(  2,1,1\right)  \right)
\]
and%
\begin{align*}
&  \left(  \text{\# of partitions of }n\text{ whose largest part is }k\right)
\\
&  =\left(  \text{\# of partitions of }4\text{ whose largest part is }3\right)
\\
&  =1\ \ \ \ \ \ \ \ \ \ \left(  \text{due to the partition }\left(
3,1\right)  \right)  .
\end{align*}
Thus, Proposition \ref{prop.pars.pkn=dual} holds for $n=4$ and $k=3$.
\end{example}

\begin{proof}
[Proof of Proposition \ref{prop.pars.pkn=dual} (sketched).]We do a
\textquotedblleft proof by picture\textquotedblright\ (it can be made rigorous
-- see Exercise \ref{exe.pars.transpose} for this). We pick $n=14$ and $k=4$
for example, and we start with the partition $\lambda=\left(  5,4,4,1\right)
$ of $n$ into $k$ parts.

We draw a table of $k$ left-aligned rows, where the length of each row equals
the corresponding part of $\lambda$ (that is, the $i$-th row from the top has
$\lambda_{i}$ boxes, where $\lambda=\left(  \lambda_{1},\lambda_{2}%
,\ldots,\lambda_{k}\right)  $):%
\[%
\begin{ytableau}
\none[5 \to] & \none& & & & & \\
\none[4 \to] & \none& & & & \\
\none[4 \to] & \none& & & & \\
\none[1 \to] & \none&
\end{ytableau}%
\]

Now, let us flip this table across the \textquotedblleft main
diagonal\textquotedblright\ (i.e., the diagonal that goes from the top-left
corner to the bottom-right corner)\footnote{This kind of flip is precisely how
you would transpose a matrix.}:%
\[%
\begin{ytableau}
\none[5] & \none[4] & \none[4] & \none[1] \\
\none[\downarrow] & \none[\downarrow] & \none[\downarrow] & \none
[\downarrow] \\
& & & \\
& & \\
& & \\
& & \\
\;
\end{ytableau}%
\]

The lengths of the rows of the resulting table again form a partition of $n$.
(In our case, this new partition is $\left(  4,3,3,3,1\right)  $.) Moreover,
the largest part of this new partition is $k$ (because the original table had
$k$ rows, so the flipped table has $k$ columns, and this means that its top
row has $k$ boxes). This procedure (i.e., turning a partition into a table,
then flipping the table across the \textquotedblleft main
diagonal\textquotedblright, and then reading the lengths of the rows of the
resulting table again as a partition) therefore gives a map from%
\[
\left\{  \text{partitions of }n\text{ into }k\text{ parts}\right\}
\]
to%
\[
\left\{  \text{partitions of }n\text{ whose largest part is }k\right\}  .
\]
Moreover, this map is a bijection (indeed, its inverse can be effected in the
exact same way, by flipping the table). This bijection is called
\emph{conjugation} of partitions, and will be studied in more detail later.

Here are some pointers to how this proof can be formalized (see Exercise
\ref{exe.pars.transpose} for much more): For any partition $\lambda=\left(
\lambda_{1},\lambda_{2},\ldots,\lambda_{k}\right)  $, we define the
\emph{Young diagram} of $\lambda$ to be the set%
\[
Y\left(  \lambda\right)  :=\left\{  \left(  i,j\right)  \in\mathbb{Z}%
^{2}\ \mid\ 1\leq i\leq k\text{ and }1\leq j\leq\lambda_{i}\right\}  .
\]
This Young diagram is precisely the table that we drew above, as long as we
agree to identify each pair $\left(  i,j\right)  \in Y\left(  \lambda\right)
$ with the box in row $i$ and column $j$. Now, the \emph{conjugate} of the
partition $\lambda$ is the partition $\lambda^{t}$ uniquely determined by%
\[
Y\left(  \lambda^{t}\right)  =\operatorname*{flip}\left(  Y\left(
\lambda\right)  \right)  =\left\{  \left(  j,i\right)  \ \mid\ \left(
i,j\right)  \in Y\left(  \lambda\right)  \right\}  .
\]
Explicitly, $\lambda^{t}$ can be defined by $\lambda^{t}=\left(  \mu_{1}%
,\mu_{2},\ldots,\mu_{p}\right)  $, where $p$ is the largest part of $\lambda$
and where%
\[
\mu_{i}=\left(  \text{\# of parts of }\lambda\text{ that are }\geq i\right)
\ \ \ \ \ \ \ \ \ \ \text{for each }i\in\left\{  1,2,\ldots,p\right\}  .
\]
(This conjugate $\lambda^{t}$ is also often called $\lambda^{\prime}$, and is
also known as the \emph{transpose} of $\lambda$.) Now, it is not hard to show
that $\left\vert \lambda^{t}\right\vert =\left\vert \lambda\right\vert $ and
$\left(  \lambda^{t}\right)  ^{t}=\lambda$ for each partition $\lambda$, and
that the largest part of $\lambda^{t}$ equals the length of $\lambda$. Using
these observations (which are proved in Exercise \ref{exe.pars.transpose}), we
see that the map%
\begin{align*}
\left\{  \text{partitions of }n\text{ into }k\text{ parts}\right\}   &
\rightarrow\left\{  \text{partitions of }n\text{ whose largest part is
}k\right\}  ,\\
\lambda &  \mapsto\lambda^{t}%
\end{align*}
is well-defined and is a bijection; thus, the above proof of Proposition
\ref{prop.pars.pkn=dual} becomes fully rigorous.
\end{proof}

The word \textquotedblleft Young\textquotedblright\ in \textquotedblleft Young
diagram\textquotedblright\ (and, later, \textquotedblleft Young
tableau\textquotedblright) does not imply any novelty (Young diagrams have
been around in some form or another since the 19th century -- if often in the
superficially different guise of \textquotedblleft Ferrers
diagrams\textquotedblright), but rather honors
\href{https://en.wikipedia.org/wiki/Alfred_Young}{Alfred Young}, who built up
the representation theory of symmetric groups (and significantly forwarded
invariant theory) using these objects.

\begin{corollary}
\label{cor.pars.p0kn=dual}Let $n\in\mathbb{N}$ and $k\in\mathbb{N}$. Then,%
\begin{align*}
&  p_{0}\left(  n\right)  +p_{1}\left(  n\right)  +\cdots+p_{k}\left(
n\right) \\
&  =\left(  \text{\# of partitions of }n\text{ whose largest part is }\leq
k\right)  .
\end{align*}

\end{corollary}

\begin{proof}
[Proof of Corollary \ref{cor.pars.p0kn=dual}.]We have%
\begin{align*}
&  p_{0}\left(  n\right)  +p_{1}\left(  n\right)  +\cdots+p_{k}\left(
n\right) \\
&  =\sum_{i=0}^{k}\underbrace{p_{i}\left(  n\right)  }_{\substack{=\left(
\text{\# of partitions of }n\text{ whose largest part is }i\right)
\\\text{(by Proposition \ref{prop.pars.pkn=dual},}\\\text{applied to }i\text{
instead of }k\text{)}}}\\
&  =\sum_{i=0}^{k}\left(  \text{\# of partitions of }n\text{ whose largest
part is }i\right) \\
&  =\left(  \text{\# of partitions of }n\text{ whose largest part is }\leq
k\right)  .
\end{align*}
This proves Corollary \ref{cor.pars.p0kn=dual}.
\end{proof}

Corollary \ref{cor.pars.p0kn=dual} leads to yet another FPS identity:

\begin{theorem}
\label{thm.pars.main-gf-0n}Let $m\in\mathbb{N}$. Then,%
\[
\sum_{n\in\mathbb{N}}\left(  p_{0}\left(  n\right)  +p_{1}\left(  n\right)
+\cdots+p_{m}\left(  n\right)  \right)  x^{n}=\prod_{k=1}^{m}\dfrac{1}%
{1-x^{k}}.
\]

\end{theorem}

\begin{proof}
[Proof of Theorem \ref{thm.pars.main-gf-0n}.]For each $n\in\mathbb{N}$, we
have%
\begin{align*}
&  p_{0}\left(  n\right)  +p_{1}\left(  n\right)  +\cdots+p_{m}\left(
n\right) \\
&  =\left(  \text{\# of partitions of }n\text{ whose largest part is }\leq
m\right)  \ \ \ \ \ \ \ \ \ \ \left(  \text{by Corollary
\ref{cor.pars.p0kn=dual}}\right) \\
&  =\left(  \text{\# of partitions of }n\text{ whose all parts are }\leq
m\right) \\
&  \ \ \ \ \ \ \ \ \ \ \ \ \ \ \ \ \ \ \ \ \left(
\begin{array}
[c]{c}%
\text{because the condition \textquotedblleft the largest part is }\leq
m\text{\textquotedblright\ for a}\\
\text{partition is clearly equivalent to \textquotedblleft all parts are }\leq
m\text{\textquotedblright}%
\end{array}
\right) \\
&  =p_{\operatorname*{parts}\leq m}\left(  n\right)  ,
\end{align*}
where $p_{\operatorname*{parts}\leq m}\left(  n\right)  $ is defined as in
Theorem \ref{thm.pars.main-gf-parts-n}. Hence,%
\[
\sum_{n\in\mathbb{N}}\underbrace{\left(  p_{0}\left(  n\right)  +p_{1}\left(
n\right)  +\cdots+p_{m}\left(  n\right)  \right)  }_{=p_{\operatorname*{parts}%
\leq m}\left(  n\right)  }x^{n}=\sum_{n\in\mathbb{N}}p_{\operatorname*{parts}%
\leq m}\left(  n\right)  x^{n}=\prod_{k=1}^{m}\dfrac{1}{1-x^{k}}%
\]
(by Theorem \ref{thm.pars.main-gf-parts-n}). This proves Theorem
\ref{thm.pars.main-gf-0n}.
\end{proof}

\subsubsection{Partition number vs. sums of divisors}

We shall now prove a curious combinatorial identity (due, I think, to Euler),
illustrating the usefulness of generating functions and logarithmic
derivatives. It connects the partition numbers $p\left(  n\right)  $ with a
further sequence, which sums the positive divisors of a given positive integer:

\begin{theorem}
\label{thm.pars.sigma1}For any positive integer $n$, let $\sigma\left(
n\right)  $ denote the sum of all positive divisors of $n$. (For example,
$\sigma\left(  6\right)  =1+2+3+6=12$ and $\sigma\left(  7\right)  =1+7=8$.)

For any $n\in\mathbb{N}$, we have%
\[
np\left(  n\right)  =\sum_{k=1}^{n}\sigma\left(  k\right)  p\left(
n-k\right)  .
\]

\end{theorem}

For example, for $n=3$, this says that%
\[
3\underbrace{p\left(  3\right)  }_{=3}=\underbrace{\sigma\left(  1\right)
}_{=1}\underbrace{p\left(  2\right)  }_{=2}+\underbrace{\sigma\left(
2\right)  }_{=3}\underbrace{p\left(  1\right)  }_{=1}+\underbrace{\sigma
\left(  3\right)  }_{=4}\underbrace{p\left(  0\right)  }_{=1}.
\]

For reference, let us give a table of the first $15$ sums $\sigma\left(
n\right)  $ defined in Theorem \ref{thm.pars.sigma1} (see the
\href{https://oeis.org/A000203}{Sequence A000203 in the OEIS} for more
values):%
\[%
\begin{tabular}
[c]{|c||c|c|c|c|c|c|c|c|c|c|c|c|c|c|c|}\hline
$n$ & $1$ & $2$ & $3$ & $4$ & $5$ & $6$ & $7$ & $8$ & $9$ & $10$ & $11$ & $12$
& $13$ & \multicolumn{1}{|c|}{$14$} & $15$\\\hline
$\sigma\left(  n\right)  $ & $1$ & $3$ & $4$ & $7$ & $6$ & $12$ & $8$ & $15$ &
$13$ & $18$ & $12$ & $28$ & $14$ & $24$ & $24$\\\hline
\end{tabular}
\ \ \ \ \ .
\]

\begin{proof}
[Proof of Theorem \ref{thm.pars.sigma1} (sketched).]Define the two FPSs%
\[
P:=\sum_{n\in\mathbb{N}}p\left(  n\right)  x^{n}\in\mathbb{Z}\left[  \left[
x\right]  \right]  \ \ \ \ \ \ \ \ \ \ \text{and}\ \ \ \ \ \ \ \ \ \ S:=\sum
_{k>0}\sigma\left(  k\right)  x^{k}\in\mathbb{Z}\left[  \left[  x\right]
\right]  .
\]
Hence,%
\begin{align}
SP  &  =\left(  \sum_{k>0}\sigma\left(  k\right)  x^{k}\right)  \left(
\sum_{n\in\mathbb{N}}p\left(  n\right)  x^{n}\right)  =\sum_{k>0}%
\ \ \sum_{n\in\mathbb{N}}\sigma\left(  k\right)  \underbrace{x^{k}p\left(
n\right)  x^{n}}_{=p\left(  n\right)  x^{k+n}}\nonumber\\
&  =\sum_{k>0}\ \ \sum_{n\in\mathbb{N}}\sigma\left(  k\right)  p\left(
n\right)  x^{k+n}=\underbrace{\sum_{k>0}\ \ \sum_{\substack{m\in
\mathbb{N};\\m\geq k}}}_{\substack{=\sum_{m\in\mathbb{N}}\ \ \sum
_{\substack{k>0;\\m\geq k}}\\=\sum_{m\in\mathbb{N}}\ \ \sum_{k=1}^{m}}%
}\sigma\left(  k\right)  p\left(  m-k\right)  \underbrace{x^{k+\left(
m-k\right)  }}_{=x^{m}}\nonumber\\
&  \ \ \ \ \ \ \ \ \ \ \ \ \ \ \ \ \ \ \ \ \left(  \text{here, we have
substituted }m-k\text{ for }n\text{ in the inner sum}\right) \nonumber\\
&  =\sum_{m\in\mathbb{N}}\ \ \sum_{k=1}^{m}\sigma\left(  k\right)  p\left(
m-k\right)  x^{m}. \label{pf.thm.pars.sigma1.SP=}%
\end{align}
From $P=\sum_{n\in\mathbb{N}}p\left(  n\right)  x^{n}$, we obtain%
\begin{equation}
P^{\prime}=\sum_{n>0}np\left(  n\right)  x^{n-1}.
\label{pf.thm.pars.sigma1.P'=}%
\end{equation}
Multiplying this equality by $x$, we obtain the somewhat nicer%
\begin{align}
xP^{\prime}  &  =x\sum_{n>0}np\left(  n\right)  x^{n-1}=\sum_{n>0}np\left(
n\right)  \underbrace{xx^{n-1}}_{=x^{n}}=\sum_{n>0}np\left(  n\right)
x^{n}\nonumber\\
&  =\sum_{n\in\mathbb{N}}np\left(  n\right)  x^{n}
\label{pf.thm.pars.sigma1.xP'=}%
\end{align}
(here, we have extended the range of the sum to include $n=0$, which caused no
change to the value of the sum, since the newly added $n=0$ addend is
$0p\left(  0\right)  x^{0}=0$).

Let us now recall the notion of a logarithmic derivative (as defined in
Definition \ref{def.fps.loder.1}). The FPS $P$ has constant coefficient
$\left[  x^{0}\right]  P=p\left(  0\right)  =1$, thus belongs to
$\mathbb{Z}\left[  \left[  x\right]  \right]  _{1}$. Hence, its logarithmic
derivative $\operatorname*{loder}P$ is well-defined.

We shall now use $\operatorname*{loder}P$ to compute $P^{\prime}$ in a
roundabout way. Namely, we have%
\[
P=\sum_{n\in\mathbb{N}}p\left(  n\right)  x^{n}=\prod_{k=1}^{\infty}\dfrac
{1}{1-x^{k}}%
\]
(by Theorem \ref{thm.pars.main-gf}). In other words,%
\begin{equation}
P=\prod_{k>0}\dfrac{1}{1-x^{k}} \label{pf.thm.pars.sigma1.P=prod}%
\end{equation}
(since the product sign $\prod\limits_{k=1}^{\infty}$ is synonymous with
$\prod\limits_{k>0}$).

But Corollary \ref{cor.fps.loder.prodk} says that any $k$ FPSs $f_{1}%
,f_{2},\ldots,f_{k}\in K\left[  \left[  x\right]  \right]  _{1}$ (for any
commutative ring $K$) satisfy%
\[
\operatorname*{loder}\left(  f_{1}f_{2}\cdots f_{k}\right)
=\operatorname*{loder}\left(  f_{1}\right)  +\operatorname*{loder}\left(
f_{2}\right)  +\cdots+\operatorname*{loder}\left(  f_{k}\right)  .
\]
In other words, the logarithmic derivative of a product of $k$ FPSs (which
belong to $K\left[  \left[  x\right]  \right]  _{1}$) equals the sum of the
logarithmic derivatives of these $k$ FPSs. By a simple limit argument (using
Theorem \ref{thm.fps.lim.prod-lim}, Theorem \ref{thm.fps.lim.sum-lim},
Proposition \ref{prop.fps.lim.sum-quot} and Proposition
\ref{prop.fps.lim.deriv-lim}), we can extend this fact to an infinite product
(as long as it is multipliable). Thus, if $\left(  f_{1},f_{2},f_{3}%
,\ldots\right)  $ is a multipliable sequence of FPSs in $K\left[  \left[
x\right]  \right]  _{1}$ (for any commutative ring $K$), then%
\[
\operatorname*{loder}\left(  f_{1}f_{2}f_{3}\cdots\right)
=\operatorname*{loder}\left(  f_{1}\right)  +\operatorname*{loder}\left(
f_{2}\right)  +\operatorname*{loder}\left(  f_{3}\right)  +\cdots;
\]
equivalently,%
\[
\operatorname*{loder}\left(  \prod_{k>0}f_{k}\right)  =\sum_{k>0}%
\operatorname*{loder}\left(  f_{k}\right)  .
\]
Applying this to $K=\mathbb{Z}$ and $f_{k}=\dfrac{1}{1-x^{k}}$, we obtain%
\[
\operatorname*{loder}\left(  \prod_{k>0}\dfrac{1}{1-x^{k}}\right)  =\sum
_{k>0}\operatorname*{loder}\dfrac{1}{1-x^{k}}.
\]
In view of (\ref{pf.thm.pars.sigma1.P=prod}), we can rewrite this as%
\begin{equation}
\operatorname*{loder}P=\sum_{k>0}\operatorname*{loder}\dfrac{1}{1-x^{k}}.
\label{pf.thm.pars.sigma1.loderP=sumk}%
\end{equation}
However, every positive integer $k$ satisfies%
\begin{align*}
\operatorname*{loder}\underbrace{\dfrac{1}{1-x^{k}}}_{=\left(  1-x^{k}\right)
^{-1}}  &  =\operatorname*{loder}\left(  \left(  1-x^{k}\right)  ^{-1}\right)
\\
&  =-\underbrace{\operatorname*{loder}\left(  1-x^{k}\right)  }%
_{\substack{=\dfrac{\left(  1-x^{k}\right)  ^{\prime}}{1-x^{k}}\\\text{(by the
definition of a}\\\text{logarithmic derivative)}}}\ \ \ \ \ \ \ \ \ \ \left(
\begin{array}
[c]{c}%
\text{by Corollary \ref{cor.fps.loder.inv},}\\
\text{applied to }f=1-x^{k}%
\end{array}
\right) \\
&  =-\dfrac{\left(  1-x^{k}\right)  ^{\prime}}{1-x^{k}}=-\dfrac{-kx^{k-1}%
}{1-x^{k}}\ \ \ \ \ \ \ \ \ \ \left(  \text{since }\left(  1-x^{k}\right)
^{\prime}=-kx^{k-1}\right) \\
&  =\dfrac{kx^{k-1}}{1-x^{k}}=kx^{k-1}\underbrace{\left(  1-x^{k}\right)
^{-1}}_{\substack{=\sum_{m\in\mathbb{N}}\left(  x^{k}\right)  ^{m}\\\text{(by
the geometric series formula)}}}\\
&  =kx^{k-1}\sum_{m\in\mathbb{N}}\left(  x^{k}\right)  ^{m}=\sum
_{m\in\mathbb{N}}k\underbrace{x^{k-1}\left(  x^{k}\right)  ^{m}}%
_{=x^{k-1+km}=x^{k\left(  m+1\right)  -1}}\\
&  =\sum_{m\in\mathbb{N}}kx^{k\left(  m+1\right)  -1}\\
&  =\sum_{i>0}kx^{ki-1}\ \ \ \ \ \ \ \ \ \ \left(
\begin{array}
[c]{c}%
\text{here, we have substituted }i\\
\text{for }m+1\text{ in the sum}%
\end{array}
\right) \\
&  =\sum_{\substack{n>0;\\k\mid n}}kx^{n-1}\ \ \ \ \ \ \ \ \ \ \left(
\begin{array}
[c]{c}%
\text{here, we have substituted }n\\
\text{for }ki\text{ in the sum}%
\end{array}
\right)  .
\end{align*}
Hence, we can rewrite (\ref{pf.thm.pars.sigma1.loderP=sumk}) as%
\begin{align*}
\operatorname*{loder}P  &  =\underbrace{\sum_{k>0}\ \ \sum
_{\substack{n>0;\\k\mid n}}}_{=\sum_{n>0}\ \ \sum_{\substack{k>0;\\k\mid n}%
}}kx^{n-1}=\sum_{n>0}\ \ \sum_{\substack{k>0;\\k\mid n}}kx^{n-1}\\
&  =\sum_{n>0}\underbrace{\left(  \sum_{\substack{k>0;\\k\mid n}}k\right)
}_{\substack{=\left(  \text{sum of all positive divisors of }n\right)
\\=\sigma\left(  n\right)  \\\text{(by the definition of }\sigma\left(
n\right)  \text{)}}}x^{n-1}=\sum_{n>0}\sigma\left(  n\right)  x^{n-1}.
\end{align*}
Comparing this with%
\[
\operatorname*{loder}P=\dfrac{P^{\prime}}{P}\ \ \ \ \ \ \ \ \ \ \left(
\text{by Definition \ref{def.fps.loder.1}}\right)  ,
\]
we obtain
\[
\dfrac{P^{\prime}}{P}=\sum_{n>0}\sigma\left(  n\right)  x^{n-1}.
\]
Multiplying both sides of this equality by $xP$, we find%
\begin{align*}
xP^{\prime}  &  =xP\cdot\sum_{n>0}\sigma\left(  n\right)  x^{n-1}=P\cdot
\sum_{n>0}\sigma\left(  n\right)  \underbrace{xx^{n-1}}_{=x^{n}}%
=P\cdot\underbrace{\sum_{n>0}\sigma\left(  n\right)  x^{n}}_{=\sum_{k>0}%
\sigma\left(  k\right)  x^{k}=S}=P\cdot S=SP\\
&  =\sum_{m\in\mathbb{N}}\ \ \sum_{k=1}^{m}\sigma\left(  k\right)  p\left(
m-k\right)  x^{m}\ \ \ \ \ \ \ \ \ \ \left(  \text{by
(\ref{pf.thm.pars.sigma1.SP=})}\right) \\
&  =\sum_{n\in\mathbb{N}}\ \ \sum_{k=1}^{n}\sigma\left(  k\right)  p\left(
n-k\right)  x^{n}\ \ \ \ \ \ \ \ \ \ \left(
\begin{array}
[c]{c}%
\text{here, we have renamed the}\\
\text{summation index }m\text{ as }n
\end{array}
\right)  .
\end{align*}
Comparing this with (\ref{pf.thm.pars.sigma1.xP'=}), we find%
\[
\sum_{n\in\mathbb{N}}np\left(  n\right)  x^{n}=\sum_{n\in\mathbb{N}}%
\ \ \sum_{k=1}^{n}\sigma\left(  k\right)  p\left(  n-k\right)  x^{n}.
\]
Comparing coefficients in front of $x^{n}$ on both sides of this equality, we
find that%
\[
np\left(  n\right)  =\sum_{k=1}^{n}\sigma\left(  k\right)  p\left(
n-k\right)  \ \ \ \ \ \ \ \ \ \ \text{for each }n\in\mathbb{N}\text{.}%
\]
This proves Theorem \ref{thm.pars.sigma1}.
\end{proof}

More generally, the following holds:

\begin{theorem}
\label{thm.pars.sigma1-I}Let $I$ be a subset of $\left\{  1,2,3,\ldots
\right\}  $. For each $n\in\mathbb{N}$, let $p_{I}\left(  n\right)  $ be the
\# of partitions $\lambda$ of $n$ such that all parts of $\lambda$ belong to
$I$.

For any positive integer $n$, let $\sigma_{I}\left(  n\right)  $ denote the
sum of all positive divisors of $n$ that belong to $I$. (For example, if
$O=\left\{  \text{all odd positive integers}\right\}  $ and $E=\left\{
\text{all even positive integers}\right\}  $, then $\sigma_{O}\left(
6\right)  =1+3=4$ and $\sigma_{E}\left(  6\right)  =2+6=8$.)

For any $n\in\mathbb{N}$, we have%
\[
np_{I}\left(  n\right)  =\sum_{k=1}^{n}\sigma_{I}\left(  k\right)
p_{I}\left(  n-k\right)  .
\]

\end{theorem}

Of course, Theorem \ref{thm.pars.sigma1} is the particular case of Theorem
\ref{thm.pars.sigma1-I} for $I=\left\{  1,2,3,\ldots\right\}  $.

\begin{proof}
[Proof of Theorem \ref{thm.pars.sigma1-I} (sketched).]Argue as in our above
proof of Theorem \ref{thm.pars.sigma1}, making some replacements: In
particular, symbols like $\sum\limits_{k>0}$ and $\prod\limits_{k>0}$ must be
replaced by $\sum\limits_{k\in I}$ and $\prod\limits_{k\in I}$ (respectively),
and Theorem \ref{thm.pars.main-gf-parts-I} must be used instead of Theorem
\ref{thm.pars.main-gf}. We leave the details to the reader.
\end{proof}

We note that Theorem \ref{thm.pars.sigma1} and even the more general Theorem
\ref{thm.pars.sigma1-I} are not hard to prove in completely elementary ways
(see, e.g., Exercise \ref{exe.pars.sigma-I-elem}). Nevertheless, our above
proof using generating functions has some interesting consequences, one of
which we might soon explore.

\subsection{Euler's pentagonal number theorem}

The following definition looks somewhat quaint; why define a notation for a
specific quadratic function?

\begin{definition}
\label{def.pars.pent-num}For any $k\in\mathbb{Z}$, define a nonnegative
integer $w_{k}\in\mathbb{N}$ by%
\[
w_{k}=\dfrac{\left(  3k-1\right)  k}{2}.
\]
This is called the $k$\emph{-th pentagonal number}.
\end{definition}

Here is a table of these pentagonal numbers:%
\[%
\begin{tabular}
[c]{|c||c|c|c|c|c|c|c|c|c|c|c|c|c|}\hline
$k$ & $\cdots$ & $-5$ & $-4$ & $-3$ & $-2$ & $-1$ & $0$ & $1$ & $2$ & $3$ &
$4$ & $5$ & $\cdots$\\\hline
$w_{k}$ & $\cdots$ & $40$ & $26$ & $15$ & $7$ & $2$ & $0$ & $1$ & $5$ & $12$ &
$22$ & $35$ & $\cdots$\\\hline
\end{tabular}
\ \ \ \ \ .
\]
Note that $w_{k}$ really is a nonnegative integer for any $k\in\mathbb{Z}$
(check this!). The name \textquotedblleft pentagonal numbers\textquotedblright%
\ is historically motivated (see
\href{https://en.wikipedia.org/wiki/Pentagonal_number}{the Wikipedia page} for
details); the only thing we need to know about them (beside their definition)
is the fact that they are nonnegative integers and grow quadratically with $k$
in both directions (i.e., when $k\rightarrow\infty$ and when $k\rightarrow
-\infty$). The latter fact ensures that the infinite sum $\sum_{k\in
\mathbb{Z}}\left(  -1\right)  ^{k}x^{w_{k}}$ is a well-defined FPS in
$\mathbb{Z}\left[  \left[  x\right]  \right]  $. Rather surprisingly, this
infinite sum coincides with a particularly simple infinite product:

\begin{theorem}
[Euler's pentagonal number theorem]\label{thm.pars.pent}We have%
\[
\prod_{k=1}^{\infty}\left(  1-x^{k}\right)  =\sum_{k\in\mathbb{Z}}\left(
-1\right)  ^{k}x^{w_{k}}.
\]

\end{theorem}

Let us write this out concretely:%
\begin{align*}
&  \prod_{k=1}^{\infty}\left(  1-x^{k}\right) \\
&  =\sum_{k\in\mathbb{Z}}\left(  -1\right)  ^{k}x^{w_{k}}\\
&  =\cdots+x^{w_{-4}}-x^{w_{-3}}+x^{w_{-2}}-x^{w_{-1}}+x^{w_{0}}-x^{w_{1}%
}+x^{w_{2}}-x^{w_{3}}+x^{w_{4}}-x^{w_{5}}\pm\cdots\\
&  =\cdots+x^{26}-x^{15}+x^{7}-x^{2}+1-x+x^{5}-x^{12}+x^{22}-x^{35}\pm\cdots\\
&  =1-x-x^{2}+x^{5}+x^{7}-x^{12}-x^{15}+x^{22}+x^{26}\pm\cdots.
\end{align*}

We will prove Theorem \ref{thm.pars.pent} in the next section (as a particular
case of Jacobi's Triple Product Identity).\footnote{See \cite{Bell06} for the
history of Theorem \ref{thm.pars.pent}.} First, let us use it to derive the
following recursive formula for the partition numbers $p\left(  n\right)  $:

\begin{corollary}
\label{cor.pars.pn-rec}For each positive integer $n$, we have%
\begin{align*}
p\left(  n\right)   &  =\sum_{\substack{k\in\mathbb{Z};\\k\neq0}}\left(
-1\right)  ^{k-1}p\left(  n-w_{k}\right) \\
&  =p\left(  n-1\right)  +p\left(  n-2\right)  -p\left(  n-5\right)  -p\left(
n-7\right) \\
&  \ \ \ \ \ \ \ \ \ \ +p\left(  n-12\right)  +p\left(  n-15\right)  -p\left(
n-22\right)  -p\left(  n-26\right)  \pm\cdots.
\end{align*}

\end{corollary}

\begin{proof}
[Proof of Corollary \ref{cor.pars.pn-rec} using Theorem \ref{thm.pars.pent}%
.]We have%
\[
\sum_{m\in\mathbb{N}}p\left(  m\right)  x^{m}=\sum_{n\in\mathbb{N}}p\left(
n\right)  x^{n}=\prod\limits_{k=1}^{\infty}\dfrac{1}{1-x^{k}}%
\]
(by Theorem \ref{thm.pars.main-gf}) and%
\[
\sum_{k\in\mathbb{Z}}\left(  -1\right)  ^{k}x^{w_{k}}=\prod\limits_{k=1}%
^{\infty}\left(  1-x^{k}\right)
\]
(by Theorem \ref{thm.pars.pent}). Multiplying these two equalities, we obtain%
\begin{align}
\left(  \sum_{m\in\mathbb{N}}p\left(  m\right)  x^{m}\right)  \cdot\left(
\sum_{k\in\mathbb{Z}}\left(  -1\right)  ^{k}x^{w_{k}}\right)   &  =\left(
\prod\limits_{k=1}^{\infty}\dfrac{1}{1-x^{k}}\right)  \cdot\left(
\prod\limits_{k=1}^{\infty}\left(  1-x^{k}\right)  \right) \nonumber\\
&  =\prod\limits_{k=1}^{\infty}\underbrace{\left(  \dfrac{1}{1-x^{k}}%
\cdot\left(  1-x^{k}\right)  \right)  }_{=1}\nonumber\\
&  =1. \label{pf.cor.pars.pn-rec.1}%
\end{align}

Now, let us fix a positive integer $n$. We shall compare the $x^{n}%
$-coefficients on both sides of (\ref{pf.cor.pars.pn-rec.1}).

The $x^{n}$-coefficient on the left hand side of (\ref{pf.cor.pars.pn-rec.1})
is%
\begin{align*}
\sum_{\substack{m\in\mathbb{N};\\k\in\mathbb{Z};\\m+w_{k}=n}}p\left(
m\right)  \cdot\left(  -1\right)  ^{k}  &  =\sum_{\substack{k\in
\mathbb{Z};\\n-w_{k}\geq0}}p\left(  n-w_{k}\right)  \cdot\left(  -1\right)
^{k}\\
&  \ \ \ \ \ \ \ \ \ \ \ \ \ \ \ \ \ \ \ \ \left(
\begin{array}
[c]{c}%
\text{here, we have replaced }m\text{ by }n-w_{k}\\
\text{in the sum, since the condition }m+w_{k}=n\\
\text{forces }m\text{ to be }n-w_{k}%
\end{array}
\right) \\
&  =\sum_{k\in\mathbb{Z}}p\left(  n-w_{k}\right)  \cdot\left(  -1\right)
^{k}\\
&  \ \ \ \ \ \ \ \ \ \ \ \ \ \ \ \ \ \ \ \ \left(
\begin{array}
[c]{c}%
\text{here, we have extended the range of}\\
\text{summation; this does not change the sum,}\\
\text{since }p\left(  n-w_{k}\right)  =0\text{ whenever }n-w_{k}<0
\end{array}
\right) \\
&  =\sum_{k\in\mathbb{Z}}\left(  -1\right)  ^{k}p\left(  n-w_{k}\right) \\
&  =\underbrace{\left(  -1\right)  ^{0}}_{=1}p\left(  n-\underbrace{w_{0}%
}_{=0}\right)  +\sum_{\substack{k\in\mathbb{Z};\\k\neq0}}\left(  -1\right)
^{k}p\left(  n-w_{k}\right) \\
&  =p\left(  n\right)  +\sum_{\substack{k\in\mathbb{Z};\\k\neq0}}\left(
-1\right)  ^{k}p\left(  n-w_{k}\right)  .
\end{align*}
But the $x^{n}$-coefficient on the right hand side of
(\ref{pf.cor.pars.pn-rec.1}) is $0$ (since $n$ is positive). Hence, comparing
the coefficients yields
\[
p\left(  n\right)  +\sum_{\substack{k\in\mathbb{Z};\\k\neq0}}\left(
-1\right)  ^{k}p\left(  n-w_{k}\right)  =0.
\]
Solving this for $p\left(  n\right)  $, we find%
\[
p\left(  n\right)  =-\sum_{\substack{k\in\mathbb{Z};\\k\neq0}}\left(
-1\right)  ^{k}p\left(  n-w_{k}\right)  =\sum_{\substack{k\in\mathbb{Z}%
;\\k\neq0}}\left(  -1\right)  ^{k-1}p\left(  n-w_{k}\right)  .
\]
Corollary \ref{cor.pars.pn-rec} is thus proved.
\end{proof}

\subsection{Jacobi's triple product identity}

\subsubsection{The identity}

Instead of proving Theorem \ref{thm.pars.pent} directly, we shall prove a
stronger result: \emph{Jacobi's triple product identity}. This identity can be
stated as follows:%
\begin{equation}
\prod_{n>0}\left(  \left(  1+q^{2n-1}z\right)  \left(  1+q^{2n-1}%
z^{-1}\right)  \left(  1-q^{2n}\right)  \right)  =\sum_{\ell\in\mathbb{Z}%
}q^{\ell^{2}}z^{\ell}. \label{eq.pars.jtp.jtp}%
\end{equation}
What are $q$ and $z$ here? It appears that (\ref{eq.pars.jtp.jtp}) should be
an identity between multivariate Laurent series (in the indeterminates $q$ and
$z$), but we have never defined such a concept. Multivariate Laurent series
can indeed be defined, but this is not as easy as the univariate case and
involves some choices (see \cite{ApaKau13} for details).

\begin{noncompile}
However, it is not hard to see that any negative power of $z$ in
(\ref{eq.pars.jtp.jtp}) is \textquotedblleft compensated\textquotedblright\ by
a positive power of $q$. For example, the $q^{2n-1}z^{-1}$ term on the left
hand side will always have total degree $\geq0$, and so will the $q^{\ell^{2}%
}z^{\ell}$ term on the right hand side. So this is not as bad as Laurent
series in general. (Commented this out because it is irrelevant and confusing.
The same could be said about $\sum_{\ell\in\mathbb{Z}}q^{-\ell}z^{\ell}$,
which is a bad series.)
\end{noncompile}

A simpler ring in which the identity (\ref{eq.pars.jtp.jtp}) can be placed is
$\left(  \mathbb{Z}\left[  z^{\pm}\right]  \right)  \left[  \left[  q\right]
\right]  $ (that is, the ring of FPSs in the indeterminate $q$ whose
coefficients are Laurent polynomials over $\mathbb{Z}$ in the indeterminate
$z$). In other words, we state the following:

\begin{theorem}
[Jacobi's triple product identity, take 1]\label{thm.pars.jtp1}In the ring
$\left(  \mathbb{Z}\left[  z^{\pm}\right]  \right)  \left[  \left[  q\right]
\right]  $, we have%
\[
\prod_{n>0}\left(  \left(  1+q^{2n-1}z\right)  \left(  1+q^{2n-1}%
z^{-1}\right)  \left(  1-q^{2n}\right)  \right)  =\sum_{\ell\in\mathbb{Z}%
}q^{\ell^{2}}z^{\ell}.
\]

\end{theorem}

However, we aren't just planning to view this identity as a formal identity
between power series; instead, we will later evaluate both sides at certain
powers of another indeterminate $x$ (i.e., we will set $q=x^{a}$ and $z=x^{b}$
for some positive integers $a$ and $b$). Alas, this is not an operation
defined on the whole ring $\left(  \mathbb{Z}\left[  z^{\pm}\right]  \right)
\left[  \left[  q\right]  \right]  $. For example, setting $q=x$ and $z=x$ in
the sum $\sum_{\ell\in\mathbb{N}}q^{\ell}z^{-\ell}\in\left(  \mathbb{Z}\left[
z^{\pm}\right]  \right)  \left[  \left[  q\right]  \right]  $ yields the
nonsensical sum $\sum_{\ell\in\mathbb{N}}x^{\ell}x^{-\ell}=\sum_{\ell
\in\mathbb{N}}1$.

Thus, Theorem \ref{thm.pars.jtp1} is not a version of Jacobi's triple product
identity that we can use for our purposes. Instead, let us interpret the
identity (\ref{eq.pars.jtp.jtp}) in a different way: Instead of treating $q$
and $z$ as indeterminates, I will set them to be powers of a single
indeterminate $x$ (more precisely, scalar multiples of such powers). This
leads us to the following version of the identity:

\begin{theorem}
[Jacobi's triple product identity, take 2]\label{thm.pars.jtp2}Let $a$ and $b$
be two integers such that $a>0$ and $a\geq\left\vert b\right\vert $. Let
$u,v\in\mathbb{Q}$ be rational numbers with $v\neq0$. In the ring
$\mathbb{Q}\left(  \left(  x\right)  \right)  $, set $q=ux^{a}$ and $z=vx^{b}%
$. Then,%
\[
\prod_{n>0}\left(  \left(  1+q^{2n-1}z\right)  \left(  1+q^{2n-1}%
z^{-1}\right)  \left(  1-q^{2n}\right)  \right)  =\sum_{\ell\in\mathbb{Z}%
}q^{\ell^{2}}z^{\ell}.
\]

\end{theorem}

Before we start proving this theorem, let us check that the infinite product
on its left hand side and the infinite sum on its right are well-defined:

\begin{itemize}
\item The infinite product is%
\begin{align*}
&  \prod_{n>0}\left(  \left(  1+q^{2n-1}z\right)  \left(  1+q^{2n-1}%
z^{-1}\right)  \left(  1-q^{2n}\right)  \right) \\
&  =\prod_{n>0}\left(  \left(  1+\left(  ux^{a}\right)  ^{2n-1}\left(
vx^{b}\right)  \right)  \left(  1+\left(  ux^{a}\right)  ^{2n-1}\left(
vx^{b}\right)  ^{-1}\right)  \left(  1-\left(  ux^{a}\right)  ^{2n}\right)
\right) \\
&  =\prod_{n>0}\left(  \left(  1+u^{2n-1}vx^{\left(  2n-1\right)  a+b}\right)
\left(  1+u^{2n-1}v^{-1}x^{\left(  2n-1\right)  a-b}\right)  \left(
1-u^{2n}x^{2na}\right)  \right)  .
\end{align*}
All factors in this product belong to the ring $\mathbb{Q}\left[  \left[
x\right]  \right]  $ (not just to $\mathbb{Q}\left(  \left(  x\right)
\right)  $), since the exponents $\left(  2n-1\right)  a+b$ and $\left(
2n-1\right)  a-b$ and $2na$ are always nonnegative for any $n>0$ (indeed, for
any $n>0$, we have $\underbrace{\left(  2n-1\right)  }_{\geq1}\underbrace{a}%
_{\geq\left\vert b\right\vert }+\,b\geq\left\vert b\right\vert +b\geq0$ and
$\underbrace{\left(  2n-1\right)  }_{\geq1}\underbrace{a}_{\geq\left\vert
b\right\vert }-\,b\geq\left\vert b\right\vert -b\geq0$ and $2n\underbrace{a}%
_{>0}>0$). Moreover, this product is multipliable, because

\begin{itemize}
\item the number $\left(  2n-1\right)  a+b$ grows linearly when $n\rightarrow
\infty$ (since $a>0$);

\item the number $\left(  2n-1\right)  a-b$ grows linearly when $n\rightarrow
\infty$ (since $a>0$);

\item the number $2na$ grows linearly when $n\rightarrow\infty$ (since $a>0$).
\end{itemize}

Thus, the infinite product is well-defined.

This argument also shows that the three subproducts $\prod_{n>0}\left(
1+q^{2n-1}z\right)  $ and $\prod_{n>0}\left(  1+q^{2n-1}z^{-1}\right)  $ and
$\prod_{n>0}\left(  1-q^{2n}\right)  $ are well-defined.

\item The infinite sum is%
\[
\sum_{\ell\in\mathbb{Z}}q^{\ell^{2}}z^{\ell}=\sum_{\ell\in\mathbb{Z}}\left(
ux^{a}\right)  ^{\ell^{2}}\left(  vx^{b}\right)  ^{\ell}=\sum_{\ell
\in\mathbb{Z}}u^{\ell^{2}}v^{\ell}x^{a\ell^{2}+b\ell}.
\]
All addends in this sum belong to the ring $\mathbb{Q}\left[  \left[
x\right]  \right]  $ (not just to $\mathbb{Q}\left(  \left(  x\right)
\right)  $), since the exponent $a\ell^{2}+b\ell$ is always nonnegative for
any $\ell\in\mathbb{Z}$ (indeed, we have $\underbrace{a}_{\geq\left\vert
b\right\vert }\underbrace{\ell^{2}}_{=\left\vert \ell\right\vert ^{2}%
}+\underbrace{b\ell}_{\geq-\left\vert b\ell\right\vert =-\left\vert
b\right\vert \cdot\left\vert \ell\right\vert }\geq\left\vert b\right\vert
\cdot\left\vert \ell\right\vert ^{2}-\left\vert b\right\vert \cdot\left\vert
\ell\right\vert =\underbrace{\left\vert b\right\vert }_{\geq0}\cdot
\underbrace{\left\vert \ell\right\vert \cdot\left(  \left\vert \ell\right\vert
-1\right)  }_{\geq0}\geq0$ for any $\ell\in\mathbb{Z}$). Moreover, this sum is
summable, because $a\ell^{2}+b\ell$ grows quadratically when $\ell
\rightarrow+\infty$ or $\ell\rightarrow-\infty$ (since $a>0$).
\end{itemize}

\subsubsection{Jacobi implies Euler}

We will give a proof of Jacobi's triple product identity that works equally
for both versions of it (Theorem \ref{thm.pars.jtp1} and Theorem
\ref{thm.pars.jtp2}). But first, let us see how it yields Euler's pentagonal
number theorem as a particular case.

\begin{proof}
[Proof of Theorem \ref{thm.pars.pent} using Theorem \ref{thm.pars.jtp2}.]Set
$q=x^{3}$ and $z=-x$ in Theorem \ref{thm.pars.jtp2}. (This means that we apply
Theorem \ref{thm.pars.jtp2} to $a=3$ and $b=1$ and $u=1$ and $v=-1$). We get%
\begin{align}
&  \prod_{n>0}\left(  \left(  1+\left(  x^{3}\right)  ^{2n-1}\left(
-x\right)  \right)  \left(  1+\left(  x^{3}\right)  ^{2n-1}\left(  -x\right)
^{-1}\right)  \left(  1-\left(  x^{3}\right)  ^{2n}\right)  \right)
\nonumber\\
&  =\sum_{\ell\in\mathbb{Z}}\left(  x^{3}\right)  ^{\ell^{2}}\left(
-x\right)  ^{\ell}. \label{pf.thm.pars.pent.1}%
\end{align}
The left hand side of this equality simplifies as follows:%
\begin{align*}
&  \prod_{n>0}\left(  \left(  1+\underbrace{\left(  x^{3}\right)
^{2n-1}\left(  -x\right)  }_{\substack{=-x^{3\left(  2n-1\right)
+1}\\=-x^{6n-2}}}\right)  \left(  1+\underbrace{\left(  x^{3}\right)
^{2n-1}\left(  -x\right)  ^{-1}}_{\substack{=-x^{3\left(  2n-1\right)
-1}\\=-x^{6n-4}}}\right)  \left(  1-\underbrace{\left(  x^{3}\right)  ^{2n}%
}_{=x^{6n}}\right)  \right) \\
&  =\prod_{n>0}\left(  \left(  1-\underbrace{x^{6n-2}}_{=\left(  x^{2}\right)
^{3n-1}}\right)  \left(  1-\underbrace{x^{6n-4}}_{=\left(  x^{2}\right)
^{3n-2}}\right)  \left(  1-\underbrace{x^{6n}}_{=\left(  x^{2}\right)  ^{3n}%
}\right)  \right) \\
&  =\prod_{n>0}\left(  \left(  1-\left(  x^{2}\right)  ^{3n-1}\right)  \left(
1-\left(  x^{2}\right)  ^{3n-2}\right)  \left(  1-\left(  x^{2}\right)
^{3n}\right)  \right) \\
&  =\prod_{k>0}\left(  1-\left(  x^{2}\right)  ^{k}\right)  ,
\end{align*}
since each positive integer $k$ can be uniquely represented as $3n-1$ or
$3n-2$ or $3n$ for some positive integer $n$.

Comparing this with (\ref{pf.thm.pars.pent.1}), we obtain%
\begin{align}
&  \prod_{k>0}\left(  1-\left(  x^{2}\right)  ^{k}\right) \nonumber\\
&  =\sum_{\ell\in\mathbb{Z}}\underbrace{\left(  x^{3}\right)  ^{\ell^{2}%
}\left(  -x\right)  ^{\ell}}_{=\left(  -1\right)  ^{\ell}x^{3\ell^{2}+\ell}%
}=\sum_{\ell\in\mathbb{Z}}\left(  -1\right)  ^{\ell}\underbrace{x^{3\ell
^{2}+\ell}}_{\substack{=x^{2w_{-\ell}}\\\text{(since }3\ell^{2}+\ell=\left(
3\ell+1\right)  \ell=\left(  3\left(  -\ell\right)  -1\right)  \left(
-\ell\right)  =2w_{-\ell}\\\text{(because }w_{-\ell}\text{ is defined as
}\left(  3\left(  -\ell\right)  -1\right)  \left(  -\ell\right)  /2\text{))}%
}}\nonumber\\
&  =\sum_{\ell\in\mathbb{Z}}\underbrace{\left(  -1\right)  ^{\ell}}_{=\left(
-1\right)  ^{-\ell}}\underbrace{x^{2w_{-\ell}}}_{=\left(  x^{2}\right)
^{w_{-\ell}}}=\sum_{\ell\in\mathbb{Z}}\left(  -1\right)  ^{-\ell}\left(
x^{2}\right)  ^{w_{-\ell}}\nonumber\\
&  =\sum_{k\in\mathbb{Z}}\left(  -1\right)  ^{k}\left(  x^{2}\right)  ^{w_{k}}
\label{pf.thm.pars.pent.2}%
\end{align}
(here, we have substituted $k$ for $-\ell$ in the sum).

Now, let us \textquotedblleft substitute $x$ for $x^{2}$\textquotedblright\ in
this equality (see below for how this works). As a result, we obtain%
\[
\prod_{k>0}\left(  1-x^{k}\right)  =\sum_{k\in\mathbb{Z}}\left(  -1\right)
^{k}x^{w_{k}}.
\]
This is Euler's pentagonal number theorem (Theorem \ref{thm.pars.pent}).
\end{proof}

What did I mean by \textquotedblleft substituting $x$ for $x^{2}%
$\textquotedblright? I meant using the following simple fact:

\begin{lemma}
\label{lem.fps.fxx=gxx}Let $K$ be a commutative ring. Let $f$ and $g$ be two
FPSs in $K\left[  \left[  x\right]  \right]  $. Assume that $f\left[
x^{2}\right]  =g\left[  x^{2}\right]  $. Then, $f=g$.
\end{lemma}

\begin{proof}
This is easy: Write $f$ and $g$ as $f=\sum_{n\in\mathbb{N}}f_{n}x^{n}$ and
$g=\sum_{n\in\mathbb{N}}g_{n}x^{n}$ where $f_{0},f_{1},f_{2},\ldots\in K$ and
$g_{0},g_{1},g_{2},\ldots\in K$. Then, $f\left[  x^{2}\right]  =\sum
_{n\in\mathbb{N}}f_{n}\left(  x^{2}\right)  ^{n}=\sum_{n\in\mathbb{N}}%
f_{n}x^{2n}$ and similarly $g\left[  x^{2}\right]  =\sum_{n\in\mathbb{N}}%
g_{n}x^{2n}$. Thus, our assumption $f\left[  x^{2}\right]  =g\left[
x^{2}\right]  $ rewrites as $\sum_{n\in\mathbb{N}}f_{n}x^{2n}=\sum
_{n\in\mathbb{N}}g_{n}x^{2n}$. Comparing $x^{2n}$-coefficients in this
equality, we conclude that $f_{n}=g_{n}$ for each $n\in\mathbb{N}$. Hence,
$\sum_{n\in\mathbb{N}}\underbrace{f_{n}}_{=g_{n}}x^{n}=\sum_{n\in\mathbb{N}%
}g_{n}x^{n}$. In other words, $f=g$ (since $f=\sum_{n\in\mathbb{N}}f_{n}x^{n}$
and $g=\sum_{n\in\mathbb{N}}g_{n}x^{n}$). This proves Lemma
\ref{lem.fps.fxx=gxx}.
\end{proof}

Lemma \ref{lem.fps.fxx=gxx} justifies our \textquotedblleft substituting $x$
for $x^{2}$\textquotedblright\ in the above proof; indeed, we can apply Lemma
\ref{lem.fps.fxx=gxx} to $K=\mathbb{Q}$ and $f=\prod_{k>0}\left(
1-x^{k}\right)  $ and $g=\sum_{k\in\mathbb{Z}}\left(  -1\right)  ^{k}x^{w_{k}%
}$ (because (\ref{pf.thm.pars.pent.2}) says that these two FPSs $f$ and $g$
satisfy $f\left[  x^{2}\right]  =g\left[  x^{2}\right]  $), and consequently
obtain $\prod_{k>0}\left(  1-x^{k}\right)  =\sum_{k\in\mathbb{Z}}\left(
-1\right)  ^{k}x^{w_{k}}$. Thus, Theorem \ref{thm.pars.pent} is proved using
Theorem \ref{thm.pars.jtp2}. It therefore remains to prove the latter.

\subsubsection{Proof of Jacobi's triple product identity}

The following proof is due to Borcherds, and I have taken it from
\cite[\S 8.3]{Camero16} (note that \cite[\S 11.2]{Loehr-BC} gives essentially
the same proof, albeit in a different language).

\begin{proof}
[Proof of Theorem \ref{thm.pars.jtp1} and Theorem \ref{thm.pars.jtp2}.]The
following argument applies equally to Theorem \ref{thm.pars.jtp1} and to
Theorem \ref{thm.pars.jtp2}. (The meanings of $q$ and $z$ differ between these
two theorems, but all the infinite sums and products considered below are
well-defined in either case.)

We will use a somewhat physics-inspired language:

\begin{itemize}
\item A \emph{level} will mean a number of the form $p+\dfrac{1}{2}$ with
$p\in\mathbb{Z}$. (Thus, there is exactly one level midway between any two
consecutive integers.)

\item A \emph{state} will mean a set of levels that contains

\begin{itemize}
\item all but finitely many negative levels, and

\item only finitely many positive levels.
\end{itemize}
\end{itemize}

Here is an example of a state:%
\[
\left\{  \dfrac{-5}{2},\ \ \dfrac{-1}{2},\ \ \dfrac{1}{2},\ \ \dfrac{3}%
{2},\ \ \dfrac{7}{2},\ \ \dfrac{13}{2}\right\}  \cup\left\{  \text{all levels
}\leq\dfrac{-9}{2}\right\}  .
\]
Visually, it can be represented as follows:%
\[%
\begin{tikzpicture}[scale=0.6]
\tikzset{myptr/.style={decoration={markings,mark=at position 1 with {\arrow
[scale=3,>=stealth]{>}}},postaction={decorate}}}
\draw[myptr] (-12, 0) -- (15.5, 0);
\filldraw[draw=black, fill=white] (-0.5, 0) circle [radius=0.3];
\filldraw[draw=black, fill=black] (-1.5, 0) circle [radius=0.3];
\filldraw[draw=black, fill=white] (-2.5, 0) circle [radius=0.3];
\filldraw[draw=black, fill=black] (-3.5, 0) circle [radius=0.3];
\filldraw[draw=black, fill=white] (-4.5, 0) circle [radius=0.3];
\filldraw[draw=black, fill=white] (-5.5, 0) circle [radius=0.3];
\filldraw[draw=black, fill=white] (-6.5, 0) circle [radius=0.3];
\filldraw[draw=black, fill=white] (-7.5, 0) circle [radius=0.3];
\filldraw[draw=black, fill=white] (-8.5, 0) circle [radius=0.3];
\filldraw[draw=black, fill=white] (-9.5, 0) circle [radius=0.3];
\filldraw[draw=black, fill=white] (-10.5, 0) circle [radius=0.3];
\filldraw[draw=black, fill=white] (-11.5, 0) circle [radius=0.3];
\node at (0, 1) [label=left:$\cdots\cdots\cdots\text{negative levels}
\cdots\cdots\cdots$] {};
\node at (-4, -1) [label=left:$\text{only \raisebox{-0.5ex}{\tikz
\draw[black,fill=white] (0,0) circle (1ex);}\,s left of here}$] {};
\draw(0, 0.7) -- (0, -1.7);
\node at (0, -2) {$0$};
\filldraw[draw=black, fill=white] (0.5, 0) circle [radius=0.3];
\filldraw[draw=black, fill=white] (1.5, 0) circle [radius=0.3];
\filldraw[draw=black, fill=black] (2.5, 0) circle [radius=0.3];
\filldraw[draw=black, fill=white] (3.5, 0) circle [radius=0.3];
\filldraw[draw=black, fill=black] (4.5, 0) circle [radius=0.3];
\filldraw[draw=black, fill=black] (5.5, 0) circle [radius=0.3];
\filldraw[draw=black, fill=white] (6.5, 0) circle [radius=0.3];
\filldraw[draw=black, fill=black] (7.5, 0) circle [radius=0.3];
\filldraw[draw=black, fill=black] (8.5, 0) circle [radius=0.3];
\filldraw[draw=black, fill=black] (9.5, 0) circle [radius=0.3];
\filldraw[draw=black, fill=black] (10.5, 0) circle [radius=0.3];
\filldraw[draw=black, fill=black] (11.5, 0) circle [radius=0.3];
\filldraw[draw=black, fill=black] (12.5, 0) circle [radius=0.3];
\filldraw[draw=black, fill=black] (13.5, 0) circle [radius=0.3];
\filldraw[draw=black, fill=black] (14.5, 0) circle [radius=0.3];
\node at (0, 1) [label=right:$\cdots\cdots\cdots\text{positive levels}
\cdots\cdots\cdots\cdots$] {};
\node at (7, -1) [label=right:$\text{only \raisebox{-0.5ex}{\tikz
\draw[black,fill=black] (0,0) circle (1ex);}\,s right of here}$] {};
\end{tikzpicture}%
\]
where

\begin{itemize}
\item A white (=hollow) circle
\raisebox{-0.5ex}{\tikz\draw[black,fill=white] (0,0) circle (1ex);}%
\ means a level that is contained in the state (you can think of it as an
\textquotedblleft electron\textquotedblright).

\item A black (=filled) circle
\raisebox{-0.5ex}{\tikz\draw[black,fill=black] (0,0) circle (1ex);}
means a level that is not contained in the state (think of it as a
\textquotedblleft hole\textquotedblright).
\end{itemize}

For any state $S$,

\begin{itemize}
\item we define the \emph{energy} of $S$ to be%
\[
\operatorname*{energy}S:=\sum_{\substack{p>0;\\p\in S}}\underbrace{2p}%
_{>0}-\sum_{\substack{p<0;\\p\notin S}}\underbrace{2p}_{<0}\in\mathbb{N}%
\]
(where the summation index $p$ in the first sum runs over the finitely many
positive levels contained in $S$, while the summation index $p$ in the second
sum runs over the finitely many negative levels not contained in $S$).

\item we define the \emph{particle number} of $S$ to be%
\begin{align*}
\operatorname*{parnum}S:=  &  \left(  \text{\# of levels }p>0\text{ such that
}p\in S\right) \\
&  \ \ \ \ \ \ \ \ \ \ -\left(  \text{\# of levels }p<0\text{ such that
}p\notin S\right)  \in\mathbb{Z}.
\end{align*}

\end{itemize}

For instance, in the above example, we have%
\[
\operatorname*{energy}S=1+3+7+13-\left(  -3\right)  -\left(  -7\right)  =34
\]
and%
\[
\operatorname*{parnum}S=4-2=2.
\]

We want to prove the identity%
\[
\prod_{n>0}\left(  \left(  1+q^{2n-1}z\right)  \left(  1+q^{2n-1}%
z^{-1}\right)  \left(  1-q^{2n}\right)  \right)  =\sum_{\ell\in\mathbb{Z}%
}q^{\ell^{2}}z^{\ell}.
\]
We will first transform this identity into an equivalent one: Namely, we move
the $\left(  1-q^{2n}\right)  $ factors from the left hand side to the right
hand side by multiplying both sides with $\prod_{n>0}\left(  1-q^{2n}\right)
^{-1}$. Thus, we can rewrite our identity as%
\[
\prod_{n>0}\left(  \left(  1+q^{2n-1}z\right)  \left(  1+q^{2n-1}%
z^{-1}\right)  \right)  =\left(  \sum_{\ell\in\mathbb{Z}}q^{\ell^{2}}z^{\ell
}\right)  \prod_{n>0}\left(  1-q^{2n}\right)  ^{-1}.
\]
We will prove this new identity by showing that both of its sides are%
\[
\sum_{S\text{ is a state}}q^{\operatorname*{energy}S}z^{\operatorname*{parnum}%
S}.
\]

\textit{Left hand side:} We have%
\begin{align}
&  \prod_{n>0}\left(  \left(  1+q^{2n-1}z\right)  \left(  1+q^{2n-1}%
z^{-1}\right)  \right) \nonumber\\
&  =\left(  \prod_{n>0}\left(  1+q^{2n-1}z\right)  \right)  \left(
\prod_{n>0}\left(  1+q^{2n-1}z^{-1}\right)  \right) \nonumber\\
&  =\left(  \prod_{p\text{ is a positive level}}\left(  1+q^{2p}z\right)
\right)  \left(  \prod_{p\text{ is a negative level}}\left(  1+q^{-2p}%
z^{-1}\right)  \right) \nonumber\\
&  \ \ \ \ \ \ \ \ \ \ \ \ \ \ \ \ \ \ \ \ \left(
\begin{array}
[c]{c}%
\text{here, we have substituted }p+\dfrac{1}{2}\text{ for }n\text{ in the
first product,}\\
\text{and have substituted }-p+\dfrac{1}{2}\text{ for }n\text{ in the second
product}%
\end{array}
\right) \nonumber\\
&  =\left(  \sum_{\substack{P\text{ is a finite set}\\\text{of positive
levels}}}\ \ \prod_{p\in P}\left(  q^{2p}z\right)  \right)  \left(
\sum_{\substack{N\text{ is a finite set}\\\text{of negative levels}}%
}\ \ \prod_{p\in N}\left(  q^{-2p}z^{-1}\right)  \right) \nonumber\\
&  \ \ \ \ \ \ \ \ \ \ \ \ \ \ \ \ \ \ \ \ \left(  \text{here, we have
expanded both products using (\ref{eq.fps.prod.binary.prod-inf2})}\right)
\nonumber\\
&  =\sum_{\substack{P\text{ is a finite set}\\\text{of positive levels}%
}}\ \ \sum_{\substack{N\text{ is a finite set}\\\text{of negative levels}%
}}\ \ \underbrace{\prod_{p\in P}\left(  q^{2p}z\right)  \ \ \prod_{p\in
N}\left(  q^{-2p}z^{-1}\right)  }_{=q^{2\left(  \text{sum of elements of
}P\right)  -2\left(  \text{sum of elements of }N\right)  }z^{\left\vert
P\right\vert -\left\vert N\right\vert }}\nonumber\\
&  =\sum_{\substack{P\text{ is a finite set}\\\text{of positive levels}%
}}\ \ \sum_{\substack{N\text{ is a finite set}\\\text{of negative levels}%
}}q^{2\left(  \text{sum of elements of }P\right)  -2\left(  \text{sum of
elements of }N\right)  }z^{\left\vert P\right\vert -\left\vert N\right\vert
}\nonumber\\
&  =\sum_{S\text{ is a state}}\underbrace{q^{2\left(  \text{sum of positive
levels in }S\right)  -2\left(  \text{sum of negative levels not in }S\right)
}}_{\substack{=q^{\operatorname*{energy}S}\\\text{(by the definition of
}\operatorname*{energy}S\text{)}}}\nonumber\\
&  \ \ \ \ \ \ \ \ \ \ \ \ \ \ \ \ \ \ \ \ \underbrace{z^{\left(  \text{number
of positive levels in }S\right)  -\left(  \text{number of negative levels not
in }S\right)  }}_{\substack{=z^{\operatorname*{parnum}S}\\\text{(by the
definition of }\operatorname*{parnum}S\text{)}}}\nonumber\\
&  \ \ \ \ \ \ \ \ \ \ \ \ \ \ \ \ \ \ \ \ \ \ \ \ \ \ \ \ \ \ \left(
\begin{array}
[c]{c}%
\text{here, we have combined }P\text{ and }N\\
\text{into a single state }S:=P\cup\overline{N}\text{,}\\
\text{where }\overline{N}=\left\{  \text{all negative levels}\right\}
\setminus N
\end{array}
\right) \nonumber\\
&  =\sum_{S\text{ is a state}}q^{\operatorname*{energy}S}%
z^{\operatorname*{parnum}S}. \label{pf.thm.pars.jtp2.lhs}%
\end{align}

\textit{Right hand side:} Recall that%
\begin{align*}
\prod_{n>0}\left(  1-x^{n}\right)  ^{-1}  &  =\prod_{n>0}\dfrac{1}{1-x^{n}%
}=\sum_{n\in\mathbb{N}}p\left(  n\right)  x^{n}\ \ \ \ \ \ \ \ \ \ \left(
\text{by Theorem \ref{thm.pars.main-gf}}\right) \\
&  =\sum_{\substack{\lambda\text{ is a}\\\text{partition}}}x^{\left\vert
\lambda\right\vert }%
\end{align*}
(because the sum $\sum_{\substack{\lambda\text{ is a}\\\text{partition}%
}}x^{\left\vert \lambda\right\vert }$ contains each monomial $x^{n}$ precisely
$p\left(  n\right)  $ times). Substituting $q^{2}$ for $x$ in this equality,
we find%
\[
\prod_{n>0}\left(  1-\left(  q^{2}\right)  ^{n}\right)  ^{-1}=\sum
_{\substack{\lambda\text{ is a}\\\text{partition}}}\left(  q^{2}\right)
^{\left\vert \lambda\right\vert }.
\]
In other words,%
\[
\prod_{n>0}\left(  1-q^{2n}\right)  ^{-1}=\sum_{\substack{\lambda\text{ is
a}\\\text{partition}}}q^{2\left\vert \lambda\right\vert }.
\]
Multiplying both sides of this equality by $\sum_{\ell\in\mathbb{Z}}%
q^{\ell^{2}}z^{\ell}$, we obtain%
\begin{align}
\left(  \sum_{\ell\in\mathbb{Z}}q^{\ell^{2}}z^{\ell}\right)  \prod
_{n>0}\left(  1-q^{2n}\right)  ^{-1}  &  =\left(  \sum_{\ell\in\mathbb{Z}%
}q^{\ell^{2}}z^{\ell}\right)  \sum_{\substack{\lambda\text{ is a}%
\\\text{partition}}}q^{2\left\vert \lambda\right\vert }\nonumber\\
&  =\sum_{\ell\in\mathbb{Z}}\ \ \sum_{\substack{\lambda\text{ is
a}\\\text{partition}}}q^{\ell^{2}+2\left\vert \lambda\right\vert }z^{\ell}.
\label{pf.thm.pars.jtp2.rhs1}%
\end{align}
We want to show that this equals $\sum\limits_{S\text{ is a state}%
}q^{\operatorname*{energy}S}z^{\operatorname*{parnum}S}$. In order to do this,
we will find a bijection%
\[
\Phi_{\ell}:\left\{  \text{partitions}\right\}  \rightarrow\left\{
\text{states with particle number }\ell\right\}
\]
for each $\ell\in\mathbb{Z}$, and we will show that this bijection satisfies
\[
\operatorname*{energy}\left(  \Phi_{\ell}\left(  \lambda\right)  \right)
=\ell^{2}+2\left\vert \lambda\right\vert \ \ \ \ \ \ \ \ \ \ \text{for each
}\ell\in\mathbb{Z}\text{ and }\lambda\in\left\{  \text{partitions}\right\}  .
\]

Let us do this. Fix $\ell\in\mathbb{Z}$. We define the state $G_{\ell}$
(called the \textquotedblleft$\ell$-ground state\textquotedblright) by%
\[
G_{\ell}:=\left\{  \text{all levels }<\ell\right\}  =\left\{  \ell-\dfrac
{1}{2},\ \ \ell-\dfrac{3}{2},\ \ \ell-\dfrac{5}{2},\ \ \ldots\right\}  .
\]
Here is how it looks like (for $\ell$ positive):%
\[%
\begin{tikzpicture}[scale=0.6]
\tikzset{myptr/.style={decoration={markings,mark=at position 1 with {\arrow
[scale=3,>=stealth]{>}}},postaction={decorate}}}
\draw[myptr] (-12, 0) -- (15.5, 0);
\filldraw[draw=black, fill=white] (-0.5, 0) circle [radius=0.3];
\filldraw[draw=black, fill=white] (-1.5, 0) circle [radius=0.3];
\filldraw[draw=black, fill=white] (-2.5, 0) circle [radius=0.3];
\filldraw[draw=black, fill=white] (-3.5, 0) circle [radius=0.3];
\filldraw[draw=black, fill=white] (-4.5, 0) circle [radius=0.3];
\filldraw[draw=black, fill=white] (-5.5, 0) circle [radius=0.3];
\filldraw[draw=black, fill=white] (-6.5, 0) circle [radius=0.3];
\filldraw[draw=black, fill=white] (-7.5, 0) circle [radius=0.3];
\filldraw[draw=black, fill=white] (-8.5, 0) circle [radius=0.3];
\filldraw[draw=black, fill=white] (-9.5, 0) circle [radius=0.3];
\filldraw[draw=black, fill=white] (-10.5, 0) circle [radius=0.3];
\filldraw[draw=black, fill=white] (-11.5, 0) circle [radius=0.3];
\node at (0, 1) [label=left:$\cdots\cdots\cdots\text{negative levels}
\cdots\cdots\cdots$] {};
\node at (3, -1) [label=left:$\text{only \raisebox{-0.5ex}{\tikz
\draw[black,fill=white] (0,0) circle (1ex);}\,s left of here}$] {};
\draw(0, 0.7) -- (0, -1.7);
\node at (0, -2) {$0$};
\draw(3, 0.7) -- (3, -1.7);
\node at (3, -2) {$\ell$};
\filldraw[draw=black, fill=white] (0.5, 0) circle [radius=0.3];
\filldraw[draw=black, fill=white] (1.5, 0) circle [radius=0.3];
\filldraw[draw=black, fill=white] (2.5, 0) circle [radius=0.3];
\filldraw[draw=black, fill=black] (3.5, 0) circle [radius=0.3];
\filldraw[draw=black, fill=black] (4.5, 0) circle [radius=0.3];
\filldraw[draw=black, fill=black] (5.5, 0) circle [radius=0.3];
\filldraw[draw=black, fill=black] (6.5, 0) circle [radius=0.3];
\filldraw[draw=black, fill=black] (7.5, 0) circle [radius=0.3];
\filldraw[draw=black, fill=black] (8.5, 0) circle [radius=0.3];
\filldraw[draw=black, fill=black] (9.5, 0) circle [radius=0.3];
\filldraw[draw=black, fill=black] (10.5, 0) circle [radius=0.3];
\filldraw[draw=black, fill=black] (11.5, 0) circle [radius=0.3];
\filldraw[draw=black, fill=black] (12.5, 0) circle [radius=0.3];
\filldraw[draw=black, fill=black] (13.5, 0) circle [radius=0.3];
\filldraw[draw=black, fill=black] (14.5, 0) circle [radius=0.3];
\node at (0, 1) [label=right:$\cdots\cdots\cdots\text{positive levels}
\cdots\cdots\cdots\cdots$] {};
\node at (3, -1) [label=right:$\text{only \raisebox{-0.5ex}{\tikz
\draw[black,fill=black] (0,0) circle (1ex);}\,s right of here}$] {};
\end{tikzpicture}%
\]

If $\ell\geq0$, then this state $G_{\ell}$ has energy
\[
\operatorname*{energy}G_{\ell}=1+3+5+\cdots+\left(  2\ell-1\right)  =\ell^{2}%
\]
and particle number%
\[
\operatorname*{parnum}G_{\ell}=\ell-0=\ell.
\]
If $\ell<0$, then it has energy%
\[
\operatorname*{energy}G_{\ell}=-\left(  -1\right)  -\left(  -3\right)
-\left(  -5\right)  -\cdots-\left(  2\ell+1\right)  =\ell^{2}%
\]
and particle number%
\[
\operatorname*{parnum}G_{\ell}=0-\left(  -\ell\right)  =\ell.
\]
Note that the answers are the same in both cases. Thus, whatever sign $\ell$
has, we have%
\[
\operatorname*{energy}G_{\ell}=\ell^{2}\ \ \ \ \ \ \ \ \ \ \text{and}%
\ \ \ \ \ \ \ \ \ \ \operatorname*{parnum}G_{\ell}=\ell.
\]

If $S$ is a state, and if $p\in S$, and if $q$ is a positive integer such that
$p+q\notin S$, then we define a new state%
\[
\operatorname*{jump}\nolimits_{p,q}S:=\left(  S\setminus\left\{  p\right\}
\right)  \cup\left\{  p+q\right\}  .
\]
We say that this state $\operatorname*{jump}\nolimits_{p,q}S$ is obtained from
$S$ by letting the electron at level $p$ \emph{jump} $q$ steps to the right.
Note that $\operatorname*{jump}\nolimits_{p,q}S$ has the same particle number
as $S$ (check this!\footnote{There are three cases:
\par
\textit{Case 1:} We have $p>0$ (that is, the particle jumps from a positive
level to a positive level).
\par
\textit{Case 2:} We have $p<0$ and $p+q>0$ (that is, the particle jumps from a
negative level to a positive level).
\par
\textit{Case 3:} We have $p+q<0$ (that is, the particle jumps from a negative
level to a negative level).
\par
Each case is easy to check.}), whereas its energy is $2q$ higher than that of
$S$ (check this!\footnote{Again, the same three cases are to be considered.}).
Thus, a jumping electron raises the energy but keeps the particle number unchanged.

For any partition $\lambda=\left(  \lambda_{1},\lambda_{2},\ldots,\lambda
_{k}\right)  $, we define the state $E_{\ell,\lambda}$ (called an
\textquotedblleft excited state\textquotedblright) by starting with the $\ell
$-ground state $G_{\ell}$, and then successively letting the $k$ electrons at
the highest levels (which are -- from highest to lowest -- the levels
$\ell-1+\dfrac{1}{2},\ \ \ell-2+\dfrac{1}{2},\ \ \ldots,\ \ \ell-k+\dfrac
{1}{2}$) jump $\lambda_{1},\lambda_{2},\ldots,\lambda_{k}$ steps to the right,
respectively (starting with the rightmost electron). In other words,%
\begin{align*}
E_{\ell,\lambda}  &  :=\operatorname*{jump}\nolimits_{\ell-k+1/2,\ \ \lambda
_{k}}\left(  \cdots\left(  \operatorname*{jump}\nolimits_{\ell
-2+1/2,\ \ \lambda_{2}}\left(  \operatorname*{jump}\nolimits_{\ell
-1+1/2,\ \ \lambda_{1}}\left(  G_{\ell}\right)  \right)  \right)  \right) \\
&  =\left\{  \text{all levels }<\ell-k\right\}  \cup\left\{  \ell-i+\dfrac
{1}{2}+\lambda_{i}\ \mid\ i\in\left\{  1,2,\ldots,k\right\}  \right\}  .
\end{align*}
(Check that these jumps are well-defined -- i.e., that each electron jumps to
an unoccupied state.)

[Example: Let $\ell=3$ and $k=4$ and $\lambda=\left(  \lambda_{1},\lambda
_{2},\lambda_{3},\lambda_{4}\right)  =\left(  4,2,2,1\right)  $. Then,%
\[
E_{\ell,\lambda}=\left\{  \text{all levels }<-1\right\}  \cup\left\{
\dfrac{1}{2},\ \ \dfrac{5}{2},\ \ \dfrac{7}{2},\ \ \dfrac{13}{2}\right\}  .
\]
Here is a picture of this state $E_{\ell,\lambda}$ and how it is constructed
by a sequence of electron jumps (the top state is $G_{\ell}$; the bottom state
is $E_{\ell,\lambda}$):%
\[%
\begin{tikzpicture}[scale=0.6]
\tikzset{myptr/.style={decoration={markings,mark=at position 1 with {\arrow
[scale=3,>=stealth]{>}}},postaction={decorate}}}
\draw[myptr] (-12, 0) -- (15.5, 0);
\filldraw[draw=black, fill=white] (-0.5, 0) circle [radius=0.3];
\filldraw[draw=black, fill=white] (-1.5, 0) circle [radius=0.3];
\filldraw[draw=black, fill=white] (-2.5, 0) circle [radius=0.3];
\filldraw[draw=black, fill=white] (-3.5, 0) circle [radius=0.3];
\filldraw[draw=black, fill=white] (-4.5, 0) circle [radius=0.3];
\filldraw[draw=black, fill=white] (-5.5, 0) circle [radius=0.3];
\filldraw[draw=black, fill=white] (-6.5, 0) circle [radius=0.3];
\filldraw[draw=black, fill=white] (-7.5, 0) circle [radius=0.3];
\filldraw[draw=black, fill=white] (-8.5, 0) circle [radius=0.3];
\filldraw[draw=black, fill=white] (-9.5, 0) circle [radius=0.3];
\filldraw[draw=black, fill=white] (-10.5, 0) circle [radius=0.3];
\filldraw[draw=black, fill=white] (-11.5, 0) circle [radius=0.3];
\node at (0, 1) [label=left:$\cdots\cdots\cdots\text{negative levels}
\cdots\cdots\cdots$] {};
\draw(0, 1.7) -- (0, -8.7);
\node at (0, 2) {$0$};
\draw(3, 1.7) -- (3, -8.7);
\node at (3, 2) {$\ell$};
\filldraw[draw=black, fill=white] (0.5, 0) circle [radius=0.3];
\filldraw[draw=black, fill=white] (1.5, 0) circle [radius=0.3];
\filldraw[draw=black, fill=white] (2.5, 0) circle [radius=0.3];
\filldraw[draw=black, fill=black] (3.5, 0) circle [radius=0.3];
\filldraw[draw=black, fill=black] (4.5, 0) circle [radius=0.3];
\filldraw[draw=black, fill=black] (5.5, 0) circle [radius=0.3];
\filldraw[draw=black, fill=black] (6.5, 0) circle [radius=0.3];
\filldraw[draw=black, fill=black] (7.5, 0) circle [radius=0.3];
\filldraw[draw=black, fill=black] (8.5, 0) circle [radius=0.3];
\filldraw[draw=black, fill=black] (9.5, 0) circle [radius=0.3];
\filldraw[draw=black, fill=black] (10.5, 0) circle [radius=0.3];
\filldraw[draw=black, fill=black] (11.5, 0) circle [radius=0.3];
\filldraw[draw=black, fill=black] (12.5, 0) circle [radius=0.3];
\filldraw[draw=black, fill=black] (13.5, 0) circle [radius=0.3];
\filldraw[draw=black, fill=black] (14.5, 0) circle [radius=0.3];
\node at (0, 1) [label=right:$\cdots\cdots\cdots\text{positive levels}
\cdots\cdots\cdots\cdots$] {};
\draw
[myptr, blue] (2.8, -0.3) -- (6.3, -1.7) node[midway, left, scale=0.7] {$\qquad
\operatorname{jump}_{\ell-1+1/2,\ \ \lambda_1}$};
\draw[myptr] (-12, -2) -- (15.5, -2);
\filldraw[draw=black, fill=white] (-0.5, -2) circle [radius=0.3];
\filldraw[draw=black, fill=white] (-1.5, -2) circle [radius=0.3];
\filldraw[draw=black, fill=white] (-2.5, -2) circle [radius=0.3];
\filldraw[draw=black, fill=white] (-3.5, -2) circle [radius=0.3];
\filldraw[draw=black, fill=white] (-4.5, -2) circle [radius=0.3];
\filldraw[draw=black, fill=white] (-5.5, -2) circle [radius=0.3];
\filldraw[draw=black, fill=white] (-6.5, -2) circle [radius=0.3];
\filldraw[draw=black, fill=white] (-7.5, -2) circle [radius=0.3];
\filldraw[draw=black, fill=white] (-8.5, -2) circle [radius=0.3];
\filldraw[draw=black, fill=white] (-9.5, -2) circle [radius=0.3];
\filldraw[draw=black, fill=white] (-10.5, -2) circle [radius=0.3];
\filldraw[draw=black, fill=white] (-11.5, -2) circle [radius=0.3];
\filldraw[draw=black, fill=white] (0.5, -2) circle [radius=0.3];
\filldraw[draw=black, fill=white] (1.5, -2) circle [radius=0.3];
\filldraw[draw=black, fill=black] (2.5, -2) circle [radius=0.3];
\filldraw[draw=black, fill=black] (3.5, -2) circle [radius=0.3];
\filldraw[draw=black, fill=black] (4.5, -2) circle [radius=0.3];
\filldraw[draw=black, fill=black] (5.5, -2) circle [radius=0.3];
\filldraw[draw=black, fill=white] (6.5, -2) circle [radius=0.3];
\filldraw[draw=black, fill=black] (7.5, -2) circle [radius=0.3];
\filldraw[draw=black, fill=black] (8.5, -2) circle [radius=0.3];
\filldraw[draw=black, fill=black] (9.5, -2) circle [radius=0.3];
\filldraw[draw=black, fill=black] (10.5, -2) circle [radius=0.3];
\filldraw[draw=black, fill=black] (11.5, -2) circle [radius=0.3];
\filldraw[draw=black, fill=black] (12.5, -2) circle [radius=0.3];
\filldraw[draw=black, fill=black] (13.5, -2) circle [radius=0.3];
\filldraw[draw=black, fill=black] (14.5, -2) circle [radius=0.3];
\draw
[myptr, blue] (1.8, -2.3) -- (3.3, -3.7) node[midway, left, scale=0.7] {$\qquad
\operatorname{jump}_{\ell-2+1/2,\ \ \lambda_2}$};
\draw[myptr] (-12, -4) -- (15.5, -4);
\filldraw[draw=black, fill=white] (-0.5, -4) circle [radius=0.3];
\filldraw[draw=black, fill=white] (-1.5, -4) circle [radius=0.3];
\filldraw[draw=black, fill=white] (-2.5, -4) circle [radius=0.3];
\filldraw[draw=black, fill=white] (-3.5, -4) circle [radius=0.3];
\filldraw[draw=black, fill=white] (-4.5, -4) circle [radius=0.3];
\filldraw[draw=black, fill=white] (-5.5, -4) circle [radius=0.3];
\filldraw[draw=black, fill=white] (-6.5, -4) circle [radius=0.3];
\filldraw[draw=black, fill=white] (-7.5, -4) circle [radius=0.3];
\filldraw[draw=black, fill=white] (-8.5, -4) circle [radius=0.3];
\filldraw[draw=black, fill=white] (-9.5, -4) circle [radius=0.3];
\filldraw[draw=black, fill=white] (-10.5, -4) circle [radius=0.3];
\filldraw[draw=black, fill=white] (-11.5, -4) circle [radius=0.3];
\filldraw[draw=black, fill=white] (0.5, -4) circle [radius=0.3];
\filldraw[draw=black, fill=black] (1.5, -4) circle [radius=0.3];
\filldraw[draw=black, fill=black] (2.5, -4) circle [radius=0.3];
\filldraw[draw=black, fill=white] (3.5, -4) circle [radius=0.3];
\filldraw[draw=black, fill=black] (4.5, -4) circle [radius=0.3];
\filldraw[draw=black, fill=black] (5.5, -4) circle [radius=0.3];
\filldraw[draw=black, fill=white] (6.5, -4) circle [radius=0.3];
\filldraw[draw=black, fill=black] (7.5, -4) circle [radius=0.3];
\filldraw[draw=black, fill=black] (8.5, -4) circle [radius=0.3];
\filldraw[draw=black, fill=black] (9.5, -4) circle [radius=0.3];
\filldraw[draw=black, fill=black] (10.5, -4) circle [radius=0.3];
\filldraw[draw=black, fill=black] (11.5, -4) circle [radius=0.3];
\filldraw[draw=black, fill=black] (12.5, -4) circle [radius=0.3];
\filldraw[draw=black, fill=black] (13.5, -4) circle [radius=0.3];
\filldraw[draw=black, fill=black] (14.5, -4) circle [radius=0.3];
\draw
[myptr, blue] (0.8, -4.3) -- (2.3, -5.7) node[midway, left, scale=0.7] {$\qquad
\operatorname{jump}_{\ell-3+1/2,\ \ \lambda_3}$};
\draw[myptr] (-12, -6) -- (15.5, -6);
\filldraw[draw=black, fill=white] (-0.5, -6) circle [radius=0.3];
\filldraw[draw=black, fill=white] (-1.5, -6) circle [radius=0.3];
\filldraw[draw=black, fill=white] (-2.5, -6) circle [radius=0.3];
\filldraw[draw=black, fill=white] (-3.5, -6) circle [radius=0.3];
\filldraw[draw=black, fill=white] (-4.5, -6) circle [radius=0.3];
\filldraw[draw=black, fill=white] (-5.5, -6) circle [radius=0.3];
\filldraw[draw=black, fill=white] (-6.5, -6) circle [radius=0.3];
\filldraw[draw=black, fill=white] (-7.5, -6) circle [radius=0.3];
\filldraw[draw=black, fill=white] (-8.5, -6) circle [radius=0.3];
\filldraw[draw=black, fill=white] (-9.5, -6) circle [radius=0.3];
\filldraw[draw=black, fill=white] (-10.5, -6) circle [radius=0.3];
\filldraw[draw=black, fill=white] (-11.5, -6) circle [radius=0.3];
\filldraw[draw=black, fill=black] (0.5, -6) circle [radius=0.3];
\filldraw[draw=black, fill=black] (1.5, -6) circle [radius=0.3];
\filldraw[draw=black, fill=white] (2.5, -6) circle [radius=0.3];
\filldraw[draw=black, fill=white] (3.5, -6) circle [radius=0.3];
\filldraw[draw=black, fill=black] (4.5, -6) circle [radius=0.3];
\filldraw[draw=black, fill=black] (5.5, -6) circle [radius=0.3];
\filldraw[draw=black, fill=white] (6.5, -6) circle [radius=0.3];
\filldraw[draw=black, fill=black] (7.5, -6) circle [radius=0.3];
\filldraw[draw=black, fill=black] (8.5, -6) circle [radius=0.3];
\filldraw[draw=black, fill=black] (9.5, -6) circle [radius=0.3];
\filldraw[draw=black, fill=black] (10.5, -6) circle [radius=0.3];
\filldraw[draw=black, fill=black] (11.5, -6) circle [radius=0.3];
\filldraw[draw=black, fill=black] (12.5, -6) circle [radius=0.3];
\filldraw[draw=black, fill=black] (13.5, -6) circle [radius=0.3];
\filldraw[draw=black, fill=black] (14.5, -6) circle [radius=0.3];
\draw
[myptr, blue] (-0.2, -6.3) -- (0.3, -7.7) node[midway, left, scale=0.7] {$\qquad
\operatorname{jump}_{\ell-4+1/2,\ \ \lambda_4}$};
\draw[myptr] (-12, -8) -- (15.5, -8);
\filldraw[draw=black, fill=black] (-0.5, -8) circle [radius=0.3];
\filldraw[draw=black, fill=white] (-1.5, -8) circle [radius=0.3];
\filldraw[draw=black, fill=white] (-2.5, -8) circle [radius=0.3];
\filldraw[draw=black, fill=white] (-3.5, -8) circle [radius=0.3];
\filldraw[draw=black, fill=white] (-4.5, -8) circle [radius=0.3];
\filldraw[draw=black, fill=white] (-5.5, -8) circle [radius=0.3];
\filldraw[draw=black, fill=white] (-6.5, -8) circle [radius=0.3];
\filldraw[draw=black, fill=white] (-7.5, -8) circle [radius=0.3];
\filldraw[draw=black, fill=white] (-8.5, -8) circle [radius=0.3];
\filldraw[draw=black, fill=white] (-9.5, -8) circle [radius=0.3];
\filldraw[draw=black, fill=white] (-10.5, -8) circle [radius=0.3];
\filldraw[draw=black, fill=white] (-11.5, -8) circle [radius=0.3];
\filldraw[draw=black, fill=white] (0.5, -8) circle [radius=0.3];
\filldraw[draw=black, fill=black] (1.5, -8) circle [radius=0.3];
\filldraw[draw=black, fill=white] (2.5, -8) circle [radius=0.3];
\filldraw[draw=black, fill=white] (3.5, -8) circle [radius=0.3];
\filldraw[draw=black, fill=black] (4.5, -8) circle [radius=0.3];
\filldraw[draw=black, fill=black] (5.5, -8) circle [radius=0.3];
\filldraw[draw=black, fill=white] (6.5, -8) circle [radius=0.3];
\filldraw[draw=black, fill=black] (7.5, -8) circle [radius=0.3];
\filldraw[draw=black, fill=black] (8.5, -8) circle [radius=0.3];
\filldraw[draw=black, fill=black] (9.5, -8) circle [radius=0.3];
\filldraw[draw=black, fill=black] (10.5, -8) circle [radius=0.3];
\filldraw[draw=black, fill=black] (11.5, -8) circle [radius=0.3];
\filldraw[draw=black, fill=black] (12.5, -8) circle [radius=0.3];
\filldraw[draw=black, fill=black] (13.5, -8) circle [radius=0.3];
\filldraw[draw=black, fill=black] (14.5, -8) circle [radius=0.3];
\end{tikzpicture}%
.
\]
Note how the order of the electrons after the jumps is the same as before --
i.e., the electron in the rightmost position before the jumps is still the
rightmost one after the jumps, etc.]

Recall that the $\ell$-ground state $G_{\ell}$ has energy $\ell^{2}$ and
particle number $\ell$. The state $E_{\ell,\lambda}$ that we have just defined
is obtained from this $\ell$-ground state $G_{\ell}$ by $k$ jumps, which are
jumps by $\lambda_{1},\lambda_{2},\ldots,\lambda_{k}$ steps respectively.
Recall that a jump by $q$ steps raises the energy by $2q$ but keeps the
particle number unchanged. Thus, the \textquotedblleft excited
state\textquotedblright\ $E_{\ell,\lambda}$ has energy $\ell^{2}%
+\underbrace{2\lambda_{1}+2\lambda_{2}+\cdots+2\lambda_{k}}_{=2\left\vert
\lambda\right\vert }=\ell^{2}+2\left\vert \lambda\right\vert $ and particle
number $\ell$. Furthermore, every state with particle number $\ell$ can be
written as $E_{\ell,\lambda}$ for a unique partition $\lambda$ (check
this!\footnote{Here is how to obtain this $\lambda$: For each $i\in\left\{
1,2,3,\ldots\right\}  $, let $u_{i}$ be the $i$-th largest level in the state,
and let $\lambda_{i}=u_{i}-\ell+i-\dfrac{1}{2}$. Then, $\left(  \lambda
_{1},\lambda_{2},\lambda_{3},\ldots\right)  $ is a weakly decreasing sequence
of nonnegative integers, and all but finitely many of its entries are $0$
(since the state has particle number $\ell$, so that it is not hard to see
that $u_{i}=\ell-i+\dfrac{1}{2}$ for any sufficiently large $i$). Removing all
$0$s from this sequence $\left(  \lambda_{1},\lambda_{2},\lambda_{3}%
,\ldots\right)  $ thus results in a finite tuple $\left(  \lambda_{1}%
,\lambda_{2},\ldots,\lambda_{k}\right)  $, which is precisely the partition
$\lambda$ whose corresponding $E_{\ell,\lambda}$ is our state.}).

Thus, we obtain a bijection%
\begin{align*}
\Phi_{\ell}:\left\{  \text{partitions}\right\}   &  \rightarrow\left\{
\text{states with particle number }\ell\right\}  ,\\
\lambda &  \mapsto E_{\ell,\lambda}.
\end{align*}
This bijection satisfies%
\[
\operatorname*{energy}\left(  \Phi_{\ell}\left(  \lambda\right)  \right)
=\operatorname*{energy}E_{\ell,\lambda}=\ell^{2}+2\left\vert \lambda
\right\vert \ \ \ \ \ \ \ \ \ \ \text{for every partition }\lambda.
\]
Hence,%
\begin{align}
\sum_{\substack{\lambda\text{ is a}\\\text{partition}}}\underbrace{q^{\ell
^{2}+2\left\vert \lambda\right\vert }}_{\substack{=q^{\operatorname*{energy}%
\left(  \Phi_{\ell}\left(  \lambda\right)  \right)  }\\\text{(since
}\operatorname*{energy}\left(  \Phi_{\ell}\left(  \lambda\right)  \right)
=\ell^{2}+2\left\vert \lambda\right\vert \text{)}}}  &  =\sum
_{\substack{\lambda\text{ is a}\\\text{partition}}}q^{\operatorname*{energy}%
\left(  \Phi_{\ell}\left(  \lambda\right)  \right)  }\nonumber\\
&  =\sum_{\substack{S\text{ is a state with}\\\text{particle number }\ell
}}q^{\operatorname*{energy}S} \label{pf.thm.pars.jtp2.rhs5}%
\end{align}
(here, we have substituted $S$ for $\Phi_{\ell}\left(  \lambda\right)  $ in
the sum, since the map $\Phi_{\ell}$ is a bijection).

Forget that we fixed $\ell$. We thus have proved (\ref{pf.thm.pars.jtp2.rhs5})
for each $\ell\in\mathbb{Z}$.

Now, (\ref{pf.thm.pars.jtp2.rhs1}) becomes%
\begin{align*}
\left(  \sum_{\ell\in\mathbb{Z}}q^{\ell^{2}}z^{\ell}\right)  \prod
_{n>0}\left(  1-q^{2n}\right)  ^{-1}  &  =\sum_{\ell\in\mathbb{Z}%
}\ \ \underbrace{\sum_{\substack{\lambda\text{ is a}\\\text{partition}%
}}q^{\ell^{2}+2\left\vert \lambda\right\vert }}_{\substack{=\sum
_{\substack{S\text{ is a state with}\\\text{particle number }\ell
}}q^{\operatorname*{energy}S}\\\text{(by (\ref{pf.thm.pars.jtp2.rhs5}))}%
}}z^{\ell}\\
&  =\sum_{\ell\in\mathbb{Z}}\ \ \sum_{\substack{S\text{ is a state
with}\\\text{particle number }\ell}}q^{\operatorname*{energy}S}%
\underbrace{z^{\ell}}_{\substack{=z^{\operatorname*{parnum}S}\\\text{(since
}\ell=\operatorname*{parnum}S\text{)}}}\\
&  =\underbrace{\sum_{\ell\in\mathbb{Z}}\ \ \sum_{\substack{S\text{ is a state
with}\\\text{particle number }\ell}}}_{=\sum_{S\text{ is a state}}%
}q^{\operatorname*{energy}S}z^{\operatorname*{parnum}S}\\
&  =\sum_{S\text{ is a state}}q^{\operatorname*{energy}S}%
z^{\operatorname*{parnum}S}.
\end{align*}

Comparing this with (\ref{pf.thm.pars.jtp2.lhs}), we obtain%
\[
\prod_{n>0}\left(  \left(  1+q^{2n-1}z\right)  \left(  1+q^{2n-1}%
z^{-1}\right)  \right)  =\left(  \sum_{\ell\in\mathbb{Z}}q^{\ell^{2}}z^{\ell
}\right)  \prod_{n>0}\left(  1-q^{2n}\right)  ^{-1}.
\]
Multiplying both sides of this identity with $\prod_{n>0}\left(
1-q^{2n}\right)  $, we find%
\[
\prod_{n>0}\left(  \left(  1+q^{2n-1}z\right)  \left(  1+q^{2n-1}%
z^{-1}\right)  \left(  1-q^{2n}\right)  \right)  =\sum_{\ell\in\mathbb{Z}%
}q^{\ell^{2}}z^{\ell}.
\]
This proves Jacobi's Triple Product Identity (Theorem \ref{thm.pars.jtp1} and
Theorem \ref{thm.pars.jtp2}).
\end{proof}

Other proofs of Jacobi's Triple Product Identity can be found in
\cite[\S 3.4]{Aigner07}, \cite[Theorem 10.2]{Wagner08}, \cite[Theorem
8.2b]{Henric77}, \cite[Theorem 5.10]{Sambal22} and \cite[\S 1.3 and
\S 1.4]{Hirsch17}.

We note that proving Theorem \ref{thm.pars.pent} using Theorem
\ref{thm.pars.jtp2} can be viewed as a kind of overkill; there are more direct
proofs of Theorem \ref{thm.pars.pent} as well (see, e.g., \cite{Zabroc03},
\cite[\S 10]{Koch16}, \cite[Exercise 5]{18f-mt3s}, \cite[\S 3]{Bell06},
\cite[Theorem VI.4]{Hulpke25} or \cite[Theorem 2.28 and \S 3.9 Solution
15]{Bona25}).

\subsubsection{Application: A recursion for the sum of divisors}

Now that Theorem \ref{thm.pars.pent} is proven, we use the occasion to
spotlight a truly unexpected application of it: Euler's recursion for the sum
of divisors. We state it in the following form:

\begin{theorem}
\label{thm.pars.euler-sum-div-rec}For any positive integer $n$, let
$\sigma\left(  n\right)  $ denote the sum of all positive divisors of $n$.
(For example, $\sigma\left(  6\right)  =1+2+3+6=12$ and $\sigma\left(
7\right)  =1+7=8$.)

Then, for each positive integer $n$, we have%
\[
\sum_{\substack{k\in\mathbb{Z};\\w_{k}<n}}\left(  -1\right)  ^{k}\sigma\left(
n-w_{k}\right)  =%
\begin{cases}
\left(  -1\right)  ^{k-1}n, & \text{if }n=w_{k}\text{ for some }k\in
\mathbb{Z};\\
0, & \text{if not.}%
\end{cases}
\]
(Here, $w_{k}$ is the $k$-th pentagonal number, defined in Definition
\ref{def.pars.pent-num}.)
\end{theorem}

The left hand side of the equality in Theorem \ref{thm.pars.euler-sum-div-rec}
looks as follows:%
\[
\sigma\left(  n\right)  -\sigma\left(  n-1\right)  -\sigma\left(  n-2\right)
+\sigma\left(  n-5\right)  +\sigma\left(  n-7\right)  -\sigma\left(
n-12\right)  -\sigma\left(  n-15\right)  \pm\cdots
\]
(the sum goes only as far as the arguments remain positive integers, so it is
a finite sum). Thus, by solving this equality for $\sigma\left(  n\right)  $,
we obtain a recursive formula for $\sigma\left(  n\right)  $. (In this form,
Theorem \ref{thm.pars.euler-sum-div-rec} appears in many places, such as
\cite[\S 23]{Ness61} and \cite[Theorem 37]{Johnso20}.) This does not give an
efficient algorithm for computing $\sigma\left(  n\right)  $ (it is much
slower than the formula in \cite[Exercise 2.18.1 \textbf{(b)}]{19s}, even
though the latter requires computing the prime factorization of $n$), but its
beauty and element of surprise (why would you expect the pentagonal numbers to
have anything to do with sums of divisors?) make up for this practical
uselessness. Even more surprisingly, it can be proved using Euler's pentagonal
number theorem, even though it has seemingly nothing to do with partitions!

\begin{proof}
[Proof of Theorem \ref{thm.pars.euler-sum-div-rec} (sketched).]First, let us
show that the right hand side in Theorem \ref{thm.pars.euler-sum-div-rec} is
well-defined. Indeed, it is easy to see that
\[
w_{0}<w_{1}<w_{-1}<w_{2}<w_{-2}<w_{3}<w_{-3}<\cdots;
\]
thus, the pentagonal numbers $w_{k}$ for all $k\in\mathbb{Z}$ are distinct.
Hence, if a positive integer $n$ satisfies $n=w_{k}$ for some $k\in\mathbb{Z}%
$, then this latter $k$ is uniquely determined by $n$, and therefore the
expression%
\[%
\begin{cases}
\left(  -1\right)  ^{k-1}n, & \text{if }n=w_{k}\text{ for some }k\in
\mathbb{Z};\\
0, & \text{if not}%
\end{cases}
\]
is well-defined. In other words, the right hand side in Theorem
\ref{thm.pars.euler-sum-div-rec} is well-defined.

Now, define the two FPSs $P$ and $S$ as in our above proof of Theorem
\ref{thm.pars.sigma1}. In that very proof, we have shown that these FPSs
satisfy
\begin{equation}
xP^{\prime}=SP. \label{pf.thm.pars.euler-sum-div-rec.xP'}%
\end{equation}

Next, we define a further FPS%
\[
Q:=\sum_{k\in\mathbb{Z}}\left(  -1\right)  ^{k}x^{w_{k}}\in\mathbb{Z}\left[
\left[  x\right]  \right]  .
\]

Multiplying the equalities%
\[
P=\sum_{n\in\mathbb{N}}p\left(  n\right)  x^{n}=\prod_{k=1}^{\infty}\dfrac
{1}{1-x^{k}}\ \ \ \ \ \ \ \ \ \ \left(  \text{by Theorem
\ref{thm.pars.main-gf}}\right)
\]
and%
\[
Q=\sum_{k\in\mathbb{Z}}\left(  -1\right)  ^{k}x^{w_{k}}=\prod\limits_{k=1}%
^{\infty}\left(  1-x^{k}\right)  \ \ \ \ \ \ \ \ \ \ \left(  \text{by Theorem
\ref{thm.pars.pent}}\right)  ,
\]
we obtain%
\[
PQ=\left(  \prod_{k=1}^{\infty}\dfrac{1}{1-x^{k}}\right)  \left(
\prod\limits_{k=1}^{\infty}\left(  1-x^{k}\right)  \right)  =\prod
\limits_{k=1}^{\infty}\underbrace{\left(  \dfrac{1}{1-x^{k}}\cdot\left(
1-x^{k}\right)  \right)  }_{=1}=1.
\]

Combining this equality with (\ref{pf.thm.pars.euler-sum-div-rec.xP'}), we can
easily see that $xQ^{\prime}=-QS$. Indeed, here are two (essentially
equivalent) ways to show this:

\begin{itemize}
\item From $PQ=1$, we obtain $Q=P^{-1}$. But $P\in\mathbb{Z}\left[  \left[
x\right]  \right]  _{1}$ (since the constant coefficient of $P$ is $p\left(
0\right)  =1$). Hence, Corollary \ref{cor.fps.loder.inv} (applied to $f=P$)
yields $\operatorname*{loder}\left(  P^{-1}\right)  =-\operatorname*{loder}P$.
In view of $P^{-1}=Q$, this rewrites as $\operatorname*{loder}%
Q=-\operatorname*{loder}P$. But the definition of a logarithmic derivative
yields $\operatorname*{loder}Q=\dfrac{Q^{\prime}}{Q}$ and
$\operatorname*{loder}P=\dfrac{P^{\prime}}{P}$. Hence,
\[
\dfrac{Q^{\prime}}{Q}=\operatorname*{loder}Q=-\operatorname*{loder}%
P=-\dfrac{P^{\prime}}{P}\ \ \ \ \ \ \ \ \ \ \left(  \text{since }%
\operatorname*{loder}P=\dfrac{P^{\prime}}{P}\right)  .
\]
Multiplying this equality by $x$, we find%
\begin{align*}
\dfrac{xQ^{\prime}}{Q}  &  =-\dfrac{xP^{\prime}}{P}=-\dfrac{SP}{P}%
\ \ \ \ \ \ \ \ \ \ \left(  \text{by (\ref{pf.thm.pars.euler-sum-div-rec.xP'}%
)}\right) \\
&  =-S,
\end{align*}
so that $xQ^{\prime}=-QS$.

\item A less conceptual but slicker way to draw the same conclusion is the
following: Taking derivatives in the equality $PQ=1$, we find $\left(
PQ\right)  ^{\prime}=1^{\prime}=0$, so that $0=\left(  PQ\right)  ^{\prime
}=P^{\prime}Q+PQ^{\prime}$ (by Theorem \ref{thm.fps.deriv.rules}
\textbf{(d)}). Thus, $PQ^{\prime}=-P^{\prime}Q$. Multiplying this equality by
$x$, we find $PQ^{\prime}\cdot x=-P^{\prime}Q\cdot x=-\underbrace{xP^{\prime}%
}_{\substack{=SP\\\text{(by (\ref{pf.thm.pars.euler-sum-div-rec.xP'}))}%
}}Q=-S\underbrace{PQ}_{=1}=-S$. Multiplying this further by $Q$, we obtain
$PQ^{\prime}\cdot xQ=-SQ=-QS$. In view of $PQ^{\prime}\cdot xQ=xQ^{\prime
}\cdot\underbrace{PQ}_{=1}=xQ^{\prime}$, we can rewrite this as $xQ^{\prime
}=-QS$.
\end{itemize}

Either way, we now know that $xQ^{\prime}=-QS$. However, from $Q=\sum
_{k\in\mathbb{Z}}\left(  -1\right)  ^{k}x^{w_{k}}$, we obtain%
\[
Q^{\prime}=\left(  \sum_{k\in\mathbb{Z}}\left(  -1\right)  ^{k}x^{w_{k}%
}\right)  ^{\prime}=\sum_{k\in\mathbb{Z}}\left(  -1\right)  ^{k}%
\underbrace{\left(  x^{w_{k}}\right)  ^{\prime}}_{\substack{=w_{k}x^{w_{k}%
-1}\\\text{(where we understand this}\\\text{expression to mean }0\text{ if
}w_{k}=0\text{)}}}=\sum_{k\in\mathbb{Z}}\left(  -1\right)  ^{k}w_{k}%
x^{w_{k}-1}.
\]
Upon multiplication by $x$, this becomes\footnote{We are using the fact that
$\left(  x^{i}\right)  ^{\prime}= ix^{i-1}$ for each $i \in\mathbb{N}$. This
holds even for $i=0$, using the convention that $0x^{-1}$ stands for $0$ (or
simply treating $x^{-1}$ as a Laurent series).}
\[
xQ^{\prime}=x\sum_{k\in\mathbb{Z}}\left(  -1\right)  ^{k}w_{k}x^{w_{k}-1}%
=\sum_{k\in\mathbb{Z}}\left(  -1\right)  ^{k}w_{k}\underbrace{xx^{w_{k}-1}%
}_{=x^{w_{k}}}=\sum_{k\in\mathbb{Z}}\left(  -1\right)  ^{k}w_{k}x^{w_{k}}.
\]
Hence,%
\begin{align}
\sum_{k\in\mathbb{Z}}\left(  -1\right)  ^{k}w_{k}x^{w_{k}}  &  =xQ^{\prime
}=-\underbrace{Q}_{=\sum_{k\in\mathbb{Z}}\left(  -1\right)  ^{k}x^{w_{k}}%
}\ \ \underbrace{S}_{\substack{=\sum_{k>0}\sigma\left(  k\right)  x^{k}%
\\=\sum_{i>0}\sigma\left(  i\right)  x^{i}}}\nonumber\\
&  =-\left(  \sum_{k\in\mathbb{Z}}\left(  -1\right)  ^{k}x^{w_{k}}\right)
\left(  \sum_{i>0}\sigma\left(  i\right)  x^{i}\right) \nonumber\\
&  =\sum_{k\in\mathbb{Z}}\ \ \sum_{i>0}\underbrace{\left(  -\left(  -1\right)
^{k}\right)  }_{=\left(  -1\right)  ^{k-1}}\underbrace{x^{w_{k}}\sigma\left(
i\right)  x^{i}}_{=\sigma\left(  i\right)  x^{w_{k}+i}}\nonumber\\
&  =\sum_{k\in\mathbb{Z}}\ \ \sum_{i>0}\left(  -1\right)  ^{k-1}\sigma\left(
i\right)  x^{w_{k}+i}\nonumber\\
&  =\underbrace{\sum_{k\in\mathbb{Z}}\ \ \sum_{m>w_{k}}}_{=\sum_{m>0}%
\ \ \sum_{\substack{k\in\mathbb{Z};\\m>w_{k}}}}\left(  -1\right)  ^{k-1}%
\sigma\left(  m-w_{k}\right)  \underbrace{x^{w_{k}+\left(  m-w_{k}\right)  }%
}_{=x^{m}}\nonumber\\
&  \ \ \ \ \ \ \ \ \ \ \ \ \ \ \ \ \ \ \ \ \left(
\begin{array}
[c]{c}%
\text{here, we have substituted }m-w_{k}\text{ for }i\\
\text{in the inner sum}%
\end{array}
\right) \nonumber\\
&  =\sum_{m>0}\ \ \sum_{\substack{k\in\mathbb{Z};\\m>w_{k}}}\left(  -1\right)
^{k-1}\sigma\left(  m-w_{k}\right)  x^{m}.
\label{pf.thm.pars.euler-sum-div-rec.4}%
\end{align}

Now, let $n$ be a positive integer. Let us compare the coefficients of $x^{n}$
on the left and right hand sides of (\ref{pf.thm.pars.euler-sum-div-rec.4}).
On the left hand side, the coefficient of $x^{n}$ is%
\[%
\begin{cases}
\left(  -1\right)  ^{k}w_{k}, & \text{if }n=w_{k}\text{ for some }%
k\in\mathbb{Z};\\
0, & \text{if not}%
\end{cases}
\]
(since the pentagonal numbers $w_{k}$ for all $k\in\mathbb{Z}$ are distinct,
and thus the different addends on the left hand side of
(\ref{pf.thm.pars.euler-sum-div-rec.4}) contribute to different monomials). On
the right hand side, the coefficient of $x^{n}$ is obviously%
\[
\sum_{\substack{k\in\mathbb{Z};\\n>w_{k}}}\left(  -1\right)  ^{k-1}%
\sigma\left(  n-w_{k}\right)  .
\]
Since these two coefficients are equal (because
(\ref{pf.thm.pars.euler-sum-div-rec.4}) is an identity), we thus conclude that%
\begin{align*}%
\begin{cases}
\left(  -1\right)  ^{k}w_{k}, & \text{if }n=w_{k}\text{ for some }%
k\in\mathbb{Z};\\
0, & \text{if not}%
\end{cases}
&  =\underbrace{\sum_{\substack{k\in\mathbb{Z};\\n>w_{k}}}}_{=\sum
_{\substack{k\in\mathbb{Z};\\w_{k}<n}}}\underbrace{\left(  -1\right)  ^{k-1}%
}_{=-\left(  -1\right)  ^{k}}\sigma\left(  n-w_{k}\right) \\
&  =-\sum_{\substack{k\in\mathbb{Z};\\w_{k}<n}}\left(  -1\right)  ^{k}%
\sigma\left(  n-w_{k}\right)  .
\end{align*}
Thus,%
\begin{align*}
\sum_{\substack{k\in\mathbb{Z};\\w_{k}<n}}\left(  -1\right)  ^{k}\sigma\left(
n-w_{k}\right)   &  =-%
\begin{cases}
\left(  -1\right)  ^{k}w_{k}, & \text{if }n=w_{k}\text{ for some }%
k\in\mathbb{Z};\\
0, & \text{if not}%
\end{cases}
\\
&  =%
\begin{cases}
-\left(  -1\right)  ^{k}w_{k}, & \text{if }n=w_{k}\text{ for some }%
k\in\mathbb{Z};\\
0, & \text{if not}%
\end{cases}
\\
&  =%
\begin{cases}
\left(  -1\right)  ^{k-1}w_{k}, & \text{if }n=w_{k}\text{ for some }%
k\in\mathbb{Z};\\
0, & \text{if not}%
\end{cases}
\\
&  \ \ \ \ \ \ \ \ \ \ \ \ \ \ \ \ \ \ \ \ \left(  \text{since }-\left(
-1\right)  ^{k}=\left(  -1\right)  ^{k-1}\right) \\
&  =%
\begin{cases}
\left(  -1\right)  ^{k-1}n, & \text{if }n=w_{k}\text{ for some }k\in
\mathbb{Z};\\
0, & \text{if not}%
\end{cases}
\end{align*}
(because if $n=w_{k}$ for some $k\in\mathbb{Z}$, then $w_{k}=n$). This proves
Theorem \ref{thm.pars.euler-sum-div-rec}.
\end{proof}

\subsection{$q$-binomial coefficients}

Next, we shall study $q$\emph{-binomial coefficients} (also known as
\emph{Gaussian binomial coefficients}, due to their origins in Gauss's
number-theoretical research \cite[\S 5]{Gauss08}). While we will define them
as generating functions for certain kinds of partitions, they are sufficiently
elementary to have relevance to various other subjects. We will scratch the
surface; more can be found in \cite[Chapters 5--7]{KacChe02}, \cite[Chapter
2]{KliSch97}, \cite[Chapter 5]{Wagner08}, \cite[\S 11.3--11.11]{Wagner20},
\cite[spread across the text]{Stanley-EC1}, \cite[\S 2.6]{GouJac83},
\cite{Johnso20} and other sources. (The book \cite{Johnso20} is particularly
recommended as a leisurely introduction to $q$-binomial coefficients and
related power series.)

\subsubsection{\label{subsec.qbin.motiv}Motivation}

For any $n\in\mathbb{N}$, we have%
\[
p\left(  n\right)  =\left(  \text{\# of partitions of }n\right)  .
\]
For any $n,k\in\mathbb{N}$, we have%
\begin{align*}
p_{k}\left(  n\right)   &  =\left(  \text{\# of partitions of }n\text{ into
}k\text{ parts}\right) \\
&  =\left(  \text{\# of partitions of }n\text{ with largest part }k\right)
\ \ \ \ \ \ \ \ \ \ \left(  \text{by Proposition \ref{prop.pars.pkn=dual}%
}\right)  .
\end{align*}
Thus, for any $n,k\in\mathbb{N}$, we have%
\begin{align*}
p_{0}\left(  n\right)  +p_{1}\left(  n\right)  +\cdots+p_{k}\left(  n\right)
&  =\left(  \text{\# of partitions of }n\text{ into \textbf{at most} }k\text{
parts}\right) \\
&  =\left(  \text{\# of partitions of }n\text{ with largest part }\leq
k\right)
\end{align*}
(by Corollary \ref{cor.pars.p0kn=dual}).

So far, so good. But how to count partitions of $n$ that both have a fixed \#
of parts (say, $k$ parts) and a fixed largest part (say, largest part $\ell$) ?

Let us first drop the size requirement -- i.e., we replace \textquotedblleft
partitions of $n$\textquotedblright\ by just \textquotedblleft
partitions\textquotedblright.

For example, how many partitions have $4$ parts and largest part $6$ ?

As in the proof of Proposition \ref{prop.pars.pkn=dual}, let us draw the Young
diagram of such a partition: For example, the partition $\left(
6,3,3,2\right)  $ has Young diagram%
\[
\ydiagram{6,3,3,2}\ \ .
\]

Consider the lower boundary of this Young diagram -- i.e., the
\textquotedblleft irregular\textquotedblright\ southeastern border between
what is in the diagram and what is outside of it. Let me mark it in thick red:%
\[%
\begin{tikzpicture}
\draw(0, 0) rectangle (1, 1);
\draw(0, 1) rectangle (1, 2);
\draw(0, 2) rectangle (1, 3);
\draw(0, 3) rectangle (1, 4);
\draw(1, 0) rectangle (2, 1);
\draw(1, 1) rectangle (2, 2);
\draw(1, 2) rectangle (2, 3);
\draw(1, 3) rectangle (2, 4);
\draw(2, 1) rectangle (3, 2);
\draw(2, 2) rectangle (3, 3);
\draw(2, 3) rectangle (3, 4);
\draw(3, 3) rectangle (4, 4);
\draw(4, 3) rectangle (5, 4);
\draw(5, 3) rectangle (6, 4);
\draw
[ultra thick, red] (0, 0) -- (2, 0) -- (2, 1) -- (3, 1) -- (3, 3) -- (6, 3) -- (6, 4);
\end{tikzpicture}%
\ \ .
\]
This lower boundary can be viewed as a lattice path from the point $\left(
0,0\right)  $ to the point $\left(  6,4\right)  $ (where we are using
Cartesian coordinates to label the intersections of grid lines, so that the
southwesternmost point in our diagram is $\left(  0,0\right)  $; note that
this is completely unrelated to our labeling of cells used in defining the
Young diagram!\footnote{For additional clarity, here are the Cartesian
coordinates of all grid points on our lattice path:%
\[%
\begin{tikzpicture}
\draw(0, 0) rectangle (1, 1);
\draw(0, 1) rectangle (1, 2);
\draw(0, 2) rectangle (1, 3);
\draw(0, 3) rectangle (1, 4);
\draw(1, 0) rectangle (2, 1);
\draw(1, 1) rectangle (2, 2);
\draw(1, 2) rectangle (2, 3);
\draw(1, 3) rectangle (2, 4);
\draw(2, 1) rectangle (3, 2);
\draw(2, 2) rectangle (3, 3);
\draw(2, 3) rectangle (3, 4);
\draw(3, 3) rectangle (4, 4);
\draw(4, 3) rectangle (5, 4);
\draw(5, 3) rectangle (6, 4);
\draw
[ultra thick, red] (0, 0) -- (2, 0) -- (2, 1) -- (3, 1) -- (3, 3) -- (6, 3) -- (6, 4);
\fill(0, 0) circle (2pt) node[below] {\scriptsize$\left(0, 0\right)$};
\fill(1, 0) circle (2pt) node[below] {\scriptsize$\left(1, 0\right)$};
\fill(2, 0) circle (2pt) node[below] {\scriptsize$\left(2, 0\right)$};
\fill(2, 1) circle (2pt) node[below right] {\scriptsize$\left(2, 1\right)$};
\fill(3, 1) circle (2pt) node[right] {\scriptsize$\left(3, 1\right)$};
\fill(3, 2) circle (2pt) node[right] {\scriptsize$\left(3, 2\right)$};
\fill(3, 3) circle (2pt) node[below left] {\scriptsize$\left(3, 3\right)$};
\fill(4, 3) circle (2pt) node[below] {\scriptsize$\left(4, 3\right)$};
\fill(5, 3) circle (2pt) node[below] {\scriptsize$\left(5, 3\right)$};
\fill(6, 3) circle (2pt) node[right] {\scriptsize$\left(6, 3\right)$};
\fill(6, 4) circle (2pt) node[right] {\scriptsize$\left(6, 4\right)$};
\end{tikzpicture}%
\ \ .
\]
}). This lattice path consists of east-steps (i.e., steps $\left(  i,j\right)
\rightarrow\left(  i+1,j\right)  $) and north-steps (i.e., steps $\left(
i,j\right)  \rightarrow\left(  i,j+1\right)  $); moreover, it begins with an
east-step (since otherwise, our partition would have fewer than $4$ parts) and
ends with a north-step (since otherwise, our partition would have largest part
$<6$). Moreover, the Young diagram (and thus the partition) is uniquely
determined by this lattice path, since its cells are precisely the cells
\textquotedblleft northwest\textquotedblright\ of this lattice path.
Conversely, any lattice path from $\left(  0,0\right)  $ to $\left(
6,4\right)  $ that consists of east-steps and north-steps and begins with an
east step and ends with a north-step uniquely determines a Young diagram and
therefore a partition. Therefore, in order to count the partitions that have
$4$ parts and largest part $6$, we only need to count such lattice paths.

To count them, we notice that any such lattice path has precisely $10$ steps
(since any step increases the sum of the coordinates by $1$; but this sum must
increase from $0+0=0$ to $6+4=10$). The first and the last steps are
predetermined; it remains to decide which of the remaining $8$ steps are
north-steps. The \# of ways to do this is $\dbinom{8}{3}$, because we want
precisely $3$ of our $8$ non-predetermined steps to be north-steps (in order
to end up at $\left(  6,4\right)  $ rather than some other point).

As a consequence of this all, we find%
\[
\left(  \text{\# of partitions with }4\text{ parts and largest part }6\right)
=\dbinom{8}{3}.
\]
More generally, by the same argument, we obtain the following:

\begin{proposition}
\label{prop.pars.qbinom.intro-count-binom}For any positive integers $k$ and
$\ell$, we have%
\[
\left(  \text{\# of partitions with }k\text{ parts and largest part }%
\ell\right)  =\dbinom{k+\ell-2}{k-1}.
\]

\end{proposition}

Note two things:

\begin{itemize}
\item This is a finite number, even without fixing the size of the partition.
This is not surprising, since you have only finitely many parts and only
finitely many options for each part.

\item The number is symmetric in $k$ and $\ell$. This, too, is not surprising,
because conjugation (as defined in the proof of Proposition
\ref{prop.pars.pkn=dual}) gives a bijection%
\begin{align*}
&  \text{from }\left\{  \text{partitions with }k\text{ parts and largest part
}\ell\right\} \\
&  \text{to }\left\{  \text{partitions with }\ell\text{ parts and largest part
}k\right\}  .
\end{align*}

\end{itemize}

Now, let us integrate the size of the partition back into our count -- i.e.,
let us try to count the partitions of a given $n\in\mathbb{N}$ with $k$ parts
and largest part $\ell$. No simple formula (like Proposition
\ref{prop.pars.qbinom.intro-count-binom}) exists for this number any more, so
we switch our focus to the generating function of such numbers (for fixed $k$
and $\ell$). In other words, we try to compute the FPS%
\begin{align*}
&  \sum_{n\in\mathbb{N}}\left(  \text{\# of partitions of }n\text{ with
}k\text{ parts and largest part }\ell\right)  x^{n}\\
&  =\sum_{\substack{\lambda\text{ is a partition}\\\text{with largest part
}\ell\\\text{and length }k}}x^{\left\vert \lambda\right\vert }.
\end{align*}

For reasons of convenience, history and simplicity, we modify this problem
slightly (without changing its essence). To wit,

\begin{itemize}
\item we rename $\ell$ as $n-k$ (note that $n$ will no longer stand for the
size of the partition);

\item we replace \textquotedblleft largest part $n-k$ and length
$k$\textquotedblright\ by \textquotedblleft largest part $\leq n-k$ and length
$\leq k$\textquotedblright\ (this changes the results of our counts, but we
can easily recover the answer to the original question from an answer to the
new question; e.g., in order to count the length-$k$ partitions, it suffices
to subtract the \# of length-$\left(  \leq k-1\right)  $-partitions from the
\# of length-$\left(  \leq k\right)  $ partitions);

\item we rename the indeterminate $x$ as $q$.
\end{itemize}

\subsubsection{Definition}

\begin{convention}
\label{conv.pars.qbinom.q}In this section, we will mostly be using FPSs in the
indeterminate $q$. That is, we call the indeterminate $q$ rather than $x$.
Thus, e.g., our formula%
\[
\prod_{n>0}\left(  1-x^{n}\right)  ^{-1}=\prod_{n>0}\dfrac{1}{1-x^{n}}%
=\sum_{n\in\mathbb{N}}p\left(  n\right)  x^{n}=\sum_{\substack{\lambda\text{
is a}\\\text{partition}}}x^{\left\vert \lambda\right\vert }%
\]
becomes%
\[
\prod_{n>0}\left(  1-q^{n}\right)  ^{-1}=\prod_{n>0}\dfrac{1}{1-q^{n}}%
=\sum_{n\in\mathbb{N}}p\left(  n\right)  q^{n}=\sum_{\substack{\lambda\text{
is a}\\\text{partition}}}q^{\left\vert \lambda\right\vert }.
\]
The ring of FPSs in the indeterminate $q$ over a commutative ring $K$ will be
denoted by $K\left[  \left[  q\right]  \right]  $. The ring of polynomials in
the indeterminate $q$ over $K$ will be denoted by $K\left[  q\right]  $.
\end{convention}

\begin{definition}
\label{def.pars.qbinom.qbinom}Let $n\in\mathbb{N}$ and $k\in\mathbb{Z}$.
\medskip

\textbf{(a)} The $q$\emph{-binomial coefficient} (or \emph{Gaussian binomial
coefficient}) $\dbinom{n}{k}_{q}$ is defined to be the polynomial%
\[
\sum_{\substack{\lambda\text{ is a partition}\\\text{with largest part }\leq
n-k\\\text{and length }\leq k}}q^{\left\vert \lambda\right\vert }\in
\mathbb{Z}\left[  q\right]  .
\]

This is also denoted by $%
\genfrac{[}{]}{0pt}{0}{n}{k}%
$ (but this notation has other meanings, too, and suppresses $q$). \medskip

\textbf{(b)} If $a$ is any element of a ring $A$, then we set%
\begin{align*}
\dbinom{n}{k}_{a}:=  &  \ \dbinom{n}{k}_{q}\left[  a\right] \\
&  \ \ \ \ \ \ \ \ \ \ \left(  \text{this means the result of substituting
}a\text{ for }q\text{ in }\dbinom{n}{k}_{q}\right) \\
=  &  \ \sum_{\substack{\lambda\text{ is a partition}\\\text{with largest part
}\leq n-k\\\text{and length }\leq k}}a^{\left\vert \lambda\right\vert }\in A.
\end{align*}

\end{definition}

\begin{remark}
\label{rmk.pars.qbinom.poly}The $\dbinom{n}{k}_{q}$ we defined in Definition
\ref{def.pars.qbinom.qbinom} \textbf{(a)} is really a polynomial, not merely a
FPS, because (for any given $n$ and $k$) there are only finitely many
partitions with largest part $\leq n-k$ and length $\leq k$.
\end{remark}

\begin{remark}
\label{rmk.pars.qbinom.neg}The notation $\dbinom{n}{k}_{q}$ (and the name
\textquotedblleft$q$-binomial coefficient\textquotedblright) suggests a
similarity to the usual binomial coefficient $\dbinom{n}{k}$. And indeed, we
will soon see that $\dbinom{n}{k}_{1}=\dbinom{n}{k}$.

Note, however, that $\dbinom{n}{k}_{q}$ is only defined for $n\in\mathbb{N}$
and $k\in\mathbb{Z}$ (unlike $\dbinom{n}{k}$, which we defined for arbitrary
$n,k\in\mathbb{C}$). It is possible to extend it to negative integers $n$, but
this will result in a Laurent polynomial. (See Exercise
\ref{exe.pars.qbinom.upneg} for this extension.)
\end{remark}

\begin{example}
\label{exa.pars.qbinom.exa1}\textbf{(a)} We have%
\begin{align*}
\dbinom{3}{2}_{q}  &  =\sum_{\substack{\lambda\text{ is a partition}%
\\\text{with largest part }\leq1\\\text{and length }\leq2}}q^{\left\vert
\lambda\right\vert }=q^{\left\vert \left(  1,1\right)  \right\vert
}+q^{\left\vert \left(  1\right)  \right\vert }+q^{\left\vert \left(
{}\right)  \right\vert }\\
&  \ \ \ \ \ \ \ \ \ \ \ \ \ \ \ \ \ \ \ \ \left(
\begin{array}
[c]{c}%
\text{since the partitions with largest part }\leq1\\
\text{and length }\leq2\text{ are }\left(  1,1\right)  \text{, }\left(
1\right)  \text{ and }\left(  {}\right)
\end{array}
\right) \\
&  =q^{2}+q^{1}+q^{0}=q^{2}+q+1
\end{align*}
and%
\begin{align*}
\dbinom{4}{2}_{q}  &  =\sum_{\substack{\lambda\text{ is a partition}%
\\\text{with largest part }\leq2\\\text{and length }\leq2}}q^{\left\vert
\lambda\right\vert }\\
&  =q^{\left\vert \left(  2,2\right)  \right\vert }+q^{\left\vert \left(
2,1\right)  \right\vert }+q^{\left\vert \left(  2\right)  \right\vert
}+q^{\left\vert \left(  1,1\right)  \right\vert }+q^{\left\vert \left(
1\right)  \right\vert }+q^{\left\vert \left(  {}\right)  \right\vert }\\
&  =q^{4}+q^{3}+q^{2}+q^{2}+q^{1}+q^{0}=q^{4}+q^{3}+2q^{2}+q+1.
\end{align*}

\textbf{(b)} For any $n\in\mathbb{N}$ and any negative $k\in\mathbb{Z}$, we
have $\dbinom{n}{k}_{q}=0$, since there exist no partitions $\lambda$ with
length $\leq k$.
\end{example}

\subsubsection{Basic properties}

Let us show two slightly different (but equivalent) ways to express
$q$-binomial coefficients:

\begin{proposition}
\label{prop.pars.qbinom.alt-defs}Let $n\in\mathbb{N}$ and $k\in\mathbb{Z}$.
\medskip

\textbf{(a)} We have%
\[
\dbinom{n}{k}_{q}=\sum_{0\leq i_{1}\leq i_{2}\leq\cdots\leq i_{k}\leq
n-k}q^{i_{1}+i_{2}+\cdots+i_{k}}.
\]
Here, the sum ranges over all weakly increasing $k$-tuples $\left(
i_{1},i_{2},\ldots,i_{k}\right)  \in\left\{  0,1,\ldots,n-k\right\}  ^{k}$. If
$k>n$, then this is an empty sum (since the set $\left\{  0,1,\ldots
,n-k\right\}  $ is empty in this case, and thus its $k$-th power $\left\{
0,1,\ldots,n-k\right\}  ^{k}$ is also empty because $k>n\geq0$). If $k<0$,
then this is also an empty sum (since there are no $k$-tuples when $k<0$).
\medskip

\textbf{(b)} Set $\operatorname*{sum}S=\sum_{s\in S}s$ for any finite set $S$
of integers. (For example, $\operatorname*{sum}\left\{  2,4,5\right\}
=2+4+5=11$.) Then, we have%
\[
\dbinom{n}{k}_{q}=\sum_{\substack{S\subseteq\left\{  1,2,\ldots,n\right\}
;\\\left\vert S\right\vert =k}}q^{\operatorname*{sum}S-\left(  1+2+\cdots
+k\right)  }.
\]

\textbf{(c)} We have%
\[
\dbinom{n}{k}_{1}=\dbinom{n}{k}.
\]

\end{proposition}

\begin{example}
For example, let us compute $\dbinom{5}{2}_{q}$ using Proposition
\ref{prop.pars.qbinom.alt-defs} \textbf{(b)}. Namely, applying Proposition
\ref{prop.pars.qbinom.alt-defs} \textbf{(b)} to $n=5$ and $k=2$, we obtain%
\begin{align*}
\dbinom{5}{2}_{q}  &  =\sum_{\substack{S\subseteq\left\{  1,2,\ldots
,5\right\}  ;\\\left\vert S\right\vert =2}}q^{\operatorname*{sum}S-\left(
1+2\right)  }\\
&  =q^{\left(  1+2\right)  -\left(  1+2\right)  }+q^{\left(  1+3\right)
-\left(  1+2\right)  }+q^{\left(  1+4\right)  -\left(  1+2\right)
}+q^{\left(  1+5\right)  -\left(  1+2\right)  }\\
&  \ \ \ \ \ \ \ \ \ \ +q^{\left(  2+3\right)  -\left(  1+2\right)
}+q^{\left(  2+4\right)  -\left(  1+2\right)  }+q^{\left(  2+5\right)
-\left(  1+2\right)  }\\
&  \ \ \ \ \ \ \ \ \ \ +q^{\left(  3+4\right)  -\left(  1+2\right)
}+q^{\left(  3+5\right)  -\left(  1+2\right)  }+q^{\left(  4+5\right)
-\left(  1+2\right)  }\\
&  \ \ \ \ \ \ \ \ \ \ \ \ \ \ \ \ \ \ \ \ \left(
\begin{array}
[c]{c}%
\text{since the }2\text{-element subsets of }\left\{  1,2,\ldots,5\right\}
\text{ are}\\
\left\{  1,2\right\}  ,\ \left\{  1,3\right\}  ,\ \left\{  1,4\right\}
,\ \left\{  1,5\right\}  ,\ \left\{  2,3\right\}  ,\ \left\{  2,4\right\}  ,\\
\left\{  2,5\right\}  ,\ \left\{  3,4\right\}  ,\ \left\{  3,5\right\}
,\ \left\{  4,5\right\}
\end{array}
\right) \\
&  =q^{0}+q^{1}+q^{2}+q^{3}+q^{2}+q^{3}+q^{4}+q^{4}+q^{5}+q^{6}\\
&  =1+q+2q^{2}+2q^{3}+2q^{4}+q^{5}+q^{6}.
\end{align*}

\end{example}

\begin{proof}
[Proof of Proposition \ref{prop.pars.qbinom.alt-defs}.]\textbf{(a)} The
definition of $\dbinom{n}{k}_{q}$ yields%
\begin{align}
\dbinom{n}{k}_{q}  &  =\sum_{\substack{\lambda\text{ is a partition}%
\\\text{with largest part }\leq n-k\\\text{and length }\leq k}}q^{\left\vert
\lambda\right\vert }\nonumber\\
&  =\sum_{\ell=0}^{k}\ \ \sum_{\substack{\lambda\text{ is a partition}%
\\\text{with largest part }\leq n-k\\\text{and length }\ell}}q^{\left\vert
\lambda\right\vert }. \label{pf.prop.pars.qbinom.alt-defs.a.1}%
\end{align}
Now, let us simplify the inner sum on the right hand side.

Fix $\ell\in\left\{  0,1,\ldots,k\right\}  $. Then, any partition $\lambda$
with length $\ell$ has the form $\left(  \lambda_{1},\lambda_{2}%
,\ldots,\lambda_{\ell}\right)  $ for some nonnegative integers $\lambda
_{1},\lambda_{2},\ldots,\lambda_{\ell}$ satisfying $\lambda_{1}\geq\lambda
_{2}\geq\cdots\geq\lambda_{\ell}>0$ (by the definitions of \textquotedblleft
partition\textquotedblright\ and \textquotedblleft length\textquotedblright).
Moreover, this partition $\lambda$ has largest part $\leq n-k$ if and only if
its entries satisfy $n-k\geq\lambda_{1}\geq\lambda_{2}\geq\cdots\geq
\lambda_{\ell}>0$. Finally, the size $\left\vert \lambda\right\vert $ of this
partition equals $\lambda_{1}+\lambda_{2}+\cdots+\lambda_{\ell}$. Hence, we
can rewrite the sum%
\[
\sum_{\substack{\lambda\text{ is a partition}\\\text{with largest part }\leq
n-k\\\text{and length }\ell}}q^{\left\vert \lambda\right\vert }%
\ \ \ \ \ \ \ \ \ \ \text{as}\ \ \ \ \ \ \ \ \ \ \sum_{\substack{\left(
\lambda_{1},\lambda_{2},\ldots,\lambda_{\ell}\right)  \in\mathbb{N}^{\ell
};\\n-k\geq\lambda_{1}\geq\lambda_{2}\geq\cdots\geq\lambda_{\ell}%
>0}}q^{\lambda_{1}+\lambda_{2}+\cdots+\lambda_{\ell}}.
\]
In other words, we have%
\begin{equation}
\sum_{\substack{\lambda\text{ is a partition}\\\text{with largest part }\leq
n-k\\\text{and length }\ell}}q^{\left\vert \lambda\right\vert }=\sum
_{\substack{\left(  \lambda_{1},\lambda_{2},\ldots,\lambda_{\ell}\right)
\in\mathbb{N}^{\ell};\\n-k\geq\lambda_{1}\geq\lambda_{2}\geq\cdots\geq
\lambda_{\ell}>0}}q^{\lambda_{1}+\lambda_{2}+\cdots+\lambda_{\ell}}.
\label{pf.prop.pars.qbinom.alt-defs.a.2}%
\end{equation}

Next, for any $k$-tuple $\left(  \lambda_{1},\lambda_{2},\ldots,\lambda
_{k}\right)  \in\mathbb{N}^{k}$, let us define $\operatorname*{numpos}\left(
\lambda_{1},\lambda_{2},\ldots,\lambda_{k}\right)  $ to be the number of
positive entries of this $k$-tuple (i.e., the number of $i\in\left\{
1,2,\ldots,k\right\}  $ satisfying $\lambda_{i}>0$). For example,
$\operatorname*{numpos}\left(  4,2,2,0\right)  =3$ and \newline%
$\operatorname*{numpos}\left(  4,2,2,1\right)  =4$ and $\operatorname*{numpos}%
\left(  0,0,0,0\right)  =0$. The following is easy but important:

\begin{statement}
\textit{Observation 1:} Let $\left(  \lambda_{1},\lambda_{2},\ldots
,\lambda_{k}\right)  \in\mathbb{N}^{k}$ be any $k$-tuple satisfying
$\lambda_{1}\geq\lambda_{2}\geq\cdots\geq\lambda_{k}\geq0$ and
$\operatorname*{numpos}\left(  \lambda_{1},\lambda_{2},\ldots,\lambda
_{k}\right)  =\ell$. Then:

\textbf{(a)} The first $\ell$ entries of $\left(  \lambda_{1},\lambda
_{2},\ldots,\lambda_{k}\right)  $ are positive (i.e., we have $\lambda_{i}>0$
for all $i\in\left\{  1,2,\ldots,\ell\right\}  $).

\textbf{(b)} The last $k-\ell$ entries of $\left(  \lambda_{1},\lambda
_{2},\ldots,\lambda_{k}\right)  $ are $0$ (i.e., we have $\lambda_{i}=0$ for
all $i\in\left\{  \ell+1,\ell+2,\ldots,k\right\}  $).

\textbf{(c)} We have $\lambda_{1}+\lambda_{2}+\cdots+\lambda_{\ell}%
=\lambda_{1}+\lambda_{2}+\cdots+\lambda_{k}$.

\textbf{(d)} The $\ell$-tuple $\left(  \lambda_{1},\lambda_{2},\ldots
,\lambda_{\ell}\right)  $ is a partition.
\end{statement}

\begin{fineprint}
[\textit{Proof of Observation 1:} We have $\operatorname*{numpos}\left(
\lambda_{1},\lambda_{2},\ldots,\lambda_{k}\right)  =\ell$. In other words, the
$k$-tuple $\left(  \lambda_{1},\lambda_{2},\ldots,\lambda_{k}\right)  $ has
exactly $\ell$ positive entries. Since this $k$-tuple is weakly decreasing
(because $\lambda_{1}\geq\lambda_{2}\geq\cdots\geq\lambda_{k}$), these $\ell$
positive entries must be concentrated at the beginning of the $k$-tuple; i.e.,
they must be the first $\ell$ entries of the $k$-tuple. Hence, the first
$\ell$ entries of $\left(  \lambda_{1},\lambda_{2},\ldots,\lambda_{k}\right)
$ are positive. This proves Observation 1 \textbf{(a)}.

We have shown that the $k$-tuple $\left(  \lambda_{1},\lambda_{2}%
,\ldots,\lambda_{k}\right)  $ has exactly $\ell$ positive entries, and they
are the first $\ell$ entries of this $k$-tuple. Hence, the remaining $k-\ell$
entries of this $k$-tuple are nonpositive. Since these entries are nonnegative
as well (because $\lambda_{1}\geq\lambda_{2}\geq\cdots\geq\lambda_{k}\geq0$),
we thus conclude that they are $0$. In other words, the last $k-\ell$ entries
of $\left(  \lambda_{1},\lambda_{2},\ldots,\lambda_{k}\right)  $ are $0$. This
proves Observation 1 \textbf{(b)}.

Furthermore,%
\begin{align*}
\lambda_{1}+\lambda_{2}+\cdots+\lambda_{k}  &  =\left(  \lambda_{1}%
+\lambda_{2}+\cdots+\lambda_{\ell}\right)  +\underbrace{\left(  \lambda
_{\ell+1}+\lambda_{\ell+2}+\cdots+\lambda_{k}\right)  }_{\substack{=0+0+\cdots
+0\\\text{(by Observation 1 \textbf{(b)})}}}\\
&  =\left(  \lambda_{1}+\lambda_{2}+\cdots+\lambda_{\ell}\right)  +\left(
0+0+\cdots+0\right)  =\lambda_{1}+\lambda_{2}+\cdots+\lambda_{\ell}.
\end{align*}
This proves Observation 1 \textbf{(c)}.

\textbf{(d)} The $\ell$-tuple $\left(  \lambda_{1},\lambda_{2},\ldots
,\lambda_{\ell}\right)  $ is weakly decreasing (since $\lambda_{1}\geq
\lambda_{2}\geq\cdots\geq\lambda_{k}$ entails $\lambda_{1}\geq\lambda_{2}%
\geq\cdots\geq\lambda_{\ell}$) and consists of positive integers (since
Observation 1 \textbf{(a)} says that we have $\lambda_{i}>0$ for all
$i\in\left\{  1,2,\ldots,\ell\right\}  $). Thus, it is a weakly decreasing
tuple of positive integers, i.e., a partition. This proves Observation 1
\textbf{(d)}.] \medskip
\end{fineprint}

Let us recall that $\ell\in\left\{  0,1,\ldots,k\right\}  $, so that $\ell\leq
k$. Hence, any $\ell$-tuple $\left(  \lambda_{1},\lambda_{2},\ldots
,\lambda_{\ell}\right)  \in\mathbb{N}^{\ell}$ can be extended to a $k$-tuple
$\left(  \lambda_{1},\lambda_{2},\ldots,\lambda_{k}\right)  \in\mathbb{N}^{k}$
by inserting $k-\ell$ zeroes at the end (i.e., by setting $\lambda_{\ell
+1}=\lambda_{\ell+2}=\cdots=\lambda_{k}=0$).\ \ \ \ \footnote{For example, if
$\ell=3$ and $k=5$, then the $\ell$-tuple $\left(  4,2,2\right)  $ gets
extended to the $k$-tuple $\left(  4,2,2,0,0\right)  $.} If the original
$\ell$-tuple $\left(  \lambda_{1},\lambda_{2},\ldots,\lambda_{\ell}\right)
\in\mathbb{N}^{\ell}$ was a partition with largest part $\leq n-k$, then the
extended $k$-tuple $\left(  \lambda_{1},\lambda_{2},\ldots,\lambda_{k}\right)
=\left(  \lambda_{1},\lambda_{2},\ldots,\lambda_{\ell},\underbrace{0,0,\ldots
,0}_{k-\ell\text{ zeroes}}\right)  $ will satisfy $n-k\geq\lambda_{1}%
\geq\lambda_{2}\geq\cdots\geq\lambda_{k}\geq0$ (since $n-k\geq\lambda_{1}%
\geq\lambda_{2}\geq\cdots\geq\lambda_{\ell}$ and $\lambda_{\ell}\geq
0=\lambda_{\ell+1}=\lambda_{\ell+2}=\cdots=\lambda_{k}\geq0$) and
$\operatorname*{numpos}\left(  \lambda_{1},\lambda_{2},\ldots,\lambda
_{k}\right)  =\ell$ (since its first $\ell$ entries are positive whereas its
remaining $k-\ell$ entries are $0$).

Conversely, if $\left(  \lambda_{1},\lambda_{2},\ldots,\lambda_{k}\right)
\in\mathbb{N}^{k}$ is a $k$-tuple satisfying $n-k\geq\lambda_{1}\geq
\lambda_{2}\geq\cdots\geq\lambda_{k}\geq0$ and $\operatorname*{numpos}\left(
\lambda_{1},\lambda_{2},\ldots,\lambda_{k}\right)  =\ell$, then $\left(
\lambda_{1},\lambda_{2},\ldots,\lambda_{\ell}\right)  $ is a
partition\footnote{because of Observation 1 \textbf{(d)}} with largest part
$\leq n-k$\ \ \ \ \footnote{since all parts $\lambda_{1},\lambda_{2}%
,\ldots,\lambda_{\ell}$ of this partition are $\leq n-k$ (because
$n-k\geq\lambda_{1}\geq\lambda_{2}\geq\cdots\geq\lambda_{k}$)} and length
$\ell$. Thus, we obtain a map%
\begin{align*}
&  \text{from }\left\{  k\text{-tuples }\left(  \lambda_{1},\lambda_{2}%
,\ldots,\lambda_{k}\right)  \in\mathbb{N}^{k}\text{ satisfying }n-k\geq
\lambda_{1}\geq\lambda_{2}\geq\cdots\geq\lambda_{k}\geq0\right. \\
&  \ \ \ \ \ \ \ \ \ \ \ \ \ \ \ \ \ \ \ \ \left.  \text{ and }%
\operatorname*{numpos}\left(  \lambda_{1},\lambda_{2},\ldots,\lambda
_{k}\right)  =\ell\right\} \\
&  \text{to }\left\{  \text{partitions with largest part }\leq n-k\text{ and
length }\ell\right\}
\end{align*}
which sends any $k$-tuple $\left(  \lambda_{1},\lambda_{2},\ldots,\lambda
_{k}\right)  $ to the partition $\left(  \lambda_{1},\lambda_{2}%
,\ldots,\lambda_{\ell}\right)  $. On the other hand, we have a map%
\begin{align*}
&  \text{from }\left\{  \text{partitions with largest part }\leq n-k\text{ and
length }\ell\right\} \\
&  \text{to }\left\{  k\text{-tuples }\left(  \lambda_{1},\lambda_{2}%
,\ldots,\lambda_{k}\right)  \in\mathbb{N}^{k}\text{ satisfying }n-k\geq
\lambda_{1}\geq\lambda_{2}\geq\cdots\geq\lambda_{k}\geq0\right. \\
&  \ \ \ \ \ \ \ \ \ \ \ \ \ \ \ \ \ \ \ \ \left.  \text{ and }%
\operatorname*{numpos}\left(  \lambda_{1},\lambda_{2},\ldots,\lambda
_{k}\right)  =\ell\right\}
\end{align*}
which sends any partition $\left(  \lambda_{1},\lambda_{2},\ldots
,\lambda_{\ell}\right)  $ to the $k$-tuple $\left(  \lambda_{1},\lambda
_{2},\ldots,\lambda_{\ell},\underbrace{0,0,\ldots,0}_{k-\ell\text{ zeroes}%
}\right)  $\ \ \ \ \footnote{because we have shown in the previous paragraph
that if $\left(  \lambda_{1},\lambda_{2},\ldots,\lambda_{\ell}\right)  $ is a
partition with largest part $\leq n-k$, then we can extend it to a $k$-tuple
$\left(  \lambda_{1},\lambda_{2},\ldots,\lambda_{k}\right)  =\left(
\lambda_{1},\lambda_{2},\ldots,\lambda_{\ell},\underbrace{0,0,\ldots
,0}_{k-\ell\text{ zeroes}}\right)  $ by inserting $k-\ell$ zeroes at the end,
and this extended $k$-tuple will satisfy $n-k\geq\lambda_{1}\geq\lambda
_{2}\geq\cdots\geq\lambda_{k}\geq0$ and $\operatorname*{numpos}\left(
\lambda_{1},\lambda_{2},\ldots,\lambda_{k}\right)  =\ell$}. These two maps are
mutually inverse\footnote{Indeed, the first map removes the last $k-\ell$
entries from a $k$-tuple, whereas the second map inserts $k-\ell$ zeroes at
the end of the partition. Thus, if we apply the first map after the second
map, we clearly recover the partition that we started with. If we apply the
second map after the first map, then we end up replacing the last $k-\ell$
entries of our $k$-tuple by zeroes. However, if $\left(  \lambda_{1}%
,\lambda_{2},\ldots,\lambda_{k}\right)  \in\mathbb{N}^{k}$ is any $k$-tuple
satisfying $n-k\geq\lambda_{1}\geq\lambda_{2}\geq\cdots\geq\lambda_{k}\geq0$
and $\operatorname*{numpos}\left(  \lambda_{1},\lambda_{2},\ldots,\lambda
_{k}\right)  =\ell$, then the last $k-\ell$ entries of this $k$-tuple are $0$
(by Observation 1 \textbf{(b)}), and therefore the $k$-tuple does not change
if we replace these $k-\ell$ entries by zeroes.}, and therefore are
bijections. In particular, this shows that the first map is a bijection. This
bijection allows us to replace our partitions $\left(  \lambda_{1},\lambda
_{2},\ldots,\lambda_{\ell}\right)  \in\mathbb{N}^{\ell}$ by $k$-tuples
$\left(  \lambda_{1},\lambda_{2},\ldots,\lambda_{k}\right)  \in\mathbb{N}^{k}$
in the sum $\sum_{\substack{\left(  \lambda_{1},\lambda_{2},\ldots
,\lambda_{\ell}\right)  \in\mathbb{N}^{\ell};\\n-k\geq\lambda_{1}\geq
\lambda_{2}\geq\cdots\geq\lambda_{\ell}>0}}q^{\lambda_{1}+\lambda_{2}%
+\cdots+\lambda_{\ell}}$; we thus find%
\begin{align*}
\sum_{\substack{\left(  \lambda_{1},\lambda_{2},\ldots,\lambda_{\ell}\right)
\in\mathbb{N}^{\ell};\\n-k\geq\lambda_{1}\geq\lambda_{2}\geq\cdots\geq
\lambda_{\ell}>0}}q^{\lambda_{1}+\lambda_{2}+\cdots+\lambda_{\ell}}  &
=\sum_{\substack{\left(  \lambda_{1},\lambda_{2},\ldots,\lambda_{k}\right)
\in\mathbb{N}^{k};\\n-k\geq\lambda_{1}\geq\lambda_{2}\geq\cdots\geq\lambda
_{k}\geq0;\\\operatorname*{numpos}\left(  \lambda_{1},\lambda_{2}%
,\ldots,\lambda_{k}\right)  =\ell}}\ \ \underbrace{q^{\lambda_{1}+\lambda
_{2}+\cdots+\lambda_{\ell}}}_{\substack{=q^{\lambda_{1}+\lambda_{2}%
+\cdots+\lambda_{k}}\\\text{(by Observation 1 \textbf{(c)})}}}\\
&  =\sum_{\substack{\left(  \lambda_{1},\lambda_{2},\ldots,\lambda_{k}\right)
\in\mathbb{N}^{k};\\n-k\geq\lambda_{1}\geq\lambda_{2}\geq\cdots\geq\lambda
_{k}\geq0;\\\operatorname*{numpos}\left(  \lambda_{1},\lambda_{2}%
,\ldots,\lambda_{k}\right)  =\ell}}q^{\lambda_{1}+\lambda_{2}+\cdots
+\lambda_{k}}.
\end{align*}

Now, (\ref{pf.prop.pars.qbinom.alt-defs.a.2}) becomes%
\begin{align}
\sum_{\substack{\lambda\text{ is a partition}\\\text{with largest part }\leq
n-k\\\text{and length }\ell}}q^{\left\vert \lambda\right\vert }  &
=\sum_{\substack{\left(  \lambda_{1},\lambda_{2},\ldots,\lambda_{\ell}\right)
\in\mathbb{N}^{\ell};\\n-k\geq\lambda_{1}\geq\lambda_{2}\geq\cdots\geq
\lambda_{\ell}>0}}q^{\lambda_{1}+\lambda_{2}+\cdots+\lambda_{\ell}}\nonumber\\
&  =\sum_{\substack{\left(  \lambda_{1},\lambda_{2},\ldots,\lambda_{k}\right)
\in\mathbb{N}^{k};\\n-k\geq\lambda_{1}\geq\lambda_{2}\geq\cdots\geq\lambda
_{k}\geq0;\\\operatorname*{numpos}\left(  \lambda_{1},\lambda_{2}%
,\ldots,\lambda_{k}\right)  =\ell}}q^{\lambda_{1}+\lambda_{2}+\cdots
+\lambda_{k}}. \label{pf.prop.pars.qbinom.alt-defs.a.5}%
\end{align}

Now, forget that we fixed $\ell$. We thus have proved
(\ref{pf.prop.pars.qbinom.alt-defs.a.5}) for each $\ell\in\left\{
0,1,\ldots,k\right\}  $. Now, (\ref{pf.prop.pars.qbinom.alt-defs.a.1}) becomes%
\begin{align*}
\dbinom{n}{k}_{q}  &  =\sum_{\ell=0}^{k}\ \ \underbrace{\sum
_{\substack{\lambda\text{ is a partition}\\\text{with largest part }\leq
n-k\\\text{and length }\ell}}q^{\left\vert \lambda\right\vert }}%
_{\substack{=\sum_{\substack{\left(  \lambda_{1},\lambda_{2},\ldots
,\lambda_{k}\right)  \in\mathbb{N}^{k};\\n-k\geq\lambda_{1}\geq\lambda_{2}%
\geq\cdots\geq\lambda_{k}\geq0;\\\operatorname*{numpos}\left(  \lambda
_{1},\lambda_{2},\ldots,\lambda_{k}\right)  =\ell}}q^{\lambda_{1}+\lambda
_{2}+\cdots+\lambda_{k}}\\\text{(by (\ref{pf.prop.pars.qbinom.alt-defs.a.5}%
))}}}\\
&  =\underbrace{\sum_{\ell=0}^{k}\ \ \sum_{\substack{\left(  \lambda
_{1},\lambda_{2},\ldots,\lambda_{k}\right)  \in\mathbb{N}^{k};\\n-k\geq
\lambda_{1}\geq\lambda_{2}\geq\cdots\geq\lambda_{k}\geq
0;\\\operatorname*{numpos}\left(  \lambda_{1},\lambda_{2},\ldots,\lambda
_{k}\right)  =\ell}}}_{=\sum_{\substack{\left(  \lambda_{1},\lambda_{2}%
,\ldots,\lambda_{k}\right)  \in\mathbb{N}^{k};\\n-k\geq\lambda_{1}\geq
\lambda_{2}\geq\cdots\geq\lambda_{k}\geq0}}}q^{\lambda_{1}+\lambda_{2}%
+\cdots+\lambda_{k}}\\
&  =\sum_{\substack{\left(  \lambda_{1},\lambda_{2},\ldots,\lambda_{k}\right)
\in\mathbb{N}^{k};\\n-k\geq\lambda_{1}\geq\lambda_{2}\geq\cdots\geq\lambda
_{k}\geq0}}q^{\lambda_{1}+\lambda_{2}+\cdots+\lambda_{k}}=\sum
_{\substack{\left(  i_{1},i_{2},\ldots,i_{k}\right)  \in\mathbb{N}^{k};\\0\leq
i_{1}\leq i_{2}\leq\cdots\leq i_{k}\leq n-k}}\underbrace{q^{i_{k}%
+i_{k-1}+\cdots+i_{1}}}_{=q^{i_{1}+i_{2}+\cdots+i_{k}}}\\
&  \ \ \ \ \ \ \ \ \ \ \left(
\begin{array}
[c]{c}%
\text{here, we have reversed the }k\text{-tuple }\left(  \lambda_{1}%
,\lambda_{2},\ldots,\lambda_{k}\right)  \text{,}\\
\text{i.e., we have substituted }\left(  i_{k},i_{k-1},\ldots,i_{1}\right)
\text{ for }\left(  \lambda_{1},\lambda_{2},\ldots,\lambda_{k}\right) \\
\text{in our sum}%
\end{array}
\right) \\
&  =\sum_{\substack{\left(  i_{1},i_{2},\ldots,i_{k}\right)  \in\mathbb{N}%
^{k};\\0\leq i_{1}\leq i_{2}\leq\cdots\leq i_{k}\leq n-k}}q^{i_{1}%
+i_{2}+\cdots+i_{k}}=\sum_{0\leq i_{1}\leq i_{2}\leq\cdots\leq i_{k}\leq
n-k}q^{i_{1}+i_{2}+\cdots+i_{k}}.
\end{align*}
This proves Proposition \ref{prop.pars.qbinom.alt-defs} \textbf{(a)}. \medskip

\textbf{(b)} There is a bijection%
\begin{align*}
&  \text{from }\left\{  \text{weakly increasing }k\text{-tuples }\left(
i_{1},i_{2},\ldots,i_{k}\right)  \in\left\{  0,1,\ldots,n-k\right\}
^{k}\right\} \\
&  \text{to }\left\{  \text{strictly increasing }k\text{-tuples }\left(
s_{1},s_{2},\ldots,s_{k}\right)  \in\left\{  1,2,\ldots,n\right\}
^{k}\right\}
\end{align*}
that sends each weakly increasing $k$-tuple $\left(  i_{1},i_{2},\ldots
,i_{k}\right)  $ to $\left(  i_{1}+1,\ i_{2}+2,\ \ldots,\ i_{k}+k\right)  $
(you can think of it as \textquotedblleft spacing the $i_{j}$s
apart\textquotedblright, i.e., increasing the distance between any two
consecutive $i_{j}$'s by $1$ and also increasing $i_{1}$ by $1$). The inverse
of this bijection sends each strictly increasing $k$-tuple $\left(
s_{1},s_{2},\ldots,s_{k}\right)  $ to $\left(  s_{1}-1,\ s_{2}-2,\ \ldots
,\ s_{k}-k\right)  $. Thus, we can substitute $\left(  s_{1}-1,\ s_{2}%
-2,\ \ldots,\ s_{k}-k\right)  $ for $\left(  i_{1},i_{2},\ldots,i_{k}\right)
$ in the sum%
\[
\sum_{0\leq i_{1}\leq i_{2}\leq\cdots\leq i_{k}\leq n-k}q^{i_{1}+i_{2}%
+\cdots+i_{k}}.
\]
Hence we obtain%
\begin{align*}
\sum_{0\leq i_{1}\leq i_{2}\leq\cdots\leq i_{k}\leq n-k}q^{i_{1}+i_{2}%
+\cdots+i_{k}}  &  =\sum_{1\leq s_{1}<s_{2}<\cdots<s_{k}\leq n}%
\underbrace{q^{\left(  s_{1}-1\right)  +\left(  s_{2}-2\right)  +\cdots
+\left(  s_{k}-k\right)  }}_{=q^{\left(  s_{1}+s_{2}+\cdots+s_{k}\right)
-\left(  1+2+\cdots+k\right)  }}\\
&  =\sum_{1\leq s_{1}<s_{2}<\cdots<s_{k}\leq n}q^{\left(  s_{1}+s_{2}%
+\cdots+s_{k}\right)  -\left(  1+2+\cdots+k\right)  }.
\end{align*}

On the other hand, there is a bijection%
\begin{align*}
&  \text{from }\left\{  \text{strictly increasing }k\text{-tuples }\left(
s_{1},s_{2},\ldots,s_{k}\right)  \in\left\{  1,2,\ldots,n\right\}
^{k}\right\} \\
&  \text{to }\left\{  k\text{-element subsets of }\left\{  1,2,\ldots
,n\right\}  \right\}
\end{align*}
that sends each $k$-tuple $\left(  s_{1},s_{2},\ldots,s_{k}\right)  $ to the
subset $\left\{  s_{1},s_{2},\ldots,s_{k}\right\}  $. (This map is indeed a
bijection, because any $k$-element subset of $\left\{  1,2,\ldots,n\right\}  $
can be written as $\left\{  s_{1},s_{2},\ldots,s_{k}\right\}  $ for a unique
strictly increasing $k$-tuple $\left(  s_{1},s_{2},\ldots,s_{k}\right)
\in\left\{  1,2,\ldots,n\right\}  ^{k}$; in fact, this is simply saying that
there is a unique way of listing the elements of this subset in increasing order.)

Because of this bijection, we have%
\[
\sum_{1\leq s_{1}<s_{2}<\cdots<s_{k}\leq n}q^{\left(  s_{1}+s_{2}+\cdots
+s_{k}\right)  -\left(  1+2+\cdots+k\right)  }=\sum_{\substack{S\subseteq
\left\{  1,2,\ldots,n\right\}  ;\\\left\vert S\right\vert =k}%
}q^{\operatorname*{sum}S-\left(  1+2+\cdots+k\right)  }%
\]
(because $s_{1}+s_{2}+\cdots+s_{k}=\operatorname*{sum}\left\{  s_{1}%
,s_{2},\ldots,s_{k}\right\}  $ for any strictly increasing $k$-tuple $\left(
s_{1},s_{2},\ldots,s_{k}\right)  \in\left\{  1,2,\ldots,n\right\}  ^{k}$).

Now, Proposition \ref{prop.pars.qbinom.alt-defs} \textbf{(a)} yields%
\begin{align*}
\dbinom{n}{k}_{q}  &  =\sum_{0\leq i_{1}\leq i_{2}\leq\cdots\leq i_{k}\leq
n-k}q^{i_{1}+i_{2}+\cdots+i_{k}}=\sum_{1\leq s_{1}<s_{2}<\cdots<s_{k}\leq
n}q^{\left(  s_{1}+s_{2}+\cdots+s_{k}\right)  -\left(  1+2+\cdots+k\right)
}\\
&  =\sum_{\substack{S\subseteq\left\{  1,2,\ldots,n\right\}  ;\\\left\vert
S\right\vert =k}}q^{\operatorname*{sum}S-\left(  1+2+\cdots+k\right)  }.
\end{align*}
This proves Proposition \ref{prop.pars.qbinom.alt-defs} \textbf{(b)}. \medskip

\textbf{(c)} Proposition \ref{prop.pars.qbinom.alt-defs} \textbf{(b)} yields
\[
\dbinom{n}{k}_{q}=\sum_{\substack{S\subseteq\left\{  1,2,\ldots,n\right\}
;\\\left\vert S\right\vert =k}}q^{\operatorname*{sum}S-\left(  1+2+\cdots
+k\right)  }.
\]
Substituting $1$ for $q$ in this equality, we find%
\begin{align*}
\dbinom{n}{k}_{1}  &  =\sum_{\substack{S\subseteq\left\{  1,2,\ldots
,n\right\}  ;\\\left\vert S\right\vert =k}}\underbrace{1^{\operatorname*{sum}%
S-\left(  1+2+\cdots+k\right)  }}_{=1}=\sum_{\substack{S\subseteq\left\{
1,2,\ldots,n\right\}  ;\\\left\vert S\right\vert =k}}1\\
&  =\left(  \text{\# of subsets }S\subseteq\left\{  1,2,\ldots,n\right\}
\text{ satisfying }\left\vert S\right\vert =k\right) \\
&  =\left(  \text{\# of }k\text{-element subsets of }\left\{  1,2,\ldots
,n\right\}  \right)  =\dbinom{n}{k}.
\end{align*}
This proves Proposition \ref{prop.pars.qbinom.alt-defs} \textbf{(c)}.
\end{proof}

The following property of $q$-binomial coefficients generalizes Proposition
\ref{prop.binom.0}:

\begin{proposition}
\label{prop.pars.qbinom.0}Let $n,k\in\mathbb{N}$ satisfy $k>n$. Then,
$\dbinom{n}{k}_{q}=0$.
\end{proposition}

\begin{proof}
From $k>n$, we obtain $n-k<0$. The definition of $\dbinom{n}{k}_{q}$ yields%
\begin{equation}
\dbinom{n}{k}_{q}=\sum_{\substack{\lambda\text{ is a partition}\\\text{with
largest part }\leq n-k\\\text{and length }\leq k}}q^{\left\vert \lambda
\right\vert }. \label{pf.prop.pars.qbinom.0.1}%
\end{equation}
The sum on the right hand side is an empty sum, since there exists no
partition with largest part $\leq n-k$ (because $n-k<0$). Thus,
(\ref{pf.prop.pars.qbinom.0.1}) rewrites as $\dbinom{n}{k}_{q}=\left(
\text{empty sum}\right)  =0$, and this proves Proposition
\ref{prop.pars.qbinom.0}.
\end{proof}

\begin{proposition}
\label{prop.pars.qbinom.n0}We have $\dbinom{n}{0}_{q}=\dbinom{n}{n}_{q}=1$ for
each $n\in\mathbb{N}$.
\end{proposition}

\begin{proof}
This is easy and left as a homework exercise (Exercise
\ref{exe.pars.qbinom.basics} \textbf{(a)}).
\end{proof}

The next convention mirrors a convention we made for the (usual) binomial coefficients:

\begin{convention}
\label{conv.pars.qbinom.neg-k}Let $n\in\mathbb{N}$. For any $k\notin%
\mathbb{Z}$, we set $\dbinom{n}{k}_{q}:=0$.
\end{convention}

The following theorem gives not one, but two analogues (\textquotedblleft%
$q$-analogues\textquotedblright) of the recurrence relation $\dbinom{n}%
{k}=\dbinom{n-1}{k-1}+\dbinom{n-1}{k}$ (from Proposition \ref{prop.binom.rec}):

\begin{theorem}
\label{thm.pars.qbinom.rec}Let $n$ be a positive integer. Let $k\in\mathbb{N}%
$. Then: \medskip

\textbf{(a)} We have%
\[
\dbinom{n}{k}_{q}=q^{n-k}\dbinom{n-1}{k-1}_{q}+\dbinom{n-1}{k}_{q}.
\]
(If $n-k<0$, then this is -- strictly speaking -- an equality in
$\mathbb{Z}\left[  q,q^{-1}\right]  $ because $q^{n-k}$ is a negative power of
$q$; but in fact Proposition \ref{prop.pars.qbinom.0} yields that
$\dbinom{n-1}{k-1}_{q}=0$ in this case, so this negative power disappears.)
\medskip

\textbf{(b)} We have%
\[
\dbinom{n}{k}_{q}=\dbinom{n-1}{k-1}_{q}+q^{k}\dbinom{n-1}{k}_{q}.
\]

\end{theorem}

\begin{proof}
\textbf{(a)} This is similar to the combinatorial proof of the recurrence
relation for binomial coefficients.

If $k=0$, then the claim we are proving boils down to $1=q^{n-k}0+1$ (because
Proposition \ref{prop.pars.qbinom.n0} yields $\dbinom{n}{0}_{q}=1$ and
$\dbinom{n-1}{0}_{q}=1$, and because Example \ref{exa.pars.qbinom.exa1}
\textbf{(b)} shows that $\dbinom{n-1}{-1}_{q}=0$). Hence, we WLOG assume that
$k>0$. Thus, $k-1\in\mathbb{N}$.

Proposition \ref{prop.pars.qbinom.alt-defs} \textbf{(b)} says that%
\begin{equation}
\dbinom{n}{k}_{q}=\sum_{\substack{S\subseteq\left\{  1,2,\ldots,n\right\}
;\\\left\vert S\right\vert =k}}q^{\operatorname*{sum}S-\left(  1+2+\cdots
+k\right)  }. \label{pf.thm.pars.qbinom.rec.nk}%
\end{equation}
Proposition \ref{prop.pars.qbinom.alt-defs} \textbf{(b)} (applied to $n-1$ and
$k-1$ instead of $n$ and $k$) yields%
\begin{equation}
\dbinom{n-1}{k-1}_{q}=\sum_{\substack{S\subseteq\left\{  1,2,\ldots
,n-1\right\}  ;\\\left\vert S\right\vert =k-1}}q^{\operatorname*{sum}S-\left(
1+2+\cdots+\left(  k-1\right)  \right)  }.
\label{pf.thm.pars.qbinom.rec.n-1k-1}%
\end{equation}
Proposition \ref{prop.pars.qbinom.alt-defs} \textbf{(b)} (applied to $n-1$
instead of $n$) yields%
\begin{equation}
\dbinom{n-1}{k}_{q}=\sum_{\substack{S\subseteq\left\{  1,2,\ldots,n-1\right\}
;\\\left\vert S\right\vert =k}}q^{\operatorname*{sum}S-\left(  1+2+\cdots
+k\right)  }. \label{pf.thm.pars.qbinom.rec.n-1k}%
\end{equation}

Let us now make two definitions:

\begin{itemize}
\item A \emph{type-1 subset} will mean a $k$-element subset of $\left\{
1,2,\ldots,n\right\}  $ that contains $n$;

\item A \emph{type-2 subset} will mean a $k$-element subset of $\left\{
1,2,\ldots,n\right\}  $ that does not contain $n$.
\end{itemize}

Each $k$-element subset of $\left\{  1,2,\ldots,n\right\}  $ is either type-1
or type-2 (but not both at the same time). Thus,%
\begin{align*}
&  \sum_{\substack{S\subseteq\left\{  1,2,\ldots,n\right\}  ;\\\left\vert
S\right\vert =k}}q^{\operatorname*{sum}S-\left(  1+2+\cdots+k\right)  }\\
&  =\sum_{\substack{S\subseteq\left\{  1,2,\ldots,n\right\}  ;\\\left\vert
S\right\vert =k;\\S\text{ is type-1}}}q^{\operatorname*{sum}S-\left(
1+2+\cdots+k\right)  }+\sum_{\substack{S\subseteq\left\{  1,2,\ldots
,n\right\}  ;\\\left\vert S\right\vert =k;\\S\text{ is type-2}}%
}q^{\operatorname*{sum}S-\left(  1+2+\cdots+k\right)  }.
\end{align*}

The type-2 subsets are precisely the $k$-element subsets of $\left\{
1,2,\ldots,n-1\right\}  $. Hence,%
\[
\sum_{\substack{S\subseteq\left\{  1,2,\ldots,n\right\}  ;\\\left\vert
S\right\vert =k;\\S\text{ is type-2}}}q^{\operatorname*{sum}S-\left(
1+2+\cdots+k\right)  }=\sum_{\substack{S\subseteq\left\{  1,2,\ldots
,n-1\right\}  ;\\\left\vert S\right\vert =k}}q^{\operatorname*{sum}S-\left(
1+2+\cdots+k\right)  }=\dbinom{n-1}{k}_{q}%
\]
(by (\ref{pf.thm.pars.qbinom.rec.n-1k})).

The type-1 subsets are just the $\left(  k-1\right)  $-element subsets of
$\left\{  1,2,\ldots,n-1\right\}  $ with an $n$ inserted into them; i.e., the
map%
\begin{align*}
\left\{  \left(  k-1\right)  \text{-element subsets of }\left\{
1,2,\ldots,n-1\right\}  \right\}   &  \rightarrow\left\{  \text{type-1
subsets}\right\}  ,\\
S  &  \mapsto S\cup\left\{  n\right\}
\end{align*}
is a bijection. Hence, substituting $S\cup\left\{  n\right\}  $ for $S$ in the
sum, we find%
\begin{align*}
\sum_{\substack{S\subseteq\left\{  1,2,\ldots,n\right\}  ;\\\left\vert
S\right\vert =k;\\S\text{ is type-1}}}q^{\operatorname*{sum}S-\left(
1+2+\cdots+k\right)  }  &  =\sum_{\substack{S\subseteq\left\{  1,2,\ldots
,n-1\right\}  ;\\\left\vert S\right\vert =k-1}%
}\underbrace{q^{\operatorname*{sum}\left(  S\cup\left\{  n\right\}  \right)
-\left(  1+2+\cdots+k\right)  }}_{\substack{=q^{\operatorname*{sum}S+n-\left(
1+2+\cdots+k\right)  }\\\text{(since }S\subseteq\left\{  1,2,\ldots
,n-1\right\}  \text{ entails }n\notin S\\\text{and thus }\operatorname*{sum}%
\left(  S\cup\left\{  n\right\}  \right)  =\operatorname*{sum}S+n\text{)}}}\\
&  =\sum_{\substack{S\subseteq\left\{  1,2,\ldots,n-1\right\}  ;\\\left\vert
S\right\vert =k-1}}\underbrace{q^{\operatorname*{sum}S+n-\left(
1+2+\cdots+k\right)  }}_{\substack{=q^{\operatorname*{sum}S+n-\left(
1+2+\cdots+\left(  k-1\right)  \right)  -k}\\=q^{n-k}q^{\operatorname*{sum}%
S-\left(  1+2+\cdots+\left(  k-1\right)  \right)  }}}\\
&  =\sum_{\substack{S\subseteq\left\{  1,2,\ldots,n-1\right\}  ;\\\left\vert
S\right\vert =k-1}}q^{n-k}q^{\operatorname*{sum}S-\left(  1+2+\cdots+\left(
k-1\right)  \right)  }\\
&  =q^{n-k}\underbrace{\sum_{\substack{S\subseteq\left\{  1,2,\ldots
,n-1\right\}  ;\\\left\vert S\right\vert =k-1}}q^{\operatorname*{sum}S-\left(
1+2+\cdots+\left(  k-1\right)  \right)  }}_{\substack{=\dbinom{n-1}{k-1}%
_{q}\\\text{(by (\ref{pf.thm.pars.qbinom.rec.n-1k-1}))}}}\\
&  =q^{n-k}\dbinom{n-1}{k-1}_{q}.
\end{align*}

All that's left to do is combining what we have found:%
\begin{align*}
\dbinom{n}{k}_{q}  &  =\sum_{\substack{S\subseteq\left\{  1,2,\ldots
,n\right\}  ;\\\left\vert S\right\vert =k}}q^{\operatorname*{sum}S-\left(
1+2+\cdots+k\right)  }\ \ \ \ \ \ \ \ \ \ \left(  \text{by
(\ref{pf.thm.pars.qbinom.rec.nk})}\right) \\
&  =\underbrace{\sum_{\substack{S\subseteq\left\{  1,2,\ldots,n\right\}
;\\\left\vert S\right\vert =k;\\S\text{ is type-1}}}q^{\operatorname*{sum}%
S-\left(  1+2+\cdots+k\right)  }}_{=q^{n-k}\dbinom{n-1}{k-1}_{q}%
}+\underbrace{\sum_{\substack{S\subseteq\left\{  1,2,\ldots,n\right\}
;\\\left\vert S\right\vert =k;\\S\text{ is type-2}}}q^{\operatorname*{sum}%
S-\left(  1+2+\cdots+k\right)  }}_{=\dbinom{n-1}{k}_{q}}\\
&  =q^{n-k}\dbinom{n-1}{k-1}_{q}+\dbinom{n-1}{k}_{q}.
\end{align*}
This proves Theorem \ref{thm.pars.qbinom.rec} \textbf{(a)}. \medskip

\textbf{(b)} This is somewhat similar to Theorem \ref{thm.pars.qbinom.rec}
\textbf{(a)} (but a little bit more complicated). It is left as a homework
exercise (Exercise \ref{exe.pars.qbinom.basics} \textbf{(b)}).
\end{proof}

Next, we shall derive a $q$-analogue of the formula $\dbinom{n}{k}%
=\dfrac{n\left(  n-1\right)  \cdots\left(  n-k+1\right)  }{k!}=\dfrac{n\left(
n-1\right)  \cdots\left(  n-k+1\right)  }{k\left(  k-1\right)  \cdots1}$:

\begin{theorem}
\label{thm.pars.qbinom.quot1}Let $n,k\in\mathbb{N}$ satisfy $n\geq k$. Then:
\medskip

\textbf{(a)} We have%
\begin{align*}
&  \left(  1-q^{k}\right)  \left(  1-q^{k-1}\right)  \cdots\left(
1-q^{1}\right)  \cdot\dbinom{n}{k}_{q}\\
&  =\left(  1-q^{n}\right)  \left(  1-q^{n-1}\right)  \cdots\left(
1-q^{n-k+1}\right)  .
\end{align*}

\textbf{(b)} We have%
\[
\dbinom{n}{k}_{q}=\dfrac{\left(  1-q^{n}\right)  \left(  1-q^{n-1}\right)
\cdots\left(  1-q^{n-k+1}\right)  }{\left(  1-q^{k}\right)  \left(
1-q^{k-1}\right)  \cdots\left(  1-q^{1}\right)  }%
\]
(in the ring $\mathbb{Z}\left[  \left[  q\right]  \right]  $ or in the field
of rational functions over $\mathbb{Q}$).
\end{theorem}

Note that part \textbf{(b)} of Theorem \ref{thm.pars.qbinom.quot1} is the more
intuitive statement, but part \textbf{(a)} is easier to substitute things in
(because substituting something for $q$ in part \textbf{(b)} requires showing
that the denominator remains invertible, whereas part \textbf{(a)} has no
denominators and thus requires no such diligence).

\begin{proof}
[Proof of Theorem \ref{thm.pars.qbinom.quot1}.]This is left as a homework
exercise (Exercise \ref{exe.pars.qbinom.basics} \textbf{(c)}). (Use induction
on $n$ and Theorem \ref{thm.pars.qbinom.rec}.)
\end{proof}

\begin{remark}
If I just gave you the fraction $\dfrac{\left(  1-q^{n}\right)  \left(
1-q^{n-1}\right)  \cdots\left(  1-q^{n-k+1}\right)  }{\left(  1-q^{k}\right)
\left(  1-q^{k-1}\right)  \cdots\left(  1-q^{1}\right)  }$, you would be
surprised to hear that it is a polynomial (i.e., that the denominator divides
the numerator) and has nonnegative coefficients. But given the way we defined
$\dbinom{n}{k}_{q}$, you are now getting this for free from Theorem
\ref{thm.pars.qbinom.quot1}.
\end{remark}

Theorem \ref{thm.pars.qbinom.quot1} \textbf{(b)} can be rewritten in a
somewhat simpler way using the following notations:

\begin{definition}
\label{def.pars.qbinom.qint}\textbf{(a)} For any $n\in\mathbb{N}$, define the
$q$\emph{-integer }$\left[  n\right]  _{q}$ to be
\[
\left[  n\right]  _{q}:=q^{0}+q^{1}+\cdots+q^{n-1}\in\mathbb{Z}\left[
q\right]  .
\]

\textbf{(b)} For any $n\in\mathbb{N}$, define the $q$\emph{-factorial
}$\left[  n\right]  _{q}!$ to be%
\[
\left[  n\right]  _{q}!\ :=\left[  1\right]  _{q}\left[  2\right]  _{q}%
\cdots\left[  n\right]  _{q}\in\mathbb{Z}\left[  q\right]  .
\]

\textbf{(c)} As usual, if $a$ is an element of a ring $A$, then $\left[
n\right]  _{a}$ and $\left[  n\right]  _{a}!$ will mean the results of
substituting $a$ for $q$ in $\left[  n\right]  _{q}$ and $\left[  n\right]
_{q}!$, respectively. Thus, explicitly, $\left[  n\right]  _{a}=a^{0}%
+a^{1}+\cdots+a^{n-1}$ and $\left[  n\right]  _{a}!=\left[  1\right]
_{a}\left[  2\right]  _{a}\cdots\left[  n\right]  _{a}$.
\end{definition}

Alternative notations for $\left[  n\right]  _{q}!$ are $\left[  n\right]
!_{q}$ (used in \cite{Loehr-BC}) and the boldface $\mathbf{(}\boldsymbol{n}%
\mathbf{)!}$ (used in \cite{Stanley-EC1}).

\begin{remark}
\label{rmk.pars.qbinom.qint.frac}For any $n\in\mathbb{N}$, we have%
\[
\left[  n\right]  _{q}=\dfrac{1-q^{n}}{1-q}\ \ \ \ \ \ \ \ \ \ \left(
\text{in }\mathbb{Z}\left[  \left[  q\right]  \right]  \text{ or in the ring
of rational functions over }\mathbb{Q}\right)
\]
and%
\[
\left[  n\right]  _{1}=n\ \ \ \ \ \ \ \ \ \ \text{and}%
\ \ \ \ \ \ \ \ \ \ \left[  n\right]  _{1}!=n!.
\]

\end{remark}

\begin{proof}
[Proof of Remark \ref{rmk.pars.qbinom.qint.frac}.]Let $n\in\mathbb{N}$. We
have%
\[
\left[  n\right]  _{q}:=q^{0}+q^{1}+\cdots+q^{n-1}=\dfrac{1-q^{n}}{1-q},
\]
since%
\begin{align*}
\left(  1-q\right)  \left(  q^{0}+q^{1}+\cdots+q^{n-1}\right)   &  =\left(
q^{0}+q^{1}+\cdots+q^{n-1}\right)  -q\left(  q^{0}+q^{1}+\cdots+q^{n-1}\right)
\\
&  =\left(  q^{0}+q^{1}+\cdots+q^{n-1}\right)  -\left(  q^{1}+q^{2}%
+\cdots+q^{n}\right) \\
&  =\underbrace{q^{0}}_{=1}-\,q^{n}=1-q^{n}.
\end{align*}

Furthermore, substituting $1$ for $q$ in the equality $\left[  n\right]
_{q}=q^{0}+q^{1}+\cdots+q^{n-1}$, we obtain%
\begin{equation}
\left[  n\right]  _{1}=1^{0}+1^{1}+\cdots+1^{n-1}=\underbrace{1+1+\cdots
+1}_{n\text{ times}}=n. \label{pf.rmk.pars.qbinom.qint.frac.4}%
\end{equation}
Substituting $1$ for $q$ in the equality $\left[  n\right]  _{q}!=\left[
1\right]  _{q}\left[  2\right]  _{q}\cdots\left[  n\right]  _{q}$, we obtain%
\[
\left[  n\right]  _{1}!=\underbrace{\left[  1\right]  _{1}}%
_{\substack{=1\\\text{(by (\ref{pf.rmk.pars.qbinom.qint.frac.4}))}%
}}\underbrace{\left[  2\right]  _{1}}_{\substack{=2\\\text{(by
(\ref{pf.rmk.pars.qbinom.qint.frac.4}))}}}\cdots\underbrace{\left[  n\right]
_{1}}_{\substack{=n\\\text{(by (\ref{pf.rmk.pars.qbinom.qint.frac.4}))}%
}}=1\cdot2\cdot\cdots\cdot n=n!.
\]

\end{proof}

\begin{theorem}
\label{thm.pars.qbinom.quot2}Let $n,k\in\mathbb{N}$ with $n\geq k$. Then,%
\[
\dbinom{n}{k}_{q}=\dfrac{\left[  n\right]  _{q}\left[  n-1\right]  _{q}%
\cdots\left[  n-k+1\right]  _{q}}{\left[  k\right]  _{q}!}=\dfrac{\left[
n\right]  _{q}!}{\left[  k\right]  _{q}!\cdot\left[  n-k\right]  _{q}!}%
\]
(in the ring $\mathbb{Z}\left[  \left[  q\right]  \right]  $ or in the ring of
rational functions over $\mathbb{Q}$).
\end{theorem}

\begin{proof}
This is left as a homework exercise (Exercise \ref{exe.pars.qbinom.basics}
\textbf{(d)}).
\end{proof}

A consequence of this theorem is the following symmetry property of
$q$-binomial coefficients:

\begin{proposition}
\label{prop.pars.qbinom.symm}Let $n,k\in\mathbb{N}$. Then,%
\[
\dbinom{n}{k}_{q}=\dbinom{n}{n-k}_{q}.
\]

\end{proposition}

\begin{proof}
This is left as a homework exercise (Exercise \ref{exe.pars.qbinom.basics}
\textbf{(e)}).
\end{proof}

\subsubsection{$q$-binomial formulas}

The properties we have seen so far are suggesting that $q$-binomial
coefficients not only generalize binomial coefficients, but also share most of
their properties in a somewhat modified form. In other words, we start
expecting most properties of binomial coefficients to generalize to
$q$-binomial coefficients, often in several ways (e.g., the recurrence of the
binomial coefficients generalized in two ways).

Let us see how this expectation holds up for the most famous property of
binomial coefficients: the binomial formula%
\[
\left(  a+b\right)  ^{n}=\sum_{k=0}^{n}\dbinom{n}{k}a^{k}b^{n-k}.
\]
This formula holds whenever $a$ and $b$ are two elements of a commutative
ring, or even more generally, whenever $a$ and $b$ are two commuting elements
of an arbitrary ring. If we want to integrate a $q$ into this formula, we need to

\begin{itemize}
\item either change the structure of the formula,

\item or modify the commutativity assumption.
\end{itemize}

This gives rise to two different \textquotedblleft$q$%
-analogues\textquotedblright\ of the binomial formula. Both are important (one
for the theory of partitions, and another for the theory of
\href{https://www.ams.org/notices/200601/what-is.pdf}{quantum groups}). Here
is the first one (first stated by Rothe in 1811, proved shortly thereafter by
Gruson \cite{Gruson18}):

\begin{theorem}
[1st $q$-binomial theorem]\label{thm.pars.qbinom.binom1}Let $K$ be a
commutative ring. Let $a,b\in K$ and $n\in\mathbb{N}$. In the polynomial ring
$K\left[  q\right]  $, we have%
\[
\left(  aq^{0}+b\right)  \left(  aq^{1}+b\right)  \cdots\left(  aq^{n-1}%
+b\right)  =\sum_{k=0}^{n}q^{k\left(  k-1\right)  /2}\dbinom{n}{k}_{q}%
a^{k}b^{n-k}.
\]

\end{theorem}

Note that setting $q=1$ in Theorem \ref{thm.pars.qbinom.binom1} (i.e.,
substituting $1$ for $q$) recovers the good old binomial formula, since all
the $n$ factors on the left hand side become $a+b$.

There is a straightforward way to prove Theorem \ref{thm.pars.qbinom.binom1}
by induction on $n$ (see Exercise \ref{exe.pars.qbinom.qbinom} \textbf{(a)}).
Let us instead give a nicer argument. This argument will rely on the following
general fact:

\begin{lemma}
\label{lem.prodrule.sum-ai-plus-bi}Let $L$ be a commutative ring. Let
$n\in\mathbb{N}$. Let $\left[  n\right]  $ denote the set $\left\{
1,2,\ldots,n\right\}  $. Let $a_{1},a_{2},\ldots,a_{n}$ be $n$ elements of
$L$. Let $b_{1},b_{2},\ldots,b_{n}$ be $n$ further elements of $L$. Then,%
\begin{equation}
\prod_{i=1}^{n}\left(  a_{i}+b_{i}\right)  =\sum_{S\subseteq\left[  n\right]
}\left(  \prod_{i\in S}a_{i}\right)  \left(  \prod_{i\in\left[  n\right]
\setminus S}b_{i}\right)  . \label{eq.lem.prodrule.sum-ai-plus-bi.eq}%
\end{equation}

\end{lemma}

Lemma \ref{lem.prodrule.sum-ai-plus-bi} is well-known and intuitively clear:
When expanding the product $\prod_{i=1}^{n}\left(  a_{i}+b_{i}\right)
=\left(  a_{1}+b_{1}\right)  \left(  a_{2}+b_{2}\right)  \cdots\left(
a_{n}+b_{n}\right)  $, you obtain a sum of $2^{n}$ terms, each of which is a
product of one addend chosen from each of the $n$ sums $a_{1}+b_{1}%
,a_{2}+b_{2},\ldots,a_{n}+b_{n}$. This is precisely what the right hand side
of (\ref{eq.lem.prodrule.sum-ai-plus-bi.eq}) is. A rigorous proof of Lemma
\ref{lem.prodrule.sum-ai-plus-bi} can be found in \cite[Exercise 6.1
\textbf{(a)}]{detnotes}.

\begin{proof}
[Proof of Theorem \ref{thm.pars.qbinom.binom1}.]Let $\left[  n\right]  $
denote the set $\left\{  1,2,\ldots,n\right\}  $. We have%
\begin{align*}
&  \left(  aq^{0}+b\right)  \left(  aq^{1}+b\right)  \cdots\left(
aq^{n-1}+b\right)  =\prod_{i=1}^{n}\left(  aq^{i-1}+b\right) \\
&  =\sum_{S\subseteq\left[  n\right]  }\underbrace{\left(  \prod_{i\in
S}\left(  aq^{i-1}\right)  \right)  }_{=a^{\left\vert S\right\vert }%
\prod_{i\in S}q^{i-1}}\underbrace{\left(  \prod_{i\in\left[  n\right]
\setminus S}b\right)  }_{=b^{\left\vert \left[  n\right]  \setminus
S\right\vert }}\\
&  \ \ \ \ \ \ \ \ \ \ \ \ \ \ \ \ \ \ \ \ \left(  \text{by Lemma
\ref{lem.prodrule.sum-ai-plus-bi}, applied to }L=K\left[  q\right]  \text{,
}a_{i}=aq^{i-1}\text{ and }b_{i}=b\right) \\
&  =\underbrace{\sum_{S\subseteq\left[  n\right]  }}_{=\sum\limits_{k=0}%
^{n}\ \ \sum\limits_{\substack{S\subseteq\left[  n\right]  ;\\\left\vert
S\right\vert =k}}}a^{\left\vert S\right\vert }\underbrace{\left(  \prod_{i\in
S}q^{i-1}\right)  }_{\substack{=q^{\operatorname*{sum}S-\left\vert
S\right\vert }\\\text{(since the sum of the exponents }i-1\\\text{over all
}i\in S\text{ is precisely }\operatorname*{sum}S-\left\vert S\right\vert
\text{)}}}b^{\left\vert \left[  n\right]  \setminus S\right\vert }\\
&  =\sum_{k=0}^{n}\ \ \sum_{\substack{S\subseteq\left[  n\right]
;\\\left\vert S\right\vert =k}}\ \ \underbrace{a^{\left\vert S\right\vert }%
}_{\substack{=a^{k}\\\text{(since }\left\vert S\right\vert =k\text{)}%
}}\ \ \underbrace{q^{\operatorname*{sum}S-\left\vert S\right\vert }%
}_{\substack{=q^{\operatorname*{sum}S-k}\\\text{(since }\left\vert
S\right\vert =k\text{)}}}\ \ \underbrace{b^{\left\vert \left[  n\right]
\setminus S\right\vert }}_{\substack{=b^{n-k}\\\text{(since }S\text{ is a
}k\text{-element}\\\text{subset of the }n\text{-element}\\\text{set }\left[
n\right]  \text{, and thus we}\\\text{have }\left\vert \left[  n\right]
\setminus S\right\vert =n-k\text{)}}}\\
&  =\sum_{k=0}^{n}\ \ \sum_{\substack{S\subseteq\left[  n\right]
;\\\left\vert S\right\vert =k}}a^{k}\underbrace{q^{\operatorname*{sum}S-k}%
}_{\substack{=q^{\operatorname*{sum}S-\left(  1+2+\cdots+k\right)
}q^{1+2+\cdots+\left(  k-1\right)  }\\=q^{\operatorname*{sum}S-\left(
1+2+\cdots+k\right)  }q^{k\left(  k-1\right)  /2}\\\text{(since }%
1+2+\cdots+\left(  k-1\right)  =k\left(  k-1\right)  /2\text{)}}}b^{n-k}\\
&  =\sum_{k=0}^{n}\ \ \sum_{\substack{S\subseteq\left[  n\right]
;\\\left\vert S\right\vert =k}}a^{k}q^{\operatorname*{sum}S-\left(
1+2+\cdots+k\right)  }q^{k\left(  k-1\right)  /2}b^{n-k}\\
&  =\sum_{k=0}^{n}q^{k\left(  k-1\right)  /2}\underbrace{\left(
\sum\limits_{\substack{S\subseteq\left[  n\right]  ;\\\left\vert S\right\vert
=k}}q^{\operatorname*{sum}S-\left(  1+2+\cdots+k\right)  }\right)
}_{\substack{=\dbinom{n}{k}_{q}\\\text{(by Proposition
\ref{prop.pars.qbinom.alt-defs} \textbf{(b)},}\\\text{since }\left[  n\right]
=\left\{  1,2,\ldots,n\right\}  \text{)}}}a^{k}b^{n-k}=\sum_{k=0}%
^{n}q^{k\left(  k-1\right)  /2}\dbinom{n}{k}_{q}a^{k}b^{n-k}.
\end{align*}
This proves Theorem \ref{thm.pars.qbinom.binom1}.
\end{proof}

The 2nd $q$-binomial theorem grows out of
noncommutativity\footnote{Historically, it arose from a study of matrices
\cite{Potter50}.}:

\begin{theorem}
[2nd $q$-binomial theorem, aka Potter's binomial theorem]%
\label{thm.pars.qbinom.binom2}Let $L$ be a commutative ring. Let $\omega\in
L$. Let $A$ be a noncommutative $L$-algebra. Let $a,b\in A$ be such that
$ba=\omega ab$. Let $n\in\mathbb{N}$. Then,%
\[
\left(  a+b\right)  ^{n}=\sum_{k=0}^{n}\dbinom{n}{k}_{\omega}a^{k}b^{n-k}.
\]

\end{theorem}

The condition $ba=\omega ab$ looks somewhat artificial -- do such elements
$a,b$ actually exist in the wild? Indeed they do, as the following examples show:

\begin{example}
Let $L=\mathbb{Z}$ and $\omega=-1$ and $A=\mathbb{Z}^{2\times2}$ (the ring of
$2\times2$-matrices with integer entries). Let%
\[
a=\left(
\begin{array}
[c]{cc}%
0 & 1\\
1 & 0
\end{array}
\right)  \ \ \ \ \ \ \ \ \ \ \text{and}\ \ \ \ \ \ \ \ \ \ b=\left(
\begin{array}
[c]{cc}%
1 & 0\\
0 & -1
\end{array}
\right)  .
\]
It is easy to check that these two matrices satisfy $ba=-ab$, that is,
$ba=\omega ab$. Thus, Theorem \ref{thm.pars.qbinom.binom2} predicts that%
\[
\left(  a+b\right)  ^{n}=\sum_{k=0}^{n}\dbinom{n}{k}_{\omega}a^{k}b^{n-k}.
\]
And this is indeed true (check it for $n=3$).
\end{example}

\begin{example}
Let $L=\mathbb{R}$. Let $A$ be the ring of $\mathbb{R}$-linear operators on
$C^{\infty}\left(  \mathbb{R}\right)  =\left\{  \text{smooth functions from
}\mathbb{R}\text{ to }\mathbb{R}\right\}  $. Let $\omega$ be any real number.

Let $b\in A$ be the differentiation operator (sending each $f\in C^{\infty
}\left(  \mathbb{R}\right)  $ to $f^{\prime}$).

Let $a\in A$ be the operator that substitutes $\omega x$ for $x$ in the
function (in other words, it shrinks the plot of the function by $\omega$ in
the $x$-direction).

Then, you can check that $ba=\omega ab$. (Indeed, $\left(  f\left(  \omega
x\right)  \right)  ^{\prime}=\omega f^{\prime}\left(  \omega x\right)  $.)
\end{example}

The proof of Theorem \ref{thm.pars.qbinom.binom2} is again a homework exercise
(Exercise \ref{exe.pars.qbinom.qbinom} \textbf{(b)}).

The two $q$-binomial theorems are not entirely unrelated: Theorem
\ref{thm.pars.qbinom.binom1} can be obtained from Theorem
\ref{thm.pars.qbinom.binom2}. (See Exercise \ref{exe.pars.qbinom.2to1} for the details.)

\subsubsection{Counting subspaces of vector spaces}

We have introduced the $q$-binomial coefficient $\dbinom{n}{k}_{q}$ as a
generating function for a certain sort of partitions -- i.e., a
\textquotedblleft weighted number\textquotedblright\ of partitions, where each
partition $\lambda$ has weight $q^{\left\vert \lambda\right\vert }$. However,
for certain integers $a$, the number $\dbinom{n}{k}_{a}$ has other
interpretations, too. A particularly striking one can be found when $a$ is the
size of a finite field.

Let us recall a few things about finite fields (see, e.g., \cite[Chapter
5]{23wa} or \cite[Theorem 15.7.3]{Artin}):

\begin{itemize}
\item For any prime power $p^{k}$, there is a finite field of size $p^{k}$; it
is unique up to isomorphism, and is therefore often called the
\textquotedblleft Galois field of size $p^{k}$\textquotedblright\ and denoted
by $\mathbb{F}_{p^{k}}$. The finite fields $\mathbb{F}_{p^{1}}$ are easiest to
construct -- they are just the quotient rings $\mathbb{Z}/p\mathbb{Z}%
=\mathbb{Z}/p$ (that is, the rings of integers modulo $p$). Higher prime
powers are more complicated. For example, the finite field $\mathbb{F}_{p^{2}%
}$ for any odd prime $p$ can be obtained by starting with $\mathbb{Z}/p$ and
adjoining a square root of an element that is not a square. It is not
$\mathbb{Z}/p^{2}$, since $\mathbb{Z}/p^{2}$ is not a field!

\item Linear algebra (i.e., the notions of vector spaces, subspaces, linear
independence, bases, matrices, Gaussian elimination, etc.) can be done over
any field. In fact, many of its concepts can be defined over any commutative
ring, but only over fields do they behave as nicely as they do over the real
numbers. Thus, much of the linear algebra that you have learned over the real
numbers remains valid over any field. (Exceptions are some properties that
rely on positivity or on characteristic $0$.)
\end{itemize}

Thus, it makes sense to talk about finite-dimensional vector spaces over
finite fields. Such spaces are finite as sets, and thus can be viewed as
combinatorial objects. An $n$-dimensional vector space over a finite field $F$
has size $\left\vert F\right\vert ^{n}$.

Now, we might wonder how many $k$-dimensional subspaces such an $n$%
-dimensional vector space has. The answer is given by the following theorem:

\begin{theorem}
\label{thm.pars.qbinom.subsp-count}Let $F$ be a finite field. Let
$n,k\in\mathbb{N}$. Let $V$ be an $n$-dimensional $F$-vector space. Then,%
\[
\dbinom{n}{k}_{\left\vert F\right\vert }=\left(  \text{\# of }%
k\text{-dimensional vector subspaces of }V\right)  .
\]

\end{theorem}

Compare this with the classical fact that if $S$ is an $n$-element set, then%
\[
\dbinom{n}{k}=\left(  \text{\# of }k\text{-element subsets of }S\right)  .
\]
This hints at an analogy between finite sets and finite-dimensional vector
spaces. Such an analogy does indeed exist; the expository paper \cite{Cohn04}
gives a great overview over its reach.

The easiest proof of Theorem \ref{thm.pars.qbinom.subsp-count} uses three
lemmas. The first one is a classical fact from linear algebra, which holds for
any vector space (not necessarily finite-dimensional) over any field (not
necessarily finite):

\begin{lemma}
\label{lem.linalg.lin-ind-via-span}Let $F$ be a field. Let $V$ be an
$F$-vector space. Let $\left(  v_{1},v_{2},\ldots,v_{k}\right)  $ be a
$k$-tuple of vectors in $V$. Then, $\left(  v_{1},v_{2},\ldots,v_{k}\right)  $
is linearly independent if and only if each $i\in\left\{  1,2,\ldots
,k\right\}  $ satisfies $v_{i}\notin\operatorname*{span}\left(  v_{1}%
,v_{2},\ldots,v_{i-1}\right)  $ (where the span $\operatorname*{span}\left(
{}\right)  $ of an empty list is understood to be the set $\left\{
\mathbf{0}\right\}  $ consisting only of the zero vector $\mathbf{0}$). In
other words, $\left(  v_{1},v_{2},\ldots,v_{k}\right)  $ is linearly
independent if and only if we have%
\begin{align*}
v_{1}  &  \notin\underbrace{\operatorname*{span}\left(  {}\right)
}_{=\left\{  \mathbf{0}\right\}  }\ \ \ \ \ \ \ \ \ \ \text{and}\\
v_{2}  &  \notin\operatorname*{span}\left(  v_{1}\right)
\ \ \ \ \ \ \ \ \ \ \text{and}\\
v_{3}  &  \notin\operatorname*{span}\left(  v_{1},v_{2}\right)
\ \ \ \ \ \ \ \ \ \ \text{and}\\
&  \cdots\ \ \ \ \ \ \ \ \ \ \text{and}\\
v_{k}  &  \notin\operatorname*{span}\left(  v_{1},v_{2},\ldots,v_{k-1}\right)
.
\end{align*}

\end{lemma}

\begin{proof}
[Proof of Lemma \ref{lem.linalg.lin-ind-via-span}.]We must prove that $\left(
v_{1},v_{2},\ldots,v_{k}\right)  $ is linearly independent if and only if each
$i\in\left\{  1,2,\ldots,k\right\}  $ satisfies $v_{i}\notin%
\operatorname*{span}\left(  v_{1},v_{2},\ldots,v_{i-1}\right)  $. This is an
\textquotedblleft if and only if\textquotedblright\ statement; we shall prove
its \textquotedblleft only if\textquotedblright\ (i.e., \textquotedblleft%
$\Longrightarrow$\textquotedblright) and \textquotedblleft
if\textquotedblright\ (i.e., \textquotedblleft$\Longleftarrow$%
\textquotedblright) directions separately:

\begin{enumerate}
\item[$\Longrightarrow:$] Assume that the $k$-tuple $\left(  v_{1}%
,v_{2},\ldots,v_{k}\right)  $ is linearly independent. Let $i\in\left\{
1,2,\ldots,k\right\}  $. If we had $v_{i}\in\operatorname*{span}\left(
v_{1},v_{2},\ldots,v_{i-1}\right)  $, then we could write $v_{i}$ in the form
$v_{i}=\alpha_{1}v_{1}+\alpha_{2}v_{2}+\cdots+\alpha_{i-1}v_{i-1}$ for some
coefficients $\alpha_{1},\alpha_{2},\ldots,\alpha_{i-1}\in F$, and therefore
these coefficients $\alpha_{1},\alpha_{2},\ldots,\alpha_{i-1}$ would satisfy
\begin{align*}
&  \underbrace{\alpha_{1}v_{1}+\alpha_{2}v_{2}+\cdots+\alpha_{i-1}v_{i-1}%
}_{=v_{i}}+\underbrace{\left(  -1\right)  v_{i}}_{=-v_{i}}%
+\underbrace{0v_{i+1}+0v_{i+2}+\cdots+0v_{k}}_{=\mathbf{0}}\\
&  =v_{i}+\left(  -v_{i}\right)  +\mathbf{0}=\mathbf{0},
\end{align*}
which would be a nontrivial linear dependence relation between $\left(
v_{1},v_{2},\ldots,v_{k}\right)  $ (nontrivial because $v_{i}$ appears in it
with the nonzero coefficient $-1$); this would contradict the linear
independence of $\left(  v_{1},v_{2},\ldots,v_{k}\right)  $. Hence, we cannot
have $v_{i}\in\operatorname*{span}\left(  v_{1},v_{2},\ldots,v_{i-1}\right)
$. In other words, we have $v_{i}\notin\operatorname*{span}\left(  v_{1}%
,v_{2},\ldots,v_{i-1}\right)  $.

Forget that we fixed $i$. We thus have shown that each $i\in\left\{
1,2,\ldots,k\right\}  $ satisfies $v_{i}\notin\operatorname*{span}\left(
v_{1},v_{2},\ldots,v_{i-1}\right)  $. This proves the \textquotedblleft%
$\Longrightarrow$\textquotedblright\ direction of our claim.

\item[$\Longleftarrow:$] Assume that each $i\in\left\{  1,2,\ldots,k\right\}
$ satisfies $v_{i}\notin\operatorname*{span}\left(  v_{1},v_{2},\ldots
,v_{i-1}\right)  $. We must prove that the $k$-tuple $\left(  v_{1}%
,v_{2},\ldots,v_{k}\right)  $ is linearly independent. Indeed, assume the
contrary. Thus, this $k$-tuple is linearly dependent. In other words, there
exist coefficients $\beta_{1},\beta_{2},\ldots,\beta_{k}\in F$ that satisfy
$\beta_{1}v_{1}+\beta_{2}v_{2}+\cdots+\beta_{k}v_{k}=\mathbf{0}$ and that are
not all zero. Consider these $\beta_{1},\beta_{2},\ldots,\beta_{k}$. At least
one $i\in\left\{  1,2,\ldots,k\right\}  $ satisfies $\beta_{i}\neq0$ (since
the coefficients $\beta_{1},\beta_{2},\ldots,\beta_{k}$ are not all zero).
Pick the \textbf{largest} such $i$. Thus, $\beta_{i}\neq0$ but $\beta
_{i+1}=\beta_{i+2}=\cdots=\beta_{k}=0$. Hence,
\begin{align*}
&  \beta_{1}v_{1}+\beta_{2}v_{2}+\cdots+\beta_{k}v_{k}\\
&  =\left(  \beta_{1}v_{1}+\beta_{2}v_{2}+\cdots+\beta_{i-1}v_{i-1}\right)
+\beta_{i}v_{i}+\underbrace{\left(  \beta_{i+1}v_{i+1}+\beta_{i+2}%
v_{i+2}+\cdots+\beta_{k}v_{k}\right)  }_{\substack{=0v_{i+1}+0v_{i+2}%
+\cdots+0v_{k}\\=\mathbf{0}}}\\
&  =\left(  \beta_{1}v_{1}+\beta_{2}v_{2}+\cdots+\beta_{i-1}v_{i-1}\right)
+\beta_{i}v_{i},
\end{align*}
so that%
\begin{align*}
\beta_{i}v_{i}  &  =\underbrace{\left(  \beta_{1}v_{1}+\beta_{2}v_{2}%
+\cdots+\beta_{k}v_{k}\right)  }_{=\mathbf{0}}-\left(  \beta_{1}v_{1}%
+\beta_{2}v_{2}+\cdots+\beta_{i-1}v_{i-1}\right) \\
&  =-\left(  \beta_{1}v_{1}+\beta_{2}v_{2}+\cdots+\beta_{i-1}v_{i-1}\right) \\
&  =\left(  -\beta_{1}\right)  v_{1}+\left(  -\beta_{2}\right)  v_{2}%
+\cdots+\left(  -\beta_{i-1}\right)  v_{i-1}\\
&  \in\operatorname*{span}\left(  v_{1},v_{2},\ldots,v_{i-1}\right)  .
\end{align*}
Since $\beta_{i}\neq0$, we thus obtain $v_{i}\in\operatorname*{span}\left(
v_{1},v_{2},\ldots,v_{i-1}\right)  $ (since the set $\operatorname*{span}%
\left(  v_{1},v_{2},\ldots,v_{i-1}\right)  $ is an $F$-vector subspace of $V$
and thus preserved under scaling). This contradicts our assumption that
$v_{i}\notin\operatorname*{span}\left(  v_{1},v_{2},\ldots,v_{i-1}\right)  $.
This contradiction shows that our assumption was wrong, and thus completes our
proof of the \textquotedblleft$\Longleftarrow$\textquotedblright\ direction of
our claim.
\end{enumerate}

Thus, both directions of our claim are proved. This concludes the proof of
Lemma \ref{lem.linalg.lin-ind-via-span}.
\end{proof}

The next lemma we are going to use is itself an answer to a rather natural
counting problem. Indeed, it is well-known that (see, e.g., \cite[Proposition
2.7.2]{19fco}) if $X$ is an $n$-element set, and if $k\in\mathbb{N}$, then
the
\begin{align}
&  \left(  \text{\# of }k\text{-tuples of distinct elements of }X\right)
\nonumber\\
&  =n\left(  n-1\right)  \left(  n-2\right)  \cdots\left(  n-k+1\right)
=\prod_{i=0}^{k-1}\left(  n-i\right)  . \label{eq.count.k-tups-dist-elts}%
\end{align}
The following lemma is a \textquotedblleft linear analogue\textquotedblright%
\ of this combinatorial fact: The $n$-element set $X$ is replaced by an
$n$-dimensional vector space $V$, and \textquotedblleft distinct
elements\textquotedblright\ are replaced by \textquotedblleft linearly
independent vectors\textquotedblright. The answer is rather similar:

\begin{lemma}
\label{lem.pars.qbinom.lin-ind-count}Let $F$ be a finite field. Let
$n,k\in\mathbb{N}$. Let $V$ be an $n$-dimensional $F$-vector space. Then,%
\begin{align*}
&  \left(  \text{\# of linearly independent }k\text{-tuples of vectors in
}V\right) \\
&  =\left(  \left\vert F\right\vert ^{n}-\left\vert F\right\vert ^{0}\right)
\left(  \left\vert F\right\vert ^{n}-\left\vert F\right\vert ^{1}\right)
\cdots\left(  \left\vert F\right\vert ^{n}-\left\vert F\right\vert
^{k-1}\right)  =\prod_{i=0}^{k-1}\left(  \left\vert F\right\vert
^{n}-\left\vert F\right\vert ^{i}\right)  .
\end{align*}

\end{lemma}

\begin{proof}
[Proof of Lemma \ref{lem.pars.qbinom.lin-ind-count}.]We have $\left\vert
V\right\vert =\left\vert F\right\vert ^{n}$ (since $V$ is an $n$-dimensional
$F$-vector space).

Lemma \ref{lem.linalg.lin-ind-via-span} says that a $k$-tuple $\left(
v_{1},v_{2},\ldots,v_{k}\right)  $ of vectors in $V$ is linearly independent
if and only if it satisfies%
\begin{align*}
v_{1}  &  \notin\underbrace{\operatorname*{span}\left(  {}\right)
}_{=\left\{  \mathbf{0}\right\}  }\ \ \ \ \ \ \ \ \ \ \text{and}\\
v_{2}  &  \notin\operatorname*{span}\left(  v_{1}\right)
\ \ \ \ \ \ \ \ \ \ \text{and}\\
v_{3}  &  \notin\operatorname*{span}\left(  v_{1},v_{2}\right)
\ \ \ \ \ \ \ \ \ \ \text{and}\\
&  \cdots\ \ \ \ \ \ \ \ \ \ \text{and}\\
v_{k}  &  \notin\operatorname*{span}\left(  v_{1},v_{2},\ldots,v_{k-1}\right)
.
\end{align*}

Thus, we can construct a linearly independent $k$-tuple $\left(  v_{1}%
,v_{2},\ldots,v_{k}\right)  $ of vectors in $V$ as follows, proceeding entry
by entry:

\begin{itemize}
\item First, we choose $v_{1}$. This has to be a vector in $V\setminus
\underbrace{\operatorname*{span}\left(  {}\right)  }_{=\left\{  \mathbf{0}%
\right\}  }$ (because it has to satisfy $v_{1}\notin\operatorname*{span}%
\left(  {}\right)  $); thus, there are $\left\vert V\setminus
\operatorname*{span}\left(  {}\right)  \right\vert =\left\vert V\right\vert
-\underbrace{\left\vert \operatorname*{span}\left(  {}\right)  \right\vert
}_{=1}=\left\vert V\right\vert -1$ options for it.

Once $v_{1}$ has been chosen, we have obtained a linearly independent
singleton list $\left(  v_{1}\right)  $. Hence, its span $\operatorname*{span}%
\left(  v_{1}\right)  $ has dimension $1$ (as an $F$-vector space) and thus
size $\left\vert F\right\vert ^{1}$. In other words, $\left\vert
\operatorname*{span}\left(  v_{1}\right)  \right\vert =\left\vert F\right\vert
^{1}$.

\item Next, we choose $v_{2}$. This has to be a vector in $V\setminus
\operatorname*{span}\left(  v_{1}\right)  $ (because it has to satisfy
$v_{2}\notin\operatorname*{span}\left(  v_{1}\right)  $); thus, there are
$\left\vert V\right\vert -\underbrace{\left\vert \operatorname*{span}\left(
v_{1}\right)  \right\vert }_{=\left\vert F\right\vert ^{1}}=\left\vert
V\right\vert -\left\vert F\right\vert ^{1}$ options for it.

Once $v_{2}$ has been chosen, we have obtained a linearly independent list
$\left(  v_{1},v_{2}\right)  $. Hence, its span $\operatorname*{span}\left(
v_{1},v_{2}\right)  $ has dimension $2$ (as an $F$-vector space) and thus size
$\left\vert F\right\vert ^{2}$. In other words, $\left\vert
\operatorname*{span}\left(  v_{1},v_{2}\right)  \right\vert =\left\vert
F\right\vert ^{2}$.

\item Next, we choose $v_{3}$. This has to be a vector in $V\setminus
\operatorname*{span}\left(  v_{1},v_{2}\right)  $ (because it has to satisfy
$v_{3}\notin\operatorname*{span}\left(  v_{1},v_{2}\right)  $); thus, there
are $\left\vert V\right\vert -\underbrace{\left\vert \operatorname*{span}%
\left(  v_{1},v_{2}\right)  \right\vert }_{=\left\vert F\right\vert ^{2}%
}=\left\vert V\right\vert -\left\vert F\right\vert ^{2}$ options for it.

Once $v_{3}$ has been chosen, we have obtained a linearly independent list
$\left(  v_{1},v_{2},v_{3}\right)  $. Hence, its span $\operatorname*{span}%
\left(  v_{1},v_{2},v_{3}\right)  $ has dimension $3$ (as an $F$-vector space)
and thus size $\left\vert F\right\vert ^{3}$. In other words, $\left\vert
\operatorname*{span}\left(  v_{1},v_{2},v_{3}\right)  \right\vert =\left\vert
F\right\vert ^{3}$.

\item And so on, until the last vector $v_{k}$ in our list has been chosen.
\end{itemize}

The total \# of ways to perform this construction is%
\[
\left(  \left\vert V\right\vert -1\right)  \left(  \left\vert V\right\vert
-\left\vert F\right\vert ^{1}\right)  \left(  \left\vert V\right\vert
-\left\vert F\right\vert ^{2}\right)  \cdots\left(  \left\vert V\right\vert
-\left\vert F\right\vert ^{k-1}\right)  .
\]

Hence,
\begin{align*}
&  \left(  \text{\# of linearly independent }k\text{-tuples of vectors in
}V\right) \\
&  =\left(  \left\vert V\right\vert -1\right)  \left(  \left\vert V\right\vert
-\left\vert F\right\vert ^{1}\right)  \left(  \left\vert V\right\vert
-\left\vert F\right\vert ^{2}\right)  \cdots\left(  \left\vert V\right\vert
-\left\vert F\right\vert ^{k-1}\right) \\
&  =\prod_{i=0}^{k-1}\left(  \left\vert V\right\vert -\left\vert F\right\vert
^{i}\right)  =\prod_{i=0}^{k-1}\left(  \left\vert F\right\vert ^{n}-\left\vert
F\right\vert ^{i}\right)
\end{align*}
(since $\left\vert V\right\vert =\left\vert F\right\vert ^{n}$). This proves
Lemma \ref{lem.pars.qbinom.lin-ind-count}.
\end{proof}

Another lemma we will need is a basic combinatorial principle (often
illustrated by the saying \textquotedblleft to count a flock of sheep, count
the legs and divide by $4$\textquotedblright):

\begin{lemma}
[Multijection principle]\label{lem.count.multijection}Let $A$ and $B$ be two
finite sets. Let $m\in\mathbb{N}$. Let $f:A\rightarrow B$ be any map. Assume
that each $b\in B$ has exactly $m$ preimages under $f$ (that is, for each
$b\in B$, there are exactly $m$ many elements $a\in A$ such that $f\left(
a\right)  =b$). Then,%
\[
\left\vert A\right\vert =m\cdot\left\vert B\right\vert .
\]

\end{lemma}

\begin{proof}
Easy and LTTR.
\end{proof}

We note that a map $f:A\rightarrow B$ satisfying the assumption of Lemma
\ref{lem.count.multijection} is often called an $m$\emph{-to-}$1$\emph{ map}.

\begin{proof}
[Proof of Theorem \ref{thm.pars.qbinom.subsp-count}.]First of all, we notice
that if $k>n$, then $\dbinom{n}{k}_{\left\vert F\right\vert }=0$ (by
Proposition \ref{prop.pars.qbinom.0}) and $\left(  \text{\# of }%
k\text{-dimensional vector subspaces of }V\right)  =0$ (since the dimension of
a subspace of $V$ is never larger than the dimension of $V$). Thus, Theorem
\ref{thm.pars.qbinom.subsp-count} is true when $k>n$. Hence, for the rest of
this proof, we WLOG assume that $k\leq n$.

We will use the shorthand \textquotedblleft linind\textquotedblright\ for the
words \textquotedblleft linearly independent\textquotedblright.

If $\left(  v_{1},v_{2},\ldots,v_{k}\right)  $ is a linind $k$-tuple of
vectors in $V$, then $\operatorname*{span}\left(  v_{1},v_{2},\ldots
,v_{k}\right)  $ is a $k$-dimensional vector subspace of $V$. Hence, we can
define a map%
\begin{align*}
f:\left\{  \text{linind }k\text{-tuples of vectors in }V\right\}   &
\rightarrow\left\{  k\text{-dimensional vector subspaces of }V\right\}  ,\\
\left(  v_{1},v_{2},\ldots,v_{k}\right)   &  \mapsto\operatorname*{span}%
\left(  v_{1},v_{2},\ldots,v_{k}\right)  .
\end{align*}
Consider this map $f$. We claim the following:

\begin{statement}
\textit{Observation 1:} Each $k$-dimensional vector subspace of $V$ has
exactly%
\[
\left(  \left\vert F\right\vert ^{k}-\left\vert F\right\vert ^{0}\right)
\left(  \left\vert F\right\vert ^{k}-\left\vert F\right\vert ^{1}\right)
\cdots\left(  \left\vert F\right\vert ^{k}-\left\vert F\right\vert
^{k-1}\right)
\]
preimages under $f$.
\end{statement}

[\textit{Proof of Observation 1:} Let $W$ be a $k$-dimensional vector subspace
of $V$. We must prove that
\[
\left(  \text{\# of preimages of }W\text{ under }f\right)  =\left(  \left\vert
F\right\vert ^{k}-\left\vert F\right\vert ^{0}\right)  \left(  \left\vert
F\right\vert ^{k}-\left\vert F\right\vert ^{1}\right)  \cdots\left(
\left\vert F\right\vert ^{k}-\left\vert F\right\vert ^{k-1}\right)  .
\]

A preimage of $W$ under $f$ is a $k$-tuple of vectors in $V$ that spans $W$
(by the very definition of $f$). Obviously, all the $k$ vectors in such a
$k$-tuple must belong to $W$. Thus, a preimage of $W$ under $f$ is a $k$-tuple
of vectors in $W$ that spans $W$. However, the vector space $W$ is
$k$-dimensional. Thus, a $k$-tuple of vectors in $W$ spans $W$ if and only if
this $k$-tuple is linind\footnote{Again, we are using a simple fact from
linear algebra here, which is true over any field (not necessarily finite): A
$k$-tuple of vectors in a $k$-dimensional vector space spans the space if and
only if it is linind.}. Therefore, a preimage of $W$ under $f$ is a linind
$k$-tuple of vectors in $W$. Hence,%
\begin{align*}
&  \left(  \text{\# of preimages of }W\text{ under }f\right) \\
&  =\left(  \text{\# of linearly independent }k\text{-tuples of vectors in
}W\right) \\
&  =\left(  \left\vert F\right\vert ^{k}-\left\vert F\right\vert ^{0}\right)
\left(  \left\vert F\right\vert ^{k}-\left\vert F\right\vert ^{1}\right)
\cdots\left(  \left\vert F\right\vert ^{k}-\left\vert F\right\vert
^{k-1}\right)
\end{align*}
(by Lemma \ref{lem.pars.qbinom.lin-ind-count}, applied to $W$ and $k$ instead
of $V$ and $n$). This proves Observation 1.]

Now, Observation 1 shows that each $k$-dimensional vector subspace of $V$ has
exactly $\left(  \left\vert F\right\vert ^{k}-\left\vert F\right\vert
^{0}\right)  \left(  \left\vert F\right\vert ^{k}-\left\vert F\right\vert
^{1}\right)  \cdots\left(  \left\vert F\right\vert ^{k}-\left\vert
F\right\vert ^{k-1}\right)  $ preimages under $f$. Hence, Lemma
\ref{lem.count.multijection} (applied to $A=\left\{  \text{linind
}k\text{-tuples of vectors in }V\right\}  $ and $B=\left\{
k\text{-dimensional vector subspaces of }V\right\}  $ and \newline$m=\left(
\left\vert F\right\vert ^{k}-\left\vert F\right\vert ^{0}\right)  \left(
\left\vert F\right\vert ^{k}-\left\vert F\right\vert ^{1}\right)
\cdots\left(  \left\vert F\right\vert ^{k}-\left\vert F\right\vert
^{k-1}\right)  $) shows that%
\begin{align*}
&  \left(  \text{\# of linind }k\text{-tuples of vectors in }V\right) \\
&  =\left(  \left\vert F\right\vert ^{k}-\left\vert F\right\vert ^{0}\right)
\left(  \left\vert F\right\vert ^{k}-\left\vert F\right\vert ^{1}\right)
\cdots\left(  \left\vert F\right\vert ^{k}-\left\vert F\right\vert
^{k-1}\right) \\
&  \ \ \ \ \ \ \ \ \ \ \cdot\left(  \text{\# of }k\text{-dimensional vector
subspaces of }V\right)  .
\end{align*}
However, Lemma \ref{lem.pars.qbinom.lin-ind-count} yields%
\begin{align*}
&  \left(  \text{\# of linind }k\text{-tuples of vectors in }V\right) \\
&  =\left(  \left\vert F\right\vert ^{n}-\left\vert F\right\vert ^{0}\right)
\left(  \left\vert F\right\vert ^{n}-\left\vert F\right\vert ^{1}\right)
\cdots\left(  \left\vert F\right\vert ^{n}-\left\vert F\right\vert
^{k-1}\right)  .
\end{align*}
Comparing these two equalities, we obtain%
\begin{align*}
&  \left(  \left\vert F\right\vert ^{k}-\left\vert F\right\vert ^{0}\right)
\left(  \left\vert F\right\vert ^{k}-\left\vert F\right\vert ^{1}\right)
\cdots\left(  \left\vert F\right\vert ^{k}-\left\vert F\right\vert
^{k-1}\right) \\
&  \ \ \ \ \ \ \ \ \ \ \cdot\left(  \text{\# of }k\text{-dimensional vector
subspaces of }V\right) \\
&  =\left(  \left\vert F\right\vert ^{n}-\left\vert F\right\vert ^{0}\right)
\left(  \left\vert F\right\vert ^{n}-\left\vert F\right\vert ^{1}\right)
\cdots\left(  \left\vert F\right\vert ^{n}-\left\vert F\right\vert
^{k-1}\right)  .
\end{align*}
Therefore,%
\begin{align*}
&  \left(  \text{\# of }k\text{-dimensional vector subspaces of }V\right) \\
&  =\dfrac{\left(  \left\vert F\right\vert ^{n}-\left\vert F\right\vert
^{0}\right)  \left(  \left\vert F\right\vert ^{n}-\left\vert F\right\vert
^{1}\right)  \cdots\left(  \left\vert F\right\vert ^{n}-\left\vert
F\right\vert ^{k-1}\right)  }{\left(  \left\vert F\right\vert ^{k}-\left\vert
F\right\vert ^{0}\right)  \left(  \left\vert F\right\vert ^{k}-\left\vert
F\right\vert ^{1}\right)  \cdots\left(  \left\vert F\right\vert ^{k}%
-\left\vert F\right\vert ^{k-1}\right)  }\\
&  =\dfrac{\prod_{i=0}^{k-1}\left(  \left\vert F\right\vert ^{n}-\left\vert
F\right\vert ^{i}\right)  }{\prod_{i=0}^{k-1}\left(  \left\vert F\right\vert
^{k}-\left\vert F\right\vert ^{i}\right)  }=\prod_{i=0}^{k-1}%
\underbrace{\dfrac{\left\vert F\right\vert ^{n}-\left\vert F\right\vert ^{i}%
}{\left\vert F\right\vert ^{k}-\left\vert F\right\vert ^{i}}}%
_{\substack{=\dfrac{\left\vert F\right\vert ^{i}\left(  \left\vert
F\right\vert ^{n-i}-1\right)  }{\left\vert F\right\vert ^{i}\left(  \left\vert
F\right\vert ^{k-i}-1\right)  }=\dfrac{\left\vert F\right\vert ^{n-i}%
-1}{\left\vert F\right\vert ^{k-i}-1}\\=\dfrac{1-\left\vert F\right\vert
^{n-i}}{1-\left\vert F\right\vert ^{k-i}}}}=\prod_{i=0}^{k-1}\dfrac
{1-\left\vert F\right\vert ^{n-i}}{1-\left\vert F\right\vert ^{k-i}}\\
&  =\dfrac{\prod_{i=0}^{k-1}\left(  1-\left\vert F\right\vert ^{n-i}\right)
}{\prod_{i=0}^{k-1}\left(  1-\left\vert F\right\vert ^{k-i}\right)  }%
=\dfrac{\left(  1-\left\vert F\right\vert ^{n}\right)  \left(  1-\left\vert
F\right\vert ^{n-1}\right)  \cdots\left(  1-\left\vert F\right\vert
^{n-k+1}\right)  }{\left(  1-\left\vert F\right\vert ^{k}\right)  \left(
1-\left\vert F\right\vert ^{k-1}\right)  \cdots\left(  1-\left\vert
F\right\vert ^{1}\right)  }\\
&  =\dbinom{n}{k}_{\left\vert F\right\vert }%
\end{align*}
(since substituting $\left\vert F\right\vert $ for $q$ in Theorem
\ref{thm.pars.qbinom.quot1} \textbf{(b)} yields \newline$\dbinom{n}%
{k}_{\left\vert F\right\vert }=\dfrac{\left(  1-\left\vert F\right\vert
^{n}\right)  \left(  1-\left\vert F\right\vert ^{n-1}\right)  \cdots\left(
1-\left\vert F\right\vert ^{n-k+1}\right)  }{\left(  1-\left\vert F\right\vert
^{k}\right)  \left(  1-\left\vert F\right\vert ^{k-1}\right)  \cdots\left(
1-\left\vert F\right\vert ^{1}\right)  }$). This proves Theorem
\ref{thm.pars.qbinom.subsp-count}.
\end{proof}

\subsubsection{Limits of $q$-binomial coefficients}

There is much more to say about $q$-binomial coefficients, but let us just
briefly focus on their limiting behavior. This is not analogous to anything
known from usual binomial coefficients; indeed, the limit $\lim
\limits_{n\rightarrow\infty}\dbinom{n}{k}$ does not exist for any positive
integer $k$. However, $q$-binomial coefficients behave much better in this regard.

Indeed, consider the $q$-binomial coefficients $\dbinom{n}{2}_{q}$ for various
values of $n$:%
\begin{align*}
\dbinom{0}{2}_{q}  &  =0,\\
\dbinom{1}{2}_{q}  &  =0,\\
\dbinom{2}{2}_{q}  &  =1,\\
\dbinom{3}{2}_{q}  &  =1+q+q^{2},\\
\dbinom{4}{2}_{q}  &  =1+q+2q^{2}+q^{3}+q^{4},\\
\dbinom{5}{2}_{q}  &  =1+q+2q^{2}+2q^{3}+2q^{4}+q^{5}+q^{6},\\
\dbinom{6}{2}_{q}  &  =1+q+2q^{2}+2q^{3}+3q^{4}+2q^{5}+2q^{6}+q^{7}+q^{8}.
\end{align*}
It appears from these examples that the sequence $\left(  \dbinom{n}{2}%
_{q}\right)  _{n\in\mathbb{N}}$ coefficientwise stabilizes\footnote{Recall
Definition \ref{def.fps.lim.coeff-stab} for the notion of \textquotedblleft
coefficientwise stabilizing\textquotedblright.} to%
\[
1+q+2q^{2}+2q^{3}+3q^{4}+3q^{5}+\cdots=\sum_{n\in\mathbb{N}}\left(
1+\left\lfloor \dfrac{n}{2}\right\rfloor \right)  q^{n}.
\]
And this is indeed the case:

\begin{proposition}
\label{prop.pars.qbinom.lim1}Let $k\in\mathbb{N}$ be fixed. Then,%
\[
\lim\limits_{n\rightarrow\infty}\dbinom{n}{k}_{q}=\sum_{n\in\mathbb{N}}\left(
p_{0}\left(  n\right)  +p_{1}\left(  n\right)  +\cdots+p_{k}\left(  n\right)
\right)  q^{n}=\prod_{i=1}^{k}\dfrac{1}{1-q^{i}}.
\]
(See Definition \ref{def.pars.pn-pkn} \textbf{(a)} for the meaning of
$p_{i}\left(  n\right)  $.)
\end{proposition}

\begin{proof}
[First proof of Proposition \ref{prop.pars.qbinom.lim1} (sketched).]For each
integer $n\geq k$, we have%
\begin{align*}
\dbinom{n}{k}_{q}  &  =\dfrac{\left(  1-q^{n}\right)  \left(  1-q^{n-1}%
\right)  \cdots\left(  1-q^{n-k+1}\right)  }{\left(  1-q^{k}\right)  \left(
1-q^{k-1}\right)  \cdots\left(  1-q^{1}\right)  }\ \ \ \ \ \ \ \ \ \ \left(
\text{by Theorem \ref{thm.pars.qbinom.quot1} \textbf{(b)}}\right) \\
&  =\dfrac{\left(  1-q^{n-k+1}\right)  \left(  1-q^{n-k+2}\right)
\cdots\left(  1-q^{n}\right)  }{\left(  1-q^{1}\right)  \left(  1-q^{2}%
\right)  \cdots\left(  1-q^{k}\right)  }\\
&  \ \ \ \ \ \ \ \ \ \ \ \ \ \ \ \ \ \ \ \ \left(  \text{here, we have turned
both products upside down}\right) \\
&  =\dfrac{\prod_{i=1}^{k}\left(  1-q^{n-k+i}\right)  }{\prod_{i=1}^{k}\left(
1-q^{i}\right)  }=\dfrac{1}{\prod_{i=1}^{k}\left(  1-q^{i}\right)  }\cdot
\prod_{i=1}^{k}\left(  1-q^{n-k+i}\right)  .
\end{align*}

However, we have $\lim\limits_{n\rightarrow\infty}q^{n}=0$ (check
this!\footnote{This is a matter of understanding Definition
\ref{def.fps.lim.coeff-stab}.}). Thus, for each $i\in\left\{  1,2,\ldots
,k\right\}  $, we have%
\[
\lim\limits_{n\rightarrow\infty}q^{n-k+i}=0
\]
(since the family $\left(  q^{n-k+i}\right)  _{n\geq k-i}$ is just a
reindexing of the family $\left(  q^{n}\right)  _{n\geq0}$), and therefore%
\[
\lim\limits_{n\rightarrow\infty}\left(  1-q^{n-k+i}\right)  =1.
\]
Hence, Corollary \ref{cor.fps.lim.sum-prod-k} (applied to $f_{i,n}%
=1-q^{n-k+i}$ and $f_{i}=1$) yields that%
\[
\lim\limits_{n\rightarrow\infty}\sum_{i=1}^{k}\left(  1-q^{n-k+i}\right)
=\sum_{i=1}^{k}1\ \ \ \ \ \ \ \ \ \ \text{and}\ \ \ \ \ \ \ \ \ \ \lim
\limits_{n\rightarrow\infty}\prod_{i=1}^{k}\left(  1-q^{n-k+i}\right)
=\prod_{i=1}^{k}1.
\]
Thus, in particular,%
\[
\lim\limits_{n\rightarrow\infty}\prod_{i=1}^{k}\left(  1-q^{n-k+i}\right)
=\prod_{i=1}^{k}1=1.
\]

Now, recall that each integer $n\geq k$ satisfies%
\[
\dbinom{n}{k}_{q}=\dfrac{1}{\prod_{i=1}^{k}\left(  1-q^{i}\right)  }\cdot
\prod_{i=1}^{k}\left(  1-q^{n-k+i}\right)  .
\]
Hence,%
\begin{align}
\lim\limits_{n\rightarrow\infty}\dbinom{n}{k}_{q}  &  =\lim
\limits_{n\rightarrow\infty}\left(  \dfrac{1}{\prod_{i=1}^{k}\left(
1-q^{i}\right)  }\cdot\prod_{i=1}^{k}\left(  1-q^{n-k+i}\right)  \right)
\nonumber\\
&  =\dfrac{1}{\prod_{i=1}^{k}\left(  1-q^{i}\right)  }\cdot\underbrace{\lim
\limits_{n\rightarrow\infty}\prod_{i=1}^{k}\left(  1-q^{n-k+i}\right)  }%
_{=1}=\dfrac{1}{\prod_{i=1}^{k}\left(  1-q^{i}\right)  }\nonumber\\
&  =\prod_{i=1}^{k}\dfrac{1}{1-q^{i}}. \label{pf.prop.pars.qbinom.lim1.rhs1}%
\end{align}

Finally, Theorem \ref{thm.pars.main-gf-0n} (with the letters $x$, $m$ and $k$
renamed as $q$, $k$ and $i$) says that%
\[
\sum_{n\in\mathbb{N}}\left(  p_{0}\left(  n\right)  +p_{1}\left(  n\right)
+\cdots+p_{k}\left(  n\right)  \right)  q^{n}=\prod_{i=1}^{k}\dfrac{1}%
{1-q^{i}}.
\]
Combining this with (\ref{pf.prop.pars.qbinom.lim1.rhs1}), we obtain%
\[
\lim\limits_{n\rightarrow\infty}\dbinom{n}{k}_{q}=\sum_{n\in\mathbb{N}}\left(
p_{0}\left(  n\right)  +p_{1}\left(  n\right)  +\cdots+p_{k}\left(  n\right)
\right)  q^{n}=\prod_{i=1}^{k}\dfrac{1}{1-q^{i}}.
\]
This proves Proposition \ref{prop.pars.qbinom.lim1}.
\end{proof}

\begin{proof}
[Second proof of Proposition \ref{prop.pars.qbinom.lim1} (sketched).]For each
$n\in\mathbb{N}$, we have%
\begin{equation}
\dbinom{n}{k}_{q}=\sum_{\substack{\lambda\text{ is a partition}\\\text{with
largest part }\leq n-k\\\text{and length }\leq k}}q^{\left\vert \lambda
\right\vert } \label{pf.prop.pars.qbinom.lim1.2nd.1}%
\end{equation}
(by the definition of $\dbinom{n}{k}_{q}$).

However, for each $n\in\mathbb{N}$, the sum $\sum_{\substack{\lambda\text{ is
a partition}\\\text{with largest part }\leq n-k\\\text{and length }\leq
k}}q^{\left\vert \lambda\right\vert }$ is a partial sum of the sum
$\sum_{\substack{\lambda\text{ is a partition}\\\text{with length }\leq
k}}q^{\left\vert \lambda\right\vert }$, and this partial sum grows by more and
more addends as $n$ increases; each addend of the sum $\sum_{\substack{\lambda
\text{ is a partition}\\\text{with length }\leq k}}q^{\left\vert
\lambda\right\vert }$ gets eventually included in this partial sum (for
sufficiently large $n$). From these observations, it is easy to obtain that%
\[
\lim\limits_{n\rightarrow\infty}\sum_{\substack{\lambda\text{ is a
partition}\\\text{with largest part }\leq n-k\\\text{and length }\leq
k}}q^{\left\vert \lambda\right\vert }=\sum_{\substack{\lambda\text{ is a
partition}\\\text{with length }\leq k}}q^{\left\vert \lambda\right\vert }.
\]
In view of (\ref{pf.prop.pars.qbinom.lim1.2nd.1}), this rewrites as%
\begin{align*}
\lim\limits_{n\rightarrow\infty}\dbinom{n}{k}_{q}  &  =\sum_{\substack{\lambda
\text{ is a partition}\\\text{with length }\leq k}}q^{\left\vert
\lambda\right\vert }=\sum_{n\in\mathbb{N}}\underbrace{\left(  \text{\# of
partitions of }n\text{ having length }\leq k\right)  }_{\substack{=p_{0}%
\left(  n\right)  +p_{1}\left(  n\right)  +\cdots+p_{k}\left(  n\right)
\\\text{(by Definition \ref{def.pars.pn-pkn} \textbf{(a)})}}}q^{n}\\
&  =\sum_{n\in\mathbb{N}}\left(  p_{0}\left(  n\right)  +p_{1}\left(
n\right)  +\cdots+p_{k}\left(  n\right)  \right)  q^{n}=\prod_{i=1}^{k}%
\dfrac{1}{1-q^{i}}%
\end{align*}
(by Theorem \ref{thm.pars.main-gf-0n}, with the letters $x$, $m$ and $k$
renamed as $q$, $k$ and $i$). Thus, Proposition \ref{prop.pars.qbinom.lim1} is
proved again.
\end{proof}

\subsection{\label{sec.par.refs}References}

Thus ends our foray into integer partitions and related FPSs. We will
partially revisit this topic later, as we discuss symmetric functions. Here
are just a few things we are omitting:

\begin{itemize}
\item In 1919, Ramanujan discovered the following three congruences for
$p\left(  n\right)  $:%
\begin{align*}
p\left(  n\right)   &  \equiv0\operatorname{mod}5\ \ \ \ \ \ \ \ \ \ \text{if
}n\equiv4\operatorname{mod}5;\\
p\left(  n\right)   &  \equiv0\operatorname{mod}7\ \ \ \ \ \ \ \ \ \ \text{if
}n\equiv5\operatorname{mod}7;\\
p\left(  n\right)   &  \equiv0\operatorname{mod}11\ \ \ \ \ \ \ \ \ \ \text{if
}n\equiv6\operatorname{mod}11.
\end{align*}
The first of these follows from the FPS equality%
\[
\sum_{n\in\mathbb{N}}p\left(  5n+4\right)  x^{n}=5\prod_{i=1}^{\infty}%
\dfrac{\left(  1-x^{5i}\right)  ^{5}}{\left(  1-x^{i}\right)  ^{6}},
\]
whose proof is far from straightforward. All of these results (and some rather
subtle generalizations) are shown in \cite[Chapter 2]{Berndt06} and
\cite[Chapters 3 and 5]{Hirsch17}; see also \cite[Chapter 3, Highlight]%
{Aigner07} for a proof of the latter equality.

\item An asymptotic expansion for $p\left(  n\right)  $ (found by Hardy and
Ramanujan in 1918) is
\begin{equation}
p\left(  n\right)  \sim\dfrac{1}{4n\sqrt{3}}\exp\left(  \pi\sqrt{\dfrac{2n}%
{3}}\right)  \ \ \ \ \ \ \ \ \ \ \text{as }n\rightarrow\infty.
\label{eq.pars.asymptotic-pn}%
\end{equation}
See \cite{Erdos42} or \cite[Proposition VIII.6]{FlaSed09} for a proof.

\item In 1770, Lagrange proved that every nonnegative integer $n$ can be
written as a sum of four perfect squares. In 1829, Jacobi strengthened this to
a counting formula: If $n$ is a positive integer, then the number of
quadruples $\left(  a,b,c,d\right)  $ of integers satisfying $n=a^{2}%
+b^{2}+c^{2}+d^{2}$ is $8$ times the sum of positive divisors of $n$ that are
not divisible by $4$. The most elementary proofs of this striking result use
partition-related FPSs and the Jacobi Triple Product Identity. (See
\cite{Hirsch87} or \cite[Theorem 7.20]{Sambal22} for self-contained proofs;
see also \cite[Chapter 2]{Hirsch17} and \cite[Chapter 3]{Berndt06} for various
related results.)

\item The
\href{https://en.wikipedia.org/wiki/Rogers-Ramanujan_identities}{Rogers--Ramanujan
identities}%
\begin{align*}
\sum_{k\in\mathbb{N}}\dfrac{x^{k^{2}}}{\left(  1-x^{1}\right)  \left(
1-x^{2}\right)  \cdots\left(  1-x^{k}\right)  }  &  =\prod_{i\in\mathbb{N}%
}\dfrac{1}{\left(  1-x^{5i+1}\right)  \left(  1-x^{5i+4}\right)
}\ \ \ \ \ \ \ \ \ \ \text{and}\\
\sum_{k\in\mathbb{N}}\dfrac{x^{k\left(  k+1\right)  }}{\left(  1-x^{1}\right)
\left(  1-x^{2}\right)  \cdots\left(  1-x^{k}\right)  }  &  =\prod
_{i\in\mathbb{N}}\dfrac{1}{\left(  1-x^{5i+2}\right)  \left(  1-x^{5i+3}%
\right)  }%
\end{align*}
can be used to count partitions into parts that are congruent to
$\pm1\operatorname{mod}5$ or congruent to $\pm2\operatorname{mod}5$,
respectively. These surprising identities can be proved using Proposition
\ref{prop.pars.qbinom.lim1} and the Jacobi Triple Product Identity; see
\cite{Doyle19} for a self-contained writeup of this proof. A whole book
\cite{Sills18} has been written about these two identities and their many variants.
\end{itemize}

Here is a list of references for further reading on partitions:

\begin{itemize}
\item The book \cite{AndEri04} by Andrews and Eriksson is a beautiful (if not
always fully precise) introduction to integer partitions and related topics.
Andrews's monograph \cite{Andrew76} goes deeper.

\item Pak's \cite{Pak06} is a survey of identities between partition numbers
(and related FPSs) with occasionally outlined proofs. (Beware: the writing is
very terse and teems with typos.)

\item Hirschhorn's \cite{Hirsch17} (subtitled \textquotedblleft a personal
journey\textquotedblright, not meant to be comprehensive) studies partitions
through the lens of (mostly purely algebraic) manipulation of FPSs.

\item Berndt's \cite{Berndt06} is another (more analytic and
number-theoretical) study of partition-related FPSs, with applications to
number theory.
\end{itemize}

\begin{noncompile}
also use \cite{Ness61} (e.g. \S 9)
\end{noncompile}

\newpage

\section{\label{chap.perm}Permutations}

We now come back to the foundations of combinatorics: We will study
permutations of finite sets. I will assume that you know their most basic
properties (see, e.g., \cite[Appendix B]{Strick13} and \cite[\S 1.5]{Goodman}
for refreshers; see also \cite[Chapter 5]{detnotes} for many more details on
inversions), and will show some more advanced results. For deeper treatments,
see \cite{Bona22}, \cite{Sagan01} and \cite[Chapter 1]{Stanley-EC1}.

\subsection{Basic definitions}

\begin{definition}
\label{def.perm.perm}Let $X$ be a set. \medskip

\textbf{(a)} A \emph{permutation} of $X$ means a bijection from $X$ to $X$.
\medskip

\textbf{(b)} It is known that the set of all permutations of $X$ is a group
under composition. This group is called the \emph{symmetric group} of $X$, and
is denoted by $S_{X}$. Its neutral element is the identity map
$\operatorname*{id}\nolimits_{X}:X\rightarrow X$. Its size is $\left\vert
X\right\vert !$ when $X$ is finite.

(Alternative notations for $S_{X}$ include $\operatorname*{Sym}\left(
X\right)  $ and $\Sigma_{X}$ and $\mathfrak{S}_{X}$ and $\mathcal{S}_{X}$.)
\medskip

\textbf{(c)} As usual in group theory, we will write $\alpha\beta$ for the
composition $\alpha\circ\beta$ when $\alpha,\beta\in S_{X}$. This is the map
that sends each $x\in X$ to $\alpha\left(  \beta\left(  x\right)  \right)  $.
\medskip

\textbf{(d)} If $\alpha\in S_{X}$ and $i\in\mathbb{Z}$, then $\alpha^{i}$
shall denote the $i$-th power of $\alpha$ in the group $S_{X}$. Note that
$\alpha^{i}=\underbrace{\alpha\circ\alpha\circ\cdots\circ\alpha}_{i\text{
times}}$ if $i\geq0$. Also, $\alpha^{0}=\operatorname*{id}\nolimits_{X}$.
Also, $\alpha^{-1}$ is the inverse of $\alpha$ in the group $S_{X}$, that is,
the inverse of the map $\alpha$.
\end{definition}

\begin{definition}
\label{def.perm.Sn-iven}Let $n\in\mathbb{Z}$. Then, $\left[  n\right]  $ shall
mean the set $\left\{  1,2,\ldots,n\right\}  $. This is an $n$-element set if
$n\geq0$, and is an empty set if $n\leq0$.

The symmetric group $S_{\left[  n\right]  }$ (consisting of all permutations
of $\left[  n\right]  $) will be denoted $S_{n}$ and called the $n$\emph{-th
symmetric group}. Its size is $n!$ (when $n\geq0$).
\end{definition}

For instance, $S_{3}$ is the group of all $6$ permutations of the set $\left[
3\right]  =\left\{  1,2,3\right\}  $.

If two sets $X$ and $Y$ are in bijection, then their symmetric groups $S_{X}$
and $S_{Y}$ are isomorphic. Intuitively, this is clear (just think of $Y$ as a
\textquotedblleft copy\textquotedblright\ of $X$ with all elements relabelled,
and use this to reinterpret each permutation of $X$ as a permutation of $Y$).
We can formalize this as the following proposition:

\begin{proposition}
\label{prop.perm.Sf}Let $X$ and $Y$ be two sets, and let $f:X\rightarrow Y$ be
a bijection. Then, for each permutation $\sigma$ of $X$, the map $f\circ
\sigma\circ f^{-1}:Y\rightarrow Y$ is a permutation of $Y$. Furthermore, the
map%
\begin{align*}
S_{f}:S_{X}  &  \rightarrow S_{Y},\\
\sigma &  \mapsto f\circ\sigma\circ f^{-1}%
\end{align*}
is a group isomorphism; thus, we obtain $S_{X}\cong S_{Y}$.
\end{proposition}

\begin{proof}
Easy and LTTR.
\end{proof}

Because of Proposition \ref{prop.perm.Sf}, if you want to understand the
symmetric groups of finite sets, you only need to understand $S_{n}$ for all
$n\in\mathbb{N}$ (because if $X$ is a finite set of size $n$, then there is a
bijection $f:X\rightarrow\left[  n\right]  $ and therefore a group isomorphism
$S_{f}:S_{X}\rightarrow S_{\left[  n\right]  }$). Thus, we will focus mostly
on $S_{n}$ in this chapter.

\begin{remark}
If $Y=X$ in Proposition \ref{prop.perm.Sf}, then the group isomorphism $S_{f}$
is conjugation by $f$ in the group $S_{X}$.
\end{remark}

Next, let us define three ways to represent a permutation:

\begin{definition}
\label{def.perm.notations}Let $n\in\mathbb{N}$ and $\sigma\in S_{n}$. We
introduce three notations for $\sigma$: \medskip

\textbf{(a)} A \emph{two-line notation} of $\sigma$ means a $2\times n$-array
\[
\left(
\begin{array}
[c]{cccc}%
p_{1} & p_{2} & \cdots & p_{n}\\
\sigma\left(  p_{1}\right)  & \sigma\left(  p_{2}\right)  & \cdots &
\sigma\left(  p_{n}\right)
\end{array}
\right)  ,
\]
where the entries $p_{1},p_{2},\ldots,p_{n}$ of the top row are the $n$
elements of $\left[  n\right]  $ in some order. Note that this is a standard
notation for any kind of map from a finite set. Commonly, we pick $p_{i}=i$,
so we get the array
\[
\left(
\begin{array}
[c]{cccc}%
1 & 2 & \cdots & n\\
\sigma\left(  1\right)  & \sigma\left(  2\right)  & \cdots & \sigma\left(
n\right)
\end{array}
\right)  .
\]

\textbf{(b)} The \emph{one-line notation} (short, \emph{OLN}) of $\sigma$
means the $n$-tuple $\left(  \sigma\left(  1\right)  ,\sigma\left(  2\right)
,\ldots,\sigma\left(  n\right)  \right)  $.

It is common to omit the commas and the parentheses when writing down the OLN
of $\sigma$. Thus, one simply writes $\sigma\left(  1\right)  \ \sigma\left(
2\right)  \ \cdots\ \sigma\left(  n\right)  $ instead of $\left(
\sigma\left(  1\right)  ,\sigma\left(  2\right)  ,\ldots,\sigma\left(
n\right)  \right)  $. Note that this omission can make the notation ambiguous
if some of the $\sigma\left(  i\right)  $ have more than one digit (for
example, the OLN $1112345678910$ can mean two different permutations of
$\left[  11\right]  $, depending on whether you read the \textquotedblleft%
$111$\textquotedblright\ part as \textquotedblleft$1,11$\textquotedblright\ or
as \textquotedblleft$11,1$\textquotedblright). However, if $n\leq10$, then
this ambiguity does not occur, and the notation is unproblematic (even without
commas and parentheses). \medskip

\textbf{(c)} The \emph{cycle digraph} of $\sigma$ is defined (informally) as follows:

\begin{itemize}
\item For each $i\in\left[  n\right]  $, draw a point (\textquotedblleft
node\textquotedblright) labelled $i$.

\item For each $i\in\left[  n\right]  $, draw an arrow (\textquotedblleft
arc\textquotedblright) from the node labelled $i$ to the node labelled
$\sigma\left(  i\right)  $.
\end{itemize}

The resulting picture is called the cycle digraph of $\sigma$.

Using the concept of \emph{digraphs} (= directed graphs), this definition can
be restated formally as follows: The \emph{cycle digraph} of $\sigma$ is the
directed graph with vertices $1,2,\ldots,n$ and arcs $i\rightarrow
\sigma\left(  i\right)  $ for all $i\in\left[  n\right]  $.
\end{definition}

\begin{example}
Let $\sigma:\left[  4\right]  \rightarrow\left[  4\right]  $ be the map that
sends the elements $1,2,3,4$ to $2,4,3,1$, respectively. Then, $\sigma$ is a
bijection, thus a permutation of $\left[  4\right]  $. \medskip

\textbf{(a)} A two-line notation of $\sigma$ is $\left(
\begin{array}
[c]{cccc}%
1 & 2 & 3 & 4\\
2 & 4 & 3 & 1
\end{array}
\right)  $. Another is $\left(
\begin{array}
[c]{cccc}%
3 & 1 & 4 & 2\\
3 & 2 & 1 & 4
\end{array}
\right)  $. Another is $\left(
\begin{array}
[c]{cccc}%
4 & 1 & 3 & 2\\
1 & 2 & 3 & 4
\end{array}
\right)  $. There are $24$ two-line notations of $\sigma$ in total, since we
can freely choose the order in which the elements of $\left[  4\right]  $
appear in the top row. \medskip

\textbf{(b)} The one-line notation of $\sigma$ is $\left(  2,4,3,1\right)  $.
Omitting the commas and the parentheses, we can rewrite this as $2431$.
\medskip

\textbf{(c)} One way to draw the cycle digraph of $\sigma$ is
\[%
\begin{tikzpicture}%
[->,shorten >=1pt,auto,node distance=3cm, thick,main node/.style={circle,fill=blue!20,draw}%
]
\node[main node] (1) {1};
\node[main node] (2) [right of=1] {2};
\node[main node] (3) [right of=2] {3};
\node[main node] (4) [right of=3] {4};
\path[-{Stealth[length=4mm]}]
(1) edge [bend left] (2)
(2) edge [bend left] (4)
(4) edge [bend left] (1)
(3) edge [loop right] (3);
\end{tikzpicture}%
\ \ .
\]
Another is
\[%
\begin{tikzpicture}%
[->,shorten >=1pt,auto,node distance=3cm, thick,main node/.style={circle,fill=blue!20,draw}%
]
\node[main node] (1) {1};
\node[main node] (2) [below left of=1] {2};
\node[main node] (3) [right of=1] {3};
\node[main node] (4) [left of=1] {4};
\path[-{Stealth[length=4mm]}]
(1) edge [bend left] (2)
(2) edge [bend left] (4)
(4) edge [bend left] (1)
(3) edge [loop right] (3);
\end{tikzpicture}%
\ \ .
\]
(When drawing cycle digraphs, one commonly tries to place the nodes in such a
way as to make the arcs as short as possible. Thus, it is natural to keep the
cycles separate in the picture. But formally speaking, any picture is fine, as
long as the nodes and arcs don't overlap.)
\end{example}

\begin{example}
Let $\sigma:\left[  10\right]  \rightarrow\left[  10\right]  $ be the map that
sends the elements $1,2,3,4,5,6,7,8,9,10$ to $5,4,3,2,6,10,1,9,8,7$,
respectively. Then, $\sigma$ is a bijection, hence a permutation of $\left[
10\right]  $. The one-line notation of $\sigma$ is $\left(
5,4,3,2,6,10,1,9,8,7\right)  $. If we omit the commas and the parentheses,
then this becomes%
\[
5\ 4\ 3\ 2\ 6\ \left(  10\right)  \ 1\ 9\ 8\ 7.
\]
(We have put the $10$ in parentheses to make its place clearer.) The cycle
digraph of $\sigma$ is%
\[%
\begin{tikzpicture}%
[->,shorten >=1pt,auto,node distance=3cm, thick,main node/.style={circle,fill=blue!20,draw}%
]
\node[main node] (0) {10};
\node[main node] (6) [above of=0] {6};
\node[main node] (5) [above left of=6] {5};
\node[main node] (1) [below left of=5] {1};
\node[main node] (7) [below of=1] {7};
\node[main node] (9) [above right of=6] {9};
\node[main node] (8) [right of=9] {8};
\node[main node] (2) [right of=0] {2};
\node[main node] (4) [right of=2] {4};
\node[main node] (3) [below of=8] {3};
\path[-{Stealth[length=4mm]}]
(1) edge [bend left] (5)
(5) edge [bend left] (6)
(6) edge [bend left] (0)
(0) edge [bend left] (7)
(7) edge [bend left] (1)
(9) edge [bend right] (8)
(8) edge [bend right] (9)
(2) edge [bend right] (4)
(4) edge [bend right] (2)
(3) edge [loop right] (3);
\end{tikzpicture}%
\ \ .
\]

\end{example}

Note that two-line notations and cycle digraphs can be defined (just as we did
in Definition \ref{def.perm.notations}) not just for permutations $\sigma\in
S_{n}$, but more generally for permutations $\sigma\in S_{X}$ of any finite
set $X$. But the one-line notation makes sense for $\sigma\in S_{n}$ only.

\subsection{Transpositions, cycles and involutions}

We shall now define some important families of permutations.

\subsubsection{Transpositions}

\begin{definition}
\label{def.perm.tij}Let $i$ and $j$ be two distinct elements of a set $X$.

Then, the \emph{transposition} $t_{i,j}$ is the permutation of $X$ that sends
$i$ to $j$, sends $j$ to $i$, and leaves all other elements of $X$ unchanged.
\end{definition}

Strictly speaking, the notation $t_{i,j}$ is somewhat ambiguous, since it
suppresses $X$. However, most of the times we will use it, the set $X$ will be
either clear from the context or irrelevant.

\begin{example}
The permutation $t_{2,4}$ of the set $\left[  7\right]  $ sends the elements
$1,2,3,4,5,6,7$ to $1,4,3,2,5,6,7$, respectively. Its one-line notation (with
commas and parentheses omitted) is therefore $1432567$.
\end{example}

Note that $t_{i,j}=t_{j,i}$ for any two distinct elements $i$ and $j$ of a set
$X$.

\begin{definition}
\label{def.perm.si}Let $n\in\mathbb{N}$ and $i\in\left[  n-1\right]  $. Then,
the \emph{simple transposition} $s_{i}$ is defined by%
\[
s_{i}:=t_{i,i+1}\in S_{n}.
\]

\end{definition}

Thus, a simple transposition is a transposition that swaps two consecutive
integers. Again, the notation $s_{i}$ suppresses $n$, but this won't usually
be a problem.

\begin{example}
The permutation $s_{2}$ of the set $\left[  7\right]  $ sends the elements
$1,2,3,4,5,6,7$ to $1,3,2,4,5,6,7$, respectively. Its one-line notation is
therefore $1324567$.
\end{example}

Here are some very basic properties of simple transpositions:\footnote{Recall
Definition \ref{def.perm.perm}. Thus, for example, $s_{i}^{2}$ means
$s_{i}s_{i}=s_{i}\circ s_{i}$, whereas $s_{i}s_{j}$ means $s_{i}\circ s_{j}$.}

\begin{proposition}
\label{prop.perm.si.rules}Let $n\in\mathbb{N}$. \medskip

\textbf{(a)} We have $s_{i}^{2}=\operatorname*{id}$ for all $i\in\left[
n-1\right]  $. In other words, we have $s_{i}=s_{i}^{-1}$ for all $i\in\left[
n-1\right]  $. \medskip

\textbf{(b)} We have $s_{i}s_{j}=s_{j}s_{i}$ for any $i,j\in\left[
n-1\right]  $ with $\left\vert i-j\right\vert >1$. \medskip

\textbf{(c)} We have $s_{i}s_{i+1}s_{i}=s_{i+1}s_{i}s_{i+1}$ for any
$i\in\left[  n-2\right]  $.
\end{proposition}

\begin{proof}
To prove that two permutations $\alpha$ and $\beta$ of $\left[  n\right]  $
are identical, it suffices to show that $\alpha\left(  k\right)  =\beta\left(
k\right)  $ for each $k\in\left[  n\right]  $. Using this strategy, we can
prove all three parts of Proposition \ref{prop.perm.si.rules}
straightforwardly (distinguishing cases corresponding to the relative
positions of $k$, $i$, $i+1$, $j$ and $j+1$). This is done in detail for
Proposition \ref{prop.perm.si.rules} \textbf{(c)} in \cite[solution to
Exercise 5.1 \textbf{(a)}]{detnotes}; the proofs of parts \textbf{(a)} and
\textbf{(b)} are easier and LTTR.
\end{proof}

\subsubsection{Cycles}

The following definition can be viewed as a generalization of Definition
\ref{def.perm.tij}:

\begin{definition}
\label{def.perm.cycs}Let $X$ be a set. Let $i_{1},i_{2},\ldots,i_{k}$ be $k$
distinct elements of $X$, where $k\geq1$. Then,%
\[
\operatorname*{cyc}\nolimits_{i_{1},i_{2},\ldots,i_{k}}%
\]
means the permutation of $X$ that sends%
\begin{align*}
&  i_{1}\text{ to }i_{2},\\
&  i_{2}\text{ to }i_{3},\\
&  i_{3}\text{ to }i_{4},\\
&  \ldots,\\
&  i_{k-1}\text{ to }i_{k},\\
&  i_{k}\text{ to }i_{1}%
\end{align*}
and leaves all other elements of $X$ unchanged. In other words,
$\operatorname*{cyc}\nolimits_{i_{1},i_{2},\ldots,i_{k}}$ means the
permutation of $X$ that satisfies%
\begin{align*}
\operatorname*{cyc}\nolimits_{i_{1},i_{2},\ldots,i_{k}}\left(  p\right)   &  =%
\begin{cases}
i_{j+1}, & \text{if }p=i_{j}\text{ for some }j\in\left\{  1,2,\ldots
,k\right\}  ;\\
p, & \text{otherwise}%
\end{cases}
\\
&
\ \ \ \ \ \ \ \ \ \ \ \ \ \ \ \ \ \ \ \ \ \ \ \ \ \ \ \ \ \ \ \ \ \ \ \ \ \ \ \ \text{for
every }p\in X,
\end{align*}
where $i_{k+1}$ means $i_{1}$.

This permutation is called a $k$\emph{-cycle}.
\end{definition}

The name \textquotedblleft$k$-cycle\textquotedblright\ harkens back to the
cycle digraph of $\operatorname*{cyc}\nolimits_{i_{1},i_{2},\ldots,i_{k}}$,
which consists of a cycle of length $k$ (containing the nodes $i_{1}%
,i_{2},\ldots,i_{k}$ in this order) along with $\left\vert X\right\vert -k$
isolated nodes (more precisely, each of the elements of $X\setminus\left\{
i_{1},i_{2},\ldots,i_{k}\right\}  $ has an arrow from itself to itself in the
cycle digraph). Here is an example:

\begin{example}
Let $X=\left[  8\right]  $. Then, the permutation $\operatorname*{cyc}%
\nolimits_{2,6,5}$ of $X$ sends
\begin{align*}
&  2\text{ to }6,\\
&  6\text{ to }5,\\
&  5\text{ to }2
\end{align*}
and leaves all other elements of $X$ unchanged. Thus, this permutation has OLN
$16342578$ and cycle digraph
\[%
\begin{tikzpicture}%
[->,shorten >=1pt,auto,node distance=3cm, thick,main node/.style={circle,fill=blue!20,draw}%
]
\node[main node] (1) {1};
\node[main node] (3) [right of=1] {3};
\node[main node] (4) [below of=1] {4};
\node[main node] (7) [right of=4] {7};
\node[main node] (8) [right of=7] {8};
\node[main node] (6) [left of=1] {6};
\node[main node] (5) [left of=6] {5};
\node[main node] (2) [below left of=6] {2};
\path[-{Stealth[length=4mm]}]
(2) edge [bend right] (6)
(6) edge [bend right] (5)
(5) edge [bend right] (2)
(3) edge [loop right] (3)
(4) edge [loop right] (4)
(7) edge [loop right] (7)
(8) edge [loop right] (8)
(1) edge [loop right] (1);
\end{tikzpicture}%
\]

\end{example}

\begin{example}
Let $X$ be a set. If $i$ and $j$ are two distinct elements of $X$, then
$\operatorname*{cyc}\nolimits_{i,j}=t_{i,j}$. Thus, the $2$-cycles in $S_{X}$
are precisely the transpositions in $S_{X}$, so there are $\dbinom{\left\vert
X\right\vert }{2}$ many of them (since any $2$-element subset $\left\{
i,j\right\}  $ of $X$ gives rise to a transposition $t_{i,j}$, and this
assignment of transpositions to $2$-element subsets is bijective).
\end{example}

Note that the $k$-cycle $\operatorname*{cyc}\nolimits_{i_{1},i_{2}%
,\ldots,i_{k}}$ is often denoted by $\left(  i_{1},i_{2},\ldots,i_{k}\right)
$, but I will not use this notation here, since it clashes with the standard
notation for $k$-tuples.

\begin{exercise}
\label{exe.perm.cyc.how-many-kcyc}Let $n\in\mathbb{N}$ and let $k\in\left[
n\right]  $. Let $X$ be an $n$-element set. How many $k$-cycles exist in
$S_{X}$ ?
\end{exercise}

\begin{proof}
[Solution to Exercise \ref{exe.perm.cyc.how-many-kcyc} (sketched).]First, we
note that there is exactly one $1$-cycle in $S_{X}$ (for $n>0$), since a
$1$-cycle is just the identity map. This should be viewed as a degenerate
case; thus, we WLOG assume that $k>1$.

For any $k$ distinct elements $i_{1},i_{2},\ldots,i_{k}$ of $X$, we have%
\[
\operatorname*{cyc}\nolimits_{i_{1},i_{2},\ldots,i_{k}}=\operatorname*{cyc}%
\nolimits_{i_{2},i_{3},\ldots,i_{k},i_{1}}=\operatorname*{cyc}\nolimits_{i_{3}%
,i_{4},\ldots,i_{k},i_{1},i_{2}}=\cdots=\operatorname*{cyc}\nolimits_{i_{k}%
,i_{1},i_{2},\ldots,i_{k-1}}.
\]
That is, $\operatorname*{cyc}\nolimits_{i_{1},i_{2},\ldots,i_{k}}$ does not
change if we cyclically rotate the list $\left(  i_{1},i_{2},\ldots
,i_{k}\right)  $.

Any $k$-cycle $\operatorname*{cyc}\nolimits_{i_{1},i_{2},\ldots,i_{k}}$
uniquely determines the elements $i_{1},i_{2},\ldots,i_{k}$ up to cyclic
rotation (since $k>1$). Indeed, if $\sigma=\operatorname*{cyc}\nolimits_{i_{1}%
,i_{2},\ldots,i_{k}}$ is a $k$-cycle, then the elements $i_{1},i_{2}%
,\ldots,i_{k}$ are precisely the elements of $X$ that are not fixed by
$\sigma$ (it is here that we use our assumption $k>1$), and furthermore, if we
know which of these elements is $i_{1}$, then we can reconstruct the remaining
elements $i_{2},i_{3},\ldots,i_{k}$ recursively by%
\[
i_{2}=\sigma\left(  i_{1}\right)  ,\ \ \ \ \ \ \ \ \ \ i_{3}=\sigma\left(
i_{2}\right)  ,\ \ \ \ \ \ \ \ \ \ i_{4}=\sigma\left(  i_{3}\right)
,\ \ \ \ \ \ \ \ \ \ \ldots,\ \ \ \ \ \ \ \ \ \ i_{k}=\sigma\left(
i_{k-1}\right)
\]
(that is, $i_{2},i_{3},\ldots,i_{k}$ are obtained by iteratively applying
$\sigma$ to $i_{1}$). Therefore, if $\sigma\in S_{X}$ is a $k$-cycle, then
there are precisely $k$ lists $\left(  i_{1},i_{2},\ldots,i_{k}\right)  $ for
which $\sigma=\operatorname*{cyc}\nolimits_{i_{1},i_{2},\ldots,i_{k}}$ (coming
from the $k$ possibilities for which of the $k$ non-fixed points of $\sigma$
should be $i_{1}$).

Hence, the map%
\begin{align*}
f:\left\{  k\text{-tuples of distinct elements of }X\right\}   &
\rightarrow\left\{  k\text{-cycles in }S_{X}\right\}  ,\\
\left(  i_{1},i_{2},\ldots,i_{k}\right)   &  \mapsto\operatorname*{cyc}%
\nolimits_{i_{1},i_{2},\ldots,i_{k}}%
\end{align*}
is a $k$-to-$1$ map (i.e., each $k$-cycle in $S_{X}$ has precisely $k$
preimages under this map). Therefore, Lemma \ref{lem.count.multijection}
(applied to $m=k$ and \newline$A=\left\{  k\text{-tuples of distinct elements
of }X\right\}  $ and $B=\left\{  k\text{-cycles in }S_{X}\right\}  $) yields%
\[
\left(  \text{\# of }k\text{-tuples of distinct elements of }X\right)
=k\cdot\left(  \text{\# of }k\text{-cycles in }S_{X}\right)  .
\]
Therefore,
\begin{align*}
\left(  \text{\# of }k\text{-cycles in }S_{X}\right)   &  =\dfrac{1}{k}%
\cdot\underbrace{\left(  \text{\# of }k\text{-tuples of distinct elements of
}X\right)  }_{\substack{=n\left(  n-1\right)  \left(  n-2\right)
\cdots\left(  n-k+1\right)  \\\text{(by (\ref{eq.count.k-tups-dist-elts}),
since }X\text{ is an }n\text{-element set)}}}\\
&  =\dfrac{1}{k}\cdot n\left(  n-1\right)  \left(  n-2\right)  \cdots\left(
n-k+1\right) \\
&  =\dbinom{n}{k}\cdot\left(  k-1\right)  !\ \ \ \ \ \ \ \ \ \ \left(
\text{by a bit of simple algebra}\right)  .
\end{align*}
This is the answer to Exercise \ref{exe.perm.cyc.how-many-kcyc} in the case
$k>1$. Hence, Exercise \ref{exe.perm.cyc.how-many-kcyc} is solved.
\end{proof}

\subsubsection{Involutions}

\emph{Involutions} are maps that are inverse to themselves:

\begin{definition}
\label{def.perm.invol}Let $X$ be a set. An \emph{involution} of $X$ means a
map $f:X\rightarrow X$ that satisfies $f\circ f=\operatorname*{id}$. Clearly,
an involution is always a permutation, and equals its own inverse.
\end{definition}

For example, the identity map $\operatorname*{id}\nolimits_{X}$ is an
involution, and any transposition $t_{i,j}\in S_{X}$ is an involution, whereas
a $k$-cycle $\operatorname*{cyc}\nolimits_{i_{1},i_{2},\ldots,i_{k}}$ with
$k>2$ is never an involution. Further examples of involutions are products of
disjoint transpositions -- i.e., permutations of the form $t_{i_{1},j_{1}%
}t_{i_{2},j_{2}}\cdots t_{i_{k},j_{k}}$, where all $2k$ elements $i_{1}%
,j_{1},i_{2},j_{2},\ldots,i_{k},j_{k}$ are distinct. In fact, any involution
of a finite set can be expressed in this form.

An involution $\sigma\in S_{n}$ has a cycle digraph consisting of $1$-cycles
and $2$-cycles. For example, here is the cycle digraph of the involution
$t_{1,5}t_{2,3}\in S_{7}$:%
\[%
\begin{tikzpicture}%
[->,shorten >=1pt,auto,node distance=3cm, thick,main node/.style={circle,fill=blue!20,draw}%
]
\node[main node] (1) {1};
\node[main node] (5) [above of=1] {5};
\node[main node] (2) [right of=1] {2};
\node[main node] (3) [above of=2] {3};
\node[main node] (4) [right of=2] {4};
\node[main node] (6) [right of=4] {6};
\node[main node] (7) [right of=6] {7};
\path[-{Stealth[length=4mm]}]
(1) edge [bend right] (5)
(5) edge [bend right] (1)
(2) edge [bend right] (3)
(3) edge [bend right] (2)
(4) edge [loop right] (4)
(7) edge [loop right] (7)
(6) edge [loop right] (6);
\end{tikzpicture}%
\]

\subsection{\label{sec.perm.inversions}Inversions, length and Lehmer codes}

\subsubsection{Inversions and lengths}

Let us define some features of arbitrary permutations of $\left[  n\right]  $:

\begin{definition}
\label{def.perm.invs}Let $n\in\mathbb{N}$ and $\sigma\in S_{n}$. \medskip

\textbf{(a)} An \emph{inversion} of $\sigma$ means a pair $\left(  i,j\right)
$ of elements of $\left[  n\right]  $ such that $i<j$ and $\sigma\left(
i\right)  >\sigma\left(  j\right)  $. \medskip

\textbf{(b)} The \emph{length} (also known as the \emph{Coxeter length}) of
$\sigma$ is the \# of inversions of $\sigma$. It is called $\ell\left(
\sigma\right)  $. (Some authors call it $\operatorname*{inv}\sigma$ instead.)
\end{definition}

(In LaTeX, the symbol \textquotedblleft$\ell$\textquotedblright\ is obtained
by typing \textquotedblleft\texttt{%
$\backslash$%
ell}\textquotedblright. If you just type \textquotedblleft\texttt{l}%
\textquotedblright, you will get \textquotedblleft$l$\textquotedblright.)

An inversion of a permutation $\sigma$ can thus be viewed as a pair of
elements of $\left[  n\right]  $ whose relative order changes when $\sigma$ is
applied to them. (We require this pair $\left(  i,j\right)  $ to satisfy $i<j$
in order not to count each such pair doubly.)

\begin{example}
Let $\pi\in S_{4}$ be the permutation with OLN $3142$. The inversions of $\pi$
are%
\begin{align*}
&  \left(  1,2\right)  \ \ \ \ \ \ \ \ \ \ \left(  \text{since }1<2\text{ and
}\underbrace{\pi\left(  1\right)  }_{=3}>\underbrace{\pi\left(  2\right)
}_{=1}\right)  \ \ \ \ \ \ \ \ \ \ \text{and}\\
&  \left(  1,4\right)  \ \ \ \ \ \ \ \ \ \ \left(  \text{since }1<4\text{ and
}\underbrace{\pi\left(  1\right)  }_{=3}>\underbrace{\pi\left(  4\right)
}_{=2}\right)  \ \ \ \ \ \ \ \ \ \ \text{and}\\
&  \left(  3,4\right)  \ \ \ \ \ \ \ \ \ \ \left(  \text{since }3<4\text{ and
}\underbrace{\pi\left(  3\right)  }_{=4}>\underbrace{\pi\left(  4\right)
}_{=2}\right)  .
\end{align*}
Thus, the length of $\pi$ is $3$.
\end{example}

For a given $n\in\mathbb{N}$ and a given $k\in\mathbb{N}$, how many
permutations $\sigma\in S_{n}$ have length $k$ ? The following proposition
gives a partial answer:

\begin{proposition}
\label{prop.perm.lengths-k-small-k}Let $n\in\mathbb{N}$. \medskip

\textbf{(a)} For any $\sigma\in S_{n}$, we have $\ell\left(  \sigma\right)
\in\left\{  0,1,\ldots,\dbinom{n}{2}\right\}  $. \medskip

\textbf{(b)} We have%
\[
\left(  \text{\# of }\sigma\in S_{n}\text{ with }\ell\left(  \sigma\right)
=0\right)  =1.
\]
Indeed, the only permutation $\sigma\in S_{n}$ with $\ell\left(
\sigma\right)  =0$ is the identity map $\operatorname*{id}$. \medskip

\textbf{(c)} We have%
\[
\left(  \text{\# of }\sigma\in S_{n}\text{ with }\ell\left(  \sigma\right)
=\dbinom{n}{2}\right)  =1.
\]
Indeed, the only permutation $\sigma\in S_{n}$ with $\ell\left(
\sigma\right)  =\dbinom{n}{2}$ is the permutation with OLN $n\left(
n-1\right)  \left(  n-2\right)  \cdots21$. (This permutation is often called
$w_{0}$.) \medskip

\textbf{(d)} If $n\geq1$, then%
\[
\left(  \text{\# of }\sigma\in S_{n}\text{ with }\ell\left(  \sigma\right)
=1\right)  =n-1.
\]
Indeed, the only permutations $\sigma\in S_{n}$ with $\ell\left(
\sigma\right)  =1$ are the simple transpositions $s_{i}$ with $i\in\left[
n-1\right]  $. \medskip

\textbf{(e)} If $n\geq2$, then%
\[
\left(  \text{\# of }\sigma\in S_{n}\text{ with }\ell\left(  \sigma\right)
=2\right)  =\dfrac{\left(  n-2\right)  \left(  n+1\right)  }{2}.
\]
Indeed, the only permutations $\sigma\in S_{n}$ with $\ell\left(
\sigma\right)  =2$ are the products $s_{i}s_{j}$ with $1\leq i<j<n$ as well as
the products $s_{i}s_{i-1}$ with $i\in\left\{  2,3,\ldots,n-1\right\}  $. If
$n\geq2$, then there are $\dfrac{\left(  n-2\right)  \left(  n+1\right)  }{2}$
such products (and they are all distinct). \medskip

\textbf{(f)} For any $k\in\mathbb{Z}$, we have%
\[
\left(  \text{\# of }\sigma\in S_{n}\text{ with }\ell\left(  \sigma\right)
=k\right)  =\left(  \text{\# of }\sigma\in S_{n}\text{ with }\ell\left(
\sigma\right)  =\dbinom{n}{2}-k\right)  .
\]

\end{proposition}

\begin{proof}
Exercise \ref{exe.perm.lengths-k-small-k}.
\end{proof}

What about the general case? Alas, there is no explicit formula for the \# of
$\sigma\in S_{n}$ with $\ell\left(  \sigma\right)  =k$. However, there is a
nice formula for the generating function%
\[
\sum_{k\in\mathbb{N}}\left(  \text{\# of }\sigma\in S_{n}\text{ with }%
\ell\left(  \sigma\right)  =k\right)  x^{k}=\sum_{\sigma\in S_{n}}%
x^{\ell\left(  \sigma\right)  }.
\]
Let us first sound it out on the case $n=3$:

\begin{example}
Written in one-line notation, the permutations of the set $\left[  3\right]  $
are $123$, $132$, $213$, $231$,\ $312$, and $321$. Their lengths are%
\begin{align*}
\ell\left(  123\right)   &  =0,\ \ \ \ \ \ \ \ \ \ \ell\left(  132\right)
=1,\ \ \ \ \ \ \ \ \ \ \ell\left(  213\right)  =1,\\
\ell\left(  231\right)   &  =2,\ \ \ \ \ \ \ \ \ \ \ell\left(  312\right)
=2,\ \ \ \ \ \ \ \ \ \ \ell\left(  321\right)  =3.
\end{align*}
Thus,
\begin{align*}
\sum_{\sigma\in S_{3}}x^{\ell\left(  \sigma\right)  }  &  =x^{\ell\left(
123\right)  }+x^{\ell\left(  132\right)  }+x^{\ell\left(  213\right)
}+x^{\ell\left(  231\right)  }+x^{\ell\left(  312\right)  }+x^{\ell\left(
321\right)  }\\
&  \ \ \ \ \ \ \ \ \ \ \ \ \ \ \ \ \ \ \ \ \left(  \text{where we are writing
each }\sigma\in S_{3}\text{ in OLN}\right) \\
&  =x^{0}+x^{1}+x^{1}+x^{2}+x^{2}+x^{3}=1+2x+2x^{2}+x^{3}\\
&  =\left(  1+x\right)  \left(  1+x+x^{2}\right)  .
\end{align*}

\end{example}

This suggests the following general result:

\begin{proposition}
\label{prop.perm.length.gf}Let $n\in\mathbb{N}$. Then,%
\begin{align*}
&  \sum_{\sigma\in S_{n}}x^{\ell\left(  \sigma\right)  }\\
&  =\prod_{i=1}^{n-1}\left(  1+x+x^{2}+\cdots+x^{i}\right) \\
&  =\left(  1+x\right)  \left(  1+x+x^{2}\right)  \left(  1+x+x^{2}%
+x^{3}\right)  \cdots\left(  1+x+x^{2}+\cdots+x^{n-1}\right) \\
&  =\left[  n\right]  _{x}!.
\end{align*}
(Here, we are using Definition \ref{def.pars.qbinom.qint}, so that $\left[
n\right]  _{x}!$ means the result of substituting $x$ for $q$ in the
$q$-factorial $\left[  n\right]  _{q}!$.)
\end{proposition}

\subsubsection{Lehmer codes}

We will prove this proposition using the so-called \emph{Lehmer code} of a
permutation, which is defined as follows:

\begin{definition}
\label{def.perm.lehmer1}Let $n\in\mathbb{N}$. The following notations will be
used throughout Section \ref{sec.perm.inversions}: \medskip

\textbf{(a)} For each $\sigma\in S_{n}$ and $i\in\left[  n\right]  $, we set%
\begin{align*}
\ell_{i}\left(  \sigma\right)  :=  &  \left(  \text{\# of all }j\in\left[
n\right]  \text{ satisfying }i<j\text{ and }\sigma\left(  i\right)
>\sigma\left(  j\right)  \right) \\
=  &  \left(  \text{\# of all }j\in\left\{  i+1,i+2,\ldots,n\right\}  \text{
such that }\sigma\left(  i\right)  >\sigma\left(  j\right)  \right)  .
\end{align*}
(The last equality sign here is clear, since the $j\in\left[  n\right]  $
satisfying $i<j$ are precisely the $j\in\left\{  i+1,i+2,\ldots,n\right\}  $.)
\medskip

\textbf{(b)} For each $m\in\mathbb{Z}$, we let $\left[  m\right]  _{0}$ denote
the set $\left\{  0,1,\ldots,m\right\}  $. (This is an empty set when $m<0$.)
\medskip

\textbf{(c)} We let $H_{n}$ denote the set
\begin{align*}
&  \left[  n-1\right]  _{0}\times\left[  n-2\right]  _{0}\times\cdots
\times\left[  n-n\right]  _{0}\\
&  =\left\{  \left(  j_{1},j_{2},\ldots,j_{n}\right)  \in\mathbb{N}^{n}%
\ \mid\ j_{i}\leq n-i\text{ for each }i\in\left[  n\right]  \right\}  .
\end{align*}
This set $H_{n}$ has size%
\begin{align*}
\left\vert H_{n}\right\vert  &  =\left\vert \left[  n-1\right]  _{0}%
\times\left[  n-2\right]  _{0}\times\cdots\times\left[  n-n\right]
_{0}\right\vert \\
&  =\underbrace{\left\vert \left[  n-1\right]  _{0}\right\vert }_{=n}%
\cdot\underbrace{\left\vert \left[  n-2\right]  _{0}\right\vert }_{=n-1}%
\cdot\cdots\cdot\underbrace{\left\vert \left[  n-n\right]  _{0}\right\vert
}_{=1}\\
&  =n\left(  n-1\right)  \cdots1=n!.
\end{align*}

\textbf{(d)} We define the map%
\begin{align*}
L:S_{n}  &  \rightarrow H_{n},\\
\sigma &  \mapsto\left(  \ell_{1}\left(  \sigma\right)  ,\ell_{2}\left(
\sigma\right)  ,\ldots,\ell_{n}\left(  \sigma\right)  \right)  .
\end{align*}
(This map is well-defined, since each $\sigma\in S_{n}$ and each $i\in\left[
n\right]  $ satisfy $\ell_{i}\left(  \sigma\right)  \in\left\{  0,1,\ldots
,n-i\right\}  =\left[  n-i\right]  _{0}$.) \medskip

\textbf{(e)} If $\sigma\in S_{n}$ is a permutation, then the $n$-tuple
$L\left(  \sigma\right)  =\left(  \ell_{1}\left(  \sigma\right)  ,\ell
_{2}\left(  \sigma\right)  ,\ldots,\ell_{n}\left(  \sigma\right)  \right)  $
is called the \emph{Lehmer code} (or just the \emph{code}) of $\sigma$.
\end{definition}

\begin{example}
\label{exa.perm.lehmer1}\textbf{(a)} If $n=6$, and if $\sigma\in S_{6}$ is the
permutation with one-line notation $451263$, then $\ell_{2}\left(
\sigma\right)  =3$ (because the numbers $j\in\left[  n\right]  $ satisfying
$2<j$ and $\sigma\left(  2\right)  >\sigma\left(  j\right)  $ are precisely
$3$, $4$ and $6$, so that there are $3$ of them) and likewise $\ell_{1}\left(
\sigma\right)  =3$ and $\ell_{3}\left(  \sigma\right)  =0$ and $\ell
_{4}\left(  \sigma\right)  =0$ and $\ell_{5}\left(  \sigma\right)  =1$ and
$\ell_{6}\left(  \sigma\right)  =0$, and thus $L\left(  \sigma\right)
=\left(  3,3,0,0,1,0\right)  \in H_{6}$. \medskip

\textbf{(b)} Here is a table of all $6$ permutations $\sigma\in S_{3}$
(written in one-line notation) and their respective Lehmer codes $L\left(
\sigma\right)  $:%
\[%
\begin{tabular}
[c]{|c|c|}\hline
$\sigma$ & $L\left(  \sigma\right)  $\\\hline\hline
$123$ & $\left(  0,0,0\right)  $\\\hline
$132$ & $\left(  0,1,0\right)  $\\\hline
$213$ & $\left(  1,0,0\right)  $\\\hline
$231$ & $\left(  1,1,0\right)  $\\\hline
$312$ & $\left(  2,0,0\right)  $\\\hline
$321$ & $\left(  2,1,0\right)  $\\\hline
\end{tabular}
\ \ .
\]

\end{example}

The Lehmer code $L\left(  \sigma\right)  $ of a permutation $\sigma$ is a
refinement of its length $\ell\left(  \sigma\right)  $, in the sense that it
gives finer information (i.e., we can reconstruct $\ell\left(  \sigma\right)
$ from $L\left(  \sigma\right)  $):

\begin{proposition}
\label{prop.perm.lehmer.l}Let $n\in\mathbb{N}$ and $\sigma\in S_{n}$. Then,
$\ell\left(  \sigma\right)  =\ell_{1}\left(  \sigma\right)  +\ell_{2}\left(
\sigma\right)  +\cdots+\ell_{n}\left(  \sigma\right)  $.
\end{proposition}

\begin{proof}
This follows from the definitions of $\ell\left(  \sigma\right)  $ and
$\ell_{i}\left(  \sigma\right)  $.
\end{proof}

The main property of Lehmer codes is that they uniquely determine
permutations, and in fact are in bijection with them (cf. Example
\ref{exa.perm.lehmer1} \textbf{(b)}):

\begin{theorem}
\label{thm.perm.lehmer.bij}Let $n\in\mathbb{N}$. Then, the map $L:S_{n}%
\rightarrow H_{n}$ is a bijection.
\end{theorem}

We shall sketch two ways of proving this theorem.

\begin{proof}
[First proof of Theorem \ref{thm.perm.lehmer.bij} (sketched).]Let $\sigma\in
S_{n}$. Let $i\in\left[  n\right]  $. Recall that the OLN of $\sigma$ is the
$n$-tuple $\sigma\left(  1\right)  \ \sigma\left(  2\right)  \ \cdots
\ \sigma\left(  n\right)  $. The definition of $\ell_{i}\left(  \sigma\right)
$ can be rewritten as follows:%
\begin{align*}
\ell_{i}\left(  \sigma\right)   &  =\left(  \text{\# of all }j\in\left\{
i+1,i+2,\ldots,n\right\}  \text{ such that }\sigma\left(  i\right)
>\sigma\left(  j\right)  \right) \\
&  =\left(  \text{\# of all entries in the OLN of }\sigma\text{ that
appear}\right. \\
&  \ \ \ \ \ \ \ \ \ \ \ \ \ \ \ \ \ \ \ \ \left.  \text{after }\sigma\left(
i\right)  \text{ but are smaller than }\sigma\left(  i\right)  \right) \\
&  =\left(  \text{\# of elements of }\left[  n\right]  \setminus\left\{
\sigma\left(  1\right)  ,\sigma\left(  2\right)  ,\ldots,\sigma\left(
i\right)  \right\}  \right. \\
&  \ \ \ \ \ \ \ \ \ \ \ \ \ \ \ \ \ \ \ \ \left.  \text{that are smaller than
}\sigma\left(  i\right)  \right)
\end{align*}
(since the elements of $\left[  n\right]  \setminus\left\{  \sigma\left(
1\right)  ,\sigma\left(  2\right)  ,\ldots,\sigma\left(  i\right)  \right\}  $
are precisely the entries that appear after $\sigma\left(  i\right)  $ in the
OLN of $\sigma$). We can replace the set $\left[  n\right]  \setminus\left\{
\sigma\left(  1\right)  ,\sigma\left(  2\right)  ,\ldots,\sigma\left(
i\right)  \right\}  $ by $\left[  n\right]  \setminus\left\{  \sigma\left(
1\right)  ,\sigma\left(  2\right)  ,\ldots,\sigma\left(  i-1\right)  \right\}
$ in this equality\footnote{If $i=1$, then the set $\left\{  \sigma\left(
1\right)  ,\sigma\left(  2\right)  ,\ldots,\sigma\left(  i-1\right)  \right\}
$ is empty, so that we have $\left[  n\right]  \setminus\left\{  \sigma\left(
1\right)  ,\sigma\left(  2\right)  ,\ldots,\sigma\left(  i-1\right)  \right\}
=\left[  n\right]  $ in this case.} (this will not change the \# of elements
of this set that are smaller than $\sigma\left(  i\right)  $, because it only
inserts the element $\sigma\left(  i\right)  $ into the set, and obviously
this element $\sigma\left(  i\right)  $ is not smaller than $\sigma\left(
i\right)  $). Thus, we obtain%
\begin{align}
\ell_{i}\left(  \sigma\right)   &  =\left(  \text{\# of elements of }\left[
n\right]  \setminus\left\{  \sigma\left(  1\right)  ,\sigma\left(  2\right)
,\ldots,\sigma\left(  i-1\right)  \right\}  \right. \nonumber\\
&  \ \ \ \ \ \ \ \ \ \ \ \ \ \ \ \ \ \ \ \ \left.  \text{that are smaller than
}\sigma\left(  i\right)  \right)  . \label{pf.thm.perm.lehmer.bij.1st.1}%
\end{align}
Therefore,%
\begin{align}
\sigma\left(  i\right)   &  =\left(  \text{the }\left(  \ell_{i}\left(
\sigma\right)  +1\right)  \text{-st smallest element of}\right. \nonumber\\
&  \ \ \ \ \ \ \ \ \ \ \ \ \ \ \ \ \ \ \ \ \left.  \text{the set }\left[
n\right]  \setminus\left\{  \sigma\left(  1\right)  ,\sigma\left(  2\right)
,\ldots,\sigma\left(  i-1\right)  \right\}  \right)
\label{pf.thm.perm.lehmer.bij.1st.2}%
\end{align}
(since $\sigma\left(  i\right)  $ is an element of $\left[  n\right]
\setminus\left\{  \sigma\left(  1\right)  ,\sigma\left(  2\right)
,\ldots,\sigma\left(  i-1\right)  \right\}  $).

Forget that we fixed $i$. We thus have proved the equality
(\ref{pf.thm.perm.lehmer.bij.1st.2}) for each $i\in\left[  n\right]  $. This
equality (\ref{pf.thm.perm.lehmer.bij.1st.2}) allows us to find $\sigma\left(
i\right)  $ if $\ell_{i}\left(  \sigma\right)  $ and $\sigma\left(  1\right)
,\sigma\left(  2\right)  ,\ldots,\sigma\left(  i-1\right)  $ are known. This
can be used to recover $\sigma$ from $L\left(  \sigma\right)  $. Let us
perform this recovery in an example: Let $n=5$, and let $\sigma\in S_{5}$
satisfy $L\left(  \sigma\right)  =\left(  3,1,2,1,0\right)  $. What is
$\sigma$ ?

From $L\left(  \sigma\right)  =\left(  3,1,2,1,0\right)  $, we obtain
$\ell_{1}\left(  \sigma\right)  =3$ and $\ell_{2}\left(  \sigma\right)  =1$
and $\ell_{3}\left(  \sigma\right)  =2$ and $\ell_{4}\left(  \sigma\right)
=1$ and $\ell_{5}\left(  \sigma\right)  =0$.

Applying (\ref{pf.thm.perm.lehmer.bij.1st.2}) to $i=1$, we find%
\begin{align*}
\sigma\left(  1\right)   &  =\left(  \text{the }\left(  \ell_{1}\left(
\sigma\right)  +1\right)  \text{-st smallest element of}\right. \\
&  \ \ \ \ \ \ \ \ \ \ \ \ \ \ \ \ \ \ \ \ \left.  \text{the set }\left[
n\right]  \setminus\underbrace{\left\{  \sigma\left(  1\right)  ,\sigma\left(
2\right)  ,\ldots,\sigma\left(  1-1\right)  \right\}  }_{=\varnothing}\right)
\\
&  =\left(  \text{the }\left(  3+1\right)  \text{-st smallest element of the
set }\left[  5\right]  \right) \\
&  \ \ \ \ \ \ \ \ \ \ \ \ \ \ \ \ \ \ \ \ \left(  \text{since }\ell
_{1}\left(  \sigma\right)  =3\text{ and }n=5\right) \\
&  =\left(  \text{the }4\text{-th smallest element of the set }\left[
5\right]  \right)  =4.
\end{align*}

Applying (\ref{pf.thm.perm.lehmer.bij.1st.2}) to $i=2$, we find%
\begin{align*}
\sigma\left(  2\right)   &  =\left(  \text{the }\left(  \ell_{2}\left(
\sigma\right)  +1\right)  \text{-st smallest element of}\right. \\
&  \ \ \ \ \ \ \ \ \ \ \ \ \ \ \ \ \ \ \ \ \left.  \text{the set }\left[
n\right]  \setminus\underbrace{\left\{  \sigma\left(  1\right)  ,\sigma\left(
2\right)  ,\ldots,\sigma\left(  2-1\right)  \right\}  }_{=\left\{
\sigma\left(  1\right)  \right\}  }\right) \\
&  =\left(  \text{the }\left(  1+1\right)  \text{-st smallest element of the
set }\left[  5\right]  \setminus\left\{  4\right\}  \right) \\
&  \ \ \ \ \ \ \ \ \ \ \ \ \ \ \ \ \ \ \ \ \left(  \text{since }\ell
_{2}\left(  \sigma\right)  =1\text{ and }n=5\text{ and }\sigma\left(
1\right)  =4\right) \\
&  =\left(  \text{the }2\text{-nd smallest element of the set }\left[
5\right]  \setminus\left\{  4\right\}  \right)  =2.
\end{align*}

Applying (\ref{pf.thm.perm.lehmer.bij.1st.2}) to $i=3$, we find%
\begin{align*}
\sigma\left(  3\right)   &  =\left(  \text{the }\left(  \ell_{3}\left(
\sigma\right)  +1\right)  \text{-st smallest element of}\right. \\
&  \ \ \ \ \ \ \ \ \ \ \ \ \ \ \ \ \ \ \ \ \left.  \text{the set }\left[
n\right]  \setminus\underbrace{\left\{  \sigma\left(  1\right)  ,\sigma\left(
2\right)  ,\ldots,\sigma\left(  3-1\right)  \right\}  }_{=\left\{
\sigma\left(  1\right)  ,\sigma\left(  2\right)  \right\}  }\right) \\
&  =\left(  \text{the }\left(  2+1\right)  \text{-st smallest element of the
set }\left[  5\right]  \setminus\left\{  4,2\right\}  \right) \\
&  \ \ \ \ \ \ \ \ \ \ \ \ \ \ \ \ \ \ \ \ \left(  \text{since }\ell
_{3}\left(  \sigma\right)  =2\text{ and }n=5\text{ and }\sigma\left(
1\right)  =4\text{ and }\sigma\left(  2\right)  =2\right) \\
&  =\left(  \text{the }3\text{-rd smallest element of the set }\left[
5\right]  \setminus\left\{  4,2\right\}  \right)  =5
\end{align*}
(since $\left[  5\right]  \setminus\left\{  4,2\right\}  =\left\{
1,3,5\right\}  $).

Continuing like this, we find $\sigma\left(  4\right)  =3$ and $\sigma\left(
5\right)  =1$. Thus, the OLN of $\sigma$ is $\sigma\left(  1\right)
\ \sigma\left(  2\right)  \ \sigma\left(  3\right)  \ \sigma\left(  4\right)
\ \sigma\left(  5\right)  =42531$.

This method allows us to reconstruct any $\sigma\in S_{n}$ from $L\left(
\sigma\right)  $ (and thus shows that the map $L$ is injective). We shall now
see what happens if we apply it to an \textbf{arbitrary} $n$-tuple
$\mathbf{j}=\left(  j_{1},j_{2},\ldots,j_{n}\right)  \in H_{n}$ instead of
$L\left(  \sigma\right)  $ (that is, we replace $\ell_{i}\left(
\sigma\right)  $ by $j_{i}$).

Thus, we define a map%
\[
M:H_{n}\rightarrow S_{n}%
\]
as follows: If $\mathbf{j}=\left(  j_{1},j_{2},\ldots,j_{n}\right)  \in H_{n}%
$, then $M\left(  \mathbf{j}\right)  $ is the map $\sigma:\left[  n\right]
\rightarrow\left[  n\right]  $ whose values $\sigma\left(  1\right)
,\sigma\left(  2\right)  ,\ldots,\sigma\left(  n\right)  $ are defined
recursively by the rule%
\begin{align}
\sigma\left(  i\right)   &  =\left(  \text{the }\left(  j_{i}+1\right)
\text{-st smallest element of}\right. \nonumber\\
&  \ \ \ \ \ \ \ \ \ \ \ \ \ \ \ \ \ \ \ \ \left.  \text{the set }\left[
n\right]  \setminus\left\{  \sigma\left(  1\right)  ,\sigma\left(  2\right)
,\ldots,\sigma\left(  i-1\right)  \right\}  \right)  .
\label{pf.thm.perm.lehmer.bij.1st.6}%
\end{align}
This map $\sigma$

\begin{itemize}
\item is always well-defined (indeed, we never run out of values in the
process of constructing $\sigma$, because of the following argument: each
$i\in\left[  n\right]  $ satisfies $j_{i}\leq n-i$\ \ \ \ \footnote{because
$\left(  j_{1},j_{2},\ldots,j_{n}\right)  \in H_{n}=\left[  n-1\right]
_{0}\times\left[  n-2\right]  _{0}\times\cdots\times\left[  n-n\right]  _{0}$
entails $j_{i}\in\left[  n-i\right]  _{0}$}, and thus the $\left(
n-i+1\right)  $-element set $\left[  n\right]  \setminus\left\{  \sigma\left(
1\right)  ,\sigma\left(  2\right)  ,\ldots,\sigma\left(  i-1\right)  \right\}
$ has a $\left(  j_{i}+1\right)  $-st smallest element\footnote{To be fully
precise: We don't know yet that $\left[  n\right]  \setminus\left\{
\sigma\left(  1\right)  ,\sigma\left(  2\right)  ,\ldots,\sigma\left(
i-1\right)  \right\}  $ is an $\left(  n-i+1\right)  $-element set, since we
haven't yet shown that the $i-1$ elements $\sigma\left(  1\right)
,\sigma\left(  2\right)  ,\ldots,\sigma\left(  i-1\right)  $ are distinct.
However, if they are not distinct, then the set $\left[  n\right]
\setminus\left\{  \sigma\left(  1\right)  ,\sigma\left(  2\right)
,\ldots,\sigma\left(  i-1\right)  \right\}  $ has more than $n-i+1$ elements,
which is just as good for our argument.}), and

\item always is a permutation of $\left[  n\right]  $ (indeed, our definition
of $\sigma\left(  i\right)  $ ensures that $\sigma\left(  i\right)  \in\left[
n\right]  \setminus\left\{  \sigma\left(  1\right)  ,\sigma\left(  2\right)
,\ldots,\sigma\left(  i-1\right)  \right\}  $, so that $\sigma\left(
i\right)  \notin\left\{  \sigma\left(  1\right)  ,\sigma\left(  2\right)
,\ldots,\sigma\left(  i-1\right)  \right\}  $, and therefore the map $\sigma$
has no two equal values; thus, $\sigma$ is injective; therefore, by the
Pigeonhole Principle for Injections\footnote{The \emph{Pigeonhole Principle
for Injections} says the following two things:
\par
\begin{itemize}
\item If $f:X\rightarrow Y$ is an injective map between two finite sets $X$
and $Y$, then $\left\vert X\right\vert \leq\left\vert Y\right\vert $.
\par
\item If $f:X\rightarrow Y$ is an injective map between two finite sets $X$
and $Y$ of the same size, then $f$ is bijective.
\end{itemize}
\par
We are here using the second statement.}, the map $\sigma$ must also be
bijective, i.e., a permutation of $\left[  n\right]  $).
\end{itemize}

\noindent Thus, the map $M$ is well-defined.

We claim that the maps $L$ and $M$ are mutually inverse. Indeed, we already
know that $M$ undoes $L$ (since applying $M$ to $L\left(  \sigma\right)  $
produces precisely our above-discussed algorithm for reconstructing $\sigma$
from $L\left(  \sigma\right)  $); in other words, we have $M\circ
L=\operatorname*{id}$. It is also easy to see that $L\circ
M=\operatorname*{id}$. Thus, the maps $L$ and $M$ are inverse, so that $L$ is
bijective. This proves Theorem \ref{thm.perm.lehmer.bij}.
\end{proof}

Our second proof of Theorem \ref{thm.perm.lehmer.bij} will be less
algorithmic, but it provides a good illustration for the use of total orders.
We will only outline it; the details can be found in \cite[solution to
Exercise 5.18]{detnotes}.

This second proof relies on a total order that can be defined on the set
$\mathbb{Z}^{n}$ of $n$-tuples of integers:

\begin{definition}
\label{def.perm.lehmer.lex-ord}Let $\left(  a_{1},a_{2},\ldots,a_{n}\right)  $
and $\left(  b_{1},b_{2},\ldots,b_{n}\right)  $ be two $n$-tuples of integers.
We say that%
\begin{equation}
\left(  a_{1},a_{2},\ldots,a_{n}\right)  <_{\operatorname*{lex}}\left(
b_{1},b_{2},\ldots,b_{n}\right)  \label{eq.def.perm.lehmer.lex-ord.ineq}%
\end{equation}
if and only if

\begin{itemize}
\item there exists some $k\in\left[  n\right]  $ such that $a_{k}\neq b_{k}$, and

\item the \textbf{smallest} such $k$ satisfies $a_{k}<b_{k}$.
\end{itemize}
\end{definition}

For example, $\left(  4,1,2,5\right)  <_{\operatorname*{lex}}\left(
4,1,3,0\right)  $ and $\left(  1,1,0,1\right)  <_{\operatorname*{lex}}\left(
2,0,0,0\right)  $. The relation (\ref{eq.def.perm.lehmer.lex-ord.ineq}) is
usually pronounced as \textquotedblleft$\left(  a_{1},a_{2},\ldots
,a_{n}\right)  $ is \emph{lexicographically smaller} than $\left(  b_{1}%
,b_{2},\ldots,b_{n}\right)  $\textquotedblright; the word \textquotedblleft
lexicographic\textquotedblright\ comes from the idea that if numbers were
letters, then a \textquotedblleft word\textquotedblright\ $a_{1}a_{2}\cdots
a_{n}$ would appear earlier in a dictionary than $b_{1}b_{2}\cdots b_{n}$ if
and only if $\left(  a_{1},a_{2},\ldots,a_{n}\right)  <_{\operatorname*{lex}%
}\left(  b_{1},b_{2},\ldots,b_{n}\right)  $.

Now, the following is easy to see:

\begin{proposition}
\label{prop.perm.lehmer.lex-ord.total}If $\mathbf{a}$ and $\mathbf{b}$ are two
distinct $n$-tuples of integers, then we have either $\mathbf{a}%
<_{\operatorname*{lex}}\mathbf{b}$ or $\mathbf{b}<_{\operatorname*{lex}%
}\mathbf{a}$.
\end{proposition}

Actually, it is not hard to show that the relation $<_{\operatorname*{lex}}$
is a total order on the set $\mathbb{Z}^{n}$ (known as the \emph{lexicographic
order}); however, Proposition \ref{prop.perm.lehmer.lex-ord.total} is the only
part of this statement that we will need.

\begin{proposition}
\label{prop.perm.lehmer.lex}Let $\sigma\in S_{n}$ and $\tau\in S_{n}$ be such
that
\begin{equation}
\left(  \sigma\left(  1\right)  ,\sigma\left(  2\right)  ,\ldots,\sigma\left(
n\right)  \right)  <_{\operatorname*{lex}}\left(  \tau\left(  1\right)
,\tau\left(  2\right)  ,\ldots,\tau\left(  n\right)  \right)  .
\label{eq.prop.perm.lehmer.lex.ass}%
\end{equation}
Then,
\begin{equation}
\left(  \ell_{1}\left(  \sigma\right)  ,\ell_{2}\left(  \sigma\right)
,\ldots,\ell_{n}\left(  \sigma\right)  \right)  <_{\operatorname*{lex}}\left(
\ell_{1}\left(  \tau\right)  ,\ell_{2}\left(  \tau\right)  ,\ldots,\ell
_{n}\left(  \tau\right)  \right)  . \label{eq.prop.perm.lehmer.lex.clm}%
\end{equation}
(In other words, $L\left(  \sigma\right)  <_{\operatorname*{lex}}L\left(
\tau\right)  $.)
\end{proposition}

\begin{proof}
[Proof of Proposition \ref{prop.perm.lehmer.lex} (sketched).](See
\cite[solution to Exercise 5.18, proof of Proposition 5.50]{detnotes} for
details.) The assumption (\ref{eq.prop.perm.lehmer.lex.ass}) shows that there
exists some $k\in\left[  n\right]  $ such that $\sigma\left(  k\right)
\neq\tau\left(  k\right)  $, and that the \textbf{smallest} such $k$ satisfies
$\sigma\left(  k\right)  <\tau\left(  k\right)  $. Consider this smallest $k$.
Thus, $\sigma\left(  i\right)  =\tau\left(  i\right)  $ for each $i\in\left[
k-1\right]  $ (since $k$ is smallest with $\sigma\left(  k\right)  \neq
\tau\left(  k\right)  $). Hence, using (\ref{pf.thm.perm.lehmer.bij.1st.1}),
we can easily see that%
\begin{equation}
\ell_{i}\left(  \sigma\right)  =\ell_{i}\left(  \tau\right)
\ \ \ \ \ \ \ \ \ \ \text{for each }i\in\left[  k-1\right]  .
\label{pf.prop.perm.lehmer.lex.1}%
\end{equation}
Let $Z$ be the set%
\[
\left[  n\right]  \setminus\left\{  \sigma\left(  1\right)  ,\sigma\left(
2\right)  ,\ldots,\sigma\left(  k-1\right)  \right\}  =\left[  n\right]
\setminus\left\{  \tau\left(  1\right)  ,\tau\left(  2\right)  ,\ldots
,\tau\left(  k-1\right)  \right\}  .
\]
Then, (\ref{pf.thm.perm.lehmer.bij.1st.1}) (applied to $i=k$) yields that
\[
\ell_{k}\left(  \sigma\right)  =\left(  \text{\# of elements of }Z\text{ that
are smaller than }\sigma\left(  k\right)  \right)  ,
\]
and similarly we have%
\[
\ell_{k}\left(  \tau\right)  =\left(  \text{\# of elements of }Z\text{ that
are smaller than }\tau\left(  k\right)  \right)  .
\]
From these two equalities, we can easily see that $\ell_{k}\left(
\sigma\right)  <\ell_{k}\left(  \tau\right)  $. (In fact, any element of $Z$
that is smaller than $\sigma\left(  k\right)  $ must also be smaller than
$\tau\left(  k\right)  $ (since $\sigma\left(  k\right)  <\tau\left(
k\right)  $), but there is at least one element of $Z$ that is smaller than
$\tau\left(  k\right)  $ but not smaller than $\sigma\left(  k\right)  $
(namely, the element $\sigma\left(  k\right)  $). Hence, there are fewer
elements of $Z$ that are smaller than $\sigma\left(  k\right)  $ than there
are elements of $Z$ that are smaller than $\tau\left(  k\right)  $.)

Combining (\ref{pf.prop.perm.lehmer.lex.1}) with $\ell_{k}\left(
\sigma\right)  <\ell_{k}\left(  \tau\right)  $, we obtain $\left(  \ell
_{1}\left(  \sigma\right)  ,\ell_{2}\left(  \sigma\right)  ,\ldots,\ell
_{n}\left(  \sigma\right)  \right)  <_{\operatorname*{lex}}\left(  \ell
_{1}\left(  \tau\right)  ,\ell_{2}\left(  \tau\right)  ,\ldots,\ell_{n}\left(
\tau\right)  \right)  $ (by Definition \ref{def.perm.lehmer.lex-ord}). This
proves Proposition \ref{prop.perm.lehmer.lex}.
\end{proof}

Now, we can easily finish our second proof of Theorem
\ref{thm.perm.lehmer.bij}:

\begin{proof}
[Second proof of Theorem \ref{thm.perm.lehmer.bij} (sketched).]We shall first
show that the map $L$ is injective.

Indeed, let $\sigma$ and $\tau$ be two distinct permutations in $S_{n}$. Then,
the two $n$-tuples $\left(  \sigma\left(  1\right)  ,\sigma\left(  2\right)
,\ldots,\sigma\left(  n\right)  \right)  $ and $\left(  \tau\left(  1\right)
,\tau\left(  2\right)  ,\ldots,\tau\left(  n\right)  \right)  $ (that is, the
OLNs of $\sigma$ and $\tau$) are distinct as well. Hence, Proposition
\ref{prop.perm.lehmer.lex-ord.total} yields that we have either $\left(
\sigma\left(  1\right)  ,\sigma\left(  2\right)  ,\ldots,\sigma\left(
n\right)  \right)  <_{\operatorname*{lex}}\left(  \tau\left(  1\right)
,\tau\left(  2\right)  ,\ldots,\tau\left(  n\right)  \right)  $ or $\left(
\tau\left(  1\right)  ,\tau\left(  2\right)  ,\ldots,\tau\left(  n\right)
\right)  <_{\operatorname*{lex}}\left(  \sigma\left(  1\right)  ,\sigma\left(
2\right)  ,\ldots,\sigma\left(  n\right)  \right)  $. In the first case, we
obtain $L\left(  \sigma\right)  <_{\operatorname*{lex}}L\left(  \tau\right)  $
(by Proposition \ref{prop.perm.lehmer.lex}); in the second case, we likewise
obtain $L\left(  \tau\right)  <_{\operatorname*{lex}}L\left(  \sigma\right)
$. In either case, we thus conclude that $L\left(  \sigma\right)  \neq
L\left(  \tau\right)  $.

Forget that we fixed $\sigma$ and $\tau$. We thus have shown that if $\sigma$
and $\tau$ are two distinct permutations in $S_{n}$, then $L\left(
\sigma\right)  \neq L\left(  \tau\right)  $. In other words, the map
$L:S_{n}\rightarrow H_{n}$ is injective. However, $L$ is a map between two
finite sets of the same size (indeed, $\left\vert S_{n}\right\vert
=n!=\left\vert H_{n}\right\vert $). Thus, the Pigeonhole Principle for
Injections shows that $L$ is bijective (since $L$ is injective). This proves
Theorem \ref{thm.perm.lehmer.bij} again.
\end{proof}

Now, we can prove Proposition \ref{prop.perm.length.gf} at last:

\begin{proof}
[Proof of Proposition \ref{prop.perm.length.gf}.]Each $\sigma\in S_{n}$
satisfies
\begin{align*}
\ell\left(  \sigma\right)   &  =\ell_{1}\left(  \sigma\right)  +\ell
_{2}\left(  \sigma\right)  +\cdots+\ell_{n}\left(  \sigma\right)
\ \ \ \ \ \ \ \ \ \ \left(  \text{by Proposition \ref{prop.perm.lehmer.l}%
}\right) \\
&  =\left(  \text{sum of the entries of }L\left(  \sigma\right)  \right)
\end{align*}
(since $L\left(  \sigma\right)  =\left(  \ell_{1}\left(  \sigma\right)
,\ell_{2}\left(  \sigma\right)  ,\ldots,\ell_{n}\left(  \sigma\right)
\right)  $). Thus,%
\begin{align*}
\sum_{\sigma\in S_{n}}x^{\ell\left(  \sigma\right)  }  &  =\sum_{\sigma\in
S_{n}}x^{\left(  \text{sum of the entries of }L\left(  \sigma\right)  \right)
}=\sum_{\left(  j_{1},j_{2},\ldots,j_{n}\right)  \in H_{n}}%
\underbrace{x^{\left(  \text{sum of the entries of }\left(  j_{1},j_{2}%
,\ldots,j_{n}\right)  \right)  }}_{=x^{j_{1}+j_{2}+\cdots+j_{n}}=x^{j_{1}%
}x^{j_{2}}\cdots x^{j_{n}}}\\
&  \ \ \ \ \ \ \ \ \ \ \ \ \ \ \ \ \ \ \ \ \left(
\begin{array}
[c]{c}%
\text{here, we have substituted }\left(  j_{1},j_{2},\ldots,j_{n}\right)
\text{ for }L\left(  \sigma\right)  \text{ in}\\
\text{the sum, since the map }L:S_{n}\rightarrow H_{n}\text{ is a bijection}%
\end{array}
\right) \\
&  =\sum_{\left(  j_{1},j_{2},\ldots,j_{n}\right)  \in H_{n}}x^{j_{1}}%
x^{j_{2}}\cdots x^{j_{n}}=\sum_{\left(  j_{1},j_{2},\ldots,j_{n}\right)
\in\left[  n-1\right]  _{0}\times\left[  n-2\right]  _{0}\times\cdots
\times\left[  n-n\right]  _{0}}x^{j_{1}}x^{j_{2}}\cdots x^{j_{n}}\\
&  \ \ \ \ \ \ \ \ \ \ \ \ \ \ \ \ \ \ \ \ \left(  \text{since }H_{n}=\left[
n-1\right]  _{0}\times\left[  n-2\right]  _{0}\times\cdots\times\left[
n-n\right]  _{0}\right) \\
&  =\left(  \sum_{j_{1}\in\left[  n-1\right]  _{0}}x^{j_{1}}\right)  \left(
\sum_{j_{2}\in\left[  n-2\right]  _{0}}x^{j_{2}}\right)  \cdots\left(
\sum_{j_{n}\in\left[  n-n\right]  _{0}}x^{j_{n}}\right) \\
&  \ \ \ \ \ \ \ \ \ \ \ \ \ \ \ \ \ \ \ \ \left(  \text{by the product rule
(\ref{eq.lem.fps.prodrule-fin-fin.eq})}\right) \\
&  =\left(  \sum_{j_{1}=0}^{n-1}x^{j_{1}}\right)  \left(  \sum_{j_{2}=0}%
^{n-2}x^{j_{2}}\right)  \cdots\left(  \sum_{j_{n}=0}^{n-n}x^{j_{n}}\right) \\
&  \ \ \ \ \ \ \ \ \ \ \ \ \ \ \ \ \ \ \ \ \left(  \text{since }\left[
m\right]  _{0}=\left\{  0,1,\ldots,m\right\}  \text{ for any }m\in
\mathbb{Z}\right) \\
&  =\left(  1+x+x^{2}+\cdots+x^{n-1}\right)  \left(  1+x+x^{2}+\cdots
+x^{n-2}\right)  \cdots\left(  1+x\right)  \left(  1\right) \\
&  =\left(  1+x+x^{2}+\cdots+x^{n-1}\right)  \left(  1+x+x^{2}+\cdots
+x^{n-2}\right)  \cdots\left(  1+x\right) \\
&  =\left(  1+x\right)  \left(  1+x+x^{2}\right)  \left(  1+x+x^{2}%
+x^{3}\right)  \cdots\left(  1+x+x^{2}+\cdots+x^{n-1}\right) \\
&  =\prod_{i=1}^{n-1}\left(  1+x+x^{2}+\cdots+x^{i}\right)  =\left[  n\right]
_{x}!
\end{align*}
(the last equality sign here is easy to check). This proves Proposition
\ref{prop.perm.length.gf}.
\end{proof}

\subsubsection{More about lengths and simples}

Let us continue studying lengths of permutations.

\begin{proposition}
\label{prop.perm.len.inv}Let $n\in\mathbb{N}$ and $\sigma\in S_{n}$. Then,
$\ell\left(  \sigma^{-1}\right)  =\ell\left(  \sigma\right)  $.
\end{proposition}

\begin{proof}
[Proof of Proposition \ref{prop.perm.len.inv} (sketched).](See \cite[Exercise
5.2 \textbf{(f)}]{detnotes} for details.)

Recall that $\ell\left(  \sigma\right)  $ is the \# of inversions of $\sigma$,
while $\ell\left(  \sigma^{-1}\right)  $ is the \# of inversions of
$\sigma^{-1}$.

Recall also that an inversion of $\sigma$ is a pair $\left(  i,j\right)
\in\left[  n\right]  \times\left[  n\right]  $ such that $i<j$ and
$\sigma\left(  i\right)  >\sigma\left(  j\right)  $. Likewise, an inversion of
$\sigma^{-1}$ is a pair $\left(  u,v\right)  \in\left[  n\right]
\times\left[  n\right]  $ such that $u<v$ and $\sigma^{-1}\left(  u\right)
>\sigma^{-1}\left(  v\right)  $.

Thus, if $\left(  i,j\right)  $ is an inversion of $\sigma$, then $\left(
\sigma\left(  j\right)  ,\sigma\left(  i\right)  \right)  $ is an inversion of
$\sigma^{-1}$. Hence, we obtain a map%
\begin{align*}
\left\{  \text{inversions of }\sigma\right\}   &  \rightarrow\left\{
\text{inversions of }\sigma^{-1}\right\}  ,\\
\left(  i,j\right)   &  \mapsto\left(  \sigma\left(  j\right)  ,\sigma\left(
i\right)  \right)  .
\end{align*}
This map is furthermore bijective (indeed, it has an inverse map, which sends
each $\left(  u,v\right)  \in\left\{  \text{inversions of }\sigma
^{-1}\right\}  $ to $\left(  \sigma^{-1}\left(  v\right)  ,\sigma^{-1}\left(
u\right)  \right)  $). Thus, the bijection principle yields%
\[
\left(  \text{\# of inversions of }\sigma\right)  =\left(  \text{\# of
inversions of }\sigma^{-1}\right)  .
\]
In other words, $\ell\left(  \sigma\right)  =\ell\left(  \sigma^{-1}\right)
$. This proves Proposition \ref{prop.perm.len.inv}.
\end{proof}

The following lemma is crucial for understanding lengths of permutations:

\begin{lemma}
[single swap lemma]\label{lem.perm.len.ssl}Let $n\in\mathbb{N}$, $\sigma\in
S_{n}$ and $k\in\left[  n-1\right]  $. Then: \medskip

\textbf{(a)} We have%
\[
\ell\left(  \sigma s_{k}\right)  =%
\begin{cases}
\ell\left(  \sigma\right)  +1, & \text{if }\sigma\left(  k\right)
<\sigma\left(  k+1\right)  ;\\
\ell\left(  \sigma\right)  -1, & \text{if }\sigma\left(  k\right)
>\sigma\left(  k+1\right)  .
\end{cases}
\]

\textbf{(b)} We have%
\[
\ell\left(  s_{k}\sigma\right)  =%
\begin{cases}
\ell\left(  \sigma\right)  +1, & \text{if }\sigma^{-1}\left(  k\right)
<\sigma^{-1}\left(  k+1\right)  ;\\
\ell\left(  \sigma\right)  -1, & \text{if }\sigma^{-1}\left(  k\right)
>\sigma^{-1}\left(  k+1\right)  .
\end{cases}
\]

[\textbf{Note:} If $i\in\left[  n\right]  $, then $\sigma\left(  i\right)  $
is the \textbf{entry} in position $i$ of the one-line notation of $\sigma$,
whereas $\sigma^{-1}\left(  i\right)  $ is the \textbf{position} in which the
number $i$ appears in the one-line notation of $\sigma$. For example, if
$\sigma=512364$ in one-line notation, then $\sigma\left(  6\right)  =4$ and
$\sigma^{-1}\left(  6\right)  =5$.]
\end{lemma}

We will only outline the proof of Lemma \ref{lem.perm.len.ssl}; a detailed
proof can be found in \cite[Exercise 5.2 \textbf{(a)}]{detnotes} (although it
is a slightly different proof).

\begin{proof}
[Proof of Lemma \ref{lem.perm.len.ssl} (sketched).]\textbf{(b)} The
OLN\footnote{Recall that \textquotedblleft OLN\textquotedblright\ means
\textquotedblleft one-line notation\textquotedblright.} of $s_{k}\sigma$ is
obtained from the OLN of $\sigma$ by swapping the two entries $k$ and $k+1$.
This is best seen on an example: For example, if $\sigma=512364$ (in OLN),
then $s_{3}\sigma=512463$. In general, this follows by observing that%
\[
\left(  s_{k}\sigma\right)  \left(  i\right)  =s_{k}\left(  \sigma\left(
i\right)  \right)  =%
\begin{cases}
k+1, & \text{if }\sigma\left(  i\right)  =k;\\
k, & \text{if }\sigma\left(  i\right)  =k+1;\\
\sigma\left(  i\right)  , & \text{otherwise}%
\end{cases}
\ \ \ \ \ \ \ \ \ \ \text{for each }i\in\left[  n\right]  .
\]

Let us now use this observation to see how the inversions of $s_{k}\sigma$
differ from the inversions of $\sigma$. Indeed, let us call an inversion
$\left(  i,j\right)  $ of a permutation $\tau$ \emph{exceptional} if we have
$\tau\left(  i\right)  =k+1$ and $\tau\left(  j\right)  =k$. All other
inversions of $\tau$ will be called \emph{non-exceptional}.

Now, we make the following observation:

\begin{statement}
\textit{Observation 1:} If $\left(  i,j\right)  $ is any non-exceptional
inversion of $\sigma$, then $\left(  i,j\right)  $ is still a non-exceptional
inversion of $s_{k}\sigma$.
\end{statement}

[\textit{Proof of Observation 1:} Let $\left(  i,j\right)  $ be a
non-exceptional inversion of $\sigma$. Thus, we have the inequality
$\sigma\left(  i\right)  >\sigma\left(  j\right)  $ (since $\left(
i,j\right)  $ is an inversion of $\sigma$). This inequality cannot get
reversed by applying $s_{k}$ to both its sides (i.e., we cannot have
$s_{k}\left(  \sigma\left(  i\right)  \right)  <s_{k}\left(  \sigma\left(
j\right)  \right)  $), since the only pair $\left(  u,v\right)  \in\left[
n\right]  \times\left[  n\right]  $ satisfying $u>v$ and $s_{k}\left(
u\right)  <s_{k}\left(  v\right)  $ is the pair $\left(  k+1,k\right)  $ (but
our pair $\left(  \sigma\left(  i\right)  ,\sigma\left(  j\right)  \right)  $
cannot equal this pair $\left(  k+1,k\right)  $, since the inversion $\left(
i,j\right)  $ of $\sigma$ is non-exceptional). Hence, we have $s_{k}\left(
\sigma\left(  i\right)  \right)  \geq s_{k}\left(  \sigma\left(  j\right)
\right)  $. In other words, $\left(  s_{k}\sigma\right)  \left(  i\right)
\geq\left(  s_{k}\sigma\right)  \left(  j\right)  $. Since the map
$s_{k}\sigma$ is injective (being a permutation of $\left[  n\right]  $), we
thus obtain $\left(  s_{k}\sigma\right)  \left(  i\right)  >\left(
s_{k}\sigma\right)  \left(  j\right)  $ (since $i\neq j$). Thus, the pair
$\left(  i,j\right)  $ is an inversion of $s_{k}\sigma$. Moreover, this
inversion $\left(  i,j\right)  $ is non-exceptional (since otherwise we would
have $\left(  s_{k}\sigma\right)  \left(  i\right)  =k+1$ and $\left(
s_{k}\sigma\right)  \left(  j\right)  =k$, which would lead to $\sigma\left(
i\right)  =k$ and $\sigma\left(  j\right)  =k+1$, which would contradict
$\sigma\left(  i\right)  >\sigma\left(  j\right)  $). Thus, we have shown that
$\left(  i,j\right)  $ is still a non-exceptional inversion of $s_{k}\sigma$.
This proves Observation 1.]

Similarly to Observation 1, we can prove the following:

\begin{statement}
\textit{Observation 2:} If $\left(  i,j\right)  $ is any non-exceptional
inversion of $s_{k}\sigma$, then $\left(  i,j\right)  $ is still a
non-exceptional inversion of $\sigma$.
\end{statement}

(Alternatively, we can obtain Observation 2 by applying Observation 1 to
$s_{k}\sigma$ instead of $\sigma$, since we have $\underbrace{s_{k}s_{k}%
}_{=s_{k}^{2}=\operatorname*{id}}\sigma=\sigma$.)

Combining Observation 1 with Observation 2, we see that the non-exceptional
inversions of $s_{k}\sigma$ are precisely the non-exceptional inversions of
$\sigma$. Hence,%
\begin{align}
&  \left(  \text{\# of non-exceptional inversions of }s_{k}\sigma\right)
\nonumber\\
&  =\left(  \text{\# of non-exceptional inversions of }\sigma\right)  .
\label{pf.lem.perm.len.ssl.eq-nonexc}%
\end{align}

What about the exceptional inversions? A permutation $\tau\in S_{n}$ has a
unique exceptional inversion if $k$ appears after $k+1$ in the OLN of $\tau$
(that is, if we have $\tau^{-1}\left(  k\right)  >\tau^{-1}\left(  k+1\right)
$); otherwise, it has none. Thus:

\begin{itemize}
\item If $\sigma^{-1}\left(  k\right)  <\sigma^{-1}\left(  k+1\right)  $, then
the permutation $s_{k}\sigma$ has a unique exceptional inversion, whereas the
permutation $\sigma$ has none.

\item If $\sigma^{-1}\left(  k\right)  >\sigma^{-1}\left(  k+1\right)  $, then
the permutation $\sigma$ has a unique exceptional inversion, whereas the
permutation $s_{k}\sigma$ has none.
\end{itemize}

Thus,%
\begin{align}
&  \left(  \text{\# of exceptional inversions of }s_{k}\sigma\right)
\nonumber\\
&  =\left(  \text{\# of exceptional inversions of }\sigma\right) \nonumber\\
&  \ \ \ \ \ \ \ \ \ \ +%
\begin{cases}
1, & \text{if }\sigma^{-1}\left(  k\right)  <\sigma^{-1}\left(  k+1\right)
;\\
-1, & \text{if }\sigma^{-1}\left(  k\right)  >\sigma^{-1}\left(  k+1\right)  .
\end{cases}
\label{pf.lem.perm.len.ssl.eq-exc}%
\end{align}

Now, recall that each inversion of a permutation $\tau\in S_{n}$ is either
exceptional or non-exceptional (and cannot be both at the same time). Thus,
adding together the two equalities (\ref{pf.lem.perm.len.ssl.eq-exc}) and
(\ref{pf.lem.perm.len.ssl.eq-nonexc}), we obtain%
\begin{align*}
&  \left(  \text{\# of inversions of }s_{k}\sigma\right) \\
&  =\left(  \text{\# of inversions of }\sigma\right)  +%
\begin{cases}
1, & \text{if }\sigma^{-1}\left(  k\right)  <\sigma^{-1}\left(  k+1\right)
;\\
-1, & \text{if }\sigma^{-1}\left(  k\right)  >\sigma^{-1}\left(  k+1\right)
\end{cases}
\\
&  =%
\begin{cases}
\left(  \text{\# of inversions of }\sigma\right)  +1, & \text{if }\sigma
^{-1}\left(  k\right)  <\sigma^{-1}\left(  k+1\right)  ;\\
\left(  \text{\# of inversions of }\sigma\right)  -1, & \text{if }\sigma
^{-1}\left(  k\right)  >\sigma^{-1}\left(  k+1\right)  .
\end{cases}
\end{align*}
In other words,%
\[
\ell\left(  s_{k}\sigma\right)  =%
\begin{cases}
\ell\left(  \sigma\right)  +1, & \text{if }\sigma^{-1}\left(  k\right)
<\sigma^{-1}\left(  k+1\right)  ;\\
\ell\left(  \sigma\right)  -1, & \text{if }\sigma^{-1}\left(  k\right)
>\sigma^{-1}\left(  k+1\right)
\end{cases}
\]
(since $\ell\left(  \tau\right)  $ denotes the \# of inversions of any
permutation $\tau$). This proves Lemma \ref{lem.perm.len.ssl} \textbf{(b)}.
\medskip

\textbf{(a)} Applying Lemma \ref{lem.perm.len.ssl} \textbf{(b)} to
$\sigma^{-1}$ instead of $\sigma$, we obtain%
\begin{align}
\ell\left(  s_{k}\sigma^{-1}\right)   &  =%
\begin{cases}
\ell\left(  \sigma^{-1}\right)  +1, & \text{if }\left(  \sigma^{-1}\right)
^{-1}\left(  k\right)  <\left(  \sigma^{-1}\right)  ^{-1}\left(  k+1\right)
;\\
\ell\left(  \sigma^{-1}\right)  -1, & \text{if }\left(  \sigma^{-1}\right)
^{-1}\left(  k\right)  >\left(  \sigma^{-1}\right)  ^{-1}\left(  k+1\right)
\end{cases}
\nonumber\\
&  =%
\begin{cases}
\ell\left(  \sigma\right)  +1, & \text{if }\sigma\left(  k\right)
<\sigma\left(  k+1\right)  ;\\
\ell\left(  \sigma\right)  -1, & \text{if }\sigma\left(  k\right)
>\sigma\left(  k+1\right)
\end{cases}
\label{pf.lem.perm.len.ssl.a.1}%
\end{align}
(since $\left(  \sigma^{-1}\right)  ^{-1}=\sigma$, and since Proposition
\ref{prop.perm.len.inv} yields $\ell\left(  \sigma^{-1}\right)  =\ell\left(
\sigma\right)  $). However, Proposition \ref{prop.perm.len.inv} (applied to
$\sigma s_{k}$ instead of $\sigma$) yields $\ell\left(  \left(  \sigma
s_{k}\right)  ^{-1}\right)  =\ell\left(  \sigma s_{k}\right)  $. In view of
$\left(  \sigma s_{k}\right)  ^{-1}=\underbrace{s_{k}^{-1}}_{=s_{k}}%
\sigma^{-1}=s_{k}\sigma^{-1}$, this rewrites as $\ell\left(  s_{k}\sigma
^{-1}\right)  =\ell\left(  \sigma s_{k}\right)  $. Comparing this with
(\ref{pf.lem.perm.len.ssl.a.1}), we obtain%
\[
\ell\left(  \sigma s_{k}\right)  =%
\begin{cases}
\ell\left(  \sigma\right)  +1, & \text{if }\sigma\left(  k\right)
<\sigma\left(  k+1\right)  ;\\
\ell\left(  \sigma\right)  -1, & \text{if }\sigma\left(  k\right)
>\sigma\left(  k+1\right)  .
\end{cases}
\]
This proves Lemma \ref{lem.perm.len.ssl} \textbf{(a)}.
\end{proof}

Lemma \ref{lem.perm.len.ssl} answers what happens to the length of a
permutation when we compose it (from the left or the right) with a simple
transposition $s_{k}$. What happens when we compose it with a non-simple
transposition? The situation is more complicated, but it is still true that
the length decreases or increases depending on whether the two entries that
are being swapped formed an inversion or not. Here is the exact answer (stated
only for $\sigma t_{i,j}$, but a version for $t_{i,j}\sigma$ can easily be
derived from it):

\begin{proposition}
\label{prop.perm.lisitij}Let $n\in\mathbb{N}$ and $\sigma\in S_{n}$. Let $i$
and $j$ be two elements of $\left[  n\right]  $ such that $i<j$. Then:
\medskip

\textbf{(a)} We have $\ell\left(  \sigma t_{i,j}\right)  <\ell\left(
\sigma\right)  $ if $\sigma\left(  i\right)  >\sigma\left(  j\right)  $. We
have $\ell\left(  \sigma t_{i,j}\right)  >\ell\left(  \sigma\right)  $ if
$\sigma\left(  i\right)  <\sigma\left(  j\right)  $. \medskip

\textbf{(b)} We have%
\[
\ell\left(  \sigma t_{i,j}\right)  =%
\begin{cases}
\ell\left(  \sigma\right)  -2\left\vert Q\right\vert -1, & \text{if }%
\sigma\left(  i\right)  >\sigma\left(  j\right)  ;\\
\ell\left(  \sigma\right)  +2\left\vert R\right\vert +1, & \text{if }%
\sigma\left(  i\right)  <\sigma\left(  j\right)  ,
\end{cases}
\]
where
\begin{align*}
Q  &  =\left\{  k\in\left\{  i+1,i+2,\ldots,j-1\right\}  \ \mid\ \sigma\left(
i\right)  >\sigma\left(  k\right)  >\sigma\left(  j\right)  \right\}
\ \ \ \ \ \ \ \ \ \ \text{and}\\
R  &  =\left\{  k\in\left\{  i+1,i+2,\ldots,j-1\right\}  \ \mid\ \sigma\left(
i\right)  <\sigma\left(  k\right)  <\sigma\left(  j\right)  \right\}  .
\end{align*}

\end{proposition}

\begin{proof}
[Proof of Proposition \ref{prop.perm.lisitij} (sketched).]\textbf{(b)} This
follows by a diligent analysis of the possible interactions between an
inversion and composition by $t_{i,j}$. To be more concrete:

\begin{itemize}
\item The fact that $\ell\left(  \sigma t_{i,j}\right)  =\ell\left(
\sigma\right)  -2\left\vert Q\right\vert -1$ when $\sigma\left(  i\right)
>\sigma\left(  j\right)  $ is \cite[Exercise 5.20]{detnotes}. A
straightforward proof was given by Elafandi in \cite[solution to Exercise
1]{18f-hw4se}. (The proof given in \cite[solution to Exercise 5.20]{detnotes}
is more circuitous, as it uses summation tricks to bypass case distinctions.)

\item The fact that $\ell\left(  \sigma t_{i,j}\right)  =\ell\left(
\sigma\right)  +2\left\vert R\right\vert +1$ when $\sigma\left(  i\right)
<\sigma\left(  j\right)  $ follows by applying the previous fact to $\sigma
t_{i,j}$ instead of $\sigma$. (Indeed, if $\sigma\left(  i\right)
<\sigma\left(  j\right)  $, then $\left(  \sigma t_{i,j}\right)  \left(
i\right)  =\sigma\left(  j\right)  >\sigma\left(  i\right)  =\left(  \sigma
t_{i,j}\right)  \left(  j\right)  $ and $\sigma\underbrace{t_{i,j}t_{i,j}%
}_{=t_{i,j}^{2}=\operatorname*{id}}=\sigma$. Moreover, when we replace
$\sigma$ by $\sigma t_{i,j}$, the sets $Q$ and $R$ trade places.)
\end{itemize}

\textbf{(a)} This follows immediately from part \textbf{(b)}.
\end{proof}

Now we come to one of the main facts about permutations of a finite set.

\begin{convention}
\label{conv.perm.simple}We recall that a \emph{simple transposition} in
$S_{n}$ means one of the $n-1$ transpositions $s_{1},s_{2},\ldots,s_{n-1}$. We
shall occasionally abbreviate \textquotedblleft simple
transposition\textquotedblright\ as \textquotedblleft\emph{simple}%
\textquotedblright.
\end{convention}

\begin{theorem}
[1st reduced word theorem for the symmetric group]%
\label{thm.perm.len.redword1}Let $n\in\mathbb{N}$ and $\sigma\in S_{n}$. Then:
\medskip

\textbf{(a)} We can write $\sigma$ as a composition (i.e., product) of
$\ell\left(  \sigma\right)  $ simples. \medskip

\textbf{(b)} The number $\ell\left(  \sigma\right)  $ is the smallest
$p\in\mathbb{N}$ such that we can write $\sigma$ as a composition of $p$
simples. \medskip

[\textbf{Keep in mind:} The composition of $0$ simples is $\operatorname*{id}%
$, since $\operatorname*{id}$ is the neutral element of the group $S_{n}$.]
\end{theorem}

\begin{example}
Let $\sigma\in S_{4}$ be the permutation $4132$ (in OLN). How can we write
$\sigma$ as a composition of simples? There are several ways to do this; for
example,%
\[
\sigma=\underbrace{s_{2}s_{3}s_{2}}_{=s_{3}s_{2}s_{3}}s_{1}=s_{3}%
s_{2}\underbrace{s_{3}s_{1}}_{=s_{1}s_{3}}=s_{3}s_{2}s_{1}s_{3}=s_{2}%
s_{1}s_{1}s_{3}s_{2}s_{1}=s_{2}s_{1}s_{3}s_{1}s_{2}s_{1}=\cdots.
\]
The shortest of these representations involve $4$ simples, precisely as
predicted by Theorem \ref{thm.perm.len.redword1} (since $\ell\left(
\sigma\right)  =4$).
\end{example}

Before we prove Theorem \ref{thm.perm.len.redword1}, let me mention a
geometric visualization of the symmetric group that will not be used in what
follows, but sheds some light on the theorem and on the role of simple transpositions:

\begin{remark}
Let $n>0$. Then, the permutations of $\left[  n\right]  $ can be represented
as the vertices of an $\left(  n-1\right)  $-dimensional polytope in
$n$-dimensional space.

Namely, each permutation $\sigma$ of $\left[  n\right]  $ gives rise to the
point%
\[
V_{\sigma}:=\left(  \sigma\left(  1\right)  ,\sigma\left(  2\right)
,\ldots,\sigma\left(  n\right)  \right)  \in\mathbb{R}^{n}.
\]
The convex hull of all these $n!$ many points $V_{\sigma}$ (for $\sigma\in
S_{n}$) is a polytope (i.e., a bounded convex polyhedron in $\mathbb{R}^{n}$).
This polytope is known as the
\emph{\href{https://en.wikipedia.org/wiki/Permutohedron}{\emph{permutahedron}%
}} (corresponding to $n$). It is actually $\left(  n-1\right)  $-dimensional,
since all its vertices lie on the hyperplane with equation $x_{1}+x_{2}%
+\cdots+x_{n}=1+2+\cdots+n$. It can be shown (see, e.g., \cite{GaiGup77}) that:

\begin{itemize}
\item The vertices of this polytope are precisely the $n!$ points $V_{\sigma}$
with $\sigma\in S_{n}$.

\item Two vertices $V_{\sigma}$ and $V_{\tau}$ are joined by an edge if and
only if $\sigma=s_{k}\tau$ for some $k\in\left[  n-1\right]  $.
\end{itemize}

\noindent The (intuitively obvious) fact that any two vertices of a polytope
can be connected by a sequence of edges therefore yields that any $\sigma\in
S_{n}$ can be written as a product of simples. This is a weaker version of
Theorem \ref{thm.perm.len.redword1} \textbf{(a)}.

We refer to textbooks on discrete geometry and geometric combinatorics for
more about polytopes and permutahedra in particular. Let me here just show the
permutahedra for $n=3$ and for $n=4$ (note that the permutahedron for $n=2$ is
a boring line segment in $\mathbb{R}^{2}$):

\begin{itemize}
\item The permutahedron for $n=3$ is a regular hexagon:%
\[%
\begin{tabular}
[c]{ccc}%
$%
\begin{tikzpicture}
\begin{axis}[
view={110}{30},
xmin=1, xmax=3,
ymin=1, ymax=3,
zmin=1, zmax=3,
axis equal image
]
\addplot3 [fill=cyan, fill opacity=0.5] table {
x y z
1 2 3
2 1 3
3 1 2
3 2 1
2 3 1
1 3 2
}--cycle;
\end{axis}
\end{tikzpicture}%
$ & $\ \ \ \ \ \ \ \ \ \ $ & $%
\begin{tikzpicture}
\filldraw
[very thick, draw=blue, fill=blue!20!white] (30 : 2) -- (90 : 2) -- (150 : 2) -- (-150 : 2) -- (-90 : 2) -- (-30 : 2) -- cycle;
\fill(  30 : 2) circle (2pt) node[right] {$(3,1,2)$};
\fill(  90 : 2) circle (2pt) node[above] {$(3,2,1)$};
\fill( 150 : 2) circle (2pt) node[left] {$(2,3,1)$};
\fill( -30 : 2) circle (2pt) node[right] {$(2,1,3)$};
\fill( -90 : 2) circle (2pt) node[below] {$(1,2,3)$};
\fill(-150 : 2) circle (2pt) node[left] {$(1,3,2)$};
\end{tikzpicture}%
$%
\end{tabular}
\
\]
The picture on the left
(\href{https://tex.stackexchange.com/a/321238/}{courtesy of tex.stackexchange
user Jake}, released under the MIT license) shows the permutahedron embedded
in $\mathbb{R}^{3}$; the picture on the right is a view from an orthogonal direction.

\item The permutahedron for $n=4$ is a truncated octahedron:%
\[%
{\def\svgwidth{\columnwidth} 
\begingroup%
  \makeatletter%
  \providecommand\color[2][]{%
    \errmessage{(Inkscape) Color is used for the text in Inkscape, but the package 'color.sty' is not loaded}%
    \renewcommand\color[2][]{}%
  }%
  \providecommand\transparent[1]{%
    \errmessage{(Inkscape) Transparency is used (non-zero) for the text in Inkscape, but the package 'transparent.sty' is not loaded}%
    \renewcommand\transparent[1]{}%
  }%
  \providecommand\rotatebox[2]{#2}%
  \ifx\svgwidth\undefined%
    \setlength{\unitlength}{515.25bp}%
    \ifx\svgscale\undefined%
      \relax%
    \else%
      \setlength{\unitlength}{\unitlength * \real{\svgscale}}%
    \fi%
  \else%
    \setlength{\unitlength}{\svgwidth}%
  \fi%
  \global\let\svgwidth\undefined%
  \global\let\svgscale\undefined%
  \makeatother%
  \begin{picture}(1,0.80349345)%
    \put(0,0){\includegraphics[width=\unitlength,page=1]{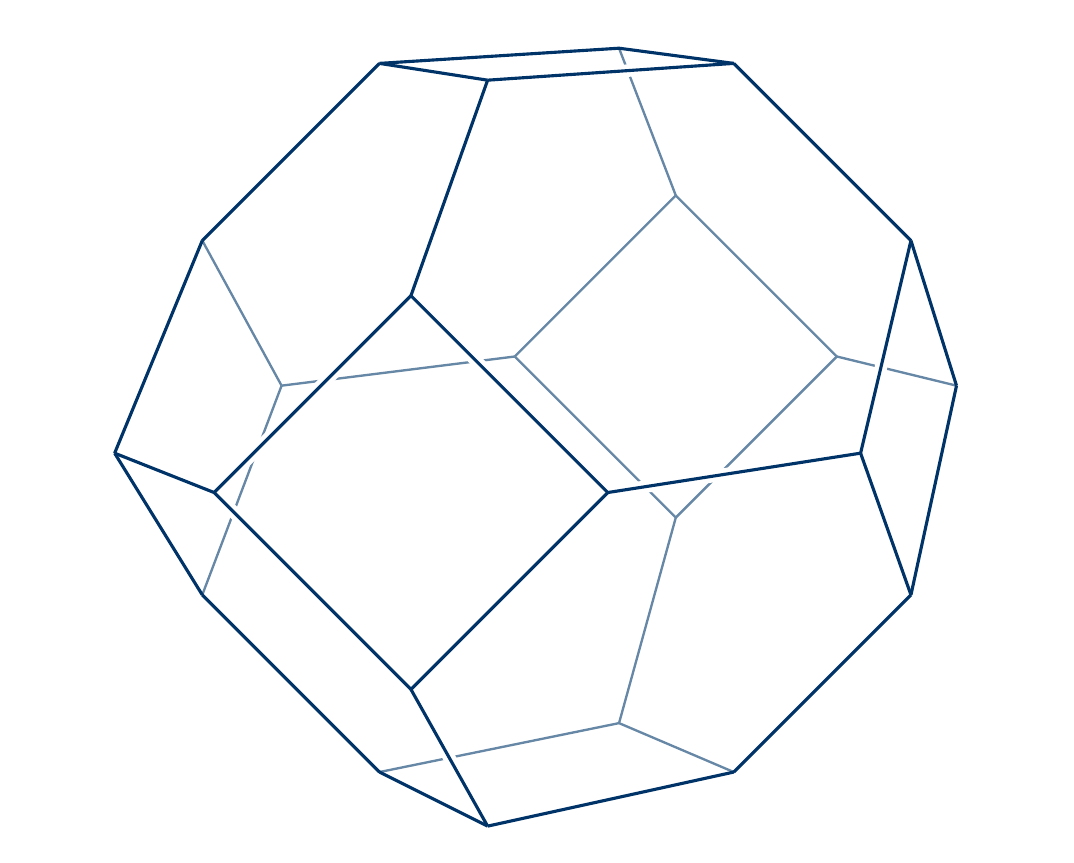}}%
    \put(0.49600568,0.77420655){\makebox(0,0)[lb]{\smash{(4,1,2,3)}}}%
    \put(0.69105662,0.74800568){\makebox(0,0)[lb]{\smash{(4,2,1,3)}}}%
    \put(0.45815997,0.70288195){\makebox(0,0)[lb]{\smash{(3,2,1,4)}}}%
    \put(0.27329825,0.7552837){\makebox(0,0)[lb]{\smash{(3,1,2,4)}}}%
    \put(0.09280335,0.58352242){\makebox(0,0)[lb]{\smash{(2,1,3,4)}}}%
    \put(0.00228952,0.36809301){\makebox(0,0)[lb]{\smash{(1,2,3,4)}}}%
    \put(0.08080335,0.23272183){\makebox(0,0)[lb]{\smash{(1,2,4,3)}}}%
    \put(0.22089651,0.33898093){\makebox(0,0)[lb]{\smash{(1,3,2,4)}}}%
    \put(0.15912795,0.4350508){\makebox(0,0)[lb]{\smash{(2,1,4,3)}}}%
    \put(0.27894309,0.53112067){\makebox(0,0)[lb]{\smash{(2,3,1,4)}}}%
    \put(0.49891689,0.46561849){\makebox(0,0)[lb]{\smash{(3,1,4,2)}}}%
    \put(0.51802897,0.61991252){\makebox(0,0)[lb]{\smash{(4,1,3,2)}}}%
    \put(0.66504498,0.46416288){\makebox(0,0)[lb]{\smash{(4,2,3,1)}}}%
    \put(0.64156608,0.30259083){\makebox(0,0)[lb]{\smash{(3,2,4,1)}}}%
    \put(0.44961557,0.33752533){\makebox(0,0)[lb]{\smash{(2,4,1,3)}}}%
    \put(0.39830262,0.15066361){\makebox(0,0)[lb]{\smash{(1,4,2,3)}}}%
    \put(0.2459198,0.06823857){\makebox(0,0)[lb]{\smash{(1,3,4,2)}}}%
    \put(0.58916434,0.13374076){\makebox(0,0)[lb]{\smash{(2,3,4,1)}}}%
    \put(0.42759229,0.01292562){\makebox(0,0)[lb]{\smash{(1,4,3,2)}}}%
    \put(0.69105662,0.06532737){\makebox(0,0)[lb]{\smash{(2,4,3,1)}}}%
    \put(0.85117307,0.22981063){\makebox(0,0)[lb]{\smash{(3,4,2,1)}}}%
    \put(0.89920801,0.44232882){\makebox(0,0)[lb]{\smash{(4,3,2,1)}}}%
    \put(0.8599067,0.5762444){\makebox(0,0)[lb]{\smash{(4,3,1,2)}}}%
    \put(0.81332737,0.37391543){\makebox(0,0)[lb]{\smash{(3,4,1,2)}}}%
    \put(0,0){\includegraphics[width=\unitlength,page=2]{Permutohedron.pdf}}%
  \end{picture}%
\endgroup%
 }%
\]
(\href{https://en.wikipedia.org/wiki/File:Permutohedron.svg}{picture courtesy
of David Eppstein on the Wikipedia}).
\end{itemize}
\end{remark}

Let us now outline a proof of Theorem \ref{thm.perm.len.redword1}.

\begin{proof}
[Proof of Theorem \ref{thm.perm.len.redword1} (sketched).]\textbf{(a)} (See
\cite[Exercise 5.2 \textbf{(e)}]{detnotes} for details.)

We proceed by induction on $\ell\left(  \sigma\right)  $:

\textit{Induction base:} If $\ell\left(  \sigma\right)  =0$, then
$\sigma=\operatorname*{id}$, so that we can write $\sigma$ as a composition of
$0$ simples.

\textit{Induction step:} Fix $h\in\mathbb{N}$. Assume (as the
IH\footnote{\textquotedblleft IH\textquotedblright\ means \textquotedblleft
induction hypothesis\textquotedblright.}) that Theorem
\ref{thm.perm.len.redword1} \textbf{(a)} holds for $\ell\left(  \sigma\right)
=h$.

Now, let $\sigma\in S_{n}$ be such that $\ell\left(  \sigma\right)  =h+1$. We
must prove that we can write $\sigma$ as a composition of $\ell\left(
\sigma\right)  $ simples.

We have $\ell\left(  \sigma\right)  =h+1>h\geq0$; hence, $\sigma$ has at least
one inversion. Thus, we cannot have $\sigma\left(  1\right)  \leq\sigma\left(
2\right)  \leq\cdots\leq\sigma\left(  n\right)  $. In other words, there
exists some $k\in\left[  n-1\right]  $ such that $\sigma\left(  k\right)
>\sigma\left(  k+1\right)  $. Let us fix such a $k$.

Lemma \ref{lem.perm.len.ssl} \textbf{(a)} yields%
\begin{align*}
\ell\left(  \sigma s_{k}\right)   &  =%
\begin{cases}
\ell\left(  \sigma\right)  +1, & \text{if }\sigma\left(  k\right)
<\sigma\left(  k+1\right)  ;\\
\ell\left(  \sigma\right)  -1, & \text{if }\sigma\left(  k\right)
>\sigma\left(  k+1\right)
\end{cases}
\\
&  =\ell\left(  \sigma\right)  -1\ \ \ \ \ \ \ \ \ \ \left(  \text{since
}\sigma\left(  k\right)  >\sigma\left(  k+1\right)  \right) \\
&  =h\ \ \ \ \ \ \ \ \ \ \left(  \text{since }\ell\left(  \sigma\right)
=h+1\right)  .
\end{align*}
Thus, the IH tells us that we can apply Theorem \ref{thm.perm.len.redword1}
\textbf{(a)} to $\sigma s_{k}$ instead of $\sigma$. We conclude that we can
write $\sigma s_{k}$ as a composition of $\ell\left(  \sigma s_{k}\right)  $
simples. In other words, we can write $\sigma s_{k}$ as a composition of $h$
simples (since $\ell\left(  \sigma s_{k}\right)  =h$). That is, we have%
\[
\sigma s_{k}=s_{i_{1}}s_{i_{2}}\cdots s_{i_{h}}\ \ \ \ \ \ \ \ \ \ \text{for
some }i_{1},i_{2},\ldots,i_{h}\in\left[  n-1\right]  .
\]
Consider these $i_{1},i_{2},\ldots,i_{h}$. Now,%
\[
\sigma=\underbrace{\sigma s_{k}}_{=s_{i_{1}}s_{i_{2}}\cdots s_{i_{h}}%
}\underbrace{s_{k}^{-1}}_{=s_{k}}=s_{i_{1}}s_{i_{2}}\cdots s_{i_{h}}s_{k}.
\]
This shows that we can write $\sigma$ as a composition of $h+1$ simples. In
other words, we can write $\sigma$ as a composition of $\ell\left(
\sigma\right)  $ simples (since $\ell\left(  \sigma\right)  =h+1$). Thus,
Theorem \ref{thm.perm.len.redword1} \textbf{(a)} holds for $\ell\left(
\sigma\right)  =h+1$. This completes the induction step, and so Theorem
\ref{thm.perm.len.redword1} \textbf{(a)} is proved by induction. \medskip

\textbf{(b)} (See \cite[Exercise 5.2 \textbf{(g)}]{detnotes} for details.)

We already know from Theorem \ref{thm.perm.len.redword1} \textbf{(a)} that we
can write $\sigma$ as a composition of $\ell\left(  \sigma\right)  $ simples.
It thus remains to show that we cannot write $\sigma$ as a composition of
fewer than $\ell\left(  \sigma\right)  $ simples.

This will clearly follow if we can show that%
\begin{equation}
\ell\left(  s_{i_{1}}s_{i_{2}}\cdots s_{i_{g}}\right)  \leq g
\label{pf.thm.perm.len.redword1.b.1}%
\end{equation}
for any $i_{1},i_{2},\ldots,i_{g}\in\left[  n-1\right]  $.

In order to prove (\ref{pf.thm.perm.len.redword1.b.1}), we first make a simple
observation: For any $\sigma\in S_{n}$ and any $k\in\left[  n-1\right]  $, we
have%
\begin{equation}
\ell\left(  \sigma s_{k}\right)  \leq\ell\left(  \sigma\right)  +1.
\label{pf.thm.perm.len.redword1.b.2}%
\end{equation}
(This follows from Lemma \ref{lem.perm.len.ssl} \textbf{(a)}, since both
numbers $\ell\left(  \sigma\right)  +1$ and $\ell\left(  \sigma\right)  -1$
are $\leq\ell\left(  \sigma\right)  +1$.) Now, for any $i_{1},i_{2}%
,\ldots,i_{g}\in\left[  n-1\right]  $, we have%
\begin{align*}
&  \ell\left(  s_{i_{1}}s_{i_{2}}\cdots s_{i_{g}}\right) \\
&  \leq\ell\left(  s_{i_{1}}s_{i_{2}}\cdots s_{i_{g-1}}\right)
+1\ \ \ \ \ \ \ \ \ \ \left(  \text{by (\ref{pf.thm.perm.len.redword1.b.2}%
)}\right) \\
&  \leq\ell\left(  s_{i_{1}}s_{i_{2}}\cdots s_{i_{g-2}}\right)  +2\\
&  \ \ \ \ \ \ \ \ \ \ \left(  \text{since (\ref{pf.thm.perm.len.redword1.b.2}%
) yields }\ell\left(  s_{i_{1}}s_{i_{2}}\cdots s_{i_{g-1}}\right)  \leq
\ell\left(  s_{i_{1}}s_{i_{2}}\cdots s_{i_{g-2}}\right)  +1\right) \\
&  \leq\ell\left(  s_{i_{1}}s_{i_{2}}\cdots s_{i_{g-3}}\right)  +3\\
&  \ \ \ \ \ \ \ \ \ \ \left(  \text{since (\ref{pf.thm.perm.len.redword1.b.2}%
) yields }\ell\left(  s_{i_{1}}s_{i_{2}}\cdots s_{i_{g-2}}\right)  \leq
\ell\left(  s_{i_{1}}s_{i_{2}}\cdots s_{i_{g-3}}\right)  +1\right) \\
&  \leq\cdots\\
&  \leq\underbrace{\ell\left(  \operatorname*{id}\right)  }_{=0}%
+\,g\ \ \ \ \ \ \ \ \ \ \left(  \text{since the product of }0\text{ simples is
}\operatorname*{id}\right) \\
&  =g.
\end{align*}
This proves (\ref{pf.thm.perm.len.redword1.b.1}), and thus concludes our proof
of Theorem \ref{thm.perm.len.redword1} \textbf{(b)}.
\end{proof}

\begin{corollary}
\label{cor.perm.red.sigtau}Let $n\in\mathbb{N}$. \medskip

\textbf{(a)} We have $\ell\left(  \sigma\tau\right)  \equiv\ell\left(
\sigma\right)  +\ell\left(  \tau\right)  \operatorname{mod}2$ for all
$\sigma\in S_{n}$ and $\tau\in S_{n}$. \medskip

\textbf{(b)} We have $\ell\left(  \sigma\tau\right)  \leq\ell\left(
\sigma\right)  +\ell\left(  \tau\right)  $ for all $\sigma\in S_{n}$ and
$\tau\in S_{n}$. \medskip

\textbf{(c)} Let $k_{1},k_{2},\ldots,k_{q}\in\left[  n-1\right]  $, and let
$\sigma=s_{k_{1}}s_{k_{2}}\cdots s_{k_{q}}$. Then, $q\geq\ell\left(
\sigma\right)  $ and $q\equiv\ell\left(  \sigma\right)  \operatorname{mod}2$.
\end{corollary}

\begin{example}
Let $n=4$. Consider the two permutations $\sigma=3214$ and $\tau=3142$ (both
written in one-line notation). Then, $\ell\left(  \sigma\right)  =3$ and
$\ell\left(  \tau\right)  =3$. Now, the permutation $\sigma\tau=1342$ has
length $\ell\left(  \sigma\tau\right)  =2$.

Corollary \ref{cor.perm.red.sigtau} \textbf{(a)} says $\ell\left(  \sigma
\tau\right)  \equiv\ell\left(  \sigma\right)  +\ell\left(  \tau\right)
\operatorname{mod}2$. In other words, $2\equiv3+3\operatorname{mod}2$.

Corollary \ref{cor.perm.red.sigtau} \textbf{(b)} says $\ell\left(  \sigma
\tau\right)  \leq\ell\left(  \sigma\right)  +\ell\left(  \tau\right)  $. In
other words, $2\leq3+3$.
\end{example}

\begin{proof}
[Proof of Corollary \ref{cor.perm.red.sigtau} (sketched).]\textbf{(a)} (See
\cite[Exercise 5.2 \textbf{(b)}]{detnotes} for details.)

For any $\sigma\in S_{n}$ and any $k\in\left[  n-1\right]  $, we have%
\begin{equation}
\ell\left(  \sigma s_{k}\right)  \equiv\ell\left(  \sigma\right)
+1\operatorname{mod}2. \label{pf.cor.perm.red.sigtau.a.1}%
\end{equation}
(This follows from Lemma \ref{lem.perm.len.ssl} \textbf{(a)}, since both
numbers $\ell\left(  \sigma\right)  +1$ and $\ell\left(  \sigma\right)  -1$
are congruent to $\ell\left(  \sigma\right)  +1$ modulo $2$.)

Now, let $\sigma\in S_{n}$ and $\tau\in S_{n}$. Theorem
\ref{thm.perm.len.redword1} \textbf{(a)} yields that we can write $\tau$ as a
composition of $\ell\left(  \tau\right)  $ simples. In other words, we can
write $\tau$ as $\tau=s_{k_{1}}s_{k_{2}}\cdots s_{k_{q}}$ for some
$k_{1},k_{2},\ldots,k_{q}\in\left[  n-1\right]  $, where $q=\ell\left(
\tau\right)  $. Consider these $k_{1},k_{2},\ldots,k_{q}$. Now,%
\begin{align*}
\ell\left(  \sigma\tau\right)   &  =\ell\left(  \sigma s_{k_{1}}s_{k_{2}%
}\cdots s_{k_{q}}\right)  \ \ \ \ \ \ \ \ \ \ \left(  \text{since }%
\tau=s_{k_{1}}s_{k_{2}}\cdots s_{k_{q}}\right) \\
&  \equiv\ell\left(  \sigma s_{k_{1}}s_{k_{2}}\cdots s_{k_{q-1}}\right)
+1\ \ \ \ \ \ \ \ \ \ \left(  \text{by (\ref{pf.cor.perm.red.sigtau.a.1}%
)}\right) \\
&  \equiv\ell\left(  \sigma s_{k_{1}}s_{k_{2}}\cdots s_{k_{q-2}}\right)  +2\\
&  \ \ \ \ \ \ \ \ \ \ \left(  \text{since (\ref{pf.cor.perm.red.sigtau.a.1})
yields }\ell\left(  \sigma s_{k_{1}}s_{k_{2}}\cdots s_{k_{q-1}}\right)
\equiv\ell\left(  \sigma s_{k_{1}}s_{k_{2}}\cdots s_{k_{q-2}}\right)
+1\operatorname{mod}2\right) \\
&  \equiv\ell\left(  \sigma s_{k_{1}}s_{k_{2}}\cdots s_{k_{q-3}}\right)  +3\\
&  \ \ \ \ \ \ \ \ \ \ \left(  \text{since (\ref{pf.cor.perm.red.sigtau.a.1})
yields }\ell\left(  \sigma s_{k_{1}}s_{k_{2}}\cdots s_{k_{q-2}}\right)
\equiv\ell\left(  \sigma s_{k_{1}}s_{k_{2}}\cdots s_{k_{q-3}}\right)
+1\operatorname{mod}2\right) \\
&  \equiv\cdots\\
&  \equiv\ell\left(  \sigma\right)  +\underbrace{q}_{=\ell\left(  \tau\right)
}=\ell\left(  \sigma\right)  +\ell\left(  \tau\right)  \operatorname{mod}2.
\end{align*}
This proves Corollary \ref{cor.perm.red.sigtau} \textbf{(a)}. \medskip

\textbf{(b)} (See \cite[Exercise 5.2 \textbf{(c)}]{detnotes} for details.)

This is analogous to the proof of Corollary \ref{cor.perm.red.sigtau}
\textbf{(a)} (but using inequalities instead of congruences, and using
(\ref{pf.thm.perm.len.redword1.b.2}) instead of
(\ref{pf.cor.perm.red.sigtau.a.1})). \medskip

\textbf{(c)} Let $k_{1},k_{2},\ldots,k_{q}\in\left[  n-1\right]  $, and let
$\sigma=s_{k_{1}}s_{k_{2}}\cdots s_{k_{q}}$. We must prove that $q\geq
\ell\left(  \sigma\right)  $ and $q\equiv\ell\left(  \sigma\right)
\operatorname{mod}2$.

From $\sigma=s_{k_{1}}s_{k_{2}}\cdots s_{k_{q}}$, we obtain%
\begin{align*}
\ell\left(  \sigma\right)   &  =\ell\left(  s_{k_{1}}s_{k_{2}}\cdots s_{k_{q}%
}\right) \\
&  \equiv\ell\left(  s_{k_{1}}s_{k_{2}}\cdots s_{k_{q-1}}\right)
+1\ \ \ \ \ \ \ \ \ \ \left(  \text{by (\ref{pf.cor.perm.red.sigtau.a.1}%
)}\right) \\
&  \equiv\ell\left(  s_{k_{1}}s_{k_{2}}\cdots s_{k_{q-2}}\right)  +2\\
&  \ \ \ \ \ \ \ \ \ \ \left(  \text{since (\ref{pf.cor.perm.red.sigtau.a.1})
yields }\ell\left(  s_{k_{1}}s_{k_{2}}\cdots s_{k_{q-1}}\right)  \equiv
\ell\left(  s_{k_{1}}s_{k_{2}}\cdots s_{k_{q-2}}\right)  +1\operatorname{mod}%
2\right) \\
&  \equiv\ell\left(  s_{k_{1}}s_{k_{2}}\cdots s_{k_{q-3}}\right)  +3\\
&  \ \ \ \ \ \ \ \ \ \ \left(  \text{since (\ref{pf.cor.perm.red.sigtau.a.1})
yields }\ell\left(  s_{k_{1}}s_{k_{2}}\cdots s_{k_{q-2}}\right)  \equiv
\ell\left(  s_{k_{1}}s_{k_{2}}\cdots s_{k_{q-3}}\right)  +1\operatorname{mod}%
2\right) \\
&  \equiv\cdots\\
&  \equiv\underbrace{\ell\left(  \operatorname*{id}\right)  }_{=0}%
+\,q\ \ \ \ \ \ \ \ \ \ \left(  \text{since the product of }0\text{ simples is
}\operatorname*{id}\right) \\
&  =q\operatorname{mod}2.
\end{align*}
In other words, $q\equiv\ell\left(  \sigma\right)  \operatorname{mod}2$. A
similar argument (but using inequalities instead of congruences, and using
(\ref{pf.thm.perm.len.redword1.b.2}) instead of
(\ref{pf.cor.perm.red.sigtau.a.1})) shows that $q\geq\ell\left(
\sigma\right)  $. Thus, Corollary \ref{cor.perm.red.sigtau} \textbf{(c)} is proved.
\end{proof}

\begin{corollary}
\label{cor.perm.generated}Let $n\in\mathbb{N}$. Then, the group $S_{n}$ is
generated by the simples $s_{1},s_{2},\ldots,s_{n-1}$.
\end{corollary}

\begin{proof}
This follows directly from Theorem \ref{thm.perm.len.redword1} \textbf{(a)}.
\end{proof}

Theorem \ref{thm.perm.len.redword1} \textbf{(a)} shows that every permutation
$\sigma\in S_{n}$ can be represented as a product of $\ell\left(
\sigma\right)  $ simples (and in most cases, this can be done in many
different ways). It turns out that there is a rather explicit way to find such
a representation:

\begin{remark}
\label{rmk.perm.redword-lehmer}Let $n\in\mathbb{N}$ and $\sigma\in S_{n}$. Let
us represent the Lehmer code of $\sigma$ visually as follows:

We draw an (empty) $n\times n$-matrix.

For each $i\in\left[  n\right]  $, we put a cross $\times$ into the cell
$\left(  i,\sigma\left(  i\right)  \right)  $ of the matrix.

In the following, I will use the case $n=6$ and $\sigma=513462$ (in one-line
notation) as a running example. In this case, the matrix looks as follows:%
\[%
\begin{tikzpicture}[scale=0.8]
\draw[step=1cm] (0, 0) grid (6, 6);
\draw(0.2, 4.8) -- (0.8, 4.2)
(0.8, 4.8) -- (0.2, 4.2);
\draw(1.2, 0.8) -- (1.8, 0.2)
(1.8, 0.8) -- (1.2, 0.2);
\draw(2.2, 3.8) -- (2.8, 3.2)
(2.8, 3.8) -- (2.2, 3.2);
\draw(3.2, 2.8) -- (3.8, 2.2)
(3.8, 2.8) -- (3.2, 2.2);
\draw(4.2, 5.8) -- (4.8, 5.2)
(4.8, 5.8) -- (4.2, 5.2);
\draw(5.2, 1.8) -- (5.8, 1.2)
(5.8, 1.8) -- (5.2, 1.2);
\end{tikzpicture}%
\]

Now, starting from each $\times$, we draw a vertical ray downwards and a
horizontal ray eastwards. I will call these two rays the \emph{Lehmer lasers}.
Here is how the rays look like in our running example:%
\[%
\begin{tikzpicture}[scale=0.8]
\draw[step=1cm] (0, 0) grid (6, 6);
\draw(0.2, 4.8) -- (0.8, 4.2)
(0.8, 4.8) -- (0.2, 4.2);
\draw(1.2, 0.8) -- (1.8, 0.2)
(1.8, 0.8) -- (1.2, 0.2);
\draw(2.2, 3.8) -- (2.8, 3.2)
(2.8, 3.8) -- (2.2, 3.2);
\draw(3.2, 2.8) -- (3.8, 2.2)
(3.8, 2.8) -- (3.2, 2.2);
\draw(4.2, 5.8) -- (4.8, 5.2)
(4.8, 5.8) -- (4.2, 5.2);
\draw(5.2, 1.8) -- (5.8, 1.2)
(5.8, 1.8) -- (5.2, 1.2);
\draw[very thick, red, ->] (0.5, 4.5) -- (6.5, 4.5);
\draw[very thick, red, ->] (0.5, 4.5) -- (0.5, -0.5);
\draw[very thick, red, ->] (1.5, 0.5) -- (6.5, 0.5);
\draw[very thick, red, ->] (1.5, 0.5) -- (1.5, -0.5);
\draw[very thick, red, ->] (2.5, 3.5) -- (6.5, 3.5);
\draw[very thick, red, ->] (2.5, 3.5) -- (2.5, -0.5);
\draw[very thick, red, ->] (3.5, 2.5) -- (6.5, 2.5);
\draw[very thick, red, ->] (3.5, 2.5) -- (3.5, -0.5);
\draw[very thick, red, ->] (4.5, 5.5) -- (6.5, 5.5);
\draw[very thick, red, ->] (4.5, 5.5) -- (4.5, -0.5);
\draw[very thick, red, ->] (5.5, 1.5) -- (6.5, 1.5);
\draw[very thick, red, ->] (5.5, 1.5) -- (5.5, -0.5);
\end{tikzpicture}
\]

Now, we draw a little circle $\circ$ into each cell that is not hit by any
laser. Here is where the circles end up in our example:%
\[%
\begin{tikzpicture}[scale=0.8]
\draw[step=1cm] (0, 0) grid (6, 6);
\draw(0.2, 4.8) -- (0.8, 4.2)
(0.8, 4.8) -- (0.2, 4.2);
\draw(1.2, 0.8) -- (1.8, 0.2)
(1.8, 0.8) -- (1.2, 0.2);
\draw(2.2, 3.8) -- (2.8, 3.2)
(2.8, 3.8) -- (2.2, 3.2);
\draw(3.2, 2.8) -- (3.8, 2.2)
(3.8, 2.8) -- (3.2, 2.2);
\draw(4.2, 5.8) -- (4.8, 5.2)
(4.8, 5.8) -- (4.2, 5.2);
\draw(5.2, 1.8) -- (5.8, 1.2)
(5.8, 1.8) -- (5.2, 1.2);
\draw[very thick, red, ->] (0.5, 4.5) -- (6.5, 4.5);
\draw[very thick, red, ->] (0.5, 4.5) -- (0.5, -0.5);
\draw[very thick, red, ->] (1.5, 0.5) -- (6.5, 0.5);
\draw[very thick, red, ->] (1.5, 0.5) -- (1.5, -0.5);
\draw[very thick, red, ->] (2.5, 3.5) -- (6.5, 3.5);
\draw[very thick, red, ->] (2.5, 3.5) -- (2.5, -0.5);
\draw[very thick, red, ->] (3.5, 2.5) -- (6.5, 2.5);
\draw[very thick, red, ->] (3.5, 2.5) -- (3.5, -0.5);
\draw[very thick, red, ->] (4.5, 5.5) -- (6.5, 5.5);
\draw[very thick, red, ->] (4.5, 5.5) -- (4.5, -0.5);
\draw[very thick, red, ->] (5.5, 1.5) -- (6.5, 1.5);
\draw[very thick, red, ->] (5.5, 1.5) -- (5.5, -0.5);
\draw[very thick, blue!70!black] (0.5, 5.5) circle (0.2);
\draw[very thick, blue!70!black] (1.5, 5.5) circle (0.2);
\draw[very thick, blue!70!black] (2.5, 5.5) circle (0.2);
\draw[very thick, blue!70!black] (3.5, 5.5) circle (0.2);
\draw[very thick, blue!70!black] (1.5, 3.5) circle (0.2);
\draw[very thick, blue!70!black] (1.5, 2.5) circle (0.2);
\draw[very thick, blue!70!black] (1.5, 1.5) circle (0.2);
\end{tikzpicture}%
\]

This picture is called the \emph{Rothe diagram} of $\sigma$.

Explicitly, a cell $\left(  i,j\right)  $ has a $\circ$ in it if and only if%
\[
\sigma\left(  i\right)  >j\text{ and }\sigma^{-1}\left(  j\right)  >i
\]
(indeed, the vertical laser in column $j$ hits cell $\left(  i,j\right)  $ if
and only if $\sigma^{-1}\left(  j\right)  \leq i$, whereas the horizontal
laser in row $i$ hits cell $\left(  i,j\right)  $ if and only if
$\sigma\left(  i\right)  \leq j$).

If we substitute $\sigma\left(  j\right)  $ for $j$ in this statement, then we
obtain the following: A cell $\left(  i,\sigma\left(  j\right)  \right)  $ has
a $\circ$ in it if and only if%
\[
\sigma\left(  i\right)  >\sigma\left(  j\right)  \text{ and }j>i.
\]
In other words, a cell $\left(  i,\sigma\left(  j\right)  \right)  $ has a
$\circ$ in it if and only if $\left(  i,j\right)  $ is an inversion of
$\sigma$.

Thus,%
\[
\ell\left(  \sigma\right)  =\left(  \text{\# of inversions of }\sigma\right)
=\left(  \text{\# of }\circ\text{'s}\right)  ,
\]
and%
\[
\ell_{i}\left(  \sigma\right)  =\left(  \text{\# of }\circ\text{'s in row
}i\right)  \ \ \ \ \ \ \ \ \ \ \text{for each }i\in\left[  n\right]  .
\]

Finally, let us label the $\circ$'s as follows: For each $i\in\left[
n\right]  $, we label the $\circ$'s in row $i$ from right to left by the
simple transpositions $s_{i},s_{i+1},s_{i+2},\ldots,s_{i^{\prime}-1}$ where
$i^{\prime}=i+\ell_{i}\left(  \sigma\right)  $. (This works, since there are
precisely $\ell_{i}\left(  \sigma\right)  =i^{\prime}-i$ many $\circ$'s in row
$i$.) Here is how this labeling looks in our running example:%
\[%
\begin{tikzpicture}[scale=0.8]
\draw[step=1cm] (0, 0) grid (6, 6);
\draw(0.2, 4.8) -- (0.8, 4.2)
(0.8, 4.8) -- (0.2, 4.2);
\draw(1.2, 0.8) -- (1.8, 0.2)
(1.8, 0.8) -- (1.2, 0.2);
\draw(2.2, 3.8) -- (2.8, 3.2)
(2.8, 3.8) -- (2.2, 3.2);
\draw(3.2, 2.8) -- (3.8, 2.2)
(3.8, 2.8) -- (3.2, 2.2);
\draw(4.2, 5.8) -- (4.8, 5.2)
(4.8, 5.8) -- (4.2, 5.2);
\draw(5.2, 1.8) -- (5.8, 1.2)
(5.8, 1.8) -- (5.2, 1.2);
\draw[very thick, red, ->] (0.5, 4.5) -- (6.5, 4.5);
\draw[very thick, red, ->] (0.5, 4.5) -- (0.5, -0.5);
\draw[very thick, red, ->] (1.5, 0.5) -- (6.5, 0.5);
\draw[very thick, red, ->] (1.5, 0.5) -- (1.5, -0.5);
\draw[very thick, red, ->] (2.5, 3.5) -- (6.5, 3.5);
\draw[very thick, red, ->] (2.5, 3.5) -- (2.5, -0.5);
\draw[very thick, red, ->] (3.5, 2.5) -- (6.5, 2.5);
\draw[very thick, red, ->] (3.5, 2.5) -- (3.5, -0.5);
\draw[very thick, red, ->] (4.5, 5.5) -- (6.5, 5.5);
\draw[very thick, red, ->] (4.5, 5.5) -- (4.5, -0.5);
\draw[very thick, red, ->] (5.5, 1.5) -- (6.5, 1.5);
\draw[very thick, red, ->] (5.5, 1.5) -- (5.5, -0.5);
\node[blue!70!black] at (0.5, 5.5) {$s_4$};
\node[blue!70!black] at (1.5, 5.5) {$s_3$};
\node[blue!70!black] at (2.5, 5.5) {$s_2$};
\node[blue!70!black] at (3.5, 5.5) {$s_1$};
\node[blue!70!black] at (1.5, 3.5) {$s_3$};
\node[blue!70!black] at (1.5, 2.5) {$s_4$};
\node[blue!70!black] at (1.5, 1.5) {$s_5$};
\end{tikzpicture}%
\]

Finally, read the matrix row by row, starting with the top row, and reading
each row from left to right. The result, in our running example, is%
\[
s_{4}s_{3}s_{2}s_{1}s_{3}s_{4}s_{5}.
\]
Reading this as a product, we obtain a product of $\ell\left(  \sigma\right)
$ simples that equals $\sigma$.
\end{remark}

The claim of Remark \ref{rmk.perm.redword-lehmer} can be restated in a more
direct (if less visual) fashion:

\begin{proposition}
\label{prop.perm.redword-lehmer}Let $n\in\mathbb{N}$. Let $\sigma\in S_{n}$.
For each $i\in\left[  n\right]  $, we set%
\begin{equation}
a_{i}:=\operatorname*{cyc}\nolimits_{i^{\prime},i^{\prime}-1,i^{\prime
}-2,\ldots,i}=s_{i^{\prime}-1}s_{i^{\prime}-2}s_{i^{\prime}-3}\cdots s_{i},
\label{eq.prop.perm.redword-lehmer.1}%
\end{equation}
where $i^{\prime}=i+\ell_{i}\left(  \sigma\right)  $. Then, $\sigma=a_{1}%
a_{2}\cdots a_{n}$. (The second equality sign in
(\ref{eq.prop.perm.redword-lehmer.1}) is not hard to check. Note that
$a_{n}=\operatorname*{id}$.)
\end{proposition}

We refer to \cite[Exercise 5.21 parts \textbf{(b)} and \textbf{(c)}]{detnotes}
for a (detailed, but annoyingly long) proof of Proposition
\ref{prop.perm.redword-lehmer}. (You are probably better off proving it yourself.)

\subsection{\label{sec.perm.sign}Signs of permutations}

The notion of the \emph{sign} (aka \emph{signature}) of a permutation is a
simple consequence of that of its length; moreover, it is rather well-known,
due to its role in the definition of a determinant. Thus we will survey its
properties quickly and without proofs. More can be found in \cite[\S 5.3 and
\S 5.6]{detnotes} and \cite[Appendix B]{Strick13}.

\begin{definition}
\label{def.perm.sign}Let $n\in\mathbb{N}$. The \emph{sign} of a permutation
$\sigma\in S_{n}$ is defined to be the integer $\left(  -1\right)
^{\ell\left(  \sigma\right)  }$.

It is denoted by $\left(  -1\right)  ^{\sigma}$ or $\operatorname*{sgn}\left(
\sigma\right)  $ or $\operatorname*{sign}\left(  \sigma\right)  $ or
$\varepsilon\left(  \sigma\right)  $. It is also known as the \emph{signature}
of $\sigma$.
\end{definition}

\begin{proposition}
\label{prop.perm.sign.props}Let $n\in\mathbb{N}$. \medskip

\textbf{(a)} The sign of the permutation $\operatorname*{id}\in S_{n}$ is
$\left(  -1\right)  ^{\operatorname*{id}}=1$. \medskip

\textbf{(b)} For any two distinct elements $i$ and $j$ of $\left[  n\right]
$, the transposition $t_{i,j}\in S_{n}$ has sign $\left(  -1\right)
^{t_{i,j}}=-1$. \medskip

\textbf{(c)} For any positive integer $k$ and any distinct elements
$i_{1},i_{2},\ldots,i_{k}\in\left[  n\right]  $, the $k$-cycle
$\operatorname*{cyc}\nolimits_{i_{1},i_{2},\ldots,i_{k}}$ has sign $\left(
-1\right)  ^{\operatorname*{cyc}\nolimits_{i_{1},i_{2},\ldots,i_{k}}}=\left(
-1\right)  ^{k-1}$. \medskip

\textbf{(d)} We have $\left(  -1\right)  ^{\sigma\tau}=\left(  -1\right)
^{\sigma}\cdot\left(  -1\right)  ^{\tau}$ for any $\sigma\in S_{n}$ and
$\tau\in S_{n}$. \medskip

\textbf{(e)} We have $\left(  -1\right)  ^{\sigma_{1}\sigma_{2}\cdots
\sigma_{p}}=\left(  -1\right)  ^{\sigma_{1}}\left(  -1\right)  ^{\sigma_{2}%
}\cdots\left(  -1\right)  ^{\sigma_{p}}$ for any $\sigma_{1},\sigma_{2}%
,\ldots,\sigma_{p}\in S_{n}$. \medskip

\textbf{(f)} We have $\left(  -1\right)  ^{\sigma^{-1}}=\left(  -1\right)
^{\sigma}$ for any $\sigma\in S_{n}$. (The left hand side here has to be
understood as $\left(  -1\right)  ^{\left(  \sigma^{-1}\right)  }$.) \medskip

\textbf{(g)} We have%
\[
\left(  -1\right)  ^{\sigma}=\prod_{1\leq i<j\leq n}\dfrac{\sigma\left(
i\right)  -\sigma\left(  j\right)  }{i-j}\ \ \ \ \ \ \ \ \ \ \text{for each
}\sigma\in S_{n}.
\]
(The product sign \textquotedblleft$\prod_{1\leq i<j\leq n}$\textquotedblright%
\ means a product over all pairs $\left(  i,j\right)  $ of integers satisfying
$1\leq i<j\leq n$. There are $\dbinom{n}{2}$ such pairs.) \medskip

\textbf{(h)} If $x_{1},x_{2},\ldots,x_{n}$ are any elements of some
commutative ring, and if $\sigma\in S_{n}$, then%
\[
\prod_{1\leq i<j\leq n}\left(  x_{\sigma\left(  i\right)  }-x_{\sigma\left(
j\right)  }\right)  =\left(  -1\right)  ^{\sigma}\prod_{1\leq i<j\leq
n}\left(  x_{i}-x_{j}\right)  .
\]

\end{proposition}

\begin{proof}
[Proof of Proposition \ref{prop.perm.sign.props} (sketched).]Most of this
follows easily from what we have proved above, but here are references to
complete proofs: \medskip

\textbf{(a)} This is \cite[Proposition 5.15 \textbf{(a)}]{detnotes}, and
follows easily from $\ell\left(  \operatorname*{id}\right)  =0$. \medskip

\textbf{(d)} This is \cite[Proposition 5.15 \textbf{(c)}]{detnotes}, and
follows easily from Corollary \ref{cor.perm.red.sigtau} \textbf{(a)}. A
different proof appears in \cite[Proposition B.13]{Strick13}. \medskip

\textbf{(b)} This is \cite[Exercise 5.10 \textbf{(b)}]{detnotes}, and follows
easily from Exercise \ref{exe.perm.len.cycles} \textbf{(a)}. \medskip

\textbf{(c)} This is \cite[Exercise 5.17 \textbf{(d)}]{detnotes}, and follows
easily from Exercise \ref{exe.perm.cycles.c-t} \textbf{(a)} and Exercise
\ref{exe.perm.len.cycles} \textbf{(b)} using Proposition
\ref{prop.perm.sign.props} \textbf{(d)}. \medskip

\textbf{(e)} This is \cite[Proposition 5.28]{detnotes}, and follows by
induction from Proposition \ref{prop.perm.sign.props} \textbf{(d)}. \medskip

\textbf{(f)} This is \cite[Proposition 5.15 \textbf{(d)}]{detnotes}, and
follows easily from Proposition \ref{prop.perm.sign.props} \textbf{(d)} or
from Proposition \ref{prop.perm.len.inv}. \medskip

\textbf{(h)} This is \cite[Exercise 5.13 \textbf{(a)}]{detnotes} (or, rather,
the straightforward generalization of \cite[Exercise 5.13 \textbf{(a)}%
]{detnotes} to arbitrary commutative rings). The proof is fairly easy: Each
factor $x_{\sigma\left(  i\right)  }-x_{\sigma\left(  j\right)  }$ on the left
hand side appears also on the right hand side, albeit with a different sign if
$\left(  i,j\right)  $ is an inversion of $\sigma$. Thus, the products on both
sides agree up to a sign, which is precisely $\left(  -1\right)  ^{\ell\left(
\sigma\right)  }=\left(  -1\right)  ^{\sigma}$. \medskip

\textbf{(g)} This is \cite[Exercise 5.13 \textbf{(c)}]{detnotes}, and is a
particular case of Proposition \ref{prop.perm.sign.props} \textbf{(h)}.
\end{proof}

\begin{corollary}
\label{cor.perm.sign.hom}Let $n\in\mathbb{N}$. The map%
\begin{align*}
S_{n}  &  \rightarrow\left\{  1,-1\right\}  ,\\
\sigma &  \mapsto\left(  -1\right)  ^{\sigma}%
\end{align*}
is a group homomorphism from the symmetric group $S_{n}$ to the order-$2$
group $\left\{  1,-1\right\}  $. (Of course, $\left\{  1,-1\right\}  $ is a
group with respect to multiplication.)
\end{corollary}

This map is known as the \emph{sign homomorphism}.

\begin{proof}
[Proof of Corollary \ref{cor.perm.sign.hom} (sketched).]Proposition
\ref{prop.perm.sign.props} \textbf{(d)} shows that this map respects
multiplication (i.e., sends products to products). However, if a map between
two groups respects multiplication, then it is automatically a group
homomorphism. Thus, Corollary \ref{cor.perm.sign.hom} follows.
\end{proof}

\begin{definition}
\label{def.perm.even-odd}Let $n\in\mathbb{N}$. A permutation $\sigma\in S_{n}$
is said to be

\begin{itemize}
\item \emph{even} if $\left(  -1\right)  ^{\sigma}=1$ (that is, if
$\ell\left(  \sigma\right)  $ is even);

\item \emph{odd} if $\left(  -1\right)  ^{\sigma}=-1$ (that is, if
$\ell\left(  \sigma\right)  $ is odd).
\end{itemize}
\end{definition}

The sign and the \textquotedblleft parity\textquotedblright\ (i.e.,
evenness/oddness) of a permutation have applications throughout mathematics
(in the definition of determinants and the construction of exterior powers) as
well as in the solution of \emph{permutation puzzles} (such as Rubik's cube
and the $15$-game; see \cite[Chapters 7--8 and Theorem 20.2.1]{Mulhol21} for
example). Even permutations are also crucial in group theory, as they form a group:

\begin{corollary}
\label{cor.perm.altgp}Let $n\in\mathbb{N}$. The set of all even permutations
in $S_{n}$ is a normal subgroup of $S_{n}$.
\end{corollary}

This subgroup is known as the $n$\emph{-th alternating group} (commonly called
$A_{n}$). If $n\geq5$, then this group is a simple group (meaning that it has
no normal subgroups besides itself and the trivial group)\footnote{See, e.g.,
\url{https://groupprops.subwiki.org/wiki/Alternating_groups_are_simple} or
\cite[Theorem 10.3.4]{Goodman} for a proof.}, and this fact has been used by
Galois to prove that the general $5$-th degree polynomial equation cannot be
solved using radicals.

\begin{proof}
[Proof of Corollary \ref{cor.perm.altgp} (sketched).]The set of all even
permutations in $S_{n}$ is the kernel of the group homomorphism $S_{n}%
\rightarrow\left\{  1,-1\right\}  $ from Corollary \ref{cor.perm.sign.hom}.
Thus, it is a normal subgroup of $S_{n}$ (since any kernel is a normal subgroup).
\end{proof}

\begin{corollary}
\label{cor.perm.num-even}Let $n\geq2$. Then,%
\[
\left(  \text{\# of even permutations }\sigma\in S_{n}\right)  =\left(
\text{\# of odd permutations }\sigma\in S_{n}\right)  =n!/2.
\]

\end{corollary}

\begin{proof}
[Proof of Corollary \ref{cor.perm.num-even} (sketched).]The symmetric group
$S_{n}$ contains the simple transposition $s_{1}$ (since $n\geq2$). If
$\sigma\in S_{n}$, then
\begin{align*}
\left(  -1\right)  ^{\sigma s_{1}}  &  =\left(  -1\right)  ^{\sigma}%
\cdot\underbrace{\left(  -1\right)  ^{s_{1}}}_{\substack{=-1\\\text{(by
Proposition \ref{prop.perm.sign.props} \textbf{(b)},}\\\text{since }%
s_{1}=t_{1,2}\text{)}}}\ \ \ \ \ \ \ \ \ \ \left(  \text{by Proposition
\ref{prop.perm.sign.props} \textbf{(d)}}\right) \\
&  =-\left(  -1\right)  ^{\sigma}.
\end{align*}
Hence, a permutation $\sigma\in S_{n}$ is even if and only if the permutation
$\sigma s_{1}$ is odd. Hence, the map%
\begin{align*}
\left\{  \text{even permutations }\sigma\in S_{n}\right\}   &  \rightarrow
\left\{  \text{odd permutations }\sigma\in S_{n}\right\}  ,\\
\sigma &  \mapsto\sigma s_{1}%
\end{align*}
is well-defined. This map is furthermore a bijection (since $S_{n}$ is a
group). Thus, the bijection principle yields
\[
\left(  \text{\# of even permutations }\sigma\in S_{n}\right)  =\left(
\text{\# of odd permutations }\sigma\in S_{n}\right)  .
\]
Both sides of this equality must furthermore equal $n!/2$, since they add up
to $\left\vert S_{n}\right\vert =n!$. This proves Corollary
\ref{cor.perm.num-even}. (See \cite[Exercise 5.4]{detnotes} for details.)
\end{proof}

As a consequence of Corollary \ref{cor.perm.num-even}, we see that%
\begin{equation}
\sum_{\sigma\in S_{n}}\left(  -1\right)  ^{\sigma}%
=0\ \ \ \ \ \ \ \ \ \ \text{for each }n\geq2.
\label{eq.cor.perm.num-even.sum-sign}%
\end{equation}
(Indeed, the sum $\sum_{\sigma\in S_{n}}\left(  -1\right)  ^{\sigma}$ can be
rewritten as
\[
\left(  \text{\# of even permutations }\sigma\in S_{n}\right)  -\left(
\text{\# of odd permutations }\sigma\in S_{n}\right)  ,
\]
since the addends corresponding to the even permutations $\sigma\in S_{n}$ are
equal to $1$ whereas the addends corresponding to the odd permutations
$\sigma\in S_{n}$ are equal to $-1$.)

We note that the sign can be defined not only for a permutation $\sigma\in
S_{n}$, but also for any permutation of any finite set $X$ (even if the set
$X$ has no chosen total order on it, as the set $\left[  n\right]  $ has).
Here is one way to do so:

\begin{proposition}
\label{prop.perm.sign.X}Let $X$ be a finite set. We want to define the sign of
any permutation of $X$.

Fix a bijection $\phi:X\rightarrow\left[  n\right]  $ for some $n\in
\mathbb{N}$. (Such a bijection always exists, since $X$ is finite.) For every
permutation $\sigma$ of $X$, set%
\[
\left(  -1\right)  _{\phi}^{\sigma}:=\left(  -1\right)  ^{\phi\circ\sigma
\circ\phi^{-1}}.
\]
Here, the right hand side is well-defined, since $\phi\circ\sigma\circ
\phi^{-1}$ is a permutation of $\left[  n\right]  $. Now: \medskip

\textbf{(a)} This number $\left(  -1\right)  _{\phi}^{\sigma}$ depends only on
the permutation $\sigma$, but not on the bijection $\phi$. (In other words, if
$\phi_{1}$ and $\phi_{2}$ are two bijections from $X$ to $\left[  n\right]  $,
then $\left(  -1\right)  _{\phi_{1}}^{\sigma}=\left(  -1\right)  _{\phi_{2}%
}^{\sigma}$.)

Thus, we shall denote $\left(  -1\right)  _{\phi}^{\sigma}$ by $\left(
-1\right)  ^{\sigma}$ from now on. We refer to this number $\left(  -1\right)
^{\sigma}$ as the \emph{sign} of the permutation $\sigma\in S_{X}$. (When
$X=\left[  n\right]  $, this notation does not clash with Definition
\ref{def.perm.sign}, since we can pick the bijection $\phi=\operatorname*{id}$
and obtain $\left(  -1\right)  _{\phi}^{\sigma}=\left(  -1\right)
^{\operatorname*{id}\circ\sigma\circ\operatorname*{id}\nolimits^{-1}}=\left(
-1\right)  ^{\sigma}$.) \medskip

\textbf{(b)} The identity permutation $\operatorname*{id}:X\rightarrow X$
satisfies $\left(  -1\right)  ^{\operatorname*{id}}=1$. \medskip

\textbf{(c)} We have $\left(  -1\right)  ^{\sigma\tau}=\left(  -1\right)
^{\sigma}\cdot\left(  -1\right)  ^{\tau}$ for any two permutations $\sigma$
and $\tau$ of $X$.
\end{proposition}

\begin{proof}
[Proof of Proposition \ref{prop.perm.sign.X} (sketched).]This all follows
quite easily from Proposition \ref{prop.perm.sign.props} \textbf{(d)}. See
\cite[Exercise 5.12]{detnotes} for a detailed proof.
\end{proof}

\subsection{The cycle decomposition}

Next, we shall discuss the \emph{cycle decomposition} (or \emph{disjoint cycle
decomposition}) of a permutation. Again, this is a fairly well-known and
elementary tool, so we will restrict ourselves to the basic properties and
omit the details.

We begin with an introductory example:

\begin{example}
\label{exa.perm.dcd.1}Let $\sigma\in S_{9}$ be the permutation with one-line
notation $461352987$. Here is its cycle digraph:%
\[%
\begin{tikzpicture}%
[->,shorten >=1pt,auto,node distance=3cm, thick,main node/.style={circle,fill=blue!20,draw}%
]
\node[main node] (1) at (-60 : 1.5) {1};
\node[main node] (3) at (60 : 1.5) {3};
\node[main node] (4) at (180 : 1.5) {4};
\path[-{Stealth[length=4mm]}]
(1) edge [bend left] (4)
(4) edge [bend left] (3)
(3) edge [bend left] (1);
\end{tikzpicture}%
\qquad\qquad%
\begin{tikzpicture}%
[->,shorten >=1pt,auto,node distance=3cm, thick,main node/.style={circle,fill=blue!20,draw}%
]
\node[main node] (6) at (-90 : 1.4) {6};
\node[main node] (2) at (90 : 1.4) {2};
\path[-{Stealth[length=4mm]}]
(2) edge [bend left] (6)
(6) edge [bend left] (2);
\end{tikzpicture}%
\qquad%
\raisebox{8ex}{
\begin{tikzpicture}%
[->,shorten >=1pt,auto,node distance=3cm, thick,main node/.style={circle,fill=blue!20,draw}%
]
\node[main node] (5) at (-90 : 0) {5};
\path[-{Stealth[length=4mm]}]
(5) edge [loop right] (5);
\end{tikzpicture}
}%
\qquad%
\begin{tikzpicture}%
[->,shorten >=1pt,auto,node distance=3cm, thick,main node/.style={circle,fill=blue!20,draw}%
]
\node[main node] (7) at (-90 : 1.4) {7};
\node[main node] (9) at (90 : 1.4) {9};
\path[-{Stealth[length=4mm]}]
(7) edge [bend left] (9)
(9) edge [bend left] (7);
\end{tikzpicture}%
\qquad%
\raisebox{8ex}{
\begin{tikzpicture}%
[->,shorten >=1pt,auto,node distance=3cm, thick,main node/.style={circle,fill=blue!20,draw}%
]
\node[main node] (8) at (-90 : 0) {8};
\path[-{Stealth[length=4mm]}]
(8) edge [loop right] (8);
\end{tikzpicture}
}%
\]
(where we have strategically arranged the cycles apart horizontally). This
digraph consists of five node-disjoint cycles (i.e., cycles that share no
nodes); thus, we can view the permutation $\sigma$ as acting on the five
subsets $\left\{  1,4,3\right\}  $, $\left\{  2,6\right\}  $, $\left\{
5\right\}  $, $\left\{  7,9\right\}  $ and $\left\{  8\right\}  $ of $\left[
9\right]  $ separately. On the first of these five subsets, $\sigma$ acts as
the $3$-cycle $\operatorname*{cyc}\nolimits_{1,4,3}$ (in the sense that
$\sigma\left(  k\right)  =\operatorname*{cyc}\nolimits_{1,4,3}\left(
k\right)  $ for each $k\in\left\{  1,4,3\right\}  $). On the second, it acts
as the $2$-cycle $\operatorname*{cyc}\nolimits_{2,6}$. On the third, it acts
as the $1$-cycle $\operatorname*{cyc}\nolimits_{5}$ (which, of course, is the
identity map). On the fourth, it acts as the $2$-cycle $\operatorname*{cyc}%
\nolimits_{7,9}$. On the fifth, it acts as the $1$-cycle $\operatorname*{cyc}%
\nolimits_{8}$ (which, again, is just the identity map). Combining these
observations, we conclude that%
\[
\sigma\left(  k\right)  =\left(  \operatorname*{cyc}\nolimits_{1,4,3}%
\circ\operatorname*{cyc}\nolimits_{2,6}\circ\operatorname*{cyc}\nolimits_{5}%
\circ\operatorname*{cyc}\nolimits_{7,9}\circ\operatorname*{cyc}\nolimits_{8}%
\right)  \left(  k\right)  \ \ \ \ \ \ \ \ \ \ \text{for each }k\in\left[
9\right]
\]
(because, when we apply the composed permutation $\operatorname*{cyc}%
\nolimits_{1,4,3}\circ\operatorname*{cyc}\nolimits_{2,6}\circ
\operatorname*{cyc}\nolimits_{5}\circ\operatorname*{cyc}\nolimits_{7,9}%
\circ\operatorname*{cyc}\nolimits_{8}$ to an element of $\left[  9\right]  $,
then four of the five cycles $\operatorname*{cyc}\nolimits_{1,4,3}%
,\operatorname*{cyc}\nolimits_{2,6},\operatorname*{cyc}\nolimits_{5}%
,\operatorname*{cyc}\nolimits_{7,9},\operatorname*{cyc}\nolimits_{8}$ will
leave this element unchanged, whereas the remaining one will move it one step
forward along the appropriate cycle of the above digraph -- which is precisely
what the permutation $\sigma$ does to our element). This entails that%
\begin{equation}
\sigma=\operatorname*{cyc}\nolimits_{1,4,3}\circ\operatorname*{cyc}%
\nolimits_{2,6}\circ\operatorname*{cyc}\nolimits_{5}\circ\operatorname*{cyc}%
\nolimits_{7,9}\circ\operatorname*{cyc}\nolimits_{8}.
\label{eq.exa.perm.dcd.1.1}%
\end{equation}

\end{example}

In Example \ref{exa.perm.dcd.1}, we have represented our permutation
$\sigma\in S_{9}$ as a composition of five cycles with the property that each
element of $\left[  9\right]  $ appears in exactly one of these cycles. This
is not specific to the permutation $\sigma$ chosen in Example
\ref{exa.perm.dcd.1}. Indeed, for any finite set $X$, any permutation
$\sigma\in S_{X}$ can be written as a composition of finitely many cycles
$\operatorname*{cyc}\nolimits_{i_{1},i_{2},\ldots,i_{k}}$ with the property
that each element of $X$ appears in exactly one of these cycles. Moreover,
this representation of $\sigma$ is unique up to

\begin{itemize}
\item swapping the cycles (for example, we could have replaced
$\operatorname*{cyc}\nolimits_{1,4,3}\circ\operatorname*{cyc}\nolimits_{2,6}$
by $\operatorname*{cyc}\nolimits_{2,6}\circ\operatorname*{cyc}%
\nolimits_{1,4,3}$ in (\ref{eq.exa.perm.dcd.1.1})), and

\item rotating each cycle (for example, we could have replaced
$\operatorname*{cyc}\nolimits_{1,4,3}$ by $\operatorname*{cyc}%
\nolimits_{4,3,1}$ or by $\operatorname*{cyc}\nolimits_{3,1,4}$ in
(\ref{eq.exa.perm.dcd.1.1})).
\end{itemize}

Let us state this in a more rigorous fashion:

\begin{theorem}
[disjoint cycle decomposition of permutations]\label{thm.perm.dcd.main}Let $X$
be a finite set. Let $\sigma$ be a permutation of $X$. Then: \medskip

\textbf{(a)} There is a list%
\begin{align*}
&  \Big(\left(  a_{1,1},a_{1,2},\ldots,a_{1,n_{1}}\right)  ,\\
&  \ \ \ \left(  a_{2,1},a_{2,2},\ldots,a_{2,n_{2}}\right)  ,\\
&  \ \ \ \ldots,\\
&  \ \ \ \left(  a_{k,1},a_{k,2},\ldots,a_{k,n_{k}}\right)  \Big)
\end{align*}
of nonempty lists of elements of $X$ such that:

\begin{itemize}
\item each element of $X$ appears exactly once in the composite list%
\begin{align*}
&  (a_{1,1},a_{1,2},\ldots,a_{1,n_{1}},\\
&  \ \ \ a_{2,1},a_{2,2},\ldots,a_{2,n_{2}},\\
&  \ \ \ \ldots,\\
&  \ \ \ a_{k,1},a_{k,2},\ldots,a_{k,n_{k}}),
\end{align*}
and

\item we have%
\[
\sigma=\operatorname*{cyc}\nolimits_{a_{1,1},a_{1,2},\ldots,a_{1,n_{1}}}%
\circ\operatorname*{cyc}\nolimits_{a_{2,1},a_{2,2},\ldots,a_{2,n_{2}}}%
\circ\cdots\circ\operatorname*{cyc}\nolimits_{a_{k,1},a_{k,2},\ldots
,a_{k,n_{k}}}.
\]

\end{itemize}

Such a list is called a \emph{disjoint cycle decomposition} (or short
\emph{DCD}) of $\sigma$. Its entries (which themselves are lists of elements
of $X$) are called the \emph{cycles} of $\sigma$. \medskip

\textbf{(b)} Any two DCDs of $\sigma$ can be obtained from each other by
(repeatedly) swapping the cycles with each other, and rotating each cycle
(i.e., replacing $\left(  a_{i,1},a_{i,2},\ldots,a_{i,n_{i}}\right)  $ by
$\left(  a_{i,2},a_{i,3},\ldots,a_{i,n_{i}},a_{i,1}\right)  $). \medskip

\textbf{(c)} Now assume that $X$ is a set of integers (or, more generally, any
totally ordered finite set). Then, there is a unique DCD
\begin{align*}
&  \Big(\left(  a_{1,1},a_{1,2},\ldots,a_{1,n_{1}}\right)  ,\\
&  \ \ \ \left(  a_{2,1},a_{2,2},\ldots,a_{2,n_{2}}\right)  ,\\
&  \ \ \ \ldots,\\
&  \ \ \ \left(  a_{k,1},a_{k,2},\ldots,a_{k,n_{k}}\right)  \Big)
\end{align*}
of $\sigma$ that satisfies the additional requirements that

\begin{itemize}
\item we have $a_{i,1}\leq a_{i,p}$ for each $i\in\left[  k\right]  $ and each
$p\in\left[  n_{i}\right]  $ (that is, each cycle in this DCD is written with
its smallest entry first), and

\item we have $a_{1,1}>a_{2,1}>\cdots>a_{k,1}$ (that is, the cycles appear in
this DCD in the order of decreasing first entries).
\end{itemize}
\end{theorem}

\begin{example}
Let $\sigma\in S_{9}$ be the permutation from Example \ref{exa.perm.dcd.1}.
Then, the representation (\ref{eq.exa.perm.dcd.1.1}) shows that%
\[
\Big(\left(  1,4,3\right)  ,\ \left(  2,6\right)  ,\ \left(  5\right)
,\ \left(  7,9\right)  ,\ \left(  8\right)  \Big)
\]
is a DCD of $\sigma$. By swapping the five cycles of this DCD, and by rotating
each cycle, we can produce various other DCDs of $\sigma$, such as
\[
\Big(\left(  7,9\right)  ,\ \left(  6,2\right)  ,\ \left(  3,1,4\right)
,\ \left(  8\right)  ,\ \left(  5\right)  \Big).
\]
The unique DCD of $\sigma$ that satisfies the two additional requirements of
Theorem \ref{thm.perm.dcd.main} \textbf{(c)} is
\[
\Big(\left(  8\right)  ,\ \left(  7,9\right)  ,\ \left(  5\right)  ,\ \left(
2,6\right)  ,\ \left(  1,4,3\right)  \Big).
\]

\end{example}

\begin{proof}
[Proof of Theorem \ref{thm.perm.dcd.main} (sketched).]This is a classical
result with an easy proof; sadly, this easy proof does not present well in
writing. I will try to be as clear as the situation allows. Some familiarity
with digraphs (= directed graphs) is recommended\footnote{See, e.g.,
\cite[\S 5.11]{Guicha20} or \cite[\S 3.5]{Loehr-BC} for brief introductions to
digraphs.}. \medskip

\textbf{(a)} Let $\mathcal{D}$ be the cycle digraph of $\sigma$, as in Example
\ref{exa.perm.dcd.1}. This cycle digraph $\mathcal{D}$ has the following two properties:

\begin{itemize}
\item \textit{Outbound uniqueness:} For each node $i$, there is exactly one
arc outgoing from $i$. (Indeed, this is the arc from $i$ to $\sigma\left(
i\right)  $, as should be clear from the construction of $\mathcal{D}$.)

\item \textit{Inbound uniqueness:} For each node $i$, there is exactly one arc
incoming into $i$. (Indeed, this is the arc from $\sigma^{-1}\left(  i\right)
$ to $i$, since $\sigma$ is a permutation and therefore invertible.)
\end{itemize}

Using these two properties, we will now show that the cycle digraph
$\mathcal{D}$ consists of several node-disjoint cycles (i.e., several cycles
that pairwise share no nodes).

Indeed, let us first observe the following: If two cycles $C$ and $D$ of
$\mathcal{D}$ have a node in common, then they are identical\footnote{Here, we
identify any cycle with its cyclic rotations. For example, if $a\rightarrow
b\rightarrow c\rightarrow a$ is a cycle, then we consider $b\rightarrow
c\rightarrow a\rightarrow b$ to be the same cycle.} (because the outbound
uniqueness property prevents these cycles from ever separating after meeting
at the common node). In other words, any two cycles $C$ and $D$ of
$\mathcal{D}$ are either identical or node-disjoint (i.e., share no nodes with
each other).

Now, let $i$ be any node of $\mathcal{D}$. Then, if we start at $i$ and follow
the outgoing arcs, then we obtain an infinite walk%
\[
\sigma^{0}\left(  i\right)  \rightarrow\sigma^{1}\left(  i\right)
\rightarrow\sigma^{2}\left(  i\right)  \rightarrow\sigma^{3}\left(  i\right)
\rightarrow\cdots
\]
along our digraph $\mathcal{D}$. Since $X$ is finite, the Pigeonhole Principle
guarantees that this walk will eventually revisit a node it has already been
to; i.e., there exist two integers $u,v\in\mathbb{N}$ with $u<v$ and
$\sigma^{u}\left(  i\right)  =\sigma^{v}\left(  i\right)  $. Let us pick two
such integers $u,v$ with the \textbf{smallest} possible $v$. Thus,
\begin{equation}
\text{the }v\text{ nodes }\sigma^{0}\left(  i\right)  ,\sigma^{1}\left(
i\right)  ,\ldots,\sigma^{v-1}\left(  i\right)  \text{ are distinct}
\label{pf.thm.perm.dcd.main.a.dist}%
\end{equation}
(since otherwise, $v$ would not be smallest possible). Now, $\sigma$ is a
permutation and thus has an inverse $\sigma^{-1}$. Applying the map
$\sigma^{-1}$ to both sides of the equality $\sigma^{u}\left(  i\right)
=\sigma^{v}\left(  i\right)  $, we obtain $\sigma^{u-1}\left(  i\right)
=\sigma^{v-1}\left(  i\right)  $. However, if we had $u\geq1$, then
(\ref{pf.thm.perm.dcd.main.a.dist}) would entail $\sigma^{u-1}\left(
i\right)  \neq\sigma^{v-1}\left(  i\right)  $ (because $0\leq\underbrace{u}%
_{<v}-\,1<v-1$), which would contradict $\sigma^{u-1}\left(  i\right)
=\sigma^{v-1}\left(  i\right)  $. Thus, we cannot have $u\geq1$. Hence, $u<1$,
so that $u=0$. Therefore, $\sigma^{u}\left(  i\right)  =\sigma^{0}\left(
i\right)  $, so that $\sigma^{0}\left(  i\right)  =\sigma^{u}\left(  i\right)
=\sigma^{v}\left(  i\right)  $. This shows that our walk $\sigma^{0}\left(
i\right)  \rightarrow\sigma^{1}\left(  i\right)  \rightarrow\sigma^{2}\left(
i\right)  \rightarrow\sigma^{3}\left(  i\right)  \rightarrow\cdots$ is
circular: it comes back to its starting node $\sigma^{0}\left(  i\right)  =i$
after $v$ steps. We have thus found a cycle in our digraph $\mathcal{D}$:%
\[
\sigma^{0}\left(  i\right)  \rightarrow\sigma^{1}\left(  i\right)
\rightarrow\sigma^{2}\left(  i\right)  \rightarrow\cdots\rightarrow\sigma
^{v}\left(  i\right)  =\sigma^{0}\left(  i\right)  .
\]
(This is indeed a cycle, since (\ref{pf.thm.perm.dcd.main.a.dist}) shows that
its first $v$ nodes are distinct.) This shows that the node $i$ lies on a
cycle $C_{i}$ of $\mathcal{D}$ (namely, the cycle that we just found).

Now, forget that we fixed $i$. We thus have shown that each node $i$ of
$\mathcal{D}$ lies on a cycle $C_{i}$. The cycles $C_{i}$ for all nodes $i\in
X$ will be called the \emph{chosen cycles}.

Any arc of our digraph $\mathcal{D}$ must belong to one of these chosen
cycles. Indeed, if $a$ is an arc from a node $i$ to a node $j$, then $a$ must
be the only arc outgoing from $i$ (by the outbound uniqueness property); but
this means that this arc $a$ belongs to the chosen cycle $C_{i}$.

Now, let us look back at what we have shown:

\begin{itemize}
\item Any node $i$ of $\mathcal{D}$ lies on one of our chosen cycles (namely,
on $C_{i}$).

\item Some of the chosen cycles may be identical, but apart from that, the
chosen cycles are pairwise node-disjoint (since any two cycles of
$\mathcal{D}$ are either identical or node-disjoint).

\item Any arc of $\mathcal{D}$ must belong to one of these chosen cycles.
\end{itemize}

Combining these facts, we conclude that $\mathcal{D}$ consists of several
node-disjoint cycles. Let us label these cycles as%
\begin{align*}
&  \left(  a_{1,1}\rightarrow a_{1,2}\rightarrow\cdots\rightarrow a_{1,n_{1}%
}\rightarrow a_{1,1}\right)  ,\\
&  \left(  a_{2,1}\rightarrow a_{2,2}\rightarrow\cdots\rightarrow a_{2,n_{2}%
}\rightarrow a_{2,1}\right)  ,\\
&  \ldots,\\
&  \left(  a_{k,1}\rightarrow a_{k,2}\rightarrow\cdots\rightarrow a_{k,n_{k}%
}\rightarrow a_{k,1}\right)
\end{align*}
(making sure to label each cycle only once). Then, each element of $X$ appears
exactly once in the composite list%
\begin{align*}
&  (a_{1,1},a_{1,2},\ldots,a_{1,n_{1}},\\
&  \ \ \ a_{2,1},a_{2,2},\ldots,a_{2,n_{2}},\\
&  \ \ \ \ldots,\\
&  \ \ \ a_{k,1},a_{k,2},\ldots,a_{k,n_{k}}),
\end{align*}
and we have%
\[
\sigma=\operatorname*{cyc}\nolimits_{a_{1,1},a_{1,2},\ldots,a_{1,n_{1}}}%
\circ\operatorname*{cyc}\nolimits_{a_{2,1},a_{2,2},\ldots,a_{2,n_{2}}}%
\circ\cdots\circ\operatorname*{cyc}\nolimits_{a_{k,1},a_{k,2},\ldots
,a_{k,n_{k}}}%
\]
(since $\sigma$ moves any node $i\in X$ one step forward along its chosen
cycle). This proves Theorem \ref{thm.perm.dcd.main} \textbf{(a)}.

Alternative proofs of Theorem \ref{thm.perm.dcd.main} \textbf{(a)} can be
found (e.g.) in \cite[Theorem 1.5.3]{Goodman} or in \cite[\S I.4, Proposition
1.21]{Knapp16} or in \cite[Chapter I, \S 5.7, Proposition 7]{Bourba74} or in
\cite[\S 1.9, proof of Theorem 1.5.1]{Sagan19} (this is essentially our proof)
or in
\url{https://proofwiki.org/wiki/Existence_and_Uniqueness_of_Cycle_Decomposition}
(see also \cite[Exercise 7 \textbf{(e)} and \textbf{(d)}]{17f-hw7s} for a
rather formalized proof). Note that some of these sources work with a slightly
modified concept of a DCD, in which they throw away the $1$-cycles (i.e., they
replace \textquotedblleft appears exactly once\textquotedblright\ by
\textquotedblleft appears at most once\textquotedblright, and require all
cycle lengths $n_{1},n_{2},\ldots,n_{k}$ to be $>1$). For instance, the DCD
(\ref{eq.exa.perm.dcd.1.1}) becomes%
\[
\sigma=\operatorname*{cyc}\nolimits_{1,4,3}\circ\operatorname*{cyc}%
\nolimits_{2,6}\circ\operatorname*{cyc}\nolimits_{7,9}%
\]
if we use this modified notion of a DCD. \medskip

\textbf{(b)} See \cite[Theorem 1.5.3]{Goodman} or \cite[Chapter I, \S 5.7,
Proposition 7]{Bourba74}. The idea is fairly simple: Let%
\begin{align*}
&  \Big(\left(  a_{1,1},a_{1,2},\ldots,a_{1,n_{1}}\right)  ,\\
&  \ \ \ \left(  a_{2,1},a_{2,2},\ldots,a_{2,n_{2}}\right)  ,\\
&  \ \ \ \ldots,\\
&  \ \ \ \left(  a_{k,1},a_{k,2},\ldots,a_{k,n_{k}}\right)  \Big)
\end{align*}
be a DCD of $\sigma$. Then, for each $i\in X$, the cycle of this DCD that
contains $i$ is uniquely determined by $\sigma$ and $i$ up to cyclic rotation
(indeed, it is a rotated version of the list $\left(  i,\sigma\left(
i\right)  ,\sigma^{2}\left(  i\right)  ,\ldots,\sigma^{r-1}\left(  i\right)
\right)  $, where $r$ is the smallest positive integer satisfying $\sigma
^{r}\left(  i\right)  =i$). Therefore, all cycles of this DCD are uniquely
determined by $\sigma$ up to cyclic rotation and up to the relative order in
which these cycles appear in the DCD. But this is precisely the claim of
Theorem \ref{thm.perm.dcd.main} \textbf{(b)}. \medskip

\textbf{(c)} In order to obtain a DCD of $\sigma$ that satisfies these two
requirements, it suffices to

\begin{itemize}
\item start with an arbitrary DCD of $\sigma$,

\item then rotate each cycle of this DCD so that it begins with its smallest
entry, and

\item then repeatedly swap these cycles so they appear in the order of
decreasing first entries.
\end{itemize}

It is clear that the result of this procedure is uniquely determined (a
consequence of Theorem \ref{thm.perm.dcd.main} \textbf{(b)}). Thus, Theorem
\ref{thm.perm.dcd.main} \textbf{(c)} is proven.
\end{proof}

\begin{definition}
\label{def.perm.cycs.cycs}Let $X$ be a finite set. Let $\sigma$ be a
permutation of $X$. \medskip

\textbf{(a)} The \emph{cycles} of $\sigma$ are defined to be the cycles in the
DCD of $\sigma$ (as defined in Theorem \ref{thm.perm.dcd.main} \textbf{(a)}).
(This includes $1$-cycles, if there are any in the DCD of $\sigma$.)

We shall equate a cycle of $\sigma$ with any of its cyclic rotations; thus,
for example, $\left(  3,1,4\right)  $ and $\left(  1,4,3\right)  $ shall be
regarded as being the same cycle (but $\left(  3,1,4\right)  $ and $\left(
3,4,1\right)  $ shall not). \medskip

\textbf{(b)} The \emph{cycle lengths partition} (aka \emph{cycle type}) of
$\sigma$ shall denote the partition of $\left\vert X\right\vert $ obtained by
writing down the lengths of the cycles of $\sigma$ in weakly decreasing order.
\end{definition}

\begin{example}
Let $\sigma\in S_{9}$ be the permutation from Example \ref{exa.perm.dcd.1}.
Then, the cycles of $\sigma$ are%
\[
\left(  1,4,3\right)  ,\ \left(  2,6\right)  ,\ \left(  5\right)  ,\ \left(
7,9\right)  ,\ \left(  8\right)  .
\]
Their lengths are $3,2,1,2,1$. Hence, the cycle lengths partition of $\sigma$
is $\left(  3,2,2,1,1\right)  $.
\end{example}

The following is obvious:

\begin{proposition}
\label{prop.perm.cycs.same}Let $X$ be a finite set. Let $i,j\in X$. Let
$\sigma$ be a permutation of $X$. Then, we have the following chain of
equivalences:%
\begin{align*}
\  &  \left(  i\text{ and }j\text{ belong to the same cycle of }\sigma\right)
\\
\Longleftrightarrow\  &  \left(  i=\sigma^{p}\left(  j\right)  \text{ for some
}p\in\mathbb{N}\right) \\
\Longleftrightarrow\  &  \left(  j=\sigma^{p}\left(  i\right)  \text{ for some
}p\in\mathbb{N}\right)  .
\end{align*}

\end{proposition}

The number of cycles of a permutation determines its sign. Let us state this
for permutations of $\left[  n\right]  $ in particular (the reader can easily
extend this to the general case using Proposition \ref{prop.perm.sign.X}):

\begin{proposition}
\label{prop.perm.cycs.sign}Let $n\in\mathbb{N}$. Let $\sigma\in S_{n}$. Let
$k\in\mathbb{N}$ be such that $\sigma$ has exactly $k$ cycles (including the
$1$-cycles). Then, $\left(  -1\right)  ^{\sigma}=\left(  -1\right)  ^{n-k}$.
\end{proposition}

\begin{proof}
[Proof of Proposition \ref{prop.perm.cycs.sign} (sketched).]Let
\begin{align*}
&  \left(  a_{1,1},a_{1,2},\ldots,a_{1,n_{1}}\right)  ,\\
&  \left(  a_{2,1},a_{2,2},\ldots,a_{2,n_{2}}\right)  ,\\
&  \ldots,\\
&  \left(  a_{k,1},a_{k,2},\ldots,a_{k,n_{k}}\right)
\end{align*}
be the $k$ cycles of $\sigma$. Thus,
\begin{align}
&  \Big(\left(  a_{1,1},a_{1,2},\ldots,a_{1,n_{1}}\right)  ,\nonumber\\
&  \ \ \ \left(  a_{2,1},a_{2,2},\ldots,a_{2,n_{2}}\right)  ,\nonumber\\
&  \ \ \ \ldots,\nonumber\\
&  \ \ \ \left(  a_{k,1},a_{k,2},\ldots,a_{k,n_{k}}\right)
\Big) \label{pf.prop.perm.cycs.sign.3}%
\end{align}
is the DCD of $\sigma$. Therefore,
\begin{align*}
\sigma &  =\operatorname*{cyc}\nolimits_{a_{1,1},a_{1,2},\ldots,a_{1,n_{1}}%
}\circ\operatorname*{cyc}\nolimits_{a_{2,1},a_{2,2},\ldots,a_{2,n_{2}}}%
\circ\cdots\circ\operatorname*{cyc}\nolimits_{a_{k,1},a_{k,2},\ldots
,a_{k,n_{k}}}\\
&  =\left(  \text{an }n_{1}\text{-cycle}\right)  \circ\left(  \text{an }%
n_{2}\text{-cycle}\right)  \circ\cdots\circ\left(  \text{an }n_{k}%
\text{-cycle}\right)  ,
\end{align*}
so that%
\begin{align}
\left(  -1\right)  ^{\sigma}  &  =\left(  -1\right)  ^{\left(  \text{an }%
n_{1}\text{-cycle}\right)  \circ\left(  \text{an }n_{2}\text{-cycle}\right)
\circ\cdots\circ\left(  \text{an }n_{k}\text{-cycle}\right)  }\nonumber\\
&  =\left(  -1\right)  ^{\left(  \text{an }n_{1}\text{-cycle}\right)  }%
\cdot\left(  -1\right)  ^{\left(  \text{an }n_{2}\text{-cycle}\right)  }%
\cdot\cdots\cdot\left(  -1\right)  ^{\left(  \text{an }n_{k}\text{-cycle}%
\right)  }\nonumber\\
&  \ \ \ \ \ \ \ \ \ \ \ \ \ \ \ \ \ \ \ \ \left(  \text{by Proposition
\ref{prop.perm.sign.props} \textbf{(e)}}\right) \nonumber\\
&  =\left(  -1\right)  ^{n_{1}-1}\cdot\left(  -1\right)  ^{n_{2}-1}\cdot
\cdots\cdot\left(  -1\right)  ^{n_{k}-1}\nonumber\\
&  \ \ \ \ \ \ \ \ \ \ \ \ \ \ \ \ \ \ \ \ \left(
\begin{array}
[c]{c}%
\text{since Proposition \ref{prop.perm.sign.props} \textbf{(c)} yields
}\left(  -1\right)  ^{\left(  \text{a }p\text{-cycle}\right)  }=\left(
-1\right)  ^{p-1}\\
\text{for any }p>0
\end{array}
\right) \nonumber\\
&  =\left(  -1\right)  ^{\left(  n_{1}-1\right)  +\left(  n_{2}-1\right)
+\cdots+\left(  n_{k}-1\right)  }. \label{pf.prop.perm.cycs.sign.5}%
\end{align}
However, recall that (\ref{pf.prop.perm.cycs.sign.3}) is a DCD of $\sigma$.
Thus, each element of $\left[  n\right]  $ appears exactly once in the
composite list%
\begin{align*}
&  (a_{1,1},a_{1,2},\ldots,a_{1,n_{1}},\\
&  \ \ \ a_{2,1},a_{2,2},\ldots,a_{2,n_{2}},\\
&  \ \ \ \ldots,\\
&  \ \ \ a_{k,1},a_{k,2},\ldots,a_{k,n_{k}}).
\end{align*}
Therefore, the length $n_{1}+n_{2}+\cdots+n_{k}$ of this composite list equals
the size $\left\vert \left[  n\right]  \right\vert =n$ of the set $\left[
n\right]  $. In other words, $n_{1}+n_{2}+\cdots+n_{k}=n$. Hence,%
\[
\left(  n_{1}-1\right)  +\left(  n_{2}-1\right)  +\cdots+\left(
n_{k}-1\right)  =\underbrace{\left(  n_{1}+n_{2}+\cdots+n_{k}\right)  }%
_{=n}-\,k=n-k.
\]
Thus, (\ref{pf.prop.perm.cycs.sign.5}) rewrites as $\left(  -1\right)
^{\sigma}=\left(  -1\right)  ^{n-k}$. This proves Proposition
\ref{prop.perm.cycs.sign}.
\end{proof}

\subsection{References}

We end our study of permutations here, although we will see them again several
times. Much more about permutations can be found in \cite{Bona22},
\cite{Kitaev11} (focussing on permutation patterns), \cite{Sagan01} (focussing
on the representation theory of the symmetric group) and various other texts.

It is worth mentioning that the symmetric groups $S_{n}$ are a particular case
of \emph{Coxeter groups} -- a class of groups highly significant to algebra,
combinatorics and geometry. One of the most combinatorial introductions to
this subject (which sheds new light on the combinatorics of symmetric groups)
is the highly readable text \cite{BjoBre05}. Other texts include
\cite{Cohen08} and (for the particularly resolute) \cite[Chapter IV]{Bourba02}.

\newpage

\section{\label{chap.sign}Alternating sums, signed counting and determinants}

This chapter is not concerned with any specific combinatorial objects like
partitions or permutations, but rather with a set of simple ideas that appear
often (and not just in combinatorics). The main objects of study here are
\emph{alternating sums} -- i.e., sums with a $\left(  -1\right)
^{\text{something}}$ factor in them. A poster child is the determinant of a
matrix. Such sums are often simplified by the cancellation that occurs in
them, with positive addends cancelling negative addends. Frequently,
understanding this cancellation is key to computing the sums. As a rule of
thumb, alternating sums are more likely to have simple closed-form answers
than non-alternating sums. For example, each $n\in\mathbb{N}$ satisfies%
\[
\sum_{k=0}^{n}\left(  -1\right)  ^{k}\dbinom{n}{k}^{3}=%
\begin{cases}
\left(  -1\right)  ^{n/2}\dfrac{\left(  3n/2\right)  !}{\left(  n/2\right)
!^{3}}, & \text{if }n\text{ is even};\\
0, & \text{if }n\text{ is odd}%
\end{cases}
\]
(see Exercise \ref{exe.fps.dixon} \textbf{(g)}), but there is no closed form
for%
\[
\sum_{k=0}^{n}\dbinom{n}{k}^{3}.
\]

The use of cancellations to simplify alternating sums is old, but systematic
surveys of applications of this technique have not appeared until recently
(\cite[Chapter 2]{Stanley-EC1}, \cite[Chapter 5]{Aigner07}, \cite[Chapter
6]{BenQui03}, \cite{BenQui08}, \cite[Chapter 2]{Sagan19}, \cite{Grinbe20}).

\subsection{\label{sec.sign.intro}Cancellations in alternating sums}

We begin with a simple binomial identity:

\begin{proposition}
[Negative hockey-stick identity]\label{prop.binom.nhs}Let $n\in\mathbb{C}$ and
$m\in\mathbb{N}$. Then,%
\begin{equation}
\sum_{k=0}^{m}\left(  -1\right)  ^{k}\dbinom{n}{k}=\left(  -1\right)
^{m}\dbinom{n-1}{m}. \label{eq.prop.binom.nhs.eq}%
\end{equation}

\end{proposition}

There are easy proofs of this proposition by induction on $m$ or by the
telescope principle (see, e.g., \cite[Exercise 4]{18f-hw2s} and \cite[Exercise
2.1.1 \textbf{(a)} and \S 7.27]{19fco}). However, let us try to prove the
proposition combinatorially. For a bijective proof, the $\left(  -1\right)
^{k}$ and $\left(  -1\right)  ^{m}$ factors would be gamestoppers, as there is
no way to get a negative number (let alone a number of variable sign) by
counting something. However, if we think of the $\left(  -1\right)  ^{k}$ as
an opportunity for cancelling addends, then we can use combinatorics pretty well:

\begin{proof}
[Combinatorial proof of Proposition \ref{prop.binom.nhs} (sketched).]We need
to prove the equality (\ref{eq.prop.binom.nhs.eq}). Both sides of this
equality are polynomial functions in $n$. Thus, we can WLOG assume that $n$ is
a positive integer (because of the polynomial identity trick that we saw in
Subsection \ref{subsec.gf.defs.cvi}). Assume this.

Set $\left[  n\right]  :=\left\{  1,2,\ldots,n\right\}  $. We shall now
introduce some notations tailored for this particular proof.

Define an \emph{acceptable set} to be a subset of $\left[  n\right]  $ that
has size $\leq m$. Clearly,%
\begin{equation}
\left(  \text{\# of acceptable sets}\right)  =\sum_{k=0}^{m}\dbinom{n}{k}
\label{pf.prop.binom.nhs.1}%
\end{equation}
(since the \# of $k$-element subsets of $\left[  n\right]  $ equals
$\dbinom{n}{k}$ for each $k\in\mathbb{Z}$). Incidentally, we note that there
is no closed form for this sum -- another instance of the phenomenon in which
alternating sums are simpler than non-alternating ones.

Define the \emph{sign} of a finite set $I$ to be $\left(  -1\right)
^{\left\vert I\right\vert }$. Then,%
\begin{align}
&  \left(  \text{the sum of the signs of all acceptable sets}\right)
\nonumber\\
&  =\sum_{k=0}^{m}\left(  -1\right)  ^{k}\dbinom{n}{k}.
\label{pf.prop.binom.nhs.2}%
\end{align}
(This can be shown just like (\ref{pf.prop.binom.nhs.1}).)

However, the sum of the signs of all acceptable sets is a sum of $1$s and
$-1$s. Let us try to cancel as many of these $1$s and $-1$s against each other
as we can, hoping that what remains will be precisely $\left(  -1\right)
^{m}\dbinom{n-1}{m}$.

In order to cancel two addends, we need to pair up two finite sets $I$ that
have opposite signs. How do we find such pairs? One way to do so is to pick
some set $I$ that does not contain the element $1$, and pair it up with
$I\cup\left\{  1\right\}  $. Alternatively, we can pick some set $I$ that
contains the element $1$, and pair it up with $I\setminus\left\{  1\right\}
$. In other words, we pair up a finite set $I$ with either $I\setminus\left\{
1\right\}  $ or $I\cup\left\{  1\right\}  $, depending on whether $1\in I$ or
$1\notin I$.

Let us do this systematically for all finite sets: For each finite set $I$, we
define the \emph{partner} of $I$ to be the set%
\[
I^{\prime}:=%
\begin{cases}
I\setminus\left\{  1\right\}  , & \text{if }1\in I;\\
I\cup\left\{  1\right\}  , & \text{if }1\notin I
\end{cases}
\ \ \ \ =I\bigtriangleup\left\{  1\right\}  ,
\]
where the notation $X\bigtriangleup Y$ means the symmetric difference $\left(
X\cup Y\right)  \setminus\left(  X\cap Y\right)  $ of two sets $X$ and $Y$ (as
in Subsection \ref{subsec.gf.defs.commrings}). It is easy to see that each
finite set $I$ satisfies $I^{\prime\prime}=I$ and $\left\vert I^{\prime
}\right\vert =\left\vert I\right\vert \pm1$, so that $\left(  -1\right)
^{\left\vert I^{\prime}\right\vert }=-\left(  -1\right)  ^{\left\vert
I\right\vert }$. Thus, if both $I$ and $I^{\prime}$ are acceptable sets, then
their contributions to the sum of the signs of all acceptable sets cancel each
other out.

This does not mean that all addends in this sum cancel. Indeed, while each
finite set has a partner, it may happen that an acceptable set has a
non-acceptable partner, and then the contribution of the former to the sum
does not get cancelled (since the partner does not contribute to the sum).
Thus, in order to see what remains of the sum after the cancellations, we need
to study the acceptable sets that have non-acceptable partners.

Fortunately, this is easy: In order for an acceptable set $I$ to have a
non-acceptable partner, it needs to satisfy $1\notin I$ and $\left\vert
I\right\vert =m$. Better yet, this is an \textquotedblleft if and only
if\textquotedblright\ statement:

\begin{statement}
\textit{Claim 1:} Let $I$ be an acceptable set. Then, the partner $I^{\prime}$
of $I$ is non-acceptable if and only if $\left(  1\notin I\text{ and
}\left\vert I\right\vert =m\right)  $.
\end{statement}

[\textit{Proof of Claim 1:} The \textquotedblleft if\textquotedblright%
\ direction is easy: If $1\notin I$ and $\left\vert I\right\vert =m$, then the
partner $I^{\prime}$ of $I$ is defined by $I^{\prime}=I\cup\left\{  1\right\}
$ and thus has size $\left\vert I^{\prime}\right\vert =\left\vert I\right\vert
+1>\left\vert I\right\vert =m$, which shows that it is non-acceptable.

It remains to prove the converse, i.e., the \textquotedblleft only
if\textquotedblright\ direction. Thus, we assume that the partner $I^{\prime}$
of $I$ is non-acceptable. We must show that $1\notin I$ and $\left\vert
I\right\vert =m$.

Since $I$ is acceptable, we have $I\subseteq\left[  n\right]  $ and
$\left\vert I\right\vert \leq m$. If we had $1\in I$, then we would have
$I^{\prime}=I\setminus\left\{  1\right\}  \subseteq I\subseteq\left[
n\right]  $ and furthermore $\left\vert I^{\prime}\right\vert \leq\left\vert
I\right\vert $ (since $I^{\prime}\subseteq I$), so that $\left\vert I^{\prime
}\right\vert \leq\left\vert I\right\vert \leq m$. This would entail that
$I^{\prime}$ is acceptable (since $I^{\prime}\subseteq\left[  n\right]  $ and
$\left\vert I^{\prime}\right\vert \leq m$), which would contradict our
assumption that $I^{\prime}$ be non-acceptable. Hence, we cannot have $1\in
I$. Thus, $1\notin I$ is proved. Hence, $I^{\prime}=I\cup\left\{  1\right\}
$. However, $1\in\left[  n\right]  $ (since $n\geq1$), and $I^{\prime}%
=I\cup\left\{  1\right\}  \subseteq\left[  n\right]  $ (since $I\subseteq
\left[  n\right]  $ and $1\in\left[  n\right]  $). Moreover, from $I^{\prime
}=I\cup\left\{  1\right\}  $, we obtain $\left\vert I^{\prime}\right\vert
=\left\vert I\cup\left\{  1\right\}  \right\vert =\left\vert I\right\vert +1$.
Now, if we had $\left\vert I\right\vert \leq m-1$, then we would obtain
$\left\vert I^{\prime}\right\vert =\left\vert I\right\vert +1\leq m$ (since
$\left\vert I\right\vert \leq m-1$), which would entail that $I^{\prime}$ is
acceptable (since $I^{\prime}\subseteq\left[  n\right]  $), which again would
contradict our assumption. Thus, we cannot have $\left\vert I\right\vert \leq
m-1$. Hence, $\left\vert I\right\vert >m-1$, so that $\left\vert I\right\vert
\geq m$ and therefore $\left\vert I\right\vert =m$ (since $\left\vert
I\right\vert \leq m$). Thus, we have shown that $1\notin I$ and $\left\vert
I\right\vert =m$. This proves the \textquotedblleft only if\textquotedblright%
\ direction and thus completes the proof of Claim 1.]

Now, recall our line of reasoning: We start with the sum of the signs of all
acceptable sets, and we cancel any two addends that correspond to an
acceptable set and its acceptable partner. What remains are the addends
corresponding to the acceptable sets that have non-acceptable partners.
According to Claim 1, these are precisely the acceptable sets $I$ that satisfy
$\left(  1\notin I\text{ and }\left\vert I\right\vert =m\right)  $. In other
words, these are precisely the $m$-element subsets of $\left[  n\right]  $
that do not contain $1$. In other words, these are precisely the $m$-element
subsets of $\left[  n\right]  \setminus\left\{  1\right\}  $ (since a subset
of $\left[  n\right]  $ that does not contain $1$ is the same as a subset of
$\left[  n\right]  \setminus\left\{  1\right\}  $). Thus, there are precisely
$\dbinom{n-1}{m}$ of these subsets (since $\left[  n\right]  \setminus\left\{
1\right\}  $ is an $\left(  n-1\right)  $-element set), and each of them has
sign $\left(  -1\right)  ^{m}$. Hence, there are precisely $\dbinom{n-1}{m}$
addends left in the sum after our cancellations, and each of these addends is
$\left(  -1\right)  ^{m}$. Hence,
\[
\left(  \text{the sum of the signs of all acceptable sets}\right)  =\left(
-1\right)  ^{m}\dbinom{n-1}{m}.
\]
Comparing this with (\ref{pf.prop.binom.nhs.2}), we obtain
\[
\sum_{k=0}^{m}\left(  -1\right)  ^{k}\dbinom{n}{k}=\left(  -1\right)
^{m}\dbinom{n-1}{m}.
\]
This proves Proposition \ref{prop.binom.nhs}.
\end{proof}

Let me outline how to formalize this argument without using vague notions like
\textquotedblleft cancelling\textquotedblright\ and \textquotedblleft pairing
up\textquotedblright. We let%
\[
\mathcal{A}:=\left\{  \text{acceptable sets}\right\}
\]
and%
\begin{align*}
\mathcal{X}:=  &  \left\{  \text{acceptable sets whose partner is
acceptable}\right\} \\
=  &  \left\{  I\subseteq\left[  n\right]  \ \mid\ \left\vert I\right\vert
\leq m\text{ but not }\left(  \left\vert I\right\vert =m\text{ and }1\notin
I\right)  \right\}
\end{align*}
(by Claim 1 in the above proof). Now, we define a map%
\begin{align*}
f:\mathcal{X}  &  \rightarrow\mathcal{X},\\
I  &  \mapsto I^{\prime}.
\end{align*}
This map $f$ is a bijection, since each $I\in\mathcal{X}$ satisfies
$I^{\prime\prime}=I$ and thus $I^{\prime}\in\mathcal{X}$. This bijection $f$
is furthermore \emph{sign-reversing}, meaning that $\left(  -1\right)
^{\left\vert f\left(  I\right)  \right\vert }=-\left(  -1\right)  ^{\left\vert
I\right\vert }$ for all $I\in\mathcal{X}$. We claim that this automatically
guarantees%
\begin{align*}
&  \left(  \text{the sum of the signs of all acceptable sets}\right) \\
&  =\left(  \text{the sum of the signs of all acceptable sets \textbf{not}
belonging to }\mathcal{X}\right)  .
\end{align*}
The reason for this equality is that in the sum of the signs of all acceptable
sets, the contributions of the sets that belong to $\mathcal{X}$ (that is, of
the acceptable sets that have acceptable partners) cancel each other out. This
principle is worth generalizing and stating as a lemma:

\begin{lemma}
[Cancellation principle, take 1]\label{lem.sign.cancel1}Let $\mathcal{A}$ be a
finite set. Let $\mathcal{X}$ be a subset of $\mathcal{A}$.

For each $I\in\mathcal{A}$, let $\operatorname*{sign}I$ be a real number (not
necessarily $1$ or $-1$). Let $f:\mathcal{X}\rightarrow\mathcal{X}$ be a
bijection with the property that%
\begin{equation}
\operatorname*{sign}\left(  f\left(  I\right)  \right)  =-\operatorname*{sign}%
I\ \ \ \ \ \ \ \ \ \ \text{for all }I\in\mathcal{X}.
\label{eq.lem.sign.cancel1.sr}%
\end{equation}
(Such a bijection $f$ is called \emph{sign-reversing}.) Then,%
\[
\sum_{I\in\mathcal{A}}\operatorname*{sign}I=\sum_{I\in\mathcal{A}%
\setminus\mathcal{X}}\operatorname*{sign}I.
\]

\end{lemma}

Note that we did \textbf{not} require that $f\circ f=\operatorname*{id}$ in
Lemma \ref{lem.sign.cancel1}; we only required that $f$ is a bijection. That
said, most examples that I know do have $f\circ f=\operatorname*{id}$.

\begin{proof}
[Proof of Lemma \ref{lem.sign.cancel1}.]Intuitively, this is clear: The
contributions of all $I\in\mathcal{X}$ to the sum $\sum_{I\in\mathcal{A}%
}\operatorname*{sign}I$ cancel out, to the extent they are not already zero.
However, rather than formalize this cancellation, let us give an even slicker argument:

We have%
\begin{align*}
\sum_{I\in\mathcal{X}}\operatorname*{sign}I  &  =\sum_{I\in\mathcal{X}%
}\underbrace{\operatorname*{sign}\left(  f\left(  I\right)  \right)
}_{\substack{=-\operatorname*{sign}I\\\text{(by (\ref{eq.lem.sign.cancel1.sr}%
))}}}\\
&  \ \ \ \ \ \ \ \ \ \ \ \ \ \ \ \ \ \ \ \ \left(
\begin{array}
[c]{c}%
\text{here, we have substituted }f\left(  I\right)  \text{ for }I\\
\text{in the sum, since }f:\mathcal{X}\rightarrow\mathcal{X}\text{ is a
bijection}%
\end{array}
\right) \\
&  =\sum_{I\in\mathcal{X}}\left(  -\operatorname*{sign}I\right)  =-\sum
_{I\in\mathcal{X}}\operatorname*{sign}I.
\end{align*}
Adding $\sum_{I\in\mathcal{X}}\operatorname*{sign}I$ to both sides of this
equality, we obtain $2\cdot\sum_{I\in\mathcal{X}}\operatorname*{sign}I=0$.
Hence, $\sum_{I\in\mathcal{X}}\operatorname*{sign}I=0$ (since any real number
$a$ satisfying $2a=0$ must satisfy $a=0$).

Now, $\mathcal{X}\subseteq\mathcal{A}$; hence, we can split the sum
$\sum_{I\in\mathcal{A}}\operatorname*{sign}I$ as follows:%
\[
\sum_{I\in\mathcal{A}}\operatorname*{sign}I=\underbrace{\sum_{I\in\mathcal{X}%
}\operatorname*{sign}I}_{=0}+\sum_{I\in\mathcal{A}\setminus\mathcal{X}%
}\operatorname*{sign}I=\sum_{I\in\mathcal{A}\setminus\mathcal{X}%
}\operatorname*{sign}I.
\]
This proves Lemma \ref{lem.sign.cancel1}.
\end{proof}

In the proof of Proposition \ref{prop.binom.nhs}, we applied Lemma
\ref{lem.sign.cancel1} to%
\begin{align*}
\mathcal{A}  &  =\left\{  \text{acceptable sets}\right\}
\ \ \ \ \ \ \ \ \ \ \text{and}\\
\mathcal{X}  &  =\left\{  \text{acceptable sets whose partner is
acceptable}\right\}  \ \ \ \ \ \ \ \ \ \ \text{and}\\
\operatorname*{sign}I  &  =\left(  -1\right)  ^{\left\vert I\right\vert }.
\end{align*}
But there are many other situations in which Lemma \ref{lem.sign.cancel1} can
be applied. For example, $\mathcal{A}$ can be some set of permutations, and
$\operatorname*{sign}\sigma$ can be the sign of $\sigma$ (as in Definition
\ref{def.perm.sign}).

Let us observe that Lemma \ref{lem.sign.cancel1} can be generalized. Indeed,
in Lemma \ref{lem.sign.cancel1}, we can replace \textquotedblleft real
number\textquotedblright\ by \textquotedblleft element of any $\mathbb{Q}%
$-vector space\textquotedblright\ or even by \textquotedblleft element of any
additive abelian group with the property that $2a=0$ implies $a=0$%
\textquotedblright. We cannot, however, remove this requirement entirely.
Indeed, if all the signs $\operatorname*{sign}I$ were the element
$\overline{1}$ of $\mathbb{Z}/2$, then the sign-reversing condition
(\ref{eq.lem.sign.cancel1.sr}) would hold automatically (since $\overline
{1}=-\overline{1}$ in $\mathbb{Z}/2$), but the claim of Lemma
\ref{lem.sign.cancel1} would not necessarily be true.

However, if we replace the word \textquotedblleft bijection\textquotedblright%
\ by \textquotedblleft involution with no fixed points\textquotedblright, then
Lemma \ref{lem.sign.cancel1} holds even without any requirements on the group:

\begin{lemma}
[Cancellation principle, take 2]\label{lem.sign.cancel2}Let $\mathcal{A}$ be a
finite set. Let $\mathcal{X}$ be a subset of $\mathcal{A}$.

For each $I\in\mathcal{A}$, let $\operatorname*{sign}I$ be an element of some
additive abelian group. Let $f:\mathcal{X}\rightarrow\mathcal{X}$ be an
involution (i.e., a map satisfying $f\circ f=\operatorname*{id}$) that has no
fixed points. Assume that%
\[
\operatorname*{sign}\left(  f\left(  I\right)  \right)  =-\operatorname*{sign}%
I\ \ \ \ \ \ \ \ \ \ \text{for all }I\in\mathcal{X}.
\]
Then,%
\[
\sum_{I\in\mathcal{A}}\operatorname*{sign}I=\sum_{I\in\mathcal{A}%
\setminus\mathcal{X}}\operatorname*{sign}I.
\]

\end{lemma}

\begin{proof}
The idea is that all addends corresponding to the $I\in\mathcal{X}$ cancel out
from the sum $\sum_{I\in\mathcal{A}}\operatorname*{sign}I$ (because they come
in pairs of addends with opposite signs). See Section
\ref{sec.details.sign.intro} for a detailed proof.
\end{proof}

A more general version of Lemma \ref{lem.sign.cancel2} allows for $f$ to have
fixed points, as long as these fixed points have sign $0$:

\begin{lemma}
[Cancellation principle, take 3]\label{lem.sign.cancel3}Let $\mathcal{A}$ be a
finite set. Let $\mathcal{X}$ be a subset of $\mathcal{A}$.

For each $I\in\mathcal{A}$, let $\operatorname*{sign}I$ be an element of some
additive abelian group. Let $f:\mathcal{X}\rightarrow\mathcal{X}$ be an
involution (i.e., a map satisfying $f\circ f=\operatorname*{id}$). Assume that%
\[
\operatorname*{sign}\left(  f\left(  I\right)  \right)  =-\operatorname*{sign}%
I\ \ \ \ \ \ \ \ \ \ \text{for all }I\in\mathcal{X}.
\]
Assume furthermore that%
\[
\operatorname*{sign}I=0\ \ \ \ \ \ \ \ \ \ \text{for all }I\in\mathcal{X}%
\text{ satisfying }f\left(  I\right)  =I.
\]
Then,%
\[
\sum_{I\in\mathcal{A}}\operatorname*{sign}I=\sum_{I\in\mathcal{A}%
\setminus\mathcal{X}}\operatorname*{sign}I.
\]

\end{lemma}

\begin{proof}
This is similar to Lemma \ref{lem.sign.cancel2}, except that the addends
corresponding to the $I\in\mathcal{X}$ satisfying $f\left(  I\right)  =I$
don't cancel (but are already zero and thus can be removed right away). See
Section \ref{sec.details.sign.intro} for a detailed proof.
\end{proof}

Let us try to use this idea in another setting. Recall the notion of
$q$-binomial coefficients, and specifically their values (Definition
\ref{def.pars.qbinom.qbinom} \textbf{(b)}).

\begin{exercise}
\label{exe.sign.-1inom}Let $n,k\in\mathbb{N}$. Simplify $\dbinom{n}{k}_{-1}$.
\end{exercise}

\begin{example}
Let us compute $\dbinom{4}{2}_{-1}$. Theorem \ref{thm.pars.qbinom.quot1}
\textbf{(b)} yields%
\[
\dbinom{4}{2}_{q}=\dfrac{\left(  1-q^{4}\right)  \left(  1-q^{3}\right)
}{\left(  1-q^{2}\right)  \left(  1-q^{1}\right)  }.
\]
We cannot substitute $-1$ for $q$ in this formula directly, since both
numerator and denominator would become $0$ if we did. However, we can first
simplify the fraction and then substitute $-1$ for $q$: We have%
\[
\dbinom{4}{2}_{q}=\dfrac{\left(  1-q^{4}\right)  \left(  1-q^{3}\right)
}{\left(  1-q^{2}\right)  \left(  1-q^{1}\right)  }=q^{4}+q^{3}+2q^{2}+q+1,
\]
so that (by substituting $-1$ for $q$) we obtain%
\[
\dbinom{4}{2}_{-1}=\left(  -1\right)  ^{4}+\left(  -1\right)  ^{3}+2\left(
-1\right)  ^{2}+\left(  -1\right)  +1=2.
\]

\end{example}

\begin{proof}
[Solution of Exercise \ref{exe.sign.-1inom} (sketched).]Proposition
\ref{prop.pars.qbinom.alt-defs} \textbf{(b)} yields%
\[
\dbinom{n}{k}_{q}=\sum_{\substack{S\subseteq\left\{  1,2,\ldots,n\right\}
;\\\left\vert S\right\vert =k}}q^{\operatorname*{sum}S-\left(  1+2+\cdots
+k\right)  },
\]
where $\operatorname*{sum}S$ denotes the sum of the elements of a finite set
$S$ of integers. Substituting $-1$ for $q$ in this equality, we find%
\[
\dbinom{n}{k}_{-1}=\sum_{\substack{S\subseteq\left\{  1,2,\ldots,n\right\}
;\\\left\vert S\right\vert =k}}\left(  -1\right)  ^{\operatorname*{sum}%
S-\left(  1+2+\cdots+k\right)  }.
\]
Using the shorthand $\left[  n\right]  $ for the set $\left\{  1,2,\ldots
,n\right\}  $, we can rewrite this as%
\begin{equation}
\dbinom{n}{k}_{-1}=\sum_{\substack{S\subseteq\left[  n\right]  ;\\\left\vert
S\right\vert =k}}\left(  -1\right)  ^{\operatorname*{sum}S-\left(
1+2+\cdots+k\right)  }. \label{sol.sign.-1inom.1}%
\end{equation}
Let us analyze the sum on the right hand side using sign-reversing
involutions. Thus, we set%
\[
\mathcal{A}:=\left\{  S\subseteq\left[  n\right]  \ \mid\ \left\vert
S\right\vert =k\right\}  =\left\{  k\text{-element subsets of }\left[
n\right]  \right\}
\]
and%
\[
\operatorname*{sign}S:=\left(  -1\right)  ^{\operatorname*{sum}S-\left(
1+2+\cdots+k\right)  }\ \ \ \ \ \ \ \ \ \ \text{for every }S\in\mathcal{A}.
\]
Hence, (\ref{sol.sign.-1inom.1}) rewrites as%
\begin{equation}
\dbinom{n}{k}_{-1}=\sum_{S\in\mathcal{A}}\operatorname*{sign}S.
\label{sol.sign.-1inom.2}%
\end{equation}
Now, we seek a reasonable subset $\mathcal{X}\subseteq\mathcal{A}$ and a
sign-reversing bijection $f:\mathcal{X}\rightarrow\mathcal{X}$ in order to
cancel addends in the sum $\sum_{S\in\mathcal{A}}\operatorname*{sign}S$.

To wit, let us try to construct $f$ as a partial map first, and then (as an
afterthought) define $\mathcal{X}$ to be the set of all $S\in\mathcal{A}$ for
which $f\left(  S\right)  $ is defined.

Consider a $k$-element subset $S$ of $\left[  n\right]  $. What is a way to
transform $S$ that leaves its size $\left\vert S\right\vert =k$ unchanged, but
flips its sign (i.e., flips the parity of $\operatorname*{sum}S$) ? One thing
we can do is \emph{switching }$1$\emph{ with }$2$. By this I mean the
following operation:

\begin{itemize}
\item If $1\in S$ and $2\notin S$, then we replace $1$ by $2$ in $S$.

\item Otherwise, if $2\in S$ and $1\notin S$, then we replace $2$ by $1$ in
$S$.

\item If none of $1$ and $2$ is in $S$, or if both are in $S$, then we leave
$S$ unchanged for now.
\end{itemize}

\noindent Thus, switching $1$ with $2$ means replacing $S$ by%
\[
\operatorname*{switch}\nolimits_{1,2}\left(  S\right)  :=%
\begin{cases}
\left(  S\setminus\left\{  1\right\}  \right)  \cup\left\{  2\right\}  , &
\text{if }1\in S\text{ and }2\notin S;\\
\left(  S\setminus\left\{  2\right\}  \right)  \cup\left\{  1\right\}  , &
\text{if }1\notin S\text{ and }2\in S;\\
S, & \text{otherwise.}%
\end{cases}
\]
For example,%
\begin{align*}
\operatorname*{switch}\nolimits_{1,2}\left(  \left\{  1,3,5\right\}  \right)
&  =\left\{  2,3,5\right\}  ;\\
\operatorname*{switch}\nolimits_{1,2}\left(  \left\{  2,3,5\right\}  \right)
&  =\left\{  1,3,5\right\}  ;\\
\operatorname*{switch}\nolimits_{1,2}\left(  \left\{  1,2,5\right\}  \right)
&  =\left\{  1,2,5\right\}  ;\\
\operatorname*{switch}\nolimits_{1,2}\left(  \left\{  3,4,5\right\}  \right)
&  =\left\{  3,4,5\right\}  .
\end{align*}

Notice that the definition of switching $1$ with $2$ can be restated in a
somewhat slicker way using symmetric differences (see Subsection
\ref{subsec.gf.defs.commrings} for the definition of symmetric differences):%
\[
\operatorname*{switch}\nolimits_{1,2}\left(  S\right)  :=%
\begin{cases}
S\bigtriangleup\left\{  1,2\right\}  , & \text{if }\left\vert S\cap\left\{
1,2\right\}  \right\vert =1;\\
S, & \text{otherwise.}%
\end{cases}
\]
Indeed, the condition \textquotedblleft$\left\vert S\cap\left\{  1,2\right\}
\right\vert =1$\textquotedblright\ is equivalent to having%
\[
\text{either }\left(  1\in S\text{ and }2\notin S\right)  \text{ or }\left(
1\notin S\text{ and }2\in S\right)  ;
\]
and in this case, the symmetric difference $S\bigtriangleup\left\{
1,2\right\}  $ is precisely the set we need (i.e., the set $\left(
S\setminus\left\{  1\right\}  \right)  \cup\left\{  2\right\}  $ if we have
$\left(  1\in S\text{ and }2\notin S\right)  $, and the set $\left(
S\setminus\left\{  2\right\}  \right)  \cup\left\{  1\right\}  $ if we have
$\left(  1\notin S\text{ and }2\in S\right)  $).

This map $\operatorname*{switch}\nolimits_{1,2}:\mathcal{A}\rightarrow
\mathcal{A}$ is certainly a bijection (and, in fact, an involution). It is not
sign-reversing on the entire set $\mathcal{A}$; however, it has the property
that the sign of $\operatorname*{switch}\nolimits_{1,2}\left(  S\right)  $ is
opposite to the sign of $S$ whenever we have $\left(  1\in S\text{ and
}2\notin S\right)  $ or $\left(  1\notin S\text{ and }2\in S\right)  $
(because in these two cases, $\operatorname*{sum}S$ either increases by $1$ or
decreases by $1$, respectively). We can restate this property as follows: The
sign of $\operatorname*{switch}\nolimits_{1,2}\left(  S\right)  $ is opposite
to the sign of $S$ whenever we have $\left\vert S\cap\left\{  1,2\right\}
\right\vert =1$. Thus, we can use $\operatorname*{switch}\nolimits_{1,2}$ to
cancel many addends from our sum $\sum_{S\in\mathcal{A}}\operatorname*{sign}%
S$. Still, many other addends (of different signs) remain, and the result is
far from simple.

Thus, we need a \textquotedblleft Plan B\textquotedblright\ if the map
$\operatorname*{switch}\nolimits_{1,2}$ does not succeed. Assuming that
$\left\vert S\cap\left\{  1,2\right\}  \right\vert \neq1$ (that is, the set
$S\in\mathcal{A}$ contains none or both of $1$ and $2$), we gain nothing by
switching $1$ with $2$ in $S$, but maybe we get lucky switching $2$ with $3$
in $S$ (which is defined in the same way as switching $1$ with $2$, but with
the obvious changes)? If that, too, fails, we can try to switch $3$ with $4$.
If that fails as well, we can try to switch $4$ with $5$, and so on, until we
get to the end of the set $\left[  n\right]  $.

In other words, we try to define a bijection $f:\mathcal{A}\rightarrow
\mathcal{A}$ as follows: For any $S\in\mathcal{A}$, we pick the
\textbf{smallest} $i\in\left[  n-1\right]  $ such that $\left\vert
S\cap\left\{  i,i+1\right\}  \right\vert =1$ (in other words, the
\textbf{smallest} $i\in\left[  n-1\right]  $ such that exactly one of the two
elements $i$ and $i+1$ belongs to $S$); and we switch $i$ with $i+1$ in $S$
(that is, we replace $S$ by $S\bigtriangleup\left\{  i,i+1\right\}  $). This
produces a new subset $S^{\prime}$ of $\left[  n\right]  $ that has the same
size as $S$ but has the opposite sign (actually, we have $\operatorname*{sum}%
S^{\prime}=\operatorname*{sum}S\pm1$), except for the two cases when
$S=\varnothing$ and when $S=\left[  n\right]  $ (these are the cases where we
cannot find any $i\in\left[  n-1\right]  $ such that $\left\vert S\cap\left\{
i,i+1\right\}  \right\vert =1$). We set $f\left(  S\right)  :=S\bigtriangleup
\left\{  i,i+1\right\}  $.

Here are some examples (for $n=4$ and $k=2$):%
\begin{align*}
f\left(  \left\{  1,3\right\}  \right)   &  =\left\{  2,3\right\}
\ \ \ \ \ \ \ \ \ \ \left(  \text{here, the smallest }i\text{ is }1\right)
;\\
f\left(  \left\{  1,4\right\}  \right)   &  =\left\{  2,4\right\}
\ \ \ \ \ \ \ \ \ \ \left(  \text{here, the smallest }i\text{ is }1\right)
;\\
f\left(  \left\{  3,4\right\}  \right)   &  =\left\{  2,4\right\}
\ \ \ \ \ \ \ \ \ \ \left(  \text{here, the smallest }i\text{ is }2\right)  .
\end{align*}
Alas, the last two of these examples show that $f$ is not injective (as
$f\left(  \left\{  1,4\right\}  \right)  =\left\{  2,4\right\}  =f\left(
\left\{  3,4\right\}  \right)  $). Thus, $f$ is not a bijection. The
underlying problem is that the $i$ that was picked in the construction of
$f\left(  S\right)  $ is not uniquely recoverable from $f\left(  S\right)  $.
Hence, our map $f$ does not work for us -- we cannot use it to cancel addends,
since we cannot cancel (e.g.) a single $1$ against multiple $-1$s.

How can we salvage this argument? We change our map $f$ to \textquotedblleft
space\textquotedblright\ our switches apart. That is, we again start by trying
to switch $1$ with $2$; if this fails, we jump straight to trying to switch
$3$ with $4$; if this fails too, we jump further to trying to switch $5$ with
$6$; and so on, until we either succeed at some switch or run out of pairs to
switch. For the explicit description of $f$, this means that instead of
picking the \textbf{smallest} $i\in\left[  n-1\right]  $ such that $\left\vert
S\cap\left\{  i,i+1\right\}  \right\vert =1$, we pick the \textbf{smallest
odd} $i\in\left[  n-1\right]  $ such that $\left\vert S\cap\left\{
i,i+1\right\}  \right\vert =1$; and then we set $f\left(  S\right)
:=S\bigtriangleup\left\{  i,i+1\right\}  $ as before.

In other words, we define our new map $f:\mathcal{A}\rightarrow\mathcal{A}$ as
follows: For any $S\in\mathcal{A}$, we set%
\[
f\left(  S\right)  :=S\bigtriangleup\left\{  i,i+1\right\}  ,
\]
where $i$ is the \textbf{smallest odd} element of $\left[  n-1\right]  $ such
that $\left\vert S\cap\left\{  i,i+1\right\}  \right\vert =1$. If no such $i$
exists, we just set $f\left(  S\right)  :=S$. (We will soon see when this happens.)

Here are some examples (for $n=8$ and $k=3$):%
\begin{align*}
f\left(  \left\{  1,3,4\right\}  \right)   &  =\left\{  2,3,4\right\}
\ \ \ \ \ \ \ \ \ \ \left(  \text{here, the smallest odd }i\text{ is
}1\right)  ;\\
f\left(  \left\{  2,4,5\right\}  \right)   &  =\left\{  1,4,5\right\}
\ \ \ \ \ \ \ \ \ \ \left(  \text{here, the smallest odd }i\text{ is
}1\right)  ;\\
f\left(  \left\{  1,2,3\right\}  \right)   &  =\left\{  1,2,4\right\}
\ \ \ \ \ \ \ \ \ \ \left(  \text{here, the smallest odd }i\text{ is
}3\right)  ;\\
f\left(  \left\{  3,5,7\right\}  \right)   &  =\left\{  4,5,7\right\}
\ \ \ \ \ \ \ \ \ \ \left(  \text{here, the smallest odd }i\text{ is
}3\right)  ;\\
f\left(  \left\{  3,4,5\right\}  \right)   &  =\left\{  3,4,6\right\}
\ \ \ \ \ \ \ \ \ \ \left(  \text{here, the smallest odd }i\text{ is
}5\right)  ;\\
f\left(  \left\{  5,6,7\right\}  \right)   &  =\left\{  5,6,8\right\}
\ \ \ \ \ \ \ \ \ \ \left(  \text{here, the smallest odd }i\text{ is
}7\right)  .
\end{align*}
And here are two more examples (for $n=8$ and $k=4$):%
\begin{align*}
f\left(  \left\{  1,2,5,7\right\}  \right)   &  =\left\{  1,2,6,7\right\}
\ \ \ \ \ \ \ \ \ \ \left(  \text{here, the smallest odd }i\text{ is
}5\right)  ;\\
f\left(  \left\{  1,2,5,6\right\}  \right)   &  =\left\{  1,2,5,6\right\}
\ \ \ \ \ \ \ \ \ \ \left(  \text{here, there is no appropriate odd }i\right)
.
\end{align*}

Once again, it is clear that the set $f\left(  S\right)  $ has size $k$
whenever $S$ does. Hence, $f:\mathcal{A}\rightarrow\mathcal{A}$ is at least a
well-defined map. This time, the map $f$ is furthermore an involution (that
is, $f\circ f=\operatorname*{id}$). Here is a quick argument for this (details
are left to the reader): Since we have \textquotedblleft
spaced\textquotedblright\ the switches apart, they don't interfere with each
other. Thus, the $i$ that gets chosen in the construction of $f\left(
S\right)  $ will again get chosen in the construction of $f\left(  f\left(
S\right)  \right)  $ (since the elements of $f\left(  S\right)  $ that are
smaller than this $i$ will not have changed from $S$). Thus, the switch that
happens in the construction of $f\left(  f\left(  S\right)  \right)  $ undoes
the switch made in the construction of $f\left(  S\right)  $, and as a result,
the set $f\left(  f\left(  S\right)  \right)  $ will be $S$ again. This shows
that $f\circ f=\operatorname*{id}$.

Thus, $f$ is an involution, hence a bijection. Moreover, $\operatorname*{sign}%
\left(  f\left(  S\right)  \right)  =-\operatorname*{sign}S$ holds whenever
$f\left(  S\right)  \neq S$ (because $f\left(  S\right)  \neq S$ implies that
$f\left(  S\right)  =S\bigtriangleup\left\{  i,i+1\right\}  $ for some $i$
satisfying $\left\vert S\cap\left\{  i,i+1\right\}  \right\vert =1$, and
therefore $\operatorname*{sum}\left(  f\left(  S\right)  \right)
=\operatorname*{sum}S\pm1$). Thus, we set%
\[
\mathcal{X}:=\left\{  S\in\mathcal{A}\ \mid\ f\left(  S\right)  \neq
S\right\}  ,
\]
and we restrict $f$ to a map $\mathcal{X}\rightarrow\mathcal{X}$ (this is
well-defined, since it is easy to see from $f\circ f=\operatorname*{id}$ that
$f\left(  S\right)  \in\mathcal{X}$ for each $S\in\mathcal{X}$). Then, the map
$f$ becomes a sign-reversing bijection from $\mathcal{X}$ to $\mathcal{X}$.
Hence, Lemma \ref{lem.sign.cancel1} yields%
\[
\sum_{I\in\mathcal{A}}\operatorname*{sign}I=\sum_{I\in\mathcal{A}%
\setminus\mathcal{X}}\operatorname*{sign}I.
\]
Renaming the index $I$ as $S$, we can rewrite this equality as%
\[
\sum_{S\in\mathcal{A}}\operatorname*{sign}S=\sum_{S\in\mathcal{A}%
\setminus\mathcal{X}}\operatorname*{sign}S.
\]
Hence, (\ref{sol.sign.-1inom.2}) becomes%
\begin{equation}
\dbinom{n}{k}_{-1}=\sum_{S\in\mathcal{A}}\operatorname*{sign}S=\sum
_{S\in\mathcal{A}\setminus\mathcal{X}}\operatorname*{sign}S.
\label{sol.sign.-1inom.3}%
\end{equation}

Now, what is $\mathcal{A}\setminus\mathcal{X}$ ? In other words, what addends
are left behind uncancelled?

In order to answer this question, we need to consider the case when $n$ is
even and the case when $n$ is odd separately. We begin with the case when $n$
is even.

A $k$-element subset $S$ of $\left[  n\right]  $ belongs to $\mathcal{A}%
\setminus\mathcal{X}$ if and only if it satisfies $f\left(  S\right)  =S$. In
other words, $S$ belongs to $\mathcal{A}\setminus\mathcal{X}$ if and only if
there exists no odd $i\in\left[  n-1\right]  $ such that $\left\vert
S\cap\left\{  i,i+1\right\}  \right\vert =1$ (because $f$ has been defined in
such a way that $f\left(  S\right)  =S$ in this case, while $f\left(
S\right)  =S\bigtriangleup\left\{  i,i+1\right\}  \neq S$ in the other case).
In other words, $S$ belongs to $\mathcal{A}\setminus\mathcal{X}$ if and only
if for each odd $i\in\left[  n-1\right]  $, the size $\left\vert S\cap\left\{
i,i+1\right\}  \right\vert $ is either $0$ or $2$. This is equivalent to
saying that if we break up the $n$ elements $1,2,\ldots,n$ into $n/2$
\textquotedblleft blocks\textquotedblright\
\[
\left\{  1,2\right\}  ,\ \ \left\{  3,4\right\}  ,\ \ \left\{  5,6\right\}
,\ \ \ldots,\ \ \left\{  n-1,n\right\}
\]
(this can be done, since $n$ is even), then the intersection of $S$ with each
block has size $0$ or $2$. In other words, this is saying that the set $S$
consists of entire blocks (i.e., each block is either fully included in $S$ or
is disjoint from $S$). In other words, this is saying that the set $S$ is a
union (possibly empty) of blocks. We call a subset $S$ of $\left[  n\right]  $
\emph{blocky} if it satisfies this condition.\footnote{For example, $\left\{
3,4,5,6,9,10\right\}  $ is a blocky subset of $\left[  10\right]  $, whereas
$\left\{  2,3,5,6\right\}  $ is not (since it neither fully includes nor is
disjoint from the block $\left\{  1,2\right\}  $).} Thus, a $k$-element subset
$S$ of $\left[  n\right]  $ belongs to $\mathcal{A}\setminus\mathcal{X}$ if
and only if it is blocky. In other words, $\mathcal{A}\setminus\mathcal{X}$ is
the set of all blocky $k$-element subsets of $\left[  n\right]  $.

How many blocky $k$-element subsets does $\left[  n\right]  $ have, and what
are their signs? Any blocky subset of $\left[  n\right]  $ has the
form\footnote{The symbol \textquotedblleft$\sqcup$\textquotedblright\ means
\textquotedblleft disjoint union\textquotedblright\ (in our case, a union of
disjoint sets).}%
\[
\left\{  i_{1},i_{1}+1\right\}  \sqcup\left\{  i_{2},i_{2}+1\right\}
\sqcup\cdots\sqcup\left\{  i_{p},i_{p}+1\right\}
\]
for some odd elements $i_{1}<i_{2}<\cdots<i_{p}$ of $\left[  n-1\right]  $.
Thus, this subset has size $2p$, which entails that its size is even. Hence,
if $k$ is odd, there are no blocky $k$-element subsets of $\left[  n\right]  $
at all. In other words, if $k$ is odd, then there are no $S\in\mathcal{A}%
\setminus\mathcal{X}$ (since $\mathcal{A}\setminus\mathcal{X}$ is the set of
all blocky $k$-element subsets of $\left[  n\right]  $). Therefore, if $k$ is
odd, then the sum $\sum_{S\in\mathcal{A}\setminus\mathcal{X}}%
\operatorname*{sign}S$ is empty, and thus (\ref{sol.sign.-1inom.3}) simplifies
to%
\begin{equation}
\dbinom{n}{k}_{-1}=\sum_{S\in\mathcal{A}\setminus\mathcal{X}}%
\operatorname*{sign}S=\left(  \text{empty sum}\right)  =0.
\label{sol.sign.-1inom.ev/od}%
\end{equation}

Let us now consider the case when $k$ is even. In this case, again, any blocky
subset $S$ of $\left[  n\right]  $ has the form%
\[
\left\{  i_{1},i_{1}+1\right\}  \sqcup\left\{  i_{2},i_{2}+1\right\}
\sqcup\cdots\sqcup\left\{  i_{p},i_{p}+1\right\}
\]
for some odd elements $i_{1}<i_{2}<\cdots<i_{p}$ of $\left[  n-1\right]  $.
Moreover, again, this subset has size $2p$. Thus, if $S$ has size $k$, then we
must have $2p=k$, so that $p=k/2$. Thus, any blocky $k$-element subset $S$ of
$\left[  n\right]  $ has the form
\[
\left\{  i_{1},i_{1}+1\right\}  \sqcup\left\{  i_{2},i_{2}+1\right\}
\sqcup\cdots\sqcup\left\{  i_{k/2},i_{k/2}+1\right\}
\]
for some odd elements $i_{1}<i_{2}<\cdots<i_{k/2}$ of $\left[  n-1\right]  $.
Hence, there are $\dbinom{n/2}{k/2}$ such subsets $S$ (since there are
precisely $\dbinom{n/2}{k/2}$ choices for these odd elements $i_{1}%
<i_{2}<\cdots<i_{k/2}$\ \ \ \ \footnote{because the set $\left[  n-1\right]  $
has $n/2$ odd elements}). Moreover, any such subset $S$ satisfies%
\begin{align*}
\operatorname*{sum}S  &  =\underbrace{i_{1}+\left(  i_{1}+1\right)  }%
_{\equiv1\operatorname{mod}2}+\underbrace{i_{2}+\left(  i_{2}+1\right)
}_{\equiv1\operatorname{mod}2}+\cdots+\underbrace{i_{k/2}+\left(
i_{k/2}+1\right)  }_{\equiv1\operatorname{mod}2}\\
&  \equiv\underbrace{1+1+\cdots+1}_{k/2\text{ times}}=k/2\operatorname{mod}2
\end{align*}
and therefore%
\[
\underbrace{\operatorname*{sum}S}_{\equiv k/2\operatorname{mod}2}%
-\underbrace{\left(  1+2+\cdots+k\right)  }_{=\dfrac{k\left(  k+1\right)  }%
{2}}\equiv k/2-\dfrac{k\left(  k+1\right)  }{2}=-k^{2}/2\equiv
0\operatorname{mod}2
\]
(since $k$ is even, so that $-k^{2}/2$ is even), and%
\begin{equation}
\operatorname*{sign}S=\left(  -1\right)  ^{\operatorname*{sum}S-\left(
1+2+\cdots+k\right)  }=1 \label{sol.sign.-1inom.od/ev.sign}%
\end{equation}
(since $\operatorname*{sum}S-\left(  1+2+\cdots+k\right)  \equiv
0\operatorname{mod}2$). Hence, the sum $\sum_{S\in\mathcal{A}\setminus
\mathcal{X}}\operatorname*{sign}S$ has $\dbinom{n/2}{k/2}$ addends (because
$\mathcal{A}\setminus\mathcal{X}$ is the set of all blocky $k$-element subsets
of $\left[  n\right]  $, and we have just shown that there are $\dbinom
{n/2}{k/2}$ such subsets), and each of these addends is $1$ (by
(\ref{sol.sign.-1inom.od/ev.sign})). Thus, this sum simplifies to%
\[
\sum_{S\in\mathcal{A}\setminus\mathcal{X}}\operatorname*{sign}S=\dbinom
{n/2}{k/2}\cdot1=\dbinom{n/2}{k/2}.
\]
Hence, (\ref{sol.sign.-1inom.3}) becomes%
\begin{equation}
\dbinom{n}{k}_{-1}=\sum_{S\in\mathcal{A}\setminus\mathcal{X}}%
\operatorname*{sign}S=\dbinom{n/2}{k/2}. \label{sol.sign.-1inom.ev/ev}%
\end{equation}

Thus, we have computed $\dbinom{n}{k}_{-1}$

\begin{itemize}
\item in the case when $n$ is even and $k$ is odd (obtaining
(\ref{sol.sign.-1inom.ev/od})), and

\item in the case when $n$ is even and $k$ is even (obtaining
(\ref{sol.sign.-1inom.ev/ev})).
\end{itemize}

It remains to handle the case when $n$ is odd. This case is different in that
the $n$ elements $1,2,\ldots,n$ are now subdivided into $\left(  n+1\right)
/2$ \textquotedblleft blocks\textquotedblright%
\[
\left\{  1,2\right\}  ,\ \ \left\{  3,4\right\}  ,\ \ \left\{  5,6\right\}
,\ \ \ldots,\ \ \left\{  n-2,n-1\right\}  ,\ \ \left\{  n\right\}  ,
\]
with the last of these blocks having size $1$. As a consequence, this time, a
blocky subset of $\left[  n\right]  $ can have odd size. Moreover, the parity
of $k$ determines whether a blocky $k$-element subset of $\left[  n\right]  $
will contain $n$:

\begin{itemize}
\item If $k$ is even, then no blocky $k$-element subset of $\left[  n\right]
$ can contain $n$ (because if it did, then it would have odd size, since all
non-$\left\{  n\right\}  $ blocks have even size).

\item If $k$ is odd, then every blocky $k$-element subset of $\left[
n\right]  $ must contain $n$ (because if it didn't, then it would have even
size, since all non-$\left\{  n\right\}  $ blocks have even size).
\end{itemize}

Thus, when classifying the blocky $k$-element subsets of $\left[  n\right]  $,
we can either dismiss $n$ immediately (if $k$ is even) or take $n$ for granted
(if $k$ is odd); in either case, the problem gets reduced to classifying the
blocky $k$-element or $\left(  k-1\right)  $-element subsets of $\left[
n-1\right]  $, which we already know how to do (since $n-1$ is even). The
result is that the \# of blocky $k$-element subsets of $\left[  n\right]  $
(in the case when $n$ is odd) is%
\[%
\begin{cases}
\dbinom{\left(  n-1\right)  /2}{k/2}, & \text{if }k\text{ is even;}\\
\dbinom{\left(  n-1\right)  /2}{\left(  k-1\right)  /2}, & \text{if }k\text{
is odd}%
\end{cases}
\ \ \ \ \ =\dbinom{\left(  n-1\right)  /2}{\left\lfloor k/2\right\rfloor },
\]
and their signs are always $1$. Hence, we obtain%
\[
\sum_{S\in\mathcal{A}\setminus\mathcal{X}}\operatorname*{sign}S=\dbinom
{\left(  n-1\right)  /2}{\left\lfloor k/2\right\rfloor }\cdot1=\dbinom{\left(
n-1\right)  /2}{\left\lfloor k/2\right\rfloor }.
\]
Thus, (\ref{sol.sign.-1inom.3}) becomes%
\begin{equation}
\dbinom{n}{k}_{-1}=\sum_{S\in\mathcal{A}\setminus\mathcal{X}}%
\operatorname*{sign}S=\dbinom{\left(  n-1\right)  /2}{\left\lfloor
k/2\right\rfloor }. \label{sol.sign.-1inom.od/*}%
\end{equation}
This is the answer to our exercise in the case when $n$ is odd.

Returning to the general case, we can now combine the formulas
(\ref{sol.sign.-1inom.ev/ev}), (\ref{sol.sign.-1inom.ev/od}) and
(\ref{sol.sign.-1inom.od/*}) into a single equality that holds in all cases:%
\begin{equation}
\dbinom{n}{k}_{-1}=%
\begin{cases}
0, & \text{if }n\text{ is even and }k\text{ is odd;}\\
\dbinom{\left\lfloor n/2\right\rfloor }{\left\lfloor k/2\right\rfloor }, &
\text{otherwise.}%
\end{cases}
\label{sol.sign.-1inom.res}%
\end{equation}

\end{proof}

This formula (\ref{sol.sign.-1inom.res}) can be generalized. Indeed, here is a
generalization of the number $-1$:

\begin{definition}
\label{def.root-of-unity.prim}Let $K$ be a field. Let $d$ be a positive
integer. \medskip

\textbf{(a)} A $d$\emph{-th root of unity} in $K$ means an element $\omega$ of
$K$ satisfying $\omega^{d}=1$. In other words, a $d$-th root of unity in $K$
means an element of $K$ whose $d$-th power is $1$. \medskip

\textbf{(b)} A \emph{primitive }$d$\emph{-th root of unity} in $K$ means an
element $\omega$ of $K$ satisfying
\[
\omega^{d}=1
\]
but
\[
\omega^{i}\neq1\ \ \ \ \ \ \ \ \ \ \text{for each }i\in\left\{  1,2,\ldots
,d-1\right\}  .
\]

In other words, a primitive $d$-th root of unity in $K$ means an element of
the multiplicative group $K^{\times}$ whose order is $d$.
\end{definition}

For $K=\mathbb{C}$, the $d$-th roots of unity are the $d$ complex numbers
\newline$e^{2\pi i0/d},e^{2\pi i1/d},e^{2\pi i2/d},\ldots,e^{2\pi i\left(
d-1\right)  /d}$ (which are the vertices of a regular $d$-gon inscribed in the
unit circle, with one vertex at $1$), whereas the primitive $d$-th roots of
unity are the numbers $e^{2\pi ig/d}$ for all $g\in\left[  d\right]  $
satisfying $\gcd\left(  g,d\right)  =1$. In particular, the $2$-nd roots of
unity in $\mathbb{C}$ are $1$ and $-1$, and the only primitive $2$-nd root of
unity is $-1$. The following picture shows the six $6$-th roots of unity in
$\mathbb{C}$ (with the two primitive $6$-th roots colored blue, and the
remaining four roots colored red):%
\[%
\begin{tikzpicture}[scale=2]
\draw[densely dotted] (-1.2,-1.2) grid (1.2,1.2);
\draw[->] (-1.2,0) -- (1.2,0);
\draw[->] (0,-1.2) -- (0,1.2);
\draw(0, 0) circle[radius=1, black];
\fill[darkred] (60*0 : 1) circle (1pt) node[above right] {$e^{2\pi i 0/6}$};
\fill[dbluecolor] (60*1 : 1) circle (1pt) node[above right] {$e^{2\pi i 1/6}%
$};
\fill[darkred] (60*2 : 1) circle (1pt) node[above left] {$e^{2\pi i 2/6}$};
\fill[darkred] (60*3 : 1) circle (1pt) node[above left] {$e^{2\pi i 3/6}$};
\fill[darkred] (60*4 : 1) circle (1pt) node[below left] {$e^{2\pi i 4/6}$};
\fill[dbluecolor] (60*5 : 1) circle (1pt) node[right] {$e^{2\pi i 5/6}$};
\end{tikzpicture}%
\ \ \ .
\]

We can now generalize (\ref{sol.sign.-1inom.res}) by replacing $-1$ by
primitive roots of unity:

\begin{theorem}
[$q$-Lucas theorem]\label{thm.sign.q-lucas}Let $K$ be a field. Let $d$ be a
positive integer. Let $\omega$ be a primitive $d$-th root of unity in $K$. Let
$n,k\in\mathbb{N}$. Let $q$ and $u$ be the quotients obtained when dividing
$n$ and $k$ by $d$ with remainder, and let $r$ and $v$ be the respective
remainders. Then,%
\begin{equation}
\dbinom{n}{k}_{\omega}=\dbinom{q}{u}\cdot\dbinom{r}{v}_{\omega}.
\label{eq.thm.sign.q-lucas.eq}%
\end{equation}

\end{theorem}

Note that the equality (\ref{eq.thm.sign.q-lucas.eq}) contains two $\omega
$-binomial coefficients and one regular binomial coefficient.

It is not hard to check that (\ref{sol.sign.-1inom.res}) is the particular
case of Theorem \ref{thm.sign.q-lucas} for $d=2$ and $\omega=-1$. Indeed, the
only possible remainders of an integer upon division by $2$ are $0$ and $1$,
and the $\omega$-binomial coefficients $\dbinom{r}{v}_{\omega}$ for
$r,v\in\left\{  0,1\right\}  $ are%
\[
\dbinom{0}{0}_{\omega}=1,\ \ \ \ \ \ \ \ \ \ \dbinom{0}{1}_{\omega
}=0,\ \ \ \ \ \ \ \ \ \ \dbinom{1}{0}_{\omega}=1,\ \ \ \ \ \ \ \ \ \ \dbinom
{1}{1}_{\omega}=1.
\]

Theorem \ref{thm.sign.q-lucas} can be proved using a generalization of
sign-reversing involutions to $d$-cycles instead of pairs\footnote{In our
proofs so far, we have been using the equality $1+\left(  -1\right)  =0$ to
cancel addends in alternating sums. This equality can be generalized as
follows: If $\omega$ is a primitive $d$-th root of unity for $d>1$, then
$1+\omega+\omega^{2}+\cdots+\omega^{d-1}=0$. This equality can be used to
cancel addends in sums that involve powers of $\omega$. The \emph{discrete
Fourier transform} (see, e.g., \cite[\S 5.6]{OlvSha}) and the
\emph{roots-of-unity filter} (see, e.g., \cite[\S 1.2.9.D]{Knuth-TAoCP1}) are
applications of this idea.}. A more algebraic proof -- using \textquotedblleft
noncommutative generating functions\textquotedblright\ -- is given in Exercise
\ref{exe.sign.qbinom.qlucas}.

\begin{noncompile}
TODO: Do the combinatorial proof of Theorem \ref{thm.sign.q-lucas} using a
Fourier generalization of sign-reversing involutions.
\end{noncompile}

\subsection{\label{sec.sign.pie}The principles of inclusion and exclusion}

We have so far been applying the \textquotedblleft yoga\textquotedblright\ of
sign-reversing involutions directly to alternating sums. However, some of its
essence can also be crystallized into useful theorems. Most famous among such
theorems are the \emph{principles of inclusion and exclusion} (also known as
\emph{Sylvester sieve theorems} or \emph{Poincar\'{e}'s theorems}%
)\footnote{People typically use the singular forms (\textquotedblleft
principle\textquotedblright\ and \textquotedblleft theorem\textquotedblright%
\ rather than \textquotedblleft principles\textquotedblright\ and
\textquotedblleft theorems\textquotedblright), but often mean different
things.}.

\subsubsection{\label{subsec.sign.pie.size}The size version}

The simplest such principle answers the following question: Assume that you
have a finite set $U$, and some subsets $A_{1},A_{2},\ldots,A_{n}$ of $U$.
Assume that, for each selection of some of these subsets, you know how many
elements of $U$ belong to all selected sets at the same time. (For instance,
you know how many elements of $U$ belong to $A_{2}$, $A_{3}$ and $A_{5}$ at
the same time.) Does this help you count the elements of $U$ that belong to
\textbf{none} of the $n$ subsets (i.e., that don't belong to $A_{1}\cup
A_{2}\cup\cdots\cup A_{n}$) ?

The answer is \textquotedblleft yes\textquotedblright, and there is an
explicit formula for this count:

\begin{theorem}
[size version of the PIE]\label{thm.pie.1}Let $n\in\mathbb{N}$. Let $U$ be a
finite set. Let $A_{1},A_{2},\ldots,A_{n}$ be $n$ subsets of $U$. Then,%
\begin{align*}
&  \left(  \text{\# of }u\in U\text{ that satisfy }u\notin A_{i}\text{ for all
}i\in\left[  n\right]  \right) \\
&  =\sum_{I\subseteq\left[  n\right]  }\left(  -1\right)  ^{\left\vert
I\right\vert }\left(  \text{\# of }u\in U\text{ that satisfy }u\in A_{i}\text{
for all }i\in I\right)  .
\end{align*}

\end{theorem}

Some explanations are in order:

\begin{itemize}
\item Here and in the following, we are using the notation $\left[  n\right]
$ for the set $\left\{  1,2,\ldots,n\right\}  $, as defined in Definition
\ref{def.perm.Sn-iven}.

\item The summation sign \textquotedblleft$\sum_{I\subseteq\left[  n\right]
}$\textquotedblright\ means a sum over all subsets $I$ of $\left[  n\right]
$. More generally, if $S$ is a given set, then the summation sign
\textquotedblleft$\sum_{I\subseteq S}$\textquotedblright\ shall always mean a
sum over all subsets $I$ of $S$.

\item The shorthand \textquotedblleft\emph{PIE}\textquotedblright\ in the name
of Theorem \ref{thm.pie.1} is short for \textquotedblleft\emph{Principle of
Inclusion and Exclusion}\textquotedblright.
\end{itemize}

In one form or another, Theorem \ref{thm.pie.1} appears in almost any text on
combinatorics (e.g., in \cite[Theorem 2.1.1]{Sagan19}, in \cite[\S 4.11]%
{Loehr-BC}, in \cite[Theorem 5.3]{Strick20}, in \cite[Theorem 2.9.7]{19fco},
or -- in almost the same form as above -- in \cite[Theorem 7.8.6]{20f}). Most
commonly, its claim is stated in the shorter (if less transparent) form%
\[
\left\vert U\setminus\left(  A_{1}\cup A_{2}\cup\cdots\cup A_{n}\right)
\right\vert =\sum_{I\subseteq\left[  n\right]  }\left(  -1\right)
^{\left\vert I\right\vert }\left\vert \bigcap\limits_{i\in I}A_{i}\right\vert
,
\]
where $\bigcap\limits_{i\in I}A_{i}$ denotes the set $\left\{  u\in
U\ \mid\ u\in A_{i}\text{ for all }i\in I\right\}  $ whenever $I$ is a subset
of $\left[  n\right]  $\ \ \ \ \footnote{This notation should be taken with a
grain of salt. When $I$ is a \textbf{nonempty} subset of $\left[  n\right]  $,
it is indeed true that the set $\bigcap\limits_{i\in I}A_{i}$ as we just
defined it (i.e., the set $\left\{  u\in U\ \mid\ u\in A_{i}\text{ for all
}i\in I\right\}  $) is the intersection of the sets $A_{i}$ over all $i\in I$
(that is, if $I=\left\{  i_{1},i_{2},\ldots,i_{k}\right\}  $, then
$\bigcap\limits_{i\in I}A_{i}=A_{i_{1}}\cap A_{i_{2}}\cap\cdots\cap A_{i_{k}}%
$). However, if $I$ is the \textbf{empty} set, then the literal intersection
of the sets $A_{i}$ over all $i\in I$ is not well-defined (indeed, by common
sense, such an intersection should contain every object whatsoever; but there
is no set that does this), whereas the set we just called $\bigcap
\limits_{i\in I}A_{i}$ is simply the set $U$. This is still justified, since
we should think of the sets $A_{i}$ not as arbitrary sets but rather as
subsets of $U$ (so that any intersections are to be taken within $U$). This,
incidentally, is the reason for the choice of the letter \textquotedblleft%
$U$\textquotedblright: We think of $U$ as the \textquotedblleft
universe\textquotedblright\ in which our objects live.}. This form is indeed
equivalent to the claim of Theorem \ref{thm.pie.1}, since we have%
\begin{align*}
&  \left\vert U\setminus\left(  A_{1}\cup A_{2}\cup\cdots\cup A_{n}\right)
\right\vert \\
&  =\left(  \text{\# of }u\in U\text{ that satisfy }u\notin A_{1}\cup
A_{2}\cup\cdots\cup A_{n}\right) \\
&  =\left(  \text{\# of }u\in U\text{ that satisfy }u\notin A_{i}\text{ for
all }i\in\left[  n\right]  \right)
\end{align*}
and since (for each subset $I$ of $\left[  n\right]  $) we have%
\[
\left\vert \bigcap\limits_{i\in I}A_{i}\right\vert =\left(  \text{\# of }u\in
U\text{ that satisfy }u\in A_{i}\text{ for all }i\in I\right)
\]
(by the definition of $\bigcap\limits_{i\in I}A_{i}$ that we just gave).

Rather than prove Theorem \ref{thm.pie.1} directly, we shall soon derive it
from more general results (in order to avoid duplicating arguments). First,
however, let us give an interpretation that makes Theorem \ref{thm.pie.1} a
little bit more intuitive, and sketch four applications (more can be found in
textbooks -- e.g., \cite[\S 2.1]{Sagan19}, \cite[Chapter 2]{Stanley-EC1},
\cite[Chapter 3]{Wildon19}, \cite[Chapter 6]{AndFen04}, \cite[\S 2.9]{19fco}).

\begin{statement}
\textbf{\textquotedblleft Rule-breaking\textquotedblright\ interpretation of
Theorem \ref{thm.pie.1}.} Assume that we are given a finite set $U$, and we
are given $n$ rules (labelled $1,2,\ldots,n$) that each element of $U$ may or
may not satisfy. (For instance, one element of $U$ might satisfy all of these
rules; another might satisfy none; yet another might satisfy rules $1$ and $3$
only. A rule can be something like \textquotedblleft thou shalt be divisible
by $5$\textquotedblright\ (if the elements of $U$ are numbers) or
\textquotedblleft thou shalt be a nonempty set\textquotedblright\ (if the
elements of $U$ are sets).)

Assume that, for each $I\subseteq\left[  n\right]  $, we know how many
elements $u\in U$ satisfy all rules in $I$ (but may or may not satisfy the
remaining rules). For example, this means that we know how many elements $u\in
U$ satisfy rules $2,3,5$ (simultaneously). Then, we can compute the \# of
elements $u\in U$ that violate all $n$ rules $1,2,\ldots,n$ by the following
formula:%
\begin{align*}
&  \left(  \text{\# of elements }u\in U\text{ that violate all }n\text{ rules
}1,2,\ldots,n\right) \\
&  =\sum_{I\subseteq\left[  n\right]  }\left(  -1\right)  ^{\left\vert
I\right\vert }\left(  \text{\# of elements }u\in U\text{ that satisfy all
rules in }I\right)  .
\end{align*}

Indeed, this formula is precisely what we obtain if we apply Theorem
\ref{thm.pie.1} to the $n$ subsets $A_{1},A_{2},\ldots,A_{n}$ defined by
setting%
\[
A_{i}:=\left\{  u\in U\ \mid\ u\text{ satisfies rule }i\right\}
\ \ \ \ \ \ \ \ \ \ \text{for each }i\in\left[  n\right]  .
\]

\end{statement}

Thus, if you have a counting problem that can be restated as \textquotedblleft
count things that violate a bunch of rules\textquotedblright, then you can
apply Theorem \ref{thm.pie.1} (in the interpretation we just gave) to
\textquotedblleft turn the problem positive\textquotedblright, i.e., to make
it about counting rule-followers instead of rule-violators. If the
\textquotedblleft positive\textquotedblright\ problem is easier, then this is
a useful technique. We will now witness this on four examples.

\subsubsection{\label{subsec.sign.pie.exas}Examples}

\textbf{Example 1.} Let $n,m\in\mathbb{N}$. Let us compute the \# of
surjective maps from $\left[  m\right]  $ to $\left[  n\right]  $. (We will
outline the argument here; details can be found in \cite[\S 7.8.2]{20f} or
\cite[\S 2.9.4]{19fco}.)

What are surjective maps? They are maps that take each element of the target
set as a value. Thus, in particular, a map $f:\left[  m\right]  \rightarrow
\left[  n\right]  $ is surjective if and only if it takes each $i\in\left[
n\right]  $ as a value.

Hence, if we impose $n$ rules $1,2,\ldots,n$ on a map $f:\left[  m\right]
\rightarrow\left[  n\right]  $, where rule $i$ says \textquotedblleft thou
shalt not take $i$ as a value\textquotedblright, then the surjective maps
$f:\left[  m\right]  \rightarrow\left[  n\right]  $ are precisely the maps
$f:\left[  m\right]  \rightarrow\left[  n\right]  $ that violate all $n$
rules. Hence,%
\begin{align*}
&  \left(  \text{\# of surjective maps }f:\left[  m\right]  \rightarrow\left[
n\right]  \right) \\
&  =\left(  \text{\# of maps }f:\left[  m\right]  \rightarrow\left[  n\right]
\text{ that violate all }n\text{ rules }1,2,\ldots,n\right) \\
&  =\sum_{I\subseteq\left[  n\right]  }\left(  -1\right)  ^{\left\vert
I\right\vert }\left(  \text{\# of maps }f:\left[  m\right]  \rightarrow\left[
n\right]  \text{ that satisfy all rules in }I\right)
\end{align*}
(by the \textquotedblleft rule-breaking\textquotedblright\ interpretation of
Theorem \ref{thm.pie.1}).

Now, fix some subset $I$ of $\left[  n\right]  $. What is the \# of maps
$f:\left[  m\right]  \rightarrow\left[  n\right]  $ that satisfy all rules in
$I$ ? A map $f:\left[  m\right]  \rightarrow\left[  n\right]  $ satisfies all
rules in $I$ if and only if it takes none of the $i\in I$ as a value, i.e., if
all its values belong to $\left[  n\right]  \setminus I$. Thus, the maps
$f:\left[  m\right]  \rightarrow\left[  n\right]  $ that satisfy all rules in
$I$ are nothing but the maps from $\left[  m\right]  $ to $\left[  n\right]
\setminus I$. The \# of such maps is therefore $\left\vert \left[  n\right]
\setminus I\right\vert ^{\left\vert \left[  m\right]  \right\vert }=\left(
n-\left\vert I\right\vert \right)  ^{m}$ (since $I\subseteq\left[  n\right]  $
entails $\left\vert \left[  n\right]  \setminus I\right\vert =\left\vert
\left[  n\right]  \right\vert -\left\vert I\right\vert =n-\left\vert
I\right\vert $, and since $\left\vert \left[  m\right]  \right\vert =m$).

Forget that we fixed $I$. We thus have shown that for each subset $I$ of
$\left[  n\right]  $, we have%
\begin{align}
&  \left(  \text{\# of maps }f:\left[  m\right]  \rightarrow\left[  n\right]
\text{ that satisfy all rules in }I\right) \nonumber\\
&  =\left(  n-\left\vert I\right\vert \right)  ^{m}. \label{eq.exa.pie1.1.3}%
\end{align}
Substituting this into the above computation, we find%
\begin{align*}
&  \left(  \text{\# of surjective maps }f:\left[  m\right]  \rightarrow\left[
n\right]  \right) \\
&  =\sum_{I\subseteq\left[  n\right]  }\left(  -1\right)  ^{\left\vert
I\right\vert }\underbrace{\left(  \text{\# of maps }f:\left[  m\right]
\rightarrow\left[  n\right]  \text{ that satisfy all rules in }I\right)
}_{\substack{=\left(  n-\left\vert I\right\vert \right)  ^{m}\\\text{(by
(\ref{eq.exa.pie1.1.3}))}}}\\
&  =\sum_{I\subseteq\left[  n\right]  }\left(  -1\right)  ^{\left\vert
I\right\vert }\left(  n-\left\vert I\right\vert \right)  ^{m}=\sum_{k=0}%
^{n}\ \ \sum_{\substack{I\subseteq\left[  n\right]  ;\\\left\vert I\right\vert
=k}}\underbrace{\left(  -1\right)  ^{\left\vert I\right\vert }\left(
n-\left\vert I\right\vert \right)  ^{m}}_{\substack{=\left(  -1\right)
^{k}\left(  n-k\right)  ^{m}\\\text{(since }\left\vert I\right\vert
=k\text{)}}}\\
&  \ \ \ \ \ \ \ \ \ \ \ \ \ \ \ \ \ \ \ \ \left(
\begin{array}
[c]{c}%
\text{here, we have split the sum according to}\\
\text{the value of }\left\vert I\right\vert
\end{array}
\right) \\
&  =\sum_{k=0}^{n}\ \ \underbrace{\sum_{\substack{I\subseteq\left[  n\right]
;\\\left\vert I\right\vert =k}}\left(  -1\right)  ^{k}\left(  n-k\right)
^{m}}_{\substack{=\dbinom{n}{k}\left(  -1\right)  ^{k}\left(  n-k\right)
^{m}\\\text{(since this is a sum of }\dbinom{n}{k}\\\text{many equal
addends)}}}=\sum_{k=0}^{n}\dbinom{n}{k}\left(  -1\right)  ^{k}\left(
n-k\right)  ^{m}\\
&  =\sum_{k=0}^{n}\left(  -1\right)  ^{k}\dbinom{n}{k}\left(  n-k\right)
^{m}.
\end{align*}

Thus, we have proved the following theorem:

\begin{theorem}
\label{thm.pie.count-sur}Let $n,m\in\mathbb{N}$. Then,%
\[
\left(  \text{\# of surjective maps }f:\left[  m\right]  \rightarrow\left[
n\right]  \right)  =\sum\limits_{k=0}^{n}\left(  -1\right)  ^{k}\dbinom{n}%
{k}\left(  n-k\right)  ^{m}.
\]

\end{theorem}

This is the simplest expression for this number. It has no product formula
(unlike the \# of injective maps $f:\left[  m\right]  \rightarrow\left[
n\right]  $, which is $n\left(  n-1\right)  \left(  n-2\right)  \cdots\left(
n-m+1\right)  $).

Before we move on to the next example, let us draw a few consequences from
Theorem \ref{thm.pie.count-sur}:

\begin{corollary}
\label{cor.pie.count-sur.cors}Let $n\in\mathbb{N}$. Then: \medskip

\textbf{(a)} We have $\sum\limits_{k=0}^{n}\left(  -1\right)  ^{k}\dbinom
{n}{k}\left(  n-k\right)  ^{m}=0$ for any $m\in\mathbb{N}$ satisfying $m<n$.
\medskip

\textbf{(b)} We have $\sum\limits_{k=0}^{n}\left(  -1\right)  ^{k}\dbinom
{n}{k}\left(  n-k\right)  ^{n}=n!$. \medskip

\textbf{(c)} We have $\sum\limits_{k=0}^{n}\left(  -1\right)  ^{k}\dbinom
{n}{k}\left(  n-k\right)  ^{m}\geq0$ for each $m\in\mathbb{N}$. \medskip

\textbf{(d)} We have $n!\mid\sum\limits_{k=0}^{n}\left(  -1\right)
^{k}\dbinom{n}{k}\left(  n-k\right)  ^{m}$ for each $m\in\mathbb{N}$.
\end{corollary}

\begin{proof}
[Proof of Corollary \ref{cor.pie.count-sur.cors} (sketched).]\textbf{(a)} Let
$m\in\mathbb{N}$ satisfy $m<n$. Theorem \ref{thm.pie.count-sur} shows that the
sum $\sum\limits_{k=0}^{n}\left(  -1\right)  ^{k}\dbinom{n}{k}\left(
n-k\right)  ^{m}$ equals the \# of surjective maps $f:\left[  m\right]
\rightarrow\left[  n\right]  $. However, there are no such maps (by the
Pigeonhole Principle for Surjections)\footnote{Indeed, the Pigeonhole
Principle for Surjections shows that there are no surjective maps from a
smaller set to a larger set. Thus, in particular, there are no surjective maps
from $\left[  m\right]  $ to $\left[  n\right]  $ (since $m<n$).}; hence, this
\# is $0$. Therefore, the sum is $0$. In other words, $\sum\limits_{k=0}%
^{n}\left(  -1\right)  ^{k}\dbinom{n}{k}\left(  n-k\right)  ^{m}=0$. This
proves Corollary \ref{cor.pie.count-sur.cors} \textbf{(a)}. \medskip

\textbf{(b)} Theorem \ref{thm.pie.count-sur} (applied to $m=n$) shows that the
sum $\sum\limits_{k=0}^{n}\left(  -1\right)  ^{k}\dbinom{n}{k}\left(
n-k\right)  ^{n}$ equals the \# of surjective maps $f:\left[  n\right]
\rightarrow\left[  n\right]  $. However, these maps are precisely the
permutations of $\left[  n\right]  $ (by the Pigeonhole Principle for
Surjections)\footnote{Indeed, the Pigeonhole Principle for Surjections shows
that any surjective map between two finite sets of the same size is bijective.
Thus, any surjective map $f:\left[  n\right]  \rightarrow\left[  n\right]  $
is bijective, and hence is a permutation of $\left[  n\right]  $. The converse
is true as well (i.e., any permutation of $\left[  n\right]  $ is a surjective
map $f:\left[  n\right]  \rightarrow\left[  n\right]  $). Thus, the surjective
maps $f:\left[  n\right]  \rightarrow\left[  n\right]  $ are precisely the
permutations of $\left[  n\right]  $.}; therefore, this \# is $n!$. Therefore,
this sum is $n!$. This proves Corollary \ref{cor.pie.count-sur.cors}
\textbf{(b)}. \medskip

\textbf{(c)} Let $m\in\mathbb{N}$. Theorem \ref{thm.pie.count-sur} shows that
the sum $\sum\limits_{k=0}^{n}\left(  -1\right)  ^{k}\dbinom{n}{k}\left(
n-k\right)  ^{m}$ equals the \# of surjective maps $f:\left[  m\right]
\rightarrow\left[  n\right]  $. Hence, this sum is $\geq0$. This proves
Corollary \ref{cor.pie.count-sur.cors} \textbf{(c)}. \medskip

\textbf{(d)} Let $m\in\mathbb{N}$. Theorem \ref{thm.pie.count-sur} shows that
the sum $\sum\limits_{k=0}^{n}\left(  -1\right)  ^{k}\dbinom{n}{k}\left(
n-k\right)  ^{m}$ equals the \# of surjective maps $f:\left[  m\right]
\rightarrow\left[  n\right]  $. Let
\[
U:=\left\{  \text{surjective maps }f:\left[  m\right]  \rightarrow\left[
n\right]  \right\}  .
\]
Thus, this sum equals $\left\vert U\right\vert $. Hence, it remains to show
that $n!\mid\left\vert U\right\vert $.

The argument that follows will use the language of group actions (although it
could be restated combinatorially)\footnote{For a refresher on group actions,
see (e.g.) \cite[\S 3.2]{Quinlan21} or \cite[\S 9.11--\S 9.15]{Loehr-BC} or
\cite[\S 6.7--\S 6.9]{Artin} or \cite[Fall 2018, Weeks 8--10]%
{Armstrong-561fa18sp19} or \cite[\S 6.1]{Aigner07}. When $G$ is a group, we
use the word \textquotedblleft$G$\emph{-set}\textquotedblright\ to mean a set
on which $G$ acts.}. If $\sigma\in S_{n}$ is any permutation and $f:\left[
m\right]  \rightarrow\left[  n\right]  $ is a surjective map, then the
composition $\sigma\circ f:\left[  m\right]  \rightarrow\left[  n\right]  $ is
again a surjective map. In other words, any $\sigma\in S_{n}$ and $f\in U$
satisfy $\sigma\circ f\in U$. Hence, the symmetric group $S_{n}$ acts on the
set $U$ by the rule\footnote{This is called \textquotedblleft acting by
post-composition\textquotedblright\ (since $\sigma\circ f$ is obtained from
$f$ by composing with $\sigma$ \textquotedblleft after\textquotedblright%
\ $f$).}%
\[
\sigma\cdot f=\sigma\circ f\ \ \ \ \ \ \ \ \ \ \text{for all }\sigma\in
S_{n}\text{ and }f\in U.
\]
(This is a well-defined group action, since composition of maps is associative.)

This turns $U$ into an $S_{n}$-set. It is not hard to see that this action of
$S_{n}$ on $U$ is free -- i.e., the stabilizer of each $f\in U$ is the trivial
subgroup $\left\{  \operatorname*{id}\right\}  $ of $S_{n}$%
\ \ \ \ \footnote{\textit{Proof.} Let $f\in U$. We must prove that the
stabilizer of $f$ is $\left\{  \operatorname*{id}\right\}  $.
\par
Let $\sigma$ belong to the stabilizer of $f$. Thus, $\sigma\circ f=f$. Now,
let $j\in\left[  n\right]  $. Recall that $f\in U$; hence, $f$ is surjective.
Thus, there exists some $i\in\left[  m\right]  $ such that $j=f\left(
i\right)  $. Consider this $i$. Now, from $j=f\left(  i\right)  $, we obtain
$\sigma\left(  j\right)  =\sigma\left(  f\left(  i\right)  \right)
=\underbrace{\left(  \sigma\circ f\right)  }_{=f}\left(  i\right)  =f\left(
i\right)  =j$. Now, forget that we fixed $j$. We thus have shown that
$\sigma\left(  j\right)  =j$ for each $j\in\left[  n\right]  $. In other
words, $\sigma=\operatorname*{id}$.
\par
Forget that we fixed $\sigma$. We thus have shown that any $\sigma$ in the
stabilizer of $f$ satisfies $\sigma=\operatorname*{id}$. In other words, the
stabilizer of $f$ is a subset of $\left\{  \operatorname*{id}\right\}  $.
Therefore, the stabilizer of $f$ must be $\left\{  \operatorname*{id}\right\}
$ (since $\operatorname*{id}$ is clearly in the stabilizer of $f$).}. Thus,
the Orbit-Stabilizer Theorem (see, e.g., \cite[Theorem 3.2.5]{Quinlan21})
shows that each orbit of this action has size%
\[
\left[  S_{n}:\left\{  \operatorname*{id}\right\}  \right]
=\underbrace{\left\vert S_{n}\right\vert }_{=n!}/\underbrace{\left\vert
\left\{  \operatorname*{id}\right\}  \right\vert }_{=1}=n!/1=n!.
\]
However, the orbits of this action form a partition of $U$ (since each element
of $U$ belongs to exactly one orbit). Thus, we have%
\[
\left\vert U\right\vert =\left(  \text{the sum of the sizes of the
orbits}\right)  =\left(  \text{\# of orbits}\right)  \cdot n!
\]
(since each orbit has size $n!$). This entails $n!\mid\left\vert U\right\vert
$, which is exactly what we wanted to show. This proves Corollary
\ref{cor.pie.count-sur.cors} \textbf{(d)}.
\end{proof}

We note that parts \textbf{(a)} and \textbf{(b)} of Corollary
\ref{cor.pie.count-sur.cors} can also be proved algebraically (see, e.g.,
\cite[Exercise 5.4.2 \textbf{(d)}]{20f} for an algebraic generalization of
Corollary \ref{cor.pie.count-sur.cors} \textbf{(a)}); but this is harder for
Corollary \ref{cor.pie.count-sur.cors} \textbf{(d)} and (to my knowledge)
impossible for Corollary \ref{cor.pie.count-sur.cors} \textbf{(c)}.

\bigskip

\noindent\textbf{Example 2.} (See \cite[\S 2.9.5]{19fco} for details.) We will
use the following definition:

\begin{definition}
\label{def.pie.dera}A \emph{derangement} of a set $X$ means a permutation of
$X$ that has no fixed points.
\end{definition}

Now, let $n\in\mathbb{N}$. How many derangements does $\left[  n\right]  $ have?

Before answering this question, we establish notations and a few examples.

Let $D_{n}$ be the \# of derangements of $\left[  n\right]  $. For example:

\begin{itemize}
\item The identity permutation $\operatorname*{id}\in S_{0}$ is a derangement,
since it has no fixed points (since $\left[  0\right]  =\varnothing$ has no
elements to begin with). Thus, $D_{0}=1$.

\item The identity permutation $\operatorname*{id}\in S_{1}$ is \textbf{not} a
derangement, since $\operatorname*{id}\left(  1\right)  =1$. Thus, $D_{1}=0$.

\item In the symmetric group $S_{2}$, the identity is \textbf{not} a
derangement, but the transposition $s_{1}=t_{1,2}$ is one. Thus, $D_{2}=1$.

\item In the symmetric group $S_{3}$, the derangements are the $3$-cycles
$\operatorname*{cyc}\nolimits_{1,2,3}$ and $\operatorname*{cyc}%
\nolimits_{1,3,2}$. Thus, $D_{3}=2$.
\end{itemize}

Here is a table of early values of $D_{n}$:%
\[%
\begin{tabular}
[c]{|c||c|c|c|c|c|c|c|c|c|c|c|}\hline
$n$ & $0$ & $1$ & $2$ & $3$ & $4$ & $5$ & $6$ & $7$ & $8$ & $9$ & $10$\\\hline
$D_{n}$ & $1$ & $0$ & $1$ & $2$ & $9$ & $44$ & $265$ & $1854$ & $14\ 833$ &
$133\ 496$ & $1\ 334\ 961$\\\hline
\end{tabular}
\ \ \ \ \ \ \ .
\]

\textbf{Note:} The number $D_{n}$ is also called the \emph{subfactorial} of
$n$ and is sometimes denoted by $!n$. Be careful with that notation: what is
$2!2$ ?

\bigskip

Let us now try to compute $D_{n}$ in general. Fix $n\in\mathbb{N}$, and set
$U:=S_{n}=\left\{  \text{permutations of }\left[  n\right]  \right\}  $. We
impose $n$ rules $1,2,\ldots,n$ on a permutation $\sigma\in U$, with rule $i$
being \textquotedblleft thou shalt leave the element $i$
fixed\textquotedblright\ (in other words, rule $i$ requires a permutation
$\sigma\in U$ to satisfy $\sigma\left(  i\right)  =i$). Now,%
\begin{align*}
D_{n}  &  =\left(  \text{\# of derangements of }\left[  n\right]  \right) \\
&  =\left(  \text{\# of permutations }u\in U\text{ that violate all }n\text{
rules }1,2,\ldots,n\right) \\
&  =\sum_{I\subseteq\left[  n\right]  }\left(  -1\right)  ^{\left\vert
I\right\vert }\left(  \text{\# of permutations }u\in U\text{ that satisfy all
rules in }I\right)
\end{align*}
(by the \textquotedblleft rule-breaking\textquotedblright\ interpretation of
Theorem \ref{thm.pie.1}).

Now, let $I$ be a subset of $\left[  n\right]  $. What is the \# of
permutations $u\in U$ that satisfy all rules in $I$ ? These permutations are
the permutations of $\left[  n\right]  $ that leave each $i\in I$ fixed (but
can do whatever they want with the remaining elements of $\left[  n\right]
$). Clearly, there are $\left(  n-\left\vert I\right\vert \right)  !$ such
permutations, since they are essentially the permutations of the $\left(
n-\left\vert I\right\vert \right)  $-element set $\left[  n\right]  \setminus
I$ (with the elements of $I$ tacked on as fixed points). (See \cite[Corollary
2.9.16]{19fco} for the technical details of this intuitively obvious argument.)

Forget that we fixed $I$. We thus have shown that%
\begin{align}
&  \left(  \text{\# of permutations }u\in U\text{ that satisfy all rules in
}I\right) \nonumber\\
&  =\left(  n-\left\vert I\right\vert \right)  ! \label{eq.exa.pie1.2.3}%
\end{align}
for any subset $I$ of $\left[  n\right]  $.

Thus, the above computation becomes%
\begin{align*}
D_{n}  &  =\sum_{I\subseteq\left[  n\right]  }\left(  -1\right)  ^{\left\vert
I\right\vert }\underbrace{\left(  \text{\# of permutations }u\in U\text{ that
satisfy all rules in }I\right)  }_{\substack{=\left(  n-\left\vert
I\right\vert \right)  !\\\text{(by (\ref{eq.exa.pie1.2.3}))}}}\\
&  =\underbrace{\sum_{I\subseteq\left[  n\right]  }}_{=\sum_{k=0}^{n}%
\ \ \sum_{\substack{I\subseteq\left[  n\right]  ;\\\left\vert I\right\vert
=k}}}\left(  -1\right)  ^{\left\vert I\right\vert }\left(  n-\left\vert
I\right\vert \right)  !=\sum_{k=0}^{n}\ \ \sum_{\substack{I\subseteq\left[
n\right]  ;\\\left\vert I\right\vert =k}}\underbrace{\left(  -1\right)
^{\left\vert I\right\vert }\left(  n-\left\vert I\right\vert \right)
!}_{\substack{=\left(  -1\right)  ^{k}\left(  n-k\right)  !\\\text{(since
}\left\vert I\right\vert =k\text{)}}}\\
&  =\sum_{k=0}^{n}\ \ \underbrace{\sum_{\substack{I\subseteq\left[  n\right]
;\\\left\vert I\right\vert =k}}\left(  -1\right)  ^{k}\left(  n-k\right)
!}_{=\dbinom{n}{k}\left(  -1\right)  ^{k}\left(  n-k\right)  !}=\sum_{k=0}%
^{n}\dbinom{n}{k}\left(  -1\right)  ^{k}\left(  n-k\right)  !\\
&  =\sum_{k=0}^{n}\left(  -1\right)  ^{k}\underbrace{\dbinom{n}{k}\left(
n-k\right)  !}_{\substack{=\dfrac{n!}{k!}\\\text{(by (\ref{eq.binom.fac-form}%
))}}}=\sum_{k=0}^{n}\left(  -1\right)  ^{k}\dfrac{n!}{k!}=n!\cdot\sum
_{k=0}^{n}\dfrac{\left(  -1\right)  ^{k}}{k!}.
\end{align*}

Let us summarize our results as a theorem:

\begin{theorem}
\label{thm.pie.count-der}Let $n\in\mathbb{N}$. Then, the \# of derangements of
$\left[  n\right]  $ is
\[
D_{n}=\sum_{k=0}^{n}\left(  -1\right)  ^{k}\dbinom{n}{k}\left(  n-k\right)
!=n!\cdot\sum_{k=0}^{n}\dfrac{\left(  -1\right)  ^{k}}{k!}.
\]

\end{theorem}

\begin{remark}
The sum on the right hand side is a partial sum of the well-known infinite
series $\sum\limits_{k=0}^{\infty}\dfrac{\left(  -1\right)  ^{k}}{k!}=e^{-1}$
(where $e=2.718\ldots$ is Euler's number). This is quite helpful in
approximating $D_{n}$; indeed, it is easy to see (using some simple estimates)
that%
\[
D_{n}=\operatorname*{round}\dfrac{n!}{e}\ \ \ \ \ \ \ \ \ \ \text{for each
}n>0,
\]
where $\operatorname*{round}x$ means the result of rounding a real number $x$
to the nearest integer (fortunately, since $e$ is irrational, we never get a tie).
\end{remark}

See Exercise \ref{exe.pie.dera.recs} and \cite[Chapter 1]{Wildon19} for more
about derangements.

\bigskip

\noindent\textbf{Example 3.} Like many other things in these lectures, the
following elementary number-theoretical result is due to Euler:

\begin{theorem}
\label{thm.pie.euler-tot}Let $c$ be a positive integer with prime
factorization $p_{1}^{a_{1}}p_{2}^{a_{2}}\cdots p_{n}^{a_{n}}$, where
$p_{1},p_{2},\ldots,p_{n}$ are distinct primes, and where $a_{1},a_{2}%
,\ldots,a_{n}$ are positive integers. Then,%
\[
\left(  \text{\# of all }u\in\left[  c\right]  \text{ that are coprime to
}c\right)  =c\cdot\prod_{i=1}^{n}\left(  1-\dfrac{1}{p_{i}}\right)
=\prod_{i=1}^{n}\left(  p_{i}^{a_{i}}-p_{i}^{a_{i}-1}\right)  .
\]

\end{theorem}

Note that the \# of all $u\in\left[  c\right]  $ that are coprime to $c$ is
usually denoted by $\phi\left(  c\right)  $ in number theory, and the map
$\phi:\left\{  1,2,3,\ldots\right\}  \rightarrow\mathbb{N}$ that sends each
$c$ to $\phi\left(  c\right)  $ is called \emph{Euler's totient function}.

Theorem \ref{thm.pie.euler-tot} can be proved in many ways (and a proof can be
found in almost any text on elementary number theory). Probably the most
transparent proof relies on the PIE:

\begin{proof}
[Proof of Theorem \ref{thm.pie.euler-tot} (sketched).](See \cite[proof of
Theorem 2.9.19]{19fco} for the details of this argument.) Let $U=\left[
c\right]  $. A number $u\in U$ is coprime to $c$ if and only if it is not
divisible by any of the prime factors $p_{1},p_{2},\ldots,p_{n}$ of $c$.
Again, this means that $u$ breaks all $n$ rules $1,2,\ldots,n$, where rule $i$
says \textquotedblleft thou shalt be divisible by $p_{i}$\textquotedblright.
Thus, by the \textquotedblleft rule-breaking\textquotedblright\ interpretation
of Theorem \ref{thm.pie.1}, we obtain%
\begin{align*}
&  \left(  \text{\# of all }u\in\left[  c\right]  \text{ that are coprime to
}c\right) \\
&  =\sum_{I\subseteq\left[  n\right]  }\left(  -1\right)  ^{\left\vert
I\right\vert }\underbrace{\left(  \text{\# of all }u\in\left[  c\right]
\text{ that are divisible by all }p_{i}\text{ with }i\in I\right)
}_{\substack{=\dfrac{c}{\prod_{i\in I}p_{i}}\\\text{(this is not hard to
prove)}}}\\
&  =\sum_{I\subseteq\left[  n\right]  }\left(  -1\right)  ^{\left\vert
I\right\vert }\dfrac{c}{\prod_{i\in I}p_{i}}=c\cdot\underbrace{\sum
_{I\subseteq\left[  n\right]  }\left(  -1\right)  ^{\left\vert I\right\vert
}\prod_{i\in I}\dfrac{1}{p_{i}}}_{\substack{=\prod_{i=1}^{n}\left(
1-\dfrac{1}{p_{i}}\right)  \\\text{(an easy consequence of
(\ref{eq.lem.prodrule.sum-ai-plus-bi.eq}))}}}\\
&  =\underbrace{c}_{\substack{=p_{1}^{a_{1}}p_{2}^{a_{2}}\cdots p_{n}^{a_{n}%
}\\=\prod_{i=1}^{n}p_{i}^{a_{i}}}}\cdot\prod_{i=1}^{n}\left(  1-\dfrac
{1}{p_{i}}\right)  =\left(  \prod_{i=1}^{n}p_{i}^{a_{i}}\right)  \cdot
\prod_{i=1}^{n}\left(  1-\dfrac{1}{p_{i}}\right) \\
&  =\prod_{i=1}^{n}\underbrace{\left(  p_{i}^{a_{i}}\left(  1-\dfrac{1}{p_{i}%
}\right)  \right)  }_{=p_{i}^{a_{i}}-p_{i}^{a_{i}-1}}=\prod_{i=1}^{n}\left(
p_{i}^{a_{i}}-p_{i}^{a_{i}-1}\right)  .
\end{align*}
This proves Theorem \ref{thm.pie.euler-tot}.
\end{proof}

\bigskip

\noindent\textbf{Example 4.} (This one is taken from \cite[Theorem
2.3.3]{Sagan19}.) Recall Theorem \ref{thm.pars.odd-dist-equal}, which states
that each $n\in\mathbb{N}$ satisfies%
\[
p_{\operatorname*{odd}}\left(  n\right)  =p_{\operatorname*{dist}}\left(
n\right)  ,
\]
where%
\begin{align*}
p_{\operatorname*{odd}}\left(  n\right)   &  :=\left(  \text{\# of partitions
of }n\text{ into odd parts}\right)  \ \ \ \ \ \ \ \ \ \ \text{and}\\
p_{\operatorname*{dist}}\left(  n\right)   &  :=\left(  \text{\# of partitions
of }n\text{ into distinct parts}\right)  .
\end{align*}

We have already proved this twice, but let us prove it again.

\begin{proof}
[Third proof of Theorem \ref{thm.pars.odd-dist-equal} (sketched).]Let
$n\in\mathbb{N}$. We set $U:=\left\{  \text{partitions of }n\right\}  $.

In this proof, the word \textquotedblleft partition\textquotedblright\ will
always mean \textquotedblleft partition of $n$\textquotedblright. Thus, a
partition cannot contain any of the entries $n+1,n+2,n+3,\ldots$ (because any
of these entries would cause the partition to have size $>n$).

We want to frame the partitions of $n$ into distinct parts as rule-breakers.
We observe that%
\begin{align*}
&  \left\{  \text{partitions of }n\text{ into distinct parts}\right\} \\
&  =\left\{  \text{partitions that contain none of the entries }%
1,2,3,\ldots\text{ twice}\right\} \\
&  \ \ \ \ \ \ \ \ \ \ \ \ \ \ \ \ \ \ \ \ \left(
\begin{array}
[c]{c}%
\text{here and in the following, the word \textquotedblleft
twice\textquotedblright}\\
\text{means \textquotedblleft at least twice\textquotedblright}%
\end{array}
\right) \\
&  =\left\{  \text{partitions that contain none of the entries }%
1,2,\ldots,n\text{ twice}\right\}
\end{align*}
(since a partition cannot contain any of the entries $n+1,n+2,n+3,\ldots$). In
other words, the partitions of $n$ into distinct parts are precisely the
partitions that break all rules $1,2,\ldots,n$, where rule $i$ says
\textquotedblleft thou shalt contain the entry $i$ twice\textquotedblright%
\footnote{Once again, \textquotedblleft twice\textquotedblright\ means
\textquotedblleft at least twice\textquotedblright.}.

Thus, applying the PIE (specifically, the \textquotedblleft
rule-breaking\textquotedblright\ interpretation of Theorem \ref{thm.pie.1}),
we obtain%
\begin{align}
&  p_{\operatorname*{dist}}\left(  n\right) \nonumber\\
&  =\sum_{I\subseteq\left[  n\right]  }\left(  -1\right)  ^{\left\vert
I\right\vert }\left(  \text{\# of partitions that satisfy all rules in
}I\right) \nonumber\\
&  =\sum_{I\subseteq\left[  n\right]  }\left(  -1\right)  ^{\left\vert
I\right\vert }\left(  \text{\# of partitions that contain each of the entries
}i\in I\text{ twice}\right)  . \label{pf.thm.pars.odd-dist-equal.3rd.3}%
\end{align}

We can play the same game with $p_{\operatorname*{odd}}\left(  n\right)  $.
This time, rule $i$ says \textquotedblleft thou shalt contain the entry
$2i$\textquotedblright. Thus, again applying the PIE, we obtain%
\begin{equation}
p_{\operatorname*{odd}}\left(  n\right)  =\sum_{I\subseteq\left[  n\right]
}\left(  -1\right)  ^{\left\vert I\right\vert }\left(  \text{\# of partitions
that contain the entry }2i\text{ for each }i\in I\right)  .
\label{pf.thm.pars.odd-dist-equal.3rd.4}%
\end{equation}

Now, comparing these two equalities, we see that in order to prove that
$p_{\operatorname*{dist}}\left(  n\right)  =p_{\operatorname*{odd}}\left(
n\right)  $, it will suffice to show that%
\begin{align*}
&  \left(  \text{\# of partitions that contain each of the entries }i\in
I\text{ twice}\right) \\
&  =\left(  \text{\# of partitions that contain the entry }2i\text{ for each
}i\in I\right)
\end{align*}
for any subset $I$ of $\left[  n\right]  $.

So let $I$ be a subset of $\left[  n\right]  $. We are looking for a
bijection
\begin{align*}
&  \text{from }\left\{  \text{partitions that contain each of the entries
}i\in I\text{ twice}\right\} \\
&  \text{to }\left\{  \text{partitions that contain the entry }2i\text{ for
each }i\in I\right\}  .
\end{align*}
Such a bijection can be obtained as follows: For each $i\in I$ (from highest
to lowest\footnote{Actually, a bit of thought reveals that the order in which
we go through the $i\in I$ does not affect the result; this becomes
particularly clear if we identify each partition with the multiset of its
entries. Thus, me saying \textquotedblleft from highest to
lowest\textquotedblright\ is unnecessary.}), we remove two copies of $i$ from
the partition, and insert a $2i$ into the partition in their stead. For
example, if $I=\left\{  2,4,5\right\}  $ and $n=33$, then our bijection sends
the partition
\[
\left(  5,5,4,4,3,3,2,2,2,2,1\right)  \text{ to }\left(
10,8,4,3,3,2,2,1\right)  .
\]
(Note that the $4$ in the resulting partition is not one of the original two
$4$s, but rather a new $4$ that was inserted when we removed two copies of
$2$. On the other hand, the two $2$s in the resulting partition are inherited
from the original partition, because (unlike the bijection $A$ in our Second
proof of Theorem \ref{thm.pars.odd-dist-equal} above) our bijection only
removes two copies of each $i\in I$.)

It is easy to see that this purported bijection really is a
bijection\footnote{Its inverse, of course, does what you would expect: For
each $i\in I$, we remove a $2i$ from the partition, and insert two copies of
$i$ in its stead.}. Thus, we have found our bijection. The bijection principle
therefore yields%
\begin{align*}
&  \left(  \text{\# of partitions that contain each of the entries }i\in
I\text{ twice}\right) \\
&  =\left(  \text{\# of partitions that contain the entry }2i\text{ for each
}i\in I\right)  .
\end{align*}

We have proved this equality for all $I\subseteq\left[  n\right]  $. Hence,
the right hand sides of the equalities (\ref{pf.thm.pars.odd-dist-equal.3rd.3}%
) and (\ref{pf.thm.pars.odd-dist-equal.3rd.4}) are equal. Thus, their left
hand sides are equal as well. In other words, $p_{\operatorname*{dist}}\left(
n\right)  =p_{\operatorname*{odd}}\left(  n\right)  $. This proves Theorem
\ref{thm.pars.odd-dist-equal} again.
\end{proof}

\begin{remark}
It is worth contrasting the above four examples (in which we applied the PIE
to solve a counting problem, obtaining an alternating sum as a result) with
our arguments in Section \ref{sec.sign.intro} (in which we computed
alternating sums using sign-reversing involutions). Sign-reversing involutions
help turn alternating sums into combinatorial problems, while the PIE moves us
in the opposite direction. The two techniques are thus, in some way, inverse
to each other. This will become less mysterious once we prove the PIE itself
using a sign-reversing involution. The PIE can also be used backwards, to turn
an alternating sum into a counting problem, which is how we proved Corollary
\ref{cor.pie.count-sur.cors} above.
\end{remark}

\subsubsection{\label{subsec.sign.pie.wt}The weighted version}

The main rule of algebra is to never turn down nature's gifts. The PIE (in the
shape of Theorem \ref{thm.pie.1}) can be generalized with zero effort, so let
us do it (\cite[Theorem 7.8.9]{20f}):

\begin{theorem}
[weighted version of the PIE]\label{thm.pie.2}Let $n\in\mathbb{N}$, and let
$U$ be a finite set. Let $A_{1},A_{2},\ldots,A_{n}$ be $n$ subsets of $U$. Let
$A$ be any additive abelian group (such as $\mathbb{R}$, or any vector space,
or any ring). Let $w:U\rightarrow A$ be any map (i.e., let $w\left(  u\right)
$ be an element of $A$ for each $u\in U$). Then,%
\begin{equation}
\sum_{\substack{u\in U;\\u\notin A_{i}\text{ for all }i\in\left[  n\right]
}}w\left(  u\right)  =\sum_{I\subseteq\left[  n\right]  }\left(  -1\right)
^{\left\vert I\right\vert }\sum_{\substack{u\in U;\\u\in A_{i}\text{ for all
}i\in I}}w\left(  u\right)  . \label{eq.thm.pie.2.eq}%
\end{equation}

\end{theorem}

We can think of each value $w\left(  u\right)  $ in Theorem \ref{thm.pie.2} as
a kind of \textquotedblleft weight\textquotedblright\ of the respective
element $u$. Thus, the left hand side of the equality (\ref{eq.thm.pie.2.eq})
is the total weight of all \textquotedblleft rule-breaking\textquotedblright%
\ $u\in U$ (that is, of all $u\in U$ that satisfy $\left(  u\notin A_{i}\text{
for all }i\in\left[  n\right]  \right)  $), whereas the inner sum
$\sum_{\substack{u\in U;\\u\in A_{i}\text{ for all }i\in I}}w\left(  u\right)
$ on the right hand side is the total weight of all $u\in U$ that satisfy
$\left(  u\in A_{i}\text{ for all }i\in I\right)  $. This is why we call
Theorem \ref{thm.pie.2} the \emph{weighted version of the PIE} (or just the
\emph{weighted PIE}).

Theorem \ref{thm.pie.1} can be obtained from Theorem \ref{thm.pie.2} by taking
$w$ to be constantly $1$ (that is, by setting $w\left(  u\right)  =1$ for each
$u\in U$). Indeed, a sum of a bunch of $1$s equals the \# of $1$s being
summed, so that sums generalize cardinalities.

With Theorem \ref{thm.pie.2}, we can take any of our above four examples from
Subsection \ref{subsec.sign.pie.exas}, and introduce weights into them --
i.e., instead of asking for a number, we sum certain \textquotedblleft
weights\textquotedblright. Little question (cf. the homework): What do you get?

But we have no time for this now. Another generalization is calling!

\subsubsection{\label{subsec.sign.pie.moeb-bool}Boolean M\"{o}bius inversion}

This generalization (we will soon see why it is a generalization) is our third
\textquotedblleft Principle of Inclusion and Exclusion\textquotedblright, but
is probably best referred to as the \emph{Boolean M\"{o}bius inversion
formula} (or the \emph{M\"{o}bius inversion formula for the Boolean lattice}):

\begin{theorem}
[Boolean M\"{o}bius inversion]\label{thm.pie.moeb}Let $S$ be a finite set. Let
$A$ be any additive abelian group.

For each subset $I$ of $S$, let $a_{I}$ and $b_{I}$ be two elements of $A$.

Assume that%
\begin{equation}
b_{I}=\sum_{J\subseteq I}a_{J}\ \ \ \ \ \ \ \ \ \ \text{for all }I\subseteq S.
\label{eq.thm.pie.moeb.ass}%
\end{equation}

Then, we also have%
\begin{equation}
a_{I}=\sum_{J\subseteq I}\left(  -1\right)  ^{\left\vert I\setminus
J\right\vert }b_{J}\ \ \ \ \ \ \ \ \ \ \text{for all }I\subseteq S.
\label{eq.thm.pie.moeb.claim}%
\end{equation}

\end{theorem}

\begin{example}
Let $S=\left[  2\right]  =\left\{  1,2\right\}  $. Then, the assumptions of
Theorem \ref{thm.pie.moeb} state that%
\begin{align*}
b_{\varnothing}  &  =a_{\varnothing};\\
b_{\left\{  1\right\}  }  &  =a_{\varnothing}+a_{\left\{  1\right\}  };\\
b_{\left\{  2\right\}  }  &  =a_{\varnothing}+a_{\left\{  2\right\}  };\\
b_{\left\{  1,2\right\}  }  &  =a_{\varnothing}+a_{\left\{  1\right\}
}+a_{\left\{  2\right\}  }+a_{\left\{  1,2\right\}  }.
\end{align*}
The claim of Theorem \ref{thm.pie.moeb} then states that%
\begin{align*}
a_{\varnothing}  &  =b_{\varnothing};\\
a_{\left\{  1\right\}  }  &  =-b_{\varnothing}+b_{\left\{  1\right\}  };\\
a_{\left\{  2\right\}  }  &  =-b_{\varnothing}+b_{\left\{  2\right\}  };\\
a_{\left\{  1,2\right\}  }  &  =b_{\varnothing}-b_{\left\{  1\right\}
}-b_{\left\{  2\right\}  }+b_{\left\{  1,2\right\}  }.
\end{align*}
These four equalities can be verified easily. For instance, let us check the
last of them:%
\begin{align*}
&  \underbrace{b_{\varnothing}}_{=a_{\varnothing}}-\underbrace{b_{\left\{
1\right\}  }}_{=a_{\varnothing}+a_{\left\{  1\right\}  }}%
-\underbrace{b_{\left\{  2\right\}  }}_{=a_{\varnothing}+a_{\left\{
2\right\}  }}+\underbrace{b_{\left\{  1,2\right\}  }}_{=a_{\varnothing
}+a_{\left\{  1\right\}  }+a_{\left\{  2\right\}  }+a_{\left\{  1,2\right\}
}}\\
&  =a_{\varnothing}-\left(  a_{\varnothing}+a_{\left\{  1\right\}  }\right)
-\left(  a_{\varnothing}+a_{\left\{  2\right\}  }\right)  +\left(
a_{\varnothing}+a_{\left\{  1\right\}  }+a_{\left\{  2\right\}  }+a_{\left\{
1,2\right\}  }\right) \\
&  =a_{\left\{  1,2\right\}  }.
\end{align*}

\end{example}

Before we prove Theorem \ref{thm.pie.moeb}, let us show that the weighted PIE
(Theorem \ref{thm.pie.2}) is a particular case of it:

\begin{proof}
[Proof of Theorem \ref{thm.pie.2} using Theorem \ref{thm.pie.moeb}.]Let
$S=\left[  n\right]  $. We note that the map
\begin{align*}
\left\{  \text{subsets of }S\right\}   &  \rightarrow\left\{  \text{subsets of
}S\right\}  ,\\
J  &  \mapsto S\setminus J
\end{align*}
is a bijection. (Indeed, this map is an involution, since each subset $J$ of
$S$ satisfies $S\setminus\left(  S\setminus J\right)  =J$.)

For each $u\in U$, define a subset $\operatorname*{Viol}u$ of $S$ by%
\[
\operatorname*{Viol}u:=\left\{  i\in S\ \mid\ u\notin A_{i}\right\}  .
\]
(In terms of the \textquotedblleft rule-breaking\textquotedblright%
\ interpretation, $\operatorname*{Viol}u$ is the set of all rules that $u$
violates.) Now, for each subset $I$ of $\left[  n\right]  $, we set%
\[
a_{I}:=\sum_{\substack{u\in U;\\\operatorname*{Viol}u=I}}w\left(  u\right)
\ \ \ \ \ \ \ \ \ \ \text{and}\ \ \ \ \ \ \ \ \ \ b_{I}:=\sum_{\substack{u\in
U;\\\operatorname*{Viol}u\subseteq I}}w\left(  u\right)  .
\]
Then, for each subset $I$ of $S$, we have%
\begin{align*}
b_{I}  &  =\sum_{\substack{u\in U;\\\operatorname*{Viol}u\subseteq I}}w\left(
u\right)  =\sum_{J\subseteq I}\ \ \underbrace{\sum_{\substack{u\in
U;\\\operatorname*{Viol}u=J}}w\left(  u\right)  }_{\substack{=a_{J}\\\text{(by
the definition of }a_{J}\text{)}}}\ \ \ \ \ \ \ \ \ \ \left(
\begin{array}
[c]{c}%
\text{here, we have split}\\
\text{the sum according to}\\
\text{the value of }\operatorname*{Viol}u
\end{array}
\right) \\
&  =\sum_{J\subseteq I}a_{J}.
\end{align*}
Thus, we can apply Theorem \ref{thm.pie.moeb}. This gives us%
\begin{equation}
a_{I}=\sum_{J\subseteq I}\left(  -1\right)  ^{\left\vert I\setminus
J\right\vert }b_{J}\ \ \ \ \ \ \ \ \ \ \text{for all }I\subseteq S.
\label{pf.thm.pie.2.4}%
\end{equation}

However, if $J$ is any subset of $S$, then the definition of $b_{J}$ yields%
\begin{equation}
b_{J}=\sum_{\substack{u\in U;\\\operatorname*{Viol}u\subseteq J}}w\left(
u\right)  =\sum_{\substack{u\in U;\\u\in A_{i}\text{ for all }i\in S\setminus
J}}w\left(  u\right)  \label{pf.thm.pie.2.5}%
\end{equation}
(because the elements $u\in U$ that satisfy $\operatorname*{Viol}u\subseteq J$
are precisely the elements $u\in U$ that satisfy $\left(  u\in A_{i}\text{ for
all }i\in S\setminus J\right)  $\ \ \ \ \footnote{\textit{Proof.} Let $u\in
U$. We must prove that the condition \textquotedblleft$\operatorname*{Viol}%
u\subseteq J$\textquotedblright\ is equivalent to \textquotedblleft$u\in
A_{i}$ for all $i\in S\setminus J$\textquotedblright.
\par
We have $\operatorname*{Viol}u=\left\{  i\in S\ \mid\ u\notin A_{i}\right\}  $
(by the definition of $\operatorname*{Viol}u$). Hence, we have the following
chain of equivalences:%
\begin{align*}
\left(  \operatorname*{Viol}u\subseteq J\right)  \  &  \Longleftrightarrow
\ \left(  \left\{  i\in S\ \mid\ u\notin A_{i}\right\}  \subseteq J\right) \\
&  \Longleftrightarrow\ \left(  \text{each }i\in S\text{ satisfying }u\notin
A_{i}\text{ belongs to }J\right) \\
&  \Longleftrightarrow\ \left(  \text{each }i\in S\text{ that does not belong
to }J\text{ must satisfy }u\in A_{i}\right) \\
&  \ \ \ \ \ \ \ \ \ \ \ \ \ \ \ \ \ \ \ \ \left(  \text{by contraposition}%
\right) \\
&  \Longleftrightarrow\ \left(  \text{each }i\in S\setminus J\text{ must
satisfy }u\in A_{i}\right) \\
&  \Longleftrightarrow\ \left(  u\in A_{i}\text{ for all }i\in S\setminus
J\right)  .
\end{align*}
Hence, the condition \textquotedblleft$\operatorname*{Viol}u\subseteq
J$\textquotedblright\ is equivalent to \textquotedblleft$u\in A_{i}$ for all
$i\in S\setminus J$\textquotedblright.}).

Now, applying (\ref{pf.thm.pie.2.4}) to $I=S$, we obtain%
\begin{align*}
a_{S}  &  =\sum_{J\subseteq S}\left(  -1\right)  ^{\left\vert S\setminus
J\right\vert }b_{J}=\sum_{J\subseteq S}\left(  -1\right)  ^{\left\vert
S\setminus J\right\vert }\sum_{\substack{u\in U;\\u\in A_{i}\text{ for all
}i\in S\setminus J}}w\left(  u\right)  \ \ \ \ \ \ \ \ \ \ \left(  \text{by
(\ref{pf.thm.pie.2.5})}\right) \\
&  =\sum_{I\subseteq S}\left(  -1\right)  ^{\left\vert I\right\vert }%
\sum_{\substack{u\in U;\\u\in A_{i}\text{ for all }i\in I}}w\left(  u\right)
\\
&  \ \ \ \ \ \ \ \ \ \ \ \ \ \ \ \ \ \ \ \ \left(
\begin{array}
[c]{c}%
\text{here, we have substituted }I\text{ for }S\setminus J\text{ in the sum,
since}\\
\text{the map }\left\{  \text{subsets of }S\right\}  \rightarrow\left\{
\text{subsets of }S\right\}  ,\ J\mapsto S\setminus J\\
\text{is a bijection}%
\end{array}
\right) \\
&  =\sum_{I\subseteq\left[  n\right]  }\left(  -1\right)  ^{\left\vert
I\right\vert }\sum_{\substack{u\in U;\\u\in A_{i}\text{ for all }i\in
I}}w\left(  u\right)  \ \ \ \ \ \ \ \ \ \ \left(  \text{since }S=\left[
n\right]  \right)  .
\end{align*}

On the other hand, the definition of $a_{S}$ yields%
\[
a_{S}=\sum_{\substack{u\in U;\\\operatorname*{Viol}u=S}}w\left(  u\right)
=\sum_{\substack{u\in U;\\u\notin A_{i}\text{ for all }i\in\left[  n\right]
}}w\left(  u\right)
\]
(since the elements $u\in U$ that satisfy $\operatorname*{Viol}u=S$ are
precisely the elements $u\in U$ that satisfy $\left(  u\notin A_{i}\text{ for
all }i\in\left[  n\right]  \right)  $\ \ \ \ \footnote{\textit{Proof.} Let
$u\in U$. We must prove that the condition \textquotedblleft%
$\operatorname*{Viol}u=S$\textquotedblright\ is equivalent to
\textquotedblleft$u\notin A_{i}$ for all $i\in\left[  n\right]  $%
\textquotedblright.
\par
We have $\operatorname*{Viol}u=\left\{  i\in S\ \mid\ u\notin A_{i}\right\}  $
(by the definition of $\operatorname*{Viol}u$). Hence, we have the following
chain of equivalences:%
\begin{align*}
\left(  \operatorname*{Viol}u=S\right)  \  &  \Longleftrightarrow\ \left(
\left\{  i\in S\ \mid\ u\notin A_{i}\right\}  =S\right) \\
&  \Longleftrightarrow\ \left(  \text{each }i\in S\text{ satisfies }u\notin
A_{i}\right) \\
&  \Longleftrightarrow\ \left(  u\notin A_{i}\text{ for all }i\in S\right) \\
&  \Longleftrightarrow\ \left(  u\notin A_{i}\text{ for all }i\in\left[
n\right]  \right)  \ \ \ \ \ \ \ \ \ \ \left(  \text{since }S=\left[
n\right]  \right)  .
\end{align*}
Hence, the condition \textquotedblleft$\operatorname*{Viol}u=S$%
\textquotedblright\ is equivalent to \textquotedblleft$u\notin A_{i}$ for all
$i\in\left[  n\right]  $\textquotedblright.}).

Comparing these two equalities, we obtain%
\[
\sum_{\substack{u\in U;\\u\notin A_{i}\text{ for all }i\in\left[  n\right]
}}w\left(  u\right)  =\sum_{I\subseteq\left[  n\right]  }\left(  -1\right)
^{\left\vert I\right\vert }\sum_{\substack{u\in U;\\u\in A_{i}\text{ for all
}i\in I}}w\left(  u\right)  .
\]
Thus, Theorem \ref{thm.pie.2} has been proved, assuming Theorem
\ref{thm.pie.moeb}.
\end{proof}

It remains to prove the latter:

\begin{proof}
[Proof of Theorem \ref{thm.pie.moeb}.]Fix a subset $Q$ of $S$. We shall prove
that%
\[
a_{Q}=\sum_{I\subseteq Q}\left(  -1\right)  ^{\left\vert Q\setminus
I\right\vert }b_{I}.
\]

We begin by rewriting the right hand side:%
\begin{align}
\sum_{I\subseteq Q}\left(  -1\right)  ^{\left\vert Q\setminus I\right\vert
}\underbrace{b_{I}}_{\substack{=\sum_{J\subseteq I}a_{J}\\\text{(by
(\ref{eq.thm.pie.moeb.ass}))}}}  &  =\sum_{I\subseteq Q}\left(  -1\right)
^{\left\vert Q\setminus I\right\vert }\underbrace{\sum_{J\subseteq I}a_{J}%
}_{=\sum_{P\subseteq I}a_{P}}=\sum_{I\subseteq Q}\left(  -1\right)
^{\left\vert Q\setminus I\right\vert }\sum_{P\subseteq I}a_{P}\nonumber\\
&  =\sum_{I\subseteq Q}\ \ \sum_{P\subseteq I}\left(  -1\right)  ^{\left\vert
Q\setminus I\right\vert }a_{P}. \label{pf.thm.pie.moeb.1}%
\end{align}

The two summation signs \textquotedblleft$\sum_{I\subseteq Q}\ \ \sum
_{P\subseteq I}$\textquotedblright\ on the right hand side of this equality
result in a sum over all pairs $\left(  I,P\right)  $ of subsets of $Q$
satisfying $P\subseteq I\subseteq Q$. The same result can be obtained by the
two summation signs \textquotedblleft$\sum_{P\subseteq Q}\ \ \sum
_{\substack{I\subseteq Q;\\P\subseteq I}}$\textquotedblright\ (indeed, the
only difference between \textquotedblleft$\sum_{I\subseteq Q}\ \ \sum
_{P\subseteq I}$\textquotedblright\ and \textquotedblleft$\sum_{P\subseteq
Q}\ \ \sum_{\substack{I\subseteq Q;\\P\subseteq I}}$\textquotedblright\ is the
order in which the two subsets $I$ and $P$ are chosen). Thus, we can replace
\textquotedblleft$\sum_{I\subseteq Q}\ \ \sum_{P\subseteq I}$%
\textquotedblright\ by \textquotedblleft$\sum_{P\subseteq Q}\ \ \sum
_{\substack{I\subseteq Q;\\P\subseteq I}}$\textquotedblright\ on the right
hand side of (\ref{pf.thm.pie.moeb.1}). Hence, (\ref{pf.thm.pie.moeb.1})
rewrites as follows:%
\begin{align}
\sum_{I\subseteq Q}\left(  -1\right)  ^{\left\vert Q\setminus I\right\vert
}b_{I}  &  =\sum_{P\subseteq Q}\ \ \sum_{\substack{I\subseteq Q;\\P\subseteq
I}}\left(  -1\right)  ^{\left\vert Q\setminus I\right\vert }a_{P}\nonumber\\
&  =\sum_{P\subseteq Q}\left(  \sum_{\substack{I\subseteq Q;\\P\subseteq
I}}\left(  -1\right)  ^{\left\vert Q\setminus I\right\vert }\right)  a_{P}.
\label{pf.thm.pie.moeb.2}%
\end{align}

We want to prove that this equals $a_{Q}$. Since the $a_{P}$'s are arbitrary
elements of an abelian group, the only way this can possibly be achieved is by
showing that the sum on the right hand side simplifies to $a_{Q}$
\textbf{formally} -- i.e., that the coefficient $\sum_{\substack{I\subseteq
Q;\\P\subseteq I}}\left(  -1\right)  ^{\left\vert Q\setminus I\right\vert }$
in front of $a_{P}$ is $0$ whenever $P\neq Q$, and is $1$ whenever $P=Q$.
Thus, we now set out to prove this. Using Definition \ref{def.iverson}, we can
restate this goal as follows: We want to prove that every subset $P$ of $Q$
satisfies%
\begin{equation}
\sum_{\substack{I\subseteq Q;\\P\subseteq I}}\left(  -1\right)  ^{\left\vert
Q\setminus I\right\vert }=\left[  P=Q\right]  . \label{pf.thm.pie.moeb.4}%
\end{equation}

We shall prove this soon (in Lemma \ref{lem.pie.two-sets-altsum} \textbf{(b)}
below). For now, let us explain how the proof of Theorem \ref{thm.pie.moeb}
can be completed if (\ref{pf.thm.pie.moeb.4}) is known to be true. Indeed,
(\ref{pf.thm.pie.moeb.2}) becomes%
\begin{align*}
\sum_{I\subseteq Q}\left(  -1\right)  ^{\left\vert Q\setminus I\right\vert
}b_{I}  &  =\sum_{P\subseteq Q}\underbrace{\left(  \sum_{\substack{I\subseteq
Q;\\P\subseteq I}}\left(  -1\right)  ^{\left\vert Q\setminus I\right\vert
}\right)  }_{\substack{=\left[  P=Q\right]  \\\text{(by
(\ref{pf.thm.pie.moeb.4}))}}}a_{P}=\sum_{P\subseteq Q}\left[  P=Q\right]
a_{P}\\
&  =\underbrace{\left[  Q=Q\right]  }_{\substack{=1\\\text{(since }%
Q=Q\text{)}}}a_{Q}+\sum_{\substack{P\subseteq Q;\\P\neq Q}}\underbrace{\left[
P=Q\right]  }_{\substack{=0\\\text{(since }P\neq Q\text{)}}}a_{P}\\
&  \ \ \ \ \ \ \ \ \ \ \ \ \ \ \ \ \ \ \ \ \left(
\begin{array}
[c]{c}%
\text{here, we have split off the addend for }P=Q\\
\text{from the sum (since }Q\subseteq Q\text{)}%
\end{array}
\right) \\
&  =a_{Q}+\underbrace{\sum_{\substack{P\subseteq Q;\\P\neq Q}}0a_{P}}%
_{=0}=a_{Q}.
\end{align*}
In other words, $a_{Q}=\sum_{I\subseteq Q}\left(  -1\right)  ^{\left\vert
Q\setminus I\right\vert }b_{I}$.

Forget that we fixed $Q$. We thus have shown that%
\[
a_{Q}=\sum_{I\subseteq Q}\left(  -1\right)  ^{\left\vert Q\setminus
I\right\vert }b_{I}\ \ \ \ \ \ \ \ \ \ \text{for all }Q\subseteq S.
\]
Renaming the indices $Q$ and $I$ as $I$ and $J$ in this statement, we obtain
the following:%
\[
a_{I}=\sum_{J\subseteq I}\left(  -1\right)  ^{\left\vert I\setminus
J\right\vert }b_{J}\ \ \ \ \ \ \ \ \ \ \text{for all }I\subseteq S.
\]
Thus, Theorem \ref{thm.pie.moeb} is proved (assuming that
(\ref{pf.thm.pie.moeb.4}) is known to be true).
\end{proof}

It now remains to prove (\ref{pf.thm.pie.moeb.4}). We shall do this as part of
the following lemma:

\begin{lemma}
\label{lem.pie.two-sets-altsum}Let $Q$ be a finite set. Let $P$ be a subset of
$Q$. Then: \medskip

\textbf{(a)} We have%
\begin{equation}
\sum_{\substack{I\subseteq Q;\\P\subseteq I}}\left(  -1\right)  ^{\left\vert
I\right\vert }=\left(  -1\right)  ^{\left\vert P\right\vert }\left[
P=Q\right]  . \label{eq.lem.pie.two-sets-altsum.a}%
\end{equation}

\textbf{(b)} We have%
\begin{equation}
\sum_{\substack{I\subseteq Q;\\P\subseteq I}}\left(  -1\right)  ^{\left\vert
Q\setminus I\right\vert }=\left[  P=Q\right]  .
\label{eq.lem.pie.two-sets-altsum.b}%
\end{equation}

\end{lemma}

As promised, Lemma \ref{lem.pie.two-sets-altsum} \textbf{(b)} (once proved)
will yield (\ref{pf.thm.pie.moeb.4}) and thus will complete our above proof of
Theorem \ref{thm.pie.moeb} (and, with it, the proofs of Theorem
\ref{thm.pie.2} and Theorem \ref{thm.pie.1}).

\begin{proof}
[Proof of Lemma \ref{lem.pie.two-sets-altsum}.]\textbf{(a)} There are many
ways to prove this (in particular, a simple one using the binomial theorem --
do you see it?); but staying true to the spirit of this chapter, we pick one
using a sign-reversing involution. (A variant of this proof can be found in
\cite[solution to Exercise 2.9.1]{19fco}.\footnote{Our sets $Q$ and $P$ are
called $S$ and $T$ in \cite[solution to Exercise 2.9.1]{19fco}.})

We must prove the equality (\ref{eq.lem.pie.two-sets-altsum.a}). If $P=Q$,
then this equality is easily seen to hold\footnote{\textit{Proof.} Assume that
$P=Q$. Then, the only subset $I$ of $Q$ that satisfies $P\subseteq I$ is the
set $Q$ itself (since any such subset $I$ has to satisfy both $Q=P\subseteq I$
and $I\subseteq Q$, which in combination entail $I=Q$). Thus, the sum
$\sum_{\substack{I\subseteq Q;\\P\subseteq I}}\left(  -1\right)  ^{\left\vert
I\right\vert }$ has only one addend, namely the addend for $I=Q$.
Consequently, this sum simplifies as follows:%
\begin{equation}
\sum_{\substack{I\subseteq Q;\\P\subseteq I}}\left(  -1\right)  ^{\left\vert
I\right\vert }=\left(  -1\right)  ^{\left\vert Q\right\vert }.
\label{pf.lem.pie.two-sets-altsum.a.fn1.1}%
\end{equation}
On the other hand, from $P=Q$, we obtain%
\[
\left(  -1\right)  ^{\left\vert P\right\vert }\left[  P=Q\right]  =\left(
-1\right)  ^{\left\vert Q\right\vert }\underbrace{\left[  Q=Q\right]
}_{\substack{=1\\\text{(since }Q=Q\text{)}}}=\left(  -1\right)  ^{\left\vert
Q\right\vert }.
\]
Comparing this with (\ref{pf.lem.pie.two-sets-altsum.a.fn1.1}), we obtain
$\sum_{\substack{I\subseteq Q;\\P\subseteq I}}\left(  -1\right)  ^{\left\vert
I\right\vert }=\left(  -1\right)  ^{\left\vert P\right\vert }\left[
P=Q\right]  $. Thus, we have shown that (\ref{eq.lem.pie.two-sets-altsum.a})
holds under the assumption that $P=Q$.}. Hence, for the rest of this proof, we
WLOG assume that $P\neq Q$. Thus, $\left[  P=Q\right]  =0$.

Now, $P$ is a \textbf{proper} subset of $Q$ (since $P$ is a subset of $Q$ and
satisfies $P\neq Q$). Hence, there exists some $q\in Q$ such that $q\notin P$.
Fix such a $q$.

Let
\[
\mathcal{A}:=\left\{  I\subseteq Q\ \mid\ P\subseteq I\right\}  ,
\]
and let
\[
\operatorname*{sign}I:=\left(  -1\right)  ^{\left\vert I\right\vert
}\ \ \ \ \ \ \ \ \ \ \text{for each }I\in\mathcal{A}.
\]
Then,%
\begin{equation}
\sum_{I\in\mathcal{A}}\operatorname*{sign}I=\sum_{\substack{I\subseteq
Q;\\P\subseteq I}}\left(  -1\right)  ^{\left\vert I\right\vert }.
\label{pf.lem.pie.two-sets-altsum.a.3}%
\end{equation}
We shall now construct an involution $f:\mathcal{A}\rightarrow\mathcal{A}$ on
the set $\mathcal{A}$; this will allow us to apply Lemma
\ref{lem.sign.cancel2} (to $\mathcal{X}=\mathcal{A}$), and easily conclude
that $\sum_{I\in\mathcal{A}}\operatorname*{sign}I=0$ (see below for the details).

Indeed, if $I$ is a subset of $Q$, then $I\cup\left\{  q\right\}  $ is also a
subset of $Q$ (because $q\in Q$). Thus, if $I\in\mathcal{A}$, then
$I\cup\left\{  q\right\}  \in\mathcal{A}$ (because $P\subseteq I$ implies
$P\subseteq I\subseteq I\cup\left\{  q\right\}  $).

On the other hand, if $I$ is a set satisfying $P\subseteq I$, then the set
$I\setminus\left\{  q\right\}  $ also satisfies $P\subseteq I\setminus\left\{
q\right\}  $ (since $q\notin P$). Thus, if $I\in\mathcal{A}$, then
$I\setminus\left\{  q\right\}  \in\mathcal{A}$ (since $I\subseteq Q$ implies
$I\setminus\left\{  q\right\}  \subseteq I\subseteq Q$).

Now, we define a map $f:\mathcal{A}\rightarrow\mathcal{A}$ by
setting\footnote{Here, the notation $X\bigtriangleup Y$ means the symmetric
difference $\left(  X\cup Y\right)  \setminus\left(  X\cap Y\right)  $ of two
sets $X$ and $Y$ (as in Subsection \ref{subsec.gf.defs.commrings}).}%
\[
f\left(  I\right)  :=I\bigtriangleup\left\{  q\right\}  =%
\begin{cases}
I\setminus\left\{  q\right\}  , & \text{if }q\in I;\\
I\cup\left\{  q\right\}  , & \text{if }q\notin I
\end{cases}
\ \ \ \ \ \ \ \ \ \ \text{for each }I\in\mathcal{A}.
\]
This map $f$ is well-defined, because (as we have just shown in the two
paragraphs above) every $I\in\mathcal{A}$ satisfies $I\cup\left\{  q\right\}
\in\mathcal{A}$ and $I\setminus\left\{  q\right\}  \in\mathcal{A}$. Moreover,
this map $f$ is an involution\footnote{Indeed, the map $f$ merely removes $q$
from a set $I$ if $q$ is contained in $I$, and inserts it into $I$ otherwise;
but this is clearly an operation that undoes itself when performed a second
time.}. This involution $f$ has no fixed points (because if $I\in\mathcal{A}$,
then $f\left(  I\right)  =I\bigtriangleup\left\{  q\right\}  \neq I$).
Furthermore, if $I\in\mathcal{A}$, then the set $f\left(  I\right)
=I\bigtriangleup\left\{  q\right\}  $ differs from $I$ in exactly one element
(namely, $q$), and thus satisfies $\left\vert f\left(  I\right)  \right\vert
=\left\vert I\right\vert \pm1$, so that%
\[
\left(  -1\right)  ^{\left\vert f\left(  I\right)  \right\vert }=-\left(
-1\right)  ^{\left\vert I\right\vert },
\]
or, equivalently,%
\[
\operatorname*{sign}\left(  f\left(  I\right)  \right)  =-\operatorname*{sign}%
I
\]
(since the definition of $\operatorname*{sign}I$ yields $\operatorname*{sign}%
I=\left(  -1\right)  ^{\left\vert I\right\vert }$, and similarly
$\operatorname*{sign}\left(  f\left(  I\right)  \right)  =\left(  -1\right)
^{\left\vert f\left(  I\right)  \right\vert }$). Thus, Lemma
\ref{lem.sign.cancel2} (applied to $\mathcal{X}=\mathcal{A}$) shows that%
\begin{align*}
\sum_{I\in\mathcal{A}}\operatorname*{sign}I  &  =\sum_{I\in\mathcal{A}%
\setminus\mathcal{A}}\operatorname*{sign}I=\left(  \text{empty sum}\right)
\ \ \ \ \ \ \ \ \ \ \left(  \text{since }\mathcal{A}\setminus\mathcal{A}%
=\varnothing\right) \\
&  =0.
\end{align*}
Comparing this with (\ref{pf.lem.pie.two-sets-altsum.a.3}), we find%
\[
\sum_{\substack{I\subseteq Q;\\P\subseteq I}}\left(  -1\right)  ^{\left\vert
I\right\vert }=0=\left(  -1\right)  ^{\left\vert P\right\vert }\left[
P=Q\right]  \ \ \ \ \ \ \ \ \ \ \left(  \text{since }\left(  -1\right)
^{\left\vert P\right\vert }\underbrace{\left[  P=Q\right]  }_{=0}=0\right)  .
\]
Thus, (\ref{eq.lem.pie.two-sets-altsum.a}) is proved. This proves Lemma
\ref{lem.pie.two-sets-altsum} \textbf{(a)}. \medskip

\textbf{(b)} If $I$ is any subset of $Q$, then $\left\vert Q\setminus
I\right\vert =\left\vert Q\right\vert -\left\vert I\right\vert \equiv
\left\vert Q\right\vert +\left\vert I\right\vert \operatorname{mod}2$ and thus%
\[
\left(  -1\right)  ^{\left\vert Q\setminus I\right\vert }=\left(  -1\right)
^{\left\vert Q\right\vert +\left\vert I\right\vert }=\left(  -1\right)
^{\left\vert Q\right\vert }\left(  -1\right)  ^{\left\vert I\right\vert }.
\]
Hence,%
\begin{align*}
\sum_{\substack{I\subseteq Q;\\P\subseteq I}}\underbrace{\left(  -1\right)
^{\left\vert Q\setminus I\right\vert }}_{=\left(  -1\right)  ^{\left\vert
Q\right\vert }\left(  -1\right)  ^{\left\vert I\right\vert }}  &
=\sum_{\substack{I\subseteq Q;\\P\subseteq I}}\left(  -1\right)  ^{\left\vert
Q\right\vert }\left(  -1\right)  ^{\left\vert I\right\vert }=\left(
-1\right)  ^{\left\vert Q\right\vert }\underbrace{\sum_{\substack{I\subseteq
Q;\\P\subseteq I}}\left(  -1\right)  ^{\left\vert I\right\vert }%
}_{\substack{=\left(  -1\right)  ^{\left\vert P\right\vert }\left[
P=Q\right]  \\\text{(by Lemma \ref{lem.pie.two-sets-altsum} \textbf{(a)})}}}\\
&  =\left(  -1\right)  ^{\left\vert Q\right\vert }\left(  -1\right)
^{\left\vert P\right\vert }\left[  P=Q\right]  .
\end{align*}
However, it is easy to see that%
\begin{equation}
\left(  -1\right)  ^{\left\vert Q\right\vert }\left(  -1\right)  ^{\left\vert
P\right\vert }\left[  P=Q\right]  =\left[  P=Q\right]
\label{pf.lem.pie.two-sets-altsum.b.1}%
\end{equation}
\footnote{\textit{Proof of (\ref{pf.lem.pie.two-sets-altsum.b.1}):} If $P\neq
Q$, then the equality (\ref{pf.lem.pie.two-sets-altsum.b.1}) boils down to
$\left(  -1\right)  ^{\left\vert Q\right\vert }\left(  -1\right)  ^{\left\vert
P\right\vert }\cdot0=0$ (since $P\neq Q$ entails $\left[  P=Q\right]  =0$),
which is obviously true. Hence, (\ref{pf.lem.pie.two-sets-altsum.b.1}) is
proved if $P\neq Q$. Thus, for the rest of this proof, we WLOG assume that
$P=Q$. Hence, $\left(  -1\right)  ^{\left\vert Q\right\vert }\left(
-1\right)  ^{\left\vert P\right\vert }=\left(  -1\right)  ^{\left\vert
Q\right\vert }\left(  -1\right)  ^{\left\vert Q\right\vert }=\left(
-1\right)  ^{\left\vert Q\right\vert +\left\vert Q\right\vert }=1$ (since
$\left\vert Q\right\vert +\left\vert Q\right\vert =2\left\vert Q\right\vert $
is even). Therefore, $\underbrace{\left(  -1\right)  ^{\left\vert Q\right\vert
}\left(  -1\right)  ^{\left\vert P\right\vert }}_{=1}\left[  P=Q\right]
=\left[  P=Q\right]  $. This proves (\ref{pf.lem.pie.two-sets-altsum.b.1}).}.
Thus,%
\[
\sum_{\substack{I\subseteq Q;\\P\subseteq I}}\left(  -1\right)  ^{\left\vert
Q\setminus I\right\vert }=\left(  -1\right)  ^{\left\vert Q\right\vert
}\left(  -1\right)  ^{\left\vert P\right\vert }\left[  P=Q\right]  =\left[
P=Q\right]  .
\]
This proves Lemma \ref{lem.pie.two-sets-altsum} \textbf{(b)}.
\end{proof}

As said above, this completes the proofs of Theorem \ref{thm.pie.moeb}, of
Theorem \ref{thm.pie.2} and of Theorem \ref{thm.pie.1}.

While Theorem \ref{thm.pie.moeb} has played the part of the ultimate
generalization to us, it can be generalized further. Indeed, it is merely a
particular case of \emph{M\"{o}bius inversion for arbitrary posets} (see,
e.g., \cite[Proposition 3.7.1]{Stanley-EC1} or \cite[Theorem 2.3.1]{Martin21}
or \cite[Theorem 5.5.5]{Sagan19} or \cite[Theorem 6.1.10]{Sam21} or
\cite[Theorem 14.6.4]{Wagner20}).

\subsection{\label{sec.cancel.moresub}More subtractive methods}

\subsubsection{\label{subsec.cancel.moresub.all-even}Sums with varying signs}

The ideas of subtraction and cancellation have several other applications in
combinatorics. One way to achieve a lot of cancellation is by summing the same
sums but with varying signs. For instance, when $P\left(  x\right)  $ is a
univariate polynomial, the sum $P\left(  x\right)  +P\left(  -x\right)  $
contains no odd powers of $x$, since all these odd powers cancel out (for
instance, $\left(  1+x\right)  ^{3}+\left(  1-x\right)  ^{3}=6x^{2}+2$).
Meanwhile, all even powers of $x$ that appear in $P\left(  x\right)  $ will
appear in $P\left(  x\right)  +P\left(  -x\right)  $ with their coefficients
doubled. Likewise, for a polynomial $P\left(  x,y\right)  $ in two variables,
we can \textquotedblleft disable\textquotedblright\ all odd powers of $x$ by
taking $P\left(  x,y\right)  +P\left(  -x,y\right)  $; moreover, by taking the
sum $P\left(  x,y\right)  +P\left(  -x,y\right)  +P\left(  x,-y\right)
+P\left(  -x,-y\right)  $, we can \textquotedblleft disable\textquotedblright%
\ all monomials $x^{i}y^{j}$ for which at least one of $i$ and $j$ is odd,
leaving only the monomials $x^{i}y^{j}$ with $i\equiv j\equiv
0\operatorname{mod}2$ intact.

An impressive example of how this trick can be used is the proof of the
following enumerative result:

\begin{theorem}
\label{thm.cancel.all-even}Let $n\in\mathbb{N}$ and $d\in\mathbb{N}$. An
$n$-tuple $\left(  x_{1},x_{2},\ldots,x_{n}\right)  \in\left[  d\right]  ^{n}$
is said to be \emph{all-even} if each element of $\left[  d\right]  $ occurs
an even number of times in this $n$-tuple (i.e., if for each $k\in\left[
d\right]  $, the number of all $i\in\left[  n\right]  $ satisfying $x_{i}=k$
is even). For example, the $4$-tuple $\left(  1,4,4,1\right)  $ and the
$6$-tuples $\left(  1,3,3,5,1,5\right)  $ and $\left(  2,4,2,4,3,3\right)  $
are all-even, while the $4$-tuples $\left(  1,2,2,4\right)  $ and $\left(
2,4,6,4\right)  $ are not.

Then, the number of all all-even $n$-tuples $\left(  x_{1},x_{2},\ldots
,x_{n}\right)  \in\left[  d\right]  ^{n}$ is
\[
\dfrac{1}{2^{d}}\sum_{k=0}^{d}\dbinom{d}{k}\left(  d-2k\right)  ^{n}.
\]

\end{theorem}

This formula entails that the sum $\sum_{k=0}^{d}\dbinom{d}{k}\left(
d-2k\right)  ^{n}$ is nonnegative and divisible by $2^{d}$. This would not be
obvious without the theorem! The presence of addends with different signs (at
least for $n$ odd) appears to make a direct combinatorial proof impossible. We
will instead use the cancellation trick we mentioned above. (Our proof is
taken from \cite[solution to Exercise 7]{18f-hw4s}.)

As an introductory example, let us take a look at the case $n=4$ and $d=2$.
Then, Theorem \ref{thm.cancel.all-even} is saying that the number of all
all-even $4$-tuples $\left(  x_{1},x_{2},x_{3},x_{4}\right)  \in\left[
2\right]  ^{4}$ is
\begin{align*}
\dfrac{1}{2^{2}}\sum_{k=0}^{2}\dbinom{2}{k}\left(  2-2k\right)  ^{4}  &
=\dfrac{1}{4}\left(  \dbinom{2}{0}\cdot2^{4}+\dbinom{2}{1}\cdot0^{4}%
+\dbinom{2}{2}\cdot\left(  -2\right)  ^{4}\right) \\
&  =\dfrac{1}{4}\left(  2^{4}+2\cdot0^{4}+\left(  -2\right)  ^{4}\right)  .
\end{align*}
We claim that this can be seen as follows: For any integers $e_{1},e_{2}$, the
distributive law yields%
\begin{align*}
\left(  e_{1}+e_{2}\right)  ^{4}  &  =e_{1}e_{1}e_{1}e_{1}+e_{1}e_{1}%
e_{1}e_{2}+e_{1}e_{1}e_{2}e_{1}+\cdots\ \ \ \ \ \ \ \ \ \ \left(  \text{a
total of }16\text{ addends}\right) \\
&  =\sum_{\left(  x_{1},x_{2},x_{3},x_{4}\right)  \in\left[  2\right]  ^{4}%
}e_{x_{1}}e_{x_{2}}e_{x_{3}}e_{x_{4}}.
\end{align*}
Summing this over all four pairs $\left(  e_{1},e_{2}\right)  \in\left\{
1,-1\right\}  ^{2}$, we find%
\begin{align}
&  \sum_{\left(  e_{1},e_{2}\right)  \in\left\{  1,-1\right\}  ^{2}}\left(
e_{1}+e_{2}\right)  ^{4}\nonumber\\
&  =\sum_{\left(  e_{1},e_{2}\right)  \in\left\{  1,-1\right\}  ^{2}}%
\ \ \sum_{\left(  x_{1},x_{2},x_{3},x_{4}\right)  \in\left[  2\right]  ^{4}%
}e_{x_{1}}e_{x_{2}}e_{x_{3}}e_{x_{4}}\nonumber\\
&  =\sum_{\left(  x_{1},x_{2},x_{3},x_{4}\right)  \in\left[  2\right]  ^{4}%
}\ \ \sum_{\left(  e_{1},e_{2}\right)  \in\left\{  1,-1\right\}  ^{2}}%
e_{x_{1}}e_{x_{2}}e_{x_{3}}e_{x_{4}} \label{eq.thm.cancel.all-even.ex1}%
\end{align}
(here, we interchanged the summation signs). But now we can simplify the inner
sum $\sum_{\left(  e_{1},e_{2}\right)  \in\left\{  1,-1\right\}  ^{2}}%
e_{x_{1}}e_{x_{2}}e_{x_{3}}e_{x_{4}}$ rather conveniently:

\begin{itemize}
\item If the $4$-tuple $\left(  x_{1},x_{2},x_{3},x_{4}\right)  $ is all-even,
then all addends of this sum are $1$ (for example, if $\left(  x_{1}%
,x_{2},x_{3},x_{4}\right)  =\left(  1,2,2,1\right)  $, then $e_{x_{1}}%
e_{x_{2}}e_{x_{3}}e_{x_{4}}=e_{1}e_{2}e_{2}e_{1}=e_{1}^{2}e_{2}^{2}=1$ for all
$\left(  e_{1},e_{2}\right)  \in\left\{  1,-1\right\}  ^{2}$), and thus the
sum is equal to $4$.

\item If the $4$-tuple $\left(  x_{1},x_{2},x_{3},x_{4}\right)  $ is not
all-even, then the addends of this sum cancel out entirely (for example, if
$\left(  x_{1},x_{2},x_{3},x_{4}\right)  =\left(  1,2,1,1\right)  $, then
$e_{x_{1}}e_{x_{2}}e_{x_{3}}e_{x_{4}}=e_{1}e_{2}e_{1}e_{1}=e_{1}^{3}e_{2}$
changes sign whenever we flip the sign of $e_{1}$, so that the addends with
$e_{1}=1$ cancel out the addends with $e_{1}=-1$), and thus the sum is just
$0$.
\end{itemize}

Altogether, we thus see that each $\left(  x_{1},x_{2},x_{3},x_{4}\right)
\in\left[  2\right]  ^{4}$ satisfies
\[
\sum_{\left(  e_{1},e_{2}\right)  \in\left\{  1,-1\right\}  ^{2}}e_{x_{1}%
}e_{x_{2}}e_{x_{3}}e_{x_{4}}=%
\begin{cases}
4, & \text{if }\left(  x_{1},x_{2},x_{3},x_{4}\right)  \text{ is all-even};\\
0, & \text{if not.}%
\end{cases}
\]
Hence, (\ref{eq.thm.cancel.all-even.ex1}) simplifies to%
\begin{align*}
\sum_{\left(  e_{1},e_{2}\right)  \in\left\{  1,-1\right\}  ^{2}}\left(
e_{1}+e_{2}\right)  ^{4}  &  =\sum_{\left(  x_{1},x_{2},x_{3},x_{4}\right)
\in\left[  2\right]  ^{4}}%
\begin{cases}
4, & \text{if }\left(  x_{1},x_{2},x_{3},x_{4}\right)  \text{ is all-even};\\
0, & \text{if not}%
\end{cases}
\\
&  =4\cdot\left(  \text{\# of all all-even }4\text{-tuples }\left(
x_{1},x_{2},x_{3},x_{4}\right)  \in\left[  2\right]  ^{4}\right)  .
\end{align*}
We thus conclude that%
\begin{align*}
\left(  \text{\# of all all-even }4\text{-tuples }\left(  x_{1},x_{2}%
,x_{3},x_{4}\right)  \in\left[  2\right]  ^{4}\right)   &  =\dfrac{1}{4}%
\sum_{\left(  e_{1},e_{2}\right)  \in\left\{  1,-1\right\}  ^{2}}\left(
e_{1}+e_{2}\right)  ^{4}\\
&  =\dfrac{1}{4}\left(  2^{4}+2\cdot0^{4}+\left(  -2\right)  ^{4}\right)
\end{align*}
(since the sum $\sum_{\left(  e_{1},e_{2}\right)  \in\left\{  1,-1\right\}
^{2}}\left(  e_{1}+e_{2}\right)  ^{4}$ contains one addend with $e_{1}%
+e_{2}=2$, two addends with $e_{1}+e_{2}=0$, and one addend with $e_{1}%
+e_{2}=-2$). This is precisely the answer that Theorem
\ref{thm.cancel.all-even} gives us.

We shall now prove Theorem \ref{thm.cancel.all-even} by formalizing and
generalizing this argument. We begin with the obvious generalization of
(\ref{eq.thm.cancel.all-even.ex1}):

\begin{lemma}
\label{lem.cancel.all-even.l1}Let $n,d\in\mathbb{N}$. Then,%
\begin{align*}
&  \sum_{\left(  e_{1},e_{2},\ldots,e_{d}\right)  \in\left\{  1,-1\right\}
^{d}}\left(  e_{1}+e_{2}+\cdots+e_{d}\right)  ^{n}\\
&  =\sum_{\left(  x_{1},x_{2},\ldots,x_{n}\right)  \in\left[  d\right]  ^{n}%
}\ \ \sum_{\left(  e_{1},e_{2},\ldots,e_{d}\right)  \in\left\{  1,-1\right\}
^{d}}e_{x_{1}}e_{x_{2}}\cdots e_{x_{n}}.
\end{align*}

\end{lemma}

\begin{proof}
For each $\left(  e_{1},e_{2},\ldots,e_{d}\right)  \in\left\{  1,-1\right\}
^{d}$, we have $e_{1}+e_{2}+\cdots+e_{d}=\sum_{k=1}^{d}e_{k}$ and thus%
\begin{align*}
&  \left(  e_{1}+e_{2}+\cdots+e_{d}\right)  ^{n}\\
&  =\left(  \sum_{k=1}^{d}e_{k}\right)  ^{n}=\prod_{i=1}^{n}\ \ \sum_{k=1}%
^{d}e_{k}=\underbrace{\sum_{\left(  k_{1},k_{2},\ldots,k_{n}\right)
\in\left[  d\right]  \times\left[  d\right]  \times\cdots\times\left[
d\right]  }}_{=\sum_{\left(  k_{1},k_{2},\ldots,k_{n}\right)  \in\left[
d\right]  ^{n}}}\ \ \underbrace{\prod_{i=1}^{n}e_{k_{i}}}_{=e_{k_{1}}e_{k_{2}%
}\cdots e_{k_{n}}}\\
&  \ \ \ \ \ \ \ \ \ \ \ \ \ \ \ \ \ \ \ \ \left(  \text{by the product rule
(\ref{eq.lem.fps.prodrule-fin-fin.eq}), applied to }m_{i}=d\text{ and }%
p_{i,k}=e_{k}\right) \\
&  =\sum_{\left(  k_{1},k_{2},\ldots,k_{n}\right)  \in\left[  d\right]  ^{n}%
}e_{k_{1}}e_{k_{2}}\cdots e_{k_{n}}=\sum_{\left(  x_{1},x_{2},\ldots
,x_{n}\right)  \in\left[  d\right]  ^{n}}e_{x_{1}}e_{x_{2}}\cdots e_{x_{n}}%
\end{align*}
(here, we have renamed the summation index). Summing this over all $\left(
e_{1},e_{2},\ldots,e_{d}\right)  \in\left\{  1,-1\right\}  ^{d}$, we find%
\begin{align*}
&  \sum_{\left(  e_{1},e_{2},\ldots,e_{d}\right)  \in\left\{  1,-1\right\}
^{d}}\left(  e_{1}+e_{2}+\cdots+e_{d}\right)  ^{n}\\
&  =\sum_{\left(  e_{1},e_{2},\ldots,e_{d}\right)  \in\left\{  1,-1\right\}
^{d}}\ \ \sum_{\left(  x_{1},x_{2},\ldots,x_{n}\right)  \in\left[  d\right]
^{n}}e_{x_{1}}e_{x_{2}}\cdots e_{x_{n}}\\
&  =\sum_{\left(  x_{1},x_{2},\ldots,x_{n}\right)  \in\left[  d\right]  ^{n}%
}\ \ \sum_{\left(  e_{1},e_{2},\ldots,e_{d}\right)  \in\left\{  1,-1\right\}
^{d}}e_{x_{1}}e_{x_{2}}\cdots e_{x_{n}}%
\end{align*}
(here, we have interchanged the summation signs). This proves Lemma
\ref{lem.cancel.all-even.l1}.
\end{proof}

Now we shall compute the inner sums on the right hand side of Lemma
\ref{lem.cancel.all-even.l1}:

\begin{lemma}
\label{lem.cancel.all-even.l2}Let $n,d\in\mathbb{N}$. Let $\left(  x_{1}%
,x_{2},\ldots,x_{n}\right)  \in\left[  d\right]  ^{n}$. Then: \medskip

\textbf{(a)} If the $n$-tuple $\left(  x_{1},x_{2},\ldots,x_{n}\right)  $ is
not all-even\footnotemark, then%
\[
\sum_{\left(  e_{1},e_{2},\ldots,e_{d}\right)  \in\left\{  1,-1\right\}  ^{d}%
}e_{x_{1}}e_{x_{2}}\cdots e_{x_{n}}=0.
\]

\textbf{(b)} If the $n$-tuple $\left(  x_{1},x_{2},\ldots,x_{n}\right)  $ is
all-even, then%
\[
\sum_{\left(  e_{1},e_{2},\ldots,e_{d}\right)  \in\left\{  1,-1\right\}  ^{d}%
}e_{x_{1}}e_{x_{2}}\cdots e_{x_{n}}=2^{d}.
\]

\end{lemma}

\footnotetext{See Theorem \ref{thm.cancel.all-even} for the definition of
\textquotedblleft all-even\textquotedblright.}

\begin{proof}
For each $k\in\left[  d\right]  $, let $m_{k}$ be the number of all
$i\in\left[  n\right]  $ satisfying $x_{i}=k$. Then, the $n$-tuple $\left(
x_{1},x_{2},\ldots,x_{n}\right)  $ is all-even if and only if all these
numbers $m_{1},m_{2},\ldots,m_{d}$ are even (by the definition of
\textquotedblleft all-even\textquotedblright).

On the other hand, every $d$-tuple $\left(  e_{1},e_{2},\ldots,e_{d}\right)
\in\mathbb{Z}^{d}$ satisfies%
\begin{equation}
e_{x_{1}}e_{x_{2}}\cdots e_{x_{n}}=e_{1}^{m_{1}}e_{2}^{m_{2}}\cdots
e_{d}^{m_{d}}. \label{pf.lem.cancel.all-even.l2.e=e}%
\end{equation}

\begin{proof}
[Proof of (\ref{pf.lem.cancel.all-even.l2.e=e}):]Let $\left(  e_{1}%
,e_{2},\ldots,e_{d}\right)  \in\mathbb{Z}^{d}$ be any $d$-tuple. Then, each
factor of the product $e_{x_{1}}e_{x_{2}}\cdots e_{x_{n}}$ is one of the $d$
entries $e_{1},e_{2},\ldots,e_{d}$ of this $d$-tuple. Moreover, each given
entry $e_{k}$ appears exactly $m_{k}$ times in the product $e_{x_{1}}e_{x_{2}%
}\cdots e_{x_{n}}$, because $k$ appears exactly $m_{k}$ times among the
subscripts $x_{1},x_{2},\ldots,x_{n}$ (since $m_{k}$ is defined as the number
of all $i\in\left[  n\right]  $ satisfying $x_{i}=k$). Thus, in total, the
product $e_{x_{1}}e_{x_{2}}\cdots e_{x_{n}}$ contains the factor $e_{1}$
exactly $m_{1}$ times, the factor $e_{2}$ exactly $m_{2}$ times, and so on;
therefore, it equals $e_{1}^{m_{1}}e_{2}^{m_{2}}\cdots e_{d}^{m_{d}}$. This
proves (\ref{pf.lem.cancel.all-even.l2.e=e}).
\end{proof}

\textbf{(a)} Assume that the $n$-tuple $\left(  x_{1},x_{2},\ldots
,x_{n}\right)  $ is not all-even. Thus, not all of the numbers $m_{1}%
,m_{2},\ldots,m_{d}$ are even (because the $n$-tuple $\left(  x_{1}%
,x_{2},\ldots,x_{n}\right)  $ is all-even if and only if all these numbers
$m_{1},m_{2},\ldots,m_{d}$ are even). In other words, there exists some
$k\in\left[  d\right]  $ such that $m_{k}$ is odd. Consider this $k$. Hence,
$\left(  -1\right)  ^{m_{k}}=-1$ (since $m_{k}$ is odd).

For each $d$-tuple $\overrightarrow{e}=\left(  e_{1},e_{2},\ldots
,e_{d}\right)  \in\left\{  1,-1\right\}  ^{d}$, we define an integer
\[
\operatorname*{sign}\overrightarrow{e}:=e_{1}^{m_{1}}e_{2}^{m_{2}}\cdots
e_{d}^{m_{d}}.
\]

Let $f:\left\{  1,-1\right\}  ^{d}\rightarrow\left\{  1,-1\right\}  ^{d}$ be
the map that sends each $d$-tuple $\left(  e_{1},e_{2},\ldots,e_{d}\right)  $
to $\left(  e_{1},e_{2},\ldots,e_{k-1},-e_{k},e_{k+1},e_{k+2},\ldots
,e_{d}\right)  $ (in other words, $f$ flips the sign of the $k$-th entry of
the $d$-tuple). This map $f$ is an involution (since flipping the sign of the
$k$-th entry twice results in the restoration of the original sign), and has
no fixed points (since our $d$-tuples $\left(  e_{1},e_{2},\ldots
,e_{d}\right)  $ come from $\left\{  1,-1\right\}  ^{d}$ and thus contain no
zeroes, so that they cannot remain unchanged when we flip the sign of an
entry). Moreover, we have%
\[
\operatorname*{sign}\left(  f\left(  \overrightarrow{e}\right)  \right)
=-\operatorname*{sign}\overrightarrow{e}\ \ \ \ \ \ \ \ \ \ \text{for all
}\overrightarrow{e}\in\left\{  1,-1\right\}  ^{d},
\]
because when we flip the sign of the $k$-th entry $e_{k}$ of a $d$-tuple
$\overrightarrow{e}=\left(  e_{1},e_{2},\ldots,e_{d}\right)  \in\left\{
1,-1\right\}  ^{d}$, the integer $\operatorname*{sign}\overrightarrow{e}%
=e_{1}^{m_{1}}e_{2}^{m_{2}}\cdots e_{d}^{m_{d}}$ is multiplied by $\left(
-1\right)  ^{m_{k}}=-1$. Renaming the variable $\overrightarrow{e}$ as $I$ in
this claim, we can rewrite it as follows:%
\[
\operatorname*{sign}\left(  f\left(  I\right)  \right)  =-\operatorname*{sign}%
I\ \ \ \ \ \ \ \ \ \ \text{for all }I\in\left\{  1,-1\right\}  ^{d}.
\]

Hence, Lemma \ref{lem.sign.cancel2} (applied to $\mathcal{A}=\left\{
1,-1\right\}  ^{d}$ and $\mathcal{X}=\left\{  1,-1\right\}  ^{d}$) yields%
\[
\sum_{I\in\left\{  1,-1\right\}  ^{d}}\operatorname*{sign}I=\sum_{I\in\left\{
1,-1\right\}  ^{d}\setminus\left\{  1,-1\right\}  ^{d}}\operatorname*{sign}%
I=\left(  \text{empty sum}\right)  =0.
\]
But%
\begin{align*}
\sum_{I\in\left\{  1,-1\right\}  ^{d}}\operatorname*{sign}I  &  =\sum_{\left(
e_{1},e_{2},\ldots,e_{d}\right)  \in\left\{  1,-1\right\}  ^{d}}%
\underbrace{\operatorname*{sign}\left(  e_{1},e_{2},\ldots,e_{d}\right)
}_{\substack{=e_{1}^{m_{1}}e_{2}^{m_{2}}\cdots e_{d}^{m_{d}}\\\text{(by the
definition of }\operatorname*{sign}\overrightarrow{e}\text{)}}}\\
&  =\sum_{\left(  e_{1},e_{2},\ldots,e_{d}\right)  \in\left\{  1,-1\right\}
^{d}}\ \ \underbrace{e_{1}^{m_{1}}e_{2}^{m_{2}}\cdots e_{d}^{m_{d}}%
}_{\substack{=e_{x_{1}}e_{x_{2}}\cdots e_{x_{n}}\\\text{(by
(\ref{pf.lem.cancel.all-even.l2.e=e}))}}}\\
&  =\sum_{\left(  e_{1},e_{2},\ldots,e_{d}\right)  \in\left\{  1,-1\right\}
^{d}}e_{x_{1}}e_{x_{2}}\cdots e_{x_{n}}.
\end{align*}
Comparing these two equalities, we find
\[
\sum_{\left(  e_{1},e_{2},\ldots,e_{d}\right)  \in\left\{  1,-1\right\}  ^{d}%
}e_{x_{1}}e_{x_{2}}\cdots e_{x_{n}}=0.
\]
This proves Lemma \ref{lem.cancel.all-even.l2} \textbf{(a)}. \medskip

\textbf{(b)} Assume that the $n$-tuple $\left(  x_{1},x_{2},\ldots
,x_{n}\right)  $ is all-even. Thus, all the numbers $m_{1},m_{2},\ldots,m_{d}$
are even (because the $n$-tuple $\left(  x_{1},x_{2},\ldots,x_{n}\right)  $ is
all-even if and only if all these numbers $m_{1},m_{2},\ldots,m_{d}$ are even).

Now, let $\left(  e_{1},e_{2},\ldots,e_{d}\right)  \in\left\{  1,-1\right\}
^{d}$ be any $d$-tuple consisting of $1$'s and $-1$'s. Then, for each
$k\in\left[  d\right]  $, we have $e_{k}\in\left\{  1,-1\right\}  $ whereas
$m_{k}$ is even (because all the numbers $m_{1},m_{2},\ldots,m_{d}$ are even);
thus we obtain%
\begin{equation}
e_{k}^{m_{k}}=1 \label{pf.lem.cancel.all-even.l2.b.1}%
\end{equation}
(since an even power of $1$ or of $-1$ is always $1$). Therefore, all the $d$
factors of the product $e_{1}^{m_{1}}e_{2}^{m_{2}}\cdots e_{d}^{m_{d}}$ equal
$1$. Hence, this whole product is a product of $1$'s and thus equals $1$. In
other words, $e_{1}^{m_{1}}e_{2}^{m_{2}}\cdots e_{d}^{m_{d}}=1$. In view of
(\ref{pf.lem.cancel.all-even.l2.e=e}), this rewrites as $e_{x_{1}}e_{x_{2}%
}\cdots e_{x_{n}}=1$.

Forget that we fixed $\left(  e_{1},e_{2},\ldots,e_{d}\right)  $. We thus have
shown that each $\left(  e_{1},e_{2},\ldots,e_{d}\right)  \in\left\{
1,-1\right\}  ^{d}$ satisfies $e_{x_{1}}e_{x_{2}}\cdots e_{x_{n}}=1$.
Therefore,%
\begin{align*}
\sum_{\left(  e_{1},e_{2},\ldots,e_{d}\right)  \in\left\{  1,-1\right\}  ^{d}%
}\underbrace{e_{x_{1}}e_{x_{2}}\cdots e_{x_{n}}}_{=1}  &  =\sum_{\left(
e_{1},e_{2},\ldots,e_{d}\right)  \in\left\{  1,-1\right\}  ^{d}}1=\left\vert
\left\{  1,-1\right\}  ^{d}\right\vert \cdot1\\
&  =\left\vert \left\{  1,-1\right\}  ^{d}\right\vert =\left\vert \left\{
1,-1\right\}  \right\vert ^{d}=2^{d}.
\end{align*}
This proves Lemma \ref{lem.cancel.all-even.l2} \textbf{(b)}.
\end{proof}

We are now ready to prove Theorem \ref{thm.cancel.all-even}:

\begin{proof}
[Proof of Theorem \ref{thm.cancel.all-even}.]Lemma
\ref{lem.cancel.all-even.l1} yields%
\begin{align}
&  \sum_{\left(  e_{1},e_{2},\ldots,e_{d}\right)  \in\left\{  1,-1\right\}
^{d}}\left(  e_{1}+e_{2}+\cdots+e_{d}\right)  ^{n}\nonumber\\
&  =\sum_{\left(  x_{1},x_{2},\ldots,x_{n}\right)  \in\left[  d\right]  ^{n}%
}\ \ \sum_{\left(  e_{1},e_{2},\ldots,e_{d}\right)  \in\left\{  1,-1\right\}
^{d}}e_{x_{1}}e_{x_{2}}\cdots e_{x_{n}}\nonumber\\
&  =\sum_{\substack{\left(  x_{1},x_{2},\ldots,x_{n}\right)  \in\left[
d\right]  ^{n};\\\left(  x_{1},x_{2},\ldots,x_{n}\right)  \text{ is all-even}%
}}\ \ \underbrace{\sum_{\left(  e_{1},e_{2},\ldots,e_{d}\right)  \in\left\{
1,-1\right\}  ^{d}}e_{x_{1}}e_{x_{2}}\cdots e_{x_{n}}}_{\substack{=2^{d}%
\\\text{(by Lemma \ref{lem.cancel.all-even.l2} \textbf{(b)})}}}\nonumber\\
&  \qquad+\sum_{\substack{\left(  x_{1},x_{2},\ldots,x_{n}\right)  \in\left[
d\right]  ^{n};\\\left(  x_{1},x_{2},\ldots,x_{n}\right)  \text{ is not
all-even}}}\ \ \underbrace{\sum_{\left(  e_{1},e_{2},\ldots,e_{d}\right)
\in\left\{  1,-1\right\}  ^{d}}e_{x_{1}}e_{x_{2}}\cdots e_{x_{n}}%
}_{\substack{=0\\\text{(by Lemma \ref{lem.cancel.all-even.l2} \textbf{(a)})}%
}}\nonumber\\
&  =\sum_{\substack{\left(  x_{1},x_{2},\ldots,x_{n}\right)  \in\left[
d\right]  ^{n};\\\left(  x_{1},x_{2},\ldots,x_{n}\right)  \text{ is all-even}%
}}2^{d}+\underbrace{\sum_{\substack{\left(  x_{1},x_{2},\ldots,x_{n}\right)
\in\left[  d\right]  ^{n};\\\left(  x_{1},x_{2},\ldots,x_{n}\right)  \text{ is
not all-even}}}0}_{=0}=\sum_{\substack{\left(  x_{1},x_{2},\ldots
,x_{n}\right)  \in\left[  d\right]  ^{n};\\\left(  x_{1},x_{2},\ldots
,x_{n}\right)  \text{ is all-even}}}2^{d}\nonumber\\
&  =\left(  \text{\# of all all-even }\left(  x_{1},x_{2},\ldots,x_{n}\right)
\in\left[  d\right]  ^{n}\right)  \cdot2^{d}.
\label{sol.count.all-even-tups.c1.pf.b.lhs}%
\end{align}

Now, we shall compute the left hand side of this equality in a different way.
Clearly, a $d$-tuple $\left(  e_{1},e_{2},\ldots,e_{d}\right)  \in\left\{
1,-1\right\}  ^{d}$ is uniquely determined by the set of all its positions
occupied by $-1$'s (that is, by the set $\left\{  i\in\left[  d\right]
\ \mid\ e_{i}=-1\right\}  $). Conversely, for each subset $S$ of $\left[
d\right]  $, there exists a $d$-tuple $\left(  e_{1},e_{2},\ldots
,e_{d}\right)  \in\left\{  1,-1\right\}  ^{d}$ whose set $\left\{  i\in\left[
d\right]  \ \mid\ e_{i}=-1\right\}  $ is $S$. Thus, the map
\begin{align*}
\Phi:\left\{  1,-1\right\}  ^{d}  &  \rightarrow\left\{  \text{subsets of
}\left[  d\right]  \right\}  ,\\
\left(  e_{1},e_{2},\ldots,e_{d}\right)   &  \mapsto\left\{  i\in\left[
d\right]  \ \mid\ e_{i}=-1\right\}
\end{align*}
is a bijection. Moreover, if this bijection $\Phi$ sends a $d$-tuple $\left(
e_{1},e_{2},\ldots,e_{d}\right)  \in\left\{  1,-1\right\}  ^{d}$ to some
subset $S$ of $\left[  d\right]  $, then%
\[
e_{1}+e_{2}+\cdots+e_{d}=d-2\cdot\left\vert S\right\vert
\]
(in fact, each entry $e_{i}$ of the $d$-tuple $\left(  e_{1},e_{2}%
,\ldots,e_{d}\right)  $ equals $-1$ if $i\in S$ and equals $1$ if $i\notin S$;
thus, the sum $e_{1}+e_{2}+\cdots+e_{d}$ has $\left\vert S\right\vert $ many
addends equal to $-1$ while its remaining $d-\left\vert S\right\vert $ addends
are equal to $1$; therefore, this sum totals to $\left\vert S\right\vert
\cdot\left(  -1\right)  +\left(  d-\left\vert S\right\vert \right)
\cdot1=-\left\vert S\right\vert +\left(  d-\left\vert S\right\vert \right)
=d-2\cdot\left\vert S\right\vert $). Thus,%
\[
\sum_{\left(  e_{1},e_{2},\ldots,e_{d}\right)  \in\left\{  1,-1\right\}  ^{d}%
}\left(  e_{1}+e_{2}+\cdots+e_{d}\right)  ^{n}=\sum_{S\subseteq\left[
d\right]  }\left(  d-2\cdot\left\vert S\right\vert \right)  ^{n}%
\]
(since the sum on the left hand side is obtained from the sum on the right
hand side by substituting $\Phi\left(  e_{1},e_{2},\ldots,e_{d}\right)  $ for
$S$). In view of%
\begin{align*}
&  \sum_{S\subseteq\left[  d\right]  }\left(  d-2\cdot\left\vert S\right\vert
\right)  ^{n}\\
&  =\sum_{k=0}^{d}\ \ \sum_{\substack{S\subseteq\left[  d\right]
;\\\left\vert S\right\vert =k}}\left(  d-2\cdot\underbrace{\left\vert
S\right\vert }_{=k}\right)  ^{n}\ \ \ \ \ \ \ \ \ \ \left(
\begin{array}
[c]{c}%
\text{here, we have split the sum}\\
\text{according to the value }\left\vert S\right\vert
\end{array}
\right) \\
&  =\sum_{k=0}^{d}\ \ \underbrace{\sum_{\substack{S\subseteq\left[  d\right]
;\\\left\vert S\right\vert =k}}\left(  d-2k\right)  ^{n}}_{=\left(  \text{\#
of all }S\subseteq\left[  d\right]  \text{ satisfying }\left\vert S\right\vert
=k\right)  \cdot\left(  d-2k\right)  ^{n}}\\
&  =\sum_{k=0}^{d}\underbrace{\left(  \text{\# of all }S\subseteq\left[
d\right]  \text{ satisfying }\left\vert S\right\vert =k\right)  }_{=\left(
\text{\# of all }k\text{-element subsets of }\left[  d\right]  \right)
=\dbinom{d}{k}}\cdot\left(  d-2k\right)  ^{n}\\
&  =\sum_{k=0}^{d}\dbinom{d}{k}\left(  d-2k\right)  ^{n},
\end{align*}
we can rewrite this as%
\[
\sum_{\left(  e_{1},e_{2},\ldots,e_{d}\right)  \in\left\{  1,-1\right\}  ^{d}%
}\left(  e_{1}+e_{2}+\cdots+e_{d}\right)  ^{n}=\sum_{k=0}^{d}\dbinom{d}%
{k}\left(  d-2k\right)  ^{n}.
\]
Comparing this with \eqref{sol.count.all-even-tups.c1.pf.b.lhs}, we obtain%
\[
\left(  \text{\# of all all-even }\left(  x_{1},x_{2},\ldots,x_{n}\right)
\in\left[  d\right]  ^{n}\right)  \cdot2^{d}=\sum_{k=0}^{d}\dbinom{d}%
{k}\left(  d-2k\right)  ^{n}.
\]
Dividing this equality by $2^{d}$, we obtain%
\[
\left(  \text{\# of all all-even }\left(  x_{1},x_{2},\ldots,x_{n}\right)
\in\left[  d\right]  ^{n}\right)  =\dfrac{1}{2^{d}}\sum_{k=0}^{d}\dbinom{d}%
{k}\left(  d-2k\right)  ^{n}.
\]
This proves Theorem \ref{thm.cancel.all-even}.
\end{proof}

\begin{noncompile}
Compare with \cite[Corollary 2.5]{Stanle18}.
\end{noncompile}

\subsection{\label{sec.sign.det}Determinants}

Determinants were introduced by Leibniz in the 17th century, and quickly
became one of the most powerful tools in mathematics. They remained so until
the early 20th century. There is a 5-volume book by Thomas Muir \cite{Muir}
that merely summarizes the results found on determinants... until 1920.

Most of these old results are still interesting and nontrivial. The relative
role of determinants in mathematics has declined mainly because other parts of
mathematics have \textquotedblleft caught up\textquotedblright\ and have
produced easier ways to reach many of the places that were previously only
accessible through the study of determinants.

As with anything else, we will just present some of the most basic results and
methods related to determinants. For more, see \cite{MuiMet60}, \cite{Zeilbe},
\cite[Chapter 6]{detnotes}, \cite[Chapter I]{Prasolov}, \cite[Chapter
9]{BruRys91}, \cite{Krattenthaler} and various other sources. A good
introduction to the most fundamental properties is \cite{Strick13}.

\begin{convention}
\label{conv.det.K}For the rest of Section \ref{sec.sign.det}, we fix a
commutative ring $K$. In most examples, $K$ will be $\mathbb{Z}$ or
$\mathbb{Q}$ or a polynomial ring.
\end{convention}

\begin{convention}
\label{conv.det.matrices}Let $n,m\in\mathbb{N}$. \medskip

\textbf{(a)} If $A$ is an $n\times m$-matrix, then $A_{i,j}$ shall mean the
$\left(  i,j\right)  $-th entry of $A$, that is, the entry of $A$ in row $i$
and column $j$. \medskip

\textbf{(b)} If $a_{i,j}$ is an element of $K$ for each $i\in\left[  n\right]
$ and each $j\in\left[  m\right]  $, then
\[
\left(  a_{i,j}\right)  _{1\leq i\leq n,\ 1\leq j\leq m}%
\]
shall denote the $n\times m$-matrix whose $\left(  i,j\right)  $-th entry is
$a_{i,j}$ for all $i\in\left[  n\right]  $ and $j\in\left[  m\right]  $.
Explicitly:%
\[
\left(  a_{i,j}\right)  _{1\leq i\leq n,\ 1\leq j\leq m}=\left(
\begin{array}
[c]{cccc}%
a_{1,1} & a_{1,2} & \cdots & a_{1,m}\\
a_{2,1} & a_{2,2} & \cdots & a_{2,m}\\
\vdots & \vdots & \ddots & \vdots\\
a_{n,1} & a_{n,2} & \cdots & a_{n,m}%
\end{array}
\right)  .
\]

Note that the letters \textquotedblleft$i$\textquotedblright\ and
\textquotedblleft$j$\textquotedblright\ in the notation \textquotedblleft%
$\left(  a_{i,j}\right)  _{1\leq i\leq n,\ 1\leq j\leq m}$\textquotedblright%
\ are not carved in stone. We could just as well use any other letters
instead, and write $\left(  a_{x,y}\right)  _{1\leq x\leq n,\ 1\leq y\leq m}$
or (somewhat misleadingly, but technically correctly) $\left(  a_{j,i}\right)
_{1\leq j\leq n,\ 1\leq i\leq m}$ for the exact same matrix. (However,
$\left(  a_{j,i}\right)  _{1\leq i\leq m,\ 1\leq j\leq n}$ is a different
matrix. Whichever index is mentioned first in the subscript after the closing
parenthesis is used to index rows; the other index is used to index columns.)
\medskip

\textbf{(c)} We let $K^{n\times m}$ denote the set of all $n\times m$-matrices
with entries in $K$. This is a $K$-module. If $n=m$, this is also a
$K$-algebra. \medskip

\textbf{(d)} Let $A\in K^{n\times m}$ be an $n\times m$-matrix. The
\emph{transpose} $A^{T}$ of $A$ is defined to be the $m\times n$-matrix whose
entries are given by%
\[
\left(  A^{T}\right)  _{i,j}=A_{j,i}\ \ \ \ \ \ \ \ \ \ \text{for all }%
i\in\left[  m\right]  \text{ and }j\in\left[  n\right]  .
\]

\end{convention}

\subsubsection{Definition}

There are several ways to define the determinant of a square matrix. The
following is the most direct one:

\begin{definition}
\label{def.det.det}Let $n\in\mathbb{N}$. Let $A\in K^{n\times n}$ be an
$n\times n$-matrix. The \emph{determinant} $\det A$ of $A$ is defined to be
the element%
\[
\sum_{\sigma\in S_{n}}\left(  -1\right)  ^{\sigma}\underbrace{A_{1,\sigma
\left(  1\right)  }A_{2,\sigma\left(  2\right)  }\cdots A_{n,\sigma\left(
n\right)  }}_{=\prod\limits_{i=1}^{n}A_{i,\sigma\left(  i\right)  }}%
\]
of $K$. Here, as before:

\begin{itemize}
\item we let $S_{n}$ denote the $n$-th symmetric group (i.e., the group of
permutations of $\left[  n\right]  =\left\{  1,2,\ldots,n\right\}  $);

\item we let $\left(  -1\right)  ^{\sigma}$ denote the sign of the permutation
$\sigma$ (as defined in Definition \ref{def.perm.sign}).
\end{itemize}
\end{definition}

\begin{example}
For $n=2$, we have%
\begin{align*}
\det A  &  =\sum_{\sigma\in S_{2}}\left(  -1\right)  ^{\sigma}A_{1,\sigma
\left(  1\right)  }A_{2,\sigma\left(  2\right)  }\\
&  =\underbrace{\left(  -1\right)  ^{\operatorname*{id}}}_{=1}%
\underbrace{A_{1,\operatorname*{id}\left(  1\right)  }}_{=A_{1,1}%
}\underbrace{A_{2,\operatorname*{id}\left(  2\right)  }}_{=A_{2,2}%
}+\underbrace{\left(  -1\right)  ^{s_{1}}}_{=-1}\underbrace{A_{1,s_{1}\left(
1\right)  }}_{=A_{1,2}}\underbrace{A_{2,s_{1}\left(  2\right)  }}_{=A_{2,1}}\\
&  \ \ \ \ \ \ \ \ \ \ \ \ \ \ \ \ \ \ \ \ \left(
\begin{array}
[c]{c}%
\text{since the two elements of }S_{2}\text{ are the identity}\\
\text{map }\operatorname*{id}\text{ and the simple transposition }%
s_{1}=t_{1,2}%
\end{array}
\right) \\
&  =A_{1,1}A_{2,2}-A_{1,2}A_{2,1}.
\end{align*}
Using less cumbersome notations, we can rewrite this as follows: For any
$a,b,a^{\prime},b^{\prime}\in K$, we have%
\[
\det\left(
\begin{array}
[c]{cc}%
a & b\\
a^{\prime} & b^{\prime}%
\end{array}
\right)  =ab^{\prime}-ba^{\prime}.
\]

Similarly, for $n=3$, we obtain%
\[
\det\left(
\begin{array}
[c]{ccc}%
a & b & c\\
a^{\prime} & b^{\prime} & c^{\prime}\\
a^{\prime\prime} & b^{\prime\prime} & c^{\prime\prime}%
\end{array}
\right)  =ab^{\prime}c^{\prime\prime}-ac^{\prime}b^{\prime\prime}-ba^{\prime
}c^{\prime\prime}+bc^{\prime}a^{\prime\prime}+ca^{\prime}b^{\prime\prime
}-cb^{\prime}a^{\prime\prime}.
\]
(The six addends on the right hand side here correspond to the six
permutations in $S_{3}$, which in one-line notation are $123$, $132$, $213$,
$231$, $312$, and $321$, respectively.)

Similarly, for $n=1$, we obtain that the determinant of the $1\times1$-matrix
$\left(
\begin{array}
[c]{c}%
a
\end{array}
\right)  $ is%
\[
\det\left(
\begin{array}
[c]{c}%
a
\end{array}
\right)  =a.
\]
Here, the \textquotedblleft$\left(
\begin{array}
[c]{c}%
a
\end{array}
\right)  $\textquotedblright\ on the left hand side is a $1\times1$-matrix.

Finally, for $n=0$, we obtain that the determinant of the $0\times0$-matrix
$\left(  {}\right)  $ (this is an empty matrix, with no rows and no columns)
is%
\[
\det\left(  {}\right)  =\left(  \text{empty product}\right)  =1.
\]

\end{example}

Some (particularly, older) texts use the notation $\left\vert A\right\vert $
instead of $\det A$ for the determinant of the matrix $A$.

The above definition of the determinant is purely combinatorial: it is an
alternating sum over the $n$-th symmetric group $S_{n}$. Typically, when
computing determinants, this definition is not in itself very useful (e.g.,
because $S_{n}$ gets rather large when $n$ is large). However, in some cases,
it suffices. Here are a few examples:

\begin{example}
\label{exa.det.hollow5x5}Let $a,b,c,d,e,\ldots,p$ be $16$ elements of $K$.
Prove that%
\[
\det\left(
\begin{array}
[c]{ccccc}%
a & b & c & d & e\\
p & 0 & 0 & 0 & f\\
o & 0 & 0 & 0 & g\\
n & 0 & 0 & 0 & h\\
m & l & k & j & i
\end{array}
\right)  =0.
\]
(The \textquotedblleft$o$\textquotedblright\ is a letter \textquotedblleft
oh\textquotedblright, not a zero. Not that it matters much...)
\end{example}

\begin{proof}
[Proof of Example \ref{exa.det.hollow5x5}.]Let $A$ be the $5\times5$-matrix
whose determinant we are trying to identify as $0$; thus, $A_{1,1}=a$ and
$A_{1,2}=b$ and $A_{3,2}=0$ and so on. Notice that%
\begin{equation}
A_{i,j}=0\ \ \ \ \ \ \ \ \ \ \text{whenever }i,j\in\left\{  2,3,4\right\}
\label{pf.exa.det.hollow5x5.0}%
\end{equation}
(since $A$ has a \textquotedblleft hollow core\textquotedblright, i.e., a
$3\times3$-square consisting entirely of zeroes in its middle). We must prove
that $\det A=0$.

Our definition of $\det A$ yields%
\begin{equation}
\det A=\sum_{\sigma\in S_{5}}\left(  -1\right)  ^{\sigma}\prod\limits_{i=1}%
^{5}A_{i,\sigma\left(  i\right)  }. \label{pf.exa.det.hollow5x5.detA=}%
\end{equation}

Now, I claim that each of the addends in the sum on the right hand side is
$0$. In other words, I claim that $\prod\limits_{i=1}^{5}A_{i,\sigma\left(
i\right)  }=0$ for each $\sigma\in S_{5}$.

To prove this, fix $\sigma\in S_{5}$. The three numbers $\sigma\left(
2\right)  ,\sigma\left(  3\right)  ,\sigma\left(  4\right)  $ are three
distinct elements of $\left[  5\right]  $ (distinct because $\sigma$ is
injective), so they cannot all belong to the $2$-element set $\left\{
1,5\right\}  $ (since there are no three distinct elements in a $2$-element
set). Hence, at least one of them must belong to the complement $\left\{
2,3,4\right\}  $ of this set. In other words, there exists some $i\in\left\{
2,3,4\right\}  $ such that $\sigma\left(  i\right)  \in\left\{  2,3,4\right\}
$. This $i$ must then satisfy $A_{i,\sigma\left(  i\right)  }=0$ (by
(\ref{pf.exa.det.hollow5x5.0}), applied to $j=\sigma\left(  i\right)  $).
Thus, we have shown that there exists some $i\in\left\{  2,3,4\right\}  $ such
that $A_{i,\sigma\left(  i\right)  }=0$.

This shows that at least one factor of the product $\prod\limits_{i=1}%
^{5}A_{i,\sigma\left(  i\right)  }$ is $0$. Thus, the entire product is $0$.

Forget that we fixed $\sigma$. We thus have proved that $\prod\limits_{i=1}%
^{5}A_{i,\sigma\left(  i\right)  }=0$ for each $\sigma\in S_{5}$. Hence,%
\[
\det A=\sum_{\sigma\in S_{5}}\left(  -1\right)  ^{\sigma}\underbrace{\prod
\limits_{i=1}^{5}A_{i,\sigma\left(  i\right)  }}_{=0}=0,
\]
and thus Example \ref{exa.det.hollow5x5} is proved.

[We note that there are various alternative proofs, e.g., using Laplace
expansion. Also, if $K$ is a field, you can argue that $\det A=0$ using rank
arguments. See \cite[Exercise 6.47 \textbf{(a)}]{detnotes} for a
generalization of Example \ref{exa.det.hollow5x5}.]
\end{proof}

\begin{example}
Let $n\in\mathbb{N}$, and let $x_{1},x_{2},\ldots,x_{n}\in K$ and $y_{1}%
,y_{2},\ldots,y_{n}\in K$. Compute%
\[
\det\left(  \left(  x_{i}y_{j}\right)  _{1\leq i\leq n,\ 1\leq j\leq
n}\right)  =\det\left(
\begin{array}
[c]{cccc}%
x_{1}y_{1} & x_{1}y_{2} & \cdots & x_{1}y_{n}\\
x_{2}y_{1} & x_{2}y_{2} & \cdots & x_{2}y_{n}\\
\vdots & \vdots & \ddots & \vdots\\
x_{n}y_{1} & x_{n}y_{2} & \cdots & x_{n}y_{n}%
\end{array}
\right)  .
\]

\end{example}

Let us experiment with small $n$'s:%
\begin{align*}
\det\left(  {}\right)   &  =1;\\
\det\left(
\begin{array}
[c]{c}%
x_{1}y_{1}%
\end{array}
\right)   &  =x_{1}y_{1};\\
\det\left(
\begin{array}
[c]{cc}%
x_{1}y_{1} & x_{1}y_{2}\\
x_{2}y_{1} & x_{2}y_{2}%
\end{array}
\right)   &  =0;\\
\det\left(
\begin{array}
[c]{ccc}%
x_{1}y_{1} & x_{1}y_{2} & x_{1}y_{3}\\
x_{2}y_{1} & x_{2}y_{2} & x_{2}y_{3}\\
x_{3}y_{1} & x_{3}y_{2} & x_{3}y_{3}%
\end{array}
\right)   &  =0.
\end{align*}
This makes us suspect the following:

\begin{proposition}
\label{prop.det.xiyj}Let $n\in\mathbb{N}$ be such that $n\geq2$. Let
$x_{1},x_{2},\ldots,x_{n}\in K$ and $y_{1},y_{2},\ldots,y_{n}\in K$. Then,%
\[
\det\left(  \left(  x_{i}y_{j}\right)  _{1\leq i\leq n,\ 1\leq j\leq
n}\right)  =0.
\]

\end{proposition}

\begin{proof}
[Proof of Proposition \ref{prop.det.xiyj} (sketched).]The definition of the
determinant yields%
\begin{align*}
\det\left(  \left(  x_{i}y_{j}\right)  _{1\leq i\leq n,\ 1\leq j\leq
n}\right)   &  =\sum_{\sigma\in S_{n}}\left(  -1\right)  ^{\sigma
}\underbrace{\left(  x_{1}y_{\sigma\left(  1\right)  }\right)  \left(
x_{2}y_{\sigma\left(  2\right)  }\right)  \cdots\left(  x_{n}y_{\sigma\left(
n\right)  }\right)  }_{=\left(  x_{1}x_{2}\cdots x_{n}\right)  \left(
y_{\sigma\left(  1\right)  }y_{\sigma\left(  2\right)  }\cdots y_{\sigma
\left(  n\right)  }\right)  }\\
&  =\sum_{\sigma\in S_{n}}\left(  -1\right)  ^{\sigma}\left(  x_{1}x_{2}\cdots
x_{n}\right)  \underbrace{\left(  y_{\sigma\left(  1\right)  }y_{\sigma\left(
2\right)  }\cdots y_{\sigma\left(  n\right)  }\right)  }_{\substack{=y_{1}%
y_{2}\cdots y_{n}\\\text{(since }\sigma\text{ is a bijection }\left[
n\right]  \rightarrow\left[  n\right]  \text{)}}}\\
&  =\sum_{\sigma\in S_{n}}\left(  -1\right)  ^{\sigma}\left(  x_{1}x_{2}\cdots
x_{n}\right)  \left(  y_{1}y_{2}\cdots y_{n}\right) \\
&  =\left(  x_{1}x_{2}\cdots x_{n}\right)  \left(  y_{1}y_{2}\cdots
y_{n}\right)  \underbrace{\sum_{\sigma\in S_{n}}\left(  -1\right)  ^{\sigma}%
}_{\substack{=0\\\text{(by (\ref{eq.cor.perm.num-even.sum-sign}))}}}=0.
\end{align*}
Thus, Proposition \ref{prop.det.xiyj} is proved.
\end{proof}

As a consequence of Proposition \ref{prop.det.xiyj} (applied to $x_{i}=x$ and
$y_{j}=1$), we see that%
\begin{equation}
\det\underbrace{\left(
\begin{array}
[c]{cccc}%
x & x & \cdots & x\\
x & x & \cdots & x\\
\vdots & \vdots & \ddots & \vdots\\
x & x & \cdots & x
\end{array}
\right)  }_{\text{an }n\times n\text{-matrix with }n\geq2}=0
\label{eq.prop.det.xiyj.cor-x}%
\end{equation}
for any $x\in K$. In other words, if all entries of a square matrix of size
$\geq2$ are equal, then the determinant of this matrix is $0$.

\begin{example}
Let $n\in\mathbb{N}$, and let $x_{1},x_{2},\ldots,x_{n}\in K$ and $y_{1}%
,y_{2},\ldots,y_{n}\in K$. Compute%
\[
\det\left(  \left(  x_{i}+y_{j}\right)  _{1\leq i\leq n,\ 1\leq j\leq
n}\right)  =\det\left(
\begin{array}
[c]{cccc}%
x_{1}+y_{1} & x_{1}+y_{2} & \cdots & x_{1}+y_{n}\\
x_{2}+y_{1} & x_{2}+y_{2} & \cdots & x_{2}+y_{n}\\
\vdots & \vdots & \ddots & \vdots\\
x_{n}+y_{1} & x_{n}+y_{2} & \cdots & x_{n}+y_{n}%
\end{array}
\right)  .
\]

\end{example}

Let us experiment with small $n$'s:%
\begin{align*}
\det\left(  {}\right)   &  =1;\\
\det\left(
\begin{array}
[c]{c}%
x_{1}+y_{1}%
\end{array}
\right)   &  =x_{1}+y_{1};\\
\det\left(
\begin{array}
[c]{cc}%
x_{1}+y_{1} & x_{1}+y_{2}\\
x_{2}+y_{1} & x_{2}+y_{2}%
\end{array}
\right)   &  =-\left(  x_{1}-x_{2}\right)  \left(  y_{1}-y_{2}\right)  ;\\
\det\left(
\begin{array}
[c]{ccc}%
x_{1}+y_{1} & x_{1}+y_{2} & x_{1}+y_{3}\\
x_{2}+y_{1} & x_{2}+y_{2} & x_{2}+y_{3}\\
x_{3}+y_{1} & x_{3}+y_{2} & x_{3}+y_{3}%
\end{array}
\right)   &  =0.
\end{align*}
So we suspect the following:

\begin{proposition}
\label{prop.det.xi+yj}Let $n\in\mathbb{N}$ be such that $n\geq3$. Let
$x_{1},x_{2},\ldots,x_{n}\in K$ and $y_{1},y_{2},\ldots,y_{n}\in K$. Then,%
\[
\det\left(  \left(  x_{i}+y_{j}\right)  _{1\leq i\leq n,\ 1\leq j\leq
n}\right)  =0.
\]

\end{proposition}

\begin{proof}
[Proof of Proposition \ref{prop.det.xi+yj} (sketched).](See \cite[Example
6.7]{detnotes} for more details.) The definition of the determinant yields%
\begin{align*}
&  \det\left(  \left(  x_{i}+y_{j}\right)  _{1\leq i\leq n,\ 1\leq j\leq
n}\right) \\
&  =\sum_{\sigma\in S_{n}}\left(  -1\right)  ^{\sigma}\underbrace{\left(
x_{1}+y_{\sigma\left(  1\right)  }\right)  \left(  x_{2}+y_{\sigma\left(
2\right)  }\right)  \cdots\left(  x_{n}+y_{\sigma\left(  n\right)  }\right)
}_{\substack{=\prod_{i=1}^{n}\left(  x_{i}+y_{\sigma\left(  i\right)
}\right)  =\sum\limits_{I\subseteq\left[  n\right]  }\left(  \prod
\limits_{i\in I}x_{i}\right)  \left(  \prod\limits_{i\in\left[  n\right]
\setminus I}y_{\sigma\left(  i\right)  }\right)  \\\text{(by
(\ref{eq.lem.prodrule.sum-ai-plus-bi.eq}), applied to }a_{i}=x_{i}\text{ and
}b_{i}=y_{\sigma\left(  i\right)  }\text{)}}}\\
&  =\sum_{\sigma\in S_{n}}\left(  -1\right)  ^{\sigma}\sum\limits_{I\subseteq
\left[  n\right]  }\left(  \prod\limits_{i\in I}x_{i}\right)  \left(
\prod\limits_{i\in\left[  n\right]  \setminus I}y_{\sigma\left(  i\right)
}\right) \\
&  =\sum\limits_{I\subseteq\left[  n\right]  }\ \ \sum_{\sigma\in S_{n}%
}\left(  -1\right)  ^{\sigma}\left(  \prod\limits_{i\in I}x_{i}\right)
\left(  \prod\limits_{i\in\left[  n\right]  \setminus I}y_{\sigma\left(
i\right)  }\right) \\
&  =\sum\limits_{I\subseteq\left[  n\right]  }\left(  \prod\limits_{i\in
I}x_{i}\right)  \sum_{\sigma\in S_{n}}\left(  -1\right)  ^{\sigma}%
\prod\limits_{i\in\left[  n\right]  \setminus I}y_{\sigma\left(  i\right)  }.
\end{align*}

Now, I claim that the inner sum is $0$ for each $I$. In other words, I claim
that%
\begin{equation}
\sum_{\sigma\in S_{n}}\left(  -1\right)  ^{\sigma}\prod\limits_{i\in\left[
n\right]  \setminus I}y_{\sigma\left(  i\right)  }%
=0\ \ \ \ \ \ \ \ \ \ \text{for each }I\subseteq\left[  n\right]  .
\label{pf.prop.det.xi+yj.innersum}%
\end{equation}
For instance, for $I=\left\{  1,2\right\}  $, this is claiming that%
\[
\sum_{\sigma\in S_{n}}\left(  -1\right)  ^{\sigma}y_{\sigma\left(  3\right)
}y_{\sigma\left(  4\right)  }\cdots y_{\sigma\left(  n\right)  }=0.
\]

[\textit{Proof of (\ref{pf.prop.det.xi+yj.innersum}):} Fix a subset $I$ of
$\left[  n\right]  $. We shall show that all addends in the sum $\sum
_{\sigma\in S_{n}}\left(  -1\right)  ^{\sigma}\prod\limits_{i\in\left[
n\right]  \setminus I}y_{\sigma\left(  i\right)  }$ cancel each other -- i.e.,
that for each addend in this sum, there is a different addend with the same
product of $y_{j}$'s but a different sign $\left(  -1\right)  ^{\sigma}$. To
achieve this, we need to pair up each $\sigma\in S_{n}$ with a different
permutation $\sigma^{\prime}=\sigma t_{u,v}\in S_{n}$ that satisfies
$\prod\limits_{i\in\left[  n\right]  \setminus I}y_{\sigma^{\prime}\left(
i\right)  }=\prod\limits_{i\in\left[  n\right]  \setminus I}y_{\sigma\left(
i\right)  }$ but $\left(  -1\right)  ^{\sigma^{\prime}}=-\left(  -1\right)
^{\sigma}$. (Indeed, this pairing will then produce the required
cancellations: the addend for each $\sigma$ will cancel the addend for the
corresponding $\sigma^{\prime}$. To be more rigorous, we are here applying
Lemma \ref{lem.sign.cancel2} to $\mathcal{A}=S_{n}$, $\mathcal{X}=S_{n}$ and
$\operatorname*{sign}\sigma=\left(  -1\right)  ^{\sigma}\prod\limits_{i\in
\left[  n\right]  \setminus I}y_{\sigma\left(  i\right)  }$ (of course, this
should not be confused for the notation $\operatorname*{sign}\sigma$ for
$\left(  -1\right)  ^{\sigma}$) and \newline$f=\left(  \text{the map }%
S_{n}\rightarrow S_{n}\text{ that sends each }\sigma\in S_{n}\text{ to the
corresponding }\sigma^{\prime}\right)  $.)

So let us construct our pairing. Indeed, from $I\subseteq\left[  n\right]  $,
we obtain $\left\vert I\right\vert +\left\vert \left[  n\right]  \setminus
I\right\vert =\left\vert \left[  n\right]  \right\vert =n\geq3$; hence, at
least one of the two sets $I$ and $\left[  n\right]  \setminus I$ has size
$>1$. In other words, we must be in one of the following two cases:

\textit{Case 1:} We have $\left\vert I\right\vert >1$.

\textit{Case 2:} We have $\left\vert \left[  n\right]  \setminus I\right\vert
>1$.

Let us first consider Case 1. In this case, we have $\left\vert I\right\vert
>1$. Thus, $\left\vert I\right\vert \geq2$. Pick two distinct elements $u$ and
$v$ of $I$. (These exist, since $\left\vert I\right\vert \geq2$.) Now, for
each permutation $\sigma\in S_{n}$, we set $\sigma^{\prime}:=\sigma t_{u,v}\in
S_{n}$. Then, each $\sigma\in S_{n}$ satisfies $\sigma^{\prime\prime
}=\underbrace{\sigma^{\prime}}_{=\sigma t_{u,v}}t_{u,v}=\sigma
\underbrace{t_{u,v}t_{u,v}}_{=t_{u,v}^{2}=\operatorname*{id}}=\sigma$. Hence,
we can pair up each $\sigma\in S_{n}$ with $\sigma^{\prime}\in S_{n}$. Any two
permutations $\sigma$ and $\sigma^{\prime}$ that are paired with each other
have different signs (indeed, if $\sigma\in S_{n}$, then $\sigma^{\prime
}=\sigma t_{u,v}$ and thus $\left(  -1\right)  ^{\sigma^{\prime}}=\left(
-1\right)  ^{\sigma t_{u,v}}=\left(  -1\right)  ^{\sigma}\underbrace{\left(
-1\right)  ^{t_{u,v}}}_{=-1}=-\left(  -1\right)  ^{\sigma}$), but the
corresponding products $\prod\limits_{i\in\left[  n\right]  \setminus
I}y_{\sigma\left(  i\right)  }$ and $\prod\limits_{i\in\left[  n\right]
\setminus I}y_{\sigma^{\prime}\left(  i\right)  }$ are equal (indeed, the
permutations $\sigma$ and $\sigma^{\prime}=\sigma t_{u,v}$ differ only in
their values at $u$ and $v$, but neither of these two values appears in any of
our two products, since $u,v\in I$). This shows that our pairing has precisely
the properties we want: Each $\sigma\in S_{n}$ satisfies $\prod\limits_{i\in
\left[  n\right]  \setminus I}y_{\sigma^{\prime}\left(  i\right)  }%
=\prod\limits_{i\in\left[  n\right]  \setminus I}y_{\sigma\left(  i\right)  }$
but $\left(  -1\right)  ^{\sigma^{\prime}}=-\left(  -1\right)  ^{\sigma}$ (and
thus $\sigma^{\prime}\neq\sigma$). As explained above, this completes the
proof of (\ref{pf.prop.det.xi+yj.innersum}) in Case 1.

Let us now consider Case 2. In this case, we have $\left\vert \left[
n\right]  \setminus I\right\vert >1$. Thus, $\left\vert \left[  n\right]
\setminus I\right\vert \geq2$. Pick two distinct elements $u$ and $v$ of
$\left[  n\right]  \setminus I$. (These exist, since $\left\vert \left[
n\right]  \setminus I\right\vert \geq2$.) We now proceed just as we did in
Case 1: For each permutation $\sigma\in S_{n}$, we set $\sigma^{\prime
}:=\sigma t_{u,v}\in S_{n}$. Then, each $\sigma\in S_{n}$ satisfies
$\sigma^{\prime\prime}=\underbrace{\sigma^{\prime}}_{=\sigma t_{u,v}}%
t_{u,v}=\sigma\underbrace{t_{u,v}t_{u,v}}_{=t_{u,v}^{2}=\operatorname*{id}%
}=\sigma$. Hence, we can pair up each $\sigma\in S_{n}$ with $\sigma^{\prime
}\in S_{n}$. Any two permutations $\sigma$ and $\sigma^{\prime}$ that are
paired with each other have different signs (this can be seen as in Case 1),
but the corresponding products $\prod\limits_{i\in\left[  n\right]  \setminus
I}y_{\sigma\left(  i\right)  }$ and $\prod\limits_{i\in\left[  n\right]
\setminus I}y_{\sigma^{\prime}\left(  i\right)  }$ are equal (indeed, the
permutation $\sigma^{\prime}=\sigma t_{u,v}$ can be obtained from $\sigma$ by
swapping the values at $u$ and $v$ (so that $\sigma^{\prime}\left(  u\right)
=\sigma\left(  v\right)  $ and $\sigma^{\prime}\left(  v\right)
=\sigma\left(  u\right)  $); thus, the products $\prod\limits_{i\in\left[
n\right]  \setminus I}y_{\sigma\left(  i\right)  }$ and $\prod\limits_{i\in
\left[  n\right]  \setminus I}y_{\sigma^{\prime}\left(  i\right)  }$ differ
only in the order in which their factors $y_{\sigma\left(  u\right)  }$ and
$y_{\sigma\left(  v\right)  }$ appear in them\footnote{Both factors
$y_{\sigma\left(  u\right)  }$ and $y_{\sigma\left(  v\right)  }$ do indeed
appear in these products, since $u$ and $v$ belong to $\left[  n\right]
\setminus I$.}; hence, these products are equal, since $K$ is commutative).
Again, this shows that our pairing has the properties we want, and thus the
proof of (\ref{pf.prop.det.xi+yj.innersum}) is complete in Case 2.

Thus, (\ref{pf.prop.det.xi+yj.innersum}) is proved in both cases.]

Now, we can finish our computation of the original determinant:%
\[
\det\left(  \left(  x_{i}+y_{j}\right)  _{1\leq i\leq n,\ 1\leq j\leq
n}\right)  =\sum\limits_{I\subseteq\left[  n\right]  }\left(  \prod
\limits_{i\in I}x_{i}\right)  \underbrace{\sum_{\sigma\in S_{n}}\left(
-1\right)  ^{\sigma}\prod\limits_{i\in\left[  n\right]  \setminus I}%
y_{\sigma\left(  i\right)  }}_{\substack{=0\\\text{(by
(\ref{pf.prop.det.xi+yj.innersum}))}}}=0.
\]
This proves Proposition \ref{prop.det.xi+yj}.
\end{proof}

\subsubsection{Basic properties}

The examples above show that pedestrian proofs of facts about determinants
(using just the definition) are possible, but become cumbersome fairly
quickly. Fortunately, determinants have a lot of properties that, once proved,
provide new methods for computing determinants.

Let us first recall a few basic facts that should (ideally) be known from a
good course on linear algebra. Recall that the transpose of a matrix $A$ is
denoted by $A^{T}$.

\begin{theorem}
[Transposes preserve determinants]\label{thm.det.transp}Let $n\in\mathbb{N}$.
If $A\in K^{n\times n}$ is any $n\times n$-matrix, then $\det\left(
A^{T}\right)  =\det A$.
\end{theorem}

\begin{proof}
See \cite[Corollary B.16]{Strick13} or \cite[Exercise 6.4]{detnotes} or
\cite[\S 5.3.2]{Laue-det}.
\end{proof}

\begin{theorem}
[Determinants of triangular matrices]\label{thm.det.triang}Let $n\in
\mathbb{N}$. Let $A\in K^{n\times n}$ be a triangular (i.e., lower-triangular
or upper-triangular) $n\times n$-matrix. Then, the determinant of the matrix
$A$ is the product of its diagonal entries. That is,%
\[
\det A=A_{1,1}A_{2,2}\cdots A_{n,n}.
\]

\end{theorem}

\begin{proof}
See \cite[Proposition B.11]{Strick13} or \cite[Exercise 6.3 and the paragraph
after Exercise 6.4]{detnotes}.
\end{proof}

As a consequence of Theorem \ref{thm.det.triang}, we see that the determinant
of a diagonal matrix is the product of its diagonal entries (since any
diagonal matrix is triangular).

\begin{theorem}
[Row operation properties]\label{thm.det.rowop}Let $n\in\mathbb{N}$. Let $A\in
K^{n\times n}$ be an $n\times n$-matrix. Then: \medskip

\textbf{(a)} If we swap two rows of $A$, then $\det A$ gets multiplied by
$-1$. \medskip

\textbf{(b)} If $A$ has a zero row (i.e., a row that consists entirely of
zeroes), then $\det A=0$. \medskip

\textbf{(c)} If $A$ has two equal rows, then $\det A=0$. \medskip

\textbf{(d)} Let $\lambda\in K$. If we multiply a row of $A$ by $\lambda$
(that is, we multiply all entries of this one row by $\lambda$, while leaving
all other entries of $A$ unchanged), then $\det A$ gets multiplied by
$\lambda$. \medskip

\textbf{(e)} If we add a row of $A$ to another row of $A$ (that is, we add
each entry of the former row to the corresponding entry of the latter), then
$\det A$ stays unchanged. \medskip

\textbf{(f)} Let $\lambda\in K$. If we add $\lambda$ times a row of $A$ to
another row of $A$ (that is, we add $\lambda$ times each entry of the former
row to the corresponding entry of the latter), then $\det A$ stays unchanged.
\medskip

\textbf{(g)} Let $B,C\in K^{n\times n}$ be two further $n\times n$-matrices.
Let $k\in\left[  n\right]  $. Assume that%
\[
\left(  \text{the }k\text{-th row of }C\right)  =\left(  \text{the }k\text{-th
row of }A\right)  +\left(  \text{the }k\text{-th row of }B\right)  ,
\]
whereas each $i\neq k$ satisfies%
\[
\left(  \text{the }i\text{-th row of }C\right)  =\left(  \text{the }i\text{-th
row of }A\right)  =\left(  \text{the }i\text{-th row of }B\right)  .
\]
Then,%
\[
\det C=\det A+\det B.
\]

\end{theorem}

\begin{example}
Let us see what Theorem \ref{thm.det.rowop} is saying in some particular cases
(specifically, for $3\times3$-matrices): \medskip

\textbf{(a)} One instance of Theorem \ref{thm.det.rowop} \textbf{(a)} is
\[
\det\left(
\begin{array}
[c]{ccc}%
a & b & c\\
a^{\prime\prime} & b^{\prime\prime} & c^{\prime\prime}\\
a^{\prime} & b^{\prime} & c^{\prime}%
\end{array}
\right)  =-\det\left(
\begin{array}
[c]{ccc}%
a & b & c\\
a^{\prime} & b^{\prime} & c^{\prime}\\
a^{\prime\prime} & b^{\prime\prime} & c^{\prime\prime}%
\end{array}
\right)  .
\]

\textbf{(b)} One instance of Theorem \ref{thm.det.rowop} \textbf{(b)} is
\[
\det\left(
\begin{array}
[c]{ccc}%
a & b & c\\
0 & 0 & 0\\
a^{\prime\prime} & b^{\prime\prime} & c^{\prime\prime}%
\end{array}
\right)  =0.
\]

\textbf{(c)} One instance of Theorem \ref{thm.det.rowop} \textbf{(c)} is
\[
\det\left(
\begin{array}
[c]{ccc}%
a & b & c\\
a^{\prime} & b^{\prime} & c^{\prime}\\
a & b & c
\end{array}
\right)  =0.
\]

\textbf{(d)} One instance of Theorem \ref{thm.det.rowop} \textbf{(d)} is
\[
\det\left(
\begin{array}
[c]{ccc}%
a & b & c\\
\lambda a^{\prime} & \lambda b^{\prime} & \lambda c^{\prime}\\
a^{\prime\prime} & b^{\prime\prime} & c^{\prime\prime}%
\end{array}
\right)  =\lambda\det\left(
\begin{array}
[c]{ccc}%
a & b & c\\
a^{\prime} & b^{\prime} & c^{\prime}\\
a^{\prime\prime} & b^{\prime\prime} & c^{\prime\prime}%
\end{array}
\right)  .
\]

\textbf{(e)} One instance of Theorem \ref{thm.det.rowop} \textbf{(e)} is
\[
\det\left(
\begin{array}
[c]{ccc}%
a & b & c\\
a^{\prime}+a^{\prime\prime} & b^{\prime}+b^{\prime\prime} & c^{\prime
}+c^{\prime\prime}\\
a^{\prime\prime} & b^{\prime\prime} & c^{\prime\prime}%
\end{array}
\right)  =\det\left(
\begin{array}
[c]{ccc}%
a & b & c\\
a^{\prime} & b^{\prime} & c^{\prime}\\
a^{\prime\prime} & b^{\prime\prime} & c^{\prime\prime}%
\end{array}
\right)  .
\]

\textbf{(f)} One instance of Theorem \ref{thm.det.rowop} \textbf{(f)} is
\[
\det\left(
\begin{array}
[c]{ccc}%
a & b & c\\
a^{\prime}+\lambda a^{\prime\prime} & b^{\prime}+\lambda b^{\prime\prime} &
c^{\prime}+\lambda c^{\prime\prime}\\
a^{\prime\prime} & b^{\prime\prime} & c^{\prime\prime}%
\end{array}
\right)  =\det\left(
\begin{array}
[c]{ccc}%
a & b & c\\
a^{\prime} & b^{\prime} & c^{\prime}\\
a^{\prime\prime} & b^{\prime\prime} & c^{\prime\prime}%
\end{array}
\right)  .
\]

\textbf{(g)} One instance of Theorem \ref{thm.det.rowop} \textbf{(g)} is%
\[
\det\underbrace{\left(
\begin{array}
[c]{ccc}%
a & b & c\\
d+d^{\prime} & e+e^{\prime} & f+f^{\prime}\\
g & h & i
\end{array}
\right)  }_{\text{this is }C}=\det\underbrace{\left(
\begin{array}
[c]{ccc}%
a & b & c\\
d & e & f\\
g & h & i
\end{array}
\right)  }_{\text{this is }A}+\det\underbrace{\left(
\begin{array}
[c]{ccc}%
a & b & c\\
d^{\prime} & e^{\prime} & f^{\prime}\\
g & h & i
\end{array}
\right)  }_{\text{this is }B}.
\]
(Specifically, this is the particular case of Theorem \ref{thm.det.rowop}
\textbf{(g)} for $n=3$ and $k=2$.)
\end{example}

Parts \textbf{(b)}, \textbf{(d)} and \textbf{(g)} of Theorem
\ref{thm.det.rowop} are commonly summarized under the mantle of
\textquotedblleft\emph{multilinearity of the determinant}\textquotedblright%
\ or \textquotedblleft\emph{linearity of the determinant in the }$k$\emph{-th
row}\textquotedblright. In fact, they say that (for any given $n\in\mathbb{N}$
and $k\in\left[  n\right]  $) if we hold all rows other than the $k$-th row of
an $n\times n$-matrix $A$ fixed, then $\det A$ depends $K$-linearly on the
$k$-th row of $A$.

\begin{proof}
[Proof of Theorem \ref{thm.det.rowop}.]\textbf{(a)} See \cite[Exercise 6.7
\textbf{(a)}]{detnotes}. This is also a particular case of \cite[Corollary
B.19]{Strick13}. \medskip

\textbf{(b)} See \cite[Exercise 6.7 \textbf{(c)}]{detnotes}. This is also
near-obvious from Definition \ref{def.det.det}. \medskip

\textbf{(c)} See \cite[Exercise 6.7 \textbf{(e)}]{detnotes} or \cite[\S 5.3.3,
property (iii)]{Laue-det} or \cite[2019-10-23 blackboard notes, Theorem
1.3.3]{19fla}.\footnote{\textit{Warning:} Several authors claim to give an
easy proof of part \textbf{(c)} using part \textbf{(a)}. This
\textquotedblleft proof\textquotedblright\ goes as follows: If $A$ has two
equal rows, then swapping these rows leaves $A$ unchanged, but (because of
Theorem \ref{thm.det.rowop}) flips the sign of $\det A$. Hence, in this case,
we have $\det A=-\det A$, so that $2\det A=0$ and therefore $\det A=0$, right?
Not so fast! In order to obtain $\det A=0$ from $2\det A=0$, we need the
element $2$ of $K$ to be invertible or at least be a non-zero-divisor (since
we have to divide by $2$). This is true when $K$ is one of the
\textquotedblleft high-school rings\textquotedblright\ $\mathbb{Z}$,
$\mathbb{Q}$, $\mathbb{R}$ and $\mathbb{C}$, but it is not true when $K$ is
the field $\mathbb{F}_{2}$ with $2$ elements (or, more generally, any field of
characteristic $2$). This slick argument can be salvaged, but in the form just
given it is incomplete.} \medskip

\textbf{(d)} See \cite[Exercise 6.7 \textbf{(g)}]{detnotes} or \cite[\S 5.3.3,
property (ii)]{Laue-det}. This is also a particular case of \cite[Corollary
B.19]{Strick13}. \medskip

\textbf{(f)} See \cite[Exercise 6.8 \textbf{(a)}]{detnotes}. This is also a
particular case of \cite[Corollary B.19]{Strick13}. \medskip

\textbf{(e)} This is the particular case of part \textbf{(f)} for $\lambda=1$.
\medskip

\textbf{(g)} See \cite[Exercise 6.7 \textbf{(i)}]{detnotes} or \cite[\S 5.3.3,
property (i)]{Laue-det} or \cite[2019-10-30 blackboard notes, Theorem
1.2.3]{19fla}.
\end{proof}

\begin{theorem}
[Column operation properties]\label{thm.det.colop}Theorem \ref{thm.det.rowop}
also holds if we replace \textquotedblleft row\textquotedblright\ by
\textquotedblleft column\textquotedblright\ throughout it.
\end{theorem}

\begin{proof}
Theorem \ref{thm.det.transp} shows that the determinant of a matrix does not
change when we replace it by its transpose; however, the rows of this
transpose $A^{T}$ are the transposes of the columns of $A$. Thus, Theorem
\ref{thm.det.colop} follows by applying Theorem \ref{thm.det.rowop} to the
transposes of all the matrices involved. (See \cite[Exercises 6.7 and
6.8]{detnotes} for the details.)
\end{proof}

\begin{corollary}
\label{cor.det.sig-row-col}Let $n\in\mathbb{N}$. Let $A\in K^{n\times n}$ and
$\tau\in S_{n}$. Then,%
\begin{equation}
\det\left(  \left(  A_{\tau\left(  i\right)  ,j}\right)  _{1\leq i\leq
n,\ 1\leq j\leq n}\right)  =\left(  -1\right)  ^{\tau}\cdot\det A
\label{eq.cor.det.sig-row-col.row}%
\end{equation}
and%
\begin{equation}
\det\left(  \left(  A_{i,\tau\left(  j\right)  }\right)  _{1\leq i\leq
n,\ 1\leq j\leq n}\right)  =\left(  -1\right)  ^{\tau}\cdot\det A.
\label{eq.cor.det.sig-row-col.col}%
\end{equation}

\end{corollary}

In words: When we permute the rows or the columns of a matrix, its determinant
gets multiplied by the sign of the permutation.

\begin{proof}
[Proof of Corollary \ref{cor.det.sig-row-col}.]Let us first prove
(\ref{eq.cor.det.sig-row-col.col}).

The definition of $\det A$ yields%
\begin{align}
\det A  &  =\sum_{\sigma\in S_{n}}\left(  -1\right)  ^{\sigma}%
\underbrace{A_{1,\sigma\left(  1\right)  }A_{2,\sigma\left(  2\right)  }\cdots
A_{n,\sigma\left(  n\right)  }}_{=\prod_{i=1}^{n}A_{i,\sigma\left(  i\right)
}}\nonumber\\
&  =\sum_{\sigma\in S_{n}}\left(  -1\right)  ^{\sigma}\prod_{i=1}%
^{n}A_{i,\sigma\left(  i\right)  }. \label{pf.cor.det.sig-row-col.0}%
\end{align}

The definition of $\det\left(  \left(  A_{i,\tau\left(  j\right)  }\right)
_{1\leq i\leq n,\ 1\leq j\leq n}\right)  $ yields%
\begin{align}
\det\left(  \left(  A_{i,\tau\left(  j\right)  }\right)  _{1\leq i\leq
n,\ 1\leq j\leq n}\right)   &  =\sum_{\sigma\in S_{n}}\left(  -1\right)
^{\sigma}\underbrace{A_{1,\tau\left(  \sigma\left(  1\right)  \right)
}A_{2,\tau\left(  \sigma\left(  2\right)  \right)  }\cdots A_{n,\tau\left(
\sigma\left(  n\right)  \right)  }}_{=\prod_{i=1}^{n}A_{i,\tau\left(
\sigma\left(  i\right)  \right)  }}\nonumber\\
&  =\sum_{\sigma\in S_{n}}\left(  -1\right)  ^{\sigma}\prod_{i=1}^{n}%
A_{i,\tau\left(  \sigma\left(  i\right)  \right)  }.
\label{pf.cor.det.sig-row-col.1}%
\end{align}

However, $S_{n}$ is a group. Thus, the map%
\begin{align*}
S_{n}  &  \rightarrow S_{n},\\
\sigma &  \mapsto\tau^{-1}\sigma
\end{align*}
is a bijection. Hence, we can substitute $\tau^{-1}\sigma$ for $\sigma$ in the
sum $\sum_{\sigma\in S_{n}}\left(  -1\right)  ^{\sigma}\prod_{i=1}%
^{n}A_{i,\tau\left(  \sigma\left(  i\right)  \right)  }$. Thus, we obtain%
\begin{align*}
&  \sum_{\sigma\in S_{n}}\left(  -1\right)  ^{\sigma}\prod_{i=1}^{n}%
A_{i,\tau\left(  \sigma\left(  i\right)  \right)  }\\
&  =\sum_{\sigma\in S_{n}}\underbrace{\left(  -1\right)  ^{\tau^{-1}\sigma}%
}_{\substack{=\left(  -1\right)  ^{\tau^{-1}}\cdot\left(  -1\right)  ^{\sigma
}\\\text{(by Proposition \ref{prop.perm.sign.props} \textbf{(d)}%
,}\\\text{applied to }\tau^{-1}\text{ and }\sigma\text{ instead of }%
\sigma\text{ and }\tau\text{)}}}\prod_{i=1}^{n}\underbrace{A_{i,\tau\left(
\left(  \tau^{-1}\sigma\right)  \left(  i\right)  \right)  }}%
_{\substack{=A_{i,\sigma\left(  i\right)  }\\\text{(since }\tau\left(  \left(
\tau^{-1}\sigma\right)  \left(  i\right)  \right)  =\left(  \tau\tau
^{-1}\sigma\right)  \left(  i\right)  =\sigma\left(  i\right)
\\\text{(because }\tau\tau^{-1}\sigma=\sigma\text{))}}}\\
&  =\sum_{\sigma\in S_{n}}\left(  -1\right)  ^{\tau^{-1}}\cdot\left(
-1\right)  ^{\sigma}\prod_{i=1}^{n}A_{i,\sigma\left(  i\right)  }%
=\underbrace{\left(  -1\right)  ^{\tau^{-1}}}_{\substack{=\left(  -1\right)
^{\tau}\\\text{(by Proposition \ref{prop.perm.sign.props} \textbf{(f)}%
,}\\\text{applied to }\tau^{-1}\text{ instead of }\sigma\text{)}}%
}\cdot\underbrace{\sum_{\sigma\in S_{n}}\left(  -1\right)  ^{\sigma}%
\prod_{i=1}^{n}A_{i,\sigma\left(  i\right)  }}_{\substack{=\det A\\\text{(by
(\ref{pf.cor.det.sig-row-col.0}))}}}\\
&  =\left(  -1\right)  ^{\tau}\cdot\det A.
\end{align*}
In view of this, we can rewrite (\ref{pf.cor.det.sig-row-col.1}) as%
\[
\det\left(  \left(  A_{i,\tau\left(  j\right)  }\right)  _{1\leq i\leq
n,\ 1\leq j\leq n}\right)  =\left(  -1\right)  ^{\tau}\cdot\det A.
\]
This proves (\ref{eq.cor.det.sig-row-col.col}).

Now, we can easily derive (\ref{eq.cor.det.sig-row-col.row}) by using the
transpose. Indeed, applying (\ref{eq.cor.det.sig-row-col.col}) to $A^{T}$
instead of $A$, we obtain%
\begin{align}
\det\left(  \left(  \left(  A^{T}\right)  _{i,\tau\left(  j\right)  }\right)
_{1\leq i\leq n,\ 1\leq j\leq n}\right)   &  =\left(  -1\right)  ^{\tau}%
\cdot\underbrace{\det\left(  A^{T}\right)  }_{\substack{=\det A\\\text{(by
Theorem \ref{thm.det.transp})}}}\nonumber\\
&  =\left(  -1\right)  ^{\tau}\cdot\det A. \label{pf.cor.det.sig-row-col.5}%
\end{align}
However, each $\left(  i,j\right)  \in\left[  n\right]  ^{2}$ satisfies%
\[
\left(  A^{T}\right)  _{i,\tau\left(  j\right)  }=A_{\tau\left(  j\right)
,i}\ \ \ \ \ \ \ \ \ \ \left(  \text{by the definition of }A^{T}\right)  .
\]
Thus,%
\[
\left(  \left(  A^{T}\right)  _{i,\tau\left(  j\right)  }\right)  _{1\leq
i\leq n,\ 1\leq j\leq n}=\left(  A_{\tau\left(  j\right)  ,i}\right)  _{1\leq
i\leq n,\ 1\leq j\leq n}=\left(  \left(  A_{\tau\left(  i\right)  ,j}\right)
_{1\leq i\leq n,\ 1\leq j\leq n}\right)  ^{T}%
\]
(again by the definition of the transpose). Therefore,%
\begin{align*}
\det\left(  \left(  \left(  A^{T}\right)  _{i,\tau\left(  j\right)  }\right)
_{1\leq i\leq n,\ 1\leq j\leq n}\right)   &  =\det\left(  \left(  \left(
A_{\tau\left(  i\right)  ,j}\right)  _{1\leq i\leq n,\ 1\leq j\leq n}\right)
^{T}\right) \\
&  =\det\left(  \left(  A_{\tau\left(  i\right)  ,j}\right)  _{1\leq i\leq
n,\ 1\leq j\leq n}\right)
\end{align*}
(by Theorem \ref{thm.det.transp}, applied to $\left(  A_{\tau\left(  i\right)
,j}\right)  _{1\leq i\leq n,\ 1\leq j\leq n}$ instead of $A$). Hence,%
\[
\det\left(  \left(  A_{\tau\left(  i\right)  ,j}\right)  _{1\leq i\leq
n,\ 1\leq j\leq n}\right)  =\det\left(  \left(  \left(  A^{T}\right)
_{i,\tau\left(  j\right)  }\right)  _{1\leq i\leq n,\ 1\leq j\leq n}\right)
=\left(  -1\right)  ^{\tau}\cdot\det A
\]
(by (\ref{pf.cor.det.sig-row-col.5})). This proves
(\ref{eq.cor.det.sig-row-col.row}). Thus, the proof of Corollary
\ref{cor.det.sig-row-col} is complete.
\end{proof}

The following is probably the most remarkable property of determinants:

\begin{theorem}
[Multiplicativity of the determinant]\label{thm.det.detAB}Let $n\in\mathbb{N}%
$. Let $A,B\in K^{n\times n}$ be two $n\times n$-matrices. Then,%
\[
\det\left(  AB\right)  =\det A\cdot\det B.
\]

\end{theorem}

\begin{proof}
See \cite[Theorem B.17]{Strick13} or \cite[Theorem 6.23]{detnotes} or
\cite[\S 5]{Zeilbe} or \cite[Lemma 4.7.5 (1)]{Ford21} or \cite[Theorem
5.7]{Laue-det}. (Many of these proofs are at least partly combinatorial, but
the one in \cite[\S 5]{Zeilbe} is fully so, constructing quite explicitly a
sign-reversing involution.)
\end{proof}

As an application of these properties, let us reprove Proposition
\ref{prop.det.xiyj} and Proposition \ref{prop.det.xi+yj}:

\begin{proof}
[Second proof of Proposition \ref{prop.det.xiyj}.]Define two $n\times
n$-matrices%
\[
A:=\left(
\begin{array}
[c]{ccccc}%
x_{1} & 0 & 0 & \cdots & 0\\
x_{2} & 0 & 0 & \cdots & 0\\
x_{3} & 0 & 0 & \cdots & 0\\
\vdots & \vdots & \vdots & \ddots & \vdots\\
x_{n} & 0 & 0 & \cdots & 0
\end{array}
\right)  \ \ \ \ \ \ \ \ \ \ \text{and}\ \ \ \ \ \ \ \ \ \ B:=\left(
\begin{array}
[c]{ccccc}%
y_{1} & y_{2} & y_{3} & \cdots & y_{n}\\
0 & 0 & 0 & \cdots & 0\\
0 & 0 & 0 & \cdots & 0\\
\vdots & \vdots & \vdots & \ddots & \vdots\\
0 & 0 & 0 & \cdots & 0
\end{array}
\right)  .
\]
(Only the first column of $A$ and the first row of $B$ have any nonzero entries.)

Now,
\begin{align*}
\left(  x_{i}y_{j}\right)  _{1\leq i\leq n,\ 1\leq j\leq n}  &  =\left(
\begin{array}
[c]{cccc}%
x_{1}y_{1} & x_{1}y_{2} & \cdots & x_{1}y_{n}\\
x_{2}y_{1} & x_{2}y_{2} & \cdots & x_{2}y_{n}\\
\vdots & \vdots & \ddots & \vdots\\
x_{n}y_{1} & x_{n}y_{2} & \cdots & x_{n}y_{n}%
\end{array}
\right) \\
&  =\underbrace{\left(
\begin{array}
[c]{ccccc}%
x_{1} & 0 & 0 & \cdots & 0\\
x_{2} & 0 & 0 & \cdots & 0\\
x_{3} & 0 & 0 & \cdots & 0\\
\vdots & \vdots & \vdots & \ddots & \vdots\\
x_{n} & 0 & 0 & \cdots & 0
\end{array}
\right)  }_{=A}\underbrace{\left(
\begin{array}
[c]{ccccc}%
y_{1} & y_{2} & y_{3} & \cdots & y_{n}\\
0 & 0 & 0 & \cdots & 0\\
0 & 0 & 0 & \cdots & 0\\
\vdots & \vdots & \vdots & \ddots & \vdots\\
0 & 0 & 0 & \cdots & 0
\end{array}
\right)  }_{=B}=AB,
\end{align*}
so that%
\[
\det\left(  \left(  x_{i}y_{j}\right)  _{1\leq i\leq n,\ 1\leq j\leq
n}\right)  =\det\left(  AB\right)  =\det A\cdot\det B
\]
(by Theorem \ref{thm.det.detAB}). However, the matrix $A$ has a zero column
(since $n\geq2$), and thus satisfies $\det A=0$ (by Theorem
\ref{thm.det.colop} \textbf{(b)}\footnote{Of course, by \textquotedblleft
Theorem \ref{thm.det.colop} \textbf{(b)}\textquotedblright, we mean
\textquotedblleft the analogue of Theorem \ref{thm.det.rowop} \textbf{(b)} for
columns instead of rows\textquotedblright.}). Hence,%
\[
\det\left(  \left(  x_{i}y_{j}\right)  _{1\leq i\leq n,\ 1\leq j\leq
n}\right)  =\underbrace{\det A}_{=0}\cdot\,\det B=0.
\]
Thus, Proposition \ref{prop.det.xiyj} is proven again.
\end{proof}

\begin{proof}
[Second proof of Proposition \ref{prop.det.xi+yj}.]Define two $n\times
n$-matrices%
\[
A:=\left(
\begin{array}
[c]{ccccc}%
x_{1} & 1 & 0 & \cdots & 0\\
x_{2} & 1 & 0 & \cdots & 0\\
x_{3} & 1 & 0 & \cdots & 0\\
\vdots & \vdots & \vdots & \ddots & \vdots\\
x_{n} & 1 & 0 & \cdots & 0
\end{array}
\right)  \ \ \ \ \ \ \ \ \ \ \text{and}\ \ \ \ \ \ \ \ \ \ B:=\left(
\begin{array}
[c]{ccccc}%
1 & 1 & 1 & \cdots & 1\\
y_{1} & y_{2} & y_{3} & \cdots & y_{n}\\
0 & 0 & 0 & \cdots & 0\\
\vdots & \vdots & \vdots & \ddots & \vdots\\
0 & 0 & 0 & \cdots & 0
\end{array}
\right)  .
\]
(Only the first two columns of $A$ and the first two rows of $B$ have any
nonzero entries.)

Now,
\begin{align*}
\left(  x_{i}+y_{j}\right)  _{1\leq i\leq n,\ 1\leq j\leq n}  &  =\left(
\begin{array}
[c]{cccc}%
x_{1}+y_{1} & x_{1}+y_{2} & \cdots & x_{1}+y_{n}\\
x_{2}+y_{1} & x_{2}+y_{2} & \cdots & x_{2}+y_{n}\\
\vdots & \vdots & \ddots & \vdots\\
x_{n}+y_{1} & x_{n}+y_{2} & \cdots & x_{n}+y_{n}%
\end{array}
\right) \\
&  =\underbrace{\left(
\begin{array}
[c]{ccccc}%
x_{1} & 1 & 0 & \cdots & 0\\
x_{2} & 1 & 0 & \cdots & 0\\
x_{3} & 1 & 0 & \cdots & 0\\
\vdots & \vdots & \vdots & \ddots & \vdots\\
x_{n} & 1 & 0 & \cdots & 0
\end{array}
\right)  }_{=A}\underbrace{\left(
\begin{array}
[c]{ccccc}%
1 & 1 & 1 & \cdots & 1\\
y_{1} & y_{2} & y_{3} & \cdots & y_{n}\\
0 & 0 & 0 & \cdots & 0\\
\vdots & \vdots & \vdots & \ddots & \vdots\\
0 & 0 & 0 & \cdots & 0
\end{array}
\right)  }_{=B}=AB,
\end{align*}
so that%
\[
\det\left(  \left(  x_{i}+y_{j}\right)  _{1\leq i\leq n,\ 1\leq j\leq
n}\right)  =\det\left(  AB\right)  =\det A\cdot\det B
\]
(by Theorem \ref{thm.det.detAB}). However, the matrix $A$ has a zero column
(since $n\geq3$), and thus satisfies $\det A=0$ (by Theorem
\ref{thm.det.colop} \textbf{(b)}). Hence,%
\[
\det\left(  \left(  x_{i}+y_{j}\right)  _{1\leq i\leq n,\ 1\leq j\leq
n}\right)  =\underbrace{\det A}_{=0}\cdot\,\det B=0.
\]
Thus, Proposition \ref{prop.det.xi+yj} is proven again.
\end{proof}

The following fact follows equally easily from everything we know so far about determinants:

\begin{corollary}
\label{cor.det.scale-row-col}Let $n\in\mathbb{N}$. Let $A\in K^{n\times n}$
and $d_{1},d_{2},\ldots,d_{n}\in K$. Then,%
\begin{equation}
\det\left(  \left(  d_{i}A_{i,j}\right)  _{1\leq i\leq n,\ 1\leq j\leq
n}\right)  =d_{1}d_{2}\cdots d_{n}\cdot\det A
\label{eq.cor.det.scale-row-col.row}%
\end{equation}
and%
\begin{equation}
\det\left(  \left(  d_{j}A_{i,j}\right)  _{1\leq i\leq n,\ 1\leq j\leq
n}\right)  =d_{1}d_{2}\cdots d_{n}\cdot\det A.
\label{eq.cor.det.scale-row-col.col}%
\end{equation}

\end{corollary}

\begin{proof}
[First proof of Corollary \ref{cor.det.scale-row-col} (sketched).]The matrix
$\left(  d_{i}A_{i,j}\right)  _{1\leq i\leq n,\ 1\leq j\leq n}$ is obtained
from the matrix $A$ by multiplying the $1$-st row by $d_{1}$, multiplying the
$2$-nd row by $d_{2}$, multiplying the $3$-rd row by $d_{3}$, and so on.
Theorem \ref{thm.det.rowop} \textbf{(d)} (applied repeatedly -- once for each
row) shows that these multiplications result in the determinant of $A$ getting
multiplied by $d_{1},d_{2},\ldots,d_{n}$ (in succession). Hence,%
\[
\det\left(  \left(  d_{i}A_{i,j}\right)  _{1\leq i\leq n,\ 1\leq j\leq
n}\right)  =d_{1}d_{2}\cdots d_{n}\cdot\det A.
\]
Thus, (\ref{eq.cor.det.scale-row-col.row}) is proved. The proof of
(\ref{eq.cor.det.scale-row-col.col}) is analogous (using columns instead of
rows). Corollary \ref{cor.det.scale-row-col} is proven.
\end{proof}

\begin{proof}
[Second proof of Corollary \ref{cor.det.scale-row-col} (sketched).]Let $D$ be
the diagonal matrix
\[
\left(
\begin{array}
[c]{cccc}%
d_{1} & 0 & \cdots & 0\\
0 & d_{2} & \cdots & 0\\
\vdots & \vdots & \ddots & \vdots\\
0 & 0 & \cdots & d_{n}%
\end{array}
\right)  \in K^{n\times n}.
\]
Then, $D$ is upper-triangular; hence, Theorem \ref{thm.det.triang} shows that
its determinant is $\det D=d_{1}d_{2}\cdots d_{n}$.

However, it is easy to see that $\left(  d_{i}A_{i,j}\right)  _{1\leq i\leq
n,\ 1\leq j\leq n}=DA$. Hence,%
\begin{align*}
\det\left(  \left(  d_{i}A_{i,j}\right)  _{1\leq i\leq n,\ 1\leq j\leq
n}\right)   &  =\det\left(  DA\right)  =\det D\cdot\det
A\ \ \ \ \ \ \ \ \ \ \left(  \text{by Theorem \ref{thm.det.detAB}}\right) \\
&  =d_{1}d_{2}\cdots d_{n}\cdot\det A\ \ \ \ \ \ \ \ \ \ \left(  \text{since
}\det D=d_{1}d_{2}\cdots d_{n}\right)  .
\end{align*}
This proves (\ref{eq.cor.det.scale-row-col.row}). A similar argument using
$\left(  d_{j}A_{i,j}\right)  _{1\leq i\leq n,\ 1\leq j\leq n}=AD$ proves
(\ref{eq.cor.det.scale-row-col.col}). Thus, Corollary
\ref{cor.det.scale-row-col} is proven again.
\end{proof}

\begin{proof}
[Third proof of Corollary \ref{cor.det.scale-row-col}.]Let us now proceed
completely elementarily. The definition of a determinant yields%
\begin{equation}
\det A=\sum_{\sigma\in S_{n}}\left(  -1\right)  ^{\sigma}A_{1,\sigma\left(
1\right)  }A_{2,\sigma\left(  2\right)  }\cdots A_{n,\sigma\left(  n\right)  }
\label{pf.cor.det.scale-row-col.3rd.1}%
\end{equation}
and%
\begin{align*}
\det\left(  \left(  d_{i}A_{i,j}\right)  _{1\leq i\leq n,\ 1\leq j\leq
n}\right)   &  =\sum_{\sigma\in S_{n}}\left(  -1\right)  ^{\sigma
}\underbrace{\left(  d_{1}A_{1,\sigma\left(  1\right)  }\right)  \left(
d_{2}A_{2,\sigma\left(  2\right)  }\right)  \cdots\left(  d_{n}A_{n,\sigma
\left(  n\right)  }\right)  }_{=\left(  d_{1}d_{2}\cdots d_{n}\right)  \left(
A_{1,\sigma\left(  1\right)  }A_{2,\sigma\left(  2\right)  }\cdots
A_{n,\sigma\left(  n\right)  }\right)  }\\
&  =\sum_{\sigma\in S_{n}}\left(  -1\right)  ^{\sigma}\left(  d_{1}d_{2}\cdots
d_{n}\right)  \left(  A_{1,\sigma\left(  1\right)  }A_{2,\sigma\left(
2\right)  }\cdots A_{n,\sigma\left(  n\right)  }\right) \\
&  =d_{1}d_{2}\cdots d_{n}\cdot\underbrace{\sum_{\sigma\in S_{n}}\left(
-1\right)  ^{\sigma}A_{1,\sigma\left(  1\right)  }A_{2,\sigma\left(  2\right)
}\cdots A_{n,\sigma\left(  n\right)  }}_{\substack{=\det A\\\text{(by
(\ref{pf.cor.det.scale-row-col.3rd.1}))}}}\\
&  =d_{1}d_{2}\cdots d_{n}\cdot\det A
\end{align*}
and%
\begin{align*}
&  \det\left(  \left(  d_{j}A_{i,j}\right)  _{1\leq i\leq n,\ 1\leq j\leq
n}\right) \\
&  =\sum_{\sigma\in S_{n}}\left(  -1\right)  ^{\sigma}\underbrace{\left(
d_{\sigma\left(  1\right)  }A_{1,\sigma\left(  1\right)  }\right)  \left(
d_{\sigma\left(  2\right)  }A_{2,\sigma\left(  2\right)  }\right)
\cdots\left(  d_{\sigma\left(  n\right)  }A_{n,\sigma\left(  n\right)
}\right)  }_{=\left(  d_{\sigma\left(  1\right)  }d_{\sigma\left(  2\right)
}\cdots d_{\sigma\left(  n\right)  }\right)  \left(  A_{1,\sigma\left(
1\right)  }A_{2,\sigma\left(  2\right)  }\cdots A_{n,\sigma\left(  n\right)
}\right)  }\\
&  =\sum_{\sigma\in S_{n}}\left(  -1\right)  ^{\sigma}\underbrace{\left(
d_{\sigma\left(  1\right)  }d_{\sigma\left(  2\right)  }\cdots d_{\sigma
\left(  n\right)  }\right)  }_{\substack{=d_{1}d_{2}\cdots d_{n}\\\text{(since
}\sigma\text{ is a permutation of the set }\left[  n\right]  \text{)}}}\left(
A_{1,\sigma\left(  1\right)  }A_{2,\sigma\left(  2\right)  }\cdots
A_{n,\sigma\left(  n\right)  }\right) \\
&  =\sum_{\sigma\in S_{n}}\left(  -1\right)  ^{\sigma}\left(  d_{1}d_{2}\cdots
d_{n}\right)  \left(  A_{1,\sigma\left(  1\right)  }A_{2,\sigma\left(
2\right)  }\cdots A_{n,\sigma\left(  n\right)  }\right) \\
&  =d_{1}d_{2}\cdots d_{n}\cdot\det A\ \ \ \ \ \ \ \ \ \ \left(  \text{as we
have seen above}\right)  .
\end{align*}
Once again, this proves Corollary \ref{cor.det.scale-row-col}.
\end{proof}

\subsubsection{Cauchy--Binet}

The multiplicativity of the determinant generalizes to non-square matrices $A$
and $B$, but the general statement is subtler and less famous:

\begin{theorem}
[Cauchy--Binet formula]\label{thm.det.CB}Let $n,m\in\mathbb{N}$. Let $A\in
K^{n\times m}$ be an $n\times m$-matrix, and let $B\in K^{m\times n}$ be an
$m\times n$-matrix. Then,%
\[
\det\left(  AB\right)  =\sum_{\substack{\left(  g_{1},g_{2},\ldots
,g_{n}\right)  \in\left[  m\right]  ^{n};\\g_{1}<g_{2}<\cdots<g_{n}}%
}\det\left(  \operatorname*{cols}\nolimits_{g_{1},g_{2},\ldots,g_{n}}A\right)
\cdot\det\left(  \operatorname*{rows}\nolimits_{g_{1},g_{2},\ldots,g_{n}%
}B\right)  .
\]
Here, we are using the following notations:

\begin{itemize}
\item We let $\operatorname*{cols}\nolimits_{g_{1},g_{2},\ldots,g_{n}}A$ be
the $n\times n$-matrix obtained from $A$ by removing all columns other than
the $g_{1}$-st, the $g_{2}$-nd, the $g_{3}$-rd, etc.. In other words,%
\[
\operatorname*{cols}\nolimits_{g_{1},g_{2},\ldots,g_{n}}A:=\left(  A_{i,g_{j}%
}\right)  _{1\leq i\leq n,\ 1\leq j\leq n}.
\]

\item We let $\operatorname*{rows}\nolimits_{g_{1},g_{2},\ldots,g_{n}}B$ be
the $n\times n$-matrix obtained from $B$ by removing all rows other than the
$g_{1}$-st, the $g_{2}$-nd, the $g_{3}$-rd, etc.. In other words,%
\[
\operatorname*{rows}\nolimits_{g_{1},g_{2},\ldots,g_{n}}B:=\left(  B_{g_{i}%
,j}\right)  _{1\leq i\leq n,\ 1\leq j\leq n}.
\]

\end{itemize}
\end{theorem}

Informally, we can rewrite the claim of Theorem \ref{thm.det.CB} as follows:%
\[
\det\left(  AB\right)  =\sum\det\left(  \text{some }n\text{ columns of
}A\right)  \cdot\det\left(  \text{the corresponding }n\text{ rows of
}B\right)  .
\]
The sum runs over all ways to form an $n\times n$-matrix by picking $n$
columns of $A$ (in increasing order, with no repetitions). The corresponding
$n$ rows of $B$ form an $n\times n$-matrix as well.

\begin{example}
Let $n=2$ and $m=3$, and let $A=\left(
\begin{array}
[c]{ccc}%
a & b & c\\
a^{\prime} & b^{\prime} & c^{\prime}%
\end{array}
\right)  $ and $B=\left(
\begin{array}
[c]{cc}%
x & x^{\prime}\\
y & y^{\prime}\\
z & z^{\prime}%
\end{array}
\right)  $. Then, Theorem \ref{thm.det.CB} yields%
\begin{align*}
\det\left(  AB\right)   &  =\sum_{\substack{\left(  g_{1},g_{2}\right)
\in\left[  3\right]  ^{2};\\g_{1}<g_{2}}}\det\left(  \operatorname*{cols}%
\nolimits_{g_{1},g_{2}}A\right)  \cdot\det\left(  \operatorname*{rows}%
\nolimits_{g_{1},g_{2}}B\right) \\
&  =\det\left(  \operatorname*{cols}\nolimits_{1,2}A\right)  \cdot\det\left(
\operatorname*{rows}\nolimits_{1,2}B\right) \\
&  \ \ \ \ \ \ \ \ \ \ +\det\left(  \operatorname*{cols}\nolimits_{1,3}%
A\right)  \cdot\det\left(  \operatorname*{rows}\nolimits_{1,3}B\right) \\
&  \ \ \ \ \ \ \ \ \ \ +\det\left(  \operatorname*{cols}\nolimits_{2,3}%
A\right)  \cdot\det\left(  \operatorname*{rows}\nolimits_{2,3}B\right) \\
&  \ \ \ \ \ \ \ \ \ \ \ \ \ \ \ \ \ \ \ \ \left(
\begin{array}
[c]{c}%
\text{since the only }2\text{-tuples }\left(  g_{1},g_{2}\right)  \in\left[
3\right]  ^{2}\\
\text{satisfying }g_{1}<g_{2}\text{ are }\left(  1,2\right)  \text{ and
}\left(  1,3\right)  \text{ and }\left(  2,3\right)
\end{array}
\right) \\
&  =\det\left(
\begin{array}
[c]{cc}%
a & b\\
a^{\prime} & b^{\prime}%
\end{array}
\right)  \cdot\det\left(
\begin{array}
[c]{cc}%
x & x^{\prime}\\
y & y^{\prime}%
\end{array}
\right) \\
&  \ \ \ \ \ \ \ \ \ \ +\det\left(
\begin{array}
[c]{cc}%
a & c\\
a^{\prime} & c^{\prime}%
\end{array}
\right)  \cdot\det\left(
\begin{array}
[c]{cc}%
x & x^{\prime}\\
z & z^{\prime}%
\end{array}
\right) \\
&  \ \ \ \ \ \ \ \ \ \ +\det\left(
\begin{array}
[c]{cc}%
b & c\\
b^{\prime} & c^{\prime}%
\end{array}
\right)  \cdot\det\left(
\begin{array}
[c]{cc}%
y & y^{\prime}\\
z & z^{\prime}%
\end{array}
\right) \\
&  =\left(  ab^{\prime}-ba^{\prime}\right)  \left(  xy^{\prime}-yx^{\prime
}\right)  +\left(  ac^{\prime}-ca^{\prime}\right)  \left(  xz^{\prime
}-zx^{\prime}\right) \\
&  \ \ \ \ \ \ \ \ \ \ +\left(  bc^{\prime}-cb^{\prime}\right)  \left(
yz^{\prime}-zy^{\prime}\right)  .
\end{align*}
You can check this against the equality%
\begin{align*}
\det\left(  AB\right)   &  =\left(  ax+by+cz\right)  \left(  a^{\prime
}x^{\prime}+b^{\prime}y^{\prime}+c^{\prime}z^{\prime}\right) \\
&  \ \ \ \ \ \ \ \ \ \ -\left(  ax^{\prime}+by^{\prime}+cz^{\prime}\right)
\left(  a^{\prime}x+b^{\prime}y+c^{\prime}z\right)  ,
\end{align*}
which is obtained by directly computing%
\[
AB=\left(
\begin{array}
[c]{ccc}%
a & b & c\\
a^{\prime} & b^{\prime} & c^{\prime}%
\end{array}
\right)  \left(
\begin{array}
[c]{cc}%
x & x^{\prime}\\
y & y^{\prime}\\
z & z^{\prime}%
\end{array}
\right)  =\left(
\begin{array}
[c]{cc}%
ax+by+cz & ax^{\prime}+by^{\prime}+cz^{\prime}\\
a^{\prime}x+b^{\prime}y+c^{\prime}z & a^{\prime}x^{\prime}+b^{\prime}%
y^{\prime}+c^{\prime}z^{\prime}%
\end{array}
\right)  .
\]

\end{example}

\begin{remark}
If $m<n$, then the claim of Theorem \ref{thm.det.CB} becomes%
\begin{align*}
\det\left(  AB\right)   &  =\sum_{\substack{\left(  g_{1},g_{2},\ldots
,g_{n}\right)  \in\left[  m\right]  ^{n};\\g_{1}<g_{2}<\cdots<g_{n}}%
}\det\left(  \operatorname*{cols}\nolimits_{g_{1},g_{2},\ldots,g_{n}}A\right)
\cdot\det\left(  \operatorname*{rows}\nolimits_{g_{1},g_{2},\ldots,g_{n}%
}B\right) \\
&  =\left(  \text{empty sum}\right) \\
&  \ \ \ \ \ \ \ \ \ \ \ \ \ \ \ \ \ \ \ \ \left(  \text{since the set
}\left[  m\right]  \text{ has no }n\text{ distinct elements}\right) \\
&  =0.
\end{align*}
When $K$ is a field, this can also be seen trivially from rank considerations
(to wit, the matrix $A$ has rank $\leq m<n$, and thus the product $AB$ has
rank $<n$ as well). When $K$ is not a field, the notion of a rank is not
available, so we do need Theorem \ref{thm.det.CB} to obtain this (although
there are ways around this).

If $m=n$, then the claim of Theorem \ref{thm.det.CB} becomes%
\begin{align*}
\det\left(  AB\right)   &  =\sum_{\substack{\left(  g_{1},g_{2},\ldots
,g_{n}\right)  \in\left[  n\right]  ^{n};\\g_{1}<g_{2}<\cdots<g_{n}}%
}\det\left(  \operatorname*{cols}\nolimits_{g_{1},g_{2},\ldots,g_{n}}A\right)
\cdot\det\left(  \operatorname*{rows}\nolimits_{g_{1},g_{2},\ldots,g_{n}%
}B\right) \\
&  =\det\underbrace{\left(  \operatorname*{cols}\nolimits_{1,2,\ldots
,n}A\right)  }_{=A}\cdot\det\underbrace{\left(  \operatorname*{rows}%
\nolimits_{1,2,\ldots,n}B\right)  }_{=B}\\
&  \ \ \ \ \ \ \ \ \ \ \ \ \ \ \ \ \ \ \ \ \left(
\begin{array}
[c]{c}%
\text{since the only }n\text{-tuple }\left(  g_{1},g_{2},\ldots,g_{n}\right)
\in\left[  n\right]  ^{n}\\
\text{satisfying }g_{1}<g_{2}<\cdots<g_{n}\text{ is }\left(  1,2,\ldots
,n\right)
\end{array}
\right) \\
&  =\det A\cdot\det B,
\end{align*}
so that we recover the multiplicativity of the determinant (Theorem
\ref{thm.det.detAB}).
\end{remark}

\begin{proof}
[Proof of Theorem \ref{thm.det.CB}.]See \cite[Theorem 6.32]{detnotes} or
\cite{Schwar16} or \cite[\S 1.2.3, Exercise 46]{Knuth-TAoCP1} or \cite[Theorem
9.53]{Loehr-BC} or \cite[Theorem 9.4]{Stanle18} or \cite[\S 2]{Zeng93} (a
fully combinatorial proof) or
\url{https://math.stackexchange.com/questions/3243063/} .
\end{proof}

\subsubsection{$\det\left(  A+B\right)  $}

So much for $\det\left(  AB\right)  $. What can we say about $\det\left(
A+B\right)  $ ? The answer is somewhat cumbersome, but still rather useful. We
need some notation to state it:

\begin{definition}
\label{def.det.sub}Let $n,m\in\mathbb{N}$. Let $A$ be an $n\times m$-matrix.

Let $U$ be a subset of $\left[  n\right]  $. Let $V$ be a subset of $\left[
m\right]  $.

Then, $\operatorname*{sub}\nolimits_{U}^{V}A$ is the $\left\vert U\right\vert
\times\left\vert V\right\vert $-matrix defined as follows:

Writing the two sets $U$ and $V$ as%
\[
U=\left\{  u_{1},u_{2},\ldots,u_{p}\right\}  \ \ \ \ \ \ \ \ \ \ \text{and}%
\ \ \ \ \ \ \ \ \ \ V=\left\{  v_{1},v_{2},\ldots,v_{q}\right\}
\]
with
\[
u_{1}<u_{2}<\cdots<u_{p}\ \ \ \ \ \ \ \ \ \ \text{and}%
\ \ \ \ \ \ \ \ \ \ v_{1}<v_{2}<\cdots<v_{q},
\]
we set%
\[
\operatorname*{sub}\nolimits_{U}^{V}A:=\left(  A_{u_{i},v_{j}}\right)  _{1\leq
i\leq p,\ 1\leq j\leq q}.
\]

Roughly speaking, $\operatorname*{sub}\nolimits_{U}^{V}A$ is the matrix
obtained from $A$ by focusing only on the $i$-th rows for $i\in U$ (that is,
removing all the other rows) and only on the $j$-th columns for $j\in V$ (that
is, removing all the other columns).

This matrix $\operatorname*{sub}\nolimits_{U}^{V}A$ is called the
\emph{submatrix of }$A$\emph{ obtained by restricting to the }$U$\emph{-rows
and the }$V$\emph{-columns}. If this matrix is square (i.e., if $\left\vert
U\right\vert =\left\vert V\right\vert $), then its determinant $\det\left(
\operatorname*{sub}\nolimits_{U}^{V}A\right)  $ is called a \emph{minor} of
$A$.
\end{definition}

\begin{example}
We have%
\[
\operatorname*{sub}\nolimits_{\left\{  1,2\right\}  }^{\left\{  1,3\right\}
}\left(
\begin{array}
[c]{ccc}%
a & b & c\\
a^{\prime} & b^{\prime} & c^{\prime}\\
a^{\prime\prime} & b^{\prime\prime} & c^{\prime\prime}%
\end{array}
\right)  =\left(
\begin{array}
[c]{cc}%
a & c\\
a^{\prime} & c^{\prime}%
\end{array}
\right)
\]
and%
\[
\operatorname*{sub}\nolimits_{\left\{  2\right\}  }^{\left\{  3\right\}
}\left(
\begin{array}
[c]{ccc}%
a & b & c\\
a^{\prime} & b^{\prime} & c^{\prime}\\
a^{\prime\prime} & b^{\prime\prime} & c^{\prime\prime}%
\end{array}
\right)  =\left(
\begin{array}
[c]{c}%
c^{\prime}%
\end{array}
\right)  .
\]

\end{example}

\begin{theorem}
\label{thm.det.det(A+B)}Let $n\in\mathbb{N}$. For any subset $I$ of $\left[
n\right]  $, we let $\widetilde{I}$ be the complement $\left[  n\right]
\setminus I$ of $I$. (For example, if $n=4$ and $I=\left\{  1,4\right\}  $,
then $\widetilde{I}=\left\{  2,3\right\}  $.)

For any finite set $S$ of integers, define $\operatorname*{sum}S:=\sum_{s\in
S}s$.

Let $A$ and $B$ be two $n\times n$-matrices in $K^{n\times n}$. Then,%
\[
\det\left(  A+B\right)  =\sum_{P\subseteq\left[  n\right]  }\ \ \sum
_{\substack{Q\subseteq\left[  n\right]  ;\\\left\vert P\right\vert =\left\vert
Q\right\vert }}\left(  -1\right)  ^{\operatorname*{sum}P+\operatorname*{sum}%
Q}\det\left(  \operatorname*{sub}\nolimits_{P}^{Q}A\right)  \cdot\det\left(
\operatorname*{sub}\nolimits_{\widetilde{P}}^{\widetilde{Q}}B\right)  .
\]

\end{theorem}

\begin{example}
For $n=2$, this is saying that%
\begin{align*}
\det\left(  A+B\right)   &  =\underbrace{\left(  -1\right)
^{\operatorname*{sum}\varnothing+\operatorname*{sum}\varnothing}}_{=\left(
-1\right)  ^{0+0}=1}\underbrace{\det\left(  \operatorname*{sub}%
\nolimits_{\varnothing}^{\varnothing}A\right)  }_{=\det\left(  {}\right)
=1}\cdot\underbrace{\det\left(  \operatorname*{sub}\nolimits_{\left\{
1,2\right\}  }^{\left\{  1,2\right\}  }B\right)  }_{=\det B}\\
&  \ \ \ \ \ \ \ \ \ \ +\underbrace{\left(  -1\right)  ^{\operatorname*{sum}%
\left\{  1\right\}  +\operatorname*{sum}\left\{  1\right\}  }}_{=\left(
-1\right)  ^{1+1}=1}\underbrace{\det\left(  \operatorname*{sub}%
\nolimits_{\left\{  1\right\}  }^{\left\{  1\right\}  }A\right)
}_{\substack{=A_{1,1}\\\text{(since }\operatorname*{sub}\nolimits_{\left\{
1\right\}  }^{\left\{  1\right\}  }A\text{ is}\\\text{the }1\times
1\text{-matrix }\left(
\begin{array}
[c]{c}%
A_{1,1}%
\end{array}
\right)  \text{)}}}\cdot\underbrace{\det\left(  \operatorname*{sub}%
\nolimits_{\left\{  2\right\}  }^{\left\{  2\right\}  }B\right)  }_{=B_{2,2}%
}\\
&  \ \ \ \ \ \ \ \ \ \ +\underbrace{\left(  -1\right)  ^{\operatorname*{sum}%
\left\{  1\right\}  +\operatorname*{sum}\left\{  2\right\}  }}_{=\left(
-1\right)  ^{1+2}=-1}\underbrace{\det\left(  \operatorname*{sub}%
\nolimits_{\left\{  1\right\}  }^{\left\{  2\right\}  }A\right)  }_{=A_{1,2}%
}\cdot\underbrace{\det\left(  \operatorname*{sub}\nolimits_{\left\{
2\right\}  }^{\left\{  1\right\}  }B\right)  }_{=B_{2,1}}\\
&  \ \ \ \ \ \ \ \ \ \ +\underbrace{\left(  -1\right)  ^{\operatorname*{sum}%
\left\{  2\right\}  +\operatorname*{sum}\left\{  1\right\}  }}_{=\left(
-1\right)  ^{2+1}=-1}\underbrace{\det\left(  \operatorname*{sub}%
\nolimits_{\left\{  2\right\}  }^{\left\{  1\right\}  }A\right)  }_{=A_{2,1}%
}\cdot\underbrace{\det\left(  \operatorname*{sub}\nolimits_{\left\{
1\right\}  }^{\left\{  2\right\}  }B\right)  }_{=B_{1,2}}\\
&  \ \ \ \ \ \ \ \ \ \ +\underbrace{\left(  -1\right)  ^{\operatorname*{sum}%
\left\{  2\right\}  +\operatorname*{sum}\left\{  2\right\}  }}_{=\left(
-1\right)  ^{2+2}=1}\underbrace{\det\left(  \operatorname*{sub}%
\nolimits_{\left\{  2\right\}  }^{\left\{  2\right\}  }A\right)  }_{=A_{2,2}%
}\cdot\underbrace{\det\left(  \operatorname*{sub}\nolimits_{\left\{
1\right\}  }^{\left\{  1\right\}  }B\right)  }_{=B_{1,1}}\\
&  \ \ \ \ \ \ \ \ \ \ +\underbrace{\left(  -1\right)  ^{\operatorname*{sum}%
\left\{  1,2\right\}  +\operatorname*{sum}\left\{  1,2\right\}  }}_{=\left(
-1\right)  ^{3+3}=1}\underbrace{\det\left(  \operatorname*{sub}%
\nolimits_{\left\{  1,2\right\}  }^{\left\{  1,2\right\}  }A\right)  }_{=\det
A}\cdot\underbrace{\det\left(  \operatorname*{sub}\nolimits_{\varnothing
}^{\varnothing}B\right)  }_{=\det\left(  {}\right)  =1}\\
&  =\det B+A_{1,1}B_{2,2}-A_{1,2}B_{2,1}-A_{2,1}B_{1,2}+A_{2,2}B_{1,1}+\det
A\\
&  =\det A+\det B-A_{1,2}B_{2,1}-A_{2,1}B_{1,2}+A_{1,1}B_{2,2}+A_{2,2}B_{1,1}.
\end{align*}

\end{example}

Theorem \ref{thm.det.det(A+B)} can be thought of as a kind of
\textquotedblleft binomial theorem\textquotedblright\ for determinants: On its
right hand side (for $n>0$) is a sum that contains both $\det A$ and $\det B$
as addends (in fact, $\det A$ is the addend for $P=Q=\left[  n\right]  $,
whereas $\det B$ is the addend for $P=Q=\varnothing$) as well as many
\textquotedblleft mixed\textquotedblright\ addends that contain both a part of
$A$ and a part of $B$.

\begin{proof}
[Proof of Theorem \ref{thm.det.det(A+B)}.]See \cite[Theorem 6.160]{detnotes}.
(Note that $\operatorname*{sub}\nolimits_{P}^{Q}A$ is called
$\operatorname*{sub}\nolimits_{w\left(  P\right)  }^{w\left(  Q\right)  }A$ in
\cite{detnotes}.) The main difficulty of the proof is bookkeeping; the
underlying idea is simple (expand everything and regroup). \medskip

\begin{fineprint}
Here is a rough outline of the argument. If $\sigma\in S_{n}$, and if $P$ is a
subset of $\left[  n\right]  $, then $\sigma\left(  P\right)  =\left\{
\sigma\left(  i\right)  \ \mid\ i\in P\right\}  $ is a subset of $\left[
n\right]  $ as well, and has the same size as $P$ (since $\sigma$ is a
permutation and therefore injective); thus, it satisfies $\left\vert
P\right\vert =\left\vert \sigma\left(  P\right)  \right\vert $.

The definition of $\det\left(  A+B\right)  $ yields%
\begin{align*}
\det\left(  A+B\right)   &  =\sum_{\sigma\in S_{n}}\left(  -1\right)
^{\sigma}\underbrace{\left(  A_{1,\sigma\left(  1\right)  }+B_{1,\sigma\left(
1\right)  }\right)  \left(  A_{2,\sigma\left(  2\right)  }+B_{2,\sigma\left(
2\right)  }\right)  \cdots\left(  A_{n,\sigma\left(  n\right)  }%
+B_{n,\sigma\left(  n\right)  }\right)  }_{\substack{=\prod_{i=1}^{n}\left(
A_{i,\sigma\left(  i\right)  }+B_{i,\sigma\left(  i\right)  }\right)
=\sum\limits_{P\subseteq\left[  n\right]  }\left(  \prod\limits_{i\in
P}A_{i,\sigma\left(  i\right)  }\right)  \left(  \prod\limits_{i\in
\widetilde{P}}B_{i,\sigma\left(  i\right)  }\right)  \\\text{(by
(\ref{eq.lem.prodrule.sum-ai-plus-bi.eq}), applied to }a_{i}=A_{i,\sigma
\left(  i\right)  }\text{ and }b_{i}=B_{i,\sigma\left(  i\right)  }\text{)}%
}}\\
&  =\sum_{\sigma\in S_{n}}\left(  -1\right)  ^{\sigma}\sum\limits_{P\subseteq
\left[  n\right]  }\left(  \prod\limits_{i\in P}A_{i,\sigma\left(  i\right)
}\right)  \left(  \prod\limits_{i\in\widetilde{P}}B_{i,\sigma\left(  i\right)
}\right) \\
&  =\sum\limits_{P\subseteq\left[  n\right]  }\ \ \sum_{\sigma\in S_{n}%
}\left(  -1\right)  ^{\sigma}\left(  \prod\limits_{i\in P}A_{i,\sigma\left(
i\right)  }\right)  \left(  \prod\limits_{i\in\widetilde{P}}B_{i,\sigma\left(
i\right)  }\right) \\
&  =\sum\limits_{P\subseteq\left[  n\right]  }\ \ \sum_{\substack{Q\subseteq
\left[  n\right]  ;\\\left\vert P\right\vert =\left\vert Q\right\vert
}}\ \ \sum_{\substack{\sigma\in S_{n};\\\sigma\left(  P\right)  =Q}}\left(
-1\right)  ^{\sigma}\left(  \prod\limits_{i\in P}A_{i,\sigma\left(  i\right)
}\right)  \left(  \prod\limits_{i\in\widetilde{P}}B_{i,\sigma\left(  i\right)
}\right)
\end{align*}
(here, we have split the inner sum according to the value of the subset
$\sigma\left(  P\right)  =\left\{  \sigma\left(  i\right)  \ \mid\ i\in
P\right\}  $, recalling that it satisfies $\left\vert P\right\vert =\left\vert
\sigma\left(  P\right)  \right\vert $).

Now, fix two subsets $P$ and $Q$ of $\left[  n\right]  $ satisfying
$\left\vert P\right\vert =\left\vert Q\right\vert $. Thus, $\left\vert
\widetilde{P}\right\vert =\left\vert \widetilde{Q}\right\vert $ as well. Write
the sets $P$, $Q$, $\widetilde{P}$ and $\widetilde{Q}$ as%
\begin{align*}
P  &  =\left\{  p_{1}<p_{2}<\cdots<p_{k}\right\}
\ \ \ \ \ \ \ \ \ \ \text{and}\ \ \ \ \ \ \ \ \ \ Q=\left\{  q_{1}%
<q_{2}<\cdots<q_{k}\right\}  \ \ \ \ \ \ \ \ \ \ \text{and}\\
\widetilde{P}  &  =\left\{  p_{1}^{\prime}<p_{2}^{\prime}<\cdots<p_{\ell
}^{\prime}\right\}  \ \ \ \ \ \ \ \ \ \ \text{and}%
\ \ \ \ \ \ \ \ \ \ \widetilde{Q}=\left\{  q_{1}^{\prime}<q_{2}^{\prime
}<\cdots<q_{\ell}^{\prime}\right\}  ,
\end{align*}
where the notation \textquotedblleft$U=\left\{  u_{1}<u_{2}<\cdots
<u_{a}\right\}  $\textquotedblright\ is just a shorthand way to say
\textquotedblleft$U=\left\{  u_{1},u_{2},\ldots,u_{a}\right\}  $ and
$u_{1}<u_{2}<\cdots<u_{a}$\textquotedblright\ (or, equivalently,
\textquotedblleft the elements of $U$ in strictly increasing order are
$u_{1},u_{2},\ldots,u_{a}$\textquotedblright). Now, for each permutation
$\sigma\in S_{n}$ satisfying $\sigma\left(  P\right)  =Q$, we see that:

\begin{itemize}
\item The elements $\sigma\left(  p_{1}\right)  ,\sigma\left(  p_{2}\right)
,\ldots,\sigma\left(  p_{k}\right)  $ are the elements $q_{1},q_{2}%
,\ldots,q_{k}$ in some order (since $\sigma\left(  P\right)  =Q$), and thus
there exists a unique permutation $\alpha\in S_{k}$ such that
\[
\sigma\left(  p_{i}\right)  =q_{\alpha\left(  i\right)  }%
\ \ \ \ \ \ \ \ \ \ \text{for each }i\in\left[  k\right]  .
\]
We denote this $\alpha$ by $\alpha_{\sigma}$.

\item The elements $\sigma\left(  p_{1}^{\prime}\right)  ,\sigma\left(
p_{2}^{\prime}\right)  ,\ldots,\sigma\left(  p_{\ell}^{\prime}\right)  $ are
the elements $q_{1}^{\prime},q_{2}^{\prime},\ldots,q_{\ell}^{\prime}$ in some
order (since $\sigma\left(  P\right)  =Q$ entails $\sigma\left(
\widetilde{P}\right)  =\widetilde{Q}$, because $\sigma$ is a permutation of
$\left[  n\right]  $), and thus there exists a unique permutation $\beta\in
S_{\ell}$ such that%
\[
\sigma\left(  p_{i}^{\prime}\right)  =q_{\beta\left(  i\right)  }^{\prime
}\ \ \ \ \ \ \ \ \ \ \text{for each }i\in\left[  \ell\right]  .
\]
We denote this $\beta$ by $\beta_{\sigma}$.
\end{itemize}

Thus, for each permutation $\sigma\in S_{n}$ satisfying $\sigma\left(
P\right)  =Q$, we have defined two permutations $\alpha_{\sigma}\in S_{k}$ and
$\beta_{\sigma}\in S_{\ell}$ that (roughly speaking) describe the actions of
$\sigma$ on the subsets $P$ and $\widetilde{P}$, respectively. It is not hard
to see that the map%
\begin{align*}
\left\{  \text{permutations }\sigma\in S_{n}\text{ satisfying }\sigma\left(
P\right)  =Q\right\}   &  \rightarrow S_{k}\times S_{\ell},\\
\sigma &  \mapsto\left(  \alpha_{\sigma},\beta_{\sigma}\right)
\end{align*}
is a bijection. Moreover, for any permutation $\sigma\in S_{n}$ satisfying
$\sigma\left(  P\right)  =Q$, we have
\begin{equation}
\left(  -1\right)  ^{\sigma}=\left(  -1\right)  ^{\operatorname*{sum}%
P+\operatorname*{sum}Q}\left(  -1\right)  ^{\alpha_{\sigma}}\left(  -1\right)
^{\beta_{\sigma}}. \label{pf.thm.det.det(A+B).signs-equal}%
\end{equation}
(Proving this is perhaps the least pleasant part of this proof, but it is pure
combinatorics. It is probably easiest to reduce this to the case when
$P=\left[  k\right]  $ and $Q=\left[  k\right]  $ by a reduction procedure
that involves multiplying $\sigma$ by $\operatorname*{sum}%
P+\operatorname*{sum}Q-2\cdot\left(  1+2+\cdots+k\right)  $ many
transpositions\footnote{Specifically: Recall the simple transpositions $s_{i}$
defined for all $i\in\left[  n-1\right]  $ in Definition \ref{def.perm.si}. By
replacing $\sigma$ by $\sigma s_{i}$ (for some $i\in\left[  n-1\right]  $), we
can swap two adjacent values of $\sigma$ (namely, $\sigma\left(  i\right)  $
and $\sigma\left(  i+1\right)  $). Furthermore, $\sigma\left(  P\right)  =Q$
implies $\left(  \sigma s_{i}\right)  \left(  s_{i}\left(  P\right)  \right)
=Q$ (since $s_{i}s_{i}=s_{i}^{2}=\operatorname*{id}$). Thus, the equality
$\sigma\left(  P\right)  =Q$ is preserved if we simultaneously replace
$\sigma$ by $\sigma s_{i}$ and replace $P$ by $s_{i}\left(  P\right)  $. Such
a simultaneous replacement will be called a \emph{shift}. Furthermore, if $i$
is chosen in such a way that $i\notin P$ and $i+1\in P$, then this shift will
be called a \emph{left shift}.
\par
Let us see what happens to $\sigma$, $P$, $\alpha_{\sigma}$ and $\beta
_{\sigma}$ when we perform a left shift. Indeed, consider a left shift which
replaces $\sigma$ by $\sigma s_{i}$ and replaces $P$ by $s_{i}\left(
P\right)  $, where $i\in\left[  n-1\right]  $ is chosen in such a way that
$i\notin P$ and $i+1\in P$. The set $s_{i}\left(  P\right)  $ is the set $P$
with the element $i+1$ replaced by $i$. Thus, $\operatorname*{sum}\left(
s_{i}\left(  P\right)  \right)  =\operatorname*{sum}P-1$. In other words, our
left shift has decremented $\operatorname*{sum}P$ by $1$. Thus, our left shift
has flipped the sign of $\left(  -1\right)  ^{\operatorname*{sum}%
P+\operatorname*{sum}Q}$ (since $\operatorname*{sum}Q$ obviously stays
unchanged). The permutation $\sigma$ has been replaced by $\sigma s_{i}$,
which is the same permutation as $\sigma$ but with the values $\sigma\left(
i\right)  $ and $\sigma\left(  i+1\right)  $ swapped. Since $i\notin P$ and
$i+1\in P$, this swap has not disturbed the relative order of the elements of
$P$ and $\widetilde{P}$ (but merely replaced $i+1$ by $i$ in $P$ and replaced
$i$ by $i+1$ in $\widetilde{P}$), so that the permutations $\alpha_{\sigma}$
and $\beta_{\sigma}$ have not changed. Our left shift thus has left
$\alpha_{\sigma}$ and $\beta_{\sigma}$ unchanged. It has, however, flipped the
sign of $\left(  -1\right)  ^{\sigma}$ (because $\left(  -1\right)  ^{\sigma
s_{i}}=\left(  -1\right)  ^{\sigma}\underbrace{\left(  -1\right)  ^{s_{i}}%
}_{=-1}=-\left(  -1\right)  ^{\sigma}$).
\par
Let us summarize: Each left shift leaves $\alpha_{\sigma}$ and $\beta_{\sigma
}$ unchanged, while flipping the signs of $\left(  -1\right)
^{\operatorname*{sum}P+\operatorname*{sum}Q}$ and $\left(  -1\right)
^{\sigma}$. However, by performing left shifts, we can move the smallest
element of $P$ one step to the left, or (if the smallest element of $P$ is
already $1$) we can move the second-smallest element of $P$ one step to the
left, or (if the two smallest elements of $P$ are already $1$ and $2$) we can
move the third-smallest element of $P$ one step to the left, and so on. After
sufficiently many such left shifts, the set $P$ will have become $\left[
k\right]  $, whereas $\alpha_{\sigma}$ and $\beta_{\sigma}$ will not have
changed (because we have seen above that each left shift leaves $\alpha
_{\sigma}$ and $\beta_{\sigma}$ unchanged). The total number of left shifts we
need for this is $\left(  p_{1}-1\right)  +\left(  p_{2}-2\right)
+\cdots+\left(  p_{k}-k\right)  =\operatorname*{sum}P-\left(  1+2+\cdots
+k\right)  $.
\par
Likewise, we can define \emph{left co-shifts}, which are operations similar to
left shifts but acting on the values rather than positions of $\sigma$ and
acting on $Q$ rather than $P$. Explicitly, a left co-shift replaces $\sigma$
by $s_{i}\sigma$ and replaces $Q$ by $s_{i}\left(  Q\right)  $, where
$i\in\left[  n-1\right]  $ is chosen such that $i\notin Q$ and $i+1\in Q$.
Again, we can see that each left co-shift leaves $\alpha_{\sigma}$ and
$\beta_{\sigma}$ unchanged, while flipping the signs of $\left(  -1\right)
^{\operatorname*{sum}P+\operatorname*{sum}Q}$ and $\left(  -1\right)
^{\sigma}$. After a sequence of $\operatorname*{sum}Q-\left(  1+2+\cdots
+k\right)  $ left co-shifts, the set $Q$ will have become $\left[  k\right]
$.
\par
Each left shift and each left co-shift multiplies the permutation $\sigma$ by
a transposition. Hence, after our $\operatorname*{sum}P-\left(  1+2+\cdots
+k\right)  $ left shifts and our $\operatorname*{sum}Q-\left(  1+2+\cdots
+k\right)  $ left co-shifts, we will have multiplied $\sigma$ by altogether%
\begin{align*}
&  \left(  \operatorname*{sum}P-\left(  1+2+\cdots+k\right)  \right)  +\left(
\operatorname*{sum}Q-\left(  1+2+\cdots+k\right)  \right) \\
&  =\operatorname*{sum}P+\operatorname*{sum}Q-2\cdot\left(  1+2+\cdots
+k\right)
\end{align*}
many transpositions. At the end of this process, we have $P=\left[  k\right]
$ and $Q=\left[  k\right]  $.}. Once we are in the case $P=\left[  k\right]  $
and $Q=\left[  k\right]  $, we can prove
(\ref{pf.thm.det.det(A+B).signs-equal}) by directly counting the inversions of
$\sigma$ (showing that $\ell\left(  \sigma\right)  =\ell\left(  \alpha
_{\sigma}\right)  +\ell\left(  \beta_{\sigma}\right)  $). Note that the proof
I give in \cite[Theorem 6.160]{detnotes} avoids proving
(\ref{pf.thm.det.det(A+B).signs-equal}), opting instead for an argument using
row permutations; see \cite[solution to Exercise 6.44]{detnotes} for the details.)

Now,%
\begin{align}
&  \sum_{\substack{\sigma\in S_{n};\\\sigma\left(  P\right)  =Q}%
}\ \ \underbrace{\left(  -1\right)  ^{\sigma}}_{\substack{=\left(  -1\right)
^{\operatorname*{sum}P+\operatorname*{sum}Q}\left(  -1\right)  ^{\alpha
_{\sigma}}\left(  -1\right)  ^{\beta_{\sigma}}\\\text{(by
(\ref{pf.thm.det.det(A+B).signs-equal}))}}}\underbrace{\left(  \prod
\limits_{i\in P}A_{i,\sigma\left(  i\right)  }\right)  }_{\substack{=\prod
_{i=1}^{k}A_{p_{i},q_{\alpha_{\sigma}\left(  i\right)  }}\\\text{(by the
definition of }\alpha_{\sigma}\text{)}}}\ \ \underbrace{\left(  \prod
\limits_{i\in\widetilde{P}}B_{i,\sigma\left(  i\right)  }\right)
}_{\substack{=\prod_{i=1}^{\ell}B_{p_{i}^{\prime},q_{\beta_{\sigma}\left(
i\right)  }^{\prime}}\\\text{(by the definition of }\beta_{\sigma}\text{)}%
}}\nonumber\\
&  =\sum_{\substack{\sigma\in S_{n};\\\sigma\left(  P\right)  =Q}}\left(
-1\right)  ^{\operatorname*{sum}P+\operatorname*{sum}Q}\left(  -1\right)
^{\alpha_{\sigma}}\left(  -1\right)  ^{\beta_{\sigma}}\left(  \prod_{i=1}%
^{k}A_{p_{i},q_{\alpha_{\sigma}\left(  i\right)  }}\right)  \left(
\prod_{i=1}^{\ell}B_{p_{i}^{\prime},q_{\beta_{\sigma}\left(  i\right)
}^{\prime}}\right) \nonumber\\
&  =\sum_{\left(  \alpha,\beta\right)  \in S_{k}\times S_{\ell}}\left(
-1\right)  ^{\operatorname*{sum}P+\operatorname*{sum}Q}\left(  -1\right)
^{\alpha}\left(  -1\right)  ^{\beta}\left(  \prod_{i=1}^{k}A_{p_{i}%
,q_{\alpha\left(  i\right)  }}\right)  \left(  \prod_{i=1}^{\ell}%
B_{p_{i}^{\prime},q_{\beta\left(  i\right)  }^{\prime}}\right) \nonumber\\
&  \ \ \ \ \ \ \ \ \ \ \ \ \ \ \ \ \ \ \ \ \left(
\begin{array}
[c]{c}%
\text{here, we have substituted }\left(  \alpha,\beta\right)  \text{ for
}\left(  \alpha_{\sigma},\beta_{\sigma}\right)  \text{ in the sum,}\\
\text{since the map }\sigma\mapsto\left(  \alpha_{\sigma},\beta_{\sigma
}\right)  \text{ is a bijection}%
\end{array}
\right) \nonumber\\
&  =\left(  -1\right)  ^{\operatorname*{sum}P+\operatorname*{sum}%
Q}\underbrace{\left(  \sum_{\alpha\in S_{k}}\left(  -1\right)  ^{\alpha}%
\prod_{i=1}^{k}A_{p_{i},q_{\alpha\left(  i\right)  }}\right)  }%
_{\substack{=\det\left(  \operatorname*{sub}\nolimits_{P}^{Q}A\right)
\\\text{(by the definition of }\operatorname*{sub}\nolimits_{P}^{Q}%
A\\\text{and its determinant)}}}\underbrace{\left(  \sum_{\beta\in S_{\ell}%
}\left(  -1\right)  ^{\beta}\prod_{i=1}^{\ell}B_{p_{i}^{\prime},q_{\beta
\left(  i\right)  }^{\prime}}\right)  }_{\substack{=\det\left(
\operatorname*{sub}\nolimits_{\widetilde{P}}^{\widetilde{Q}}B\right)
\\\text{(by the definition of }\operatorname*{sub}\nolimits_{\widetilde{P}%
}^{\widetilde{Q}}B\\\text{and its determinant)}}}\nonumber\\
&  =\left(  -1\right)  ^{\operatorname*{sum}P+\operatorname*{sum}Q}\det\left(
\operatorname*{sub}\nolimits_{P}^{Q}A\right)  \cdot\det\left(
\operatorname*{sub}\nolimits_{\widetilde{P}}^{\widetilde{Q}}B\right)  .
\label{pf.thm.det.det(A+B).sumPQ}%
\end{align}

Forget that we fixed $P$ and $Q$. We thus have proved
(\ref{pf.thm.det.det(A+B).sumPQ}) for any two subsets $P$ and $Q$ of $\left[
n\right]  $ satisfying $\left\vert P\right\vert =\left\vert Q\right\vert $.
Thus, our original computation of $\det\left(  A+B\right)  $ becomes%
\begin{align*}
\det\left(  A+B\right)   &  =\sum\limits_{P\subseteq\left[  n\right]
}\ \ \sum_{\substack{Q\subseteq\left[  n\right]  ;\\\left\vert P\right\vert
=\left\vert Q\right\vert }}\ \ \underbrace{\sum_{\substack{\sigma\in
S_{n};\\\sigma\left(  P\right)  =Q}}\left(  -1\right)  ^{\sigma}\left(
\prod\limits_{i\in P}A_{i,\sigma\left(  i\right)  }\right)  \left(
\prod\limits_{i\in\widetilde{P}}B_{i,\sigma\left(  i\right)  }\right)
}_{\substack{=\left(  -1\right)  ^{\operatorname*{sum}P+\operatorname*{sum}%
Q}\det\left(  \operatorname*{sub}\nolimits_{P}^{Q}A\right)  \cdot\det\left(
\operatorname*{sub}\nolimits_{\widetilde{P}}^{\widetilde{Q}}B\right)
\\\text{(by (\ref{pf.thm.det.det(A+B).sumPQ}))}}}\\
&  =\sum_{P\subseteq\left[  n\right]  }\ \ \sum_{\substack{Q\subseteq\left[
n\right]  ;\\\left\vert P\right\vert =\left\vert Q\right\vert }}\left(
-1\right)  ^{\operatorname*{sum}P+\operatorname*{sum}Q}\det\left(
\operatorname*{sub}\nolimits_{P}^{Q}A\right)  \cdot\det\left(
\operatorname*{sub}\nolimits_{\widetilde{P}}^{\widetilde{Q}}B\right)  .
\end{align*}
This proves Theorem \ref{thm.det.det(A+B)}.
\end{fineprint}
\end{proof}

We shall soon see a few applications of Theorem \ref{thm.det.det(A+B)}. First,
let us observe a simple property of diagonal matrices:\footnote{We are using
the Iverson bracket notation (see Definition \ref{def.iverson}) again.}

\begin{lemma}
\label{lem.det.minors-diag}Let $n\in\mathbb{N}$. Let $d_{1},d_{2},\ldots
,d_{n}\in K$. Let%
\[
D:=\left(  d_{i}\left[  i=j\right]  \right)  _{1\leq i\leq n,\ 1\leq j\leq
n}=\left(
\begin{array}
[c]{cccc}%
d_{1} & 0 & \cdots & 0\\
0 & d_{2} & \cdots & 0\\
\vdots & \vdots & \ddots & \vdots\\
0 & 0 & \cdots & d_{n}%
\end{array}
\right)  \in K^{n\times n}%
\]
be the diagonal $n\times n$-matrix with diagonal entries $d_{1},d_{2}%
,\ldots,d_{n}$. Then: \medskip

\textbf{(a)} We have $\det\left(  \operatorname*{sub}\nolimits_{P}%
^{P}D\right)  =\prod_{i\in P}d_{i}$ for any subset $P$ of $\left[  n\right]
$. \medskip

\textbf{(b)} Let $P$ and $Q$ be two distinct subsets of $\left[  n\right]  $
satisfying $\left\vert P\right\vert =\left\vert Q\right\vert $. Then,
$\det\left(  \operatorname*{sub}\nolimits_{P}^{Q}D\right)  =0$.
\end{lemma}

\begin{proof}
[Proof of Lemma \ref{lem.det.minors-diag}.]This is \cite[Lemma 6.163]%
{detnotes} (slightly rewritten); see \cite[Exercise 6.49]{detnotes} for a
detailed proof. (That said, it is almost obvious: In part \textbf{(a)}, the
submatrix $\operatorname*{sub}\nolimits_{P}^{P}D$ is itself diagonal, and its
diagonal entries are precisely the $d_{i}$ for $i\in P$. In part \textbf{(b)},
the submatrix $\operatorname*{sub}\nolimits_{P}^{Q}D$ has a zero row (indeed,
from $\left\vert P\right\vert =\left\vert Q\right\vert $ and $P\neq Q$, we see
that there exists some $i\in P\setminus Q$, and then the corresponding row of
$\operatorname*{sub}\nolimits_{P}^{Q}D$ is zero) and thus has determinant $0$.)
\end{proof}

Lemma \ref{lem.det.minors-diag} gives very simple formulas for minors of
diagonal matrices. Thus, the formula of Theorem \ref{thm.det.det(A+B)} becomes
simpler when the matrix $B$ is diagonal:

\begin{theorem}
\label{thm.det.det(A+D)}Let $n\in\mathbb{N}$. Let $A$ and $D$ be two $n\times
n$-matrices in $K^{n\times n}$ such that the matrix $D$ is diagonal. Let
$d_{1},d_{2},\ldots,d_{n}$ be the diagonal entries of the diagonal matrix $D$,
so that%
\[
D=\left(  d_{i}\left[  i=j\right]  \right)  _{1\leq i\leq n,\ 1\leq j\leq
n}=\left(
\begin{array}
[c]{cccc}%
d_{1} & 0 & \cdots & 0\\
0 & d_{2} & \cdots & 0\\
\vdots & \vdots & \ddots & \vdots\\
0 & 0 & \cdots & d_{n}%
\end{array}
\right)  \in K^{n\times n}.
\]
Then,%
\[
\det\left(  A+D\right)  =\sum_{P\subseteq\left[  n\right]  }\det\left(
\operatorname*{sub}\nolimits_{P}^{P}A\right)  \cdot\prod_{i\in\left[
n\right]  \setminus P}d_{i}.
\]

\end{theorem}

The minors $\det\left(  \operatorname*{sub}\nolimits_{P}^{P}A\right)  $ of an
$n\times n$-matrix $A$ are known as the \emph{principal minors} of $A$.

\begin{example}
For $n=3$, the claim of Theorem \ref{thm.det.det(A+D)} is%
\begin{align*}
&  \det\left(  A+D\right) \\
&  =\sum_{P\subseteq\left[  3\right]  }\det\left(  \operatorname*{sub}%
\nolimits_{P}^{P}A\right)  \cdot\prod_{i\in\left[  3\right]  \setminus P}%
d_{i}\\
&  =\underbrace{\det\left(  \operatorname*{sub}\nolimits_{\varnothing
}^{\varnothing}A\right)  }_{=1}\cdot\underbrace{\prod_{i\in\left[  3\right]
\setminus\varnothing}d_{i}}_{=d_{1}d_{2}d_{3}}+\underbrace{\det\left(
\operatorname*{sub}\nolimits_{\left\{  1\right\}  }^{\left\{  1\right\}
}A\right)  }_{=A_{1,1}}\cdot\underbrace{\prod_{i\in\left[  3\right]
\setminus\left\{  1\right\}  }d_{i}}_{=d_{2}d_{3}}\\
&  \ \ \ \ \ \ \ \ \ \ +\underbrace{\det\left(  \operatorname*{sub}%
\nolimits_{\left\{  2\right\}  }^{\left\{  2\right\}  }A\right)  }_{=A_{2,2}%
}\cdot\underbrace{\prod_{i\in\left[  3\right]  \setminus\left\{  2\right\}
}d_{i}}_{=d_{1}d_{3}}+\underbrace{\det\left(  \operatorname*{sub}%
\nolimits_{\left\{  3\right\}  }^{\left\{  3\right\}  }A\right)  }_{=A_{3,3}%
}\cdot\underbrace{\prod_{i\in\left[  3\right]  \setminus\left\{  3\right\}
}d_{i}}_{=d_{1}d_{2}}\\
&  \ \ \ \ \ \ \ \ \ \ +\underbrace{\det\left(  \operatorname*{sub}%
\nolimits_{\left\{  1,2\right\}  }^{\left\{  1,2\right\}  }A\right)  }%
_{=\det\left(
\begin{array}
[c]{cc}%
A_{1,1} & A_{1,2}\\
A_{2,1} & A_{2,2}%
\end{array}
\right)  }\cdot\underbrace{\prod_{i\in\left[  3\right]  \setminus\left\{
1,2\right\}  }d_{i}}_{=d_{3}}+\underbrace{\det\left(  \operatorname*{sub}%
\nolimits_{\left\{  1,3\right\}  }^{\left\{  1,3\right\}  }A\right)  }%
_{=\det\left(
\begin{array}
[c]{cc}%
A_{1,1} & A_{1,3}\\
A_{3,1} & A_{3,3}%
\end{array}
\right)  }\cdot\underbrace{\prod_{i\in\left[  3\right]  \setminus\left\{
1,3\right\}  }d_{i}}_{=d_{2}}\\
&  \ \ \ \ \ \ \ \ \ \ +\underbrace{\det\left(  \operatorname*{sub}%
\nolimits_{\left\{  2,3\right\}  }^{\left\{  2,3\right\}  }A\right)  }%
_{=\det\left(
\begin{array}
[c]{cc}%
A_{2,2} & A_{2,3}\\
A_{3,2} & A_{3,3}%
\end{array}
\right)  }\cdot\underbrace{\prod_{i\in\left[  3\right]  \setminus\left\{
2,3\right\}  }d_{i}}_{=d_{1}}+\underbrace{\det\left(  \operatorname*{sub}%
\nolimits_{\left\{  1,2,3\right\}  }^{\left\{  1,2,3\right\}  }A\right)
}_{=\det A}\cdot\underbrace{\prod_{i\in\left[  3\right]  \setminus\left\{
1,2,3\right\}  }d_{i}}_{\substack{=\left(  \text{empty product}\right)
\\=1}}\\
&  =d_{1}d_{2}d_{3}+A_{1,1}d_{2}d_{3}+A_{2,2}d_{1}d_{3}+A_{3,3}d_{1}d_{2}\\
&  \ \ \ \ \ \ \ \ \ \ +\det\left(
\begin{array}
[c]{cc}%
A_{1,1} & A_{1,2}\\
A_{2,1} & A_{2,2}%
\end{array}
\right)  \cdot d_{3}+\det\left(
\begin{array}
[c]{cc}%
A_{1,1} & A_{1,3}\\
A_{3,1} & A_{3,3}%
\end{array}
\right)  \cdot d_{2}\\
&  \ \ \ \ \ \ \ \ \ \ +\det\left(
\begin{array}
[c]{cc}%
A_{2,2} & A_{2,3}\\
A_{3,2} & A_{3,3}%
\end{array}
\right)  \cdot d_{1}+\det A.
\end{align*}

\end{example}

\begin{proof}
[Proof of Theorem \ref{thm.det.det(A+D)} (sketched).](See \cite[Corollary
6.162]{detnotes} for details.) We shall use the notations $\widetilde{I}$ and
$\operatorname*{sum}S$ as defined in Theorem \ref{thm.det.det(A+B)}. If $P$
and $Q$ are two subsets of $\left[  n\right]  $ satisfying $\left\vert
P\right\vert =\left\vert Q\right\vert $ but $P\neq Q$, then their complements
$\widetilde{P}$ and $\widetilde{Q}$ are also distinct (since $P\neq Q$) and
satisfy $\left\vert \widetilde{P}\right\vert =\left\vert \widetilde{Q}%
\right\vert $ (since $\left\vert P\right\vert =\left\vert Q\right\vert $), and
therefore Lemma \ref{lem.det.minors-diag} \textbf{(b)} (applied to
$\widetilde{P}$ and $\widetilde{Q}$ instead of $P$ and $Q$) yields%
\begin{equation}
\det\left(  \operatorname*{sub}\nolimits_{\widetilde{P}}^{\widetilde{Q}%
}D\right)  =0. \label{pf.thm.det.det(A+D).ex0}%
\end{equation}

Now, Theorem \ref{thm.det.det(A+B)} (applied to $B=D$) yields%
\begin{align*}
\det\left(  A+D\right)   &  =\sum_{P\subseteq\left[  n\right]  }%
\ \ \sum_{\substack{Q\subseteq\left[  n\right]  ;\\\left\vert P\right\vert
=\left\vert Q\right\vert }}\left(  -1\right)  ^{\operatorname*{sum}%
P+\operatorname*{sum}Q}\det\left(  \operatorname*{sub}\nolimits_{P}%
^{Q}A\right)  \cdot\underbrace{\det\left(  \operatorname*{sub}%
\nolimits_{\widetilde{P}}^{\widetilde{Q}}D\right)  }_{\substack{\text{This is
}0\text{ if }P\neq Q\\\text{(by (\ref{pf.thm.det.det(A+D).ex0}))}}}\\
&  =\sum_{P\subseteq\left[  n\right]  }\underbrace{\left(  -1\right)
^{\operatorname*{sum}P+\operatorname*{sum}P}}_{=1}\det\left(
\operatorname*{sub}\nolimits_{P}^{P}A\right)  \cdot\underbrace{\det\left(
\operatorname*{sub}\nolimits_{\widetilde{P}}^{\widetilde{P}}D\right)
}_{\substack{=\prod_{i\in\widetilde{P}}d_{i}\\\text{(by Lemma
\ref{lem.det.minors-diag} \textbf{(a)})}}}\\
&  \ \ \ \ \ \ \ \ \ \ \ \ \ \ \ \ \ \ \ \ \left(
\begin{array}
[c]{c}%
\text{here, we have removed all addends with }P\neq Q\\
\text{from the double sum, since these addends are }0
\end{array}
\right) \\
&  =\sum_{P\subseteq\left[  n\right]  }\det\left(  \operatorname*{sub}%
\nolimits_{P}^{P}A\right)  \cdot\prod_{i\in\widetilde{P}}d_{i}.
\end{align*}
This proves Theorem \ref{thm.det.det(A+D)} (since $\widetilde{P}=\left[
n\right]  \setminus P$).
\end{proof}

As a particular case of Theorem \ref{thm.det.det(A+D)}, we quickly obtain the
following formula for a class of determinants that frequently appear in graph theory:

\begin{proposition}
\label{prop.det.x+ai}Let $n\in\mathbb{N}$. Let $d_{1},d_{2},\ldots,d_{n}\in K$
and $x\in K$. Let $F$ be the $n\times n$-matrix
\[
\left(  x+d_{i}\left[  i=j\right]  \right)  _{1\leq i\leq n,\ 1\leq j\leq
n}=\left(
\begin{array}
[c]{cccc}%
x+d_{1} & x & \cdots & x\\
x & x+d_{2} & \cdots & x\\
\vdots & \vdots & \ddots & \vdots\\
x & x & \cdots & x+d_{n}%
\end{array}
\right)  \in K^{n\times n}.
\]
Then,%
\[
\det F=d_{1}d_{2}\cdots d_{n}+x\sum_{i=1}^{n}d_{1}d_{2}\cdots\widehat{d_{i}%
}\cdots d_{n},
\]
where the hat over the \textquotedblleft$d_{i}$\textquotedblright\ means
\textquotedblleft omit the $d_{i}$ factor\textquotedblright\ (that is, the
expression \textquotedblleft$d_{1}d_{2}\cdots\widehat{d_{i}}\cdots d_{n}%
$\textquotedblright\ is to be understood as \textquotedblleft$d_{1}d_{2}\cdots
d_{i-1}d_{i+1}d_{i+2}\cdots d_{n}$\textquotedblright).
\end{proposition}

\begin{proof}
[Proof of Proposition \ref{prop.det.x+ai} (sketched).]Define the two $n\times
n$-matrices%
\begin{align*}
A  &  :=\left(  x\right)  _{1\leq i\leq n,\ 1\leq j\leq n}=\left(
\begin{array}
[c]{cccc}%
x & x & \cdots & x\\
x & x & \cdots & x\\
\vdots & \vdots & \ddots & \vdots\\
x & x & \cdots & x
\end{array}
\right)  \in K^{n\times n}\ \ \ \ \ \ \ \ \ \ \text{and}\\
D  &  :=\left(  d_{i}\left[  i=j\right]  \right)  _{1\leq i\leq n,\ 1\leq
j\leq n}=\left(
\begin{array}
[c]{cccc}%
d_{1} & 0 & \cdots & 0\\
0 & d_{2} & \cdots & 0\\
\vdots & \vdots & \ddots & \vdots\\
0 & 0 & \cdots & d_{n}%
\end{array}
\right)  \in K^{n\times n}.
\end{align*}
Then, it is clear that $F=A+D$. Moreover, the matrix $D$ is diagonal, and its
diagonal entries are $d_{1},d_{2},\ldots,d_{n}$. Hence, Theorem
\ref{thm.det.det(A+D)} yields%
\begin{equation}
\det\left(  A+D\right)  =\sum_{P\subseteq\left[  n\right]  }\det\left(
\operatorname*{sub}\nolimits_{P}^{P}A\right)  \cdot\prod_{i\in\left[
n\right]  \setminus P}d_{i}. \label{pf.prop.det.x+ai.1}%
\end{equation}

However, if $P$ is a subset of $\left[  n\right]  $, then $\operatorname*{sub}%
\nolimits_{P}^{P}A$ is a submatrix of $A$ and thus has the form $\left(
\begin{array}
[c]{cccc}%
x & x & \cdots & x\\
x & x & \cdots & x\\
\vdots & \vdots & \ddots & \vdots\\
x & x & \cdots & x
\end{array}
\right)  $ (since all entries of $A$ equal $x$). If the subset $P$ of $\left[
n\right]  $ has size $\geq2$, then this submatrix $\operatorname*{sub}%
\nolimits_{P}^{P}A=\left(
\begin{array}
[c]{cccc}%
x & x & \cdots & x\\
x & x & \cdots & x\\
\vdots & \vdots & \ddots & \vdots\\
x & x & \cdots & x
\end{array}
\right)  $ has size $\left\vert P\right\vert \geq2$ and therefore has
determinant $0$ (by (\ref{eq.prop.det.xiyj.cor-x}), applied to $\left\vert
P\right\vert $ instead of $n$). In other words, we have%
\begin{equation}
\det\left(  \operatorname*{sub}\nolimits_{P}^{P}A\right)  =0
\label{pf.prop.det.x+ai.2}%
\end{equation}
whenever $P\subseteq\left[  n\right]  $ satisfies $\left\vert P\right\vert
\geq2$.

Now, (\ref{pf.prop.det.x+ai.1}) becomes%
\begin{align*}
\det\left(  A+D\right)   &  =\sum_{P\subseteq\left[  n\right]  }\det\left(
\operatorname*{sub}\nolimits_{P}^{P}A\right)  \cdot\prod_{i\in\left[
n\right]  \setminus P}d_{i}\\
&  =\sum_{\substack{P\subseteq\left[  n\right]  ;\\\left\vert P\right\vert
\leq1}}\det\left(  \operatorname*{sub}\nolimits_{P}^{P}A\right)  \cdot
\prod_{i\in\left[  n\right]  \setminus P}d_{i}+\sum_{\substack{P\subseteq
\left[  n\right]  ;\\\left\vert P\right\vert \geq2}}\underbrace{\det\left(
\operatorname*{sub}\nolimits_{P}^{P}A\right)  }_{\substack{=0\\\text{(by
(\ref{pf.prop.det.x+ai.2}))}}}\cdot\prod_{i\in\left[  n\right]  \setminus
P}d_{i}\\
&  \ \ \ \ \ \ \ \ \ \ \ \ \ \ \ \ \ \ \ \ \left(  \text{since each
}P\subseteq\left[  n\right]  \text{ satisfies either }\left\vert P\right\vert
\leq1\text{ or }\left\vert P\right\vert \geq2\right) \\
&  =\sum_{\substack{P\subseteq\left[  n\right]  ;\\\left\vert P\right\vert
\leq1}}\det\left(  \operatorname*{sub}\nolimits_{P}^{P}A\right)  \cdot
\prod_{i\in\left[  n\right]  \setminus P}d_{i}\\
&  =\underbrace{\det\left(  \operatorname*{sub}\nolimits_{\varnothing
}^{\varnothing}A\right)  }_{\substack{=1\\\text{(since }\operatorname*{sub}%
\nolimits_{\varnothing}^{\varnothing}A\text{ is}\\\text{a }0\times
0\text{-matrix)}}}\cdot\underbrace{\prod_{i\in\left[  n\right]  \setminus
\varnothing}d_{i}}_{\substack{=\prod_{i\in\left[  n\right]  }d_{i}%
\\=d_{1}d_{2}\cdots d_{n}}}+\sum_{p=1}^{n}\underbrace{\det\left(
\operatorname*{sub}\nolimits_{\left\{  p\right\}  }^{\left\{  p\right\}
}A\right)  }_{\substack{=x\\\text{(since }\operatorname*{sub}%
\nolimits_{\left\{  p\right\}  }^{\left\{  p\right\}  }A=\left(
\begin{array}
[c]{c}%
x
\end{array}
\right)  \text{)}}}\cdot\underbrace{\prod_{i\in\left[  n\right]
\setminus\left\{  p\right\}  }d_{i}}_{=d_{1}d_{2}\cdots\widehat{d_{p}}\cdots
d_{n}}\\
&  \ \ \ \ \ \ \ \ \ \ \ \ \ \ \ \ \ \ \ \ \left(
\begin{array}
[c]{c}%
\text{since the subsets }P\text{ of }\left[  n\right]  \text{ satisfying
}\left\vert P\right\vert \leq1\\
\text{are the }n+1\text{ subsets }\varnothing\text{, }\left\{  1\right\}
\text{, }\left\{  2\right\}  \text{, }\ldots\text{, }\left\{  n\right\}
\end{array}
\right) \\
&  =d_{1}d_{2}\cdots d_{n}+\sum_{p=1}^{n}x\cdot d_{1}d_{2}\cdots
\widehat{d_{p}}\cdots d_{n}\\
&  =d_{1}d_{2}\cdots d_{n}+x\sum_{p=1}^{n}d_{1}d_{2}\cdots\widehat{d_{p}%
}\cdots d_{n}\\
&  =d_{1}d_{2}\cdots d_{n}+x\sum_{i=1}^{n}d_{1}d_{2}\cdots\widehat{d_{i}%
}\cdots d_{n}%
\end{align*}
(here, we have renamed the summation index $p$ as $i$). In view of $F=A+D$,
this rewrites as%
\[
\det F=d_{1}d_{2}\cdots d_{n}+x\sum_{i=1}^{n}d_{1}d_{2}\cdots\widehat{d_{i}%
}\cdots d_{n}.
\]
Thus, Proposition \ref{prop.det.x+ai} is proven.

(See \url{https://math.stackexchange.com/questions/2110766/} for some
different proofs of Proposition \ref{prop.det.x+ai}.)
\end{proof}

Another application of Theorem \ref{thm.det.det(A+D)} is an explicit formula
for the characteristic polynomial of a matrix. We recall that the
characteristic polynomial of an $n\times n$-matrix $A\in K^{n\times n}$ is
defined to be the polynomial\footnote{Here, $I_{n}$ denotes the $n\times n$
identity matrix.} $\det\left(  xI_{n}-A\right)  \in K\left[  x\right]  $ (some
authors define it to be $\det\left(  A-xI_{n}\right)  $ instead, but this is
the same up to sign). This is a polynomial of degree $n$, whose leading term
is $x^{n}$, whose next-highest term is $-\operatorname*{Tr}A\cdot x^{n-1}$
where $\operatorname*{Tr}A:=\sum_{i=1}^{n}A_{i,i}$ is the \emph{trace} of $A$,
and whose constant term is $\left(  -1\right)  ^{n}\det A$. We shall extend
this by explicitly computing all coefficients of this polynomial. For the sake
of simplicity, we will compute $\det\left(  A+xI_{n}\right)  $ instead of
$\det\left(  xI_{n}-A\right)  $ (this is tantamount to replacing $x$ by $-x$
and multiplying the determinant by $\left(  -1\right)  ^{n}$), and we will
take $x$ to be an element of $K$ rather than an indeterminate (but this
setting is more general, since we can take $K$ itself to be a polynomial ring
and then choose $x$ to be its indeterminate). Here is our formula:

\begin{proposition}
\label{prop.det.charpol-explicit}Let $n\in\mathbb{N}$. Let $A\in K^{n\times
n}$ be an $n\times n$-matrix. Let $x\in K$. Let $I_{n}$ denote the $n\times n$
identity matrix. Then,%
\[
\det\left(  A+xI_{n}\right)  =\sum_{P\subseteq\left[  n\right]  }\det\left(
\operatorname*{sub}\nolimits_{P}^{P}A\right)  \cdot x^{n-\left\vert
P\right\vert }=\sum_{k=0}^{n}\left(  \sum_{\substack{P\subseteq\left[
n\right]  ;\\\left\vert P\right\vert =n-k}}\det\left(  \operatorname*{sub}%
\nolimits_{P}^{P}A\right)  \right)  x^{k}.
\]

\end{proposition}

\begin{proof}
[Proof of Proposition \ref{prop.det.charpol-explicit} (sketched).](See
\cite[Corollary 6.164]{detnotes} for details.) The matrix $xI_{n}$ is
diagonal, and its diagonal entries are $x,x,\ldots,x$; in fact,%
\[
xI_{n}=\left(  x\left[  i=j\right]  \right)  _{1\leq i\leq n,\ 1\leq j\leq
n}=\left(
\begin{array}
[c]{cccc}%
x & 0 & \cdots & 0\\
0 & x & \cdots & 0\\
\vdots & \vdots & \ddots & \vdots\\
0 & 0 & \cdots & x
\end{array}
\right)  .
\]
Hence, Theorem \ref{thm.det.det(A+D)} (applied to $D=xI_{n}$ and $d_{i}=x$)
yields%
\begin{align*}
&  \det\left(  A+xI_{n}\right) \\
&  =\sum_{P\subseteq\left[  n\right]  }\det\left(  \operatorname*{sub}%
\nolimits_{P}^{P}A\right)  \cdot\underbrace{\prod_{i\in\left[  n\right]
\setminus P}x}_{\substack{=x^{\left\vert \left[  n\right]  \setminus
P\right\vert }=x^{n-\left\vert P\right\vert }\\\text{(since }\left\vert
\left[  n\right]  \setminus P\right\vert =\left\vert \left[  n\right]
\right\vert -\left\vert P\right\vert =n-\left\vert P\right\vert \text{)}%
}}=\sum_{P\subseteq\left[  n\right]  }\det\left(  \operatorname*{sub}%
\nolimits_{P}^{P}A\right)  \cdot x^{n-\left\vert P\right\vert }\\
&  =\sum_{k=0}^{n}\ \ \sum_{\substack{P\subseteq\left[  n\right]
;\\\left\vert P\right\vert =k}}\det\left(  \operatorname*{sub}\nolimits_{P}%
^{P}A\right)  x^{n-k}\ \ \ \ \ \ \ \ \ \ \left(
\begin{array}
[c]{c}%
\text{here, we have split the sum}\\
\text{according to the value of }\left\vert P\right\vert
\end{array}
\right) \\
&  =\sum_{k=0}^{n}\ \ \sum_{\substack{P\subseteq\left[  n\right]
;\\\left\vert P\right\vert =n-k}}\det\left(  \operatorname*{sub}%
\nolimits_{P}^{P}A\right)  x^{k}\ \ \ \ \ \ \ \ \ \ \left(
\begin{array}
[c]{c}%
\text{here, we have substituted }n-k\\
\text{for }k\text{ in the sum}%
\end{array}
\right) \\
&  =\sum_{k=0}^{n}\left(  \sum_{\substack{P\subseteq\left[  n\right]
;\\\left\vert P\right\vert =n-k}}\det\left(  \operatorname*{sub}%
\nolimits_{P}^{P}A\right)  \right)  x^{k}.
\end{align*}
Proposition \ref{prop.det.charpol-explicit} is proven.
\end{proof}

\subsubsection{Factoring the matrix}

Next, we will see some tricks for computing determinants.

Let us compute a determinant that recently went viral on the internet after
Timothy Gowers livestreamed himself computing it\footnote{See
\url{https://www.youtube.com/watch?v=byjhpzEoXFs} and
\url{https://www.youtube.com/watch?v=frvBdaqLgLo} and
\url{https://www.youtube.com/watch?v=m8R9rVb0M5o} .} (\cite[Exercise
6.11]{detnotes}, \cite{EdelStrang}):

\begin{proposition}
\label{prop.det.pascal-LU}Let $n\in\mathbb{N}$. Let $A$ be the $n\times
n$-matrix%
\[
\left(  \dbinom{i+j-2}{i-1}\right)  _{1\leq i\leq n,\ 1\leq j\leq n}=\left(
\begin{array}
[c]{cccc}%
\dbinom{0}{0} & \dbinom{1}{0} & \cdots & \dbinom{n-1}{0}\\
\dbinom{1}{1} & \dbinom{2}{1} & \cdots & \dbinom{n}{1}\\
\vdots & \vdots & \ddots & \vdots\\
\dbinom{n-1}{n-1} & \dbinom{n}{n-1} & \cdots & \dbinom{2n-2}{n-1}%
\end{array}
\right)  .
\]
(For example, for $n=4$, we have $A=\left(
\begin{array}
[c]{cccc}%
1 & 1 & 1 & 1\\
1 & 2 & 3 & 4\\
1 & 3 & 6 & 10\\
1 & 4 & 10 & 20
\end{array}
\right)  $.)

Then, $\det A=1$.
\end{proposition}

There are many ways to prove Proposition \ref{prop.det.pascal-LU}, but here is
a particularly simple one:

\begin{proof}
[Proof of Proposition \ref{prop.det.pascal-LU} (sketched).](See \cite[Exercise
6.11]{detnotes} for details.) For any $i,j\in\left[  n\right]  $, we have%
\begin{align*}
A_{i,j}  &  =\dbinom{i+j-2}{i-1}\ \ \ \ \ \ \ \ \ \ \left(  \text{by the
definition of }A\right) \\
&  =\dbinom{\left(  i-1\right)  +\left(  j-1\right)  }{i-1}=\dbinom{\left(
i-1\right)  +\left(  j-1\right)  }{j-1}\ \ \ \ \ \ \ \ \ \ \left(  \text{by
Theorem \ref{thm.binom.sym}}\right) \\
&  =\sum_{k=0}^{i-1}\dbinom{i-1}{k}\underbrace{\dbinom{j-1}{j-1-k}%
}_{\substack{=\dbinom{j-1}{k}\\\text{(by Theorem \ref{thm.binom.sym})}}}\\
&  \ \ \ \ \ \ \ \ \ \ \ \ \ \ \ \ \ \ \ \ \left(
\begin{array}
[c]{c}%
\text{by Proposition \ref{prop.binom.vandermonde.NN}, applied to }i-1\text{,
}j-1\text{ and }j-1\\
\text{instead of }a\text{, }b\text{ and }n
\end{array}
\right) \\
&  =\sum_{k=0}^{i-1}\dbinom{i-1}{k}\dbinom{j-1}{k}=\sum_{k=0}^{n-1}%
\dbinom{i-1}{k}\dbinom{j-1}{k}\\
&  \ \ \ \ \ \ \ \ \ \ \ \ \ \ \ \ \ \ \ \ \left(
\begin{array}
[c]{c}%
\text{here, we extended the sum upwards to }k=n-1\text{,}\\
\text{but this has not changed the value of the sum,}\\
\text{since all newly introduced addends are }0\\
\text{(since }\dbinom{i-1}{k}=0\text{ whenever }k>i-1\text{)}%
\end{array}
\right) \\
&  =\sum_{k=1}^{n}\dbinom{i-1}{k-1}\dbinom{j-1}{k-1}\\
&  \ \ \ \ \ \ \ \ \ \ \ \ \ \ \ \ \ \ \ \ \left(  \text{here, we have
substituted }k-1\text{ for }k\text{ in the sum}\right)  .
\end{align*}
If we define two $n\times n$-matrices%
\[
L:=\left(  \dbinom{i-1}{k-1}\right)  _{1\leq i\leq n,\ 1\leq k\leq
n}\ \ \ \ \ \ \ \ \ \ \text{and}\ \ \ \ \ \ \ \ \ \ U:=\left(  \dbinom
{j-1}{k-1}\right)  _{1\leq k\leq n,\ 1\leq j\leq n},
\]
then this rewrites as%
\[
A_{i,j}=\sum_{k=1}^{n}L_{i,k}U_{k,j}=\left(  LU\right)  _{i,j}%
\]
(by the definition of the matrix product $LU$). Since this equality holds for
all $i,j\in\left[  n\right]  $, we thus conclude that $A=LU$.

Notice, however, that the matrix $L$ is lower-triangular (because if $i<k$,
then $i-1<k-1$ and therefore $L_{i,k}=\dbinom{i-1}{k-1}=0$), and thus (by
Theorem \ref{thm.det.triang}) its determinant is the product of its diagonal
entries. In other words,%
\[
\det L=\dbinom{1-1}{1-1}\dbinom{2-1}{2-1}\cdots\dbinom{n-1}{n-1}=\prod
_{k=1}^{n}\underbrace{\dbinom{k-1}{k-1}}_{=1}=1.
\]

Similarly, the matrix $U$ is upper-triangular, and its determinant is $\det
U=1$ as well.

Now, from $A=LU$, we obtain%
\begin{align*}
\det A  &  =\det\left(  LU\right)  =\underbrace{\det L}_{=1}\cdot
\,\underbrace{\det U}_{=1}\ \ \ \ \ \ \ \ \ \ \left(  \text{by Theorem
\ref{thm.det.detAB}}\right) \\
&  =1;
\end{align*}
this proves Proposition \ref{prop.det.pascal-LU}.
\end{proof}

How can you discover such a proof? Our serendipitous factorization of $A$ as
$LU$ might appear unmotivated, but from the viewpoint of linear algebra it is
an instance of a well-known and well-understood kind of factorization, known
as the \emph{LU-decomposition} or the \emph{LU-factorization}. Over a field,
almost every square matrix\footnote{I will not go into details as to what
\textquotedblleft almost every\textquotedblright\ means here.} has an
LU-decomposition (i.e., a factorization as a product of a lower-triangular
matrix with an upper-triangular matrix). This LU-decomposition is unique if
you require (e.g.) the diagonal entries of the lower-triangular factor to all
equal $1$. It can furthermore be algorithmically computed using Gaussian
elimination (see, e.g., \cite[\S 1.3, Theorem 1.3]{OlvSha}). Now, computing
the LU-decomposition of the matrix $A$ from Proposition
\ref{prop.det.pascal-LU} for $n=4$, we find%
\[
\left(
\begin{array}
[c]{cccc}%
1 & 1 & 1 & 1\\
1 & 2 & 3 & 4\\
1 & 3 & 6 & 10\\
1 & 4 & 10 & 20
\end{array}
\right)  =\underbrace{\left(
\begin{array}
[c]{cccc}%
1 & 0 & 0 & 0\\
1 & 1 & 0 & 0\\
1 & 2 & 1 & 0\\
1 & 3 & 3 & 1
\end{array}
\right)  }_{\text{this is the lower-triangular factor}}%
\ \ \ \underbrace{\left(
\begin{array}
[c]{cccc}%
1 & 1 & 1 & 1\\
0 & 1 & 2 & 3\\
0 & 0 & 1 & 3\\
0 & 0 & 0 & 1
\end{array}
\right)  }_{\text{this is the upper-triangular factor}}.
\]
The entries of both factors appear to be the binomial coefficients familiar
from Pascal's triangle. This suggests that we might have%
\[
L=\left(  \dbinom{i-1}{k-1}\right)  _{1\leq i\leq n,\ 1\leq k\leq
n}\ \ \ \ \ \ \ \ \ \ \text{and}\ \ \ \ \ \ \ \ \ \ U=\left(  \dbinom
{j-1}{k-1}\right)  _{1\leq k\leq n,\ 1\leq j\leq n},
\]
not just for $n=4$ but also for arbitrary $n\in\mathbb{N}$. And once this
guess has been made, it is easy to prove that $A=LU$ (our proof above is not
the only one possible; four proofs appear in \cite{EdelStrang}).

This is not the only example where LU-decomposition helps compute a
determinant (see, e.g., \cite[\S 2.6]{Krattenthaler} for examples). Sometimes
it is helpful to transpose a matrix, or to permute its rows or columns to
obtain a matrix with a good LU-decomposition.

\subsubsection{Factor hunting}

The next trick -- known as \emph{factor hunting} -- works not only for
determinants; however, determinants provide some of the simplest examples.

\begin{theorem}
[Vandermonde determinant]\label{thm.det.vander}Let $n\in\mathbb{N}$. Let
$a_{1},a_{2},\ldots,a_{n}$ be $n$ elements of $K$. Then: \medskip

\textbf{(a)} We have%
\[
\det\left(  \left(  a_{i}^{n-j}\right)  _{1\leq i\leq n,\ 1\leq j\leq
n}\right)  =\prod_{1\leq i<j\leq n}\left(  a_{i}-a_{j}\right)  .
\]

\textbf{(b)} We have%
\[
\det\left(  \left(  a_{j}^{n-i}\right)  _{1\leq i\leq n,\ 1\leq j\leq
n}\right)  =\prod_{1\leq i<j\leq n}\left(  a_{i}-a_{j}\right)  .
\]

\textbf{(c)} We have%
\[
\det\left(  \left(  a_{i}^{j-1}\right)  _{1\leq i\leq n,\ 1\leq j\leq
n}\right)  =\prod_{1\leq j<i\leq n}\left(  a_{i}-a_{j}\right)  .
\]

\textbf{(d)} We have%
\[
\det\left(  \left(  a_{j}^{i-1}\right)  _{1\leq i\leq n,\ 1\leq j\leq
n}\right)  =\prod_{1\leq j<i\leq n}\left(  a_{i}-a_{j}\right)  .
\]

\end{theorem}

\begin{example}
Here is what the four parts of Theorem \ref{thm.det.vander} say for $n=3$:
\medskip

\textbf{(a)} We have $\det\left(
\begin{array}
[c]{ccc}%
a_{1}^{2} & a_{1} & 1\\
a_{2}^{2} & a_{2} & 1\\
a_{3}^{2} & a_{3} & 1
\end{array}
\right)  =\left(  a_{1}-a_{2}\right)  \left(  a_{1}-a_{3}\right)  \left(
a_{2}-a_{3}\right)  $. \medskip

\textbf{(b)} We have $\det\left(
\begin{array}
[c]{ccc}%
a_{1}^{2} & a_{2}^{2} & a_{3}^{2}\\
a_{1} & a_{2} & a_{3}\\
1 & 1 & 1
\end{array}
\right)  =\left(  a_{1}-a_{2}\right)  \left(  a_{1}-a_{3}\right)  \left(
a_{2}-a_{3}\right)  $. \medskip

\textbf{(c)} We have $\det\left(
\begin{array}
[c]{ccc}%
1 & a_{1} & a_{1}^{2}\\
1 & a_{2} & a_{2}^{2}\\
1 & a_{3} & a_{3}^{2}%
\end{array}
\right)  =\left(  a_{2}-a_{1}\right)  \left(  a_{3}-a_{1}\right)  \left(
a_{3}-a_{2}\right)  $. \medskip

\textbf{(d)} We have $\det\left(
\begin{array}
[c]{ccc}%
1 & 1 & 1\\
a_{1} & a_{2} & a_{3}\\
a_{1}^{2} & a_{2}^{2} & a_{3}^{2}%
\end{array}
\right)  =\left(  a_{2}-a_{1}\right)  \left(  a_{3}-a_{1}\right)  \left(
a_{3}-a_{2}\right)  $.
\end{example}

Theorem \ref{thm.det.vander} is a classical and important result, known as the
\emph{Vandermonde determinant}. Many different proofs are known (see, e.g.,
\cite[Theorem 6.46]{detnotes} or \cite[\S 5.3]{Aigner07} or \cite[Section
III.8.6, Example (1)]{Bourba74} or \cite[Theorem 1]{GriHyp}; a combinatorial
proof can also be found in Exercise \ref{exe.det.vander-tour}; two more proofs
are obtained in Exercise \ref{exe.det.vander.kratt-lem6} and Exercise
\ref{exe.det.binom-by-vander}). We will now sketch a proof using factor
hunting and polynomials. We will first focus on proving part \textbf{(a)} of
Theorem \ref{thm.det.vander}, and afterwards derive the other parts from it.

The first step of our proof is reducing Theorem \ref{thm.det.vander}
\textbf{(a)} to the \textquotedblleft particular case\textquotedblright\ in
which $K$ is the polynomial ring $\mathbb{Z}\left[  x_{1},x_{2},\ldots
,x_{n}\right]  $ and the elements $a_{1},a_{2},\ldots,a_{n}$ are the
indeterminates $x_{1},x_{2},\ldots,x_{n}$. This is merely a particular case
(one possible choice of $K$ and $a_{1},a_{2},\ldots,a_{n}$ among many);
however, as we will soon see, proving Theorem \ref{thm.det.vander}
\textbf{(a)} in this particular case will quickly entail that Theorem
\ref{thm.det.vander} \textbf{(a)} holds in the general case. Let us elaborate
on this argument. First, let us state Theorem \ref{thm.det.vander}
\textbf{(a)} in this particular case as a lemma:

\begin{lemma}
\label{lem.det.vander.a.pol}Let $n\in\mathbb{N}$. Consider the polynomial ring
$\mathbb{Z}\left[  x_{1},x_{2},\ldots,x_{n}\right]  $ in $n$ indeterminates
$x_{1},x_{2},\ldots,x_{n}$ with integer coefficients. In this ring, we have%
\begin{equation}
\det\left(  \left(  x_{i}^{n-j}\right)  _{1\leq i\leq n,\ 1\leq j\leq
n}\right)  =\prod_{1\leq i<j\leq n}\left(  x_{i}-x_{j}\right)  .
\label{eq.lem.det.vander.a.pol.1}%
\end{equation}

\end{lemma}

We can derive Theorem \ref{thm.det.vander} \textbf{(a)} from Lemma
\ref{lem.det.vander.a.pol} as follows:

\begin{proof}
[Proof of Theorem \ref{thm.det.vander} \textbf{(a)} using Lemma
\ref{lem.det.vander.a.pol}.]The equality%
\begin{equation}
\det\left(  \left(  a_{i}^{n-j}\right)  _{1\leq i\leq n,\ 1\leq j\leq
n}\right)  =\prod_{1\leq i<j\leq n}\left(  a_{i}-a_{j}\right)
\label{pf.thm.det.vander.a.1}%
\end{equation}
follows from the equality (\ref{eq.lem.det.vander.a.pol.1}) by substituting
$a_{1},a_{2},\ldots,a_{n}$ for $x_{1},x_{2},\ldots,x_{n}$. \medskip

\begin{fineprint}
This is sufficiently clear to be considered a complete proof, but just in
case, here are a few details.

We can substitute $a_{1},a_{2},\ldots,a_{n}$ for $x_{1},x_{2},\ldots,x_{n}$ in
any polynomial $f\in\mathbb{Z}\left[  x_{1},x_{2},\ldots,x_{n}\right]  $,
since $a_{1},a_{2},\ldots,a_{n}$ are $n$ elements of a commutative ring
(namely, of $K$). It is obvious that $\prod_{1\leq i<j\leq n}\left(
x_{i}-x_{j}\right)  $ becomes $\prod_{1\leq i<j\leq n}\left(  a_{i}%
-a_{j}\right)  $ when we substitute $a_{1},a_{2},\ldots,a_{n}$ for
$x_{1},x_{2},\ldots,x_{n}$. It is perhaps a bit less obvious that $\det\left(
\left(  x_{i}^{n-j}\right)  _{1\leq i\leq n,\ 1\leq j\leq n}\right)  $ becomes
$\det\left(  \left(  a_{i}^{n-j}\right)  _{1\leq i\leq n,\ 1\leq j\leq
n}\right)  $ when we substitute $a_{1},a_{2},\ldots,a_{n}$ for $x_{1}%
,x_{2},\ldots,x_{n}$. To convince our skeptical selves of this, we expand both
determinants: The definition of the determinant yields%
\begin{align*}
\det\left(  \left(  x_{i}^{n-j}\right)  _{1\leq i\leq n,\ 1\leq j\leq
n}\right)   &  =\sum_{\sigma\in S_{n}}\left(  -1\right)  ^{\sigma}%
x_{1}^{n-\sigma\left(  1\right)  }x_{2}^{n-\sigma\left(  2\right)  }\cdots
x_{n}^{n-\sigma\left(  n\right)  }\ \ \ \ \ \ \ \ \ \ \text{and}\\
\det\left(  \left(  a_{i}^{n-j}\right)  _{1\leq i\leq n,\ 1\leq j\leq
n}\right)   &  =\sum_{\sigma\in S_{n}}\left(  -1\right)  ^{\sigma}%
a_{1}^{n-\sigma\left(  1\right)  }a_{2}^{n-\sigma\left(  2\right)  }\cdots
a_{n}^{n-\sigma\left(  n\right)  },
\end{align*}
and it is clear that substituting $a_{1},a_{2},\ldots,a_{n}$ for $x_{1}%
,x_{2},\ldots,x_{n}$ transforms \newline$\sum_{\sigma\in S_{n}}\left(
-1\right)  ^{\sigma}x_{1}^{n-\sigma\left(  1\right)  }x_{2}^{n-\sigma\left(
2\right)  }\cdots x_{n}^{n-\sigma\left(  n\right)  }$ into $\sum_{\sigma\in
S_{n}}\left(  -1\right)  ^{\sigma}a_{1}^{n-\sigma\left(  1\right)  }%
a_{2}^{n-\sigma\left(  2\right)  }\cdots a_{n}^{n-\sigma\left(  n\right)  }$.
\end{fineprint}

Thus, (\ref{pf.thm.det.vander.a.1}) follows from
(\ref{eq.lem.det.vander.a.pol.1}). In other words, Theorem
\ref{thm.det.vander} \textbf{(a)} follows from Lemma
\ref{lem.det.vander.a.pol}.
\end{proof}

Arguments like the one we just used are frequently applied in algebra; see
\cite{Conrad-ui} for some more examples.

It now remains to prove Lemma \ref{lem.det.vander.a.pol}.

\begin{proof}
[Proof of Lemma \ref{lem.det.vander.a.pol} (sketched).]We set%
\[
f:=\det\left(  \left(  x_{i}^{n-j}\right)  _{1\leq i\leq n,\ 1\leq j\leq
n}\right)  \ \ \ \ \ \ \ \ \ \ \text{and}\ \ \ \ \ \ \ \ \ \ g:=\prod_{1\leq
i<j\leq n}\left(  x_{i}-x_{j}\right)  .
\]
Thus, we must prove that $f=g$.

We have%
\begin{align}
f  &  =\det\left(  \left(  x_{i}^{n-j}\right)  _{1\leq i\leq n,\ 1\leq j\leq
n}\right) \nonumber\\
&  =\sum_{\sigma\in S_{n}}\left(  -1\right)  ^{\sigma}x_{1}^{n-\sigma\left(
1\right)  }x_{2}^{n-\sigma\left(  2\right)  }\cdots x_{n}^{n-\sigma\left(
n\right)  } \label{pf.lem.det.vander.a.pol.f=1}%
\end{align}
(by the definition of a determinant). The right hand side of this equality is
a homogeneous polynomial in $x_{1},x_{2},\ldots,x_{n}$ of degree
$\dfrac{n\left(  n-1\right)  }{2}$\ \ \ \ \footnote{because it is a
$\mathbb{Z}$-linear combination of the monomials $x_{1}^{n-\sigma\left(
1\right)  }x_{2}^{n-\sigma\left(  2\right)  }\cdots x_{n}^{n-\sigma\left(
n\right)  }$, each of which has degree%
\begin{align*}
\left(  n-\sigma\left(  1\right)  \right)  +\left(  n-\sigma\left(  2\right)
\right)  +\cdots+\left(  n-\sigma\left(  n\right)  \right)   &  =n\cdot
n-\underbrace{\left(  \sigma\left(  1\right)  +\sigma\left(  2\right)
+\cdots+\sigma\left(  n\right)  \right)  }_{\substack{=1+2+\cdots
+n\\\text{(since }\sigma\text{ is a permutation of }\left[  n\right]
\text{)}}}\\
&  =n\cdot n-\underbrace{\left(  1+2+\cdots+n\right)  }_{=\dfrac{n\left(
n+1\right)  }{2}}=n\cdot n-\dfrac{n\left(  n+1\right)  }{2}\\
&  =\dfrac{n\left(  n-1\right)  }{2}%
\end{align*}
}. Thus, $f$ is a homogeneous polynomial in $x_{1},x_{2},\ldots,x_{n}$ of
degree $\dfrac{n\left(  n-1\right)  }{2}$. Furthermore, the monomial
$x_{1}^{n-1}x_{2}^{n-2}\cdots x_{n}^{n-n}$ appears with coefficient $1$ on the
right hand side of (\ref{pf.lem.det.vander.a.pol.f=1}) (indeed, all the $n!$
addends in the sum on this right hand side contain distinct monomials, and
thus only the addend for $\sigma=\operatorname*{id}$ makes any contribution to
the coefficient of the monomial $x_{1}^{n-1}x_{2}^{n-2}\cdots x_{n}^{n-n}$).
Hence,%
\begin{equation}
\left[  x_{1}^{n-1}x_{2}^{n-2}\cdots x_{n}^{n-n}\right]  f=1.
\label{pf.lem.det.vander.a.pol.f-LC}%
\end{equation}
Therefore, $f\neq0$.

Now, let $u$ and $v$ be two elements of $\left[  n\right]  $ satisfying $u<v$.
Then, the polynomial $f$ becomes $0$ when we set $x_{u}$ equal to $x_{v}$
(that is, when we substitute $x_{v}$ for $x_{u}$). Indeed, when we set $x_{u}$
equal to $x_{v}$, the matrix $\left(  x_{i}^{n-j}\right)  _{1\leq i\leq
n,\ 1\leq j\leq n}$ becomes a matrix that has two equal rows (namely, its
$u$-th and its $v$-th row both become equal to $\left(  x_{v}^{n-1}%
,x_{v}^{n-2},\ldots,x_{v}^{n-n}\right)  $), and thus its determinant becomes
$0$ (by Theorem \ref{thm.det.rowop} \textbf{(c)}); but this means precisely
that $f$ becomes $0$ (since $f$ is the determinant of this matrix).

Now, we recall the following well-known property of univariate polynomials:

\begin{statement}
\textit{Root factoring-off theorem:} Let $R$ be a commutative ring. Let $p\in
R\left[  t\right]  $ be a univariate polynomial that has a root $r\in R$.
Then, this polynomial $p$ is divisible by $t-r$ (in the ring $R\left[
t\right]  $).
\end{statement}

Using this property, we can easily see that if a polynomial $p\in
\mathbb{Z}\left[  x_{1},x_{2},\ldots,x_{n}\right]  $ becomes $0$ when we set
$x_{u}$ equal to $x_{v}$, then this polynomial $p$ is a multiple of
$x_{u}-x_{v}$ (in the ring $\mathbb{Z}\left[  x_{1},x_{2},\ldots,x_{n}\right]
$). (Indeed, viewing $p$ as a polynomial in the single indeterminate $x_{u}$
over the ring $\mathbb{Z}\left[  x_{1},x_{2},\ldots,\widehat{x_{u}}%
,\ldots,x_{n}\right]  $\ \ \ \ \footnote{The meaning of the hat over the
\textquotedblleft$x_{u}$\textquotedblright\ is as in Proposition
\ref{prop.det.x+ai}: It signifies that the entry $x_{u}$ is omitted from the
list (that is, \textquotedblleft$x_{1},x_{2},\ldots,\widehat{x_{u}}%
,\ldots,x_{n}$\textquotedblright\ means \textquotedblleft$x_{1},x_{2}%
,\ldots,x_{u-1},x_{u+1},x_{u+2},\ldots,x_{n}$\textquotedblright).}, we see
that $x_{v}$ is a root of $p$, and therefore the root factoring-off theorem
yields that $p$ is divisible by $x_{u}-x_{v}$.\ \ \ \ \footnote{Here, we are
tacitly using the canonical isomorphism between the polynomial rings
\[
\mathbb{Z}\left[  x_{1},x_{2},\ldots,x_{n}\right]
\ \ \ \ \ \ \ \ \ \ \text{and}\ \ \ \ \ \ \ \ \ \ \left(  \mathbb{Z}\left[
x_{1},x_{2},\ldots,\widehat{x_{u}},\ldots,x_{n}\right]  \right)  \left[
x_{u}\right]  .
\]
This isomorphism allows us to treat multivariate polynomials in $\mathbb{Z}%
\left[  x_{1},x_{2},\ldots,x_{n}\right]  $ as univariate polynomials in the
indeterminate $x_{u}$ over the ring $\mathbb{Z}\left[  x_{1},x_{2}%
,\ldots,\widehat{x_{u}},\ldots,x_{n}\right]  $, and vice versa (and ensures
that the notion of divisibility does not change between the former and the
latter).})

Applying this observation to $p=f$, we conclude that $f$ is a multiple of
$x_{u}-x_{v}$ (in the ring $\mathbb{Z}\left[  x_{1},x_{2},\ldots,x_{n}\right]
$), since we know that $f$ becomes $0$ when we set $x_{u}$ equal to $x_{v}$.

Forget that we fixed $u$ and $v$. We thus have shown that $f$ is a multiple of
$x_{u}-x_{v}$ (in the ring $\mathbb{Z}\left[  x_{1},x_{2},\ldots,x_{n}\right]
$) whenever $u$ and $v$ are two elements of $\left[  n\right]  $ satisfying
$u<v$. In other words, $f$ is a multiple of each of the linear polynomials%
\[%
\begin{array}
[c]{cccc}%
x_{1}-x_{2}, & x_{1}-x_{3}, & \ldots, & x_{1}-x_{n},\\
& x_{2}-x_{3}, & \ldots, & x_{2}-x_{n},\\
&  & \ddots & \vdots\\
&  &  & x_{n-1}-x_{n}%
\end{array}
\]
(in the ring $\mathbb{Z}\left[  x_{1},x_{2},\ldots,x_{n}\right]  $). However,
these linear polynomials are irreducible\footnote{Check this! (Actually, any
linear polynomial over $\mathbb{Z}$ is irreducible if the gcd of its
coefficients is $1$.)} and mutually non-associate (i.e., no two of them are
associate\footnote{Recall that two elements $a$ and $b$ of a commutative ring
$R$ are said to be \emph{associate} if there exists some unit $u$ of $R$ such
that $a=ub$. Being associate is known (and easily verified) to be an
equivalence relation.})\footnote{Check this! (Recall that the only units of
the polynomial ring $\mathbb{Z}\left[  x_{1},x_{2},\ldots,x_{n}\right]  $ are
$1$ and $-1$.)}. Since the ring $\mathbb{Z}\left[  x_{1},x_{2},\ldots
,x_{n}\right]  $ is a unique factorization domain\footnote{This is a
nontrivial, but rather well-known result in abstract algebra. Proofs can be
found in \cite[Theorem 3.7.4]{Ford21}, \cite[Remark after Corollary
8.21]{Knapp16}, \cite[Chapter IV, Theorems 4.8 and 4.9]{MiRiRu87} and
\cite[Essay 1.4, Corollary of Theorem 1 and Corollary 1 of Theorem
2]{Edward22}.}, we thus conclude that any polynomial $p\in\mathbb{Z}\left[
x_{1},x_{2},\ldots,x_{n}\right]  $ that is a multiple of each of these linear
polynomials must necessarily be a multiple of their product $\prod_{1\leq
i<j\leq n}\left(  x_{i}-x_{j}\right)  $. Hence, the polynomial $f$ must be a
multiple of this product $\prod_{1\leq i<j\leq n}\left(  x_{i}-x_{j}\right)  $
(since $f$ is a multiple of each of the linear polynomials above). In other
words, the polynomial $f$ must be a multiple of $g$ (since $g=\prod_{1\leq
i<j\leq n}\left(  x_{i}-x_{j}\right)  $). In other words, $f=gq$ for some
$q\in\mathbb{Z}\left[  x_{1},x_{2},\ldots,x_{n}\right]  $. Consider this $q$.
Clearly, $gq=f\neq0$ and thus $q\neq0$ and $g\neq0$.

The ring $\mathbb{Z}$ is an integral domain. Thus, any two nonzero polynomials
$a$ and $b$ over $\mathbb{Z}$ satisfy $\deg\left(  ab\right)  =\deg a+\deg b$.
Applying this to $a=g$ and $b=q$, we find $\deg\left(  gq\right)  =\deg g+\deg
q$. In other words, $\deg f=\deg g+\deg q$ (since $gq=f$).

Now, $g=\prod_{1\leq i<j\leq n}\left(  x_{i}-x_{j}\right)  $ and thus%
\begin{align*}
\deg g  &  =\deg\left(  \prod_{1\leq i<j\leq n}\left(  x_{i}-x_{j}\right)
\right)  =\sum_{1\leq i<j\leq n}\underbrace{\deg\left(  x_{i}-x_{j}\right)
}_{=1}=\sum_{1\leq i<j\leq n}1\\
&  =\left(  \text{\# of pairs }\left(  i,j\right)  \in\left[  n\right]
^{2}\text{ satisfying }i<j\right)  =\dbinom{n}{2}=\dfrac{n\left(  n-1\right)
}{2}.
\end{align*}
However, we have $\deg f\leq\dfrac{n\left(  n-1\right)  }{2}$ (since $f$ is a
homogeneous polynomial of degree $\dfrac{n\left(  n-1\right)  }{2}$). In view
of $\deg g=\dfrac{n\left(  n-1\right)  }{2}$, this rewrites as $\deg f\leq\deg
g$. Hence,%
\[
\deg g\geq\deg f=\deg g+\deg q.
\]
Therefore, $0\geq\deg q$. This shows that the polynomial $q$ is constant. In
other words, $q\in\mathbb{Z}$.

It remains to show that this constant $q$ is $1$. However, this can be done by
comparing some coefficients of $f$ and $g$. Indeed, let us look at the
coefficient of $x_{1}^{n-1}x_{2}^{n-2}\cdots x_{n}^{n-n}$ (as we already know
this coefficient for $f$ to be $1$).

Expanding the product $\prod_{1\leq i<j\leq n}\left(  x_{i}-x_{j}\right)  $,
we obtain a sum of several (in fact, $2^{n\left(  n-1\right)  /2}$ many)
monomials with $+$ and $-$ signs. I claim that among these monomials, the
monomial $x_{1}^{n-1}x_{2}^{n-2}\cdots x_{n}^{n-n}$ will appear exactly once,
and with a $+$ sign. Indeed, in order to obtain $x_{1}^{n-1}x_{2}^{n-2}\cdots
x_{n}^{n-n}$ when expanding the product
\[%
\begin{array}
[c]{ccccc}%
\prod_{1\leq i<j\leq n}\left(  x_{i}-x_{j}\right)  = & \left(  x_{1}%
-x_{2}\right)  & \left(  x_{1}-x_{3}\right)  & \cdots & \left(  x_{1}%
-x_{n}\right) \\
&  & \left(  x_{2}-x_{3}\right)  & \cdots & \left(  x_{2}-x_{n}\right) \\
&  &  & \ddots & \vdots\\
&  &  &  & \left(  x_{n-1}-x_{n}\right)  ,
\end{array}
\]
it is necessary to pick the $x_{1}$ minuends from all $n-1$ factors in the
first row (since none of the other factors contain any $x_{1}$), then to pick
the $x_{2}$ minuends from all $n-2$ factors in the second row (since none of
the remaining factors contain any $x_{2}$), and so on -- i.e., to take the
minuend (rather than the subtrahend) from each factor. Thus, only one of the
monomials obtained by the expansion will be $x_{1}^{n-1}x_{2}^{n-2}\cdots
x_{n}^{n-n}$, and it will appear with a $+$ sign. Hence, the coefficient of
$x_{1}^{n-1}x_{2}^{n-2}\cdots x_{n}^{n-n}$ in the product $\prod_{1\leq
i<j\leq n}\left(  x_{i}-x_{j}\right)  $ is $1$. In other words, the
coefficient of $x_{1}^{n-1}x_{2}^{n-2}\cdots x_{n}^{n-n}$ in $g$ is $1$ (since
$g=\prod_{1\leq i<j\leq n}\left(  x_{i}-x_{j}\right)  $). In other words,%
\[
\left[  x_{1}^{n-1}x_{2}^{n-2}\cdots x_{n}^{n-n}\right]  g=1.
\]

Now, recall that $f=gq$. Hence,%
\begin{align*}
\left[  x_{1}^{n-1}x_{2}^{n-2}\cdots x_{n}^{n-n}\right]  f  &  =\left[
x_{1}^{n-1}x_{2}^{n-2}\cdots x_{n}^{n-n}\right]  \left(  gq\right) \\
&  =q\cdot\underbrace{\left[  x_{1}^{n-1}x_{2}^{n-2}\cdots x_{n}^{n-n}\right]
g}_{=1}\ \ \ \ \ \ \ \ \ \ \left(  \text{since }q\in\mathbb{Z}\right) \\
&  =q,
\end{align*}
so that $q=\left[  x_{1}^{n-1}x_{2}^{n-2}\cdots x_{n}^{n-n}\right]  f=1$ (by
(\ref{pf.lem.det.vander.a.pol.f-LC})). Hence, $f=g\underbrace{q}_{=1}=g$. As
we said, this completes the proof of Lemma \ref{lem.det.vander.a.pol}.
\end{proof}

The technique used in the above proof of Lemma \ref{lem.det.vander.a.pol} may
appear somewhat underhanded: Instead of computing our polynomial $f$ upfront,
we have kept lopping off linear factors from it until a constant polynomial
remained (for degree reasons). This technique is called \emph{identification
of factors} or \emph{factor hunting}, and is used in various different places,
but particularly often in the computation of determinants (multiple examples
are given in \cite[\S 2.4 and further below]{Krattenthaler}). While I consider
it to be aesthetically inferior to sufficiently direct approaches, it has
proved itself useful in situations in which no direct approaches are known to work.

\begin{proof}
[Proof of Theorem \ref{thm.det.vander} (sketched).]\textbf{(a)} We have
already given a proof of Theorem \ref{thm.det.vander} \textbf{(a)} (and with
Lemma \ref{lem.det.vander.a.pol} established, this proof is now complete).
\medskip

\textbf{(b)} The matrix $\left(  a_{j}^{n-i}\right)  _{1\leq i\leq n,\ 1\leq
j\leq n}$ is the transpose of the matrix $\left(  a_{i}^{n-j}\right)  _{1\leq
i\leq n,\ 1\leq j\leq n}$. Thus, Theorem \ref{thm.det.vander} \textbf{(b)}
follows from Theorem \ref{thm.det.vander} \textbf{(a)} using Theorem
\ref{thm.det.transp}. \medskip

\textbf{(c)} Let $A$ be the $n\times n$-matrix $\left(  a_{i}^{n-j}\right)
_{1\leq i\leq n,\ 1\leq j\leq n}$. Then, Theorem \ref{thm.det.vander}
\textbf{(a)} says that $\det A=\prod_{1\leq i<j\leq n}\left(  a_{i}%
-a_{j}\right)  $.

Let $\tau\in S_{n}$ be the permutation of $\left[  n\right]  $ that sends each
$j\in\left[  n\right]  $ to $n+1-j$. Then, each $i,j\in\left[  n\right]  $
satisfy $j-1=n-\tau\left(  j\right)  $ and therefore
\[
a_{i}^{j-1}=a_{i}^{n-\tau\left(  j\right)  }=A_{i,\tau\left(  j\right)
}\ \ \ \ \ \ \ \ \ \ \left(  \text{by the definition of }A\right)  .
\]
Hence, $\left(  a_{i}^{j-1}\right)  _{1\leq i\leq n,\ 1\leq j\leq n}=\left(
A_{i,\tau\left(  j\right)  }\right)  _{1\leq i\leq n,\ 1\leq j\leq n}$ and
thus%
\[
\det\left(  \left(  a_{i}^{j-1}\right)  _{1\leq i\leq n,\ 1\leq j\leq
n}\right)  =\det\left(  \left(  A_{i,\tau\left(  j\right)  }\right)  _{1\leq
i\leq n,\ 1\leq j\leq n}\right)  =\left(  -1\right)  ^{\tau}\cdot\det A
\]
(by (\ref{eq.cor.det.sig-row-col.col})).

However, the permutation $\tau$ has OLN $\left(  n,n-1,n-2,\ldots,1\right)  $.
Thus, each pair $\left(  i,j\right)  \in\left[  n\right]  ^{2}$ satisfying
$i<j$ is an inversion of $\tau$. Therefore, the Coxeter length $\ell\left(
\tau\right)  $ of $\tau$ is the \# of all pairs $\left(  i,j\right)
\in\left[  n\right]  ^{2}$ satisfying $i<j$. Thus, the sign of $\tau$ is%
\begin{align*}
\left(  -1\right)  ^{\tau}  &  =\left(  -1\right)  ^{\ell\left(  \tau\right)
}\ \ \ \ \ \ \ \ \ \ \left(  \text{by the definition of the sign of a
permutation}\right) \\
&  =\left(  -1\right)  ^{\left(  \text{\# of all pairs }\left(  i,j\right)
\in\left[  n\right]  ^{2}\text{ satisfying }i<j\right)  }\\
&  \ \ \ \ \ \ \ \ \ \ \ \ \ \ \ \ \ \ \ \ \left(  \text{since }\ell\left(
\tau\right)  \text{ is the \# of all pairs }\left(  i,j\right)  \in\left[
n\right]  ^{2}\text{ satisfying }i<j\right) \\
&  =\prod_{1\leq i<j\leq n}\left(  -1\right)  .
\end{align*}

Now,
\begin{align*}
&  \det\left(  \left(  a_{i}^{j-1}\right)  _{1\leq i\leq n,\ 1\leq j\leq
n}\right) \\
&  =\underbrace{\left(  -1\right)  ^{\tau}}_{=\prod_{1\leq i<j\leq n}\left(
-1\right)  }\cdot\underbrace{\det A}_{=\prod_{1\leq i<j\leq n}\left(
a_{i}-a_{j}\right)  }=\left(  \prod_{1\leq i<j\leq n}\left(  -1\right)
\right)  \cdot\left(  \prod_{1\leq i<j\leq n}\left(  a_{i}-a_{j}\right)
\right) \\
&  =\prod_{1\leq i<j\leq n}\underbrace{\left(  \left(  -1\right)  \left(
a_{i}-a_{j}\right)  \right)  }_{=a_{j}-a_{i}}=\prod_{1\leq i<j\leq n}\left(
a_{j}-a_{i}\right)  =\prod_{1\leq j<i\leq n}\left(  a_{i}-a_{j}\right)
\end{align*}
(here, we have renamed the index $\left(  i,j\right)  $ as $\left(
j,i\right)  $ in the product). This proves Theorem \ref{thm.det.vander}
\textbf{(c)}. \medskip

\textbf{(d)} The matrix $\left(  a_{j}^{i-1}\right)  _{1\leq i\leq n,\ 1\leq
j\leq n}$ is the transpose of the matrix $\left(  a_{i}^{j-1}\right)  _{1\leq
i\leq n,\ 1\leq j\leq n}$. Thus, Theorem \ref{thm.det.vander} \textbf{(d)}
follows from Theorem \ref{thm.det.vander} \textbf{(c)} using Theorem
\ref{thm.det.transp}.
\end{proof}

The Vandermonde determinant is itself a useful tool in the computation of
various other determinants. Here is an example:

\begin{proposition}
\label{prop.det.(xi+yj)n-1}Let $n\in\mathbb{N}$. Let $x_{1},x_{2},\ldots
,x_{n}\in K$ and $y_{1},y_{2},\ldots,y_{n}\in K$. Then,%
\begin{align*}
&  \det\left(  \left(  \left(  x_{i}+y_{j}\right)  ^{n-1}\right)  _{1\leq
i\leq n,\ 1\leq j\leq n}\right) \\
&  =\det\left(
\begin{array}
[c]{cccc}%
\left(  x_{1}+y_{1}\right)  ^{n-1} & \left(  x_{1}+y_{2}\right)  ^{n-1} &
\cdots & \left(  x_{1}+y_{n}\right)  ^{n-1}\\
\left(  x_{2}+y_{1}\right)  ^{n-1} & \left(  x_{2}+y_{2}\right)  ^{n-1} &
\cdots & \left(  x_{2}+y_{n}\right)  ^{n-1}\\
\vdots & \vdots & \ddots & \vdots\\
\left(  x_{n}+y_{1}\right)  ^{n-1} & \left(  x_{n}+y_{2}\right)  ^{n-1} &
\cdots & \left(  x_{n}+y_{n}\right)  ^{n-1}%
\end{array}
\right) \\
&  =\left(  \prod_{k=0}^{n-1}\dbinom{n-1}{k}\right)  \left(  \prod_{1\leq
i<j\leq n}\left(  x_{i}-x_{j}\right)  \right)  \left(  \prod_{1\leq i<j\leq
n}\left(  y_{j}-y_{i}\right)  \right)  .
\end{align*}

\end{proposition}

\begin{proof}
[First proof of Proposition \ref{prop.det.(xi+yj)n-1} (sketched).]Here is a
rough outline of a proof that uses factor hunting (in the same way as in our
above proofs of Theorem \ref{thm.det.vander} \textbf{(a)} and Lemma
\ref{lem.det.vander.a.pol}). We WLOG assume that $x_{1},x_{2},\ldots,x_{n}$
and $y_{1},y_{2},\ldots,y_{n}$ are distinct indeterminates in a polynomial
ring over $\mathbb{Z}$. (This is again an assumption that we can make, because
the argument that we used to derive Theorem \ref{thm.det.vander} \textbf{(a)}
from Lemma \ref{lem.det.vander.a.pol} can be applied here as well.) Then, we
can easily see that%
\[
\det\left(  \left(  \left(  x_{i}+y_{j}\right)  ^{n-1}\right)  _{1\leq i\leq
n,\ 1\leq j\leq n}\right)
\]
is a homogeneous polynomial of degree $n\left(  n-1\right)  $. This polynomial
vanishes if we set any $x_{u}$ equal to any $x_{v}$ (for $u<v$), and also
vanishes if we set any $y_{u}$ equal to any $y_{v}$ (for $u<v$). Thus we have
identified $n\left(  n-1\right)  $ linear factors of this polynomial (namely,
the differences $x_{u}-x_{v}$ and $y_{v}-y_{u}$ for $u<v$), and we can again
conclude (since any polynomial ring over $\mathbb{Z}$ is a unique
factorization domain) that
\[
\det\left(  \left(  \left(  x_{i}+y_{j}\right)  ^{n-1}\right)  _{1\leq i\leq
n,\ 1\leq j\leq n}\right)  =\left(  \prod_{1\leq i<j\leq n}\left(  x_{i}%
-x_{j}\right)  \right)  \left(  \prod_{1\leq i<j\leq n}\left(  y_{j}%
-y_{i}\right)  \right)  \cdot q
\]
for a constant $q\in\mathbb{Z}$. It remains to prove that this constant $q$
equals $\prod_{k=0}^{n-1}\dbinom{n-1}{k}$. This can be done by studying the
coefficients of the monomial%
\[
x_{1}^{n-1}x_{2}^{n-2}\cdots x_{n}^{n-n}y_{1}^{0}y_{2}^{1}\cdots y_{n}^{n-1}%
\]
in $\det\left(  \left(  \left(  x_{i}+y_{j}\right)  ^{n-1}\right)  _{1\leq
i\leq n,\ 1\leq j\leq n}\right)  $ and in $\left(  \prod_{1\leq i<j\leq
n}\left(  x_{i}-x_{j}\right)  \right)  \left(  \prod_{1\leq i<j\leq n}\left(
y_{j}-y_{i}\right)  \right)  $. We leave the details to the reader.
\end{proof}

\begin{proof}
[Second proof of Proposition \ref{prop.det.(xi+yj)n-1}.](See \cite[Exercise
6.17 \textbf{(b)}]{detnotes} for details.) Let
\[
C:=\left(  \left(  x_{i}+y_{j}\right)  ^{n-1}\right)  _{1\leq i\leq n,\ 1\leq
j\leq n}\in K^{n\times n}.
\]
For any $i,j\in\left[  n\right]  $, we have%
\begin{align*}
C_{i,j}  &  =\left(  x_{i}+y_{j}\right)  ^{n-1}=\sum_{k=0}^{n-1}\dbinom
{n-1}{k}x_{i}^{k}y_{j}^{n-1-k}\ \ \ \ \ \ \ \ \ \ \left(  \text{by the
binomial formula}\right) \\
&  =\sum_{k=1}^{n}\dbinom{n-1}{k-1}x_{i}^{k-1}y_{j}^{n-k}%
\end{align*}
(here, we have substituted $k-1$ for $k$ in the sum). If we define two
$n\times n$-matrices $P$ and $Q$ by
\[
P:=\left(  \dbinom{n-1}{k-1}x_{i}^{k-1}\right)  _{1\leq i\leq n,\ 1\leq k\leq
n}\ \ \ \ \ \ \ \ \ \ \text{and}\ \ \ \ \ \ \ \ \ \ Q:=\left(  y_{j}%
^{n-k}\right)  _{1\leq k\leq n,\ 1\leq j\leq n},
\]
then we can rewrite this as%
\[
C_{i,j}=\sum_{k=1}^{n}P_{i,k}Q_{k,j}=\left(  PQ\right)  _{i,j}%
\]
(by the definition of the matrix product). Since this holds for all
$i,j\in\left[  n\right]  $, we thus obtain $C=PQ$. Hence,
\begin{equation}
\det C=\det\left(  PQ\right)  =\det P\cdot\det Q
\label{pf.prop.det.(xi+yj)n-1.2nd.detA=1}%
\end{equation}
(by Theorem \ref{thm.det.detAB}). It now remains to compute $\det P$ and $\det
Q$.

From $Q=\left(  y_{j}^{n-k}\right)  _{1\leq k\leq n,\ 1\leq j\leq n}=\left(
y_{j}^{n-i}\right)  _{1\leq i\leq n,\ 1\leq j\leq n}$, we obtain%
\begin{equation}
\det Q=\det\left(  \left(  y_{j}^{n-i}\right)  _{1\leq i\leq n,\ 1\leq j\leq
n}\right)  =\prod_{1\leq i<j\leq n}\left(  y_{i}-y_{j}\right)
\label{pf.prop.det.(xi+yj)n-1.2nd.detQ=}%
\end{equation}
(by Theorem \ref{thm.det.vander} \textbf{(b)}, applied to $a_{i}=y_{i}$).

From $P=\left(  \dbinom{n-1}{k-1}x_{i}^{k-1}\right)  _{1\leq i\leq n,\ 1\leq
k\leq n}=\left(  \dbinom{n-1}{j-1}x_{i}^{j-1}\right)  _{1\leq i\leq n,\ 1\leq
j\leq n}$, we obtain%
\begin{align}
\det P  &  =\det\left(  \left(  \dbinom{n-1}{j-1}x_{i}^{j-1}\right)  _{1\leq
i\leq n,\ 1\leq j\leq n}\right) \nonumber\\
&  =\underbrace{\dbinom{n-1}{1-1}\dbinom{n-1}{2-1}\cdots\dbinom{n-1}{n-1}%
}_{=\prod_{k=0}^{n-1}\dbinom{n-1}{k}}\cdot\underbrace{\det\left(  \left(
x_{i}^{j-1}\right)  _{1\leq i\leq n,\ 1\leq j\leq n}\right)  }%
_{\substack{=\prod_{1\leq j<i\leq n}\left(  x_{i}-x_{j}\right)  \\\text{(by
Theorem \ref{thm.det.vander} \textbf{(c)},}\\\text{applied to }a_{i}%
=x_{i}\text{)}}}\nonumber\\
&  \ \ \ \ \ \ \ \ \ \ \ \ \ \ \ \ \ \ \ \ \left(
\begin{array}
[c]{c}%
\text{by (\ref{eq.cor.det.scale-row-col.col}), applied to }A=\left(
x_{i}^{j-1}\right)  _{1\leq i\leq n,\ 1\leq j\leq n}\\
\text{and }d_{i}=\dbinom{n-1}{i-1}%
\end{array}
\right) \nonumber\\
&  =\left(  \prod_{k=0}^{n-1}\dbinom{n-1}{k}\right)  \cdot\prod_{1\leq j<i\leq
n}\left(  x_{i}-x_{j}\right) \nonumber\\
&  =\left(  \prod_{k=0}^{n-1}\dbinom{n-1}{k}\right)  \cdot\prod_{1\leq i<j\leq
n}\left(  x_{j}-x_{i}\right)  \label{pf.prop.det.(xi+yj)n-1.2nd.detP=}%
\end{align}
(here, we have renamed the index $\left(  j,i\right)  $ as $\left(
i,j\right)  $ in the second product).

Now, (\ref{pf.prop.det.(xi+yj)n-1.2nd.detA=1}) becomes%
\begin{align*}
\det C  &  =\det P\cdot\det Q\\
&  =\left(  \prod_{k=0}^{n-1}\dbinom{n-1}{k}\right)  \cdot\left(  \prod_{1\leq
i<j\leq n}\left(  x_{j}-x_{i}\right)  \right)  \cdot\prod_{1\leq i<j\leq
n}\left(  y_{i}-y_{j}\right) \\
&  \ \ \ \ \ \ \ \ \ \ \ \ \ \ \ \ \ \ \ \ \left(  \text{by
(\ref{pf.prop.det.(xi+yj)n-1.2nd.detP=}) and
(\ref{pf.prop.det.(xi+yj)n-1.2nd.detQ=})}\right) \\
&  =\left(  \prod_{k=0}^{n-1}\dbinom{n-1}{k}\right)  \cdot\prod_{1\leq i<j\leq
n}\underbrace{\left(  \left(  x_{j}-x_{i}\right)  \left(  y_{i}-y_{j}\right)
\right)  }_{=\left(  x_{i}-x_{j}\right)  \left(  y_{j}-y_{i}\right)  }\\
&  =\left(  \prod_{k=0}^{n-1}\dbinom{n-1}{k}\right)  \cdot\prod_{1\leq i<j\leq
n}\left(  \left(  x_{i}-x_{j}\right)  \left(  y_{j}-y_{i}\right)  \right) \\
&  =\left(  \prod_{k=0}^{n-1}\dbinom{n-1}{k}\right)  \left(  \prod_{1\leq
i<j\leq n}\left(  x_{i}-x_{j}\right)  \right)  \left(  \prod_{1\leq i<j\leq
n}\left(  y_{j}-y_{i}\right)  \right)  .
\end{align*}
This proves Proposition \ref{prop.det.(xi+yj)n-1} (since $C=\left(  \left(
x_{i}+y_{j}\right)  ^{n-1}\right)  _{1\leq i\leq n,\ 1\leq j\leq n}$).
\end{proof}

\subsubsection{Laplace expansion}

Let us next recall another fundamental property of determinants: \emph{Laplace
expansion}. We will use the following notation:

\begin{convention}
\label{conv.mat.tilde}Let $n\in\mathbb{N}$. Let $A$ be an $n\times n$-matrix.
Let $i,j\in\left[  n\right]  $. Then, we set%
\[
A_{\sim i,\sim j}:=\operatorname*{sub}\nolimits_{\left[  n\right]
\setminus\left\{  i\right\}  }^{\left[  n\right]  \setminus\left\{  j\right\}
}A\ \ \ \ \ \ \ \ \ \ \left(  \text{using the notation from Definition
\ref{def.det.sub}}\right)  .
\]
This is the $\left(  n-1\right)  \times\left(  n-1\right)  $-matrix obtained
from $A$ by removing its $i$-th row and its $j$-th column.
\end{convention}

\begin{example}
If $A=\left(
\begin{array}
[c]{ccc}%
a & b & c\\
a^{\prime} & b^{\prime} & c^{\prime}\\
a^{\prime\prime} & b^{\prime\prime} & c^{\prime\prime}%
\end{array}
\right)  $, then $A_{\sim1,\sim2}=\left(
\begin{array}
[c]{cc}%
a^{\prime} & c^{\prime}\\
a^{\prime\prime} & c^{\prime\prime}%
\end{array}
\right)  $.
\end{example}

Now, we can state the theorem underlying Laplace expansion:

\begin{theorem}
\label{thm.det.laplace}Let $n\in\mathbb{N}$. Let $A\in K^{n\times n}$ be an
$n\times n$-matrix. \medskip

\textbf{(a)} For every $p\in\left[  n\right]  $, we have%
\[
\det A=\sum_{q=1}^{n}\left(  -1\right)  ^{p+q}A_{p,q}\det\left(  A_{\sim
p,\sim q}\right)  .
\]

\textbf{(b)} For every $q\in\left[  n\right]  $, we have%
\[
\det A=\sum_{p=1}^{n}\left(  -1\right)  ^{p+q}A_{p,q}\det\left(  A_{\sim
p,\sim q}\right)  .
\]

\end{theorem}

Note that some authors denote the minors $\det\left(  A_{\sim p,\sim
q}\right)  $ in Theorem \ref{thm.det.laplace} by $A_{p,q}$. This is, of
course, totally incompatible with our notations.

\begin{proof}
[Proof of Theorem \ref{thm.det.laplace}.]See \cite[Theorem 6.82]{detnotes} or
\cite[Lemma 4.7.7]{Ford21} or \cite[5.8 and 5.8']{Laue-det} or
\cite[Proposition B.24 and Proposition B.25]{Strick13} or \cite[Theorem
9.48]{Loehr-BC}.
\end{proof}

Theorem \ref{thm.det.laplace} yields several ways to compute the determinant
of a matrix\footnote{In fact, some authors use Theorem \ref{thm.det.laplace}
as a \textbf{definition} of the determinant. (However, this is somewhat
tricky, as it requires proving that all the values obtained for $\det A$ by
applying Theorem \ref{thm.det.laplace} are actually equal.)}. When we compute
a determinant $\det A$ using Theorem \ref{thm.det.laplace} \textbf{(a)}, we
say that we \emph{expand this determinant along the }$p$\emph{-th row}. When
we compute a determinant $\det A$ using Theorem \ref{thm.det.laplace}
\textbf{(b)}, we say that we \emph{expand this determinant along the }%
$q$\emph{-th column}.

Theorem \ref{thm.det.laplace} has many applications, some of which you have
probably seen in your course on linear algebra. (A few might appear on the
homework.) The main theoretical application of Theorem \ref{thm.det.laplace}
is the concept of the \emph{adjugate} (or \emph{classical adjoint}) of a
matrix, which we shall introduce in a few moments.

First, let us see what happens if we replace the $A_{p,q}$ in Theorem
\ref{thm.det.laplace} by entries from another row or another column. In fact,
the respective sums become $0$ (instead of $\det A$), as the following
proposition shows:

\begin{proposition}
\label{prop.det.laplace.0}Let $n\in\mathbb{N}$. Let $A\in K^{n\times n}$ be an
$n\times n$-matrix. Let $r\in\left[  n\right]  $. \medskip

\textbf{(a)} For every $p\in\left[  n\right]  $ satisfying $p\neq r$, we have%
\[
0=\sum_{q=1}^{n}\left(  -1\right)  ^{p+q}A_{r,q}\det\left(  A_{\sim p,\sim
q}\right)  .
\]

\textbf{(b)} For every $q\in\left[  n\right]  $ satisfying $q\neq r$, we have%
\[
0=\sum_{p=1}^{n}\left(  -1\right)  ^{p+q}A_{p,r}\det\left(  A_{\sim p,\sim
q}\right)  .
\]

\end{proposition}

\begin{proof}
See \cite[Proposition 6.96]{detnotes}. This is also implicit in \cite[proof of
Proposition B.28]{Strick13} and \cite[proof of Theorem 9.50]{Loehr-BC}.
\medskip

\begin{fineprint}
Here is a sketch of the proof:

\textbf{(a)} Let $p\in\left[  n\right]  $ satisfy $p\neq r$. Let $C$ be the
matrix $A$ with its $p$-th row replaced by its $r$-th row. Then, the matrix
$C$ has two equal rows, so that $\det C=0$ (by Theorem \ref{thm.det.rowop}
\textbf{(c)}). On the other hand, expanding $\det C$ along the $p$-th row
(i.e., applying Theorem \ref{thm.det.laplace} \textbf{(a)} to $C$ instead of
$A$) yields%
\[
\det C=\sum_{q=1}^{n}\left(  -1\right)  ^{p+q}\underbrace{C_{p,q}}_{=A_{r,q}%
}\det\underbrace{\left(  C_{\sim p,\sim q}\right)  }_{=A_{\sim p,\sim q}}%
=\sum_{q=1}^{n}\left(  -1\right)  ^{p+q}A_{r,q}\det\left(  A_{\sim p,\sim
q}\right)  .
\]
Comparing these two equalities, we obtain the claim of Proposition
\ref{prop.det.laplace.0} \textbf{(a)}. A similar argument proves Proposition
\ref{prop.det.laplace.0} \textbf{(b)}.
\end{fineprint}
\end{proof}

We can now define the adjugate of a matrix:

\begin{definition}
\label{def.det.adj}Let $n\in\mathbb{N}$. Let $A\in K^{n\times n}$ be an
$n\times n$-matrix. We define a new $n\times n$-matrix $\operatorname*{adj}%
A\in K^{n\times n}$ by%
\[
\operatorname*{adj}A=\left(  \left(  -1\right)  ^{i+j}\det\left(  A_{\sim
j,\sim i}\right)  \right)  _{1\leq i\leq n,\ 1\leq j\leq n}.
\]
This matrix $\operatorname*{adj}A$ is called the \emph{adjugate} of the matrix
$A$. (Some older texts call it the \emph{adjoint}, but this name has since
been conquered by a different notion. As a compromise, some still call
$\operatorname*{adj}A$ the \emph{classical adjoint} of $A$.)
\end{definition}

\begin{example}
The adjugate of the $0\times0$-matrix is the $0\times0$-matrix.

The adjugate of a $1\times1$-matrix $\left(
\begin{array}
[c]{c}%
a
\end{array}
\right)  $ is $\operatorname*{adj}\left(
\begin{array}
[c]{c}%
a
\end{array}
\right)  =\left(
\begin{array}
[c]{c}%
1
\end{array}
\right)  $.

The adjugate of a $2\times2$-matrix $\left(
\begin{array}
[c]{cc}%
a & b\\
c & d
\end{array}
\right)  $ is
\[
\operatorname*{adj}\left(
\begin{array}
[c]{cc}%
a & b\\
c & d
\end{array}
\right)  =\left(
\begin{array}
[c]{cc}%
d & -b\\
-c & a
\end{array}
\right)  .
\]

The adjugate of a $3\times3$-matrix $\left(
\begin{array}
[c]{ccc}%
a & b & c\\
d & e & f\\
g & h & i
\end{array}
\right)  $ is
\[
\operatorname*{adj}\left(
\begin{array}
[c]{ccc}%
a & b & c\\
d & e & f\\
g & h & i
\end{array}
\right)  =\left(
\begin{array}
[c]{ccc}%
ei-fh & ch-bi & bf-ce\\
fg-di & ai-cg & cd-af\\
dh-ge & bg-ah & ae-bd
\end{array}
\right)  .
\]

\end{example}

The main property of the adjugate $\operatorname*{adj}A$ is its connection to
the inverse $A^{-1}$ of a matrix $A$. Indeed, if an $n\times n$-matrix $A$ is
invertible, then its inverse $A^{-1}$ is $\dfrac{1}{\det A}\cdot
\operatorname*{adj}A$. More generally, even if $A$ is not invertible, the
product of $\operatorname*{adj}A$ with $A$ (in either order) equals $\left(
\det A\right)  \cdot I_{n}$ (where $I_{n}$ is the $n\times n$ identity
matrix). Let us state this as a theorem:

\begin{theorem}
\label{thm.det.adj.inverse}Let $n\in\mathbb{N}$. Let $A\in K^{n\times n}$ be
an $n\times n$-matrix. Let $I_{n}$ denote the $n\times n$ identity matrix.
Then,%
\[
A\cdot\left(  \operatorname*{adj}A\right)  =\left(  \operatorname*{adj}%
A\right)  \cdot A=\left(  \det A\right)  \cdot I_{n}.
\]

\end{theorem}

\begin{proof}
See \cite[Theorem 6.100]{detnotes} or \cite[Lemma 4.7.9]{Ford21} or
\cite[Theorem 9.50]{Loehr-BC} or \cite[Proposition B.28]{Strick13}. \medskip

\begin{fineprint}
Here is a sketch of the argument: In order to show that $A\cdot\left(
\operatorname*{adj}A\right)  =\left(  \det A\right)  \cdot I_{n}$, it suffices
to check that the $\left(  i,j\right)  $-th entry of $A\cdot\left(
\operatorname*{adj}A\right)  $ equals $\det A$ whenever $i=j$ and equals $0$
otherwise. However, this follows from Theorem \ref{thm.det.laplace}
\textbf{(a)} (in the case $i=j$) and Proposition \ref{prop.det.laplace.0}
\textbf{(a)} (in the case $i\neq j$). Thus, $A\cdot\left(  \operatorname*{adj}%
A\right)  =\left(  \det A\right)  \cdot I_{n}$ is proved. Similarly, $\left(
\operatorname*{adj}A\right)  \cdot A=\left(  \det A\right)  \cdot I_{n}$ can
be shown.
\end{fineprint}
\end{proof}

More about the adjugate matrix can be found in \cite[\S 6.15]{detnotes} and
\cite[\S 5.4--\S 5.6]{trach}; see also \cite{Robins05} for some applications.

There is also a generalization of Theorem \ref{thm.det.laplace}, called
\emph{Laplace expansion along multiple rows (or columns)}:

\begin{theorem}
\label{thm.det.laplace-multi}Let $n\in\mathbb{N}$. Let $A\in K^{n\times n}$ be
an $n\times n$-matrix. We shall use the notations $\widetilde{I}$ and
$\operatorname*{sum}S$ as defined in Theorem \ref{thm.det.det(A+B)}. \medskip

\textbf{(a)} For every subset $P$ of $\left[  n\right]  $, we have%
\[
\det A=\sum_{\substack{Q\subseteq\left[  n\right]  ;\\\left\vert Q\right\vert
=\left\vert P\right\vert }}\left(  -1\right)  ^{\operatorname*{sum}%
P+\operatorname*{sum}Q}\det\left(  \operatorname*{sub}\nolimits_{P}%
^{Q}A\right)  \det\left(  \operatorname*{sub}\nolimits_{\widetilde{P}%
}^{\widetilde{Q}}A\right)  .
\]

\textbf{(b)} For every subset $Q$ of $\left[  n\right]  $, we have%
\[
\det A=\sum_{\substack{P\subseteq\left[  n\right]  ;\\\left\vert P\right\vert
=\left\vert Q\right\vert }}\left(  -1\right)  ^{\operatorname*{sum}%
P+\operatorname*{sum}Q}\det\left(  \operatorname*{sub}\nolimits_{P}%
^{Q}A\right)  \det\left(  \operatorname*{sub}\nolimits_{\widetilde{P}%
}^{\widetilde{Q}}A\right)  .
\]

\end{theorem}

\begin{example}
Let $n=4$ and $A=\left(
\begin{array}
[c]{cccc}%
a & b & c & d\\
a^{\prime} & b^{\prime} & c^{\prime} & d^{\prime}\\
a^{\prime\prime} & b^{\prime\prime} & c^{\prime\prime} & d^{\prime\prime}\\
a^{\prime\prime\prime} & b^{\prime\prime\prime} & c^{\prime\prime\prime} &
d^{\prime\prime\prime}%
\end{array}
\right)  $ and $P=\left\{  3,4\right\}  $. Then, Theorem
\ref{thm.det.laplace-multi} \textbf{(a)} says that%
\begin{align*}
&  \det A\\
&  =\sum_{\substack{Q\subseteq\left[  n\right]  ;\\\left\vert Q\right\vert
=\left\vert P\right\vert }}\left(  -1\right)  ^{\operatorname*{sum}%
P+\operatorname*{sum}Q}\det\left(  \operatorname*{sub}\nolimits_{P}%
^{Q}A\right)  \det\left(  \operatorname*{sub}\nolimits_{\widetilde{P}%
}^{\widetilde{Q}}A\right) \\
&  =\left(  -1\right)  ^{\operatorname*{sum}\left\{  3,4\right\}
+\operatorname*{sum}\left\{  1,2\right\}  }\det\left(  \operatorname*{sub}%
\nolimits_{\left\{  3,4\right\}  }^{\left\{  1,2\right\}  }A\right)
\det\left(  \operatorname*{sub}\nolimits_{\widetilde{\left\{  3,4\right\}  }%
}^{\widetilde{\left\{  1,2\right\}  }}A\right) \\
&  \ \ \ \ \ \ \ \ \ \ +\left(  -1\right)  ^{\operatorname*{sum}\left\{
3,4\right\}  +\operatorname*{sum}\left\{  1,3\right\}  }\det\left(
\operatorname*{sub}\nolimits_{\left\{  3,4\right\}  }^{\left\{  1,3\right\}
}A\right)  \det\left(  \operatorname*{sub}\nolimits_{\widetilde{\left\{
3,4\right\}  }}^{\widetilde{\left\{  1,3\right\}  }}A\right) \\
&  \ \ \ \ \ \ \ \ \ \ +\left(  -1\right)  ^{\operatorname*{sum}\left\{
3,4\right\}  +\operatorname*{sum}\left\{  1,4\right\}  }\det\left(
\operatorname*{sub}\nolimits_{\left\{  3,4\right\}  }^{\left\{  1,4\right\}
}A\right)  \det\left(  \operatorname*{sub}\nolimits_{\widetilde{\left\{
3,4\right\}  }}^{\widetilde{\left\{  1,4\right\}  }}A\right) \\
&  \ \ \ \ \ \ \ \ \ \ +\left(  -1\right)  ^{\operatorname*{sum}\left\{
3,4\right\}  +\operatorname*{sum}\left\{  2,3\right\}  }\det\left(
\operatorname*{sub}\nolimits_{\left\{  3,4\right\}  }^{\left\{  2,3\right\}
}A\right)  \det\left(  \operatorname*{sub}\nolimits_{\widetilde{\left\{
3,4\right\}  }}^{\widetilde{\left\{  2,3\right\}  }}A\right) \\
&  \ \ \ \ \ \ \ \ \ \ +\left(  -1\right)  ^{\operatorname*{sum}\left\{
3,4\right\}  +\operatorname*{sum}\left\{  2,4\right\}  }\det\left(
\operatorname*{sub}\nolimits_{\left\{  3,4\right\}  }^{\left\{  2,4\right\}
}A\right)  \det\left(  \operatorname*{sub}\nolimits_{\widetilde{\left\{
3,4\right\}  }}^{\widetilde{\left\{  2,4\right\}  }}A\right) \\
&  \ \ \ \ \ \ \ \ \ \ +\left(  -1\right)  ^{\operatorname*{sum}\left\{
3,4\right\}  +\operatorname*{sum}\left\{  3,4\right\}  }\det\left(
\operatorname*{sub}\nolimits_{\left\{  3,4\right\}  }^{\left\{  3,4\right\}
}A\right)  \det\left(  \operatorname*{sub}\nolimits_{\widetilde{\left\{
3,4\right\}  }}^{\widetilde{\left\{  3,4\right\}  }}A\right) \\
&  =\det\left(
\begin{array}
[c]{cc}%
a^{\prime\prime} & b^{\prime\prime}\\
a^{\prime\prime\prime} & b^{\prime\prime\prime}%
\end{array}
\right)  \det\left(
\begin{array}
[c]{cc}%
c & d\\
c^{\prime} & d^{\prime}%
\end{array}
\right)  -\det\left(
\begin{array}
[c]{cc}%
a^{\prime\prime} & c^{\prime\prime}\\
a^{\prime\prime\prime} & c^{\prime\prime\prime}%
\end{array}
\right)  \det\left(
\begin{array}
[c]{cc}%
b & d\\
b^{\prime} & d^{\prime}%
\end{array}
\right) \\
&  \ \ \ \ \ \ \ \ \ \ +\det\left(
\begin{array}
[c]{cc}%
a^{\prime\prime} & d^{\prime\prime}\\
a^{\prime\prime\prime} & d^{\prime\prime\prime}%
\end{array}
\right)  \det\left(
\begin{array}
[c]{cc}%
b & c\\
b^{\prime} & c^{\prime}%
\end{array}
\right)  +\det\left(
\begin{array}
[c]{cc}%
b^{\prime\prime} & c^{\prime\prime}\\
b^{\prime\prime\prime} & c^{\prime\prime\prime}%
\end{array}
\right)  \det\left(
\begin{array}
[c]{cc}%
a & d\\
a^{\prime} & d^{\prime}%
\end{array}
\right) \\
&  \ \ \ \ \ \ \ \ \ \ -\det\left(
\begin{array}
[c]{cc}%
b^{\prime\prime} & d^{\prime\prime}\\
b^{\prime\prime\prime} & d^{\prime\prime\prime}%
\end{array}
\right)  \det\left(
\begin{array}
[c]{cc}%
a & c\\
a^{\prime} & c^{\prime}%
\end{array}
\right)  +\det\left(
\begin{array}
[c]{cc}%
c^{\prime\prime} & d^{\prime\prime}\\
c^{\prime\prime\prime} & d^{\prime\prime\prime}%
\end{array}
\right)  \det\left(
\begin{array}
[c]{cc}%
a & b\\
a^{\prime} & b^{\prime}%
\end{array}
\right)  .
\end{align*}

\end{example}

\begin{proof}
[Proof of Theorem \ref{thm.det.laplace-multi}.]See \cite[Theorem
6.156]{detnotes}.
\end{proof}

\subsubsection{Desnanot--Jacobi and Dodgson condensation}

We come to more exotic results. The following theorem is one of several
versions of the \emph{Desnanot--Jacobi formula}:

\begin{theorem}
[Desnanot--Jacobi formula, take 1]\label{thm.det.des-jac-1}Let $n\in
\mathbb{N}$ be such that $n\geq2$. Let $A\in K^{n\times n}$ be an $n\times n$-matrix.

Let $A^{\prime}$ be the $\left(  n-2\right)  \times\left(  n-2\right)
$-matrix%
\[
\operatorname*{sub}\nolimits_{\left\{  2,3,\ldots,n-1\right\}  }^{\left\{
2,3,\ldots,n-1\right\}  }A=\left(  A_{i+1,j+1}\right)  _{1\leq i\leq
n-2,\ 1\leq j\leq n-2}.
\]
(This is precisely what remains of the matrix $A$ when we remove the first
row, the last row, the first column and the last column.) Then,%
\begin{align*}
\det A\cdot\det\left(  A^{\prime}\right)   &  =\det\left(  A_{\sim1,\sim
1}\right)  \cdot\det\left(  A_{\sim n,\sim n}\right)  -\det\left(
A_{\sim1,\sim n}\right)  \cdot\det\left(  A_{\sim n,\sim1}\right) \\
&  =\det\left(
\begin{array}
[c]{cc}%
\det\left(  A_{\sim1,\sim1}\right)  & \det\left(  A_{\sim1,\sim n}\right) \\
\det\left(  A_{\sim n,\sim1}\right)  & \det\left(  A_{\sim n,\sim n}\right)
\end{array}
\right)  .
\end{align*}

\end{theorem}

\begin{example}
For $n=3$, this is saying that%
\begin{align*}
&  \det\left(
\begin{array}
[c]{ccc}%
a & b & c\\
a^{\prime} & b^{\prime} & c^{\prime}\\
a^{\prime\prime} & b^{\prime\prime} & c^{\prime\prime}%
\end{array}
\right)  \cdot\det\left(
\begin{array}
[c]{c}%
b^{\prime}%
\end{array}
\right) \\
&  =\det\left(
\begin{array}
[c]{cc}%
b^{\prime} & c^{\prime}\\
b^{\prime\prime} & c^{\prime\prime}%
\end{array}
\right)  \cdot\det\left(
\begin{array}
[c]{cc}%
a & b\\
a^{\prime} & b^{\prime}%
\end{array}
\right)  -\det\left(
\begin{array}
[c]{cc}%
a^{\prime} & b^{\prime}\\
a^{\prime\prime} & b^{\prime\prime}%
\end{array}
\right)  \cdot\det\left(
\begin{array}
[c]{cc}%
b & c\\
b^{\prime} & c^{\prime}%
\end{array}
\right)  .
\end{align*}

\end{example}

\begin{proof}
[Proof of Theorem \ref{thm.det.des-jac-1}.]See \cite[Proposition
6.122]{detnotes} or \cite[\S 3.5, proof of Theorem 3.12]{Bresso99} or
\cite{zeilberger-twotime}.
\end{proof}

Theorem \ref{thm.det.des-jac-1} provides a recursive way of computing
determinants: Indeed, in the setting of Theorem \ref{thm.det.des-jac-1}, if
$\det\left(  A^{\prime}\right)  $ is invertible (which, when $K$ is a field,
simply means that $\det\left(  A^{\prime}\right)  \neq0$), then Theorem
\ref{thm.det.des-jac-1} yields%
\begin{equation}
\det A=\dfrac{\det\left(  A_{\sim1,\sim1}\right)  \cdot\det\left(  A_{\sim
n,\sim n}\right)  -\det\left(  A_{\sim1,\sim n}\right)  \cdot\det\left(
A_{\sim n,\sim1}\right)  }{\det\left(  A^{\prime}\right)  }.
\label{eq.thm.det.des-jac-1.frac}%
\end{equation}
The five matrices appearing on the right hand side of this are smaller than
$A$, so their determinants are often easier to compute than $\det A$. In
particular, if you are proving something by strong induction on $n$, you will
occasionally be able to use the induction hypothesis to compute these
determinants. This method of recursively simplifying determinants is often
known as \emph{Dodgson condensation}, as it was popularized (perhaps even
discovered) by Charles Lutwidge Dodgson (aka Lewis Carroll) in \cite[Appendix
II]{Dodgso67}. We outline a sample application:

\begin{theorem}
[Cauchy determinant]\label{thm.det.cauchy}Let $n\in\mathbb{N}$. Let
$x_{1},x_{2},\ldots,x_{n}$ be $n$ elements of $K$. Let $y_{1},y_{2}%
,\ldots,y_{n}$ be $n$ elements of $K$. Assume that $x_{i}+y_{j}$ is invertible
in $K$ for each $\left(  i,j\right)  \in\left[  n\right]  ^{2}$. Then,%
\[
\det\left(  \left(  \dfrac{1}{x_{i}+y_{j}}\right)  _{1\leq i\leq n,\ 1\leq
j\leq n}\right)  =\dfrac{\prod\limits_{1\leq i<j\leq n}\left(  \left(
x_{i}-x_{j}\right)  \left(  y_{i}-y_{j}\right)  \right)  }{\prod
\limits_{\left(  i,j\right)  \in\left[  n\right]  ^{2}}\left(  x_{i}%
+y_{j}\right)  }.
\]

\end{theorem}

Once again, there are many ways to prove this (see, e.g., \cite[Exercise
6.18]{detnotes}, \cite[\S 1.3]{Prasolov}, \cite[Theorem 2]{Gri-19.9}, or
\url{https://proofwiki.org/wiki/Value_of_Cauchy_Determinant} ). But using the
Desnanot--Jacobi identity, there is a rather straightforward proof of Theorem
\ref{thm.det.cauchy}. Indeed, if $A$ is a \emph{Cauchy matrix} (i.e., a matrix
of the form $\left(  \dfrac{1}{x_{i}+y_{j}}\right)  _{1\leq i\leq n,\ 1\leq
j\leq n}$), then so is each submatrix of $A$. Thus, if we proceed by strong
induction on $n$, we can use the induction hypothesis to compute all five
determinants on the right hand side of (\ref{eq.thm.det.des-jac-1.frac}). The
only difficulty is making sure that $\det\left(  A^{\prime}\right)  $ is
invertible. To achieve this, we again have to WLOG assume that our $x$'s and
$y$'s are indeterminates in a polynomial ring, and we have to rewrite the
claim of Theorem \ref{thm.det.cauchy} in the form%
\[
\det\left(  \left(  \prod_{k\neq j}\left(  x_{i}+y_{k}\right)  \right)
_{1\leq i\leq n,\ 1\leq j\leq n}\right)  =\prod\limits_{1\leq i<j\leq
n}\left(  \left(  x_{i}-x_{j}\right)  \left(  y_{i}-y_{j}\right)  \right)
\]
in order for both sides to actually be polynomials in $\mathbb{Z}\left[
x_{1},x_{2},\ldots,x_{n},y_{1},y_{2},\ldots,y_{n}\right]  $. The details are
left to the reader.

Theorem \ref{thm.det.des-jac-1} can be generalized significantly. Here is the
simplest generalization, in which the special role played by the first and
last rows and the first and last columns is instead given to any two rows and
any two columns:

\begin{theorem}
\label{thm.det.des-jac-2}Let $n\in\mathbb{N}$ be such that $n\geq2$. Let $p$,
$q$, $u$ and $v$ be four elements of $\left[  n\right]  $ such that $p<q$ and
$u<v$. Let $A$ be an $n\times n$-matrix. Then,%
\begin{align*}
&  \det A\cdot\det\left(  \operatorname*{sub}\nolimits_{\left[  n\right]
\setminus\left\{  p,q\right\}  }^{\left[  n\right]  \setminus\left\{
u,v\right\}  }A\right) \\
&  =\det\left(  A_{\sim p,\sim u}\right)  \cdot\det\left(  A_{\sim q,\sim
v}\right)  -\det\left(  A_{\sim p,\sim v}\right)  \cdot\det\left(  A_{\sim
q,\sim u}\right)  .
\end{align*}

\end{theorem}

\begin{proof}
See \cite[Theorem 6.126]{detnotes}.
\end{proof}

Even more generally, \emph{Jacobi's complementary minor theorem for adjugates}
(appearing, e.g., in \cite[Theorem 5.22]{trach}, or in equivalent forms in
\cite[Theorem 2.5.2]{Prasolov} and \cite[(13)]{BruSch83}) says the following:

\begin{theorem}
[Jacobi's complementary minor theorem for adjugates]%
\label{thm.det.jacobi-complement}Let $n\in\mathbb{N}$. For any subset $I$ of
$\left[  n\right]  $, we let $\widetilde{I}$ denote the complement $\left[
n\right]  \setminus I$ of $I$. Set $\operatorname*{sum}S=\sum_{s\in S}s$ for
any finite set $S$ of integers. (For example, $\operatorname*{sum}\left\{
2,4,5\right\}  =2+4+5=11$.)

Let $A\in K^{n\times n}$ be an $n\times n$-matrix. Let $P$ and $Q$ be two
subsets of $\left[  n\right]  $ such that $\left\vert P\right\vert =\left\vert
Q\right\vert \geq1$. Then,%
\[
\det\left(  \operatorname*{sub}\nolimits_{P}^{Q}\left(  \operatorname*{adj}%
A\right)  \right)  =\left(  -1\right)  ^{\operatorname*{sum}%
P+\operatorname*{sum}Q}\left(  \det A\right)  ^{\left\vert Q\right\vert
-1}\det\left(  \operatorname*{sub}\nolimits_{\widetilde{Q}}^{\widetilde{P}%
}A\right)  .
\]

\end{theorem}

Theorem \ref{thm.det.des-jac-2} is the particular case of Theorem
\ref{thm.det.jacobi-complement} for $P=\left\{  u,v\right\}  $ and $Q=\left\{
p,q\right\}  $.\ \ \ \ \footnote{Indeed, if we set $P=\left\{  u,v\right\}  $
and $Q=\left\{  p,q\right\}  $ in the situation of Theorem
\ref{thm.det.des-jac-2}, then%
\begin{align*}
\operatorname*{sub}\nolimits_{P}^{Q}\left(  \operatorname*{adj}A\right)   &
=\left(
\begin{array}
[c]{cc}%
\left(  -1\right)  ^{u+p}\det\left(  A_{\sim p,\sim u}\right)  & \left(
-1\right)  ^{u+q}\det\left(  A_{\sim q,\sim u}\right) \\
\left(  -1\right)  ^{v+p}\det\left(  A_{\sim p,\sim v}\right)  & \left(
-1\right)  ^{v+q}\det\left(  A_{\sim q,\sim v}\right)
\end{array}
\right)  \ \ \ \ \ \ \ \ \ \ \text{and}\\
\operatorname*{sub}\nolimits_{\widetilde{Q}}^{\widetilde{P}}A  &
=\operatorname*{sub}\nolimits_{\left[  n\right]  \setminus\left\{
p,q\right\}  }^{\left[  n\right]  \setminus\left\{  u,v\right\}  }A;
\end{align*}
thus, Theorem \ref{thm.det.jacobi-complement} is easily seen to boil down to
Theorem \ref{thm.det.des-jac-2} in this case (the powers of $-1$ all cancel).}
Theorem \ref{thm.det.des-jac-1} is, of course, the particular case of Theorem
\ref{thm.det.des-jac-2} for $p=1$ and $q=n$ and $u=1$ and $v=n$.

\begin{noncompile}
TODO: Here I should eventually add a proof of Theorem
\ref{thm.det.jacobi-complement} a la Prasolov. However, it needs the
determinant of a block-triangular matrix.
\end{noncompile}

\subsection{\label{sec.det.comb.lgv}The Lindstr\"{o}m--Gessel--Viennot lemma}

We have so far mostly been discussing algebraic properties of determinants, if
often from a combinatorial point of view. We shall now see a situation where
determinants naturally appear in combinatorics.

Namely, we will see an application of determinants to lattice path enumeration
-- i.e., to the counting of paths on a certain infinite digraph called the
\emph{integer lattice}. A survey of this subject can be found in
\cite{Kratte17}; we will restrict ourselves to one of the most accessible
highlights: the \emph{Lindstr\"{o}m--Gessel--Viennot lemma} (for short:
\emph{LGV lemma}). This \textquotedblleft lemma\textquotedblright\ is by now
so classical and popular that it is often deservedly called the
\emph{Lindstr\"{o}m--Gessel--Viennot theorem}. Alternate treatments of this
theorem can be found in \cite[\S 2.5]{Sagan19}, in \cite[\S 2.7]{Stanley-EC1}
and in \cite[\S 2]{GesVie89}. Some applications predating the general
statement of the theorem can be found in \cite{GesVie85}, and a surprising
recent generalization in \cite[Theorem 2.5]{Talask12}.

\subsubsection{Definitions}

We have already seen lattice paths (and even counted them in Subsection
\ref{subsec.qbin.motiv}). We shall now introduce them formally and study them
in greater depth.

\Needspace{20pc}

\begin{convention}
\label{conv.lgv.digraph}\textbf{(a)} Recall that \textquotedblleft%
\emph{digraph}\textquotedblright\ means \textquotedblleft directed
graph\textquotedblright, i.e., a graph whose edges are directed (and are
called \emph{arcs}). Against a widespread convention, we will allow our
digraphs to be infinite (i.e., to have infinitely many vertices and arcs).
\medskip

\textbf{(b)} A digraph $D$ will be called \emph{path-finite} if it has the
property that for any two vertices $u$ and $v$, there are only finitely many
paths from $u$ to $v$. (Thus, in particular, such paths can be counted.)
\medskip

\textbf{(c)} A digraph $D$ will be called \emph{acyclic} if it has no directed
cycles. For example, the digraph $%
\raisebox{-3.3pc}{
\begin{tikzpicture}%
[->,shorten >=1pt, thick,main node/.style={circle,fill=blue!20,draw}]
\node[main node] (1) at (0 : 1) {1};
\node[main node] (2) at (90 : 1) {2};
\node[main node] (3) at (180 : 1) {3};
\node[main node] (4) at (270 : 1) {4};
\path[-{Stealth[length=4mm]}]
(1) edge (2)
(2) edge (3)
(4) edge (3)
(1) edge (4);
\end{tikzpicture}
}%
$ is acyclic, whereas the digraph $%
\raisebox{-3.3pc}{
\begin{tikzpicture}%
[->,shorten >=1pt, thick,main node/.style={circle,fill=blue!20,draw}]
\node[main node] (1) at (0 : 1) {1};
\node[main node] (2) at (90 : 1) {2};
\node[main node] (3) at (180 : 1) {3};
\node[main node] (4) at (270 : 1) {4};
\path[-{Stealth[length=4mm]}]
(1) edge (2)
(2) edge (3)
(3) edge (4)
(4) edge (1);
\end{tikzpicture}
}%
$ is not. \medskip

\textbf{(d)} A \emph{simple digraph} $D$ means a digraph whose arcs are merely
pairs of distinct vertices (i.e., each arc is a pair $\left(  u,v\right)  $ of
two vertices $u$ and $v$ with $u\neq v$).
\end{convention}

We note that a path may contain $0$ arcs (in which case its starting and
ending point are identical).

\begin{definition}
\label{def.lgv.lattice}We consider the infinite simple digraph with vertex set
$\mathbb{Z}^{2}$ (so the vertices are pairs of integers) and arcs%
\begin{equation}
\left(  i,j\right)  \rightarrow\left(  i+1,j\right)
\ \ \ \ \ \ \ \ \ \ \text{for all }\left(  i,j\right)  \in\mathbb{Z}^{2}
\label{eq.def.lgv.lattice.east}%
\end{equation}
and%
\begin{equation}
\left(  i,j\right)  \rightarrow\left(  i,j+1\right)
\ \ \ \ \ \ \ \ \ \ \text{for all }\left(  i,j\right)  \in\mathbb{Z}^{2}.
\label{eq.def.lgv.lattice.north}%
\end{equation}
The arcs of the form (\ref{eq.def.lgv.lattice.east}) are called
\textquotedblleft\emph{east-steps}\textquotedblright\ or \textquotedblleft%
\emph{right-steps}\textquotedblright; the arcs of the form
(\ref{eq.def.lgv.lattice.north}) are called \textquotedblleft%
\emph{north-steps}\textquotedblright\ or \textquotedblleft\emph{up-steps}%
\textquotedblright.

The vertices of this digraph will be called \emph{lattice points} or
\emph{grid points} or simply \emph{points}. They will be represented as points
in the Cartesian plane (in the usual way: the vertex $\left(  i,j\right)
\in\mathbb{Z}^{2}$ is represented as the point with x-coordinate $i$ and
y-coordinate $j$).

The entire digraph will be denoted by $\mathbb{Z}^{2}$ and called the
\emph{integer lattice} or \emph{integer grid} (or, to be short, just
\emph{lattice} or \emph{grid}). Here is a picture of a small part of this
digraph $\mathbb{Z}^{2}$:%
\[%
\begin{tikzpicture}
\foreach\x/\xtext in {-1, 0, 1, 2, 3, 4, 5}
{
\foreach\y/\ytext in {-1, 0, 1, 2, 3, 4, 5}
{
\ifnum\x>-1
\draw(\x,\y+0.12) edge[->, thick, >=stealth, darkred] (\x,\y+0.88);
\fi\ifnum\y>-1
\draw(\x+0.12,\y) edge[->, thick, >=stealth, dbluecolor] (\x+0.88,\y);
\fi\ifnum\x>-1
\ifnum\y>-1
\draw(\x,\y) circle[radius=0.1, black];
\fi\fi}
};
\foreach\x/\xtext in {0, 1, 2, 3, 4, 5}
\draw(\x cm,1pt) node[anchor=north west] {$\xtext$};
\foreach\y/\ytext in {0, 1, 2, 3, 4, 5}
\draw(1pt,\y cm) node[anchor=south east] {$\ytext$};
\end{tikzpicture}
\]
(with east-steps colored blue and north-steps colored dark-red). Of course,
the digraph continues indefinitely in all directions. In the following, we
will not draw the vertices as circles, nor will we draw the arcs as arrows; we
will simply draw the grid lines in order to avoid crowding our pictures.

However, $\mathbb{Z}^{2}$ is also an abelian group under addition. Thus,
points can be added and subtracted entrywise; e.g., for any $\left(
a,b\right)  \in\mathbb{Z}^{2}$ and $\left(  c,d\right)  \in\mathbb{Z}^{2}$, we
have%
\begin{align*}
\left(  a,b\right)  +\left(  c,d\right)   &  =\left(  a+c,b+d\right)
\ \ \ \ \ \ \ \ \ \ \text{and }\\
\left(  a,b\right)  -\left(  c,d\right)   &  =\left(  a-c,b-d\right)  .
\end{align*}
Thus, the digraph $\mathbb{Z}^{2}$ has an arc from a vertex $u$ to a vertex
$v$ if and only if $v-u\in\left\{  \left(  0,1\right)  ,\ \left(  1,0\right)
\right\}  $.

The digraph $\mathbb{Z}^{2}$ is acyclic (i.e., it has no directed cycles).
Thus, its paths are the same as its walks. We call these paths the
\emph{lattice paths} (or just \emph{paths}). Thus, a lattice path is a finite
tuple $\left(  v_{0},v_{1},\ldots,v_{n}\right)  $ of points $v_{i}%
\in\mathbb{Z}^{2}$ with the property that%
\begin{equation}
v_{i}-v_{i-1}\in\left\{  \left(  0,1\right)  ,\ \left(  1,0\right)  \right\}
\ \ \ \ \ \ \ \ \ \ \text{for each }i\in\left[  n\right]  .
\label{eq.def.lgv.lattice.path.arcs}%
\end{equation}

The \emph{step sequence} of a path $\left(  v_{0},v_{1},\ldots,v_{n}\right)  $
is defined to be the $n$-tuple $\left(  v_{1}-v_{0},\ v_{2}-v_{1}%
,\ \ldots,\ v_{n}-v_{n-1}\right)  $. Each entry of this $n$-tuple is either
$\left(  0,1\right)  $ or $\left(  1,0\right)  $ (because of
(\ref{eq.def.lgv.lattice.path.arcs})). We will often write $U$ and $R$ for the
pairs $\left(  0,1\right)  $ and $\left(  1,0\right)  $, respectively (as they
correspond to an \textbf{u}p-step and a \textbf{r}ight-step). Informally
speaking, the step sequence of a path records the directions (i.e., east or
north) of all steps of the path.
\end{definition}

\begin{example}
\label{exa.lgv.lattice.path.53}Here is a path from $\left(  0,0\right)  $ to
$\left(  5,3\right)  $:%
\[%
\begin{tikzpicture}
\draw[densely dotted] (0,0) grid (5.2,3.2);
\draw[->] (0,0) -- (0,3.2);
\draw[->] (0,0) -- (5.2,0);
\foreach\x/\xtext in {0, 1, 2, 3, 4, 5}
\draw(\x cm,1pt) -- (\x cm,-1pt) node[anchor=north] {$\xtext$};
\foreach\y/\ytext in {0, 1, 2, 3}
\draw(1pt,\y cm) -- (-1pt,\y cm) node[anchor=east] {$\ytext$};
\begin{scope}[thick,>=stealth,darkred]
\draw(0,0) edge[->] (0,1);
\draw(0,1) edge[->] (1,1);
\draw(1,1) edge[->] (2,1);
\draw(2,1) edge[->] (3,1);
\draw(3,1) edge[->] (3,2);
\draw(3,2) edge[->] (4,2);
\draw(4,2) edge[->] (4,3);
\draw(4,3) edge[->] (5,3);
\end{scope}
\end{tikzpicture}%
\ \ .
\]
Formally speaking, this path is the $9$-tuple%
\[
\left(  \left(  0,0\right)  ,\left(  0,1\right)  ,\left(  1,1\right)  ,\left(
2,1\right)  ,\left(  3,1\right)  ,\left(  3,2\right)  ,\left(  4,2\right)
,\left(  4,3\right)  ,\left(  5,3\right)  \right)  .
\]
Its step sequence\ (i.e., the sequence of the directions of its steps) is
$URRRURUR$ (meaning that its first step is an up-step, its second step is a
right-step, its third step is a right-step, and so on).
\end{example}

Clearly, any path is uniquely determined by its starting point and its step sequence.

Note that we are considering one of the simplest possible notions of a lattice
path here. In more advanced texts, the word \textquotedblleft lattice
path\textquotedblright\ is often used for paths in digraphs more complicated
than $\mathbb{Z}^{2}$ (for instance, a digraph with the same vertex set
$\mathbb{Z}^{2}$ but allowing steps in all four directions). However, the
digraph we are considering is perhaps the most useful for algebraic combinatorics.

\subsubsection{Counting paths from $\left(  a,b\right)  $ to $\left(
c,d\right)  $}

In Subsection \ref{subsec.qbin.motiv}, we have counted the lattice paths from
$\left(  0,0\right)  $ to $\left(  6,4\right)  $ that begin with an east-step
and end with a north-step. These are, of course, in bijection with the paths
from $\left(  1,0\right)  $ to $\left(  6,3\right)  $ (since the first and
last step are predetermined and thus can be ignored). Let us now generalize
this by counting paths between any two lattice points:

\begin{proposition}
\label{prop.lgv.1-paths.ct}Let $\left(  a,b\right)  \in\mathbb{Z}^{2}$ and
$\left(  c,d\right)  \in\mathbb{Z}^{2}$ be two points. Then,%
\[
\left(  \text{\# of paths from }\left(  a,b\right)  \text{ to }\left(
c,d\right)  \right)  =%
\begin{cases}
\dbinom{c+d-a-b}{c-a}, & \text{if }c+d\geq a+b;\\
0, & \text{if }c+d<a+b.
\end{cases}
\]

\end{proposition}

\begin{proof}
[Proof of Proposition \ref{prop.lgv.1-paths.ct}.]This is just a formalization
(and generalization) of the reasoning we used in Subsection
\ref{subsec.qbin.motiv}.

We shall first show the following two observations:

\begin{statement}
\textit{Observation 1:} Each path from $\left(  a,b\right)  $ to $\left(
c,d\right)  $ has exactly $c+d-a-b$ steps\footnote{A \textquotedblleft
step\textquotedblright\ of a path means an arc of this path.}.
\end{statement}

\begin{statement}
\textit{Observation 2:} Each path from $\left(  a,b\right)  $ to $\left(
c,d\right)  $ has exactly $c-a$ east-steps.
\end{statement}

\begin{fineprint}
[\textit{Proof of Observation 1:} We define the \emph{coordinate sum} of a
point $\left(  x,y\right)  \in\mathbb{Z}^{2}$ to be $x+y$. We shall denote
this coordinate sum by $\operatorname*{cs}\left(  x,y\right)  $. We observe
that the coordinate sum of a point increases by exactly $1$ along each arc of
$\mathbb{Z}^{2}$: That is, if $u\rightarrow v$ is an arc of $\mathbb{Z}^{2}$,
then%
\begin{equation}
\operatorname*{cs}\left(  v\right)  -\operatorname*{cs}\left(  u\right)  =1.
\label{pf.prop.lgv.1-paths.ct.o1.pf.1}%
\end{equation}
(This is because we can write $u$ in the form $u=\left(  i,j\right)  $ and
then must have either $v=\left(  i+1,j\right)  $ or $v=\left(  i,j+1\right)
$; but this entails $\operatorname*{cs}\left(  v\right)  =\underbrace{i+j}%
_{=\operatorname*{cs}\left(  u\right)  }+\,1=\operatorname*{cs}\left(
u\right)  +1$ in either case.)

Let $\left(  v_{0},v_{1},\ldots,v_{n}\right)  $ be a path from $\left(
a,b\right)  $ to $\left(  c,d\right)  $. Thus, $v_{0}=\left(  a,b\right)  $
and $v_{n}=\left(  c,d\right)  $. Moreover, for each $i\in\left[  n\right]  $,
we know that $v_{i-1}\rightarrow v_{i}$ is an arc of $\mathbb{Z}^{2}$, and
thus we have%
\[
\operatorname*{cs}\left(  v_{i}\right)  -\operatorname*{cs}\left(
v_{i-1}\right)  =1
\]
(by (\ref{pf.prop.lgv.1-paths.ct.o1.pf.1}), applied to $u=v_{i-1}$ and
$v=v_{i}$). Summing these equalities over all $i\in\left[  n\right]  $, we
obtain%
\[
\sum_{i=1}^{n}\left(  \operatorname*{cs}\left(  v_{i}\right)
-\operatorname*{cs}\left(  v_{i-1}\right)  \right)  =\sum_{i=1}^{n}1=n.
\]
Hence,%
\begin{align*}
n  &  =\sum_{i=1}^{n}\left(  \operatorname*{cs}\left(  v_{i}\right)
-\operatorname*{cs}\left(  v_{i-1}\right)  \right)  =\operatorname*{cs}%
\underbrace{\left(  v_{n}\right)  }_{=\left(  c,d\right)  }-\operatorname*{cs}%
\underbrace{\left(  v_{0}\right)  }_{=\left(  a,b\right)  }%
\ \ \ \ \ \ \ \ \ \ \left(  \text{by the telescope principle}\right) \\
&  =\underbrace{\operatorname*{cs}\left(  c,d\right)  }_{=c+d}%
-\underbrace{\operatorname*{cs}\left(  a,b\right)  }_{=a+b}=c+d-\left(
a+b\right)  =c+d-a-b.
\end{align*}
In other words, the path $\left(  v_{0},v_{1},\ldots,v_{n}\right)  $ has
exactly $c+d-a-b$ steps (since this path clearly has $n$ steps).

Forget that we fixed $\left(  v_{0},v_{1},\ldots,v_{n}\right)  $. We thus have
shown that each path $\left(  v_{0},v_{1},\ldots,v_{n}\right)  $ from $\left(
a,b\right)  $ to $\left(  c,d\right)  $ has exactly $c+d-a-b$ steps. This
proves Observation 1.] \medskip
\end{fineprint}

\begin{fineprint}
[\textit{Proof of Observation 2:} For any point $v\in\mathbb{Z}^{2}$, we
define $\operatorname*{x}\left(  v\right)  $ to be the x-coordinate of $v$.
(Thus, $\operatorname*{x}\left(  x,y\right)  =x$ for each $\left(  x,y\right)
\in\mathbb{Z}^{2}$.)

Obviously, the x-coordinate of a point increases by exactly $1$ along each
east-step and stays unchanged along each north-step: That is, if $u\rightarrow
v$ is an arc of $\mathbb{Z}^{2}$, then%
\begin{equation}
\operatorname*{x}\left(  v\right)  -\operatorname*{x}\left(  u\right)  =%
\begin{cases}
1, & \text{if }u\rightarrow v\text{ is an east-step;}\\
0, & \text{if }u\rightarrow v\text{ is a north-step.}%
\end{cases}
\label{pf.prop.lgv.1-paths.ct.o2.pf.1}%
\end{equation}

Let $\left(  v_{0},v_{1},\ldots,v_{n}\right)  $ be a path from $\left(
a,b\right)  $ to $\left(  c,d\right)  $. Thus, $v_{0}=\left(  a,b\right)  $
and $v_{n}=\left(  c,d\right)  $. Moreover, for each $i\in\left[  n\right]  $,
we know that $v_{i-1}\rightarrow v_{i}$ is an arc of $\mathbb{Z}^{2}$, and
thus we have%
\[
\operatorname*{x}\left(  v_{i}\right)  -\operatorname*{x}\left(
v_{i-1}\right)  =%
\begin{cases}
1, & \text{if }v_{i-1}\rightarrow v_{i}\text{ is an east-step;}\\
0, & \text{if }v_{i-1}\rightarrow v_{i}\text{ is a north-step}%
\end{cases}
\]
(by (\ref{pf.prop.lgv.1-paths.ct.o2.pf.1}), applied to $u=v_{i-1}$ and
$v=v_{i}$). Summing these equalities over all $i\in\left[  n\right]  $, we
obtain%
\begin{align*}
\sum_{i=1}^{n}\left(  \operatorname*{x}\left(  v_{i}\right)
-\operatorname*{x}\left(  v_{i-1}\right)  \right)   &  =\sum_{i=1}^{n}%
\begin{cases}
1, & \text{if }v_{i-1}\rightarrow v_{i}\text{ is an east-step;}\\
0, & \text{if }v_{i-1}\rightarrow v_{i}\text{ is a north-step}%
\end{cases}
\\
&  =\left(  \text{\# of }i\in\left[  n\right]  \text{ such that }%
v_{i-1}\rightarrow v_{i}\text{ is an east-step}\right) \\
&  =\left(  \text{\# of east-steps in the path }\left(  v_{0},v_{1}%
,\ldots,v_{n}\right)  \right)  .
\end{align*}
Hence,%
\begin{align*}
&  \left(  \text{\# of east-steps in the path }\left(  v_{0},v_{1}%
,\ldots,v_{n}\right)  \right) \\
&  =\sum_{i=1}^{n}\left(  \operatorname*{x}\left(  v_{i}\right)
-\operatorname*{x}\left(  v_{i-1}\right)  \right)  =\operatorname*{x}%
\underbrace{\left(  v_{n}\right)  }_{=\left(  c,d\right)  }-\operatorname*{x}%
\underbrace{\left(  v_{0}\right)  }_{=\left(  a,b\right)  }%
\ \ \ \ \ \ \ \ \ \ \left(  \text{by the telescope principle}\right) \\
&  =\underbrace{\operatorname*{x}\left(  c,d\right)  }_{=c}%
-\underbrace{\operatorname*{x}\left(  a,b\right)  }_{=a}=c-a.
\end{align*}
In other words, the path $\left(  v_{0},v_{1},\ldots,v_{n}\right)  $ has
exactly $c-a$ east-steps.

Forget that we fixed $\left(  v_{0},v_{1},\ldots,v_{n}\right)  $. We thus have
shown that each path $\left(  v_{0},v_{1},\ldots,v_{n}\right)  $ from $\left(
a,b\right)  $ to $\left(  c,d\right)  $ has exactly $c-a$ east-steps. This
proves Observation 2.] \medskip
\end{fineprint}

Observation 1 immediately shows that no path from $\left(  a,b\right)  $ to
$\left(  c,d\right)  $ exists when $c+d-a-b<0$. In other words, no path from
$\left(  a,b\right)  $ to $\left(  c,d\right)  $ exists when $c+d<a+b$. In
other words, $\left(  \text{\# of paths from }\left(  a,b\right)  \text{ to
}\left(  c,d\right)  \right)  =0$ when $c+d<a+b$. This proves Proposition
\ref{prop.lgv.1-paths.ct} in the case when $c+d<a+b$. Hence, for the rest of
this proof of Proposition \ref{prop.lgv.1-paths.ct}, we WLOG assume that
$c+d\geq a+b$. Thus, $c+d-a-b\geq0$, so that $c+d-a-b\in\mathbb{N}$.

Observations 1 and 2 have a sort of (common) converse:

\begin{statement}
\textit{Observation 3:} Let $p$ be a path that starts at the point $\left(
a,b\right)  $ and has exactly $c+d-a-b$ steps. Assume that exactly $c-a$ of
these steps are east-steps. Then, this path $p$ ends at $\left(  c,d\right)  $.
\end{statement}

\begin{fineprint}
[\textit{Proof of Observation 3:} Let $\left(  c^{\prime},d^{\prime}\right)  $
be the point at which this path $p$ ends. Then, we can apply Observation 1 to
$\left(  c^{\prime},d^{\prime}\right)  $ instead of $\left(  c,d\right)  $,
and hence conclude that this path has exactly $c^{\prime}+d^{\prime}-a-b$
steps. Since we already know that this path has exactly $c+d-a-b$ steps, we
therefore conclude that $c^{\prime}+d^{\prime}-a-b=c+d-a-b$. In other words,
$c^{\prime}+d^{\prime}=c+d$. Similarly, using Observation 2, we can find that
$c^{\prime}=c$. Subtracting this equality from $c^{\prime}+d^{\prime}=c+d$, we
obtain $d^{\prime}=d$. Combining $c^{\prime}=c$ with $d^{\prime}=d$, we find
$\left(  c^{\prime},d^{\prime}\right)  =\left(  c,d\right)  $. Thus, our path
$p$ ends at $\left(  c,d\right)  $ (since we know that it ends at $\left(
c^{\prime},d^{\prime}\right)  $). This proves Observation 3.] \medskip
\end{fineprint}

Now, combining Observations 1, 2 and 3, we see that the paths from $\left(
a,b\right)  $ to $\left(  c,d\right)  $ are precisely the paths that start at
$\left(  a,b\right)  $ and have exactly $c+d-a-b$ steps and exactly $c-a$
east-steps (among these $c+d-a-b$ steps). Such a path is therefore uniquely
determined if we know \textbf{which} $c-a$ of its $c+d-a-b$ steps are
east-steps. Thus, specifying such a path is equivalent to specifying a
$\left(  c-a\right)  $-element subset of the $\left(  c+d-a-b\right)
$-element set\footnote{Here we are tacitly using $c+d-a-b\geq0$.} $\left[
c+d-a-b\right]  $. The bijection principle thus yields\footnote{To make this
more formal: We are saying that the map
\begin{align*}
\left\{  \text{paths from }\left(  a,b\right)  \text{ to }\left(  c,d\right)
\right\}   &  \rightarrow\left\{  \left(  c-a\right)  \text{-element subsets
of }\left[  c+d-a-b\right]  \right\}  ,\\
\left(  v_{0},v_{1},\ldots,v_{c+d-a-b}\right)   &  \mapsto\left\{  i\in\left[
c+d-a-b\right]  \ \mid\ \text{the arc }v_{i-1}\rightarrow v_{i}\text{ is an
east-step}\right\}
\end{align*}
is a bijection, and we are applying the bijection principle to this
bijection.}%
\begin{align*}
&  \left(  \text{\# of paths from }\left(  a,b\right)  \text{ to }\left(
c,d\right)  \right) \\
&  =\left(  \text{\# of }\left(  c-a\right)  \text{-element subsets of
}\left[  c+d-a-b\right]  \right) \\
&  =\dbinom{c+d-a-b}{c-a}\ \ \ \ \ \ \ \ \ \ \left(  \text{since }%
c+d-a-b\in\mathbb{N}\right)  .
\end{align*}
This proves Proposition \ref{prop.lgv.1-paths.ct} (since we have assumed that
$c+d\geq a+b$).
\end{proof}

\subsubsection{Path tuples, nipats and ipats}

Now, let us try to count something more interesting: tuples of
non-intersecting paths.

\begin{definition}
\label{def.lgv.path-tups}Let $k\in\mathbb{N}$. \medskip

\textbf{(a)} A $k$\emph{-vertex} means a $k$-tuple of lattice points (i.e., a
$k$-tuple of vertices). For example, $\left(  \left(  1,2\right)  ,\ \left(
4,5\right)  ,\ \left(  7,4\right)  \right)  $ is a $3$-vertex. \medskip

\textbf{(b)} If $\mathbf{A}=\left(  A_{1},A_{2},\ldots,A_{k}\right)  $ is a
$k$-vertex, and if $\sigma\in S_{k}$ is a permutation, then $\sigma\left(
\mathbf{A}\right)  $ shall denote the $k$-vertex $\left(  A_{\sigma\left(
1\right)  },A_{\sigma\left(  2\right)  },\ldots,A_{\sigma\left(  k\right)
}\right)  $. For instance, for the simple transposition $s_{1}\in S_{3}$, we
have $s_{1}\left(  A,B,C\right)  =\left(  B,A,C\right)  $ for any $3$-vertex
$\left(  A,B,C\right)  $. \medskip

\textbf{(c)} If $\mathbf{A}=\left(  A_{1},A_{2},\ldots,A_{k}\right)  $ and
$\mathbf{B}=\left(  B_{1},B_{2},\ldots,B_{k}\right)  $ are two $k$-vertices,
then a \emph{path tuple} from $\mathbf{A}$ to $\mathbf{B}$ means a $k$-tuple
$\left(  p_{1},p_{2},\ldots,p_{k}\right)  $, where each $p_{i}$ is a path from
$A_{i}$ to $B_{i}$. \medskip

\textbf{(d)} A path tuple $\left(  p_{1},p_{2},\ldots,p_{k}\right)  $ is said
to be \emph{non-intersecting} if no two of the paths $p_{1},p_{2},\ldots
,p_{k}$ have any vertex in common. (Visually speaking, this not only forbids
them from crossing each other, but also forbids them from touching or bouncing
off each other, or starting or ending at the same point.)

We shall abbreviate \textquotedblleft non-intersecting path
tuple\textquotedblright\ as \textquotedblleft\emph{nipat}\textquotedblright.
(Historically, the more common abbreviation is \textquotedblleft%
\emph{NILP}\textquotedblright, for \textquotedblleft non-intersecting lattice
paths\textquotedblright, but I prefer \textquotedblleft
nipat\textquotedblright\ as it stresses the tupleness.) \medskip

\textbf{(e)} A path tuple $\left(  p_{1},p_{2},\ldots,p_{k}\right)  $ is said
to be \emph{intersecting} if it is not non-intersecting (i.e., if two of its
paths have a vertex in common).

We shall abbreviate \textquotedblleft intersecting path
tuple\textquotedblright\ as \textquotedblleft\emph{ipat}\textquotedblright.
\end{definition}

\begin{example}
Here are some path tuples for $k=3$: \medskip

\textbf{(a)} The following path tuple is a nipat:%
\[%
\begin{tikzpicture}
\draw[densely dotted] (0,0) grid (7.2,7.2);
\draw[->] (0,0) -- (0,7.2);
\draw[->] (0,0) -- (7.2,0);
\foreach\x/\xtext in {0, 1, 2, 3, 4, 5, 6, 7}
\draw(\x cm,1pt) -- (\x cm,-1pt) node[anchor=north] {$\xtext$};
\foreach\y/\ytext in {0, 1, 2, 3, 4, 5, 6, 7}
\draw(1pt,\y cm) -- (-1pt,\y cm) node[anchor=east] {$\ytext$};
\node
[circle,fill=white,draw=black,text=black,inner sep=1pt] (A3) at (2,3) {$A_3$};
\node
[circle,fill=white,draw=black,text=black,inner sep=1pt] (A2) at (3,1) {$A_2$};
\node
[circle,fill=white,draw=black,text=black,inner sep=1pt] (A1) at (5,1) {$A_1$};
\node
[circle,fill=white,draw=black,text=black,inner sep=1pt] (B3) at (4,6) {$B_3$};
\node
[circle,fill=white,draw=black,text=black,inner sep=1pt] (B2) at (5,6) {$B_2$};
\node
[circle,fill=white,draw=black,text=black,inner sep=1pt] (B1) at (6,4) {$B_1$};
\begin{scope}[thick,>=stealth,darkred]
\draw(A3) edge[->] (2,4);
\draw(2,4) edge[->] (3,4);
\draw(3,4) edge[->] (3,5);
\draw(3,5) edge[->] (4,5);
\draw(4,5) edge[->] (B3);
\draw(3.4,5) node[anchor=south east] {$p_3$};
\end{scope}
\begin{scope}[thick,>=stealth,dbluecolor]
\draw(A2) edge[->] (4,1);
\draw(4,1) edge[->] (4,2);
\draw(4,2) edge[->] (4,3);
\draw(4,3) edge[->] (4,4);
\draw(4,4) edge[->] (5,4);
\draw(5,4) edge[->] (5,5);
\draw(5,5) edge[->] (B2);
\draw(5,4.4) node[anchor=east] {$p_2$};
\end{scope}
\begin{scope}[thick,>=stealth,dgreencolor]
\draw(A1) edge[->] (5,2);
\draw(5,2) edge[->] (6,2);
\draw(6,2) edge[->] (6,3);
\draw(6,3) edge[->] (B1);
\draw(6,2.4) node[anchor=east] {$p_1$};
\end{scope}
\end{tikzpicture}
\ \ .
\]

\textbf{(b)} The following path tuple is an ipat:%
\[%
\begin{tikzpicture}
\draw[densely dotted] (0,0) grid (7.2,7.2);
\draw[->] (0,0) -- (0,7.2);
\draw[->] (0,0) -- (7.2,0);
\foreach\x/\xtext in {0, 1, 2, 3, 4, 5, 6, 7}
\draw(\x cm,1pt) -- (\x cm,-1pt) node[anchor=north] {$\xtext$};
\foreach\y/\ytext in {0, 1, 2, 3, 4, 5, 6, 7}
\draw(1pt,\y cm) -- (-1pt,\y cm) node[anchor=east] {$\ytext$};
\node
[circle,fill=white,draw=black,text=black,inner sep=1pt] (A3) at (2,3) {$A_3$};
\node
[circle,fill=white,draw=black,text=black,inner sep=1pt] (A2) at (3,1) {$A_2$};
\node
[circle,fill=white,draw=black,text=black,inner sep=1pt] (A1) at (5,1) {$A_1$};
\node
[circle,fill=white,draw=black,text=black,inner sep=1pt] (B3) at (4,6) {$B_3$};
\node
[circle,fill=white,draw=black,text=black,inner sep=1pt] (B2) at (5,6) {$B_2$};
\node
[circle,fill=white,draw=black,text=black,inner sep=1pt] (B1) at (6,4) {$B_1$};
\begin{scope}[thick,>=stealth,darkred]
\draw(A3) edge[->] (2,4);
\draw(2,4) edge[->] (3,4);
\draw(3,4) edge[->] (3,5);
\draw(3,5) edge[->] (4,5);
\draw(4,5) edge[->] (B3);
\draw(3.4,5) node[anchor=south east] {$p_3$};
\end{scope}
\begin{scope}[thick,>=stealth,dbluecolor]
\draw(A2) edge[->] (4,1);
\draw(4,1) edge[->] (4,2);
\draw(4,2) edge[->] (4,3);
\draw(4,3) edge[->] (5,3);
\draw(5,3) edge[->] (5,4);
\draw(5,4) edge[->] (5,5);
\draw(5,5) edge[->] (B2);
\draw(5,4.4) node[anchor=east] {$p_2$};
\end{scope}
\begin{scope}[thick,>=stealth,dgreencolor]
\draw(A1) edge[->] (5,2);
\draw(5,2) edge[->] (5,3);
\draw(5,3) edge[->] (6,3);
\draw(6,3) edge[->] (B1);
\draw(5,2.4) node[anchor=west] {$p_1$};
\end{scope}
\end{tikzpicture}%
\ \ .
\]

\textbf{(c)} The following path tuple is an ipat, too (for several reasons):%
\[%
\begin{tikzpicture}
\draw[densely dotted] (0,0) grid (7.2,7.2);
\draw[->] (0,0) -- (0,7.2);
\draw[->] (0,0) -- (7.2,0);
\foreach\x/\xtext in {0, 1, 2, 3, 4, 5, 6, 7}
\draw(\x cm,1pt) -- (\x cm,-1pt) node[anchor=north] {$\xtext$};
\foreach\y/\ytext in {0, 1, 2, 3, 4, 5, 6, 7}
\draw(1pt,\y cm) -- (-1pt,\y cm) node[anchor=east] {$\ytext$};
\node
[circle,fill=white,draw=black,text=black,inner sep=1pt] (A3) at (2,3) {$A_3$};
\node
[circle,fill=white,draw=black,text=black,inner sep=1pt] (A2) at (3,1) {$A_2$};
\node
[circle,fill=white,draw=black,text=black,inner sep=1pt] (A1) at (1,1) {$A_1$};
\node
[circle,fill=white,draw=black,text=black,inner sep=1pt] (B3) at (4,6) {$B_3$};
\node
[circle,fill=white,draw=black,text=black,inner sep=1pt] (B2) at (5,6) {$B_2$};
\node
[circle,fill=white,draw=black,text=black,inner sep=1pt] (B1) at (4,3) {$B_1$};
\begin{scope}[thick,>=stealth,darkred]
\draw(A3) edge[->] (2,4);
\draw(2,4) edge[->] (3,4);
\draw(3,4) edge[->] (4,4);
\draw(4,4) edge[->] (4,5);
\draw(4,5) edge[->] (B3);
\draw(4,4) node[anchor=west] {$p_3$};
\end{scope}
\begin{scope}[thick,>=stealth,dbluecolor]
\draw(A2) edge[->] (3,2);
\draw(3,2) edge[->, bend left] (3,3);
\draw(3,3) edge[->] (3,4);
\draw(3,4) edge[->] (3,5);
\draw(3,5) edge[->] (4,5);
\draw(4,5) edge[->] (5,5);
\draw(5,5) edge[->] (B2);
\draw(3,5) node[anchor=east] {$p_2$};
\end{scope}
\begin{scope}[thick,>=stealth,dgreencolor]
\draw(A1) edge[->] (1,2);
\draw(1,2) edge[->] (2,2);
\draw(2,2) edge[->] (3,2);
\draw(3,2) edge[->, bend right] (3,3);
\draw(3,3) edge[->] (B1);
\draw(2,2) node[anchor=north] {$p_1$};
\end{scope}
\end{tikzpicture}%
\ \ .
\]
(In this tuple, the paths $p_{1}$ and $p_{2}$ even have an arc in common.
Don't let the picture confuse you: The two curved arcs are actually one and
the same arc of $\mathbb{Z}^{2}$ appearing in two paths, not two different arcs.)
\end{example}

\subsubsection{The LGV lemma for two paths}

As I mentioned, we want to count nipats. Here is a first result on the case
$k=2$:

\begin{proposition}
[LGV lemma for two paths]\label{prop.lgv.2paths.count}Let $\left(
A,A^{\prime}\right)  $ and $\left(  B,B^{\prime}\right)  $ be two $2$-vertices
(i.e., let $A,A^{\prime},B,B^{\prime}$ be four lattice points). Then,%
\begin{align*}
&  \det\left(
\begin{array}
[c]{cc}%
\left(  \text{\# of paths from }A\text{ to }B\right)  & \left(  \text{\# of
paths from }A\text{ to }B^{\prime}\right) \\
\left(  \text{\# of paths from }A^{\prime}\text{ to }B\right)  & \left(
\text{\# of paths from }A^{\prime}\text{ to }B^{\prime}\right)
\end{array}
\right) \\
&  =\left(  \text{\# of nipats from }\left(  A,A^{\prime}\right)  \text{ to
}\left(  B,B^{\prime}\right)  \right) \\
&  \ \ \ \ \ \ \ \ \ \ -\left(  \text{\# of nipats from }\left(  A,A^{\prime
}\right)  \text{ to }\left(  B^{\prime},B\right)  \right)  .
\end{align*}

\end{proposition}

\begin{example}
Let $A=\left(  0,0\right)  $ and $A^{\prime}=\left(  1,1\right)  $ and
$B=\left(  2,2\right)  $ and $B^{\prime}=\left(  3,3\right)  $. Then,
Proposition \ref{prop.lgv.2paths.count} says that%
\begin{align*}
\det\left(
\begin{array}
[c]{cc}%
\dbinom{4}{2} & \dbinom{6}{3}\\
\dbinom{2}{1} & \dbinom{4}{2}%
\end{array}
\right)   &  =\left(  \text{\# of nipats from }\left(  A,A^{\prime}\right)
\text{ to }\left(  B,B^{\prime}\right)  \right) \\
&  \ \ \ \ \ \ \ \ \ \ -\left(  \text{\# of nipats from }\left(  A,A^{\prime
}\right)  \text{ to }\left(  B^{\prime},B\right)  \right)  .
\end{align*}
(The matrix entries on the left hand side have been computed using Proposition
\ref{prop.lgv.1-paths.ct}.)

And indeed, this equality is easily verified. There are $2$ nipats from
$\left(  A,A^{\prime}\right)  $ to $\left(  B,B^{\prime}\right)  $, one of
which is
\[%
\begin{tikzpicture}
\draw[densely dotted] (0,0) grid (3.2,3.2);
\draw[->] (0,0) -- (0,3.2);
\draw[->] (0,0) -- (3.2,0);
\foreach\x/\xtext in {0, 1, 2, 3}
\draw(\x cm,1pt) -- (\x cm,-1pt) node[anchor=north] {$\xtext$};
\foreach\y/\ytext in {0, 1, 2, 3}
\draw(1pt,\y cm) -- (-1pt,\y cm) node[anchor=east] {$\ytext$};
\node
[circle,fill=white,draw=black,text=black,inner sep=1pt] (A') at (1,1) {$A'$};
\node
[circle,fill=white,draw=black,text=black,inner sep=1pt] (A)  at (0,0) {$A$};
\node
[circle,fill=white,draw=black,text=black,inner sep=1pt] (B') at (3,3) {$B'$};
\node
[circle,fill=white,draw=black,text=black,inner sep=1pt] (B)  at (2,2) {$B$};
\begin{scope}[thick,>=stealth,darkred]
\draw(A) edge[->] (0,1);
\draw(0,1) edge[->] (0,2);
\draw(0,2) edge[->] (1,2);
\draw(1,2) edge[->] (B);
\end{scope}
\begin{scope}[thick,>=stealth,dbluecolor]
\draw(A') edge[->] (2,1);
\draw(2,1) edge[->] (3,1);
\draw(3,1) edge[->] (3,2);
\draw(3,2) edge[->] (B');
\end{scope}
\end{tikzpicture}%
\]
while the other is its reflection in the $x=y$ diagonal. There are $6$ nipats
from $\left(  A,A^{\prime}\right)  $ to $\left(  B^{\prime},B\right)  $, three
of which are
\[%
\begin{tikzpicture}
\draw[densely dotted] (0,0) grid (3.2,3.2);
\draw[->] (0,0) -- (0,3.2);
\draw[->] (0,0) -- (3.2,0);
\foreach\x/\xtext in {0, 1, 2, 3}
\draw(\x cm,1pt) -- (\x cm,-1pt) node[anchor=north] {$\xtext$};
\foreach\y/\ytext in {0, 1, 2, 3}
\draw(1pt,\y cm) -- (-1pt,\y cm) node[anchor=east] {$\ytext$};
\node
[circle,fill=white,draw=black,text=black,inner sep=1pt] (A') at (1,1) {$A'$};
\node
[circle,fill=white,draw=black,text=black,inner sep=1pt] (A)  at (0,0) {$A$};
\node
[circle,fill=white,draw=black,text=black,inner sep=1pt] (B') at (3,3) {$B'$};
\node
[circle,fill=white,draw=black,text=black,inner sep=1pt] (B)  at (2,2) {$B$};
\begin{scope}[thick,>=stealth,darkred]
\draw(A) edge[->] (0,1);
\draw(0,1) edge[->] (0,2);
\draw(0,2) edge[->] (1,2);
\draw(1,2) edge[->] (1,3);
\draw(1,3) edge[->] (2,3);
\draw(2,3) edge[->] (B');
\end{scope}
\begin{scope}[thick,>=stealth,dbluecolor]
\draw(A') edge[->] (2,1);
\draw(2,1) edge[->] (B);
\end{scope}
\end{tikzpicture}%
\ \ \ \ \ \ \ \ \ \
\begin{tikzpicture}
\draw[densely dotted] (0,0) grid (3.2,3.2);
\draw[->] (0,0) -- (0,3.2);
\draw[->] (0,0) -- (3.2,0);
\foreach\x/\xtext in {0, 1, 2, 3}
\draw(\x cm,1pt) -- (\x cm,-1pt) node[anchor=north] {$\xtext$};
\foreach\y/\ytext in {0, 1, 2, 3}
\draw(1pt,\y cm) -- (-1pt,\y cm) node[anchor=east] {$\ytext$};
\node
[circle,fill=white,draw=black,text=black,inner sep=1pt] (A') at (1,1) {$A'$};
\node
[circle,fill=white,draw=black,text=black,inner sep=1pt] (A)  at (0,0) {$A$};
\node
[circle,fill=white,draw=black,text=black,inner sep=1pt] (B') at (3,3) {$B'$};
\node
[circle,fill=white,draw=black,text=black,inner sep=1pt] (B)  at (2,2) {$B$};
\begin{scope}[thick,>=stealth,darkred]
\draw(A) edge[->] (0,1);
\draw(0,1) edge[->] (0,2);
\draw(0,2) edge[->] (0,3);
\draw(0,3) edge[->] (1,3);
\draw(1,3) edge[->] (2,3);
\draw(2,3) edge[->] (B');
\end{scope}
\begin{scope}[thick,>=stealth,dbluecolor]
\draw(A') edge[->] (2,1);
\draw(2,1) edge[->] (B);
\end{scope}
\end{tikzpicture}%
\ \ \ \ \ \ \ \ \ \
\begin{tikzpicture}
\draw[densely dotted] (0,0) grid (3.2,3.2);
\draw[->] (0,0) -- (0,3.2);
\draw[->] (0,0) -- (3.2,0);
\foreach\x/\xtext in {0, 1, 2, 3}
\draw(\x cm,1pt) -- (\x cm,-1pt) node[anchor=north] {$\xtext$};
\foreach\y/\ytext in {0, 1, 2, 3}
\draw(1pt,\y cm) -- (-1pt,\y cm) node[anchor=east] {$\ytext$};
\node
[circle,fill=white,draw=black,text=black,inner sep=1pt] (A') at (1,1) {$A'$};
\node
[circle,fill=white,draw=black,text=black,inner sep=1pt] (A)  at (0,0) {$A$};
\node
[circle,fill=white,draw=black,text=black,inner sep=1pt] (B') at (3,3) {$B'$};
\node
[circle,fill=white,draw=black,text=black,inner sep=1pt] (B)  at (2,2) {$B$};
\begin{scope}[thick,>=stealth,darkred]
\draw(A) edge[->] (0,1);
\draw(0,1) edge[->] (0,2);
\draw(0,2) edge[->] (0,3);
\draw(0,3) edge[->] (1,3);
\draw(1,3) edge[->] (2,3);
\draw(2,3) edge[->] (B');
\end{scope}
\begin{scope}[thick,>=stealth,dbluecolor]
\draw(A') edge[->] (1,2);
\draw(1,2) edge[->] (B);
\end{scope}
\end{tikzpicture}%
\]
while the other three are their reflections in the $x=y$ diagonal. The right
hand side of the above equality is thus $2-6=-4$, which is also the left hand side.
\end{example}

\begin{example}
\label{exa.lgv.2paths.2}Let $A=\left(  0,0\right)  $ and $A^{\prime}=\left(
-1,1\right)  $ and $B=\left(  2,2\right)  $ and $B^{\prime}=\left(
0,3\right)  $. Then, the claim of Proposition \ref{prop.lgv.2paths.count}
simplifies, since $\left(  \text{\# of nipats from }\left(  A,A^{\prime
}\right)  \text{ to }\left(  B^{\prime},B\right)  \right)  =0$ in this case.
Here is a picture of the four points that should make this visually clear:
\[%
\begin{tikzpicture}
\draw[densely dotted] (-1.2, -0.2) grid (2.2,3.2);
\node
[circle,fill=white,draw=black,text=black,inner sep=1pt] (A') at (-1,1) {$A'$};
\node
[circle,fill=white,draw=black,text=black,inner sep=1pt] (A)  at (0,0) {$A$};
\node
[circle,fill=white,draw=black,text=black,inner sep=1pt] (B') at (0,3) {$B'$};
\node
[circle,fill=white,draw=black,text=black,inner sep=1pt] (B)  at (2,2) {$B$};
\end{tikzpicture}%
\ \ .
\]

\end{example}

\begin{proof}
[Proof of Proposition \ref{prop.lgv.2paths.count}.]We have%
\begin{align}
&  \det\left(
\begin{array}
[c]{cc}%
\left(  \text{\# of paths from }A\text{ to }B\right)  & \left(  \text{\# of
paths from }A\text{ to }B^{\prime}\right) \\
\left(  \text{\# of paths from }A^{\prime}\text{ to }B\right)  & \left(
\text{\# of paths from }A^{\prime}\text{ to }B^{\prime}\right)
\end{array}
\right) \nonumber\\
&  =\underbrace{\left(  \text{\# of paths from }A\text{ to }B\right)
\cdot\left(  \text{\# of paths from }A^{\prime}\text{ to }B^{\prime}\right)
}_{\substack{=\left(  \text{\# of path tuples from }\left(  A,A^{\prime
}\right)  \text{ to }\left(  B,B^{\prime}\right)  \right)  \\\text{(by the
product rule, since a path tuple from }\left(  A,A^{\prime}\right)  \text{ to
}\left(  B,B^{\prime}\right)  \text{ is just}\\\text{a pair consisting of a
path from }A\text{ to }B\text{ and a path from }A^{\prime}\text{ to }%
B^{\prime}\text{)}}}\nonumber\\
&  \ \ \ \ \ \ \ \ \ \ -\underbrace{\left(  \text{\# of paths from }A\text{ to
}B^{\prime}\right)  \cdot\left(  \text{\# of paths from }A^{\prime}\text{ to
}B\right)  }_{\substack{=\left(  \text{\# of path tuples from }\left(
A,A^{\prime}\right)  \text{ to }\left(  B^{\prime},B\right)  \right)
\\\text{(by the product rule, since a path tuple from }\left(  A,A^{\prime
}\right)  \text{ to }\left(  B^{\prime},B\right)  \text{ is just}\\\text{a
pair consisting of a path from }A\text{ to }B^{\prime}\text{ and a path from
}A^{\prime}\text{ to }B\text{)}}}\nonumber\\
&  =\left(  \text{\# of path tuples from }\left(  A,A^{\prime}\right)  \text{
to }\left(  B,B^{\prime}\right)  \right) \nonumber\\
&  \ \ \ \ \ \ \ \ \ \ -\left(  \text{\# of path tuples from }\left(
A,A^{\prime}\right)  \text{ to }\left(  B^{\prime},B\right)  \right)  .
\label{pf.prop.lgv.2paths.count.1}%
\end{align}
We need to show that on the right hand side, all the intersecting path tuples
cancel each other out (so that only the nipats remain).

Our $k$-vertices are $2$-vertices; thus, our path tuples are pairs. Hence,
such a path tuple $\left(  p,p^{\prime}\right)  $ is intersecting if and only
if $p$ and $p^{\prime}$ have a vertex in common. We shall use these common
vertices to define a sign-reversing involution on the intersecting path
tuples. Specifically, we do the following:

Define the set%
\begin{align*}
\mathcal{A}  &  :=\left\{  \text{path tuples from }\left(  A,A^{\prime
}\right)  \text{ to }\left(  B,B^{\prime}\right)  \right\} \\
&  \ \ \ \ \ \ \ \ \ \ \sqcup\left\{  \text{path tuples from }\left(
A,A^{\prime}\right)  \text{ to }\left(  B^{\prime},B\right)  \right\}  .
\end{align*}
Here, the symbol \textquotedblleft$\sqcup$\textquotedblright\ means
\textquotedblleft disjoint union (of sets)\textquotedblright, which is a way
of uniting two sets without removing duplicates (i.e., even if the sets are
not disjoint, we treat them as disjoint for the purpose of the union, and
therefore include two copies of each common element). As a consequence of us
taking the disjoint union, each path tuple in $\mathcal{A}$ \textquotedblleft
remembers\textquotedblright\ whether it comes from the set \newline$\left\{
\text{path tuples from }\left(  A,A^{\prime}\right)  \text{ to }\left(
B,B^{\prime}\right)  \right\}  $ or from the set \newline$\left\{  \text{path
tuples from }\left(  A,A^{\prime}\right)  \text{ to }\left(  B^{\prime
},B\right)  \right\}  $ (and if these two sets have a path tuple in common,
then $\mathcal{A}$ has two copies of it, each of which remembers from which
set it comes). However, in practice, this is barely relevant: Indeed, the only
case in which the sets \newline$\left\{  \text{path tuples from }\left(
A,A^{\prime}\right)  \text{ to }\left(  B,B^{\prime}\right)  \right\}  $ and
$\left\{  \text{path tuples from }\left(  A,A^{\prime}\right)  \text{ to
}\left(  B^{\prime},B\right)  \right\}  $ can fail to be disjoint is the case
when $B=B^{\prime}$; however, in this case, the claim we are proving is
trivial anyway, since there are no nipats\footnote{Indeed, if $p$ and
$p^{\prime}$ are two paths with the same destination, then $p$ and $p^{\prime
}$ automatically have a vertex in common.}, and our matrix has determinant $0$
(since it has two equal columns).

Define a subset $\mathcal{X}$ of $\mathcal{A}$ by%
\[
\mathcal{X}:=\left\{  \text{ipats in }\mathcal{A}\right\}  =\left\{  \left(
p,p^{\prime}\right)  \in\mathcal{A}\ \mid\ p\text{ and }p^{\prime}\text{ have
a vertex in common}\right\}  .
\]
Hence, $\mathcal{A}\setminus\mathcal{X}=\left\{  \text{nipats in }%
\mathcal{A}\right\}  $.

For each $\left(  p,p^{\prime}\right)  \in\mathcal{A}$, we set%
\[
\operatorname*{sign}\left(  p,p^{\prime}\right)  :=%
\begin{cases}
1, & \text{if }\left(  p,p^{\prime}\right)  \text{ is a path tuple from
}\left(  A,A^{\prime}\right)  \text{ to }\left(  B,B^{\prime}\right)  ;\\
-1, & \text{if }\left(  p,p^{\prime}\right)  \text{ is a path tuple from
}\left(  A,A^{\prime}\right)  \text{ to }\left(  B^{\prime},B\right)  .
\end{cases}
\]
(This is well-defined, because each $\left(  p,p^{\prime}\right)
\in\mathcal{A}$ is either a path tuple from $\left(  A,A^{\prime}\right)  $ to
$\left(  B,B^{\prime}\right)  $ or a path tuple from $\left(  A,A^{\prime
}\right)  $ to $\left(  B^{\prime},B\right)  $ but never both at the same
time\footnote{because we took the \textbf{disjoint} union of $\left\{
\text{path tuples from }\left(  A,A^{\prime}\right)  \text{ to }\left(
B,B^{\prime}\right)  \right\}  $ and $\left\{  \text{path tuples from }\left(
A,A^{\prime}\right)  \text{ to }\left(  B^{\prime},B\right)  \right\}  $}.)
Thus,
\begin{align*}
&  \left(  \text{\# of path tuples from }\left(  A,A^{\prime}\right)  \text{
to }\left(  B,B^{\prime}\right)  \right) \\
&  \ \ \ \ \ \ \ \ \ \ -\left(  \text{\# of path tuples from }\left(
A,A^{\prime}\right)  \text{ to }\left(  B^{\prime},B\right)  \right) \\
&  =\sum_{\left(  p,p^{\prime}\right)  \in\mathcal{A}}\operatorname*{sign}%
\left(  p,p^{\prime}\right)
\end{align*}
and%
\begin{align*}
&  \left(  \text{\# of nipats from }\left(  A,A^{\prime}\right)  \text{ to
}\left(  B,B^{\prime}\right)  \right) \\
&  \ \ \ \ \ \ \ \ \ \ -\left(  \text{\# of nipats from }\left(  A,A^{\prime
}\right)  \text{ to }\left(  B^{\prime},B\right)  \right) \\
&  =\sum_{\left(  p,p^{\prime}\right)  \in\mathcal{A}\setminus\mathcal{X}%
}\operatorname*{sign}\left(  p,p^{\prime}\right)  \ \ \ \ \ \ \ \ \ \ \left(
\text{since }\mathcal{A}\setminus\mathcal{X}=\left\{  \text{nipats in
}\mathcal{A}\right\}  \right)  .
\end{align*}
We want to prove that the left hand sides of these two equalities are equal.
Thus, it clearly suffices to show that the right hand sides are equal. By
Lemma \ref{lem.sign.cancel2}, it suffices to find a sign-reversing involution
$f:\mathcal{X}\rightarrow\mathcal{X}$ that has no fixed points.

So let us define our sign-reversing involution $f:\mathcal{X}\rightarrow
\mathcal{X}$. For each path tuple $\left(  p,p^{\prime}\right)  \in
\mathcal{X}$, we define $f\left(  p,p^{\prime}\right)  $ as follows:

\begin{itemize}
\item Since $\left(  p,p^{\prime}\right)  \in\mathcal{X}$, the paths $p$ and
$p^{\prime}$ have a vertex in common. There might be several; let $v$ be the
first one. (The first one on $p$ or the first one on $p^{\prime}$ ? Doesn't
matter, because these are the same thing. Indeed, if the first vertex on $p$
that is contained in $p^{\prime}$ was different from the first vertex on
$p^{\prime}$ that is contained in $p$, then we could obtain a nontrivial
circuit of our digraph by walking from the former vertex to the latter vertex
along $p$ and then back along $p^{\prime}$; but this is impossible, since our
digraph is acyclic.)

We call this vertex $v$ the \emph{first intersection} of $\left(  p,p^{\prime
}\right)  $.

\item Call the part of $p$ that comes after $v$ the \emph{tail} of $p$. Call
the part of $p$ that comes before $v$ the \emph{head} of $p$.

Call the part of $p^{\prime}$ that comes after $v$ the \emph{tail} of
$p^{\prime}$. Call the part of $p^{\prime}$ that comes before $v$ the
\emph{head} of $p^{\prime}$.

\item Now, we exchange the tails of the paths $p$ and $p^{\prime}$. That is,
we set%
\begin{align*}
q  &  :=\left(  \text{head of }p\right)  \cup\left(  \text{tail of }p^{\prime
}\right)  \ \ \ \ \ \ \ \ \ \ \text{and}\\
q^{\prime}  &  :=\left(  \text{head of }p^{\prime}\right)  \cup\left(
\text{tail of }p\right)
\end{align*}
(where the symbol \textquotedblleft$\cup$\textquotedblright\ means combining a
path ending at $v$ with a path starting at $v$ in the obvious way), and set
$f\left(  p,p^{\prime}\right)  :=\left(  q,q^{\prime}\right)  $.
\end{itemize}

Thus, we have defined a map $f:\mathcal{X}\rightarrow\mathcal{X}$ (in a
moment, we will explain why it is well-defined). Here is an example:
\[
\text{If }\left(  p,p^{\prime}\right)  =%
\raisebox{-6pc}{
\begin{tikzpicture}
\draw[densely dotted] (0.8,0.8) grid (5.2,6.2);
\node
[circle,fill=white,draw=black,text=black,inner sep=1pt] (A') at (1,3) {$A^{\prime
}$};
\node
[circle,fill=white,draw=black,text=black,inner sep=1pt] (A) at (2,1) {$A$};
\node
[circle,fill=white,draw=black,text=black,inner sep=1pt] (B') at (5,4) {$B^{\prime
}$};
\node
[circle,fill=white,draw=black,text=black,inner sep=1pt] (B) at (5,6) {$B$};
\begin{scope}[thick,>=stealth,darkred]
\draw(A) edge[->] (2,2);
\draw(2,2) node[anchor=west] {$p$};
\draw(2,2) edge[->] (2,3);
\draw(2,3) edge[->] (3,3);
\draw(3,3) edge[->] (4,3);
\draw(4,3) edge[->] (4,4);
\draw(4,4) edge[->] (4,5);
\draw(4,5) edge[->] (5,5);
\draw(5,5) edge[->] (B);
\draw(4,5) node[anchor=south] {$p$};
\end{scope}
\begin{scope}[thick,>=stealth,dbluecolor]
\draw(A') edge[->] (2,3);
\draw(1.6,3) node[anchor=north] {$p^{\prime}$};
\draw(2,3) edge[->] (2,4);
\draw(2,4) node[anchor=east] {$p^{\prime}$};
\draw(2,4) edge[->] (3,4);
\draw(3,4) edge[->] (4,4);
\draw(4,4) edge[->] (B');
\draw(4.4,4) node[anchor=north] {$p^{\prime}$};
\end{scope}
\end{tikzpicture}
}%
\text{, then }\left(  q,q^{\prime}\right)  =%
\raisebox{-6pc}{\begin{tikzpicture}
\draw[densely dotted] (0.8,0.8) grid (5.2,6.2);
\node
[circle,fill=white,draw=black,text=black,inner sep=1pt] (A') at (1,3) {$A^{\prime
}$};
\node
[circle,fill=white,draw=black,text=black,inner sep=1pt] (A) at (2,1) {$A$};
\node
[circle,fill=white,draw=black,text=black,inner sep=1pt] (B') at (5,4) {$B^{\prime
}$};
\node
[circle,fill=white,draw=black,text=black,inner sep=1pt] (B) at (5,6) {$B$};
\begin{scope}[thick,>=stealth,darkred]
\draw(A) edge[->] (2,2);
\draw(2,2) node[anchor=west] {$q$};
\draw(2,2) edge[->] (2,3);
\draw(2,3) edge[->] (2,4);
\draw(2,4) edge[->] (3,4);
\draw(2,4) node[anchor=east] {$q$};
\draw(3,4) edge[->] (4,4);
\draw(4,4) edge[->] (B');
\draw(4.4,4) node[anchor=north] {$q$};
\end{scope}
\begin{scope}[thick,>=stealth,dbluecolor]
\draw(A') edge[->] (2,3);
\draw(1.6,3) node[anchor=north] {$q^{\prime}$};
\draw(2,3) edge[->] (3,3);
\draw(3,3) edge[->] (4,3);
\draw(4,3) edge[->] (4,4);
\draw(4,4) edge[->] (4,5);
\draw(4,5) edge[->] (5,5);
\draw(5,5) edge[->] (B);
\draw(4,5) node[anchor=south] {$q^{\prime}$};
\end{scope}
\end{tikzpicture}
}%
\text{.}%
\]
Here is the same configuration, with the point $v$ marked and with the tails
of the two paths drawn extra-thick:%
\[%
\raisebox{-6pc}{
\begin{tikzpicture}
\draw[densely dotted] (0.8,0.8) grid (5.2,6.2);
\node
[circle,fill=white,draw=black,text=black,inner sep=1pt] (A') at (1,3) {$A^{\prime
}$};
\node
[circle,fill=white,draw=black,text=black,inner sep=1pt] (A) at (2,1) {$A$};
\node
[circle,fill=white,draw=black,text=black,inner sep=1pt] (B') at (5,4) {$B^{\prime
}$};
\node
[circle,fill=white,draw=black,text=black,inner sep=1pt] (B) at (5,6) {$B$};
\begin{scope}[thick,>=stealth,darkred]
\draw(A) edge[->] (2,2);
\draw(2,2) node[anchor=west] {$p$};
\draw(2,2) edge[->] (2,3);
\draw(2,3) edge[ultra thick, ->] (3,3);
\draw(3,3) edge[ultra thick, ->] (4,3);
\draw(4,3) edge[ultra thick, ->] (4,4);
\draw(4,4) edge[ultra thick, ->] (4,5);
\draw(4,5) edge[ultra thick, ->] (5,5);
\draw(5,5) edge[ultra thick, ->] (B);
\draw(4,5) node[anchor=south] {$p$};
\end{scope}
\begin{scope}[thick,>=stealth,dbluecolor]
\draw(A') edge[->] (2,3);
\draw(1.6,3) node[anchor=north] {$p^{\prime}$};
\draw(2,3) edge[ultra thick, ->] (2,4);
\draw(2,4) edge[ultra thick, ->] (3,4);
\draw(2,4) node[anchor=east] {$p^{\prime}$};
\draw(3,4) edge[ultra thick, ->] (4,4);
\draw(4,4) edge[ultra thick, ->] (B');
\draw(4.4,4) node[anchor=north] {$p^{\prime}$};
\end{scope}
\node
[circle,fill=green,draw=black,text=black,inner sep=1pt] (v) at (2,3) {$v$};
\end{tikzpicture}
}%
\ \ .
\]

We note that if $\left(  p,p^{\prime}\right)  \in\mathcal{X}$ is an ipat from
$\left(  A,A^{\prime}\right)  $ to $\left(  B,B^{\prime}\right)  $, then
$f\left(  p,p^{\prime}\right)  =\left(  q,q^{\prime}\right)  $ is an ipat from
$\left(  A,A^{\prime}\right)  $ to $\left(  B^{\prime},B\right)  $ (because by
exchanging the tails of $p$ and $p^{\prime}$, we have caused the two paths to
exchange their destinations as well), and vice versa. Thus, $f\left(
p,p^{\prime}\right)  \in\mathcal{X}$ whenever $\left(  p,p^{\prime}\right)
\in\mathcal{X}$. This shows that the map $f:\mathcal{X}\rightarrow\mathcal{X}$
is well-defined.

Furthermore, let $\left(  p,p^{\prime}\right)  \in\mathcal{X}$ be any ipat,
and let $\left(  q,q^{\prime}\right)  =f\left(  p,p^{\prime}\right)  $. Then,
the vertex $v$ chosen in the definition of $f\left(  p,p^{\prime}\right)  $
(that is, the first intersection of $\left(  p,p^{\prime}\right)  $) is still
the first intersection of $\left(  q,q^{\prime}\right)  $ (because when we
exchange the tails of $p$ and $p^{\prime}$, we do not change their heads, and
thus the resulting paths $q$ and $q^{\prime}$ still do not intersect until
$v$), and therefore this vertex gets chosen again if we apply our map $f$ to
$\left(  q,q^{\prime}\right)  $. As a consequence, $f\left(  q,q^{\prime
}\right)  $ is again $\left(  p,p^{\prime}\right)  $ (since exchanging the
tails of $q$ and $q^{\prime}$ simply undoes the changes incurred when we
exchanged the tails of $p$ and $p^{\prime}$). Thus,
\[
\left(  f\circ f\right)  \left(  p,p^{\prime}\right)  =f\left(
\underbrace{f\left(  p,p^{\prime}\right)  }_{=\left(  q,q^{\prime}\right)
}\right)  =f\left(  q,q^{\prime}\right)  =\left(  p,p^{\prime}\right)  .
\]

Forget that we fixed $\left(  p,p^{\prime}\right)  $. We thus have shown that
$\left(  f\circ f\right)  \left(  p,p^{\prime}\right)  =\left(  p,p^{\prime
}\right)  $ for each $\left(  p,p^{\prime}\right)  \in\mathcal{X}$. Hence,
$f\circ f=\operatorname*{id}$. In other words, $f$ is an involution on
$\mathcal{X}$. Moreover, this involution $f$ is sign-reversing (i.e.,
satisfies $\operatorname*{sign}\left(  f\left(  p,p^{\prime}\right)  \right)
=-\operatorname*{sign}\left(  p,p^{\prime}\right)  $ for any $\left(
p,p^{\prime}\right)  \in\mathcal{X}$)\ \ \ \ \footnote{This is because we have
observed above that if $\left(  p,p^{\prime}\right)  \in\mathcal{X}$ is an
ipat from $\left(  A,A^{\prime}\right)  $ to $\left(  B,B^{\prime}\right)  $,
then $f\left(  p,p^{\prime}\right)  =\left(  q,q^{\prime}\right)  $ is an ipat
from $\left(  A,A^{\prime}\right)  $ to $\left(  B^{\prime},B\right)  $, and
vice versa.}. As a consequence of the latter fact, we see that $f$ has no
fixed points (i.e., that we have $f\left(  p,p^{\prime}\right)  \neq\left(
p,p^{\prime}\right)  $ for any $\left(  p,p^{\prime}\right)  \in\mathcal{X}$).
Hence, Lemma \ref{lem.sign.cancel2} yields%
\begin{equation}
\sum_{\left(  p,p^{\prime}\right)  \in\mathcal{A}}\operatorname*{sign}\left(
p,p^{\prime}\right)  =\sum_{\left(  p,p^{\prime}\right)  \in\mathcal{A}%
\setminus\mathcal{X}}\operatorname*{sign}\left(  p,p^{\prime}\right)  .
\label{pf.prop.lgv.2paths.count.at}%
\end{equation}
As we have explained above, the left hand side of this equality is%
\begin{align*}
&  \left(  \text{\# of path tuples from }\left(  A,A^{\prime}\right)  \text{
to }\left(  B,B^{\prime}\right)  \right) \\
&  \ \ \ \ \ \ \ \ \ \ -\left(  \text{\# of path tuples from }\left(
A,A^{\prime}\right)  \text{ to }\left(  B^{\prime},B\right)  \right)  ,
\end{align*}
whereas its right hand side is%
\begin{align*}
&  \left(  \text{\# of nipats from }\left(  A,A^{\prime}\right)  \text{ to
}\left(  B,B^{\prime}\right)  \right) \\
&  \ \ \ \ \ \ \ \ \ \ -\left(  \text{\# of nipats from }\left(  A,A^{\prime
}\right)  \text{ to }\left(  B^{\prime},B\right)  \right)
\end{align*}
(since $\mathcal{A}\setminus\mathcal{X}=\left\{  \text{nipats in }%
\mathcal{A}\right\}  $). Hence, (\ref{pf.prop.lgv.2paths.count.at}) rewrites
as%
\begin{align*}
&  \left(  \text{\# of path tuples from }\left(  A,A^{\prime}\right)  \text{
to }\left(  B,B^{\prime}\right)  \right) \\
&  \ \ \ \ \ \ \ \ \ \ -\left(  \text{\# of path tuples from }\left(
A,A^{\prime}\right)  \text{ to }\left(  B^{\prime},B\right)  \right) \\
&  =\left(  \text{\# of nipats from }\left(  A,A^{\prime}\right)  \text{ to
}\left(  B,B^{\prime}\right)  \right) \\
&  \ \ \ \ \ \ \ \ \ \ -\left(  \text{\# of nipats from }\left(  A,A^{\prime
}\right)  \text{ to }\left(  B^{\prime},B\right)  \right)  .
\end{align*}
In view of (\ref{pf.prop.lgv.2paths.count.1}), this rewrites as%
\begin{align*}
&  \det\left(
\begin{array}
[c]{cc}%
\left(  \text{\# of paths from }A\text{ to }B\right)  & \left(  \text{\# of
paths from }A\text{ to }B^{\prime}\right) \\
\left(  \text{\# of paths from }A^{\prime}\text{ to }B\right)  & \left(
\text{\# of paths from }A^{\prime}\text{ to }B^{\prime}\right)
\end{array}
\right) \\
&  =\left(  \text{\# of nipats from }\left(  A,A^{\prime}\right)  \text{ to
}\left(  B,B^{\prime}\right)  \right) \\
&  \ \ \ \ \ \ \ \ \ \ -\left(  \text{\# of nipats from }\left(  A,A^{\prime
}\right)  \text{ to }\left(  B^{\prime},B\right)  \right)  .
\end{align*}
This completes our proof of Proposition \ref{prop.lgv.2paths.count}.
\end{proof}

As in Example \ref{exa.lgv.2paths.2}, the proposition becomes particularly
nice when we have $\left(  \text{\# of nipats from }\left(  A,A^{\prime
}\right)  \text{ to }\left(  B^{\prime},B\right)  \right)  =0$. Here is a
sufficient criterion for when this happens:

\begin{proposition}
[baby Jordan curve theorem]\label{prop.lgv.jordan-2}Let $A$, $B$, $A^{\prime}$
and $B^{\prime}$ be four lattice points satisfying%
\begin{align}
\operatorname*{x}\left(  A^{\prime}\right)   &  \leq\operatorname*{x}\left(
A\right)  ,\ \ \ \ \ \ \ \ \ \ \operatorname*{y}\left(  A^{\prime}\right)
\geq\operatorname*{y}\left(  A\right)  ,\label{eq.prop.lgv.jordan-2.1}\\
\operatorname*{x}\left(  B^{\prime}\right)   &  \leq\operatorname*{x}\left(
B\right)  ,\ \ \ \ \ \ \ \ \ \ \operatorname*{y}\left(  B^{\prime}\right)
\geq\operatorname*{y}\left(  B\right)  . \label{eq.prop.lgv.jordan-2.2}%
\end{align}
Here, $\operatorname*{x}\left(  P\right)  $ and $\operatorname*{y}\left(
P\right)  $ denote the two coordinates of any point $P\in\mathbb{Z}^{2}$.

Let $p$ be any path from $A$ to $B^{\prime}$. Let $p^{\prime}$ be any path
from $A^{\prime}$ to $B$. Then, $p$ and $p^{\prime}$ have a vertex in common.
\end{proposition}

Note that the condition (\ref{eq.prop.lgv.jordan-2.1}) can be restated as
\textquotedblleft the point $A^{\prime}$ lies weakly northwest of
$A$\textquotedblright, where \textquotedblleft weakly northwest of
$A$\textquotedblright\ allows for the options \textquotedblleft due north of
$A$\textquotedblright, \textquotedblleft due west of $A$\textquotedblright%
\ and \textquotedblleft at $A$\textquotedblright. Likewise,
(\ref{eq.prop.lgv.jordan-2.2}) can be restated as \textquotedblleft the point
$B^{\prime}$ lies weakly northwest of $B$\textquotedblright. The following
picture illustrates the situation of Proposition \ref{prop.lgv.jordan-2}:%
\[%
\begin{tikzpicture}
\draw[densely dotted] (0,0) grid (7.2,7.2);
\draw[->] (0,0) -- (0,7.2);
\draw[->] (0,0) -- (7.2,0);
\foreach\x/\xtext in {0, 1, 2, 3, 4, 5, 6, 7}
\draw(\x cm,1pt) -- (\x cm,-1pt) node[anchor=north] {$\xtext$};
\foreach\y/\ytext in {0, 1, 2, 3, 4, 5, 6, 7}
\draw(1pt,\y cm) -- (-1pt,\y cm) node[anchor=east] {$\ytext$};
\node
[circle,fill=white,draw=black,text=black,inner sep=1pt] (A') at (1,3) {$A^{\prime
}$};
\node
[circle,fill=white,draw=black,text=black,inner sep=1pt] (A) at (2,1) {$A$};
\node
[circle,fill=white,draw=black,text=black,inner sep=1pt] (B') at (5,6) {$B^{\prime
}$};
\node
[circle,fill=white,draw=black,text=black,inner sep=1pt] (B) at (6,4) {$B$};
\begin{scope}[thick,>=stealth,darkred]
\draw(A) edge[->] (3,1);
\draw(2.6,1) node[anchor=north] {$p$};
\draw(3,1) edge[->] (3,2);
\draw(3,2) edge[->] (4,2);
\draw(4,2) edge[->] (4,3);
\draw(4,3) edge[->] (4,4);
\draw(4,4) edge[->] (4,5);
\draw(4,5) edge[->] (5,5);
\draw(5,5) edge[->] (B');
\draw(5,5.4) node[anchor=east] {$p$};
\end{scope}
\begin{scope}[thick,>=stealth,dbluecolor]
\draw(A') edge[->] (2,3);
\draw(1.6,3) node[anchor=north] {$p^{\prime}$};
\draw(2,3) edge[->] (2,4);
\draw(2,4) edge[->] (3,4);
\draw(3,4) edge[->] (4,4);
\draw(4,4) edge[->] (5,4);
\draw(5,4) edge[->] (B);
\draw(5.4,4) node[anchor=north] {$p^{\prime}$};
\end{scope}
\end{tikzpicture}%
\]

Proposition \ref{prop.lgv.jordan-2} has an intuitive plausibility to it (one
can think of the path $p$ as creating a \textquotedblleft
river\textquotedblright\ that the path $p^{\prime}$ must necessarily cross
somewhere), but it is not obvious from a mathematical perspective. We give a
rigorous proof in Section \ref{sec.details.det.comb}.

Proposition \ref{prop.lgv.2paths.count} is just the $k=2$ case of a more
general theorem that we will soon derive; however, it already has a nice application:

\begin{corollary}
\label{cor.lgv.binom-unimod}Let $n,k\in\mathbb{N}$. Then, $\dbinom{n}{k}%
^{2}\geq\dbinom{n}{k-1}\cdot\dbinom{n}{k+1}$.
\end{corollary}

\begin{proof}
[Proof of Corollary \ref{cor.lgv.binom-unimod} (sketched).]This is easy to see
algebraically, but here is a combinatorial proof: Define four lattice points
$A=\left(  1,0\right)  $ and $A^{\prime}=\left(  0,1\right)  $ and $B=\left(
k+1,n-k\right)  $ and $B^{\prime}=\left(  k,n-k+1\right)  $. Then, Proposition
\ref{prop.lgv.2paths.count} yields%
\begin{align*}
&  \det\left(
\begin{array}
[c]{cc}%
\left(  \text{\# of paths from }A\text{ to }B\right)  & \left(  \text{\# of
paths from }A\text{ to }B^{\prime}\right) \\
\left(  \text{\# of paths from }A^{\prime}\text{ to }B\right)  & \left(
\text{\# of paths from }A^{\prime}\text{ to }B^{\prime}\right)
\end{array}
\right) \\
&  =\left(  \text{\# of nipats from }\left(  A,A^{\prime}\right)  \text{ to
}\left(  B,B^{\prime}\right)  \right)  -\underbrace{\left(  \text{\# of nipats
from }\left(  A,A^{\prime}\right)  \text{ to }\left(  B^{\prime},B\right)
\right)  }_{\substack{=0\\\text{(since Proposition \ref{prop.lgv.jordan-2}
yields that any path}\\\text{from }A\text{ to }B^{\prime}\text{ and any path
from }A^{\prime}\text{ to }B\text{ have a}\\\text{vertex in common)}}}\\
&  =\left(  \text{\# of nipats from }\left(  A,A^{\prime}\right)  \text{ to
}\left(  B,B^{\prime}\right)  \right)  \geq0.
\end{align*}
However, Proposition \ref{prop.lgv.1-paths.ct} yields
\begin{align*}
&  \left(
\begin{array}
[c]{cc}%
\left(  \text{\# of paths from }A\text{ to }B\right)  & \left(  \text{\# of
paths from }A\text{ to }B^{\prime}\right) \\
\left(  \text{\# of paths from }A^{\prime}\text{ to }B\right)  & \left(
\text{\# of paths from }A^{\prime}\text{ to }B^{\prime}\right)
\end{array}
\right) \\
&  =\left(
\begin{array}
[c]{cc}%
\dbinom{n}{k} & \dbinom{n}{k-1}\\
\dbinom{n}{k+1} & \dbinom{n}{k}%
\end{array}
\right)  ,
\end{align*}
so the determinant on the left hand side is $\dbinom{n}{k}^{2}-\dbinom{n}%
{k-1}\cdot\dbinom{n}{k+1}$. Thus, we have obtained%
\[
\dbinom{n}{k}^{2}-\dbinom{n}{k-1}\cdot\dbinom{n}{k+1}\geq0,
\]
and this proves Corollary \ref{cor.lgv.binom-unimod}.
\end{proof}

Corollary \ref{cor.lgv.binom-unimod} is often stated as \textquotedblleft the
sequence $\dbinom{n}{0},\dbinom{n}{1},\ldots,\dbinom{n}{n}$ is
log-concave\textquotedblright. There are many more log-concave sequences in
combinatorics (see, e.g., \cite{Sagan19}, \cite{Stanle89} and \cite{Brande14}
for more).

\subsubsection{The LGV lemma for $k$ paths}

The following proposition -- which is one of the weakest forms of the
\emph{LGV lemma} (short for \emph{Lindstr\"{o}m--Gessel--Viennot lemma}) --
extends the logic of Proposition \ref{prop.lgv.2paths.count} to nipats between
$k$-vertices for general $k$).

\begin{proposition}
[LGV lemma, lattice counting version]\label{prop.lgv.kpaths.count}Let
$k\in\mathbb{N}$. Let $\mathbf{A}=\left(  A_{1},A_{2},\ldots,A_{k}\right)  $
and $\mathbf{B}=\left(  B_{1},B_{2},\ldots,B_{k}\right)  $ be two
$k$-vertices. Then,%
\begin{align*}
&  \det\left(  \left(  \text{\# of paths from }A_{i}\text{ to }B_{j}\right)
_{1\leq i\leq k,\ 1\leq j\leq k}\right) \\
&  =\sum_{\sigma\in S_{k}}\left(  -1\right)  ^{\sigma}\left(  \text{\# of
nipats from }\mathbf{A}\text{ to }\sigma\left(  \mathbf{B}\right)  \right)  .
\end{align*}

\end{proposition}

The right hand side of this equality can be viewed as a signed count of all
$k$-tuples of paths that start at the points $A_{1},A_{2},\ldots,A_{k}$ in
\textbf{this} order, but end at the points $B_{1},B_{2},\ldots,B_{k}$ in
\textbf{some} order. For example, for $k=3$, the claim of Proposition
\ref{prop.lgv.kpaths.count} takes the form%
\begin{align*}
&  \det\left(  \left(  \text{\# of paths from }A_{i}\text{ to }B_{j}\right)
_{1\leq i\leq3,\ 1\leq j\leq3}\right) \\
&  =\left(  \text{\# of nipats from }\left(  A_{1},A_{2},A_{3}\right)  \text{
to }\left(  B_{1},B_{2},B_{3}\right)  \right) \\
&  \ \ \ \ \ \ \ \ \ \ -\left(  \text{\# of nipats from }\left(  A_{1}%
,A_{2},A_{3}\right)  \text{ to }\left(  B_{1},B_{3},B_{2}\right)  \right) \\
&  \ \ \ \ \ \ \ \ \ \ -\left(  \text{\# of nipats from }\left(  A_{1}%
,A_{2},A_{3}\right)  \text{ to }\left(  B_{2},B_{1},B_{3}\right)  \right) \\
&  \ \ \ \ \ \ \ \ \ \ +\left(  \text{\# of nipats from }\left(  A_{1}%
,A_{2},A_{3}\right)  \text{ to }\left(  B_{2},B_{3},B_{1}\right)  \right) \\
&  \ \ \ \ \ \ \ \ \ \ +\left(  \text{\# of nipats from }\left(  A_{1}%
,A_{2},A_{3}\right)  \text{ to }\left(  B_{3},B_{1},B_{2}\right)  \right) \\
&  \ \ \ \ \ \ \ \ \ \ -\left(  \text{\# of nipats from }\left(  A_{1}%
,A_{2},A_{3}\right)  \text{ to }\left(  B_{3},B_{2},B_{1}\right)  \right)  .
\end{align*}

\begin{proof}
[Proof of Proposition \ref{prop.lgv.kpaths.count} (sketched).]We adapt the
idea of our proof of Proposition \ref{prop.lgv.2paths.count}, but we have to
be more systematic now. Define a set%
\[
\mathcal{A}:=\left\{  \left(  \sigma,\mathbf{p}\right)  \ \mid\ \sigma\in
S_{k}\text{, and}\ \mathbf{p}\text{ is a path tuple from }\mathbf{A}\text{ to
}\sigma\left(  \mathbf{B}\right)  \right\}
\]
\footnote{This set $\mathcal{A}$ is meant to generalize the set $\mathcal{A}$
that was used in the proof of Proposition \ref{prop.lgv.2paths.count}. The
reason why we are defining it to be
\begin{align*}
&  \left\{  \left(  \sigma,\mathbf{p}\right)  \ \mid\ \sigma\in S_{k}\text{,
and}\ \mathbf{p}\text{ is a path tuple from }\mathbf{A}\text{ to }%
\sigma\left(  \mathbf{B}\right)  \right\} \\
&  \text{instead of }\left\{  \mathbf{p}\ \mid\ \mathbf{p}\text{ is a path
tuple from }\mathbf{A}\text{ to }\sigma\left(  \mathbf{B}\right)  \text{ for
some }\sigma\in S_{k}\right\}
\end{align*}
is to make sure that each path tuple in $\mathcal{A}$ \textquotedblleft
remembers\textquotedblright\ which permutation $\sigma\in S_{k}$ it comes
from. (This is the same rationale that caused us to take the disjoint union in
the proof of Proposition \ref{prop.lgv.2paths.count}; but now we are
explicitly inserting the $\sigma$ into the elements of $\mathcal{A}$ rather
than handwaving about disjoint unions.)}. Define a subset $\mathcal{X}$ of
$\mathcal{A}$ by\footnote{Here we are saying that a pair $\left(
\sigma,\mathbf{p}\right)  \in\mathcal{A}$ is an \emph{ipat} if the path tuple
$\mathbf{p}$ is an ipat, and we are saying that a pair $\left(  \sigma
,\mathbf{p}\right)  \in\mathcal{A}$ is a \emph{nipat} if the path tuple
$\mathbf{p}$ is a nipat. This is a bit sloppy ($\sigma$ has nothing to do with
whether $\mathbf{p}$ is an ipat or a nipat), but we hope that no confusion
will ensue.}%
\[
\mathcal{X}:=\left\{  \text{ipats in }\mathcal{A}\right\}  =\left\{  \left(
\sigma,\mathbf{p}\right)  \in\mathcal{A}\ \mid\ \mathbf{p}\text{ is
intersecting}\right\}  .
\]
Set%
\[
\operatorname*{sign}\left(  \sigma,\mathbf{p}\right)  :=\left(  -1\right)
^{\sigma}\ \ \ \ \ \ \ \ \ \ \text{for each }\left(  \sigma,\mathbf{p}\right)
\in\mathcal{A}.
\]

Again, we need to find a sign-reversing involution $f:\mathcal{X}%
\rightarrow\mathcal{X}$ that has no fixed points.

Again, we construct this involution by exchanging the tails of two
intersecting paths\footnote{This is probably obvious, but just in case: We say
that two paths \emph{intersect} if they have a vertex in common.} in our path
tuple. There is a complication now, due to the fact that there might be
several pairs of intersecting paths. We have to come up with a rule for
picking one such pair so that when we apply $f$ again to the result of the
exchange, then we again pick the same pair. Otherwise, $f$ won't be an involution!

There are different ways to do this. Here is one: If $\left(  \sigma,\left(
p_{1},p_{2},\ldots,p_{k}\right)  \right)  \in\mathcal{X}$, then we construct
$f\left(  \sigma,\left(  p_{1},p_{2},\ldots,p_{k}\right)  \right)
\in\mathcal{X}$ as follows:

\begin{itemize}
\item We say that a point $u$ is \emph{crowded} if it is contained in at least
two of our paths $p_{1},p_{2},\ldots,p_{k}$. Since $\left(  p_{1},p_{2}%
,\ldots,p_{k}\right)  $ is intersecting, there exists at least one crowded point.

\item We pick the smallest $i\in\left[  k\right]  $ such that $p_{i}$ contains
a crowded point.

\item Then, we pick the first crowded point $v$ on $p_{i}$.

\item Then, we pick the largest $j\in\left[  k\right]  $ such that $v$ belongs
to $p_{j}$. (Note that $j>i$, since $v$ is crowded.)

\item Call the part of $p_{i}$ that comes after $v$ the \emph{tail} of $p_{i}%
$. Call the part of $p_{i}$ that comes before $v$ the \emph{head} of $p_{i}$.

Call the part of $p_{j}$ that comes after $v$ the \emph{tail} of $p_{j}$. Call
the part of $p_{j}$ that comes before $v$ the \emph{head} of $p_{j}$.

\item Then, we exchange the tails of the paths $p_{i}$ and $p_{j}$ (while
leaving all other paths unchanged).

\item We let $\left(  q_{1},q_{2},\ldots,q_{k}\right)  $ be the resulting path
tuple\footnote{Thus,%
\[
q_{i}=\left(  \text{head of }p_{i}\right)  \cup\left(  \text{tail of }%
p_{j}\right)  \ \ \ \ \ \ \ \ \ \ \text{and}\ \ \ \ \ \ \ \ \ \ q_{j}=\left(
\text{head of }p_{j}\right)  \cup\left(  \text{tail of }p_{i}\right)
\]
(where \textquotedblleft head\textquotedblright\ means \textquotedblleft part
until $v$\textquotedblright, and \textquotedblleft tail\textquotedblright%
\ means \textquotedblleft part after $v$\textquotedblright). Furthermore, we
have $q_{m}=p_{m}$ for any $m\in\left[  k\right]  \setminus\left\{
i,j\right\}  $, since we have left all paths other than $p_{i}$ and $p_{j}$
unchanged.}.

\item We set $\sigma^{\prime}:=\sigma t_{i,j}$, where $t_{i,j}$ is the
transposition in $S_{k}$ that swaps $i$ with $j$. Thus, $\left(  q_{1}%
,q_{2},\ldots,q_{k}\right)  $ is a path tuple from $\mathbf{A}$ to
$\sigma^{\prime}\left(  \mathbf{B}\right)  $ (because exchanging the tails of
the paths $p_{i}$ and $p_{j}$ has switched their ending points $B_{\sigma
\left(  i\right)  }$ and $B_{\sigma\left(  j\right)  }$ to $B_{\sigma\left(
j\right)  }=B_{\sigma^{\prime}\left(  i\right)  }$ and $B_{\sigma\left(
i\right)  }=B_{\sigma^{\prime}\left(  j\right)  }$, respectively).

\item Finally, we set%
\[
f\left(  \sigma,\left(  p_{1},p_{2},\ldots,p_{k}\right)  \right)  :=\left(
\sigma^{\prime},\left(  q_{1},q_{2},\ldots,q_{k}\right)  \right)  .
\]

\end{itemize}

This defines a map $f:\mathcal{X}\rightarrow\mathcal{X}$ (again, it is not
hard to see that it is well-defined). Here is an example: If $k=5$ and if
$\left(  p_{1},p_{2},\ldots,p_{k}\right)  $ is%
\[%
\begin{tikzpicture}
\draw[densely dotted] (-0.2,-0.2) grid (6.2,7.2);
\node
[circle,fill=white,draw=black,text=black,inner sep=1pt] (A1) at (1,0) {$A_1$};
\node
[circle,fill=white,draw=black,text=black,inner sep=1pt] (A2) at (2,1) {$A_2$};
\node
[circle,fill=white,draw=black,text=black,inner sep=1pt] (A3) at (1,2) {$A_3$};
\node
[circle,fill=white,draw=black,text=black,inner sep=1pt] (A4) at (0,3) {$A_4$};
\node
[circle,fill=white,draw=black,text=black,inner sep=1pt] (A5) at (0,4) {$A_5$};
\node
[circle,fill=white,draw=black,text=black,inner sep=1pt] (B1) at (3,0) {$B'_1$}%
;
\node
[circle,fill=white,draw=black,text=black,inner sep=1pt] (B2) at (6,6) {$B'_2$}%
;
\node
[circle,fill=white,draw=black,text=black,inner sep=1pt] (B3) at (5,7) {$B'_3$}%
;
\node
[circle,fill=white,draw=black,text=black,inner sep=1pt] (B4) at (4,3) {$B'_4$}%
;
\node
[circle,fill=white,draw=black,text=black,inner sep=1pt] (B5) at (6,7) {$B'_5$}%
;
\begin{scope}[thick,>=stealth,darkred]
\draw(A1) edge[->] (2,0);
\draw(2,0) node[anchor=north] {$p_1$};
\draw(2,0) edge[->] (B1);
\end{scope}
\begin{scope}[thick,>=stealth,dbluecolor]
\draw(A2) edge[->] (2,2);
\draw(2,1.9) node[anchor=west] {$p_2$};
\draw(2,2) edge[->] (2,3);
\draw(2,3) edge[->] (2,4);
\draw(2,4) edge[->] (3,4);
\draw(3,4) edge[->] (4,4);
\draw(4,4) edge[->] (4,5);
\draw(4,4.5) node[anchor=west] {$p_2$};
\draw(4,5) edge[->] (5,5);
\draw(5,5) edge[->] (5,6);
\draw(5,6) edge[->] (B2);
\draw(5.3,6) node[anchor=north] {$p_2$};
\end{scope}
\begin{scope}[thick,>=stealth,dgreencolor]
\draw(A3) edge[->] (1,3);
\draw(1,2.5) node[anchor=west] {$p_3$};
\draw(1,3) edge[->] (1,4);
\draw(1,4) node[anchor=west] {$p_3$};
\draw(1,4) edge[->] (1,5);
\draw(1,5) edge[->] (2,5);
\draw(2,5) edge[->] (3,5);
\draw(3,5) edge[->] (3,6);
\draw(3,5.5) node[anchor=south east] {$p_3$};
\draw(3,6) edge[->] (4,6);
\draw(4,6) edge[->] (5,6);
\draw(5,6) edge[->] (B3);
\draw(5,6) node[anchor=south east] {$p_3$};
\end{scope}
\begin{scope}[thick,>=stealth,red!90!yellow]
\draw(A4) edge[->] (1,3);
\draw(3,3) node[anchor=north] {$p_4$};
\draw(1,3) edge[->] (2,3);
\draw(2,3) edge[->] (3,3);
\draw(3,3) edge[->] (B4);
\end{scope}
\begin{scope}[thick,>=stealth,cyan!80!black]
\draw(A5) edge[->] (0,5);
\draw(0,4.5) node[anchor=east] {$p_5$};
\draw(0,5) edge[->] (1,5);
\draw(1,5) edge[->] (1,6);
\draw(1,6) edge[->] (1,7);
\draw(1,7) edge[->] (2,7);
\draw(1,6.5) node[anchor=east] {$p_5$};
\draw(2,7) edge[->] (3,7);
\draw(3,7) edge[->] (4,7);
\draw(4,7) edge[->] (B3);
\draw(B3) edge[->] (B5);
\end{scope}
\end{tikzpicture}%
\]
(where $B_{m}^{\prime}$ is shorthand for $B_{\sigma\left(  m\right)  }$), then
$v$ is the point where paths $p_{2}$ and $p_{4}$ intersect, and we have $i=2$
and $j=4$, and therefore $\left(  q_{1},q_{2},\ldots,q_{k}\right)  $ is%
\[%
\begin{tikzpicture}
\draw[densely dotted] (-0.2,-0.2) grid (6.2,7.2);
\node
[circle,fill=white,draw=black,text=black,inner sep=1pt] (A1) at (1,0) {$A_1$};
\node
[circle,fill=white,draw=black,text=black,inner sep=1pt] (A2) at (2,1) {$A_2$};
\node
[circle,fill=white,draw=black,text=black,inner sep=1pt] (A3) at (1,2) {$A_3$};
\node
[circle,fill=white,draw=black,text=black,inner sep=1pt] (A4) at (0,3) {$A_4$};
\node
[circle,fill=white,draw=black,text=black,inner sep=1pt] (A5) at (0,4) {$A_5$};
\node
[circle,fill=white,draw=black,text=black,inner sep=1pt] (B1) at (3,0) {$B'_1$}%
;
\node
[circle,fill=white,draw=black,text=black,inner sep=1pt] (B2) at (6,6) {$B'_2$}%
;
\node
[circle,fill=white,draw=black,text=black,inner sep=1pt] (B3) at (5,7) {$B'_3$}%
;
\node
[circle,fill=white,draw=black,text=black,inner sep=1pt] (B4) at (4,3) {$B'_4$}%
;
\node
[circle,fill=white,draw=black,text=black,inner sep=1pt] (B5) at (6,7) {$B'_5$}%
;
\begin{scope}[thick,>=stealth,darkred]
\draw(A1) edge[->] (2,0);
\draw(2,0) node[anchor=north] {$q_1$};
\draw(2,0) edge[->] (B1);
\end{scope}
\begin{scope}[thick,>=stealth,dbluecolor]
\draw(A2) edge[->] (2,2);
\draw(2,1.9) node[anchor=west] {$q_2$};
\draw(2,2) edge[->] (2,3);
\draw(2,3) edge[ultra thick, ->] (3,3);
\draw(3,3) node[anchor=north] {$q_2$};
\draw(3,3) edge[ultra thick, ->] (B4);
\end{scope}
\begin{scope}[thick,>=stealth,dgreencolor]
\draw(A3) edge[->] (1,3);
\draw(1,2.5) node[anchor=west] {$q_3$};
\draw(1,3) edge[->] (1,4);
\draw(1,4) node[anchor=west] {$q_3$};
\draw(1,4) edge[->] (1,5);
\draw(1,5) edge[->] (2,5);
\draw(2,5) edge[->] (3,5);
\draw(3,5) edge[->] (3,6);
\draw(3,5.5) node[anchor=south east] {$q_3$};
\draw(3,6) edge[->] (4,6);
\draw(4,6) edge[->] (5,6);
\draw(5,6) edge[->] (B3);
\draw(5,6) node[anchor=south east] {$q_3$};
\end{scope}
\begin{scope}[thick,>=stealth,red!90!yellow]
\draw(A4) edge[->] (1,3);
\draw(1,3) edge[->] (2,3);
\draw(2,3) edge[ultra thick, ->] (2,4);
\draw(2,4) edge[ultra thick, ->] (3,4);
\draw(3,4) edge[ultra thick, ->] (4,4);
\draw(4,4) edge[ultra thick, ->] (4,5);
\draw(4,4.5) node[anchor=west] {$q_4$};
\draw(4,5) edge[ultra thick, ->] (5,5);
\draw(5,5) edge[ultra thick, ->] (5,6);
\draw(5,6) edge[ultra thick, ->] (B2);
\draw(5.3,6) node[anchor=north] {$q_4$};
\end{scope}
\begin{scope}[thick,>=stealth,cyan!80!black]
\draw(A5) edge[->] (0,5);
\draw(0,4.5) node[anchor=east] {$q_5$};
\draw(0,5) edge[->] (1,5);
\draw(1,5) edge[->] (1,6);
\draw(1,6) edge[->] (1,7);
\draw(1,7) edge[->] (2,7);
\draw(1,6.5) node[anchor=east] {$q_5$};
\draw(2,7) edge[->] (3,7);
\draw(3,7) edge[->] (4,7);
\draw(4,7) edge[->] (B3);
\draw(B3) edge[->] (B5);
\end{scope}
\node
[circle,fill=green,draw=black,text=black,inner sep=1pt] (v) at (2,3) {$v$};
\end{tikzpicture}%
\]
(where we again have drawn the exchanged tails extra-thick).

Convince yourself that the map $f$ defined above really is a sign-reversing
involution from $\mathcal{X}$ to $\mathcal{X}$. (This means showing that
applying $f$ twice in succession to any given $\left(  \sigma,\mathbf{p}%
\right)  \in\mathcal{X}$ returns $\left(  \sigma,\mathbf{p}\right)  $, and
that we have $\operatorname*{sign}\left(  f\left(  \sigma,\mathbf{p}\right)
\right)  =-\operatorname*{sign}\left(  \sigma,\mathbf{p}\right)  $ for any
$\left(  \sigma,\mathbf{p}\right)  \in\mathcal{X}$. The proof of the first
claim is straightforward but somewhat lengthy when written out; we have
banished it into the Appendix (Section \ref{sec.details.det.comb}). The proof
of the second claim, of course, relies on parts \textbf{(b)} and \textbf{(d)}
of Proposition \ref{prop.perm.sign.props}.)

Thus, we have defined a sign-reversing involution $f:\mathcal{X}%
\rightarrow\mathcal{X}$. This involution $f$ has no fixed points (since it is
sign-reversing). It is now easy to complete the proof: Lemma
\ref{lem.sign.cancel2} yields%
\[
\sum_{\left(  \sigma,\mathbf{p}\right)  \in\mathcal{A}}\operatorname*{sign}%
\left(  \sigma,\mathbf{p}\right)  =\sum_{\left(  \sigma,\mathbf{p}\right)
\in\mathcal{A}\setminus\mathcal{X}}\operatorname*{sign}\left(  \sigma
,\mathbf{p}\right)  .
\]
In view of our definition of $\operatorname*{sign}\left(  \sigma
,\mathbf{p}\right)  $, this rewrites as%
\begin{equation}
\sum_{\left(  \sigma,\mathbf{p}\right)  \in\mathcal{A}}\left(  -1\right)
^{\sigma}=\sum_{\left(  \sigma,\mathbf{p}\right)  \in\mathcal{A}%
\setminus\mathcal{X}}\left(  -1\right)  ^{\sigma}.
\label{pf.prop.lgv.kpaths.count.at}%
\end{equation}
The left hand side of this equality is%
\begin{align*}
\sum_{\left(  \sigma,\mathbf{p}\right)  \in\mathcal{A}}\left(  -1\right)
^{\sigma}  &  =\sum_{\sigma\in S_{k}}\left(  -1\right)  ^{\sigma
}\underbrace{\left(  \text{\# of path tuples from }\mathbf{A}\text{ to }%
\sigma\left(  \mathbf{B}\right)  \right)  }_{\substack{=\prod_{i=1}^{k}\left(
\text{\# of paths from }A_{i}\text{ to }B_{\sigma\left(  i\right)  }\right)
\\\text{(by the product rule, since a path tuple from }\mathbf{A}\text{ to
}\sigma\left(  \mathbf{B}\right)  \text{ is just}\\\text{a tuple }\left(
p_{1},p_{2},\ldots,p_{k}\right)  \text{, where each }p_{i}\text{ is a path
from }A_{i}\text{ to }B_{\sigma\left(  i\right)  }\text{)}}}\\
&  \ \ \ \ \ \ \ \ \ \ \ \ \ \ \ \ \ \ \ \ \left(  \text{by the definition of
}\mathcal{A}\right) \\
&  =\sum_{\sigma\in S_{k}}\left(  -1\right)  ^{\sigma}\prod_{i=1}^{k}\left(
\text{\# of paths from }A_{i}\text{ to }B_{\sigma\left(  i\right)  }\right) \\
&  =\det\left(  \left(  \text{\# of paths from }A_{i}\text{ to }B_{j}\right)
_{1\leq i\leq k,\ 1\leq j\leq k}\right)
\end{align*}
(by the definition of the determinant), whereas the right hand side is%
\begin{align*}
\sum_{\left(  \sigma,\mathbf{p}\right)  \in\mathcal{A}\setminus\mathcal{X}%
}\left(  -1\right)  ^{\sigma}  &  =\sum_{\left(  \sigma,\mathbf{p}\right)
\in\mathcal{A}\text{ is a nipat}}\left(  -1\right)  ^{\sigma}%
\ \ \ \ \ \ \ \ \ \ \left(
\begin{array}
[c]{c}%
\text{since }\mathcal{X}=\left\{  \text{ipats in }\mathcal{A}\right\} \\
\text{entails }\mathcal{A}\setminus\mathcal{X}=\left\{  \text{nipats in
}\mathcal{A}\right\}
\end{array}
\right) \\
&  =\sum_{\sigma\in S_{k}}\left(  -1\right)  ^{\sigma}\left(  \text{\# of
nipats from }\mathbf{A}\text{ to }\sigma\left(  \mathbf{B}\right)  \right)
\end{align*}
(by the definition of $\mathcal{A}$). Thus, (\ref{pf.prop.lgv.kpaths.count.at}%
) rewrites as
\begin{align*}
&  \det\left(  \left(  \text{\# of paths from }A_{i}\text{ to }B_{j}\right)
_{1\leq i\leq k,\ 1\leq j\leq k}\right) \\
&  =\sum_{\sigma\in S_{k}}\left(  -1\right)  ^{\sigma}\left(  \text{\# of
nipats from }\mathbf{A}\text{ to }\sigma\left(  \mathbf{B}\right)  \right)  .
\end{align*}
This proves Proposition \ref{prop.lgv.kpaths.count}.
\end{proof}

\subsubsection{The weighted version}

So far we have just been counting paths; but we can easily introduce weights
to obtain a more general result:

\begin{theorem}
[LGV lemma, lattice weight version]\label{thm.lgv.kpaths.wt}Let $K$ be a
commutative ring.

For each arc $a$ of the digraph $\mathbb{Z}^{2}$, let $w\left(  a\right)  $ be
an element of $K$. We call this element $w\left(  a\right)  $ the
\emph{weight} of $a$.

For each path $p$ of $\mathbb{Z}^{2}$, define the \emph{weight} $w\left(
p\right)  $ of $p$ by%
\[
w\left(  p\right)  :=\prod_{a\text{ is an arc of }p}w\left(  a\right)  .
\]

For each path tuple $\mathbf{p}=\left(  p_{1},p_{2},\ldots,p_{k}\right)  $,
define the \emph{weight} $w\left(  \mathbf{p}\right)  $ of $\mathbf{p}$ by%
\[
w\left(  \mathbf{p}\right)  :=w\left(  p_{1}\right)  w\left(  p_{2}\right)
\cdots w\left(  p_{k}\right)  .
\]

Let $k\in\mathbb{N}$. Let $\mathbf{A}=\left(  A_{1},A_{2},\ldots,A_{k}\right)
$ and $\mathbf{B}=\left(  B_{1},B_{2},\ldots,B_{k}\right)  $ be two
$k$-vertices. Then,%
\[
\det\left(  \left(  \sum_{p:A_{i}\rightarrow B_{j}}w\left(  p\right)  \right)
_{1\leq i\leq k,\ 1\leq j\leq k}\right)  =\sum_{\sigma\in S_{k}}\left(
-1\right)  ^{\sigma}\sum_{\substack{\mathbf{p}\text{ is a nipat}\\\text{from
}\mathbf{A}\text{ to }\sigma\left(  \mathbf{B}\right)  }}w\left(
\mathbf{p}\right)  .
\]
Here, \textquotedblleft$p:A_{i}\rightarrow B_{j}$\textquotedblright\ means
\textquotedblleft$p$ is a path from $A_{i}$ to $B_{j}$\textquotedblright.
\end{theorem}

Note that an empty path (i.e., a path with no arcs) has weight $1$.

Clearly, Proposition \ref{prop.lgv.kpaths.count} is the particular case of
Theorem \ref{thm.lgv.kpaths.wt} when $K=\mathbb{Z}$ and $w\left(  a\right)
=1$ for all arcs $a$ (because in this case, all the weights $w\left(
p\right)  $ and $w\left(  \mathbf{p}\right)  $ of paths and path tuples equal
$1$, and therefore the sums over paths or nipats become the \#s of paths or nipats).

\begin{proof}
[Proof of Theorem \ref{thm.lgv.kpaths.wt}.]The same argument as for
Proposition \ref{prop.lgv.kpaths.count} can be used here; just replace
$\operatorname*{sign}\left(  \sigma,\mathbf{p}\right)  :=\left(  -1\right)
^{\sigma}$ by $\operatorname*{sign}\left(  \sigma,\mathbf{p}\right)  =\left(
-1\right)  ^{\sigma}\cdot w\left(  \mathbf{p}\right)  $. (The only new
observation required is that when we exchange the tails of two paths in our
path tuple, the weight of the path tuple does not change. This is rather
clear: The weight of a path tuple is the product of the weights of all arcs in
all paths of the tuple\footnote{An arc will appear multiple times in the
product if it appears in multiple paths.}. When we exchange the tails of two
paths, some arcs get moved from one path to the other, but the total product
stays unchanged.)
\end{proof}

\subsubsection{Generalization to acyclic digraphs}

We can generalize Theorem \ref{thm.lgv.kpaths.wt} further. Indeed, we have
barely used anything specific to $\mathbb{Z}^{2}$ in our proofs; all we used
is that $\mathbb{Z}^{2}$ is a path-finite acyclic digraph. Thus, Theorem
\ref{thm.lgv.kpaths.wt} remains true if we replace $\mathbb{Z}^{2}$ by an
arbitrary such digraph. We thus obtain the following more general result:

\begin{theorem}
[LGV lemma, digraph weight version]\label{thm.lgv.kpaths.wt-dg}Let $K$ be a
commutative ring.

Let $D$ be a path-finite (but possibly infinite) acyclic digraph. We extend
Definition \ref{def.lgv.path-tups} to $D$ instead of $\mathbb{Z}^{2}$ (with
the obvious changes: \textquotedblleft lattice points\textquotedblright%
\ becomes \textquotedblleft vertices of $D$\textquotedblright).

For each arc $a$ of the digraph $D$, let $w\left(  a\right)  $ be an element
of $K$. We call this element $w\left(  a\right)  $ the \emph{weight} of $a$.

For each path $p$ of $D$, define the \emph{weight} $w\left(  p\right)  $ of
$p$ by%
\[
w\left(  p\right)  :=\prod_{a\text{ is an arc of }p}w\left(  a\right)  .
\]

For each path tuple $\mathbf{p}=\left(  p_{1},p_{2},\ldots,p_{k}\right)  $,
define the \emph{weight} $w\left(  \mathbf{p}\right)  $ of $\mathbf{p}$ by%
\[
w\left(  \mathbf{p}\right)  :=w\left(  p_{1}\right)  w\left(  p_{2}\right)
\cdots w\left(  p_{k}\right)  .
\]

Let $k\in\mathbb{N}$. Let $\mathbf{A}=\left(  A_{1},A_{2},\ldots,A_{k}\right)
$ and $\mathbf{B}=\left(  B_{1},B_{2},\ldots,B_{k}\right)  $ be two
$k$-vertices (i.e., two $k$-tuples of vertices of $D$). Then,%
\[
\det\left(  \left(  \sum_{p:A_{i}\rightarrow B_{j}}w\left(  p\right)  \right)
_{1\leq i\leq k,\ 1\leq j\leq k}\right)  =\sum_{\sigma\in S_{k}}\left(
-1\right)  ^{\sigma}\sum_{\substack{\mathbf{p}\text{ is a nipat}\\\text{from
}\mathbf{A}\text{ to }\sigma\left(  \mathbf{B}\right)  }}w\left(
\mathbf{p}\right)  .
\]
Here, \textquotedblleft$p:A_{i}\rightarrow B_{j}$\textquotedblright\ means
\textquotedblleft$p$ is a path from $A_{i}$ to $B_{j}$\textquotedblright.
\end{theorem}

\begin{proof}
Completely analogous to the proof of Theorem \ref{thm.lgv.kpaths.wt}.
\end{proof}

\subsubsection{The nonpermutable case}

One nice thing about the digraph $\mathbb{Z}^{2}$, however, is that in many
cases, the sum%
\[
\sum_{\sigma\in S_{k}}\left(  -1\right)  ^{\sigma}\sum_{\substack{\mathbf{p}%
\text{ is a nipat}\\\text{from }\mathbf{A}\text{ to }\sigma\left(
\mathbf{B}\right)  }}w\left(  \mathbf{p}\right)
\]
has only one nonzero addend. We have already seen this happen often in the
$k=2$ case (thanks to Proposition \ref{prop.lgv.jordan-2}). Here is the
analogous statement for general $k$:

\begin{corollary}
[LGV lemma, nonpermutable lattice weight version]\label{cor.lgv.kpaths.wt-np}%
Consider the setting of Theorem \ref{thm.lgv.kpaths.wt}, but additionally
assume that%
\begin{align}
\operatorname*{x}\left(  A_{1}\right)   &  \geq\operatorname*{x}\left(
A_{2}\right)  \geq\cdots\geq\operatorname*{x}\left(  A_{k}\right)
;\label{eq.cor.lgv.kpaths.wt-np.xA}\\
\operatorname*{y}\left(  A_{1}\right)   &  \leq\operatorname*{y}\left(
A_{2}\right)  \leq\cdots\leq\operatorname*{y}\left(  A_{k}\right)
;\label{eq.cor.lgv.kpaths.wt-np.yA}\\
\operatorname*{x}\left(  B_{1}\right)   &  \geq\operatorname*{x}\left(
B_{2}\right)  \geq\cdots\geq\operatorname*{x}\left(  B_{k}\right)
;\label{eq.cor.lgv.kpaths.wt-np.xB}\\
\operatorname*{y}\left(  B_{1}\right)   &  \leq\operatorname*{y}\left(
B_{2}\right)  \leq\cdots\leq\operatorname*{y}\left(  B_{k}\right)  .
\label{eq.cor.lgv.kpaths.wt-np.yB}%
\end{align}
Here, $\operatorname*{x}\left(  P\right)  $ and $\operatorname*{y}\left(
P\right)  $ denote the two coordinates of any point $P\in\mathbb{Z}^{2}$.

Then, there are no nipats from $\mathbf{A}$ to $\sigma\left(  \mathbf{B}%
\right)  $ when $\sigma\in S_{k}$ is not the identity permutation
$\operatorname*{id}\in S_{k}$. Therefore, the claim of Theorem
\ref{thm.lgv.kpaths.wt} simplifies to
\begin{equation}
\det\left(  \left(  \sum_{p:A_{i}\rightarrow B_{j}}w\left(  p\right)  \right)
_{1\leq i\leq k,\ 1\leq j\leq k}\right)  =\sum_{\substack{\mathbf{p}\text{ is
a nipat}\\\text{from }\mathbf{A}\text{ to }\mathbf{B}}}w\left(  \mathbf{p}%
\right)  .\ \ \ \ \ \ \ \ \ \ \label{eq.cor.lgv.kpaths.wt-np.claim}%
\end{equation}

\end{corollary}

\begin{proof}
[Proof of Corollary \ref{cor.lgv.kpaths.wt-np} (sketched).]This is easy using
Proposition \ref{prop.lgv.jordan-2}. Here are the details: \medskip

\begin{fineprint}
Let $\sigma\in S_{k}$ be a permutation that is not the identity permutation
$\operatorname*{id}\in S_{k}$. Then, we don't have $\sigma\left(  1\right)
\leq\sigma\left(  2\right)  \leq\cdots\leq\sigma\left(  k\right)  $ (since
$\sigma$ is not $\operatorname*{id}$). In other words, there exists some
$i\in\left[  k-1\right]  $ such that $\sigma\left(  i\right)  >\sigma\left(
i+1\right)  $. Consider this $i$.

Now, let $\mathbf{p}$ be a nipat from $\mathbf{A}$ to $\sigma\left(
\mathbf{B}\right)  $. Write $\mathbf{p}$ in the form $\mathbf{p}=\left(
p_{1},p_{2},\ldots,p_{k}\right)  $. Thus, $p_{i}$ is a path from $A_{i}$ to
$B_{\sigma\left(  i\right)  }$, whereas $p_{i+1}$ is a path from $A_{i+1}$ to
$B_{\sigma\left(  i+1\right)  }$. Moreover, $p_{i}$ and $p_{i+1}$ have no
vertex in common (since $\mathbf{p}$ is a nipat).

The sequence $\left(  \operatorname*{x}\left(  B_{1}\right)
,\operatorname*{x}\left(  B_{2}\right)  ,\ldots,\operatorname*{x}\left(
B_{k}\right)  \right)  $ is weakly decreasing (by
(\ref{eq.cor.lgv.kpaths.wt-np.xB})). In other words, if $m$ and $n$ are two
elements of $\left[  k\right]  $ satisfying $m>n$, then $\operatorname*{x}%
\left(  B_{m}\right)  \leq\operatorname*{x}\left(  B_{n}\right)  $. Applying
this to $m=\sigma\left(  i\right)  $ and $n=\sigma\left(  i+1\right)  $, we
obtain $\operatorname*{x}\left(  B_{\sigma\left(  i\right)  }\right)
\leq\operatorname*{x}\left(  B_{\sigma\left(  i+1\right)  }\right)  $ (since
$\sigma\left(  i\right)  >\sigma\left(  i+1\right)  $). Likewise, using
(\ref{eq.cor.lgv.kpaths.wt-np.yB}), we can obtain $\operatorname*{y}\left(
B_{\sigma\left(  i\right)  }\right)  \geq\operatorname*{y}\left(
B_{\sigma\left(  i+1\right)  }\right)  $. However,
(\ref{eq.cor.lgv.kpaths.wt-np.xA}) shows that $\operatorname*{x}\left(
A_{i}\right)  \geq\operatorname*{x}\left(  A_{i+1}\right)  $. In other words,
$\operatorname*{x}\left(  A_{i+1}\right)  \leq\operatorname*{x}\left(
A_{i}\right)  $. Furthermore, (\ref{eq.cor.lgv.kpaths.wt-np.yA}) shows that
$\operatorname*{y}\left(  A_{i}\right)  \leq\operatorname*{y}\left(
A_{i+1}\right)  $. In other words, $\operatorname*{y}\left(  A_{i+1}\right)
\geq\operatorname*{y}\left(  A_{i}\right)  $.

Hence, Proposition \ref{prop.lgv.jordan-2} (applied to $A=A_{i}$,
$B=B_{\sigma\left(  i+1\right)  }$, $A^{\prime}=A_{i+1}$, $B^{\prime
}=B_{\sigma\left(  i\right)  }$, $p=p_{i}$ and $p^{\prime}=p_{i+1}$) yields
that $p_{i}$ and $p_{i+1}$ have a vertex in common. This contradicts the fact
that $p_{i}$ and $p_{i+1}$ have no vertex in common.

Forget that we fixed $\mathbf{p}$. We thus have found a contradiction for each
nipat $\mathbf{p}$ from $\mathbf{A}$ to $\sigma\left(  \mathbf{B}\right)  $.
Hence, there are no nipats from $\mathbf{A}$ to $\sigma\left(  \mathbf{B}%
\right)  $.

Forget that we fixed $\sigma$. We thus have proved that there are no nipats
from $\mathbf{A}$ to $\sigma\left(  \mathbf{B}\right)  $ when $\sigma\in
S_{k}$ is not the identity permutation $\operatorname*{id}\in S_{k}$. Hence,
if $\sigma\in S_{k}$ is not the identity permutation $\operatorname*{id}\in
S_{k}$, then%
\begin{equation}
\sum_{\substack{\mathbf{p}\text{ is a nipat}\\\text{from }\mathbf{A}\text{ to
}\sigma\left(  \mathbf{B}\right)  }}w\left(  \mathbf{p}\right)  =\left(
\text{empty sum}\right)  =0. \label{pf.cor.lgv.kpaths.wt-np.at}%
\end{equation}

Now, Theorem \ref{thm.lgv.kpaths.wt} yields%
\begin{align*}
&  \det\left(  \left(  \sum_{p:A_{i}\rightarrow B_{j}}w\left(  p\right)
\right)  _{1\leq i\leq k,\ 1\leq j\leq k}\right) \\
&  =\sum_{\sigma\in S_{k}}\left(  -1\right)  ^{\sigma}\sum
_{\substack{\mathbf{p}\text{ is a nipat}\\\text{from }\mathbf{A}\text{ to
}\sigma\left(  \mathbf{B}\right)  }}w\left(  \mathbf{p}\right) \\
&  =\underbrace{\left(  -1\right)  ^{\operatorname*{id}}}_{=1}\sum
_{\substack{\mathbf{p}\text{ is a nipat}\\\text{from }\mathbf{A}\text{ to
}\operatorname*{id}\left(  \mathbf{B}\right)  }}w\left(  \mathbf{p}\right)
+\sum_{\substack{\sigma\in S_{k};\\\sigma\neq\operatorname*{id}}}\left(
-1\right)  ^{\sigma}\underbrace{\sum_{\substack{\mathbf{p}\text{ is a
nipat}\\\text{from }\mathbf{A}\text{ to }\sigma\left(  \mathbf{B}\right)
}}w\left(  \mathbf{p}\right)  }_{\substack{=0\\\text{(by
(\ref{pf.cor.lgv.kpaths.wt-np.at}))}}}\\
&  \ \ \ \ \ \ \ \ \ \ \ \ \ \ \ \ \ \ \ \ \left(  \text{here, we have split
off the addend for }\sigma=\operatorname*{id}\text{ from the sum}\right) \\
&  =\sum_{\substack{\mathbf{p}\text{ is a nipat}\\\text{from }\mathbf{A}\text{
to }\operatorname*{id}\left(  \mathbf{B}\right)  }}w\left(  \mathbf{p}\right)
+\underbrace{\sum_{\substack{\sigma\in S_{k};\\\sigma\neq\operatorname*{id}%
}}\left(  -1\right)  ^{\sigma}0}_{=0}=\sum_{\substack{\mathbf{p}\text{ is a
nipat}\\\text{from }\mathbf{A}\text{ to }\operatorname*{id}\left(
\mathbf{B}\right)  }}w\left(  \mathbf{p}\right)  =\sum_{\substack{\mathbf{p}%
\text{ is a nipat}\\\text{from }\mathbf{A}\text{ to }\mathbf{B}}}w\left(
\mathbf{p}\right)
\end{align*}
(since $\operatorname*{id}\left(  \mathbf{B}\right)  =\mathbf{B}$). The proof
of Corollary \ref{cor.lgv.kpaths.wt-np} is now complete.
\end{fineprint}
\end{proof}

\begin{corollary}
\label{cor.lgv.binom-det-nonneg}Let $k\in\mathbb{N}$. Let $a_{1},a_{2}%
,\ldots,a_{k}$ and $b_{1},b_{2},\ldots,b_{k}$ be nonnegative integers such
that
\[
a_{1}\geq a_{2}\geq\cdots\geq a_{k}\ \ \ \ \ \ \ \ \ \ \text{and}%
\ \ \ \ \ \ \ \ \ \ b_{1}\geq b_{2}\geq\cdots\geq b_{k}.
\]
Then,%
\[
\det\left(  \left(  \dbinom{a_{i}}{b_{j}}\right)  _{1\leq i\leq k,\ 1\leq
j\leq k}\right)  \geq0.
\]

\end{corollary}

For example, if $a_{1}\geq a_{2}\geq a_{3}\geq0$ and $b_{1}\geq b_{2}\geq
b_{3}\geq0$, then%
\[
\det\left(
\begin{array}
[c]{ccc}%
\dbinom{a_{1}}{b_{1}} & \dbinom{a_{1}}{b_{2}} & \dbinom{a_{1}}{b_{3}}\\
\dbinom{a_{2}}{b_{1}} & \dbinom{a_{2}}{b_{2}} & \dbinom{a_{2}}{b_{3}}\\
\dbinom{a_{3}}{b_{1}} & \dbinom{a_{3}}{b_{2}} & \dbinom{a_{3}}{b_{3}}%
\end{array}
\right)  \geq0.
\]

\begin{proof}
[Proof of Corollary \ref{cor.lgv.binom-det-nonneg} (sketched).]Set
$K=\mathbb{Z}$, and set $w\left(  a\right)  :=1$ for each arc $a$ of
$\mathbb{Z}^{2}$. Define the lattice points%
\[
A_{i}:=\left(  0,-a_{i}\right)  \ \ \ \ \ \ \ \ \ \ \text{and}%
\ \ \ \ \ \ \ \ \ \ B_{i}:=\left(  b_{i},-b_{i}\right)
\]
for all $i\in\left[  k\right]  $. These lattice points satisfy the assumptions
of Corollary \ref{cor.lgv.kpaths.wt-np}. Hence,
(\ref{eq.cor.lgv.kpaths.wt-np.claim}) entails%
\begin{equation}
\det\left(  \left(  \sum_{p:A_{i}\rightarrow B_{j}}w\left(  p\right)  \right)
_{1\leq i\leq k,\ 1\leq j\leq k}\right)  =\sum_{\substack{\mathbf{p}\text{ is
a nipat}\\\text{from }\mathbf{A}\text{ to }\mathbf{B}}}w\left(  \mathbf{p}%
\right)  .\nonumber
\end{equation}
Since all the weights $w\left(  p\right)  $ and $w\left(  \mathbf{p}\right)  $
are $1$ in our situation, we can rewrite this as%
\begin{equation}
\det\left(  \left(  \text{\# of paths from }A_{i}\text{ to }B_{j}\right)
_{1\leq i\leq k,\ 1\leq j\leq k}\right)  =\left(  \text{\# of nipats from
}\mathbf{A}\text{ to }\mathbf{B}\right)  .\nonumber
\end{equation}

Using Proposition \ref{prop.lgv.1-paths.ct}, we can easily see that the matrix
on the left hand side of this equality is $\left(  \dbinom{a_{i}}{b_{j}%
}\right)  _{1\leq i\leq k,\ 1\leq j\leq k}$. Thus, this equality rewrites as%
\begin{equation}
\det\left(  \left(  \dbinom{a_{i}}{b_{j}}\right)  _{1\leq i\leq k,\ 1\leq
j\leq k}\right)  =\left(  \text{\# of nipats from }\mathbf{A}\text{ to
}\mathbf{B}\right)  .\nonumber
\end{equation}
Its left hand side is therefore $\geq0$ (since its right hand side is $\geq
0$). This proves Corollary \ref{cor.lgv.binom-det-nonneg}.
\end{proof}

\begin{corollary}
\label{cor.lgv.catalan-hankel-det-0}Let $k\in\mathbb{N}$. Recall the Catalan
numbers $c_{n}=\dfrac{1}{n+1}\dbinom{2n}{n}$ for all $n\in\mathbb{N}$. Then,%
\[
\det\left(  \left(  c_{i+j-2}\right)  _{1\leq i\leq k,\ 1\leq j\leq k}\right)
=\det\left(
\begin{array}
[c]{cccc}%
c_{0} & c_{1} & \cdots & c_{k-1}\\
c_{1} & c_{2} & \cdots & c_{k}\\
\vdots & \vdots & \ddots & \vdots\\
c_{k-1} & c_{k} & \cdots & c_{2k-2}%
\end{array}
\right)  =1.
\]

\end{corollary}

\begin{proof}
[Proof of Corollary \ref{cor.lgv.catalan-hankel-det-0} (sketched).]We will use
not the lattice $\mathbb{Z}^{2}$, but a different digraph. Namely, we use the
simple digraph with vertex set $\mathbb{Z}\times\mathbb{N}$ (that is, the
vertices are the lattice points that lie on the x-axis or above it) and arcs%
\[
\left(  i,j\right)  \rightarrow\left(  i+1,\ j+1\right)
\ \ \ \ \ \ \ \ \ \ \text{for all }\left(  i,j\right)  \in\mathbb{Z}%
\times\mathbb{N}%
\]
and%
\[
\left(  i,j\right)  \rightarrow\left(  i+1,\ j-1\right)
\ \ \ \ \ \ \ \ \ \ \text{for all }\left(  i,j\right)  \in\mathbb{Z}%
\times\mathbb{P},
\]
where $\mathbb{P}:=\left\{  1,2,3,\ldots\right\}  $. Here is a picture of a
small part of this digraph:%
\[%
\begin{tikzpicture}
\foreach\x/\xtext in {-1, 0, 1, 2, 3, 4, 5}
{
\foreach\y/\ytext in {0, 1, 2, 3, 4, 5, 6}
{
\ifnum\y<6
\draw(\x+0.12,\y+0.12) edge[->, thick, >=stealth, darkred] (\x+0.88,\y+0.88);
\fi\ifnum\y>0
\draw(\x+0.12,\y-0.12) edge[->, thick, >=stealth, dbluecolor] (\x
+0.88,\y-0.88);
\fi\ifnum\x>-1
\ifnum\y<6
\draw(\x,\y) circle[radius=0.1, black];
\fi\fi}
};
\foreach\x/\xtext in {0, 1, 2, 3, 4, 5}
\draw(\x cm,1pt) node[anchor=north west] {$\xtext$};
\foreach\y/\ytext in {0, 1, 2, 3, 4, 5}
\draw(1pt,\y cm) node[left=4pt] {$\ytext$};
\end{tikzpicture}%
\ \ .
\]
As we know, the Catalan number $c_{n}$ counts the paths from $\left(
0,0\right)  $ to $\left(  2n,0\right)  $ on this digraph (indeed, these are
just the Dyck paths\footnote{See Example 2 in Section \ref{sec.gf.exas} for
the definition of a Dyck path.}). Hence, $c_{n}$ also counts the paths from
$\left(  i,0\right)  $ to $\left(  2n+i,0\right)  $ whenever $i\in\mathbb{Z}$
(because these are just the Dyck paths shifted by $i$ in the x-direction). It
is easy to see that this digraph is acyclic and path-finite.

Now, define two $k$-vertices $\mathbf{A}=\left(  A_{1},A_{2},\ldots
,A_{k}\right)  $ and $\mathbf{B}=\left(  B_{1},B_{2},\ldots,B_{k}\right)  $ by
setting%
\[
A_{i}:=\left(  -2\left(  i-1\right)  ,0\right)  \ \ \ \ \ \ \ \ \ \ \text{and}%
\ \ \ \ \ \ \ \ \ \ B_{i}:=\left(  2\left(  i-1\right)  ,0\right)
\]
for all $i\in\left[  k\right]  $. It is not hard to show (see Exercise
\ref{exe.lgv.catalan-hankel-det-0} \textbf{(a)}) that there is only one nipat
from $\mathbf{A}$ to $\mathbf{B}$, which is shown in the case $k=4$ on the
following picture:\footnote{In this picture, we are drawing only
\textquotedblleft half\textquotedblright\ of the grid. Indeed, the vertices
$\left(  i,j\right)  $ of our digraph $\mathbb{Z}\times\mathbb{N}$ can be
classified into \emph{even vertices} (i.e., the ones for which $i+j$ is even)
and \emph{odd vertices} (i.e., the ones for which $i+j$ is odd). Any arc
either connects two even vertices or connects two odd vertices. Hence, a path
starting at an even vertex cannot contain any odd vertex (and vice versa).
Since all our vertices $A_{1},A_{2},\ldots,A_{k}$ and $B_{1},B_{2}%
,\ldots,B_{k}$ are even, we thus don't have to bother even drawing the odd
vertices (as they have no chance to appear in any paths between our vertices).
As a consequence, we are drawing only the grid lines containing the even
vertices.}%
\[%
\begin{tikzpicture}
\draw[densely dotted] (-6.5,5.5) -- (-5.5,6.5);
\draw[densely dotted] (-6.5,3.5) -- (-3.5,6.5);
\draw[densely dotted] (-6.5,1.5) -- (-1.5,6.5);
\draw[densely dotted] (-6,0) -- (0.5,6.5);
\draw[densely dotted] (-4,0) -- (2.5,6.5);
\draw[densely dotted] (-2,0) -- (4.5,6.5);
\draw[densely dotted] ( 0,0) -- (6.5,6.5);
\draw[densely dotted] ( 2,0) -- (6.5,4.5);
\draw[densely dotted] ( 4,0) -- (6.5,2.5);
\draw[densely dotted] ( 6,0) -- (6.5,0.5);
\draw[densely dotted] ( 6,0) -- (-0.5,6.5);
\draw[densely dotted] ( 4,0) -- (-2.5,6.5);
\draw[densely dotted] ( 2,0) -- (-4.5,6.5);
\draw[densely dotted] ( 0,0) -- (-6.5,6.5);
\draw[densely dotted] (-2,0) -- (-6.5,4.5);
\draw[densely dotted] (-4,0) -- (-6.5,2.5);
\draw[densely dotted] (-6,0) -- (-6.5,0.5);
\draw[densely dotted] (6.5,5.5) -- (5.5,6.5);
\draw[densely dotted] (6.5,3.5) -- (3.5,6.5);
\draw[densely dotted] (6.5,1.5) -- (1.5,6.5);
\node
[circle,fill=white,draw=black,text=black,inner sep=1pt] (A1) at (0,0) {$A_1$};
\node
[circle,fill=white,draw=black,text=black,inner sep=1pt] (A2) at (-2,0) {$A_2$}%
;
\node
[circle,fill=white,draw=black,text=black,inner sep=1pt] (A3) at (-4,0) {$A_3$}%
;
\node
[circle,fill=white,draw=black,text=black,inner sep=1pt] (A4) at (-6,0) {$A_4$}%
;
\node
[circle,fill=white,draw=black,text=black,inner sep=1pt] (B1) at (0,0) {$B_1$};
\node
[circle,fill=white,draw=black,text=black,inner sep=1pt] (B2) at (2,0) {$B_2$};
\node
[circle,fill=white,draw=black,text=black,inner sep=1pt] (B3) at (4,0) {$B_3$};
\node
[circle,fill=white,draw=black,text=black,inner sep=1pt] (B4) at (6,0) {$B_4$};
\begin{scope}[thick,>=stealth,dbluecolor]
\draw(A2) edge[->] (-1,1);
\draw(-1,1) edge[->] (0,2);
\draw(0,2) edge[->] (1,1);
\draw(1,1) edge[->] (B2);
\end{scope}
\begin{scope}[thick,>=stealth,dgreencolor]
\draw(A3) edge[->] (-3,1);
\draw(-3,1) edge[->] (-2,2);
\draw(-2,2) edge[->] (-1,3);
\draw(-1,3) edge[->] (0,4);
\draw(0,4) edge[->] (1,3);
\draw(1,3) edge[->] (2,2);
\draw(2,2) edge[->] (3,1);
\draw(3,1) edge[->] (B3);
\end{scope}
\begin{scope}[thick,>=stealth,darkred]
\draw(A4) edge[->] (-5,1);
\draw(-5,1) edge[->] (-4,2);
\draw(-4,2) edge[->] (-3,3);
\draw(-3,3) edge[->] (-2,4);
\draw(-2,4) edge[->] (-1,5);
\draw(-1,5) edge[->] (0,6);
\draw(0,6) edge[->] (1,5);
\draw(1,5) edge[->] (2,4);
\draw(2,4) edge[->] (3,3);
\draw(3,3) edge[->] (4,2);
\draw(4,2) edge[->] (5,1);
\draw(5,1) edge[->] (B4);
\end{scope}
\end{tikzpicture}%
\]
(the point $A_{1}$ coincides with $B_{1}$, and the path from $A_{1}$ to
$B_{1}$ is invisible, since it has no arcs). Moreover, it can be shown (see
Exercise \ref{exe.lgv.catalan-hankel-det-0} \textbf{(b)}) that there are no
nipats from $\mathbf{A}$ to $\sigma\left(  \mathbf{B}\right)  $ when
$\sigma\in S_{k}$ is not the identity permutation $\operatorname*{id}\in
S_{k}$. (This is analogous to Corollary \ref{cor.lgv.kpaths.wt-np}.) Hence, if
we set $K=\mathbb{Z}$ and $w\left(  a\right)  =1$ for each arc $a$ of our
digraph, then the relevant analogue of (\ref{eq.cor.lgv.kpaths.wt-np.claim})
entails%
\begin{equation}
\det\left(  \left(  \sum_{p:A_{i}\rightarrow B_{j}}w\left(  p\right)  \right)
_{1\leq i\leq k,\ 1\leq j\leq k}\right)  =\sum_{\substack{\mathbf{p}\text{ is
a nipat}\\\text{from }\mathbf{A}\text{ to }\mathbf{B}}}w\left(  \mathbf{p}%
\right) \nonumber
\end{equation}
(by the same reasoning as in the proof of Corollary \ref{cor.lgv.kpaths.wt-np}%
). The right hand side of this equality is $1$ (since there is only one nipat
from $\mathbf{A}$ to $\mathbf{B}$), while the matrix on the left hand side is
easily seen to be $\left(  c_{i+j-2}\right)  _{1\leq i\leq k,\ 1\leq j\leq k}$
(since the \# of paths from $A_{i}$ to $B_{j}$ is the Catalan number
$c_{i+j-2}$). This yields the claim of Corollary
\ref{cor.lgv.catalan-hankel-det-0}. The details are LTTR.
\end{proof}

The LGV lemma in all its variants is one major place in which combinatorial
questions reduce to the computation of determinants. Other such places are the
\emph{matrix-tree theorem} (see, e.g., \cite[\S 4]{Zeilbe}, \cite[\S 3.17]%
{Loehr-BC}, \cite[Theorems 9.8 and 10.4]{Stanle18}, \cite[\S 5.14]{23s}) and
the enumeration of perfect matchings or domino tilings (see, e.g.,
\cite{Stucky15}, \cite[\S 12.12--\S 12.13]{Loehr-BC}, \cite[\S 10.1]%
{Aigner07}). Soon, we will also encounter determinants in the study of
symmetric functions.

See also \cite{BruRys91} for a selection of other intersections between
combinatorics and linear algebra.

\begin{noncompile}
TODO: Cover the matrix-tree theorem (OH8).
\end{noncompile}

\newpage

\section{\label{chap.sf}Symmetric functions}

This final chapter is devoted to the theory of \emph{symmetric functions}.
Specifically, we will restrict ourselves to \emph{symmetric polynomials} (the
\textquotedblleft functions\textquotedblright\ part is a technical tweak that
makes the theory neater but we won't have time to introduce). Serious
treatments of the subject can be found in \cite{Wildon20}, \cite[Chapters
10--11]{Loehr-BC}, \cite{Egge19}, \cite{MenRem15}, \cite{Macdon95},
\cite[Chapter 8]{Aigner07}, \cite[Chapter 7]{Stanley-EC2}, \cite[Chapter
7]{Sagan19}, \cite[Chapter 4]{Sagan01}, \cite{Krishn86}, \cite[Chapters
14--19]{FoaHan04}, \cite{Savage22}, \cite{SmiTut24}, \cite{Anders19},
\cite[Chapter 2]{GriRei} and \cite{LLPT95}.

We begin with some oversimplified historical motivation.

Symmetric polynomials first(?) appeared in the study of roots of polynomials.
Consider a monic univariate polynomial%
\[
f=x^{n}+a_{1}x^{n-1}+a_{2}x^{n-2}+\cdots+a_{n}x^{0}\in\mathbb{C}\left[
x\right]
\]
(note the nonstandard labeling of coefficients). Let $r_{1},r_{2},\ldots
,r_{n}\in\mathbb{C}$ be the roots of this polynomial (listed with
multiplicities). Then, Fran\c{c}ois Vi\`{e}te (aka Franciscus Vieta) noticed
that%
\[
f=\left(  x-r_{1}\right)  \left(  x-r_{2}\right)  \cdots\left(  x-r_{n}%
\right)  ,
\]
so that (by comparing coefficients) we see that%
\begin{align*}
r_{1}+r_{2}+\cdots+r_{n}  &  =-a_{1};\\
\sum_{i<j}r_{i}r_{j}  &  =a_{2};\\
\sum_{i<j<k}r_{i}r_{j}r_{k}  &  =-a_{3};\\
&  \ldots;\\
r_{1}r_{2}\cdots r_{n}  &  =\left(  -1\right)  ^{n}a_{n}.
\end{align*}
These equalities are now known as \emph{Viete's formulas}. They allow
computing certain expressions in the $r_{i}$'s without having to compute the
$r_{i}$'s themselves. For instance, we can compute $r_{1}^{2}+r_{2}^{2}%
+\cdots+r_{n}^{2}$ (that is, the sum of the squares of all roots of $f$) by
observing that%
\[
\left(  r_{1}+r_{2}+\cdots+r_{n}\right)  ^{2}=\left(  r_{1}^{2}+r_{2}%
^{2}+\cdots+r_{n}^{2}\right)  +2\sum_{i<j}r_{i}r_{j},
\]
so that%
\[
r_{1}^{2}+r_{2}^{2}+\cdots+r_{n}^{2}=\left(  \underbrace{r_{1}+r_{2}%
+\cdots+r_{n}}_{=-a_{1}}\right)  ^{2}-2\underbrace{\sum_{i<j}r_{i}r_{j}%
}_{=a_{2}}=a_{1}^{2}-2a_{2}.
\]
This shows, among other things, that $r_{1}^{2}+r_{2}^{2}+\cdots+r_{n}^{2}$ is
an integer if the coefficients of $f$ are integers. Newton and others found
similar formulas for $r_{1}^{3}+r_{2}^{3}+\cdots+r_{n}^{3}$ and other such
polynomials. (These formulas are now known as the \emph{Newton--Girard
identities} -- see Theorem \ref{thm.sf.NG} below.) Gauss extended this to
arbitrary symmetric polynomials in $r_{1},r_{2},\ldots,r_{n}$ (by
algorithmically expressing them as polynomials in $a_{1},a_{2},\ldots,a_{n}$),
and used it in one of his proofs of the Fundamental Theorem of Algebra
\cite{Gauss16}; Galois used this to build what is now known as Galois theory
(even though modern treatments of Galois theory often avoid symmetric
polynomials); some harbingers of this can be seen in Cardano's solution of the
cubic equation. See \cite{Tignol16} and \cite{Armstrong-561fa18sp19} for the
real history.

Here is a simple modern application of the same ideas: Let $A\in
\mathbb{C}^{n\times n}$ be a matrix with eigenvalues $\lambda_{1},\lambda
_{2},\ldots,\lambda_{n}$ (listed with algebraic multiplicities). Let
$f\in\mathbb{C}\left[  x\right]  $ be a univariate polynomial. The
\emph{spectral mapping theorem} says that the eigenvalues of the matrix
$f\left[  A\right]  $ are $f\left[  \lambda_{1}\right]  ,f\left[  \lambda
_{2}\right]  ,\ldots,f\left[  \lambda_{n}\right]  $ (here, I am using the
notation $f\left[  a\right]  $ for the value of $f$ at some element $a$; this
is usually written $f\left(  a\right)  $). Thus, the characteristic polynomial
of $f\left[  A\right]  $ is
\begin{align*}
\chi_{f\left[  A\right]  }  &  =\left(  x-f\left[  \lambda_{1}\right]
\right)  \left(  x-f\left[  \lambda_{2}\right]  \right)  \cdots\left(
x-f\left[  \lambda_{n}\right]  \right) \\
&  =x^{n}-\left(  f\left[  \lambda_{1}\right]  +f\left[  \lambda_{2}\right]
+\cdots+f\left[  \lambda_{n}\right]  \right)  x^{n-1}+\left(  \sum
_{i<j}f\left[  \lambda_{i}\right]  f\left[  \lambda_{j}\right]  \right)
x^{n-2}\pm\cdots.
\end{align*}
Hence, all coefficients of $\chi_{f\left[  A\right]  }$ are symmetric
polynomials in the $\lambda_{i}$s that depend only on $f$ (not on $A$). In
particular, this shows that $\chi_{f\left[  A\right]  }$ is uniquely
determined by $f$ and $\chi_{A}$. But can you compute $\chi_{f\left[
A\right]  }$ exactly in terms of $f$ and $\chi_{A}$ without computing the
roots $\lambda_{1},\lambda_{2},\ldots,\lambda_{n}$ ? Yes, if you know how to
express \textbf{any} symmetric polynomial in the $\lambda_{i}$s in terms of
the coefficients of $\chi_{A}$. This is the same problem that Gauss solved
with his algorithm for expressing an arbitrary symmetric polynomial in the
roots of a polynomial in terms of the coefficients of the polynomial.
Incidentally, this algorithm also becomes helpful when one tries to generalize
the spectral mapping theorem to matrices over arbitrary commutative rings.
Here, eigenvalues don't always exist (let alone $n$ of them), so it becomes
necessary to restate the theorem in a language that does not rely on them.
Again, symmetric polynomials provide the way to do this.

\subsection{\label{sec.sf.sp}Definitions and examples of symmetric
polynomials}

\subsubsection{Defining symmetric polynomials}

\begin{convention}
\label{conv.sf.KN}Fix a commutative ring $K$. Fix an $N\in\mathbb{N}$.
(Perhaps $n$ would be more conventional, but lowercase letters are chronically
in short supply in this subject.)

Throughout this chapter, we will keep $K$ and $N$ fixed. We will use
Definition \ref{def.perm.Sn-iven}.
\end{convention}

Recall that $S_{N}$ denotes the $N$-th symmetric group, i.e., the group of all
permutations of the set $\left[  N\right]  :=\left\{  1,2,\ldots,N\right\}  $.

\begin{definition}
\label{def.sf.PS}\textbf{(a)} Let $\mathcal{P}$ be the polynomial ring
$K\left[  x_{1},x_{2},\ldots,x_{N}\right]  $ in $N$ variables over $K$. This
is not just a ring; it is a commutative $K$-algebra. \medskip

\textbf{(b)} The symmetric group $S_{N}$ acts on the set $\mathcal{P}$
according to the formula%
\[
\sigma\cdot f=f\left[  x_{\sigma\left(  1\right)  },x_{\sigma\left(  2\right)
},\ldots,x_{\sigma\left(  N\right)  }\right]  \ \ \ \ \ \ \ \ \ \ \text{for
any }\sigma\in S_{N}\text{ and any }f\in\mathcal{P}.
\]
Here, $f\left[  a_{1},a_{2},\ldots,a_{N}\right]  $ means the result of
substituting $a_{1},a_{2},\ldots,a_{N}$ for the indeterminates $x_{1}%
,x_{2},\ldots,x_{N}$ in a polynomial $f\in\mathcal{P}$.

(For example, if $N=4$ and $\sigma=\operatorname*{cyc}\nolimits_{1,2,3}\in
S_{4}$, then $\sigma\cdot f=f\left[  x_{\sigma\left(  1\right)  }%
,x_{\sigma\left(  2\right)  },x_{\sigma\left(  3\right)  },x_{\sigma\left(
4\right)  }\right]  =f\left[  x_{2},x_{3},x_{1},x_{4}\right]  $ for any
$f\in\mathcal{P}$, so that, for example,%
\[
\sigma\cdot\left(  2x_{1}+3x_{2}^{2}+4x_{3}-x_{4}^{15}\right)  =2x_{2}%
+3x_{3}^{2}+4x_{1}-x_{4}^{15}%
\]
and $\sigma\cdot\left(  x_{1}-x_{3}x_{4}\right)  =x_{2}-x_{1}x_{4}$.)

Roughly speaking, the group $S_{N}$ is thus acting on $\mathcal{P}$ by
permuting variables: A permutation $\sigma\in S_{N}$ transforms a polynomial
$f$ by substituting $x_{\sigma\left(  i\right)  }$ for each $x_{i}$.

Note that this action of $S_{N}$ on $\mathcal{P}$ is a well-defined group
action (as we will see in Proposition \ref{prop.sf.SN-acts} below). \medskip

\textbf{(c)} A polynomial $f\in\mathcal{P}$ is said to be \emph{symmetric} if
it satisfies%
\[
\sigma\cdot f=f\ \ \ \ \ \ \ \ \ \ \text{for all }\sigma\in S_{N}.
\]

\textbf{(d)} We let $\mathcal{S}$ be the set of all symmetric polynomials
$f\in\mathcal{P}$.
\end{definition}

\begin{example}
\label{exa.sf.PS1}Let $N=3$ and $K=\mathbb{Q}$, and let us rename the
indeterminates $x_{1},x_{2},x_{3}$ as $x,y,z$. Then: \medskip

\textbf{(a)} We have $x+y+z\in\mathcal{S}$ (since, for example, the simple
transposition $s_{1}\in S_{3}$ satisfies $s_{1}\cdot\left(  x+y+z\right)
=y+x+z=x+y+z$, and similarly any other $\sigma\in S_{3}$ also satisfies
$\sigma\cdot\left(  x+y+z\right)  =x+y+z$). \medskip

\textbf{(b)} We have $x+y\notin\mathcal{S}$ (since the transposition
$t_{1,3}\in S_{3}$ satisfies $t_{1,3}\cdot\left(  x+y\right)  =z+y\neq x+y$).
\medskip

\textbf{(c)} We have $\left(  x-y\right)  \left(  y-z\right)  \left(
z-x\right)  \notin\mathcal{S}$ (since the simple transposition $s_{1}\in
S_{3}$ transforms $\left(  x-y\right)  \left(  y-z\right)  \left(  z-x\right)
$ into
\begin{align*}
s_{1}\cdot\left(  \left(  x-y\right)  \left(  y-z\right)  \left(  z-x\right)
\right)   &  =\left(  y-x\right)  \left(  x-z\right)  \left(  z-y\right) \\
&  =-\left(  x-y\right)  \left(  y-z\right)  \left(  z-x\right) \\
&  \neq\left(  x-y\right)  \left(  y-z\right)  \left(  z-x\right)  ,
\end{align*}
because $-1\neq1$ in $\mathbb{Q}$). Actually, the polynomial $\left(
x-y\right)  \left(  y-z\right)  \left(  z-x\right)  $ is an example of an
\emph{antisymmetric} polynomial -- i.e., a polynomial $f\in\mathcal{P}$ such
that $\sigma\cdot f=\left(  -1\right)  ^{\sigma}f$ for all $\sigma\in S_{N}$.
However, if $K=\mathbb{Z}/2$ (or if $K$ is a $\mathbb{Z}/2$-algebra), then
antisymmetric polynomials and symmetric polynomials are the same thing.
\medskip

\textbf{(d)} We have $\left(  \left(  x-y\right)  \left(  y-z\right)  \left(
z-x\right)  \right)  ^{2}\in\mathcal{S}$. More generally, if $f\in\mathcal{P}$
is antisymmetric, then $f^{2}$ is symmetric, so that $f^{2}\in\mathcal{S}$.
\medskip

\textbf{(e)} We have $37\in\mathcal{S}$. More generally, any constant
polynomial $f\in\mathcal{P}$ is symmetric. \medskip

\textbf{(f)} We have $\left(  1-x\right)  \left(  1-y\right)  \left(
1-z\right)  \in\mathcal{S}$. \medskip

\textbf{(g)} We have $\dfrac{1}{\left(  1-x\right)  \left(  1-y\right)
\left(  1-z\right)  }\notin\mathcal{S}$, because this is not a polynomial. It
is an example of a \emph{symmetric power series}.
\end{example}

Some basic properties of our current setup are worth mentioning:

\begin{proposition}
\label{prop.sf.SN-acts}The action of $S_{N}$ on $\mathcal{P}$ is a
well-defined group action. In other words, the following holds: \medskip

\textbf{(a)} We have $\operatorname*{id}\nolimits_{\left[  N\right]  }\cdot
f=f$ for every $f\in\mathcal{P}$. \medskip

\textbf{(b)} We have $\left(  \sigma\tau\right)  \cdot f=\sigma\cdot\left(
\tau\cdot f\right)  $ for every $\sigma,\tau\in S_{N}$ and $f\in\mathcal{P}$.
\end{proposition}

The proof of this proposition is straightforward, but due to its somewhat
slippery nature (the two substitutions in part \textbf{(b)} are a particularly
frequent source of confusion), we present it in full:

\begin{fineprint}
\begin{proof}
[Proof of Proposition \ref{prop.sf.SN-acts}.]\textbf{(a)} If $f\in\mathcal{P}%
$, then the definition of $\operatorname*{id}\nolimits_{\left[  N\right]
}\cdot f$ yields
\[
\operatorname*{id}\nolimits_{\left[  N\right]  }\cdot f=f\left[
x_{\operatorname*{id}\left(  1\right)  },x_{\operatorname*{id}\left(
2\right)  },\ldots,x_{\operatorname*{id}\left(  N\right)  }\right]  =f\left[
x_{1},x_{2},\ldots,x_{N}\right]  =f.
\]
This proves Proposition \ref{prop.sf.SN-acts} \textbf{(a)}. \medskip

\textbf{(b)} Let $\sigma,\tau\in S_{N}$ and $f\in\mathcal{P}$. The definition
of $\sigma\cdot\left(  \tau\cdot f\right)  $ yields%
\begin{equation}
\sigma\cdot\left(  \tau\cdot f\right)  =\left(  \tau\cdot f\right)  \left[
x_{\sigma\left(  1\right)  },x_{\sigma\left(  2\right)  },\ldots
,x_{\sigma\left(  N\right)  }\right]  . \label{pf.prop.sf.SN-acts.b.1}%
\end{equation}

Write the polynomial $f$ in the form
\begin{equation}
f=\sum_{\left(  a_{1},a_{2},\ldots,a_{N}\right)  \in\mathbb{N}^{N}}f_{\left(
a_{1},a_{2},\ldots,a_{N}\right)  }x_{1}^{a_{1}}x_{2}^{a_{2}}\cdots
x_{N}^{a_{N}}, \label{pf.prop.sf.SN-acts.b.f=}%
\end{equation}
where $f_{\left(  a_{1},a_{2},\ldots,a_{N}\right)  }\in K$ are its
coefficients. The definition of $\tau\cdot f$ yields%
\[
\tau\cdot f=f\left[  x_{\tau\left(  1\right)  },x_{\tau\left(  2\right)
},\ldots,x_{\tau\left(  N\right)  }\right]  =\sum_{\left(  a_{1},a_{2}%
,\ldots,a_{N}\right)  \in\mathbb{N}^{N}}f_{\left(  a_{1},a_{2},\ldots
,a_{N}\right)  }x_{\tau\left(  1\right)  }^{a_{1}}x_{\tau\left(  2\right)
}^{a_{2}}\cdots x_{\tau\left(  N\right)  }^{a_{N}}%
\]
(here, we have substituted $x_{\tau\left(  1\right)  },x_{\tau\left(
2\right)  },\ldots,x_{\tau\left(  N\right)  }$ for $x_{1},x_{2},\ldots,x_{N}$
on both sides of (\ref{pf.prop.sf.SN-acts.b.f=})). Substituting $x_{\sigma
\left(  1\right)  },x_{\sigma\left(  2\right)  },\ldots,x_{\sigma\left(
N\right)  }$ for $x_{1},x_{2},\ldots,x_{N}$ on both sides of this equality, we
obtain%
\begin{align*}
&  \left(  \tau\cdot f\right)  \left[  x_{\sigma\left(  1\right)  }%
,x_{\sigma\left(  2\right)  },\ldots,x_{\sigma\left(  N\right)  }\right] \\
&  =\sum_{\left(  a_{1},a_{2},\ldots,a_{N}\right)  \in\mathbb{N}^{N}%
}f_{\left(  a_{1},a_{2},\ldots,a_{N}\right)  }\underbrace{x_{\sigma\left(
\tau\left(  1\right)  \right)  }^{a_{1}}x_{\sigma\left(  \tau\left(  2\right)
\right)  }^{a_{2}}\cdots x_{\sigma\left(  \tau\left(  N\right)  \right)
}^{a_{N}}}_{\substack{=x_{\left(  \sigma\tau\right)  \left(  1\right)
}^{a_{1}}x_{\left(  \sigma\tau\right)  \left(  2\right)  }^{a_{2}}\cdots
x_{\left(  \sigma\tau\right)  \left(  N\right)  }^{a_{N}}\\\text{(since
}\sigma\left(  \tau\left(  i\right)  \right)  =\left(  \sigma\tau\right)
\left(  i\right)  \text{ for all }i\in\left[  N\right]  \text{)}}}\\
&  \ \ \ \ \ \ \ \ \ \ \ \ \ \ \ \ \ \ \ \ \left(  \text{since our
substitution replaces each }x_{i}\text{ by }x_{\sigma\left(  i\right)
}\right) \\
&  =\sum_{\left(  a_{1},a_{2},\ldots,a_{N}\right)  \in\mathbb{N}^{N}%
}f_{\left(  a_{1},a_{2},\ldots,a_{N}\right)  }x_{\left(  \sigma\tau\right)
\left(  1\right)  }^{a_{1}}x_{\left(  \sigma\tau\right)  \left(  2\right)
}^{a_{2}}\cdots x_{\left(  \sigma\tau\right)  \left(  N\right)  }^{a_{N}}.
\end{align*}
On the other hand, the definition of the action of $S_{N}$ on $\mathcal{P}$
yields%
\begin{align*}
\left(  \sigma\tau\right)  \cdot f  &  =f\left[  x_{\left(  \sigma\tau\right)
\left(  1\right)  },x_{\left(  \sigma\tau\right)  \left(  2\right)  }%
,\ldots,x_{\left(  \sigma\tau\right)  \left(  N\right)  }\right] \\
&  =\sum_{\left(  a_{1},a_{2},\ldots,a_{N}\right)  \in\mathbb{N}^{N}%
}f_{\left(  a_{1},a_{2},\ldots,a_{N}\right)  }x_{\left(  \sigma\tau\right)
\left(  1\right)  }^{a_{1}}x_{\left(  \sigma\tau\right)  \left(  2\right)
}^{a_{2}}\cdots x_{\left(  \sigma\tau\right)  \left(  N\right)  }^{a_{N}}%
\end{align*}
(here, we have substituted $x_{\left(  \sigma\tau\right)  \left(  1\right)
},x_{\left(  \sigma\tau\right)  \left(  2\right)  },\ldots,x_{\left(
\sigma\tau\right)  \left(  N\right)  }$ for $x_{1},x_{2},\ldots,x_{N}$ on both
sides of (\ref{pf.prop.sf.SN-acts.b.f=})). Comparing these two equalities, we
obtain%
\[
\left(  \sigma\tau\right)  \cdot f=\left(  \tau\cdot f\right)  \left[
x_{\sigma\left(  1\right)  },x_{\sigma\left(  2\right)  },\ldots
,x_{\sigma\left(  N\right)  }\right]  =\sigma\cdot\left(  \tau\cdot f\right)
\]
(by (\ref{pf.prop.sf.SN-acts.b.1})). This proves Proposition
\ref{prop.sf.SN-acts} \textbf{(b)}.
\end{proof}
\end{fineprint}

We recall the following notions from abstract algebra:

\begin{itemize}
\item A $K$\emph{-algebra isomorphism} from a $K$-algebra $A$ to a $K$-algebra
$B$ means an invertible $K$-algebra morphism $f:A\rightarrow B$ such that its
inverse $f^{-1}:B\rightarrow A$ is a $K$-algebra morphism from $B$ to $A$.
(Actually, any invertible $K$-algebra morphism is an isomorphism.)

\item A $K$\emph{-algebra automorphism} of a $K$-algebra $A$ means a
$K$-algebra isomorphism from $A$ to $A$.
\end{itemize}

\begin{proposition}
\label{prop.sf.SN-acts-by-alg-auts}The group $S_{N}$ acts on $\mathcal{P}$ by
$K$-algebra automorphisms. In other words, for each $\sigma\in S_{N}$, the map%
\begin{align*}
\mathcal{P}  &  \rightarrow\mathcal{P},\\
f  &  \mapsto\sigma\cdot f
\end{align*}
is a $K$-algebra automorphism of $\mathcal{P}$ (that is, a $K$-algebra
isomorphism from $\mathcal{P}$ to $\mathcal{P}$).
\end{proposition}

\begin{fineprint}
\begin{proof}
[Proof of Proposition \ref{prop.sf.SN-acts-by-alg-auts} (sketched).]Fix
$\sigma\in S_{N}$. For any $f,g\in\mathcal{P}$, we have%
\begin{align}
\sigma\cdot\left(  fg\right)   &  =\left(  fg\right)  \left[  x_{\sigma\left(
1\right)  },x_{\sigma\left(  2\right)  },\ldots,x_{\sigma\left(  N\right)
}\right]  \ \ \ \ \ \ \ \ \ \ \left(  \text{by the definition of the action of
}S_{N}\text{ on }\mathcal{P}\right) \nonumber\\
&  =\underbrace{f\left[  x_{\sigma\left(  1\right)  },x_{\sigma\left(
2\right)  },\ldots,x_{\sigma\left(  N\right)  }\right]  }_{\substack{=\sigma
\cdot f\\\text{(by the definition of the action of }S_{N}\text{ on
}\mathcal{P}\text{)}}}\cdot\underbrace{g\left[  x_{\sigma\left(  1\right)
},x_{\sigma\left(  2\right)  },\ldots,x_{\sigma\left(  N\right)  }\right]
}_{\substack{=\sigma\cdot g\\\text{(by the definition of the action of }%
S_{N}\text{ on }\mathcal{P}\text{)}}}\nonumber\\
&  =\left(  \sigma\cdot f\right)  \cdot\left(  \sigma\cdot g\right)  .
\label{pf.prop.sf.SN-acts-by-alg-auts.1}%
\end{align}
Thus, the map%
\begin{align*}
\mathcal{P}  &  \rightarrow\mathcal{P},\\
f  &  \mapsto\sigma\cdot f
\end{align*}
respects multiplication. Similarly, this map respects addition, respects
scaling, respects the zero and respects the unity. Hence, this map is a
$K$-algebra morphism from $\mathcal{P}$ to $\mathcal{P}$. Furthermore, this
map is invertible, since its inverse is the map%
\begin{align*}
\mathcal{P}  &  \rightarrow\mathcal{P},\\
f  &  \mapsto\sigma^{-1}\cdot f.
\end{align*}
Thus, this map is an invertible $K$-algebra morphism from $\mathcal{P}$ to
$\mathcal{P}$, and therefore a $K$-algebra isomorphism from $\mathcal{P}$ to
$\mathcal{P}$. In other words, this map is a $K$-algebra automorphism of
$\mathcal{P}$. This proves Proposition \ref{prop.sf.SN-acts-by-alg-auts}.
\end{proof}
\end{fineprint}

\begin{theorem}
\label{thm.sf.S-subalg}The subset $\mathcal{S}$ is a $K$-subalgebra of
$\mathcal{P}$.
\end{theorem}

\begin{fineprint}
\begin{proof}
[Proof of Theorem \ref{thm.sf.S-subalg} (sketched).]We need to show that
$\mathcal{S}$ is closed under addition, multiplication and scaling, and that
$\mathcal{S}$ contains the zero and the unity of $\mathcal{P}$. Let me just
show that $\mathcal{S}$ is closed under multiplication (since all the other
claims are equally easy): Let $f,g\in\mathcal{S}$. We must show that
$fg\in\mathcal{S}$.

The polynomial $f$ is symmetric (since $f\in\mathcal{S}$); in other words,
$\sigma\cdot f=f$ for each $\sigma\in S_{N}$. Similarly, $\sigma\cdot g=g$ for
each $\sigma\in S_{N}$. Now, from (\ref{pf.prop.sf.SN-acts-by-alg-auts.1}), we
see that%
\[
\text{for each }\sigma\in S_{N}\text{, we have}\ \ \ \ \ \ \ \ \ \ \sigma
\cdot\left(  fg\right)  =\underbrace{\left(  \sigma\cdot f\right)
}_{\substack{=f\\\text{(as we have seen)}}}\cdot\underbrace{\left(
\sigma\cdot g\right)  }_{\substack{=g\\\text{(as we have seen)}}}=fg.
\]
This shows that $fg$ is symmetric, i.e., we have $fg\in\mathcal{S}$.

Forget that we fixed $f,g$. We thus have shown that $fg\in\mathcal{S}$ for any
$f,g\in\mathcal{S}$. This shows that $\mathcal{S}$ is closed under
multiplication. As explained above, this concludes our proof of Theorem
\ref{thm.sf.S-subalg}.
\end{proof}
\end{fineprint}

\begin{definition}
\label{def.sf.ring-of-symm}The $K$-subalgebra $\mathcal{S}$ of $\mathcal{P}$
is called the \emph{ring of symmetric polynomials in }$N$\emph{ variables over
}$K$.
\end{definition}

\subsubsection{The polynomials $e_{n}$, $h_{n}$ and $p_{n}$: definitions}

Some more terminology is worth defining:

\begin{definition}
\label{def.sf.monomial}\textbf{(a)} A \emph{monomial} is an expression of the
form $x_{1}^{a_{1}}x_{2}^{a_{2}}\cdots x_{N}^{a_{N}}$ with $a_{1},a_{2}%
,\ldots,a_{N}\in\mathbb{N}$. \medskip

\textbf{(b)} The \emph{degree} $\deg\mathfrak{m}$ of a monomial $\mathfrak{m}%
=x_{1}^{a_{1}}x_{2}^{a_{2}}\cdots x_{N}^{a_{N}}$ is defined to be $a_{1}%
+a_{2}+\cdots+a_{N}\in\mathbb{N}$. \medskip

\textbf{(c)} A monomial $\mathfrak{m}=x_{1}^{a_{1}}x_{2}^{a_{2}}\cdots
x_{N}^{a_{N}}$ is said to be \emph{squarefree} if $a_{1},a_{2},\ldots,a_{N}%
\in\left\{  0,1\right\}  $. (This is saying that no square or higher power of
an indeterminate appears in $\mathfrak{m}$; thus the name \textquotedblleft
squarefree\textquotedblright.) \medskip

\textbf{(d)} A monomial $\mathfrak{m}=x_{1}^{a_{1}}x_{2}^{a_{2}}\cdots
x_{N}^{a_{N}}$ is said to be \emph{primal} if there is at most one
$i\in\left[  N\right]  $ satisfying $a_{i}>0$. (This is saying that the
monomial $\mathfrak{m}$ contains no two distinct indeterminates. Thus, a
primal monomial is just $1$ or a power of an indeterminate.)
\end{definition}

Now we can define some specific symmetric polynomials:\footnote{We are
following the convention that any set $S$ and any negative integer $n$ satisfy
$S^{n}=\varnothing$, since there are no $n$-tuples for negative $n$.}

\begin{definition}
\label{def.sf.ehp}\textbf{(a)} For each $n\in\mathbb{Z}$, define a symmetric
polynomial $e_{n}\in\mathcal{S}$ by%
\[
e_{n}=\sum_{\substack{\left(  i_{1},i_{2},\ldots,i_{n}\right)  \in\left[
N\right]  ^{n};\\i_{1}<i_{2}<\cdots<i_{n}}}x_{i_{1}}x_{i_{2}}\cdots x_{i_{n}%
}=\left(  \text{sum of all squarefree monomials of degree }n\right)  .
\]
This $e_{n}$ is called the $n$\emph{-th elementary symmetric polynomial} in
$x_{1},x_{2},\ldots,x_{N}$. \medskip

\textbf{(b)} For each $n\in\mathbb{Z}$, define a symmetric polynomial
$h_{n}\in\mathcal{S}$ by%
\[
h_{n}=\sum_{\substack{\left(  i_{1},i_{2},\ldots,i_{n}\right)  \in\left[
N\right]  ^{n};\\i_{1}\leq i_{2}\leq\cdots\leq i_{n}}}x_{i_{1}}x_{i_{2}}\cdots
x_{i_{n}}=\left(  \text{sum of all monomials of degree }n\right)  .
\]
This $h_{n}$ is called the $n$\emph{-th complete homogeneous symmetric
polynomial} in $x_{1},x_{2},\ldots,x_{N}$. \medskip

\textbf{(c)} For each $n\in\mathbb{Z}$, define a symmetric polynomial
$p_{n}\in\mathcal{S}$ by%
\begin{align*}
p_{n}  &  =%
\begin{cases}
x_{1}^{n}+x_{2}^{n}+\cdots+x_{N}^{n}, & \text{if }n>0;\\
1, & \text{if }n=0;\\
0, & \text{if }n<0
\end{cases}
\\
&  =\left(  \text{sum of all primal monomials of degree }n\right)  .
\end{align*}
This $p_{n}$ is called the $n$\emph{-th power sum} in $x_{1},x_{2}%
,\ldots,x_{N}$.
\end{definition}

\begin{example}
\label{exa.sf.ehp.1}\textbf{(a)} The $2$-nd elementary symmetric polynomial is%
\begin{align*}
e_{2}  &  =\sum_{\substack{\left(  i_{1},i_{2}\right)  \in\left[  N\right]
^{2};\\i_{1}<i_{2}}}x_{i_{1}}x_{i_{2}}=\sum_{\substack{\left(  i,j\right)
\in\left[  N\right]  ^{2};\\i<j}}x_{i}x_{j}\\
&
\begin{array}
[c]{cccccccc}%
= & x_{1}x_{2} & + & x_{1}x_{3} & + & \cdots & + & x_{1}x_{N}\\
&  & + & x_{2}x_{3} & + & \cdots & + & x_{2}x_{N}\\
&  &  &  & \ddots & \cdots & \cdots & \cdots\\
&  &  &  &  &  & + & x_{N-1}x_{N}.
\end{array}
\end{align*}

\textbf{(b)} The $2$-nd complete homogeneous symmetric polynomial is%
\begin{align*}
h_{2}  &  =\sum_{\substack{\left(  i_{1},i_{2}\right)  \in\left[  N\right]
^{2};\\i_{1}\leq i_{2}}}x_{i_{1}}x_{i_{2}}=\sum_{\substack{\left(  i,j\right)
\in\left[  N\right]  ^{2};\\i\leq j}}x_{i}x_{j}\\
&
\begin{array}
[c]{cccccccccc}%
= & x_{1}^{2} & + & x_{1}x_{2} & + & x_{1}x_{3} & + & \cdots & + & x_{1}%
x_{N}\\
&  & + & x_{2}^{2} & + & x_{2}x_{3} & + & \cdots & + & x_{2}x_{N}\\
&  &  &  &  &  & \ddots & \cdots & \cdots & \cdots\\
&  &  &  &  &  & + & x_{N-1}^{2} & + & x_{N-1}x_{N}\\
&  &  &  &  &  &  &  & + & x_{N}^{2}.
\end{array}
\end{align*}

\textbf{(c)} The $2$-nd power sum is%
\[
p_{2}=x_{1}^{2}+x_{2}^{2}+\cdots+x_{N}^{2}.
\]
In light of the above two formulas for $e_{2}$ and $h_{2}$, we thus have
$h_{2}=p_{2}+e_{2}$. \medskip

\textbf{(d)} We have%
\[
e_{1}=h_{1}=p_{1}=x_{1}+x_{2}+\cdots+x_{N}.
\]

\textbf{(e)} We have%
\[
e_{0}=h_{0}=p_{0}=1.
\]

\textbf{(f)} If $n<0$, then $e_{n}=h_{n}=p_{n}=0$. \medskip

\textbf{(g)} If $N=0$ (that is, if $\mathcal{P}$ is a polynomial ring in $0$
variables), then $e_{n}=h_{n}=p_{n}=0$ for all $n>0$ (because there are no
monomials of positive degree in this case), but $e_{0}=h_{0}=p_{0}=1$ still
holds. This is a boring border case which, however, is important to get right.
\end{example}

\begin{proposition}
\label{prop.sf.en=0}For each integer $n>N$, we have $e_{n}=0$.
\end{proposition}

\begin{proof}
Let $n>N$ be an integer. Then, the set $\left[  N\right]  $ has no $n$
distinct elements. Thus, there exists no $n$-tuple $\left(  i_{1},i_{2}%
,\ldots,i_{n}\right)  \in\left[  N\right]  ^{n}$ satisfying $i_{1}%
<i_{2}<\cdots<i_{n}$ (because if $\left(  i_{1},i_{2},\ldots,i_{n}\right)  $
was such an $n$-tuple, then its $n$ entries $i_{1},i_{2},\ldots,i_{n}$ would
be $n$ distinct elements of $\left[  N\right]  $).

Now, the definition of $e_{n}$ yields%
\begin{align*}
e_{n}  &  =\sum_{\substack{\left(  i_{1},i_{2},\ldots,i_{n}\right)  \in\left[
N\right]  ^{n};\\i_{1}<i_{2}<\cdots<i_{n}}}x_{i_{1}}x_{i_{2}}\cdots x_{i_{n}%
}=\left(  \text{empty sum}\right) \\
&  \ \ \ \ \ \ \ \ \ \ \ \ \ \ \ \ \ \ \ \ \left(
\begin{array}
[c]{c}%
\text{since there exists no }n\text{-tuple }\left(  i_{1},i_{2},\ldots
,i_{n}\right)  \in\left[  N\right]  ^{n}\\
\text{satisfying }i_{1}<i_{2}<\cdots<i_{n}%
\end{array}
\right) \\
&  =0.
\end{align*}
This proves Proposition \ref{prop.sf.en=0}.
\end{proof}

Thus, there are only $N$ \textquotedblleft interesting\textquotedblright%
\ elementary symmetric polynomials: namely, $e_{1},e_{2},\ldots,e_{N}$. All
other $e_{n}$'s are either $1$ or $0$.

In contrast, there are infinitely many \textquotedblleft
interesting\textquotedblright\ complete homogeneous symmetric polynomials and
power sums (provided that $N>0$). For example, for $N=2$, we have $h_{5}%
=x_{1}^{5}+x_{1}^{4}x_{2}+x_{1}^{3}x_{2}^{2}+x_{1}^{2}x_{2}^{3}+x_{1}x_{2}%
^{4}+x_{2}^{5}$ and $p_{5}=x_{1}^{5}+x_{2}^{5}$.

\subsubsection{The polynomials $e_{n}$, $h_{n}$ and $p_{n}$: Newton--Girard
formulas}

We have so far defined three sequences $\left(  e_{0},e_{1},e_{2}%
,\ldots\right)  $, $\left(  h_{0},h_{1},h_{2},\ldots\right)  $ and $\left(
p_{0},p_{1},p_{2},\ldots\right)  $ of symmetric polynomials. The following
theorem (known as the \emph{Newton--Girard formulas} or the
\emph{Newton--Girard identities}) relates these three sequences:

\begin{theorem}
[Newton--Girard formulas]\label{thm.sf.NG}For any positive integer $n$, we
have%
\begin{align}
\sum_{j=0}^{n}\left(  -1\right)  ^{j}e_{j}h_{n-j}  &
=0;\label{eq.thm.sf.NG.eh}\\
\sum_{j=1}^{n}\left(  -1\right)  ^{j-1}e_{n-j}p_{j}  &  =ne_{n}%
;\label{eq.thm.sf.NG.ep}\\
\sum_{j=1}^{n}h_{n-j}p_{j}  &  =nh_{n}. \label{eq.thm.sf.NG.hp}%
\end{align}

\end{theorem}

\begin{example}
The formula (\ref{eq.thm.sf.NG.ep}), applied to $n=2$, says that $e_{1}%
p_{1}-p_{2}=2e_{2}$. Therefore,%
\[
p_{2}=e_{1}p_{1}-2e_{2}=\left(  x_{1}+x_{2}+\cdots+x_{N}\right)  \left(
x_{1}+x_{2}+\cdots+x_{N}\right)  -2\sum_{i<j}x_{i}x_{j}.
\]

\end{example}

Before we prove Theorem \ref{thm.sf.NG}, we establish some equalities in the
polynomial rings $\mathcal{P}\left[  t\right]  $ and $\mathcal{P}\left[
u,v\right]  $ (here, $t,u,v$ are three new indeterminates) and in the FPS ring
$\mathcal{P}\left[  \left[  t\right]  \right]  $:

\begin{proposition}
\label{prop.sf.e-h-FPS}\textbf{(a)} In the polynomial ring $\mathcal{P}\left[
t\right]  $, we have%
\[
\prod_{i=1}^{N}\left(  1-tx_{i}\right)  =\sum_{n\in\mathbb{N}}\left(
-1\right)  ^{n}t^{n}e_{n}.
\]

\textbf{(b)} In the polynomial ring $\mathcal{P}\left[  u,v\right]  $, we have%
\[
\prod_{i=1}^{N}\left(  u-vx_{i}\right)  =\sum_{n=0}^{N}\left(  -1\right)
^{n}u^{N-n}v^{n}e_{n}.
\]

\textbf{(c)} In the FPS ring $\mathcal{P}\left[  \left[  t\right]  \right]  $,
we have%
\[
\prod_{i=1}^{N}\dfrac{1}{1-tx_{i}}=\sum_{n\in\mathbb{N}}t^{n}h_{n}.
\]

\textbf{(d)} Recall the logarithmic derivative $\operatorname*{loder}f$ of any
FPS $f$ with constant term $1$ (see Definition \ref{def.fps.loder.1}). Then,%
\[
\operatorname*{loder}\prod_{i=1}^{N}\left(  1-tx_{i}\right)  =-\sum
_{n\in\mathbb{N}}t^{n}p_{n+1}.
\]
(Note that the logarithmic derivative is taken in the FPS ring $\mathcal{P}%
\left[  \left[  t\right]  \right]  $; thus, $x_{1},x_{2},\ldots,x_{N}$ are
viewed as constants.)
\end{proposition}

\begin{proof}
[Proof of Proposition \ref{prop.sf.e-h-FPS} (sketched).]\textbf{(a)} For each
$i\in\left\{  1,2,\ldots,N\right\}  $, we have%
\[
1+tx_{i}=\sum_{a\in\left\{  0,1\right\}  }\left(  tx_{i}\right)  ^{a}.
\]
Multiplying these equalities over all $i\in\left\{  1,2,\ldots,N\right\}  $,
we obtain%
\begin{align*}
\prod_{i=1}^{N}\left(  1+tx_{i}\right)   &  =\prod_{i=1}^{N}\ \ \sum
_{a\in\left\{  0,1\right\}  }\left(  tx_{i}\right)  ^{a}=\sum_{\left(
a_{1},a_{2},\ldots,a_{N}\right)  \in\left\{  0,1\right\}  ^{N}}%
\underbrace{\left(  tx_{1}\right)  ^{a_{1}}\left(  tx_{2}\right)  ^{a_{2}%
}\cdots\left(  tx_{N}\right)  ^{a_{N}}}_{=t^{a_{1}+a_{2}+\cdots+a_{N}}%
x_{1}^{a_{1}}x_{2}^{a_{2}}\cdots x_{N}^{a_{N}}}\\
&  \ \ \ \ \ \ \ \ \ \ \ \ \ \ \ \ \ \ \ \ \left(  \text{by Proposition
\ref{prop.fps.prodrule-fin-fin}}\right) \\
&  =\sum_{\left(  a_{1},a_{2},\ldots,a_{N}\right)  \in\left\{  0,1\right\}
^{N}}t^{a_{1}+a_{2}+\cdots+a_{N}}x_{1}^{a_{1}}x_{2}^{a_{2}}\cdots x_{N}%
^{a_{N}}=\sum_{\substack{\mathfrak{m}\text{ is a squarefree}\\\text{monomial}%
}}t^{\deg\mathfrak{m}}\mathfrak{m}\\
&  \ \ \ \ \ \ \ \ \ \ \ \ \ \ \ \ \ \ \ \ \left(
\begin{array}
[c]{c}%
\text{since the squarefree monomials are precisely}\\
\text{the }x_{1}^{a_{1}}x_{2}^{a_{2}}\cdots x_{N}^{a_{N}}\text{ with }\left(
a_{1},a_{2},\ldots,a_{N}\right)  \in\left\{  0,1\right\}  ^{N}\text{,}\\
\text{and since the degree of such a monomial}\\
\text{is precisely }a_{1}+a_{2}+\cdots+a_{N}%
\end{array}
\right) \\
&  =\sum_{n\in\mathbb{N}}\ \ \sum_{\substack{\mathfrak{m}\text{ is a
squarefree}\\\text{monomial of degree }n}}t^{n}\mathfrak{m}%
\ \ \ \ \ \ \ \ \ \ \left(
\begin{array}
[c]{c}%
\text{here, we have split the sum}\\
\text{according to the value of }\deg\mathfrak{m}%
\end{array}
\right) \\
&  =\sum_{n\in\mathbb{N}}t^{n}\underbrace{\sum_{\substack{\mathfrak{m}\text{
is a squarefree}\\\text{monomial of degree }n}}\mathfrak{m}}%
_{\substack{=\left(  \text{sum of all squarefree monomials of degree
}n\right)  \\=e_{n}\\\text{(by the definition of }e_{n}\text{)}}}\\
&  =\sum_{n\in\mathbb{N}}t^{n}e_{n}.
\end{align*}
Substituting $-t$ for $t$ on both sides of this equality, we obtain%
\[
\prod_{i=1}^{N}\left(  1-tx_{i}\right)  =\sum_{n\in\mathbb{N}}%
\underbrace{\left(  -t\right)  ^{n}}_{=\left(  -1\right)  ^{n}t^{n}}e_{n}%
=\sum_{n\in\mathbb{N}}\left(  -1\right)  ^{n}t^{n}e_{n}.
\]
This proves Proposition \ref{prop.sf.e-h-FPS} \textbf{(a)}. \medskip

\textbf{(b)} This is similar to part \textbf{(a)}. (See Exercise
\ref{exe.sf.e-h-FPS.b} for details.) \medskip

\textbf{(c)} For each $i\in\left\{  1,2,\ldots,N\right\}  $, we have%
\begin{align*}
\dfrac{1}{1-tx_{i}}  &  =1+tx_{i}+\left(  tx_{i}\right)  ^{2}+\left(
tx_{i}\right)  ^{3}+\cdots\\
&  \ \ \ \ \ \ \ \ \ \ \ \ \ \ \ \ \ \ \ \ \left(
\begin{array}
[c]{c}%
\text{by substituting }tx_{i}\text{ for }x\text{ in the}\\
\text{geometric series formula (\ref{eq.sec.gf.exas.1.1/1-x})}%
\end{array}
\right) \\
&  =\sum_{a\in\mathbb{N}}\left(  tx_{i}\right)  ^{a}.
\end{align*}
Multiplying these equalities over all $i\in\left\{  1,2,\ldots,N\right\}  $,
we obtain%
\begin{align*}
\prod_{i=1}^{N}\dfrac{1}{1-tx_{i}}  &  =\prod_{i=1}^{N}\ \ \sum_{a\in
\mathbb{N}}\left(  tx_{i}\right)  ^{a}\\
&  =\sum_{\left(  a_{1},a_{2},\ldots,a_{N}\right)  \in\mathbb{N}^{N}%
}\underbrace{\left(  tx_{1}\right)  ^{a_{1}}\left(  tx_{2}\right)  ^{a_{2}%
}\cdots\left(  tx_{N}\right)  ^{a_{N}}}_{=t^{a_{1}+a_{2}+\cdots+a_{N}}%
x_{1}^{a_{1}}x_{2}^{a_{2}}\cdots x_{N}^{a_{N}}}\\
&  \ \ \ \ \ \ \ \ \ \ \ \ \ \ \ \ \ \ \ \ \left(  \text{by Proposition
\ref{prop.fps.prodrule-fin-inf}}\right) \\
&  =\sum_{\left(  a_{1},a_{2},\ldots,a_{N}\right)  \in\mathbb{N}^{N}}%
t^{a_{1}+a_{2}+\cdots+a_{N}}x_{1}^{a_{1}}x_{2}^{a_{2}}\cdots x_{N}^{a_{N}%
}=\sum_{\mathfrak{m}\text{ is a monomial}}t^{\deg\mathfrak{m}}\mathfrak{m}\\
&  \ \ \ \ \ \ \ \ \ \ \ \ \ \ \ \ \ \ \ \ \left(
\begin{array}
[c]{c}%
\text{since the monomials are precisely}\\
\text{the }x_{1}^{a_{1}}x_{2}^{a_{2}}\cdots x_{N}^{a_{N}}\text{ with }\left(
a_{1},a_{2},\ldots,a_{N}\right)  \in\mathbb{N}^{N}\text{,}\\
\text{and since the degree of such a monomial}\\
\text{is precisely }a_{1}+a_{2}+\cdots+a_{N}%
\end{array}
\right) \\
&  =\sum_{n\in\mathbb{N}}\ \ \sum_{\substack{\mathfrak{m}\text{ is a
monomial}\\\text{of degree }n}}t^{n}\mathfrak{m}\ \ \ \ \ \ \ \ \ \ \left(
\begin{array}
[c]{c}%
\text{here, we have split the sum}\\
\text{according to the value of }\deg\mathfrak{m}%
\end{array}
\right) \\
&  =\sum_{n\in\mathbb{N}}t^{n}\underbrace{\sum_{\substack{\mathfrak{m}\text{
is a monomial}\\\text{of degree }n}}\mathfrak{m}}_{\substack{=\left(
\text{sum of all monomials of degree }n\right)  =h_{n}\\\text{(by the
definition of }h_{n}\text{)}}}=\sum_{n\in\mathbb{N}}t^{n}h_{n}.
\end{align*}
This proves Proposition \ref{prop.sf.e-h-FPS} \textbf{(c)}. \medskip

\textbf{(d)} For each $i\in\left\{  1,2,\ldots,N\right\}  $, we have%
\begin{align*}
\operatorname*{loder}\left(  1-tx_{i}\right)   &  =\dfrac{\left(
1-tx_{i}\right)  ^{\prime}}{1-tx_{i}}\ \ \ \ \ \ \ \ \ \ \left(  \text{by
Definition \ref{def.fps.loder.1}}\right) \\
&  =\dfrac{-x_{i}}{1-tx_{i}}\ \ \ \ \ \ \ \ \ \ \left(  \text{since }\left(
1-tx_{i}\right)  ^{\prime}=-x_{i}\right) \\
&  =-x_{i}\dfrac{1}{1-tx_{i}}\\
&  =-x_{i}\sum_{k\in\mathbb{N}}\left(  tx_{i}\right)  ^{k}%
\ \ \ \ \ \ \ \ \ \ \left(
\begin{array}
[c]{c}%
\text{since substituting }tx_{i}\text{ for }x\\
\text{in the formula }\dfrac{1}{1-x}=\sum_{k\in\mathbb{N}}x^{k}\\
\text{yields }\dfrac{1}{1-tx_{i}}=\sum_{k\in\mathbb{N}}\left(  tx_{i}\right)
^{k}%
\end{array}
\right) \\
&  =-\sum_{k\in\mathbb{N}}\underbrace{x_{i}\left(  tx_{i}\right)  ^{k}%
}_{=t^{k}x_{i}^{k+1}}=-\sum_{k\in\mathbb{N}}t^{k}x_{i}^{k+1}.
\end{align*}

But Corollary \ref{cor.fps.loder.prodk} (applied to $K=\mathcal{P}$ and $k=N$
and $f_{i}=1-tx_{i}$) yields%
\begin{align*}
\operatorname*{loder}\prod_{i=1}^{N}\left(  1-tx_{i}\right)   &  =\sum
_{i=1}^{N}\underbrace{\operatorname*{loder}\left(  1-tx_{i}\right)
}_{\substack{=-\sum_{k\in\mathbb{N}}t^{k}x_{i}^{k+1}\\\text{(as we just
showed)}}}=\sum_{i=1}^{N}\left(  -\sum_{k\in\mathbb{N}}t^{k}x_{i}^{k+1}\right)
\\
&  =-\sum_{k\in\mathbb{N}}t^{k}\underbrace{\sum_{i=1}^{N}x_{i}^{k+1}%
}_{\substack{=x_{1}^{k+1}+x_{2}^{k+1}+\cdots+x_{N}^{k+1}\\=p_{k+1}\\\text{(by
the definition of }p_{k+1}\text{)}}}=-\sum_{k\in\mathbb{N}}t^{k}p_{k+1}\\
&  =-\sum_{n\in\mathbb{N}}t^{n}p_{n+1}.
\end{align*}
This proves Proposition \ref{prop.sf.e-h-FPS} \textbf{(d)}.
\end{proof}

Let us now prove the Newton--Girard formulas:

\begin{proof}
[Proof of the 1st Newton--Girard formula (\ref{eq.thm.sf.NG.eh}).]In the FPS
ring $\mathcal{P}\left[  \left[  t\right]  \right]  $, we have%
\[
\prod_{i=1}^{N}\left(  1-tx_{i}\right)  =\sum_{n\in\mathbb{N}}\left(
-1\right)  ^{n}t^{n}e_{n}\ \ \ \ \ \ \ \ \ \ \left(  \text{by Proposition
\ref{prop.sf.e-h-FPS} \textbf{(a)}}\right)
\]
and%
\[
\prod_{i=1}^{N}\dfrac{1}{1-tx_{i}}=\sum_{n\in\mathbb{N}}t^{n}h_{n}%
\ \ \ \ \ \ \ \ \ \ \left(  \text{by Proposition \ref{prop.sf.e-h-FPS}
\textbf{(c)}}\right)  .
\]
Multiplying these two equalities, we obtain%
\begin{align*}
&  \left(  \prod_{i=1}^{N}\left(  1-tx_{i}\right)  \right)  \left(
\prod_{i=1}^{N}\dfrac{1}{1-tx_{i}}\right) \\
&  =\left(  \sum_{n\in\mathbb{N}}\left(  -1\right)  ^{n}t^{n}e_{n}\right)
\left(  \sum_{n\in\mathbb{N}}t^{n}h_{n}\right)  =\left(  \sum_{j\in\mathbb{N}%
}\left(  -1\right)  ^{j}t^{j}e_{j}\right)  \left(  \sum_{k\in\mathbb{N}}%
t^{k}h_{k}\right) \\
&  \ \ \ \ \ \ \ \ \ \ \ \ \ \ \ \ \ \ \ \ \left(
\begin{array}
[c]{c}%
\text{here, we have renamed the summation}\\
\text{indices }n\text{ and }n\text{ (in the two sums) as }j\text{ and }k
\end{array}
\right) \\
&  =\underbrace{\sum_{j\in\mathbb{N}}\ \ \sum_{k\in\mathbb{N}}}_{=\sum
_{\left(  j,k\right)  \in\mathbb{N}^{2}}}\left(  -1\right)  ^{j}%
\underbrace{t^{j}e_{j}t^{k}h_{k}}_{=e_{j}h_{k}t^{j+k}}=\sum_{\left(
j,k\right)  \in\mathbb{N}^{2}}\left(  -1\right)  ^{j}e_{j}h_{k}t^{j+k}\\
&  =\sum_{n\in\mathbb{N}}\left(  \sum_{\substack{\left(  j,k\right)
\in\mathbb{N}^{2};\\j+k=n}}\left(  -1\right)  ^{j}e_{j}h_{k}\right)  t^{n}%
\end{align*}
(here, we have split the sum according to the value of $j+k$). Comparing this
with%
\[
\left(  \prod_{i=1}^{N}\left(  1-tx_{i}\right)  \right)  \left(  \prod
_{i=1}^{N}\dfrac{1}{1-tx_{i}}\right)  =\prod_{i=1}^{N}\underbrace{\left(
\left(  1-tx_{i}\right)  \cdot\dfrac{1}{1-tx_{i}}\right)  }_{=1}=\prod
_{i=1}^{N}1=1,
\]
we obtain%
\[
1=\sum_{n\in\mathbb{N}}\left(  \sum_{\substack{\left(  j,k\right)
\in\mathbb{N}^{2};\\j+k=n}}\left(  -1\right)  ^{j}e_{j}h_{k}\right)  t^{n}.
\]
This is an equality between two FPSs in $\mathcal{P}\left[  \left[  t\right]
\right]  $. Comparing coefficients in front of $t^{n}$, we conclude that each
positive integer $n$ satisfies%
\[
0=\sum_{\substack{\left(  j,k\right)  \in\mathbb{N}^{2};\\j+k=n}}\left(
-1\right)  ^{j}e_{j}h_{k}=\sum_{j=0}^{n}\left(  -1\right)  ^{j}e_{j}h_{n-j}%
\]
(here, we have substituted $\left(  j,n-j\right)  $ for $\left(  j,k\right)  $
in the sum, since the map $\left\{  0,1,\ldots,n\right\}  \rightarrow\left\{
\left(  j,k\right)  \in\mathbb{N}^{2}\ \mid\ j+k=n\right\}  $ that sends each
$j\in\left\{  0,1,\ldots,n\right\}  $ to the pair $\left(  j,n-j\right)  $ is
a bijection). This proves the 1st Newton--Girard formula
(\ref{eq.thm.sf.NG.eh}).
\end{proof}

\begin{proof}
[Proof of the 2nd Newton--Girard formula (\ref{eq.thm.sf.NG.ep}).]Proposition
\ref{prop.sf.e-h-FPS} \textbf{(d)} yields
\[
\operatorname*{loder}\prod_{i=1}^{N}\left(  1-tx_{i}\right)  =-\sum
_{n\in\mathbb{N}}t^{n}p_{n+1}.
\]
Thus,%
\begin{align*}
-\sum_{n\in\mathbb{N}}t^{n}p_{n+1}  &  =\operatorname*{loder}\underbrace{\prod
_{i=1}^{N}\left(  1-tx_{i}\right)  }_{\substack{=\sum_{n\in\mathbb{N}}\left(
-1\right)  ^{n}t^{n}e_{n}\\\text{(by Proposition \ref{prop.sf.e-h-FPS}
\textbf{(a)})}}}=\operatorname*{loder}\left(  \sum_{n\in\mathbb{N}}\left(
-1\right)  ^{n}t^{n}e_{n}\right) \\
&  =\dfrac{\left(  \sum_{n\in\mathbb{N}}\left(  -1\right)  ^{n}t^{n}%
e_{n}\right)  ^{\prime}}{\sum_{n\in\mathbb{N}}\left(  -1\right)  ^{n}%
t^{n}e_{n}}\ \ \ \ \ \ \ \ \ \ \left(  \text{by Definition
\ref{def.fps.loder.1}}\right) \\
&  =\dfrac{\sum_{n>0}\left(  -1\right)  ^{n}nt^{n-1}e_{n}}{\sum_{n\in
\mathbb{N}}\left(  -1\right)  ^{n}t^{n}e_{n}}%
\end{align*}
(since the definition of a derivative yields $\left(  \sum_{n\in\mathbb{N}%
}\left(  -1\right)  ^{n}t^{n}e_{n}\right)  ^{\prime}=\sum_{n>0}\left(
-1\right)  ^{n}nt^{n-1}e_{n}$). Multiplying this by $-\sum_{n\in\mathbb{N}%
}\left(  -1\right)  ^{n}t^{n}e_{n}$ (so that the minus sign is moved to the
right-hand side and the denominator to the left-hand side), we obtain%
\begin{align*}
\left(  \sum_{n\in\mathbb{N}}\left(  -1\right)  ^{n}t^{n}e_{n}\right)  \left(
\sum_{n\in\mathbb{N}}t^{n}p_{n+1}\right)   &  =-\sum_{n>0}\left(  -1\right)
^{n}nt^{n-1}e_{n}=\sum_{n>0}\underbrace{\left(  -\left(  -1\right)
^{n}\right)  }_{=\left(  -1\right)  ^{n-1}}ne_{n}t^{n-1}\\
&  =\sum_{n>0}\left(  -1\right)  ^{n-1}ne_{n}t^{n-1}.
\end{align*}
Hence,%
\begin{align*}
&  \sum_{n>0}\left(  -1\right)  ^{n-1}ne_{n}t^{n-1}\\
&  =\left(  \sum_{n\in\mathbb{N}}\left(  -1\right)  ^{n}t^{n}e_{n}\right)
\left(  \sum_{n\in\mathbb{N}}t^{n}p_{n+1}\right) \\
&  =\left(  \sum_{n\in\mathbb{N}}\underbrace{t^{n}p_{n+1}}_{=p_{n+1}t^{n}%
}\right)  \left(  \sum_{n\in\mathbb{N}}\left(  -1\right)  ^{n}%
\underbrace{t^{n}e_{n}}_{=e_{n}t^{n}}\right) \\
&  =\left(  \sum_{n\in\mathbb{N}}p_{n+1}t^{n}\right)  \left(  \sum
_{n\in\mathbb{N}}\left(  -1\right)  ^{n}e_{n}t^{n}\right) \\
&  =\sum_{n\in\mathbb{N}}\left(  \sum_{i=0}^{n}p_{i+1}\left(  -1\right)
^{n-i}e_{n-i}\right)  t^{n}\\
&  \ \ \ \ \ \ \ \ \ \ \ \ \ \ \ \ \ \ \ \ \left(  \text{by the definition of
the product of two FPSs}\right) \\
&  =\sum_{n>0}\left(  \sum_{i=0}^{n-1}p_{i+1}\left(  -1\right)  ^{n-1-i}%
e_{n-1-i}\right)  t^{n-1}%
\end{align*}
(here, we have substituted $n-1$ for $n$ in the sum). Comparing $t^{n-1}%
$-coefficients on both sides of this equality, we conclude that every positive
integer $n$ satisfies%
\begin{align*}
\left(  -1\right)  ^{n-1}ne_{n}  &  =\sum_{i=0}^{n-1}p_{i+1}\left(  -1\right)
^{n-1-i}e_{n-1-i}=\sum_{j=1}^{n}\underbrace{p_{j-1+1}}_{=p_{j}}%
\underbrace{\left(  -1\right)  ^{n-1-\left(  j-1\right)  }}_{=\left(
-1\right)  ^{n-j}}\underbrace{e_{n-1-\left(  j-1\right)  }}_{=e_{n-j}}\\
&  \ \ \ \ \ \ \ \ \ \ \ \ \ \ \ \ \ \ \ \ \left(  \text{here, we have
substituted }j-1\text{ for }i\text{ in the sum}\right) \\
&  =\sum_{j=1}^{n}p_{j}\left(  -1\right)  ^{n-j}e_{n-j}=\sum_{j=1}%
^{n}\underbrace{\left(  -1\right)  ^{n-j}}_{=\left(  -1\right)  ^{n-1}\left(
-1\right)  ^{j-1}}e_{n-j}p_{j}\\
&  =\left(  -1\right)  ^{n-1}\sum_{j=1}^{n}\left(  -1\right)  ^{j-1}%
e_{n-j}p_{j}%
\end{align*}
and therefore (upon cancelling the element $\left(  -1\right)  ^{n-1}$, which
is clearly invertible)%
\[
ne_{n}=\sum_{j=1}^{n}\left(  -1\right)  ^{j-1}e_{n-j}p_{j}.
\]
This proves (\ref{eq.thm.sf.NG.ep}).
\end{proof}

\begin{proof}
[Proof of the 3rd Newton--Girard formula (\ref{eq.thm.sf.NG.hp}).]Proposition
\ref{prop.sf.e-h-FPS} \textbf{(d)} yields
\[
\operatorname*{loder}\prod_{i=1}^{N}\left(  1-tx_{i}\right)  =-\sum
_{n\in\mathbb{N}}t^{n}p_{n+1}.
\]
However,%
\[
\prod_{i=1}^{N}\underbrace{\left(  1-tx_{i}\right)  }_{=\left(  \dfrac
{1}{1-tx_{i}}\right)  ^{-1}}=\prod_{i=1}^{N}\left(  \dfrac{1}{1-tx_{i}%
}\right)  ^{-1}=\left(  \prod_{i=1}^{N}\dfrac{1}{1-tx_{i}}\right)
^{-1}=\left(  \sum_{n\in\mathbb{N}}t^{n}h_{n}\right)  ^{-1}%
\]
(since Proposition \ref{prop.sf.e-h-FPS} \textbf{(c)} yields $\prod_{i=1}%
^{N}\dfrac{1}{1-tx_{i}}=\sum_{n\in\mathbb{N}}t^{n}h_{n}$), so that%
\[
\operatorname*{loder}\prod_{i=1}^{N}\left(  1-tx_{i}\right)
=\operatorname*{loder}\left(  \left(  \sum_{n\in\mathbb{N}}t^{n}h_{n}\right)
^{-1}\right)  =-\operatorname*{loder}\left(  \sum_{n\in\mathbb{N}}t^{n}%
h_{n}\right)
\]
(by Corollary \ref{cor.fps.loder.inv}, applied to $K=\mathcal{P}$ and
$f=\sum_{n\in\mathbb{N}}t^{n}h_{n}$). Thus,%
\[
-\operatorname*{loder}\left(  \sum_{n\in\mathbb{N}}t^{n}h_{n}\right)
=\operatorname*{loder}\prod_{i=1}^{N}\left(  1-tx_{i}\right)  =-\sum
_{n\in\mathbb{N}}t^{n}p_{n+1}.
\]
Cancelling the minus sign on both sides, we thus obtain%
\[
\operatorname*{loder}\left(  \sum_{n\in\mathbb{N}}t^{n}h_{n}\right)
=\sum_{n\in\mathbb{N}}t^{n}p_{n+1}.
\]
Hence,%
\begin{align*}
\sum_{n\in\mathbb{N}}t^{n}p_{n+1}  &  =\operatorname*{loder}\left(  \sum
_{n\in\mathbb{N}}t^{n}h_{n}\right) \\
&  =\dfrac{\left(  \sum_{n\in\mathbb{N}}t^{n}h_{n}\right)  ^{\prime}}%
{\sum_{n\in\mathbb{N}}t^{n}h_{n}}\ \ \ \ \ \ \ \ \ \ \left(  \text{by
Definition \ref{def.fps.loder.1}}\right) \\
&  =\dfrac{\sum_{n>0}nt^{n-1}h_{n}}{\sum_{n\in\mathbb{N}}t^{n}h_{n}}%
\end{align*}
(since the definition of a derivative yields $\left(  \sum_{n\in\mathbb{N}%
}t^{n}h_{n}\right)  ^{\prime}=\sum_{n>0}nt^{n-1}h_{n}$). Multiplying this by
$\sum_{n\in\mathbb{N}}t^{n}h_{n}$, we obtain%
\[
\left(  \sum_{n\in\mathbb{N}}t^{n}h_{n}\right)  \left(  \sum_{n\in\mathbb{N}%
}t^{n}p_{n+1}\right)  =\sum_{n>0}nt^{n-1}h_{n}=\sum_{n>0}nh_{n}t^{n-1}.
\]
Hence,%
\begin{align*}
&  \sum_{n>0}nh_{n}t^{n-1}\\
&  =\left(  \sum_{n\in\mathbb{N}}t^{n}h_{n}\right)  \left(  \sum
_{n\in\mathbb{N}}t^{n}p_{n+1}\right) \\
&  =\left(  \sum_{n\in\mathbb{N}}\underbrace{t^{n}p_{n+1}}_{=p_{n+1}t^{n}%
}\right)  \left(  \sum_{n\in\mathbb{N}}\underbrace{t^{n}h_{n}}_{=h_{n}t^{n}%
}\right) \\
&  =\left(  \sum_{n\in\mathbb{N}}p_{n+1}t^{n}\right)  \left(  \sum
_{n\in\mathbb{N}}h_{n}t^{n}\right) \\
&  =\sum_{n\in\mathbb{N}}\left(  \sum_{i=0}^{n}p_{i+1}h_{n-i}\right)  t^{n}\\
&  \ \ \ \ \ \ \ \ \ \ \ \ \ \ \ \ \ \ \ \ \left(  \text{by the definition of
the product of two FPSs}\right) \\
&  =\sum_{n>0}\left(  \sum_{i=0}^{n-1}p_{i+1}h_{n-1-i}\right)  t^{n-1}%
\end{align*}
(here, we have substituted $n-1$ for $n$ in the sum). Comparing $t^{n-1}%
$-coefficients on both sides of this equality, we conclude that every positive
integer $n$ satisfies%
\begin{align*}
nh_{n}  &  =\sum_{i=0}^{n-1}p_{i+1}h_{n-1-i}=\sum_{j=1}^{n}%
\underbrace{p_{j-1+1}}_{=p_{j}}\underbrace{h_{n-1-\left(  j-1\right)  }%
}_{=h_{n-j}}\\
&  \ \ \ \ \ \ \ \ \ \ \ \ \ \ \ \ \ \ \ \ \left(  \text{here, we have
substituted }j-1\text{ for }i\text{ in the sum}\right) \\
&  =\sum_{j=1}^{n}p_{j}h_{n-j}.
\end{align*}
This proves (\ref{eq.thm.sf.NG.hp}).
\end{proof}

Thus, the proof of Theorem \ref{thm.sf.NG} is complete. There exist proofs of
other kinds as well (Exercise \ref{exe.sf.NG.p}).

Note that our above proof of Theorem \ref{thm.sf.NG} attests to the usefulness
of generating functions: Even though the polynomials $f\in\mathcal{P}$ already
involve $N$ variables $x_{1},x_{2},\ldots,x_{N}$, the proof proceeds by
adjoining yet another variable $t$ (to form the ring $\mathcal{P}\left[
\left[  t\right]  \right]  $).

\subsubsection{The Fundamental Theorem of Symmetric Polynomials}

The Newton--Girard formulas can be used to express the $e_{i}$'s and the
$h_{i}$'s in terms of each other, and the $p_{i}$'s in terms of the $e_{i}$'s
and the $h_{i}$'s, and finally the $e_{i}$'s and the $h_{i}$'s in terms of the
$p_{i}$'s when $K$ is a commutative $\mathbb{Q}$-algebra (i.e., when the
numbers $1,2,3,\ldots$ have inverses in $K$). More generally, it turns out
that we can express any symmetric polynomial in terms of $e_{i}$'s or of
$h_{i}$'s or (if $K$ is a commutative $\mathbb{Q}$-algebra) of $p_{i}$'s:

\begin{theorem}
[Fundamental Theorem of Symmetric Polynomials, due to Gauss et al.]%
\label{thm.sf.ftsf}\textbf{(a)} The elementary symmetric polynomials
$e_{1},e_{2},\ldots,e_{N}$ are algebraically independent (over $K$) and
generate the $K$-algebra $\mathcal{S}$.

In other words, each $f\in\mathcal{S}$ can be uniquely written as a polynomial
in $e_{1},e_{2},\ldots,e_{N}$.

In yet other words, the map%
\begin{align*}
\underbrace{K\left[  y_{1},y_{2},\ldots,y_{N}\right]  }_{\substack{\text{a
polynomial ring}\\\text{in }N\text{ variables}}}  &  \rightarrow\mathcal{S},\\
g  &  \mapsto g\left[  e_{1},e_{2},\ldots,e_{N}\right]
\end{align*}
is a $K$-algebra isomorphism. \medskip

\textbf{(b)} The complete homogeneous symmetric polynomials $h_{1}%
,h_{2},\ldots,h_{N}$ are algebraically independent (over $K$) and generate the
$K$-algebra $\mathcal{S}$.

In other words, each $f\in\mathcal{S}$ can be uniquely written as a polynomial
in $h_{1},h_{2},\ldots,h_{N}$.

In yet other words, the map%
\begin{align*}
\underbrace{K\left[  y_{1},y_{2},\ldots,y_{N}\right]  }_{\substack{\text{a
polynomial ring}\\\text{in }N\text{ variables}}}  &  \rightarrow\mathcal{S},\\
g  &  \mapsto g\left[  h_{1},h_{2},\ldots,h_{N}\right]
\end{align*}
is a $K$-algebra isomorphism. \medskip

\textbf{(c)} Now assume that $K$ is a commutative $\mathbb{Q}$-algebra (e.g.,
a field of characteristic $0$). Then, the power sums $p_{1},p_{2},\ldots
,p_{N}$ are algebraically independent (over $K$) and generate the $K$-algebra
$\mathcal{S}$.

In other words, each $f\in\mathcal{S}$ can be uniquely written as a polynomial
in $p_{1},p_{2},\ldots,p_{N}$.

In yet other words, the map%
\begin{align*}
\underbrace{K\left[  y_{1},y_{2},\ldots,y_{N}\right]  }_{\substack{\text{a
polynomial ring}\\\text{in }N\text{ variables}}}  &  \rightarrow\mathcal{S},\\
g  &  \mapsto g\left[  p_{1},p_{2},\ldots,p_{N}\right]
\end{align*}
is a $K$-algebra isomorphism.
\end{theorem}

\begin{example}
\textbf{(a)} Theorem \ref{thm.sf.ftsf} \textbf{(a)} yields that $p_{3}$ can be
uniquely written as a polynomial in $e_{1},e_{2},\ldots,e_{N}$. How to write
it this way?

Here is a method that (more generally) can be used to express $p_{n}$ (for any
given $n>0$) as a polynomial in $e_{1},e_{2},\ldots,e_{n}$. This method is
recursive, so we assume that all the \textquotedblleft
smaller\textquotedblright\ power sums $p_{1},p_{2},\ldots,p_{n-1}$ have
already been expressed in this way. Now, the 2nd Newton--Girard formula
(\ref{eq.thm.sf.NG.ep}) yields
\begin{align*}
ne_{n}  &  =\sum\limits_{j=1}^{n}\left(  -1\right)  ^{j-1}e_{n-j}p_{j}%
=\sum\limits_{j=1}^{n-1}\left(  -1\right)  ^{j-1}e_{n-j}p_{j}+\left(
-1\right)  ^{n-1}\underbrace{e_{n-n}}_{=e_{0}=1}p_{n}\\
&  \ \ \ \ \ \ \ \ \ \ \ \ \ \ \ \ \ \ \ \ \left(
\begin{array}
[c]{c}%
\text{here, we have split off the addend}\\
\text{for }j=n\text{ from the sum}%
\end{array}
\right) \\
&  =\sum\limits_{j=1}^{n-1}\left(  -1\right)  ^{j-1}e_{n-j}p_{j}+\left(
-1\right)  ^{n-1}p_{n}.
\end{align*}
Solving this equality for $p_{n}$, we obtain
\[
p_{n}=\left(  -1\right)  ^{n-1}\left(  ne_{n}-\sum\limits_{j=1}^{n-1}\left(
-1\right)  ^{j-1}e_{n-j}p_{j}\right)  .
\]
The right hand side can now be expressed in terms of $e_{1},e_{2},\ldots
,e_{n}$ (since the only power sums appearing in it are $p_{1},p_{2}%
,\ldots,p_{n-1}$, which we already know how to express in these terms);
therefore, we obtain an expression of $p_{n}$ as a polynomial in $e_{1}%
,e_{2},\ldots,e_{n}$.

For example, here is what we obtain for $n\in\left[  4\right]  $ by following
this method:
\begin{align*}
p_{1}  &  =e_{1};\\
p_{2}  &  =e_{1}^{2}-2e_{2};\\
p_{3}  &  =e_{1}^{3}-3e_{2}e_{1}+3e_{3};\\
p_{4}  &  =e_{1}^{4}-4e_{2}e_{1}^{2}+2e_{2}^{2}+4e_{3}e_{1}-4e_{4}.
\end{align*}
If $N<n$, then this expression of $p_{n}$ as a polynomial in $e_{1}%
,e_{2},\ldots,e_{n}$ becomes an expression as a polynomial in $e_{1}%
,e_{2},\ldots,e_{N}$ if we throw away all addends that contain one of
$e_{N+1},e_{N+2},\ldots,e_{n}$ as factor (we are allowed to do this, because
Proposition \ref{prop.sf.en=0} shows that all these addends are $0$). For
example, if $N=2$, then the expression $p_{4}=e_{1}^{4}-4e_{2}e_{1}^{2}%
+2e_{2}^{2}+4e_{3}e_{1}-4e_{4}$ becomes $p_{4}=e_{1}^{4}-4e_{2}e_{1}%
^{2}+2e_{2}^{2}$ this way. \medskip

\textbf{(b)} Theorem \ref{thm.sf.ftsf} \textbf{(b)} yields that $p_{3}$ can be
uniquely written as a polynomial in $h_{1},h_{2},\ldots,h_{N}$. How to write
it this way?

In part \textbf{(a)}, we have given a method to express $p_{n}$ (for any given
$n>0$) as a polynomial in $e_{1},e_{2},\ldots,e_{n}$. A similar method (but
using (\ref{eq.thm.sf.NG.hp}) instead of (\ref{eq.thm.sf.NG.ep})) can be used
to express $p_{n}$ (for any given $n>0$) as a polynomial in $h_{1}%
,h_{2},\ldots,h_{n}$. For example, for $n\in\left[  4\right]  $, this method
produces%
\begin{align*}
p_{1}  &  =h_{1};\\
p_{2}  &  =-h_{1}^{2}+2h_{2};\\
p_{3}  &  =h_{1}^{3}-3h_{2}h_{1}+3h_{3};\\
p_{4}  &  =-h_{1}^{4}+4h_{2}h_{1}^{2}-2h_{2}^{2}-4h_{3}h_{1}+4h_{4}.
\end{align*}
(The similarity with the analogous formulas expressing $p_{n}$ in terms of
$e_{1},e_{2},\ldots,e_{n}$ is not accidental -- the formulas are indeed
identical when $n$ is odd and differ in all signs when $n$ is even. Proving
this is Exercise \ref{exe.sf.eh-pol-same} \textbf{(b)}.)

So we can express $p_{n}$ as a polynomial in $h_{1},h_{2},\ldots,h_{n}$.
However, expressing $p_{n}$ as a polynomial in $h_{1},h_{2},\ldots,h_{N}$ is
harder when $N<n$. For example, if $N=2$, then the former expression is
$p_{4}=-h_{1}^{4}+4h_{2}h_{1}^{2}-2h_{2}^{2}-4h_{3}h_{1}+4h_{4}$, while the
latter is $p_{4}=-h_{1}^{4}+2h_{2}^{2}$; there is no easy way to get the
latter from the former. \medskip

\textbf{(c)} Assume that $K$ is a $\mathbb{Q}$-algebra. Theorem
\ref{thm.sf.ftsf} \textbf{(c)} yields that $e_{3}$ can be uniquely written as
a polynomial in $p_{1},p_{2},\ldots,p_{N}$. How to write it this way?

In part \textbf{(a)}, we have given a method to express $p_{n}$ (for any given
$n>0$) as a polynomial in $e_{1},e_{2},\ldots,e_{n}$. The crux of this method
was to solve the equation (\ref{eq.thm.sf.NG.ep}) for $p_{n}$. If we instead
solve it for $e_{n}$ (which is almost immediate: it gives $e_{n}=\dfrac{1}%
{n}\sum_{j=1}^{n}\left(  -1\right)  ^{j-1}e_{n-j}p_{j}$), then we obtain a
method for expressing $e_{n}$ (for any given $n>0$) as a polynomial in
$p_{1},p_{2},\ldots,p_{n}$. Applied to all $n\in\left[  4\right]  $, this
method produces%
\begin{align*}
e_{1}  &  =p_{1};\\
e_{2}  &  =\dfrac{1}{2}p_{1}^{2}-\dfrac{1}{2}p_{2};\\
e_{3}  &  =\dfrac{1}{6}p_{1}^{3}-\dfrac{1}{2}p_{2}p_{1}+\dfrac{1}{3}p_{3};\\
e_{4}  &  =\dfrac{1}{24}p_{1}^{4}-\dfrac{1}{4}p_{2}p_{1}^{2}+\dfrac{1}{8}%
p_{2}^{2}+\dfrac{1}{3}p_{3}p_{1}-\dfrac{1}{4}p_{4}.
\end{align*}
Note the fractions on the right hand sides! This is why we required $K$ to be
a $\mathbb{Q}$-algebra in Theorem \ref{thm.sf.ftsf} \textbf{(c)}. In general,
we cannot express $e_{n}$ in terms of $p_{1},p_{2},\ldots,p_{n}$ if the
integer $n$ is not invertible in $K$.

The question of expressing $e_{n}$ as a polynomial in $p_{1},p_{2}%
,\ldots,p_{N}$ (as opposed to $p_{1},p_{2},\ldots,p_{n}$) is easily reduced to
what we just have done: If $n\leq N$, then we have answered it already; if
$n>N$, then the answer is $e_{n}=0$ (by Proposition \ref{prop.sf.en=0}).
\medskip

\textbf{(d)} Not even algebraic independence of $p_{1},p_{2},\ldots,p_{N}$ is
true in general (if we don't assume that $K$ is a $\mathbb{Q}$-algebra)!
Indeed, if $K=\mathbb{Z}/2$, then%
\[
p_{1}^{2}=\left(  x_{1}+x_{2}+\cdots+x_{N}\right)  ^{2}=\underbrace{x_{1}%
^{2}+x_{2}^{2}+\cdots+x_{N}^{2}}_{=p_{2}}+\underbrace{2\sum_{i<j}x_{i}x_{j}%
}_{\substack{=0\\\text{if }K=\mathbb{Z}/2}}=p_{2}.
\]
More generally, if $K=\mathbb{Z}/p$ for some prime $p$, then the Idiot's
Binomial Formula (i.e., the formula $\left(  x+y\right)  ^{p}=x^{p}+y^{p}$
that holds in any commutative $\mathbb{Z}/p$-algebra) yields $p_{1}^{p}=p_{p}%
$. (Did I mention that lowercase letters are in short supply in the theory of
symmetric polynomials?) \medskip

\textbf{(e)} If $N=3$, then Theorem \ref{thm.sf.ftsf} \textbf{(b)} yields that
$h_{4}$ can be written as a polynomial in $h_{1},h_{2},h_{3}$. Here is how
this looks like:%
\[
h_{4}=h_{1}^{4}-3h_{2}h_{1}^{2}+h_{2}^{2}+2h_{3}h_{1}.
\]

\textbf{(f)} If $N=3$, then
\[
\left(  \left(  x-y\right)  \left(  y-z\right)  \left(  z-x\right)  \right)
^{2}=e_{2}^{2}e_{1}^{2}-4e_{2}^{3}-4e_{3}e_{1}^{3}+18e_{3}e_{2}e_{1}%
-27e_{3}^{2}%
\]
(where we are again denoting $x_{1},x_{2},x_{3}$ by $x,y,z$).
\end{example}

We omit the proof of Theorem \ref{thm.sf.ftsf} for now; various proofs can be
found in \cite[Theorem 171]{Benede25}, \cite[proof of Theorem 1]{BluCos16},
\cite[Theorem 1.2.1]{Dumas08}, \cite[Theorem 9.6.6]{Goodman}, \cite[\S 1.1]%
{Smith95} and many other places.\footnote{Some of these sources only state the
theorem in the case when $K$ is a field, or even only when $K=\mathbb{C}$; but
the proofs they give can easily be generalized to any commutative ring $K$.}

\subsubsection{A faster symmetry check}

Let us record a useful criterion for showing that a polynomial is symmetric:

\begin{lemma}
\label{lem.sf.simples-enough}For each $i\in\left[  N-1\right]  $, we consider
the simple transposition $s_{i}\in S_{N}$ defined in Definition
\ref{def.perm.si} (applied to $n=N$).

Let $f\in\mathcal{P}$. Assume that
\begin{equation}
s_{k}\cdot f=f\ \ \ \ \ \ \ \ \ \ \text{for each }k\in\left[  N-1\right]  .
\label{eq.lem.sf.simples-enough.ass}%
\end{equation}
Then, the polynomial $f$ is symmetric.
\end{lemma}

In plainer terms, Lemma \ref{lem.sf.simples-enough} says that if a polynomial
$f\in\mathcal{P}$ remains unchanged whenever we swap two adjacent
indeterminates (i.e., it remains unchanged if we swap $x_{1}$ with $x_{2}$; it
remains unchanged if we swap $x_{2}$ with $x_{3}$; it remains unchanged if we
swap $x_{3}$ with $x_{4}$; etc.), then this polynomial $f$ is symmetric. For
example, for $N=3$, it says that if a polynomial $f\in K\left[  x_{1}%
,x_{2},x_{3}\right]  $ satisfies $f\left[  x_{2},x_{1},x_{3}\right]  =f$ and
$f\left[  x_{1},x_{3},x_{2}\right]  =f$, then $f$ is symmetric.

\begin{proof}
[Proof of Lemma \ref{lem.sf.simples-enough}.]This follows from Corollary
\ref{cor.perm.generated} or from Theorem \ref{thm.perm.len.redword1}
\textbf{(a)}. See Section \ref{sec.details.sf.sp} for the details of this proof.
\end{proof}

\subsection{\label{sec.sf.m}$N$-partitions and monomial symmetric polynomials}

Recall that an \emph{(integer) partition} means a weakly decreasing finite
tuple of positive integers -- such as $\left(  5,3,3,2,1\right)  $.

Let us define a variant of this notion:

\begin{definition}
\label{def.sf.Npar}An $N$\emph{-partition} will mean a weakly decreasing
$N$-tuple of nonnegative integers. In other words, an $N$-partition means an
$N$-tuple $\left(  \lambda_{1},\lambda_{2},\ldots,\lambda_{N}\right)
\in\mathbb{N}^{N}$ with $\lambda_{1}\geq\lambda_{2}\geq\cdots\geq\lambda_{N}$.
\end{definition}

For example, $\left(  5,3,3,2,1,0,0\right)  $ is a $7$-partition.

Per se, an $N$-partition can contain zeroes and thus is not always a
partition. However, the $N$-partitions are \textquotedblleft more or less the
same\textquotedblright\ as the partitions of length $\leq N$. Indeed:

\begin{proposition}
\label{prop.sf.Npar-as-par}There is a bijection%
\begin{align*}
\left\{  \text{partitions of length }\leq N\right\}   &  \rightarrow\left\{
N\text{-partitions}\right\}  ,\\
\left(  \lambda_{1},\lambda_{2},\ldots,\lambda_{\ell}\right)   &
\mapsto\left(  \lambda_{1},\lambda_{2},\ldots,\lambda_{\ell}%
,\underbrace{0,0,\ldots,0}_{N-\ell\text{ zeroes}}\right)  .
\end{align*}

\end{proposition}

\begin{proof}
Straightforward. (We essentially did this back in our proof of Proposition
\ref{prop.pars.qbinom.alt-defs} \textbf{(a)}, although we used the letter $k$
instead of $N$ back then.)
\end{proof}

The $N$-partitions turn out to be closely connected to the ring $\mathcal{S}$.
Indeed, we will soon see various bases of the $K$-module $\mathcal{S}$, all of
which are indexed by the $N$-partitions. We shall construct the simplest one
in a moment. First, we define some auxiliary notations:

\begin{definition}
\label{def.sf.sort}Let $a=\left(  a_{1},a_{2},\ldots,a_{N}\right)
\in\mathbb{N}^{N}$. Then: \medskip

\textbf{(a)} We let $x^{a}$ denote the monomial $x_{1}^{a_{1}}x_{2}^{a_{2}%
}\cdots x_{N}^{a_{N}}$. \medskip

\textbf{(b)} We let $\operatorname*{sort}a$ mean the $N$-partition obtained
from $a$ by sorting the entries of $a$ in weakly decreasing order.
\end{definition}

For example, if $N=5$, then $x^{\left(  1,5,0,4,4\right)  }=x_{1}^{1}x_{2}%
^{5}x_{3}^{0}x_{4}^{4}x_{5}^{4}=x_{1}x_{2}^{5}x_{4}^{4}x_{5}^{4}$ and
\newline$\operatorname*{sort}\left(  1,5,0,4,4\right)  =\left(
5,4,4,1,0\right)  $.

\begin{definition}
\label{def.sf.m}Let $\lambda$ be any $N$-partition. Then, we define a
symmetric polynomial $m_{\lambda}\in\mathcal{S}$ by%
\[
m_{\lambda}:=\sum_{\substack{a\in\mathbb{N}^{N};\\\operatorname*{sort}%
a=\lambda}}x^{a}.
\]
This is called the \emph{monomial symmetric polynomial corresponding to
}$\lambda$.
\end{definition}

\begin{example}
Let $N=3$. Then,%
\begin{align*}
m_{\left(  2,1,0\right)  }  &  =\sum_{\substack{a\in\mathbb{N}^{3}%
;\\\operatorname*{sort}a=\left(  2,1,0\right)  }}x^{a}=x^{\left(
2,1,0\right)  }+x^{\left(  2,0,1\right)  }+x^{\left(  1,2,0\right)
}+x^{\left(  1,0,2\right)  }+x^{\left(  0,2,1\right)  }+x^{\left(
0,1,2\right)  }\\
&  =x_{1}^{2}x_{2}+x_{1}^{2}x_{3}+x_{1}x_{2}^{2}+x_{1}x_{3}^{2}+x_{2}^{2}%
x_{3}+x_{2}x_{3}^{2}%
\end{align*}
and%
\begin{align*}
m_{\left(  3,2,1\right)  }  &  =\sum_{\substack{a\in\mathbb{N}^{3}%
;\\\operatorname*{sort}a=\left(  3,2,1\right)  }}x^{a}=x^{\left(
3,2,1\right)  }+x^{\left(  3,1,2\right)  }+x^{\left(  2,3,1\right)
}+x^{\left(  2,1,3\right)  }+x^{\left(  1,3,2\right)  }+x^{\left(
1,2,3\right)  }\\
&  =x_{1}^{3}x_{2}^{2}x_{3}+x_{1}^{3}x_{2}x_{3}^{2}+x_{1}^{2}x_{2}^{3}%
x_{3}+x_{1}^{2}x_{2}x_{3}^{3}+x_{1}x_{2}^{3}x_{3}^{2}+x_{1}x_{2}^{2}x_{3}%
^{3}\\
&  =x_{1}^{3}x_{2}^{2}x_{3}+\left(  \text{all other }5\text{ permutations of
this monomial}\right)
\end{align*}
and%
\begin{align*}
m_{\left(  2,2,1\right)  }  &  =\sum_{\substack{a\in\mathbb{N}^{3}%
;\\\operatorname*{sort}a=\left(  2,2,1\right)  }}x^{a}=x^{\left(
2,2,1\right)  }+x^{\left(  2,1,2\right)  }+x^{\left(  1,2,2\right)  }\\
&  =x_{1}^{2}x_{2}^{2}x_{3}+x_{1}^{2}x_{2}x_{3}^{2}+x_{1}x_{2}^{2}x_{3}^{2}%
\end{align*}
and%
\[
m_{\left(  2,2,2\right)  }=\sum_{\substack{a\in\mathbb{N}^{3}%
;\\\operatorname*{sort}a=\left(  2,2,2\right)  }}x^{a}=x_{1}^{2}x_{2}^{2}%
x_{3}^{2}.
\]

\end{example}

Our symmetric polynomials $e_{n}$, $h_{n}$ and $p_{n}$ so far can be easily
expressed in terms of monomial symmetric polynomials:

\begin{proposition}
\label{prop.sf.ehp-through-m}\textbf{(a)} For each $n\in\left\{
0,1,\ldots,N\right\}  $, we have%
\[
e_{n}=m_{\left(  1,1,\ldots,1,0,0,\ldots,0\right)  },
\]
where $\left(  1,1,\ldots,1,0,0,\ldots,0\right)  $ is the $N$-tuple that
begins with $n$ many $1$'s and ends with $N-n$ many $0$'s. \medskip

\textbf{(b)} For each $n\in\mathbb{N}$, we have%
\[
h_{n}=\sum_{\substack{\lambda\text{ is an }N\text{-partition;}\\\left\vert
\lambda\right\vert =n}}m_{\lambda},
\]
where the \emph{size} $\left\vert \lambda\right\vert $ of an $N$-partition
$\lambda$ is defined to be the sum of its entries (i.e., if $\lambda=\left(
\lambda_{1},\lambda_{2},\ldots,\lambda_{N}\right)  $, then $\left\vert
\lambda\right\vert :=\lambda_{1}+\lambda_{2}+\cdots+\lambda_{N}$). \medskip

\textbf{(c)} Assume that $N>0$. For each $n\in\mathbb{N}$, we have%
\[
p_{n}=m_{\left(  n,0,0,\ldots,0\right)  },
\]
where $\left(  n,0,0,\ldots,0\right)  $ is the $N$-tuple that begins with an
$n$ and ends with $N-1$ zeroes.
\end{proposition}

\begin{proof}
Easy and LTTR.
\end{proof}

The monomial symmetric polynomials $m_{\lambda}$ are the \textquotedblleft
building blocks\textquotedblright\ of symmetric polynomials, in the same way
as the monomials are the \textquotedblleft building blocks\textquotedblright%
\ of polynomials. Here is a way to make this precise:

\begin{theorem}
\label{thm.sf.m-basis}\textbf{(a)} The family $\left(  m_{\lambda}\right)
_{\lambda\text{ is an }N\text{-partition}}$ is a basis of the $K$-module
$\mathcal{S}$. \medskip

\textbf{(b)} Each symmetric polynomial $f\in\mathcal{S}$ satisfies%
\[
f=\sum_{\substack{\lambda=\left(  \lambda_{1},\lambda_{2},\ldots,\lambda
_{N}\right)  \\\text{is an }N\text{-partition}}}\left(  \left[  x_{1}%
^{\lambda_{1}}x_{2}^{\lambda_{2}}\cdots x_{N}^{\lambda_{N}}\right]  f\right)
m_{\lambda}.
\]

\textbf{(c)} Let $n\in\mathbb{N}$. Let%
\[
\mathcal{S}_{n}:=\left\{  \text{homogeneous symmetric polynomials }%
f\in\mathcal{P}\text{ of degree }n\right\}
\]
(where we understand the zero polynomial $0\in\mathcal{P}$ to be homogeneous
of every degree). Then, $\mathcal{S}_{n}$ is a $K$-submodule of $\mathcal{S}$.
\medskip

\textbf{(d)} Define the \emph{size} of any $N$-partition $\lambda=\left(
\lambda_{1},\lambda_{2},\ldots,\lambda_{N}\right)  $ to be the number
$\lambda_{1}+\lambda_{2}+\cdots+\lambda_{N}\in\mathbb{N}$. Then, the family
$\left(  m_{\lambda}\right)  _{\lambda\text{ is an }N\text{-partition of size
}n}$ is a basis of the $K$-module $\mathcal{S}_{n}$.
\end{theorem}

\begin{example}
\label{exa.sf.m-basis.1}Let $N=3$, and let us rename the indeterminates
$x_{1},x_{2},x_{3}$ as $x,y,z$. The polynomial $\left(  x+y\right)  \left(
y+z\right)  \left(  z+x\right)  $ is symmetric, thus belongs to $\mathcal{S}$.
Expanding it, we find%
\begin{align*}
\left(  x+y\right)  \left(  y+z\right)  \left(  z+x\right)   &
=\underbrace{x^{2}y+x^{2}z+y^{2}x+y^{2}z+z^{2}x+z^{2}y}_{=m_{\left(
2,1,0\right)  }}+\,2\underbrace{xyz}_{=m_{\left(  1,1,1\right)  }}\\
&  =m_{\left(  2,1,0\right)  }+2m_{\left(  1,1,1\right)  }.
\end{align*}
Thus, we have written $\left(  x+y\right)  \left(  y+z\right)  \left(
z+x\right)  $ as a $K$-linear combination of $m_{\lambda}$'s for various
$N$-partitions $\lambda$. The same procedure (i.e., expanding, and then
collecting monomials that differ only in the order of their exponents, such as
the monomials $x^{2}y,x^{2}z,y^{2}x,y^{2}z,z^{2}x,z^{2}y$ in our example) can
be applied to any symmetric polynomial $f\in\mathcal{S}$, and always results
in a representation of $f$ as a $K$-linear combination of $m_{\lambda}$'s
(because the symmetry of $f$ ensures that monomials that differ only in the
order of their exponents appear in $f$ with equal coefficients).
\end{example}

The proof of Theorem \ref{thm.sf.m-basis} will rely on a simple proposition
that expresses how a permutation $\sigma\in S_{N}$ transforms the coefficients
of a polynomial $f\in\mathcal{P}$ (guess what: it permutes these coefficients):

\begin{proposition}
\label{prop.sf.sigma-pol-coeff}Let $\sigma\in S_{N}$ and $f\in\mathcal{P}$.
Then,%
\[
\left[  x_{1}^{a_{1}}x_{2}^{a_{2}}\cdots x_{N}^{a_{N}}\right]  \left(
\sigma\cdot f\right)  =\left[  x_{1}^{a_{\sigma\left(  1\right)  }}%
x_{2}^{a_{\sigma\left(  2\right)  }}\cdots x_{N}^{a_{\sigma\left(  N\right)
}}\right]  f
\]
for any $\left(  a_{1},a_{2},\ldots,a_{N}\right)  \in\mathbb{N}^{N}$.
\end{proposition}

Here, as in Section \ref{sec.gf.multivar}, we are using the notation $\left[
x_{1}^{a_{1}}x_{2}^{a_{2}}\cdots x_{N}^{a_{N}}\right]  g$ for the coefficient
of a monomial $x_{1}^{a_{1}}x_{2}^{a_{2}}\cdots x_{N}^{a_{N}}$ in a polynomial
$g\in\mathcal{P}$.

The proof of Proposition \ref{prop.sf.sigma-pol-coeff} is quite easy and can
be found in Section \ref{sec.details.sf.m}.

\begin{proof}
[Proof of Theorem \ref{thm.sf.m-basis} (sketched).]Here is a rough outline of
the proof; we leave the details to the reader. \medskip

\textbf{(a)} The method shown in Example \ref{exa.sf.m-basis.1} shows that
each $f\in\mathcal{S}$ is a $K$-linear combination of the family $\left(
m_{\lambda}\right)  _{\lambda\text{ is an }N\text{-partition}}$. Thus, this
family spans $\mathcal{S}$. It remains to show that this family is
$K$-linearly independent.

To show this, we observe the following: If you expand a linear combination
$\sum\limits_{\lambda\text{ is an }N\text{-partition}}a_{\lambda}m_{\lambda}$
(where $a_{\lambda}\in K$) as a sum of monomials, then none of the addends can
cancel (since no two $m_{\lambda}$'s have any monomial in common\footnote{and
since each $m_{\lambda}$ contains at least one monomial (trivial observation,
but should not be forgotten)}); thus, the linear combination cannot be $0$
unless all the $a_{\lambda}$'s are $0$. This proves the $K$-linear
independence of the family $\left(  m_{\lambda}\right)  _{\lambda\text{ is an
}N\text{-partition}}$. Thus, the proof of Theorem \ref{thm.sf.m-basis}
\textbf{(a)} is complete. \medskip

\textbf{(b)} This should be clear from Example \ref{exa.sf.m-basis.1} as well.
\medskip

\textbf{(c)} This is rather obvious: Any $K$-linear combination of homogeneous
polynomials of degree $n$ is again homogeneous of degree $n$. \medskip

\textbf{(d)} This follows by the same argument as part \textbf{(a)}, except
that we now need to observe that homogeneous polynomials of degree $n$ are
$K$-linear combinations of degree-$n$ monomials (rather than arbitrary monomials).
\end{proof}

\subsection{\label{sec.sf.schur}Schur polynomials}

\subsubsection{Alternants}

Here is one way to generate symmetric polynomials:

\begin{example}
Let $N=3$, and let us again abbreviate the indeterminates $x_{1},x_{2},x_{3}$
as $x,y,z$. For simplicity, we assume that $K$ is a field. As we know (from
the Vandermonde determinant -- specifically, Theorem \ref{thm.det.vander}
\textbf{(a)}), we have%
\[
\det\left(
\begin{array}
[c]{ccc}%
x^{2} & x & 1\\
y^{2} & y & 1\\
z^{2} & z & 1
\end{array}
\right)  =\prod_{i<j}\left(  x_{i}-x_{j}\right)  =\left(  x-y\right)  \left(
x-z\right)  \left(  y-z\right)  .
\]

What about similar determinants, such as $\det\left(
\begin{array}
[c]{ccc}%
x^{5} & x^{3} & 1\\
y^{5} & y^{3} & 1\\
z^{5} & z^{3} & 1
\end{array}
\right)  $ ? Just as in the proof of Lemma \ref{lem.det.vander.a.pol} (in
which we computed the original Vandermonde determinant), we can argue that
this is a polynomial in $x,y,z$ that is divisible by each of $x-y$ and $x-z$
and $y-z$ (since it becomes $0$ if we set one of $x,y,z$ equal to another).
Hence,%
\[
\det\left(
\begin{array}
[c]{ccc}%
x^{5} & x^{3} & 1\\
y^{5} & y^{3} & 1\\
z^{5} & z^{3} & 1
\end{array}
\right)  =\left(  x-y\right)  \left(  x-z\right)  \left(  y-z\right)  \cdot q
\]
for some $q\in K\left[  x,y,z\right]  $. However, this time, degree
considerations yield $\deg q=8-3=5$, so $q$ is no longer just a constant. What
is $q$ ? Using computer algebra, we see that%
\begin{align*}
q  &  =\dfrac{\det\left(
\begin{array}
[c]{ccc}%
x^{5} & x^{3} & 1\\
y^{5} & y^{3} & 1\\
z^{5} & z^{3} & 1
\end{array}
\right)  }{\left(  x-y\right)  \left(  x-z\right)  \left(  y-z\right)
}=\dfrac{-x^{5}y^{3}+x^{5}z^{3}+x^{3}y^{5}-x^{3}z^{5}-y^{5}z^{3}+y^{3}z^{5}%
}{-x^{2}y+x^{2}z+xy^{2}-xz^{2}-y^{2}z+yz^{2}}\\
&  =x^{2}y^{3}+x^{3}y^{2}+x^{2}z^{3}+x^{3}z^{2}+y^{2}z^{3}+y^{3}z^{2}%
+xyz^{3}\\
&  \ \ \ \ \ \ \ \ \ \ +xy^{3}z+x^{3}yz+2xy^{2}z^{2}+2x^{2}yz^{2}+2x^{2}%
y^{2}z\\
&  =m_{\left(  3,2,0\right)  }+m_{\left(  3,1,1\right)  }+2m_{\left(
2,2,1\right)  }\in\mathcal{S}.
\end{align*}
Note that $q\in\mathcal{S}$ can be easily seen without computing $q$. Indeed,
if we swap two of our variables $x,y,z$, then both $\det\left(
\begin{array}
[c]{ccc}%
x^{5} & x^{3} & 1\\
y^{5} & y^{3} & 1\\
z^{5} & z^{3} & 1
\end{array}
\right)  $ and $\left(  x-y\right)  \left(  x-z\right)  \left(  y-z\right)  $
get multiplied by $-1$, so their ratio $q$ stays unchanged. This shows that
$\sigma\cdot q=q$ whenever $\sigma\in S_{3}$ is a transposition. Since the
transpositions generate the group $S_{3}$ (indeed, Corollary
\ref{cor.perm.generated} yields that the simple transpositions $s_{1},s_{2}$
generate $S_{3}$), this entails that $\sigma\cdot q=q$ for any $\sigma\in
S_{3}$ (not just for transpositions). This means that $q$ is symmetric.

There is nothing special about the exponents $5$ and $3$ and $0$ in the above
determinant. More generally, for any $a,b,c\in\mathbb{N}$, we can define the
so-called \emph{alternant}%
\[
\det\left(
\begin{array}
[c]{ccc}%
x^{a} & x^{b} & x^{c}\\
y^{a} & y^{b} & y^{c}\\
z^{a} & z^{b} & z^{c}%
\end{array}
\right)  \in\mathcal{P}.
\]
When studying this alternant, we can WLOG assume that $a,b,c$ are distinct
(since otherwise, the alternant is just $0$) and furthermore assume that
$a>b>c$ (since the general case is reduced to this one by swapping the columns
around). The alternant is then a polynomial divisible by $x-y$ and $x-z$ and
$y-z$ (since it becomes $0$ if we set one of $x,y,z$ equal to another), and
thus divisible by $\left(  x-y\right)  \left(  x-z\right)  \left(  y-z\right)
=\det\left(
\begin{array}
[c]{ccc}%
x^{2} & x & 1\\
y^{2} & y & 1\\
z^{2} & z & 1
\end{array}
\right)  $ (the simplest nonzero alternant). Moreover, the ratio%
\[
\dfrac{\det\left(
\begin{array}
[c]{ccc}%
x^{a} & x^{b} & x^{c}\\
y^{a} & y^{b} & y^{c}\\
z^{a} & z^{b} & z^{c}%
\end{array}
\right)  }{\left(  x-y\right)  \left(  x-z\right)  \left(  y-z\right)
}=\dfrac{\det\left(
\begin{array}
[c]{ccc}%
x^{a} & x^{b} & x^{c}\\
y^{a} & y^{b} & y^{c}\\
z^{a} & z^{b} & z^{c}%
\end{array}
\right)  }{\det\left(
\begin{array}
[c]{ccc}%
x^{2} & x & 1\\
y^{2} & y & 1\\
z^{2} & z & 1
\end{array}
\right)  }%
\]
is a symmetric polynomial in $x,y,z$ (by the same argument that we used
before). Some experimentation suggests that all coefficients of this
polynomial are positive integers (to be rigorous, nonnegative integers). There
is probably no way of showing this without explicitly finding this polynomial
-- and the best way to do so is by defining this polynomial combinatorially.
\end{example}

Let us prepare for doing this.

\begin{definition}
\label{def.sf.alternants}\textbf{(a)} We let $\rho$ be the $N$-tuple $\left(
N-1,N-2,\ldots,N-N\right)  \in\mathbb{N}^{N}$. \medskip

\textbf{(b)} For any $N$-tuple $\alpha=\left(  \alpha_{1},\alpha_{2}%
,\ldots,\alpha_{N}\right)  \in\mathbb{N}^{N}$, we define%
\[
a_{\alpha}:=\det\left(  \underbrace{\left(  x_{i}^{\alpha_{j}}\right)  _{1\leq
i\leq N,\ 1\leq j\leq N}}_{\in\mathcal{P}^{N\times N}}\right)  \in
\mathcal{P}.
\]
This is called the $\alpha$\emph{-alternant} (of $x_{1},x_{2},\ldots,x_{N}$).
\end{definition}

For example, for $N=3$, we have%
\[
a_{\left(  5,3,0\right)  }=\det\left(
\begin{array}
[c]{ccc}%
x^{5} & x^{3} & 1\\
y^{5} & y^{3} & 1\\
z^{5} & z^{3} & 1
\end{array}
\right)  \ \ \ \ \ \ \ \ \ \ \left(  \text{where }\left(  x,y,z\right)
=\left(  x_{1},x_{2},x_{3}\right)  \right)  .
\]

Note that the definition of $a_{\rho}$ yields%
\begin{align}
a_{\rho}  &  =\det\left(  \left(  x_{i}^{\rho_{j}}\right)  _{1\leq i\leq
N,\ 1\leq j\leq N}\right)  =\det\left(  \left(  x_{i}^{N-j}\right)  _{1\leq
i\leq N,\ 1\leq j\leq N}\right) \nonumber\\
&  \ \ \ \ \ \ \ \ \ \ \ \ \ \ \ \ \ \ \ \ \left(  \text{since }\rho
_{j}=N-j\text{ for each }j\in\left[  N\right]  \right) \nonumber\\
&  =\prod_{1\leq i<j\leq N}\left(  x_{i}-x_{j}\right)
\label{eq.def.sf.alternants.arho=vdm}%
\end{align}
(by Theorem \ref{thm.det.vander} \textbf{(a)}, applied to $N$, $\mathcal{P}$
and $x_{i}$ instead of $n$, $K$ and $a_{i}$).

Thus, we suspect:

\begin{conjecture}
\label{conj.sf.schur-alt}For every $\alpha=\left(  \alpha_{1},\alpha
_{2},\ldots,\alpha_{N}\right)  \in\mathbb{N}^{N}$, the alternant $a_{\alpha}$
is a multiple of $a_{\rho}$ in the polynomial ring $\mathcal{P}$. Furthermore,
if $\alpha_{1}>\alpha_{2}>\cdots>\alpha_{N}$, then the polynomial $a_{\alpha
}/a_{\rho}$ has positive (more precisely, nonnegative) integer coefficients.
\end{conjecture}

We will prove this by explicitly constructing $a_{\alpha}/a_{\rho}$ combinatorially.

\subsubsection{Young diagrams and Schur polynomials}

Let us first define the \emph{Young diagram} of an $N$-partition. This is
analogous to the definition of the Young diagram of a partition (which we did
back in the proof of Proposition \ref{prop.pars.pkn=dual}):

\begin{definition}
\label{def.sf.ydiag}Let $\lambda$ be an $N$-partition.

The \emph{Young diagram} of $\lambda$ is defined to be a table of $N$
left-aligned rows, with the $i$-th row (counted from the top, as always)
having $\lambda_{i}$ boxes. Formally, the Young diagram of $\lambda$ is
defined as the set%
\[
\left\{  \left(  i,j\right)  \ \mid\ i\in\left[  N\right]  \text{ and }%
j\in\left[  \lambda_{i}\right]  \right\}  \subseteq\left\{  1,2,3,\ldots
\right\}  ^{2}.
\]
We visually represent each element $\left(  i,j\right)  $ of this Young
diagram as a box in row $i$ and column $j$; thus we obtain a table with $N$
left-aligned rows (some of which might be empty).

We denote the Young diagram of $\lambda$ by $Y\left(  \lambda\right)  $.
\end{definition}

For example, the $3$-partition $\left(  4,1,0\right)  $ has Young diagram%
\[
\ydiagram{4,1,0}\ \ .
\]
(The $3$-rd row is invisible since it has length $0$.) The four boxes in the
$1$-st (i.e., topmost) row of this diagram are $\left(  1,1\right)  $,
$\left(  1,2\right)  $, $\left(  1,3\right)  $ and $\left(  1,4\right)  $
(from left to right), while the single box in its $2$-nd row is $\left(
2,1\right)  $.

Now, we are going to fill our Young diagrams -- i.e., to put numbers in the boxes:

\begin{definition}
\label{def.sf.ytab}Let $\lambda$ be an $N$-partition.

A \emph{Young tableau} of shape $\lambda$ means a way of filling the boxes of
$Y\left(  \lambda\right)  $ with elements of $\left[  N\right]  $ (one element
per box). Formally speaking, it is defined as a map $T:Y\left(  \lambda
\right)  \rightarrow\left[  N\right]  $. We visually represent such a map by
filling in the number $T\left(  i,j\right)  $ into each box $\left(
i,j\right)  $.

We often abbreviate \textquotedblleft Young tableau\textquotedblright\ as
\textquotedblleft\emph{tableau}\textquotedblright. The plural of
\textquotedblleft tableau\textquotedblright\ is \textquotedblleft
tableaux\textquotedblright.
\end{definition}

For instance, here is a Young tableau of shape $\left(  4,3,3,0,0,0,0,\ldots
,0\right)  $ (defined for any $N\geq7$):%
\[
\ytableaushort{1724,336,261}\ \ .
\]
Formally speaking, this is a map $T:Y\left(  4,3,3,0,0,0,0,\ldots,0\right)
\rightarrow\left[  N\right]  $ that sends the pairs $\left(  1,1\right)
,\left(  1,2\right)  ,\left(  1,3\right)  ,\left(  1,4\right)  ,\left(
2,1\right)  ,\ldots$ to $1,7,2,4,3,\ldots$, respectively.

We will use some visually inspired language when talking about Young diagrams
and tableaux:

\begin{itemize}
\item For instance, the \emph{entry} of a tableau $T$ in box $\left(
i,j\right)  $ will mean the value $T\left(  i,j\right)  $.

\item Also, the $u$\emph{-th row} of a tableau $T$ (for a given $u\geq1$) will
mean the sequence of all entries of $T$ in the boxes $\left(  i,j\right)  $
with $i=u$.

\item Likewise, the $v$\emph{-th column} of a tableau $T$ (for a given
$v\geq1$) will mean the sequence of all entries of $T$ in the boxes $\left(
i,j\right)  $ with $j=v$.

\item If $T$ is a Young tableau of shape $\lambda$, then the boxes of
$Y\left(  \lambda\right)  $ will also be called the boxes of $T$.

\item Two boxes of a Young diagram (or of a tableau) are said to be
\emph{adjacent} if they have an edge in common when drawn on the picture
(i.e., when one of them has the form $\left(  i,j\right)  $, while the other
has the form $\left(  i,j+1\right)  $ or $\left(  i+1,j\right)  $).

\item The words \textquotedblleft north\textquotedblright, \textquotedblleft
west\textquotedblright, \textquotedblleft south\textquotedblright\ and
\textquotedblleft east\textquotedblright\ are to be understood according to
the picture of a Young diagram: e.g., the box $\left(  2,4\right)  $ lies one
step north and three steps west of the box $\left(  3,7\right)  $.
\end{itemize}

Some tableaux are better than others:

\begin{definition}
\label{def.sf.ssyt}Let $\lambda$ be an $N$-partition.

A Young tableau $T$ of shape $\lambda$ is said to be \emph{semistandard} if
its entries

\begin{itemize}
\item increase weakly along each row (from left to right);

\item increase strictly down each column (from top to bottom).
\end{itemize}

Formally speaking, this means that a Young tableau $T:Y\left(  \lambda\right)
\rightarrow\left[  N\right]  $ is semistandard if and only if

\begin{itemize}
\item we have $T\left(  i,j\right)  \leq T\left(  i,j+1\right)  $ for any
$\left(  i,j\right)  \in Y\left(  \lambda\right)  $ satisfying $\left(
i,j+1\right)  \in Y\left(  \lambda\right)  $;

\item we have $T\left(  i,j\right)  <T\left(  i+1,j\right)  $ for any $\left(
i,j\right)  \in Y\left(  \lambda\right)  $ satisfying $\left(  i+1,j\right)
\in Y\left(  \lambda\right)  $.
\end{itemize}

We let $\operatorname*{SSYT}\left(  \lambda\right)  $ denote the set of all
semistandard Young tableaux of shape $\lambda$. (This depends on $N$ as well,
but $N$ is fixed, so we omit it from our notation.) We will usually say
\textquotedblleft\emph{semistandard tableau}\textquotedblright\ instead of
\textquotedblleft semistandard Young tableau\textquotedblright.
\end{definition}

\begin{example}
Consider the following $6$ Young tableaux of shape $\left(  4,3,1,0,0,0,\ldots
,0\right)  $:%
\begin{align}
&  \ytableaushort{1334,235,4}\qquad\ytableaushort{2134,345,6}\qquad
\ytableaushort{1123,245,6}\label{eq.def.sf.ssyt.1.1}\\
& \nonumber\\
&  \ytableaushort{1234,567,8}\qquad\ytableaushort{1111,222,3}\qquad
\ytableaushort{1234,123,1}\nonumber
\end{align}
Which of these $6$ tableaux are semistandard? The first one is not
semistandard, since the entries in its second column do not strictly increase
down the column. The second one is not semistandard, since the entries in its
first row do not weakly increase along the row. The third one is semistandard.
The fourth one is semistandard, too. The fifth one is semistandard, too
(actually it has a special property: each of its entries is the smallest
possible value that an entry of a semistandard tableau could have in its box).
The sixth one is not semistandard, again because of the columns.
\end{example}

\begin{definition}
\label{def.sf.ytab.xT}Let $\lambda$ be an $N$-partition. If $T$ is any Young
tableau of shape $\lambda$, then we define the corresponding monomial%
\[
x_{T}:=\prod_{c\text{ is a box of }Y\left(  \lambda\right)  }x_{T\left(
c\right)  }=\prod_{\left(  i,j\right)  \in Y\left(  \lambda\right)
}x_{T\left(  i,j\right)  }=\prod_{k=1}^{N}x_{k}^{\left(  \text{\# of times
}k\text{ appears in }T\right)  }.
\]

\end{definition}

For example, the three Young tableaux in (\ref{eq.def.sf.ssyt.1.1}) have
corresponding monomials%
\begin{align*}
&  x_{1}x_{3}x_{3}x_{4}x_{2}x_{3}x_{5}x_{4}=x_{1}x_{2}x_{3}^{3}x_{4}^{2}%
x_{5},\\
&  x_{2}x_{1}x_{3}x_{4}x_{3}x_{4}x_{5}x_{6},\\
&  x_{1}x_{1}x_{2}x_{3}x_{2}x_{4}x_{5}x_{6}.
\end{align*}

\begin{definition}
\label{def.sf.schur}Let $\lambda$ be an $N$-partition. We define the
\emph{Schur polynomial} $s_{\lambda}\in\mathcal{P}$ by%
\[
s_{\lambda}:=\sum_{T\in\operatorname*{SSYT}\left(  \lambda\right)  }x_{T}.
\]

\end{definition}

\begin{example}
\label{exa.sf.schur-h-e}\textbf{(a)} Let $n\in\mathbb{N}$. Consider the
$N$-partition $\left(  n,0,0,\ldots,0\right)  $. The semistandard tableaux $T$
of shape $\left(  n,0,0,\ldots,0\right)  $ are simply the fillings of a single
row with $n$ elements of $\left[  N\right]  $ that weakly increase from left
to right:%
\[
T=\begin{ytableau}
i_1 & i_2 & \cdots & i_n
\end{ytableau}\ \ \ \ \ \ \ \ \ \ \text{with }i_{1}\leq i_{2}\leq\cdots\leq
i_{n}.
\]
Thus,
\[
s_{\left(  n,0,0,\ldots,0\right)  }=\sum_{i_{1}\leq i_{2}\leq\cdots\leq i_{n}%
}x_{i_{1}}x_{i_{2}}\cdots x_{i_{n}}=h_{n}.
\]

\textbf{(b)} Let $n\in\left\{  0,1,\ldots,N\right\}  $. Consider the
$N$-partition $\left(  1,1,\ldots,1,0,0,\ldots,0\right)  $ (with $n$ ones and
$N-n$ zeroes). The semistandard tableaux $T$ of shape $\left(  1,1,\ldots
,1,0,0,\ldots,0\right)  $ are simply the fillings of a single column with $n$
elements of $\left[  N\right]  $ that strictly increase from top to bottom:%
\[
T=\begin{ytableau}
i_1 \\ i_2 \\ \vdots \\ i_n
\end{ytableau}\ \ \ \ \ \ \ \ \ \ \text{with }i_{1}<i_{2}<\cdots<i_{n}.
\]
Hence,%
\[
s_{\left(  1,1,\ldots,1,0,0,\ldots,0\right)  \text{ (with }n\text{ ones and
}N-n\text{ zeroes)}}=\sum_{i_{1}<i_{2}<\cdots<i_{n}}x_{i_{1}}x_{i_{2}}\cdots
x_{i_{n}}=e_{n}.
\]

\textbf{(c)} Assume that $N\geq2$. Consider the $N$-partition $\left(
2,1,0,0,0,\ldots,0\right)  $ (with all entries from the third on being $0$).
The semistandard tableaux $T$ of shape $\left(  2,1,0,0,0,\ldots,0\right)  $
all have the form%
\[
T=\ytableaushort{ij,k}\ \ \ \ \ \ \ \ \ \ \text{with }i\leq j\text{ and }i<k.
\]
Hence,%
\begin{align*}
s_{\left(  2,1,0,0,0,\ldots,0\right)  }  &  =\sum_{\substack{i\leq
j;\\i<k}}x_{i}x_{j}x_{k}=\underbrace{\sum_{\substack{i<k;\\j=i}}x_{i}%
x_{j}x_{k}}_{=\sum_{i<k}x_{i}x_{i}x_{k}}+\underbrace{\sum_{i<j<k}x_{i}%
x_{j}x_{k}}_{=e_{3}}+\underbrace{\sum_{\substack{i<k;\\j=k}}x_{i}x_{j}x_{k}%
}_{=\sum_{i<k}x_{i}x_{k}x_{k}}+\underbrace{\sum_{i<k<j}x_{i}x_{j}x_{k}%
}_{=e_{3}}\\
&  \ \ \ \ \ \ \ \ \ \ \ \ \ \ \ \ \ \ \ \ \left(
\begin{array}
[c]{c}%
\text{since each triple }\left(  i,j,k\right)  \text{ of elements of }\left[
N\right] \\
\text{that satisfies }i\leq j\text{ and }i<k\text{ must satisfy}\\
\text{\textbf{exactly one} of the four}\\
\text{conditions }\left(  i<k\text{ and }j=i\right)  \text{ and }i<j<k\\
\text{and }\left(  i<k\text{ and }j=k\right)  \text{ and }i<k<j\text{,}\\
\text{and conversely, each triple satisfying one}\\
\text{of the latter four conditions must}\\
\text{satisfy }i\leq j\text{ and }i<k
\end{array}
\right) \\
&  =\sum_{i<k}x_{i}x_{i}x_{k}+e_{3}+\sum_{i<k}x_{i}x_{k}x_{k}+e_{3}%
=2e_{3}+\sum_{i<k}\underbrace{x_{i}x_{i}}_{=x_{i}^{2}}x_{k}+\sum_{i<k}%
x_{i}\underbrace{x_{k}x_{k}}_{=x_{k}^{2}}\\
&  =2e_{3}+\underbrace{\sum_{i<k}x_{i}^{2}x_{k}+\sum_{i<k}x_{i}x_{k}^{2}%
}_{=m_{\left(  2,1,0,0,\ldots,0\right)  }}=2e_{3}+\underbrace{m_{\left(
2,1,0,0,\ldots,0\right)  }}_{\substack{=e_{2}e_{1}-3e_{3}\\\text{(check
this!)}}}\\
&  =2e_{3}+e_{2}e_{1}-3e_{3}=e_{2}e_{1}-e_{3}.
\end{align*}

\end{example}

\begin{theorem}
\label{thm.sf.schur-symm}Let $\lambda$ be an $N$-partition. Then:\medskip

\textbf{(a)} The polynomial $s_{\lambda}$ is symmetric. \medskip

\textbf{(b)} We have
\[
a_{\lambda+\rho}=a_{\rho}\cdot s_{\lambda}.
\]
Here, the addition on $\mathbb{N}^{N}$ is defined entrywise: that is, if
$\alpha=\left(  \alpha_{1},\alpha_{2},\ldots,\alpha_{N}\right)  $ and
$\beta=\left(  \beta_{1},\beta_{2},\ldots,\beta_{N}\right)  $ are two
$N$-tuples in $\mathbb{N}^{N}$, then we set%
\[
\alpha+\beta:=\left(  \alpha_{1}+\beta_{1},\alpha_{2}+\beta_{2},\ldots
,\alpha_{N}+\beta_{N}\right)  .
\]

\end{theorem}

This theorem (once it will be proved) will yield Conjecture
\ref{conj.sf.schur-alt} in the case when $\alpha_{1}>\alpha_{2}>\cdots
>\alpha_{N}$. Indeed, if $\alpha=\left(  \alpha_{1},\alpha_{2},\ldots
,\alpha_{N}\right)  \in\mathbb{N}^{N}$ is an $N$-tuple satisfying $\alpha
_{1}>\alpha_{2}>\cdots>\alpha_{N}$, then we can write $\alpha$ as
$\alpha=\lambda+\rho$ for some $N$-partition $\lambda$ (namely, for
$\lambda=\alpha-\rho=\left(  \alpha_{1}-\left(  N-1\right)  ,\alpha
_{2}-\left(  N-2\right)  ,\ldots,\alpha_{N}-\left(  N-N\right)  \right)  $),
and then Theorem \ref{thm.sf.schur-symm} \textbf{(b)} will yield $a_{\alpha
}=a_{\rho}\cdot s_{\lambda}$, so that $a_{\alpha}/a_{\rho}=s_{\lambda}$, which
is a symmetric polynomial (by Theorem \ref{thm.sf.schur-symm} \textbf{(a)})
and furthermore has positive coefficients (by its combinatorial definition).
Once Conjecture \ref{conj.sf.schur-alt} is proved in the case when $\alpha
_{1}>\alpha_{2}>\cdots>\alpha_{N}$, the validity of its first claim in the
general case easily follows (because the alternant $a_{\alpha}$ is $0$ when
two of $\alpha_{1},\alpha_{2},\ldots,\alpha_{N}$ are equal\footnote{This is
easy to prove (but see Lemma \ref{lem.sf.alternant-0} \textbf{(a)} below for a
proof).}, and otherwise can be reduced to the case $\alpha_{1}>\alpha
_{2}>\cdots>\alpha_{N}$ by swapping the columns around\footnote{See Lemma
\ref{lem.sf.alternant-0} \textbf{(b)} below for the details of this.}).

We will prove Theorem \ref{thm.sf.schur-symm} \textbf{(a)} soon and Theorem
\ref{thm.sf.schur-symm} \textbf{(b)} later.

\subsubsection{Skew Young diagrams and skew Schur polynomials}

Before we get to the proof, let us generalize the situation a bit. First, we
will replace our Young diagrams by so-called \emph{skew Young diagrams}, which
are set differences of two Young diagrams; here are four
examples:\footnote{These are four skew Young diagrams, not five, since the
second from the right is disconnected. The rightmost skew Young diagram is the
regular Young diagram $Y\left(  \left(  2,1\right)  \right)  $; this is okay,
as any regular Young diagram is a skew one.}%
\[%
\begin{tabular}
[c]{|cccccccccccc|}\hline
&  & \multicolumn{1}{c|}{} &  &  & \multicolumn{1}{c|}{} &  &  &
\multicolumn{1}{c|}{} &  &  & \\
& $\ydiagram{1+2,0+2}$ & \multicolumn{1}{c|}{} &  & $\ydiagram{2+1,1+1,0+2}$ &
\multicolumn{1}{c|}{} &  & $\ydiagram{2+1,0+1}$ & \multicolumn{1}{c|}{} &  &
$\ydiagram{0+2,0+1}$ & \\
&  & \multicolumn{1}{c|}{} &  &  & \multicolumn{1}{c|}{} &  &  &
\multicolumn{1}{c|}{} &  &  & \\\hline
\end{tabular}
\]
Second, we unshackle the $N$ that bounds the length of our partitions (i.e.,
the $N$ in \textquotedblleft$N$-partition\textquotedblright) from the $N$ that
bounds the entries of our tableaux (i.e., the $N$ in \textquotedblleft%
$\mathcal{P}=K\left[  x_{1},x_{2},\ldots,x_{N}\right]  $\textquotedblright%
\ and in \textquotedblleft a tableau has entries in $\left[  N\right]
$\textquotedblright). The latter $N$ will still be $N$, while the former $N$
will become a new nonnegative integer $M$, which can be arbitrarily small or
large. So our partitions will be $M$-partitions, but our tableaux will still
have their entries in $\left[  N\right]  $.

\begin{definition}
\label{def.sf.par-subset}Let $M\in\mathbb{N}$. Let $\lambda$ and $\mu$ be two
$M$-partitions.

We say that $\mu\subseteq\lambda$ if and only if $Y\left(  \mu\right)
\subseteq Y\left(  \lambda\right)  $ (where $Y\left(  \mu\right)  $ and
$Y\left(  \lambda\right)  $ are defined as in Definition \ref{def.sf.ydiag},
with $N$ replaced by $M$). Equivalently, $\mu\subseteq\lambda$ if and only if%
\[
\text{each }i\in\left[  M\right]  \text{ satisfies }\mu_{i}\leq\lambda_{i}%
\]
(where we write the $M$-partitions $\mu$ and $\lambda$ as $\mu=\left(  \mu
_{1},\mu_{2},\ldots,\mu_{M}\right)  $ and $\lambda=\left(  \lambda_{1}%
,\lambda_{2},\ldots,\lambda_{M}\right)  $). Thus we have defined a partial
order $\subseteq$ on the set of all $M$-partitions.
\end{definition}

For example, we have $\left(  3,2,0\right)  \subseteq\left(  4,2,1\right)  $
(since $3\leq4$ and $2\leq2$ and $0\leq1$), but we don't have $\left(
2,2,0\right)  \subseteq\left(  3,1,0\right)  $ (since $2>1$). We can see this
on the Young diagrams:%
\begin{align*}
&  \text{we have }\ydiagram{3,2,0} \subseteq\ydiagram{4,2,1}\text{\ \ ,}\\
& \\
&  \text{but we don't have }\ydiagram{2,2,0} \subseteq
\ydiagram{3,1,0}\text{\ \ .}%
\end{align*}

\begin{definition}
\label{def.sf.skew-diag}Let $M\in\mathbb{N}$. Let $\lambda$ and $\mu$ be two
$M$-partitions such that $\mu\subseteq\lambda$. Then, we define the \emph{skew
Young diagram} $Y\left(  \lambda/\mu\right)  $ to be the set difference%
\begin{align*}
Y\left(  \lambda\right)  \setminus Y\left(  \mu\right)   &  =\left\{  \left(
i,j\right)  \ \mid\ i\in\left[  M\right]  \text{ and }j\in\left[  \lambda
_{i}\right]  \setminus\left[  \mu_{i}\right]  \right\} \\
&  =\left\{  \left(  i,j\right)  \ \mid\ i\in\left[  M\right]  \text{ and
}j\in\mathbb{Z}\text{ and }\mu_{i}<j\leq\lambda_{i}\right\}  .
\end{align*}

\end{definition}

For example,%
\[
Y\left(  \left(  4,3,1\right)  /\left(  2,1,0\right)  \right)
=\ydiagram{2+2,1+2,0+1}\ \ .
\]
We note that any row or column of a skew Young diagram $Y\left(  \lambda
/\mu\right)  $ is contiguous, i.e., has no holes between boxes. Better yet, if
$\left(  a,b\right)  $ and $\left(  e,f\right)  $ are two boxes of $Y\left(
\lambda/\mu\right)  $, then any box $\left(  c,d\right)  $ that lies
\textquotedblleft between\textquotedblright\ them (i.e., that satisfies $a\leq
c\leq e$ and $b\leq d\leq f$) must also belong to $Y\left(  \lambda
/\mu\right)  $. Let us state this as a lemma:

\begin{lemma}
[Convexity of skew Young diagrams]\label{lem.sf.skew-diag.convexity}Let
$M\in\mathbb{N}$. Let $\lambda$ and $\mu$ be two $M$-partitions such that
$\mu\subseteq\lambda$. Let $\left(  a,b\right)  $ and $\left(  e,f\right)  $
be two elements of $Y\left(  \lambda/\mu\right)  $. Let $\left(  c,d\right)
\in\mathbb{Z}^{2}$ satisfy $a\leq c\leq e$ and $b\leq d\leq f$. Then, $\left(
c,d\right)  \in Y\left(  \lambda/\mu\right)  $.
\end{lemma}

\begin{proof}
[Hints to the proof of Lemma \ref{lem.sf.skew-diag.convexity}.]This follows
easily from the definition of $Y\left(  \lambda/\mu\right)  $ using the fact
that $\lambda$ and $\mu$ are weakly decreasing $M$-tuples. A detailed proof
can be found in Section \ref{sec.details.sf.schur}.
\end{proof}

Lemma \ref{lem.sf.skew-diag.convexity} is known as the \emph{convexity} of
$Y\left(  \lambda/\mu\right)  $ (albeit in a very specific sense of the word
\textquotedblleft convexity\textquotedblright).

Next, we can define Young tableaux of shape $\lambda/\mu$ whenever $\lambda$
and $\mu$ are two $M$-partitions satisfying $\mu\subseteq\lambda$. The
definition is analogous to Definition \ref{def.sf.ytab}, except that we are
only filling the boxes of $Y\left(  \lambda/\mu\right)  $ (rather than all
boxes of $Y\left(  \lambda\right)  $) this time:

\begin{definition}
\label{def.sf.skew-tab}Let $M\in\mathbb{N}$. Let $\lambda$ and $\mu$ be two
$M$-partitions such that $\mu\subseteq\lambda$. A \emph{Young tableau} of
shape $\lambda/\mu$ means a way of filling the boxes of $Y\left(  \lambda
/\mu\right)  $ with elements of $\left[  N\right]  $ (one element per box).
Formally speaking, it is defined as a map $T:Y\left(  \lambda/\mu\right)
\rightarrow\left[  N\right]  $. We visually represent such a map by filling in
the number $T\left(  i,j\right)  $ into each box $\left(  i,j\right)  $.

Young tableaux of shape $\lambda/\mu$ are often called \emph{skew Young
tableaux}.

If we don't have $\mu\subseteq\lambda$, then we agree that there are no Young
tableaux of shape $\lambda/\mu$.
\end{definition}

The notion of a semistandard tableau of shape $\lambda/\mu$ is, again, defined
in the same way as for shape $\lambda$:

\begin{definition}
\label{def.sf.skew-ssyt}Let $M\in\mathbb{N}$. Let $\lambda$ and $\mu$ be two
$M$-partitions.

A Young tableau $T$ of shape $\lambda/\mu$ is said to be \emph{semistandard}
if its entries

\begin{itemize}
\item increase weakly along each row (from left to right);

\item increase strictly down each column (from top to bottom).
\end{itemize}

Formally speaking, this means that a Young tableau $T:Y\left(  \lambda
/\mu\right)  \rightarrow\left[  N\right]  $ is semistandard if and only if

\begin{itemize}
\item we have $T\left(  i,j\right)  \leq T\left(  i,j+1\right)  $ for any
$\left(  i,j\right)  \in Y\left(  \lambda/\mu\right)  $ satisfying $\left(
i,j+1\right)  \in Y\left(  \lambda/\mu\right)  $;

\item we have $T\left(  i,j\right)  <T\left(  i+1,j\right)  $ for any $\left(
i,j\right)  \in Y\left(  \lambda/\mu\right)  $ satisfying $\left(
i+1,j\right)  \in Y\left(  \lambda/\mu\right)  $.
\end{itemize}

We let $\operatorname*{SSYT}\left(  \lambda/\mu\right)  $ denote the set of
all semistandard Young tableaux of shape $\lambda/\mu$. We will usually say
\textquotedblleft\emph{semistandard tableau}\textquotedblright\ instead of
\textquotedblleft semistandard Young tableau\textquotedblright.
\end{definition}

For example, here is a semistandard Young tableau of shape $\left(
4,3,1\right)  /\left(  2,1,0\right)  $:%
\[
\ytableaushort{\none\none 13,\none 22,1}\ \ .
\]
Meanwhile, there are no Young tableaux of shape $\left(  3,2,1\right)
/\left(  2,2,2\right)  $ (since we don't have $\left(  2,2,2\right)
\subseteq\left(  3,2,1\right)  $), and thus the set $\operatorname*{SSYT}%
\left(  \left(  3,2,1\right)  /\left(  2,2,2\right)  \right)  $ is empty.

The phrases \textquotedblleft increase weakly along each row\textquotedblright%
\ and \textquotedblleft increase strictly down each column\textquotedblright%
\ in Definition \ref{def.sf.skew-ssyt} have been formalized in terms of
adjacent entries: e.g., we have declared \textquotedblleft increase weakly
along each row\textquotedblright\ to mean \textquotedblleft$T\left(
i,j\right)  \leq T\left(  i,j+1\right)  $\textquotedblright\ rather than
\textquotedblleft$T\left(  i,j_{1}\right)  \leq T\left(  i,j_{2}\right)  $
whenever $j_{1}\leq j_{2}$\textquotedblright. However, since any row or column
of $Y\left(  \lambda/\mu\right)  $ is contiguous, the latter stronger meaning
actually follows from the former. To wit:

\begin{lemma}
\label{lem.sf.skew-ssyt.increase}Let $M\in\mathbb{N}$. Let $\lambda$ and $\mu$
be two $M$-partitions. Let $T$ be a semistandard Young tableau of shape
$\lambda/\mu$. Then: \medskip

\textbf{(a)} If $\left(  i,j_{1}\right)  $ and $\left(  i,j_{2}\right)  $ are
two elements of $Y\left(  \lambda/\mu\right)  $ satisfying $j_{1}\leq j_{2}$,
then $T\left(  i,j_{1}\right)  \leq T\left(  i,j_{2}\right)  $. \medskip

\textbf{(b)} If $\left(  i_{1},j\right)  $ and $\left(  i_{2},j\right)  $ are
two elements of $Y\left(  \lambda/\mu\right)  $ satisfying $i_{1}\leq i_{2}$,
then $T\left(  i_{1},j\right)  \leq T\left(  i_{2},j\right)  $. \medskip

\textbf{(c)} If $\left(  i_{1},j\right)  $ and $\left(  i_{2},j\right)  $ are
two elements of $Y\left(  \lambda/\mu\right)  $ satisfying $i_{1}<i_{2}$, then
$T\left(  i_{1},j\right)  <T\left(  i_{2},j\right)  $. \medskip

\textbf{(d)} If $\left(  i_{1},j_{1}\right)  $ and $\left(  i_{2}%
,j_{2}\right)  $ are two elements of $Y\left(  \lambda/\mu\right)  $
satisfying $i_{1}\leq i_{2}$ and $j_{1}\leq j_{2}$, then $T\left(  i_{1}%
,j_{1}\right)  \leq T\left(  i_{2},j_{2}\right)  $. \medskip

\textbf{(e)} If $\left(  i_{1},j_{1}\right)  $ and $\left(  i_{2}%
,j_{2}\right)  $ are two elements of $Y\left(  \lambda/\mu\right)  $
satisfying $i_{1}<i_{2}$ and $j_{1}\leq j_{2}$, then $T\left(  i_{1}%
,j_{1}\right)  <T\left(  i_{2},j_{2}\right)  $.
\end{lemma}

\begin{proof}
[Hints to the proof of Lemma \ref{lem.sf.skew-ssyt.increase}.]This is easy
using Lemma \ref{lem.sf.skew-diag.convexity}. A detailed proof can be found in
Section \ref{sec.details.sf.schur}.
\end{proof}

We extend Definition \ref{def.sf.ytab.xT} to skew tableaux:

\begin{definition}
\label{def.sf.ytab.skew-xT}Let $M\in\mathbb{N}$. Let $\lambda$ and $\mu$ be
two $M$-partitions. If $T$ is any Young tableau of shape $\lambda/\mu$, then
we define the corresponding monomial%
\[
x_{T}:=\prod_{c\text{ is a box of }Y\left(  \lambda/\mu\right)  }x_{T\left(
c\right)  }=\prod_{\left(  i,j\right)  \in Y\left(  \lambda/\mu\right)
}x_{T\left(  i,j\right)  }=\prod_{k=1}^{N}x_{k}^{\left(  \text{\# of times
}k\text{ appears in }T\right)  }.
\]

\end{definition}

For example,%
\[
\text{if }T=\ytableaushort{\none\none 13,\none 22,344}\ \ \text{, then }%
x_{T}=x_{1}x_{3}x_{2}x_{2}x_{3}x_{4}x_{4}=x_{1}x_{2}^{2}x_{3}^{2}x_{4}^{2}.
\]

Finally, we generalize Definition \ref{def.sf.schur} to $\lambda/\mu$:

\begin{definition}
\label{def.sf.skew-schur}Let $M\in\mathbb{N}$. Let $\lambda$ and $\mu$ be two
$M$-partitions. We define the \emph{skew Schur polynomial} $s_{\lambda/\mu}%
\in\mathcal{P}$ by%
\[
s_{\lambda/\mu}:=\sum_{T\in\operatorname*{SSYT}\left(  \lambda/\mu\right)
}x_{T}.
\]

\end{definition}

\begin{example}
\textbf{(a)} We have%
\[
s_{\left(  3,3,2,0,0,0,\ldots,0\right)  /\left(  2,1,0,0,0,\ldots,0\right)
}=\sum_{i<j\geq k<\ell\geq m}x_{i}x_{j}x_{k}x_{\ell}x_{m},
\]
since the semistandard tableaux of shape $\left(  3,3,2,0,0,0,\ldots,0\right)
/\left(  2,1,0,0,0,\ldots,0\right)  $ have the form
$\ytableaushort{\none\none i,\none kj,m\ell}$ for $i,j,k,\ell,m\in\left[
N\right]  $ satisfying $i<j\geq k<\ell\geq m$. \medskip

\textbf{(b)} We have%
\[
s_{\left(  3,2,1,0,0,0,\ldots,0\right)  /\left(  2,1,0,0,0,\ldots,0\right)
}=\sum_{i,j,k}x_{i}x_{j}x_{k},
\]
since the semistandard tableaux of shape $\left(  3,2,1,0,0,0,\ldots,0\right)
/\left(  2,1,0,0,0,\ldots,0\right)  $ have the form
$\ytableaushort{\none\none i,\none j,k}$ for $i,j,k\in\left[  N\right]  $
satisfying no special requirements (as there are never two entries in the same
row or in the same column of the tableau). Thus,%
\[
s_{\left(  3,2,1,0,0,0,\ldots,0\right)  /\left(  2,1,0,0,0,\ldots,0\right)
}=\sum_{i,j,k}x_{i}x_{j}x_{k}=\left(  \sum_{i}x_{i}\right)  ^{3}=\left(
x_{1}+x_{2}+\cdots+x_{N}\right)  ^{3}.
\]

\end{example}

\begin{theorem}
\label{thm.sf.skew-schur-symm}Let $M\in\mathbb{N}$. Let $\lambda$ and $\mu$ be
any two $M$-partitions. Then, the polynomial $s_{\lambda/\mu}$ is symmetric.
\end{theorem}

\begin{remark}
\label{rmk.sf.skew-0}If we set $M=N$ and $\mathbf{0}=\left(  0,0,\ldots
,0\right)  \in\mathbb{N}^{N}$, then $s_{\lambda/\mathbf{0}}=s_{\lambda}$
(since the semistandard tableaux of shape $\lambda/\mathbf{0}$ are precisely
the semistandard tableaux of shape $\lambda$). Hence, Theorem
\ref{thm.sf.skew-schur-symm} generalizes Theorem \ref{thm.sf.schur-symm}
\textbf{(a)}.
\end{remark}

We will now prove Theorem \ref{thm.sf.skew-schur-symm} bijectively, using a
beautiful set of combinatorial bijections known as the \emph{Bender--Knuth
involutions}.

\subsubsection{The Bender--Knuth involutions}

\begin{proof}
[Proof of Theorem \ref{thm.sf.skew-schur-symm}.]For each $i\in\left[
N-1\right]  $, we consider the simple transposition $s_{i}\in S_{N}$ defined
in Definition \ref{def.perm.si} (applied to $n=N$). As we recall, this is the
transposition that swaps $i$ with $i+1$. (The notation $s_{i}$ for this
transposition is unfortunately similar to the notation $s_{\lambda}$ for Schur
polynomials; however, this should not cause any confusion, since the only
Schur polynomial that will appear in this proof is $s_{\lambda/\mu}$, which
cannot be mistaken for a transposition.)

We must show that the polynomial $s_{\lambda/\mu}$ is symmetric. According to
Lemma \ref{lem.sf.simples-enough}, it will suffice to show that $s_{k}\cdot
s_{\lambda/\mu}=s_{\lambda/\mu}$ for each $k\in\left[  N-1\right]  $. \medskip

So let us fix some $k\in\left[  N-1\right]  $. In order to prove $s_{k}\cdot
s_{\lambda/\mu}=s_{\lambda/\mu}$, we will construct a bijection%
\[
\beta_{k}:\operatorname*{SSYT}\left(  \lambda/\mu\right)  \rightarrow
\operatorname*{SSYT}\left(  \lambda/\mu\right)
\]
(called the $k$-th \emph{Bender--Knuth involution}) that

\begin{itemize}
\item interchanges the \# of $k$'s with the \# of $\left(  k+1\right)  $'s in
a tableau (that is, if a tableau $T$ has $a$ many $k$'s and $b$ many $\left(
k+1\right)  $'s, then $\beta_{k}\left(  T\right)  $ will have $b$ many $k$'s
and $a$ many $\left(  k+1\right)  $'s);

\item leaves all the other entries of the tableau unchanged.
\end{itemize}

Indeed, once such a bijection $\beta_{k}$ is constructed, we can easily see
that
\[
x_{\beta_{k}\left(  T\right)  }=s_{k}\cdot x_{T}\ \ \ \ \ \ \ \ \ \ \text{for
each }T\in\operatorname*{SSYT}\left(  \lambda/\mu\right)  ,
\]
so that applying $s_{k}$ to $s_{\lambda/\mu}=\sum\limits_{T\in
\operatorname*{SSYT}\left(  \lambda/\mu\right)  }x_{T}$ will amount to
permuting the addends in the sum.\footnote{We will explain this argument in
more detail at the end of this proof.} \medskip

We construct the map $\beta_{k}$ as follows: Let $T\in\operatorname*{SSYT}%
\left(  \lambda/\mu\right)  $. We focus on the $k$'s and the $\left(
k+1\right)  $'s in $T$ (that is, on the entries of $T$ that are equal to $k$
or to $k+1$). An entry $k$ in $T$ will be called \emph{matched} if there is a
$k+1$ directly underneath it (in the same column). An entry $k+1$ in $T$ will
be called \emph{matched} if there is a $k$ directly above it (in the same
column). All other $k$'s and $\left(  k+1\right)  $'s in $T$ will be called
\emph{free}. Let us see an example of this first:

\begin{example}
\label{exa.pf.thm.sf.skew-schur-symm.1}For this example, let $k=2$. Here is a
semistandard tableau $T\in\operatorname*{SSYT}\left(  \left(
9,8,5,2,0,0\right)  /\left(  3,1,1,0,0,0\right)  \right)  $ with the free
entries printed in boldface and the matched entries printed on a grey
background:%
\[
T=\ytableaushort{
\none\none\none 1{*(lightgray) \textcolor{red}2}{*(lightgray) \textcolor{red}2}{\textcolor{red}{\mathbf{2}}}{\textcolor{red}{\mathbf{2}}}{\textcolor{blue}{\mathbf{3}}},
\none 11{*(lightgray) \textcolor{red}2}{*(lightgray) \textcolor{blue}3}{*(lightgray) \textcolor{blue}3}46,
\none {\textcolor{red}{\mathbf{2}}}{\textcolor{blue}{\mathbf{3}}}{*(lightgray) \textcolor{blue}3}5,
{\textcolor{red}{\mathbf{2}}}4
}\ \ .
\]
(We have color-coded the entries so that $2$'s are red, $3$'s are blue, and
all other entries are black. You can mostly forget about the black entries,
since our construction of $\beta_{k}\left(  T\right)  $ will neither change
them nor depend on them.)
\end{example}

We note that the entries of $T$ increase weakly along each row (since $T$ is
semistandard), and increase strictly down each column (for the same reason).
Now, an entry $k$ in $T$ is matched if and only if there is a $k+1$ anywhere
in its column (because if there is a $k+1$ anywhere in its column, then this
$k+1$ must be directly underneath the $k$\ \ \ \ \footnote{since the entries
of each column of $T$ increase down this column}, and therefore the $k$ is
matched). Likewise, an entry $k+1$ in $T$ is matched if and only if there is a
$k$ anywhere in its column. Thus, matched entries come in pairs: If a $k$ in
$T$ is matched, then the $k+1$ directly underneath it is also matched, and
conversely, if a $k+1$ in $T$ is matched, then the $k$ directly above it is
also matched. Hence, there is an obvious bijection between the sets $\left\{
\text{matched }k\text{'s in }T\right\}  $ and $\left\{  \text{matched }\left(
k+1\right)  \text{'s in }T\right\}  $\ \ \ \ \footnote{Strictly speaking,
these sets should consist not of the entries, but rather of the boxes in which
they are located.}. Thus, by the bijection principle, we have%
\begin{align}
&  \left(  \text{\# of matched }k\text{'s in }T\right) \nonumber\\
&  =\left(  \text{\# of matched }\left(  k+1\right)  \text{'s in }T\right)  .
\label{pf.thm.sf.skew-schur-symm.num-matched-1}%
\end{align}

Each column of $T$ that contains a matched entry must contain exactly two
matched entries (one $k$ and one $k+1$); we shall refer to these two entries
as each other's \textquotedblleft\emph{partners}\textquotedblright.

Our goal is to modify some of the entries $k$ and $k+1$ in such a way that we
obtain a new semistandard tableau that has as many $k$'s as our original
tableau $T$ had $\left(  k+1\right)  $'s, and has as many $\left(  k+1\right)
$'s as our original tableau $T$ had $k$'s. We do not want to change any
entries other than $k$'s and $\left(  k+1\right)  $'s; nor do we want to
replace any $k$'s or $\left(  k+1\right)  $'s by entries other than $k$ and
$k+1$.

These requirements force us to leave all matched entries (both $k$'s and
$\left(  k+1\right)  $'s) unchanged. Indeed, if an entry is matched, then its
column contains both a $k$ and a $k+1$, and thus neither of these two entries
can be changed without breaking the \textquotedblleft entries increase
strictly down each column\textquotedblright\ condition in Definition
\ref{def.sf.skew-ssyt}. Thus, the matched entries will have to stay unchanged.

On the other hand, we can arbitrarily replace the free entries by $k$'s or
$\left(  k+1\right)  $'s, as long as we make sure to keep the rows weakly
increasing; the columns will stay strictly increasing no matter what we do
(because a column containing a free $k$ does not contain any $k+1$, and a
column containing a free $k+1$ does not contain any $k$), so our tableau will
remain semistandard.

In view of these observations, let us perform the following procedure:

\begin{itemize}
\item For each row of $T$, if there are $a$ free $k$'s and $b$ free $\left(
k+1\right)  $'s in this row, we replace them by $b$ many $k$'s and $a$ many
$\left(  k+1\right)  $'s (placed in this order, from left to right).
\end{itemize}

We define $\beta_{k}\left(  T\right)  $ to be the result of this procedure.

\begin{example}
If $k$ and $T$ are as in Example \ref{exa.pf.thm.sf.skew-schur-symm.1}, then%
\[
\beta_{k}\left(  T\right)  =\ytableaushort{
\none\none\none 1{*(lightgray) \textcolor{red}2}{*(lightgray) \textcolor{red}2}{\textcolor{red}{\mathbf{2}}}{\textcolor{blue}{\mathbf{3}}}{\textcolor{blue}{\mathbf{3}}},
\none 11{*(lightgray) \textcolor{red}2}{*(lightgray) \textcolor{blue}3}{*(lightgray) \textcolor{blue}3}46,
\none {\textcolor{red}{\mathbf{2}}}{\textcolor{blue}{\mathbf{3}}}{*(lightgray) \textcolor{blue}3}5,
{\textcolor{blue}{\mathbf{3}}}4
}\ \ .
\]
Indeed:

\begin{itemize}
\item The $1$-st row of $T$ had two free $2$'s and one free $3$, so we
replaced them by one $2$ and two $3$'s.

\item The $2$-nd row had zero free $2$'s and zero free $3$'s, so we replaced
them by zero $2$'s and zero $3$'s. (Of course, this did not change anything.)

\item The $3$-rd row had one free $2$ and one free $3$, so we replaced them by
one $2$ and one $3$. (Of course, this did not change anything.)

\item The $4$-th row had one free $2$ and zero free $3$'s, so we replaced them
by zero $2$'s and one $3$.
\end{itemize}
\end{example}

Thus, $\beta_{k}\left(  T\right)  $ is obtained from $T$ by \textquotedblleft
flipping the imbalance between free $k$'s and free $\left(  k+1\right)
$'s\textquotedblright\ in each row of $T$ (so that a row that was heavy on
free $k$'s becomes equally heavy on free $\left(  k+1\right)  $'s, and vice
versa). Rows that have equally many free $k$'s and free $\left(  k+1\right)
$'s stay unchanged.

In order to make sure that $\beta_{k}$ is a well-defined map from
$\operatorname*{SSYT}\left(  \lambda/\mu\right)  $ to $\operatorname*{SSYT}%
\left(  \lambda/\mu\right)  $, we need to show that the tableau $\beta
_{k}\left(  T\right)  $ is semistandard. As we have already explained, the
columns are not at issue (all the entries that were changed were free, so the
columns remain strictly increasing), but we need to convince ourselves that
the rows are still weakly increasing. It is clear that the \textbf{new}
entries in each row (i.e., the entries of $\beta_{k}\left(  T\right)  $ in the
cells that used to contain the free entries of $T$) are in the right order
(i.e., any new $k$ stands further left than any new $k+1$), and it is also
clear that the \textbf{matched} entries in each row are in the right order
(since they are unchanged from $T$); however, it is imaginable that the order
between new and matched entries has gotten messed up (e.g., a larger new entry
stands further left than a smaller matched entry in $\beta_{k}\left(
T\right)  $). We need to prove that this does not happen.

To prove this, we make the following observation:

\begin{statement}
\textit{Observation 1:} Each row of $T$ can be subdivided into the following
six blocks:%
\[
\hspace{-3pc}%
\begin{tikzpicture}[scale=0.8]
\draw(0, 0) rectangle (3, 1);
\draw(3, 0) rectangle (6, 1);
\draw(6, 0) rectangle (8, 1);
\draw(8, 0) rectangle (11.3, 1);
\draw(11.3, 0) rectangle (15.8, 1);
\draw(15.8, 0) rectangle (19.3, 1);
\node at (1.5, 0.5) {\small entries $<k$};
\node at (4.5, 0.5) {\small matched $k$'s};
\node at (7, 0.5) {\small free $k$'s};
\node at (9.6, 0.5) {\small free $\left(k+1\right)$'s};
\node at (13.6, 0.5) {\small matched $\left(k+1\right)$'s};
\node at (17.6, 0.5) {\small entries $>k+1$};
\end{tikzpicture}
\ \ \ ,
\]
which appear in this order from left to right. (Each of these blocks can be empty.)
\end{statement}

[\textit{Proof of Observation 1:} Consider some row of $T$. The entries in
this row increase weakly from left to right (since $T$ is semistandard), so it
is clear that all entries $<k$ stand further left than all $k$'s, which in
turn stand further left than all $\left(  k+1\right)  $'s, which in turn stand
further left than all entries $>k+1$. It therefore remains to show that all
matched $k$'s stand further left than all free $k$'s, and that all free
$\left(  k+1\right)  $'s stand further left than all matched $\left(
k+1\right)  $'s. We will prove the first of these two statements only (since
the proof of the second is largely analogous).

So we need to show that all matched $k$'s in our row stand further left than
all free $k$'s. Assume the contrary. Thus, there is some matched $k$ that
stands further right than some free $k$. Let the former matched $k$ stand in
box $\left(  u,v\right)  $, and let the latter free $k$ stand in box $\left(
u,v^{\prime}\right)  $; thus, $v>v^{\prime}$ (since the matched $k$ stands
further right than the free $k$) and $T\left(  u,v^{\prime}\right)  =k$ (since
there is a $k$ in box $\left(  u,v^{\prime}\right)  $). Since the $k$ in box
$\left(  u,v\right)  $ is matched, there is a $k+1$ directly underneath it; in
other words, the box $\left(  u+1,v\right)  $ belongs to $Y\left(  \lambda
/\mu\right)  $ and we have $T\left(  u+1,v\right)  =k+1$. Here is an
illustration of this situation:%
\[%
\begin{tikzpicture}[scale=1.3]
\draw(0, 1) rectangle (1, 2);
\draw(4, 0) rectangle (5, 1);
\draw(4, 1) rectangle (5, 2);
\node(a2) at (0.5, 1.5) {$k$};
\node(a3) at (4.5, 0.5) {$k+1$};
\node(a4) at (4.5, 1.5) {$k$};
\node at (2.5, 1.5) {$\cdots$};
\node at (3.5, 1.5) {$\cdots$};
\node at (1.5, 1.5) {$\cdots$};
\node at (-1.5, 0.5) {$\text{row } u+1 \to$};
\node at (-1.5, 1.5) {$\text{row } u \to$};
\node at (0.5, 3) {$\substack{\displaystyle\text{column} \vspace
{0.3pc}\\ \displaystyle v' \vspace{0.3pc}\\ \displaystyle\downarrow}$};
\node at (4.5, 3) {$\substack{\displaystyle\text{column} \vspace
{0.3pc}\\ \displaystyle v \vspace{0.3pc}\\ \displaystyle\downarrow}$};
\end{tikzpicture}%
\ \ .
\]

Now, the Young diagram of $\lambda/\mu$ contains the box $\left(  u,v^{\prime
}\right)  $, but also contains the box $\left(  u+1,v\right)  $, which lies
one row south and some number of columns east of the former box (since
$v>v^{\prime}$). Hence, the Young diagram of $\lambda/\mu$ must have a box
$\left(  u+1,v^{\prime}\right)  $ directly underneath the box $\left(
u,v^{\prime}\right)  $ (that is, the box $\left(  u+1,v^{\prime}\right)  $
must belong to $Y\left(  \lambda/\mu\right)  $\ \ \ \ \footnote{Here is a
rigorous \textit{proof:} We know that $\left(  u,v^{\prime}\right)  $ and
$\left(  u+1,v\right)  $ are two elements of $Y\left(  \lambda/\mu\right)  $.
Moreover, $\left(  u+1,v^{\prime}\right)  \in\mathbb{Z}^{2}$ satisfies $u\leq
u+1\leq u+1$ and $v^{\prime}\leq v^{\prime}\leq v$ (since $v>v^{\prime}$).
Thus, Lemma \ref{lem.sf.skew-diag.convexity} (applied to $\left(  u,v^{\prime
}\right)  $, $\left(  u+1,v\right)  $ and $\left(  u+1,v^{\prime}\right)  $
instead of $\left(  a,b\right)  $, $\left(  e,f\right)  $ and $\left(
c,d\right)  $) yields $\left(  u+1,v^{\prime}\right)  \in Y\left(  \lambda
/\mu\right)  $. Qed.}). Here is an illustration of this:%
\[%
\begin{tikzpicture}[scale=1.3]
\draw(0, 0) rectangle (1, 1);
\draw(0, 1) rectangle (1, 2);
\draw(4, 0) rectangle (5, 1);
\draw(4, 1) rectangle (5, 2);
\node(a1) at (0.5, 0.5) {$?$};
\node(a2) at (0.5, 1.5) {$k$};
\node(a3) at (4.5, 0.5) {$k+1$};
\node(a4) at (4.5, 1.5) {$k$};
\node at (2.5, 0.5) {$\cdots$};
\node at (2.5, 1.5) {$\cdots$};
\node at (3.5, 0.5) {$\cdots$};
\node at (3.5, 1.5) {$\cdots$};
\node at (1.5, 0.5) {$\cdots$};
\node at (1.5, 1.5) {$\cdots$};
\node at (-1.5, 0.5) {$\text{row } u+1 \to$};
\node at (-1.5, 1.5) {$\text{row } u \to$};
\node at (0.5, 3) {$\substack{\displaystyle\text{column} \vspace
{0.3pc}\\ \displaystyle v' \vspace{0.3pc}\\ \displaystyle\downarrow}$};
\node at (4.5, 3) {$\substack{\displaystyle\text{column} \vspace
{0.3pc}\\ \displaystyle v \vspace{0.3pc}\\ \displaystyle\downarrow}$};
\end{tikzpicture}%
\ \ .
\]

Thus, we know that there is a box $\left(  u+1,v^{\prime}\right)  $ in the
Young diagram of $\lambda/\mu$. The entry of the tableau $T$ in this box
$\left(  u+1,v^{\prime}\right)  $ must satisfy $T\left(  u+1,v^{\prime
}\right)  >T\left(  u,v^{\prime}\right)  $ (because the entries of $T$
increase strictly down each column) and $T\left(  u+1,v^{\prime}\right)  \leq
T\left(  u+1,v\right)  $ (because the entries of $T$ increase weakly along
each row, and because $v>v^{\prime}$ shows that the entry $T\left(
u+1,v^{\prime}\right)  $ lies further left than $T\left(  u+1,v\right)
$\ \ \ \ \footnote{Strictly speaking, we are using Lemma
\ref{lem.sf.skew-ssyt.increase} \textbf{(a)} here (applying it to $\left(
u+1,v^{\prime}\right)  $ and $\left(  u+1,v\right)  $ instead of $\left(
i,j_{1}\right)  $ and $\left(  i,j_{2}\right)  $).}). These two inequalities
rewrite as $T\left(  u+1,v^{\prime}\right)  >k$ and $T\left(  u+1,v^{\prime
}\right)  \leq k+1$ (since $T\left(  u,v^{\prime}\right)  =k$ and $T\left(
u+1,v\right)  =k+1$). Thus, the number $T\left(  u+1,v^{\prime}\right)  $ lies
in the half-open interval $\left(  k,k+1\right]  $. Since $T\left(
u+1,v^{\prime}\right)  $ is an integer, we must thus have $T\left(
u+1,v^{\prime}\right)  =k+1$. In other words, the entry $k$ in box $\left(
u,v^{\prime}\right)  $ of our tableau $T$ has a $k+1$ directly underneath it.
Therefore, this entry $k$ is matched. This contradicts our assumption that
this $k$ is free. This contradiction shows that our assumption was false.
Thus, we have shown that all matched $k$'s in our row stand further left than
all free $k$'s. Similarly, we can show that all free $\left(  k+1\right)  $'s
stand further left than all matched $\left(  k+1\right)  $'s (the argument is
analogous, but it uses the $\left(  u-1\right)  $-st row of $T$ rather than
the $\left(  u+1\right)  $-st one). As explained above, this completes the
proof of Observation 1.] \medskip

Observation 1 entails that the free entries in each row of $T$ are stuck
together between the rightmost matched $k$ and the leftmost matched $k+1$.
Hence, replacing these free entries as in our above definition of $\beta
_{k}\left(  T\right)  $ does not mess up the weakly increasing order of the
entries in this row. This completes our proof that $\beta_{k}\left(  T\right)
$ is a semistandard tableau. \medskip

Hence, the map $\beta_{k}:\operatorname*{SSYT}\left(  \lambda/\mu\right)
\rightarrow\operatorname*{SSYT}\left(  \lambda/\mu\right)  $ is well-defined.
This map $\beta_{k}$ is called the $k$-th \emph{Bender--Knuth involution}.
\medskip

We shall now show that this map $\beta_{k}$ is a bijection. Better yet, we
will show that it is an involution (i.e., that $\beta_{k}\circ\beta
_{k}=\operatorname*{id}$):

\begin{statement}
\textit{Observation 2:} We have $\beta_{k}\circ\beta_{k}=\operatorname*{id}$.
\end{statement}

[\textit{Proof of Observation 2:} We need to check that $\beta_{k}\left(
\beta_{k}\left(  T\right)  \right)  =T$ for each $T\in\operatorname*{SSYT}%
\left(  \lambda/\mu\right)  $. So let $T\in\operatorname*{SSYT}\left(
\lambda/\mu\right)  $ be arbitrary. Recall our above definition of free and
matched $k$'s and $\left(  k+1\right)  $'s in $T$. The construction of
$\beta_{k}\left(  T\right)  $ replaced some free entries while leaving the
matched ones unchanged.

Now, we claim that the matched entries of $\beta_{k}\left(  T\right)  $ stand
in the exact same boxes as the matched entries of $T$, whereas the free
entries of $\beta_{k}\left(  T\right)  $ stand in the exact same boxes as the
free entries of $T$. Indeed, the matched entries of $T$ remain matched in
$\beta_{k}\left(  T\right)  $ (since neither these entries themselves, nor
their \textquotedblleft partners\textquotedblright\ have changed in the
construction of $\beta_{k}\left(  T\right)  $), whereas the free entries of
$T$ remain free in $\beta_{k}\left(  T\right)  $ (since the construction of
$\beta_{k}\left(  T\right)  $ cannot have produced any \textquotedblleft
partners\textquotedblright\ for them\footnote{Why not? Let's take, for
example, a free $k$. This free $k$ is the only $k$ in its column (since the
entries of $T$ increase strictly down each column), and there is no $k+1$ in
its column (since otherwise, the $k$ would be matched, not free). Thus, there
is no entry in its column that could become a \textquotedblleft
partner\textquotedblright\ for it in $\beta_{k}\left(  T\right)  $. An
analogous argument applies to a free $k+1$.}).

This has the consequence that if we apply our definition of $\beta_{k}$ to the
semistandard tableau $\beta_{k}\left(  T\right)  $ (to construct $\beta
_{k}\left(  \beta_{k}\left(  T\right)  \right)  $), then we end up undoing the
very changes that transformed $T$ into $\beta_{k}\left(  T\right)  $ (indeed,
in each row, the original imbalance between free $k$'s and free $\left(
k+1\right)  $'s that was flipped in the construction of $\beta_{k}\left(
T\right)  $ gets flipped again, and thus gets restored to its original
state\footnote{In more detail: If some row of $T$ had $a$ free $k$'s and $b$
free $\left(  k+1\right)  $'s, then the same row of $\beta_{k}\left(
T\right)  $ has $b$ free $k$'s and $a$ free $\left(  k+1\right)  $'s, and
therefore the same row of $\beta_{k}\left(  \beta_{k}\left(  T\right)
\right)  $ will, in turn, have $a$ free $k$'s and $b$ free $\left(
k+1\right)  $'s again; but this means that its free entries are the same as in
$T$.}). Hence, $\beta_{k}\left(  \beta_{k}\left(  T\right)  \right)  =T$.

Forget that we fixed $T$. We thus have shown that $\beta_{k}\left(  \beta
_{k}\left(  T\right)  \right)  =T$ for each $T\in\operatorname*{SSYT}\left(
\lambda/\mu\right)  $. In other words, $\beta_{k}\circ\beta_{k}%
=\operatorname*{id}$. This proves Observation 2.] \medskip

Observation 2 shows that $\beta_{k}$ is an involution. Hence, $\beta_{k}$ is a
bijection. \medskip

The map $\beta_{k}$ leaves all entries of the tableau other than $k$'s and
$\left(  k+1\right)  $'s unchanged (because of how we defined $\beta_{k}$).
Let us now show that $\beta_{k}$ interchanges the \# of $k$'s with the \# of
$\left(  k+1\right)  $'s in a tableau (that is, if a tableau $T$ has $a$ many
$k$'s and $b$ many $\left(  k+1\right)  $'s, then $\beta_{k}\left(  T\right)
$ will have $b$ many $k$'s and $a$ many $\left(  k+1\right)  $'s). More
precisely, we shall show the following:

\begin{statement}
\textit{Observation 3:} Let $T\in\operatorname*{SSYT}\left(  \lambda
/\mu\right)  $. Then,%
\begin{equation}
\left(  \text{\# of }k\text{'s in }\beta_{k}\left(  T\right)  \right)
=\left(  \text{\# of }\left(  k+1\right)  \text{'s in }T\right)
\label{pf.thm.sf.skew-schur-symm.num-k}%
\end{equation}
and
\begin{equation}
\left(  \text{\# of }\left(  k+1\right)  \text{'s in }\beta_{k}\left(
T\right)  \right)  =\left(  \text{\# of }k\text{'s in }T\right)  .
\label{pf.thm.sf.skew-schur-symm.num-k+1}%
\end{equation}
Moreover, if $i\in\left[  N\right]  $ satisfies $i\neq k$ and $i\neq k+1$,
then
\begin{equation}
\left(  \text{\# of }i\text{'s in }\beta_{k}\left(  T\right)  \right)
=\left(  \text{\# of }i\text{'s in }T\right)  .
\label{pf.thm.sf.skew-schur-symm.num-i}%
\end{equation}

\end{statement}

[\textit{Proof of Observation 3:} We recall a simple fact we noticed during
our proof of Observation 2 above: The matched entries of $\beta_{k}\left(
T\right)  $ stand in the exact same boxes as the matched entries of $T$,
whereas the free entries of $\beta_{k}\left(  T\right)  $ stand in the exact
same boxes as the free entries of $T$. Thus, the matched entries of $T$ remain
matched in $\beta_{k}\left(  T\right)  $, whereas the free entries of $T$
remain free in $\beta_{k}\left(  T\right)  $ (even if some of them change
their values).

Since the matched entries of $T$ remain unchanged under the map $\beta_{k}$,
we therefore have%
\begin{equation}
\left(  \text{\# of matched }k\text{'s in }\beta_{k}\left(  T\right)  \right)
=\left(  \text{\# of matched }k\text{'s in }T\right)
\ \ \ \ \ \ \ \ \ \ \label{pf.thm.sf.skew-schur-symm.num-matched-2}%
\end{equation}
and%
\begin{equation}
\left(  \text{\# of matched }\left(  k+1\right)  \text{'s in }\beta_{k}\left(
T\right)  \right)  =\left(  \text{\# of matched }\left(  k+1\right)  \text{'s
in }T\right)
.\ \ \ \ \ \ \ \ \ \ \label{pf.thm.sf.skew-schur-symm.num-matched-3}%
\end{equation}

On the other hand, the map $\beta_{k}$ flips the imbalance between free $k$'s
and free $\left(  k+1\right)  $'s in each row (but all these free entries
remain free, whereas the matched entries of $T$ remain matched in $\beta
_{k}\left(  T\right)  $); therefore, it also flips the total imbalance between
free $k$'s and free $\left(  k+1\right)  $'s in the entire
tableau\footnote{since the total \# of free $k$'s in a tableau equals the sum
of the \#s of free $k$'s in all rows (and the same holds for free $\left(
k+1\right)  $'s)}. Thus,%
\begin{equation}
\left(  \text{\# of free }k\text{'s in }\beta_{k}\left(  T\right)  \right)
=\left(  \text{\# of free }\left(  k+1\right)  \text{'s in }T\right)
\ \ \ \ \ \ \ \ \ \ \label{pf.thm.sf.skew-schur-symm.num-free-2}%
\end{equation}
and%
\begin{equation}
\left(  \text{\# of free }\left(  k+1\right)  \text{'s in }\beta_{k}\left(
T\right)  \right)  =\left(  \text{\# of free }k\text{'s in }T\right)
.\ \ \ \ \ \ \ \ \ \ \label{pf.thm.sf.skew-schur-symm.num-free-3}%
\end{equation}

Now, each $k$ in $\beta_{k}\left(  T\right)  $ is either free or matched (but
not both at the same time). Hence,%
\begin{align*}
&  \left(  \text{\# of }k\text{'s in }\beta_{k}\left(  T\right)  \right) \\
&  =\underbrace{\left(  \text{\# of free }k\text{'s in }\beta_{k}\left(
T\right)  \right)  }_{\substack{=\left(  \text{\# of free }\left(  k+1\right)
\text{'s in }T\right)  \\\text{(by (\ref{pf.thm.sf.skew-schur-symm.num-free-2}%
))}}}+\underbrace{\left(  \text{\# of matched }k\text{'s in }\beta_{k}\left(
T\right)  \right)  }_{\substack{=\left(  \text{\# of matched }k\text{'s in
}T\right)  \\\text{(by (\ref{pf.thm.sf.skew-schur-symm.num-matched-2}))}}}\\
&  =\left(  \text{\# of free }\left(  k+1\right)  \text{'s in }T\right)
+\underbrace{\left(  \text{\# of matched }k\text{'s in }T\right)
}_{\substack{=\left(  \text{\# of matched }\left(  k+1\right)  \text{'s in
}T\right)  \\\text{(by (\ref{pf.thm.sf.skew-schur-symm.num-matched-1}))}}}\\
&  =\left(  \text{\# of free }\left(  k+1\right)  \text{'s in }T\right)
+\left(  \text{\# of matched }\left(  k+1\right)  \text{'s in }T\right) \\
&  =\left(  \text{\# of }\left(  k+1\right)  \text{'s in }T\right)
\end{align*}
(since each $k+1$ in $T$ is either free or matched, but not both at the same time).

Also, each $k+1$ in $\beta_{k}\left(  T\right)  $ is either free or matched
(but not both at the same time). Hence,%
\begin{align*}
&  \left(  \text{\# of }\left(  k+1\right)  \text{'s in }\beta_{k}\left(
T\right)  \right) \\
&  =\underbrace{\left(  \text{\# of free }\left(  k+1\right)  \text{'s in
}\beta_{k}\left(  T\right)  \right)  }_{\substack{=\left(  \text{\# of free
}k\text{'s in }T\right)  \\\text{(by
(\ref{pf.thm.sf.skew-schur-symm.num-free-3}))}}}+\underbrace{\left(  \text{\#
of matched }\left(  k+1\right)  \text{'s in }\beta_{k}\left(  T\right)
\right)  }_{\substack{=\left(  \text{\# of matched }\left(  k+1\right)
\text{'s in }T\right)  \\\text{(by
(\ref{pf.thm.sf.skew-schur-symm.num-matched-3}))}}}\\
&  =\left(  \text{\# of free }k\text{'s in }T\right)  +\underbrace{\left(
\text{\# of matched }\left(  k+1\right)  \text{'s in }T\right)  }%
_{\substack{=\left(  \text{\# of matched }k\text{'s in }T\right)  \\\text{(by
(\ref{pf.thm.sf.skew-schur-symm.num-matched-1}))}}}\\
&  =\left(  \text{\# of free }k\text{'s in }T\right)  +\left(  \text{\# of
matched }k\text{'s in }T\right) \\
&  =\left(  \text{\# of }k\text{'s in }T\right)
\end{align*}
(since each $k$ in $T$ is either free or matched, but not both at the same time).

Moreover, if $i\in\left[  N\right]  $ satisfies $i\neq k$ and $i\neq k+1$,
then
\[
\left(  \text{\# of }i\text{'s in }\beta_{k}\left(  T\right)  \right)
=\left(  \text{\# of }i\text{'s in }T\right)
\]
(since the map $\beta_{k}$ leaves all $i$'s in $T$ unchanged\footnote{because
$i\neq k$ and $i\neq k+1$}, and does not replace any other entries by $i$'s).
Thus, Observation 3 is proved.] \medskip

From Observation 3, we can easily conclude the following:

\begin{statement}
\textit{Observation 4:} We have $x_{\beta_{k}\left(  T\right)  }=s_{k}\cdot
x_{T}$ for each $T\in\operatorname*{SSYT}\left(  \lambda/\mu\right)  $.
\end{statement}

\begin{fineprint}
[\textit{Proof of Observation 4:} For the sake of completeness, here is a
detailed proof. Let $T\in\operatorname*{SSYT}\left(  \lambda/\mu\right)  $.
Then,
\begin{align}
x_{T}  &  =\prod_{i=1}^{N}x_{i}^{\left(  \text{\# of times }i\text{ appears in
}T\right)  }\ \ \ \ \ \ \ \ \ \ \left(  \text{by Definition
\ref{def.sf.ytab.skew-xT}}\right) \nonumber\\
&  =\prod_{i=1}^{N}x_{i}^{\left(  \text{\# of }i\text{'s in }T\right)
}\nonumber\\
&  \ \ \ \ \ \ \ \ \ \ \ \ \ \ \ \ \ \ \ \ \left(  \text{since }\left(
\text{\# of times }i\text{ appears in }T\right)  =\left(  \text{\# of
}i\text{'s in }T\right)  \text{ for each }i\in\left[  N\right]  \right)
\nonumber\\
&  =x_{k}^{\left(  \text{\# of }k\text{'s in }T\right)  }\cdot x_{k+1}%
^{\left(  \text{\# of }\left(  k+1\right)  \text{'s in }T\right)  }\cdot
\prod_{\substack{i\in\left[  N\right]  ;\\i\neq k\text{ and }i\neq k+1}%
}x_{i}^{\left(  \text{\# of }i\text{'s in }T\right)  }
\label{pf.thm.sf.skew-schur-symm.xtx.pf.1}%
\end{align}
(here, we have split off the factors for $i=k$ and for $i=k+1$ from the
product). The same argument (applied to $\beta_{k}\left(  T\right)  $ instead
of $T$) yields%
\begin{align*}
x_{\beta_{k}\left(  T\right)  }  &  =\underbrace{x_{k}^{\left(  \text{\# of
}k\text{'s in }\beta_{k}\left(  T\right)  \right)  }}_{\substack{=x_{k}%
^{\left(  \text{\# of }\left(  k+1\right)  \text{'s in }T\right)  }\\\text{(by
(\ref{pf.thm.sf.skew-schur-symm.num-k}))}}}\cdot\underbrace{x_{k+1}^{\left(
\text{\# of }\left(  k+1\right)  \text{'s in }\beta_{k}\left(  T\right)
\right)  }}_{\substack{=x_{k+1}^{\left(  \text{\# of }k\text{'s in }T\right)
}\\\text{(by (\ref{pf.thm.sf.skew-schur-symm.num-k+1}))}}}\cdot\prod
_{\substack{i\in\left[  N\right]  ;\\i\neq k\text{ and }i\neq k+1}%
}\underbrace{x_{i}^{\left(  \text{\# of }i\text{'s in }\beta_{k}\left(
T\right)  \right)  }}_{\substack{=x_{i}^{\left(  \text{\# of }i\text{'s in
}T\right)  }\\\text{(by (\ref{pf.thm.sf.skew-schur-symm.num-i}))}}}\\
&  =x_{k}^{\left(  \text{\# of }\left(  k+1\right)  \text{'s in }T\right)
}\cdot x_{k+1}^{\left(  \text{\# of }k\text{'s in }T\right)  }\cdot
\prod_{\substack{i\in\left[  N\right]  ;\\i\neq k\text{ and }i\neq k+1}%
}x_{i}^{\left(  \text{\# of }i\text{'s in }T\right)  }\\
&  =x_{k+1}^{\left(  \text{\# of }k\text{'s in }T\right)  }\cdot
x_{k}^{\left(  \text{\# of }\left(  k+1\right)  \text{'s in }T\right)  }%
\cdot\prod_{\substack{i\in\left[  N\right]  ;\\i\neq k\text{ and }i\neq
k+1}}x_{i}^{\left(  \text{\# of }i\text{'s in }T\right)  }.
\end{align*}
On the other hand, applying the transposition $s_{k}$ (or, more precisely, the
action of this transposition $s_{k}\in S_{N}$ on the ring $\mathcal{P}$) to
both sides of the equality (\ref{pf.thm.sf.skew-schur-symm.xtx.pf.1}), we
obtain%
\begin{align*}
s_{k}\cdot x_{T}  &  =s_{k}\cdot\left(  x_{k}^{\left(  \text{\# of }k\text{'s
in }T\right)  }\cdot x_{k+1}^{\left(  \text{\# of }\left(  k+1\right)
\text{'s in }T\right)  }\cdot\prod_{\substack{i\in\left[  N\right]  ;\\i\neq
k\text{ and }i\neq k+1}}x_{i}^{\left(  \text{\# of }i\text{'s in }T\right)
}\right) \\
&  =x_{k+1}^{\left(  \text{\# of }k\text{'s in }T\right)  }\cdot
x_{k}^{\left(  \text{\# of }\left(  k+1\right)  \text{'s in }T\right)  }%
\cdot\prod_{\substack{i\in\left[  N\right]  ;\\i\neq k\text{ and }i\neq
k+1}}x_{i}^{\left(  \text{\# of }i\text{'s in }T\right)  }%
\end{align*}
(since the action of $s_{k}$ on $\mathcal{P}$ swaps the indeterminates $x_{k}$
and $x_{k+1}$ while leaving all other indeterminates $x_{i}$ unchanged).
Comparing the last two equalities, we obtain $x_{\beta_{k}\left(  T\right)
}=s_{k}\cdot x_{T}$. This proves Observation 4.] \medskip
\end{fineprint}

Now, the definition of $s_{\lambda/\mu}$ yields
\begin{equation}
s_{\lambda/\mu}=\sum\limits_{T\in\operatorname*{SSYT}\left(  \lambda
/\mu\right)  }x_{T}. \label{pf.thm.sf.skew-schur-symm.slm=}%
\end{equation}
Applying the permutation $s_{k}\in S_{N}$ (or, rather, the action of this
permutation on the ring $\mathcal{P}$) to both sides of this equality, we
obtain%
\begin{align*}
s_{k}\cdot s_{\lambda/\mu}  &  =s_{k}\cdot\sum\limits_{T\in
\operatorname*{SSYT}\left(  \lambda/\mu\right)  }x_{T}=\sum\limits_{T\in
\operatorname*{SSYT}\left(  \lambda/\mu\right)  }\underbrace{s_{k}\cdot x_{T}%
}_{\substack{=x_{\beta_{k}\left(  T\right)  }\\\text{(by Observation 4)}}}\\
&  \ \ \ \ \ \ \ \ \ \ \ \ \ \ \ \ \ \ \ \ \left(
\begin{array}
[c]{c}%
\text{since the group }S_{N}\text{ acts on the ring }\mathcal{P}\\
\text{by }K\text{-algebra automorphisms, and thus}\\
\text{the action of }s_{k}\text{ on }\mathcal{P}\text{ is }K\text{-linear}%
\end{array}
\right) \\
&  =\sum\limits_{T\in\operatorname*{SSYT}\left(  \lambda/\mu\right)  }%
x_{\beta_{k}\left(  T\right)  }=\sum\limits_{T\in\operatorname*{SSYT}\left(
\lambda/\mu\right)  }x_{T}\\
&  \ \ \ \ \ \ \ \ \ \ \ \ \ \ \ \ \ \ \ \ \left(
\begin{array}
[c]{c}%
\text{here, we have substituted }T\text{ for }\beta_{k}\left(  T\right)
\text{ in the sum,}\\
\text{since the map }\beta_{k}:\operatorname*{SSYT}\left(  \lambda/\mu\right)
\rightarrow\operatorname*{SSYT}\left(  \lambda/\mu\right) \\
\text{is a bijection}%
\end{array}
\right) \\
&  =s_{\lambda/\mu}\ \ \ \ \ \ \ \ \ \ \left(  \text{by
(\ref{pf.thm.sf.skew-schur-symm.slm=})}\right)  .
\end{align*}

Now, forget that we fixed $k$. We thus have shown that
\begin{equation}
s_{k}\cdot s_{\lambda/\mu}=s_{\lambda/\mu}\ \ \ \ \ \ \ \ \ \ \text{for each
}k\in\left[  N-1\right]  . \label{pf.thm.sf.skew-schur-symm.at}%
\end{equation}
Hence, Lemma \ref{lem.sf.simples-enough} (applied to $f=s_{\lambda/\mu}$)
shows that the polynomial $s_{\lambda/\mu}$ is symmetric. This proves Theorem
\ref{thm.sf.skew-schur-symm}.
\end{proof}

As we already mentioned, Theorem \ref{thm.sf.schur-symm} \textbf{(a)} is a
particular case of Theorem \ref{thm.sf.skew-schur-symm} (namely, the one
obtained when we set $M=N$ and $\mu=\mathbf{0}=\left(  0,0,\ldots,0\right)
$), because $s_{\lambda/\mathbf{0}}=s_{\lambda}$.

\subsubsection{The Littlewood--Richardson rule}

The Bender--Knuth involutions have served us well in the above proof of
Theorem \ref{thm.sf.skew-schur-symm}, but they have much more to offer. We
will soon see them prove one of the most famous results in the theory of
symmetric polynomials, namely the Littlewood--Richardson rule. In the process,
we will also (finally) prove Theorem \ref{thm.sf.schur-symm} \textbf{(b)}.

The Littlewood--Richardson rule has its roots in the representation theory of
the classical groups (specifically, $\operatorname*{GL}\nolimits_{N}\left(
\mathbb{C}\right)  $). We shall say a few words about this motivation before
we move on to stating the rule itself (which is purely combinatorial, as is
its proof). At the simplest level, the Littlewood--Richardson rule is about
expanding the product $s_{\nu}s_{\lambda}$ of two Schur polynomials as a sum
of other Schur polynomials. For instance, for $N=4$, we have%
\[
s_{2100}s_{1100}=s_{2111}+s_{2210}+s_{3110}+s_{3200},
\]
where we are omitting commas and parentheses for brevity (i.e., we are writing
$2100$ for the $N$-partition $\left(  2,1,0,0\right)  $, and likewise for the
other $N$-partitions). Likewise, for $N=3$, we have%
\[
s_{210}s_{210}=s_{222}+2s_{321}+s_{330}+s_{411}+s_{420}.
\]
I think it was Hermann Weyl who originally proved the existence of such an
expansion (i.e., that any product $s_{\nu}s_{\lambda}$ of two Schur
polynomials is a sum of Schur polynomials). The original proof used Lie group
representations. The idea of the proof, in a nutshell, is the following (skip
this paragraph if you are unfamiliar with representation theory): The
irreducible polynomial representations of the classical group
$\operatorname*{GL}\nolimits_{N}\left(  \mathbb{C}\right)  $ are (more or
less) in bijection with the $N$-partitions, meaning that there is an
irreducible polynomial representation $V_{\lambda}$ for each $N$-partition
$\lambda$, and all irreps (= irreducible polynomial representations) of
$\operatorname*{GL}\nolimits_{N}\left(  \mathbb{C}\right)  $ have this
form\footnote{At least if one uses the \textquotedblleft
right\textquotedblright\ definition of a polynomial representation. See
\cite[\S 5 and \S 6]{KraPro10} or \cite[\S 6.1]{Prasad-rep} for details.}.
These $V_{\lambda}$'s are known as the \emph{Weyl modules}, or in a slightly
more general form as the \emph{Schur functors}. The tensor product of two such
irreps can be decomposed as a direct sum of irreps (since polynomial
representations of $\operatorname*{GL}\nolimits_{N}\left(  \mathbb{C}\right)
$ are completely reducible):%
\[
V_{\nu}\otimes V_{\lambda}\cong\bigoplus_{\omega\text{ is an }%
N\text{-partition}}\ \ \underbrace{V_{\omega}^{c\left(  \nu,\lambda
,\omega\right)  }}_{\substack{\text{a direct sum of }c\left(  \nu
,\lambda,\omega\right)  \\\text{many }V_{\omega}\text{'s}}}.
\]
The multiplicities $c\left(  \nu,\lambda,\omega\right)  $ in this
decomposition are precisely the coefficients that you get when you decompose
the product of the Schur polynomials $s_{\nu}$ and $s_{\lambda}$ as a sum of
Schur polynomials:%
\[
s_{\nu}s_{\lambda}=\sum_{\omega\text{ is an }N\text{-partition}}c\left(
\nu,\lambda,\omega\right)  s_{\omega}.
\]
In fact, the Schur polynomials $s_{\lambda}$ are the so-called
\emph{characters} of the irreps $V_{\lambda}$, and it is known that tensor
products of representations correspond to products of their characters.

All of this, in the detail it deserves, is commonly taught in a 1st or 2nd
course on representation theory (e.g., \cite{Proces07} or \cite{EGHetc11} or
\cite[Chapter 6]{Prasad-rep}). But we are here for something else: we want to
know these $c\left(  \nu,\lambda,\omega\right)  $'s. In other words, we want a
formula that expands a product $s_{\nu}s_{\lambda}$ as a finite sum of Schur polynomials.

Such a formula was first conjectured by
\href{https://mathshistory.st-andrews.ac.uk/Biographies/Littlewood_Dudley/}{Dudley
Ernest Littlewood} and
\href{https://mathshistory.st-andrews.ac.uk/Biographies/Richardson_Archibald/}{Archibald
Read Richardson} in 1934. It remained unproven for 40 years, not least because
the statement was not very clear. In the 1970s, proofs were found
independently by
\href{https://mathshistory.st-andrews.ac.uk/Biographies/Schutzenberger/}{Marcel-Paul
Sch\"{u}tzenberger} and Glanffrwd Thomas. Since then, at least a dozen
different proofs have appeared. The proof that I will show was published by
Stembridge in 1997 (in \cite{Stembr02}, perhaps one of the most readable
papers in all of mathematics), and crystallizes decades of work by many
authors (Gasharov's somewhat similar proof \cite{Gashar98} probably being the
main harbinger). It will prove not just an expansion for $s_{\nu}s_{\lambda}$,
but also a generalization (replacing $s_{\lambda}$ by a skew Schur polynomial
$s_{\lambda/\mu}$) found by Zelevinsky in 1981 (\cite{Zelevi81}), as well as
Theorem \ref{thm.sf.schur-symm} \textbf{(b)}. My presentation of this proof
will follow \cite[\S 2.6]{GriRei} (which, in turn, elaborates on
\cite{Stembr02}).

To state the Littlewood--Richardson rule, we need some notions and notations.
We begin with the notations:

\begin{definition}
\label{def.sf.tuple-addition}\textbf{(a)} We let $\mathbf{0}$ denote the
$N$-tuple $\left(  0,0,\ldots,0\right)  \in\mathbb{N}^{N}$. \medskip

\textbf{(b)} Let $\alpha=\left(  \alpha_{1},\alpha_{2},\ldots,\alpha
_{N}\right)  $ and $\beta=\left(  \beta_{1},\beta_{2},\ldots,\beta_{N}\right)
$ be two $N$-tuples in $\mathbb{N}^{N}$. Then, we set%
\begin{align*}
\alpha+\beta &  :=\left(  \alpha_{1}+\beta_{1},\alpha_{2}+\beta_{2}%
,\ldots,\alpha_{N}+\beta_{N}\right)  \ \ \ \ \ \ \ \ \ \ \text{and}\\
\alpha-\beta &  :=\left(  \alpha_{1}-\beta_{1},\alpha_{2}-\beta_{2}%
,\ldots,\alpha_{N}-\beta_{N}\right)  .
\end{align*}
Note that $\alpha+\beta\in\mathbb{N}^{N}$, whereas $\alpha-\beta\in
\mathbb{Z}^{N}$.
\end{definition}

Of course, the addition operation $+$ that we just defined on the set
$\mathbb{N}^{N}$ is associative and commutative, and the $N$-tuple
$\mathbf{0}$ is its neutral element. The subtraction operation $-$ undoes $+$.
Note that the operation $+$ defined in Definition \ref{def.sf.tuple-addition}
is precisely the one that we used in Theorem \ref{thm.sf.schur-symm}
\textbf{(b)}. We notice that (using the notation of Definition
\ref{def.sf.sort} \textbf{(a)}) we have%
\begin{equation}
x^{\alpha}x^{\beta}=x^{\alpha+\beta} \label{eq.def.sf.tuple-addition.xab}%
\end{equation}
for any two $N$-tuples $\alpha,\beta\in\mathbb{N}^{N}$ (check this!).

Our next piece of notation is mostly a bookkeeping device:

\begin{definition}
\label{def.sf.content}Let $M\in\mathbb{N}$. Let $\lambda$ and $\mu$ be two
$M$-partitions. Let $T$ be a tableau of shape $\lambda/\mu$. We define the
\emph{content} of $T$ to be the $N$-tuple $\left(  a_{1},a_{2},\ldots
,a_{N}\right)  $, where%
\[
a_{i}:=\left(  \text{\# of }i\text{'s in }T\right)  =\left(  \text{\# of boxes
}c\text{ of }T\text{ such that }T\left(  c\right)  =i\right)  .
\]
We denote this $N$-tuple by $\operatorname*{cont}T$.
\end{definition}

For instance, if $N=5$, then $\operatorname*{cont}%
\ytableaushort{112,4}=\left(  2,1,0,1,0\right)  $.

Note that%
\begin{equation}
x_{T}=x^{\operatorname*{cont}T}\ \ \ \ \ \ \ \ \ \ \text{for any tableau }T.
\label{eq.def.sf.content.xT=}%
\end{equation}
(Indeed, both sides of this equality equal $\prod_{i=1}^{N}x_{i}^{\left(
\text{\# of }i\text{'s in }T\right)  }$.)

Another notation lets us cut certain columns out of a tableau:

\begin{definition}
\label{def.sf.col-tab}Let $M\in\mathbb{N}$. Let $\lambda$ and $\mu$ be two
$M$-partitions. Let $T$ be a tableau of shape $\lambda/\mu$. Let $j$ be a
positive integer. Then, $\operatorname{col}_{\geq j}T$ means the restriction
of $T$ to columns $j,j+1,j+2,\ldots$ (that is, the result of removing the
first $j-1$ columns from $T$). Formally speaking, this means the restriction
of the map $T$ to the set $\left\{  \left(  u,v\right)  \in Y\left(
\lambda/\mu\right)  \ \mid\ v\geq j\right\}  $.
\end{definition}

For example,%
\begin{align*}
\operatorname{col}_{\geq3}\ytableaushort{\none 112,\none 23,135,22}\ \  &
=\ \ \ \ \ \ \ \ \ \ \ \ \ \ytableaushort{12,3,5}\ \ \ \ \ \ \ \ \ \ \text{and}%
\\
& \\
\operatorname{col}_{\geq5}\ytableaushort{\none 112,\none 23,135,22}\ \  &
=\ \ \left(  \text{empty tableau}\right)  .
\end{align*}

\begin{remark}
\label{rmk.sf.col-tab.shape}What shape does the tableau $\operatorname{col}%
_{\geq j}T$ in Definition \ref{def.sf.col-tab} have?

We don't care, since we will only need this tableau for its content
$\operatorname{cont}\left(  \operatorname{col}_{\geq j}T\right)  $ (which is
defined independently of the shape). However, the answer is not hard to give:
If $\lambda=\left(  \lambda_{1},\lambda_{2},\ldots,\lambda_{M}\right)  $ and
$\mu=\left(  \mu_{1},\mu_{2},\ldots,\mu_{M}\right)  $, then
$\operatorname{col}_{\geq j}T$ is a skew Young tableau of shape
$\widetilde{\lambda}/\widetilde{\mu}$, where%
\begin{align*}
\widetilde{\lambda}  &  =\left(  \max\left\{  j-1,\lambda_{1}\right\}
,\ \max\left\{  j-1,\lambda_{2}\right\}  ,\ \ldots,\ \max\left\{
j-1,\lambda_{M}\right\}  \right)  \ \ \ \ \ \ \ \ \ \ \text{and}\\
\widetilde{\mu}  &  =\left(  \max\left\{  j-1,\mu_{1}\right\}  ,\ \max\left\{
j-1,\mu_{2}\right\}  ,\ \ldots,\ \max\left\{  j-1,\mu_{M}\right\}  \right)  .
\end{align*}
(Thus, the first $j-1$ columns of $\operatorname{col}_{\geq j}T$ are empty,
i.e., have no boxes.)
\end{remark}

Note that $\operatorname{col}_{\geq1}T=T$ for any tableau $T$.

Now we are ready to define a nontrivial notion:

\begin{definition}
\label{def.sf.yamanouchi}Let $M\in\mathbb{N}$. Let $\lambda$ and $\mu$ be two
$M$-partitions, and let $\nu$ be an $N$-partition. A semistandard tableau $T$
of shape $\lambda/\mu$ is said to be $\nu$\emph{-Yamanouchi} (this is an
adjective) if for each positive integer $j$, the $N$-tuple $\nu
+\operatorname*{cont}\left(  \operatorname{col}_{\geq j}T\right)
\in\mathbb{N}^{N}$ is an $N$-partition (i.e., weakly decreasing).
\end{definition}

This is a complex and somewhat confusing notion; before we move on, let us
thus give a metaphor that might help clarify it, and several examples.

\begin{remark}
\label{rmk.sf.yamanouchi.votes}Definition \ref{def.sf.yamanouchi} becomes
somewhat easier to conceptualize (and memorize) through a voting metaphor
(which, incidentally, is the reason why $\mathbf{0}$-Yamanouchi tableaux are
sometimes called \textquotedblleft ballot tableaux\textquotedblright):

Let $M,\lambda,\mu,\nu$ and $T$ be as in Definition \ref{def.sf.yamanouchi}.
Consider an election between $N$ candidates numbered $1,2,\ldots,N$. Regard
each entry $i$ of $T$ as a single vote for candidate $i$. Thus, for example,
the tableau $\ytableaushort{\none 12,25}$ has one vote for candidate $1$, two
votes for candidate $2$, and one for candidate $5$. Now, we count the votes by
keeping an \textquotedblleft tally board\textquotedblright, i.e., an $N$-tuple
$\left(  a_{1},a_{2},\ldots,a_{N}\right)  \in\mathbb{N}^{N}$ that records how
many votes each candidate has received (namely, candidate $i$ has received
$a_{i}$ votes). Assume that, at the beginning of our counting process, the
tally board is $\nu$ (as a consequence of ballot stuffing, or because some
votes have already been counted on the previous day). Now, we process the
votes from the tableau $T$, column by column, starting with the rightmost
column and moving left. Each time a column is processed, all the votes from
this column are simultaneously added to our tally board. Thus, after the
rightmost column is processed, our tally board is $\nu+\operatorname*{cont}%
\left(  \operatorname{col}_{\geq j}T\right)  $, where $j$ is the index of the
rightmost column (i.e., the rightmost column is the $j$-th column). Then, the
second-to-rightmost column gets processed, and the tally board becomes
$\nu+\operatorname*{cont}\left(  \operatorname{col}_{\geq j-1}T\right)  $. And
so on, until all columns have been processed.

Now, the tableau $T$ is $\nu$-Yamanouchi if and only if the tally board has
stayed weakly decreasing (i.e., candidate $1$ has at least as many votes as
candidate $2$, who in turn has at least as many votes as candidate $3$, who in
turn has at least as many votes as candidate $4$, and so on) throughout the
vote counting process. This is just a trivial restatement of the definition of
\textquotedblleft$\nu$-Yamanouchi\textquotedblright, but in my impression it
is conducive to understanding.

One takeaway from this interpretation is the following useful feature of the
vote counting process: No candidate gains more than one vote at a single time
(because no column of $T$ has two equal entries). Thus, the number of votes
for any given candidate increases only in small steps (viz., not at all or by
only $1$ vote).
\end{remark}

\begin{example}
\label{exa.sf.yamanouchi.1}\textbf{(a)} Let $N=3$ and $\nu=\mathbf{0}=\left(
0,0,0\right)  $. Which of the following six tableaux are $\mathbf{0}%
$-Yamanouchi?%
\begin{align*}
T_{1}  &  =\ytableaushort{\none 11,22}\ \ ,\qquad T_{2}%
=\ytableaushort{\none 11,23}\ \ ,\qquad T_{3}%
=\ytableaushort{\none 12,22}\ \ ,\\
& \\
T_{4}  &  =\ytableaushort{\none 1,1,2}\ \ ,\qquad T_{5}%
=\ytableaushort{\none 1,1,3}\ \ ,\qquad T_{6}%
=\ytableaushort{\none\none 11,\none 122,23}\ \ .
\end{align*}

Note that all six of these tableaux are semistandard.

Let us check whether $T_{1}$ is $\mathbf{0}$-Yamanouchi. Indeed, we compute
the $N$-tuple $\nu+\operatorname*{cont}\left(  \operatorname{col}_{\geq
j}T\right)  \in\mathbb{N}^{N}$ for each positive integer $j$, obtaining
\begin{align*}
\nu+\operatorname*{cont}\left(  \operatorname{col}_{\geq1}T_{1}\right)   &
=\mathbf{0}+\left(  2,2,0\right)  =\left(  2,2,0\right)  ;\\
\nu+\operatorname*{cont}\left(  \operatorname{col}_{\geq2}T_{1}\right)   &
=\mathbf{0}+\left(  2,1,0\right)  =\left(  2,1,0\right)  ;\\
\nu+\operatorname*{cont}\left(  \operatorname{col}_{\geq3}T_{1}\right)   &
=\mathbf{0}+\left(  1,0,0\right)  =\left(  1,0,0\right)  ;\\
\nu+\operatorname*{cont}\left(  \operatorname{col}_{\geq j}T_{1}\right)   &
=\mathbf{0}+\left(  0,0,0\right)  =\left(  0,0,0\right)  \text{ for each
}j\geq4\text{.}%
\end{align*}
All of the results $\left(  2,2,0\right)  $, $\left(  2,1,0\right)  $,
$\left(  1,0,0\right)  $ and $\left(  0,0,0\right)  $ are $N$-partitions.
Thus, $T_{1}$ is $\mathbf{0}$-Yamanouchi.

Let us check whether $T_{2}$ is $\mathbf{0}$-Yamanouchi. Indeed,%
\begin{align*}
\nu+\operatorname*{cont}\left(  \operatorname{col}_{\geq1}T_{2}\right)   &
=\mathbf{0}+\left(  2,1,1\right)  =\left(  2,1,1\right)  \text{ is an
}N\text{-partition;}\\
\nu+\operatorname*{cont}\left(  \operatorname{col}_{\geq2}T_{2}\right)   &
=\mathbf{0}+\left(  2,0,1\right)  =\left(  2,0,1\right)  \text{ is
\textbf{not} an }N\text{-partition.}%
\end{align*}
Thus, $T_{2}$ is \textbf{not} $\mathbf{0}$-Yamanouchi.

The tableau $T_{3}$ is \textbf{not} $\mathbf{0}$-Yamanouchi, since
$\nu+\operatorname*{cont}\left(  \operatorname{col}_{\geq1}T_{3}\right)
=\left(  1,3,0\right)  $ is not an $N$-partition.

The tableau $T_{4}$ is $\mathbf{0}$-Yamanouchi.

The tableau $T_{5}$ is \textbf{not} $\mathbf{0}$-Yamanouchi, since
$\nu+\operatorname*{cont}\left(  \operatorname{col}_{\geq1}T_{5}\right)
=\left(  2,0,1\right)  $ is not an $N$-partition.

The tableau $T_{6}$ is $\mathbf{0}$-Yamanouchi. \medskip

\textbf{(b)} So we know that $T_{2},T_{3},T_{5}$ are not $\mathbf{0}%
$-Yamanouchi. However, they are $\nu$-Yamanouchi for some other $N$-partitions
$\nu$. For example:

\begin{itemize}
\item The tableau $T_{2}$ becomes $\nu$-Yamanouchi for $\nu=\left(
1,1,0\right)  $.

\item The tableau $T_{3}$ becomes $\nu$-Yamanouchi for $\nu=\left(
2,0,0\right)  $.

\item The tableau $T_{5}$ becomes $\nu$-Yamanouchi for $\nu=\left(
1,1,0\right)  $.
\end{itemize}

\noindent These are, in a sense, the \textquotedblleft
minimal\textquotedblright\ choices of $\nu$ for this to happen, but of course
there are many other choices of $\nu$ that work.
\end{example}

We can now state the Littlewood--Richardson rule (in one of its many forms):

\begin{theorem}
[Zelevinsky's generalized Littlewood--Richardson rule, in Yamanouchi
form]\label{thm.sf.lr-zy}Let $M\in\mathbb{N}$. Let $\lambda$ and $\mu$ be two
$M$-partitions, and let $\nu$ be an $N$-partition. Then,%
\begin{equation}
s_{\nu}\cdot s_{\lambda/\mu}=\sum_{\substack{T\text{ is a }\nu
\text{-Yamanouchi}\\\text{semistandard tableau}\\\text{of shape }\lambda/\mu
}}s_{\nu+\operatorname*{cont}T}. \label{eq.thm.sf.lr-zy.eq}%
\end{equation}

\end{theorem}

Some comments are in order:

\begin{itemize}
\item In the sum on the right hand side of (\ref{eq.thm.sf.lr-zy.eq}), the
Schur polynomial $s_{\nu+\operatorname*{cont}T}$ is always well-defined.
Indeed, if $T$ is a $\nu$-Yamanouchi semistandard tableau of shape
$\lambda/\mu$, then $\nu+\operatorname*{cont}\underbrace{T}%
_{=\operatorname{col}_{\geq1}T}=\nu+\operatorname*{cont}\left(
\operatorname{col}_{\geq1}T\right)  $ is an $N$-partition (by the definition
of \textquotedblleft$\nu$-Yamanouchi\textquotedblright), so that
$s_{\nu+\operatorname*{cont}T}$ is a well-defined Schur polynomial.

\item Theorem \ref{thm.sf.lr-zy} expresses a product of a regular Schur
polynomial $s_{\nu}$ with a skew Schur polynomial $s_{\lambda/\mu}$ as a sum
of Schur polynomials. You can get a similar formula for the product of two
regular Schur polynomials by setting $M=N$ and $\mu=\mathbf{0}=\left(
0,0,\ldots,0\right)  $ in Theorem \ref{thm.sf.lr-zy}, so that $s_{\lambda/\mu
}$ becomes $s_{\lambda/\mathbf{0}}=s_{\lambda}$.
\end{itemize}

\begin{example}
\label{exa.sf.lr-zy.1}Let us apply Theorem \ref{thm.sf.lr-zy} to $M=3$ and
$N=3$ and $\nu=\left(  1,0,0\right)  $ and $\lambda=\left(  2,1,0\right)  $
and $\mu=\mathbf{0}=\left(  0,0,0\right)  $. Thus we get%
\begin{equation}
s_{\left(  1,0,0\right)  }\cdot s_{\left(  2,1,0\right)  }=\sum
_{\substack{T\text{ is a }\left(  1,0,0\right)  \text{-Yamanouchi}%
\\\text{semistandard tableau}\\\text{of shape }\left(  2,1,0\right)
/\mathbf{0}}}s_{\left(  1,0,0\right)  +\operatorname*{cont}T}.
\label{eq.exa.sf.lr-zy.1.1}%
\end{equation}
What are the $T$'s in the sum? The $\left(  1,0,0\right)  $-Yamanouchi
semistandard tableaux of shape $\left(  2,1,0\right)  /\mathbf{0}$ are
\[
\ytableaushort{11,2}\ \ ,\ \ \ \ \ \ \ \ \ \ \ytableaushort{12,2}\ \ ,\ \ \ \ \ \ \ \ \ \ \ytableaushort{12,3}\ \ ,
\]
and the corresponding addends of our sum are%
\begin{align*}
s_{\left(  1,0,0\right)  +\left(  2,1,0\right)  }  &  =s_{\left(
3,1,0\right)  },\\
s_{\left(  1,0,0\right)  +\left(  1,2,0\right)  }  &  =s_{\left(
2,2,0\right)  },\\
s_{\left(  1,0,0\right)  +\left(  1,1,1\right)  }  &  =s_{\left(
2,1,1\right)  }.
\end{align*}
Thus, the equality (\ref{eq.exa.sf.lr-zy.1.1}) rewrites as%
\[
s_{\left(  1,0,0\right)  }\cdot s_{\left(  2,1,0\right)  }=s_{\left(
3,1,0\right)  }+s_{\left(  2,2,0\right)  }+s_{\left(  2,1,1\right)  }.
\]
Note that $s_{\left(  1,0,0\right)  }=x_{1}+x_{2}+x_{3}$ and $s_{\left(
2,1,0\right)  }=\sum_{i\leq j\text{ and }i<k}x_{i}x_{j}x_{k}$.
\end{example}

We will prove the Littlewood--Richardson rule as a consequence of the
following lemma:

\begin{lemma}
[Stembridge's Lemma]\label{lem.sf.stemb-lem}Let $M\in\mathbb{N}$. Let
$\lambda$ and $\mu$ be two $M$-partitions, and let $\nu$ be an $N$-partition.
Then,%
\[
a_{\nu+\rho}\cdot s_{\lambda/\mu}=\sum_{\substack{T\text{ is a }%
\nu\text{-Yamanouchi}\\\text{semistandard tableau}\\\text{of shape }%
\lambda/\mu}}a_{\nu+\operatorname*{cont}T+\rho}.
\]

\end{lemma}

Before we prove this lemma, let us explore its consequences. One of them is
the Littlewood--Richardson rule; another is Theorem \ref{thm.sf.schur-symm}
\textbf{(b)}. Let us first see how the latter can be derived from the lemma.
This derivation, in turn, relies on another (simple) lemma:

\begin{lemma}
\label{lem.sf.tab-greater-i}Let $\lambda$ be any $N$-partition. Let $T$ be a
semistandard tableau of shape $\lambda$. Then, $T\left(  i,j\right)  \geq i$
for each $\left(  i,j\right)  \in Y\left(  \lambda\right)  $.
\end{lemma}

\begin{proof}
[Proof of Lemma \ref{lem.sf.tab-greater-i}.]This lemma is an easy consequence
of the fact that the entries of a semistandard tableau increase strictly down
each column. A detailed proof is given in Section \ref{sec.details.sf.schur}.
\end{proof}

\begin{proof}
[Proof of Theorem \ref{thm.sf.schur-symm} \textbf{(b)} using Lemma
\ref{lem.sf.stemb-lem}.]Set $M=N$. Recall that we have $\mathbf{0}=\left(
0,0,\ldots,0\right)  \in\mathbb{N}^{N}=\mathbb{N}^{M}$ (since $N=M$). Applying
Lemma \ref{lem.sf.stemb-lem} to $\mu=\mathbf{0}$ and $\nu=\mathbf{0}$, we
obtain%
\[
a_{\mathbf{0}+\rho}\cdot s_{\lambda/\mathbf{0}}=\sum_{\substack{T\text{ is a
}\mathbf{0}\text{-Yamanouchi}\\\text{semistandard tableau}\\\text{of shape
}\lambda/\mathbf{0}}}a_{\mathbf{0}+\operatorname*{cont}T+\rho}.
\]
This rewrites as%
\begin{equation}
a_{\rho}\cdot s_{\lambda}=\sum_{\substack{T\text{ is a }\mathbf{0}%
\text{-Yamanouchi}\\\text{semistandard tableau}\\\text{of shape }%
\lambda/\mathbf{0}}}a_{\operatorname*{cont}T+\rho}
\label{pf.thm.sf.schur-symm.b.1}%
\end{equation}
(since $\mathbf{0}+\rho=\rho$ and $s_{\lambda/\mathbf{0}}=s_{\lambda}$ and
$\mathbf{0}+\operatorname*{cont}T=\operatorname*{cont}T$).

Now, we shall analyze the sum on the right hand side. What are the
$\mathbf{0}$-Yamanouchi semistandard tableaux of shape $\lambda/\mathbf{0}$ ?
One such tableau is easy to construct: namely, the one tableau (of shape
$\lambda/\mathbf{0}$) whose all entries in the $1$-st row are $1$'s, all
entries in the $2$-nd row are $2$'s, all entries in the $3$-rd row are $3$'s,
and so on. Let us call this tableau \emph{minimalistic}, and denote it by
$T_{0}$. Formally speaking, this minimalistic tableau $T_{0}$ is defined to be
the map $Y\left(  \lambda/\mathbf{0}\right)  \rightarrow\left[  N\right]  $
that sends each $\left(  i,j\right)  \in Y\left(  \lambda/\mathbf{0}\right)  $
to $i$. For example, if $N=4$ and $\lambda=\left(  4,2,2,1\right)  $, then
this minimalistic tableau is%
\[
T_{0}=\ytableaushort{1111,22,33,4}\ \ .
\]

It turns out that this minimalistic tableau is the only $T$ on the right hand
side of (\ref{pf.thm.sf.schur-symm.b.1}). This will follow from the following
two observations:

\begin{statement}
\textit{Observation 1:} The minimalistic tableau $T_{0}$ is a $\mathbf{0}%
$-Yamanouchi semistandard tableau of shape $\lambda/\mathbf{0}$.
\end{statement}

\begin{statement}
\textit{Observation 2:} If $T$ is a $\mathbf{0}$-Yamanouchi semistandard
tableau of shape $\lambda/\mathbf{0}$, then $T=T_{0}$.
\end{statement}

\begin{fineprint}
[\textit{Proof of Observation 1:} It is clear that $T_{0}$ is a semistandard
tableau of shape $\lambda/\mathbf{0}$. Thus, we only need to show that it is
$\mathbf{0}$-Yamanouchi. In other words, we need to show that for each
positive integer $j$, the $N$-tuple $\mathbf{0}+\operatorname*{cont}\left(
\operatorname{col}_{\geq j}T_{0}\right)  \in\mathbb{N}^{N}$ is an
$N$-partition (i.e., weakly decreasing).

This can be done directly: Write $\lambda$ in the form $\lambda=\left(
\lambda_{1},\lambda_{2},\ldots,\lambda_{N}\right)  $. Thus, $\lambda_{1}%
\geq\lambda_{2}\geq\cdots\geq\lambda_{N}$ (since $\lambda$ is an
$N$-partition). Let $j$ be a positive integer. Recall that the tableau $T_{0}$
is minimalistic; hence, the restricted tableau $\operatorname{col}_{\geq
j}T_{0}$ is itself minimalistic (meaning that all its entries in the $1$-st
row are $1$'s, all entries in the $2$-nd row are $2$'s, all entries in the
$3$-rd row are $3$'s, and so on). Therefore, for each $i\in\left[  N\right]
$, we have%
\begin{align}
&  \left(  \text{\# of }i\text{'s in }\operatorname{col}_{\geq j}T_{0}\right)
\nonumber\\
&  =\left(  \text{\# of boxes in the }i\text{-th row of }\operatorname{col}%
_{\geq j}T_{0}\right) \nonumber\\
&  =\max\left\{  \lambda_{i}-\left(  j-1\right)  ,0\right\}
\label{pf.thm.sf.schur-symm.b.o1.pf.num=max}%
\end{align}
(since the $i$-th row of $T_{0}$ has $\lambda_{i}$ many boxes, and thus the
$i$-th row of $\operatorname{col}_{\geq j}T_{0}$ has $\max\left\{  \lambda
_{i}-\left(  j-1\right)  ,0\right\}  $ many boxes). Now, the $N$-tuple
\begin{align*}
&  \mathbf{0}+\operatorname*{cont}\left(  \operatorname{col}_{\geq j}%
T_{0}\right) \\
&  =\operatorname*{cont}\left(  \operatorname{col}_{\geq j}T_{0}\right) \\
&  =\left(  \text{\# of }1\text{'s in }\operatorname{col}_{\geq j}%
T_{0},\ \ \ \text{\# of }2\text{'s in }\operatorname{col}_{\geq j}%
T_{0},\ \ \ldots,\ \ \text{\# of }N\text{'s in }\operatorname{col}_{\geq
j}T_{0}\right) \\
&  =\left(  \max\left\{  \lambda_{1}-\left(  j-1\right)  ,0\right\}
,\ \ \max\left\{  \lambda_{2}-\left(  j-1\right)  ,0\right\}  ,\ \ \ldots
,\ \ \max\left\{  \lambda_{N}-\left(  j-1\right)  ,0\right\}  \right) \\
&  \ \ \ \ \ \ \ \ \ \ \ \ \ \ \ \ \ \ \ \ \left(  \text{by
(\ref{pf.thm.sf.schur-symm.b.o1.pf.num=max})}\right)
\end{align*}
is weakly decreasing (since $\lambda_{1}\geq\lambda_{2}\geq\cdots\geq
\lambda_{N}$ quickly yields $\max\left\{  \lambda_{1}-\left(  j-1\right)
,0\right\}  \geq\max\left\{  \lambda_{2}-\left(  j-1\right)  ,0\right\}
\geq\cdots\geq\max\left\{  \lambda_{N}-\left(  j-1\right)  ,0\right\}  $), and
thus is an $N$-partition. Forget that we fixed $j$. Thus, we have shown that
for each positive integer $j$, the $N$-tuple $\mathbf{0}+\operatorname*{cont}%
\left(  \operatorname{col}_{\geq j}T_{0}\right)  \in\mathbb{N}^{N}$ is an
$N$-partition. This proves that the tableau $T_{0}$ is $\mathbf{0}%
$-Yamanouchi. This completes the proof of Observation 1.] \medskip
\end{fineprint}

\begin{fineprint}
[\textit{Proof of Observation 2:} Let $T$ be a $\mathbf{0}$-Yamanouchi
semistandard tableau of shape $\lambda/\mathbf{0}$. We must prove that
$T=T_{0}$.

Let us assume the contrary. Thus, $T\neq T_{0}$. Hence, there exists some
$\left(  i,j\right)  \in Y\left(  \lambda/\mathbf{0}\right)  $ satisfying
$T\left(  i,j\right)  \neq T_{0}\left(  i,j\right)  $. Choose such an $\left(
i,j\right)  $ with maximum possible $j$. More precisely, among all such pairs
$\left(  i,j\right)  $ with maximum possible $j$, we choose one with the
minimum possible $i$.

Thus, for each $\left(  i^{\prime},j^{\prime}\right)  \in Y\left(
\lambda/\mathbf{0}\right)  $ satisfying $j^{\prime}>j$, we have%
\begin{equation}
T\left(  i^{\prime},j^{\prime}\right)  =T_{0}\left(  i^{\prime},j^{\prime
}\right)  \label{pf.thm.sf.schur-symm.b.o2.pf.maxj}%
\end{equation}
(since we have chosen $\left(  i,j\right)  $ to have maximum possible $j$
among the pairs satisfying $T\left(  i,j\right)  \neq T_{0}\left(  i,j\right)
$). Furthermore, for each $\left(  i^{\prime},j^{\prime}\right)  \in Y\left(
\lambda/\mathbf{0}\right)  $ satisfying $i^{\prime}<i$ and $j^{\prime}=j$, we
have%
\begin{equation}
T\left(  i^{\prime},j^{\prime}\right)  =T_{0}\left(  i^{\prime},j^{\prime
}\right)  \label{pf.thm.sf.schur-symm.b.o2.pf.mini}%
\end{equation}
(since we have chosen $\left(  i,j\right)  $ to have minimum possible $i$
among the maximum-$j$ pairs satisfying $T\left(  i,j\right)  \neq T_{0}\left(
i,j\right)  $).

The definition of the minimalistic tableau $T_{0}$ yields $T_{0}\left(
i,j\right)  =i$. Set $p:=T\left(  i,j\right)  $. Hence, $p=T\left(
i,j\right)  \neq T_{0}\left(  i,j\right)  =i$.\ \ \ \ \footnote{Here is an
example of what our tableau $T$ can look like at this point (for $N=6$ and
$\lambda=\left(  6,5,5,2,2,1\right)  $ and $\left(  i,j\right)  =\left(
3,2\right)  $):%
\[
\ytableaushort{?11111,?2222,?p333,??,??,?}\ \ .
\]
Here, the known entries come from (\ref{pf.thm.sf.schur-symm.b.o2.pf.maxj})
and (\ref{pf.thm.sf.schur-symm.b.o2.pf.mini}) (since the definition of the
minimalistic tableau $T_{0}$ shows that $T_{0}\left(  i^{\prime},j^{\prime
}\right)  =i^{\prime}$ for each $\left(  i^{\prime},j^{\prime}\right)  \in
Y\left(  \lambda/\mathbf{0}\right)  $).}

The number $p$ appears at least once in the $j$-th column of $T$ (since
$p=T\left(  i,j\right)  $), and thus appears at least once in the restricted
tableau $\operatorname{col}_{\geq j}T$ (since this restricted tableau contains
the $j$-th column of $T$).

The definition of $Y\left(  \lambda/\mathbf{0}\right)  $ yields $Y\left(
\lambda/\mathbf{0}\right)  =Y\left(  \lambda\right)  \setminus
\underbrace{Y\left(  \mathbf{0}\right)  }_{=\varnothing}=Y\left(
\lambda\right)  $. Hence, a tableau of shape $\lambda/\mathbf{0}$ is the same
as a tableau of shape $\lambda$. Thus, $T$ is a tableau of shape $\lambda$
(since $T$ is a tableau of shape $\lambda/\mathbf{0}$). Since $T$ is
semistandard, we can thus apply Lemma \ref{lem.sf.tab-greater-i}, and conclude
that $T\left(  i,j\right)  \geq i$. Hence, $p=T\left(  i,j\right)  \geq i$.
Combining this with $p\neq i$, we obtain $p>i$. In other words, $i<p$.

Now, recall that $T$ is $\mathbf{0}$-Yamanouchi; hence, $\mathbf{0}%
+\operatorname*{cont}\left(  \operatorname{col}_{\geq j}T\right)  $ is an
$N$-partition (by the definition of \textquotedblleft$\mathbf{0}%
$-Yamanouchi\textquotedblright). In other words, $\operatorname*{cont}\left(
\operatorname{col}_{\geq j}T\right)  $ is an $N$-partition (since
$\mathbf{0}+\operatorname*{cont}\left(  \operatorname{col}_{\geq j}T\right)
=\operatorname*{cont}\left(  \operatorname{col}_{\geq j}T\right)  $). Write
this $N$-partition $\operatorname*{cont}\left(  \operatorname{col}_{\geq
j}T\right)  $ as $\left(  a_{1},a_{2},\ldots,a_{N}\right)  $. For each
$k\in\left[  N\right]  $, its entry $a_{k}$ is the \# of $k$'s in
$\operatorname{col}_{\geq j}T$ (by the definition of $\operatorname*{cont}%
\left(  \operatorname{col}_{\geq j}T\right)  $). Applying this to $k=i$, we
see that $a_{i}$ is the \# of $i$'s in $\operatorname{col}_{\geq j}T$.

Similarly, $a_{p}$ is the \# of $p$'s in $\operatorname{col}_{\geq j}T$.
Hence, $a_{p}\geq1$ (since we know that the number $p$ appears at least once
in the restricted tableau $\operatorname{col}_{\geq j}T$). However, $a_{1}\geq
a_{2}\geq\cdots\geq a_{N}$ (since $\left(  a_{1},a_{2},\ldots,a_{N}\right)  $
is an $N$-partition), and thus $a_{i}\geq a_{p}$ (since $i<p$). Hence,
$a_{i}\geq a_{p}\geq1$. In other words, the number $i$ appears at least once
in the restricted tableau $\operatorname{col}_{\geq j}T$ (since $a_{i}$ is the
\# of $i$'s in $\operatorname{col}_{\geq j}T$). In other words, the number $i$
appears at least once in one of the columns $j,j+1,j+2,\ldots$ of the tableau
$T$. In other words, there exists some $\left(  i^{\prime},j^{\prime}\right)
\in Y\left(  \lambda/\mathbf{0}\right)  $ satisfying $j^{\prime}\geq j$ and
$T\left(  i^{\prime},j^{\prime}\right)  =i$. Consider this $\left(  i^{\prime
},j^{\prime}\right)  $.

Let us first assume (for the sake of contradiction) that $j^{\prime}>j$. Thus,
(\ref{pf.thm.sf.schur-symm.b.o2.pf.maxj}) yields $T\left(  i^{\prime
},j^{\prime}\right)  =T_{0}\left(  i^{\prime},j^{\prime}\right)  =i^{\prime}$
(by the definition of the minimalistic tableau $T_{0}$). Therefore,
$i^{\prime}=T\left(  i^{\prime},j^{\prime}\right)  =i$. Hence, we can rewrite
$\left(  i^{\prime},j^{\prime}\right)  \in Y\left(  \lambda/\mathbf{0}\right)
$ and $T\left(  i^{\prime},j^{\prime}\right)  =i$ as $\left(  i,j^{\prime
}\right)  \in Y\left(  \lambda/\mathbf{0}\right)  $ and $T\left(  i,j^{\prime
}\right)  =i$. Also, $j<j^{\prime}$ (since $j^{\prime}>j$). However, the
tableau $T$ is semistandard; thus, its entries increase weakly along each row.
Therefore, from $j<j^{\prime}$, we obtain $T\left(  i,j\right)  \leq T\left(
i,j^{\prime}\right)  $\ \ \ \ \footnote{Strictly speaking, this follows by
applying Lemma \ref{lem.sf.skew-ssyt.increase} \textbf{(a)} to $\left(
i,j\right)  $ and $\left(  i,j^{\prime}\right)  $ instead of $\left(
i,j_{1}\right)  $ and $\left(  i,j_{2}\right)  $.}. Thus, $p=T\left(
i,j\right)  \leq T\left(  i,j^{\prime}\right)  =i$. But this contradicts $p>i$.

This contradiction shows that our assumption (that $j^{\prime}>j$) was false.
Hence, we must have $j^{\prime}\leq j$. Combined with $j^{\prime}\geq j$, this
yields $j^{\prime}=j$. Thus, we can rewrite $\left(  i^{\prime},j^{\prime
}\right)  \in Y\left(  \lambda/\mathbf{0}\right)  $ and $T\left(  i^{\prime
},j^{\prime}\right)  =i$ as $\left(  i^{\prime},j\right)  \in Y\left(
\lambda/\mathbf{0}\right)  $ and $T\left(  i^{\prime},j\right)  =i$.

We assume (for the sake of contradiction) that $i^{\prime}<i$. Hence,
(\ref{pf.thm.sf.schur-symm.b.o2.pf.mini}) yields $T\left(  i^{\prime
},j^{\prime}\right)  =T_{0}\left(  i^{\prime},j^{\prime}\right)  =i^{\prime}$
(by the definition of the minimalistic tableau $T_{0}$), so that $i^{\prime
}=T\left(  i^{\prime},j^{\prime}\right)  =i$; but this contradicts $i^{\prime
}<i$.

This contradiction shows that our assumption (that $i^{\prime}<i$) was false.
Hence, we must have $i^{\prime}\geq i$. In other words, $i\leq i^{\prime}$.
However, the tableau $T$ is semistandard; thus, its entries increase strictly
down each column. Therefore, from $i\leq i^{\prime}$, we obtain $T\left(
i,j\right)  \leq T\left(  i^{\prime},j\right)  $\ \ \ \ \footnote{Strictly
speaking, this follows by applying Lemma \ref{lem.sf.skew-ssyt.increase}
\textbf{(b)} to $\left(  i,j\right)  $ and $\left(  i^{\prime},j\right)  $
instead of $\left(  i_{1},j\right)  $ and $\left(  i_{2},j\right)  $.}. Thus,
$T\left(  i^{\prime},j\right)  \geq T\left(  i,j\right)  =p>i$, so that
$i<T\left(  i^{\prime},j\right)  =i$. Thus we have obtained a contradiction
again. This contradiction shows that our assumption was false; hence,
$T=T_{0}$. This proves Observation 2.] \medskip
\end{fineprint}

\begin{noncompile}
[\textit{Old (not quite formal and not very readable) proof of Observation 2:}
Let $T$ be any $\mathbf{0}$-Yamanouchi semistandard tableau of shape
$\lambda/\mathbf{0}$. We must prove that $T=T_{0}$.

Let $c$ be the \# of (nonempty) columns of $T$. (Note that $c=\lambda_{1}$ if
$N>0$.)

We now claim that for each $m\in\left\{  0,1,\ldots,c\right\}  $,
\begin{equation}
\text{the }m\text{ rightmost columns of }T\text{ agree with }T_{0}
\label{pf.thm.sf.schur-symm.b.o1.pf.1}%
\end{equation}
(i.e., the entries in these columns are equal to the corresponding entries of
$T_{0}$).

[\textit{Proof of (\ref{pf.thm.sf.schur-symm.b.o1.pf.1}):} We shall prove this
by induction on $m$. The \textit{induction base} (i.e., the case $m=0$) is
obvious (indeed, the $0$ rightmost columns of $T$ have no entries, and thus
agree with $T_{0}$).

For the \textit{induction step}, we fix some $k\in\left[  c\right]  $, and we
assume (as the induction hypothesis) that
(\ref{pf.thm.sf.schur-symm.b.o1.pf.1}) holds for $m=k-1$ (that is, the $k-1$
rightmost columns of $T$ agree with $T_{0}$). We must then show that
(\ref{pf.thm.sf.schur-symm.b.o1.pf.1}) also holds for $m=k$ (that is, we must
show that the $k$ rightmost columns of $T$ agree with $T_{0}$). It clearly
suffices to prove that the $k$-th column of $T$ from the right has the same
entries as the corresponding column of $T_{0}$ (since the induction hypothesis
guarantees that the same holds for all columns to the right of this column).

We shall prove this by example: Here is how $T$ can look like (for $N=6$ and
$\lambda=\left(  6,5,5,2,2,1\right)  $ and $c=6$ and $k=5$):%
\begin{equation}
\ytableaushort{\ast ?1111,\ast ?222,\ast ?333,\ast ?,\ast ?,\ast}\ \ ,
\label{pf.thm.sf.schur-symm.b.o1.pf.1.pf.2}%
\end{equation}
where the question marks stand for the entries in the $k$-th column of $T$
from the right, whereas the asterisks stand for the entries to its left (which
we are not currently interested in). We want to prove that the $k$-th column
of $T$ from the right has the same entries as the corresponding column of
$T_{0}$, which are $1,2,3,4,5$ from top to bottom. In other words, we need to
prove that the question marks in (\ref{pf.thm.sf.schur-symm.b.o1.pf.1.pf.2})
are $1,2,3,4,5$ from top to bottom. For the topmost three question marks, this
follows from the semistandardness of $T$: Indeed, the topmost question mark is
$1$ (because it has to be $\leq$ to the entry directly to its right, but that
entry is a $1$). Therefore, the second question mark from the top is $2$
(since it has to be greater than the $1$ above it, but $\leq$ to the $2$
directly to its right). Thus, the third question mark from the top is $3$
(since it has to be greater than the $2$ above it, but $\leq$ to the $3$
directly to its right). It remains to show that the remaining two question
marks are $4$ and $5$ (from top to bottom). Here, we have to use the
$\mathbf{0}$-Yamanouchi condition. Indeed, let us refer to the $k$-th column
of $T$ from the right as the \emph{questionable column} (since it is the
column with the question marks). Let $j\in\left[  c\right]  $ be such that the
questionable column is the $j$-th column of $T$ from the left. Now, recall
that $T$ is $\mathbf{0}$-Yamanouchi; hence, $\mathbf{0}+\operatorname*{cont}%
\left(  \operatorname{col}_{\geq j}T\right)  $ is an $N$-partition (by the
definition of \textquotedblleft$\mathbf{0}$-Yamanouchi\textquotedblright). In
other words, $\operatorname*{cont}\left(  \operatorname{col}_{\geq j}T\right)
$ is an $N$-partition (since $\mathbf{0}+\operatorname*{cont}\left(
\operatorname{col}_{\geq j}T\right)  =\operatorname*{cont}\left(
\operatorname{col}_{\geq j}T\right)  $). Write this $N$-partition
$\operatorname*{cont}\left(  \operatorname{col}_{\geq j}T\right)  $ as
$\left(  a_{1},a_{2},\ldots,a_{N}\right)  $. For each $i\in\left[  N\right]
$, its entry $a_{i}$ is the \# of $i$'s in the $j$-th column of $T$ (that is,
the questionable column) and the columns further right. We have $a_{1}\geq
a_{2}\geq\cdots\geq a_{N}$ (since $\left(  a_{1},a_{2},\ldots,a_{N}\right)  $
is an $N$-partition). Now, we can identify the question mark in the $4$-th
row: Indeed, let $p$ be this question mark (or, to be more precise, the entry
that it stands for). Then, $p>3$ (since it must be greater than the $3$ above
it). Moreover, $a_{p}>0$ (since there is a $p$ in the questionable column --
namely, the question mark we are currently discussing). If we had $p>4$, then
the questionable column of $T$ would not contain a $4$ at all, and therefore
we would have $a_{4}=0$ (since the columns to the right of this column don't
contain any $4$ either) and thus $a_{p}>0=a_{4}$, which would contradict $p>4$
(since $a_{1}\geq a_{2}\geq\cdots\geq a_{N}$). Thus, we cannot have $p>4$.
Hence, we must have $p\leq4$ and therefore $p=4$ (since $p>3$). Thus, we have
shown that the question mark in the $4$-th row must be equal to $4$.
Similarly, the question mark directly below it must be equal to $5$. Thus, we
have proved that the five question marks in the questionable column must be
$1,2,3,4,5$ from top to bottom. Hence, the questionable column (i.e., the
$k$-th column of $T$ from the right) has the same entries as the corresponding
column of $T_{0}$. This completes the induction step. Thus,
(\ref{pf.thm.sf.schur-symm.b.o1.pf.1}) is proved.]

Applying (\ref{pf.thm.sf.schur-symm.b.o1.pf.1}) to $m=c$, we conclude that the
$c$ rightmost columns of $T$ agree with $T_{0}$. In other words, $T$ agrees
with $T_{0}$ entirely (since $T$ has only $c$ columns). In other words,
$T=T_{0}$. Thus we have proved Observation 2.] \medskip
\end{noncompile}

Combining Observation 1 with Observation 2, we see that the minimalistic
tableau $T_{0}$ is the \textbf{only} $\mathbf{0}$-Yamanouchi semistandard
tableau of shape $\lambda/\mathbf{0}$. Hence, the sum on the right hand side
of (\ref{pf.thm.sf.schur-symm.b.1}) has only one addend, namely the addend for
$T=T_{0}$. Thus, (\ref{pf.thm.sf.schur-symm.b.1}) simplifies to
\[
a_{\rho}\cdot s_{\lambda}=a_{\operatorname*{cont}\left(  T_{0}\right)  +\rho
}=a_{\lambda+\rho},
\]
since it is easy to see that $\operatorname*{cont}\left(  T_{0}\right)
=\lambda$. This proves Theorem \ref{thm.sf.schur-symm} \textbf{(b)} (using
Lemma \ref{lem.sf.stemb-lem}).
\end{proof}

Let us furthermore derive Theorem \ref{thm.sf.lr-zy} from Lemma
\ref{lem.sf.stemb-lem}. This relies on some elementary properties of certain
polynomials. The underlying notion is defined in an arbitrary commutative ring:

\begin{definition}
\label{def.cring.reg}Let $L$ be a commutative ring. Let $a\in L$. The element
$a$ of $L$ is said to be \emph{regular} if and only if every $x\in L$
satisfying $ax=0$ satisfies $x=0$.
\end{definition}

Regular elements of a commutative ring are often called \textquotedblleft%
\emph{non-zero-divisors}\textquotedblright\footnote{This name is somewhat
murky in the literature (and is best avoided). In fact, many authors prefer to
consider $0$ to be a non-zero-divisor as well (so that they can say that an
integral domain has no zero-divisors, rather than saying that the only
zero-divisor in an integral domain is $0$), even though $0$ is not a regular
element (unless the ring $L$ is trivial). This exception tends to make the
notion of a non-zero-divisor fickle and unreliable.} or \emph{cancellable}
elements. The latter word is explained by the following simple fact:

\begin{lemma}
\label{lem.cring.reg.cancel}Let $L$ be a commutative ring. Let $a,u,v\in L$ be
such that $a$ is regular. Assume that $au=av$. Then, $u=v$.
\end{lemma}

\begin{proof}
[Proof of Lemma \ref{lem.cring.reg.cancel}.]We have $a\left(  u-v\right)
=au-av=0$ (since $au=av$). However, $a$ is regular; in other words, every
$x\in L$ satisfying $ax=0$ satisfies $x=0$ (by the definition of
\textquotedblleft regular\textquotedblright). Applying this to $x=u-v$, we
obtain $u-v=0$ (since $a\left(  u-v\right)  =0$). Thus, $u=v$. This proves
Lemma \ref{lem.cring.reg.cancel}.
\end{proof}

Lemma \ref{lem.cring.reg.cancel} shows that regular elements of a commutative
ring can be cancelled when they appear as factors on both sides of an
equality. To make use of this, we need to actually find nontrivial regular
elements. Here is one:

\begin{lemma}
\label{lem.sf.arho-reg}The element $a_{\rho}$ of the polynomial ring
$\mathcal{P}$ is regular.
\end{lemma}

\begin{proof}
[Proof of Lemma \ref{lem.sf.arho-reg} (sketched).]There are different ways to
prove this. One is to define a lexicographic order on the monomials in
$\mathcal{P}$, and to argue that the leading coefficient of $a_{\rho}$ with
respect to this order is $1$ (which is a regular element of $K$); this uses a
bit of multivariate polynomial theory (see \cite[Theorem 40.7]{Warner90} or
\cite[Proposition 6.2.10]{23wa} for the properties of polynomials that are
used here).\footnote{This argument can in fact be used to show a more general
statement: Namely, for any $N$-tuple $\alpha=\left(  \alpha_{1},\alpha
_{2},\ldots,\alpha_{N}\right)  \in\mathbb{N}^{N}$ satisfying $\alpha
_{1}>\alpha_{2}>\cdots>\alpha_{N}$, the alternant $a_{\alpha}$ is regular in
$\mathcal{P}$. (But we won't need this statement.)}

Here is a more elementary proof. First, we observe that each of the
indeterminates $x_{1},x_{2},\ldots,x_{N}$ is regular (as an element of
$\mathcal{P}$). Indeed, multiplying a polynomial $f$ by an indeterminate
$x_{i}$ merely shifts the coefficients of $f$ to different monomials; thus, if
$x_{i}f=0$, then $f=0$. Next, we conclude that the polynomial $x_{i}-x_{j}%
\in\mathcal{P}$ is regular whenever $1\leq i<j\leq N$. Indeed, this polynomial
$x_{i}-x_{j}$ is the image of the indeterminate $x_{i}$ under a certain
$K$-algebra automorphism of $\mathcal{P}$ (namely, under the automorphism that
sends $x_{i}$ to $x_{i}-x_{j}$ while leaving all other indeterminates
unchanged\footnote{This is an automorphism, because its inverse is easily
constructed (namely, it sends $x_{i}$ to $x_{i}+x_{j}$ while leaving all other
indeterminates unchanged).}), and therefore is regular because $x_{i}$ is
regular (and because any ring automorphism sends regular elements to regular
elements). However, it is easy to see (see \cite[Proposition 2.3]{regpol} for
a proof) that any finite product of regular elements is again regular. Thus,
the element $\prod_{1\leq i<j\leq N}\left(  x_{i}-x_{j}\right)  \in
\mathcal{P}$ is regular (since it is the product of the regular elements
$x_{i}-x_{j}$ for $1\leq i<j\leq N$). In view of
(\ref{eq.def.sf.alternants.arho=vdm}), this rewrites as follows: The element
$a_{\rho}\in\mathcal{P}$ is regular. This proves Lemma \ref{lem.sf.arho-reg}.
\end{proof}

We can now derive Theorem \ref{thm.sf.lr-zy} from Lemma \ref{lem.sf.stemb-lem}:

\begin{proof}
[Proof of Theorem \ref{thm.sf.lr-zy} using Lemma \ref{lem.sf.stemb-lem}%
.]Theorem \ref{thm.sf.schur-symm} \textbf{(b)} (applied to $\nu$ instead of
$\lambda$) tells us that $a_{\nu+\rho}=a_{\rho}\cdot s_{\nu}$. However, Lemma
\ref{lem.sf.stemb-lem} says that%
\begin{align*}
a_{\nu+\rho}\cdot s_{\lambda/\mu}  &  =\sum_{\substack{T\text{ is a }%
\nu\text{-Yamanouchi}\\\text{semistandard tableau}\\\text{of shape }%
\lambda/\mu}}\underbrace{a_{\nu+\operatorname*{cont}T+\rho}}%
_{\substack{=a_{\rho}\cdot s_{\nu+\operatorname*{cont}T}\\\text{(by Theorem
\ref{thm.sf.schur-symm} \textbf{(b)}}\\\text{(applied to }\nu
+\operatorname*{cont}T\text{ instead of }\lambda\text{),}\\\text{since }%
\nu+\operatorname*{cont}T\text{ is an }N\text{-partition}\\\text{(as we have
shown in a comment}\\\text{after Theorem \ref{thm.sf.lr-zy}))}}}\\
&  =\sum_{\substack{T\text{ is a }\nu\text{-Yamanouchi}\\\text{semistandard
tableau}\\\text{of shape }\lambda/\mu}}a_{\rho}\cdot s_{\nu
+\operatorname*{cont}T}=a_{\rho}\cdot\sum_{\substack{T\text{ is a }%
\nu\text{-Yamanouchi}\\\text{semistandard tableau}\\\text{of shape }%
\lambda/\mu}}s_{\nu+\operatorname*{cont}T}.
\end{align*}
In view of $a_{\nu+\rho}=a_{\rho}\cdot s_{\nu}$, this equality rewrites as%
\[
a_{\rho}\cdot s_{\nu}\cdot s_{\lambda/\mu}=a_{\rho}\cdot\sum
_{\substack{T\text{ is a }\nu\text{-Yamanouchi}\\\text{semistandard
tableau}\\\text{of shape }\lambda/\mu}}s_{\nu+\operatorname*{cont}T}.
\]
Since the element $a_{\rho}\in\mathcal{P}$ is regular (by Lemma
\ref{lem.sf.arho-reg}), we can cancel $a_{\rho}$ from this equality (i.e., we
can apply Lemma \ref{lem.cring.reg.cancel} to $L=\mathcal{P}$ and $a=a_{\rho}$
and $u=s_{\nu}\cdot s_{\lambda/\mu}$ and $v=\sum_{\substack{T\text{ is a }%
\nu\text{-Yamanouchi}\\\text{semistandard tableau}\\\text{of shape }%
\lambda/\mu}}s_{\nu+\operatorname*{cont}T}$). As a result, we obtain%
\[
s_{\nu}\cdot s_{\lambda/\mu}=\sum_{\substack{T\text{ is a }\nu
\text{-Yamanouchi}\\\text{semistandard tableau}\\\text{of shape }\lambda/\mu
}}s_{\nu+\operatorname*{cont}T}.
\]
Thus, Theorem \ref{thm.sf.lr-zy} is proven.
\end{proof}

Thus it remains to prove Stembridge's lemma. Before we do so, let us spell out
two simple properties of alternants that will be used in the proof:

\begin{lemma}
\label{lem.sf.alternant-0}Let $\alpha\in\mathbb{N}^{N}$. \medskip

\textbf{(a)} If the $N$-tuple $\alpha$ has two equal entries, then $a_{\alpha
}=0$. \medskip

\textbf{(b)} Let $\beta\in\mathbb{N}^{N}$ be an $N$-tuple obtained from
$\alpha$ by swapping two entries. Then, $a_{\beta}=-a_{\alpha}$.
\end{lemma}

\begin{proof}
[Proof of Lemma \ref{lem.sf.alternant-0}.]This is an easy consequence of
Definition \ref{def.sf.alternants} \textbf{(b)}. See Section
\ref{sec.details.sf.schur} for a detailed proof.
\end{proof}

\begin{proof}
[Proof of Lemma \ref{lem.sf.stemb-lem}.]For any $\beta\in\mathbb{N}^{N}$ and
any $i\in\left[  N\right]  $, we let $\beta_{i}$ denote the $i$-th entry of
$\beta$. Thus, for example, $\rho_{k}=N-k$ for each $k\in\left[  N\right]  $
(since $\rho=\left(  N-1,N-2,\ldots,N-N\right)  $).

Since addition on $\mathbb{N}^{N}$ is defined entrywise, we have $\left(
\beta+\gamma\right)  _{i}=\beta_{i}+\gamma_{i}$ for any $\beta,\gamma
\in\mathbb{N}^{N}$ and $i\in\left[  N\right]  $.

The group $S_{N}$ acts on $\mathcal{P}$ by $K$-algebra automorphisms. Hence,
in particular, we have%
\[
\sigma\cdot\left(  fg\right)  =\left(  \sigma\cdot f\right)  \cdot\left(
\sigma\cdot g\right)
\]
for any $\sigma\in S_{N}$ and any $f,g\in\mathcal{P}$. In other words, we have%
\begin{equation}
\left(  \sigma\cdot f\right)  \cdot\left(  \sigma\cdot g\right)  =\sigma
\cdot\left(  fg\right)  \label{pf.lem.sf.stemb-lem.sigfg=}%
\end{equation}
for any $\sigma\in S_{N}$ and any $f,g\in\mathcal{P}$.

The polynomial $s_{\lambda/\mu}$ is symmetric (by Theorem
\ref{thm.sf.skew-schur-symm}). The definition of $s_{\lambda/\mu}$ yields%
\begin{equation}
s_{\lambda/\mu}=\sum_{T\in\operatorname*{SSYT}\left(  \lambda/\mu\right)
}\underbrace{x_{T}}_{\substack{=x^{\operatorname*{cont}T}\\\text{(by
(\ref{eq.def.sf.content.xT=}))}}}=\sum_{T\in\operatorname*{SSYT}\left(
\lambda/\mu\right)  }x^{\operatorname*{cont}T}.
\label{pf.lem.sf.stemb-lem.slm=}%
\end{equation}

For any $\beta\in\mathbb{N}^{N}$, we have%
\begin{align}
a_{\beta}  &  =\det\left(  \left(  x_{i}^{\beta_{j}}\right)  _{1\leq i\leq
N,\ 1\leq j\leq N}\right)  \ \ \ \ \ \ \ \ \ \ \left(  \text{by the definition
of the alternant }a_{\beta}\right) \nonumber\\
&  =\det\left(  \left(  x_{j}^{\beta_{i}}\right)  _{1\leq i\leq N,\ 1\leq
j\leq N}\right) \nonumber\\
&  \ \ \ \ \ \ \ \ \ \ \ \ \ \ \ \ \ \ \ \ \left(
\begin{array}
[c]{c}%
\text{by Theorem \ref{thm.det.transp},}\\
\text{since }\left(  x_{i}^{\beta_{j}}\right)  _{1\leq i\leq N,\ 1\leq j\leq
N}=\left(  \left(  x_{j}^{\beta_{i}}\right)  _{1\leq i\leq N,\ 1\leq j\leq
N}\right)  ^{T}%
\end{array}
\right) \nonumber\\
&  =\sum_{\sigma\in S_{N}}\left(  -1\right)  ^{\sigma}\underbrace{x_{\sigma
\left(  1\right)  }^{\beta_{1}}x_{\sigma\left(  2\right)  }^{\beta_{2}}\cdots
x_{\sigma\left(  N\right)  }^{\beta_{N}}}_{\substack{=\sigma\cdot\left(
x_{1}^{\beta_{1}}x_{2}^{\beta_{2}}\cdots x_{N}^{\beta_{N}}\right)
\\\text{(because the action of }\sigma\in S_{N}\\\text{on }\mathcal{P}\text{
substitutes }x_{\sigma\left(  i\right)  }\text{ for each }x_{i}\text{)}%
}}\nonumber\\
&  \ \ \ \ \ \ \ \ \ \ \ \ \ \ \ \ \ \ \ \ \left(  \text{by the definition of
a determinant}\right) \nonumber\\
&  =\sum_{\sigma\in S_{N}}\left(  -1\right)  ^{\sigma}\sigma\cdot
\underbrace{\left(  x_{1}^{\beta_{1}}x_{2}^{\beta_{2}}\cdots x_{N}^{\beta_{N}%
}\right)  }_{=x^{\beta}}=\sum_{\sigma\in S_{N}}\left(  -1\right)  ^{\sigma
}\sigma\cdot x^{\beta}. \label{pf.lem.sf.stemb-lem.abeta=}%
\end{align}
Applying this to $\beta=\nu+\rho$, we obtain%
\[
a_{\nu+\rho}=\sum_{\sigma\in S_{N}}\left(  -1\right)  ^{\sigma}\sigma\cdot
x^{\nu+\rho}.
\]
Multiplying both sides of this equality by $s_{\lambda/\mu}$, we find%
\begin{align}
a_{\nu+\rho}\cdot s_{\lambda/\mu}  &  =\left(  \sum_{\sigma\in S_{N}}\left(
-1\right)  ^{\sigma}\sigma\cdot x^{\nu+\rho}\right)  \cdot s_{\lambda/\mu
}\nonumber\\
&  =\sum_{\sigma\in S_{N}}\left(  -1\right)  ^{\sigma}\left(  \sigma\cdot
x^{\nu+\rho}\right)  \cdot\underbrace{s_{\lambda/\mu}}_{\substack{=\sigma\cdot
s_{\lambda/\mu}\\\text{(since }s_{\lambda/\mu}\text{ is symmetric,}\\\text{so
that }\sigma\cdot s_{\lambda/\mu}=s_{\lambda/\mu}\text{)}}}\nonumber\\
&  =\sum_{\sigma\in S_{N}}\left(  -1\right)  ^{\sigma}\underbrace{\left(
\sigma\cdot x^{\nu+\rho}\right)  \cdot\left(  \sigma\cdot s_{\lambda/\mu
}\right)  }_{\substack{=\sigma\cdot\left(  x^{\nu+\rho}s_{\lambda/\mu}\right)
\\\text{(by (\ref{pf.lem.sf.stemb-lem.sigfg=}))}}}\nonumber\\
&  =\sum_{\sigma\in S_{N}}\left(  -1\right)  ^{\sigma}\sigma\cdot\left(
x^{\nu+\rho}s_{\lambda/\mu}\right)  . \label{pf.lem.stemb-lem.anu+rslm}%
\end{align}

However, multiplying both sides of (\ref{pf.lem.sf.stemb-lem.slm=}) by
$x^{\nu+\rho}$, we find%
\begin{align*}
x^{\nu+\rho}s_{\lambda/\mu}  &  =x^{\nu+\rho}\sum_{T\in\operatorname*{SSYT}%
\left(  \lambda/\mu\right)  }x^{\operatorname*{cont}T}=\sum_{T\in
\operatorname*{SSYT}\left(  \lambda/\mu\right)  }\underbrace{x^{\nu+\rho
}x^{\operatorname*{cont}T}}_{\substack{=x^{\nu+\rho+\operatorname*{cont}%
T}\\\text{(by (\ref{eq.def.sf.tuple-addition.xab}))}}}\\
&  =\sum_{T\in\operatorname*{SSYT}\left(  \lambda/\mu\right)  }%
\underbrace{x^{\nu+\rho+\operatorname*{cont}T}}_{=x^{\nu+\operatorname*{cont}%
T+\rho}}=\sum_{T\in\operatorname*{SSYT}\left(  \lambda/\mu\right)  }%
x^{\nu+\operatorname*{cont}T+\rho}.
\end{align*}
Thus, any $\sigma\in S_{N}$ satisfies%
\[
\sigma\cdot\left(  x^{\nu+\rho}s_{\lambda/\mu}\right)  =\sigma\cdot\left(
\sum_{T\in\operatorname*{SSYT}\left(  \lambda/\mu\right)  }x^{\nu
+\operatorname*{cont}T+\rho}\right)  =\sum_{T\in\operatorname*{SSYT}\left(
\lambda/\mu\right)  }\sigma\cdot x^{\nu+\operatorname*{cont}T+\rho}%
\]
(since the action of $\sigma$ on $\mathcal{P}$ is a $K$-algebra automorphism
of $\mathcal{P}$, and thus in particular is $K$-linear). Therefore,
(\ref{pf.lem.stemb-lem.anu+rslm}) becomes%
\begin{align}
a_{\nu+\rho}\cdot s_{\lambda/\mu}  &  =\sum_{\sigma\in S_{N}}\left(
-1\right)  ^{\sigma}\underbrace{\sigma\cdot\left(  x^{\nu+\rho}s_{\lambda/\mu
}\right)  }_{=\sum_{T\in\operatorname*{SSYT}\left(  \lambda/\mu\right)
}\sigma\cdot x^{\nu+\operatorname*{cont}T+\rho}}\nonumber\\
&  =\sum_{\sigma\in S_{N}}\left(  -1\right)  ^{\sigma}\sum_{T\in
\operatorname*{SSYT}\left(  \lambda/\mu\right)  }\sigma\cdot x^{\nu
+\operatorname*{cont}T+\rho}\nonumber\\
&  =\sum_{T\in\operatorname*{SSYT}\left(  \lambda/\mu\right)  }%
\ \ \underbrace{\sum_{\sigma\in S_{N}}\left(  -1\right)  ^{\sigma}\sigma\cdot
x^{\nu+\operatorname*{cont}T+\rho}}_{\substack{=a_{\nu+\operatorname*{cont}%
T+\rho}\\\text{(by (\ref{pf.lem.sf.stemb-lem.abeta=}), applied to }\beta
=\nu+\operatorname*{cont}T+\rho\text{)}}}\nonumber\\
&  =\sum_{T\in\operatorname*{SSYT}\left(  \lambda/\mu\right)  }a_{\nu
+\operatorname*{cont}T+\rho}. \label{pf.lem.stemb-lem.before-cancel}%
\end{align}
This almost looks like the claim we want to prove, but the sum on the right
hand side has too many addends: It runs over all semistandard tableaux of
shape $\lambda/\mu$, while we only want it to run over the ones that are $\nu
$-Yamanouchi. Thus, we will now try to cancel the extraneous addends (i.e.,
the addends corresponding to the $T$'s that are not $\nu$-Yamanouchi).

Let us first make this a bit more precise. We define two sets%
\begin{align*}
\mathcal{A}  &  :=\operatorname*{SSYT}\left(  \lambda/\mu\right)
\ \ \ \ \ \ \ \ \ \ \text{and}\\
\mathcal{X}  &  :=\left\{  T\in\operatorname*{SSYT}\left(  \lambda/\mu\right)
\ \mid\ T\text{ is not }\nu\text{-Yamanouchi}\right\}  .
\end{align*}
For each $T\in\mathcal{A}$, we define an element $\operatorname*{sign}%
T\in\mathcal{P}$ by
\[
\operatorname*{sign}T:=a_{\nu+\operatorname*{cont}T+\rho}.
\]
Thus, (\ref{pf.lem.stemb-lem.before-cancel}) rewrites as%
\begin{equation}
a_{\nu+\rho}\cdot s_{\lambda/\mu}=\sum_{T\in\mathcal{A}}\operatorname*{sign}T.
\label{pf.lem.stemb-lem.before-cancel-sign}%
\end{equation}
We shall now construct a sign-reversing involution $f:\mathcal{X}%
\rightarrow\mathcal{X}$.

Indeed, let $T\in\mathcal{X}$. Thus, $T$ is a semistandard tableau of shape
$\lambda/\mu$ that is not $\nu$-Yamanouchi (by the definition of $\mathcal{X}%
$). Hence, there exists at least one $j\geq1$ such that the $N$-tuple
$\nu+\operatorname*{cont}\left(  \operatorname{col}_{\geq j}T\right)  $ is
\textbf{not} an $N$-partition (by the definition of \textquotedblleft$\nu
$-Yamanouchi\textquotedblright). Any such $j$ will be called a \emph{violator}
of $T$. Thus, there exists at least one violator of $T$. In other words, the
set of all violators of $T$ is nonempty. On the other hand, this set is
finite\footnote{\textit{Proof.} Let $j\geq1$ be larger than each entry of
$\lambda$. Then, the restricted tableau $\operatorname{col}_{\geq j}T$ is
empty and thus satisfies $\operatorname*{cont}\left(  \operatorname{col}_{\geq
j}T\right)  =\mathbf{0}$. Hence, $\nu+\underbrace{\operatorname*{cont}\left(
\operatorname{col}_{\geq j}T\right)  }_{=\mathbf{0}}=\nu+\mathbf{0}=\nu$,
which is an $N$-partition by assumption. Thus, $j$ is not a violator of $T$
(by the definition of a \textquotedblleft violator\textquotedblright).
\par
Forget that we fixed $j$. We thus have shown that if $j\geq1$ is larger than
each entry of $\lambda$, then $j$ is not a violator of $T$. Hence, if $j\geq1$
is sufficiently high, then $j$ is not a violator of $T$. Thus, the set of
violators of $T$ is bounded from above, and therefore finite (since it is a
set of positive integers).}. Hence, this set has a maximum element. In other
words, the largest violator of $T$ exists.\footnote{For example, if
$\nu=\mathbf{0}$, then the tableaux
\[
T_{2}=\ytableaushort{\none 11,23}\ \ ,\qquad T_{3}%
=\ytableaushort{\none 12,22}\ \ ,\qquad T_{5}=\ytableaushort{\none 1,1,3}
\]
from Example \ref{exa.sf.yamanouchi.1} have largest violators $2,3,1$,
respectively.}

Let $j$ be the \textbf{largest} violator of $T$. Then, $\nu
+\operatorname*{cont}\left(  \operatorname{col}_{\geq j}T\right)  $ is not an
$N$-partition, but $\nu+\operatorname*{cont}\left(  \operatorname{col}_{\geq
j+1}T\right)  $ is an $N$-partition (since $j$ is the \textbf{largest}
violator of $T$).

Define two $N$-tuples $b\in\mathbb{N}^{N}$ and $c\in\mathbb{N}^{N}$ by
$b:=\nu+\operatorname*{cont}\left(  \operatorname{col}_{\geq j}T\right)  $ and
$c:=\nu+\operatorname*{cont}\left(  \operatorname{col}_{\geq j+1}T\right)  $.
Thus, $b$ is not an $N$-partition\footnote{since $\nu+\operatorname*{cont}%
\left(  \operatorname{col}_{\geq j}T\right)  $ is not an $N$-partition}, but
$c$ is an $N$-partition\footnote{since $\nu+\operatorname*{cont}\left(
\operatorname{col}_{\geq j+1}T\right)  $ is an $N$-partition}.

Since $b$ is not an $N$-partition, there exists some $k\in\left[  N-1\right]
$ such that $b_{k}<b_{k+1}$. Such a $k$ will be called a \emph{misstep} of
$T$. Thus, there exists a misstep of $T$. Let $k$ be the \textbf{smallest}
misstep of $T$. Then, $b_{k}<b_{k+1}$ (since $k$ is a misstep of $T$).
Furthermore, $c_{k}\geq c_{k+1}$ (since $c$ is an $N$-partition).

\begin{example}
\label{exa.pf.lem.stemb-lem.exa1}For this example, let $N=7$ and $M=5$ and
$\nu=\left(  4,2,2,0,0,0,0\right)  $ and $\lambda=\left(  7,7,6,5,4\right)  $
and $\mu=\left(  6,2,2,0,0\right)  $. Let $T$ be the following semistandard
tableau of shape $\lambda/\mu$:%
\[
T=\ytableaushort{
\none\none\none\none\none\none 2,
\none\none 11223,
\none\none 2334,
13356,
3456
}\ \ .
\]
We have
\begin{align*}
\nu+\underbrace{\operatorname*{cont}\left(  \operatorname{col}_{\geq
j}T\right)  }_{=\mathbf{0}}  &  =\nu=\left(  4,2,2,0,0,0,0\right)
\ \ \ \ \ \ \ \ \ \ \text{for each }j\geq8;\\
\nu+\underbrace{\operatorname*{cont}\left(  \operatorname{col}_{\geq
7}T\right)  }_{=\left(  0,1,1,0,0,0,0\right)  }  &  =\nu+\left(
0,1,1,0,0,0,0\right)  =\left(  4,3,3,0,0,0,0\right)  ;\\
\nu+\underbrace{\operatorname*{cont}\left(  \operatorname{col}_{\geq
6}T\right)  }_{=\left(  0,2,1,1,0,0,0\right)  }  &  =\nu+\left(
0,2,1,1,0,0,0\right)  =\left(  4,4,3,1,0,0,0\right)  ;\\
\nu+\underbrace{\operatorname*{cont}\left(  \operatorname{col}_{\geq
5}T\right)  }_{=\left(  0,3,2,1,0,1,0\right)  }  &  =\nu+\left(
0,3,2,1,0,1,0\right)  =\left(  4,5,4,1,0,1,0\right)  .
\end{align*}
Thus, $\nu+\operatorname*{cont}\left(  \operatorname{col}_{\geq5}T\right)  $
is not an $N$-partition. This shows that $T$ is not $\nu$-Yamanouchi (so that
$T\in\mathcal{X}$), and in fact $5$ is the largest violator of $T$. Thus,
according to our above instructions, we set
\begin{align*}
j  &  :=5\ \ \ \ \ \ \ \ \ \ \text{and}\\
b  &  :=\nu+\operatorname*{cont}\left(  \operatorname{col}_{\geq j}T\right)
=\nu+\operatorname*{cont}\left(  \operatorname{col}_{\geq5}T\right)  =\left(
4,5,4,1,0,1,0\right)  \ \ \ \ \ \ \ \ \ \ \text{and}\\
c  &  :=\nu+\operatorname*{cont}\left(  \operatorname{col}_{\geq j+1}T\right)
=\nu+\operatorname*{cont}\left(  \operatorname{col}_{\geq6}T\right)  =\left(
4,4,3,1,0,0,0\right)  .
\end{align*}
The missteps of $T$ are the numbers $k\in\left[  N-1\right]  $ such that
$b_{k}<b_{k+1}$; these numbers are $1$ and $5$ (since $b_{1}<b_{2}$ and
$b_{5}<b_{6}$). Thus, the smallest misstep of $T$ is $1$. Hence, we set $k:=1$.
\end{example}

Let us next make a few general observations about $b$ and $c$.

The restrictions $\operatorname{col}_{\geq j}T$ and $\operatorname{col}_{\geq
j+1}T$ of $T$ are \textquotedblleft almost the same\textquotedblright: The
only difference between them is that the $j$-th column of $T$ is included in
$\operatorname{col}_{\geq j}T$ but not in $\operatorname{col}_{\geq j+1}T$.
Hence,%
\[
\operatorname*{cont}\left(  \operatorname{col}_{\geq j}T\right)
=\operatorname*{cont}\left(  \operatorname{col}_{\geq j+1}T\right)
+\operatorname*{cont}\left(  \operatorname{col}_{j}T\right)  ,
\]
where $\operatorname{col}_{j}T$ denotes the $j$-th column of $T$ (or, to be
more precise, the restriction of $T$ to the $j$-th column). Now,%
\begin{align}
b  &  =\nu+\underbrace{\operatorname*{cont}\left(  \operatorname{col}_{\geq
j}T\right)  }_{=\operatorname*{cont}\left(  \operatorname{col}_{\geq
j+1}T\right)  +\operatorname*{cont}\left(  \operatorname{col}_{j}T\right)
}=\underbrace{\nu+\operatorname*{cont}\left(  \operatorname{col}_{\geq
j+1}T\right)  }_{=c}+\operatorname*{cont}\left(  \operatorname{col}%
_{j}T\right) \nonumber\\
&  =c+\operatorname*{cont}\left(  \operatorname{col}_{j}T\right)  .
\label{pf.lem.stemb-lem.b=c+cont}%
\end{align}

Now, recall that the tableau $T$ is semistandard; thus, its entries increase
strictly down each column. Hence, in particular, the entries of the $j$-th
column of $T$ increase strictly down this column. Therefore, any given number
$i\in\left[  N\right]  $ appears at most once in this column. In other words,
any given number $i\in\left[  N\right]  $ appears at most once in
$\operatorname{col}_{j}T$. In other words, $\left(  \operatorname*{cont}%
\left(  \operatorname{col}_{j}T\right)  \right)  _{i}\leq1$ for each
$i\in\left[  N\right]  $ (because $\left(  \operatorname*{cont}\left(
\operatorname{col}_{j}T\right)  \right)  _{i}$ counts how often $i$ appears in
$\operatorname{col}_{j}T$). Applying this inequality to $i=k+1$, we obtain
$\left(  \operatorname*{cont}\left(  \operatorname{col}_{j}T\right)  \right)
_{k+1}\leq1$. Now, from (\ref{pf.lem.stemb-lem.b=c+cont}), we obtain%
\[
b_{k+1}=\left(  c+\operatorname*{cont}\left(  \operatorname{col}_{j}T\right)
\right)  _{k+1}=c_{k+1}+\underbrace{\left(  \operatorname*{cont}\left(
\operatorname{col}_{j}T\right)  \right)  _{k+1}}_{\leq1}\leq c_{k+1}+1,
\]
so that $b_{k}<b_{k+1}\leq c_{k+1}+1$. Since $b_{k}$ and $c_{k+1}+1$ are
integers, this entails $b_{k}\leq\left(  c_{k+1}+1\right)  -1=c_{k+1}$.
However, (\ref{pf.lem.stemb-lem.b=c+cont}) also yields%
\[
b_{k}=\left(  c+\operatorname*{cont}\left(  \operatorname{col}_{j}T\right)
\right)  _{k}=c_{k}+\underbrace{\left(  \operatorname*{cont}\left(
\operatorname{col}_{j}T\right)  \right)  _{k}}_{\substack{\geq0\\\text{(since
}\left(  \operatorname*{cont}\left(  \operatorname{col}_{j}T\right)  \right)
_{k}\text{ counts}\\\text{how often }k\text{ appears in }\operatorname{col}%
_{j}T\text{)}}}\geq c_{k},
\]
so that $c_{k}\leq b_{k}\leq c_{k+1}$. Combining this with $c_{k}\geq c_{k+1}%
$, we obtain $c_{k}=c_{k+1}$. Hence, $c_{k+1}=c_{k}$. Now, combining
$b_{k}\leq c_{k+1}=c_{k}$ with $c_{k}\leq b_{k}$, we obtain $b_{k}=c_{k}$.
Comparing this with $b_{k}=c_{k}+\left(  \operatorname*{cont}\left(
\operatorname{col}_{j}T\right)  \right)  _{k}$, we obtain $c_{k}+\left(
\operatorname*{cont}\left(  \operatorname{col}_{j}T\right)  \right)
_{k}=c_{k}$, so that $\left(  \operatorname*{cont}\left(  \operatorname{col}%
_{j}T\right)  \right)  _{k}=0$. In other words, the number $k$ appears $0$
times in $\operatorname{col}_{j}T$ (since $\left(  \operatorname*{cont}\left(
\operatorname{col}_{j}T\right)  \right)  _{k}$ counts how often $k$ appears in
$\operatorname{col}_{j}T$). In other words, the number $k$ does not appear in
$\operatorname{col}_{j}T$. In other words, the number $k$ does not appear in
the $j$-th column of $T$.

Let us make one more simple observation, which we will not use until later: We
have $b_{k}<b_{k+1}$, so that $b_{k}\leq b_{k+1}-1$ (since $b_{k}$ and
$b_{k+1}$ are integers). Thus, $b_{k}+1\leq b_{k+1}$. Combining this with
$b_{k+1}\leq\underbrace{c_{k+1}}_{=c_{k}\leq b_{k}}+\,1\leq b_{k}+1$, we
obtain%
\begin{equation}
b_{k}+1=b_{k+1}. \label{pf.lem.stemb-lem.ck=ck+1}%
\end{equation}

Now, let $\operatorname{col}_{<j}T$ be the restriction of $T$ to columns
$1,2,\ldots,j-1$ (that is, the result of removing all but the first $j-1$
columns from $T$). Formally speaking, this means the restriction of the map
$T$ to the set $\left\{  \left(  u,v\right)  \in Y\left(  \lambda/\mu\right)
\ \mid\ v<j\right\}  $. This restriction $\operatorname{col}_{<j}T$ is a
semistandard skew tableau of a certain (skew) shape; thus, we can apply the
Bender--Knuth involution $\beta_{k}$ (from our above proof of Theorem
\ref{thm.sf.skew-schur-symm}) to this tableau $\operatorname{col}_{<j}T$
instead of $T$. Let $T^{\ast}$ be the tableau obtained from $T$ by applying
$\beta_{k}$ \textbf{only to the columns }$1,2,\ldots,j-1$ of $T$ (that is,
replacing $\operatorname{col}_{<j}T$ by $\beta_{k}\left(  \operatorname{col}%
_{<j}T\right)  $), while leaving the columns $j,j+1,j+2,\ldots$ unchanged.
Thus, formally, $T^{\ast}$ is the tableau of shape $Y\left(  \lambda
/\mu\right)  $ defined by%
\begin{align}
\operatorname{col}_{<j}\left(  T^{\ast}\right)   &  =\beta_{k}\left(
\operatorname{col}_{<j}T\right)  \ \ \ \ \ \ \ \ \ \ \text{and}%
\label{pf.lem.stemb-lem.T*1}\\
\operatorname{col}_{\geq j}\left(  T^{\ast}\right)   &  =\operatorname{col}%
_{\geq j}T \label{pf.lem.stemb-lem.T*2}%
\end{align}
(where $\operatorname{col}_{<j}\left(  T^{\ast}\right)  $ is defined just as
$\operatorname{col}_{<j}T$ was defined, except that we are using $T^{\ast}$
instead of $T$).

\begin{example}
\label{exa.pf.lem.stemb-lem.exa2}Let $M$, $N$, $\nu$, $\lambda$, $\mu$ and $T$
be as in Example \ref{exa.pf.lem.stemb-lem.exa1}. Let us now compute $T^{\ast
}$. As we know, $j=5$ and $k=1$. Thus, in order to obtain $T^{\ast}$, we need
to apply the Bender--Knuth involution $\beta_{k}=\beta_{1}$ \textbf{only to
the columns }$1,2,\ldots,j-1$ of $T$ (that is, only to the first $j-1=4$
columns of $T$), while leaving the columns $5,6,7,\ldots$ unchanged. Here is
what this looks like:%
\[
T=\ytableaushort{
\none\none\none\none\none\none {*(lightgray)2},
\none\none 11{*(lightgray)2}{*(lightgray)2}{*(lightgray)3},
\none\none 23{*(lightgray)3}{*(lightgray)4},
1335{*(lightgray)6},
3456
}\ \ \Longrightarrow\ \ T^{\ast}=\ytableaushort{
\none\none\none\none\none\none {*(lightgray)2},
\none\none 12{*(lightgray)2}{*(lightgray)2}{*(lightgray)3},
\none\none 23{*(lightgray)3}{*(lightgray)4},
2335{*(lightgray)6},
3456
}
\]
(where we have grayed out all boxes in columns $5,6,7,\ldots$, because the
entries in these boxes stay unchanged and are ignored by the Bender--Knuth involution).
\end{example}

We shall now show that $T^{\ast}\in\mathcal{X}$. Indeed, let us first check
that the tableau $T^{\ast}$ is semistandard. We know that the tableau $T$ is
semistandard, so that its restrictions $\operatorname{col}_{<j}T$ and
$\operatorname{col}_{\geq j}T$ are semistandard; thus, $\beta_{k}\left(
\operatorname{col}_{<j}T\right)  $ is semistandard as well (since the
Bender--Knuth involution $\beta_{k}$ sends semistandard tableaux to
semistandard tableaux). Now, recall that the tableau $T^{\ast}$ is obtained
from $T$ by applying $\beta_{k}$ to columns $1,2,\ldots,j-1$ only; thus,
$T^{\ast}$ is obtained by glueing the tableaux $\beta_{k}\left(
\operatorname{col}_{<j}T\right)  $ and $\operatorname{col}_{\geq j}T$ together
(along a vertical line). Hence:

\begin{itemize}
\item Each column of $T^{\ast}$ is either a column of $\beta_{k}\left(
\operatorname{col}_{<j}T\right)  $ or a column of $\operatorname{col}_{\geq
j}T$ (depending on whether it is one of columns $1,2,\ldots,j-1$ or one of
columns $j,j+1,j+2,\ldots$). In either case, the entries of this column
increase strictly down this column (because the tableaux $\beta_{k}\left(
\operatorname{col}_{<j}T\right)  $ and $\operatorname{col}_{\geq j}T$ are
semistandard). Thus, we have shown that the entries of $T^{\ast}$ increase
strictly down each column.

\item It is not hard to see that the entries of $T^{\ast}$ increase weakly
along each row\footnote{\textit{Proof.} Let $i\in\left[  M\right]  $. We must
prove that the entries of $T^{\ast}$ increase weakly along the $i$-th row of
$T^{\ast}$. Assume the contrary. Thus, there exist two adjacent entries in the
$i$-th row of $T^{\ast}$ that are out of order (in the sense that the one
lying further left is larger than the one lying further right). In other
words, there exists some positive integer $u$ such that $\left(  i,u\right)
\in Y\left(  \lambda/\mu\right)  $ and $\left(  i,u+1\right)  \in Y\left(
\lambda/\mu\right)  $ and $T^{\ast}\left(  i,u\right)  >T^{\ast}\left(
i,u+1\right)  $. Consider this $u$.
\par
Recall that $T^{\ast}$ is obtained by glueing the tableaux $\beta_{k}\left(
\operatorname{col}_{<j}T\right)  $ and $\operatorname{col}_{\geq j}T$ together
(along a vertical line). Thus, the $i$-th row of $T^{\ast}$ is obtained by
glueing the $i$-th row of $\beta_{k}\left(  \operatorname{col}_{<j}T\right)  $
together with the $i$-th row of $\operatorname{col}_{\geq j}T$. In other
words, this row consists of two blocks, looking as follows:%
\[%
\begin{tikzpicture}
\draw(0, 0) rectangle (7, 1);
\draw(7, 0) rectangle (12, 1);
\node at (3.5, 0.5) {\small$i$-th row of $\beta_k\left(\operatorname
{col}_{< j}T\right)$};
\node at (9.5, 0.5) {\small$i$-th row of $\operatorname{col}_{\geq j}T$};
\end{tikzpicture}%
\ \ .
\]
We shall refer to these two blocks as the \emph{left block} and the
\emph{right block} (so the left block is the $i$-th row of $\beta_{k}\left(
\operatorname{col}_{<j}T\right)  $, whereas the right block is the $i$-th row
of $\operatorname{col}_{\geq j}T$). The boundary between the two blocks falls
between the $\left(  j-1\right)  $-st and $j$-th columns; the left block
covers columns $1,2,\ldots,j-1$, while the right block covers columns
$j,j+1,j+2,\ldots$.
\par
The entries of the left block increase weakly from left to right (since this
left block is a row of the tableau $\beta_{k}\left(  \operatorname{col}%
_{<j}T\right)  $, which is semistandard). Thus, if both boxes $\left(
i,u\right)  $ and $\left(  i,u+1\right)  $ belonged to the left block, then we
would have $T^{\ast}\left(  i,u\right)  \leq T^{\ast}\left(  i,u+1\right)  $,
which would contradict $T^{\ast}\left(  i,u\right)  >T^{\ast}\left(
i,u+1\right)  $. Hence, it is impossible for both boxes $\left(  i,u\right)  $
and $\left(  i,u+1\right)  $ to belong to the left block; thus, at least one
of these boxes must belong to the right block. Therefore, $\left(
i,u+1\right)  $ belongs to the right block (since the right block is further
right than the left block).
\par
The entries of the right block also increase weakly from left to right (since
this right block is a row of the tableau $\operatorname{col}_{\geq j}T$, which
is semistandard). Thus, if both boxes $\left(  i,u\right)  $ and $\left(
i,u+1\right)  $ belonged to the right block, then we would have $T^{\ast
}\left(  i,u\right)  \leq T^{\ast}\left(  i,u+1\right)  $, which would
contradict $T^{\ast}\left(  i,u\right)  >T^{\ast}\left(  i,u+1\right)  $.
Hence, it is impossible for both boxes $\left(  i,u\right)  $ and $\left(
i,u+1\right)  $ to belong to the right block; thus, at least one of these
boxes must belong to the left block. Therefore, $\left(  i,u\right)  $ belongs
to the left block (since $\left(  i,u+1\right)  $ belongs to the right block).
\par
Thus, the boxes $\left(  i,u\right)  $ and $\left(  i,u+1\right)  $ straddle
the boundary between the left block and the right block. Since this boundary
falls between the $\left(  j-1\right)  $-st and $j$-th columns, this entails
that the box $\left(  i,u\right)  $ lies on the $\left(  j-1\right)  $-st
column, while the box $\left(  i,u+1\right)  $ lies on the $j$-th column. In
other words, $u=j-1$ and $u+1=j$. Thus, the inequality $T^{\ast}\left(
i,u\right)  >T^{\ast}\left(  i,u+1\right)  $ can be rewritten as $T^{\ast
}\left(  i,j-1\right)  >T^{\ast}\left(  i,j\right)  $. Moreover, from $u=j-1$,
we obtain $j-1=u$ and thus $\left(  i,j-1\right)  =\left(  i,u\right)  \in
Y\left(  \lambda/\mu\right)  $. Furthermore, from $u+1=j$, we obtain $j=u+1$
and thus $\left(  i,j\right)  =\left(  i,u+1\right)  \in Y\left(  \lambda
/\mu\right)  $.
\par
The equality (\ref{pf.lem.stemb-lem.T*2}) shows that the entries of $T^{\ast}$
in columns $j,j+1,j+2,\ldots$ equal the corresponding entries of $T$. Thus, in
particular, we have $T^{\ast}\left(  i,j\right)  =T\left(  i,j\right)  $
(since the box $\left(  i,j\right)  $ lies in column $j$). Thus, $T^{\ast
}\left(  i,j-1\right)  >T^{\ast}\left(  i,j\right)  =T\left(  i,j\right)  $.
\par
On the other hand, the tableau $T$ is semistandard, so that its entries
increase weakly along each row. Hence, $T\left(  i,j-1\right)  \leq T\left(
i,j\right)  $. Therefore, $T\left(  i,j\right)  \geq T\left(  i,j-1\right)  $,
so that $T^{\ast}\left(  i,j-1\right)  >T\left(  i,j\right)  \geq T\left(
i,j-1\right)  $.
\par
We know that the number $k$ does not appear in the $j$-th column of $T$; thus,
$T\left(  i,j\right)  \neq k$ (since $T\left(  i,j\right)  $ is an entry in
the $j$-th column of $T$).
\par
We recall a simple property of the Bender-Knuth involution $\beta_{k}$ (which
follows directly from the construction of $\beta_{k}$): When we apply
$\beta_{k}$ to a semistandard tableau,
\par
\begin{itemize}
\item some $k$'s get replaced by $\left(  k+1\right)  $'s,
\par
\item some $\left(  k+1\right)  $'s get replaced by $k$'s, and
\par
\item all other entries remain unchanged.
\end{itemize}
\par
Thus, in particular, when we apply $\beta_{k}$ to a semistandard tableau, the
only entries that can get replaced by larger entries are $k$'s, and in that
case they can only be replaced by $\left(  k+1\right)  $'s. In other words, if
some entry of a semistandard tableau gets replaced by a larger entry when we
apply $\beta_{k}$ to the tableau, then this entry must have been $k$ before
applying $\beta_{k}$, and must get replaced by $k+1$ when $\beta_{k}$ is
applied.
\par
Since $T^{\ast}$ is obtained from $T$ by applying $\beta_{k}$ to columns
$1,2,\ldots,j-1$ (while all other columns remain unchanged), we thus conclude
that if some entry of $T$ gets replaced by a larger entry when we pass from
$T$ to $T^{\ast}$, then this entry must have been $k$ in $T$, and must get
replaced by $k+1$ in $T^{\ast}$. Let us restate this in a more formal
language: If $\left(  p,q\right)  \in Y\left(  \lambda/\mu\right)  $ satisfies
$T^{\ast}\left(  p,q\right)  >T\left(  p,q\right)  $, then%
\[
T\left(  p,q\right)  =k\ \ \ \ \ \ \ \ \ \ \text{and}%
\ \ \ \ \ \ \ \ \ \ T^{\ast}\left(  p,q\right)  =k+1.
\]
\par
We can apply this to $\left(  p,q\right)  =\left(  i,j-1\right)  $ (since
$T^{\ast}\left(  i,j-1\right)  >T\left(  i,j-1\right)  $), and thus conclude
that $T\left(  i,j-1\right)  =k$ and $T^{\ast}\left(  i,j-1\right)  =k+1$.
\par
Now, from $T\left(  i,j-1\right)  =k$, we obtain $k=T\left(  i,j-1\right)
\leq T\left(  i,j\right)  $. On the other hand, $T^{\ast}\left(  i,j-1\right)
>T\left(  i,j\right)  $, so that $T\left(  i,j\right)  <T^{\ast}\left(
i,j-1\right)  =k+1$. Since $T\left(  i,j\right)  $ and $k+1$ are integers,
this entails $T\left(  i,j\right)  \leq\left(  k+1\right)  -1=k$. Combining
this with $k\leq T\left(  i,j\right)  $, we obtain $T\left(  i,j\right)  =k$.
This contradicts $T\left(  i,j\right)  \neq k$. This contradiction shows that
our assumption was wrong. Hence, we have shown that the entries of $T^{\ast}$
increase weakly along the $i$-th row of $T^{\ast}$. Qed.}.
\end{itemize}

Combining these two conclusions, we conclude that $T^{\ast}$ is a semistandard
tableau. In other words, $T^{\ast}\in\operatorname*{SSYT}\left(  \lambda
/\mu\right)  $.

Recall that $\nu+\operatorname*{cont}\left(  \operatorname{col}_{\geq
j}T\right)  $ is not an $N$-partition. In view of (\ref{pf.lem.stemb-lem.T*2}%
), this rewrites as follows: $\nu+\operatorname*{cont}\left(
\operatorname{col}_{\geq j}\left(  T^{\ast}\right)  \right)  $ is not an
$N$-partition. Hence, the tableau $T^{\ast}$ is not $\nu$-Yamanouchi. Thus,
$T^{\ast}\in\mathcal{X}$ (by the definition of $\mathcal{X}$, since $T^{\ast
}\in\operatorname*{SSYT}\left(  \lambda/\mu\right)  $).

Forget that we fixed $T$. Thus, for each tableau $T\in\mathcal{X}$, we have
constructed a tableau $T^{\ast}\in\mathcal{X}$. Let $f:\mathcal{X}%
\rightarrow\mathcal{X}$ be the map that sends each $T\in\mathcal{X}$ to
$T^{\ast}$. We shall now show that $f$ is a sign-reversing involution. First,
we shall show the following:

\begin{statement}
\textit{Observation 1:} The map $f$ is an involution.
\end{statement}

[\textit{Proof of Observation 1:} We must show that $f\circ
f=\operatorname*{id}$. In other words, we must prove that $f\left(  f\left(
T\right)  \right)  =T$ for each $T\in\mathcal{X}$.

Let $T\in\mathcal{X}$. Then, the definition of $f$ yields $f\left(  T\right)
=T^{\ast}$ and $f\left(  T^{\ast}\right)  =\left(  T^{\ast}\right)  ^{\ast}$.
Recall how $T^{\ast}$ is constructed from $T$:

\begin{itemize}
\item We let $j$ be the \textbf{largest} violator of $T$. This is the largest
$j\geq1$ such that $\nu+\operatorname*{cont}\left(  \operatorname{col}_{\geq
j}T\right)  $ is not an $N$-partition.

\item We let $k$ be the \textbf{smallest} misstep of $T$. This is the smallest
$k\in\left[  N-1\right]  $ such that $\left(  \nu+\operatorname*{cont}\left(
\operatorname{col}_{\geq j}T\right)  \right)  _{k}<\left(  \nu
+\operatorname*{cont}\left(  \operatorname{col}_{\geq j}T\right)  \right)
_{k+1}$.\ \ \ \ \footnote{Indeed, a misstep of $T$ was defined to be a
$k\in\left[  N-1\right]  $ such that $b_{k}<b_{k+1}$, where $b=\nu
+\operatorname*{cont}\left(  \operatorname{col}_{\geq j}T\right)  $. In other
words, a misstep of $T$ means a $k\in\left[  N-1\right]  $ such that $\left(
\nu+\operatorname*{cont}\left(  \operatorname{col}_{\geq j}T\right)  \right)
_{k}<\left(  \nu+\operatorname*{cont}\left(  \operatorname{col}_{\geq
j}T\right)  \right)  _{k+1}$.}

\item We apply the Bender--Knuth involution $\beta_{k}$ \textbf{only to the
columns }$1,2,\ldots,j-1$ of $T$, while leaving the columns $j,j+1,j+2,\ldots$
unchanged. The result is $T^{\ast}$.
\end{itemize}

The construction of $\left(  T^{\ast}\right)  ^{\ast}$ from $T^{\ast}$
proceeds similarly:

\begin{itemize}
\item We let $j^{\prime}$ be the \textbf{largest} violator of $T^{\ast}$.

\item We let $k^{\prime}$ be the \textbf{smallest} misstep of $T^{\ast}$.

\item We apply the Bender--Knuth involution $\beta_{k^{\prime}}$ \textbf{only
to the columns }$1,2,\ldots,j^{\prime}-1$ of $T^{\ast}$, while leaving the
columns $j^{\prime},j^{\prime}+1,j^{\prime}+2,\ldots$ unchanged. The result is
$\left(  T^{\ast}\right)  ^{\ast}$.
\end{itemize}

We claim that this construction undoes the previous construction and recovers
$T$ (so that $\left(  T^{\ast}\right)  ^{\ast}=T$). To see this, we argue as follows:

\begin{itemize}
\item We know that $\operatorname{col}_{\geq j}\left(  T^{\ast}\right)
=\operatorname{col}_{\geq j}T$, so that $j$ is the largest violator of
$T^{\ast}$ (since $j$ is the largest violator of $T$)\ \ \ \ \footnote{Here is
the argument in some more detail:
\par
We have%
\begin{equation}
\operatorname{col}_{\geq j}\left(  T^{\ast}\right)  =\operatorname{col}_{\geq
j}T, \label{pf.lem.stemb-lem.o1.pf.fn1.1}%
\end{equation}
and therefore we also have
\begin{equation}
\operatorname{col}_{\geq p}\left(  T^{\ast}\right)  =\operatorname{col}_{\geq
p}T\ \ \ \ \ \ \ \ \ \ \text{for any integer }p>j
\label{pf.lem.stemb-lem.o1.pf.fn1.2}%
\end{equation}
(since the tableau $\operatorname{col}_{\geq p}T$ is obtained by removing some
columns from $\operatorname{col}_{\geq j}T$, whereas the tableau
$\operatorname{col}_{\geq p}\left(  T^{\ast}\right)  $ is obtained in the same
fashion from $\operatorname{col}_{\geq j}\left(  T^{\ast}\right)  $).
\par
We know that $j$ is the largest violator of $T$. In other words, the $N$-tuple
$\nu+\operatorname*{cont}\left(  \operatorname{col}_{\geq j}T\right)  $ is not
an $N$-partition, but $\nu+\operatorname*{cont}\left(  \operatorname{col}%
_{\geq p}T\right)  $ is an $N$-partition for any integer $p>j$. In view of
(\ref{pf.lem.stemb-lem.o1.pf.fn1.1}) and (\ref{pf.lem.stemb-lem.o1.pf.fn1.2}),
we can rewrite this as follows: The $N$-tuple $\nu+\operatorname*{cont}\left(
\operatorname{col}_{\geq j}\left(  T^{\ast}\right)  \right)  $ is not an
$N$-partition, but $\nu+\operatorname*{cont}\left(  \operatorname{col}_{\geq
p}\left(  T^{\ast}\right)  \right)  $ is an $N$-partition for any integer
$p>j$. In other words, $j$ is the largest violator of $T^{\ast}$.}. Therefore,
$j^{\prime}=j$.

\item Knowing that $j^{\prime}=j$ and $\operatorname{col}_{\geq j}\left(
T^{\ast}\right)  =\operatorname{col}_{\geq j}T$, we now conclude that $k$ is
the smallest misstep of $T^{\ast}$ (since $k$ is the smallest misstep of $T$).
Therefore, $k^{\prime}=k$.

\item Knowing that $j^{\prime}=j$ and $k^{\prime}=k$, we conclude that
$\left(  T^{\ast}\right)  ^{\ast}$ is obtained from $T^{\ast}$ by the exact
same operation that we used to obtain $T^{\ast}$ from $T$: namely, by applying
the Bender--Knuth involution $\beta_{k}$ only to the columns $1,2,\ldots,j-1$
(while leaving the columns $j,j+1,j+2,\ldots$ unchanged). However, this
operation undoes itself when applied a second time, because the Bender--Knuth
involution $\beta_{k}$ is an involution\footnote{This has been shown during
our proof of Theorem \ref{thm.sf.skew-schur-symm}.}. Thus, we conclude that
$\left(  T^{\ast}\right)  ^{\ast}=T$.
\end{itemize}

Thus, $f\left(  \underbrace{f\left(  T\right)  }_{=T^{\ast}}\right)  =f\left(
T^{\ast}\right)  =\left(  T^{\ast}\right)  ^{\ast}=T$.

Forget that we fixed $T$. We thus have proved that $f\left(  f\left(
T\right)  \right)  =T$ for each $T\in\mathcal{X}$. As explained, this
completes the proof of Observation 1.] \medskip

Next, we shall show two observations about the effect of the map $f$ on the
sign of a tableau:

\begin{statement}
\textit{Observation 2:} We have $\operatorname*{sign}\left(  f\left(
T\right)  \right)  =-\operatorname*{sign}T$ for all $T\in\mathcal{X}$.
\end{statement}

[\textit{Proof of Observation 2:} Let $T\in\mathcal{X}$. We must show that
$\operatorname*{sign}\left(  f\left(  T\right)  \right)
=-\operatorname*{sign}T$.

Define two $N$-tuples $\alpha\in\mathbb{N}^{N}$ and $\gamma\in\mathbb{N}^{N}$
by $\alpha:=\nu+\operatorname*{cont}T+\rho$ and $\gamma:=\nu
+\operatorname*{cont}\left(  T^{\ast}\right)  +\rho$. Then, $f\left(
T\right)  =T^{\ast}$ (by the definition of $f$), and thus
\begin{align*}
\operatorname*{sign}\left(  f\left(  T\right)  \right)   &
=\operatorname*{sign}\left(  T^{\ast}\right)  =a_{\nu+\operatorname*{cont}%
\left(  T^{\ast}\right)  +\rho}\ \ \ \ \ \ \ \ \ \ \left(  \text{by the
definition of }\operatorname*{sign}\left(  T^{\ast}\right)  \right) \\
&  =a_{\gamma}\ \ \ \ \ \ \ \ \ \ \left(  \text{since }\nu
+\operatorname*{cont}\left(  T^{\ast}\right)  +\rho=\gamma\right)  .
\end{align*}
Also, the definition of $\operatorname*{sign}T$ yields
\[
\operatorname*{sign}T=a_{\nu+\operatorname*{cont}T+\rho}=a_{\alpha
}\ \ \ \ \ \ \ \ \ \ \left(  \text{since }\nu+\operatorname*{cont}%
T+\rho=\alpha\right)  .
\]
Now, we are going to show that the $N$-tuple $\gamma$ is obtained from
$\alpha$ by swapping two entries. Once this is shown, we will easily conclude
$\operatorname*{sign}\left(  f\left(  T\right)  \right)
=-\operatorname*{sign}T$ by applying Lemma \ref{lem.sf.alternant-0}
\textbf{(b)}.

We recall the notations from the construction of $T^{\ast}$: Let $j$ be the
largest violator of $T$. Let $k$ be the smallest misstep of $T$. Define an
$N$-tuple $b\in\mathbb{N}^{N}$ by $b:=\nu+\operatorname*{cont}\left(
\operatorname{col}_{\geq j}T\right)  $. Then, $b_{k}+1=b_{k+1}$ (as we have
proved in (\ref{pf.lem.stemb-lem.ck=ck+1})). However,
\begin{align}
b_{k}  &  =\left(  \nu+\operatorname*{cont}\left(  \operatorname{col}_{\geq
j}T\right)  \right)  _{k}\ \ \ \ \ \ \ \ \ \ \left(  \text{since }%
b=\nu+\operatorname*{cont}\left(  \operatorname{col}_{\geq j}T\right)  \right)
\nonumber\\
&  =\nu_{k}+\underbrace{\left(  \operatorname*{cont}\left(  \operatorname{col}%
_{\geq j}T\right)  \right)  _{k}}_{\substack{=\left(  \text{\# of }k\text{'s
in }\operatorname{col}_{\geq j}T\right)  \\\text{(by Definition
\ref{def.sf.content})}}}\nonumber\\
&  =\nu_{k}+\left(  \text{\# of }k\text{'s in }\operatorname{col}_{\geq
j}T\right)  . \label{pf.lem.stemb-lem.o2.pf.bk=}%
\end{align}
The same argument (applied to $k+1$ instead of $k$) yields
\begin{equation}
b_{k+1}=\nu_{k+1}+\left(  \text{\# of }\left(  k+1\right)  \text{'s in
}\operatorname{col}_{\geq j}T\right)  . \label{pf.lem.stemb-lem.o2.pf.bk+1=}%
\end{equation}

We shall now show that $\gamma_{k}=\alpha_{k+1}$. Indeed, the vertical line
that separates the $\left(  j-1\right)  $-st and $j$-th columns cuts the
tableau $T$ into its two parts $\operatorname{col}_{<j}T$ and
$\operatorname{col}_{\geq j}T$. Thus, every $i\in\left[  N\right]  $ satisfies%
\begin{equation}
\left(  \text{\# of }i\text{'s in }T\right)  =\left(  \text{\# of }i\text{'s
in }\operatorname{col}_{<j}T\right)  +\left(  \text{\# of }i\text{'s in
}\operatorname{col}_{\geq j}T\right)
.\ \ \ \ \ \ \ \ \ \ \label{pf.lem.stemb-lem.o2.pf.numiinT}%
\end{equation}
The same argument (applied to $T^{\ast}$ instead of $T$) shows that every
$i\in\left[  N\right]  $ satisfies%
\begin{equation}
\left(  \text{\# of }i\text{'s in }T^{\ast}\right)  =\left(  \text{\# of
}i\text{'s in }\operatorname{col}_{<j}\left(  T^{\ast}\right)  \right)
+\left(  \text{\# of }i\text{'s in }\operatorname{col}_{\geq j}\left(
T^{\ast}\right)  \right)
.\ \ \ \ \ \ \ \ \ \ \label{pf.lem.stemb-lem.o2.pf.numiinT*}%
\end{equation}

Now, the definition of $\operatorname*{cont}\left(  T^{\ast}\right)  $ yields
\begin{align*}
\left(  \operatorname*{cont}\left(  T^{\ast}\right)  \right)  _{k}  &
=\left(  \text{\# of }k\text{'s in }T^{\ast}\right) \\
&  =\left(  \text{\# of }k\text{'s in }\underbrace{\operatorname{col}%
_{<j}\left(  T^{\ast}\right)  }_{\substack{=\beta_{k}\left(
\operatorname{col}_{<j}T\right)  \\\text{(by (\ref{pf.lem.stemb-lem.T*1}))}%
}}\right)  +\left(  \text{\# of }k\text{'s in }\underbrace{\operatorname{col}%
_{\geq j}\left(  T^{\ast}\right)  }_{\substack{=\operatorname{col}_{\geq
j}T\\\text{(by (\ref{pf.lem.stemb-lem.T*2}))}}}\right) \\
&  \ \ \ \ \ \ \ \ \ \ \ \ \ \ \ \ \ \ \ \ \left(  \text{by
(\ref{pf.lem.stemb-lem.o2.pf.numiinT*}), applied to }i=k\right) \\
&  =\underbrace{\left(  \text{\# of }k\text{'s in }\beta_{k}\left(
\operatorname{col}_{<j}T\right)  \right)  }_{\substack{=\left(  \text{\# of
}\left(  k+1\right)  \text{'s in }\operatorname{col}_{<j}T\right)  \\\text{(by
(\ref{pf.thm.sf.skew-schur-symm.num-k}),}\\\text{applied to }%
\operatorname{col}_{<j}T\text{ instead of }T\text{)}}}+\underbrace{\left(
\text{\# of }k\text{'s in }\operatorname{col}_{\geq j}T\right)  }%
_{\substack{=b_{k}-\nu_{k}\\\text{(by (\ref{pf.lem.stemb-lem.o2.pf.bk=}))}}}\\
&  =\left(  \text{\# of }\left(  k+1\right)  \text{'s in }\operatorname{col}%
_{<j}T\right)  +\underbrace{b_{k}}_{\substack{=b_{k+1}-1\\\text{(since }%
b_{k}+1=b_{k+1}\text{)}}}-\,\nu_{k}\\
&  =\left(  \text{\# of }\left(  k+1\right)  \text{'s in }\operatorname{col}%
_{<j}T\right)  +b_{k+1}-1-\nu_{k}.
\end{align*}
However, $\gamma=\nu+\operatorname*{cont}\left(  T^{\ast}\right)  +\rho$, so
that%
\begin{align}
\gamma_{k}  &  =\left(  \nu+\operatorname*{cont}\left(  T^{\ast}\right)
+\rho\right)  _{k}\nonumber\\
&  =\nu_{k}+\underbrace{\left(  \operatorname*{cont}\left(  T^{\ast}\right)
\right)  _{k}}_{=\left(  \text{\# of }\left(  k+1\right)  \text{'s in
}\operatorname{col}_{<j}T\right)  +b_{k+1}-1-\nu_{k}}+\underbrace{\rho_{k}%
}_{\substack{=N-k\\\text{(by the definition of }\rho\text{)}}}\nonumber\\
&  =\nu_{k}+\left(  \text{\# of }\left(  k+1\right)  \text{'s in
}\operatorname{col}_{<j}T\right)  +b_{k+1}-1-\nu_{k}+N-k\nonumber\\
&  =\left(  \text{\# of }\left(  k+1\right)  \text{'s in }\operatorname{col}%
_{<j}T\right)  +b_{k+1}-1+N-k. \label{pf.lem.stemb-lem.o2.pf.gammak=}%
\end{align}

On the other hand, the definition of $\operatorname*{cont}T$ yields
\begin{align*}
\left(  \operatorname*{cont}T\right)  _{k+1}  &  =\left(  \text{\# of }\left(
k+1\right)  \text{'s in }T\right) \\
&  =\left(  \text{\# of }\left(  k+1\right)  \text{'s in }\operatorname{col}%
_{<j}T\right)  +\underbrace{\left(  \text{\# of }\left(  k+1\right)  \text{'s
in }\operatorname{col}_{\geq j}T\right)  }_{\substack{=b_{k+1}-\nu
_{k+1}\\\text{(by (\ref{pf.lem.stemb-lem.o2.pf.bk+1=}))}}}\\
&  \ \ \ \ \ \ \ \ \ \ \ \ \ \ \ \ \ \ \ \ \left(  \text{by
(\ref{pf.lem.stemb-lem.o2.pf.numiinT}), applied to }i=k+1\right) \\
&  =\left(  \text{\# of }\left(  k+1\right)  \text{'s in }\operatorname{col}%
_{<j}T\right)  +b_{k+1}-\nu_{k+1}.
\end{align*}
However, $\alpha=\nu+\operatorname*{cont}T+\rho$, so that%
\begin{align*}
\alpha_{k+1}  &  =\left(  \nu+\operatorname*{cont}T+\rho\right)  _{k+1}\\
&  =\nu_{k+1}+\underbrace{\left(  \operatorname*{cont}T\right)  _{k+1}%
}_{=\left(  \text{\# of }\left(  k+1\right)  \text{'s in }\operatorname{col}%
_{<j}T\right)  +b_{k+1}-\nu_{k+1}}+\underbrace{\rho_{k+1}}%
_{\substack{=N-\left(  k+1\right)  \\\text{(by the definition of }\rho
\text{)}}}\\
&  =\nu_{k+1}+\left(  \text{\# of }\left(  k+1\right)  \text{'s in
}\operatorname{col}_{<j}T\right)  +b_{k+1}-\nu_{k+1}+N-\left(  k+1\right) \\
&  =\left(  \text{\# of }\left(  k+1\right)  \text{'s in }\operatorname{col}%
_{<j}T\right)  +b_{k+1}-1+N-k.
\end{align*}
Comparing this with (\ref{pf.lem.stemb-lem.o2.pf.gammak=}), we obtain%
\begin{equation}
\gamma_{k}=\alpha_{k+1}. \label{pf.lem.stemb-lem.o2.pf.gammak=alk+1}%
\end{equation}

A similar argument (using (\ref{pf.thm.sf.skew-schur-symm.num-k+1}) instead of
(\ref{pf.thm.sf.skew-schur-symm.num-k})) can be used to show that%
\begin{equation}
\gamma_{k+1}=\alpha_{k}. \label{pf.lem.stemb-lem.o2.pf.gammak+1=alk}%
\end{equation}
(See Section \ref{sec.details.sf.schur} for the details of this argument.)

A further argument of this form (using (\ref{pf.thm.sf.skew-schur-symm.num-i})
instead of (\ref{pf.thm.sf.skew-schur-symm.num-k})) can be used to show that%
\begin{equation}
\gamma_{i}=\alpha_{i}\ \ \ \ \ \ \ \ \ \ \text{for each }i\in\left[  N\right]
\text{ satisfying }i\neq k\text{ and }i\neq k+1.
\label{pf.lem.stemb-lem.o2.pf.gammai=ali}%
\end{equation}
(See Section \ref{sec.details.sf.schur} for the details of this argument.)

Combining the three equalities (\ref{pf.lem.stemb-lem.o2.pf.gammak=alk+1}),
(\ref{pf.lem.stemb-lem.o2.pf.gammak+1=alk}) and
(\ref{pf.lem.stemb-lem.o2.pf.gammai=ali}), we see that the $N$-tuple $\gamma$
is obtained from $\alpha$ by swapping two entries (namely, the $k$-th and the
$\left(  k+1\right)  $-st entry). Thus, Lemma \ref{lem.sf.alternant-0}
\textbf{(b)} (applied to $\gamma$ instead of $\beta$) yields that $a_{\gamma
}=-a_{\alpha}$. This rewrites as $\operatorname*{sign}\left(  f\left(
T\right)  \right)  =-\operatorname*{sign}T$ (since $\operatorname*{sign}%
\left(  f\left(  T\right)  \right)  =a_{\gamma}$ and $\operatorname*{sign}%
T=a_{\alpha}$). This proves Observation 2.]

\begin{statement}
\textit{Observation 3:} We have $\operatorname*{sign}T=0$ for all
$T\in\mathcal{X}$ satisfying $f\left(  T\right)  =T$.
\end{statement}

[\textit{Proof of Observation 3:} Let $T\in\mathcal{X}$ be such that $f\left(
T\right)  =T$.

We shall use all the notations that we have introduced in the proof of
Observation 2 above. In particular, we define two $N$-tuples $\alpha
\in\mathbb{N}^{N}$ and $\gamma\in\mathbb{N}^{N}$ by $\alpha:=\nu
+\operatorname*{cont}T+\rho$ and $\gamma:=\nu+\operatorname*{cont}\left(
T^{\ast}\right)  +\rho$, and we let $k$ be the smallest misstep of $T$.

The definition of $f$ yields $f\left(  T\right)  =T^{\ast}$, so that $T^{\ast
}=f\left(  T\right)  =T$. Thus,
\[
\gamma=\nu+\underbrace{\operatorname*{cont}\left(  T^{\ast}\right)
}_{=\operatorname*{cont}T}+\,\rho=\nu+\operatorname*{cont}T+\rho=\alpha.
\]
Hence, $\gamma_{k}=\alpha_{k}$. However, the equality
(\ref{pf.lem.stemb-lem.o2.pf.gammak=alk+1}) (which we have shown in the proof
of Observation 2) yields $\gamma_{k}=\alpha_{k+1}$. Comparing these two
equalities, we obtain $\alpha_{k}=\alpha_{k+1}$. Therefore, the $N$-tuple
$\alpha\in\mathbb{N}^{N}$ has two equal entries (namely, its $k$-th and its
$\left(  k+1\right)  $-st entry). Thus, Lemma \ref{lem.sf.alternant-0}
\textbf{(a)} yields $a_{\alpha}=0$.

However, the definition of $\operatorname*{sign}T$ yields%
\begin{align*}
\operatorname*{sign}T  &  =a_{\nu+\operatorname*{cont}T+\rho}=a_{\alpha
}\ \ \ \ \ \ \ \ \ \ \left(  \text{since }\nu+\operatorname*{cont}%
T+\rho=\alpha\right) \\
&  =0.
\end{align*}
This proves Observation 3.] \medskip

Now, let us combine what we have shown. We know that the map $f:\mathcal{X}%
\rightarrow\mathcal{X}$ is an involution (by Observation 1). Moreover, we have%
\[
\operatorname*{sign}\left(  f\left(  I\right)  \right)  =-\operatorname*{sign}%
I\ \ \ \ \ \ \ \ \ \ \text{for all }I\in\mathcal{X}%
\]
(by Observation 2, applied to $T=I$). Furthermore,%
\[
\operatorname*{sign}I=0\ \ \ \ \ \ \ \ \ \ \text{for all }I\in\mathcal{X}%
\text{ satisfying }f\left(  I\right)  =I
\]
(by Observation 3, applied to $T=I$). Therefore, Lemma \ref{lem.sign.cancel3}
(applied to the additive abelian group $\mathcal{P}$) yields%
\[
\sum_{I\in\mathcal{A}}\operatorname*{sign}I=\sum_{I\in\mathcal{A}%
\setminus\mathcal{X}}\operatorname*{sign}I.
\]
Renaming the summation index $I$ as $T$ on both sides of this equality, we
obtain%
\[
\sum_{T\in\mathcal{A}}\operatorname*{sign}T=\sum_{T\in\mathcal{A}%
\setminus\mathcal{X}}\operatorname*{sign}T.
\]
Hence, (\ref{pf.lem.stemb-lem.before-cancel-sign}) rewrites as
\begin{align*}
a_{\nu+\rho}\cdot s_{\lambda/\mu}  &  =\sum_{T\in\mathcal{A}\setminus
\mathcal{X}}\underbrace{\operatorname*{sign}T}_{\substack{=a_{\nu
+\operatorname*{cont}T+\rho}\\\text{(by the definition of }%
\operatorname*{sign}T\text{)}}}=\sum_{T\in\mathcal{A}\setminus\mathcal{X}%
}a_{\nu+\operatorname*{cont}T+\rho}\\
&  =\sum_{\substack{T\text{ is a }\nu\text{-Yamanouchi}\\\text{semistandard
tableau}\\\text{of shape }\lambda/\mu}}a_{\nu+\operatorname*{cont}T+\rho}%
\end{align*}
(since $\mathcal{A}\setminus\mathcal{X}$ is the set of all $\nu$-Yamanouchi
semistandard tableaux of shape $\lambda/\mu$\ \ \ \ \footnote{because
\[
\mathcal{A}=\operatorname*{SSYT}\left(  \lambda/\mu\right)
\ \ \ \ \ \ \ \ \ \ \text{and}\ \ \ \ \ \ \ \ \ \ \mathcal{X}=\left\{
T\in\operatorname*{SSYT}\left(  \lambda/\mu\right)  \ \mid\ T\text{ is not
}\nu\text{-Yamanouchi}\right\}
\]
}). This proves Lemma \ref{lem.sf.stemb-lem}.
\end{proof}

\subsubsection{The Pieri rules}

Having proved the Littlewood--Richardson rule, let us discuss a few more of
its consequences. As we know from Example \ref{exa.sf.schur-h-e}, the complete
homogeneous symmetric polynomials $h_{n}$ and the elementary symmetric
polynomials $e_{n}$ (for $n\in\left\{  0,1,\ldots,N\right\}  $) are instances
of Schur polynomials. Thus, by appropriately specializing Theorem
\ref{thm.sf.lr-zy}, we can obtain rules for expressing products of the form
$h_{n}s_{\mu}$ and $e_{n}s_{\mu}$ as sums of Schur polynomials. These rules
are known as the \emph{Pieri rules}. To formulate them, we need some more notations:

\begin{definition}
\label{def.sf.strips}Let $\lambda$ and $\mu$ be two $N$-partitions. \medskip

\textbf{(a)} We write $\lambda/\mu$ for the pair $\left(  \mu,\lambda\right)
$. Such a pair is called a \emph{skew partition}. \medskip

\textbf{(b)} We say that $\lambda/\mu$ is a \emph{horizontal strip} if we have
$\mu\subseteq\lambda$ and the Young diagram $Y\left(  \lambda/\mu\right)  $
has no two boxes lying in the same column. \medskip

\textbf{(c)} We say that $\lambda/\mu$ is a \emph{vertical strip} if we have
$\mu\subseteq\lambda$ and the Young diagram $Y\left(  \lambda/\mu\right)  $
has no two boxes lying in the same row. \medskip

Now, let $n\in\mathbb{N}$. \medskip

\textbf{(d)} We say that $\lambda/\mu$ is a \emph{horizontal }$n$\emph{-strip}
if $\lambda/\mu$ is a horizontal strip and satisfies $\left\vert Y\left(
\lambda/\mu\right)  \right\vert =n$. \medskip

\textbf{(e)} We say that $\lambda/\mu$ is a \emph{vertical }$n$\emph{-strip}
if $\lambda/\mu$ is a vertical strip and satisfies $\left\vert Y\left(
\lambda/\mu\right)  \right\vert =n$.
\end{definition}

\Needspace{40pc}

\begin{example}
Let $N=4$. \medskip

\textbf{(a)} If $\lambda=\left(  8,7,4,3\right)  $ and $\mu=\left(
7,4,4,1\right)  $, then we have $\mu\subseteq\lambda$, and the Young diagram
$Y\left(  \lambda/\mu\right)  $ looks as follows:%
\[
Y\left(  \lambda/\mu\right)  =\ydiagram{7+1,4+3,4+0,1+2}\ \ .
\]
From this picture, it is clear that this skew partition $\lambda/\mu$ is a
horizontal strip (and, in fact, a horizontal $6$-strip, since $\left\vert
Y\left(  \lambda/\mu\right)  \right\vert =6$), but not a vertical strip
(since, e.g., there are $3$ boxes in the second row of $Y\left(  \lambda
/\mu\right)  $). \medskip

\textbf{(b)} If $\lambda=\left(  3,3,2,1\right)  $ and $\mu=\left(
2,2,1,0\right)  $, then we have $\mu\subseteq\lambda$, and the Young diagram
$Y\left(  \lambda/\mu\right)  $ looks as follows:%
\[
Y\left(  \lambda/\mu\right)  =\ydiagram{2+1,2+1,1+1,0+1}\ \ .
\]
From this picture, it is clear that this skew partition $\lambda/\mu$ is a
vertical strip (and, in fact, a vertical $4$-strip, since $\left\vert Y\left(
\lambda/\mu\right)  \right\vert =4$), but not a horizontal strip (since there
are $2$ boxes in the third column of $Y\left(  \lambda/\mu\right)  $).
\medskip

\textbf{(c)} If $\lambda=\left(  4,3,1,1\right)  $ and $\mu=\left(
3,2,1,0\right)  $, then we have $\mu\subseteq\lambda$, and the Young diagram
$Y\left(  \lambda/\mu\right)  $ looks as follows:%
\[
Y\left(  \lambda/\mu\right)  =\ydiagram{3+1,2+1,1+0,0+1}\ \ .
\]
From this picture, it is clear that this skew partition $\lambda/\mu$ is both
a horizontal strip (and, in fact, a horizontal $3$-strip) and a vertical strip
(and, in fact, a vertical $3$-strip). \medskip

\textbf{(d)} If $\lambda=\left(  3,3,2,1\right)  $ and $\mu=\left(
1,1,1,1\right)  $, then we have $\mu\subseteq\lambda$, and the Young diagram
$Y\left(  \lambda/\mu\right)  $ looks as follows:%
\[
Y\left(  \lambda/\mu\right)  =\ydiagram{1+2,1+2,1+1,1+0}\ \ .
\]
From this picture, it is clear that this skew partition $\lambda/\mu$ is
neither a horizontal strip nor a vertical strip.
\end{example}

Horizontal and vertical strips can also be characterized in terms of the
entries of the partitions:

\begin{proposition}
\label{prop.sf.strips.entries}Let $\lambda=\left(  \lambda_{1},\lambda
_{2},\ldots,\lambda_{N}\right)  $ and $\mu=\left(  \mu_{1},\mu_{2},\ldots
,\mu_{N}\right)  $ be two $N$-partitions. \medskip

\textbf{(a)} The skew partition $\lambda/\mu$ is a horizontal strip if and
only if we have%
\[
\lambda_{1}\geq\mu_{1}\geq\lambda_{2}\geq\mu_{2}\geq\cdots\geq\lambda_{N}%
\geq\mu_{N}.
\]

\textbf{(b)} The skew partition $\lambda/\mu$ is a vertical strip if and only
if we have%
\[
\mu_{i}\leq\lambda_{i}\leq\mu_{i}+1\ \ \ \ \ \ \ \ \ \ \text{for each }%
i\in\left[  N\right]  .
\]

\end{proposition}

\begin{proof}
See Exercise \ref{exe.sf.strips.entries}.
\end{proof}

We can now state the \emph{Pieri rules}:

\begin{theorem}
[Pieri rules]\label{thm.sf.pieri}Let $n\in\mathbb{N}$. Let $\mu$ be an
$N$-partition. Then: \medskip

\textbf{(a)} We have%
\[
h_{n}s_{\mu}=\sum_{\substack{\lambda\text{ is an }N\text{-partition;}%
\\\lambda/\mu\text{ is a horizontal }n\text{-strip}}}s_{\lambda}.
\]

\textbf{(b)} We have%
\[
e_{n}s_{\mu}=\sum_{\substack{\lambda\text{ is an }N\text{-partition;}%
\\\lambda/\mu\text{ is a vertical }n\text{-strip}}}s_{\lambda}.
\]

\end{theorem}

\begin{example}
Let $N=4$ and $\mu=\left(  2,1,1,0\right)  $. Then: \medskip

\textbf{(a)} Theorem \ref{thm.sf.pieri} \textbf{(a)} (applied to $n=2$) yields%
\[
h_{2}s_{\mu}=\sum_{\substack{\lambda\text{ is an }N\text{-partition;}%
\\\lambda/\mu\text{ is a horizontal }2\text{-strip}}}s_{\lambda}=s_{\left(
2,2,1,1\right)  }+s_{\left(  3,1,1,1\right)  }+s_{\left(  3,2,1,0\right)
}+s_{\left(  4,1,1,0\right)  },
\]
since the $N$-partitions $\lambda$ for which $\lambda/\mu$ is a horizontal
$2$-strip are precisely the four $N$-partitions $\left(  2,2,1,1\right)  $,
$\left(  3,1,1,1\right)  $, $\left(  3,2,1,0\right)  $ and $\left(
4,1,1,0\right)  $. Here are the Young diagrams $Y\left(  \lambda\right)  $ of
these four $N$-partitions $\lambda$ (with the $Y\left(  \mu\right)  $
subdiagram colored red each time):%
\[
\ydiagram{2,2,1,1} *[*(red)]{2,1,1}\ \ ,\ \ \ \ \ \ \ \ \ \ \ydiagram{3,1,1,1}
*[*(red)]{2,1,1}\ \ ,\ \ \ \ \ \ \ \ \ \ \ydiagram{3,2,1,0} *[*(red)]{2,1,1}%
\ \ ,\ \ \ \ \ \ \ \ \ \ \ydiagram{4,1,1,0} *[*(red)]{2,1,1}\ \ .
\]

\textbf{(b)} Theorem \ref{thm.sf.pieri} \textbf{(b)} (applied to $n=2$)
yields
\[
e_{2}s_{\mu}=\sum_{\substack{\lambda\text{ is an }N\text{-partition;}%
\\\lambda/\mu\text{ is a vertical }2\text{-strip}}}s_{\lambda}=s_{\left(
2,2,1,1\right)  }+s_{\left(  2,2,2,0\right)  }+s_{\left(  3,1,1,1\right)
}+s_{\left(  3,2,1,0\right)  },
\]
since the $N$-partitions $\lambda$ for which $\lambda/\mu$ is a vertical
$2$-strip are precisely the four $N$-partitions $\left(  2,2,1,1\right)  $,
$\left(  2,2,2,0\right)  $, $\left(  3,1,1,1\right)  $ and $\left(
3,2,1,0\right)  $. Here are the Young diagrams $Y\left(  \lambda\right)  $ of
these four $N$-partitions $\lambda$ (with the $Y\left(  \mu\right)  $
subdiagram colored red each time):%
\[
\ydiagram{2,2,1,1} *[*(red)]{2,1,1}\ \ ,\ \ \ \ \ \ \ \ \ \ \ydiagram{2,2,2,0}
*[*(red)]{2,1,1}\ \ ,\ \ \ \ \ \ \ \ \ \ \ydiagram{3,1,1,1} *[*(red)]{2,1,1}%
\ \ ,\ \ \ \ \ \ \ \ \ \ \ydiagram{3,2,1,0} *[*(red)]{2,1,1}\ \ .
\]

\end{example}

\begin{proof}
[Proof of Theorem \ref{thm.sf.pieri}.]See Exercise \ref{exe.sf.pieri}.
\end{proof}

\subsubsection{The Jacobi--Trudi identities}

The \emph{Jacobi--Trudi identities} are determinantal formulas expressing a
skew Schur polynomial $s_{\lambda/\mu}$ in terms of complete homogeneous
symmetric polynomials $h_{n}$ or elementary symmetric polynomials $e_{n}$. We
begin with the former:

\begin{theorem}
[First Jacobi--Trudi formula]\label{thm.sf.jt-h}Let $M\in\mathbb{N}$. Let
$\lambda=\left(  \lambda_{1},\lambda_{2},\ldots,\lambda_{M}\right)  $ and
$\mu=\left(  \mu_{1},\mu_{2},\ldots,\mu_{M}\right)  $ be two $M$-partitions
(i.e., weakly decreasing $M$-tuples of nonnegative integers). Then,%
\[
s_{\lambda/\mu}=\det\left(  \left(  h_{\lambda_{i}-\mu_{j}-i+j}\right)
_{1\leq i\leq M,\ 1\leq j\leq M}\right)  .
\]

\end{theorem}

\begin{example}
If $M=3$, then Theorem \ref{thm.sf.jt-h} says that%
\[
s_{\lambda/\mu}=\det\left(  \left(  h_{\lambda_{i}-\mu_{j}-i+j}\right)
_{1\leq i\leq3,\ 1\leq j\leq3}\right)  =\det\left(
\begin{array}
[c]{ccc}%
h_{\lambda_{1}-\mu_{1}} & h_{\lambda_{1}-\mu_{2}+1} & h_{\lambda_{1}-\mu
_{3}+2}\\
h_{\lambda_{2}-\mu_{1}-1} & h_{\lambda_{2}-\mu_{2}} & h_{\lambda_{2}-\mu
_{3}+1}\\
h_{\lambda_{3}-\mu_{1}-2} & h_{\lambda_{3}-\mu_{2}-1} & h_{\lambda_{3}-\mu
_{3}}%
\end{array}
\right)  .
\]
For instance, if $\lambda=\left(  4,2,1\right)  $ and $\mu=\left(
1,0,0\right)  $, then%
\begin{align*}
s_{\lambda/\mu}  &  =\det\left(
\begin{array}
[c]{ccc}%
h_{4-1} & h_{4-0+1} & h_{4-0+2}\\
h_{2-1-1} & h_{2-0} & h_{2-0+1}\\
h_{1-1-2} & h_{1-0-1} & h_{1-0}%
\end{array}
\right)  =\det\left(
\begin{array}
[c]{ccc}%
h_{3} & h_{5} & h_{6}\\
h_{0} & h_{2} & h_{3}\\
h_{-2} & h_{0} & h_{1}%
\end{array}
\right) \\
&  =\det\left(
\begin{array}
[c]{ccc}%
h_{3} & h_{5} & h_{6}\\
1 & h_{2} & h_{3}\\
0 & 1 & h_{1}%
\end{array}
\right)  \ \ \ \ \ \ \ \ \ \ \left(  \text{since }h_{0}=1\text{ and }%
h_{-2}=0\right)  .
\end{align*}

\end{example}

The proof of Theorem \ref{thm.sf.jt-h} is not too hard using what we have
learnt about lattice paths in Section \ref{sec.det.comb.lgv}. Here is an
outline (with some details left to exercises):

\begin{proof}
[Proof of Theorem \ref{thm.sf.jt-h} (sketched).]The case $N=0$ is easy and
left to the reader. We thus WLOG assume that $N\geq1$.

Let us follow Convention \ref{conv.lgv.digraph}, Definition
\ref{def.lgv.lattice} and Definition \ref{def.lgv.path-tups}. We will work
with the digraph $\mathbb{Z}^{2}$. For each arc $a$ of the digraph
$\mathbb{Z}^{2}$, we define an element $w\left(  a\right)  \in\mathcal{P}$
(called the \emph{weight} of $a$) as follows:

\begin{itemize}
\item If $a$ is an east-step $\left(  i,j\right)  \rightarrow\left(
i+1,j\right)  $ with $j\in\left[  N\right]  $, then we set $w\left(  a\right)
:=x_{j}$.

\item If $a$ is any other arc, then we set $w\left(  a\right)  :=1$.
\end{itemize}

For each path $p$ of $\mathbb{Z}^{2}$, define the \emph{weight} $w\left(
p\right)  $ of $p$ by%
\[
w\left(  p\right)  :=\prod_{a\text{ is an arc of }p}w\left(  a\right)  .
\]

Now, it is not hard to see the following (compare with Proposition
\ref{prop.lgv.1-paths.ct}):

\begin{statement}
\textit{Observation 1:} Let $a$ and $c$ be two integers. Then,%
\[
\sum_{\substack{p\text{ is a path}\\\text{from }\left(  a,1\right)  \text{ to
}\left(  c,N\right)  }}w\left(  p\right)  =h_{c-a}.
\]

\end{statement}

See Exercise \ref{exe.sf.jt-h} \textbf{(a)} for a proof of Observation 1.

Next, for each path tuple $\mathbf{p}=\left(  p_{1},p_{2},\ldots,p_{k}\right)
$, let us define the \emph{weight} $w\left(  \mathbf{p}\right)  $ of
$\mathbf{p}$ by%
\[
w\left(  \mathbf{p}\right)  :=w\left(  p_{1}\right)  w\left(  p_{2}\right)
\cdots w\left(  p_{k}\right)  .
\]

Set $k:=M$ (for the sake of convenience). Thus, $\lambda=\left(  \lambda
_{1},\lambda_{2},\ldots,\lambda_{k}\right)  $ and $\mu=\left(  \mu_{1},\mu
_{2},\ldots,\mu_{k}\right)  $.

Define two $k$-vertices $\mathbf{A}=\left(  A_{1},A_{2},\ldots,A_{k}\right)  $
and $\mathbf{B}=\left(  B_{1},B_{2},\ldots,B_{k}\right)  $ by setting%
\[
A_{i}:=\left(  \mu_{i}-i,\ 1\right)  \ \ \ \ \ \ \ \ \ \ \text{and}%
\ \ \ \ \ \ \ \ \ \ B_{i}:=\left(  \lambda_{i}-i,\ N\right)
\ \ \ \ \ \ \ \ \ \ \text{for each }i\in\left[  k\right]  .
\]
It is easy to see that the conditions (\ref{eq.cor.lgv.kpaths.wt-np.xA}),
(\ref{eq.cor.lgv.kpaths.wt-np.yA}), (\ref{eq.cor.lgv.kpaths.wt-np.xB}) and
(\ref{eq.cor.lgv.kpaths.wt-np.yB}) of Corollary \ref{cor.lgv.kpaths.wt-np} are
satisfied. Hence, (\ref{eq.cor.lgv.kpaths.wt-np.claim}) yields%
\begin{equation}
\det\left(  \left(  \sum_{p:A_{i}\rightarrow B_{j}}w\left(  p\right)  \right)
_{1\leq i\leq k,\ 1\leq j\leq k}\right)  =\sum_{\substack{\mathbf{p}\text{ is
a nipat}\\\text{from }\mathbf{A}\text{ to }\mathbf{B}}}w\left(  \mathbf{p}%
\right)  ,\ \ \ \ \ \ \ \ \ \ \label{pf.thm.sf.jt-h.5}%
\end{equation}
where \textquotedblleft$p:A_{i}\rightarrow B_{j}$\textquotedblright\ means
\textquotedblleft$p$ is a path from $A_{i}$ to $B_{j}$\textquotedblright. The
left hand side of this equality is%
\begin{align}
&  \det\left(  \left(  \sum_{p:A_{i}\rightarrow B_{j}}w\left(  p\right)
\right)  _{1\leq i\leq k,\ 1\leq j\leq k}\right) \nonumber\\
&  =\det\left(  \left(  h_{\left(  \lambda_{j}-j\right)  -\left(  \mu
_{i}-i\right)  }\right)  _{1\leq i\leq k,\ 1\leq j\leq k}\right)
\ \ \ \ \ \ \ \ \ \ \left(  \text{by Observation 1}\right) \nonumber\\
&  =\det\left(  \left(  h_{\left(  \lambda_{i}-i\right)  -\left(  \mu
_{j}-j\right)  }\right)  _{1\leq i\leq k,\ 1\leq j\leq k}\right)
\ \ \ \ \ \ \ \ \ \ \left(  \text{by Theorem \ref{thm.det.transp}}\right)
\nonumber\\
&  =\det\left(  \left(  h_{\lambda_{i}-\mu_{j}-i+j}\right)  _{1\leq i\leq
k,\ 1\leq j\leq k}\right) \nonumber\\
&  =\det\left(  \left(  h_{\lambda_{i}-\mu_{j}-i+j}\right)  _{1\leq i\leq
M,\ 1\leq j\leq M}\right)  \label{pf.thm.sf.jt-h.6}%
\end{align}
(since $k=M$). We shall now analyze the right hand side of
(\ref{pf.thm.sf.jt-h.5}). To that purpose, we need to understand the nipats
from $\mathbf{A}$ to $\mathbf{B}$.

We define the \emph{height} of an east-step $\left(  i,j\right)
\rightarrow\left(  i+1,j\right)  $ to be the number $j$. We define the
\emph{height sequence} of a path $p$ to be the sequence of the heights of the
east-steps of $p$ (going from the starting point to the ending point of $p$).
For example, the path shown in Example \ref{exa.lgv.lattice.path.53} has
height sequence $\left(  1,1,1,2,3\right)  $. It is clear that the height
sequence of a path is always weakly increasing.

If $\mathbf{p}=\left(  p_{1},p_{2},\ldots,p_{k}\right)  $ is a nipat from
$\mathbf{A}$ to $\mathbf{B}$, we let $T\left(  \mathbf{p}\right)  $ be the
tableau of shape $Y\left(  \lambda/\mu\right)  $ such that the entries in the
$i$-th row of $T\left(  \mathbf{p}\right)  $ (for each $i\in\left[  k\right]
$) are the entries of the height sequence of $p_{i}$.

\Needspace{14pc}

\begin{example}
Let $N=6$ and $M=3$ (so that $k=M=3$) and $\lambda=\left(  4,2,1\right)  $ and
$\mu=\left(  1,0,0\right)  $. Here is a nipat $\mathbf{p}$ from $\mathbf{A}$
to $\mathbf{B}$, and the corresponding tableau $T\left(  \mathbf{p}\right)  $:%
\[
\mathbf{p}=%
\raisebox{-6.5pc}{
\begin{tikzpicture}
\draw[densely dotted] (-4.2,-0.2) grid (3.2, 5.2);
\draw[very thin] (-4.2,0) -- (3.2,0);
\draw[very thin] (-4.2,5) -- (3.2,5);
\node
[circle,fill=white,draw=black,text=black,inner sep=1pt] (A1) at (0,0) {$A_1$};
\node
[circle,fill=white,draw=black,text=black,inner sep=1pt] (A2) at (-2,0) {$A_2$}%
;
\node
[circle,fill=white,draw=black,text=black,inner sep=1pt] (A3) at (-3,0) {$A_3$}%
;
\node
[circle,fill=white,draw=black,text=black,inner sep=1pt] (B1) at (3,5) {$B_1$};
\node
[circle,fill=white,draw=black,text=black,inner sep=1pt] (B2) at (0,5) {$B_2$};
\node
[circle,fill=white,draw=black,text=black,inner sep=1pt] (B3) at (-2,5) {$B_3$}%
;
\begin{scope}[thick,>=stealth,darkred]
\draw(A1) edge[->] (1,0);
\draw(1,0) edge[->] (1,1);
\draw(1,1) edge[->] (2,1);
\draw(2,1) edge[->] (2,2);
\draw(2,2) edge[->] (2,3);
\draw(2,3) edge[->] (2,4);
\draw(2,4) edge[->] (3,4);
\draw(3,4) edge[->] (B1);
\draw(2,1.5) node[anchor=west] {$p_1$};
\end{scope}
\begin{scope}[thick,>=stealth,dbluecolor]
\draw(A2) edge[->] (-2,1);
\draw(-2,1) edge[->] (-1,1);
\draw(-1,1) edge[->] (-1,2);
\draw(-1,2) edge[->] (-1,3);
\draw(-1,3) edge[->] (-1,4);
\draw(-1,4) edge[->] (0,4);
\draw(0,4) edge[->] (B2);
\draw(-1,1.5) node[anchor=west] {$p_2$};
\end{scope}
\begin{scope}[thick,>=stealth,dgreencolor]
\draw(A3) edge[->] (-3,1);
\draw(-3,1) edge[->] (-3,2);
\draw(-3,2) edge[->] (-3,3);
\draw(-3,3) edge[->] (-2,3);
\draw(-2,3) edge[->] (-2,4);
\draw(-2,4) edge[->] (B3);
\draw(-3,2.5) node[anchor=east] {$p_3$};
\end{scope}
\end{tikzpicture}
}%
\ \ \Longrightarrow\ \ T\left(  \mathbf{p}\right)
=\raisebox{1.4pc}{\ytableaushort{\none 125,25,4}}\ \ .
\]

\end{example}

We could have defined the tableau $T\left(  \mathbf{p}\right)  $ just as
easily for any path tuple $\mathbf{p}$ from $\mathbf{A}$ to $\mathbf{B}$ (not
just for a nipat); however, the case of a nipat is particularly useful,
because it turns out that the tableau $T\left(  \mathbf{p}\right)  $ is
semistandard if and only if $\mathbf{p}$ is a nipat. Moreover, the following
stronger statement holds:

\begin{statement}
\textit{Observation 2:} There is a bijection%
\begin{align*}
\left\{  \text{nipats from }\mathbf{A}\text{ to }\mathbf{B}\right\}   &
\rightarrow\operatorname*{SSYT}\left(  \lambda/\mu\right)  ,\\
\mathbf{p}  &  \mapsto T\left(  \mathbf{p}\right)  .
\end{align*}

\end{statement}

See Exercise \ref{exe.sf.jt-h} \textbf{(b)} for a proof of Observation 2.

It is easy to see that $w\left(  \mathbf{p}\right)  =x_{T\left(
\mathbf{p}\right)  }$ for any nipat $\mathbf{p}$ from $\mathbf{A}$ to
$\mathbf{B}$. Hence,%
\begin{align*}
\sum_{\substack{\mathbf{p}\text{ is a nipat}\\\text{from }\mathbf{A}\text{ to
}\mathbf{B}}}\underbrace{w\left(  \mathbf{p}\right)  }_{=x_{T\left(
\mathbf{p}\right)  }}  &  =\sum_{\substack{\mathbf{p}\text{ is a
nipat}\\\text{from }\mathbf{A}\text{ to }\mathbf{B}}}x_{T\left(
\mathbf{p}\right)  }=\sum_{T\in\operatorname*{SSYT}\left(  \lambda/\mu\right)
}x_{T}\\
&  \ \ \ \ \ \ \ \ \ \ \ \ \ \ \ \ \ \ \ \ \left(
\begin{array}
[c]{c}%
\text{here, we have substituted }T\text{ for }T\left(  \mathbf{p}\right) \\
\text{in the sum, since the map in}\\
\text{Observation 2 is a bijection}%
\end{array}
\right) \\
&  =s_{\lambda/\mu}\ \ \ \ \ \ \ \ \ \ \left(  \text{since }s_{\lambda/\mu
}\text{ is defined to be }\sum_{T\in\operatorname*{SSYT}\left(  \lambda
/\mu\right)  }x_{T}\right)  .
\end{align*}
Thus,%
\begin{align*}
s_{\lambda/\mu}  &  =\sum_{\substack{\mathbf{p}\text{ is a nipat}\\\text{from
}\mathbf{A}\text{ to }\mathbf{B}}}w\left(  \mathbf{p}\right)  =\det\left(
\left(  \sum_{p:A_{i}\rightarrow B_{j}}w\left(  p\right)  \right)  _{1\leq
i\leq k,\ 1\leq j\leq k}\right)  \ \ \ \ \ \ \ \ \ \ \left(  \text{by
(\ref{pf.thm.sf.jt-h.5})}\right) \\
&  =\det\left(  \left(  h_{\lambda_{i}-\mu_{j}-i+j}\right)  _{1\leq i\leq
M,\ 1\leq j\leq M}\right)  \ \ \ \ \ \ \ \ \ \ \left(  \text{by
(\ref{pf.thm.sf.jt-h.6})}\right)  .
\end{align*}
This proves Theorem \ref{thm.sf.jt-h}.
\end{proof}

Our above proof of Theorem \ref{thm.sf.jt-h} is essentially taken from
\cite[First proof of Theorem 7.16.1]{Stanley-EC2}; other proofs can be found
in \cite[Exercise 2.7.13]{GriRei} (see also \cite[paragraph after Theorem
2.4.6]{GriRei} for several references).

The \emph{second Jacobi--Trudi formula} involves elementary symmetric
polynomials $e_{n}$ (instead of $h_{n}$) and transpose partitions (as in
Exercise \ref{exe.pars.transpose}):

\begin{theorem}
[Second Jacobi--Trudi formula]\label{thm.sf.jt-e}Let $\lambda$ and $\mu$ be
two partitions. Let $\lambda^{t}$ and $\mu^{t}$ be the transposes of $\lambda$
and $\mu$. Let $M\in\mathbb{N}$ be such that both $\lambda^{t}$ and $\mu^{t}$
have length $\leq M$. We extend the partitions $\lambda^{t}$ and $\mu^{t}$ to
$M$-tuples (by inserting zeroes at the end). Write these $M$-tuples
$\lambda^{t}$ and $\mu^{t}$ as $\lambda^{t}=\left(  \lambda_{1}^{t}%
,\lambda_{2}^{t},\ldots,\lambda_{M}^{t}\right)  $ and $\mu^{t}=\left(  \mu
_{1}^{t},\mu_{2}^{t},\ldots,\mu_{M}^{t}\right)  $. Then,%
\[
s_{\lambda/\mu}=\det\left(  \left(  e_{\lambda_{i}^{t}-\mu_{j}^{t}%
-i+j}\right)  _{1\leq i\leq M,\ 1\leq j\leq M}\right)  .
\]

\end{theorem}

\begin{example}
If $M=3$, then Theorem \ref{thm.sf.jt-e} says that%
\[
s_{\lambda/\mu}=\det\left(  \left(  e_{\lambda_{i}^{t}-\mu_{j}^{t}%
-i+j}\right)  _{1\leq i\leq M,\ 1\leq j\leq M}\right)  =\det\left(
\begin{array}
[c]{ccc}%
e_{\lambda_{1}^{t}-\mu_{1}^{t}} & e_{\lambda_{1}^{t}-\mu_{2}^{t}+1} &
e_{\lambda_{1}^{t}-\mu_{3}^{t}+2}\\
e_{\lambda_{2}^{t}-\mu_{1}^{t}-1} & e_{\lambda_{2}^{t}-\mu_{2}^{t}} &
e_{\lambda_{2}^{t}-\mu_{3}^{t}+1}\\
e_{\lambda_{3}^{t}-\mu_{1}^{t}-2} & e_{\lambda_{3}^{t}-\mu_{2}^{t}-1} &
e_{\lambda_{3}^{t}-\mu_{3}^{t}}%
\end{array}
\right)  .
\]
For instance, if $\lambda=\left(  3,2,2\right)  $ and $\mu=\left(
1,1,0\right)  $, then $\lambda^{t}=\left(  3,3,1\right)  $ and $\mu
^{t}=\left(  2\right)  =\left(  2,0,0\right)  $ (here, we have extended the
partition $\mu^{t}$ to an $M$-tuple by inserting zeroes at the end), so that
this becomes
\begin{align*}
s_{\lambda/\mu}  &  =\det\left(
\begin{array}
[c]{ccc}%
e_{3-2} & e_{3-0+1} & e_{3-0+2}\\
e_{3-2-1} & e_{3-0} & e_{3-0+1}\\
e_{1-2-2} & e_{1-0-1} & e_{1-0}%
\end{array}
\right)  =\det\left(
\begin{array}
[c]{ccc}%
e_{1} & e_{4} & e_{5}\\
e_{0} & e_{3} & e_{4}\\
e_{-3} & e_{0} & e_{1}%
\end{array}
\right) \\
&  =\det\left(
\begin{array}
[c]{ccc}%
e_{1} & e_{4} & e_{5}\\
1 & e_{3} & e_{4}\\
0 & 1 & e_{1}%
\end{array}
\right)  \ \ \ \ \ \ \ \ \ \ \left(  \text{since }e_{0}=1\text{ and }%
e_{-3}=0\right)  .
\end{align*}

\end{example}

\begin{proof}
[Proof of Theorem \ref{thm.sf.jt-e}.]See Exercise \ref{exe.sf.jt-e}.
\end{proof}

\newpage

\appendix

\newpage

\section{\label{chap.hw}Homework exercises}

What follows is a collection of problems (of varying difficulty) that are
meant to illuminate, expand upon and otherwise complement the above text.

The numbers in the squares (like \fbox{3}) are the experience points you gain
for solving the problems. They are a mix of difficulty rating and relevance
score: The harder or more important the problem, the larger is the number in
the square. I believe a \fbox{5} represents a good graduate-level homework
problem that requires thinking and work. A \fbox{3} usually requires some
thinking \textbf{or} work. A \fbox{1} is a warm-up question. A \fbox{7} should
be somewhat too hard for regular homework. Anything above \fbox{10} is not
really meant as homework, but I'd be excited to hear your ideas. Multi-part
exercises sometimes have points split between the parts -- i.e., if parts
\textbf{(b)} and \textbf{(c)} of an exercise are solved using the same idea,
then they may both be assigned \fbox{3} points even if each for itself would
be a \fbox{5}.

In solving an exercise, you can freely use (without proof) the claims of all
exercises above it.

Your goal (for an A grade in the 2024 iteration of Math 531) is to gain at
least 20 experience points from each of the Chapters 3--7 (counting Chapter 2
as part of Chapter 3).

\subsection{Before we start...}

\subsubsection{Binomial coefficients and elementary counting}

\begin{definition}
\label{def.floor-ceil}Let $x\in\mathbb{R}$. Then:

\begin{itemize}
\item We let $\left\lfloor x\right\rfloor $ denote the largest integer $\leq
x$. This integer $\left\lfloor x\right\rfloor $ is called the \emph{floor} of
$x$, or the result of \textquotedblleft\emph{rounding down} $x$%
\textquotedblright.

\item We let $\left\lceil x\right\rceil $ denote the smallest integer $\geq
x$. This integer $\left\lceil x\right\rceil $ is called the \emph{ceiling} of
$x$, or the result of \textquotedblleft\emph{rounding up} $x$%
\textquotedblright.
\end{itemize}
\end{definition}

\begin{example}
We have
\begin{align*}
\left\lfloor 3\right\rfloor  &  =3;\ \ \ \ \ \ \ \ \ \ \left\lfloor \sqrt
{2}\right\rfloor =1;\ \ \ \ \ \ \ \ \ \ \left\lfloor \pi\right\rfloor
=3;\ \ \ \ \ \ \ \ \ \ \left\lfloor -\pi\right\rfloor =-4;\\
\left\lceil 3\right\rceil  &  =3;\ \ \ \ \ \ \ \ \ \ \left\lceil \sqrt
{2}\right\rceil =2;\ \ \ \ \ \ \ \ \ \ \left\lceil \pi\right\rceil
=4;\ \ \ \ \ \ \ \ \ \ \left\lceil -\pi\right\rceil =-3.
\end{align*}

\end{example}

Let us note that each $x\in\mathbb{R}$ satisfies $\left\lfloor x\right\rfloor
\leq x\leq\left\lceil x\right\rceil $.

\begin{exercise}
\fbox{3} Let $n\in\mathbb{N}$. Prove that%
\[
\sum_{k=0}^{n}\dbinom{-2}{k}=\left(  -1\right)  ^{n}\left\lfloor \dfrac
{n+2}{2}\right\rfloor .
\]

\end{exercise}

\begin{noncompile}
(The notation $\left\lfloor x\right\rfloor $ stands for the \emph{floor} of a
real number $x$. This is defined as the largest integer that is $\leq x$. For
example, $\left\lfloor \pi\right\rfloor =\left\lfloor 3\right\rfloor =3$.)
\end{noncompile}

The next exercise is concerned with the notion of \emph{lacunar sets}. This
notion appears all over combinatorics (particularly in connection with
Fibonacci numbers), so chances are we will meet it again.

\begin{definition}
\label{def.lacunar.lac}A set $S$ of integers is said to be \emph{lacunar} if
it contains no two consecutive integers (i.e., there is no integer $i$ such
that both $i\in S$ and $i+1\in S$).
\end{definition}

For example, the set $\left\{  1,5,7\right\}  $ is lacunar, but $\left\{
1,5,6\right\}  $ is not. Any $1$-element subset of $\mathbb{Z}$ is lacunar,
and so is the empty set.

Some people say \textquotedblleft sparse\textquotedblright\ instead of
\textquotedblleft lacunar\textquotedblright, but the word \textquotedblleft
sparse\textquotedblright\ also has other meanings.

\begin{example}
The lacunar subsets of $\left\{  1,2,3,4,5\right\}  $ are%
\begin{align*}
&  \varnothing,\ \ \ \ \ \ \ \ \ \ \left\{  1\right\}
,\ \ \ \ \ \ \ \ \ \ \left\{  2\right\}  ,\ \ \ \ \ \ \ \ \ \ \left\{
3\right\}  ,\ \ \ \ \ \ \ \ \ \ \left\{  4\right\}
,\ \ \ \ \ \ \ \ \ \ \left\{  5\right\}  ,\ \ \ \ \ \ \ \ \ \ \left\{
1,3\right\}  ,\\
&  \left\{  1,4\right\}  ,\ \ \ \ \ \ \ \ \ \ \left\{  1,5\right\}
,\ \ \ \ \ \ \ \ \ \ \left\{  2,4\right\}  ,\ \ \ \ \ \ \ \ \ \ \left\{
2,5\right\}  ,\ \ \ \ \ \ \ \ \ \ \left\{  3,5\right\}
,\ \ \ \ \ \ \ \ \ \ \left\{  1,3,5\right\}  .
\end{align*}

\end{example}

\begin{exercise}
\label{exe.lacunar.counts}Let $n\in\mathbb{N}$. \medskip

\textbf{(a)} \fbox{1} Prove that the total \# of lacunar subsets of $\left\{
1,2,\ldots,n\right\}  $ is the Fibonacci number $f_{n+2}$. \medskip

\textbf{(b)} \fbox{1} Let $k\in\left\{  0,1,\ldots,n+1\right\}  $. Prove that
the total \# of $k$-element lacunar subsets of $\left\{  1,2,\ldots,n\right\}
$ equals $\dbinom{n+1-k}{k}$. \medskip

\textbf{(c)} \fbox{1} What goes wrong with the claim of part \textbf{(b)} if
$k>n+1$ ? \medskip

\textbf{(d)} \fbox{1} Find the largest possible size of a lacunar subset of
$\left\{  1,2,\ldots,n\right\}  $. \medskip

\textbf{(e)} \fbox{1} Prove that $f_{n+1}=\sum_{k=0}^{n}\dbinom{n-k}{k}$ for
each $n\in\left\{  -1,0,1,\ldots\right\}  $.
\end{exercise}

\begin{exercise}
\label{exe.binom.sum-n-choose-4k}Let $n$ be a positive integer. \medskip

\textbf{(a)} \fbox{2} Prove that%
\[
\dbinom{n}{0}+\dbinom{n}{2}+\dbinom{n}{4}+\dbinom{n}{6}+\cdots=\sum
_{k\in\mathbb{N}}\dbinom{n}{2k}=2^{n-1}.
\]

\textbf{(b)} \fbox{3} Prove that%
\[
\dbinom{n}{0}+\dbinom{n}{4}+\dbinom{n}{8}+\dbinom{n}{12}+\cdots=\sum
_{k\in\mathbb{N}}\dbinom{n}{4k}=2^{n-2}+2^{n/2-1}\cos\dfrac{\pi n}{4}.
\]

[\textbf{Hint:} For part \textbf{(a)}, compute $\left(  1+1\right)
^{n}+\left(  1-1\right)  ^{n}$. For part \textbf{(b)}, compute $\left(
1+1\right)  ^{n}+\left(  1+i\right)  ^{n}+\left(  1-1\right)  ^{n}+\left(
1-i\right)  ^{n}$, where $i=\sqrt{-1}\in\mathbb{C}$ is the imaginary unit.]
\end{exercise}

From now on, we shall use the so-called
\href{https://en.wikipedia.org/wiki/Iverson_bracket}{\emph{Iverson bracket
notation}}:

\begin{definition}
\label{def.iverson}If $\mathcal{A}$ is any logical statement, then we define
an integer $\left[  \mathcal{A}\right]  \in\left\{  0,1\right\}  $ by%
\[
\left[  \mathcal{A}\right]  =%
\begin{cases}
1, & \text{if }\mathcal{A}\text{ is true};\\
0, & \text{if }\mathcal{A}\text{ is false}.
\end{cases}
\]
For example, $\left[  1+1=2\right]  =1$ (since $1+1=2$ is true), whereas
$\left[  1+1=1\right]  =0$ (since $1+1=1$ is false).

If $\mathcal{A}$ is any logical statement, then the integer $\left[
\mathcal{A}\right]  $ is known as the \emph{truth value} of $\mathcal{A}$.
\end{definition}

\begin{exercise}
\label{exe.binom.pascal-matrices.LcLd}\textbf{(a)} \fbox{3} Prove that%
\[
\dbinom{n}{a}\dbinom{a}{b}=\dbinom{n}{b}\dbinom{n-b}{a-b}%
\ \ \ \ \ \ \ \ \ \ \text{for any }n,a,b\in\mathbb{C}.
\]

\textbf{(b)} \fbox{3} Let $N\in\mathbb{N}$.

For each $c\in\mathbb{C}$, let $L_{c}\in\mathbb{C}^{N\times N}$ be the
$N\times N$-matrix whose rows are indexed $0,1,\ldots,N-1$ and whose columns
are indexed $0,1,\ldots,N-1$, and whose $\left(  i,j\right)  $-th entry is
$\dbinom{i}{j}c^{i-j}$ for each $i,j\in\left\{  0,1,\ldots,N-1\right\}  $.
(The expression \textquotedblleft$\dbinom{i}{j}c^{i-j}$\textquotedblright%
\ should be understood as $0$ if $i<j$, even if $c$ itself is $0$.)

[For example, if $N=5$, then $L_{c}=\left(
\begin{array}
[c]{ccccc}%
1 & 0 & 0 & 0 & 0\\
1c & 1 & 0 & 0 & 0\\
1c^{2} & 2c & 1 & 0 & 0\\
1c^{3} & 3c^{2} & 3c & 1 & 0\\
1c^{4} & 4c^{3} & 6c^{2} & 4c & 1
\end{array}
\right)  $.]

Prove that%
\[
L_{c}L_{d}=L_{c+d}\ \ \ \ \ \ \ \ \ \ \text{for any }c,d\in\mathbb{C}.
\]

\textbf{(c)} \fbox{1} Prove that the matrices $L_{1}$ and $L_{-1}$ are
mutually inverse.
\end{exercise}

\begin{exercise}
Let $p$ be a prime number. \medskip

\textbf{(a)} \fbox{2} Prove that $p\mid\dbinom{p}{k}$ for each $k\in\left\{
1,2,\ldots,p-1\right\}  $. \medskip

\textbf{(b)} \fbox{3} Let $a,b\in\mathbb{N}$. Prove that%
\[
\dbinom{ap}{bp}\equiv\dbinom{a}{b}\operatorname{mod}p^{2}.
\]

\textbf{(c)} \fbox{3} Prove the claim of part \textbf{(b)} still holds if we
replace \textquotedblleft$a,b\in\mathbb{N}$\textquotedblright\ by
\textquotedblleft$a,b\in\mathbb{Z}$\textquotedblright. \medskip

[\textbf{Hint:} The following suggests a combinatorial solution (algebraic
solutions also exist).

Consider the cyclic group $C_{p}=\mathbb{Z}/p\mathbb{Z}$ with $p$ elements.

For part \textbf{(a)}, let $U$ be the set of all $p$-tuples of elements of
$\left\{  0,1\right\}  $ with the property that exactly $k$ entries of the
$p$-tuple are $1$. The group $C_{p}$ acts on $U$ by cyclic rotation. Argue
that each orbit of this action has size divisible by $p$.

For part \textbf{(b)}, let $W$ be the set of all $p\times a$-matrices with
entries in $\left\{  0,1\right\}  $ and having the property that the sum of
all entries of the matrix is $bp$ (that is, exactly $bp$ entries are $1$).
Construct an action of the group $C_{p}^{a}=\underbrace{C_{p}\times
C_{p}\times\cdots\times C_{p}}_{a\text{ times}}$ on $W$ in which the $k$-th
$C_{p}$ factor cyclically rotates the entries of the $k$-th column of the
matrix. Argue that all but $\dbinom{a}{b}$ orbits of this action have size
divisible by $p^{2}$, and conclude by writing $\left\vert W\right\vert $ as
the sum of the sizes of the orbits.

For part \textbf{(c)}, fix $b\in\mathbb{N}$ and $p$ (why is it enough to
consider $b\in\mathbb{N}$?), and show that the remainders of $\dbinom{ap}%
{bp}-\dbinom{a}{b}$ modulo $p^{2}$ are periodic as a function in $a$.]
\end{exercise}

\subsection{Generating functions}

The notations of Chapter \ref{chap.gf} shall be used here. In particular, we
fix a commutative ring $K$.

\subsubsection{Examples}

All the properties of generating functions that have been used without proof
in Section \ref{sec.gf.exas} can also be used in the following exercises.

\begin{exercise}
\fbox{2} The \emph{Lucas sequence} is the sequence $\left(  \ell_{0},\ell
_{1},\ell_{2},\ldots\right)  $ of integers defined recursively by%
\[
\ell_{0}=2,\ \ \ \ \ \ \ \ \ \ \ell_{1}=1,\ \ \ \ \ \ \ \ \ \ \ell_{n}%
=\ell_{n-1}+\ell_{n-2}\text{ for each }n\geq2.
\]
(Thus, $\ell_{2}=3$ and $\ell_{3}=4$ and $\ell_{4}=7$ and so on.)

Find an explicit formula for $\ell_{n}$ analogous to Binet's formula for the
Fibonacci numbers.
\end{exercise}

\begin{exercise}
\label{exe.gf.cn-easy-manip}\fbox{1} Prove that $\dfrac{1}{n+1}\dbinom{2n}%
{n}=\dbinom{2n}{n}-\dbinom{2n}{n-1}$ for any $n\in\mathbb{N}$.
\end{exercise}

\begin{exercise}
\fbox{3} Let $q$ and $d$ be any two real numbers. Let $\left(  a_{0}%
,a_{1},a_{2},\ldots\right)  $ be a sequence of real numbers such that each
$n\geq1$ satisfies $a_{n}=qa_{n-1}+d$. (This can be viewed as a common
generalization of arithmetic and geometric sequences.)

Find an explicit formula for $a_{n}$ in terms of $q$, $d$ and $a_{0}$. (The
formula may depend on whether $q$ is $1$ or not.)
\end{exercise}

\begin{exercise}
\fbox{5} Find and prove an explicit formula for the coefficient of $x^{n}$ in
the formal power series $\dfrac{1}{1-x-x^{2}+x^{3}}$.
\end{exercise}

\begin{exercise}
Recall the Fibonacci sequence $\left(  f_{0},f_{1},f_{2},\ldots\right)  $.
\medskip

\textbf{(a)} \fbox{2} Prove that%
\[
f_{0}+f_{2}x+f_{4}x^{2}+f_{6}x^{3}+\cdots=\sum_{k\geq0}f_{2k}x^{k}=\dfrac
{x}{x^{2}-3x+1}.
\]

\textbf{(b)} \fbox{2} Find a degree-$2$ linear recurrence relation for the
sequence $\left(  f_{0},f_{2},f_{4},f_{6},\ldots\right)  $. That is, find two
numbers $a$ and $b$ such that each $n\geq2$ satisfies $f_{2n}=af_{2n-2}%
+bf_{2n-4}$. \medskip

[\textbf{Hint:} For part \textbf{(a)}, start with the generating function
$F\left(  x\right)  $ from Section \ref{sec.gf.exas}, and compute the
\textquotedblleft average\textquotedblright\ $\dfrac{F\left(  x\right)
+F\left(  -x\right)  }{2}$ in two different ways: On the one hand, this
\textquotedblleft average\textquotedblright\ is $f_{0}+f_{2}x^{2}+f_{4}%
x^{4}+f_{6}x^{6}+\cdots$; on the other hand, it is a sum of two fractions.
Compare the results, and \textquotedblleft substitute $x^{1/2}$ for
$x$\textquotedblright\ (that is, replace each $x^{2n}$ by $x^{n}$).]
\end{exercise}

\begin{exercise}
\label{exe.fps.motzkin.1}A \emph{Motzkin word} of length $n$ (where
$n\in\mathbb{N}$) is an $n$-tuple whose entries belong to the set $\left\{
0,1,-1\right\}  $ and sum up to $0$ (that is, it contains equally many $1$'s
and $-1$'s), and that has the additional property that for each $k$, we have
\begin{align*}
&  \left(  \text{\# of }\left.  -1\right.  \text{'s among its first }k\text{
entries}\right) \\
&  \leq\left(  \text{\# of }1\text{'s among its first }k\text{ entries}%
\right)  .
\end{align*}
A \emph{Motzkin path} of length $n$ is a path from the point $\left(
0,0\right)  $ to the point $\left(  n,0\right)  $ in the Cartesian plane that
moves only using \textquotedblleft NE-steps\textquotedblright\ (i.e., steps of
the form $\left(  x,y\right)  \rightarrow\left(  x+1,y+1\right)  $),
\textquotedblleft SE-steps\textquotedblright\ (i.e., steps of the form
$\left(  x,y\right)  \rightarrow\left(  x+1,y-1\right)  $) and
\textquotedblleft E-steps\textquotedblright\ (i.e., steps of the form $\left(
x,y\right)  \rightarrow\left(  x+1,y\right)  $) and never falls below the
x-axis (i.e., does not contain any point $\left(  x,y\right)  $ with $y<0$).
For example, here is a Motzkin path from $\left(  0,0\right)  $ to $\left(
6,0\right)  $:%
\[
\text{%
\begin{tikzpicture}
\draw(0, 0) -- (1, 1) -- (2, 1) -- (3, 2) -- (4, 1) -- (5, 0) -- (6, 0);
\filldraw(0, 0) circle [fill=red, radius=0.1];
\filldraw(1, 1) circle [fill=red, radius=0.1];
\filldraw(2, 1) circle [fill=red, radius=0.1];
\filldraw(3, 2) circle [fill=red, radius=0.1];
\filldraw(4, 1) circle [fill=red, radius=0.1];
\filldraw(5, 0) circle [fill=red, radius=0.1];
\filldraw(6, 0) circle [fill=red, radius=0.1];
\end{tikzpicture}%
\ \ .}%
\]

For each $n\in\mathbb{N}$, we define the \emph{Motzkin number} $m_{n}$ by%
\[
m_{n}:=\left(  \text{\# of Motzkin paths from }\left(  0,0\right)  \text{ to
}\left(  n,0\right)  \right)  .
\]
Here is a table of the first $12$ Motzkin numbers $m_{n}$:%
\[%
\begin{tabular}
[c]{|c||c|c|c|c|c|c|c|c|c|c|c|c|}\hline
$n$ & $0$ & $1$ & $2$ & $3$ & $4$ & $5$ & $6$ & $7$ & $8$ & $9$ & $10$ &
$11$\\\hline
$m_{n}$ & $1$ & $1$ & $2$ & $4$ & $9$ & $21$ & $51$ & $127$ & $323$ & $835$ &
$2188$ & $5798$\\\hline
\end{tabular}
\ \ \ \ \ \ .
\]

\textbf{(a)} \fbox{1} Prove that $m_{n}$ is also the \# of Motzkin words of
length $n$ for each $n\in\mathbb{N}$. \medskip

\textbf{(b)} \fbox{2} Let $c_{n}$ be the $n$-th Catalan number (defined in
Section \ref{sec.gf.exas}) for each $n\in\mathbb{N}$. Prove that
\[
m_{n}=\sum_{k=0}^{n}\dbinom{n}{2k}c_{k}\ \ \ \ \ \ \ \ \ \ \text{for each
}n\in\mathbb{N}.
\]
(The sum here could just as well range from $k=0$ to $\left\lfloor
n/2\right\rfloor $ or range over all $k\in\mathbb{N}$, since $\dbinom{n}%
{2k}=0$ when $2k>n$.) \medskip

\textbf{(c)} \fbox{2} Prove that%
\[
m_{n}=\sum_{k=0}^{n}\dfrac{1}{k+1}\dbinom{n}{k}\dbinom{n-k}{k}=\dfrac{1}%
{n+1}\sum_{k=0}^{n}\dbinom{n+1}{k+1}\dbinom{n-k}{k}%
\]
for each $n\in\mathbb{N}$. \medskip

[Despite the analogy between the Motzkin numbers $m_{n}$ and the Catalan
numbers $c_{n}$, there is no formula for $m_{n}$ as simple as
(\ref{eq.sec.gf.exas.cn=1/(n+1)}) or (\ref{eq.sec.gf.exas.cn=diff}).]
\end{exercise}

\subsubsection{Definitions}

\begin{exercise}
\label{exe.fps.altern-quasivdm}\fbox{4} Let $n\in\mathbb{N}$ and
$m\in\mathbb{C}$. Prove that%
\[
\sum_{k=0}^{n}\left(  -1\right)  ^{k}\dbinom{m}{k}\dbinom{m}{n-k}=%
\begin{cases}
\left(  -1\right)  ^{n/2}\dbinom{m}{n/2}, & \text{if }n\text{ is even};\\
0, & \text{if }n\text{ is odd.}%
\end{cases}
\]

\end{exercise}

\begin{exercise}
\label{exe.fps.intdom}\fbox{3} Recall that a commutative ring $L$ is said to
be an \emph{integral domain} if it is nontrivial (i.e., its zero and its unity
are distinct) and has the property that if $a,b\in L$ satisfy $ab=0$, then
$a=0$ or $b=0$.

Let $K$ be an integral domain. Prove that the ring $K\left[  \left[  x\right]
\right]  $ is an integral domain. \medskip

[\textbf{Hint:} The analogous fact for the polynomial ring $K\left[  x\right]
$ is well-known. It is commonly proved by noticing that if two polynomials are
nonzero, then the leading term of their product equals the product of their
leading terms. This argument does not immediately apply to FPSs, since nonzero
FPSs usually have no leading term. What do nonzero FPSs have, though?]
\end{exercise}

\begin{exercise}
An FPS $a\in K\left[  \left[  x\right]  \right]  $ will be called \emph{even}
if it satisfies $\left[  x^{1}\right]  a=\left[  x^{3}\right]  a=\left[
x^{5}\right]  a=\cdots=0$ (that is, if it satisfies $\left[  x^{n}\right]
a=0$ for all odd $n\in\mathbb{N}$).

Let $f\in K\left[  \left[  x\right]  \right]  $ be any FPS. Write $f$ in the
form $f=\sum_{n\in\mathbb{N}}f_{n}x^{n}$ with $f_{0},f_{1},f_{2},\ldots\in K$
(so that $f_{n}=\left[  x^{n}\right]  f$ for all $n\in\mathbb{N}$). We define
$\widetilde{f}$ to be the FPS $\sum_{n\in\mathbb{N}}f_{n}\left(  -x\right)
^{n}=\sum_{n\in\mathbb{N}}\left(  -1\right)  ^{n}f_{n}x^{n}$. (Using the
notations of Definition \ref{def.fps.subs}, this FPS $\widetilde{f}$ is the
composition $f\left[  -x\right]  =f\circ\left(  -x\right)  $.) \medskip

\textbf{(a)} \fbox{1} Show that the FPS $f+\widetilde{f}$ is even. \medskip

\textbf{(b)} \fbox{2} Show that the FPS $f\cdot\widetilde{f}$ is even.
\end{exercise}

\subsubsection{Dividing FPSs}

\begin{exercise}
\textbf{(a)} \fbox{4} Prove that%
\[
\sum_{n\in\mathbb{N}}\dfrac{2^{n}x^{2^{n}}}{1+x^{2^{n}}}=\dfrac{x}{1-x}.
\]
(In particular, show that the sum on the left hand side is well-defined.)
\medskip

\textbf{(b)} \fbox{3} For each positive integer $n$, let $\nu_{2}\left(
n\right)  $ be the highest $k\in\mathbb{N}$ such that $2^{k}\mid n$.
(Equivalently, $\nu_{2}\left(  n\right)  $ is the exponent with which $2$
appears in the prime factorization of $n$; when $n$ is odd, this is understood
to be $0$. For example, $\nu_{2}\left(  40\right)  =3$ and $\nu_{2}\left(
41\right)  =0$.)

Prove that%
\[
\sum_{n\in\mathbb{N}}\dfrac{x^{2^{n}}}{1-x^{2^{n}}}=\sum_{n>0}\left(  \nu
_{2}\left(  n\right)  +1\right)  x^{n}.
\]
(In particular, show that the sum on the left hand side is well-defined.)
\end{exercise}

\subsubsection{Polynomials}

The following exercise is a generalization of Binet's formula for the
Fibonacci sequence (Example 1 in Section \ref{sec.gf.exas}):

\begin{exercise}
Let $F$ be a field of characteristic $0$ (that is, a field that is a
$\mathbb{Q}$-algebra). Let $d$ be a positive integer, and let $p_{1}%
,p_{2},\ldots,p_{d}$ be $d$ elements of $F$. Let $p\in F\left[  x\right]  $ be
the polynomial $1-\sum_{i=1}^{d}p_{i}x^{i}$.

Let $\left(  a_{0},a_{1},a_{2},\ldots\right)  $ be a sequence of elements of
$F$ with the property that each integer $n\geq d$ satisfies%
\[
a_{n}=\sum_{i=1}^{d}p_{i}a_{n-i}.
\]
(Such a sequence is said to be a \emph{linearly recursive sequence with
constant coefficients}. For example, if $d=2$ and $p_{1}=1$ and $p_{2}=1$,
then each $n\geq2$ must satisfy $a_{n}=a_{n-1}+a_{n-2}$, that is, the
recursive equation of the Fibonacci sequence. Of course, the starting values
$a_{0},a_{1},\ldots,a_{d-1}$ of the sequence can be arbitrary.) \medskip

\textbf{(a)} \fbox{3} Prove that there is some polynomial $q\in F\left[
x\right]  $ of degree $<d$ (this allows $q=0$) such that
\[
a_{0}+a_{1}x+a_{2}x^{2}+a_{3}x^{3}+\cdots=\dfrac{q}{p}%
\ \ \ \ \ \ \ \ \ \ \text{in }F\left[  \left[  x\right]  \right]  .
\]

\textbf{(b)} \fbox{2} Assume that the polynomial $p\in F\left[  x\right]  $
can be factored as%
\[
p=\left(  1-r_{1}x\right)  \left(  1-r_{2}x\right)  \cdots\left(
1-r_{d}x\right)
\]
for some distinct elements $r_{1},r_{2},\ldots,r_{d}$ of $F$. Prove that there
exist $d$ scalars $\lambda_{1},\lambda_{2},\ldots,\lambda_{d}\in F$ such that
each $n\in\mathbb{N}$ satisfies%
\[
a_{n}=\sum_{i=1}^{d}\lambda_{i}r_{i}^{n}.
\]

\textbf{(c)} \fbox{3} Now, assume instead that the polynomial $p\in F\left[
x\right]  $ satisfies $\deg p=d$ and can be factored as%
\[
p=\left(  1-r_{1}x\right)  ^{m_{1}}\left(  1-r_{2}x\right)  ^{m_{2}}%
\cdots\left(  1-r_{k}x\right)  ^{m_{k}}%
\]
for some distinct nonzero elements $r_{1},r_{2},\ldots,r_{k}$ of $F$ and some
nonnegative integers $m_{1},m_{2},\ldots,m_{k}$. Prove that there exist $k$
polynomials $u_{1},u_{2},\ldots,u_{k}\in F\left[  x\right]  $ such that $\deg
u_{i}<m_{i}$ for each $i\in\left\{  1,2,\ldots,k\right\}  $, and such that
each $n\in\mathbb{N}$ satisfies%
\[
a_{n}=\sum_{i=1}^{k}u_{i}\left(  n\right)  r_{i}^{n}.
\]

\end{exercise}

The next exercise reveals an application of FPSs to number theory (more such
applications will appear later on):

\begin{exercise}
Let $p$ and $q$ be two coprime positive integers. We define the set%
\[
S\left(  p,q\right)  :=\left\{  ap+bq\ \mid\ \left(  a,b\right)  \in
\mathbb{N}\times\mathbb{N}\right\}  .
\]
(For example, if $p=3$ and $q=5$, then $S\left(  p,q\right)  =\left\{
0,3,5,6,8,9,10,11,\ldots\right\}  $, where the \textquotedblleft$\ldots
$\textquotedblright\ is saying that all integers $\geq8$ belong to $S\left(
p,q\right)  $. The set $S\left(  p,q\right)  $ can be viewed as the set of all
denominations that can be paid with $p$-cent coins and $q$-cent coins, without
getting change.) \medskip

\textbf{(a)} \fbox{3} Prove that%
\[
\sum_{n\in S\left(  p,q\right)  }x^{n}=\dfrac{1-x^{pq}}{\left(  1-x^{p}%
\right)  \left(  1-x^{q}\right)  }.
\]

\textbf{(b)} \fbox{3} Prove that every integer $n>pq-p-q$ belongs to $S\left(
p,q\right)  $, whereas the integer $pq-p-q$ itself does not. \medskip

[\textbf{Hint:} For part \textbf{(a)}, describe the coefficient of $x^{m}$ in
\[
\left(  1-x^{p}\right)  \left(  1-x^{q}\right)  \sum_{n\in S\left(
p,q\right)  }x^{n}=\sum_{n\in S\left(  p,q\right)  }\left(  x^{n}%
-x^{n+p}-x^{n+q}+x^{n+p+q}\right)
\]
in a form revealing that it is $0$ unless $m=0$ or $m=pq$. Now, part
\textbf{(b)} can be solved as follows: First, show that every sufficiently
high $n\in\mathbb{N}$ belongs to $S\left(  p,q\right)  $. Hence,
$\sum_{\substack{n\in\mathbb{N};\\n\notin S\left(  p,q\right)  }}x^{n}%
=\dfrac{1}{1-x}-\dfrac{1-x^{pq}}{\left(  1-x^{p}\right)  \left(
1-x^{q}\right)  }$ is a polynomial. Finding the largest integer that doesn't
belong to $S\left(  p,q\right)  $ means finding the degree of this polynomial.]
\end{exercise}

\begin{exercise}
Let $N\in\mathbb{N}$. Let $P_{N}$ denote the $\mathbb{C}$-vector space of all
polynomials $f\in\mathbb{C}\left[  x\right]  $ of degree $<N$. Consider the
matrices $L_{c}$ for all $c\in\mathbb{C}$ defined in Exercise
\ref{exe.binom.pascal-matrices.LcLd} \textbf{(b)}.

For each $c\in\mathbb{C}$, let $B_{c}$ be the basis $\left(  \left(
x-c\right)  ^{0},\left(  x-c\right)  ^{1},\ldots,\left(  x-c\right)
^{N-1}\right)  $ of $P_{N}$. (This is a basis, since it is the image of the
monomial basis $\left(  x^{0},x^{1},\ldots,x^{N-1}\right)  $ under the
\textquotedblleft substitute $x-c$ for $x$\textquotedblright\ automorphism.)
\medskip

\textbf{(a)} \fbox{2} Let $c,d\in\mathbb{C}$. Prove that $\left(
L_{c}\right)  ^{T}$ (that is, the transpose of $L_{c}$) is the change-of-basis
matrix from the basis $B_{c+d}$ to the basis $B_{d}$. (This means that
\[
\left(  x-d\right)  ^{j}=\sum_{i=0}^{N-1}\left(  \left(  L_{c}\right)
^{T}\right)  _{i,j}\left(  x-\left(  c+d\right)  \right)  ^{i}%
\ \ \ \ \ \ \ \ \ \ \text{for any }j\in\left\{  0,1,\ldots,N-1\right\}  ,
\]
where $\left(  \left(  L_{c}\right)  ^{T}\right)  _{i,j}$ denotes the $\left(
i,j\right)  $-th entry of the matrix $\left(  L_{c}\right)  ^{T}$.) \medskip

\textbf{(b)} \fbox{2} Use this to give a new solution to Exercise
\ref{exe.binom.pascal-matrices.LcLd} \textbf{(b)} (without using Exercise
\ref{exe.binom.pascal-matrices.LcLd} \textbf{(a)}).
\end{exercise}

\subsubsection{Substitution and evaluation of power series}

\begin{exercise}
\label{exe.fps.comp-moebius}\fbox{1} Let $K$ be a commutative ring. Let $a\in
K\left[  \left[  x\right]  \right]  $ be the FPS $\dfrac{x}{x-1}$. Prove that
$a\circ a=x$.
\end{exercise}

\begin{exercise}
\label{exe.fps.comp-inv}Let $K$ be a commutative ring. Let $a\in K\left[
\left[  x\right]  \right]  $ be an FPS such that $\left[  x^{0}\right]  a=0$.

A \emph{compositional inverse} of $a$ shall mean a FPS $b\in K\left[  \left[
x\right]  \right]  $ such that $\left[  x^{0}\right]  b=0$ and $a\circ b=x$
and $b\circ a=x$.

Prove the following: \medskip

\textbf{(a)} \fbox{1} If a compositional inverse of $a$ exists, then it is
unique. \medskip

\textbf{(b)} \fbox{4} A compositional inverse of $a$ exists if and only if
$\left[  x^{1}\right]  a$ is invertible in $K$.
\end{exercise}

\subsubsection{Derivatives of FPSs}

\begin{exercise}
\fbox{2} Let $f\in K\left[  \left[  x\right]  \right]  $ be an FPS. Let $p$
and $q$ be two coprime nonnegative integers. Prove that the coefficient
$\left[  x^{q}\right]  \left(  f^{p}\right)  $ is a multiple of $p$ (that is,
there exists some $c\in K$ such that $\left[  x^{q}\right]  \left(
f^{p}\right)  =pc$). \medskip

[\textbf{Hint:} More generally, prove that $q\cdot\left[  x^{q}\right]
\left(  f^{p}\right)  $ is a multiple of $p$ whether or not $p$ and $q$ are
coprime. (Think about the coefficients of $\left(  f^{p}\right)  ^{\prime}$.)]
\end{exercise}

The next exercise is concerned with generalizing the two equalities%
\begin{align*}
1+x+x^{2}+x^{3}+\cdots &  =\dfrac{1}{1-x}\ \ \ \ \ \ \ \ \ \ \text{and}\\
0+1x+2x^{2}+3x^{3}+\cdots &  =\dfrac{x}{\left(  1-x\right)  ^{2}}%
\end{align*}
that we have encountered in Section \ref{sec.gf.exas} (as
(\ref{eq.sec.gf.exas.1.1/1-x}) and (\ref{eq.sec.gf.exas.4.x/(1-x)2}), respectively).

\begin{exercise}
\label{exe.gf.eulerian-pol.basics}For any $m\in\mathbb{N}$, we define an FPS%
\[
Q_{m}:=\sum_{n\in\mathbb{N}}n^{m}x^{n}=0^{m}x^{0}+1^{m}x^{1}+2^{m}x^{2}%
+\cdots\in\mathbb{Z}\left[  \left[  x\right]  \right]  .
\]
For example,
\begin{align*}
Q_{0}  &  =x^{0}+x^{1}+x^{2}+x^{3}+\cdots=\dfrac{1}{1-x};\\
Q_{1}  &  =0x^{0}+1x^{1}+2x^{2}+3x^{3}+\cdots=\dfrac{x}{\left(  1-x\right)
^{2}}\ \ \ \ \ \ \ \ \ \ \left(  \text{by (\ref{eq.sec.gf.exas.4.x/(1-x)2}%
)}\right)  ;
\end{align*}
it can furthermore be shown that%
\begin{align*}
Q_{2}  &  =0x^{0}+1x^{1}+4x^{2}+9x^{3}+\cdots=\dfrac{x\left(  1+x\right)
}{\left(  1-x\right)  ^{3}};\\
Q_{3}  &  =0x^{0}+1x^{1}+8x^{2}+27x^{3}+\cdots=\dfrac{x\left(  1+4x+x^{2}%
\right)  }{\left(  1-x\right)  ^{4}};\\
Q_{4}  &  =0x^{0}+1x^{1}+16x^{2}+81x^{3}+\cdots=\dfrac{x\left(  1+11x+11x^{2}%
+x^{3}\right)  }{\left(  1-x\right)  ^{5}}.
\end{align*}
The expressions become more complicated as $m$ increases, but one will still
notice that each $Q_{m}$ has the form $\dfrac{A_{m}}{\left(  1-x\right)
^{m+1}}$, where $A_{m}$ is a polynomial of degree $m$ that has constant term
$0$ (unless $m=0$) and whose coefficients have a \textquotedblleft
palindromic\textquotedblright\ symmetry (in the sense that the sequence of
nonzero coefficients is symmetric across its middle). Let us prove this.

For each $m\in\mathbb{N}$, we define an FPS%
\[
A_{m}:=\left(  1-x\right)  ^{m+1}Q_{m}\in\mathbb{Z}\left[  \left[  x\right]
\right]  .
\]
(Thus, $Q_{m}=\dfrac{A_{m}}{\left(  1-x\right)  ^{m+1}}$, so that the $A_{m}$
we just defined are the $A_{m}$ we are interested in -- but we don't yet know
that they are polynomials.)

Let $\vartheta:\mathbb{Z}\left[  \left[  x\right]  \right]  \rightarrow
\mathbb{Z}\left[  \left[  x\right]  \right]  $ be the $\mathbb{Z}$-linear map
that sends each FPS $f\in\mathbb{Z}\left[  \left[  x\right]  \right]  $ to
$xf^{\prime}$. (That is, $\vartheta$ takes the derivative of an FPS and then
multiplies it by $x$.) \medskip

\textbf{(a)} \fbox{1} Prove that $\vartheta\left(  fg\right)  =\vartheta
\left(  f\right)  \cdot g+f\cdot\vartheta\left(  g\right)  $ for any
$f,g\in\mathbb{Z}\left[  \left[  x\right]  \right]  $. (In the lingo of
algebraists, this is saying that $\vartheta$ is a \emph{derivation} of
$\mathbb{Z}\left[  \left[  x\right]  \right]  $.) \medskip

\textbf{(b)} \fbox{1} Prove that $\vartheta\left(  \left(  1-x\right)
^{k}\right)  =-kx\left(  1-x\right)  ^{k-1}$ for each $k\in\mathbb{Z}$.
\medskip

\textbf{(c)} \fbox{1} Prove that $Q_{m}=\vartheta\left(  Q_{m-1}\right)  $ for
each $m>0$. \medskip

\textbf{(d)} \fbox{2} Prove that $A_{m}=mxA_{m-1}+x\left(  1-x\right)
A_{m-1}^{\prime}$ for each $m>0$. \medskip

\textbf{(e)} \fbox{1} Conclude that $A_{m}$ is a polynomial of degree $\leq m$
for each $m\in\mathbb{N}$. \medskip

\textbf{(f)} \fbox{1} Show that $\left[  x^{0}\right]  \left(  A_{m}\right)
=0$ for each $m>0$. \medskip

\textbf{(g)} \fbox{2} Show that $\left[  x^{i}\right]  \left(  A_{m}\right)
=\left(  m-i+1\right)  \left[  x^{i-1}\right]  \left(  A_{m-1}\right)
+i\left[  x^{i}\right]  \left(  A_{m-1}\right)  $ for each $m>0$ and each
$i>0$. \medskip

\textbf{(h)} \fbox{3} Show that $\left[  x^{i}\right]  \left(  A_{m}\right)
=\left[  x^{m+1-i}\right]  \left(  A_{m}\right)  $ for each $m>0$ and each
$i\in\left\{  0,1,\ldots,m+1\right\}  $. \medskip

The polynomials $A_{0},A_{1},A_{2},\ldots$ are known as the \emph{Eulerian
polynomials}.
\end{exercise}

The next exercise puts some of the theory of FPSs into a functional-analytic
context. Note that no functional analysis -- beyond a few basic definitions --
is needed to understand and to solve it, since our notion of convergence for
FPSs is a very simple one compared to what analysts deal with on a regular
basis. Yet the viewpoint can be helpful, particularly if you are familiar with
some functional analysis.

\begin{exercise}
\label{exe.fps.banach-space}For any nonzero FPS $f\in K\left[  \left[
x\right]  \right]  $, define the \emph{order} $\operatorname*{ord}\left(
f\right)  $ of $f$ to be the smallest $m\in\mathbb{N}$ such that $\left[
x^{m}\right]  f\neq0$. Further define the \emph{norm} $\left\vert \left\vert
f\right\vert \right\vert $ of an FPS $f\in K\left[  \left[  x\right]  \right]
$ to be the rational number $\dfrac{1}{2^{\operatorname*{ord}\left(  f\right)
}}$ if $f$ is nonzero. If $f$ is zero, set $\left\vert \left\vert f\right\vert
\right\vert :=0$.

This norm on $K\left[  \left[  x\right]  \right]  $ gives rise to a metric
$d:K\left[  \left[  x\right]  \right]  \times K\left[  \left[  x\right]
\right]  \rightarrow\mathbb{Q}$ on $K\left[  \left[  x\right]  \right]  $,
defined by%
\[
d\left(  f,g\right)  =\left\vert \left\vert f-g\right\vert \right\vert
\ \ \ \ \ \ \ \ \ \ \text{for any }f,g\in K\left[  \left[  x\right]  \right]
.
\]
This metric, in turn, induces a topology on $K\left[  \left[  x\right]
\right]  $.

[Note that the norm we have defined is \textbf{not} a norm in the sense of
functional analysis, even if $K$ is $\mathbb{R}$ or $\mathbb{C}$, since (for
example) $\left\vert \left\vert 2f\right\vert \right\vert $ equals $\left\vert
\left\vert f\right\vert \right\vert $ rather than $2\left\vert \left\vert
f\right\vert \right\vert $. However, it is used in the same way to define a
metric and thus a topology.] \medskip

\textbf{(a)} \fbox{3} Prove that $K\left[  \left[  x\right]  \right]  $ is a
\href{https://en.wikipedia.org/wiki/Complete_metric_space}{complete metric
space} with respect to this metric. \medskip

\textbf{(b)} \fbox{3} Prove that the maps%
\begin{align*}
K\left[  \left[  x\right]  \right]  \times K\left[  \left[  x\right]  \right]
&  \rightarrow K\left[  \left[  x\right]  \right]  ,\\
\left(  f,g\right)   &  \mapsto f+g
\end{align*}
and%
\begin{align*}
K\left[  \left[  x\right]  \right]  \times K\left[  \left[  x\right]  \right]
&  \rightarrow K\left[  \left[  x\right]  \right]  ,\\
\left(  f,g\right)   &  \mapsto fg
\end{align*}
and%
\begin{align*}
K\left[  \left[  x\right]  \right]  \times K\left[  \left[  x\right]  \right]
_{0}  &  \rightarrow K\left[  \left[  x\right]  \right]  ,\\
\left(  f,g\right)   &  \mapsto f\circ g
\end{align*}
are continuous with respect to the topologies induced by this metric. (Recall
that $K\left[  \left[  x\right]  \right]  _{0}$ denotes the subset of
$K\left[  \left[  x\right]  \right]  $ consisting of all FPSs $g\in K\left[
\left[  x\right]  \right]  $ satisfying $\left[  x^{0}\right]  g=0$. This
subset becomes a topological space by inheriting a subspace topology from
$K\left[  \left[  x\right]  \right]  $.) \medskip

\textbf{(c)} \fbox{1} Prove that the map%
\begin{align*}
K\left[  \left[  x\right]  \right]   &  \rightarrow K\left[  \left[  x\right]
\right]  ,\\
f  &  \mapsto f^{\prime}%
\end{align*}
is \href{https://en.wikipedia.org/wiki/Lipschitz_continuity}{Lipschitz
continuous} with
\href{https://en.wikipedia.org/wiki/Lipschitz_continuity}{Lipschitz constant}
$2$. \medskip

\textbf{(d)} \fbox{1} Assume that $K$ is a commutative $\mathbb{Q}$-algebra.
Let $\int$ denote the $K$-linear map from $K\left[  \left[  x\right]  \right]
$ to $K\left[  \left[  x\right]  \right]  $ that sends each FPS $\sum
_{n\in\mathbb{N}}a_{n}x^{n}$ to $\sum_{n\in\mathbb{N}}\dfrac{1}{n+1}%
a_{n}x^{n+1}$. (This map $\int$ is an algebraic analogue of the
antiderivative.) Prove that this map $\int$ is Lipschitz continuous with
Lipschitz constant $\dfrac{1}{2}$. \medskip

[\textbf{Hint:} For part \textbf{(b)}, the topology on the product of two
metric spaces is induced by the sup metric, which is given by%
\[
d_{\sup}\left(  \left(  f_{1},g_{1}\right)  ,\left(  f_{2},g_{2}\right)
\right)  =\max\left\{  d\left(  f_{1},f_{2}\right)  ,d\left(  g_{1}%
,g_{2}\right)  \right\}  .
\]
Show that all three maps are Lipschitz continuous with Lipschitz constant $1$
-- i.e., that any $\left(  f_{1},g_{1}\right)  $ and $\left(  f_{2}%
,g_{2}\right)  $ in the respective product spaces satisfy%
\begin{align*}
d\left(  f_{1}+g_{1},f_{2}+g_{2}\right)   &  \leq d_{\sup}\left(  \left(
f_{1},g_{1}\right)  ,\left(  f_{2},g_{2}\right)  \right)
\ \ \ \ \ \ \ \ \ \ \text{and}\\
d\left(  f_{1}g_{1},f_{2}g_{2}\right)   &  \leq d_{\sup}\left(  \left(
f_{1},g_{1}\right)  ,\left(  f_{2},g_{2}\right)  \right)
\ \ \ \ \ \ \ \ \ \ \text{and}\\
d\left(  f_{1}\circ g_{1},f_{2}\circ g_{2}\right)   &  \leq d_{\sup}\left(
\left(  f_{1},g_{1}\right)  ,\left(  f_{2},g_{2}\right)  \right)  .
\end{align*}
]
\end{exercise}

\begin{exercise}
\fbox{3} Let $K$ be a commutative $\mathbb{Q}$-algebra. Let $f\in K\left[
\left[  x\right]  \right]  $ be any FPS. Prove that there exists a
\textbf{unique} FPS $g\in K\left[  \left[  x\right]  \right]  $ satisfying
$\left[  x^{0}\right]  g=0$ and $g^{\prime}=f\circ g$. \medskip

[\textbf{Hint:} This is an algebraic version of local existence and uniqueness
of a solution of an ODE. There is an elementary recursive way to prove this.
However, a more elegant way is to rewrite the ODE $g^{\prime}=f\circ g$ as an
integral equation $g=\int\left(  f\circ g\right)  $, where $\int$ is the map
defined in Exercise \ref{exe.fps.banach-space} \textbf{(d)}. This integral
equation says that $g$ is a fixed point of the (nonlinear) operator $K\left[
\left[  x\right]  \right]  _{0}\rightarrow K\left[  \left[  x\right]  \right]
_{0},\ h\mapsto\int\left(  f\circ h\right)  $, where $K\left[  \left[
x\right]  \right]  _{0}$ is the set of all FPSs with constant term $0$. Now,
apply \href{https://en.wikipedia.org/wiki/Banach_fixed-point_theorem}{the
Banach fixed-point theorem}.]
\end{exercise}

\begin{exercise}
\label{exe.fps.stirling2.1}A \emph{set partition} of a set $U$ means a set
$\left\{  S_{1},S_{2},\ldots,S_{k}\right\}  $ of disjoint nonempty subsets
$S_{1},S_{2},\ldots,S_{k}$ of $U$ that satisfy $S_{1}\cup S_{2}\cup\cdots\cup
S_{k}=U$. These sets $S_{1},S_{2},\ldots,S_{k}$ are called the \emph{blocks}
(or the \emph{parts}) of this set partition.

For any $n,k\in\mathbb{N}$, we define $S\left(  n,k\right)  $ to be the number
of set partitions of $\left[  n\right]  $ that have $k$ blocks. For example,
$S\left(  4,3\right)  =6$, since the set partitions of $\left[  4\right]  $
that have $3$ blocks are%
\begin{align*}
&  \left\{  \left\{  1,2\right\}  ,\ \left\{  3\right\}  ,\ \left\{
4\right\}  \right\}  ,\ \ \ \ \ \ \ \ \ \ \left\{  \left\{  1,3\right\}
,\ \left\{  2\right\}  ,\ \left\{  4\right\}  \right\}
,\ \ \ \ \ \ \ \ \ \ \left\{  \left\{  1,4\right\}  ,\ \left\{  2\right\}
,\ \left\{  3\right\}  \right\}  ,\\
&  \left\{  \left\{  2,3\right\}  ,\ \left\{  1\right\}  ,\ \left\{
4\right\}  \right\}  ,\ \ \ \ \ \ \ \ \ \ \left\{  \left\{  2,4\right\}
,\ \left\{  1\right\}  ,\ \left\{  3\right\}  \right\}
,\ \ \ \ \ \ \ \ \ \ \left\{  \left\{  3,4\right\}  ,\ \left\{  1\right\}
,\ \left\{  2\right\}  \right\}  .
\end{align*}

The number $S\left(  n,k\right)  $ is called a \emph{Stirling number of the
2nd kind}. \medskip

\textbf{(a)} \fbox{1} Show that every positive integer $n$ satisfies
\begin{align*}
S\left(  n,n\right)   &  =1;\ \ \ \ \ \ \ \ \ \ S\left(  n,0\right)
=0;\ \ \ \ \ \ \ \ \ \ S\left(  n,1\right)  =1;\\
S\left(  n,2\right)   &  =2^{n-1}-1;\ \ \ \ \ \ \ \ \ \ S\left(  n,n-1\right)
=\dbinom{n}{2};\\
S\left(  n,k\right)   &  =0\ \ \ \ \ \ \ \ \ \ \text{for all }k>n.
\end{align*}
(The first and the last equalities here hold for $n=0$ as well. However,
$S\left(  0,0\right)  =1$ and $S\left(  0,1\right)  =0$.) \medskip

\textbf{(b)} \fbox{3} Show that every $n\in\mathbb{N}$ satisfies%
\[
x^{n}=\sum_{k=0}^{n}S\left(  n,k\right)  x^{\underline{k}}%
\ \ \ \ \ \ \ \ \ \ \text{in }K\left[  x\right]  ,
\]
where $x^{\underline{k}}$ denotes the polynomial
\[
\prod_{i=0}^{k-1}\left(  x-i\right)  =x\left(  x-1\right)  \left(  x-2\right)
\cdots\left(  x-k+1\right)  .
\]

\textbf{(c)} \fbox{1} Show that every $n\in\mathbb{N}$ satisfies%
\[
\sum_{k=0}^{n}\left(  -1\right)  ^{k}k!\cdot S\left(  n,k\right)  =\left(
-1\right)  ^{n}.
\]

\textbf{(d)} \fbox{2} Show that every $n\in\mathbb{N}$ satisfies%
\[
\sum_{k=1}^{n}\left(  -1\right)  ^{k-1}\left(  k-1\right)  !\cdot S\left(
n,k\right)  =%
\begin{cases}
1, & \text{if }n=1;\\
0, & \text{if }n\neq1.
\end{cases}
\]

[\textbf{Hint:} For part \textbf{(d)}, take the derivative of part
\textbf{(b)} and evaluate at $x=0$.]
\end{exercise}

\subsubsection{Exponentials and logarithms}

Recall that $K\left[  \left[  x\right]  \right]  _{1}=\left\{  f\in K\left[
\left[  x\right]  \right]  \ \mid\ \left[  x^{0}\right]  f=1\right\}  $, and
that $\operatorname*{loder}f=\dfrac{f^{\prime}}{f}$ for any $f\in K\left[
\left[  x\right]  \right]  _{1}$.

\begin{exercise}
\label{exe.fps.loder.2}\fbox{1} For this exercise, let $K$ be any commutative
ring (not necessarily a $\mathbb{Q}$-algebra). Let $f\in K\left[  \left[
x\right]  \right]  _{1}$ and $g\in K\left[  \left[  x\right]  \right]  _{0}$
be two FPSs. Note that $f\circ g\in K\left[  \left[  x\right]  \right]  _{1}$
(by Lemma \ref{lem.fps.Exp-Log-maps-wd} \textbf{(b)}), so that
$\operatorname*{loder}\left(  f\circ g\right)  $ is well-defined. Prove that%
\[
\operatorname*{loder}\left(  f\circ g\right)  =\left(  \left(
\operatorname*{loder}f\right)  \circ g\right)  \cdot g^{\prime}.
\]

\end{exercise}

\begin{exercise}
\label{exe.fps.stirling2.3}Recall the Stirling numbers of the 2nd kind
$S\left(  n,k\right)  $ defined in Exercise \ref{exe.fps.stirling2.1}.
\medskip

\textbf{(a)} \fbox{2} Show that every positive integers $n$ and $k$ satisfy%
\[
S\left(  n,k\right)  =k\cdot S\left(  n-1,k\right)  +S\left(  n-1,k-1\right)
.
\]

\textbf{(b)} \fbox{3} Show that every $k\in\mathbb{N}$ satisfies%
\[
\sum_{n\in\mathbb{N}}\dfrac{S\left(  n,k\right)  }{n!}x^{n}=\dfrac{\left(
\exp\left[  x\right]  -1\right)  ^{k}}{k!}\ \ \ \ \ \ \ \ \ \ \text{for
}K=\mathbb{Q}.
\]

[\textbf{Hint:} For part \textbf{(b)}, denote the left hand side by $f_{k}$,
and show that $f_{k}^{\prime}=kf_{k}+f_{k-1}$ for each $k>0$.]
\end{exercise}

\subsubsection{Non-integer powers}

\begin{exercise}
\label{exe.fps.not-squares}Let $K$ be a nontrivial commutative ring. \medskip

\textbf{(a)} \fbox{2} Prove that there exists no FPS $f\in K\left[  \left[
x\right]  \right]  $ such that $f^{2}=x$. (Do not assume that $K$ is a field
or an integral domain!) \medskip

\textbf{(b)} \fbox{2} More generally: Let $f\in K\left[  \left[  x\right]
\right]  $ and $n\in\mathbb{N}$. Assume that
\[
f^{n}=ax^{m}+\sum_{i>m}a_{i}x^{i}%
\]
for some $m\in\mathbb{N}$, some \textbf{invertible} $a\in K$ and some elements
$a_{m+1},a_{m+2},a_{m+3},\ldots\in K$. Prove that $n\mid m$. [This might
require a little bit of commutative algebra -- specifically the fact that any
nontrivial commutative ring has a maximal ideal.] \medskip

\textbf{(c)} \fbox{1} Now assume that $K=\mathbb{Z}/2$ is the field with $2$
elements. Prove that there exists no FPS $f\in K\left[  \left[  x\right]
\right]  $ such that $f^{2}=1+x$.
\end{exercise}

\begin{exercise}
\label{exe.fps.power-c.rules}\fbox{1} Prove Theorem
\ref{thm.fps.power-c.rules}.
\end{exercise}

\begin{exercise}
\fbox{5} Recall the Catalan numbers $c_{0},c_{1},c_{2},\ldots$ introduced in
Example 2 in Section \ref{sec.gf.exas}. Prove that%
\[
\sum_{k=0}^{n}c_{2k}c_{2\left(  n-k\right)  }=4^{n}c_{n}%
\ \ \ \ \ \ \ \ \ \ \text{for each }n\in\mathbb{N}.
\]

\end{exercise}

\begin{exercise}
\fbox{4} \textbf{(a)} Prove that there exists a unique sequence $\left(
a_{0},a_{1},a_{2},\ldots\right)  $ of rational numbers that satisfies
$a_{0}=1$ and%
\[
\sum_{k=0}^{n}a_{k}a_{n-k}=1\ \ \ \ \ \ \ \ \ \ \text{for all }n\in
\mathbb{N}.
\]

\textbf{(b)} Find an explicit formula for the $n$-th entry $a_{n}$ of this
sequence (in terms of binomial coefficients).
\end{exercise}

\subsubsection{Integer compositions}

\begin{exercise}
\fbox{3} Let $n$ be a positive integer. Let $m\in\mathbb{N}$. Prove that%
\[
\sum_{\substack{\left(  a_{1},a_{2},\ldots,a_{n}\right)  \in\mathbb{N}%
^{n};\\a_{1}+a_{2}+\cdots+a_{n}=m}}a_{1}a_{2}\cdots a_{n}=\dbinom{n+m-1}%
{2n-1}.
\]

\end{exercise}

\begin{exercise}
\label{exe.fps.comps.fibonaccis-3}\fbox{5} Let $n$ be a positive integer.
Recall the Fibonacci sequence $\left(  f_{0},f_{1},f_{2},\ldots\right)  $.
Prove that: \medskip

\textbf{(a)} The \# of compositions $\left(  \alpha_{1},\alpha_{2}%
,\ldots,\alpha_{m}\right)  $ of $n$ such that $\alpha_{1},\alpha_{2}%
,\ldots,\alpha_{m}$ are odd is $f_{n}$. \medskip

\textbf{(b)} The \# of compositions $\left(  \alpha_{1},\alpha_{2}%
,\ldots,\alpha_{m}\right)  $ of $n$ such that $\alpha_{i}\geq2$ for each
$i\in\left\{  1,2,\ldots,m\right\}  $ is $f_{n-1}$. \medskip

\textbf{(c)} The \# of compositions $\left(  \alpha_{1},\alpha_{2}%
,\ldots,\alpha_{m}\right)  $ of $n$ such that $\alpha_{i}\leq2$ for each
$i\in\left\{  1,2,\ldots,m\right\}  $ is $f_{n+1}$. \medskip
\end{exercise}

(Note that Exercise \ref{exe.fps.comps.fibonaccis-3} is behind many
appearances of the Fibonacci numbers in the research literature, e.g., in the
theory of \textquotedblleft peak algebras\textquotedblright.)

It is surprising that one and the same sequence (the Fibonacci sequence)
answers the three different counting questions in Exercise
\ref{exe.fps.comps.fibonaccis-3}. Even more surprisingly, this generalizes:

\begin{exercise}
\fbox{4} Let $n$ and $k$ be two positive integers such that $k>1$.

Let $u$ be the \# of compositions $\alpha$ of $n$ such that each entry of
$\alpha$ is congruent to $1$ modulo $k$.

Let $v$ be the \# of compositions $\beta$ of $n+k-1$ such that each entry of
$\beta$ is $\geq k$.

Let $w$ be the \# of compositions $\gamma$ of $n-1$ such that each entry of
$\gamma$ is either $1$ or $k$.

Prove that $u=v=w$.
\end{exercise}

\subsubsection{$x^{n}$-equivalence}

\begin{exercise}
\fbox{2} Assume that $K$ is a commutative $\mathbb{Q}$-algebra. Let
$n\in\mathbb{N}$. Let $c\in K$ and $a,b\in K\left[  \left[  x\right]  \right]
_{1}$ satisfy $a\overset{x^{n}}{\equiv}b$. Prove that $a^{c}\overset{x^{n}%
}{\equiv}b^{c}$. (See Definition \ref{def.fps.Exp-Log-maps} and Definition
\ref{def.fps.power-c} for the meanings of $K\left[  \left[  x\right]  \right]
_{1}$, $a^{c}$ and $b^{c}$.)
\end{exercise}

\begin{exercise}
\fbox{3} Let $a,b\in K\left[  \left[  x\right]  \right]  $ be two FPSs that
have compositional inverses. (See Exercise \ref{exe.fps.comp-inv} for the
meaning of \textquotedblleft compositional inverse\textquotedblright.) Let
$\widetilde{a}$ and $\widetilde{b}$ be the compositional inverses of $a$ and
$b$. Let $n\in\mathbb{N}$ be such that $a\overset{x^{n}}{\equiv}b$. Prove that
$\widetilde{a}\overset{x^{n}}{\equiv}\widetilde{b}$.
\end{exercise}

\begin{exercise}
\label{exe.fps.CMJ1232}\fbox{3} Consider the sequence of polynomials
$p_{0},p_{1},p_{2},\ldots\in\mathbb{Z}\left[  x\right]  $ defined recursively
by $p_{0}=x$ and
\[
p_{k}=p_{k-1}^{2}+x\ \ \ \ \ \ \ \ \ \ \text{for each }k\geq1.
\]
Let $k\in\mathbb{N}$. Prove that the first $k+2$ coefficients of $p_{k}$ agree
with the first $k+2$ coefficients of the FPS $\sum_{n\geq0}c_{n}x^{n+1}$,
where $c_{0},c_{1},c_{2},\ldots$ are the Catalan numbers (see Example 2 in
Section \ref{sec.gf.exas}). In other words, prove that $\left[  x^{0}\right]
p_{k}=0$ and%
\[
\left[  x^{i}\right]  p_{k}=c_{i-1}\ \ \ \ \ \ \ \ \ \ \text{for each }%
i\in\left[  k+1\right]  .
\]

\end{exercise}

\subsubsection{Infinite products}

\begin{exercise}
\fbox{2} Prove that each nonnegative integer can be written uniquely in the
form $\sum_{k\geq1}a_{k}\cdot k!$, for some sequence $\left(  a_{1}%
,a_{2},a_{3},\ldots\right)  $ of integers satisfying $\left(  0\leq a_{k}\leq
k\text{ for each }k\geq1\right)  $ and $\left(  a_{k}=0\text{ for all but
finitely many }k\geq1\right)  $. \medskip

[\textbf{Hint:} Simplify the FPS $\prod_{k\geq1}\left(  1+x^{k!}%
+x^{2k!}+\cdots+x^{k\cdot k!}\right)  $.]
\end{exercise}

\begin{exercise}
\label{exe.fps.witt-fac}\textbf{(a)} \fbox{1} Prove that the family $\left(
1-a_{i}x^{i}\right)  _{i\in\left\{  1,2,3,\ldots\right\}  }$ is multipliable
whenever $\left(  a_{1},a_{2},a_{3},\ldots\right)  \in K^{\left\{
1,2,3,\ldots\right\}  }$ is a sequence of elements of $K$. \medskip

\textbf{(b)} \fbox{3} Let $f\in K\left[  \left[  x\right]  \right]  $ be an
FPS with constant term $\left[  x^{0}\right]  f=1$. Prove that there is a
unique sequence $\left(  a_{1},a_{2},a_{3},\ldots\right)  \in K^{\left\{
1,2,3,\ldots\right\}  }$ such that%
\[
f=\prod_{i=1}^{\infty}\left(  1-a_{i}x^{i}\right)  .
\]

We call this sequence $\left(  a_{1},a_{2},a_{3},\ldots\right)  $ the
\emph{Witt coordinate sequence} of $f$. \medskip

\textbf{(c)} \fbox{1} Find the Witt coordinate sequence of the FPS $\dfrac
{1}{1-x}$.
\end{exercise}

\begin{exercise}
\label{exe.fps.necklace-fac}In this exercise, we set $K=\mathbb{Z}$. \medskip

\textbf{(a)} \fbox{1} Prove that the family $\left(  \left(  1-x^{i}\right)
^{-a_{i}}\right)  _{i\in\left\{  1,2,3,\ldots\right\}  }$ is multipliable
whenever $\left(  a_{1},a_{2},a_{3},\ldots\right)  \in\mathbb{Z}^{\left\{
1,2,3,\ldots\right\}  }$ is a sequence of integers. \medskip

\textbf{(b)} \fbox{1} Let $f\in\mathbb{Z}\left[  \left[  x\right]  \right]  $
be an FPS with constant term $\left[  x^{0}\right]  f=1$. Prove that there is
a unique sequence $\left(  a_{1},a_{2},a_{3},\ldots\right)  \in\mathbb{Z}%
^{\left\{  1,2,3,\ldots\right\}  }$ such that%
\[
f=\prod_{i=1}^{\infty}\left(  1-x^{i}\right)  ^{-a_{i}}.
\]

We call this sequence $\left(  a_{1},a_{2},a_{3},\ldots\right)  $ the
\emph{necklace coordinate sequence} of $f$. \medskip

\textbf{(c)} \fbox{1} Prove that the necklace coordinate sequence of $1+x$ is
$\left(  1,-1,0,0,0,\ldots\right)  $ (with infinitely many zeroes after the
$-1$). \medskip

\textbf{(d)} \fbox{2} Let $\left(  a_{1},a_{2},a_{3},\ldots\right)  $ be the
necklace coordinate sequence of an FPS $f\in\mathbb{Z}\left[  \left[
x\right]  \right]  $, where $\left[  x^{0}\right]  f=1$. Prove that%
\[
n\cdot\left[  x^{n}\right]  f=\sum_{d=1}^{n}da_{d}\sum_{k=1}^{\left\lfloor
n/d\right\rfloor }\left[  x^{n-kd}\right]  f\ \ \ \ \ \ \ \ \ \ \text{for each
}n\in\mathbb{N}.
\]

\textbf{(e)} \fbox{1} What is the relation between the necklace coordinate
sequences of $f$ and $\dfrac{1}{f}$ ?\medskip

\textbf{(f)} \fbox{3} Let $q\in\mathbb{Z}$. Let $\left(  a_{1},a_{2}%
,a_{3},\ldots\right)  $ be the necklace coordinate sequence of the FPS
$\dfrac{1}{1-qx}$. Prove that each $n\geq1$ satisfies
\[
q^{n}=\sum_{d\mid n}da_{d}.
\]
Here and in the following, the symbol \textquotedblleft$\sum_{d\mid n}%
$\textquotedblright\ means a sum over all positive divisors $d$ of $n$.
\medskip

\textbf{(g)} \fbox{2} Let $\mu:\left\{  1,2,3,\ldots\right\}  \rightarrow
\mathbb{Z}$ be the \emph{number-theoretical M\"{o}bius function}; this is the
function given by
\[
\mu\left(  d\right)  :=%
\begin{cases}
\left(  -1\right)  ^{\left(  \text{\# of prime divisors of }d\right)  }, &
\text{if }d\text{ is squarefree};\\
0, & \text{if not}%
\end{cases}
\]
for all $d\geq1$ (where $d$ is said to be \emph{squarefree} if $d$ is a
product of \textbf{distinct} primes). Let $q$ and $\left(  a_{1},a_{2}%
,a_{3},\ldots\right)  $ be as in part \textbf{(f)}. Prove that
\[
a_{n}=\dfrac{1}{n}\sum_{d\mid n}\mu\left(  d\right)  q^{n/d}%
\ \ \ \ \ \ \ \ \ \ \text{for each }n\geq1.
\]

\textbf{(h)} \fbox{3} Now let $q\in\mathbb{N}$, and let $\left(  a_{1}%
,a_{2},a_{3},\ldots\right)  $ be as in part \textbf{(f)}. Let $n\geq1$ be a
positive integer. The cyclic group $C_{n}=\mathbb{Z}/n\mathbb{Z}$ acts on the
set $\left[  q\right]  ^{n}$ by rotating $n$-tuples (i.e., the generator of
$C_{n}$ sends any $n$-tuple $\left(  r_{1},r_{2},\ldots,r_{n}\right)
\in\left[  q\right]  ^{n}$ to $\left(  r_{2},r_{3},\ldots,r_{n},r_{1}\right)
$). We define a $q$\emph{-ary }$n$\emph{-necklace} (or $q$\emph{-ary necklace
of length }$n$) to mean an orbit of this action (i.e., an equivalence class of
$n$-tuples in $\left[  q\right]  ^{n}$ under rotation). An $n$-necklace is
said to be \emph{aperiodic} if it contains $n$ many $n$-tuples (i.e., if the
group $C_{n}$ acts freely on it). (For example, $\left\{  \left(
1,1,3\right)  ,\ \left(  1,3,1\right)  ,\ \left(  3,1,1\right)  \right\}  $ is
an aperiodic $q$-ary $3$-necklace for any $q\geq3$, while the $4$-necklace
$\left\{  \left(  1,2,1,2\right)  ,\ \left(  2,1,2,1\right)  \right\}  $ is
not aperiodic.)

Prove that $a_{n}$ is the \# of aperiodic $q$-ary $n$-necklaces for each
$n\geq1$.\medskip

\textbf{(i)} \fbox{3} Let $u,v\in\mathbb{Z}$. Let $\left(  a_{1},a_{2}%
,a_{3},\ldots\right)  $ and $\left(  b_{1},b_{2},b_{3},\ldots\right)  $ be the
necklace coordinate sequences of the FPSs $\dfrac{1}{1-ux}$ and $\dfrac
{1}{1-vx}$, respectively. Prove that%
\[
\prod_{i=1}^{\infty}\left(  1-vx^{i}\right)  ^{-a_{i}}=\prod_{i=1}^{\infty
}\left(  1-ux^{i}\right)  ^{-b_{i}}.
\]

\textbf{(j)} \fbox{3} Let $u,v\in\mathbb{Z}$. Let $\left(  a_{1},a_{2}%
,a_{3},\ldots\right)  $, $\left(  b_{1},b_{2},b_{3},\ldots\right)  $ and
$\left(  c_{1},c_{2},c_{3},\ldots\right)  $ be the necklace coordinate
sequences of the FPSs $\dfrac{1}{1-ux}$, $\dfrac{1}{1-vx}$ and $\dfrac
{1}{1-uvx}$, respectively. Prove that%
\[
c_{n}=\sum_{\substack{i,j\geq1;\\\operatorname{lcm}\left(  i,j\right)
=n}}\gcd\left(  i,j\right)  a_{i}b_{j}\ \ \ \ \ \ \ \ \ \ \text{for all }%
n\geq1.
\]

\textbf{(k)} \fbox{2} Let $f\in\mathbb{Z}\left[  \left[  x\right]  \right]  $
be an FPS with $\left[  x^{0}\right]  f=1$. Let $\left(  a_{1},a_{2}%
,a_{3},\ldots\right)  $ and $\left(  b_{1},b_{2},b_{3},\ldots\right)  $ be the
necklace coordinate sequences of the FPSs $f$ and $f\left[  -x\right]  $,
respectively. Prove that%
\[
b_{n}=%
\begin{cases}
-a_{n}, & \text{if }n\text{ is odd};\\
a_{n}, & \text{if }4\mid n;\\
a_{n}+a_{n/2}, & \text{if }n\equiv2\operatorname{mod}4
\end{cases}
\ \ \ \ \ \ \ \ \ \ \text{for all }n\geq1.
\]

[\textbf{Note:} We are working with $K=\mathbb{Z}$ in this exercise, since the
$\left(  -a_{i}\right)  $-th power $\left(  1-x^{i}\right)  ^{-a_{i}}$ would
not be defined for $a_{i}$ in an arbitrary commutative ring $K$. However,
Theorem \ref{thm.fps.gen-newton} shows that it is defined whenever $K$ is a
commutative $\mathbb{Q}$-algebra. More generally, it can be defined whenever
$K$ is a so-called \emph{binomial ring}, which essentially means that all
binomial coefficients $\dbinom{a}{k}$ with $a\in K$ and $k\in\mathbb{N}$ exist
in $K$ and satisfy the same relations as the usual binomial coefficients. See
\cite{Xantch17} for more about binomial rings.

The low points for part \textbf{(b)} are due to the analogy to Exercise
\ref{exe.fps.witt-fac}.

In part \textbf{(g)}, a crucial ingredient is the identity $\sum_{d\mid n}%
\mu\left(  d\right)  =%
\begin{cases}
1, & \text{if }n=1;\\
0, & \text{if }n>1
\end{cases}
$ for all integers $n\geq1$. This is an easy consequence of cancellation
arguments as in Section \ref{sec.sign.intro}; see, e.g., \cite[Theorem
4.29]{Loehr-BC}.

In view of part \textbf{(h)}, the original formula $\dfrac{1}{1-qx}%
=\prod_{i=1}^{\infty}\left(  1-x^{i}\right)  ^{-a_{i}}$ defining the necklace
coordinate sequence of $\dfrac{1}{1-qx}$ can be interpreted combinatorially:
It says that $q^{n}$ (for $q\in\mathbb{N}$ and $n\geq1$) equals the \# of
multisets of aperiodic $q$-ary necklaces whose lengths add up to $n$. This
equality can also be proved bijectively; two such bijections are given by the
Chen--Fox--Lyndon factorization \cite[Theorem 6.1.27]{GriRei} and the
Gessel--Reutenauer bijection \cite[Theorem 6.6.29]{GriRei}.

This exercise is a window into the theory of necklaces; see \cite[\S 6.6]%
{GriRei}, \cite[\S 17]{HW-Witt} and \cite{MetRot83} for more.]
\end{exercise}

Next come some more exercises on the technicalities of multipliability and
infinite products. The first one is a (partial) converse to Theorem
\ref{thm.fps.1+f-mulable}:

\begin{exercise}
\fbox{5} Let $\left(  \mathbf{a}_{i}\right)  _{i\in I}$ be a multipliable
family of FPSs such that each $\mathbf{a}_{i}$ is invertible (in $K\left[
\left[  x\right]  \right]  $). Prove that the family $\left(  \mathbf{a}%
_{i}-1\right)  _{i\in I}$ is summable.
\end{exercise}

\begin{exercise}
Let $\left(  \mathbf{a}_{i}\right)  _{i\in I}$ and $\left(  \mathbf{b}%
_{i}\right)  _{i\in I}$ be two families of FPSs. \medskip

\textbf{(a)} \fbox{1} If $\left(  \mathbf{a}_{i}\right)  _{i\in I}$ and
$\left(  \mathbf{b}_{i}\right)  _{i\in I}$ are multipliable, is it necessarily
true that the family $\left(  \mathbf{a}_{i}+\mathbf{b}_{i}\right)  _{i\in I}$
is multipliable? \medskip

\textbf{(b)} \fbox{1} If $\left(  \mathbf{a}_{i}\right)  _{i\in I}$ is
summable and $\left(  \mathbf{b}_{i}\right)  _{i\in I}$ is multipliable, is it
necessarily true that the family $\left(  \mathbf{a}_{i}+\mathbf{b}%
_{i}\right)  _{i\in I}$ is multipliable? \medskip

\textbf{(c)} \fbox{1} Does the answer to part \textbf{(b)} change if we
additionally assume that $\mathbf{b}_{i}$ is invertible for each $i\in I$ ?
\end{exercise}

\begin{exercise}
\label{exe.fps.prodrule-inf-inf-0}\fbox{5} Prove the following generalization
of Proposition \ref{prop.fps.prodrule-inf-inf}:

Let $I$ be a set. For any $i\in I$, let $S_{i}$ be a set. Set%
\[
\overline{S}=\left\{  \left(  i,k\right)  \ \mid\ i\in I\text{ and }k\in
S_{i}\text{ and }k\neq0\right\}  .
\]

For any $i\in I$ and any $k\in S_{i}$, let $p_{i,k}$ be an element of
$K\left[  \left[  x\right]  \right]  $. Assume that
\[
p_{i,0}=1\ \ \ \ \ \ \ \ \ \ \text{for any }i\in I\text{ satisfying }0\in
S_{i}.
\]
Assume further that the family $\left(  p_{i,k}\right)  _{\left(  i,k\right)
\in\overline{S}}$ is summable. Then, the product $\prod_{i\in I}\ \ \sum_{k\in
S_{i}}p_{i,k}$ is well-defined (i.e., the family $\left(  p_{i,k}\right)
_{k\in S_{i}}$ is summable for each $i\in I$, and the family $\left(
\sum_{k\in S_{i}}p_{i,k}\right)  _{i\in I}$ is multipliable), and we have%
\begin{equation}
\prod_{i\in I}\ \ \sum_{k\in S_{i}}p_{i,k}=\sum_{\substack{\left(
k_{i}\right)  _{i\in I}\in\prod_{i\in I}S_{i}\\\text{is essentially finite}%
}}\ \ \prod_{i\in I}p_{i,k_{i}}. \label{eq.exe.fps.prodrule-inf-inf-0.eq}%
\end{equation}

\end{exercise}

\begin{exercise}
\label{exe.fps.Exp-Log-infsum}\fbox{3} Prove Proposition
\ref{prop.fps.Exp-Log-infsum} and Proposition \ref{prop.fps.Exp-Log-infprod}.
\end{exercise}

\subsubsection{The generating function of a weighted set}

\begin{exercise}
\label{exe.fps.bitstrings-1x1}Let $G$ be the set of all nonempty bitstrings
that start and end with a $1$. (Thus, for example, it contains $\left(
1\right)  $ and $\left(  1,1\right)  $ and $\left(  1,0,1\right)  $ and
$\left(  1,0,1,1\right)  $, but not $\left(  0,1\right)  $ or $\left(
1,0\right)  $.)

Let $H$ be the set of all bitstrings in $G$ that contain no two consecutive
$1$'s. (For instance, $\left(  1,0,1\right)  $ belongs to $H$, but $\left(
1,1\right)  $ and $\left(  1,0,1,1\right)  $ do not.)

Both $G$ and $H$ become weighted sets (with the weight of a bitstring being
its length). \medskip

\textbf{(a)} \fbox{1} Prove that $\overline{G}=x+\dfrac{x^{2}}{1-2x}%
=\dfrac{x\left(  1-x\right)  }{1-2x}$. \medskip

\textbf{(b)} \fbox{3} Prove that $G\cong H^{1}+H^{2}+H^{3}+\cdots$. \medskip

\textbf{(c)} \fbox{1} Compute $\overline{H}$. \medskip

\textbf{(d)} \fbox{2} Recall the notion of a lacunar set (Definition
\ref{def.lacunar.lac}). Show that
\begin{align*}
&  \left(  \text{\# of bitstrings in }H\text{ having weight }n\right) \\
&  =\left(  \text{\# of lacunar subsets of }\left[  n-4\right]  \right)
\ \ \ \ \ \ \ \ \ \ \text{for every }n\geq3.
\end{align*}

\textbf{(e)} \fbox{2} Use this and Exercise \ref{exe.lacunar.counts}
\textbf{(a)} to recover the generating function (\ref{eq.sec.gf.exas.1.Fx=1})
of the Fibonacci sequence.
\end{exercise}

\begin{exercise}
\label{exe.fps.motzkin.2}Recall the Motzkin paths and the Motzkin numbers
$m_{n}$, defined in Exercise \ref{exe.fps.motzkin.1}. \medskip

\textbf{(a)} \fbox{3} Consider the weighted set
\[
M:=\left\{  \text{Motzkin paths from }\left(  0,0\right)  \text{ to }\left(
n,0\right)  \text{ with }n\in\mathbb{N}\right\}  ,
\]
where the weight $\left\vert p\right\vert $ of a Motzkin path $p$ from
$\left(  0,0\right)  $ to $\left(  n,0\right)  $ is defined to be $n$. Also
consider the weighted set $X=\left\{  1\right\}  $ (with weight $\left\vert
1\right\vert =1$) and the weighted set $\mathbb{N}$ (with weights given by
$\left\vert n\right\vert =n$ for each $n\in\mathbb{N}$). Note that
$\mathbb{N}\cong1+X+X^{2}+X^{3}+\cdots$ (an infinite disjoint union). Prove
that%
\[
M\cong\mathbb{N}\times\left(  1+X^{2}\times M\times M\right)  .
\]

\textbf{(b)} \fbox{2} Conclude that the weight generating function
$\overline{M}=\sum_{n\in\mathbb{N}}m_{n}x^{n}$ of $M$ is
\[
\overline{M}=\dfrac{1-x-\sqrt{\left(  1+x\right)  \left(  1-3x\right)  }%
}{2x^{2}}.
\]

\textbf{(c)} \fbox{2} Prove that%
\[
m_{n}=-\dfrac{1}{2}\sum_{k=0}^{n+2}\left(  -3\right)  ^{n+2-k}\dbinom{1/2}%
{k}\dbinom{1/2}{n+2-k}\ \ \ \ \ \ \ \ \ \ \text{for each }n\in\mathbb{N}.
\]

\end{exercise}

\begin{exercise}
\textbf{(a)} \fbox{3} Extend the analysis in Subsection
\ref{subsec.gf.weighted-set.domino} to domino tilings of height-$4$ rectangles
by defining the weighted set
\[
F:=\left\{  \text{\textbf{faultfree} domino tilings of }R_{n,4}\text{ with
}n\in\mathbb{N}\right\}
\]
and showing that $\overline{F}=x+x^{2}+\dfrac{x^{2}}{1-x^{2}}+2\cdot
\dfrac{x^{2}}{1-x}$. \medskip

\textbf{(b)} \fbox{3} Use this to find an explicit formula for $d_{n,4}$ (as
defined in Definition \ref{def.domino.shapes-and-tilings} \textbf{(e)}) that
uses only quadratic irrationalities. (Note that the formula will be rather
intricate and contain nested square roots.)
\end{exercise}

\begin{exercise}
\fbox{5} This exercise is about a variation on domino tilings. We define
\textquotedblleft shapes\textquotedblright\ and the specific shapes $R_{n,m}$
as in Definition \ref{def.domino.shapes-and-tilings}. Fix a positive integer
$k$.

A $k$\emph{-omino} means a size-$k$ shape of the form%
\begin{align*}
&  \left\{  \left(  i+1,j\right)  ,\ \left(  i+2,j\right)  ,\ \ldots,\ \left(
i+k,j\right)  \right\}  \text{ (a \textquotedblleft\emph{horizontal }%
}k\text{\emph{-omino}\textquotedblright)}\ \ \ \ \ \ \ \ \ \ \text{or}\\
&  \left\{  \left(  i,j+1\right)  ,\ \left(  i,j+2\right)  ,\ \ldots,\ \left(
i,j+k\right)  \right\}  \text{ (a \textquotedblleft\emph{vertical }%
}k\text{\emph{-omino}\textquotedblright)}%
\end{align*}
for some $\left(  i,j\right)  \in\mathbb{Z}^{2}$.

A $k$\emph{-omino tiling} of a shape $S$ is a set partition of $S$ into
$k$-ominos (i.e., a set of disjoint $k$-ominos whose union is $S$).

Prove that the shape $R_{n,m}$ has a $k$-omino tiling if and only if we have
$k\mid n$ or $k\mid m$. \medskip

[\textbf{Hint:} Consider $R_{n,m}$ as a weighted set, where the weight of a
square $\left(  i,j\right)  \in\mathbb{N}^{2}$ is defined to be $i+j$. If
$R_{n,m}$ has a $k$-omino tiling, then show that the weight generating
function $\overline{R_{n,m}}=\sum_{\left(  i,j\right)  \in R_{n,m}}x^{i+j}$
must be divisible by $1+x+x^{2}+\cdots+x^{k-1}$ (as a polynomial in
$\mathbb{Q}\left[  x\right]  $, for example). However, $\overline{R_{n,m}}$
has a simple form.]
\end{exercise}

\subsubsection{Limits of FPSs}

\begin{exercise}
\label{exe.fps.lim.sum-prod-quot}\textbf{(a)} \fbox{1} Prove Proposition
\ref{prop.fps.lim.deriv-lim}. \medskip

\textbf{(b)} \fbox{2} Prove Proposition \ref{prop.fps.lim.comp}. \medskip

\textbf{(c)} \fbox{1} Let $\left(  f_{0},f_{1},f_{2},\ldots\right)  $ be a
sequence of FPSs that have compositional inverses. (See Exercise
\ref{exe.fps.comp-inv} for the definition of a compositional inverse.) Let $f$
be a further FPS such that $\lim\limits_{i\rightarrow\infty}f_{i}=f$. Prove
that $f$ also has a compositional inverse, and furthermore that%
\[
\lim\limits_{i\rightarrow\infty}\widetilde{f_{i}}=\widetilde{f},
\]
where $\widetilde{g}$ denotes the compositional inverse of an FPS $g$.
\end{exercise}

\begin{exercise}
\textbf{(a)} \fbox{2} Let $f\in K\left[  \left[  x\right]  \right]  $ be any
FPS whose constant term $\left[  x^{0}\right]  f$ is nilpotent. (An element
$u\in K$ is said to be \emph{nilpotent} if there exists some $m\in\mathbb{N}$
such that $u^{m}=0$.) Prove that $\lim\limits_{i\rightarrow\infty}f^{i}=0$.
\medskip

\textbf{(b)} \fbox{2} Assume that $K$ is an integral domain. Let $g\in
K\left[  \left[  x\right]  \right]  $ be any FPS. Prove that $\lim
\limits_{i\rightarrow\infty}g^{i}$ exists if and only if $g=1$ or $\left[
x^{0}\right]  g=0$.
\end{exercise}

\Needspace{20pc}

\begin{exercise}
\fbox{5} Recall the Catalan numbers $c_{0},c_{1},c_{2},\ldots$ introduced in
Example 2 in Section \ref{sec.gf.exas}, and the corresponding FPS $C\left(
x\right)  =\sum_{n\in\mathbb{N}}c_{n}x^{n}$. Prove that%
\[
1-\dfrac{1}{C\left(  x\right)  }=\dfrac{x}{1-\dfrac{x}{1-\dfrac{x}{\ddots}}},
\]
where the continued fraction on the right hand side is to be understood as
\[
\lim\limits_{n\rightarrow\infty}\underbrace{\dfrac{x}{1-\dfrac{x}{1-\dfrac{x}{%
\begin{array}
[c]{ccc}%
1- &  & \\
& \ddots & \\
&  & -\dfrac{x}{1-x}%
\end{array}
}}}}_{\text{with }n\text{ layers}}.
\]

(This requires checking that the $n$-layered finite continued fractions are
well-defined and converge to a limit in $K\left[  \left[  x\right]  \right]  $.)
\end{exercise}

\begin{exercise}
\label{exe.fps.prods.telescope1}Let $\left(  a_{0},a_{1},a_{2},\ldots\right)
$ be an infinite sequence of FPSs in $K\left[  \left[  x\right]  \right]  $
such that $\lim\limits_{n\rightarrow\infty}a_{n}$ exists. Prove the following:
\medskip

\textbf{(a)} \fbox{2} We have%
\[
\sum_{i\in\mathbb{N}}\left(  a_{i}-a_{i+1}\right)  =a_{0}-\lim
\limits_{n\rightarrow\infty}a_{n}%
\]
(and, in particular, the family $\left(  a_{i}-a_{i+1}\right)  _{i\in
\mathbb{N}}$ is summable). \medskip

\textbf{(b)} \fbox{2} If all the FPSs $a_{1},a_{2},a_{3},\ldots$ are
invertible, then%
\[
\prod_{i\in\mathbb{N}}\dfrac{a_{i}}{a_{i+1}}=\dfrac{a_{0}}{\lim
\limits_{n\rightarrow\infty}a_{n}}%
\]
(and, in particular, the family $\left(  \dfrac{a_{i}}{a_{i+1}}\right)
_{i\in\mathbb{N}}$ is multipliable). \medskip

(These are infinite analogues of the telescope rules for sums and products.)
\end{exercise}

\subsubsection{Laurent power series}

\begin{exercise}
\fbox{2} While $K\left[  \left[  x^{\pm}\right]  \right]  $ is not a ring,
some elements of $K\left[  \left[  x^{\pm}\right]  \right]  $ can still be
multiplied. For instance, define three elements $a,b,c\in K\left[  \left[
x^{\pm}\right]  \right]  $ by%
\begin{align*}
a  &  =1+x^{-1}+x^{-2}+x^{-3}+\cdots,\\
b  &  =1-x,\\
c  &  =1+x+x^{2}+x^{3}+\cdots.
\end{align*}

\textbf{(a)} Find $ab$ and $bc$ and $a\left(  bc\right)  $ and $\left(
ab\right)  c$. \medskip

\textbf{(b)} Why is it not surprising that $a\left(  bc\right)  \neq\left(
ab\right)  c$ ?
\end{exercise}

\begin{exercise}
\label{exe.fps.laure.field}\fbox{2} Let $K$ be a field. Prove that the
$K$-algebra $K\left(  \left(  x\right)  \right)  $ is a field.

[\textbf{Hint:} You can take it for granted that $K\left(  \left(  x\right)
\right)  $ is a commutative $K$-algebra, as the proof is a mutatis-mutandis
variant of the analogous proof for usual FPSs.]
\end{exercise}

Here are some applications of Laurent polynomials:

\begin{exercise}
\label{exe.fps.laurent.reed-dawson}Let $n\in\mathbb{N}$. \medskip

\textbf{(a)} \fbox{5} Prove that%
\[
\sum_{k=0}^{n}\left(  -2\right)  ^{n-k}\dbinom{n}{k}\dbinom{2k}{k}=\dbinom
{n}{n/2}.
\]
(Recall that $\dbinom{n}{i}=0$ whenever $i\notin\mathbb{N}$.) \medskip

\textbf{(b)} \fbox{2} More generally, prove that%
\[
\sum_{k=0}^{n}\left(  -2\right)  ^{n-k}\dbinom{n}{k}\dbinom{2k}{k+p}%
=\dbinom{n}{\left(  n+p\right)  /2}\ \ \ \ \ \ \ \ \ \ \text{for any }%
p\in\mathbb{Z}.
\]

[\textbf{Hint:} Compute the coefficients of the Laurent polynomial $\left(
\left(  x+x^{-1}\right)  ^{2}-2\right)  ^{n}$ in two ways.]
\end{exercise}

\begin{exercise}
\label{exe.fps.laure.res}For any Laurent series $f=\sum_{n\in\mathbb{Z}}%
f_{n}x^{n}\in K\left(  \left(  x\right)  \right)  $ (with $f_{n}\in K$), we
define the \emph{residue} of $f$ to be its $x^{-1}$-coefficient $f_{-1}$. We
denote this residue by $\operatorname*{Res}f$. (This is an algebraic analogue
of the \textquotedblleft residue at $0$\textquotedblright\ from complex analysis.)

The \emph{order }$\operatorname*{ord}f$ of a nonzero Laurent series
$f=\sum_{n\in\mathbb{Z}}f_{n}x^{n}\in K\left(  \left(  x\right)  \right)  $
(with $f_{n}\in K$) shall mean the smallest $n\in\mathbb{Z}$ satisfying
$f_{n}\neq0$. (This is well-defined, since all sufficiently low $n\in
\mathbb{Z}$ satisfy $f_{n}=0$ by the definition of a Laurent series.) The
\emph{trailing coefficient} of a nonzero Laurent series $f\in K\left(  \left(
x\right)  \right)  $ means the coefficient $\left[  x^{\operatorname*{ord}%
f}\right]  f$ (that is, the $x^{\operatorname*{ord}f}$-coefficient of $f$).
For example, the Laurent series $-x^{-2}+3+7x$ has order $-2$ and trailing
coefficient $-1$.

The \emph{derivative} $f^{\prime}$ of a Laurent series $f\in K\left(  \left(
x\right)  \right)  $ is defined as follows: If $f=\sum_{n\in\mathbb{Z}}%
f_{n}x^{n}$ with $f_{n}\in K$, then $f^{\prime}:=\sum_{n\in\mathbb{Z}}%
nf_{n}x^{n-1}$.

Prove the following: \medskip

\textbf{(a)} \fbox{1} Any $f\in K\left(  \left(  x\right)  \right)  $
satisfies $\operatorname*{Res}\left(  f^{\prime}\right)  =0$. \medskip

\textbf{(b)} \fbox{3} Any $n\in\mathbb{N}$ and any $f\in K\left(  \left(
x\right)  \right)  $ satisfy $\operatorname*{Res}\left(  f^{n}f^{\prime
}\right)  =0$. \medskip

\textbf{(c)} \fbox{1} If $f\in K\left(  \left(  x\right)  \right)  $ is a
nonzero Laurent series whose trailing coefficient is invertible (in $K$), then
$f$ is invertible in $K\left(  \left(  x\right)  \right)  $. (Keep in mind
that the word \textquotedblleft invertible\textquotedblright\ refers to
multiplicative inverses, not compositional inverses.) \medskip

\textbf{(d)} \fbox{3} If $f\in K\left(  \left(  x\right)  \right)  $ is a
nonzero Laurent series whose trailing coefficient is invertible (in $K$), then
each $n\in\mathbb{Z}$ satisfies%
\[
\operatorname*{Res}\left(  f^{n}f^{\prime}\right)  =%
\begin{cases}
0, & \text{if }n\neq-1;\\
\operatorname*{ord}f, & \text{if }n=-1.
\end{cases}
\]
(To be fully precise, \textquotedblleft$\operatorname*{ord}f$%
\textquotedblright\ here means the element $\left(  \operatorname*{ord}%
f\right)  \cdot1_{K}$ of the ring $K$.) \medskip

[\textbf{Hint:} In parts of this exercise, it may be expedient to first prove
the claim under the assumption that $K$ is a $\mathbb{Q}$-algebra (so that
$1,2,3,\ldots$ can be divided by in $K$), and then to argue that the
assumption can be lifted.]
\end{exercise}

\begin{exercise}
\label{exe.fps.lagr-inv}Let $f=\sum_{n>0}f_{n}x^{n}$ (with $f_{1},f_{2}%
,f_{3},\ldots\in K$) be an FPS in $K\left[  \left[  x\right]  \right]  $ whose
constant term is $0$. Assume that $f$ has a compositional inverse
$g=\sum_{n>0}g_{n}x^{n}$ (with $g_{1},g_{2},g_{3},\ldots\in K$). \medskip

\textbf{(a)} \fbox{2} Prove that there exists a unique FPS $h\in K\left[
\left[  x\right]  \right]  $ with $x=fh$. (This FPS $h$ is usually denoted by
$\dfrac{x}{f}$, but this notation is not an instance of Definition
\ref{def.commring.fracs} \textbf{(b)}, since $f$ is not invertible.) \medskip

\textbf{(b)} \fbox{4} Prove the \emph{Lagrange inversion formula}, which says
that%
\[
n\cdot g_{n}=\left[  x^{n-1}\right]  \left(  h^{n}\right)
\ \ \ \ \ \ \ \ \ \ \text{for any positive integer }n.
\]

\textbf{(c)} \fbox{2} The \emph{Lambert W series} is defined to be the
compositional inverse of the FPS $x\cdot\exp\left[  x\right]  =\sum
_{n\in\mathbb{N}}\dfrac{x^{n+1}}{n!}$ in $\mathbb{Q}\left[  \left[  x\right]
\right]  $. Find an explicit formula for the $x^{n}$-coefficient of this
series. \medskip

\textbf{(d)} \fbox{2} Consider again the FPS $C\left(  x\right)  =c_{0}%
+c_{1}x+c_{2}x^{2}+\cdots\in\mathbb{Q}\left[  \left[  x\right]  \right]  $
from Example 2 in Section \ref{sec.gf.exas}. Let us rename it as $C$. We
proved the equality $C=1+xC^{2}$ in that example (albeit we wrote it as
$C\left(  x\right)  =1+x\left(  C\left(  x\right)  \right)  ^{2}$). Set
$f=x-x^{2}$ and $g=xC$. Show that the FPS $g$ is a compositional inverse of
$f$. Use the Lagrange inversion formula to reprove the formula $c_{n}%
=\dfrac{1}{n+1}\dbinom{2n}{n}$ without the quadratic formula. \medskip

\textbf{(e)} \fbox{2} Let $m>1$ be an integer. Let $D\in\mathbb{Q}\left[
\left[  x\right]  \right]  $ be an FPS with constant term $1$ that satisfies
$D=1+x^{m-1}D^{m}$. Find an explicit formula for the $x^{n}$-coefficient of
$D$. (This generalizes part \textbf{(d)}, which is obtained for $m=2$.)
\medskip

[\textbf{Hint:} For part \textbf{(b)}, take derivatives of both sides in
$g\circ f=x$ to obtain $\left(  g^{\prime}\circ f\right)  \cdot f^{\prime}=1$;
then divide by $f^{n}$ in the Laurent series ring $K\left(  \left(  x\right)
\right)  $, and rewrite $g^{\prime}\circ f$ as $\sum_{k>0}kg_{k}f^{k-1}$. Now
take residues and use Exercise \ref{exe.fps.laure.res} \textbf{(d)}.] \medskip

[\textbf{Remarks:} Part \textbf{(b)} is remarkable for connecting the
compositional inverse $g$ of $f$ with the multiplicative inverse $h$ of $f/x$.
Since multiplicative inverses are usually easier to compute, it is a helpful
tool for the computation of compositional inverses.

The Lambert W series in part \textbf{(c)} is the Taylor series of the
\href{https://en.wikipedia.org/wiki/Lambert_W_function}{Lambert W function}.]
\end{exercise}

The following exercise tells a cautionary tale about applying some of our
results past their stated assumptions:

\begin{exercise}
\fbox{2} Extending Definition \ref{def.fps.lim.coeff-stab} to the $K$-module
$K\left[  \left[  x^{\pm}\right]  \right]  $, we obtain the notion of a limit
of a sequence of doubly infinite power series. \medskip

\textbf{(a)} Prove that $\lim\limits_{n\rightarrow\infty}\left(  x^{n}%
+x^{-n}\right)  =0$. \medskip

\textbf{(b)} Prove that $\lim\limits_{n\rightarrow\infty}\left(  \left(
x^{n}+x^{-n}\right)  ^{2}\right)  =2$. \medskip

\textbf{(c)} Can Proposition \ref{prop.fps.lim.sum-prod} be generalized to
Laurent series instead of FPSs? \medskip

[\textbf{Note:} This does not mean that the notion of limits of sequences of
Laurent series is completely useless. They behave reasonably as long as
multiplication is not involved.]
\end{exercise}

\subsubsection{Multivariate FPSs}

\begin{exercise}
\label{exe.fps.sum-n-choose-k-xn-elementary}\fbox{2} Let $k\in\mathbb{N}$.
Prove that%
\[
\sum_{n\in\mathbb{N}}\dbinom{n}{k}x^{n}=\dfrac{x^{k}}{\left(  1-x\right)
^{k+1}}%
\]
without using multivariate power series.
\end{exercise}

\begin{exercise}
\fbox{2} Prove that%
\[
\sum_{n\in\mathbb{N}}\dfrac{x^{n}}{1-yq^{n}}=\sum_{n\in\mathbb{N}}\dfrac
{y^{n}}{1-xq^{n}}\ \ \ \ \ \ \ \ \ \ \text{in the ring }K\left[  \left[
x,y,q\right]  \right]  .
\]

\end{exercise}

The next two exercises are concerned with FPSs in two indeterminates $x$ and
$y$.

\begin{exercise}
\fbox{4} For any $n\in\mathbb{N}$ and $m\in\mathbb{N}$, we let $f\left(
m,n\right)  $ denote the \# of $n$-tuples $\left(  \alpha_{1},\alpha
_{2},\ldots,\alpha_{n}\right)  $ of integers satisfying $\left\vert \alpha
_{1}\right\vert +\left\vert \alpha_{2}\right\vert +\cdots+\left\vert
\alpha_{n}\right\vert \leq m$. \medskip

\textbf{(a)} Prove that $\sum_{\left(  n,m\right)  \in\mathbb{N}%
\times\mathbb{N}}f\left(  m,n\right)  x^{m}y^{n}=\dfrac{1}{1-x-y-xy}$ in
$K\left[  \left[  x,y\right]  \right]  $. \medskip

\textbf{(b)} Prove that $f\left(  m,n\right)  =f\left(  n,m\right)  $ for all
$n,m\in\mathbb{N}$.
\end{exercise}

\begin{exercise}
\fbox{4} Let $f\in K\left[  \left[  x,y\right]  \right]  $ be an FPS such that
each positive integer $b$ satisfies $f\left[  x,x^{b}\right]  =0$. Prove that
$f=0$.
\end{exercise}

\begin{exercise}
\label{exe.fps.stirling2.4}\fbox{1} Recall the Stirling numbers of the 2nd
kind $S\left(  n,k\right)  $ studied in Exercise \ref{exe.fps.stirling2.1} and
Exercise \ref{exe.fps.stirling2.3}. Show that in $\mathbb{Q}\left[  \left[
x,y\right]  \right]  $, we have%
\[
\sum_{n\in\mathbb{N}}\ \ \sum_{k\in\mathbb{N}}\dfrac{S\left(  n,k\right)
}{n!}x^{n}y^{k}=\exp\left[  y\cdot\left(  \exp\left[  x\right]  -1\right)
\right]  .
\]

\end{exercise}

\begin{exercise}
\fbox{3} Recall the Eulerian polynomials $A_{m}\in\mathbb{Z}\left[  x\right]
$ from Exercise \ref{exe.gf.eulerian-pol.basics}. Prove that in $\mathbb{Q}%
\left[  \left[  x,y\right]  \right]  $, we have%
\[
\sum_{n\in\mathbb{N}}A_{n}\cdot\dfrac{y^{n}}{n!}=\dfrac{1-x}{1-x\exp\left[
\left(  1-x\right)  y\right]  }.
\]

\end{exercise}

Next comes an application of multivariate polynomials to proving a famous
binomial identity:

\begin{exercise}
\label{exe.fps.dixon}In this exercise, we shall prove \emph{Dixon's identity},
which states that%
\begin{equation}
\sum_{k\in\mathbb{Z}}\left(  -1\right)  ^{k}\dbinom{b+c}{c+k}\dbinom{c+a}%
{a+k}\dbinom{a+b}{b+k}=\dfrac{\left(  a+b+c\right)  !}{a!b!c!}
\label{eq.exe.fps.dixon.eq}%
\end{equation}
for any $a,b,c\in\mathbb{N}$. \medskip

\textbf{(a)} \fbox{1} Set%
\[
F\left(  a,b,c\right)  :=\sum_{k\in\mathbb{Z}}\left(  -1\right)  ^{k}%
\dbinom{b+c}{c+k}\dbinom{c+a}{a+k}\dbinom{a+b}{b+k}%
\]
for any $a,b,c\in\mathbb{N}$. Prove that $F\left(  a,b,c\right)  $ is
well-defined (i.e., the sum in this definition is summable). \medskip

\textbf{(b)} \fbox{1} Prove that (\ref{eq.exe.fps.dixon.eq}) holds whenever
$a=0$ or $b=0$ or $c=0$. \medskip

\textbf{(c)} \fbox{3} Prove that every $a,b,c\in\mathbb{N}$ satisfy%
\[
F\left(  a,b,c\right)  =\left(  -1\right)  ^{a+b+c}\cdot\left[  x^{2a}%
y^{2b}z^{2c}\right]  \left(  \left(  y-z\right)  ^{b+c}\left(  z-x\right)
^{c+a}\left(  x-y\right)  ^{a+b}\right)  .
\]
(Here, we are using polynomials in three indeterminates $x,y,z$.) \medskip

\textbf{(d)} \fbox{3} Prove that every three positive integers $a,b,c$
satisfy
\[
F\left(  a,b,c\right)  =F\left(  a-1,b,c\right)  +F\left(  a,b-1,c\right)
+F\left(  a,b,c-1\right)  .
\]

\textbf{(e)} \fbox{1} Prove that every three positive integers $a,b,c$
satisfy
\[
\dfrac{\left(  a+b+c\right)  !}{a!b!c!}=\dfrac{\left(  a-1+b+c\right)
!}{\left(  a-1\right)  !b!c!}+\dfrac{\left(  a+b-1+c\right)  !}{a!\left(
b-1\right)  !c!}+\dfrac{\left(  a+b+c-1\right)  !}{a!b!\left(  c-1\right)
!}.
\]

\textbf{(f)} \fbox{1} Prove (\ref{eq.exe.fps.dixon.eq}). \medskip

\textbf{(g)} \fbox{2} Show that each $n\in\mathbb{N}$ satisfies%
\[
\sum_{k=0}^{n}\left(  -1\right)  ^{k}\dbinom{n}{k}^{3}=%
\begin{cases}
\left(  -1\right)  ^{n/2}\dfrac{\left(  3n/2\right)  !}{\left(  n/2\right)
!^{3}}, & \text{if }n\text{ is even};\\
0, & \text{if }n\text{ is odd.}%
\end{cases}
\]

\textbf{(h)} \fbox{2} Show that
\[
\sum_{k\in\mathbb{Z}}\left(  -1\right)  ^{k}\dbinom{2a}{a+k}\dbinom{2b}%
{b+k}\dbinom{2c}{c+k}=\dfrac{\left(  a+b+c\right)  !\cdot\left(  2a\right)
!\left(  2b\right)  !\left(  2c\right)  !}{a!b!c!\cdot\left(  b+c\right)
!\left(  c+a\right)  !\left(  a+b\right)  !}%
\]
for any $a,b,c\in\mathbb{N}$. \medskip

[\textbf{Hint:} For part \textbf{(d)}, use the fact that%
\[
x^{2}\left(  y-z\right)  +y^{2}\left(  z-x\right)  +z^{2}\left(  x-y\right)
=-\left(  y-z\right)  \left(  z-x\right)  \left(  x-y\right)  ,
\]
and keep in mind that $\left[  x^{i}y^{j}z^{k}\right]  \left(  x^{m}p\right)
=\left[  x^{i-m}y^{j}z^{k}\right]  p$ for any $i,j,k,m,p$ with $i\geq m$.]
\end{exercise}

\subsection{Integer partitions and $q$-binomial coefficients}

The notations of Chapter \ref{chap.pars} shall be used here.

\subsubsection{Partition basics}

\begin{exercise}
\label{exe.pars.transpose}The purpose of this exercise is to make the proof of
Proposition \ref{prop.pars.pkn=dual} rigorous.

For any partition $\lambda=\left(  \lambda_{1},\lambda_{2},\ldots,\lambda
_{k}\right)  $, we define the \emph{Young diagram} $Y\left(  \lambda\right)  $
of $\lambda$ to be the finite set%
\[
\left\{  \left(  i,j\right)  \ \mid\ i\in\left\{  1,2,\ldots,k\right\}  \text{
and }j\in\left\{  1,2,\ldots,\lambda_{i}\right\}  \right\}  .
\]
Visually, this set $Y\left(  \lambda\right)  $ is represented by drawing each
$\left(  i,j\right)  \in Y\left(  \lambda\right)  $ as a cell of an
(invisible) matrix, namely as the cell in row $i$ and in column $j$. The
resulting picture is a table of $k$ left-aligned rows, where the $i$-th row
(counted from the top) has exactly $\lambda_{i}$ cells. For example, if
$\lambda=\left(  4,2,1\right)  $, then the Young diagram $Y\left(
\lambda\right)  $ of $\lambda$ is%
\[%
\ydiagram{4,2,1}%
\ \ .
\]
(This is only one way to draw Young diagrams; it is known as \emph{English
notation} or \emph{matrix notation}, since our labeling of cells matches the
way the cells of a matrix are commonly labeled. If we flip our pictures across
a horizontal axis, we would get \emph{French notation} aka \emph{Cartesian
notation}, as the labeling of cells would then match the Cartesian coordinates
of their centers.) \medskip

\textbf{(a)} \fbox{1} Prove that $\left\vert Y\left(  \lambda\right)
\right\vert =\left\vert \lambda\right\vert $ for any partition $\lambda$.
\medskip

\textbf{(b)} \fbox{1} Prove that the Young diagram $Y\left(  \lambda\right)  $
uniquely determines the partition $\lambda$. \medskip

A \emph{NW-set} shall mean a subset $S$ of $\left\{  1,2,3,\ldots\right\}
^{2}$ with the following property: If $\left(  i,j\right)  \in S$ and $\left(
i^{\prime},j^{\prime}\right)  \in\left\{  1,2,3,\ldots\right\}  ^{2}$ satisfy
$i^{\prime}\leq i$ and $j^{\prime}\leq j$, then $\left(  i^{\prime},j^{\prime
}\right)  \in S$ as well. (In terms of our above visual model, this means that
walking northwest from a cell of $S$ never moves you out of $S$, unless you
walk out of the matrix. For example, the set%
\[
\ydiagram{2+4,1+3}
\]
is not a NW-set, since the left neighbor of the leftmost cell in the topmost
row is not in this set.) \medskip

\textbf{(c)} \fbox{1} Prove that $Y\left(  \lambda\right)  $ is a NW-set for
each partition $\lambda$. \medskip

\textbf{(d)} \fbox{2} Prove that any finite NW-set has the form $Y\left(
\lambda\right)  $ for a unique partition $\lambda$. \medskip

Now, let $\operatorname*{flip}:\left\{  1,2,3,\ldots\right\}  ^{2}%
\rightarrow\left\{  1,2,3,\ldots\right\}  ^{2}$ be the map that sends each
$\left(  i,j\right)  \in\left\{  1,2,3,\ldots\right\}  ^{2}$ to $\left(
j,i\right)  $. Visually, this map $\operatorname*{flip}$ is a reflection in
the \textquotedblleft main diagonal\textquotedblright\ (the diagonal going
from the northwest to the southeast). We can apply $\operatorname*{flip}$ to a
subset of $\left\{  1,2,3,\ldots\right\}  ^{2}$ by applying
$\operatorname*{flip}$ to each element of this subset. For example:%
\[
\operatorname*{flip}\text{ sends \ \ }%
\ydiagram{4,2,1}%
\ \ \text{ to \ \ }%
\ydiagram{3,2,1,1}%
\ \ \text{.}%
\]

\textbf{(e)} \fbox{1} Prove that for any partition $\lambda$, there is a
unique partition $\lambda^{t}$ such that $Y\left(  \lambda^{t}\right)
=\operatorname*{flip}\left(  Y\left(  \lambda\right)  \right)  $. \medskip

This partition $\lambda^{t}$ is called the \emph{transpose} (or
\emph{conjugate}) of $\lambda$. \medskip

\textbf{(f)} \fbox{1} Prove that if $\lambda=\left(  \lambda_{1},\lambda
_{2},\ldots,\lambda_{k}\right)  $ is a partition, then the partition
$\lambda^{t}$ has exactly $\lambda_{1}$ parts. (Here, we set $\lambda_{1}=0$
if $k=0$.) \medskip

\textbf{(g)} \fbox{1} Prove that if $\lambda=\left(  \lambda_{1},\lambda
_{2},\ldots,\lambda_{k}\right)  $, then the $i$-th part of the partition
$\lambda^{t}$ equals the \# of all $j\in\left\{  1,2,\ldots,k\right\}  $ such
that $\lambda_{j}\geq i$. \medskip

\textbf{(h)} \fbox{1} Prove that $\left\vert \lambda^{t}\right\vert
=\left\vert \lambda\right\vert $ and $\left(  \lambda^{t}\right)  ^{t}%
=\lambda$ for any partition $\lambda$.
\end{exercise}

\begin{exercise}
\fbox{5} A partition will be called \emph{binarial} if all its parts are
powers of $2$. For instance, $\left(  8,2,2,1\right)  $ is a binarial
partition of $13$. Recall that the length of a partition $\lambda$ is denoted
by $\ell\left(  \lambda\right)  $.

Let $n>1$ be an integer. Prove that%
\[
\sum_{\substack{\lambda\text{ is a binarial}\\\text{partition of }n}}\left(
-1\right)  ^{\ell\left(  \lambda\right)  }=0.
\]
In other words, prove that the \# of binarial partitions of $n$ having even
length equals the \# of binarial partitions of $n$ having odd length.
\end{exercise}

\begin{exercise}
\fbox{3} A partition will be called \emph{trapezoidal} if it has the form
$\left(  j,j-1,j-2,\ldots,i\right)  $ for some integers $i\leq j$. (For
instance, $\left(  4,3,2\right)  $ and $\left(  5\right)  $ are trapezoidal partitions.)

Let $n$ be a positive integer. Prove that%
\[
\left(  \text{\# of trapezoidal partitions of }n\right)  =\left(  \text{\# of
odd positive divisors of }n\right)  .
\]

\end{exercise}

\begin{convention}
\label{conv.pars.vdash}Let $n$ be an integer. The notation \textquotedblleft%
$\lambda\vdash n$\textquotedblright\ shall mean \textquotedblleft$\lambda$ is
a partition of $n$\textquotedblright. Thus, for example, the summation sign
\textquotedblleft$\sum_{\lambda\vdash n}$\textquotedblright\ means a sum over
all partitions $\lambda$ of $n$.
\end{convention}

\begin{exercise}
\fbox{5} If $\lambda$ is a partition, and if $i$ is a positive integer, then
$m_{i}\left(  \lambda\right)  $ shall mean the \# of parts of $\lambda$ that
are equal to $i$. For instance, $m_{3}\left(  5,3,3,2\right)  =2$ and
$m_{4}\left(  5,3,3,2\right)  =0$.

Fix an $n\in\mathbb{N}$. \medskip

\textbf{(a)} Prove that%
\[
\prod_{\lambda\vdash n}\ \ \prod_{i=1}^{\infty}\left(  m_{i}\left(
\lambda\right)  \right)  !=\prod_{\lambda\vdash n}\ \ \prod_{i=1}^{\infty
}i^{m_{i}\left(  \lambda\right)  }.
\]

\textbf{(b)} More generally, prove that%
\[
\prod_{\lambda\vdash n}\ \ \prod_{\substack{\left(  i,j\right)  \in\left\{
1,2,3,\ldots\right\}  ^{2};\\j\leq m_{i}\left(  \lambda\right)  }}x_{j}%
=\prod_{\lambda\vdash n}\ \ \prod_{i=1}^{\infty}x_{i}^{m_{i}\left(
\lambda\right)  }%
\]
as monomials in $x_{1},x_{2},x_{3},\ldots$.

(For example, for $n=3$, this is saying that%
\begin{align*}
&  \underbrace{\left(  x_{1}\right)  }_{\text{factors for }\lambda=\left(
3\right)  }\cdot\underbrace{\left(  x_{1}x_{1}\right)  }_{\text{factors for
}\lambda=\left(  2,1\right)  }\cdot\underbrace{\left(  x_{1}x_{2}x_{3}\right)
}_{\text{factors for }\lambda=\left(  1,1,1\right)  }\\
&  =\underbrace{\left(  x_{3}^{1}\right)  }_{\text{factors for }%
\lambda=\left(  3\right)  }\cdot\underbrace{\left(  x_{1}^{1}x_{2}^{1}\right)
}_{\text{factors for }\lambda=\left(  2,1\right)  }\cdot\underbrace{\left(
x_{1}^{3}\right)  }_{\text{factors for }\lambda=\left(  1,1,1\right)  }.
\end{align*}
Make sure you understand why part \textbf{(a)} is a particular case of
\textbf{(b)}.)
\end{exercise}

\begin{exercise}
Recall the notations from Exercise \ref{exe.pars.transpose}.

We say that a partition $\lambda$ is a \emph{single-cell upgrade} of a
partition $\mu$ if we have $Y\left(  \mu\right)  \subseteq Y\left(
\lambda\right)  $ and $\left\vert Y\left(  \lambda\right)  \setminus Y\left(
\mu\right)  \right\vert =1$. (This is just saying that the Young diagram of
$\lambda$ is obtained from that of $\mu$ by adding one single cell.)

For instance, the single-cell upgrades of the partition $\left(  2,2,1\right)
$ are $\left(  3,2,1\right)  $, $\left(  2,2,2\right)  $ and $\left(
2,2,1,1\right)  $. On the other hand, the partition $\left(  2,2,1\right)  $
is a single-cell upgrade of each of the partitions $\left(  2,2\right)  $ and
$\left(  2,1,1\right)  $.

For any partition $\lambda$, let $\gamma\left(  \lambda\right)  $ denote the
\# of \textbf{distinct} parts of $\lambda$. For instance, $\gamma\left(
5,5,3,2,1,1\right)  =4$.

Prove the following:

\textbf{(a)} \fbox{2} For any partition $\lambda$, we have%
\begin{align*}
&  \left(  \text{\# of partitions }\mu\text{ such that }\lambda\text{ is a
single-cell upgrade of }\mu\right) \\
&  =\gamma\left(  \lambda\right)  .
\end{align*}

\textbf{(b)} \fbox{3} For any partition $\lambda$, we have%
\begin{align*}
&  \left(  \text{\# of single-cell upgrades of }\lambda\right) \\
&  =\left(  \text{\# of partitions }\mu\text{ such that }\lambda\text{ is a
single-cell upgrade of }\mu\right)  +1.
\end{align*}

\textbf{(c)} \fbox{3} We have%
\[
\sum_{\lambda\text{ is a partition}}\gamma\left(  \lambda\right)
x^{\left\vert \lambda\right\vert }=\dfrac{x}{1-x}\prod_{i=1}^{\infty}\dfrac
{1}{1-x^{i}}\ \ \ \ \ \ \ \ \ \ \text{in }\mathbb{Z}\left[  \left[  x\right]
\right]  .
\]

\end{exercise}

\begin{exercise}
\fbox{4} Let $d$ be a positive integer. Let $n\in\mathbb{N}$. Prove that
\begin{align*}
&  \left(  \text{\# of partitions of }n\text{ that have no part divisible by
}d\right) \\
&  =\left(  \text{\# of partitions of }n\text{ that have no }d\text{ equal
parts}\right)  .
\end{align*}

(For instance, the partition $\left(  4,2,2,2,2\right)  $ has no part
divisible by $3$, but it has $3$ equal parts.) \medskip

[\textbf{Remark:} Theorem \ref{thm.pars.odd-dist-equal} is the particular case
of this exercise for $d=2$.]
\end{exercise}

\begin{exercise}
\fbox{4} Let $n\in\mathbb{N}$. Let $p_{\operatorname*{dist}\operatorname*{odd}%
}\left(  n\right)  $ be the \# of partitions $\lambda$ of $n$ such that all
parts of $\lambda$ are distinct and odd. (For example, $\left(  7,3,1\right)
$ is such a partition.) Let $p_{+}\left(  n\right)  $ be the \# of partitions
of $n$ that have an even \# of even parts. (For example, $\left(
5,4,2\right)  $ is such a partition.) Let $p_{-}\left(  n\right)  $ be the \#
of partitions of $n$ that have an odd \# of even parts. (For example, $\left(
5,4,1\right)  $ is such a partition.) Prove that%
\[
p_{+}\left(  n\right)  -p_{-}\left(  n\right)  =p_{\operatorname*{dist}%
\operatorname*{odd}}\left(  n\right)  .
\]

\end{exercise}

\begin{exercise}
For any partition $\lambda=\left(  \lambda_{1},\lambda_{2},\ldots,\lambda
_{k}\right)  $, we define the integer
\[
\operatorname*{alt}\lambda=\sum_{i=1}^{k}\left(  -1\right)  ^{i-1}\lambda
_{i}=\lambda_{1}-\lambda_{2}+\lambda_{3}-\lambda_{4}\pm\cdots+\left(
-1\right)  ^{k-1}\lambda_{k}.
\]

\textbf{(a)} \fbox{1} Prove that $\operatorname*{alt}\lambda\in\mathbb{N}$ for
each partition $\lambda$. \medskip

\textbf{(b)} \fbox{1} Define the transpose $\lambda^{t}$ of a partition
$\lambda$ as in Exercise \ref{exe.pars.transpose}. Show that
$\operatorname*{alt}\left(  \lambda^{t}\right)  =\left(  \text{\# of odd parts
of }\lambda\right)  $ for any partition $\lambda$. \medskip

\textbf{(c)} \fbox{8} Prove that each $n,k\in\mathbb{N}$ satisfy%
\begin{align*}
&  \left(  \text{\# of partitions }\lambda\text{ of }n\text{ into odd parts
such that }\ell\left(  \lambda\right)  =k\right) \\
&  =\left(  \text{\# of partitions }\lambda\text{ of }n\text{ into distinct
parts such that }\operatorname*{alt}\lambda=k\right)  .
\end{align*}

\textbf{(d)} \fbox{1} Derive a new proof of Theorem
\ref{thm.pars.odd-dist-equal} from this.
\end{exercise}

\begin{exercise}
\label{exe.pars.maj-order.1}Let $n\in\mathbb{N}$. Let $\operatorname*{Par}%
\nolimits_{n}$ denote the set of all partitions of $n$.

We define a partial order $\preccurlyeq$ on the set $\operatorname*{Par}%
\nolimits_{n}$ as follows: For two partitions $\lambda=\left(  \lambda
_{1},\lambda_{2},\ldots,\lambda_{k}\right)  $ and $\mu=\left(  \mu_{1},\mu
_{2},\ldots,\mu_{\ell}\right)  $, we set $\lambda\preccurlyeq\mu$ if and only
if each positive integer $i$ satisfies%
\begin{equation}
\lambda_{1}+\lambda_{2}+\cdots+\lambda_{i}\leq\mu_{1}+\mu_{2}+\cdots+\mu_{i}.
\label{eq.exe.pars.maj-order.1.def}%
\end{equation}
Here, we set $\lambda_{j}:=0$ for each $j>k$, and we set $\mu_{j}:=0$ for each
$j>\ell$.

(For example, for $n=5$, we have $\left(  2,1,1,1\right)  \preccurlyeq\left(
2,2,1\right)  $, since we have
\begin{align*}
2  &  \leq2,\\
2+1  &  \leq2+2,\\
2+1+1  &  \leq2+2+1,\\
2+1+1+1  &  \leq2+2+1+0,\\
2+1+1+1+0  &  \leq2+2+1+0+0,
\end{align*}
and so on. Note that there are infinitely many inequalities to be checked, but
only finitely many of them are relevant, since both sides of
(\ref{eq.exe.pars.maj-order.1.def}) are essentially finite sums that stop
growing at some point.) \medskip

\textbf{(a)} \fbox{1} For any $n\geq6$, find two partitions $\lambda$ and
$\mu$ of $n$ satisfying neither $\lambda\preccurlyeq\mu$ nor $\mu
\preccurlyeq\lambda$. (This shows that $\preccurlyeq$ is not a total order for
$n\geq6$.) \medskip

\textbf{(b)} \fbox{2} Prove that two partitions $\lambda=\left(  \lambda
_{1},\lambda_{2},\ldots,\lambda_{k}\right)  $ and $\mu=\left(  \mu_{1},\mu
_{2},\ldots,\mu_{\ell}\right)  $ of $n$ satisfy $\lambda\preccurlyeq\mu$ if
and only if each $i\in\left\{  1,2,\ldots,k\right\}  $ satisfies
(\ref{eq.exe.pars.maj-order.1.def}). \medskip

\textbf{(c)} \fbox{2} Prove that two partitions $\lambda=\left(  \lambda
_{1},\lambda_{2},\ldots,\lambda_{k}\right)  $ and $\mu=\left(  \mu_{1},\mu
_{2},\ldots,\mu_{\ell}\right)  $ of $n$ satisfy $\lambda\preccurlyeq\mu$ if
and only if each $i\in\left\{  1,2,\ldots,\ell\right\}  $ satisfies
(\ref{eq.exe.pars.maj-order.1.def}). \medskip

\textbf{(d)} \fbox{4} Prove that two partitions $\lambda$ and $\mu$ of $n$
satisfy $\lambda\preccurlyeq\mu$ if and only if they satisfy $\mu
^{t}\preccurlyeq\lambda^{t}$. (See Exercise \ref{exe.pars.transpose} for the
definition of $\lambda^{t}$ and $\mu^{t}$.) \medskip

[\textbf{Note:} The partial order $\preccurlyeq$ is called the \emph{dominance
order} or the \emph{majorization order}; it is rather important in the theory
of symmetric functions.] \medskip

[\textbf{Hint:} It is helpful to identify a partition $\lambda=\left(
\lambda_{1},\lambda_{2},\ldots,\lambda_{k}\right)  $ with the weakly
decreasing essentially finite sequence $\widehat{\lambda}=\left(  \lambda
_{1},\lambda_{2},\ldots,\lambda_{k},0,0,0,\ldots\right)  $.]
\end{exercise}

\begin{exercise}
\label{exe.pars.transpose-sum}For any two partitions $\lambda=\left(
\lambda_{1},\lambda_{2},\ldots,\lambda_{k}\right)  $ and $\mu=\left(  \mu
_{1},\mu_{2},\ldots,\mu_{\ell}\right)  $, we define two partitions
$\lambda+\mu$ and $\lambda\sqcup\mu$ as follows:

\begin{itemize}
\item We let $\lambda+\mu$ be the partition $\left(  \lambda_{1}+\mu
_{1},\ \lambda_{2}+\mu_{2},\ \ldots,\ \lambda_{m}+\mu_{m}\right)  $, where we
set $m=\max\left\{  k,\ell\right\}  $, and where we set $\lambda_{j}:=0$ for
each $j>k$, and where we set $\mu_{j}:=0$ for each $j>\ell$.

\item We let $\lambda\sqcup\mu$ be the partition obtained by sorting the
entries of the list $\left(  \lambda_{1},\lambda_{2},\ldots,\lambda_{k}%
,\mu_{1},\mu_{2},\ldots,\mu_{\ell}\right)  $ in weakly decreasing order.
\end{itemize}

For example, $\left(  3,2,1\right)  +\left(  4,2\right)  =\left(
3+4,\ 2+2,\ 1+0\right)  =\left(  7,4,1\right)  $ and $\left(  3,2,1\right)
\sqcup\left(  4,2\right)  =\left(  4,3,2,2,1\right)  $. \medskip

Let $\lambda$ and $\mu$ be two partitions. Recall the definition of the
transpose of a partition (Exercise \ref{exe.pars.transpose}). \medskip

\textbf{(a)} \fbox{2} Prove that $\left(  \lambda\sqcup\mu\right)
^{t}=\lambda^{t}+\mu^{t}$. \medskip

\textbf{(b)} \fbox{2} Prove that $\left(  \lambda+\mu\right)  ^{t}=\lambda
^{t}\sqcup\mu^{t}$. \medskip

[\textbf{Hint:} Same as for Exercise \ref{exe.pars.maj-order.1}.]
\end{exercise}

\begin{exercise}
\label{exe.pars.sigma-I-elem}\fbox{3} Prove Theorem \ref{thm.pars.sigma1-I}
without any use of FPSs. \medskip

[\textbf{Hint:} Theorem \ref{thm.pars.sigma1-I} can be easily reduced to the
case when $I$ is finite; then, induction on $\left\vert I\right\vert $ becomes
available as a technique. A bijective proof also exists.]
\end{exercise}

\subsubsection{Euler's pentagonal number theorem}

Euler's pentagonal number theorem can be used freely in the following exercises.

\begin{exercise}
\textbf{(a)} \fbox{2} Prove that $p\left(  n\right)  \leq p\left(  n-1\right)
+p\left(  n-2\right)  $ for any $n>0$. \medskip

\textbf{(b)} \fbox{3} Prove that $p\left(  n\right)  \leq p\left(  n-1\right)
+p\left(  n-2\right)  -p\left(  n-5\right)  $ for any $n>0$. \medskip

\textbf{(c)} \fbox{2} Prove that $p\left(  n\right)  \geq p\left(  n-1\right)
+p\left(  n-2\right)  -p\left(  n-5\right)  -p\left(  n-7\right)  $ for any
$n\in\mathbb{N}$. \medskip

\textbf{(d)} \fbox{2} Use part \textbf{(a)} to obtain an upper bound for
$p\left(  n\right)  $ in terms of Fibonacci numbers.
\end{exercise}

\subsubsection{Jacobi's triple product identity}

\begin{exercise}
\textbf{(a)} \fbox{3} Prove that
\[
\prod_{m=1}^{\infty}\dfrac{1-x^{m}}{1+x^{m}}=\sum_{k\in\mathbb{Z}}\left(
-1\right)  ^{k}x^{k^{2}}.
\]

\textbf{(b)} \fbox{3} Prove that
\[
\prod_{m=1}^{\infty}\dfrac{1-x^{2m}}{1-x^{2m-1}}=\sum_{k\in\mathbb{N}%
}x^{k\left(  k+1\right)  /2}.
\]

[\textbf{Hint:} Both times, start by substituting appropriate values for $q$
and $z$ in the Jacobi Triple Product Identity.]
\end{exercise}

\begin{exercise}
\fbox{4} Prove that%
\[
\prod_{n=1}^{\infty}\left(  \left(  1-u^{n}v^{n-1}\right)  \left(
1-u^{n-1}v^{n}\right)  \left(  1-u^{n}v^{n}\right)  \right)  =\sum
_{k\in\mathbb{Z}}\left(  -1\right)  ^{k}u^{k\left(  k-1\right)  /2}v^{k\left(
k+1\right)  /2}%
\]
in the FPS ring $K\left[  \left[  u,v\right]  \right]  $. \medskip

[\textbf{Hint:} This is the Jacobi Triple Product Identity after (or, better,
before) a substitution.]
\end{exercise}

\begin{exercise}
\textbf{(a)} \fbox{6} Prove that
\[
\prod_{n=1}^{\infty}\left(  1-x^{n}\right)  ^{3}=\sum_{k=0}^{\infty}\left(
-1\right)  ^{k}\left(  2k+1\right)  x^{k\left(  k+1\right)  /2}%
\ \ \ \ \ \ \ \ \ \ \text{in }K\left[  \left[  x\right]  \right]  .
\]

\textbf{(b)} \fbox{3} Conclude that
\[
\prod_{n=1}^{\infty}\left(  1-x^{n}\right)  ^{3}=\sum_{k\in\mathbb{Z}}\left(
4k+1\right)  x^{k\left(  2k+1\right)  }\ \ \ \ \ \ \ \ \ \ \text{in }K\left[
\left[  x\right]  \right]  .
\]

\end{exercise}

\subsubsection{$q$-binomial coefficients}

\begin{exercise}
\label{exe.pars.qbinom.basics}\textbf{(a)} \fbox{1} Prove Proposition
\ref{prop.pars.qbinom.n0}. \medskip

\textbf{(b)} \fbox{3} Prove Theorem \ref{thm.pars.qbinom.rec} \textbf{(b)}.
\medskip

\textbf{(c)} \fbox{3} Prove Theorem \ref{thm.pars.qbinom.quot1}. \medskip

\textbf{(d)} \fbox{2} Prove Theorem \ref{thm.pars.qbinom.quot2}. \medskip

\textbf{(e)} \fbox{1} Prove Proposition \ref{prop.pars.qbinom.symm}. (Don't
forget the case $k>n$.)
\end{exercise}

\begin{exercise}
\label{exe.pars.qbinom.qbinom}\textbf{(a)} \fbox{2} Prove Theorem
\ref{thm.pars.qbinom.binom1} by induction on $n$. \medskip

\textbf{(b)} \fbox{2} Prove Theorem \ref{thm.pars.qbinom.binom2}.
\end{exercise}

\begin{exercise}
\label{exe.pars.qbinom.palin}\fbox{4} Let $n,k\in\mathbb{N}$ satisfy $n\geq
k$. Prove the following: \medskip

\textbf{(a)} The polynomial $\dbinom{n}{k}_{q}\in\mathbb{Z}\left[  q\right]  $
has degree $k\left(  n-k\right)  $. \medskip

\textbf{(b)} For each $i\in\left\{  0,1,\ldots,k\left(  n-k\right)  \right\}
$, we have%
\[
\left[  q^{i}\right]  \dbinom{n}{k}_{q}=\left[  q^{k\left(  n-k\right)
-i}\right]  \dbinom{n}{k}_{q}.
\]
(That is, the sequence of coefficients of this polynomial $\dbinom{n}{k}_{q}$
is palindromic.) \medskip

\textbf{(c)} We have $\dbinom{n}{k}_{q^{-1}}=q^{-k\left(  n-k\right)  }%
\dbinom{n}{k}_{q}$ in the Laurent polynomial ring $\mathbb{Z}\left[  q^{\pm
}\right]  $.
\end{exercise}

\begin{exercise}
\label{exe.pars.qbinom.upneg}Let us extend the definition of the $q$-integer
$\left[  n\right]  _{q}$ (Definition \ref{def.pars.qbinom.qint} \textbf{(a)})
to the case when $n$ is a negative integer as follows: If $n$ is a negative
integer, then we set%
\[
\left[  n\right]  _{q}:=-q^{-1}-q^{-2}-\cdots-q^{n}=-\sum_{k=n}^{-1}q^{k}%
\in\mathbb{Z}\left[  q^{\pm}\right]  .
\]
(This is a Laurent polynomial, not a polynomial any more.)

Furthermore, inspired by Theorem \ref{thm.pars.qbinom.quot2}, let us extend
the definition of $\dbinom{n}{k}_{q}$ to the case when $n$ is a negative
integer by setting%
\[
\dbinom{n}{k}_{q}=\dfrac{\left[  n\right]  _{q}\left[  n-1\right]  _{q}%
\cdots\left[  n-k+1\right]  _{q}}{\left[  k\right]  _{q}!}%
\ \ \ \ \ \ \ \ \ \ \text{for all }n\in\mathbb{Z}\text{ and }k\in\mathbb{N}%
\]
and%
\[
\dbinom{n}{k}_{q}=0\ \ \ \ \ \ \ \ \ \ \text{for all }n\in\mathbb{Z}\text{ and
}k\notin\mathbb{N}.
\]
The right hand sides of these two equalities are to be understood in the ring
$\mathbb{Z}\left(  \left(  q\right)  \right)  $ of Laurent series. (From
Theorem \ref{thm.pars.qbinom.quot2} and Convention
\ref{conv.pars.qbinom.neg-k}, we know that this definition does not conflict
with our existing definitions of $\dbinom{n}{k}_{q}$ for $n\in\mathbb{N}$.)
\medskip

\textbf{(a)} \fbox{1} Prove that $\left[  n\right]  _{1}=n$ for any
$n\in\mathbb{Z}$. \medskip

\textbf{(b)} \fbox{1} Prove that $\dbinom{n}{k}_{1}=\dbinom{n}{k}$ for any
$n\in\mathbb{Z}$ and $k\in\mathbb{Z}$. \medskip

\textbf{(c)} \fbox{1} Prove that $\left[  n\right]  _{q}=\dfrac{1-q^{n}}{1-q}$
in the ring $\mathbb{Z}\left(  \left(  q\right)  \right)  $ of Laurent series
for any $n\in\mathbb{Z}$. \medskip

\textbf{(d)} \fbox{1} Prove that $\left[  -n\right]  _{q}=-q^{-n}\left[
n\right]  _{q}$ for any $n\in\mathbb{Z}$. \medskip

\textbf{(e)} \fbox{1} Prove the \textquotedblleft$\emph{q}$\emph{-upper
negation formula}\textquotedblright\ (a $q$-analogue of Theorem
\ref{thm.binom.upneg-n}): If $n\in\mathbb{Z}$ and $k\in\mathbb{Z}$, then%
\[
\dbinom{-n}{k}_{q}=\left(  -1\right)  ^{k}q^{-kn-k\left(  k-1\right)
/2}\dbinom{k+n-1}{k}_{q}.
\]

\textbf{(f)} \fbox{1} Prove that $\dbinom{n}{k}_{q}\in\mathbb{Z}\left[
q^{\pm}\right]  $ (that is, $\dbinom{n}{k}_{q}$ is a Laurent polynomial, not
just a Laurent series) for any $n\in\mathbb{Z}$ and $k\in\mathbb{Z}$. \medskip

\textbf{(g)} \fbox{2} Prove that Theorem \ref{thm.pars.qbinom.rec} holds for
any $n\in\mathbb{Z}$ (not just for positive integers $n$). \medskip

\textbf{(h)} \fbox{1} Prove that Proposition \ref{prop.pars.qbinom.symm} does
\textbf{not} hold for negative $n$ (in general).
\end{exercise}

\begin{exercise}
Consider the ring $\mathbb{Z}\left[  \left[  z,q\right]  \right]  $ of FPSs in
two indeterminates $z$ and $q$. Let $n\in\mathbb{N}$. \medskip

\textbf{(a)} \fbox{1} Prove that%
\[
\prod_{i=0}^{n-1}\left(  1+zq^{i}\right)  =\sum_{k=0}^{n}q^{k\left(
k-1\right)  /2}\dbinom{n}{k}_{q}z^{k}.
\]

\textbf{(b)} \fbox{4} Prove that%
\[
\prod_{i=0}^{n-1}\dfrac{1}{1-zq^{i}}=\sum_{k\in\mathbb{N}}\dbinom{n+k-1}%
{k}_{q}z^{k},
\]
where we set $\dbinom{-1}{0}_{q}:=1$ (this is consistent with the definition
in Exercise \ref{exe.pars.qbinom.upneg}). (Up to sign, this is a $q$-analogue
of Proposition \ref{prop.fps.anti-newton-binom}.)
\end{exercise}

\begin{exercise}
\fbox{5} Let $n\in\mathbb{N}$. Prove the \emph{Gauss formula}%
\[
\sum_{k=0}^{n}\left(  -1\right)  ^{k}\dbinom{n}{k}_{q}=%
\begin{cases}
0, & \text{if }n\text{ is odd};\\
\left(  1-q^{1}\right)  \left(  1-q^{3}\right)  \left(  1-q^{5}\right)
\cdots\left(  1-q^{n-1}\right)  , & \text{if }n\text{ is even}.
\end{cases}
\]

\end{exercise}

\begin{exercise}
\fbox{5} Let $n\in\mathbb{N}$. Prove that%
\[
\sum_{k=0}^{n}q^{k}\dbinom{n}{k}_{q^{2}}=\left(  1+q^{1}\right)  \left(
1+q^{2}\right)  \cdots\left(  1+q^{n}\right)  .
\]

\end{exercise}

\begin{exercise}
Prove the following $q$-analogues of the Vandermonde convolution identity:
\medskip

\textbf{(a)} \fbox{3} If $a,b\in\mathbb{N}$ and $n\in\mathbb{N}$, then%
\[
\dbinom{a+b}{n}_{q}=\sum_{k=0}^{n}q^{k\left(  b-n+k\right)  }\dbinom{a}{k}%
_{q}\dbinom{b}{n-k}_{q}.
\]
(Note that the exponent $k\left(  b-n+k\right)  $ might occasionally be
negative, but in those cases we have $\dbinom{b}{n-k}_{q}=0$.) \medskip

\textbf{(b)} \fbox{3} If $a,b\in\mathbb{N}$ and $n\in\mathbb{N}$, then%
\[
\dbinom{a+b}{n}_{q}=\sum_{k=0}^{n}q^{\left(  a-k\right)  \left(  n-k\right)
}\dbinom{a}{k}_{q}\dbinom{b}{n-k}_{q}.
\]
(Note that the exponent $\left(  a-k\right)  \left(  n-k\right)  $ might
occasionally be negative, but in those cases we have $\dbinom{a}{k}_{q}=0$.)
\end{exercise}

\begin{exercise}
\fbox{3} Let $r,s,u,v\in\mathbb{Z}$ with $r>s>0$. Prove that%
\[
\lim\limits_{n\rightarrow\infty}\dbinom{rn+u}{sn+v}_{q}=\prod_{k=1}^{\infty
}\dfrac{1}{1-q^{k}}\ \ \ \ \ \ \ \ \ \ \text{in }\mathbb{Z}\left[  \left[
q\right]  \right]  .
\]
(See Definition \ref{def.fps.lim.coeff-stab} for the meaning of
\textquotedblleft$\lim\limits_{n\rightarrow\infty}$\textquotedblright\ used here.)
\end{exercise}

\begin{exercise}
We shall work in the ring $\left(  \mathbb{Z}\left[  z^{\pm}\right]  \right)
\left[  \left[  q\right]  \right]  $. \medskip

\textbf{(a)} \fbox{4} Prove \emph{MacMahon's identity}, which says that%
\[
\left(  \prod_{i=1}^{m}\left(  1+q^{2i-1}z\right)  \right)  \cdot\left(
\prod_{j=1}^{n}\left(  1+q^{2j-1}z^{-1}\right)  \right)  =\sum_{k=-n}%
^{m}q^{k^{2}}\dbinom{m+n}{k+n}_{q^{2}}z^{k}%
\]
for any $m,n\in\mathbb{N}$. \medskip

\textbf{(b)} \fbox{4} Recover Theorem \ref{thm.pars.jtp1} by taking the limit
$m\rightarrow\infty$ and then $n\rightarrow\infty$. (This yields a new proof
of Theorem \ref{thm.pars.jtp1}.)
\end{exercise}

The following exercise does not explicitly involve $q$-binomial coefficients,
but they can be used profitably in its solution:

\begin{exercise}
We work in the ring $\mathbb{Z}\left[  \left[  z,q\right]  \right]  $ of FPSs
in two indeterminates $z$ and $q$. \medskip

\textbf{(a)} \fbox{2} Prove that%
\[
\prod_{i=1}^{\infty}\left(  1+zq^{i}\right)  =\sum_{k\in\mathbb{N}}%
\dfrac{z^{k}q^{k\left(  k+1\right)  /2}}{\left(  1-q^{1}\right)  \left(
1-q^{2}\right)  \cdots\left(  1-q^{k}\right)  }.
\]

\textbf{(b)} \fbox{2} Prove that%
\[
\prod_{i=1}^{\infty}\dfrac{1}{1-zq^{i}}=\sum_{k\in\mathbb{N}}\dfrac{z^{k}%
q^{k}}{\left(  1-q^{1}\right)  \left(  1-q^{2}\right)  \cdots\left(
1-q^{k}\right)  }.
\]

\textbf{(c)} \fbox{5} Prove that%
\[
\prod_{i=1}^{\infty}\dfrac{1}{1-zq^{i}}=\sum_{k\in\mathbb{N}}\dfrac
{z^{k}q^{k^{2}}}{\left(  \left(  1-q^{1}\right)  \left(  1-q^{2}\right)
\cdots\left(  1-q^{k}\right)  \right)  \cdot\left(  \left(  1-zq^{1}\right)
\left(  1-zq^{2}\right)  \cdots\left(  1-zq^{k}\right)  \right)  }.
\]

[\textbf{Hint:} For part \textbf{(c)}, define the \emph{h-index} $h\left(
\lambda\right)  $ of a partition $\lambda$ to be the largest $i\in\mathbb{N}$
such that $\lambda$ has at least $i$ parts that are $\geq i$. For instance,
$h\left(  4,3,1,1\right)  =2$ and $h\left(  4,3,3,1\right)  =3$. Note that
$h\left(  \lambda\right)  $ is the size of the largest square that fits into
the Young diagram of $\lambda$. What remains if this square is removed from
the Young diagram of $\lambda$ ? See also
\href{https://en.wikipedia.org/wiki/H-index}{the Hirsch index}.]
\end{exercise}

\begin{exercise}
\label{exe.pars.qbinom.matrices-of-rank}Let $F$ be a finite field. Let
$n\in\mathbb{N}$ and $k\in\mathbb{N}$. Prove the following: \medskip

\textbf{(a)} \fbox{3} If $W$ is a $k$-dimensional $F$-vector space, then
\begin{align*}
&  \left(  \text{\# of }n\text{-tuples }\left(  w_{1},w_{2},\ldots
,w_{n}\right)  \text{ of vectors in }W\text{ that span }W\right) \\
&  =\left(  \left\vert F\right\vert ^{n}-1\right)  \cdot\left(  \left\vert
F\right\vert ^{n}-\left\vert F\right\vert \right)  \cdot\cdots\cdot\left(
\left\vert F\right\vert ^{n}-\left\vert F\right\vert ^{k-1}\right)
=\prod_{i=0}^{k-1}\left(  \left\vert F\right\vert ^{n}-\left\vert F\right\vert
^{i}\right)  .
\end{align*}

\textbf{(b)} \fbox{3} Let $m\in\mathbb{N}$. Then,%
\begin{align*}
&  \left(  \text{\# of }m\times n\text{-matrices }A\in F^{m\times n}\text{
satisfying }\operatorname*{rank}A=k\right) \\
&  =\dbinom{m}{k}_{\left\vert F\right\vert }\cdot\underbrace{\left(
\left\vert F\right\vert ^{n}-1\right)  \cdot\left(  \left\vert F\right\vert
^{n}-\left\vert F\right\vert \right)  \cdot\cdots\cdot\left(  \left\vert
F\right\vert ^{n}-\left\vert F\right\vert ^{k-1}\right)  }_{=\prod_{i=0}%
^{k-1}\left(  \left\vert F\right\vert ^{n}-\left\vert F\right\vert
^{i}\right)  }\\
&  =\dbinom{m}{k}_{\left\vert F\right\vert }\cdot\dbinom{n}{k}_{\left\vert
F\right\vert }\cdot\left\vert F\right\vert ^{k\left(  k-1\right)  /2}\left(
\left\vert F\right\vert -1\right)  ^{k}\left[  k\right]  _{\left\vert
F\right\vert }!.
\end{align*}

[\textbf{Hint:} For part \textbf{(a)}, what is the connection between
injective linear maps and surjective linear maps (between finite-dimensional
vector spaces)?]
\end{exercise}

\begin{exercise}
\fbox{3} Let $m,n\in\mathbb{N}$. Prove that%
\[
\sum_{k=0}^{m}\dbinom{m}{k}_{q}\dbinom{n}{k}_{q}q^{k\left(  k-1\right)
/2}\left(  q-1\right)  ^{k}\left[  k\right]  _{q}!=q^{mn}.
\]

[\textbf{Hint:} This is a polynomial identity in $q$, and there are infinitely
many prime numbers. What does this suggest?]
\end{exercise}

\begin{exercise}
For any $a,b\in K$ and any $n\in\mathbb{N}$, let us set%
\[
Q_{n}\left(  a,b\right)  :=\left(  aq^{0}+b\right)  \left(  aq^{1}+b\right)
\cdots\left(  aq^{n-1}+b\right)  .
\]

\textbf{(a)} \fbox{3} Prove that any $n\in\mathbb{N}$ and any $a,b,c,d\in K$
satisfy%
\[
\sum_{k=0}^{n}\dbinom{n}{k}_{q}Q_{k}\left(  a,b\right)  \cdot Q_{n-k}\left(
c,d\right)  =\sum_{k=0}^{n}\dbinom{n}{k}_{q}Q_{k}\left(  c,b\right)  \cdot
Q_{n-k}\left(  a,d\right)  .
\]

\textbf{(b)} \fbox{1} Prove that any $n\in\mathbb{N}$ and any $a,b,c\in K$
satisfy%
\[
Q_{n}\left(  a,b\right)  =\sum_{k=0}^{n}\dbinom{n}{k}_{q}Q_{k}\left(
a,c\right)  \cdot Q_{n-k}\left(  -c,b\right)  .
\]

[Note that this generalizes Theorem \ref{thm.pars.qbinom.binom1}.]
\end{exercise}

The following exercise gives a way to derive Theorem
\ref{thm.pars.qbinom.binom1} from Theorem \ref{thm.pars.qbinom.binom2}:

\begin{exercise}
\label{exe.pars.qbinom.2to1}Let $L$ be a commutative ring, and let $a,b,q\in
L$. (Note that if we let $L$ be the polynomial ring $K\left[  q\right]  $,
then this setting becomes the setting of Theorem \ref{thm.pars.qbinom.binom1}.)

Consider the polynomial ring $L\left[  u\right]  $. Let $A$ be the
(noncommutative) $L$-algebra $\operatorname*{End}\nolimits_{L}\left(  L\left[
u\right]  \right)  $ of all endomorphisms of the $L$-module $L\left[
u\right]  $. (Its multiplication is composition of endomorphisms. Note that
$L\left[  u\right]  $ is a free $L$-module of infinite rank, with basis
$\left(  u^{0},u^{1},u^{2},\ldots\right)  $; thus, the elements of $A$ can be
viewed as $\infty\times\infty$-matrices with each column having only finitely
many nonzero entries. But it is easier to just think of the elements of $A$ as
$L$-linear maps $L\left[  u\right]  \rightarrow L\left[  u\right]  $, just as
we defined $A$.)

Let $\alpha\in A=\operatorname*{End}\nolimits_{L}\left(  L\left[  u\right]
\right)  $ be the $L$-module endomorphism of $L\left[  u\right]  $ that
satisfies%
\[
\alpha\left(  u^{i}\right)  =aq^{i}u^{i+1}\ \ \ \ \ \ \ \ \ \ \text{for any
}i\in\mathbb{N}.
\]
(Thus, $\alpha$ sends each polynomial $f\in L\left[  u\right]  $ to $au\cdot
f\left[  qu\right]  $.)

Let $\beta\in A=\operatorname*{End}\nolimits_{L}\left(  L\left[  u\right]
\right)  $ be the $L$-module endomorphism of $L\left[  u\right]  $ that
satisfies%
\[
\beta\left(  u^{i}\right)  =bu^{i+1}\ \ \ \ \ \ \ \ \ \ \text{for any }%
i\in\mathbb{N}.
\]
(Thus, $\beta$ multiplies each polynomial $f\in L\left[  u\right]  $ by $bu$.)
\medskip

\textbf{(a)} \fbox{1} Prove that $\alpha\beta=q\beta\alpha$. \medskip

\textbf{(b)} \fbox{2} Prove that $\left(  \beta+\alpha\right)  ^{k}\left(
1\right)  =\left(  aq^{0}+b\right)  \left(  aq^{1}+b\right)  \cdots\left(
aq^{k-1}+b\right)  u^{k}$ for each $k\in\mathbb{N}$. \medskip

\textbf{(c)} \fbox{2} Rederive Theorem \ref{thm.pars.qbinom.binom1} by
applying Theorem \ref{thm.pars.qbinom.binom2} to $\beta$, $\alpha$ and $q$
instead of $a$, $b$ and $\omega$.
\end{exercise}

The following exercise shows an application of $q$-binomial coefficients using
a technique called the \emph{roots of unity filter} (see also Exercise
\ref{exe.binom.sum-n-choose-4k} for a similar technique):

\begin{exercise}
\label{exe.pars.qbinom.imo1995-6}Let $p$ be a prime. Let $\Omega$ be the set
of all complex numbers $z$ satisfying $z^{p}=1$. It is well-known that%
\[
\Omega=\left\{  e^{2\pi ig/p}\ \mid\ g\in\left\{  0,1,\ldots,p-1\right\}
\right\}
\]
(where the letters $e$ and $i$ have the usual meanings they have in complex
analysis) and, in particular, $\left\vert \Omega\right\vert =p$ and
$1\in\Omega$.

Let $n$ be a positive integer. \medskip

\textbf{(a)} \fbox{1} Prove that $\dbinom{np-1}{p-1}_{\omega}=1$ for each
$\omega\in\Omega\setminus\left\{  1\right\}  $. \medskip

\textbf{(b)} \fbox{2} Prove that $\dbinom{np}{p}_{\omega}=n$ for each
$\omega\in\Omega\setminus\left\{  1\right\}  $. \medskip

\textbf{(c)} \fbox{1} Prove that $\sum_{\omega\in\Omega}\omega^{k}=%
\begin{cases}
p, & \text{if }p\mid k;\\
0, & \text{if }p\nmid k
\end{cases}
$ for any $k\in\mathbb{N}$. \medskip

\textbf{(d)} \fbox{3} Prove that
\begin{align*}
&  \left(  \text{\# of }p\text{-element subsets }S\text{ of }\left[
np\right]  \text{ satisfying }p\mid\operatorname*{sum}S\right) \\
&  =\dfrac{\dbinom{np}{p}+\left(  -1\right)  ^{p-1}\left(  p-1\right)  n}{p}.
\end{align*}
Here, $\operatorname*{sum}S$ is defined as in Proposition
\ref{prop.pars.qbinom.alt-defs} \textbf{(b)}. \medskip

[\textbf{Hint:} For part \textbf{(d)}, compute $\sum_{\omega\in\Omega}%
\omega^{p\left(  p+1\right)  /2}\dbinom{np}{p}_{\omega}$ in two ways.]
\medskip

[\textbf{Remark:} The claim of part \textbf{(d)} generalizes Problem 6 from
the International Mathematical Olympiad 1995.]
\end{exercise}

We have seen $q$-analogues of nonnegative integers and of binomial
coefficients. The following exercise shows that quite a few more things have
$q$-analogues:

\begin{exercise}
\label{exe.pars.qbinom.q-deriv}Let $R$ be the ring $K\left[  \left[
q,x\right]  \right]  $ of FPSs in the two indeterminates $q$ and $x$. In this
exercise, we shall pretend that $q$ is a constant, so that $f\left[
qx\right]  $ (for a given FPS $f\in R$) shall denote the result of
substituting $qx$ for $x$ in $f$ while leaving $q$ unchanged.

Define the map $D_{q}:R\rightarrow R$ by letting%
\[
D_{q}f=\dfrac{f-f\left[  qx\right]  }{\left(  1-q\right)  x}%
\ \ \ \ \ \ \ \ \ \ \text{for each }f\in R.
\]

\textbf{(a)} \fbox{1} Prove that this is well-defined, i.e., that $f-f\left[
qx\right]  $ is really divisible by $\left(  1-q\right)  x$ in $R$ for any
$f\in R$. \medskip

The map $D_{q}$ is known as the $q$\emph{-derivative} or the \emph{Jackson
derivative}. It is easily seen to be $K\left[  \left[  q\right]  \right]
$-linear. \medskip

\textbf{(b)} \fbox{1} Prove that $D_{q}\left(  x^{n}\right)  =\left[
n\right]  _{q}x^{n-1}$ for each $n\in\mathbb{N}$. (For $n=0$, read the right
hand side as $0$.) \medskip

\textbf{(c)} \fbox{1} Prove the \textquotedblleft$q$-Leibniz
rule\textquotedblright: For any $f,g\in R$, we have%
\[
D_{q}\left(  fg\right)  =\left(  D_{q}f\right)  \cdot g+f\left[  qx\right]
\cdot D_{q}g.
\]

\textbf{(d)} \fbox{1} Let $\exp_{q}\in R$ be the FPS $\sum_{n\in\mathbb{N}%
}\dfrac{x^{n}}{\left[  n\right]  _{q}!}$ (called the \textquotedblleft%
$q$\emph{-exponential}\textquotedblright). Prove that $D_{q}\exp_{q}=\exp_{q}%
$. (Note that, unlike the usual exponential $\exp$, this does not require that
$K$ be a $\mathbb{Q}$-algebra, since $\left[  n\right]  _{q}!$ is a FPS with
constant term $1$ and thus always invertible.) \medskip

\textbf{(e)} \fbox{2} Prove that each $m\in\mathbb{N}$ satisfies%
\[
D_{q}\left(  \dfrac{1}{\prod_{i=0}^{m-1}\left(  1-q^{i}x\right)  }\right)
=\dfrac{\left[  m\right]  _{q}}{\prod_{i=0}^{m}\left(  1-q^{i}x\right)  }.
\]

Now, we shall discuss a $q$-analogue of the Eulerian polynomials from Exercise
\ref{exe.gf.eulerian-pol.basics}.

Let $K=\mathbb{Z}$. For any $m\in\mathbb{N}$, we define an FPS%
\[
Q_{q,m}:=\sum_{n\in\mathbb{N}}\left[  n\right]  _{q}^{m}x^{n}=\left[
0\right]  _{q}^{m}x^{0}+\left[  1\right]  _{q}^{m}x^{1}+\left[  2\right]
_{q}^{m}x^{2}+\cdots\in R.
\]
For example,%
\begin{align*}
Q_{q,0}  &  =\sum_{n\in\mathbb{N}}\underbrace{\left[  n\right]  _{q}^{0}}%
_{=1}x^{n}=\sum_{n\in\mathbb{N}}x^{n}=\dfrac{1}{1-x};\\
Q_{q,1}  &  =\sum_{n\in\mathbb{N}}\underbrace{\left[  n\right]  _{q}^{1}%
}_{=\left[  n\right]  _{q}}x^{n}=\sum_{n\in\mathbb{N}}\left[  n\right]
_{q}x^{n}=\dfrac{x}{\left(  1-x\right)  \left(  1-qx\right)  }%
\ \ \ \ \ \ \ \ \ \ \left(  \text{why?}\right)  ;
\end{align*}
it can furthermore be shown that%
\[
Q_{q,2}=\sum_{n\in\mathbb{N}}\left[  n\right]  _{q}^{2}x^{n}=\dfrac{x\left(
1+qx\right)  }{\left(  1-x\right)  \left(  1-qx\right)  \left(  1-q^{2}%
x\right)  }.
\]

It appears that each $Q_{q,m}$ has the form $\dfrac{A_{q,m}}{\left(
1-x\right)  \left(  1-qx\right)  \cdots\left(  1-q^{m}x\right)  }$ for some
polynomial $A_{q,m}\in\mathbb{Z}\left[  q,x\right]  $ of degree $m$ in $x$.
Let us prove this.

For each $m\in\mathbb{N}$, we define an FPS%
\[
A_{q,m}:=\left(  \prod_{i=0}^{m}\left(  1-q^{i}x\right)  \right)  Q_{q,m}\in
R.
\]
(Thus, $Q_{q,m}=\dfrac{A_{q,m}}{\left(  1-x\right)  \left(  1-qx\right)
\cdots\left(  1-q^{m}x\right)  }$.) \medskip

Let $\vartheta_{q}:R\rightarrow R$ be the $\mathbb{Z}$-linear map that sends
each FPS $f\in R$ to $x\cdot D_{q}f$. \medskip

\textbf{(f)} \fbox{1} Prove that $Q_{q,m}=\vartheta_{q}\left(  Q_{q,m-1}%
\right)  $ for each $m>0$. \medskip

\textbf{(g)} \fbox{2} Prove that $A_{q,m}=\left[  m\right]  _{q}%
xA_{q,m-1}+x\left(  1-x\right)  \cdot D_{q}A_{q,m-1}$ for each $m>0$. \medskip

\textbf{(h)} \fbox{1} Conclude that $A_{q,m}$ is a polynomial in $x$ and $q$
for each $m\in\mathbb{N}$. \medskip

\textbf{(i)} \fbox{2} Show that this polynomial $A_{q,m}$ has the form
$A_{q,m}=x^{m}q^{m\left(  m-1\right)  /2}+R_{q,m}$, where $R_{q,m}$ is a
$\mathbb{Z}$-linear combination of $x^{i}q^{j}$ with $i<m$ and $j<m\left(
m-1\right)  /2$. \medskip

The polynomials $A_{q,m}$ are called \emph{Carlitz's }$q$\emph{-Eulerian
polynomials}.
\end{exercise}

\subsection{Permutations}

The notations of Chapter \ref{chap.perm} shall be used here. In particular, if
$X$ is a set, then $S_{X}$ shall mean the symmetric group of $X$; and if $n$
is a nonnegative integer, then $S_{n}$ shall mean the symmetric group
$S_{\left[  n\right]  }$ of the set $\left[  n\right]  :=\left\{
1,2,\ldots,n\right\}  $.

\subsubsection{Basic definitions}

\begin{exercise}
\fbox{1} Let $n\in\mathbb{N}$ and $\sigma\in S_{n}$. What is the easiest way
to obtain \medskip

\textbf{(a)} a two-line notation of $\sigma^{-1}$ from a two-line notation of
$\sigma$ ? \medskip

\textbf{(b)} the one-line notation of $\sigma^{-1}$ from the one-line notation
of $\sigma$ ? \medskip

\textbf{(c)} the cycle digraph of $\sigma^{-1}$ from the cycle digraph of
$\sigma$ ?
\end{exercise}

\subsubsection{Transpositions, cycles and involutions}

\begin{exercise}
\label{exe.perm.cycles.c-t}Let $X$ be a set. Prove the following: \medskip

\textbf{(a)} \fbox{2} For any $k$ distinct elements $i_{1},i_{2},\ldots,i_{k}$
of $X$, we have%
\[
\operatorname*{cyc}\nolimits_{i_{1},i_{2},\ldots,i_{k}}=\underbrace{t_{i_{1}%
,i_{2}}t_{i_{2},i_{3}}\cdots t_{i_{k-1},i_{k}}}_{k-1\text{ transpositions}}.
\]

\textbf{(b)} \fbox{1} For any $k$ distinct elements $i_{1},i_{2},\ldots,i_{k}$
of $X$ and any $\sigma\in S_{X}$, we have%
\[
\sigma\operatorname*{cyc}\nolimits_{i_{1},i_{2},\ldots,i_{k}}\sigma
^{-1}=\operatorname*{cyc}\nolimits_{\sigma\left(  i_{1}\right)  ,\sigma\left(
i_{2}\right)  ,\ldots,\sigma\left(  i_{k}\right)  }.
\]

\end{exercise}

\begin{exercise}
\label{exe.perm.w0.basics}Let $n\in\mathbb{N}$. Let $w_{0}\in S_{n}$ be the
permutation that sends each $k\in\left[  n\right]  $ to $n+1-k$. Thus, $w_{0}$
is the permutation that \textquotedblleft reflects\textquotedblright\ all
numbers from $1$ to $n$ across the middle of the interval $\left[  n\right]
$. It is the unique strictly decreasing permutation of $\left[  n\right]  $.
In one-line notation, $w_{0}$ is $\left(  n,n-1,n-2,\ldots,2,1\right)  $.
\medskip

\textbf{(a)} \fbox{1} Prove that $w_{0}$ is an involution of $\left[
n\right]  $. \medskip

\textbf{(b)} \fbox{1} Prove that $w_{0}=t_{1,n}t_{2,n-1}\cdots t_{k,n+1-k}$,
where $k=\left\lfloor \dfrac{n}{2}\right\rfloor $. \medskip

\textbf{(c)} \fbox{3} Prove that
\begin{align*}
w_{0}  &  =\operatorname*{cyc}\nolimits_{1,2,\ldots,n}\operatorname*{cyc}%
\nolimits_{1,2,\ldots,n-1}\operatorname*{cyc}\nolimits_{1,2,\ldots,n-2}%
\cdots\operatorname*{cyc}\nolimits_{1}\\
&  =\operatorname*{cyc}\nolimits_{1}\operatorname*{cyc}\nolimits_{2,1}%
\operatorname*{cyc}\nolimits_{3,2,1}\cdots\operatorname*{cyc}%
\nolimits_{n,n-1,\ldots,1}.
\end{align*}

\end{exercise}

\begin{exercise}
\label{exe.perm.wilson-p-cycles}Let $p$ be a prime number. Let $Z$ be the set
of all $p$-cycles in the symmetric group $S_{p}$.

Let $\zeta$ be the specific $p$-cycle $\operatorname*{cyc}%
\nolimits_{1,2,\ldots,p}\in S_{p}$. Note that $\zeta$ has order $p$ in the
group $S_{p}$, and thus generates a cyclic subgroup $\left\langle
\zeta\right\rangle $ of order $p$. \medskip

\textbf{(a)} \fbox{2} Prove that a permutation $\sigma\in S_{p}$ satisfies
$\sigma\zeta=\zeta\sigma$ if and only if $\sigma\in\left\langle \zeta
\right\rangle $ (that is, if and only if $\sigma$ is a power of $\zeta$).
\medskip

\textbf{(b)} \fbox{2} Prove that $\left\vert \left\langle \zeta\right\rangle
\cap Z\right\vert =p-1$. \medskip

\textbf{(c)} \fbox{1} Prove that the cyclic group $\left\langle \zeta
\right\rangle $ acts on the set $Z$ by conjugation:%
\[
\alpha\rightharpoonup\sigma=\alpha\sigma\alpha^{-1}%
\ \ \ \ \ \ \ \ \ \ \text{for any }\alpha\in\left\langle \zeta\right\rangle
\text{ and }\sigma\in Z
\]
(where the symbol \textquotedblleft$\rightharpoonup$\textquotedblright\ means
the action of a group $G$ on a $G$-set $X$ -- i.e., we let $g\rightharpoonup
x$ denote the result of a group element $g\in G$ acting on some $x\in X$).
\medskip

\textbf{(d)} \fbox{1} Find the fixed points of this action. \medskip

\textbf{(e)} \fbox{1} Prove \emph{Wilson's theorem} from elementary number
theory, which states that%
\[
\left(  p-1\right)  !\equiv-1\operatorname{mod}p.
\]

\end{exercise}

\subsubsection{Inversions, length and Lehmer codes}

\begin{exercise}
\label{exe.perm.lengths-k-small-k}\fbox{4} Prove Proposition
\ref{prop.perm.lengths-k-small-k}.
\end{exercise}

\begin{exercise}
\label{exe.perm.len.cycles}Let $n\in\mathbb{N}$. Prove the following: \medskip

\textbf{(a)} \fbox{1} We have $\ell\left(  t_{i,j}\right)  =2\left\vert
i-j\right\vert -1$ for any distinct $i,j\in\left[  n\right]  $. \medskip

\textbf{(b)} \fbox{2} We have $\ell\left(  \operatorname*{cyc}%
\nolimits_{i+1,i+2,\ldots,i+k}\right)  =k-1$ for any integers $i$ and $k$ with
$0\leq i<i+k\leq n$. \medskip

\textbf{(c)} \fbox{4} We have $\ell\left(  \operatorname*{cyc}\nolimits_{i_{1}%
,i_{2},\ldots,i_{k}}\right)  \geq k-1$ for any $k$ distinct elements
$i_{1},i_{2},\ldots,i_{k}\in\left[  n\right]  $. \medskip

\textbf{(d)} \fbox{1} Are the $k$-cycles of the form $\operatorname*{cyc}%
\nolimits_{i+1,i+2,\ldots,i+k}$ the only $k$-cycles whose length is $k-1$ ?
\end{exercise}

\begin{exercise}
\label{exe.perm.lehmer.some-props}Let $\sigma\in S_{n}$ and $i\in\left[
n\right]  $. Prove the following (using the notation of Definition
\ref{def.perm.lehmer1} \textbf{(a)}): \medskip

\textbf{(a)} \fbox{1} We have $\ell_{i}\left(  \sigma\right)  =\left\vert
\left[  \sigma\left(  i\right)  -1\right]  \setminus\sigma\left(  \left[
i\right]  \right)  \right\vert $. \medskip

\textbf{(b)} \fbox{1} We have $\ell_{i}\left(  \sigma\right)  =\left\vert
\left[  \sigma\left(  i\right)  -1\right]  \setminus\sigma\left(  \left[
i-1\right]  \right)  \right\vert $. \medskip

\textbf{(c)} \fbox{1} We have $\sigma\left(  i\right)  \leq i+\ell_{i}\left(
\sigma\right)  $. \medskip

\textbf{(d)} \fbox{2} Assume that $i\in\left[  n-1\right]  $. We have
$\sigma\left(  i\right)  >\sigma\left(  i+1\right)  $ if and only if $\ell
_{i}\left(  \sigma\right)  >\ell_{i+1}\left(  \sigma\right)  $.
\end{exercise}

\subsubsection{V-permutations}

\begin{exercise}
\label{exe.perm.V-perm-as-cyc}\fbox{5} Let $n\in\mathbb{N}$. For each
$r\in\left[  n\right]  $, let $c_{r}$ denote the permutation
$\operatorname*{cyc}\nolimits_{r,r-1,\ldots,2,1}\in S_{n}$. (Thus,
$c_{1}=\operatorname*{cyc}\nolimits_{1}=\operatorname*{id}$ and $c_{2}%
=\operatorname*{cyc}\nolimits_{2,1}=s_{1}$.)

Let $G=\left\{  g_{1},g_{2},\ldots,g_{p}\right\}  $ be a subset of $\left[
n\right]  $, with $g_{1}<g_{2}<\cdots<g_{p}$. Let $\sigma\in S_{n}$ be the
permutation $c_{g_{1}}c_{g_{2}}\cdots c_{g_{p}}$. \medskip

[\textbf{Example:} If $n=6$ and $p=2$ and $G=\left\{  2,5\right\}  $, then
$\sigma=c_{2}c_{5}=\operatorname*{cyc}\nolimits_{2,1}\operatorname*{cyc}%
\nolimits_{5,4,3,2,1}$. In one-line notation, this permutation $\sigma$ is
$521346$.] \medskip

Prove the following: \medskip

\textbf{(a)} We have $\sigma\left(  1\right)  >\sigma\left(  2\right)
>\cdots>\sigma\left(  p\right)  $. \medskip

\textbf{(b)} We have $\sigma\left(  \left[  p\right]  \right)  =G$. \medskip

\textbf{(c)} We have $\sigma\left(  p+1\right)  <\sigma\left(  p+2\right)
<\cdots<\sigma\left(  n\right)  $. \medskip

(Note that a chain of inequalities that involves less than two numbers is
considered to be vacuously true. For example, Exercise
\ref{exe.perm.V-perm-as-cyc} \textbf{(c)} is vacuously true when $p=n-1$ and
also when $p=n$.)
\end{exercise}

\begin{exercise}
\label{exe.perm.V-perm-equivs}\fbox{8} Let $n\in\mathbb{N}$. Define the
permutations $c_{r}$ as in Exercise \ref{exe.perm.V-perm-as-cyc}.

Let $\sigma\in S_{n}$. We will use the notations from Definition
\ref{def.perm.lehmer1}. \medskip

\textbf{(a)} Prove that the following five statements are equivalent:

\begin{itemize}
\item \textit{Statement 1:} We have $\sigma\left(  1\right)  >\sigma\left(
2\right)  >\cdots>\sigma\left(  p\right)  $ and $\sigma\left(  p+1\right)
<\sigma\left(  p+2\right)  <\cdots<\sigma\left(  n\right)  $ for some
$p\in\left\{  0,1,\ldots,n\right\}  $. (In other words, the one-line notation
of $\sigma$ is decreasing at first, then increasing.)

\item \textit{Statement 2:} We have $\sigma=c_{g_{1}}c_{g_{2}}\cdots c_{g_{p}%
}$ for some elements $g_{1}<g_{2}<\cdots<g_{p}$ of $\left[  n\right]  $.

\item \textit{Statement 3:} For each $i\in\left[  n\right]  $, the set
$\sigma^{-1}\left(  \left[  i\right]  \right)  $ is an integer interval (i.e.,
there exist integers $u$ and $v$ such that $\sigma^{-1}\left(  \left[
i\right]  \right)  =\left\{  u,u+1,u+2,\ldots,v\right\}  $).

\item \textit{Statement 4:} If $i,j,k\in\left[  n\right]  $ satisfy $i<j<k$,
then we have $\sigma\left(  i\right)  >\sigma\left(  j\right)  $ or
$\sigma\left(  k\right)  >\sigma\left(  j\right)  $.

\item \textit{Statement 5:} We have $\ell_{1}\left(  \sigma\right)  >\ell
_{2}\left(  \sigma\right)  >\cdots>\ell_{p}\left(  \sigma\right)  $ and
$\ell_{p+1}\left(  \sigma\right)  =\ell_{p+2}\left(  \sigma\right)
=\cdots=\ell_{n}\left(  \sigma\right)  =0$ for some $p\in\left\{
0,1,\ldots,n\right\}  $. (In other words, the $n$-tuple $L\left(
\sigma\right)  $ is strictly decreasing until it reaches $0$, and then remains
at $0$.)
\end{itemize}

\textbf{(b)} Permutations $\sigma\in S_{n}$ satisfying the above five
statements are known as \textquotedblleft\emph{V-permutations}%
\textquotedblright\ (as their plot looks somewhat like the letter
\textquotedblleft V\textquotedblright: decreasing at first, then increasing).

Assume that $n>0$. Prove that the \# of V-permutations in $S_{n}$ is $2^{n-1}$.

[\textbf{Example:} If $n=3$, then the V-permutations in $S_{n}$ are (in
one-line notation) $123$ and $213$ and $312$ and $321$.]
\end{exercise}

\begin{exercise}
\label{exe.perm.V-perm-qcount}\fbox{4} Let $n\in\mathbb{N}$. Let $K$ be a
commutative ring, and let $a,b\in K$. Recall the $q$-binomial coefficients
from Definition \ref{def.pars.qbinom.qbinom} \textbf{(a)}. Prove that in the
polynomial ring $K\left[  q\right]  $, we have%
\[
\sum_{p=0}^{n}a^{p}b^{n-p}\sum_{\substack{\sigma\in S_{n};\\\sigma\left(
1\right)  >\sigma\left(  2\right)  >\cdots>\sigma\left(  p\right)
;\\\sigma\left(  p+1\right)  <\sigma\left(  p+2\right)  <\cdots<\sigma\left(
n\right)  }}q^{\ell\left(  \sigma\right)  }=\left(  aq^{0}+b\right)  \left(
aq^{1}+b\right)  \cdots\left(  aq^{n-1}+b\right)
\]
and%
\[
\sum_{p=0}^{n}a^{p}b^{n-p}\sum_{\substack{\sigma\in S_{n};\\\sigma\left(
1\right)  >\sigma\left(  2\right)  >\cdots>\sigma\left(  p\right)
;\\\sigma\left(  p+1\right)  <\sigma\left(  p+2\right)  <\cdots<\sigma\left(
n\right)  }}q^{\ell\left(  \sigma\right)  }=\sum_{k=0}^{n}q^{k\left(
k-1\right)  /2}\dbinom{n}{k}_{q}a^{k}b^{n-k}.
\]

[Both equalities can be proved combinatorially, thus leading to a
combinatorial proof of Theorem \ref{thm.pars.qbinom.binom1}.]
\end{exercise}

\subsubsection{Fixed points}

The next two exercises are concerned with the fixed points of maps (not only
of permutations).

\begin{definition}
\label{def.perm.fix}Let $X$ be a set. Let $f:X\rightarrow X$ be a map. Then:
\medskip

\textbf{(a)} A \emph{fixed point} of $f$ means an $x\in X$ satisfying
$f\left(  x\right)  =x$. \medskip

\textbf{(b)} We let $\operatorname*{Fix}f$ denote the set of all fixed points
of $f$. (This is the set $\left\{  x\in X\ \mid\ f\left(  x\right)
=x\right\}  $.)
\end{definition}

\begin{exercise}
\fbox{4} Let $X$ and $Y$ be two finite sets. Let $f:X\rightarrow Y$ and
$g:Y\rightarrow X$ be two maps. Prove that%
\[
\left\vert \operatorname*{Fix}\left(  f\circ g\right)  \right\vert =\left\vert
\operatorname*{Fix}\left(  g\circ f\right)  \right\vert .
\]

\end{exercise}

\begin{exercise}
Let $X$ be a finite set. \medskip

\textbf{(a)} \fbox{4} Prove that each permutation $\sigma\in S_{X}$ is a
composition of two involutions of $X$. \medskip

\textbf{(b)} \fbox{2} Prove that each permutation $\sigma\in S_{X}$ is
conjugate to its inverse $\sigma^{-1}$ in the symmetric group $S_{X}$.
\medskip

[\textbf{Hint:} For part \textbf{(a)}, it is easier to show the following
stronger claim: If $\sigma\in S_{X}$ and $Y\subseteq\operatorname*{Fix}\sigma$
and $p\in X\setminus Y$, then there exist two involutions $\alpha,\beta\in
S_{X}$ such that $\sigma=\alpha\circ\beta$ and $Y\subseteq\operatorname*{Fix}%
\alpha$ and $Y\cup\left\{  p\right\}  \subseteq\operatorname*{Fix}\beta$.

For part \textbf{(b)}, you can use part \textbf{(a)}.]
\end{exercise}

\subsubsection{More on inversions}

The next two exercises concern the inversions of a permutation. They use the
following definition:

\begin{definition}
\label{def.perm.Inv}Let $n\in\mathbb{N}$. For every $\sigma\in S_{n}$, we let
$\operatorname*{Inv}\sigma$ denote the set of all inversions of $\sigma$.
\end{definition}

We know from Corollary \ref{cor.perm.red.sigtau} \textbf{(b)} that any
$n\in\mathbb{N}$ and any two permutations $\sigma$ and $\tau$ in $S_{n}$
satisfy the inequality $\ell\left(  \sigma\tau\right)  \leq\ell\left(
\sigma\right)  +\ell\left(  \tau\right)  $. In the following exercise, we will
see when this inequality becomes an equality:

\begin{exercise}
\label{exe.perm.Inv.sub}\fbox{6} Let $n\in\mathbb{N}$. Let $\sigma\in S_{n}$
and $\tau\in S_{n}$. \medskip

\textbf{(a)} Prove that $\ell\left(  \sigma\tau\right)  =\ell\left(
\sigma\right)  +\ell\left(  \tau\right)  $ holds if and only if
$\operatorname*{Inv}\tau\subseteq\operatorname*{Inv}\left(  \sigma\tau\right)
$. \medskip

\textbf{(b)} Prove that $\ell\left(  \sigma\tau\right)  =\ell\left(
\sigma\right)  +\ell\left(  \tau\right)  $ holds if and only if
$\operatorname*{Inv}\left(  \sigma^{-1}\right)  \subseteq\operatorname*{Inv}%
\left(  \tau^{-1}\sigma^{-1}\right)  $. \medskip

\textbf{(c)} Prove that $\operatorname*{Inv}\sigma\subseteq\operatorname*{Inv}%
\tau$ holds if and only if $\ell\left(  \tau\right)  =\ell\left(  \tau
\sigma^{-1}\right)  +\ell\left(  \sigma\right)  $. \medskip

\textbf{(d)} Prove that if $\operatorname*{Inv}\sigma=\operatorname*{Inv}\tau
$, then $\sigma=\tau$. \medskip

\textbf{(e)} Prove that $\ell\left(  \sigma\tau\right)  =\ell\left(
\sigma\right)  +\ell\left(  \tau\right)  $ holds if and only if $\left(
\operatorname*{Inv}\sigma\right)  \cap\left(  \operatorname*{Inv}\left(
\tau^{-1}\right)  \right)  =\varnothing$.
\end{exercise}

Exercise \ref{exe.perm.Inv.sub} \textbf{(d)} shows that if two permutations in
$S_{n}$ have the same set of inversions, then they are equal. In other words,
a permutation in $S_{n}$ is uniquely determined by its set of inversions. The
next exercise shows what set of inversions a permutation can have:

\begin{exercise}
\label{exe.perm.invset}\fbox{7} Let $n\in\mathbb{N}$. Let $G=\left\{  \left(
i,j\right)  \in\mathbb{Z}^{2}\ \mid\ 1\leq i<j\leq n\right\}  $.

A subset $U$ of $G$ is said to be \emph{transitive} if every $a,b,c\in\left[
n\right]  $ satisfying $\left(  a,b\right)  \in U$ and $\left(  b,c\right)
\in U$ also satisfy $\left(  a,c\right)  \in U$.

A subset $U$ of $G$ is said to be \emph{inversive} if there exists a
$\sigma\in S_{n}$ such that $U=\operatorname*{Inv}\sigma$.

Let $U$ be a subset of $G$. Prove that $U$ is inversive if and only if both
$U$ and $G\setminus U$ are transitive.
\end{exercise}

\subsubsection{When transpositions generate $S_{X}$}

The next exercise uses a tiny bit of graph theory (the notion of connectedness
of a graph):

\begin{exercise}
\label{exe.perm.transp-gen}Let $X$ be a nonempty finite set. Let $G$ be a
loopless undirected graph with vertex set $X$. For each edge $e$ of $G$, we
let $t_{e}$ denote the transposition $t_{i,j}\in S_{X}$, where $i$ and $j$ are
the two endpoints of $e$. These transpositions $t_{e}$ for all edges $e$ of
$G$ will be called the $G$\emph{-edge transpositions}. \medskip

\textbf{(a)} \fbox{5} Prove that the $G$-edge transpositions generate the
symmetric group $S_{X}$ if and only if the graph $G$ is connected. \medskip

[\textbf{Example:} If $G$ is the path graph $%
\raisebox{-0.7pc}{
\begin{tikzpicture}%
[scale=1.5,thick,main node/.style={circle,fill=blue!20,draw}]
\node[main node] (1) at (1, 0) {$1$};
\node[main node] (2) at (2, 0) {$2$};
\node[main node] (3) at (3, 0) {$3$};
\node(4) at (4, 0) {$\cdots$};
\node[main node] (n) at (5, 0) {$n$};
\draw(1) -- (2) -- (3) -- (4) -- (n);
\end{tikzpicture}
}%
$ on the vertex set $X=\left[  n\right]  $, then the $G$-edge transpositions
are precisely the simple transpositions $s_{1},s_{2},\ldots,s_{n-1}$. In this
case, the claim of this exercise becomes the claim of Corollary
\ref{cor.perm.generated}.] \medskip

\textbf{(b)} \fbox{5} Let $\sigma\in S_{X}$ be a permutation that can be
written as a product of $G$-edge transpositions. Prove that $\sigma$ can be
written as a product of at most $\dbinom{n}{2}$ many $G$-edge transpositions,
where $n=\left\vert X\right\vert $.
\end{exercise}

Note that the bound $\dbinom{n}{2}$ in Exercise \ref{exe.perm.transp-gen}
\textbf{(b)} cannot be improved when $G$ is the path graph $%
\raisebox{-0.7pc}{
\begin{tikzpicture}%
[scale=1.5,thick,main node/.style={circle,fill=blue!20,draw}]
\node[main node] (1) at (1, 0) {$1$};
\node[main node] (2) at (2, 0) {$2$};
\node[main node] (3) at (3, 0) {$3$};
\node(4) at (4, 0) {$\cdots$};
\node[main node] (n) at (5, 0) {$n$};
\draw(1) -- (2) -- (3) -- (4) -- (n);
\end{tikzpicture}
}%
$, because in this case we are dealing with the Coxeter length of $\sigma$,
which can be $\dbinom{n}{2}$. However, for other graphs $G$, there can be
better bounds. In particular, $G$-edge transpositions are particularly
well-understood when $G$ is a \textquotedblleft star graph\textquotedblright;
this is a graph in which a single chosen vertex $x\in X$ is adjacent to all
the other vertices, but the other vertices are not adjacent with each other.
For instance, here is how this graph looks for $X=\left[  6\right]  $ and
$x=1$:%
\[%
\begin{tikzpicture}[scale=1.5,thick]
\begin{scope}[every node/.style={circle,fill=blue!20,draw}]
\node(1) at (0,0) {$1$};
\node(2) at (0*360/5:1) {$2$};
\node(3) at (1*360/5:1) {$3$};
\node(4) at (2*360/5:1) {$4$};
\node(5) at (3*360/5:1) {$5$};
\node(6) at (4*360/5:1) {$6$};
\end{scope}
\begin{scope}[every edge/.style={draw=black,very thick}, every loop/.style={}]
\path[-] (1) edge (2) (1) edge (3) (1) edge (4) (1) edge (5) (1) edge (6);
\end{scope}
\end{tikzpicture}%
\ \ .
\]
For such a graph $G$, the $G$-edge transpositions are just the transpositions
that swap $x$ with some other element of $X$. The factorizations of a
permutation into such transpositions are called \emph{star factorizations}.
For this kind of factorization, the $\dbinom{n}{2}$ bound in Exercise
\ref{exe.perm.transp-gen} \textbf{(b)} can be improved to $2n-2$ and even
better. See Exercise \ref{exe.perm.star-transp-fac} for this.

\subsubsection{Pattern avoidance}

A warm-up exercise for this subsection:

\begin{exercise}
Let $n\in\mathbb{N}$. \medskip

\textbf{(a)} \fbox{1} Find the \# of permutations $\sigma\in S_{n}$ such that
each $i\in\left[  n\right]  $ satisfies $\sigma\left(  i\right)  \leq i+1$.
\medskip

\textbf{(b)} \fbox{1} Find the \# of permutations $\sigma\in S_{n}$ such that
each $i\in\left[  n\right]  $ satisfies $i-1\leq\sigma\left(  i\right)  $.
\medskip

\textbf{(c)} \fbox{2} Find the \# of permutations $\sigma\in S_{n}$ such that
each $i\in\left[  n\right]  $ satisfies $i-1\leq\sigma\left(  i\right)  \leq
i+1$.
\end{exercise}

The next few exercises cover some of the most basic results in the theory of
\emph{pattern avoidance} (see \cite[Chapter 4]{Bona22} and \cite{Kitaev11} for
much more\footnote{There is \href{https://permutationpatterns.com/}{a yearly
conference} on this subject!}). This can be viewed as one possible way of
generalizing monotonicity (i.e., increasingness and decreasingness). We begin
by defining some basic concepts:

\begin{definition}
Let $\mathbf{t}=\left(  t_{1},t_{2},\ldots,t_{n}\right)  $ be an arbitrary
(finite) tuple. A \emph{subsequence} of $\mathbf{t}$ means a tuple of the form
$\left(  t_{i_{1}},t_{i_{2}},\ldots,t_{i_{m}}\right)  $, where $i_{1}%
,i_{2},\ldots,i_{m}$ are $m$ elements of $\left[  n\right]  $ satisfying
$i_{1}<i_{2}<\cdots<i_{m}$.
\end{definition}

\begin{example}
Let $\mathbf{t}=\left(  5,1,6,2,3,4\right)  $. Then, $\left(  1,3\right)  $
and $\left(  5,6,2\right)  $ are subsequences of $\mathbf{t}$ (indeed, if we
write $\mathbf{t}$ as $\left(  t_{1},t_{2},\ldots,t_{6}\right)  $, then
$\left(  1,3\right)  =\left(  t_{2},t_{5}\right)  $ and $\left(  5,6,2\right)
=\left(  t_{1},t_{3},t_{4}\right)  $), whereas $\left(  1,5\right)  $ and
$\left(  2,1,6\right)  $ and $\left(  1,1,6\right)  $ are not.
\end{example}

\begin{remark}
Let $\mathbf{t}=\left(  t_{1},t_{2},\ldots,t_{n}\right)  $ be an arbitrary
(finite) tuple. Then, the $1$-tuple $\left(  t_{i}\right)  $ is a subsequence
of $\mathbf{t}$ for any $i\in\left[  n\right]  $. So is the empty $0$-tuple
$\left(  {}\right)  $. Also, the tuple $\mathbf{t}$ itself is a subsequence of
$\mathbf{t}$ (and is the only length-$n$ subsequence of $\mathbf{t}$).
\end{remark}

Next, we need to define the notion of \emph{equally ordered tuples}. Roughly
speaking, these are tuples of the same length that might differ in their
values, but agree in the relative order of their values (e.g., if one tuple
has a smaller value in position $2$ than in position $5$, then so does the
other tuple). Here is the formal definition:

\begin{definition}
Let $\mathbf{a}=\left(  a_{1},a_{2},\ldots,a_{k}\right)  $ and $\mathbf{b}%
=\left(  b_{1},b_{2},\ldots,b_{k}\right)  $ be two $k$-tuples of integers. We
say that $\mathbf{a}$ and $\mathbf{b}$ are \emph{equally ordered} (to each
other) if for every pair $\left(  i,j\right)  \in\left[  k\right]
\times\left[  k\right]  $, we have the logical equivalence%
\[
\left(  a_{i}<a_{j}\right)  \Longleftrightarrow\left(  b_{i}<b_{j}\right)  .
\]
This relation is clearly symmetric in $\mathbf{a}$ and $\mathbf{b}$ (that is,
$\mathbf{a}$ and $\mathbf{b}$ are equally ordered if and only if $\mathbf{b}$
and $\mathbf{a}$ are equally ordered).

We agree that a $k$-tuple and an $\ell$-tuple are never equally ordered when
$k\neq\ell$.
\end{definition}

\begin{example}
\textbf{(a)} The two triples $\left(  3,1,6\right)  $ and $\left(
1,0,2\right)  $ are equally ordered. \medskip

\textbf{(b)} The two quadruples $\left(  3,1,1,2\right)  $ and $\left(
4,1,1,3\right)  $ are equally ordered. \medskip

\textbf{(c)} The two triples $\left(  3,1,2\right)  $ and $\left(
2,1,3\right)  $ are not equally ordered (indeed, we have $2<3$, but we don't
have $3<2$).
\end{example}

Now, we can define the notion of a \emph{pattern} in a tuple:

\begin{definition}
Let $\mathbf{s}=\left(  s_{1},s_{2},\ldots,s_{n}\right)  $ and $\mathbf{u}%
=\left(  u_{1},u_{2},\ldots,u_{m}\right)  $ be two tuples of integers.

A $\mathbf{u}$\emph{-pattern} in $\mathbf{s}$ means a subsequence of
$\mathbf{s}$ that is equally ordered to the tuple $\mathbf{u}$. (In
particular, this subsequence must have the same length as $\mathbf{u}$.)

In the following, when we talk about $\mathbf{u}$-patterns, we will often
write the tuple $\mathbf{u}$ without commas and parentheses. (For example, we
shall abbreviate \textquotedblleft$\left(  2,3,1\right)  $%
-pattern\textquotedblright\ as \textquotedblleft$231$%
-pattern\textquotedblright.)
\end{definition}

\begin{example}
Let $\mathbf{s}=\left(  s_{1},s_{2},\ldots,s_{n}\right)  $ be a tuple of
integers. Let us see what $\mathbf{u}$-patterns in $\mathbf{s}$ mean for
various specific tuples $\mathbf{u}$: \medskip

\textbf{(a)} A $21$-pattern in $\mathbf{s}$ is a subsequence $\left(
s_{i},s_{j}\right)  $ of $\mathbf{s}$ with $s_{i}>s_{j}$ (and, of course,
$i<j$, by the definition of a subsequence). For example, the tuple $\left(
4,5,2,1\right)  $ has five $21$-patterns (namely, $\left(  4,2\right)  $,
$\left(  4,1\right)  $, $\left(  5,2\right)  $, $\left(  5,1\right)  $ and
$\left(  2,1\right)  $). \medskip

\textbf{(b)} A $123$-pattern in $\mathbf{s}$ is a subsequence $\left(
s_{i},s_{j},s_{k}\right)  $ with $s_{i}<s_{j}<s_{k}$ (and, of course, $i<j<k$,
by the definition of a subsequence). For example, the tuple $\left(
2,1,3,5\right)  $ has two $123$-patterns (namely, $\left(  1,3,5\right)  $ and
$\left(  2,3,5\right)  $). \medskip

\textbf{(c)} A $231$-pattern in $\mathbf{s}$ is a subsequence $\left(
s_{i},s_{j},s_{k}\right)  $ with $s_{k}<s_{i}<s_{j}$ (and, of course, $i<j<k$,
by the definition of a subsequence). For example, the tuple $\left(
1,2,5,1\right)  $ has one $231$-pattern (namely, $\left(  2,5,1\right)  $).
\end{example}

Finally, we can define \emph{pattern avoidance}:

\begin{definition}
Let $\mathbf{s}$ and $\mathbf{u}$ be two tuples of integers. We say that
$\mathbf{s}$ is $\mathbf{u}$\emph{-avoiding} if there is no $\mathbf{u}%
$-pattern in $\mathbf{s}$.
\end{definition}

\begin{example}
The tuple $\left(  2,1,5,3,4\right)  $ is

\begin{itemize}
\item not $123$-avoiding, since it contains the $123$-pattern $\left(
1,3,4\right)  $ (and also the $123$-pattern $\left(  2,3,4\right)  $);

\item not $132$-avoiding, since it contains the $132$-pattern $\left(
1,5,3\right)  $;

\item not $312$-avoiding, since it contains the $312$-pattern $\left(
5,3,4\right)  $;

\item $231$-avoiding (check this!).
\end{itemize}
\end{example}

\begin{example}
\textbf{(a)} A tuple $\mathbf{s}=\left(  s_{1},s_{2},\ldots,s_{n}\right)  $ is
$21$-avoiding if and only if it is weakly increasing (i.e., satisfies
$s_{1}\leq s_{2}\leq\cdots\leq s_{n}$). Indeed, a $21$-pattern in $\mathbf{s}$
is a subsequence $\left(  s_{i},s_{j}\right)  $ of $\mathbf{s}$ with
$s_{i}>s_{j}$; thus, the non-existence of such $21$-patterns is equivalent to
$s_{1}\leq s_{2}\leq\cdots\leq s_{n}$. \medskip

\textbf{(b)} Likewise, a tuple $\mathbf{s}=\left(  s_{1},s_{2},\ldots
,s_{n}\right)  $ is $12$-avoiding if and only if it is weakly decreasing
(i.e., satisfies $s_{1}\geq s_{2}\geq\cdots\geq s_{n}$).
\end{example}

Finally, our concept of pattern avoidance can be extended from tuples to
permutations in the most obvious manner:

\begin{definition}
Let $\mathbf{u}$ be a tuple of integers. Let $n\in\mathbb{N}$ and $\sigma\in
S_{n}$. We say that the permutation $\sigma$ is $\mathbf{u}$\emph{-avoiding}
if the OLN of $\sigma$ (that is, the $n$-tuple $\left(  \sigma\left(
1\right)  ,\sigma\left(  2\right)  ,\ldots,\sigma\left(  n\right)  \right)  $)
is $\mathbf{u}$-avoiding.
\end{definition}

\begin{example}
Let $n\in\mathbb{N}$. The only $21$-avoiding permutation $\sigma\in S_{n}$ is
the identity permutation $\operatorname*{id}\in S_{n}$ (since it is the only
permutation whose OLN is weakly increasing). Likewise, the only $12$-avoiding
permutation $\sigma\in S_{n}$ is the permutation $w_{0}\in S_{n}$ from
Exercise \ref{exe.perm.w0.basics}.
\end{example}

After all this build-up, we can now study $\mathbf{u}$-avoidance for more
complicated patterns $\mathbf{u}$ than $21$ and $12$:

\begin{exercise}
\label{exe.perm.patt.132-av}\fbox{6} Let $n\in\mathbb{N}$. Let $c_{n}$ denote
the $n$-th Catalan number (from Example 2 in Section \ref{sec.gf.exas}).
\medskip

\textbf{(a)} Prove that
\[
\left(  \text{\# of }132\text{-avoiding permutations in }S_{n}\right)
=c_{n}.
\]

\textbf{(b)} Prove that
\[
\left(  \text{\# of }231\text{-avoiding permutations in }S_{n}\right)
=c_{n}.
\]

\textbf{(c)} Prove that
\[
\left(  \text{\# of }213\text{-avoiding permutations in }S_{n}\right)
=c_{n}.
\]

\textbf{(d)} Prove that
\[
\left(  \text{\# of }312\text{-avoiding permutations in }S_{n}\right)
=c_{n}.
\]

[\textbf{Hint:} Easy bijections show that parts \textbf{(a)}, \textbf{(b)},
\textbf{(c)} and \textbf{(d)} are equivalent. For part \textbf{(b)}, proceed
recursively: Assume that $n>0$, and let $\sigma\in S_{n}$, and let
$i=\sigma^{-1}\left(  n\right)  $. Show that the permutation $\sigma$ is
$231$-avoiding if and only if the two tuples $\left(  \sigma\left(  1\right)
,\sigma\left(  2\right)  ,\ldots,\sigma\left(  i-1\right)  \right)  $ and
$\left(  \sigma\left(  i+1\right)  ,\sigma\left(  i+2\right)  ,\ldots
,\sigma\left(  n\right)  \right)  $ are $231$-avoiding and satisfy $\left\{
\sigma\left(  1\right)  ,\sigma\left(  2\right)  ,\ldots,\sigma\left(
i-1\right)  \right\}  =\left\{  1,2,\ldots,i-1\right\}  $ and $\left\{
\sigma\left(  i+1\right)  ,\sigma\left(  i+2\right)  ,\ldots,\sigma\left(
n\right)  \right\}  =\left\{  i,i+1,\ldots,n-1\right\}  $. This yields a
recursive equation for the \# of $231$-avoiding permutations in $S_{n}$.]
\end{exercise}

\begin{exercise}
\label{exe.perm.patt.123-av}\fbox{8} Let $n\in\mathbb{N}$. Let $c_{n}$ denote
the $n$-th Catalan number (from Example 2 in Section \ref{sec.gf.exas}).
\medskip

\textbf{(a)} Prove that
\[
\left(  \text{\# of }123\text{-avoiding permutations in }S_{n}\right)
=c_{n}.
\]

\textbf{(b)} Prove that
\[
\left(  \text{\# of }321\text{-avoiding permutations in }S_{n}\right)
=c_{n}.
\]

[\textbf{Hint:} Consider any $321$-avoiding permutation $\sigma\in S_{n}$. A
\emph{record} of $\sigma$ means a value $\sigma\left(  i\right)  $ for some
$i\in\left[  n\right]  $ satisfying
\[
\sigma\left(  i\right)  >\sigma\left(  j\right)  \ \ \ \ \ \ \ \ \ \ \text{for
all }j\in\left[  i-1\right]  .
\]
(Equivalently, it is an entry in the OLN of $\sigma$ that is larger than all
entries further left.) Let $b_{1},b_{2},\ldots,b_{k}$ be the records of a
permutation $\sigma\in S_{n}$, written in increasing order (or, equivalently,
in the order of their appearance in the OLN of $\sigma$). (For example, if
$\sigma=14672385$, then these are $1,4,6,7,8$.) Argue first that the OLN of
$\sigma$ becomes weakly increasing when all records are removed from it. (For
example, $14672385$ becomes $235$ this way.) Write the OLN of $\sigma$ in the
form%
\[
\left(  \sigma\left(  1\right)  ,\sigma\left(  2\right)  ,\ldots,\sigma\left(
n\right)  \right)  =\left(  b_{1},\underbrace{\ldots}_{\substack{\text{some
}i_{1}\\\text{entries}}},b_{2},\underbrace{\ldots}_{\substack{\text{some
}i_{2}\\\text{entries}}},\ldots,b_{k-1},\underbrace{\ldots}%
_{\substack{\text{some }i_{k-1}\\\text{entries}}},b_{k},\underbrace{\ldots
}_{\substack{\text{some }i_{k}\\\text{entries}}}\right)  ,
\]
where $i_{1},i_{2},\ldots,i_{k}\in\mathbb{N}$. Set $c_{i}:=b_{i}-b_{i-1}$ for
each $i\in\left[  k\right]  $, where $b_{0}:=0$. Now, define $B\left(
\sigma\right)  $ to be the Dyck path%
\[
N^{c_{1}}S^{i_{1}+1}N^{c_{2}}S^{i_{2}+1}\cdots N^{c_{k}}S^{i_{k}+1},
\]
where an \textquotedblleft$N^{j}$\textquotedblright\ means $j$ consecutive
NE-steps, and where an \textquotedblleft$S^{j}$\textquotedblright\ means $j$
consecutive SE-steps. For example, $\sigma=14672385$ leads to%
\begin{align*}
&  B\left(  \sigma\right) \\
&  =N^{1}S^{1}N^{3}S^{1}N^{2}S^{1}N^{1}S^{3}N^{1}S^{2}\\
&  =NSNNNSNNSNSSSNSS\\
&  =%
\begin{tikzpicture}[scale=0.8]
\draw
(0, 0) -- (1, 1) -- (2, 0) -- (3, 1) -- (4, 2) -- (5, 3) -- (6, 2) -- (7, 3) -- (8, 4) -- (9, 3) -- (10, 4) -- (11, 3) -- (12, 2) -- (13, 1) -- (14, 2) -- (15, 1) -- (16, 0);
\fill
[red!50!white] (0, 0) -- (1, 1) -- (2, 0) -- (3, 1) -- (4, 2) -- (5, 3) -- (6, 2) -- (7, 3) -- (8, 4) -- (9, 3) -- (10, 4) -- (11, 3) -- (12, 2) -- (13, 1) -- (14, 2) -- (15, 1) -- (16, 0);
\filldraw(0, 0) circle [fill=red, radius=0.1];
\filldraw(1, 1) circle [fill=red, radius=0.1];
\filldraw(2, 0) circle [fill=red, radius=0.1];
\filldraw(3, 1) circle [fill=red, radius=0.1];
\filldraw(4, 2) circle [fill=red, radius=0.1];
\filldraw(5, 3) circle [fill=red, radius=0.1];
\filldraw(6, 2) circle [fill=red, radius=0.1];
\filldraw(7, 3) circle [fill=red, radius=0.1];
\filldraw(8, 4) circle [fill=red, radius=0.1];
\filldraw(9, 3) circle [fill=red, radius=0.1];
\filldraw(10, 4) circle [fill=red, radius=0.1];
\filldraw(11, 3) circle [fill=red, radius=0.1];
\filldraw(12, 2) circle [fill=red, radius=0.1];
\filldraw(13, 1) circle [fill=red, radius=0.1];
\filldraw(14, 2) circle [fill=red, radius=0.1];
\filldraw(15, 1) circle [fill=red, radius=0.1];
\filldraw(16, 0) circle [fill=red, radius=0.1];
\end{tikzpicture}%
\end{align*}
Prove that the map%
\begin{align*}
B:\left\{  321\text{-avoiding permutations in }S_{n}\right\}   &
\rightarrow\left\{  \text{Dyck paths from }\left(  0,0\right)  \text{ to
}\left(  2n,0\right)  \right\}  ,\\
\sigma &  \mapsto B\left(  \sigma\right)
\end{align*}
is well-defined and a bijection. (In proving surjectivity, don't forget to
check that the $b_{i}$'s really are the records of $\sigma$ !)]
\end{exercise}

Exercises \ref{exe.perm.patt.132-av} and \ref{exe.perm.patt.123-av} can be
combined into a single statement, which says that for any $\tau\in S_{3}$ and
any $n\in\mathbb{N}$, the \# of $\left(  \tau\left(  1\right)  ,\tau\left(
2\right)  ,\tau\left(  3\right)  \right)  $-avoiding permutations in $S_{n}$
equals the Catalan number $c_{n}$, independently of $\tau$. The independence
appears almost too good to be true. Parts of this miracle survive even for
$\tau\in S_{4}$; for example, for any $n\in\mathbb{N}$, we have%
\begin{align*}
&  \left(  \text{\# of }4132\text{-avoiding permutations in }S_{n}\right) \\
&  =\left(  \text{\# of }3142\text{-avoiding permutations in }S_{n}\right)
\end{align*}
(a result of Stankova \cite[Theorem 3.1]{Stanko94}), but this number does not
equal the \# of $1324$-avoiding permutations in $S_{n}$ (in general), nor does
it have any simple formula. A
\href{https://en.wikipedia.org/wiki/Enumerations_of_specific_permutation_classes}{Wikipedia
page} collects known results about these and similar numbers.

The $321$-avoiding permutations have an important role in \emph{Schubert
calculus}, a topic that connects combinatorics with algebraic geometry
(roughly speaking, the study of varieties of matrices satisfying certain rank
inequalities). See \cite[Chapters 2--3]{Manive01} for this.

We state a few simple results about permutations avoiding several patterns:

\begin{definition}
Let $\mathbf{u}$ and $\mathbf{v}$ be two tuples of integers. A permutation in
$S_{n}$ (or a tuple of integers) is said to be $\left(  \mathbf{u}%
,\mathbf{v}\right)  $\emph{-avoiding} if and only if it is both $\mathbf{u}%
$-avoiding and $\mathbf{v}$-avoiding. Similarly we define $\left(
\mathbf{u},\mathbf{v},\mathbf{w}\right)  $\emph{-avoiding} permutations (where
$\mathbf{u}$, $\mathbf{v}$ and $\mathbf{w}$ are three tuples of integers).
\end{definition}

\begin{exercise}
Let $n$ be a positive integer. \medskip

\textbf{(a)} \fbox{3} Prove that%
\[
\left(  \text{\# of }\left(  231,321\right)  \text{-avoiding permutations in
}S_{n}\right)  =2^{n-1}.
\]

\textbf{(b)} \fbox{2} Prove that%
\[
\left(  \text{\# of }\left(  132,231\right)  \text{-avoiding permutations in
}S_{n}\right)  =2^{n-1}.
\]

\textbf{(c)} \fbox{2} Prove that%
\[
\left(  \text{\# of }\left(  123,321\right)  \text{-avoiding permutations in
}S_{n}\right)  =0\ \ \ \ \ \ \ \ \ \ \text{if }n>4.
\]

\textbf{(d)} \fbox{2} Prove that%
\[
\left(  \text{\# of }\left(  231,321,312\right)  \text{-avoiding permutations
in }S_{n}\right)  =f_{n+1},
\]
where $\left(  f_{0},f_{1},f_{2},\ldots\right)  $ is the Fibonacci sequence
(as in Section \ref{sec.gf.exas}). \medskip

\textbf{(e)} \fbox{2} Prove that%
\[
\left(  \text{\# of }\left(  123,132,231\right)  \text{-avoiding permutations
in }S_{n}\right)  =n.
\]

[\textbf{Hint:} In parts \textbf{(a)}, \textbf{(b)} and \textbf{(d)}, the
permutations you are counting have already appeared in one of the previous
problems under a different guise.]
\end{exercise}

Let us connect $132$-avoidance with Lehmer codes:

\begin{exercise}
\textbf{(a)} \fbox{3} Let $n\in\mathbb{N}$ and $\sigma\in S_{n}$. Recall
Definition \ref{def.perm.lehmer1} \textbf{(a)}. Prove that $\sigma$ is
$132$-avoiding if and only if%
\[
\ell_{1}\left(  \sigma\right)  \geq\ell_{2}\left(  \sigma\right)  \geq
\cdots\geq\ell_{n}\left(  \sigma\right)  .
\]

\textbf{(b)} \fbox{2} Let $n\in\mathbb{N}$ and $k\in\mathbb{N}$. Prove that%
\[
\left(  \text{\# of }132\text{-avoiding permutations }\sigma\in S_{n}\text{
satisfying }\ell\left(  \sigma\right)  =k\right)  \leq p\left(  k\right)  ,
\]
where $p\left(  k\right)  $ is the $k$-th partition number (see Definition
\ref{def.pars.pn-pkn} \textbf{(b)}). Moreover, show that this inequality is an
equality if and only if $k<n$ or $n=k=0$.
\end{exercise}

Avoidance of multiple patterns $\tau\in S_{4}$ has also been studied. One of
the most interesting classes are the $\left(  2413,3142\right)  $-avoiding
permutations; these are known as the \emph{separable permutations}, and have
appeared in the study of polynomial inequalities (see \cite{Ghys17}). Let us
state just a few of their properties:

\begin{exercise}
\label{exe.perm.separable-equiv}If $\alpha\in S_{n}$ and $\beta\in S_{m}$ are
two permutations, then we define a new permutation $\alpha\oplus\beta\in
S_{n+m}$ by
\[
\left(  \alpha\oplus\beta\right)  \left(  k\right)  =%
\begin{cases}
\alpha\left(  k\right)  , & \text{if }k\leq n;\\
n+\beta\left(  k-n\right)  , & \text{if }k>n
\end{cases}
\ \ \ \ \ \ \ \ \ \ \text{for each }k\in\left[  n+m\right]  ,
\]
and we define a new permutation $\alpha\ominus\beta\in S_{n+m}$ by
\[
\left(  \alpha\ominus\beta\right)  \left(  k\right)  =%
\begin{cases}
n+\beta\left(  k\right)  , & \text{if }k\leq m;\\
\alpha\left(  k-m\right)  , & \text{if }k>m
\end{cases}
\ \ \ \ \ \ \ \ \ \ \text{for each }k\in\left[  n+m\right]  .
\]
Thus we have defined two binary operations $\oplus$ and $\ominus$ on the set
$S_{1}\cup S_{2}\cup S_{3}\cup\cdots$ of all permutations of the sets $\left[
1\right]  ,\left[  2\right]  ,\left[  3\right]  ,\ldots$. \medskip

\textbf{(a)} \fbox{1} Prove that both operations $\oplus$ and $\ominus$ are
associative. \medskip

\textbf{(b)} \fbox{7} Let $n>0$ and $\sigma\in S_{n}$. Prove that $\sigma$ is
separable if and only if $\sigma$ can be obtained from the permutation
$\operatorname*{id}\nolimits_{\left[  1\right]  }\in S_{1}$ by successive
applications of the operations $\oplus$ and $\ominus$. (For instance, $\left(
\operatorname*{id}\nolimits_{\left[  1\right]  }\ominus\operatorname*{id}%
\nolimits_{\left[  1\right]  }\right)  \ominus\left(  \operatorname*{id}%
\nolimits_{\left[  1\right]  }\oplus\operatorname*{id}\nolimits_{\left[
1\right]  }\right)  $ is the permutation in $S_{4}$ with OLN $\left(
3,4,2,1\right)  $.) \medskip

\textbf{(c)} \fbox{6} Prove that each $n>0$ satisfies%
\[
\left(  \text{\# of separable permutations }\sigma\in S_{n}\right)
=\sum_{k=0}^{n-1}\dfrac{1}{k+1}\dbinom{n-1+k}{k}\dbinom{n-1}{k}.
\]
(This number is known as the $\left(  n-1\right)  $-st \emph{large
Schr\"{o}der number}.)
\end{exercise}

The following is a result from \cite{AdeGri26}, specialized to the symmetric group:

\begin{exercise}
\label{exe.perm.adenbaum-len-sym}Let $n\in\mathbb{N}$. For each permutation
$\sigma\in S_{n}$, we let $F_{2}\left(  \sigma\right)  $ denote the set of all
pairs $\left(  \alpha,\beta\right)  \in S_{n}\times S_{n}$ satisfying
$\alpha\circ\beta=\sigma$. We say that a permutation $\sigma\in S_{n}$ is
\emph{length-symmetric} if for every $p,q\in\mathbb{N}$, we have%
\begin{align*}
&  \left(  \text{\# of pairs }\left(  \alpha,\beta\right)  \in F_{2}\left(
\sigma\right)  \text{ satisfying }\ell\left(  \alpha\right)  =p\text{ and
}\ell\left(  \beta\right)  =q\right) \\
&  =\left(  \text{\# of pairs }\left(  \alpha,\beta\right)  \in F_{2}\left(
\sigma\right)  \text{ satisfying }\ell\left(  \alpha\right)  =q\text{ and
}\ell\left(  \beta\right)  =p\right)
\end{align*}
(equivalently: if the generating function $\sum_{\left(  \alpha,\beta\right)
\in F_{2}\left(  \sigma\right)  }x^{\ell\left(  \alpha\right)  }y^{\ell\left(
\beta\right)  }\in\mathbb{Z}\left[  x,y\right]  $ is symmetric in $x$ and
$y$). \medskip

\textbf{(a)} \fbox{8} Prove that if $\sigma\in S_{n}$ is separable (that is,
$\left(  2413,3142\right)  $-avoiding), then $\sigma$ is length-symmetric.
\medskip

\textbf{(b)} \fbox{1} Find an example showing that the converse is not true.
\medskip

\textbf{(c)} \fbox{?} Find more general conditions under which $\sigma$ is length-symmetric.
\end{exercise}

\subsubsection{The cycle decomposition}

\begin{exercise}
\label{exe.perm.cycles-lcm}\fbox{1} Let $X$ be a finite set. Let $\sigma$ be a
permutation of $X$. Prove that the order of $\sigma$ in the symmetric group
$S_{X}$ equals the lcm of the lengths of all cycles of $\sigma$.
\end{exercise}

\begin{exercise}
\fbox{3} Let $X$ be a finite set. Let $\sigma$ and $\tau$ be two permutations
of $X$. Prove that $\sigma$ and $\tau$ are conjugate in the symmetric group
$S_{X}$ if and only if the cycle lengths partition of $\sigma$ equals the
cycle lengths partition of $\tau$.
\end{exercise}

The disjoint cycle decomposition of a permutation allows us to define its
\emph{reflection length}, which is an analogue of the Coxeter length that we
have defined in Definition \ref{def.perm.invs} \textbf{(b)}:

\begin{exercise}
\label{exe.perm.refl-len}Let $X$ be a finite set. Let $n=\left\vert
X\right\vert $.

The \emph{reflection length} (aka \emph{absolute length}) of a permutation
$\sigma\in S_{X}$ is defined to be $n-i$, where $i$ is the number of cycles
(including the $1$-cycles) in the disjoint cycle decomposition of $\sigma$.
(For example, any $k$-cycle in $S_{X}$ has reflection length $k-1$, since it
has $1$ cycle of length $k$ and $n-k$ cycles of length $1$.) The reflection
length of a permutation $\sigma\in S_{X}$ is denoted by $\ell
_{\operatorname*{refl}}\left(  \sigma\right)  $.

Prove the following: \medskip

\textbf{(a)} \fbox{1} For any $\sigma\in S_{X}$, we have $\ell
_{\operatorname*{refl}}\left(  \sigma^{-1}\right)  =\ell_{\operatorname*{refl}%
}\left(  \sigma\right)  $. \medskip

\textbf{(b)} \fbox{1} If $\sigma$ and $\tau$ are two conjugate elements in the
group $S_{X}$, then $\ell_{\operatorname*{refl}}\left(  \sigma\right)
=\ell_{\operatorname*{refl}}\left(  \tau\right)  $. \medskip

\textbf{(c)} \fbox{3} For any $\sigma\in S_{X}$ and any two distinct elements
$i$ and $j$ of $X$, we have%
\[
\ell_{\operatorname*{refl}}\left(  \sigma t_{i,j}\right)  =\ell
_{\operatorname*{refl}}\left(  t_{i,j}\sigma\right)  =%
\begin{cases}
\ell_{\operatorname*{refl}}\left(  \sigma\right)  +1, & \text{if we don't have
}i\overset{\sigma}{\sim}j;\\
\ell_{\operatorname*{refl}}\left(  \sigma\right)  -1, & \text{if
}i\overset{\sigma}{\sim}j.
\end{cases}
\]
Here, the notation \textquotedblleft$i\overset{\sigma}{\sim}j$%
\textquotedblright\ means \textquotedblleft$i$ and $j$ belong to the same
cycle of $\sigma$\textquotedblright\ (that is, \textquotedblleft there exists
some $p\in\mathbb{N}$ such that $i=\sigma^{p}\left(  j\right)  $%
\textquotedblright). \medskip

\textbf{(d)} \fbox{2} If $\sigma\in S_{X}$, then the number $\ell
_{\operatorname*{refl}}\left(  \sigma\right)  $ is the smallest $p\in
\mathbb{N}$ such that we can write $\sigma$ as a composition of $p$
transpositions. \medskip

\textbf{(e)} \fbox{2} For any $\sigma\in S_{X}$ and $\tau\in S_{X}$, we have
$\ell_{\operatorname*{refl}}\left(  \sigma\tau\right)  \leq\ell
_{\operatorname*{refl}}\left(  \sigma\right)  +\ell_{\operatorname*{refl}%
}\left(  \tau\right)  $. \medskip

\textbf{(f)} \fbox{2} For any $\sigma\in S_{X}$ and $\tau\in S_{X}$, we have
$\ell_{\operatorname*{refl}}\left(  \sigma\tau\right)  \equiv\ell
_{\operatorname*{refl}}\left(  \sigma\right)  +\ell_{\operatorname*{refl}%
}\left(  \tau\right)  \operatorname{mod}2$. \medskip

\textbf{(g)} \fbox{1} If $X=\left[  n\right]  $ and $\sigma\in S_{X}$, then
$\ell_{\operatorname*{refl}}\left(  \sigma\right)  \leq\ell\left(
\sigma\right)  $. (Keep in mind that the Coxeter length $\ell\left(
\sigma\right)  $ is defined only when $X=\left[  n\right]  $, not for general
finite sets $X$.)
\end{exercise}

\begin{exercise}
\label{exe.perm.transpos-code}Let $n\in\mathbb{N}$. Let $\sigma\in S_{n}$.
Define $\ell_{\operatorname*{refl}}\left(  \sigma\right)  $ as in Exercise
\ref{exe.perm.refl-len} (setting $X=\left[  n\right]  $). \medskip

\textbf{(a)} \fbox{4} Prove that there is a unique $n$-tuple $\left(
i_{1},i_{2},\ldots,i_{n}\right)  \in\left[  1\right]  \times\left[  2\right]
\times\cdots\times\left[  n\right]  $ such that%
\[
\sigma=t_{1,i_{1}}\circ t_{2,i_{2}}\circ\cdots\circ t_{n,i_{n}}.
\]
Here, we define $t_{i,i}$ to be the identity permutation $\operatorname*{id}%
\in S_{n}$ for each $i\in\left[  n\right]  $. \medskip

\textbf{(b)} \fbox{3} Consider this $n$-tuple $\left(  i_{1},i_{2}%
,\ldots,i_{n}\right)  $. Prove that
\[
\ell_{\operatorname*{refl}}\left(  \sigma\right)  =\left(  \text{\# of all
}k\in\left[  n\right]  \text{ satisfying }i_{k}\neq k\right)  .
\]

\end{exercise}

\begin{exercise}
Fix a commutative ring $K$ and a nonnegative integer $n\in\mathbb{N}$.

For each $\sigma\in S_{n}$, we define the \emph{permutation matrix}
$P_{\sigma}$ to be the $n\times n$-matrix%
\[
\left(  \left[  i=\sigma\left(  j\right)  \right]  \right)  _{i,j\in\left[
n\right]  }\in K^{n\times n}.
\]
(This is the $n\times n$-matrix whose $\left(  i,j\right)  $-th entry is
$\left[  i=\sigma\left(  j\right)  \right]  $, where we are using the notation
of Definition \ref{def.pars.iverson}.) For instance, if $n=4$ and
$\sigma=3124$ in one-line notation, then%
\[
P_{\sigma}=\left(
\begin{array}
[c]{cccc}%
0 & 1 & 0 & 0\\
0 & 0 & 1 & 0\\
1 & 0 & 0 & 0\\
0 & 0 & 0 & 1
\end{array}
\right)  .
\]

\textbf{(a)} \fbox{1} Let $\operatorname*{GL}\nolimits_{n}\left(  K\right)  $
be the group of all invertible $n\times n$-matrices over $K$. Prove that the
map
\begin{align*}
S_{n}  &  \rightarrow\operatorname*{GL}\nolimits_{n}\left(  K\right)  ,\\
\sigma &  \mapsto P_{\sigma}%
\end{align*}
is a group homomorphism. \medskip

\textbf{(b)} \fbox{2} Assuming that $K$ is a field, prove that each $\sigma\in
S_{n}$ satisfies $\operatorname*{rank}\left(  P_{\sigma}-I_{n}\right)
=\ell_{\operatorname*{refl}}\left(  \sigma\right)  $, where $\ell
_{\operatorname*{refl}}\left(  \sigma\right)  $ is defined as in Exercise
\ref{exe.perm.refl-len} (setting $X=\left[  n\right]  $). (Here, $I_{n}$
denotes the $n\times n$ identity matrix.) \medskip

\textbf{(c)} \fbox{7} Assuming that $K$ is a field, prove that two
permutations $\sigma,\tau\in S_{n}$ are conjugate in the group $S_{n}$ if and
only if their permutation matrices $P_{\sigma}$ and $P_{\tau}$ are similar
(i.e., conjugate in the group $\operatorname*{GL}\nolimits_{n}\left(
K\right)  $). \medskip

\textbf{(d)} \fbox{2} Prove the claim of part \textbf{(c)} more generally if
$K$ is any nontrivial commutative ring. \medskip

[\textbf{Hint:} For part \textbf{(c)}, it helps to show that the cycle lengths
partition of a permutation $\sigma$ can be uniquely recovered if one knows the
number of cycles of each of its powers $\sigma^{1},\sigma^{2},\sigma
^{3},\ldots$.]
\end{exercise}

\begin{exercise}
For each $n\in\mathbb{N}$, let $a_{n}$ denote the \# of all permutations
$\sigma\in S_{n}$ that have odd order in the symmetric group $S_{n}$. (By
Exercise \ref{exe.perm.cycles-lcm}, these are precisely the permutations whose
all cycles have odd length.) Also set $a_{n}:=0$ for all $n<0$. Prove the
following: \medskip

\textbf{(a)} \fbox{1} If $X$ is any finite set, then the \# of all
permutations $\sigma\in S_{X}$ that have odd order in $S_{X}$ equals
$a_{\left\vert X\right\vert }$. \medskip

\textbf{(b)} \fbox{3} We have $a_{n}=\sum\limits_{k\in\mathbb{N}\text{ even}%
}k!\dbinom{n-1}{k}a_{n-1-k}$ for each $n>0$. \medskip

Set $b_{n}:=\dfrac{a_{n}}{n!}$ for each $n\in\mathbb{N}$. Also, set $b_{n}:=0$
for all $n<0$. Prove that: \medskip

\textbf{(c)} \fbox{1} We have $nb_{n}=\sum\limits_{k\in\mathbb{N}\text{ even}%
}b_{n-1-k}=b_{n-1}+b_{n-3}+b_{n-5}+\cdots$ for each $n>0$. \medskip

Let $b\in\mathbb{Q}\left[  \left[  x\right]  \right]  $ be the FPS
$\sum\limits_{n\in\mathbb{N}}b_{n}x^{n}$. Prove the following: \medskip

\textbf{(d)} \fbox{1} We have $b^{\prime}=\underbrace{\left(  1+x^{2}%
+x^{4}+x^{6}+\cdots\right)  }_{=\sum_{k\in\mathbb{N}}x^{2k}}\cdot\,b$.
\medskip

\textbf{(e)} \fbox{2} We have $\operatorname*{Log}b=\dfrac{1}{2}\left(
\operatorname*{Log}\left(  1+x\right)  -\operatorname*{Log}\left(  1-x\right)
\right)  $. \medskip

\textbf{(f)} \fbox{1} We have $b=\left(  \dfrac{1+x}{1-x}\right)
^{1/2}=\left(  1+x\right)  \cdot\left(  1-x^{2}\right)  ^{-1/2}$. \medskip

\textbf{(g)} \fbox{2} We have $b_{n}=\left(  -1\right)  ^{\left\lfloor
n/2\right\rfloor }\dbinom{-1/2}{\left\lfloor n/2\right\rfloor }$ for each
$n\in\mathbb{N}$. \medskip

\textbf{(h)} \fbox{2} For each $n\in\mathbb{N}$, we have
\[
a_{n}=n!\cdot\left(  -1\right)  ^{\left\lfloor n/2\right\rfloor }\dbinom
{-1/2}{\left\lfloor n/2\right\rfloor }=%
\begin{cases}
c_{n}^{2}, & \text{if }n\text{ is even};\\
nc_{n}^{2}, & \text{if }n\text{ is odd,}%
\end{cases}
\]
where $c_{n}=1\cdot3\cdot5\cdot\cdots\cdot\left(  2\left\lfloor
n/2\right\rfloor -1\right)  $ is the product of all odd positive integers
smaller than $n$. \medskip

\textbf{(i)} \fbox{3} For each $n\in\mathbb{N}$, let $a_{n}^{\prime}$ be the
\# of all permutations $\sigma\in S_{n}$ whose all cycles have even length.
Show that%
\[
a_{n}^{\prime}=%
\begin{cases}
a_{n}, & \text{if }n\text{ is even};\\
0, & \text{if }n\text{ is odd}%
\end{cases}
\ \ \ \ \ \ \ \ \ \ \text{for each }n\in\mathbb{N}.
\]

[\textbf{Hint:} In part \textbf{(b)}, first choose the cycle of $\sigma$ that
contains the element $1$.]
\end{exercise}

\begin{exercise}
\label{exe.perm.star-transp-fac}\fbox{3} Let $n\in\mathbb{N}$. Let $X$ be an
$n$-element set. Fix some $x\in X$. An $x$\emph{-transposition} shall mean a
transposition of the form $t_{x,y}$ with $y\in X\setminus\left\{  x\right\}  $.

Let $\sigma\in S_{X}$ be a permutation that has $m$ cycles and has $f$ fixed
points distinct from $x$. Prove that $\sigma$ can be written as a product of
$n+m-2f-2$ many $x$-transpositions.
\end{exercise}

\subsubsection{Reduced words}

In this subsection, we shall take a closer look at how permutations in the
symmetric groups $S_{n}$ can be represented as products of simple
transpositions $s_{i}$. Most exercises here are particular cases of standard
results about Coxeter groups (see, e.g., \cite{BjoBre05} and \cite{Bourba02}
for introductions), but it is worth seeing them in the special yet rather
intuitive setting of symmetric groups.

We recall Definition \ref{def.perm.si}.

\begin{definition}
\label{def.perm.redwords.redword}Let $n\in\mathbb{N}$. Let $\sigma\in S_{n}$.
\medskip

\textbf{(a)} A \emph{Coxeter word} for $\sigma$ shall mean a tuple $\left(
i_{1},i_{2},\ldots,i_{k}\right)  \in\left[  n-1\right]  ^{k}$ satisfying
$\sigma=s_{i_{1}}s_{i_{2}}\cdots s_{i_{k}}$. \medskip

\textbf{(b)} A \emph{reduced word} for $\sigma$ shall mean a Coxeter word for
$\sigma$ that has the smallest length among all Coxeter words for $\sigma$.
\end{definition}

Note that Theorem \ref{thm.perm.len.redword1} \textbf{(a)} shows that any
permutation $\sigma\in S_{n}$ has a Coxeter word. Furthermore, Theorem
\ref{thm.perm.len.redword1} \textbf{(b)} says that the length $\ell\left(
\sigma\right)  $ of a permutation $\sigma\in S_{n}$ is the smallest length of
a Coxeter word for $\sigma$. Thus, a reduced word for a permutation $\sigma\in
S_{n}$ is the same as a Coxeter word for $\sigma$ that has length $\ell\left(
\sigma\right)  $. (This is why $\ell\left(  \sigma\right)  $ is called
\textquotedblleft length\textquotedblright\ of $\sigma$.)

\begin{example}
\label{exa.perm.redwords.1}Let $n=5$, and let $\sigma\in S_{5}$ be the
permutation whose OLN is $32415$. Then, $\sigma=s_{1}s_{3}s_{2}s_{1}$; thus,
$\left(  1,3,2,1\right)  $ is a Coxeter word for $\sigma$. Other Coxeter words
for $\sigma$ are $\left(  3,1,2,1\right)  $ and $\left(  3,2,1,2\right)  $ and
$\left(  1,3,1,1,2,1\right)  $ and $\left(  3,1,2,3,1,3\right)  $. The reduced
words for $\sigma$ are $\left(  1,3,2,1\right)  $ and $\left(  3,1,2,1\right)
$ and $\left(  3,2,1,2\right)  $.
\end{example}

Note that each permutation $\sigma\in S_{n}$ has finitely many reduced words,
but infinitely many Coxeter words (unless $n\leq1$). The identity permutation
$\operatorname*{id}$ has only one reduced word -- namely, the $0$-tuple
$\left(  {}\right)  $ -- but usually many Coxeter words, such as $\left(
1,2,3,3,2,1\right)  $.

We shall now study the combinatorics of Coxeter and reduced words of a
permutation $\sigma\in S_{n}$ in more depth. First, let us view them from a
different perspective:

\begin{definition}
Let $n\in\mathbb{N}$ and $\sigma\in S_{n}$. A \emph{sorting sequence} for
$\sigma$ shall mean a sequence $\left(  \sigma_{0},\sigma_{1},\ldots
,\sigma_{k}\right)  $ of permutations $\sigma_{i}\in S_{n}$ with the property
that $\sigma_{0}=\sigma$ and $\sigma_{k}=\operatorname*{id}$ and that for each
$i\in\left[  k\right]  $, the permutation $\sigma_{i}$ is obtained from
$\sigma_{i-1}$ by swapping two consecutive entries $\sigma_{i-1}\left(
h\right)  $ and $\sigma_{i-1}\left(  h+1\right)  $ into the correct order
(i.e., $\sigma_{i-1}\left(  h\right)  >\sigma_{i-1}\left(  h+1\right)  $, but
$\sigma_{i}\left(  h\right)  <\sigma_{i}\left(  h+1\right)  $).
\end{definition}

Thus, intuitively, a sorting sequence for $\sigma$ is a way of sorting its OLN
(i.e., the list $\sigma\left(  1\right)  \ \sigma\left(  2\right)
\ \cdots\ \sigma\left(  n\right)  $) into increasing order by repeatedly
swapping two out-of-order consecutive entries. For example, if $n=5$, and if
$\sigma\in S_{5}$ is the permutation whose OLN is $32415$, then%
\[
\left(  32415,\ 32145,\ 31245,\ 13245,\ 12345\right)
\]
is a sorting sequence of $\sigma$ (one of three such sequences).

\begin{exercise}
\fbox{2} Let $n\in\mathbb{N}$ and $\sigma\in S_{n}$. Find a bijection between
$\left\{  \text{reduced words for }\sigma\right\}  $ and $\left\{
\text{sorting sequences for }\sigma\right\}  $.
\end{exercise}

\begin{exercise}
\fbox{1} Let $n\in\mathbb{N}$. Let $\sigma\in S_{n}$. Let $\left(  i_{1}%
,i_{2},\ldots,i_{k}\right)  $ be a reduced word for $\sigma$. Let
$u,v\in\left\{  0,1,\ldots,k\right\}  $ be such that $u\leq v$. Prove that
$\left(  i_{u+1},i_{u+2},\ldots,i_{v}\right)  $ is a reduced word for
$s_{i_{u+1}}s_{i_{u+2}}\cdots s_{i_{v}}$.
\end{exercise}

\begin{exercise}
\fbox{2} Let $n\in\mathbb{N}$. Let $\sigma\in S_{n}$. Let $k\in\left[
n-1\right]  $. Prove that $\sigma$ has a reduced word whose last entry is $k$
if and only if $\sigma\left(  k\right)  >\sigma\left(  k+1\right)  $.
\end{exercise}

\begin{exercise}
\fbox{4} Let $n\in\mathbb{N}$. Let $\sigma\in S_{n}$. Let $u,v\in\left[
n-1\right]  $. Assume that $\sigma$ has a reduced word whose last entry is
$u$. Assume further that $\sigma$ has a reduced word whose last entry is $v$.
Prove the following: \medskip

\textbf{(a)} If $\left\vert u-v\right\vert >1$, then $\sigma$ has a reduced
word whose last two entries are $u$ and $v$ (in this order). \medskip

\textbf{(b)} If $\left\vert u-v\right\vert =1$, then $\sigma$ has a reduced
word whose last three entries are $u,v,u$ (in this order).
\end{exercise}

Let us now discuss some ways to transform Coxeter words. For instance:

\begin{itemize}
\item If a tuple of the form $\left(  \ldots,2,5,\ldots\right)  $ (that is, a
tuple that has two adjacent entries $2$ and $5$) is a Coxeter word for some
permutation $\sigma\in S_{n}$, then the tuple $\left(  \ldots,5,2,\ldots
\right)  $ (that is, the result of swapping these two adjacent entries) is a
Coxeter word for the same permutation $\sigma$, since Proposition
\ref{prop.perm.si.rules} \textbf{(b)} yields $s_{2}s_{5}=s_{5}s_{2}$.

\item If a tuple of the form $\left(  \ldots,2,3,2,\ldots\right)  $ is a
Coxeter word for some permutation $\sigma\in S_{n}$, then the tuple $\left(
\ldots,3,2,3,\ldots\right)  $ (that is, the result of replacing the three
adjacent entries $2,3,2$ by $3,2,3$, while leaving all remaining entries
unchanged) is a Coxeter word for the same permutation $\sigma$, since
Proposition \ref{prop.perm.si.rules} \textbf{(c)} yields $s_{2}s_{3}%
s_{2}=s_{3}s_{2}s_{3}$.

\item If a Coxeter word for $\sigma\in S_{n}$ has two adjacent entries that
are equal, then we can remove these two entries and still have a Coxeter word
for $\sigma$, since Proposition \ref{prop.perm.si.rules} \textbf{(a)} yields
$s_{i}s_{i}=s_{i}^{2}=\operatorname*{id}$ for each $i\in\left[  n-1\right]  $.
\end{itemize}

Generalizing these three observations, we obtain the following ways to change
a Coxeter word:

\begin{definition}
\label{def.perm.redwords.moves}Let $n\in\mathbb{N}$. Let $\sigma\in S_{n}$.

Given a Coxeter word $\mathbf{i}=\left(  i_{1},i_{2},\ldots,i_{k}\right)  $
for $\sigma$, we can obtain other Coxeter words for $\sigma$ by the following
three kinds of transformations: \medskip

\textbf{(a)} We can pick two adjacent entries $i_{u}$ and $i_{u+1}$ of
$\mathbf{i}$ that satisfy $\left\vert i_{u}-i_{u+1}\right\vert >1$, and swap
them (that is, replace the $u$-th and $\left(  u+1\right)  $-st entries of
$\mathbf{i}$ by $i_{u+1}$ and $i_{u}$, respectively). This is called a
\emph{commutation move}, and results in a new Coxeter word for $\sigma$, since
Proposition \ref{prop.perm.si.rules} \textbf{(b)} yields $s_{i_{u}}s_{i_{u+1}%
}=s_{i_{u+1}}s_{i_{u}}$.

For example, we can use such a move to transform the Coxeter word $\left(
1,2,3,1,2\right)  $ into $\left(  1,2,1,3,2\right)  $. \medskip

\textbf{(b)} We can pick three adjacent entries $i_{u}$, $i_{u+1}$ and
$i_{u+2}$ of $\mathbf{i}$ that satisfy $i_{u}=i_{u+2}=i_{u+1}\pm1$ (by which
we mean that we have either $i_{u}=i_{u+2}=i_{u+1}+1$ or $i_{u}=i_{u+2}%
=i_{u+1}-1$), and replace these three entries by $i_{u+1}$, $i_{u}$ and
$i_{u+1}$, respectively. This is called a \emph{braid move}, and results in a
new Coxeter word for $\sigma$, since we have $s_{i_{u}}s_{i_{u+1}%
}\underbrace{s_{i_{u+2}}}_{=s_{i_{u}}}=s_{i_{u}}s_{i_{u+1}}s_{i_{u}%
}=s_{i_{u+1}}s_{i_{u}}s_{i_{u+1}}$ by Proposition \ref{prop.perm.si.rules}
\textbf{(c)}.

For example, we can use such a move to transform the Coxeter word $\left(
1,2,1,3,2\right)  $ into $\left(  2,1,2,3,2\right)  $, and we can use another
such move to transform this result further into $\left(  2,1,3,2,3\right)  $.
\medskip

\textbf{(c)} We can pick two adjacent entries $i_{u}$ and $i_{u+1}$ of
$\mathbf{i}$ that are equal, and remove both of them from $\mathbf{i}$. This
is called a \emph{contraction move}, and results in a new Coxeter word for
$\sigma$, since we have $s_{i_{u}}\underbrace{s_{i_{u+1}}}_{=s_{i_{u}}%
}=s_{i_{u}}s_{i_{u}}=s_{i_{u}}^{2}=\operatorname*{id}$ by Proposition
\ref{prop.perm.si.rules} \textbf{(a)}.

For example, we can use such a move to transform the Coxeter word $\left(
3,2,2,1\right)  $ into $\left(  3,1\right)  $.
\end{definition}

Of course, these transformations can only be applied to a Coxeter word when
the respective requirements are met. For example, none of these
transformations can be applied to the Coxeter word $\left(  1,2\right)  $,
since it has neither two adjacent entries $i_{u}$ and $i_{u+1}$ that satisfy
$\left\vert i_{u}-i_{u+1}\right\vert >1$, nor three adjacent entries $i_{u}$,
$i_{u+1}$ and $i_{u+2}$ that satisfy $i_{u}=i_{u+2}=i_{u+1}\pm1$, nor two
adjacent entries $i_{u}$ and $i_{u+1}$ that are equal. On the other hand, we
can apply any of the three kinds of transformation to the Coxeter word
$\left(  4,1,2,1,2,2,4,4\right)  $, and we even have multiple choices for each
of them (e.g., we can apply a braid move to replace the \textquotedblleft%
$1,2,1$\textquotedblright\ by \textquotedblleft$2,1,2$\textquotedblright, but
we can also apply a braid move to replace the \textquotedblleft$2,1,2$%
\textquotedblright\ by \textquotedblleft$1,2,1$\textquotedblright\ instead).
Thus, starting with a single Coxeter word, we obtain a whole tapestry of
Coxeter words by applying moves to it. Note that each commutation move and
each braid move can be undone by a move of the same kind, while contraction
moves cannot be undone.

We notice that commutation moves and braid moves don't change the length of a
Coxeter word. Thus, if they are applied to a reduced word, they result in
another reduced word.

\begin{example}
\label{exam.perm.redword.graph}Let $n=4$, and let $\sigma\in S_{4}$ be the
permutation whose OLN is $4321$. It is easy to see that $\left(
1,2,1,3,2,1\right)  $ is a reduced word for $\sigma$. Omitting commas and
parentheses, we shorten this reduced word to $121321$. By repeatedly applying
commutation moves and braid moves, we can transform this reduced word into
$212321$, then further into $213231$, then into $231231$, and so on, and in
other directions too (as there are often several moves available). Let us draw
the result as a graph, with the nodes being all the reduced words that we
obtain, and the edges signifying braid moves and commutation moves (the thick
edges stand for commutation moves):%
\[%
\begin{tikzpicture}
\fill(80 : 4.5) circle (2pt) node[above right] {$321323$};
\fill(50 : 4.5) circle (2pt) node[right] {$321232$};
\fill(20 : 4.5) circle (2pt) node[right] {$312132$};
\fill(0 : 3.5) circle (2pt) node[left] {$312312$};
\fill(0 : 5) circle (2pt) node[right] {$132132$};
\fill(-80 : 4.5) circle (2pt) node[below right] {$123121$};
\fill(-50 : 4.5) circle (2pt) node[right] {$123212$};
\fill(-20 : 4.5) circle (2pt) node[right] {$132312$};
\draw(80 : 4.5) -- (50 : 4.5) -- (20 : 4.5);
\draw(-80 : 4.5) -- (-50 : 4.5) -- (-20 : 4.5);
\draw[ultra thick] (0 : 5) -- (20 : 4.5) -- (0 : 3.5) -- (-20 : 4.5) -- cycle;
\fill(100 : 4.5) circle (2pt) node[above left] {$323123$};
\fill(130 : 4.5) circle (2pt) node[left] {$232123$};
\fill(160 : 4.5) circle (2pt) node[left] {$231213$};
\fill(180 : 3.5) circle (2pt) node[right] {$231231$};
\fill(180 : 5) circle (2pt) node[left] {$213213$};
\fill(-100 : 4.5) circle (2pt) node[below left] {$121321$};
\fill(-130 : 4.5) circle (2pt) node[left] {$212321$};
\fill(-160 : 4.5) circle (2pt) node[left] {$213231$};
\draw(100 : 4.5) -- (130 : 4.5) -- (160 : 4.5);
\draw(-100 : 4.5) -- (-130 : 4.5) -- (-160 : 4.5);
\draw
[ultra thick] (180 : 5) -- (160 : 4.5) -- (180 : 3.5) -- (-160 : 4.5) -- cycle;
\draw[ultra thick] (80 : 4.5) -- (100 : 4.5);
\draw[ultra thick] (-80 : 4.5) -- (-100 : 4.5);
\end{tikzpicture}%
\]

It turns out that \textbf{each} reduced word for $\sigma$ appears as a node on
this graph. In other words, each reduced word for $\sigma$ can be obtained
from $121321$ by a sequence of commutation moves and braid moves. This is not
a coincidence, but a general result, known as \emph{Matsumoto's theorem for
the symmetric group}:
\end{example}

\begin{exercise}
\fbox{6} Let $n\in\mathbb{N}$. Let $\sigma\in S_{n}$. Let $\mathbf{i}$ and
$\mathbf{j}$ be two reduced words for $\sigma$. Prove that $\mathbf{i}$ can be
transformed into $\mathbf{j}$ by a sequence of commutation moves and braid moves.

[\textbf{Hint:} Induct on $\ell\left(  \sigma\right)  $.]
\end{exercise}

The graph in Example \ref{exam.perm.redword.graph} is furthermore bipartite;
better yet, any cycle has an even \# of thick edges and an even \# of thin
edges. This, too, is not a coincidence:

\begin{exercise}
Let $n\in\mathbb{N}$. Let $\sigma\in S_{n}$. Let $\mathbf{i}$ be a reduced
word for $\sigma$. Assume that we have transformed $\mathbf{i}$ into itself by
a sequence of commutation moves and braid moves. \medskip

\textbf{(a)} \fbox{3} Prove that the \# of braid moves in this sequence must
be even. \medskip

\textbf{(b)} \fbox{7} Prove that the \# of commutation moves in this sequence
must be even.
\end{exercise}

By throwing contraction moves into the mix, we can furthermore reduce
non-reduced Coxeter words:

\begin{exercise}
\fbox{6} Let $n\in\mathbb{N}$. Let $\sigma\in S_{n}$. Let $\mathbf{i}$ be a
Coxeter word for $\sigma$, and let $\mathbf{j}$ be a reduced word for $\sigma
$. Prove that $\mathbf{i}$ can be transformed into $\mathbf{j}$ by a sequence
of commutation moves, braid moves and contraction moves.
\end{exercise}

How many reduced words does a given permutation $\sigma\in S_{n}$ have? For
most $\sigma$, there is no nice formula for the answer. However, in at least
one specific case, a surprising (and deep) formula exists, which I am here
mentioning less as a reasonable exercise than as a curiosity:

\begin{exercise}
\fbox{50} Let $n\in\mathbb{N}$. Define the permutation $w_{0}\in S_{n}$ as in
Exercise \ref{exe.perm.w0.basics}. (Note that the $\sigma$ in Example
\ref{exam.perm.redword.graph} is the $w_{0}$ for $n=4$.)

Prove that the \# of reduced words for $w_{0}$ is%
\[
\dfrac{\dbinom{n}{2}!}{\prod_{i=1}^{n}\left(  2n-2i+1\right)  ^{i-1}}%
=\dbinom{n}{2}!\cdot\prod_{1\leq i<j\leq n}\dfrac{2}{i+j-1}.
\]

(\fbox{2} Prove the equality sign here.)
\end{exercise}

This number, incidentally, is the largest \# of reduced words that a
permutation in $S_{n}$ can have. On the other extreme (both of the \# of
reduced words and the difficulty of the proof), here is a characterization of
permutations that have a unique reduced word:

\begin{exercise}
\fbox{3} Let $n$ be a positive integer. Let $\sigma\in S_{n}$. Prove that
there is only one reduced word for $\sigma$ if and only if $\sigma$ has the
form $\operatorname*{cyc}\nolimits_{i,i+1,i+2,\ldots,j}$ or
$\operatorname*{cyc}\nolimits_{j,j-1,j-2,\ldots,i}$ for some $i,j\in\left[
n\right]  $ satisfying $i\leq j$. (If $i=j$, then these cycles are just
$\operatorname*{id}$.)
\end{exercise}

There is a connection between braid moves and pattern avoidance (Exercise
\ref{exe.perm.patt.132-av}):

\begin{exercise}
\fbox{4} Let $n\in\mathbb{N}$. Let $\sigma\in S_{n}$. Prove that $\sigma$ is
$321$-avoiding if and only if every two reduced words for $\sigma$ can be
transformed into each other by a sequence of commutation moves (without using
any braid moves).
\end{exercise}

As an application of reduced words, the following group-theoretical
characterization of the symmetric group $S_{n}$ easily follows:

\begin{exercise}
\fbox{3} Let $n\in\mathbb{N}$. Prove that the group $S_{n}$ is isomorphic to
the group with generators $g_{1},g_{2},\ldots,g_{n-1}$ and relations%
\begin{align*}
g_{i}^{2}  &  =1\ \ \ \ \ \ \ \ \ \ \text{for all }i\in\left[  n-1\right]  ;\\
g_{i}g_{j}  &  =g_{j}g_{i}\ \ \ \ \ \ \ \ \ \ \text{for all }i,j\in\left[
n-1\right]  \text{ satisfying }\left\vert i-j\right\vert >1;\\
g_{i}g_{i+1}g_{i}  &  =g_{i+1}g_{i}g_{i+1}\ \ \ \ \ \ \ \ \ \ \text{for all
}i\in\left[  n-2\right]  .
\end{align*}
(The isomorphism sends each $g_{i}$ to $s_{i}\in S_{n}$.)
\end{exercise}

This is known as the \emph{Coxeter-Moore presentation} of $S_{n}$.

\subsubsection{Descents}

\emph{Descents} are one of the most elementary features of a permutation
$\sigma\in S_{n}$: they are just the positions at which $\sigma$ decreases
(from that position to the next). Formally, they are defined as follows:

\begin{definition}
\label{def.perm.descents}Let $n\in\mathbb{N}$. Let $\sigma\in S_{n}$ be a
permutation. \medskip

\textbf{(a)} A \emph{descent} of $\sigma$ means an $i\in\left[  n-1\right]  $
such that $\sigma\left(  i\right)  >\sigma\left(  i+1\right)  $. \medskip

\textbf{(b)} The \emph{descent set} of $\sigma$ is defined to be the set of
all descents of $\sigma$. This set is denoted by $\operatorname*{Des}\sigma$.
\end{definition}

\begin{example}
The permutation $\sigma\in S_{7}$ with OLN $3146275$ has descents $1$ (since
$\sigma\left(  1\right)  >\sigma\left(  2\right)  $) and $4$ (since
$\sigma\left(  4\right)  >\sigma\left(  5\right)  $) and $6$ (since
$\sigma\left(  6\right)  >\sigma\left(  7\right)  $). Thus, it has descent set
$\operatorname*{Des}\sigma=\left\{  1,4,6\right\}  $.
\end{example}

\begin{exercise}
\label{exe.perm.descents.easy}\fbox{2} Let $n\in\mathbb{N}$. \medskip

\textbf{(a)} How many $\sigma\in S_{n}$ have exactly $0$ descents? \medskip

\textbf{(b)} How many $\sigma\in S_{n}$ have exactly $1$ descent? \medskip

\textbf{(c)} How many $\sigma\in S_{n}$ have exactly $n-1$ descents? \medskip

\textbf{(d)} Prove that the \# of all $\sigma\in S_{n}$ satisfying
$1\in\operatorname*{Des}\sigma$ (that is, $\sigma\left(  1\right)
>\sigma\left(  2\right)  $) is $\dfrac{n!}{2}$. (Here, we assume that $n\geq
2$.) \medskip

\textbf{(e)} Prove that the \# of all $\sigma\in S_{n}$ satisfying
$1,2\in\operatorname*{Des}\sigma$ (that is, $\sigma\left(  1\right)
>\sigma\left(  2\right)  >\sigma\left(  3\right)  $) is $\dfrac{n!}{6}$.
(Here, we assume that $n\geq3$.) \medskip

\textbf{(f)} How many $\sigma\in S_{n}$ satisfy $1,3\in\operatorname*{Des}%
\sigma$ (that is, $\sigma\left(  1\right)  >\sigma\left(  2\right)  $ and
$\sigma\left(  3\right)  >\sigma\left(  4\right)  $) ? (Here, we assume that
$n\geq4$.)
\end{exercise}

The following exercise generalizes parts \textbf{(d)}, \textbf{(e)} and
\textbf{(f)} of Exercise \ref{exe.perm.descents.easy}:

\begin{exercise}
\label{exe.perm.descents.subI}Let $n\in\mathbb{N}$. Let $I$ be a subset of
$\left[  n-1\right]  $. Write $I$ in the form $I=\left\{  c_{1},c_{2}%
,\ldots,c_{k}\right\}  $ with $c_{1}<c_{2}<\cdots<c_{k}$. Set $c_{0}:=0$ and
$c_{k+1}:=n$. For each $i\in\left[  k+1\right]  $, set $d_{i}:=c_{i}-c_{i-1}$.
Note that the $k+1$ numbers $d_{1},d_{2},\ldots,d_{k+1}$ are precisely the
lengths of the intervals into which the elements of $I$ subdivide the interval
$\left[  0,n\right]  $. \medskip

\textbf{(a)} \fbox{3} Prove that%
\[
\left(  \text{\# of }\sigma\in S_{n}\text{ satisfying }\operatorname{Des}%
\sigma\subseteq I\right)  =\dfrac{n!}{d_{1}!d_{2}!\cdots d_{k+1}!}.
\]

\textbf{(b)} \fbox{5} Let us use the notations from Definition
\ref{def.pars.qbinom.qint} \textbf{(b)}. Prove that%
\[
\sum_{\substack{\sigma\in S_{n};\\\operatorname*{Des}\sigma\subseteq
I}}q^{\ell\left(  \sigma\right)  }=\dfrac{\left[  n\right]  _{q}!}{\left[
d_{1}\right]  _{q}!\left[  d_{2}\right]  _{q}!\cdots\left[  d_{k+1}\right]
_{q}!}\ \ \ \ \ \ \ \ \ \ \text{in the ring }\mathbb{Z}\left[  q\right]
\text{.}%
\]

\end{exercise}

Note that Exercise \ref{exe.perm.descents.subI} \textbf{(b)} generalizes both
Exercise \ref{exe.perm.descents.subI} \textbf{(a)} (obtained by setting $q=1$)
and Proposition \ref{prop.perm.length.gf} (obtained by setting $I=\left[
n-1\right]  $ and $q=x$).

What about counting permutations $\sigma\in S_{n}$ satisfying
$\operatorname*{Des}\sigma=I$ rather than $\operatorname*{Des}\sigma\subseteq
I$ ? See Exercise \ref{exe.lgv.descents.isI} further below for this.

\begin{exercise}
\label{exe.perm.descents.stirl}\fbox{3} Let $n,k\in\mathbb{N}$ with $n>0$.
Recall from Exercise \ref{exe.fps.stirling2.1} that $S\left(  n,k\right)  $
denotes a Stirling number of the 2nd kind. Prove that $k!\cdot S\left(
n,k\right)  $ is the \# of pairs $\left(  \sigma,I\right)  $, where $\sigma\in
S_{n}$ is a permutation and $I$ is a $\left(  k-1\right)  $-element subset of
$\left[  n-1\right]  $ satisfying $\operatorname*{Des}\sigma\subseteq I$.
\end{exercise}

Meanwhile, let us connect descents with Eulerian polynomials:

\begin{exercise}
\label{exe.perm.descents.eulerian-pols}\fbox{5} Let $n$ be a positive integer.
Consider the polynomials $A_{m}\in\mathbb{Z}\left[  x\right]  $ defined for
all $m\in\mathbb{N}$ in Exercise \ref{exe.gf.eulerian-pol.basics}. Prove that%
\[
\sum_{\sigma\in S_{n}}x^{\left\vert \operatorname*{Des}\sigma\right\vert
+1}=A_{n}\ \ \ \ \ \ \ \ \ \ \text{in the ring }\mathbb{Z}\left[  x\right]  .
\]

(For example, for $n=4$, we have $\sum_{\sigma\in S_{4}}x^{\left\vert
\operatorname*{Des}\sigma\right\vert +1}=x+11x^{2}+11x^{3}+x^{4}=A_{4}$.)
\end{exercise}

This exercise can also be $q$-generalized:

\begin{exercise}
\label{exe.perm.descents.qeulerian-pols}\fbox{7} Let $n$ be a positive
integer. Consider the polynomials $A_{q,m}\in\mathbb{Z}\left[  q,x\right]  $
defined for all $m\in\mathbb{N}$ in Exercise \ref{exe.pars.qbinom.q-deriv}.
(We know from Exercise \ref{exe.pars.qbinom.q-deriv} \textbf{(h)} that they
are polynomials.) Prove that%
\[
\sum_{\sigma\in S_{n}}x^{\left\vert \operatorname*{Des}\sigma\right\vert
+1}q^{\operatorname*{sum}\left(  \operatorname*{Des}\sigma\right)  }%
=A_{q,n}\ \ \ \ \ \ \ \ \ \ \text{in the ring }\mathbb{Z}\left[  q,x\right]
.
\]
Here, $\operatorname*{sum}S$ is defined as in Proposition
\ref{prop.pars.qbinom.alt-defs} \textbf{(b)}.
\end{exercise}

The following exercise is a celebrated result of Foata:

\begin{exercise}
Let $n\in\mathbb{N}$. For any permutation $\sigma$, we let
$\operatorname*{maj}\sigma$ denote the sum of all descents of $\sigma$. (This
is known as the \emph{major index} of $\sigma$, and is precisely the number
$\operatorname*{sum}\left(  \operatorname*{Des}\sigma\right)  $ from Exercise
\ref{exe.perm.descents.qeulerian-pols}. For example, if $\sigma\in S_{5}$ has
OLN $34152$, then $\operatorname*{Des}\sigma=\left\{  2,4\right\}  $ and
$\operatorname*{maj}\sigma=2+4=6$.) \medskip

\textbf{(a)} \fbox{2} Prove that%
\[
\sum_{\sigma\in S_{n}}x^{\ell\left(  \sigma\right)  }=\sum_{\sigma\in S_{n}%
}x^{\operatorname*{maj}\sigma}\ \ \ \ \ \ \ \ \ \ \text{in }\mathbb{Z}\left[
x\right]  .
\]
In other words, prove that each $i\in\mathbb{N}$ satisfies%
\begin{align*}
&  \left(  \text{\# of permutations }\sigma\in S_{n}\text{ with }\ell\left(
\sigma\right)  =i\right) \\
&  =\left(  \text{\# of permutations }\sigma\in S_{n}\text{ with
}\operatorname*{maj}\sigma=i\right)  .
\end{align*}

\textbf{(b)} \fbox{8} Prove that%
\[
\sum_{\sigma\in S_{n}}x^{\ell\left(  \sigma\right)  }y^{\operatorname*{maj}%
\sigma}=\sum_{\sigma\in S_{n}}x^{\operatorname*{maj}\sigma}y^{\ell\left(
\sigma\right)  }=\sum_{\sigma\in S_{n}}x^{\operatorname*{maj}\sigma
}y^{\operatorname*{maj}\left(  \sigma^{-1}\right)  }%
\ \ \ \ \ \ \ \ \ \ \text{in }\mathbb{Z}\left[  x,y\right]  .
\]
In other words, prove that each $i\in\mathbb{N}$ and $j\in\mathbb{N}$ satisfy%
\begin{align*}
&  \left(  \text{\# of permutations }\sigma\in S_{n}\text{ with }\ell\left(
\sigma\right)  =i\text{ and }\operatorname*{maj}\sigma=j\right) \\
&  =\left(  \text{\# of permutations }\sigma\in S_{n}\text{ with
}\operatorname*{maj}\sigma=i\text{ and }\ell\left(  \sigma\right)  =j\right)
\\
&  =\left(  \text{\# of permutations }\sigma\in S_{n}\text{ with
}\operatorname*{maj}\sigma=i\text{ and }\operatorname*{maj}\left(  \sigma
^{-1}\right)  =j\right)  .
\end{align*}

\end{exercise}

The following exercise refines Exercise \ref{exe.perm.patt.132-av}
\textbf{(b)}:

\begin{exercise}
\textbf{(a)} \fbox{3} Let $n,k\in\mathbb{N}$ be such that $k<n$. Prove that%
\begin{align*}
&  \left(  \text{\# of }231\text{-avoiding permutations }\sigma\in S_{n}\text{
with }\left\vert \operatorname*{Des}\sigma\right\vert =k\right) \\
&  =\dfrac{1}{n}\dbinom{n}{k}\dbinom{n}{k+1}.
\end{align*}
(This number is known as the $\left(  n,k+1\right)  $-th \emph{Narayana
number}.) \medskip

\textbf{(b)} \fbox{1} Reprove Exercise \ref{exe.perm.patt.132-av} \textbf{(b)}
using this. \medskip

[\textbf{Hint:} In part \textbf{(a)}, the bivariate generating function
$C:=\sum_{n\in\mathbb{N}}\ \ \sum_{\substack{\sigma\in S_{n}\text{
is}\\231\text{-avoiding}}}t^{\left\vert \operatorname*{Des}\sigma\right\vert
}x^{n}\in\mathbb{Z}\left[  \left[  t,x\right]  \right]  $ can help. After
setting $B:=C-1$, show that $B=x\left(  1+B\right)  \left(  1+tB\right)  $.
Lagrange inversion (Exercise \ref{exe.fps.lagr-inv}) offers a comfortable way
of solving this equation.]
\end{exercise}

In contrast, the \# of $123$-avoiding permutations $\sigma\in S_{n}$ with
$\left\vert \operatorname*{Des}\sigma\right\vert =k$ is much more complicated
(see \cite{BaBoSi09}). The miracle of equinumerous avoidance classes (as seen
in Exercise \ref{exe.perm.patt.132-av} and Exercise \ref{exe.perm.patt.123-av}%
) does not persist here as we fix $\left\vert \operatorname*{Des}%
\sigma\right\vert $.

\subsection{Alternating sums, signed counting and determinants}

The notations of Chapter \ref{chap.sign} shall be used here.

\subsubsection{Cancellations in alternating sums}

\begin{exercise}
\label{exe.sign.cancel2}\fbox{3} Prove a generalization of Lemma
\ref{lem.sign.cancel3} in which $f$ is only required to be a bijection, not an
involution, but the assumption \textquotedblleft$\operatorname*{sign}I=0$ for
all $I\in\mathcal{X}$ satisfying $f\left(  I\right)  =I$\textquotedblright\ is
replaced by the stronger assumption \textquotedblleft$\operatorname*{sign}I=0$
for all $I\in\mathcal{X}$ and all \textbf{odd} $k\in\mathbb{N}$ satisfying
$f^{k}\left(  I\right)  =I$\textquotedblright.
\end{exercise}

\begin{exercise}
\label{exe.sign.dyck-signed-sum}Recall the concepts of Dyck words and Dyck
paths defined in Example 2 in Section \ref{sec.gf.exas}.

Let $n\in\mathbb{N}$.

If $w\in\left\{  0,1\right\}  ^{2n}$ is a $2n$-tuple, and if $k\in\left\{
0,1,\ldots,2n\right\}  $, then we define the $k$\emph{-height} $h_{k}\left(
w\right)  $ of $w$ to be the number%
\begin{align*}
&  \left(  \text{\# of }1\text{'s among the first }k\text{ entries of
}w\right) \\
&  \ \ \ \ \ \ \ \ \ \ -\left(  \text{\# of }0\text{'s among the first
}k\text{ entries of }w\right)  .
\end{align*}
If $w$ is a Dyck word, then this $k$-height $h_{k}\left(  w\right)  $ is a
nonnegative integer.

[For example, if $n=2$ and $w=\left(  1,0,0,1\right)  $, then $h_{3}\left(
w\right)  =1-2=-1<0$, which shows that $w$ is not a Dyck word.]

Furthermore, if $w\in\left\{  0,1\right\}  ^{2n}$ is a $2n$-tuple, then we
define the \emph{area} $\operatorname*{area}\left(  w\right)  $ of $w$ to be
the number%
\[
\operatorname*{area}\left(  w\right)  :=\sum_{k=0}^{2n}h_{k}\left(  w\right)
,
\]
and we define the \emph{sign} $\operatorname*{sign}\left(  w\right)  $ of $w$
to be the number $\left(  -1\right)  ^{\left(  \operatorname*{area}\left(
w\right)  -n\right)  /2}$ (we will soon see that this is well-defined).

[For example, if $n=5$ and $w=\left(  1,1,0,1,1,0,0,0,1,0\right)  $, then%
\begin{align*}
\operatorname*{area}\left(  w\right)   &  =\sum_{k=0}^{10}h_{k}\left(
w\right) \\
&  =\underbrace{h_{0}\left(  w\right)  }_{=0}+\underbrace{h_{1}\left(
w\right)  }_{=1}+\underbrace{h_{2}\left(  w\right)  }_{=2}+\underbrace{h_{3}%
\left(  w\right)  }_{=1}+\underbrace{h_{4}\left(  w\right)  }_{=2}%
+\underbrace{h_{5}\left(  w\right)  }_{=3}\\
&  \ \ \ \ \ \ \ \ \ \ +\underbrace{h_{6}\left(  w\right)  }_{=2}%
+\underbrace{h_{7}\left(  w\right)  }_{=1}+\underbrace{h_{8}\left(  w\right)
}_{=0}+\underbrace{h_{9}\left(  w\right)  }_{=1}+\underbrace{h_{10}\left(
w\right)  }_{=0}\\
&  =0+1+2+1+2+3+2+1+0+1+0=13
\end{align*}
and $\operatorname*{sign}\left(  w\right)  =\left(  -1\right)  ^{\left(
\operatorname*{area}\left(  w\right)  -n\right)  /2}=\left(  -1\right)
^{\left(  13-5\right)  /2}=1$. The names \textquotedblleft$k$%
-height\textquotedblright\ and \textquotedblleft area\textquotedblright\ are
not accidental: If $w$ is a Dyck word, then the \textquotedblleft
heights\textquotedblright\ $h_{0}\left(  w\right)  ,h_{1}\left(  w\right)
,\ldots,h_{2n}\left(  w\right)  $ really are the heights (i.e., the
y-coordinates) of the points on the Dyck path corresponding to the Dyck word
$w$; furthermore, the number $\operatorname*{area}\left(  w\right)  $ really
is the area of the \textquotedblleft mountain range\textquotedblright\ under
the Dyck path.] \medskip

\textbf{(a)} \fbox{1} Prove that any $w\in\left\{  0,1\right\}  ^{2n}$
satisfies $\left(  \operatorname*{area}\left(  w\right)  -n\right)
/2\in\mathbb{Z}$ (so that $\left(  -1\right)  ^{\left(  \operatorname*{area}%
\left(  w\right)  -n\right)  /2}$ really is well-defined). \medskip

\textbf{(b)} \fbox{2} Prove that a $2n$-tuple $w\in\left\{  0,1\right\}
^{2n}$ is a Dyck word of length $2n$ if and only if it satisfies%
\[
h_{2i-1}\left(  w\right)  \geq0\ \ \ \ \ \ \ \ \ \ \text{for all }i\in\left\{
1,2,\ldots,n\right\}
\]
and $h_{2n}\left(  w\right)  =0$. \medskip

\textbf{(c)} \fbox{4} Recall the Catalan numbers $c_{0},c_{1},c_{2},\ldots$ as
introduced in Section \ref{sec.gf.exas}. Assume that $n$ is a positive
integer. Prove that%
\[
\sum_{\substack{w\text{ is a Dyck word}\\\text{of length }2n}%
}\operatorname*{sign}\left(  w\right)  =%
\begin{cases}
\left(  -1\right)  ^{\left(  n-1\right)  /2}c_{\left(  n-1\right)  /2}, &
\text{if }n\text{ is odd;}\\
0, & \text{if }n\text{ is even.}%
\end{cases}
\]

[\textbf{Hint:} In part \textbf{(c)}, find a sign-reversing involution on a
certain set of Dyck words of length $2n$ that preserves all the
\textquotedblleft odd heights\textquotedblright\ $h_{1}\left(  w\right)
,h_{3}\left(  w\right)  ,\ldots,h_{2n-1}\left(  w\right)  $ while changing one
of the \textquotedblleft even heights\textquotedblright\ $h_{k}\left(
w\right)  $ by $2$.]
\end{exercise}

\begin{noncompile}
\begin{exercise}
THE\ FOLLOWING\ NEEDS\ FIXING!

\fbox{3} Let $n,k\in\mathbb{N}$. Prove that%
\[
\sum_{i=0}^{\left\lfloor n/2\right\rfloor }\left(  -1\right)  ^{i}%
\dbinom{n+1-i}{i}\dbinom{n-2i}{k-i}=%
\begin{cases}
1, & \text{if }k=0\text{;}\\
0, & \text{if }k>0\text{.}%
\end{cases}
\]

[\textbf{Hint:} There are algebraic ways to prove this, but here is a
combinatorial one: How many pairs $\left(  I,J\right)  $ of disjoint subsets
of $\left[  n\right]  $ exist with the property that no $i\in I$ satisfies
$i+1\in I\cup J$ ? Recall Exercise \ref{exe.lacunar.counts} \textbf{(b)}.]
\end{exercise}
\end{noncompile}

\begin{exercise}
\textbf{(a)} \fbox{2} Use Exercise \ref{exe.perm.descents.stirl} to give a
new, combinatorial proof of the identity $\sum_{k=0}^{n}\left(  -1\right)
^{k}k!\cdot S\left(  n,k\right)  =\left(  -1\right)  ^{n}$ from Exercise
\ref{exe.fps.stirling2.1} \textbf{(c)}. \medskip

\textbf{(b)} \fbox{2} Reprove the identity $\sum_{k=1}^{n}\left(  -1\right)
^{k-1}\left(  k-1\right)  !\cdot S\left(  n,k\right)  =%
\begin{cases}
1, & \text{if }n=1;\\
0, & \text{if }n\neq1
\end{cases}
$ from Exercise \ref{exe.fps.stirling2.1} \textbf{(d)} in a similar
combinatorial way. \medskip

[\textbf{Hint:} For part \textbf{(b)}, focus on the permutations $\sigma$
satisfying $\sigma\left(  1\right)  =1$.]
\end{exercise}

The next exercise is not about alternating sums, but rather about proving the
$q$-Lucas theorem (Theorem \ref{thm.sign.q-lucas}):

\begin{exercise}
\label{exe.sign.qbinom.qlucas}Let $K$ be a field. Let $d$ be a positive
integer. Let $\omega$ be a primitive $d$-th root of unity in $K$. \medskip

\textbf{(a)} \fbox{2} Prove that $\dbinom{d}{k}_{\omega}=0$ for each
$k\in\left\{  1,2,\ldots,d-1\right\}  $. \medskip

Now, let $A$ be a noncommutative $K$-algebra, and let $a,b\in A$ be such that
$ba=\omega ab$. \medskip

\textbf{(b)} \fbox{1} Prove that $\left(  a+b\right)  ^{d}=a^{d}+b^{d}$.
\medskip

\textbf{(c)} \fbox{3} Prove that $a^{d}$ and $b^{d}$ commute with both $a$ and
$b$ (that is, we have $uv=vu$ for each $u\in\left\{  a^{d},b^{d}\right\}  $
and each $v\in\left\{  a,b\right\}  $). \medskip

\textbf{(d)} \fbox{5} Prove Theorem \ref{thm.sign.q-lucas}. \medskip

[\textbf{Hint:} For part \textbf{(a)}, show that $\dbinom{n}{k}_{q}%
=\dfrac{\left[  n\right]  _{q}}{\left[  k\right]  _{q}}\dbinom{n-1}{k-1}_{q}$
for all $n>0$ and $k>0$. For part \textbf{(d)}, first construct a
noncommutative $K$-algebra $A$ and two elements $a,b\in A$ satisfying
$ba=\omega ab$ and such that all the monomials $a^{i}b^{j}$ are $K$-linearly
independent. Use Exercise \ref{exe.pars.qbinom.2to1} for this. In this
$K$-algebra, expand both sides of $\left(  a+b\right)  ^{n}=\left(  \left(
a+b\right)  ^{d}\right)  ^{q}\left(  a+b\right)  ^{r}$. Alternatively, there
is a commutative approach using Theorem \ref{thm.pars.qbinom.binom1}.]
\end{exercise}

\subsubsection{The principles of inclusion and exclusion}

\begin{exercise}
\fbox{3} Let $n\in\mathbb{N}$. Prove that $\sum\limits_{k=0}^{n}\left(
-1\right)  ^{k}\dbinom{n}{k}\left(  n-k\right)  ^{n+1}=\dbinom{n+1}{2}\cdot
n!$.
\end{exercise}

\begin{exercise}
\fbox{1} Let $n,m\in\mathbb{N}$. Prove that
\[
\left(  \text{\# of surjective maps }f:\left[  m\right]  \rightarrow\left[
n\right]  \right)  =n!\cdot S\left(  m,n\right)  ,
\]
where $S\left(  m,n\right)  $ is the Stirling number of the 2nd kind (as
defined in Exercise \ref{exe.fps.stirling2.1}).
\end{exercise}

The next exercise is concerned with the derangement numbers $D_{n}$ from
Definition \ref{def.pie.dera}.

\begin{exercise}
\label{exe.pie.dera.recs}\textbf{(a)} \fbox{1} Prove that $D_{n}%
=nD_{n-1}+\left(  -1\right)  ^{n}$ for all $n\geq1$. \medskip

\textbf{(b)} \fbox{1} Prove that $D_{n}=\left(  n-1\right)  \left(
D_{n-1}+D_{n-2}\right)  $ for all $n\geq2$. \medskip

\textbf{(c)} \fbox{2} Prove that $n!=\sum_{k=0}^{n}\dbinom{n}{k}D_{n-k}$ for
all $n\in\mathbb{N}$. \medskip

\textbf{(d)} \fbox{1} Show that $\sum_{n\in\mathbb{N}}\dfrac{D_{n}}{n!}%
x^{n}=\dfrac{\exp\left[  -x\right]  }{1-x}$ in the FPS ring $\mathbb{Q}\left[
\left[  x\right]  \right]  $.
\end{exercise}

\begin{exercise}
\fbox{3} Reprove Theorem \ref{thm.fps.comps.num-wpcomps-n-k} using the PIE.
\end{exercise}

\begin{exercise}
\fbox{4} For any $n,m\in\mathbb{N}$, we define a polynomial $Z_{m,n}%
\in\mathbb{Z}\left[  x\right]  $ by%
\[
Z_{m,n}=\sum_{k=0}^{n}\left(  -1\right)  ^{k}\dbinom{n}{k}\left(
x^{n-k}-1\right)  ^{m}.
\]
Prove that $Z_{m,n}=Z_{n,m}$ for all $m,n\in\mathbb{N}$.
\end{exercise}

\begin{exercise}
Let $n$ be a positive integer. Let $a_{1},a_{2},\ldots,a_{n}$ be any $n$
integers. \medskip

\textbf{(a)} \fbox{4} Show that%
\[
\max\left\{  a_{1},a_{2},\ldots,a_{n}\right\}  =\sum_{k=1}^{n}\left(
-1\right)  ^{k-1}\sum_{1\leq i_{1}<i_{2}<\cdots<i_{k}\leq n}\min\left\{
a_{i_{1}},a_{i_{2}},\ldots,a_{i_{k}}\right\}  .
\]

\textbf{(b)} \fbox{2} More generally: Show that%
\[
F\left(  \max\left\{  a_{1},a_{2},\ldots,a_{n}\right\}  \right)  =\sum
_{k=1}^{n}\left(  -1\right)  ^{k-1}\sum_{1\leq i_{1}<i_{2}<\cdots<i_{k}\leq
n}F\left(  \min\left\{  a_{i_{1}},a_{i_{2}},\ldots,a_{i_{k}}\right\}  \right)
\]
for any function $F:\mathbb{Z}\rightarrow\mathbb{R}$.
\end{exercise}

The following exercise is about a sequence of rather useful identities,
sometimes known as the \emph{polarization identities}\footnote{See, e.g.,
\cite[proof of Lemma 2.4.1]{Smith95} for a striking application of part
\textbf{(c)} to invariant theory.}:

\begin{exercise}
Let $n\in\mathbb{N}$. Let $A$ be a commutative ring. Let $v_{1},v_{2}%
,\ldots,v_{n}\in A$ and $w\in A$. Prove the following: \medskip

\textbf{(a)} \fbox{4} For each $m\in\mathbb{N}$, we have%
\[
\sum_{I\subseteq\left[  n\right]  }\left(  -1\right)  ^{n-\left\vert
I\right\vert }\left(  w+\sum_{i\in I}v_{i}\right)  ^{m}=\sum
_{\substack{\left(  i_{1},i_{2},\ldots,i_{m}\right)  \in\left\{
0,1,\ldots,n\right\}  ^{m};\\\left[  n\right]  \subseteq\left\{  i_{1}%
,i_{2},\ldots,i_{m}\right\}  }}v_{i_{1}}v_{i_{2}}\cdots v_{i_{m}},
\]
where we set $v_{0}:=w$. \medskip

\textbf{(b)} \fbox{1} For each $m\in\left\{  0,1,\ldots,n-1\right\}  $, we
have%
\[
\sum_{I\subseteq\left[  n\right]  }\left(  -1\right)  ^{n-\left\vert
I\right\vert }\left(  w+\sum_{i\in I}v_{i}\right)  ^{m}=0.
\]

\textbf{(c)} \fbox{1} We have%
\[
\sum_{I\subseteq\left[  n\right]  }\left(  -1\right)  ^{n-\left\vert
I\right\vert }\left(  w+\sum_{i\in I}v_{i}\right)  ^{n}=n!v_{1}v_{2}\cdots
v_{n}.
\]

\textbf{(d)} \fbox{2} We have%
\[
\sum_{I\subseteq\left[  n\right]  }\left(  -1\right)  ^{n-\left\vert
I\right\vert }\left(  \sum_{i\in I}v_{i}-\sum_{i\in\left[  n\right]  \setminus
I}v_{i}\right)  ^{n}=2^{n}n!v_{1}v_{2}\cdots v_{n}.
\]

\end{exercise}

\begin{exercise}
\fbox{5} Let $A$ and $B$ be two finite sets. Let $R$ be a subset of $A\times
B$.

For any subset $X$ of $A$, we define $M\left(  X\right)  $ to be the set%
\[
\left\{  b\in B\ \mid\ \text{there exists some }x\in X\text{ such that
}\left(  x,b\right)  \in R\right\}  .
\]

For any subset $Y$ of $B$, we define $N\left(  Y\right)  $ to be the set%
\[
\left\{  a\in A\ \mid\ \text{there exists some }y\in Y\text{ such that
}\left(  a,y\right)  \in R\right\}  .
\]

Prove that%
\[
\sum_{\substack{X\subseteq A;\\M\left(  X\right)  =B}}\left(  -1\right)
^{\left\vert X\right\vert }=\sum_{\substack{Y\subseteq B;\\N\left(  Y\right)
=A}}\left(  -1\right)  ^{\left\vert Y\right\vert }.
\]

[\textbf{Remark:} Those familiar with graph theory can think of $A$, $B$ and
$R$ as forming a bipartite graph (with vertex set $A\sqcup B$ and edge set
$R$). In that case, $M\left(  X\right)  $ is the \textquotedblleft neighbor
set\textquotedblright\ of $X$ (that is, the set of all vertices that have at
least one neighbor in $X$), and likewise $N\left(  Y\right)  $ is the
\textquotedblleft neighbor set\textquotedblright\ of $Y$.]
\end{exercise}

\begin{exercise}
\label{exe.pie.screaming-toes}\fbox{5} Let $n>1$ be an integer. Consider $n$
people standing in a circle. Each of them looks down at someone else's feet
(i.e., at the feet of one of the other $n-1$ persons). A bell sounds, and
every person (simultaneously) looks up at the eyes of the person whose feet
they have been ogling. If two people make eye contact, they scream. Show that
the probability that no one screams is%
\[
\sum_{k=0}^{n}\left(  -1\right)  ^{k}\dfrac{n\left(  n-1\right)  \cdots\left(
n-2k+1\right)  }{\left(  n-1\right)  ^{2k}\cdot2^{k}\cdot k!}.
\]

Here is a combinatorial restatement of the question (if you prefer not to deal
with probabilities): A pair $\left(  i,j\right)  $ of elements of $\left[
n\right]  $ is said to \emph{scream} at a map $f:\left[  n\right]
\rightarrow\left[  n\right]  $ if it satisfies $f\left(  i\right)  =j$ and
$f\left(  j\right)  =i$. A map $f:\left[  n\right]  \rightarrow\left[
n\right]  $ is \emph{silent} if no pair $\left(  i,j\right)  \in\left[
n\right]  \times\left[  n\right]  $ screams at $f$. Prove that the \# of all
silent maps $f:\left[  n\right]  \rightarrow\left[  n\right]  $ is%
\[
\sum_{k=0}^{n}\left(  -1\right)  ^{k}\dfrac{n\left(  n-1\right)  \cdots\left(
n-2k+1\right)  }{2^{k}\cdot k!}\left(  n-1\right)  ^{n-2k}.
\]

\end{exercise}

The following two exercises show some applications of the methods of Chapter
\ref{chap.sign} to graph theory.

\begin{exercise}
\label{exe.sign.chromatic.1}Let $G$ be a finite undirected graph with vertex
set $V$ and edge set $E$. Fix $n\in\mathbb{N}$.

An $n$\emph{-coloring} of $G$ means a map $c:V\rightarrow\left[  n\right]  $.
If $c:V\rightarrow\left[  n\right]  $ is an $n$-coloring, then we regard the
values $c\left(  v\right)  $ of $c$ as the \textquotedblleft
colors\textquotedblright\ of the respective vertices $v$.

An $n$-coloring $c$ of $G$ is said to be \emph{proper} if there exists no edge
of $G$ whose two endpoints $v$ and $w$ satisfy $c\left(  v\right)  =c\left(
w\right)  $. (In other words, an $n$-coloring of $G$ is said to be proper if
and only if there is no edge whose two endpoints have the same color.)

Let $\chi_{G}\left(  n\right)  $ denote the \# of proper $n$-colorings of $G$.
\medskip

\textbf{(a)} \fbox{5} Prove that
\[
\chi_{G}\left(  n\right)  =\sum_{F\subseteq E}\left(  -1\right)  ^{\left\vert
F\right\vert }n^{\operatorname*{conn}\left(  V,F\right)  },
\]
where $\operatorname*{conn}\left(  V,F\right)  $ denotes the \# of connected
components of the graph with vertex set $V$ and edge set $F$.

This shows, in particular, that $\chi_{G}\left(  n\right)  $ is a polynomial
function in $n$. (The corresponding polynomial is known as the \emph{chromatic
polynomial} of $G$.) \medskip

\textbf{(b)} \fbox{1} Find an explicit formula for $\chi_{G}\left(  n\right)
$ if $G$ is a path graph $%
\raisebox{-0.7pc}{
\begin{tikzpicture}%
[scale=1.5,thick,main node/.style={circle,fill=blue!20,draw}]
\node[main node] (1) at (1, 0) {$1$};
\node[main node] (2) at (2, 0) {$2$};
\node[main node] (3) at (3, 0) {$3$};
\node(4) at (4, 0) {$\cdots$};
\node[main node] (n) at (5, 0) {$m$};
\draw(1) -- (2) -- (3) -- (4) -- (n);
\end{tikzpicture}
}%
$ with $m$ vertices. \medskip

\textbf{(c)} \fbox{2} Find an explicit formula for $\chi_{G}\left(  n\right)
$ if $G$ is a cycle graph with $m$ vertices.
\end{exercise}

\begin{exercise}
\fbox{7} Let $G$ be a finite undirected graph with vertex set $V$ and edge set
$E$. Fix a vertex $v\in V$.

Given any subset $F$ of $E$, we define an $F$\emph{-path} to be a path of $G$
whose edges all belong to $F$.

A subset $F$ of $E$ is said to \emph{infect} an edge $e\in E$ if there is an
$F$-path leading from $v$ to some endpoint of $e$. (Note that this is
automatically satisfied if $v$ is an endpoint of $e$, since the empty path is
always an $F$-path.)

A subset $F$ of $E$ is said to be \emph{pandemic} if it infects each edge
$e\in E$.

Prove that%
\[
\sum_{\substack{F\subseteq E\text{ is}\\\text{pandemic}}}\left(  -1\right)
^{\left\vert F\right\vert }=\left[  E=\varnothing\right]  .
\]

[\textbf{Example:} Let $G$ be the following graph:%
\[%
\begin{tikzpicture}%
[-,>=stealth',shorten >=1pt,auto,node distance=3cm, thick,main node/.style={circle,fill=blue!20,draw}%
]
\node[main node] (1) {$v$};
\node[main node] [above of=1] (2) {$p$};
\node[main node] [right of=1] (3) {$w$};
\node[main node] [above of=3] (4) {$q$};
\node[main node] [right of=3] (5) {$t$};
\node[main node] [above of=5] (6) {$r$};
\path[every node/.style={font=\sffamily\small}] (1) edge node {$1$}
(2) (2) edge node {$2$} (4) (4) edge node {$3$}
(6) (6) edge [bend left] node {$4$} (5) (5) edge node {$5$}
(3) (3) edge node {$6$} (1) (5) edge [bend left] node {$7$}
(6) (4) edge node {$8$} (3);
\end{tikzpicture}%
\]
(where the vertex $v$ is the vertex labelled $v$). Then, for example, the set
$\left\{  1,2\right\}  \subseteq E$ infects edges $1,2,3,6,8$ (but none of the
other edges). The set $\left\{  1,2,5\right\}  $ infects the same edges as
$\left\{  1,2\right\}  $ (indeed, the additional edge $5$ does not increase
its infectiousness, since it is not on any $\left\{  1,2,5\right\}  $-path
from $v$). The set $\left\{  1,2,3\right\}  $ infects every edge other than
$5$. The set $\left\{  1,2,3,4\right\}  $ infects each edge, and thus is pandemic.]
\end{exercise}

\begin{exercise}
\fbox{3} Let $K$ be a commutative ring. Let $n\in\mathbb{N}$. Let $A=\left(
a_{i,j}\right)  _{1\leq i\leq n,\ 1\leq j\leq n}\in K^{n\times n}$ be an
$n\times n$-matrix. Then, the \emph{permanent} $\operatorname*{per}A$ of $A$
is defined to be the element%
\[
\sum_{\sigma\in S_{n}}a_{1,\sigma\left(  1\right)  }a_{2,\sigma\left(
2\right)  }\cdots a_{n,\sigma\left(  n\right)  }%
\]
of $K$ (where $S_{n}$ is the $n$-th symmetric group). Prove the \emph{Ryser
formula}%
\[
\operatorname*{per}A=\left(  -1\right)  ^{n}\sum_{I\subseteq\left[  n\right]
}\left(  -1\right)  ^{\left\vert I\right\vert }\prod_{j=1}^{n}\ \ \sum_{i\in
I}a_{i,j}.
\]

\end{exercise}

The following exercise is a variant of Theorem \ref{thm.pie.moeb}:

\begin{exercise}
\label{exe.pie.moeb-sym}\fbox{2} Let $S$ be a finite set. Let $A$ be any
additive abelian group.

For each subset $I$ of $S$, let $a_{I}$ and $b_{I}$ be two elements of $A$.

Assume that%
\[
b_{I}=\sum_{J\subseteq I}\left(  -1\right)  ^{\left\vert J\right\vert }%
a_{J}\ \ \ \ \ \ \ \ \ \ \text{for all }I\subseteq S.
\]

Then, prove that we also have%
\[
a_{I}=\sum_{J\subseteq I}\left(  -1\right)  ^{\left\vert J\right\vert }%
b_{J}\ \ \ \ \ \ \ \ \ \ \text{for all }I\subseteq S.
\]

\end{exercise}

The next exercise is an analogue of Exercise \ref{exe.pie.moeb-sym} with sets
replaced by numbers:

\begin{exercise}
\label{exe.pie.bin-inv}\fbox{2} Let $A$ be any additive abelian group. Let
$\left(  a_{0},a_{1},\ldots,a_{n}\right)  $ and $\left(  b_{0},b_{1}%
,\ldots,b_{n}\right)  $ be two $\left(  n+1\right)  $-tuples of elements of
$A$. Assume that%
\[
b_{m}=\sum_{i=0}^{m}\left(  -1\right)  ^{i}\dbinom{m}{i}a_{i}%
\ \ \ \ \ \ \ \ \ \ \text{for all }m\in\left\{  0,1,\ldots,n\right\}  .
\]
Then, prove that we also have%
\[
a_{m}=\sum_{i=0}^{m}\left(  -1\right)  ^{i}\dbinom{m}{i}b_{i}%
\ \ \ \ \ \ \ \ \ \ \text{for all }m\in\left\{  0,1,\ldots,n\right\}  .
\]

[\textbf{Hint:} There is a direct proof, but it is perhaps neater to derive
this from Exercise \ref{exe.pie.moeb-sym}.] \medskip

[\textbf{Remark:} This result is known as \emph{binomial inversion} (in one of
its forms). The $\left(  n+1\right)  $-tuple $\left(  b_{0},b_{1},\ldots
,b_{n}\right)  $ is called the \emph{binomial transform} of $\left(
a_{0},a_{1},\ldots,a_{n}\right)  $.]
\end{exercise}

The next exercise is a \textquotedblleft mirror version\textquotedblright\ of
Exercise \ref{exe.pie.bin-inv}:

\begin{exercise}
\label{exe.enum.anti-bin-inv}\fbox{2} Let $A$ be any additive abelian group.
Let $\left(  a_{0},a_{1},\ldots,a_{n}\right)  $ and $\left(  b_{0}%
,b_{1},\ldots,b_{n}\right)  $ be two $\left(  n+1\right)  $-tuples of elements
of $A$. Assume that%
\[
b_{m}=\sum_{i=m}^{n}\left(  -1\right)  ^{i}\dbinom{i}{m}a_{i}%
\ \ \ \ \ \ \ \ \ \ \text{for all }m\in\left\{  0,1,\ldots,n\right\}  .
\]
Prove that we also have%
\[
a_{m}=\sum_{i=m}^{n}\left(  -1\right)  ^{i}\dbinom{i}{m}b_{i}%
\ \ \ \ \ \ \ \ \ \ \text{for all }m\in\left\{  0,1,\ldots,n\right\}  .
\]

[\textbf{Hint:} Restate Exercise \ref{exe.pie.bin-inv} as a claim that a
certain matrix is inverse to itself. Now recall that the inverse of the
transpose of a matrix is the transpose of its inverse. Alternatively, prove
the claim directly.]
\end{exercise}

\begin{noncompile}
\begin{exercise}
\label{exe.pie.cheb-bin-inv}\textbf{(a)} \fbox{1} Show that all $m\in
\mathbb{N}$ and $i\in\left\{  0,1,\ldots,m\right\}  $ satisfy $\dfrac{m}%
{m-i}\dbinom{m-i}{i}\in\mathbb{Z}$. Here and in the following, we agree to
understand $\dfrac{m}{m-i}$ as being $1$ if $m=i=0$. \medskip

\textbf{(b)} \fbox{4} Let $A$ be any additive abelian group. Let $\left(
a_{0},a_{1},\ldots,a_{n}\right)  $ and $\left(  b_{0},b_{1},\ldots
,b_{n}\right)  $ be two $\left(  n+1\right)  $-tuples of elements of $A$.
Prove that we have%
\[
b_{m}=\sum_{i=0}^{\left\lfloor m/2\right\rfloor }\dbinom{m}{i}a_{m-2i}%
\ \ \ \ \ \ \ \ \ \ \text{for all }m\in\left\{  0,1,\ldots,n\right\}
\]
if and only if we have%
\[
a_{m}=\sum_{i=0}^{\left\lfloor m/2\right\rfloor }\left(  -1\right)  ^{i}%
\dfrac{m}{m-i}\dbinom{m-i}{i}b_{m-2i}\ \ \ \ \ \ \ \ \ \ \text{for all }%
m\in\left\{  0,1,\ldots,n\right\}  .
\]

\end{exercise}
\end{noncompile}

The next few exercises show some ways of generalizing the Principle of
Inclusion and Exclusion (in its original form -- Theorem \ref{thm.pie.1}). The
first one replaces the question \textquotedblleft how many elements of $U$
belong to none of the $n$ subsets $A_{1},A_{2},\ldots,A_{n}$\textquotedblright%
\ by \textquotedblleft how many elements of $U$ belong to exactly $k$ of the
$n$ subsets $A_{1},A_{2},\ldots,A_{n}$\textquotedblright:

\begin{exercise}
\label{exe.pie.dan1}\fbox{4} Let $n\in\mathbb{N}$, and let $U$ be a finite
set. Let $A_{1},A_{2},\ldots,A_{n}$ be $n$ subsets of $U$. Let $k\in
\mathbb{N}$. Let
\[
S_{k}:=\left\{  u\in U\ \mid\ \text{the \# of all }i\in\left[  n\right]
\text{ satisfying }u\in A_{i}\text{ is }k\right\}  .
\]
Show that%
\[
\left\vert S_{k}\right\vert =\sum_{I\subseteq\left[  n\right]  }\left(
-1\right)  ^{\left\vert I\right\vert -k}\dbinom{\left\vert I\right\vert }%
{k}\left(  \text{\# of }u\in U\text{ that satisfy }u\in A_{i}\text{ for all
}i\in I\right)  .
\]

\end{exercise}

\begin{exercise}
\label{exe.pie.bon}\fbox{5} Let $n\in\mathbb{N}$, and let $U$ be a finite set.
Let $A_{1},A_{2},\ldots,A_{n}$ be $n$ subsets of $U$. Let $m\in\mathbb{N}$.
\medskip

\textbf{(a)} For each $u\in U$, let $c\left(  u\right)  $ be the \# of all
$i\in\left[  n\right]  $ satisfying $u\in A_{i}$. Show that%
\begin{align*}
&  \sum_{\substack{I\subseteq\left[  n\right]  ;\\\left\vert I\right\vert \leq
m}}\left(  -1\right)  ^{\left\vert I\right\vert }\left(  \text{\# of }u\in
U\text{ that satisfy }u\in A_{i}\text{ for all }i\in I\right) \\
&  =\left(  -1\right)  ^{m}\sum_{u\in U}\dbinom{c\left(  u\right)  -1}{m}.
\end{align*}

\textbf{(b)} Conclude the \emph{Bonferroni inequalities}, which say that%
\[
\sum_{\substack{I\subseteq\left[  n\right]  ;\\\left\vert I\right\vert \leq
m}}\left(  -1\right)  ^{m-\left\vert I\right\vert }\left(  \text{\# of }u\in
U\text{ that satisfy }u\in A_{i}\text{ for all }i\in I\right)  \geq0
\]
if $U=A_{1}\cup A_{2}\cup\cdots\cup A_{n}$.
\end{exercise}

\begin{exercise}
\textbf{(a)} \fbox{4} Find a common generalization of Exercise
\ref{exe.pie.dan1} and Exercise \ref{exe.pie.bon} that has the form%
\begin{align*}
&  \sum_{\substack{I\subseteq\left[  n\right]  ;\\\left\vert I\right\vert \leq
m}}\left(  -1\right)  ^{\left\vert I\right\vert }\dbinom{\left\vert
I\right\vert }{k}\left(  \text{\# of }u\in U\text{ that satisfy }u\in
A_{i}\text{ for all }i\in I\right) \\
&  =\sum_{u\in U}\ \ \ \underbrace{\cdots}_{\text{some expression involving
}c\left(  u\right)  }.
\end{align*}

\textbf{(b)} \fbox{2} Generalize this further by including weights on the
elements of $U$ (similarly to how Theorem \ref{thm.pie.2} generalizes Theorem
\ref{thm.pie.1}).
\end{exercise}

The next exercise generalizes Theorem \ref{thm.pie.2} in a similar way as
$q$-binomial coefficients generalize binomial coefficients:

\begin{exercise}
\label{exe.pie.1+q.1}\textbf{(a)} \fbox{4} Let $n\in\mathbb{N}$, and let $U$
be a finite set. Let $A_{1},A_{2},\ldots,A_{n}$ be $n$ subsets of $U$. Let $K$
be any commutative ring. Let $w:U\rightarrow K$ be any map (i.e., let
$w\left(  u\right)  $ be an element of $K$ for each $u\in U$). Let $q\in K$.
Prove that%
\[
\sum_{u\in U}\left(  1+q\right)  ^{\left(  \text{\# of }i\in\left[  n\right]
\text{ satisfying }u\in A_{i}\right)  }w\left(  u\right)  =\sum_{I\subseteq
\left[  n\right]  }q^{\left\vert I\right\vert }\sum_{\substack{u\in U;\\u\in
A_{i}\text{ for all }i\in I}}w\left(  u\right)
\]
in $K$. \medskip

\textbf{(b)} \fbox{1} Derive Theorem \ref{thm.pie.2} as a particular case of
part \textbf{(a)}. \medskip

\textbf{(c)} \fbox{2} Prove that each $n\in\mathbb{N}$ satisfies%
\[
\sum_{\sigma\in S_{n}}q^{\left\vert \operatorname*{Fix}\sigma\right\vert
}=\sum_{k=0}^{n}\dfrac{n!}{k!}\left(  q-1\right)  ^{k}%
\]
in the polynomial ring $\mathbb{Z}\left[  q\right]  $. (See Definition
\ref{def.perm.fix} for the definition of $\operatorname*{Fix}\sigma$.)
\end{exercise}

Next comes another counting problem that can be solved in many ways:

\begin{exercise}
\fbox{3} Let $A$ be a finite additive abelian group (with its neutral element
denoted by $0$). Let $n\in\mathbb{N}$. Show that%
\begin{align*}
&  \left(  \text{\# of }n\text{-tuples }\left(  a_{1},a_{2},\ldots
,a_{n}\right)  \in\left(  A\setminus\left\{  0\right\}  \right)  ^{n}\text{
such that }a_{1}+a_{2}+\cdots+a_{n}=0\right) \\
&  =\dfrac{\left(  \left\vert A\right\vert -1\right)  ^{n}+\left(  -1\right)
^{n}\left(  \left\vert A\right\vert -1\right)  }{\left\vert A\right\vert }.
\end{align*}

\end{exercise}

Next comes a generalization of Theorem \ref{thm.pars.odd-dist-equal}:

\begin{exercise}
\fbox{3} Let $n\in\mathbb{N}$ and $k\in\mathbb{N}$.

Let $p_{\operatorname*{odd},k}\left(  n\right)  $ be the \# of partitions of
$n$ that have exactly $k$ distinct even parts. (For instance, the partition
$\left(  7,5,4,4,3,2\right)  $ has exactly $2$ distinct even parts, namely $4$
and $2$.)

Let $p_{\operatorname*{dist},k}\left(  n\right)  $ be the \# of partitions
$\lambda$ of $n$ that have exactly $k$ numbers appear in $\lambda$ more than
once (in the sense that there are exactly $k$ distinct integers $i$ such that
$i$ appears more than once in $\lambda$). (For instance, the partition
$\left(  7,4,2,2,1,1,1\right)  $ has exactly $2$ numbers appear more than
once, namely $2$ and $1$.)

Prove that%
\[
p_{\operatorname*{odd},k}\left(  n\right)  =p_{\operatorname*{dist},k}\left(
n\right)  .
\]

\end{exercise}

\begin{exercise}
\fbox{5} Solve Exercise \ref{exe.fps.laurent.reed-dawson} \textbf{(b)} again
using the Principle of Inclusion and Exclusion.
\end{exercise}

\subsubsection{More subtractive methods}

The following exercise gives a variation and a generalization of Theorem
\ref{thm.cancel.all-even}:

\begin{exercise}
Let $n\in\mathbb{N}$ and $d\in\mathbb{N}$. \medskip

\textbf{(a)} \fbox{2} An $n$-tuple $\left(  x_{1},x_{2},\ldots,x_{n}\right)
\in\left[  d\right]  ^{n}$ is said to be \textbf{all-odd} if each element of
$\left[  d\right]  $ occurs an odd number of times in this $n$-tuple (i.e., if
for each $k\in\left[  d\right]  $, the number of all $i\in\left[  n\right]  $
satisfying $x_{i}=k$ is odd). Prove that the number of all all-odd $n$-tuples
$\left(  x_{1},x_{2},\ldots,x_{n}\right)  \in\left[  d\right]  ^{n}$ is
\[
\dfrac{1}{2^{d}}\sum_{k=0}^{d}\left(  -1\right)  ^{k}\dbinom{d}{k}\left(
d-2k\right)  ^{n}.
\]

\textbf{(b)} \fbox{2} More generally: Let $p_{1},p_{2},\ldots,p_{d}$ be $d$
integers. Let $e$ be the number of $k\in\left[  d\right]  $ for which $p_{k}$
is odd. Prove that the number of all $n$-tuples $\left(  x_{1},x_{2}%
,\ldots,x_{n}\right)  \in\left[  d\right]  ^{n}$ that satisfy%
\[
\left(  \text{\# of }i\in\left[  n\right]  \text{ satisfying }x_{i}=k\right)
\equiv p_{k}\operatorname{mod}2\ \ \ \ \ \ \ \ \ \ \text{for all }k\in\left[
d\right]
\]
is
\[
\dfrac{1}{2^{d}}\sum_{i=0}^{d}\ \ \sum_{j=0}^{e}\left(  -1\right)  ^{j}%
\dbinom{e}{j}\dbinom{d-e}{i-j}\left(  d-2i\right)  ^{n}.
\]

\end{exercise}

\subsubsection{Determinants}

We fix a commutative ring $K$.

\Needspace{10pc}

\begin{exercise}
\label{exe.det.arrowhead}\fbox{3} Let $n$ be a positive integer. Let
$a_{1},a_{2},\ldots,a_{n}\in K$ and $b_{1},b_{2},\ldots,b_{n-1}\in K$ and
$c_{1},c_{2},\ldots,c_{n-1}\in K$. Let $A$ be the $n\times n$-matrix%
\[
\left(
\begin{array}
[c]{ccccc}%
a_{1} & 0 & \cdots & 0 & c_{1}\\
0 & a_{2} & \cdots & 0 & c_{2}\\
\vdots & \vdots & \ddots & \vdots & \vdots\\
0 & 0 & \cdots & a_{n-1} & c_{n-1}\\
b_{1} & b_{2} & \cdots & b_{n-1} & a_{n}%
\end{array}
\right)  .
\]
(This is the matrix whose $\left(  i,j\right)  $-th entry is $%
\begin{cases}
a_{i}, & \text{if }i=j;\\
b_{j}, & \text{if }i=n\text{ and }j\neq n;\\
c_{i}, & \text{if }i\neq n\text{ and }j=n;\\
0, & \text{if }i\neq n\text{ and }j\neq n\text{ and }i\neq j
\end{cases}
$ for all $i\in\left[  n\right]  $ and $j\in\left[  n\right]  $.) Prove that%
\[
\det A=a_{1}a_{2}\cdots a_{n}-\sum_{i=1}^{n-1}b_{i}c_{i}\prod_{\substack{j\in
\left[  n-1\right]  ;\\j\neq i}}a_{j}.
\]

\end{exercise}

\begin{exercise}
\label{exe.det.schur-lem}\fbox{3} Let $n\in\mathbb{N}$. Let $A$ be an $n\times
n$-matrix. Let $b_{1},b_{2},\ldots,b_{n}$ be $n$ elements of $K$. Prove that%
\[
\sum\limits_{k=1}^{n}\det\left(  \left(  A_{i,j}b_{i}^{\left[  j=k\right]
}\right)  _{1\leq i\leq n,\ 1\leq j\leq n}\right)  =\left(  b_{1}+b_{2}%
+\cdots+b_{n}\right)  \det A
\]
(where we are using Definition \ref{def.iverson}). Equivalently (rewritten in
a friendlier but longer form): Prove that%
\begin{align*}
&  \det\left(
\begin{array}
[c]{cccc}%
A_{1,1}b_{1} & A_{1,2} & \cdots & A_{1,n}\\
A_{2,1}b_{2} & A_{2,2} & \cdots & A_{2,n}\\
\vdots & \vdots & \ddots & \vdots\\
A_{n,1}b_{n} & A_{n,2} & \cdots & A_{n,n}%
\end{array}
\right)  +\det\left(
\begin{array}
[c]{cccc}%
A_{1,1} & A_{1,2}b_{1} & \cdots & A_{1,n}\\
A_{2,1} & A_{2,2}b_{2} & \cdots & A_{2,n}\\
\vdots & \vdots & \ddots & \vdots\\
A_{n,1} & A_{n,2}b_{n} & \cdots & A_{n,n}%
\end{array}
\right) \\
&  \ \ \ \ \ \ \ \ \ \ +\cdots+\det\left(
\begin{array}
[c]{cccc}%
A_{1,1} & A_{1,2} & \cdots & A_{1,n}b_{1}\\
A_{2,1} & A_{2,2} & \cdots & A_{2,n}b_{2}\\
\vdots & \vdots & \ddots & \vdots\\
A_{n,1} & A_{n,2} & \cdots & A_{n,n}b_{n}%
\end{array}
\right) \\
&  =\left(  b_{1}+b_{2}+\cdots+b_{n}\right)  \det\left(
\begin{array}
[c]{cccc}%
A_{1,1} & A_{1,2} & \cdots & A_{1,n}\\
A_{2,1} & A_{2,2} & \cdots & A_{2,n}\\
\vdots & \vdots & \ddots & \vdots\\
A_{n,1} & A_{n,2} & \cdots & A_{n,n}%
\end{array}
\right)  .
\end{align*}

\end{exercise}

\begin{exercise}
\label{exe.det.chio}Let $n$ be a positive integer. Let $A\in K^{n\times n}$ be
an $n\times n$-matrix. \medskip

\textbf{(a)} \fbox{4} Prove that the equality
\[
\det\left(  \left(  A_{i,j}A_{n,n}-A_{i,n}A_{n,j}\right)  _{1\leq i\leq
n-1,\ 1\leq j\leq n-1}\right)  =A_{n,n}^{n-2}\cdot\det A
\]
holds if the element $A_{n,n}$ of $K$ is invertible. \medskip

\textbf{(b)} \fbox{2} Prove that this equality also holds if $n\geq2$ (whether
or not $A_{n,n}$ is invertible). \medskip

[\textbf{Hint:} For part \textbf{(a)}, observe that $A_{i,j}A_{n,n}%
-A_{i,n}A_{n,j}=A_{n,n}\cdot\left(  A_{i,j}-\dfrac{A_{i,n}}{A_{n,n}}%
A_{n,j}\right)  $.]
\end{exercise}

The following two exercises give some applications of determinants:

\begin{exercise}
Let $n\in\mathbb{N}$. Let $a_{1},a_{2},\ldots,a_{n}\in K$ and $b_{1}%
,b_{2},\ldots,b_{n}\in K$. \medskip

\textbf{(a)} \fbox{2} Use the Cauchy--Binet identity (Theorem \ref{thm.det.CB}%
, applied to appropriate $2\times n$- and $n\times2$-matrices) to show that
\[
\left(  \sum_{k=1}^{n}a_{k}^{2}\right)  \left(  \sum_{k=1}^{n}b_{k}%
^{2}\right)  -\left(  \sum_{k=1}^{n}a_{k}b_{k}\right)  ^{2}=\sum_{1\leq
i<j\leq n}\left(  a_{i}b_{j}-a_{j}b_{i}\right)  ^{2}.
\]

\textbf{(b)} \fbox{1} If $K=\mathbb{R}$, then conclude the
\emph{Cauchy--Schwarz inequality}%
\[
\left(  \sum_{k=1}^{n}a_{k}^{2}\right)  \left(  \sum_{k=1}^{n}b_{k}%
^{2}\right)  \geq\left(  \sum_{k=1}^{n}a_{k}b_{k}\right)  ^{2}.
\]

\end{exercise}

\begin{exercise}
Let $n$ be a positive integer. \medskip

\textbf{(a)} \fbox{2} Prove that
\[
\sum_{\substack{\sigma\in S_{n}\text{ is a}\\\text{derangement}}}\left(
-1\right)  ^{\sigma}=\left(  -1\right)  ^{n-1}\left(  n-1\right)  .
\]
(See Definition \ref{def.pie.dera} for the notion of a derangement.) \medskip

\textbf{(b)} \fbox{2} Prove that%
\[
\sum_{\sigma\in S_{n}}\left(  -1\right)  ^{\sigma}x^{\left\vert
\operatorname*{Fix}\sigma\right\vert }=\left(  x+n-1\right)  \left(
x-1\right)  ^{n-1}%
\]
for any $x\in K$. (See Definition \ref{def.perm.fix} \textbf{(b)} for the
definition of $\left\vert \operatorname*{Fix}\sigma\right\vert $.) \medskip

\textbf{(c)} \fbox{3} Prove that%
\[
\sum_{\sigma\in S_{n}}\dfrac{\left(  -1\right)  ^{\sigma}}{\left\vert
\operatorname*{Fix}\sigma\right\vert +1}=\left(  -1\right)  ^{n+1}\dfrac
{n}{n+1}.
\]
(This is Problem B6 on the Putnam competition 2005.) \medskip

\textbf{(d)} \fbox{2} Let $j\in\left[  n\right]  $. Assume that $n>1$. Prove
that%
\[
\sum_{\substack{\sigma\in S_{n};\\\sigma\left(  j\right)  =j}}\dfrac{\left(
-1\right)  ^{\sigma}}{\left\vert \operatorname*{Fix}\sigma\right\vert
}=\left(  -1\right)  ^{n}\dfrac{n-1}{n}.
\]

\end{exercise}

In the next exercise, you are asked to reconstruct a proof of the Vandermonde
determinant (specifically, of Theorem \ref{thm.det.vander} \textbf{(d)}) using
a special kind of directed graphs -- the \emph{tournaments}. This is by far
not the easiest proof of Theorem \ref{thm.det.vander} \textbf{(d)}, but is
perhaps the most combinatorial. \Needspace{20pc}

\begin{exercise}
\label{exe.det.vander-tour}We define a \emph{tournament} to be a simple
directed graph with the property that for any two distinct vertices $i$ and
$j$, exactly one of the arcs $\left(  i,j\right)  $ and $\left(  j,i\right)  $
belongs to the graph. For example, there are $8$ tournaments with vertex set
$\left[  3\right]  $, namely%
\begin{equation}%
\begin{tabular}
[c]{|c|c|c|c|}\hline
$%
\begin{tikzpicture}%
[->,shorten >=1pt, thick,main node/.style={circle,fill=blue!20,draw}]
\node[main node] (1) at (-60 : 1) {1};
\node[main node] (2) at (60 : 1) {2};
\node[main node] (3) at (180 : 1) {3};
\path[-{Stealth[length=4mm]}]
(1) edge (2)
(1) edge (3)
(2) edge (3);
\end{tikzpicture}%
$ & $%
\begin{tikzpicture}%
[->,shorten >=1pt, thick,main node/.style={circle,fill=blue!20,draw}]
\node[main node] (1) at (-60 : 1) {1};
\node[main node] (2) at (60 : 1) {2};
\node[main node] (3) at (180 : 1) {3};
\path[-{Stealth[length=4mm]}]
(1) edge (2)
(1) edge (3)
(3) edge (2);
\end{tikzpicture}%
$ & $%
\begin{tikzpicture}%
[->,shorten >=1pt, thick,main node/.style={circle,fill=blue!20,draw}]
\node[main node] (1) at (-60 : 1) {1};
\node[main node] (2) at (60 : 1) {2};
\node[main node] (3) at (180 : 1) {3};
\path[-{Stealth[length=4mm]}]
(1) edge (2)
(3) edge (1)
(2) edge (3);
\end{tikzpicture}%
$ & $%
\begin{tikzpicture}%
[->,shorten >=1pt, thick,main node/.style={circle,fill=blue!20,draw}]
\node[main node] (1) at (-60 : 1) {1};
\node[main node] (2) at (60 : 1) {2};
\node[main node] (3) at (180 : 1) {3};
\path[-{Stealth[length=4mm]}]
(1) edge (2)
(3) edge (1)
(3) edge (2);
\end{tikzpicture}%
$\\\hline
$%
\begin{tikzpicture}%
[->,shorten >=1pt, thick,main node/.style={circle,fill=blue!20,draw}]
\node[main node] (1) at (-60 : 1) {1};
\node[main node] (2) at (60 : 1) {2};
\node[main node] (3) at (180 : 1) {3};
\path[-{Stealth[length=4mm]}]
(2) edge (1)
(1) edge (3)
(2) edge (3);
\end{tikzpicture}%
$ & $%
\begin{tikzpicture}%
[->,shorten >=1pt, thick,main node/.style={circle,fill=blue!20,draw}]
\node[main node] (1) at (-60 : 1) {1};
\node[main node] (2) at (60 : 1) {2};
\node[main node] (3) at (180 : 1) {3};
\path[-{Stealth[length=4mm]}]
(2) edge (1)
(1) edge (3)
(3) edge (2);
\end{tikzpicture}%
$ & $%
\begin{tikzpicture}%
[->,shorten >=1pt, thick,main node/.style={circle,fill=blue!20,draw}]
\node[main node] (1) at (-60 : 1) {1};
\node[main node] (2) at (60 : 1) {2};
\node[main node] (3) at (180 : 1) {3};
\path[-{Stealth[length=4mm]}]
(2) edge (1)
(3) edge (1)
(2) edge (3);
\end{tikzpicture}%
$ & $%
\begin{tikzpicture}%
[->,shorten >=1pt, thick,main node/.style={circle,fill=blue!20,draw}]
\node[main node] (1) at (-60 : 1) {1};
\node[main node] (2) at (60 : 1) {2};
\node[main node] (3) at (180 : 1) {3};
\path[-{Stealth[length=4mm]}]
(2) edge (1)
(3) edge (1)
(3) edge (2);
\end{tikzpicture}%
$\\\hline
\end{tabular}
\ \ \ \ \ . \label{eq.exe.det.vander-tour.tour3}%
\end{equation}
(Note that \textquotedblleft simple directed graph\textquotedblright\ implies
that any arc is merely a pair of two distinct vertices; thus, in particular,
there are no arcs of the form $\left(  i,i\right)  $.) \medskip

Fix $n\in\mathbb{N}$. Let $T$ be the set of all tournaments with vertex set
$\left[  n\right]  $. It is easy to see that $\left\vert T\right\vert
=2^{n\left(  n-1\right)  /2}$.

For any permutation $\sigma\in S_{n}$, we define $P_{\sigma}\in T$ to be the
tournament with vertex set $\left[  n\right]  $ and with arcs%
\[
\left(  \sigma\left(  i\right)  ,\sigma\left(  j\right)  \right)
\ \ \ \ \ \ \ \ \ \ \text{for all }i\in\left[  n\right]  \text{ and }%
j\in\left[  n\right]  \text{ satisfying }i<j.
\]

(For example, in the above table (\ref{eq.exe.det.vander-tour.tour3}) of
tournaments with vertex set $\left[  3\right]  $, the first tournament is
$P_{\operatorname*{id}}$, while the second tournament is $P_{s_{2}}$.)

We define the \emph{scoreboard} $\operatorname*{scb}D$ of a tournament $D\in
T$ to be the $n$-tuple $\left(  s_{1},s_{2},\ldots,s_{n}\right)  \in
\mathbb{N}^{n}$, where
\begin{align*}
s_{j}:=  &  \left(  \text{\# of arcs of }D\text{ that end at }j\right) \\
=  &  \left(  \text{\# of }i\in\left[  n\right]  \text{ such that }\left(
i,j\right)  \text{ is an arc of }D\right)
\end{align*}
for each $j\in\left[  n\right]  $.

We say that a tournament $D\in T$ is \emph{injective} if all $n$ entries of
its scoreboard $\operatorname*{scb}D$ are distinct. \medskip

\textbf{(a)} \fbox{2} Prove that a tournament $D\in T$ is injective if and
only if it has the form $P_{\sigma}$ for some $\sigma\in S_{n}$. \medskip

\textbf{(b)} \fbox{2} Prove that the tournaments $P_{\sigma}$ for all
$\sigma\in S_{n}$ are distinct. \medskip

Now, let $a_{1},a_{2},\ldots,a_{n}$ be $n$ elements of $K$. For each
tournament $D\in T$, we define the following:

\begin{itemize}
\item For each arc $e=\left(  i,j\right)  $ of $D$, we define the
\emph{weight} $w\left(  e\right)  $ of $e$ to be $\left(  -1\right)  ^{\left[
i>j\right]  }a_{j}$ (where we are using Definition \ref{def.iverson}).

\item We define the \emph{weight} $w\left(  D\right)  $ of $D$ to be
\[
\prod_{e\text{ is an arc of }D}w\left(  e\right)  =\prod_{\left(  i,j\right)
\text{ is an arc of }D}\left(  \left(  -1\right)  ^{\left[  i>j\right]  }%
a_{j}\right)  .
\]

\end{itemize}

\textbf{(c)} \fbox{2} Prove that $\prod_{1\leq j<i\leq n}\left(  a_{i}%
-a_{j}\right)  =\sum_{D\in T}w\left(  D\right)  $. \medskip

\textbf{(d)} \fbox{2} Prove that $\det\left(  \left(  a_{j}^{i-1}\right)
_{1\leq i\leq n,\ 1\leq j\leq n}\right)  =\sum_{\substack{D\in T\text{
is}\\\text{injective}}}w\left(  D\right)  $. \medskip

\textbf{(e)} \fbox{3} Prove that $\sum_{\substack{D\in T\text{ is
not}\\\text{injective}}}w\left(  D\right)  =0$. \medskip

\textbf{(f)} \fbox{1} Conclude that Theorem \ref{thm.det.vander} \textbf{(d)}
holds. \medskip

[\textbf{Hint:} In part \textbf{(e)}, use a sign-reversing involution. Namely,
if $D\in T$ is not injective, then its scoreboard $\operatorname*{scb}%
D=\left(  s_{1},s_{2},\ldots,s_{n}\right)  $ has two equal entries -- i.e.,
there exists a pair $\left(  u,v\right)  $ of two integers $u,v\in\left[
n\right]  $ such that $u<v$ and $s_{u}=s_{v}$. Pick such a pair $\left(
u,v\right)  $ with smallest possible $v$ (the $u$ is then uniquely determined
(why?)), and relabel the vertices $u$ and $v$ of $D$ as $v$ and $u$ (so that
any arcs of the forms $\left(  u,k\right)  $, $\left(  v,k\right)  $, $\left(
k,u\right)  $ or $\left(  k,v\right)  $ become $\left(  v,k\right)  $,
$\left(  u,k\right)  $, $\left(  k,v\right)  $ or $\left(  k,u\right)  $,
respectively). Argue that the new tournament $D^{\prime}$ is still not
injective and satisfies $w\left(  D^{\prime}\right)  =-w\left(  D\right)  $.]
\end{exercise}

Another proof of the Vandermonde determinant (Theorem \ref{thm.det.vander}
\textbf{(c)} to be specific) proceeds through a generalization:

\begin{exercise}
\label{exe.det.vander.hyperfac-pf}Let $n\in\mathbb{N}$. Let $a_{1}%
,a_{2},\ldots,a_{n}\in K$. Let $p_{1},p_{2},\ldots,p_{n}$ be $n$ polynomials
in $K\left[  x\right]  $ with the property that%
\[
\deg p_{j}\leq j-1\ \ \ \ \ \ \ \ \ \ \text{for each }j\in\left[  n\right]  .
\]
(In particular, $p_{1}$ is constant.) \medskip

\textbf{(a)} \fbox{3} Prove that
\[
\det\left(  \left(  p_{j}\left(  a_{i}\right)  \right)  _{1\leq i\leq
n,\ 1\leq j\leq n}\right)  =\left(  \prod_{j=1}^{n}\left[  x^{j-1}\right]
p_{j}\right)  \cdot\det\left(  \left(  a_{i}^{j-1}\right)  _{1\leq i\leq
n,\ 1\leq j\leq n}\right)  .
\]
(Here, $p_{j}\left(  a_{i}\right)  $ means the evaluation $p_{j}\left[
a_{i}\right]  $, of course.) \medskip

\textbf{(b)} \fbox{2} By applying this to the polynomials $p_{j}:=\left(
x-a_{1}\right)  \left(  x-a_{2}\right)  \cdots\left(  x-a_{j-1}\right)  $,
obtain a new proof of Theorem \ref{thm.det.vander} \textbf{(c)}. \medskip

\textbf{(c)} \fbox{1} Conclude that%
\[
\det\left(  \left(  p_{j}\left(  a_{i}\right)  \right)  _{1\leq i\leq
n,\ 1\leq j\leq n}\right)  =\left(  \prod_{j=1}^{n}\left[  x^{j-1}\right]
p_{j}\right)  \cdot\prod_{1\leq j<i\leq n}\left(  a_{i}-a_{j}\right)  .
\]

[\textbf{Hint:} Part \textbf{(a)} can be done in many ways, but the simplest
is probably by factoring the matrix $\left(  p_{j}\left(  a_{i}\right)
\right)  _{1\leq i\leq n,\ 1\leq j\leq n}$ as a product.]
\end{exercise}

And here is yet another generalization of the Vandermonde determinant:

\begin{exercise}
\label{exe.det.vander.kratt-lem6}Let $n\in\mathbb{N}$. Let $a_{1},a_{2}%
,\ldots,a_{n}\in K$ and $b_{1},b_{2},\ldots,b_{n}\in K$. Define $n$
polynomials $q_{1},q_{2},\ldots,q_{n}\in K\left[  x\right]  $ by setting%
\[
q_{j}=\left(  x-b_{j+1}\right)  \left(  x-b_{j+2}\right)  \cdots\left(
x-b_{n}\right)  \ \ \ \ \ \ \ \ \ \ \text{for each }j\in\left[  n\right]  .
\]
(In particular, $q_{n}=\left(  \text{empty product}\right)  =1$.) Furthermore,
let $p_{1},p_{2},\ldots,p_{n}$ be $n$ polynomials in $K\left[  x\right]  $
with the property that%
\[
\deg p_{j}\leq j-1\ \ \ \ \ \ \ \ \ \ \text{for each }j\in\left[  n\right]  .
\]
(In particular, $p_{1}$ is constant.) \medskip

\textbf{(a)} \fbox{6} Prove that
\[
\det\left(  \left(  p_{j}\left(  a_{i}\right)  q_{j}\left(  a_{i}\right)
\right)  _{1\leq i\leq n,\ 1\leq j\leq n}\right)  =\left(  \prod_{j=1}%
^{n}p_{j}\left(  b_{j}\right)  \right)  \cdot\prod_{1\leq i<j\leq n}\left(
a_{i}-a_{j}\right)  .
\]
(Again, $f\left(  a_{i}\right)  $ means the evaluation $f\left[  a_{i}\right]
$.) \medskip

\textbf{(b)} \fbox{1} Use this to obtain a new proof of Theorem
\ref{thm.det.vander} \textbf{(a)}. \medskip

\textbf{(c)} \fbox{3} Use this to prove Theorem \ref{thm.det.cauchy}.
\end{exercise}

As applications of Exercise \ref{exe.det.vander.hyperfac-pf} \textbf{(c)},
several determinants consisting of binomial coefficients can be computed:

\begin{exercise}
\label{exe.det.binom-by-vander}Let $n\in\mathbb{N}$. Let $a_{1},a_{2}%
,\ldots,a_{n}\in\mathbb{C}$. Let $H\left(  n\right)  $ denote the positive
integer $\left(  n-1\right)  !\cdot\left(  n-2\right)  !\cdot\cdots
\cdot1!=\prod_{k=1}^{n-1}k!$; this is known as the \emph{hyperfactorial} of
$n$. (In particular, $H\left(  0\right)  =H\left(  1\right)  =H\left(
2\right)  =1$ and $H\left(  3\right)  =2$ and $H\left(  4\right)  =12$ and
$H\left(  5\right)  =288$.) \medskip

\textbf{(a)} \fbox{1} Prove that%
\[
\det\left(  \left(  \dbinom{a_{i}}{j-1}\right)  _{1\leq i\leq n,\ 1\leq j\leq
n}\right)  =\dfrac{\prod_{1\leq j<i\leq n}\left(  a_{i}-a_{j}\right)
}{H\left(  n\right)  }.
\]

\textbf{(b)} \fbox{1} Conclude that $H\left(  n\right)  \mid\prod_{1\leq
j<i\leq n}\left(  a_{i}-a_{j}\right)  $ for any $n$ integers $a_{1}%
,a_{2},\ldots,a_{n}$. \medskip

\textbf{(c)} \fbox{1} Prove that $H\left(  n\right)  =\prod_{1\leq j<i\leq
n}\left(  i-j\right)  $. \medskip

\textbf{(d)} \fbox{3} Prove that%
\[
\det\left(  \left(  \dbinom{a_{i}+j}{j-1}\right)  _{1\leq i\leq n,\ 1\leq
j\leq n}\right)  =\dfrac{\prod_{1\leq j<i\leq n}\left(  a_{i}-a_{j}\right)
}{H\left(  n\right)  }.
\]

\textbf{(e)} \fbox{2} Assume that $a_{1},a_{2},\ldots,a_{n}\in\mathbb{N}$.
Prove that%
\[
\det\left(  \left(  \left(  a_{i}+j\right)  !\right)  _{1\leq i\leq n,\ 1\leq
j\leq n}\right)  =\left(  \prod_{i=1}^{n}\left(  a_{i}+1\right)  !\right)
\left(  \prod_{1\leq j<i\leq n}\left(  a_{i}-a_{j}\right)  \right)  .
\]

\textbf{(f)} \fbox{3} Assume that $a_{1},a_{2},\ldots,a_{n}\in\mathbb{N}$.
Prove that%
\[
\det\left(  \left(  \dfrac{1}{\left(  a_{i}+j\right)  !}\right)  _{1\leq i\leq
n,\ 1\leq j\leq n}\right)  =\dfrac{\prod_{1\leq i<j\leq n}\left(  a_{i}%
-a_{j}\right)  }{\prod_{i=1}^{n}\left(  a_{i}+n\right)  !}.
\]

\textbf{(g)} \fbox{2} Prove that%
\[
\det\left(  \left(  \dbinom{a_{i}+j}{i-1}\right)  _{1\leq i\leq n,\ 1\leq
j\leq n}\right)  =1.
\]

\end{exercise}

Combining parts \textbf{(b)} and \textbf{(c)} of Exercise
\ref{exe.det.binom-by-vander}, we see that any $n$ integers $a_{1}%
,a_{2},\ldots,a_{n}$ satisfy%
\[
\prod_{1\leq j<i\leq n}\left(  i-j\right)  \mid\prod_{1\leq j<i\leq n}\left(
a_{i}-a_{j}\right)  .
\]
In other words, the product of the pairwise differences between $n$ given
integers $a_{1},a_{2},\ldots,a_{n}$ is always divisible by the product of the
pairwise differences between $1,2,\ldots,n$. This curious fact is one of the
beginnings of Bhargava's theory of \textit{generalized factorials}
(\cite{Bharga00}). We will not elaborate further on this theory here, but let
us point out that our curious fact has a similarly curious analogue, in which
differences are replaced by differences of squares, and the numbers
$1,2,\ldots,n$ are replaced by $0,1,\ldots,n-1$. This analogue, too, can be
proved using determinants:

\begin{exercise}
For each $n\in\mathbb{N}$, we let $H^{\prime}\left(  n\right)  $ denote the
number
\[
\prod_{k=2}^{n-1}\dfrac{\left(  2k\right)  !}{2}=\dfrac{2!}{2}\cdot\dfrac
{4!}{2}\cdot\dfrac{6!}{2}\cdot\cdots\cdot\dfrac{\left(  2\left(  n-1\right)
\right)  !}{2}.
\]

For each $a\in\mathbb{R}$ and $k\in\mathbb{N}$, we define a real number%
\[
\dbinom{a}{k}^{\prime}:=\dfrac{\left(  a\left(  a+1\right)  \left(
a+2\right)  \cdots\left(  a+k-1\right)  \right)  \cdot\left(  a\left(
a-1\right)  \left(  a-2\right)  \cdots\left(  a-k+1\right)  \right)  }{\left(
2k\right)  !/2}.
\]

Prove the following: \medskip

\textbf{(a)} \fbox{2} We have $H^{\prime}\left(  n\right)  =\prod_{0\leq
j<i\leq n-1}\left(  i^{2}-j^{2}\right)  $ for each $n\in\mathbb{N}$. \medskip

\textbf{(b)} \fbox{1} We have $\dbinom{a}{k}^{\prime}=\dbinom{a+k}{2k}%
+\dbinom{a+k-1}{2k}\in\mathbb{Z}$ for each $a\in\mathbb{Z}$ and $k\in
\mathbb{N}$. \medskip

\textbf{(c)} \fbox{1} We have $\dbinom{a}{k}^{\prime}=\dfrac{2}{\left(
2k\right)  !}\cdot\prod_{r=0}^{k-1}\left(  a^{2}-r^{2}\right)  $ for each
$a\in\mathbb{R}$ and $k\in\mathbb{N}$. \medskip

\textbf{(d)} \fbox{3} For any $n$ integers $a_{1},a_{2},\ldots,a_{n}$ with
$n>0$, we have%
\[
\det\left(  \left(  \dbinom{a_{i}}{j-1}^{\prime}\right)  _{1\leq i\leq
n,\ 1\leq j\leq n}\right)  =2\cdot\dfrac{\prod_{1\leq j<i\leq n}\left(
a_{i}^{2}-a_{j}^{2}\right)  }{H^{\prime}\left(  n\right)  }.
\]

\textbf{(e)} \fbox{1} For any $n$ integers $a_{1},a_{2},\ldots,a_{n}$, we have
$H^{\prime}\left(  n\right)  \mid\prod_{1\leq j<i\leq n}\left(  a_{i}%
^{2}-a_{j}^{2}\right)  $. \medskip

\textbf{(f)} \fbox{1} Conclude that $\prod_{0\leq j<i\leq n-1}\left(
i^{2}-j^{2}\right)  \mid\prod_{1\leq j<i\leq n}\left(  a_{i}^{2}-a_{j}%
^{2}\right)  $ for any $n$ integers $a_{1},a_{2},\ldots,a_{n}$. \medskip

\textbf{(g)} \fbox{1} Is it true that $\prod_{0\leq j<i\leq n-1}\left(
i^{3}-j^{3}\right)  \mid\prod_{1\leq j<i\leq n}\left(  a_{i}^{3}-a_{j}%
^{3}\right)  $ for any $n$ integers $a_{1},a_{2},\ldots,a_{n}$ ? \medskip

[\textbf{Hint:} In part \textbf{(d)}, use Exercise
\ref{exe.det.vander.hyperfac-pf} \textbf{(c)} again, but this time apply it to
$a_{1}^{2},a_{2}^{2},\ldots,a_{n}^{2}$ instead of $a_{1},a_{2},\ldots,a_{n}$.]
\end{exercise}

The following exercise gives some variations on Proposition
\ref{prop.det.(xi+yj)n-1}:

\begin{exercise}
\label{exe.det.(xi+yj)m}Let $n$ be a positive integer. Let $x_{1},x_{2}%
,\ldots,x_{n}$ be $n$ elements of $K$. Let $y_{1},y_{2},\ldots,y_{n}$ be $n$
elements of $K$. \medskip

\textbf{(a)} \fbox{2} For every $m\in\left\{  0,1,\ldots,n-2\right\}  $, prove
that%
\[
\det\left(  \left(  \left(  x_{i}+y_{j}\right)  ^{m}\right)  _{1\leq i\leq
n,\ 1\leq j\leq n}\right)  =0.
\]

\textbf{(b)} \fbox{3} Let $p_{0},p_{1},\ldots,p_{n-1}$ be $n$ elements of $K$.
Let $P\in K\left[  x\right]  $ be the polynomial $\sum_{k=0}^{n-1}p_{k}x^{k}$.
Prove that%
\begin{align*}
&  \det\left(  \left(  P\left(  x_{i}+y_{j}\right)  \right)  _{1\leq i\leq
n,\ 1\leq j\leq n}\right) \\
&  =p_{n-1}^{n}\left(  \prod_{k=0}^{n-1}\dbinom{n-1}{k}\right)  \left(
\prod_{1\leq i<j\leq n}\left(  x_{i}-x_{j}\right)  \right)  \left(
\prod_{1\leq i<j\leq n}\left(  y_{j}-y_{i}\right)  \right)  .
\end{align*}

\textbf{(c)} \fbox{2} Let $p_{0},p_{1},\ldots,p_{n-1}$ be $n$ elements of $K$.
Let $P\in K\left[  x\right]  $ be the polynomial $\sum_{k=0}^{n-1}p_{k}x^{k}$.
Prove that%
\[
\det\left(  \left(  P\left(  x_{i}y_{j}\right)  \right)  _{1\leq i\leq
n,\ 1\leq j\leq n}\right)  =\left(  \prod_{k=0}^{n-1}p_{k}\right)  \left(
\prod_{1\leq i<j\leq n}\left(  x_{i}-x_{j}\right)  \right)  \left(
\prod_{1\leq i<j\leq n}\left(  y_{i}-y_{j}\right)  \right)  .
\]

\textbf{(d)} \fbox{3} Prove that%
\begin{align*}
&  \det\left(  \left(  \left(  x_{i}+y_{j}\right)  ^{n}\right)  _{1\leq i\leq
n,\ 1\leq j\leq n}\right) \\
&  =\left(  \prod_{1\leq i<j\leq n}\left(  x_{i}-x_{j}\right)  \right)
\left(  \prod_{1\leq i<j\leq n}\left(  y_{j}-y_{i}\right)  \right) \\
&  \ \ \ \ \ \ \ \ \ \ \cdot\sum_{r=0}^{n}\left(  \prod_{k\in\left\{
0,1,\ldots,n\right\}  \setminus\left\{  r\right\}  }\dbinom{n}{k}\right)
\cdot e_{r}\left[  x_{1},x_{2},\ldots,x_{n}\right]  \cdot e_{n-r}\left[
y_{1},y_{2},\ldots,y_{n}\right]  ,
\end{align*}
where we are using the following notation (a particular case of Definition
\ref{def.sf.ehp} \textbf{(a)}): For any $r\in\mathbb{N}$ and any $n$ scalars
$z_{1},z_{2},\ldots,z_{n}$, we set%
\[
e_{r}\left[  z_{1},z_{2},\ldots,z_{n}\right]  :=\sum_{\substack{\left(
i_{1},i_{2},\ldots,i_{r}\right)  \in\left[  n\right]  ^{r};\\i_{1}%
<i_{2}<\cdots<i_{r}}}z_{i_{1}}z_{i_{2}}\cdots z_{i_{r}}%
\]
(so that $e_{0}\left[  z_{1},z_{2},\ldots,z_{n}\right]  =1$ and $e_{1}\left[
z_{1},z_{2},\ldots,z_{n}\right]  =z_{1}+z_{2}+\cdots+z_{n}$ and $e_{2}\left[
z_{1},z_{2},\ldots,z_{n}\right]  =\sum_{1\leq i<j\leq n}z_{i}z_{j}$ and so on).
\end{exercise}

\begin{exercise}
\label{exe.det.hankel-power-sums}\fbox{4} Let $n\in\mathbb{N}$. Let
$u_{1},u_{2},\ldots,u_{n}\in K$ and $a_{1},a_{2},\ldots,a_{n}\in K$. For each
$k\in\mathbb{N}$, we set%
\[
z_{k}:=u_{1}a_{1}^{k}+u_{2}a_{2}^{k}+\cdots+u_{n}a_{n}^{k}.
\]
Prove that%
\begin{align*}
\det\left(  \left(  z_{i+j-2}\right)  _{1\leq i\leq n,\ 1\leq j\leq n}\right)
&  =\det\left(
\begin{array}
[c]{cccc}%
z_{0} & z_{1} & \cdots & z_{n-1}\\
z_{1} & z_{2} & \cdots & z_{n}\\
\vdots & \vdots & \ddots & \vdots\\
z_{n-1} & z_{n} & \cdots & z_{2n-2}%
\end{array}
\right) \\
&  =u_{1}u_{2}\cdots u_{n}\cdot\prod_{1\leq i<j\leq n}\left(  a_{i}%
-a_{j}\right)  ^{2}.
\end{align*}

\end{exercise}

\begin{exercise}
\label{exe.det.hessenberg.normalized}\fbox{5} Let $n>0$. Let $A\in K^{n\times
n}$ be a matrix with the property that%
\[
A_{i,i+1}=1\ \ \ \ \ \ \ \ \ \ \text{for every }i\in\left[  n-1\right]
\]
and%
\[
A_{i,j}=0\ \ \ \ \ \ \ \ \ \ \text{for every }i,j\in\left[  n\right]  \text{
satisfying }j>i+1.
\]
(Such a matrix $A$ is called a \emph{normalized lower Hessenberg matrix}. For
example, for $n=4$, such a matrix has the form $\left(
\begin{array}
[c]{cccc}%
\ast & 1 & 0 & 0\\
\ast & \ast & 1 & 0\\
\ast & \ast & \ast & 1\\
\ast & \ast & \ast & \ast
\end{array}
\right)  $, where each asterisk $\ast$ stands for an arbitrary entry.)

For each subset $I$ of $\left[  n-1\right]  $, we define an element
$p_{I}\left(  A\right)  \in K$ as follows: Write the subset $I$ in the form
$\left\{  i_{1},i_{2},\ldots,i_{k}\right\}  $ with $i_{1}<i_{2}<\cdots<i_{k}$.
Additionally, set $i_{0}:=0$ and $i_{k+1}:=n$. Then, set%
\[
p_{I}\left(  A\right)  :=A_{i_{1},i_{0}+1}A_{i_{2},i_{1}+1}\cdots
A_{i_{k+1},i_{k}+1}=\prod_{u=1}^{k+1}A_{i_{u},i_{u-1}+1}.
\]

Prove that%
\[
\det A=\sum_{I\subseteq\left[  n-1\right]  }\left(  -1\right)
^{n-1-\left\vert I\right\vert }p_{I}\left(  A\right)  .
\]

\end{exercise}

Normalized lower Hessenberg matrices can be used to provide an
\textquotedblleft explicit\textquotedblright\ determinantal formula for the
coefficients of the inverse of an FPS:

\begin{exercise}
\fbox{4} Let $f=\sum_{k\in\mathbb{N}}f_{k}x^{k}$ be a FPS in $K\left[  \left[
x\right]  \right]  $ (with $f_{0},f_{1},f_{2},\ldots\in K$). Assume that
$f_{0}=1$. Set $f_{k}:=0$ for all negative $k\in\mathbb{Z}$. Let $g=f^{-1}$ be
the multiplicative inverse of $f$ in $K\left[  \left[  x\right]  \right]  $,
and let $g_{0},g_{1},g_{2},\ldots$ be the coefficients of $g$ (so that
$g=\sum_{k\in\mathbb{N}}g_{k}x^{k}$). Prove that each $n\in\mathbb{N}$
satisfies%
\begin{align*}
g_{n}  &  =\left(  -1\right)  ^{n}\det\left(  \left(  f_{i-j+1}\right)
_{1\leq i\leq n,\ 1\leq j\leq n}\right) \\
&  =\left(  -1\right)  ^{n}\det\left(
\begin{array}
[c]{cccccc}%
f_{1} & 1 & 0 & 0 & \cdots & 0\\
f_{2} & f_{1} & 1 & 0 & \cdots & 0\\
f_{3} & f_{2} & f_{1} & 1 & \cdots & 0\\
f_{4} & f_{3} & f_{2} & f_{1} & \cdots & 0\\
\vdots & \vdots & \vdots & \vdots & \ddots & \vdots\\
f_{n} & f_{n-1} & f_{n-2} & f_{n-3} & \cdots & f_{1}%
\end{array}
\right)  .
\end{align*}

\end{exercise}

Our next exercise is concerned with \emph{tridiagonal matrices}. These are
matrices whose all entries are zero except for those on the diagonal and
\textquotedblleft its neighbors\textquotedblright. For example, a $3\times
3$-matrix is tridiagonal if it has the form $\left(
\begin{array}
[c]{ccc}%
\ast & \ast & 0\\
\ast & \ast & \ast\\
0 & \ast & \ast
\end{array}
\right)  $ (where each asterisk $\ast$ means an arbitrary entry). \Needspace{9pc}

\begin{exercise}
\label{exe.det.tridiag}Let $n\in\mathbb{N}$. Let $a_{1},a_{2},\ldots,a_{n}\in
K$ and $b_{1},b_{2},\ldots,b_{n-1}\in K$ and $c_{1},c_{2},\ldots,c_{n-1}\in
K$. Set%
\[
A=\left(
\begin{array}
[c]{ccccccc}%
a_{1} & b_{1} & 0 & \cdots & 0 & 0 & 0\\
c_{1} & a_{2} & b_{2} & \cdots & 0 & 0 & 0\\
0 & c_{2} & a_{3} & \cdots & 0 & 0 & 0\\
\vdots & \vdots & \vdots & \ddots & \vdots & \vdots & \vdots\\
0 & 0 & 0 & \cdots & a_{n-2} & b_{n-2} & 0\\
0 & 0 & 0 & \cdots & c_{n-2} & a_{n-1} & b_{n-1}\\
0 & 0 & 0 & \cdots & 0 & c_{n-1} & a_{n}%
\end{array}
\right)  \in K^{n\times n}.
\]
Formally speaking, this matrix $A$ is defined to be%
\[
A=\left(
\begin{cases}
a_{i}, & \text{if }i=j;\\
b_{i}, & \text{if }i=j-1;\\
c_{j}, & \text{if }i=j+1;\\
0, & \text{otherwise}%
\end{cases}
\right)  _{1\leq i\leq n,\ 1\leq j\leq n}.
\]
The matrix $A$ is called a \emph{tridiagonal matrix}, and its determinant is
known as the \emph{continuant} of the numbers $a_{i},b_{i},c_{i}$. \medskip

\textbf{(a)} \fbox{2} Prove that every $m\in\left\{  2,3,\ldots,n\right\}  $
satisfies%
\[
\det\left(  A_{:m}\right)  =a_{m}\det\left(  A_{:m-1}\right)  -b_{m-1}%
c_{m-1}\det\left(  A_{:m-2}\right)  ,
\]
where we set $A_{:k}:=\operatorname*{sub}\nolimits_{1,2,\ldots,k}%
^{1,2,\ldots,k}A=\left(  A_{i,j}\right)  _{1\leq i\leq k,\ 1\leq j\leq k}$ for
each $k\in\left\{  0,1,\ldots,n\right\}  $. \medskip

\textbf{(b)} \fbox{3} Recall the notion of a \textquotedblleft lacunar
set\textquotedblright\ as defined in Definition \ref{def.lacunar.lac}. If $I$
is a set of integers, then we let $I^{+}:=\left\{  i+1\ \mid\ i\in I\right\}
$. Prove that%
\[
\det A=\sum_{\substack{I\subseteq\left[  n-1\right]  \text{ is}%
\\\text{lacunar}}}\left(  \prod_{i\in\left[  n\right]  \setminus\left(  I\cup
I^{+}\right)  }a_{i}\right)  \left(  \prod_{i\in I}\left(  -b_{i}c_{i}\right)
\right)  .
\]

\textbf{(c)} \fbox{1} Compute $\det A$ in the case when
\begin{align*}
a_{i}  &  =1\ \ \ \ \ \ \ \ \ \ \left(  \text{for all }i\in\left[  n\right]
\right)  \ \ \ \ \ \ \ \ \ \ \text{and}\\
b_{i}  &  =1\ \ \ \ \ \ \ \ \ \ \left(  \text{for all }i\in\left[  n-1\right]
\right)  \ \ \ \ \ \ \ \ \ \ \text{and}\\
c_{i}  &  =-1\ \ \ \ \ \ \ \ \ \ \left(  \text{for all }i\in\left[
n-1\right]  \right)  .
\end{align*}

\textbf{(d)} \fbox{1} Compute $\det A$ in the case when
\begin{align*}
a_{i}  &  =2\ \ \ \ \ \ \ \ \ \ \left(  \text{for all }i\in\left[  n\right]
\right)  \ \ \ \ \ \ \ \ \ \ \text{and}\\
b_{i}  &  =c_{i}=1\ \ \ \ \ \ \ \ \ \ \left(  \text{for all }i\in\left[
n-1\right]  \right)  .
\end{align*}

\textbf{(e)} \fbox{3} If $d_{0},d_{1},\ldots,d_{n}\in K$ are arbitrary
elements, and if we have%
\begin{align*}
a_{i}  &  =d_{i-1}+d_{i}\ \ \ \ \ \ \ \ \ \ \text{for all }i\in\left[
n\right]  ,\ \ \ \ \ \ \ \ \ \ \text{and}\\
b_{i}  &  =c_{i}=-d_{i}\ \ \ \ \ \ \ \ \ \ \text{for all }i\in\left[
n-1\right]  \text{,}%
\end{align*}
then prove that
\[
\det A=\sum_{k=0}^{n}d_{0}d_{1}\cdots\widehat{d_{k}}\cdots d_{n},
\]
where the hat over the \textquotedblleft$d_{k}$\textquotedblright\ is defined
as in Proposition \ref{prop.det.x+ai}. (The matrix $A$ in this case has a
physical interpretation as the stiffness matrix of a mass-spring chain; see
\cite[\S 6.1]{OlvSha}. Note that we could replace the condition $b_{i}%
=c_{i}=-d_{i}$ by $b_{i}=c_{i}=d_{i}$ without changing $\det A$, but the minus
signs come from the physical backstory.) \medskip

\textbf{(f)} \fbox{2} We return to the general case. Define $A_{:k}$ as in
part \textbf{(a)}. Assume that $n\geq1$. Prove that%
\[
\dfrac{\det A}{\det\left(  A_{:n-1}\right)  }=a_{n}-\dfrac{b_{n-1}c_{n-1}%
}{a_{n-1}-\dfrac{b_{n-2}c_{n-2}}{a_{n-2}-\dfrac{b_{n-3}c_{n-3}}{%
\begin{array}
[c]{ccc}%
a_{n-3}- &  & \\
& \ddots & \\
&  & -\dfrac{b_{2}c_{2}}{a_{2}-\dfrac{b_{1}c_{1}}{a_{1}}}%
\end{array}
}}},
\]
provided that all denominators in this equality are invertible. \medskip

[\textbf{Remark:} If we set $a_{i}=2x$ for all $i\in\left[  n\right]  $ and
$b_{i}=c_{i}=-1$ for all $i\in\left[  n-1\right]  $, then $\det A$ becomes the
$n$-th \emph{Chebyshev polynomial of the second kind}, commonly denoted
$U_{n}\left(  x\right)  $. Some properties of Chebyshev polynomials can be
generalized to determinants of arbitrary tridiagonal matrices.

Part \textbf{(f)} can be seen as a formula for expressing a continued fraction
as a ratio of determinants. The connections between continued fractions and
determinants run much deeper.]
\end{exercise}

\begin{exercise}
Let $a\in K$ be arbitrary. Let $n\in\mathbb{N}$. Define a tridiagonal $n\times
n$-matrix $A\in K^{n\times n}$ as in Exercise \ref{exe.det.tridiag}, after
setting%
\begin{align*}
a_{i}  &  =a\ \ \ \ \ \ \ \ \ \ \left(  \text{for all }i\in\left[  n\right]
\right)  \ \ \ \ \ \ \ \ \ \ \text{and}\\
b_{i}  &  =n-i\ \ \ \ \ \ \ \ \ \ \left(  \text{for all }i\in\left[
n-1\right]  \right)  \ \ \ \ \ \ \ \ \ \ \text{and}\\
c_{i}  &  =i\ \ \ \ \ \ \ \ \ \ \left(  \text{for all }i\in\left[  n-1\right]
\right)  .
\end{align*}
Thus,%
\[
A=\left(
\begin{array}
[c]{ccccccc}%
a & n-1 & 0 & \cdots & 0 & 0 & 0\\
1 & a & n-2 & \cdots & 0 & 0 & 0\\
0 & 2 & a & \cdots & 0 & 0 & 0\\
\vdots & \vdots & \vdots & \ddots & \vdots & \vdots & \vdots\\
0 & 0 & 0 & \cdots & a & 2 & 0\\
0 & 0 & 0 & \cdots & n-2 & a & 1\\
0 & 0 & 0 & \cdots & 0 & n-1 & a
\end{array}
\right)  .
\]

\textbf{(a)} \fbox{6} Prove that%
\[
\det A=\prod_{k=0}^{n-1}\left(  a-2k+n-1\right)  .
\]
(This particular matrix $A$, for $a=0$, is called the \emph{Kac matrix}; thus,
our formula for $\det A$ computes its characteristic polynomial.) \medskip

\textbf{(b)} \fbox{?} Can you find a bijective proof using the formula in
Exercise \ref{exe.det.tridiag} \textbf{(b)}?
\end{exercise}

Another variation on the Vandermonde determinant:

\begin{exercise}
\fbox{3} Let $n\in\mathbb{N}$. Let $a_{1},a_{2},\ldots,a_{n}\in K$ and
$b_{1},b_{2},\ldots,b_{n}\in K$. Prove that%
\[
\det\left(  \left(  a_{i}^{n-j}b_{i}^{j-1}\right)  _{1\leq i\leq n,\ 1\leq
j\leq n}\right)  =\prod_{1\leq i<j\leq n}\left(  a_{i}b_{j}-a_{j}b_{i}\right)
.
\]

\end{exercise}

\subsubsection{The Lindstr\"{o}m--Gessel--Viennot lemma}

The following exercises are concerned with Section \ref{sec.det.comb.lgv}.

\begin{exercise}
\fbox{5} Consider the following variants of the definition of the map $f$ in
the proof of Proposition \ref{prop.lgv.kpaths.count}: \medskip

\textbf{(a)} Define $v$ to be the last (rather than the first) crowded point
on $p_{i}$. \medskip

\textbf{(b)} Define $j$ to be the smallest (rather than the largest) element
of $\left[  k\right]  $ such that $v$ belongs to $p_{j}$ (not counting $i$, of
course). \medskip

\textbf{(c)} Instead of choosing $v$ first and $j$ later, choose $j$ and $v$
as follows: First, pick $j$ to be the largest element of $\left[  k\right]  $
such that the path $p_{i}$ intersects $p_{j}$; then, define $v$ to be the
first point where $p_{i}$ intersects $p_{j}$. \medskip

\textbf{(d)} Instead of defining $v$ to be the first crowded point on $p_{i}$,
choose some arbitrary total order on the vertex set of the digraph, and define
$v$ to be the largest crowded point on $p_{i}$ with respect to this order.
(The total order needs to be chosen in advance; it must not depend on
$\mathbf{p}$ or $\sigma$.) \medskip

Which of these variants \textquotedblleft work\textquotedblright\ (i.e., lead
to well-defined sign-reversing involutions $f:\mathcal{X}\rightarrow
\mathcal{X}$) ?

(You are not required to try out all combinations of these variants; just
analyze each variant for itself.)
\end{exercise}

\begin{exercise}
\label{exe.lgv.catalan-hankel-det-0}Complete the proof of Corollary
\ref{cor.lgv.catalan-hankel-det-0} sketched above: \medskip

\textbf{(a)} \fbox{2} Show that there is only one nipat from $\mathbf{A}$ to
$\mathbf{B}$. \medskip

\textbf{(b)} \fbox{3} Show that there are no nipats from $\mathbf{A}$ to
$\sigma\left(  \mathbf{B}\right)  $ when $\sigma\in S_{k}$ is not the identity
permutation $\operatorname*{id}\in S_{k}$. \medskip

Furthermore: \medskip

\textbf{(c)} \fbox{3} Prove the analogue of Corollary
\ref{cor.lgv.catalan-hankel-det-0} that says that each $k\in\mathbb{N}$
satisfies%
\[
\det\left(  \left(  c_{i+j-1}\right)  _{1\leq i\leq k,\ 1\leq j\leq k}\right)
=\det\left(
\begin{array}
[c]{cccc}%
c_{1} & c_{2} & \cdots & c_{k}\\
c_{2} & c_{3} & \cdots & c_{k+1}\\
\vdots & \vdots & \ddots & \vdots\\
c_{k} & c_{k+1} & \cdots & c_{2k-1}%
\end{array}
\right)  =1.
\]

\textbf{(d)} \fbox{2} Show that the Catalan sequence $\left(  c_{0}%
,c_{1},c_{2},\ldots\right)  $ is the only sequence $\left(  a_{0},a_{1}%
,a_{2},\ldots\right)  $ of real numbers that satisfies%
\begin{align*}
\det\left(  \left(  a_{i+j-2}\right)  _{1\leq i\leq k,\ 1\leq j\leq k}\right)
&  =\det\left(  \left(  a_{i+j-1}\right)  _{1\leq i\leq k,\ 1\leq j\leq
k}\right)  =1\\
&  \ \ \ \ \ \ \ \ \ \ \ \ \ \ \ \ \ \ \ \ \text{for all }k\in\mathbb{N}.
\end{align*}

\textbf{(e)} \fbox{3} Compute $\det\left(  \left(  c_{i+j}\right)  _{1\leq
i\leq k,\ 1\leq j\leq k}\right)  $ for each $k\in\mathbb{N}$.
\end{exercise}

The next exercise should be contrasted with Exercise
\ref{exe.perm.descents.subI}.

\begin{exercise}
\label{exe.lgv.descents.isI}Let $n\in\mathbb{N}$. Let $I$ be a subset of
$\left[  n-1\right]  $. Write $I$ in the form $I=\left\{  c_{1},c_{2}%
,\ldots,c_{k}\right\}  $ with $c_{1}<c_{2}<\cdots<c_{k}$. Set $c_{0}:=0$ and
$c_{k+1}:=n$. Recall Definition \ref{def.perm.descents}. \medskip

\textbf{(a)} \fbox{3} Prove that%
\[
\left(  \text{\# of }\sigma\in S_{n}\text{ satisfying }\operatorname{Des}%
\sigma=I\right)  =\det\left(  \left(  \dbinom{n-c_{i-1}}{c_{j}-c_{i-1}%
}\right)  _{1\leq i\leq k+1,\ 1\leq j\leq k+1}\right)  .
\]

\textbf{(b)} \fbox{5} Let us use the notations from Definition
\ref{def.pars.qbinom.qbinom}. Prove that%
\[
\sum_{\substack{\sigma\in S_{n};\\\operatorname*{Des}\sigma=I}}q^{\ell\left(
\sigma\right)  }=\det\left(  \left(  \dbinom{n-c_{i-1}}{c_{j}-c_{i-1}}%
_{q}\right)  _{1\leq i\leq k+1,\ 1\leq j\leq k+1}\right)
\ \ \ \ \ \ \ \ \ \ \text{in the ring }\mathbb{Z}\left[  q\right]  \text{.}%
\]

[\textbf{Hint:} One way to approach this is by observing that the matrices in
question are transposes of normalized lower Hessenberg matrices as in Exercise
\ref{exe.det.hessenberg.normalized}. Another is to apply the LGV lemma (in its
weighted version for part \textbf{(b)}) to the $\left(  k+1\right)  $-vertices
$\mathbf{A}=\left(  A_{1},A_{2},\ldots,A_{k+1}\right)  $ and $\mathbf{B}%
=\left(  B_{1},B_{2},\ldots,B_{k+1}\right)  $ defined by%
\[
A_{i}:=\left(  0,c_{i-1}\right)  \ \ \ \ \ \ \ \ \ \ \text{and}%
\ \ \ \ \ \ \ \ \ \ B_{i}:=\left(  n-c_{i},c_{i}\right)  .
\]
Here is an example, for $n=10$ and $I=\left\{  5,6,9\right\}  $ and $\sigma\in
S_{10}$ with OLN $\left(  2,3,5,7,10,9,1,6,8,4\right)  $:%
\[%
\begin{tikzpicture}
\draw[densely dotted] (-0.2,-0.2) grid (10.2, 10.2);
\draw[very thin] (0,0) -- (0,10.2);
\draw[very thin] (-0.2,10.2) -- (10.2,-0.2);
\node
[circle,fill=white,draw=black,text=black,inner sep=1pt] (A1) at (0,0) {$A_1$};
\node
[circle,fill=white,draw=black,text=black,inner sep=1pt] (A2) at (0,5) {$A_2$};
\node
[circle,fill=white,draw=black,text=black,inner sep=1pt] (A3) at (0,6) {$A_3$};
\node
[circle,fill=white,draw=black,text=black,inner sep=1pt] (A4) at (0,9) {$A_4$};
\node
[circle,fill=white,draw=black,text=black,inner sep=1pt] (B1) at (5,5) {$B_1$};
\node
[circle,fill=white,draw=black,text=black,inner sep=1pt] (B2) at (4,6) {$B_2$};
\node
[circle,fill=white,draw=black,text=black,inner sep=1pt] (B3) at (1,9) {$B_3$};
\node
[circle,fill=white,draw=black,text=black,inner sep=1pt] (B4) at (0,10) {$B_4$}%
;
\begin{scope}[thick,>=stealth,darkred]
\draw(A1) edge[->] (1,0);
\draw(1,0) edge[->] (1,1);
\draw(1,0.5) node[anchor=west] {$1$};
\draw(1,1) edge[->] (1,2);
\draw(1,1.5) node[anchor=west] {$1$};
\draw(1,2) edge[->] (2,2);
\draw(2,2) edge[->] (2,3);
\draw(2,2.5) node[anchor=west] {$2$};
\draw(2,3) edge[->] (3,3);
\draw(3,3) edge[->] (3,4);
\draw(3,3.5) node[anchor=west] {$3$};
\draw(3,4) edge[->] (4,4);
\draw(4,4) edge[->] (5,4);
\draw(5,4) edge[->] (B1);
\draw(5,4.5) node[anchor=east] {$5$};
\end{scope}
\begin{scope}[thick,>=stealth,dbluecolor]
\draw(A2) edge[->] (1,5);
\draw(1,5) edge[->] (2,5);
\draw(2,5) edge[->] (3,5);
\draw(3,5) edge[->] (4,5);
\draw(4,5) edge[->] (B2);
\draw(4,5.5) node[anchor=east] {$4$};
\end{scope}
\begin{scope}[thick,>=stealth,dgreencolor]
\draw(A3) edge[->] (0,7);
\draw(0,6.5) node[anchor=east] {$0$};
\draw(0,7) edge[->] (1,7);
\draw(1,7) edge[->] (1,8);
\draw(1,7.5) node[anchor=east] {$1$};
\draw(1,8) edge[->] (B3);
\draw(1,8.5) node[anchor=east] {$1$};
\end{scope}
\begin{scope}[thick,>=stealth,red!90!yellow]
\draw(A4) edge[->] (B4);
\draw(0,9.5) node[anchor=east] {$0$};
\end{scope}
\end{tikzpicture}%
\ \ .
\]
Here, the x-coordinates of the north-steps of the paths (from bottom to top:
$1,1,2,3,5,4,0,1,1,0$) are the entries of the Lehmer code $L\left(
\sigma\right)  =\left(  1,1,2,3,5,4,0,1,1,0\right)  $ of $\sigma$.]
\end{exercise}

\begin{noncompile}
This is \cite[\S 3]{GesVie85}; see also \cite[Example 2.2.4]{Stanley-EC1}.
Another formula for part \textbf{(a)} is \cite[Exercise 4 \textbf{(b)}%
]{18f-hw4s}.
\end{noncompile}

The next exercise sketches out a visual proof of the Cauchy--Binet formula
(Theorem \ref{thm.det.CB}) using the LGV lemma:

\begin{exercise}
\label{exe.lgv.cb}Let $K$ be a commutative ring. Let $n,m\in\mathbb{N}$. Let
$A\in K^{n\times m}$ be an $n\times m$-matrix, and let $B\in K^{m\times n}$ be
an $m\times n$-matrix.

Let $D$ be the digraph with $2n+m$ vertices labeled%
\[
1,2,\ldots,n,\ \ \ 1^{\prime},2^{\prime},\ldots,m^{\prime},\ \ \ 1^{\prime
\prime},2^{\prime\prime},\ldots,n^{\prime\prime},
\]
and with arcs%
\[
i\rightarrow j^{\prime}\ \ \ \ \ \ \ \ \ \ \text{for all }i\in\left[
n\right]  \text{ and }j\in\left[  m\right]
\]
and%
\[
i^{\prime}\rightarrow j^{\prime\prime}\ \ \ \ \ \ \ \ \ \ \text{for all }%
i\in\left[  m\right]  \text{ and }j\in\left[  n\right]  .
\]
Here is how $D$ looks like for $n=2$ and $m=4$:%
\[%
\begin{tikzpicture}[scale=2]
\node
[circle,fill=white,draw=black,text=black,inner sep=3pt] (A1) at (2,2) {$1$};
\node
[circle,fill=white,draw=black,text=black,inner sep=3pt] (A2) at (3,2) {$2$};
\node
[circle,fill=white,draw=black,text=black,inner sep=2pt] (B1) at (1,1) {$1'$};
\node
[circle,fill=white,draw=black,text=black,inner sep=2pt] (B2) at (2,1) {$2'$};
\node
[circle,fill=white,draw=black,text=black,inner sep=2pt] (B3) at (3,1) {$3'$};
\node
[circle,fill=white,draw=black,text=black,inner sep=2pt] (B4) at (4,1) {$4'$};
\node
[circle,fill=white,draw=black,text=black,inner sep=1pt] (C1) at (2,0) {$1''$};
\node
[circle,fill=white,draw=black,text=black,inner sep=1pt] (C2) at (3,0) {$2''$};
\begin{scope}[thick,>=stealth,darkred]
\draw(A1) edge[->] (B1);
\draw(A1) edge[->] (B2);
\draw(A1) edge[->] (B3);
\draw(A1) edge[->] (B4);
\draw(A2) edge[->] (B1);
\draw(A2) edge[->] (B2);
\draw(A2) edge[->] (B3);
\draw(A2) edge[->] (B4);
\end{scope}
\begin{scope}[thick,>=stealth,dbluecolor]
\draw(B1) edge[->] (C1);
\draw(B1) edge[->] (C2);
\draw(B2) edge[->] (C1);
\draw(B2) edge[->] (C2);
\draw(B3) edge[->] (C1);
\draw(B3) edge[->] (C2);
\draw(B4) edge[->] (C1);
\draw(B4) edge[->] (C2);
\end{scope}
\end{tikzpicture}%
\ \ .
\]

\textbf{(a)} \fbox{1} Prove that this digraph $D$ is acyclic. \medskip

Now, for each arc $a$ of $D$, we define a weight $w\left(  a\right)  \in K$ as follows:

\begin{itemize}
\item If $a$ is the arc $i\rightarrow j^{\prime}$ for some $i\in\left[
n\right]  $ and $j\in\left[  m\right]  $, then we set $w\left(  a\right)
:=A_{i,j}$.

\item If $a$ is the arc $i^{\prime}\rightarrow j^{\prime\prime}$ for some
$i\in\left[  m\right]  $ and $j\in\left[  n\right]  $, then we set $w\left(
a\right)  :=B_{i,j}$.
\end{itemize}

\textbf{(b)} \fbox{1} Prove that%
\[
\left(  \sum_{p:i\rightarrow j^{\prime\prime}}w\left(  p\right)  \right)
_{1\leq i\leq n,\ 1\leq j\leq n}=AB.
\]
Here, \textquotedblleft$p:u\rightarrow v$\textquotedblright\ means
\textquotedblleft$p$ is a path from $u$ to $v$\textquotedblright\ whenever $u$
and $v$ are two vertices of $D$. \medskip

\textbf{(c)} \fbox{3} Define two $n$-vertices $\mathbf{A}$ and $\mathbf{B}$ by
$\mathbf{A}=\left(  1,2,\ldots,n\right)  $ and $\mathbf{B}=\left(
1^{\prime\prime},2^{\prime\prime},\ldots,n^{\prime\prime}\right)  $. Prove
that%
\begin{align*}
&  \sum_{\sigma\in S_{n}}\left(  -1\right)  ^{\sigma}\sum
_{\substack{\mathbf{p}\text{ is a nipat}\\\text{from }\mathbf{A}\text{ to
}\sigma\left(  \mathbf{B}\right)  }}w\left(  \mathbf{p}\right) \\
&  =\sum_{\substack{\left(  g_{1},g_{2},\ldots,g_{n}\right)  \in\left[
m\right]  ^{n};\\g_{1}<g_{2}<\cdots<g_{n}}}\det\left(  \operatorname*{cols}%
\nolimits_{g_{1},g_{2},\ldots,g_{n}}A\right)  \cdot\det\left(
\operatorname*{rows}\nolimits_{g_{1},g_{2},\ldots,g_{n}}B\right)  ,
\end{align*}
where we are using the notations of Theorem \ref{thm.lgv.kpaths.wt-dg} (with
$k=n$) and of Theorem \ref{thm.det.CB}. \medskip

\textbf{(d)} \fbox{1} Prove Theorem \ref{thm.det.CB}.
\end{exercise}

\begin{exercise}
Consider the situation of Theorem \ref{thm.lgv.kpaths.wt-dg}. Assume that our
digraph $D$ has vertex set $\left[  n\right]  $ for some $n\in\mathbb{N}$. Let
$E$ be the set of all arcs of $D$.

Let $M\in K^{n\times n}$ be the $n\times n$-matrix whose $\left(  i,j\right)
$-th entry is given by%
\[
M_{i,j}=\sum_{\substack{a\in E\text{ is an arc}\\\text{from }i\text{ to }%
j}}w\left(  a\right)  \ \ \ \ \ \ \ \ \ \ \text{for all }i,j\in\left[
n\right]  .
\]
(Note that if $D$ is a simple digraph, then the sum on the right hand side of
this equality has at most one addend.) \medskip

\textbf{(a)} \fbox{2} Prove that each $k\in\mathbb{N}$ satisfies%
\[
M^{k}=\left(  \sum_{\substack{p:i\rightarrow j\text{ is a path}\\\text{with
}k\text{ steps}}}w\left(  p\right)  \right)  _{1\leq i\leq n,\ 1\leq j\leq
n}.
\]
Here, \textquotedblleft$p:i\rightarrow j$\textquotedblright\ means
\textquotedblleft$p$ is a path from $i$ to $j$\textquotedblright. \medskip

\textbf{(b)} \fbox{2} Prove that $M^{n}=0_{n\times n}$ (the zero matrix).
\medskip

\textbf{(c)} \fbox{2} Let $I_{n}\in K^{n\times n}$ denote the $n\times n$
identity matrix. Prove that
\[
\left(  I_{n}-M\right)  ^{-1}=\left(  \sum_{p:i\rightarrow j}w\left(
p\right)  \right)  _{1\leq i\leq n,\ 1\leq j\leq n}.
\]

\end{exercise}

\subsection{Symmetric functions}

The notations of Chapter \ref{chap.sf} shall be used here. In particular, an
integer $N\in\mathbb{N}$ and a commutative ring $K$ are fixed.

\subsubsection{Definitions and examples of symmetric polynomials}

\begin{exercise}
\label{exe.sf.e-h-FPS.b}\fbox{2} Prove Proposition \ref{prop.sf.e-h-FPS}
\textbf{(b)}.
\end{exercise}

\begin{exercise}
\label{exe.sf.rec1}\fbox{2} Let $M\in\left\{  0,1,\ldots,N\right\}  $ and
$n\in\mathbb{N}$. Prove that%
\[
e_{n}\left[  x_{1},x_{2},\ldots,x_{N}\right]  =\sum_{i=0}^{n}e_{i}\left[
x_{1},x_{2},\ldots,x_{M}\right]  \cdot e_{n-i}\left[  x_{M+1},x_{M+2}%
,\ldots,x_{N}\right]
\]
and%
\[
h_{n}\left[  x_{1},x_{2},\ldots,x_{N}\right]  =\sum_{i=0}^{n}h_{i}\left[
x_{1},x_{2},\ldots,x_{M}\right]  \cdot h_{n-i}\left[  x_{M+1},x_{M+2}%
,\ldots,x_{N}\right]
\]
and (if $n$ is positive)%
\[
p_{n}\left[  x_{1},x_{2},\ldots,x_{N}\right]  =p_{n}\left[  x_{1},x_{2}%
,\ldots,x_{M}\right]  +p_{n}\left[  x_{M+1},x_{M+2},\ldots,x_{N}\right]  .
\]

\end{exercise}

\begin{exercise}
\label{exe.sf.NG.p}\textbf{(a)} \fbox{3} Reprove Theorem \ref{thm.sf.NG} by
induction on $N$. \medskip

\textbf{(b)} \fbox{4} Reprove Theorem \ref{thm.sf.NG} combinatorially, by
expanding the left-hand sides and showing that all monomials cancel except for
those that end up constituting the right-hand sides.
\end{exercise}

\begin{exercise}
\label{exe.sf.NG45}\fbox{5} Prove that each positive integer $n$ satisfies%
\begin{align*}
\sum_{j=1}^{n}\left(  -1\right)  ^{j-1}je_{j}h_{n-j}  &  =p_{n}%
\ \ \ \ \ \ \ \ \ \ \text{and}\\
\sum_{j=1}^{n}\left(  -1\right)  ^{n-j}jh_{j}e_{n-j}  &  =p_{n}.
\end{align*}

\end{exercise}

\begin{exercise}
\label{exe.sf.qbinom}Let $n\in\mathbb{N}$. \medskip

\textbf{(a)} \fbox{1} Prove that%
\begin{align*}
e_{n}\left[  \underbrace{1,1,\ldots,1}_{N\text{ times}}\right]   &
=\dbinom{N}{n}\ \ \ \ \ \ \ \ \ \ \text{and}\\
h_{n}\left[  \underbrace{1,1,\ldots,1}_{N\text{ times}}\right]   &
=\dbinom{N+n-1}{n}.
\end{align*}

\textbf{(b)} \fbox{2} Prove that%
\begin{align*}
e_{n}\left[  q^{0},q^{1},\ldots,q^{N-1}\right]   &  =q^{n\left(  n-1\right)
/2}\dbinom{N}{n}_{q}\ \ \ \ \ \ \ \ \ \ \text{and}\\
h_{n}\left[  q^{0},q^{1},\ldots,q^{N-1}\right]   &  =\dbinom{N+n-1}{n}_{q}%
\end{align*}
in the ring $\mathbb{Z}\left[  q\right]  $, where we are using the notation of
Definition \ref{def.pars.qbinom.qbinom} (and where $\dbinom{-1}{0}_{q}:=1$).
\medskip

\textbf{(c)} \fbox{2} Recover a nontrivial identity between $q$-binomial
coefficients by substituting $q^{0},q^{1},\ldots,q^{N-1}$ into an identity
between symmetric polynomials. (There are several valid answers here.)
\medskip

\textbf{(d)} \fbox{3} For any $m,k\in\mathbb{N}$, the \emph{unsigned Stirling
number of the 1st kind} $c\left(  m,k\right)  \in\mathbb{N}$ is defined to be
the \# of all permutations $\sigma\in S_{m}$ that have exactly $k$ cycles (see
Definition \ref{def.perm.cycs.cycs} \textbf{(a)}). Prove that%
\[
e_{n}\left[  1,2,\ldots,N\right]  =c\left(  N+1,N+1-n\right)  .
\]
Here, we agree to understand $c\left(  n,m\right)  $ as $0$ when $m$ is
negative. \medskip

\textbf{(e)} \fbox{3} For any $m,k\in\mathbb{N}$, the \emph{Stirling number of
the 2nd kind} $S\left(  m,k\right)  \in\mathbb{N}$ is defined to be the \# of
all set partitions of the set $\left[  m\right]  $ into $k$ parts (i.e., the
\# of sets $\left\{  U_{1},U_{2},\ldots,U_{k}\right\}  $ consisting of $k$
disjoint nonempty subsets $U_{1},U_{2},\ldots,U_{k}$ of $\left[  m\right]  $
such that $U_{1}\cup U_{2}\cup\cdots\cup U_{k}=\left[  m\right]  $). Prove
that%
\[
h_{n}\left[  1,2,\ldots,N\right]  =S\left(  N+n,N\right)  .
\]

\end{exercise}

\begin{exercise}
\label{exe.sf.petrie.1}For each $n\in\mathbb{N}$ and each positive integer
$k$, we define the $\left(  n,k\right)  $\emph{-th Petrie symmetric polynomial
}$g_{k,n}\in\mathcal{S}$ by%
\begin{align*}
g_{k,n}  &  =\sum_{\substack{\left(  i_{1},i_{2},\ldots,i_{n}\right)
\in\left[  N\right]  ^{n};\\i_{1}\leq i_{2}\leq\cdots\leq i_{n};\\\text{no
}k\text{ of the numbers }i_{1},i_{2},\ldots,i_{n}\text{ are equal}}}x_{i_{1}%
}x_{i_{2}}\cdots x_{i_{n}}\\
&  =\sum_{\substack{\left(  a_{1},a_{2},\ldots,a_{N}\right)  \in\left\{
0,1,\ldots,k-1\right\}  ^{N};\\a_{1}+a_{2}+\cdots+a_{N}=n}}x_{1}^{a_{1}}%
x_{2}^{a_{2}}\cdots x_{N}^{a_{N}}.
\end{align*}

\textbf{(a)} \fbox{1} Prove that $g_{2,n}=e_{n}$ for each $n\in\mathbb{N}$.
\medskip

\textbf{(b)} \fbox{1} Prove that $g_{k,n}=h_{n}$ for each $n\in\mathbb{N}$ and
each $k>n$. \medskip

\textbf{(c)} \fbox{1} Prove that $g_{n,n}=h_{n}-p_{n}$ for each $n>0$.
\medskip

\textbf{(d)} \fbox{3} Prove that each $n\in\mathbb{N}$ and each $k>0$ satisfy
\[
g_{k,n}=\sum_{i\in\mathbb{N}}\left(  -1\right)  ^{i}e_{i}\left[  x_{1}%
^{k},x_{2}^{k},\ldots,x_{N}^{k}\right]  \cdot h_{n-ki}.
\]
(The sum on the right hand side is well-defined, since $h_{n-ki}=0$ whenever
$i>\dfrac{n}{k}$.) \medskip

\textbf{(e)} \fbox{3} Set
\[
c_{i,j}:=%
\begin{cases}
2, & \text{if }i\equiv j\operatorname{mod}3;\\
-1, & \text{if }i\not \equiv j\operatorname{mod}3
\end{cases}
\ \ \ \ \ \ \ \ \ \ \text{for any }i,j\in\mathbb{Z}.
\]
Prove that each even $n\in\mathbb{N}$ satisfies%
\[
g_{3,n}=e_{n/2}^{2}+\sum_{i=0}^{\left(  n-2\right)  /2}c_{i,n-i}e_{i}e_{n-i},
\]
and that each odd $n\in\mathbb{N}$ satisfies%
\[
g_{3,n}=-\sum_{i=0}^{\left(  n-1\right)  /2}c_{i,n-i}e_{i}e_{n-i}.
\]

\end{exercise}

The next few exercises are about expressing the three families $\left(
h_{n}\right)  _{n\in\mathbb{N}}$, $\left(  e_{n}\right)  _{n\in\mathbb{N}}$,
$\left(  p_{n}\right)  _{n\geq1}$ of symmetric polynomials in terms of each
other. The first exercise is concerned with the case $N=2$:

\begin{exercise}
\label{exe.sf.newton.N=2}\textbf{(a)} \fbox{3} Prove that each $n\in
\mathbb{N}$ satisfies%
\[
h_{n}=\sum_{k=0}^{\left\lfloor n/2\right\rfloor }\left(  -1\right)
^{k}\dbinom{n-k}{k}e_{1}^{n-2k}e_{2}^{k}\ \ \ \ \ \ \ \ \ \ \text{if }N=2.
\]

\textbf{(b)} \fbox{3} Prove that each $n>0$ satisfies%
\[
p_{n}=\sum_{k=0}^{\left\lfloor n/2\right\rfloor }\left(  -1\right)  ^{k}%
\dfrac{n}{n-k}\dbinom{n-k}{k}e_{1}^{n-2k}e_{2}^{k}\ \ \ \ \ \ \ \ \ \ \text{if
}N=2.
\]

[\textbf{Hint:} Exercise \ref{exe.det.tridiag} \textbf{(b)} is useful for part
\textbf{(a)}. Compare also with \cite[Proposition 4.9.18]{20f}. Part
\textbf{(b)} can be proved using inclusion/exclusion: Imagine the $n$ numbers
$1,2,\ldots,n$ as the vertices of a regular $n$-gon. Show that $\dfrac{n}%
{n-k}\dbinom{n-k}{k}$ counts the ways to choose $k$ distinct vertices of this
$n$-gon such that no two chosen vertices are adjacent. Now consider the ways
to write either $x_{1}$ or $x_{2}$ on each vertex of this $n$-gon. There are
various other proofs as well, probably easier to find.]
\end{exercise}

Let us now return to the general case where $N$ is arbitrary:

\begin{exercise}
\label{exe.sf.eh-pol-same}\textbf{(a)} \fbox{3} Prove that there exists a
family of polynomials $\left(  P_{1},P_{2},P_{3},\ldots\right)  $, with each
$P_{n}$ being a polynomial in the ring $\mathbb{Z}\left[  x_{1},x_{2}%
,\ldots,x_{n}\right]  $, such that every positive integer $n$ satisfies%
\[
e_{n}=P_{n}\left[  h_{1},h_{2},\ldots,h_{n}\right]
\ \ \ \ \ \ \ \ \ \ \text{and}\ \ \ \ \ \ \ \ \ \ h_{n}=P_{n}\left[
e_{1},e_{2},\ldots,e_{n}\right]  .
\]
(This family begins with $P_{1}=x_{1}$ and $P_{2}=x_{1}^{2}-x_{2}$ and
$P_{3}=x_{1}^{3}-2x_{1}x_{2}+x_{3}$.) \medskip

\textbf{(b)} \fbox{2} Prove that there exists a family of polynomials $\left(
Q_{1},Q_{2},Q_{3},\ldots\right)  $, with each $Q_{n}$ being a polynomial in
the ring $\mathbb{Z}\left[  x_{1},x_{2},\ldots,x_{n}\right]  $, such that
every positive integer $n$ satisfies%
\[
p_{n}=\left(  -1\right)  ^{n-1}Q_{n}\left[  h_{1},h_{2},\ldots,h_{n}\right]
\ \ \ \ \ \ \ \ \ \ \text{and}\ \ \ \ \ \ \ \ \ \ p_{n}=Q_{n}\left[
e_{1},e_{2},\ldots,e_{n}\right]  .
\]
(This family begins with $Q_{1}=x_{1}$ and $Q_{2}=x_{1}^{2}-2x_{2}$ and
$Q_{3}=x_{1}^{3}-3x_{1}x_{2}+3x_{3}$.) \medskip

\textbf{(c)} \fbox{3} Express the polynomials $P_{n}$ explicitly as
determinants of certain matrices. \medskip

\textbf{(d)} \fbox{3} Express the polynomials $Q_{n}$ explicitly as
determinants of certain matrices. \medskip

\textbf{(e)} \fbox{1} What is the coefficient of $x_{i}$ in the polynomial
$P_{i}$ ? \medskip

\textbf{(f)} \fbox{1} What is the coefficient of $x_{i}$ in the polynomial
$Q_{i}$ ?
\end{exercise}

\begin{exercise}
Let $n\in\mathbb{N}$. For any permutation $\sigma\in S_{n}$, let
$\operatorname*{type}\sigma$ denote the cycle length partition of $\sigma$
(see Definition \ref{def.perm.cycs.cycs} \textbf{(b)}). Furthermore, for any
partition $\lambda=\left(  \lambda_{1},\lambda_{2},\ldots,\lambda_{\ell
}\right)  $ (with $\lambda_{1},\lambda_{2},\ldots,\lambda_{\ell}>0$), we
define the symmetric polynomial%
\[
p_{\lambda}:=p_{\lambda_{1}}p_{\lambda_{2}}\cdots p_{\lambda_{\ell}}%
\in\mathcal{S}.
\]

\textbf{(a)} \fbox{2} Prove that every $\sigma\in S_{n}$ satisfies%
\[
p_{\operatorname*{type}\sigma}=\sum_{\substack{f:\left[  n\right]
\rightarrow\left[  N\right]  \text{ is a map;}\\f\circ\sigma=f}}x_{f\left(
1\right)  }x_{f\left(  2\right)  }\cdots x_{f\left(  n\right)  }.
\]

\textbf{(b)} \fbox{4} Prove that%
\[
n!h_{n}=\sum_{\sigma\in S_{n}}p_{\operatorname*{type}\sigma}.
\]

[\textbf{Note:} Part \textbf{(b)} gives an expression of $n!h_{n}$ as a
polynomial (with integer coefficients) in $p_{1},p_{2},\ldots,p_{n}$.]
\end{exercise}

The following exercise is a symmetric-functions analogue of Exercise
\ref{exe.sign.chromatic.1}:

\begin{exercise}
\label{exe.sf.chromatic1}\fbox{3} Let us use the notations of Exercise
\ref{exe.sign.chromatic.1}. (Thus, $G$ is a finite undirected graph with
vertex set $V$ and edge set $E$.)

Define the \emph{chromatic symmetric polynomial} $X_{G}$ (in $N$ variables
$x_{1},x_{2},\ldots,x_{N}$) to be the polynomial%
\[
\sum_{\substack{c:V\rightarrow\left[  N\right]  \text{ is a}\\\text{proper
}N\text{-coloring}}}\ \ \prod_{v\in V}x_{c\left(  v\right)  }\in\mathcal{P}.
\]
For instance, if the graph $G$ is the length-$2$ path graph $%
\raisebox{-0.7pc}{
\begin{tikzpicture}%
[scale=1.5,thick,main node/.style={circle,fill=blue!20,draw}]
\node[main node] (1) at (1, 0) {$1$};
\node[main node] (2) at (2, 0) {$2$};
\node[main node] (3) at (3, 0) {$3$};
\draw(1) -- (2) -- (3);
\end{tikzpicture}
}%
$, then a proper $N$-coloring of $G$ is a map $c:\left[  3\right]
\rightarrow\left[  N\right]  $ satisfying $c\left(  1\right)  \neq c\left(
2\right)  $ and $c\left(  2\right)  \neq c\left(  3\right)  $, and therefore
its chromatic symmetric polynomial $X_{G}$ is%
\begin{align*}
\sum_{\substack{c:\left[  3\right]  \rightarrow\left[  N\right]  ;\\c\left(
1\right)  \neq c\left(  2\right)  ;\\c\left(  2\right)  \neq c\left(
3\right)  }}\ \ \prod_{v\in\left[  3\right]  }x_{c\left(  v\right)  }  &
=\sum_{\substack{i,j,k\in\left[  N\right]  ;\\i\neq j;\\j\neq k}}x_{i}%
x_{j}x_{k}=\underbrace{\sum_{\substack{i,j\in\left[  N\right]  ;\\i\neq
j}}x_{i}^{2}x_{j}}_{=p_{2}e_{1}-p_{3}}+\underbrace{\sum_{\substack{i,j,k\in
\left[  N\right]  ;\\i,j,k\text{ are distinct}}}x_{i}x_{j}x_{k}}_{=6e_{3}}\\
&  =p_{2}e_{1}-p_{3}+6e_{3}\\
&  =e_{2}e_{1}+3e_{3}\ \ \ \ \ \ \ \ \ \ \left(  \text{by some computation}%
\right) \\
&  =p_{1}^{3}-2p_{1}p_{2}+p_{3}\ \ \ \ \ \ \ \ \ \ \left(  \text{by some
computation}\right)  .
\end{align*}

\textbf{(a)} \fbox{1} Prove that $X_{G}\in\mathcal{S}$. \medskip

\textbf{(b)} \fbox{3} Prove that%
\[
X_{G}=\sum_{F\subseteq E}\left(  -1\right)  ^{\left\vert F\right\vert }%
\prod_{\substack{C\text{ is a connected component}\\\text{of the graph with
vertex set }V\\\text{and edge set }F}}p_{\left\vert C\right\vert }.
\]
(Here, $\left\vert C\right\vert $ denotes the number of vertices in the
connected component $C$.) \medskip

\textbf{(c)} \fbox{1} Prove that
\[
X_{G}\left[  \underbrace{1,1,\ldots,1}_{N\text{ ones}}\right]  =\chi
_{G}\left(  N\right)
\]
(where $\chi_{G}\left(  N\right)  $ is as defined in Exercise
\ref{exe.sign.chromatic.1}). \medskip

Next, let us generalize $X_{G}$: For each vertex $v\in V$, let $w\left(
v\right)  $ be a positive integer; we shall call $w\left(  v\right)  $ the
\emph{weight} of $v$. We define the \emph{weighted chromatic symmetric
polynomial} $X_{G,w}$ to be%
\[
\sum_{\substack{c:V\rightarrow\left[  N\right]  \text{ is a}\\\text{proper
}N\text{-coloring}}}\ \ \prod_{v\in V}x_{c\left(  v\right)  }^{w\left(
v\right)  }\in\mathcal{P}.
\]
Thus, if all weights $w\left(  v\right)  $ equal $1$, then $X_{G,w}=X_{G}$.
\medskip

\textbf{(d)} \fbox{2} Generalize the claim of part \textbf{(b)} to $X_{G,w}$.
\end{exercise}

\begin{exercise}
\label{exe.sf.smirnov-stagdes}Let $n$ be a positive integer. If $w\in\left[
N\right]  ^{n}$ is an $n$-tuple, then

\begin{itemize}
\item we let $w_{1},w_{2},\ldots,w_{n}$ denote the $n$ entries of $w$ (so that
$w=\left(  w_{1},w_{2},\ldots,w_{n}\right)  $);

\item we define the \emph{descent set} $\operatorname*{Des}w$ of $w$ by
\[
\operatorname*{Des}w:=\left\{  i\in\left[  n-1\right]  \ \mid\ w_{i}%
>w_{i+1}\right\}  ;
\]

\item we define the \emph{stagnation set} $\operatorname*{Stag}w$ of $w$ by
\[
\operatorname*{Stag}w:=\left\{  i\in\left[  n-1\right]  \ \mid\ w_{i}%
=w_{i+1}\right\}  ;
\]

\item we define the monomial $x_{w}$ to be $x_{w_{1}}x_{w_{2}}\cdots x_{w_{n}%
}$.
\end{itemize}

\noindent For instance, the $7$-tuple $\left(  2,2,4,1,4,4,2\right)  $ has
descent set $\left\{  3,6\right\}  $ and stagnation set $\left\{  1,5\right\}
$. \medskip

\textbf{(a)} \fbox{1} Fix $s\in\mathbb{N}$. Prove that
\[
\sum_{\substack{w\in\left[  N\right]  ^{n};\\\left\vert \operatorname*{Stag}%
w\right\vert =s}}x_{w}\in\mathcal{S}.
\]

\textbf{(b)} \fbox{1} Identify this sum as a sum of weighted chromatic
symmetric polynomials of path graphs (see Exercise \ref{exe.sf.chromatic1}).
\medskip

\textbf{(c)} \fbox{3} Fix $d\in\mathbb{N}$ and $s\in\mathbb{N}$. Prove that
\[
\sum_{\substack{w\in\left[  N\right]  ^{n};\\\left\vert \operatorname*{Des}%
w\right\vert =d;\\\left\vert \operatorname*{Stag}w\right\vert =s}}x_{w}%
\in\mathcal{S}.
\]

\textbf{(d)} \fbox{6} Fix $d\in\mathbb{N}$ and $s\in\mathbb{N}$. Prove that
the three polynomials
\[
\sum_{\substack{w\in\left[  N\right]  ^{n};\\\left\vert \operatorname*{Des}%
w\right\vert =d;\\\left\vert \operatorname*{Stag}w\right\vert =s;\\w_{1}%
<w_{n}}}x_{w},\ \ \ \ \ \ \ \ \ \ \sum_{\substack{w\in\left[  N\right]
^{n};\\\left\vert \operatorname*{Des}w\right\vert =d;\\\left\vert
\operatorname*{Stag}w\right\vert =s;\\w_{1}=w_{n}}}x_{w}%
,\ \ \ \ \ \ \ \ \ \ \sum_{\substack{w\in\left[  N\right]  ^{n};\\\left\vert
\operatorname*{Des}w\right\vert =d;\\\left\vert \operatorname*{Stag}%
w\right\vert =s;\\w_{1}>w_{n}}}x_{w}%
\]
all belong to $\mathcal{S}$.
\end{exercise}

\begin{exercise}
\label{exe.sf.hm-N+1-as-sum-of-fracs}\fbox{4} Prove that every $m\in
\mathbb{N}$ satisfies%
\[
\sum_{k=1}^{N}\dfrac{x_{k}^{m}}{\prod_{i\in\left[  N\right]  \setminus\left\{
k\right\}  }\left(  x_{k}-x_{i}\right)  }=h_{m-N+1}.
\]
(This is an equality in the localization of the polynomial ring $\mathcal{P}$
at the multiplicative subset generated by the pairwise differences
$x_{i}-x_{j}$ for all $i<j$. If $K$ is a field, you can also view it as an
equality in the field of rational functions $K\left(  x_{1},x_{2},\ldots
,x_{N}\right)  $.)

[\textbf{Example:} For instance, if $N=3$ (and if we rename $x_{1},x_{2}%
,x_{3}$ as $x,y,z$), then this exercise is claiming that
\[
\dfrac{x^{m}}{\left(  x-y\right)  \left(  x-z\right)  }+\dfrac{y^{m}}{\left(
y-z\right)  \left(  y-x\right)  }+\dfrac{z^{m}}{\left(  z-x\right)  \left(
z-y\right)  }=h_{m-2}\left[  x,y,z\right]  .
\]
Note that the right hand side is $0$ when $m<2$.]
\end{exercise}

\begin{exercise}
\fbox{2} Assume that $N>0$. Prove that
\[
\sum_{k=1}^{N}\dfrac{1}{x_{k}\prod_{i\in\left[  N\right]  \setminus\left\{
k\right\}  }\left(  x_{k}-x_{i}\right)  }=\dfrac{\left(  -1\right)  ^{N-1}%
}{x_{1}x_{2}\cdots x_{N}}.
\]
(This is an equality in the localization of the polynomial ring $\mathcal{P}$
at the multiplicative subset generated by the variables $x_{1},x_{2}%
,\ldots,x_{N}$ and their pairwise differences $x_{i}-x_{j}$ for all $i<j$. If
$K$ is a field, you can also view it as an equality in the field of rational
functions $K\left(  x_{1},x_{2},\ldots,x_{N}\right)  $.)

[\textbf{Hint:} This can be derived from Exercise
\ref{exe.sf.hm-N+1-as-sum-of-fracs} via a neat trick.]
\end{exercise}

\begin{exercise}
\label{exe.sf.nica-sum-1}\fbox{4} Prove that
\[
\sum_{k=1}^{N}\dfrac{x_{k}\prod_{i\in\left[  N\right]  \setminus\left\{
k\right\}  }\left(  x_{k}+x_{i}\right)  }{\prod_{i\in\left[  N\right]
\setminus\left\{  k\right\}  }\left(  x_{k}-x_{i}\right)  }=x_{1}+x_{2}%
+\cdots+x_{N}.
\]
(The same comments as in Exercise \ref{exe.sf.hm-N+1-as-sum-of-fracs} apply.)

[\textbf{Hint:} It suffices to prove this in the field $\mathbb{Q}\left(
x_{1},x_{2},\ldots,x_{N}\right)  $ of rational functions over $\mathbb{Q}$
(why?). Define the polynomial $P\left(  t\right)  :=\prod_{i\in\left[
N\right]  }\left(  t+x_{i}\right)  $ (in the indeterminate $t$ over this
field). Observe that the numerator $x_{k}\prod_{i\in\left[  N\right]
\setminus\left\{  k\right\}  }\left(  x_{k}+x_{i}\right)  $ can be rewritten
as $\dfrac{1}{2}P\left(  x_{k}\right)  $ (why?). On the other hand, $P\left(
t\right)  =t^{N}+e_{1}t^{N-1}+e_{2}t^{N-2}+\cdots+e_{N}t^{0}$. Now apply
Exercise \ref{exe.sf.hm-N+1-as-sum-of-fracs}.]
\end{exercise}

\begin{exercise}
\label{exe.sf.nica-sum-2}\fbox{5} Prove that
\[
\sum_{k=1}^{N}\dfrac{\prod_{i\in\left[  N\right]  \setminus\left\{  k\right\}
}\left(  x_{k}+x_{i}\right)  }{\prod_{i\in\left[  N\right]  \setminus\left\{
k\right\}  }\left(  x_{k}-x_{i}\right)  }=%
\begin{cases}
0, & \text{if }N\text{ is even};\\
1, & \text{if }N\text{ is odd.}%
\end{cases}
\]
(The same comments as in Exercise \ref{exe.sf.hm-N+1-as-sum-of-fracs} apply.)

[\textbf{Hint:} It suffices to prove this in the field $\mathbb{Q}\left(
x_{1},x_{2},\ldots,x_{N}\right)  $ of rational functions over $\mathbb{Q}$
(why?). Define the polynomial $P\left(  t\right)  :=\prod_{i\in\left[
N\right]  }\left(  x_{i}+t\right)  -\prod_{i\in\left[  N\right]  }\left(
x_{i}-t\right)  $ (in the indeterminate $t$ over this field). This polynomial
$P\left(  t\right)  $ is divisible by $2t$. Let $Q\left(  t\right)  $ be the
quotient. What then?]
\end{exercise}

\begin{exercise}
Let $n\in\mathbb{N}$. \medskip

\textbf{(a)} \fbox{3} Prove that%
\[
h_{n}\left[  x_{1}^{2},x_{2}^{2},\ldots,x_{N}^{2}\right]  =\sum_{i=0}%
^{2n}\left(  -1\right)  ^{i}h_{i}h_{2n-i}.
\]

\textbf{(b)} \fbox{3} Prove that%
\[
e_{n}\left[  x_{1}^{2},x_{2}^{2},\ldots,x_{N}^{2}\right]  =\sum_{i=0}%
^{2n}\left(  -1\right)  ^{n-i}e_{i}e_{2n-i}.
\]

\textbf{(c)} \fbox{2} Solve Exercise \ref{exe.fps.altern-quasivdm} again using
part \textbf{(b)}.
\end{exercise}

\begin{exercise}
\fbox{4} Let $i\in\left[  N+1\right]  $ and $p\in\mathbb{N}$. Prove that%
\[
h_{p}\left[  x_{i},x_{i+1},\ldots,x_{N}\right]  =\sum_{t=0}^{i-1}\left(
-1\right)  ^{t}e_{t}\left[  x_{1},x_{2},\ldots,x_{i-1}\right]  \cdot h_{p-t}.
\]
(The \textquotedblleft$h_{p-t}$\textquotedblright\ at the end of the right
hand side means $h_{p-t}\left[  x_{1},x_{2},\ldots,x_{N}\right]  $.)
\end{exercise}

\begin{exercise}
\textbf{(a)} \fbox{2} Prove that each $i\in\left[  N\right]  $ and
$j\in\left[  N\right]  $ satisfy%
\[
\dfrac{\partial e_{j}}{\partial x_{i}}=e_{j-1}\left[  x_{1},x_{2}%
,\ldots,\widehat{x_{i}},\ldots,x_{N}\right]  ,
\]
where the hat over the \textquotedblleft$x_{i}$\textquotedblright\ means
\textquotedblleft omit the $x_{i}$ entry\textquotedblright\ (that is, the
expression \textquotedblleft$x_{1},x_{2},\ldots,\widehat{x_{i}},\ldots,x_{N}%
$\textquotedblright\ is to be understood as \textquotedblleft$x_{1}%
,x_{2},\ldots,x_{i-1},x_{i+1},x_{i+2},\ldots,x_{N}$\textquotedblright).
\medskip

\textbf{(b)} \fbox{3} Prove that%
\[
\det\left(  \left(  \dfrac{\partial e_{j}}{\partial x_{i}}\right)  _{1\leq
i\leq N,\ 1\leq j\leq N}\right)  =\prod_{1\leq i<j\leq N}\left(  x_{i}%
-x_{j}\right)  .
\]
\medskip

\textbf{(c)} \fbox{5} Prove that%
\[
\det\left(  \left(  \dfrac{\partial h_{j}}{\partial x_{i}}\right)  _{1\leq
i\leq N,\ 1\leq j\leq N}\right)  =\prod_{1\leq i<j\leq N}\left(  x_{j}%
-x_{i}\right)  .
\]

[\textbf{Hint:} For part \textbf{(c)}, it is helpful to first prove the
following more general result:

Let $u_{1},u_{2},\ldots,u_{N}$ be any $N$ polynomials in $\mathcal{P}$. Let
$v_{1},v_{2},\ldots,v_{N}$ be any $N$ polynomials in the polynomial ring
$K\left[  y_{1},y_{2},\ldots,y_{N}\right]  $ (in $N$ indeterminates
$y_{1},y_{2},\ldots,y_{N}$). For each $j\in\left[  N\right]  $, let
$w_{j}:=v_{j}\left[  u_{1},u_{2},\ldots,u_{N}\right]  \in\mathcal{P}$ be the
polynomial obtained from $v_{j}$ by substituting $u_{1},u_{2},\ldots,u_{N}$
for $y_{1},y_{2},\ldots,y_{N}$. Then,%
\begin{align*}
&  \left(  \dfrac{\partial w_{j}}{\partial x_{i}}\right)  _{1\leq i\leq
N,\ 1\leq j\leq N}\\
&  =\left(  \dfrac{\partial u_{j}}{\partial x_{i}}\right)  _{1\leq i\leq
N,\ 1\leq j\leq N}\cdot\left(  \dfrac{\partial v_{j}}{\partial y_{i}}\left[
u_{1},u_{2},\ldots,u_{N}\right]  \right)  _{1\leq i\leq N,\ 1\leq j\leq N}.
\end{align*}
This fact is a generalization of the chain rule and can be proved, e.g., by
decomposing $v_{j}$ into monomials and considering each monomial separately.
Now, use parts \textbf{(a)} and \textbf{(e)} of Exercise
\ref{exe.sf.eh-pol-same} to apply this to our setting.]
\end{exercise}

\subsubsection{$N$-partitions and monomial symmetric polynomials}

\begin{exercise}
\fbox{3} Let $M\in\left\{  0,1,\ldots,N\right\}  $. For any $M$-partition
$\mu=\left(  \mu_{1},\mu_{2},\ldots,\mu_{M}\right)  $ and any $\left(
N-M\right)  $-partition $\nu=\left(  \nu_{1},\nu_{2},\ldots,\nu_{N-M}\right)
$, we let $\mu\sqcup\nu$ denote the $N$-partition obtained by sorting the
$N$-tuple $\left(  \mu_{1},\mu_{2},\ldots,\mu_{M},\nu_{1},\nu_{2},\ldots
,\nu_{N-M}\right)  $ in weakly decreasing order. (For example, $\left(
3,2,0\right)  \sqcup\left(  4,2,1,1\right)  =\left(  4,3,2,2,1,1,0\right)  $.)

Let $\lambda$ be any $N$-partition. Prove that%
\begin{align*}
&  m_{\lambda}\left[  x_{1},x_{2},\ldots,x_{N}\right] \\
&  =\sum_{\substack{\mu\text{ is an }M\text{-partition;}\\\nu\text{ is an
}\left(  N-M\right)  \text{-partition;}\\\mu\sqcup\nu=\lambda}}m_{\mu}\left[
x_{1},x_{2},\ldots,x_{M}\right]  \cdot m_{\nu}\left[  x_{M+1},x_{M+2}%
,\ldots,x_{N}\right]  .
\end{align*}

\end{exercise}

\begin{exercise}
\label{exe.sf.prodxi+xj-tour}\fbox{2} We shall use the notion of a tournament,
as defined in Exercise \ref{exe.det.vander-tour}. Let $T$ be the set of all
tournaments with vertex set $\left[  N\right]  $. We define the
\emph{scoreboard} $\operatorname*{scb}D$ of a tournament $D\in T$ to be the
$N$-tuple $\left(  s_{1},s_{2},\ldots,s_{N}\right)  \in\mathbb{N}^{N}$, where
\begin{align*}
s_{j}:=  &  \left(  \text{\# of arcs of }D\text{ that end at }j\right) \\
=  &  \left(  \text{\# of }i\in\left[  N\right]  \text{ such that }\left(
i,j\right)  \text{ is an arc of }D\right)
\end{align*}
for each $j\in\left[  N\right]  $.

Prove that%
\[
\prod_{1\leq i<j\leq N}\left(  x_{i}+x_{j}\right)  =\sum_{\substack{\lambda
\text{ is an }N\text{-partition}\\\text{of }N\left(  N-1\right)
/2}}t_{\lambda}m_{\lambda},
\]
where $t_{\lambda}$ denotes the \# of tournaments $D\in T$ with scoreboard
$\operatorname*{scb}D=\lambda$.
\end{exercise}

\subsubsection{Schur polynomials}

\begin{exercise}
\fbox{2} Assume that $N>1$. Without using the Jacobi--Trudi identities, prove
that $s_{\left(  n,1,0,0,\ldots,0\right)  }=h_{n}h_{1}-h_{n+1}$ for each
positive integer $n$. (Here, $\left(  n,1,0,0,\ldots,0\right)  $ denotes the
$N$-partition whose first two entries are $n$ and $1$ while all remaining
entries equal $0$.)
\end{exercise}

\begin{exercise}
Let $\lambda=\left(  \lambda_{1},\lambda_{2},\ldots,\lambda_{N}\right)  $ be
any $N$-partition. \medskip

\textbf{(a)} \fbox{2} Let $k\in\mathbb{N}$. Show that%
\[
s_{\left(  \lambda_{1}+k,\lambda_{2}+k,\ldots,\lambda_{N}+k\right)  }=\left(
x_{1}x_{2}\cdots x_{N}\right)  ^{k}\cdot s_{\lambda}.
\]

\textbf{(b)} \fbox{4} Let $k\in\mathbb{N}$ be such that $k\geq\lambda_{1}$.
(If $N=0$, then $\lambda_{1}$ is understood to mean $0$.) Let $\mu$ be the
$N$-partition $\left(  k-\lambda_{N},k-\lambda_{N-1},\ldots,k-\lambda
_{1}\right)  $. Prove that%
\[
s_{\mu}=\left(  x_{1}x_{2}\cdots x_{N}\right)  ^{k}\cdot s_{\lambda}\left[
x_{1}^{-1},x_{2}^{-1},\ldots,x_{N}^{-1}\right]
\]
in the Laurent polynomial ring $K\left[  x_{1}^{\pm},x_{2}^{\pm},\ldots
,x_{N}^{\pm}\right]  $. (We have never formally defined this ring, but it
should suffice to know that the elements of $K\left[  x_{1}^{\pm},x_{2}^{\pm
},\ldots,x_{N}^{\pm}\right]  $ are formal $K$-linear combinations of
\textquotedblleft Laurent monomials\textquotedblright\ $x_{1}^{a_{1}}%
x_{2}^{a_{2}}\cdots x_{N}^{a_{N}}$ with $a_{1},a_{2},\ldots,a_{N}\in
\mathbb{Z}$.)
\end{exercise}

\begin{exercise}
Let $M\in\mathbb{N}$. Let $\lambda=\left(  \lambda_{1},\lambda_{2}%
,\ldots,\lambda_{M}\right)  $ and $\mu=\left(  \mu_{1},\mu_{2},\ldots,\mu
_{M}\right)  $ be any two $M$-partitions. \medskip

\textbf{(a)} \fbox{1} Show that each $k\in\mathbb{N}$ satisfies%
\[
s_{\left(  \lambda_{1}+k,\lambda_{2}+k,\ldots,\lambda_{M}+k\right)  /\left(
\mu_{1}+k,\mu_{2}+k,\ldots,\mu_{M}+k\right)  }=s_{\lambda/\mu}.
\]

\textbf{(b)} \fbox{4} Let $p\in\mathbb{N}$ be such that $p\geq\lambda_{1}$ and
$p\geq\mu_{1}$. Show that%
\[
s_{\left(  p-\mu_{M},p-\mu_{M-1},\ldots,p-\mu_{1}\right)  /\left(
p-\lambda_{M},p-\lambda_{M-1},\ldots,p-\lambda_{1}\right)  }=s_{\lambda/\mu}.
\]

[\textbf{Note:} The Young diagram of \newline$\left(  \lambda_{1}%
+k,\lambda_{2}+k,\ldots,\lambda_{M}+k\right)  /\left(  \mu_{1}+k,\mu
_{2}+k,\ldots,\mu_{M}+k\right)  $ is obtained from $Y\left(  \lambda
/\mu\right)  $ by a parallel shift, whereas the Young diagram of $\left(
p-\mu_{M},p-\mu_{M-1},\ldots,p-\mu_{1}\right)  /\left(  p-\lambda
_{M},p-\lambda_{M-1},\ldots,p-\lambda_{1}\right)  $ is obtained from $Y\left(
\lambda/\mu\right)  $ by a $180^{\circ}$-rotation.]
\end{exercise}

\begin{exercise}
Let $M\in\mathbb{N}$. Let $\lambda$ and $\mu$ be two $M$-partitions with
$\mu\subseteq\lambda$. Recall the Bender--Knuth involutions $\beta
_{k}:\operatorname*{SSYT}\left(  \lambda/\mu\right)  \rightarrow
\operatorname*{SSYT}\left(  \lambda/\mu\right)  $ defined for all $k\in\left[
N-1\right]  $ in the proof of Theorem \ref{thm.sf.skew-schur-symm}. \medskip

\textbf{(a)} \fbox{1} Prove that $\beta_{i}\circ\beta_{j}=\beta_{j}\circ
\beta_{i}$ whenever $i$ and $j$ are two elements of $\left[  N-1\right]  $
satisfying $\left\vert i-j\right\vert >1$. \medskip

\textbf{(b)} \fbox{4} Prove that $\left(  \beta_{1}\circ\beta_{2}\right)
^{6}=\operatorname*{id}$ if $\mu=\mathbf{0}$ (where $\mathbf{0}$ means the
$M$-partition $\left(  0,0,\ldots,0\right)  $). \medskip

\textbf{(c)} \fbox{2} Find an example where $\left(  \beta_{1}\circ\beta
_{2}\right)  ^{3}\neq\operatorname*{id}$ for $\mu=\mathbf{0}$.
\end{exercise}

The next exercise answers a rather natural question about the definition of a
semistandard tableau (Definition \ref{def.sf.ssyt}): What happens when one
replaces \textquotedblleft increase strictly down each
column\textquotedblright\ by \textquotedblleft increase weakly down each
column\textquotedblright? This replacement gives rise to a more liberal notion
of semistandard tableaux (which I call \textquotedblleft semi-semistandard
tableaux\textquotedblright), and to a variant of the Schur polynomial that
correspondingly has more terms. As the following exercise shows, alas, this
new polynomial is rarely ever symmetric:

\begin{exercise}
\label{exe.sf.fools-schur-symm}Let $\lambda$ be an $N$-partition, or (more
generally) a partition of any length. A Young tableau $T$ of shape $\lambda$
will be called \emph{semi-semistandard} if its entries

\begin{itemize}
\item increase weakly along each row;

\item increase weakly down each column.
\end{itemize}

For instance, the tableau $\ytableaushort{11,12}$ is semi-semistandard but not
semistandard. (Semi-semistandard tableaux are actually known as \emph{reverse
plane partitions} for historical reasons: If you replace the words
\textquotedblleft increase\textquotedblright\ by \textquotedblleft
decrease\textquotedblright, then they become \textquotedblleft%
2-dimensional\textquotedblright\ analogues of partitions, in that they are
tables of positive integers that weakly decrease in two directions.)

We define the \textquotedblleft fool's Schur polynomial\textquotedblright%
\ $\widehat{s}_{\lambda}\in\mathcal{P}$ by%
\[
\widehat{s}_{\lambda}:=\sum_{T\in\operatorname*{SSSYT}\left(  \lambda\right)
}x_{T},
\]
where $\operatorname*{SSSYT}\left(  \lambda\right)  $ is the set of all
semi-semistandard tableaux of shape $\lambda$. \medskip

\textbf{(a)} \fbox{2} Assume that $N\geq3$ and $K=\mathbb{Z}$. Prove that the
polynomial $\widehat{s}_{\lambda}$ is symmetric if and only if the Young
diagram $Y\left(  \lambda\right)  $ either consists of a single row (i.e., we
have $\lambda=\left(  n,0,0,\ldots,0\right)  $ for some $n\in\mathbb{N}$) or
consists of a single column (i.e., we have $\lambda=\left(  1,1,\ldots
,1,0,0,\ldots,0\right)  $ for some number of $1$'s). \medskip

\textbf{(b)} \fbox{2} Now, replace the assumption $N\geq3$ by $N=2$. Prove
that the polynomial $\widehat{s}_{\lambda}$ is symmetric if and only if the
Young diagram $Y\left(  \lambda\right)  $ is a rectangle (i.e., all nonzero
entries of $\lambda$ are equal). \medskip

[\textbf{Hint:} Compare the coefficients of $x_{1}x_{2}^{k}$ and $x_{1}%
^{k}x_{2}$, as well as the coefficients of $x_{1}x_{2}x_{3}^{k}$ and
$x_{1}x_{2}^{k}x_{3}$.]
\end{exercise}

The next exercise asks you to explore some possible (or impossible) variations
on the proof of Lemma \ref{lem.sf.stemb-lem}.

\begin{exercise}
Recall the definition of the sign-reversing involution $f:\mathcal{X}%
\rightarrow\mathcal{X}$ in the proof of Lemma \ref{lem.sf.stemb-lem}. \medskip

\textbf{(a)} \fbox{2} Would $f$ still be a sign-reversing involution if we
defined $j$ to be the \textbf{smallest} (rather than largest) violator of $T$
? \medskip

\textbf{(b)} \fbox{1} Would $f$ still be a sign-reversing involution if we
defined $k$ to be the \textbf{largest} (rather than smallest) misstep of $T$ ?
\end{exercise}

Here are some more properties of Schur polynomials:

\begin{exercise}
\fbox{3} Let $\rho:=\left(  N-1,N-2,\ldots,N-N\right)  \in\mathbb{N}^{N}$.
Prove that%
\[
s_{\rho}=\prod_{1\leq i<j\leq N}\left(  x_{i}+x_{j}\right)  .
\]

\end{exercise}

\begin{exercise}
Let $\lambda=\left(  \lambda_{1},\lambda_{2},\ldots,\lambda_{N}\right)  $ be
an $N$-partition. \medskip

\textbf{(a)} \fbox{2} Prove that%
\[
s_{\lambda}\left(  q^{0},q^{1},\ldots,q^{N-1}\right)  =\prod_{1\leq i<j\leq
N}\dfrac{q^{\lambda_{i}+N-i}-q^{\lambda_{j}+N-j}}{q^{N-i}-q^{N-j}}%
\]
in the polynomial ring $\mathbb{Z}\left[  q\right]  $. \medskip

\textbf{(b)} \fbox{2} Using $q$-integers (Definition
\ref{def.pars.qbinom.qint}), rewrite this as%
\[
s_{\lambda}\left(  q^{0},q^{1},\ldots,q^{N-1}\right)  =q^{\eta\left(
\lambda\right)  }\prod_{1\leq i<j\leq N}\dfrac{\left[  \lambda_{i}-\lambda
_{j}+j-i\right]  _{q}}{\left[  j-i\right]  _{q}},
\]
where $\eta\left(  \lambda\right)  =\sum_{j=1}^{N}\left(  j-1\right)
\lambda_{j}$. \medskip

\textbf{(c)} \fbox{1} Prove that the \# of semistandard tableaux of shape
$\lambda$ with entries in $\left[  N\right]  $ is%
\[
s_{\lambda}\left(  \underbrace{1,1,\ldots,1}_{N\text{ times}}\right)
=\prod_{1\leq i<j\leq N}\dfrac{\lambda_{i}-\lambda_{j}+j-i}{j-i}.
\]

\end{exercise}

\begin{exercise}
\label{exe.sf.strips.entries}\fbox{3} Prove Proposition
\ref{prop.sf.strips.entries}.
\end{exercise}

\begin{exercise}
\label{exe.sf.pieri}\textbf{(a)} \fbox{4} Prove Theorem \ref{thm.sf.pieri}
\textbf{(a)}. \medskip

\textbf{(b)} \fbox{4} Prove Theorem \ref{thm.sf.pieri} \textbf{(b)}. \medskip

[\textbf{Hint:} For part \textbf{(a)}, apply Theorem \ref{thm.sf.lr-zy} to
$N$, $\left(  n,0,0,\ldots,0\right)  $, $\mathbf{0}$ and $\mu$ instead of $M$,
$\lambda$, $\mu$ and $\nu$, and characterize the $\mu$-Yamanouchi semistandard
tableaux of shape $\left(  n,0,0,\ldots,0\right)  /\mathbf{0}$. Proceed
likewise for part \textbf{(b)}.]
\end{exercise}

\begin{exercise}
\label{exe.sf.strips.count}\fbox{3} Let $n\in\mathbb{N}$ be such that $n\leq
N$. Prove that
\begin{align*}
&  \left(  \text{\# of horizontal strips }\lambda/\mu\text{ satisfying
}\left\vert \lambda\right\vert =n\right) \\
&  =\left(  \text{\# of vertical strips }\lambda/\mu\text{ satisfying
}\left\vert \lambda\right\vert =n\right) \\
&  =\sum_{i=0}^{n}p\left(  i\right)  \cdot p\left(  n-i\right)  .
\end{align*}
(See Definition \ref{def.pars.pn-pkn} \textbf{(b)} for the definition of
$p\left(  i\right)  $ and $p\left(  n-i\right)  $.)
\end{exercise}

\begin{exercise}
\label{exe.sf.jt-h}Complete the proof of Theorem \ref{thm.sf.jt-h} sketched
above: \medskip

\textbf{(a)} \fbox{2} Prove Observation 1. \medskip

\textbf{(b)} \fbox{3} Prove Observation 2.
\end{exercise}

\begin{exercise}
\label{exe.sf.jt-e}\fbox{6} Prove Theorem \ref{thm.sf.jt-e}.
\end{exercise}

\begin{exercise}
\fbox{5} Let $n$ be a positive integer. Prove that%
\[
p_{n}=\sum_{i=0}^{\min\left\{  n,N\right\}  -1}\left(  -1\right)
^{i}s_{Q\left(  i\right)  },
\]
where $Q\left(  i\right)  $ is the $N$-partition $\left(
n-i,\underbrace{1,1,\ldots,1}_{i\text{ many }1\text{'s}}%
,\underbrace{0,0,\ldots,0}_{N-i-1\text{ many }0\text{'s}}\right)  $.
\end{exercise}

\begin{exercise}
\label{exe.sf.ec2-exe7.34}\fbox{5} Let $\lambda=\left(  \lambda_{1}%
,\lambda_{2},\ldots,\lambda_{N}\right)  $ and $\mu=\left(  \mu_{1},\mu
_{2},\ldots,\mu_{N}\right)  $ be two $N$-partitions. Prove that%
\[
s_{\lambda}s_{\mu}=\det\left(  \left(  h_{\lambda_{i}+\mu_{N+1-j}-i+j}\right)
_{1\leq i\leq N,\ 1\leq j\leq N}\right)  .
\]

[\textbf{Hint:} Fix some $m\in\mathbb{N}$ that is larger than all of $\mu
_{1},\mu_{2},\ldots,\mu_{N}$, and consider the skew partition $\nu
/\kappa=\left(  \nu_{1},\nu_{2},\ldots,\nu_{N}\right)  /\left(  \kappa
_{1},\kappa_{2},\ldots,\kappa_{N}\right)  $, where $\nu_{i}:=\lambda_{i}+m$
and $\kappa_{i}:=m-\mu_{N+1-i}$ for all $i\in\left[  N\right]  $. What are the
semistandard tableaux of shape $\nu/\kappa$, and what does this mean for
$s_{\nu/\kappa}$ ?]
\end{exercise}

\begin{exercise}
\fbox{5} Prove the \emph{flagged first Jacobi--Trudi formula}: Let
$M\in\mathbb{N}$. Let $\lambda=\left(  \lambda_{1},\lambda_{2},\ldots
,\lambda_{M}\right)  $ and $\mu=\left(  \mu_{1},\mu_{2},\ldots,\mu_{M}\right)
$ be two $M$-partitions (i.e., weakly decreasing $M$-tuples of nonnegative
integers). Let $\alpha=\left(  \alpha_{1},\alpha_{2},\ldots,\alpha_{M}\right)
$ and $\beta=\left(  \beta_{1},\beta_{2},\ldots,\beta_{M}\right)  $ be two
weakly increasing sequences of elements of $\left\{  0,1,\ldots,N\right\}  $.
(As usual, \textquotedblleft weakly increasing\textquotedblright\ means that
$\alpha_{1}\leq\alpha_{2}\leq\cdots\leq\alpha_{M}$ and $\beta_{1}\leq\beta
_{2}\leq\cdots\leq\beta_{M}$.) Let%
\[
s_{\lambda/\mu,\ \alpha/\beta}:=\sum_{\substack{T\in\operatorname*{SSYT}%
\left(  \lambda/\mu\right)  ;\\\alpha_{i}<T\left(  i,j\right)  \leq\beta
_{i}\text{ for all }\left(  i,j\right)  \in Y\left(  \lambda/\mu\right)
}}x_{T}\in\mathcal{P}.
\]
(This is the same sum as $s_{\lambda/\mu}$, but restricted to those
semistandard tableaux whose entries in the $i$-th row belong to the half-open
interval $\left(  \alpha_{i},\beta_{i}\right]  $ for each $i\in\left[
M\right]  $. This polynomial is known as a \emph{(row-)flagged Schur
polynomial}; in general, it is not symmetric.) Then,%
\[
s_{\lambda/\mu,\ \alpha/\beta}=\det\left(  \left(  h_{\lambda_{i}-\mu_{j}%
-i+j}\left[  x_{\alpha_{j}+1},x_{\alpha_{j}+2},\ldots,x_{\beta_{i}}\right]
\right)  _{1\leq i\leq M,\ 1\leq j\leq M}\right)  .
\]
(If $\alpha_{j}\geq\beta_{i}$, then \textquotedblleft$x_{\alpha_{j}%
+1},x_{\alpha_{j}+2},\ldots,x_{\beta_{i}}$\textquotedblright\ is understood to
be an empty list, so that $h_{\lambda_{i}-\mu_{j}-i+j}\left[  x_{\alpha_{j}%
+1},x_{\alpha_{j}+2},\ldots,x_{\beta_{i}}\right]  $ is the complete
homogeneous symmetric polynomial $h_{\lambda_{i}-\mu_{j}-i+j}$ of $0$
variables; this equals $1$ if $\lambda_{i}-\mu_{j}-i+j=0$ and $0$ otherwise.)
\end{exercise}

\begin{exercise}
\label{exe.sf.differential-sum}Let $d$ denote the differential operator%
\[
\dfrac{\partial}{\partial x_{1}}+\dfrac{\partial}{\partial x_{2}}%
+\cdots+\dfrac{\partial}{\partial x_{N}}\text{ on }\mathcal{P}.
\]
Explicitly, this operator $d$ is the $K$-linear map from $\mathcal{P}$ to
$\mathcal{P}$ that sends
\[
\text{each monomial }x_{1}^{a_{1}}x_{2}^{a_{2}}\cdots x_{N}^{a_{N}}\text{ to
}\sum_{\substack{i\in\left[  N\right]  ;\\a_{i}>0}}a_{i}\underbrace{x_{1}%
^{a_{1}}x_{2}^{a_{2}}\cdots x_{i-1}^{a_{i-1}}x_{i}^{a_{i}-1}x_{i+1}^{a_{i+1}%
}x_{i+2}^{a_{i+2}}\cdots x_{N}^{a_{N}}}_{\substack{\text{This is just the
monomial }x_{1}^{a_{1}}x_{2}^{a_{2}}\cdots x_{N}^{a_{N}}\text{,}\\\text{with
the exponent on }x_{i}\text{ decreased by }1}}.
\]
As usual for linear operators, we abbreviate $d\left(  f\right)  $ by $df$
(when $f\in\mathcal{P}$). \medskip

Prove the following: \medskip

\textbf{(a)} \fbox{1} We have $d\left(  fg\right)  =\left(  df\right)  \cdot
g+f\cdot dg$ for any $f,g\in\mathcal{P}$. (In the lingo of algebraists, this
is saying that $d$ is a \emph{derivation} on $\mathcal{P}$.) \medskip

\textbf{(b)} \fbox{1} We have $d\left(  e_{n}\right)  =\left(  N-n+1\right)
e_{n-1}$ for any $n\in\mathbb{Z}$. \medskip

\textbf{(c)} \fbox{1} We have $d\left(  h_{n}\right)  =\left(  N+n-1\right)
h_{n-1}$ for any $n\in\mathbb{Z}$. \medskip

\textbf{(d)} \fbox{1} We have $d\left(  p_{n}\right)  =np_{n-1}$ for any
$n>1$. \medskip

\textbf{(e)} \fbox{1} We have $d\left(  \sigma\cdot f\right)  =\sigma
\cdot\left(  df\right)  $ for any $f\in\mathcal{P}$ and any $\sigma\in S_{N}$.
\medskip

\textbf{(f)} \fbox{1} We have $df\in\mathcal{S}$ for each $f\in\mathcal{S}$.
\medskip

\textbf{(g)} \fbox{2} We have $d\left(  a_{\rho}\right)  =0$. (Here, $a_{\rho
}$ is as in Definition \ref{def.sf.alternants}.) \medskip

\textbf{(h)} \fbox{5} Let $\lambda=\left(  \lambda_{1},\lambda_{2}%
,\ldots,\lambda_{N}\right)  $ be an $N$-partition. Set $\lambda_{N+1}:=0$. Let
$D$ be the set of all $i\in\left[  N\right]  $ such that $\lambda_{i}%
>\lambda_{i+1}$. For each $i\in D$, let $\lambda-e_{i}$ denote the
$N$-partition obtained from $\lambda$ by subtracting $1$ from the $i$-th
entry; that is, we let%
\[
\lambda-e_{i}:=\left(  \lambda_{1},\lambda_{2},\ldots,\lambda_{i-1}%
,\lambda_{i}-1,\lambda_{i+1},\lambda_{i+2},\ldots,\lambda_{N}\right)  .
\]
(This is an $N$-partition, since $i\in D$ entails $\lambda_{i}>\lambda_{i+1}$
and thus $\lambda_{i}-1\geq\lambda_{i+1}$.)

Then,%
\begin{equation}
d\left(  s_{\lambda}\right)  =\sum_{i\in D}\left(  \lambda_{i}+N-i\right)
s_{\lambda-e_{i}}. \label{eq.exe.sf.differential-sum.dsl}%
\end{equation}

[\textbf{Hint:} For part \textbf{(h)}, study how $d$ acts on the alternant
$a_{\lambda+\rho}$, and show that $d\left(  a_{\rho}f\right)  =a_{\rho
}d\left(  f\right)  $ for any $f\in\mathcal{P}$.] \medskip

[\textbf{Remark:} The equality (\ref{eq.exe.sf.differential-sum.dsl}) can be
restated in terms of Young diagrams as follows:%
\[
d\left(  s_{\lambda}\right)  =\sum k_{\lambda,\mu}s_{\mu},
\]
where

\begin{itemize}
\item the sum ranges over all $N$-partitions $\mu$ such that the Young diagram
$Y\left(  \mu\right)  $ can be obtained from $Y\left(  \lambda\right)  $ by
removing a single box (without shifting the remaining boxes);

\item the coefficient $k_{\lambda,\mu}$ is defined to be $j+N-i$ if $\left(
i,j\right)  $ is the box that must be removed from $Y\left(  \lambda\right)  $
to obtain $Y\left(  \mu\right)  $.
\end{itemize}

\noindent For example, for $N=4$ and $\lambda=\left(  3,1,1,0\right)  $, this
says that%
\begin{align*}
d\left(  s_{\left(  3,1,1,0\right)  }\right)   &  =\underbrace{\left(
\lambda_{1}+N-1\right)  }_{=3+4-1=6}s_{\left(  2,1,1,0\right)  }%
+\underbrace{\left(  \lambda_{3}+N-3\right)  }_{=1+4-3=2}s_{\left(
3,1,0,0\right)  }\\
&  =6s_{\left(  2,1,1,0\right)  }+2s_{\left(  3,1,0,0\right)  }.\text{ ]}%
\end{align*}

\end{exercise}

\begin{exercise}
\label{exe.sf.petrie.2}\fbox{7} \textbf{(a)} Let $n\in\mathbb{N}$, and let $k$
be a positive integer. Recall the $\left(  n,k\right)  $-Petrie symmetric
polynomial $g_{k,n}$ defined in Exercise \ref{exe.sf.petrie.1}. Show that
$g_{k,n}$ can be written in the form%
\[
g_{k,n}=\sum\limits_{\lambda\vdash n}d_{\lambda}s_{\lambda}%
\ \ \ \ \ \ \ \ \ \ \text{for some }d_{\lambda}\in\left\{  0,1,-1\right\}  .
\]
(The notation \textquotedblleft$\lambda\vdash n$\textquotedblright\ means
\textquotedblleft$\lambda$ is a partition of $n$\textquotedblright, so that
the sum ranges over all partitions of $n$. In this case, we furthermore
restrict it to the $N$-partitions of $n$.) \medskip

\textbf{(b)} Show that the numbers $d_{\lambda}$ in this expression are
explicitly given by%
\[
d_{\lambda}=\det\left(  \left(  \left[  0\leq\lambda_{i}-i+j<k\right]
\right)  _{1\leq i\leq M,\ 1\leq j\leq M}\right)  ,
\]
where the partition $\lambda$ has been written as $\lambda=\left(  \lambda
_{1},\lambda_{2},\ldots,\lambda_{M}\right)  $ for some $M\in\mathbb{N}$, and
where we are using the Iverson bracket notation (see Definition
\ref{def.pars.iverson}).
\end{exercise}

\begin{exercise}
Let $I$ be the ideal of the ring $\mathcal{P}$ generated by the $N$
polynomials $e_{1},e_{2},\ldots,e_{N}$. \medskip

\textbf{(a)} \fbox{2} Prove that any homogeneous symmetric polynomial
$f\in\mathcal{S}$ of positive degree is contained in $I$. \medskip

\textbf{(b)} \fbox{9} Prove that the quotient ring $\mathcal{P}/I$ is a free
$K$-module with basis
\[
\left(  \overline{x_{1}^{a_{1}}x_{2}^{a_{2}}\cdots x_{N}^{a_{N}}}\right)
_{\left(  a_{1},a_{2},\ldots,a_{N}\right)  \in H_{N}},
\]
where the $N!$-element set $H_{N}$ is defined as in Definition
\ref{def.perm.lehmer1} \textbf{(c)} (for $n=N$). (Here, the notation
$\overline{f}$ denotes the projection of a polynomial $f\in\mathcal{P}$ onto
the quotient ring $\mathcal{P}/I$.) \medskip

\textbf{(c)} \fbox{5} More generally: Let $u_{1},u_{2},\ldots,u_{N}$ be $N$
polynomials in $\mathcal{P}$ such that
\[
\deg u_{i}<i\ \ \ \ \ \ \ \ \ \ \text{for each }i\in\left[  N\right]  .
\]
(The polynomials $u_{i}$ need not be symmetric or homogeneous.) Let
$I^{\prime}$ be the ideal of the ring $\mathcal{P}$ generated by the $N$
polynomials
\[
e_{1}-u_{1},\ \ e_{2}-u_{2},\ \ \ldots,\ \ e_{N}-u_{N}.
\]
Prove that the quotient ring $\mathcal{P}/I^{\prime}$ is a free $K$-module
with basis
\[
\left(  \overline{x_{1}^{a_{1}}x_{2}^{a_{2}}\cdots x_{N}^{a_{N}}}\right)
_{\left(  a_{1},a_{2},\ldots,a_{N}\right)  \in H_{N}},
\]
where the $N!$-element set $H_{N}$ is defined as in Definition
\ref{def.perm.lehmer1} \textbf{(c)} (for $n=N$). (Here, the notation
$\overline{f}$ denotes the projection of a polynomial $f\in\mathcal{P}$ onto
the quotient ring $\mathcal{P}/I^{\prime}$.)
\end{exercise}

\newpage

\section{\label{chp.details}Omitted details and proofs}

This chapter contains some proofs (and parts of proofs) that have been omitted
from the text above -- usually because they are technical arguments or of
tangential interest only.

Each of the proofs given below uses the notations and conventions of the
chapter and section in which the respective claim appears.

\subsection{\label{sec.details.gf.xneq}$x^{n}$-equivalence}

\begin{fineprint}
\begin{proof}
[Detailed proof of Theorem \ref{thm.fps.xneq.props}.]\textbf{(a)} We claim the following:

\begin{statement}
\textit{Claim 1:} We have $f\overset{x^{n}}{\equiv}f$ for each $f\in K\left[
\left[  x\right]  \right]  $.
\end{statement}

[\textit{Proof of Claim 1:} Let $f\in K\left[  \left[  x\right]  \right]  $.
Obviously, each $m\in\left\{  0,1,\ldots,n\right\}  $ satisfies $\left[
x^{m}\right]  f=\left[  x^{m}\right]  f$. In other words, we have
$f\overset{x^{n}}{\equiv}f$ (by Definition \ref{def.fps.xneq}). This proves
Claim 1.]

\begin{statement}
\textit{Claim 2:} If three FPSs $f,g,h\in K\left[  \left[  x\right]  \right]
$ satisfy $f\overset{x^{n}}{\equiv}g$ and $g\overset{x^{n}}{\equiv}h$, then
$f\overset{x^{n}}{\equiv}h$.
\end{statement}

[\textit{Proof of Claim 2:} Let $f,g,h\in K\left[  \left[  x\right]  \right]
$ be three FPSs satisfying $f\overset{x^{n}}{\equiv}g$ and $g\overset{x^{n}%
}{\equiv}h$. We must show that $f\overset{x^{n}}{\equiv}h$.

We have $f\overset{x^{n}}{\equiv}g$. In other words,
\begin{equation}
\text{each }m\in\left\{  0,1,\ldots,n\right\}  \text{ satisfies }\left[
x^{m}\right]  f=\left[  x^{m}\right]  g
\label{pf.thm.fps.xneq.props.a.c2.pf.1}%
\end{equation}
(by Definition \ref{def.fps.xneq}).

We have $g\overset{x^{n}}{\equiv}h$. In other words,
\begin{equation}
\text{each }m\in\left\{  0,1,\ldots,n\right\}  \text{ satisfies }\left[
x^{m}\right]  g=\left[  x^{m}\right]  h
\label{pf.thm.fps.xneq.props.a.c2.pf.2}%
\end{equation}
(by Definition \ref{def.fps.xneq}).

Now, each $m\in\left\{  0,1,\ldots,n\right\}  $ satisfies%
\begin{align*}
\left[  x^{m}\right]  f  &  =\left[  x^{m}\right]
g\ \ \ \ \ \ \ \ \ \ \left(  \text{by (\ref{pf.thm.fps.xneq.props.a.c2.pf.1}%
)}\right) \\
&  =\left[  x^{m}\right]  h\ \ \ \ \ \ \ \ \ \ \left(  \text{by
(\ref{pf.thm.fps.xneq.props.a.c2.pf.2})}\right)  .
\end{align*}
In other words, we have $f\overset{x^{n}}{\equiv}h$ (by Definition
\ref{def.fps.xneq}). This proves Claim 2.]

\begin{statement}
\textit{Claim 3:} If two FPSs $f,g\in K\left[  \left[  x\right]  \right]  $
satisfy $f\overset{x^{n}}{\equiv}g$, then $g\overset{x^{n}}{\equiv}f$.
\end{statement}

[\textit{Proof of Claim 3:} Let $f,g\in K\left[  \left[  x\right]  \right]  $
be two FPSs satisfying $f\overset{x^{n}}{\equiv}g$. We must show that
$g\overset{x^{n}}{\equiv}f$.

We have $f\overset{x^{n}}{\equiv}g$. In other words,
\[
\text{each }m\in\left\{  0,1,\ldots,n\right\}  \text{ satisfies }\left[
x^{m}\right]  f=\left[  x^{m}\right]  g
\]
(by Definition \ref{def.fps.xneq}). In other words,
\[
\text{each }m\in\left\{  0,1,\ldots,n\right\}  \text{ satisfies }\left[
x^{m}\right]  g=\left[  x^{m}\right]  f.
\]
In other words, we have $g\overset{x^{n}}{\equiv}f$ (by Definition
\ref{def.fps.xneq}). This proves Claim 3.]

Now, the relation $\overset{x^{n}}{\equiv}$ is reflexive (by Claim 1),
transitive (by Claim 2) and symmetric (by Claim 3). In other words, this
relation $\overset{x^{n}}{\equiv}$ is an equivalence relation. This proves
Theorem \ref{thm.fps.xneq.props} \textbf{(a)}. \medskip

\textbf{(b)} Let $a,b,c,d\in K\left[  \left[  x\right]  \right]  $ be four
FPSs satisfying $a\overset{x^{n}}{\equiv}b$ and $c\overset{x^{n}}{\equiv}d$.

We have $a\overset{x^{n}}{\equiv}b$. In other words,
\begin{equation}
\text{each }m\in\left\{  0,1,\ldots,n\right\}  \text{ satisfies }\left[
x^{m}\right]  a=\left[  x^{m}\right]  b \label{pf.thm.fps.xneq.props.b.ab}%
\end{equation}
(by Definition \ref{def.fps.xneq}).

We have $c\overset{x^{n}}{\equiv}d$. In other words,
\begin{equation}
\text{each }m\in\left\{  0,1,\ldots,n\right\}  \text{ satisfies }\left[
x^{m}\right]  c=\left[  x^{m}\right]  d \label{pf.thm.fps.xneq.props.b.cd}%
\end{equation}
(by Definition \ref{def.fps.xneq}).

Now, every $m\in\left\{  0,1,\ldots,n\right\}  $ satisfies%
\begin{equation}
\left[  x^{m}\right]  \left(  a+c\right)  =\left[  x^{m}\right]  a+\left[
x^{m}\right]  c \label{pf.thm.fps.xneq.props.b.+1}%
\end{equation}
(by (\ref{pf.thm.fps.ring.xn(a+b)=}), applied to $m$, $a$ and $c$ instead of
$n$, $\mathbf{a}$ and $\mathbf{b}$) and%
\begin{equation}
\left[  x^{m}\right]  \left(  b+d\right)  =\left[  x^{m}\right]  b+\left[
x^{m}\right]  d \label{pf.thm.fps.xneq.props.b.+2}%
\end{equation}
(by (\ref{pf.thm.fps.ring.xn(a+b)=}), applied to $m$, $b$ and $d$ instead of
$n$, $\mathbf{a}$ and $\mathbf{b}$). Hence, every $m\in\left\{  0,1,\ldots
,n\right\}  $ satisfies%
\begin{align*}
\left[  x^{m}\right]  \left(  a+c\right)   &  =\underbrace{\left[
x^{m}\right]  a}_{\substack{=\left[  x^{m}\right]  b\\\text{(by
(\ref{pf.thm.fps.xneq.props.b.ab}))}}}+\underbrace{\left[  x^{m}\right]
c}_{\substack{=\left[  x^{m}\right]  d\\\text{(by
(\ref{pf.thm.fps.xneq.props.b.cd}))}}}\ \ \ \ \ \ \ \ \ \ \left(  \text{by
(\ref{pf.thm.fps.xneq.props.b.+1})}\right) \\
&  =\left[  x^{m}\right]  b+\left[  x^{m}\right]  d=\left[  x^{m}\right]
\left(  b+d\right)  \ \ \ \ \ \ \ \ \ \ \left(  \text{by
(\ref{pf.thm.fps.xneq.props.b.+2})}\right)  .
\end{align*}
In other words, $a+c\overset{x^{n}}{\equiv}b+d$ (by Definition
\ref{def.fps.xneq}). Thus, we have proved (\ref{eq.thm.fps.xneq.props.b.+}).
The same argument (but with all \textquotedblleft$+$\textquotedblright\ signs
replaced by \textquotedblleft$-$\textquotedblright\ signs, and with all
references to (\ref{pf.thm.fps.ring.xn(a+b)=}) replaced by references to
(\ref{pf.thm.fps.ring.xn(a-b)=})) can be used to prove
(\ref{eq.thm.fps.xneq.props.b.-}). It remains to prove
(\ref{eq.thm.fps.xneq.props.b.*}).

Let $m\in\left\{  0,1,\ldots,n\right\}  $. Then, $m\leq n$.

Now, let $i\in\left\{  0,1,\ldots,m\right\}  $. Then, $i\in\left\{
0,1,\ldots,m\right\}  \subseteq\left\{  0,1,\ldots,n\right\}  $ (since $m\leq
n$). Hence, (\ref{pf.thm.fps.xneq.props.b.ab}) (applied to $i$ instead of $m$)
yields $\left[  x^{i}\right]  a=\left[  x^{i}\right]  b$. Furthermore, from
$i\in\left\{  0,1,\ldots,m\right\}  $, we obtain $m-i\in\left\{
0,1,\ldots,m\right\}  \subseteq\left\{  0,1,\ldots,n\right\}  $ (since $m\leq
n$). Hence, (\ref{pf.thm.fps.xneq.props.b.cd}) (applied to $m-i$ instead of
$m$) yields $\left[  x^{m-i}\right]  c=\left[  x^{m-i}\right]  d$. Multiplying
the equalities $\left[  x^{i}\right]  a=\left[  x^{i}\right]  b$ and $\left[
x^{m-i}\right]  c=\left[  x^{m-i}\right]  d$, we obtain%
\begin{equation}
\left[  x^{i}\right]  a\cdot\left[  x^{m-i}\right]  c=\left[  x^{i}\right]
b\cdot\left[  x^{m-i}\right]  d. \label{pf.thm.fps.xneq.props.b.prod}%
\end{equation}

Forget that we fixed $i$. We thus have proved
(\ref{pf.thm.fps.xneq.props.b.prod}) for each $i\in\left\{  0,1,\ldots
,m\right\}  $. Now, (\ref{pf.thm.fps.ring.xn(ab)=2}) (applied to $m$, $a$ and
$c$ instead of $n$, $\mathbf{a}$ and $\mathbf{b}$) yields%
\[
\left[  x^{m}\right]  \left(  ac\right)  =\sum_{i=0}^{m}\underbrace{\left[
x^{i}\right]  a\cdot\left[  x^{m-i}\right]  c}_{\substack{=\left[
x^{i}\right]  b\cdot\left[  x^{m-i}\right]  d\\\text{(by
(\ref{pf.thm.fps.xneq.props.b.prod}))}}}=\sum_{i=0}^{m}\left[  x^{i}\right]
b\cdot\left[  x^{m-i}\right]  d.
\]
On the other hand, (\ref{pf.thm.fps.ring.xn(ab)=2}) (applied to $m$, $b$ and
$d$ instead of $n$, $\mathbf{a}$ and $\mathbf{b}$) yields%
\[
\left[  x^{m}\right]  \left(  bd\right)  =\sum_{i=0}^{m}\left[  x^{i}\right]
b\cdot\left[  x^{m-i}\right]  d.
\]
Comparing these two equalities, we find $\left[  x^{m}\right]  \left(
ac\right)  =\left[  x^{m}\right]  \left(  bd\right)  $.

Forget that we fixed $m$. We thus have shown that each $m\in\left\{
0,1,\ldots,n\right\}  $ satisfies $\left[  x^{m}\right]  \left(  ac\right)
=\left[  x^{m}\right]  \left(  bd\right)  $. In other words, $ac\overset{x^{n}%
}{\equiv}bd$ (by Definition \ref{def.fps.xneq}). Thus we have proved
(\ref{eq.thm.fps.xneq.props.b.*}), so that the proof of Theorem
\ref{thm.fps.xneq.props} \textbf{(b)} is complete. \medskip

\textbf{(c)} Let $a,b\in K\left[  \left[  x\right]  \right]  $ be two FPSs
satisfying $a\overset{x^{n}}{\equiv}b$. We have $a\overset{x^{n}}{\equiv}b$.
In other words,
\begin{equation}
\text{each }m\in\left\{  0,1,\ldots,n\right\}  \text{ satisfies }\left[
x^{m}\right]  a=\left[  x^{m}\right]  b \label{pf.thm.fps.xneq.props.c.ab}%
\end{equation}
(by Definition \ref{def.fps.xneq}). However, each $m\in\left\{  0,1,\ldots
,n\right\}  $ satisfies
\begin{equation}
\left[  x^{m}\right]  \left(  \lambda a\right)  =\lambda\cdot\left[
x^{m}\right]  a \label{pf.thm.fps.xneq.props.c.1}%
\end{equation}
(by (\ref{pf.thm.fps.ring.xn(la)=}), applied to $m$ and $a$ instead of $n$ and
$\mathbf{a}$) and%
\begin{equation}
\left[  x^{m}\right]  \left(  \lambda b\right)  =\lambda\cdot\left[
x^{m}\right]  b \label{pf.thm.fps.xneq.props.c.2}%
\end{equation}
(by (\ref{pf.thm.fps.ring.xn(la)=}), applied to $m$ and $b$ instead of $n$ and
$\mathbf{a}$). Now, each $m\in\left\{  0,1,\ldots,n\right\}  $ satisfies
\begin{align*}
\left[  x^{m}\right]  \left(  \lambda a\right)   &  =\lambda\cdot
\underbrace{\left[  x^{m}\right]  a}_{\substack{=\left[  x^{m}\right]
b\\\text{(by (\ref{pf.thm.fps.xneq.props.c.ab}))}}}\ \ \ \ \ \ \ \ \ \ \left(
\text{by (\ref{pf.thm.fps.xneq.props.c.1})}\right) \\
&  =\lambda\cdot\left[  x^{m}\right]  b=\left[  x^{m}\right]  \left(  \lambda
b\right)  \ \ \ \ \ \ \ \ \ \ \left(  \text{by
(\ref{pf.thm.fps.xneq.props.c.2})}\right)  .
\end{align*}
In other words, $\lambda a\overset{x^{n}}{\equiv}\lambda b$ (by Definition
\ref{def.fps.xneq}). This proves Theorem \ref{thm.fps.xneq.props}
\textbf{(c)}. \medskip

\textbf{(d)} Let $a,b\in K\left[  \left[  x\right]  \right]  $ be two FPSs
satisfying $a\overset{x^{n}}{\equiv}b$. We have $a\overset{x^{n}}{\equiv}b$.
In other words,
\begin{equation}
\text{each }m\in\left\{  0,1,\ldots,n\right\}  \text{ satisfies }\left[
x^{m}\right]  a=\left[  x^{m}\right]  b \label{pf.thm.fps.xneq.props.d.ab}%
\end{equation}
(by Definition \ref{def.fps.xneq}).

Now, we shall show that each $m\in\left\{  0,1,\ldots,n\right\}  $ satisfies
\begin{equation}
\left[  x^{m}\right]  \left(  a^{-1}\right)  =\left[  x^{m}\right]  \left(
b^{-1}\right)  . \label{pf.thm.fps.xneq.props.d.goal}%
\end{equation}

[\textit{Proof of (\ref{pf.thm.fps.xneq.props.d.goal}):} We shall prove
(\ref{pf.thm.fps.xneq.props.d.goal}) by strong induction on $m$:

\textit{Induction step:} Fix some $k\in\left\{  0,1,\ldots,n\right\}  $. We
assume (as an induction hypothesis) that (\ref{pf.thm.fps.xneq.props.d.goal})
is true for any $m<k$. In other words, for any $m\in\left\{  0,1,\ldots
,n\right\}  $ satisfying $m<k$, we have
\begin{equation}
\left[  x^{m}\right]  \left(  a^{-1}\right)  =\left[  x^{m}\right]  \left(
b^{-1}\right)  . \label{pf.thm.fps.xneq.props.d.goal.pf.IH}%
\end{equation}
We shall now prove that (\ref{pf.thm.fps.xneq.props.d.goal}) is true for
$m=k$. In other words, we shall prove that $\left[  x^{k}\right]  \left(
a^{-1}\right)  =\left[  x^{k}\right]  \left(  b^{-1}\right)  $.

Proposition \ref{prop.fps.invertible} shows that the FPS $a$ is invertible in
$K\left[  \left[  x\right]  \right]  $ if and only if its constant term
$\left[  x^{0}\right]  a$ is invertible in $K$. Hence, its constant term
$\left[  x^{0}\right]  a$ is invertible in $K$ (since $a$ is invertible in
$K\left[  \left[  x\right]  \right]  $). Note that $k\leq n$ (since
$k\in\left\{  0,1,\ldots,n\right\}  $).

Applying (\ref{pf.thm.fps.ring.xn(ab)=2}) to $k$, $a$ and $a^{-1}$ instead of
$n$, $\mathbf{a}$ and $\mathbf{b}$, we obtain%
\begin{align*}
\left[  x^{k}\right]  \left(  aa^{-1}\right)   &  =\sum_{i=0}^{k}\left[
x^{i}\right]  a\cdot\left[  x^{k-i}\right]  \left(  a^{-1}\right) \\
&  =\left[  x^{0}\right]  a\cdot\left[  x^{k}\right]  \left(  a^{-1}\right)
+\sum_{i=1}^{k}\left[  x^{i}\right]  a\cdot\left[  x^{k-i}\right]  \left(
a^{-1}\right)
\end{align*}
(here, we have split off the addend for $i=0$ from the sum). Thus,%
\[
\left[  x^{0}\right]  a\cdot\left[  x^{k}\right]  \left(  a^{-1}\right)
+\sum_{i=1}^{k}\left[  x^{i}\right]  a\cdot\left[  x^{k-i}\right]  \left(
a^{-1}\right)  =\left[  x^{k}\right]  \underbrace{\left(  aa^{-1}\right)
}_{=1}=\left[  x^{k}\right]  1,
\]
so that%
\[
\left[  x^{0}\right]  a\cdot\left[  x^{k}\right]  \left(  a^{-1}\right)
=\left[  x^{k}\right]  1-\sum_{i=1}^{k}\left[  x^{i}\right]  a\cdot\left[
x^{k-i}\right]  \left(  a^{-1}\right)  .
\]
We can divide this equality by $\left[  x^{0}\right]  a$ (since $\left[
x^{0}\right]  a$ is invertible in $K$), and thus obtain%
\begin{equation}
\left[  x^{k}\right]  \left(  a^{-1}\right)  =\left(  \left[  x^{0}\right]
a\right)  ^{-1}\cdot\left(  \left[  x^{k}\right]  1-\sum_{i=1}^{k}\left[
x^{i}\right]  a\cdot\left[  x^{k-i}\right]  \left(  a^{-1}\right)  \right)  .
\label{pf.thm.fps.xneq.props.d.goal.pf.xka-1}%
\end{equation}
The same argument (applied to $b$ instead of $a$) yields%
\begin{equation}
\left[  x^{k}\right]  \left(  b^{-1}\right)  =\left(  \left[  x^{0}\right]
b\right)  ^{-1}\cdot\left(  \left[  x^{k}\right]  1-\sum_{i=1}^{k}\left[
x^{i}\right]  b\cdot\left[  x^{k-i}\right]  \left(  b^{-1}\right)  \right)  .
\label{pf.thm.fps.xneq.props.d.goal.pf.xkb-1}%
\end{equation}

However, we observe the following:

\begin{itemize}
\item We have $0\in\left\{  0,1,\ldots,n\right\}  $ (since $n\in\mathbb{N}$)
and thus $\left[  x^{0}\right]  a=\left[  x^{0}\right]  b$ (by
(\ref{pf.thm.fps.xneq.props.d.ab}), applied to $m=0$).

\item Each $i\in\left\{  1,2,\ldots,k\right\}  $ satisfies $i\in\left\{
1,2,\ldots,k\right\}  \subseteq\left\{  0,1,\ldots,n\right\}  $ (since
$1\geq0$ and $k\leq n$) and therefore
\begin{equation}
\left[  x^{i}\right]  a=\left[  x^{i}\right]  b
\label{pf.thm.fps.xneq.props.d.goal.pf.obs2}%
\end{equation}
(by (\ref{pf.thm.fps.xneq.props.d.ab}), applied to $m=i$).

\item Each $i\in\left\{  1,2,\ldots,k\right\}  $ satisfies $k-i\in\left\{
0,1,\ldots,k-1\right\}  \subseteq\left\{  0,1,\ldots,n\right\}  $ (since
$k-1\leq k\leq n$) and $k-i<k$ (since $k-i\in\left\{  0,1,\ldots,k-1\right\}
$, so that $k-i\leq k-1<k$), and therefore
\begin{equation}
\left[  x^{k-i}\right]  \left(  a^{-1}\right)  =\left[  x^{k-i}\right]
\left(  b^{-1}\right)  \label{pf.thm.fps.xneq.props.d.goal.pf.obs3}%
\end{equation}
(by (\ref{pf.thm.fps.xneq.props.d.goal.pf.IH}), applied to $m=k-i$).
\end{itemize}

Hence, (\ref{pf.thm.fps.xneq.props.d.goal.pf.xka-1}) becomes%
\begin{align*}
\left[  x^{k}\right]  \left(  a^{-1}\right)   &  =\left(  \underbrace{\left[
x^{0}\right]  a}_{=\left[  x^{0}\right]  b}\right)  ^{-1}\cdot\left(  \left[
x^{k}\right]  1-\sum_{i=1}^{k}\underbrace{\left[  x^{i}\right]  a}%
_{\substack{=\left[  x^{i}\right]  b\\\text{(by
(\ref{pf.thm.fps.xneq.props.d.goal.pf.obs2}))}}}\cdot\underbrace{\left[
x^{k-i}\right]  \left(  a^{-1}\right)  }_{\substack{=\left[  x^{k-i}\right]
\left(  b^{-1}\right)  \\\text{(by (\ref{pf.thm.fps.xneq.props.d.goal.pf.obs3}%
))}}}\right) \\
&  =\left(  \left[  x^{0}\right]  b\right)  ^{-1}\cdot\left(  \left[
x^{k}\right]  1-\sum_{i=1}^{k}\left[  x^{i}\right]  b\cdot\left[
x^{k-i}\right]  \left(  b^{-1}\right)  \right) \\
&  =\left[  x^{k}\right]  \left(  b^{-1}\right)  \ \ \ \ \ \ \ \ \ \ \left(
\text{by (\ref{pf.thm.fps.xneq.props.d.goal.pf.xkb-1})}\right)  .
\end{align*}
In other words, (\ref{pf.thm.fps.xneq.props.d.goal}) is true for $m=k$. This
completes the induction step. Thus, (\ref{pf.thm.fps.xneq.props.d.goal}) is proved.]

Thus, we have shown that each $m\in\left\{  0,1,\ldots,n\right\}  $ satisfies
$\left[  x^{m}\right]  \left(  a^{-1}\right)  =\left[  x^{m}\right]  \left(
b^{-1}\right)  $. In other words, $a^{-1}\overset{x^{n}}{\equiv}b^{-1}$ (by
Definition \ref{def.fps.xneq}). This proves Theorem \ref{thm.fps.xneq.props}
\textbf{(d)}. \medskip

\textbf{(e)} Let $a,b,c,d\in K\left[  \left[  x\right]  \right]  $ be four
FPSs satisfying $a\overset{x^{n}}{\equiv}b$ and $c\overset{x^{n}}{\equiv}d$.
Assume that the FPSs $c$ and $d$ are invertible. Then, Theorem
\ref{thm.fps.xneq.props} \textbf{(d)} (applied to $c$ and $d$ instead of $a$
and $b$) yields $c^{-1}\overset{x^{n}}{\equiv}d^{-1}$. Hence,
(\ref{eq.thm.fps.xneq.props.b.*}) (applied to $c^{-1}$ and $d^{-1}$ instead of
$c$ and $d$) yields $ac^{-1}\overset{x^{n}}{\equiv}bd^{-1}$. In other words,
$\dfrac{a}{c}\overset{x^{n}}{\equiv}\dfrac{b}{d}$ (since $ac^{-1}=\dfrac{a}%
{c}$ and $bd^{-1}=\dfrac{b}{d}$). This proves Theorem \ref{thm.fps.xneq.props}
\textbf{(e)}. \medskip

\textbf{(f)} Let us first prove (\ref{eq.thm.fps.xneq.props.e.+}):

[\textit{Proof of (\ref{eq.thm.fps.xneq.props.e.+}):} We proceed by induction
on $\left\vert S\right\vert $:

\textit{Induction base:} It is easy to see that
(\ref{eq.thm.fps.xneq.props.e.+}) holds for $\left\vert S\right\vert
=0$\ \ \ \ \footnote{\textit{Proof.} Let $S$, $\left(  a_{s}\right)  _{s\in
S}$ and $\left(  b_{s}\right)  _{s\in S}$ be as in Theorem
\ref{thm.fps.xneq.props} \textbf{(f)}, and assume that $\left\vert
S\right\vert =0$. From $\left\vert S\right\vert =0$, we obtain $S=\varnothing
$. Hence,%
\begin{align*}
\sum_{s\in S}a_{s}  &  =\left(  \text{empty sum}\right)
=0\ \ \ \ \ \ \ \ \ \ \text{and}\\
\sum_{s\in S}b_{s}  &  =\left(  \text{empty sum}\right)  =0.
\end{align*}
Comparing these two equalities, we obtain $\sum_{s\in S}a_{s}=\sum_{s\in
S}b_{s}$. Hence, $\sum_{s\in S}a_{s}\overset{x^{n}}{\equiv}\sum_{s\in S}b_{s}$
(since the relation $\overset{x^{n}}{\equiv}$ is an equivalence relation and
thus is reflexive). Thus, we have proved (\ref{eq.thm.fps.xneq.props.e.+})
under the assumption that $\left\vert S\right\vert =0$.}.

\textit{Induction step:} Let $k\in\mathbb{N}$. Assume (as the induction
hypothesis) that (\ref{eq.thm.fps.xneq.props.e.+}) holds for $\left\vert
S\right\vert =k$. We must prove that (\ref{eq.thm.fps.xneq.props.e.+}) holds
for $\left\vert S\right\vert =k+1$.

So let $S$, $\left(  a_{s}\right)  _{s\in S}$ and $\left(  b_{s}\right)
_{s\in S}$ be as in Theorem \ref{thm.fps.xneq.props} \textbf{(f)}, and assume
that $\left\vert S\right\vert =k+1$. Then, $\left\vert S\right\vert
=k+1>k\geq0$, so that the set $S$ is nonempty. In other words, there exists
some $t\in S$. Consider this $t$. Each $s\in S\setminus\left\{  t\right\}  $
satisfies $s\in S\setminus\left\{  t\right\}  \subseteq S$ and thus
$a_{s}\overset{x^{n}}{\equiv}b_{s}$ (by (\ref{eq.thm.fps.xneq.props.e.ass})).
Moreover, from $t\in S$, we obtain $\left\vert S\setminus\left\{  t\right\}
\right\vert =\left\vert S\right\vert -1=k$ (since $\left\vert S\right\vert
=k+1$). Hence, we can apply (\ref{eq.thm.fps.xneq.props.e.+}) to
$S\setminus\left\{  t\right\}  $ instead of $S$ (since our induction
hypothesis says that (\ref{eq.thm.fps.xneq.props.e.+}) holds for $\left\vert
S\right\vert =k$). As a result, we obtain%
\[
\sum_{s\in S\setminus\left\{  t\right\}  }a_{s}\overset{x^{n}}{\equiv}%
\sum_{s\in S\setminus\left\{  t\right\}  }b_{s}.
\]
On the other hand, $a_{t}\overset{x^{n}}{\equiv}b_{t}$ (by
(\ref{eq.thm.fps.xneq.props.e.ass}), applied to $s=t$). Hence,
(\ref{eq.thm.fps.xneq.props.b.+}) (applied to $a=a_{t}$ and $b=b_{t}$ and
$c=\sum_{s\in S\setminus\left\{  t\right\}  }a_{s}$ and $d=\sum_{s\in
S\setminus\left\{  t\right\}  }b_{s}$) yields%
\[
a_{t}+\sum_{s\in S\setminus\left\{  t\right\}  }a_{s}\overset{x^{n}}{\equiv
}b_{t}+\sum_{s\in S\setminus\left\{  t\right\}  }b_{s}.
\]
In view of%
\[
\sum_{s\in S}a_{s}=a_{t}+\sum_{s\in S\setminus\left\{  t\right\}  }%
a_{s}\ \ \ \ \ \ \ \ \ \ \left(
\begin{array}
[c]{c}%
\text{here, we have split off the}\\
\text{addend for }s=t\text{ from the sum}%
\end{array}
\right)
\]
and%
\[
\sum_{s\in S}b_{s}=b_{t}+\sum_{s\in S\setminus\left\{  t\right\}  }%
b_{s}\ \ \ \ \ \ \ \ \ \ \left(
\begin{array}
[c]{c}%
\text{here, we have split off the}\\
\text{addend for }s=t\text{ from the sum}%
\end{array}
\right)  ,
\]
this rewrites as%
\[
\sum_{s\in S}a_{s}\overset{x^{n}}{\equiv}\sum_{s\in S}b_{s}.
\]
Hence, we have shown that (\ref{eq.thm.fps.xneq.props.e.+}) holds for
$\left\vert S\right\vert =k+1$. This completes the induction step. Thus, the
induction proof of (\ref{eq.thm.fps.xneq.props.e.+}) is complete.]

We have now proved (\ref{eq.thm.fps.xneq.props.e.+}). The exact same argument
(but with all sums replaced by products, and with the reference to
(\ref{eq.thm.fps.xneq.props.b.+}) replaced by a reference to
(\ref{eq.thm.fps.xneq.props.b.*})) can be used to prove
(\ref{eq.thm.fps.xneq.props.e.*}). Hence, the proof of Theorem
\ref{thm.fps.xneq.props} \textbf{(f)} is complete.
\end{proof}

\begin{proof}
[Detailed proof of Proposition \ref{prop.fps.xneq-multiple}.]For each
$m\in\left\{  0,1,\ldots,n\right\}  $, we have%
\begin{equation}
\left[  x^{m}\right]  \left(  f-g\right)  =\left[  x^{m}\right]  f-\left[
x^{m}\right]  g \label{pf.prop.fps.xneq-multiple.1}%
\end{equation}
(by (\ref{pf.thm.fps.ring.xn(a-b)=}), applied to $m$, $f$ and $g$ instead of
$n$, $\mathbf{a}$ and $\mathbf{b}$).

Lemma \ref{lem.fps.muls-of-xn} (applied to $f-g$ and $n+1$ instead of $f$ and
$k$) shows that the first $n+1$ coefficients of the FPS $f-g$ are $0$ if and
only if $f-g$ is a multiple of $x^{n+1}$.

Now, we have the following chain of logical equivalences:%
\begin{align*}
\left(  f\overset{x^{n}}{\equiv}g\right)  \  &  \Longleftrightarrow\ \left(
\text{each }m\in\left\{  0,1,\ldots,n\right\}  \text{ satisfies }%
\underbrace{\left[  x^{m}\right]  f=\left[  x^{m}\right]  g}%
_{\Longleftrightarrow\ \left(  \left[  x^{m}\right]  f-\left[  x^{m}\right]
g=0\right)  }\right) \\
&  \ \ \ \ \ \ \ \ \ \ \ \ \ \ \ \ \ \ \ \ \left(  \text{by Definition
\ref{def.fps.xneq}}\right) \\
&  \Longleftrightarrow\ \left(  \text{each }m\in\left\{  0,1,\ldots,n\right\}
\text{ satisfies }\underbrace{\left[  x^{m}\right]  f-\left[  x^{m}\right]
g}_{\substack{=\left[  x^{m}\right]  \left(  f-g\right)  \\\text{(by
(\ref{pf.prop.fps.xneq-multiple.1}))}}}=0\right) \\
&  \Longleftrightarrow\ \left(  \text{each }m\in\left\{  0,1,\ldots,n\right\}
\text{ satisfies }\left[  x^{m}\right]  \left(  f-g\right)  =0\right) \\
&  \Longleftrightarrow\ \left(  \text{the first }n+1\text{ coefficients of the
FPS }f-g\text{ are }0\right) \\
&  \Longleftrightarrow\ \left(  f-g\text{ is a multiple of }x^{n+1}\right)
\end{align*}
(since we have seen that the first $n+1$ coefficients of the FPS $f-g$ are $0$
if and only if $f-g$ is a multiple of $x^{n+1}$). In other words, we have
$f\overset{x^{n}}{\equiv}g$ if and only if the FPS $f-g$ is a multiple of
$x^{n+1}$. This proves Proposition \ref{prop.fps.xneq-multiple}.
\end{proof}

\begin{proof}
[Detailed proof of Proposition \ref{prop.fps.xneq.comp}.]Write the FPS $a$ in
the form $a=\sum_{i\in\mathbb{N}}a_{i}x^{i}$ with $a_{0},a_{1},a_{2},\ldots\in
K$. Thus,
\begin{equation}
a_{m}=\left[  x^{m}\right]  a\ \ \ \ \ \ \ \ \ \ \text{for each }%
m\in\mathbb{N} \label{pf.prop.fps.xneq.comp.am=xma}%
\end{equation}
(by the definition of $\left[  x^{m}\right]  a$). Furthermore, Definition
\ref{def.fps.subs} yields%
\begin{equation}
a\circ c=\sum_{i\in\mathbb{N}}a_{i}c^{i} \label{pf.prop.fps.xneq.comp.ac=}%
\end{equation}
(since $a=\sum_{i\in\mathbb{N}}a_{i}x^{i}$ with $a_{0},a_{1},a_{2},\ldots\in
K$).

Write the FPS $b$ in the form $b=\sum_{i\in\mathbb{N}}b_{i}x^{i}$ with
$b_{0},b_{1},b_{2},\ldots\in K$. Thus,
\begin{equation}
b_{m}=\left[  x^{m}\right]  b\ \ \ \ \ \ \ \ \ \ \text{for each }%
m\in\mathbb{N} \label{pf.prop.fps.xneq.comp.bm=xmb}%
\end{equation}
(by the definition of $\left[  x^{m}\right]  b$). Furthermore, Definition
\ref{def.fps.subs} yields%
\begin{equation}
b\circ d=\sum_{i\in\mathbb{N}}b_{i}d^{i} \label{pf.prop.fps.xneq.comp.bd=}%
\end{equation}
(since $b=\sum_{i\in\mathbb{N}}b_{i}x^{i}$ with $b_{0},b_{1},b_{2},\ldots\in
K$).

We have $a\overset{x^{n}}{\equiv}b$. In other words,%
\begin{equation}
\text{each }m\in\left\{  0,1,\ldots,n\right\}  \text{ satisfies }\left[
x^{m}\right]  a=\left[  x^{m}\right]  b \label{pf.prop.fps.xneq.comp.xma=xmb}%
\end{equation}
(by the definition of \textquotedblleft$a\overset{x^{n}}{\equiv}%
b$\textquotedblright). Hence, each $m\in\left\{  0,1,\ldots,n\right\}  $
satisfies%
\begin{align}
a_{m}  &  =\left[  x^{m}\right]  a\ \ \ \ \ \ \ \ \ \ \left(  \text{by
(\ref{pf.prop.fps.xneq.comp.am=xma})}\right) \nonumber\\
&  =\left[  x^{m}\right]  b\ \ \ \ \ \ \ \ \ \ \left(  \text{by
(\ref{pf.prop.fps.xneq.comp.xma=xmb})}\right) \nonumber\\
&  =b_{m}\ \ \ \ \ \ \ \ \ \ \left(  \text{by
(\ref{pf.prop.fps.xneq.comp.bm=xmb})}\right)  .
\label{pf.prop.fps.xneq.comp.am=bm}%
\end{align}

Now, we claim the following:

\begin{statement}
\textit{Claim 1:} Let $i\in\left\{  0,1,\ldots,n\right\}  $. Then, each
$m\in\left\{  0,1,\ldots,n\right\}  $ satisfies%
\[
\left[  x^{m}\right]  \left(  a_{i}c^{i}\right)  =\left[  x^{m}\right]
\left(  b_{i}d^{i}\right)  .
\]

\end{statement}

[\textit{Proof of Claim 1:} Let $S$ be the set $\left\{  1,2,\ldots,i\right\}
$. This set $S$ is finite, and satisfies $\left\vert S\right\vert =i$.
Moreover, we have $c\overset{x^{n}}{\equiv}d$ for each $s\in S$ (by
assumption). Hence, (\ref{eq.thm.fps.xneq.props.e.*}) (applied to $a_{s}=c$
and $b_{s}=d$) yields%
\[
\prod_{s\in S}c\overset{x^{n}}{\equiv}\prod_{s\in S}d.
\]
In view of
\begin{align*}
\prod_{s\in S}c  &  =c^{\left\vert S\right\vert }=c^{i}%
\ \ \ \ \ \ \ \ \ \ \left(  \text{since }\left\vert S\right\vert =i\right)
\ \ \ \ \ \ \ \ \ \ \text{and}\\
\prod_{s\in S}d  &  =d^{\left\vert S\right\vert }=d^{i}%
\ \ \ \ \ \ \ \ \ \ \left(  \text{since }\left\vert S\right\vert =i\right)  ,
\end{align*}
we can rewrite this as $c^{i}\overset{x^{n}}{\equiv}d^{i}$. In other words,
\begin{equation}
\text{each }m\in\left\{  0,1,\ldots,n\right\}  \text{ satisfies }\left[
x^{m}\right]  \left(  c^{i}\right)  =\left[  x^{m}\right]  \left(
d^{i}\right)  \label{pf.prop.fps.xneq.comp.c1.pf.3}%
\end{equation}
(by the definition of \textquotedblleft$c^{i}\overset{x^{n}}{\equiv}d^{i}%
$\textquotedblright).

Now, let $m\in\left\{  0,1,\ldots,n\right\}  $. Then,
(\ref{pf.thm.fps.ring.xn(la)=}) (applied to $m$, $a_{i}$ and $c^{i}$ instead
of $n$, $\lambda$ and $\mathbf{a}$) yields $\left[  x^{m}\right]  \left(
a_{i}c^{i}\right)  =a_{i}\cdot\left[  x^{m}\right]  \left(  c^{i}\right)  $.
Similarly, $\left[  x^{m}\right]  \left(  b_{i}d^{i}\right)  =b_{i}%
\cdot\left[  x^{m}\right]  \left(  d^{i}\right)  $. On the other hand,
(\ref{pf.prop.fps.xneq.comp.am=bm}) (applied to $i$ instead of $m$) yields
$a_{i}=b_{i}$. Hence,
\begin{align*}
\left[  x^{m}\right]  \left(  a_{i}c^{i}\right)   &  =\underbrace{a_{i}%
}_{=b_{i}}\cdot\underbrace{\left[  x^{m}\right]  \left(  c^{i}\right)
}_{\substack{=\left[  x^{m}\right]  \left(  d^{i}\right)  \\\text{(by
(\ref{pf.prop.fps.xneq.comp.c1.pf.3}))}}}=b_{i}\cdot\left[  x^{m}\right]
\left(  d^{i}\right) \\
&  =\left[  x^{m}\right]  \left(  b_{i}d^{i}\right)
\ \ \ \ \ \ \ \ \ \ \left(  \text{since }\left[  x^{m}\right]  \left(
b_{i}d^{i}\right)  =b_{i}\cdot\left[  x^{m}\right]  \left(  d^{i}\right)
\right)  .
\end{align*}
This proves Claim 1.] \medskip

Next, we claim the following:

\begin{statement}
\textit{Claim 2:} Let $m\in\mathbb{N}$. Let $i\in\mathbb{N}\setminus\left\{
0,1,\ldots,m\right\}  $. Then,%
\begin{equation}
\left[  x^{m}\right]  \left(  c^{i}\right)  =0
\label{pf.prop.fps.xneq.comp.c2.c}%
\end{equation}
and%
\begin{equation}
\left[  x^{m}\right]  \left(  d^{i}\right)  =0.
\label{pf.prop.fps.xneq.comp.c2.d}%
\end{equation}

\end{statement}

[\textit{Proof of Claim 2:} We have $i\in\mathbb{N}\setminus\left\{
0,1,\ldots,m\right\}  =\left\{  m+1,m+2,m+3,\ldots\right\}  $, so that $i\geq
m+1$ and therefore $m\leq i-1$. Hence, $m\in\left\{  0,1,\ldots,i-1\right\}  $
(since $m\in\mathbb{N}$).

By assumption, we have $\left[  x^{0}\right]  c=0$. In other words, the $0$-th
coefficient of $c$ is $0$. In other words, the first $1$ coefficient of the
FPS $c$ is $0$. However, Lemma \ref{lem.fps.muls-of-xn} (applied to $k=1$ and
$f=c$) yields that the first $1$ coefficient of the FPS $c$ is $0$ if and only
if $c$ is a multiple of $x^{1}$. Hence, $c$ is a multiple of $x^{1}$ (since
the first $1$ coefficient of the FPS $c$ is $0$). In other words, $c=x^{1}h$
for some $h\in K\left[  \left[  x\right]  \right]  $. Consider this $h$. Now,
$c=x^{1}h=xh$, so that $c^{i}=\left(  xh\right)  ^{i}=x^{i}h^{i}$. However,
Lemma \ref{lem.fps.first-n-coeffs-of-xna} (applied to $i$ and $h^{i}$ instead
of $k$ and $a$) yields that the first $i$ coefficients of the FPS $x^{i}h^{i}$
are $0$. In other words, the first $i$ coefficients of the FPS $c^{i}$ are $0$
(since $c^{i}=x^{i}h^{i}$). In other words, $\left[  x^{j}\right]  \left(
c^{i}\right)  =0$ for all $j\in\left\{  0,1,\ldots,i-1\right\}  $. We can
apply this to $j=m$ (since $m\in\left\{  0,1,\ldots,i-1\right\}  $), and thus
obtain $\left[  x^{m}\right]  \left(  c^{i}\right)  =0$. This proves
(\ref{pf.prop.fps.xneq.comp.c2.c}). The proof of
(\ref{pf.prop.fps.xneq.comp.c2.d}) is analogous (but uses the FPS $d$ instead
of $c$). Thus, Claim 2 is proven.] \medskip

Next, we generalize Claim 1 to all $i\in\mathbb{N}$:

\begin{statement}
\textit{Claim 3:} Let $i\in\mathbb{N}$. Let $m\in\left\{  0,1,\ldots
,n\right\}  $. Then,%
\[
\left[  x^{m}\right]  \left(  a_{i}c^{i}\right)  =\left[  x^{m}\right]
\left(  b_{i}d^{i}\right)  .
\]

\end{statement}

[\textit{Proof of Claim 3:} If $i\in\left\{  0,1,\ldots,m\right\}  $, then
this follows from Claim 1. Thus, for the rest of this proof, we WLOG assume
that we don't have $i\in\left\{  0,1,\ldots,m\right\}  $. Hence,
$i\notin\left\{  0,1,\ldots,m\right\}  $, so that $i\in\mathbb{N}%
\setminus\left\{  0,1,\ldots,m\right\}  $ (since $i\in\mathbb{N}$). Thus,
Claim 2 applies, and therefore (\ref{pf.prop.fps.xneq.comp.c2.c}) and
(\ref{pf.prop.fps.xneq.comp.c2.d}) hold.

Now, (\ref{pf.thm.fps.ring.xn(la)=}) (applied to $m$, $a_{i}$ and $c^{i}$
instead of $n$, $\lambda$ and $\mathbf{a}$) yields $\left[  x^{m}\right]
\left(  a_{i}c^{i}\right)  =a_{i}\cdot\underbrace{\left[  x^{m}\right]
\left(  c^{i}\right)  }_{\substack{=0\\\text{(by
(\ref{pf.prop.fps.xneq.comp.c2.c}))}}}=0$. The same argument (applied to
$b_{i}$ and $d$ instead of $a_{i}$ and $c$) yields $\left[  x^{m}\right]
\left(  b_{i}d^{i}\right)  =0$. Comparing these two equalities, we obtain
$\left[  x^{m}\right]  \left(  a_{i}c^{i}\right)  =\left[  x^{m}\right]
\left(  b_{i}d^{i}\right)  $. This proves Claim 3.] \medskip

Now, each $m\in\left\{  0,1,\ldots,n\right\}  $ satisfies%
\begin{align*}
\left[  x^{m}\right]  \left(  a\circ c\right)   &  =\left[  x^{m}\right]
\left(  \sum_{i\in\mathbb{N}}a_{i}c^{i}\right)  \ \ \ \ \ \ \ \ \ \ \left(
\text{by (\ref{pf.prop.fps.xneq.comp.ac=})}\right) \\
&  =\sum_{i\in\mathbb{N}}\underbrace{\left[  x^{m}\right]  \left(  a_{i}%
c^{i}\right)  }_{\substack{=\left[  x^{m}\right]  \left(  b_{i}d^{i}\right)
\\\text{(by Claim 3)}}}=\sum_{i\in\mathbb{N}}\left[  x^{m}\right]  \left(
b_{i}d^{i}\right)  =\left[  x^{m}\right]  \underbrace{\left(  \sum
_{i\in\mathbb{N}}b_{i}d^{i}\right)  }_{\substack{=b\circ d\\\text{(by
(\ref{pf.prop.fps.xneq.comp.bd=}))}}}=\left[  x^{m}\right]  \left(  b\circ
d\right)  .
\end{align*}
In other words, $a\circ c\overset{x^{n}}{\equiv}b\circ d$ (by the definition
of \textquotedblleft$a\circ c\overset{x^{n}}{\equiv}b\circ d$%
\textquotedblright). This proves Proposition \ref{prop.fps.xneq.comp}.
\end{proof}
\end{fineprint}

\subsection{\label{sec.details.gf.prod}Infinite products}

\begin{fineprint}
\begin{proof}
[Detailed proof of Proposition \ref{prop.fps.multipliable.prod-wd2}.]Let us
(temporarily) use two different notations for our two different definitions of
a product: We let $\prod_{i\in I}\mathbf{a}_{i}$ denote the finite product
$\prod_{i\in I}\mathbf{a}_{i}$ defined in the usual way (i.e., defined as in
any commutative ring), whereas $\widetilde{\prod_{i\in I}}\mathbf{a}_{i}$
shall mean the product $\prod_{i\in I}\mathbf{a}_{i}$ defined according to
Definition \ref{def.fps.multipliable} \textbf{(b)}. Thus, our goal is to show
that $\widetilde{\prod_{i\in I}}\mathbf{a}_{i}=\prod_{i\in I}\mathbf{a}_{i}$.

Definition \ref{def.fps.multipliable} \textbf{(b)} shows that if
$n\in\mathbb{N}$, and if $M$ is a finite subset of $I$ that determines the
$x^{n}$-coefficient in the product of $\left(  \mathbf{a}_{i}\right)  _{i\in
I}$, then
\begin{equation}
\left[  x^{n}\right]  \left(  \widetilde{\prod_{i\in I}}\mathbf{a}_{i}\right)
=\left[  x^{n}\right]  \left(  \prod_{i\in M}\mathbf{a}_{i}\right)  .
\label{pf.prop.fps.multipliable.prod-wd2.1}%
\end{equation}

Now, let $n\in\mathbb{N}$. The set $I$ is finite (since the family $\left(
\mathbf{a}_{i}\right)  _{i\in I}$ is finite), and thus is a finite subset of
$I$. Moreover, every finite subset $J$ of $I$ satisfying $I\subseteq
J\subseteq I$ satisfies%
\[
\left[  x^{n}\right]  \left(  \prod_{i\in J}\mathbf{a}_{i}\right)  =\left[
x^{n}\right]  \left(  \prod_{i\in I}\mathbf{a}_{i}\right)
\]
(because combining $I\subseteq J$ and $J\subseteq I$ yields $J=I$, and thus
$\left[  x^{n}\right]  \left(  \prod_{i\in J}\mathbf{a}_{i}\right)  =\left[
x^{n}\right]  \left(  \prod_{i\in I}\mathbf{a}_{i}\right)  $). In other words,
$I$ determines the $x^{n}$-coefficient in the product of $\left(
\mathbf{a}_{i}\right)  _{i\in I}$ (by the definition of \textquotedblleft
determining the $x^{n}$-coefficient in a product\textquotedblright). Hence, we
can apply (\ref{pf.prop.fps.multipliable.prod-wd2.1}) to $M=I$. As a
consequence, we obtain%
\[
\left[  x^{n}\right]  \left(  \widetilde{\prod_{i\in I}}\mathbf{a}_{i}\right)
=\left[  x^{n}\right]  \left(  \prod_{i\in I}\mathbf{a}_{i}\right)  .
\]

Now, forget that we fixed $n$. We thus have proved%
\[
\left[  x^{n}\right]  \left(  \widetilde{\prod_{i\in I}}\mathbf{a}_{i}\right)
=\left[  x^{n}\right]  \left(  \prod_{i\in I}\mathbf{a}_{i}\right)
\ \ \ \ \ \ \ \ \ \ \text{for each }n\in\mathbb{N}.
\]
In other words, each coefficient of the FPS $\widetilde{\prod_{i\in I}%
}\mathbf{a}_{i}$ equals the corresponding coefficient of $\prod_{i\in
I}\mathbf{a}_{i}$. Hence, $\widetilde{\prod_{i\in I}}\mathbf{a}_{i}%
=\prod_{i\in I}\mathbf{a}_{i}$. This completes the proof of Proposition
\ref{prop.fps.multipliable.prod-wd2}.
\end{proof}
\end{fineprint}

\begin{fineprint}
\begin{proof}
[Detailed proof of Lemma \ref{lem.fps.prod.irlv.fin}.]We shall prove Lemma
\ref{lem.fps.prod.irlv.fin} by induction on $\left\vert J\right\vert $ (this
is allowed, since the set $J$ is supposed to be finite):

\textit{Induction base:} Lemma \ref{lem.fps.prod.irlv.fin} is true in the case
when $\left\vert J\right\vert =0$\ \ \ \ \footnote{\textit{Proof.} Let $a$,
$\left(  f_{i}\right)  _{i\in J}$ and $n$ be as in Lemma
\ref{lem.fps.prod.irlv.fin}. Assume that $\left\vert J\right\vert =0$. Thus,
$J=\varnothing$, so that the product $\prod_{i\in J}\left(  1+f_{i}\right)  $
is empty. Hence, $\prod_{i\in J}\left(  1+f_{i}\right)  =\left(  \text{empty
product}\right)  =1$. Hence, we have%
\[
\left[  x^{m}\right]  \left(  a\underbrace{\prod_{i\in J}\left(
1+f_{i}\right)  }_{=1}\right)  =\left[  x^{m}\right]
a\ \ \ \ \ \ \ \ \ \ \text{for each }m\in\left\{  0,1,\ldots,n\right\}  .
\]
Thus, we have proved Lemma \ref{lem.fps.prod.irlv.fin} under the assumption
that $\left\vert J\right\vert =0$. Therefore, Lemma
\ref{lem.fps.prod.irlv.fin} is true in the case when $\left\vert J\right\vert
=0$.}.

\textit{Induction step:} Let $k\in\mathbb{N}$. Assume (as the induction
hypothesis) that Lemma \ref{lem.fps.prod.irlv.fin} is true in the case when
$\left\vert J\right\vert =k$. We must prove that Lemma
\ref{lem.fps.prod.irlv.fin} is true in the case when $\left\vert J\right\vert
=k+1$.

Let $a$, $\left(  f_{i}\right)  _{i\in J}$ and $n$ be as in Lemma
\ref{lem.fps.prod.irlv.fin}. Assume that $\left\vert J\right\vert =k+1$. Thus,
$\left\vert J\right\vert =k+1>k\geq0$; hence, the set $J$ is nonempty. In
other words, there exists some $j\in J$. Consider this $j$. We have $j\in J$
and therefore
\begin{equation}
\left[  x^{m}\right]  \left(  f_{j}\right)  =0\ \ \ \ \ \ \ \ \ \ \text{for
each }m\in\left\{  0,1,\ldots,n\right\}  \label{pf.lem.fps.prod.irlv.fin.fj}%
\end{equation}
(by (\ref{eq.lem.fps.prod.irlv.fin.ass}), applied to $i=j$). Hence, Lemma
\ref{lem.fps.prod.irlv.1} (applied to $a\prod_{i\in J\setminus\left\{
j\right\}  }\left(  1+f_{i}\right)  $ and $f_{j}$ instead of $a$ and $f$)
yields that%
\begin{align*}
\left[  x^{m}\right]  \left(  \left(  a\prod_{i\in J\setminus\left\{
j\right\}  }\left(  1+f_{i}\right)  \right)  \left(  1+f_{j}\right)  \right)
&  =\left[  x^{m}\right]  \left(  a\prod_{i\in J\setminus\left\{  j\right\}
}\left(  1+f_{i}\right)  \right) \\
\ \ \ \ \ \ \ \ \ \ \text{for each }m  &  \in\left\{  0,1,\ldots,n\right\}  .
\end{align*}
In view of%
\[
a\underbrace{\prod_{i\in J}\left(  1+f_{i}\right)  }_{\substack{=\left(
1+f_{j}\right)  \prod_{i\in J\setminus\left\{  j\right\}  }\left(
1+f_{i}\right)  \\\text{(here, we have split off the factor}\\\text{for
}i=j\text{ from the product)}}}=a\left(  1+f_{j}\right)  \prod_{i\in
J\setminus\left\{  j\right\}  }\left(  1+f_{i}\right)  =\left(  a\prod_{i\in
J\setminus\left\{  j\right\}  }\left(  1+f_{i}\right)  \right)  \left(
1+f_{j}\right)  ,
\]
we can rewrite this as follows: We have%
\begin{align}
\left[  x^{m}\right]  \left(  a\prod_{i\in J}\left(  1+f_{i}\right)  \right)
&  =\left[  x^{m}\right]  \left(  a\prod_{i\in J\setminus\left\{  j\right\}
}\left(  1+f_{i}\right)  \right) \label{pf.lem.fps.prod.irlv.fin.2}\\
\ \ \ \ \ \ \ \ \ \ \text{for each }m  &  \in\left\{  0,1,\ldots,n\right\}
.\nonumber
\end{align}

On the other hand, each $i\in J\setminus\left\{  j\right\}  $ satisfies $i\in
J\setminus\left\{  j\right\}  \subseteq J$ and therefore%
\[
\left[  x^{m}\right]  \left(  f_{i}\right)  =0\ \ \ \ \ \ \ \ \ \ \text{for
each }m\in\left\{  0,1,\ldots,n\right\}
\]
(by (\ref{eq.lem.fps.prod.irlv.fin.ass})). Furthermore, from $j\in J$, we
obtain $\left\vert J\setminus\left\{  j\right\}  \right\vert =\left\vert
J\right\vert -1=k$ (since $\left\vert J\right\vert =k+1$). Hence, we can apply
Lemma \ref{lem.fps.prod.irlv.fin} to $J\setminus\left\{  j\right\}  $ instead
of $J$ (because our induction hypothesis tells us that Lemma
\ref{lem.fps.prod.irlv.fin} is true in the case when $\left\vert J\right\vert
=k$). We thus conclude that%
\begin{equation}
\left[  x^{m}\right]  \left(  a\prod_{i\in J\setminus\left\{  j\right\}
}\left(  1+f_{i}\right)  \right)  =\left[  x^{m}\right]
a\ \ \ \ \ \ \ \ \ \ \text{for each }m\in\left\{  0,1,\ldots,n\right\}  .
\label{pf.lem.fps.prod.irlv.fin.4}%
\end{equation}

Now, for each $m\in\left\{  0,1,\ldots,n\right\}  $, we obtain%
\begin{align*}
\left[  x^{m}\right]  \left(  a\prod_{i\in J}\left(  1+f_{i}\right)  \right)
&  =\left[  x^{m}\right]  \left(  a\prod_{i\in J\setminus\left\{  j\right\}
}\left(  1+f_{i}\right)  \right)  \ \ \ \ \ \ \ \ \ \ \left(  \text{by
(\ref{pf.lem.fps.prod.irlv.fin.2})}\right) \\
&  =\left[  x^{m}\right]  a\ \ \ \ \ \ \ \ \ \ \left(  \text{by
(\ref{pf.lem.fps.prod.irlv.fin.4})}\right)  .
\end{align*}
This is precisely the claim of Lemma \ref{lem.fps.prod.irlv.fin}. Thus, we
have proved that Lemma \ref{lem.fps.prod.irlv.fin} is true in the case when
$\left\vert J\right\vert =k+1$. This completes the induction step. Thus, the
induction proof of Lemma \ref{lem.fps.prod.irlv.fin} is complete.
\end{proof}

\begin{proof}
[Detailed proof of Theorem \ref{thm.fps.1+f-mulable}.]The family $\left(
f_{i}\right)  _{i\in I}$ is summable. In other words,%
\[
\text{for each }n\in\mathbb{N}\text{, all but finitely many }i\in I\text{
satisfy }\left[  x^{n}\right]  f_{i}=0
\]
(by the definition of \textquotedblleft summable\textquotedblright). In other
words, for each $n\in\mathbb{N}$, there exists a finite subset $I_{n}$ of $I$
such that%
\begin{equation}
\text{all }i\in I\setminus I_{n}\text{ satisfy }\left[  x^{n}\right]  f_{i}=0.
\label{pf.thm.fps.1+f-mulable.In}%
\end{equation}
Consider this subset $I_{n}$. Thus, all the sets $I_{0},I_{1},I_{2},\ldots$
are finite subsets of $I$.

Now, let $n\in\mathbb{N}$ be arbitrary. Let $M:=I_{0}\cup I_{1}\cup\cdots\cup
I_{n}$. Then, $M$ is a union of $n+1$ finite subsets of $I$ (because all the
sets $I_{0},I_{1},I_{2},\ldots$ are finite subsets of $I$), and thus itself is
a finite subset of $I$. Moreover,
\begin{equation}
\text{all }m\in\left\{  0,1,\ldots,n\right\}  \text{ and all }i\in I\setminus
M\text{ satisfy }\left[  x^{m}\right]  f_{i}=0.
\label{pf.thm.fps.1+f-mulable.uni}%
\end{equation}

[\textit{Proof of (\ref{pf.thm.fps.1+f-mulable.uni}):} Let $m\in\left\{
0,1,\ldots,n\right\}  $ and $i\in I\setminus M$. We must show that $\left[
x^{m}\right]  f_{i}=0$.

From $m\in\left\{  0,1,\ldots,n\right\}  $, we obtain $I_{m}\subseteq
I_{0}\cup I_{1}\cup\cdots\cup I_{n}=M$, so that $M\supseteq I_{m}$ and thus
$I\setminus\underbrace{M}_{\supseteq I_{m}}\subseteq I\setminus I_{m}$. Hence,
$i\in I\setminus M\subseteq I\setminus I_{m}$. Therefore,
(\ref{pf.thm.fps.1+f-mulable.In}) (applied to $m$ instead of $n$) yields
$\left[  x^{m}\right]  f_{i}=0$. This proves (\ref{pf.thm.fps.1+f-mulable.uni}).]

Now, we shall prove that the set $M$ determines the $x^{n}$-coefficient in the
product of $\left(  1+f_{i}\right)  _{i\in I}$. Indeed, let $J$ be a finite
subset of $I$ satisfying $M\subseteq J\subseteq I$. Let $a$ be the FPS
$\prod_{i\in M}\left(  1+f_{i}\right)  $ (this is well-defined, since the set
$M$ is finite). We have $J\setminus M\subseteq J$, so that the set $J\setminus
M$ is finite (since $J$ is finite). Hence, $\left(  f_{i}\right)  _{i\in
J\setminus M}$ is a finite family of FPSs. Moreover, each $i\in J\setminus M$
satisfies%
\[
\left[  x^{m}\right]  f_{i}=0\ \ \ \ \ \ \ \ \ \ \text{for each }m\in\left\{
0,1,\ldots,n\right\}
\]
(by (\ref{pf.thm.fps.1+f-mulable.uni}), because $i\in\underbrace{J}_{\subseteq
I}\setminus M\subseteq I\setminus M$). Thus, Lemma \ref{lem.fps.prod.irlv.fin}
(applied to $J\setminus M$ instead of $J$) yields%
\[
\left[  x^{m}\right]  \left(  a\prod_{i\in J\setminus M}\left(  1+f_{i}%
\right)  \right)  =\left[  x^{m}\right]  a\ \ \ \ \ \ \ \ \ \ \text{for each
}m\in\left\{  0,1,\ldots,n\right\}  .
\]
Applying this to $m=n$, we find%
\begin{equation}
\left[  x^{n}\right]  \left(  a\prod_{i\in J\setminus M}\left(  1+f_{i}%
\right)  \right)  =\left[  x^{n}\right]  a=\left[  x^{n}\right]  \left(
\prod_{i\in M}\left(  1+f_{i}\right)  \right)
\label{pf.thm.fps.1+f-mulable.4}%
\end{equation}
(since $a=\prod_{i\in M}\left(  1+f_{i}\right)  $). However, the finite set
$J$ is the union of the two disjoint sets $M$ and $J\setminus M$ (since
$M\subseteq J$). Hence, the product $\prod_{i\in J}\left(  1+f_{i}\right)  $
can be split as follows:%
\[
\prod_{i\in J}\left(  1+f_{i}\right)  =\underbrace{\left(  \prod_{i\in
M}\left(  1+f_{i}\right)  \right)  }_{\substack{=a\\\text{(by the definition
of }a\text{)}}}\left(  \prod_{i\in J\setminus M}\left(  1+f_{i}\right)
\right)  =a\prod_{i\in J\setminus M}\left(  1+f_{i}\right)  .
\]
In view of this, we can rewrite (\ref{pf.thm.fps.1+f-mulable.4}) as
\[
\left[  x^{n}\right]  \left(  \prod_{i\in J}\left(  1+f_{i}\right)  \right)
=\left[  x^{n}\right]  \left(  \prod_{i\in M}\left(  1+f_{i}\right)  \right)
.
\]

Forget that we fixed $J$. We thus have shown that every finite subset $J$ of
$I$ satisfying $M\subseteq J\subseteq I$ satisfies%
\[
\left[  x^{n}\right]  \left(  \prod_{i\in J}\left(  1+f_{i}\right)  \right)
=\left[  x^{n}\right]  \left(  \prod_{i\in M}\left(  1+f_{i}\right)  \right)
.
\]
In other words, the set $M$ determines the $x^{n}$-coefficient in the product
of $\left(  1+f_{i}\right)  _{i\in I}$ (according to Definition
\ref{def.fps.determines-xn-coeff} \textbf{(b)}). Hence, the $x^{n}%
$-coefficient in the product of $\left(  1+f_{i}\right)  _{i\in I}$ is
finitely determined (according to Definition
\ref{def.fps.xn-coeff-fin-determined} \textbf{(b)}, since the set $M$ is finite).

Forget that we fixed $n$. Hence, we have shown that for each $n\in\mathbb{N}$,
the $x^{n}$-coefficient in the product of $\left(  1+f_{i}\right)  _{i\in I}$
is finitely determined. In other words, each coefficient in this product is
finitely determined. In other words, the family $\left(  1+f_{i}\right)
_{i\in I}$ is multipliable (by the definition of \textquotedblleft
multipliable\textquotedblright). This proves Theorem \ref{thm.fps.1+f-mulable}.
\end{proof}

\begin{proof}
[Detailed proof of Proposition \ref{prop.fps.1-mulable}.]Let $\left(
\mathbf{a}_{i}\right)  _{i\in I}\in K\left[  \left[  x\right]  \right]  ^{I}$
be a family of FPSs. Assume that all but finitely many entries of this family
$\left(  \mathbf{a}_{i}\right)  _{i\in I}$ equal $1$ (that is, all but
finitely many $i\in I$ satisfy $\mathbf{a}_{i}=1$). We must prove that this
family is multipliable.

We have assumed that all but finitely many $i\in I$ satisfy $\mathbf{a}_{i}%
=1$. In other words, there exists a finite subset $M$ of $I$ such that%
\begin{equation}
\text{all }i\in I\setminus M\text{ satisfy }\mathbf{a}_{i}=1.
\label{pf.prop.fps.1-mulable.abf}%
\end{equation}
Consider this $M$.

Let $n\in\mathbb{N}$. Every finite subset $J$ of $I$ satisfying $M\subseteq
J\subseteq I$ satisfies%
\[
\left[  x^{n}\right]  \left(  \prod_{i\in J}\mathbf{a}_{i}\right)  =\left[
x^{n}\right]  \left(  \prod_{i\in M}\mathbf{a}_{i}\right)
\]
\footnote{\textit{Proof.} Let $J$ be a finite subset of $I$ satisfying
$M\subseteq J\subseteq I$. Then, each $i\in J\setminus M$ satisfies
$i\in\underbrace{J}_{\subseteq I}\setminus M\subseteq I\setminus M$ and
therefore%
\begin{equation}
\mathbf{a}_{i}=1 \label{pf.prop.fps.1-mulable.fn2.1}%
\end{equation}
(by (\ref{pf.prop.fps.1-mulable.abf})). However, the set $J$ is the union of
the two disjoint sets $M$ and $J\setminus M$ (since $M\subseteq J$). Hence, we
can split the product $\prod_{i\in J}\mathbf{a}_{i}$ as follows:%
\[
\prod_{i\in J}\mathbf{a}_{i}=\left(  \prod_{i\in M}\mathbf{a}_{i}\right)
\left(  \prod_{i\in J\setminus M}\underbrace{\mathbf{a}_{i}}%
_{\substack{=1\\\text{(by (\ref{pf.prop.fps.1-mulable.fn2.1}))}}}\right)
=\left(  \prod_{i\in M}\mathbf{a}_{i}\right)  \underbrace{\left(  \prod_{i\in
J\setminus M}1\right)  }_{=1}=\prod_{i\in M}\mathbf{a}_{i}.
\]
Therefore, $\left[  x^{n}\right]  \left(  \prod_{i\in J}\mathbf{a}_{i}\right)
=\left[  x^{n}\right]  \left(  \prod_{i\in M}\mathbf{a}_{i}\right)  $, qed.}.
In other words, the set $M$ determines the $x^{n}$-coefficient in the product
of $\left(  \mathbf{a}_{i}\right)  _{i\in I}$ (by the definition of
\textquotedblleft determining the $x^{n}$-coefficient in a
product\textquotedblright). Hence, the $x^{n}$-coefficient in the product of
$\left(  \mathbf{a}_{i}\right)  _{i\in I}$ is finitely determined (by the
definition of \textquotedblleft finitely determined\textquotedblright).

Forget that we fixed $n$. We thus have proved that for each $n\in\mathbb{N}$,
the $x^{n}$-coefficient in the product of $\left(  \mathbf{a}_{i}\right)
_{i\in I}$ is finitely determined. In other words, each coefficient in the
product of $\left(  \mathbf{a}_{i}\right)  _{i\in I}$ is finitely determined.
In other words, the family $\left(  \mathbf{a}_{i}\right)  _{i\in I}$ is
multipliable (by the definition of \textquotedblleft
multipliable\textquotedblright). This proves Proposition
\ref{prop.fps.1-mulable}.
\end{proof}
\end{fineprint}

\begin{fineprint}
\begin{proof}
[Detailed proof of Lemma \ref{lem.fps.mulable.approx}.]Fix $m\in\left\{
0,1,\ldots,n\right\}  $. Recall that the family $\left(  \mathbf{a}%
_{i}\right)  _{i\in I}$ is multipliable. In other words, each coefficient in
its product is finitely determined. Hence, in particular, the $x^{m}%
$-coefficient in the product of $\left(  \mathbf{a}_{i}\right)  _{i\in I}$ is
finitely determined. In other words, there is a finite subset $M$ of $I$ that
determines the $x^{m}$-coefficient in the product of $\left(  \mathbf{a}%
_{i}\right)  _{i\in I}$. Consider this subset $M$, and denote it by $M_{m}$.
Thus, $M_{m}$ is a finite subset of $I$ that determines the $x^{m}%
$-coefficient in the product of $\left(  \mathbf{a}_{i}\right)  _{i\in I}$.

Forget that we fixed $m$. Thus, for each $m\in\left\{  0,1,\ldots,n\right\}
$, we have defined a finite subset $M_{m}$ of $I$. Let $M$ be the union
$M_{0}\cup M_{1}\cup\cdots\cup M_{n}$ of these (altogether $n+1$) subsets.
Thus, $M$ is a union of finitely many finite subsets of $I$; hence, $M$ itself
is a finite subset of $I$.

Now, let $m\in\left\{  0,1,\ldots,n\right\}  $. We shall show that $M$
determines the $x^{m}$-coefficient in the product of $\left(  \mathbf{a}%
_{i}\right)  _{i\in I}$.

Indeed, let $N$ be a finite subset of $I$ satisfying $M\subseteq N\subseteq
I$. We have $M_{m}\subseteq M$ (since $M$ is the union $M_{0}\cup M_{1}%
\cup\cdots\cup M_{n}$, while $M_{m}$ is one of the sets in this union). Hence,
$M_{m}\subseteq M\subseteq N$. Now, recall that the set $M_{m}$ determines the
$x^{m}$-coefficient in the product of $\left(  \mathbf{a}_{i}\right)  _{i\in
I}$ (by the definition of $M_{m}$). In other words, every finite subset $J$ of
$I$ satisfying $M_{m}\subseteq J\subseteq I$ satisfies%
\begin{equation}
\left[  x^{m}\right]  \left(  \prod_{i\in J}\mathbf{a}_{i}\right)  =\left[
x^{m}\right]  \left(  \prod_{i\in M_{m}}\mathbf{a}_{i}\right)
\label{pf.lem.fps.mulable.approx.1}%
\end{equation}
(by the definition of what it means to \textquotedblleft determine the $x^{m}%
$-coefficient in the product of $\left(  \mathbf{a}_{i}\right)  _{i\in I}%
$\textquotedblright). Applying this to $J=N$, we obtain%
\[
\left[  x^{m}\right]  \left(  \prod_{i\in N}\mathbf{a}_{i}\right)  =\left[
x^{m}\right]  \left(  \prod_{i\in M_{m}}\mathbf{a}_{i}\right)
\]
(since $N$ is a finite subset of $I$ satisfying $M_{m}\subseteq N\subseteq
I$). On the other hand, we can apply (\ref{pf.lem.fps.mulable.approx.1}) to
$J=M$ (since $M$ is a finite subset of $I$ satisfying $M_{m}\subseteq
M\subseteq I$), and thus obtain%
\[
\left[  x^{m}\right]  \left(  \prod_{i\in M}\mathbf{a}_{i}\right)  =\left[
x^{m}\right]  \left(  \prod_{i\in M_{m}}\mathbf{a}_{i}\right)  .
\]
Comparing these two equalities, we obtain%
\[
\left[  x^{m}\right]  \left(  \prod_{i\in N}\mathbf{a}_{i}\right)  =\left[
x^{m}\right]  \left(  \prod_{i\in M}\mathbf{a}_{i}\right)  .
\]

Forget that we fixed $N$. We thus have shown that every finite subset $N$ of
$I$ satisfying $M\subseteq N\subseteq I$ satisfies $\left[  x^{m}\right]
\left(  \prod_{i\in N}\mathbf{a}_{i}\right)  =\left[  x^{m}\right]  \left(
\prod_{i\in M}\mathbf{a}_{i}\right)  $. Renaming $N$ as $J$ in this statement,
we obtain the following: Every finite subset $J$ of $I$ satisfying $M\subseteq
J\subseteq I$ satisfies $\left[  x^{m}\right]  \left(  \prod_{i\in
J}\mathbf{a}_{i}\right)  =\left[  x^{m}\right]  \left(  \prod_{i\in
M}\mathbf{a}_{i}\right)  $. In other words, $M$ determines the $x^{m}%
$-coefficient in the product of $\left(  \mathbf{a}_{i}\right)  _{i\in I}$ (by
the definition of what it means to \textquotedblleft determine the $x^{m}%
$-coefficient in the product of $\left(  \mathbf{a}_{i}\right)  _{i\in I}%
$\textquotedblright).

Forget that we fixed $m$. We thus have shown that $M$ determines the $x^{m}%
$-coefficient in the product of $\left(  \mathbf{a}_{i}\right)  _{i\in I}$ for
each $m\in\left\{  0,1,\ldots,n\right\}  $. In other words, $M$ determines the
first $n+1$ coefficients in the product of $\left(  \mathbf{a}_{i}\right)
_{i\in I}$. In other words, $M$ is an $x^{n}$-approximator for $\left(
\mathbf{a}_{i}\right)  _{i\in I}$ (by the definition of an \textquotedblleft%
$x^{n}$-approximator\textquotedblright). Hence, there exists an $x^{n}%
$-approximator for $\left(  \mathbf{a}_{i}\right)  _{i\in I}$. This proves
Lemma \ref{lem.fps.mulable.approx}.
\end{proof}
\end{fineprint}

\begin{fineprint}
\begin{proof}
[Detailed proof of Proposition \ref{prop.fps.infprod-approx-xneq}.]The set $M$
is an $x^{n}$-approximator for $\left(  \mathbf{a}_{i}\right)  _{i\in I}$. In
other words, $M$ is a finite subset of $I$ that determines the first $n+1$
coefficients in the product of $\left(  \mathbf{a}_{i}\right)  _{i\in I}$ (by
the definition of an \textquotedblleft$x^{n}$-approximator\textquotedblright).
\medskip

\textbf{(a)} Let $m\in\left\{  0,1,\ldots,n\right\}  $. Recall that $M$
determines the first $n+1$ coefficients in the product of $\left(
\mathbf{a}_{i}\right)  _{i\in I}$. Thus, in particular, $M$ determines the
$x^{m}$-coefficient in the product of $\left(  \mathbf{a}_{i}\right)  _{i\in
I}$ (since $m\in\left\{  0,1,\ldots,n\right\}  $). In other words, every
finite subset $J$ of $I$ satisfying $M\subseteq J\subseteq I$ satisfies%
\begin{equation}
\left[  x^{m}\right]  \left(  \prod_{i\in J}\mathbf{a}_{i}\right)  =\left[
x^{m}\right]  \left(  \prod_{i\in M}\mathbf{a}_{i}\right)
\label{pf.prop.fps.infprod-approx-xneq.a.1}%
\end{equation}
(by the definition of \textquotedblleft determining the $x^{m}$-coefficient in
the product of $\left(  \mathbf{a}_{i}\right)  _{i\in I}$\textquotedblright).

Forget that we fixed $m$. We thus have proved that every $m\in\left\{
0,1,\ldots,n\right\}  $ and every finite subset $J$ of $I$ satisfying
$M\subseteq J\subseteq I$ satisfy (\ref{pf.prop.fps.infprod-approx-xneq.a.1}).

Now, let $J$ be a finite subset of $I$ satisfying $M\subseteq J\subseteq I$.
Then, each $m\in\left\{  0,1,\ldots,n\right\}  $ satisfies $\left[
x^{m}\right]  \left(  \prod_{i\in J}\mathbf{a}_{i}\right)  =\left[
x^{m}\right]  \left(  \prod_{i\in M}\mathbf{a}_{i}\right)  $ (by
(\ref{pf.prop.fps.infprod-approx-xneq.a.1})). In other words, we have
$\prod_{i\in J}\mathbf{a}_{i}\overset{x^{n}}{\equiv}\prod_{i\in M}%
\mathbf{a}_{i}$ (by Definition \ref{def.fps.xneq}). This proves Proposition
\ref{prop.fps.infprod-approx-xneq} \textbf{(a)}. \medskip

\textbf{(b)} Assume that the family $\left(  \mathbf{a}_{i}\right)  _{i\in I}$
is multipliable. Let $m\in\left\{  0,1,\ldots,n\right\}  $. Thus, $M$
determines the $x^{m}$-coefficient in the product of $\left(  \mathbf{a}%
_{i}\right)  _{i\in I}$ (since $M$ determines the first $n+1$ coefficients in
the product of $\left(  \mathbf{a}_{i}\right)  _{i\in I}$).

The product $\prod_{i\in I}\mathbf{a}_{i}$ is defined according to Definition
\ref{def.fps.multipliable} \textbf{(b)}. Specifically, Definition
\ref{def.fps.multipliable} \textbf{(b)} (with $n$ and $M$ renamed as $k$ and
$N$) shows that the product $\prod_{i\in I}\mathbf{a}_{i}$ is defined to be
the FPS whose $x^{k}$-coefficient (for any $k\in\mathbb{N}$) can be computed
as follows: If $k\in\mathbb{N}$, and if $N$ is a finite subset of $I$ that
determines the $x^{k}$-coefficient in the product of $\left(  \mathbf{a}%
_{i}\right)  _{i\in I}$, then
\[
\left[  x^{k}\right]  \left(  \prod_{i\in I}\mathbf{a}_{i}\right)  =\left[
x^{k}\right]  \left(  \prod_{i\in N}\mathbf{a}_{i}\right)  .
\]
We can apply this to $k=m$ and $N=M$ (since $M$ is a finite subset of $I$ that
determines the $x^{m}$-coefficient in the product of $\left(  \mathbf{a}%
_{i}\right)  _{i\in I}$), and thus obtain%
\[
\left[  x^{m}\right]  \left(  \prod_{i\in I}\mathbf{a}_{i}\right)  =\left[
x^{m}\right]  \left(  \prod_{i\in M}\mathbf{a}_{i}\right)  .
\]

Forget that we fixed $m$. We thus have proved that each $m\in\left\{
0,1,\ldots,n\right\}  $ satisfies $\left[  x^{m}\right]  \left(  \prod_{i\in
I}\mathbf{a}_{i}\right)  =\left[  x^{m}\right]  \left(  \prod_{i\in
M}\mathbf{a}_{i}\right)  $. In other words, $\prod_{i\in I}\mathbf{a}%
_{i}\overset{x^{n}}{\equiv}\prod_{i\in M}\mathbf{a}_{i}$ (by Definition
\ref{def.fps.xneq}). This proves Proposition
\ref{prop.fps.infprod-approx-xneq} \textbf{(b)}.
\end{proof}
\end{fineprint}

\begin{fineprint}
\begin{proof}
[Detailed proof of Proposition \ref{prop.fps.union-mulable}.]\textbf{(a)} Fix
$n\in\mathbb{N}$. We know that the family $\left(  \mathbf{a}_{i}\right)
_{i\in J}$ is multipliable. Hence, there exists an $x^{n}$-approximator $U$
for $\left(  \mathbf{a}_{i}\right)  _{i\in J}$ (by Lemma
\ref{lem.fps.mulable.approx}, applied to $J$ instead of $I$). Consider this
$U$.

We also know that the family $\left(  \mathbf{a}_{i}\right)  _{i\in I\setminus
J}$ is multipliable. Hence, there exists an $x^{n}$-approximator $V$ for
$\left(  \mathbf{a}_{i}\right)  _{i\in I\setminus J}$ (by Lemma
\ref{lem.fps.mulable.approx}, applied to $I\setminus J$ instead of $I$).
Consider this $V$.

We know that $U$ is an $x^{n}$-approximator for $\left(  \mathbf{a}%
_{i}\right)  _{i\in J}$. In other words, $U$ is a finite subset of $J$ that
determines the first $n+1$ coefficients in the product of $\left(
\mathbf{a}_{i}\right)  _{i\in J}$ (by the definition of an $x^{n}%
$-approximator). Hence, in particular, $U$ is finite. Similarly, $V$ is
finite. Moreover, $U\subseteq J$ (since $U$ is a subset of $J$); similarly,
$V\subseteq I\setminus J$.

Let $M=U\cup V$. This set $M$ is finite (since $U$ and $V$ are finite).
Moreover, using $U\subseteq J\subseteq I$ and $V\subseteq I\setminus
J\subseteq I$, we obtain $M=\underbrace{U}_{\subseteq I}\cup\underbrace{V}%
_{\subseteq I}\subseteq I\cup I=I$. Hence, $M$ is a finite subset of $I$. Note
that the sets $U$ and $V$ are disjoint\footnote{Indeed, $\underbrace{U}%
_{\subseteq J}\cap\underbrace{V}_{\subseteq I\setminus J}\subseteq
J\cap\left(  I\setminus J\right)  =\varnothing$, so that $U\cap V=\varnothing
$.}. Hence, the set $M$ is the union of its two disjoint subsets $U$ and $V$
(since $M=U\cup V$).

Now, let $N$ be a finite subset of $I$ satisfying $M\subseteq N\subseteq I$.
We shall show that
\[
\left[  x^{n}\right]  \left(  \prod_{i\in N}\mathbf{a}_{i}\right)  =\left[
x^{n}\right]  \left(  \prod_{i\in M}\mathbf{a}_{i}\right)  .
\]

Indeed, let $m\in\left\{  0,1,\ldots,n\right\}  $. The set $N\cap J$ is a
finite subset of $J$ (since $N$ is a finite subset of $I$) and satisfies
$U\subseteq N\cap J$ (since $U\subseteq U\cup V=M\subseteq N$ and $U\subseteq
J$). Now, recall that $U$ determines the first $n+1$ coefficients in the
product of $\left(  \mathbf{a}_{i}\right)  _{i\in J}$. Hence, $U$ determines
the $x^{m}$-coefficient in the product of $\left(  \mathbf{a}_{i}\right)
_{i\in J}$ (since $m\in\left\{  0,1,\ldots,n\right\}  $). In other words,
every finite subset $T$ of $J$ satisfying $U\subseteq T\subseteq J$ satisfies%
\[
\left[  x^{m}\right]  \left(  \prod_{i\in T}\mathbf{a}_{i}\right)  =\left[
x^{m}\right]  \left(  \prod_{i\in U}\mathbf{a}_{i}\right)
\]
(by the definition of what it means to \textquotedblleft determine the $x^{m}%
$-coefficient in the product of $\left(  \mathbf{a}_{i}\right)  _{i\in J}%
$\textquotedblright). We can apply this to $T=N\cap J$ (since $N\cap J$ is a
finite subset of $J$ satisfying $U\subseteq N\cap J\subseteq J$), and thus
obtain%
\begin{equation}
\left[  x^{m}\right]  \left(  \prod_{i\in N\cap J}\mathbf{a}_{i}\right)
=\left[  x^{m}\right]  \left(  \prod_{i\in U}\mathbf{a}_{i}\right)  .
\label{pf.prop.fps.union-mulable.5J}%
\end{equation}
The same argument (applied to $I\setminus J$ and $V$ instead of $J$ and $U$)
yields%
\[
\left[  x^{m}\right]  \left(  \prod_{i\in N\cap\left(  I\setminus J\right)
}\mathbf{a}_{i}\right)  =\left[  x^{m}\right]  \left(  \prod_{i\in
V}\mathbf{a}_{i}\right)  .
\]
In view of $N\cap\left(  I\setminus J\right)  =\underbrace{\left(  N\cap
I\right)  }_{\substack{=N\\\text{(since }N\subseteq I\text{)}}}\setminus
\,J=N\setminus J$, this rewrites as%
\begin{equation}
\left[  x^{m}\right]  \left(  \prod_{i\in N\setminus J}\mathbf{a}_{i}\right)
=\left[  x^{m}\right]  \left(  \prod_{i\in V}\mathbf{a}_{i}\right)  .
\label{pf.prop.fps.union-mulable.5IoJ}%
\end{equation}

Forget that we fixed $m$. We thus have proved the equalities
(\ref{pf.prop.fps.union-mulable.5J}) and (\ref{pf.prop.fps.union-mulable.5IoJ}%
) for each $m\in\left\{  0,1,\ldots,n\right\}  $. Hence, Lemma
\ref{lem.fps.prod.irlv.cong-mul} (applied to $a=\prod_{i\in N\cap J}%
\mathbf{a}_{i}$ and $b=\prod_{i\in U}\mathbf{a}_{i}$ and $c=\prod_{i\in
N\setminus J}\mathbf{a}_{i}$ and $d=\prod_{i\in V}\mathbf{a}_{i}$) yields that%
\[
\left[  x^{m}\right]  \left(  \left(  \prod_{i\in N\cap J}\mathbf{a}%
_{i}\right)  \left(  \prod_{i\in N\setminus J}\mathbf{a}_{i}\right)  \right)
=\left[  x^{m}\right]  \left(  \left(  \prod_{i\in U}\mathbf{a}_{i}\right)
\left(  \prod_{i\in V}\mathbf{a}_{i}\right)  \right)
\]
for each $m\in\left\{  0,1,\ldots,n\right\}  $. In view of
\[
\prod_{i\in N}\mathbf{a}_{i}=\left(  \prod_{i\in N\cap J}\mathbf{a}%
_{i}\right)  \left(  \prod_{i\in N\setminus J}\mathbf{a}_{i}\right)
\ \ \ \ \ \ \ \ \ \ \left(
\begin{array}
[c]{c}%
\text{since the set }N\text{ is the union of its}\\
\text{two disjoint subsets }N\cap J\text{ and }N\setminus J
\end{array}
\right)
\]
and%
\[
\prod_{i\in M}\mathbf{a}_{i}=\left(  \prod_{i\in U}\mathbf{a}_{i}\right)
\left(  \prod_{i\in V}\mathbf{a}_{i}\right)  \ \ \ \ \ \ \ \ \ \ \left(
\begin{array}
[c]{c}%
\text{since the set }M\text{ is the union of its}\\
\text{two disjoint subsets }U\text{ and }V
\end{array}
\right)  ,
\]
this rewrites as follows: We have%
\[
\left[  x^{m}\right]  \left(  \prod_{i\in N}\mathbf{a}_{i}\right)  =\left[
x^{m}\right]  \left(  \prod_{i\in M}\mathbf{a}_{i}\right)
\]
for each $m\in\left\{  0,1,\ldots,n\right\}  $. Applying this to $m=n$, we
obtain%
\[
\left[  x^{n}\right]  \left(  \prod_{i\in N}\mathbf{a}_{i}\right)  =\left[
x^{n}\right]  \left(  \prod_{i\in M}\mathbf{a}_{i}\right)  .
\]

Forget that we fixed $N$. We thus have shown that every finite subset $N$ of
$I$ satisfying $M\subseteq N\subseteq I$ satisfies $\left[  x^{n}\right]
\left(  \prod_{i\in N}\mathbf{a}_{i}\right)  =\left[  x^{n}\right]  \left(
\prod_{i\in M}\mathbf{a}_{i}\right)  $. In other words, $M$ determines the
$x^{n}$-coefficient in the product of $\left(  \mathbf{a}_{i}\right)  _{i\in
I}$ (by the definition of what it means to \textquotedblleft determine the
$x^{n}$-coefficient in the product of $\left(  \mathbf{a}_{i}\right)  _{i\in
I}$\textquotedblright). Hence, the $x^{n}$-coefficient in the product of
$\left(  \mathbf{a}_{i}\right)  _{i\in I}$ is finitely determined (by the
definition of \textquotedblleft finitely determined\textquotedblright, since
$M$ is a finite subset of $I$).

Forget that we fixed $n$. We thus have proved that for each $n\in\mathbb{N}$,
the $x^{n}$-coefficient in the product of $\left(  \mathbf{a}_{i}\right)
_{i\in I}$ is finitely determined. In other words, each coefficient in the
product of $\left(  \mathbf{a}_{i}\right)  _{i\in I}$ is finitely determined.
In other words, the family $\left(  \mathbf{a}_{i}\right)  _{i\in I}$ is
multipliable (by the definition of \textquotedblleft
multipliable\textquotedblright). This proves Proposition
\ref{prop.fps.union-mulable} \textbf{(a)}. \medskip

\textbf{(b)} Proposition \ref{prop.fps.union-mulable} \textbf{(a)} shows that
the family $\left(  \mathbf{a}_{i}\right)  _{i\in I}$ is multipliable.

Let $n\in\mathbb{N}$. In our above proof of Proposition
\ref{prop.fps.union-mulable} \textbf{(a)}, we have seen the following:

\begin{itemize}
\item There exists an $x^{n}$-approximator $U$ for $\left(  \mathbf{a}%
_{i}\right)  _{i\in J}$.

\item There exists an $x^{n}$-approximator $V$ for $\left(  \mathbf{a}%
_{i}\right)  _{i\in I\setminus J}$.
\end{itemize}

Consider these $U$ and $V$. Let $M=U\cup V$. In our above proof of Proposition
\ref{prop.fps.union-mulable} \textbf{(a)}, we have seen the following:

\begin{itemize}
\item The set $M$ is a finite subset of $I$.

\item The set $M$ is the union of its two disjoint subsets $U$ and $V$.

\item The set $M$ determines the $x^{n}$-coefficient in the product of
$\left(  \mathbf{a}_{i}\right)  _{i\in I}$.
\end{itemize}

Now, the definition of the infinite product $\prod_{i\in I}\mathbf{a}_{i}$
(namely, Definition \ref{def.fps.multipliable} \textbf{(b)}) yields that%
\begin{equation}
\left[  x^{n}\right]  \left(  \prod_{i\in I}\mathbf{a}_{i}\right)  =\left[
x^{n}\right]  \left(  \prod_{i\in M}\mathbf{a}_{i}\right)
\label{pf.prop.fps.union-mulable.b.2}%
\end{equation}
(since $M$ is a finite subset of $I$ that determines the $x^{n}$-coefficient
in the product of $\left(  \mathbf{a}_{i}\right)  _{i\in I}$). On the other
hand, the set $U$ is an $x^{n}$-approximator for $\left(  \mathbf{a}%
_{i}\right)  _{i\in J}$. Thus, Proposition \ref{prop.fps.infprod-approx-xneq}
\textbf{(b)} (applied to $J$ and $U$ instead of $I$ and $M$) yields
\begin{equation}
\prod_{i\in J}\mathbf{a}_{i}\overset{x^{n}}{\equiv}\prod_{i\in U}%
\mathbf{a}_{i} \label{pf.prop.fps.union-mulable.b.3}%
\end{equation}
(since the family $\left(  \mathbf{a}_{i}\right)  _{i\in J}$ is multipliable).
Furthermore, the set $V$ is an $x^{n}$-approximator for $\left(
\mathbf{a}_{i}\right)  _{i\in I\setminus J}$. Thus, Proposition
\ref{prop.fps.infprod-approx-xneq} \textbf{(b)} (applied to $I\setminus J$ and
$V$ instead of $I$ and $M$) yields
\begin{equation}
\prod_{i\in I\setminus J}\mathbf{a}_{i}\overset{x^{n}}{\equiv}\prod_{i\in
V}\mathbf{a}_{i} \label{pf.prop.fps.union-mulable.b.4}%
\end{equation}
(since the family $\left(  \mathbf{a}_{i}\right)  _{i\in I\setminus J}$ is
multipliable). From (\ref{pf.prop.fps.union-mulable.b.3}) and
(\ref{pf.prop.fps.union-mulable.b.4}), we obtain%
\[
\left(  \prod_{i\in J}\mathbf{a}_{i}\right)  \left(  \prod_{i\in I\setminus
J}\mathbf{a}_{i}\right)  \overset{x^{n}}{\equiv}\left(  \prod_{i\in
U}\mathbf{a}_{i}\right)  \left(  \prod_{i\in V}\mathbf{a}_{i}\right)
\]
(by (\ref{eq.thm.fps.xneq.props.b.*}), applied to $a=\prod_{i\in J}%
\mathbf{a}_{i}$ and $b=\prod_{i\in U}\mathbf{a}_{i}$ and $c=\prod_{i\in
I\setminus J}\mathbf{a}_{i}$ and $d=\prod_{i\in V}\mathbf{a}_{i}$). In view
of
\[
\prod_{i\in M}\mathbf{a}_{i}=\left(  \prod_{i\in U}\mathbf{a}_{i}\right)
\left(  \prod_{i\in V}\mathbf{a}_{i}\right)  \ \ \ \ \ \ \ \ \ \ \left(
\begin{array}
[c]{c}%
\text{since the set }M\text{ is the union of its}\\
\text{two disjoint subsets }U\text{ and }V
\end{array}
\right)  ,
\]
this rewrites as
\[
\left(  \prod_{i\in J}\mathbf{a}_{i}\right)  \left(  \prod_{i\in I\setminus
J}\mathbf{a}_{i}\right)  \overset{x^{n}}{\equiv}\prod_{i\in M}\mathbf{a}_{i}.
\]
In other words,
\[
\text{each }m\in\left\{  0,1,\ldots,n\right\}  \text{ satisfies }\left[
x^{m}\right]  \left(  \left(  \prod_{i\in J}\mathbf{a}_{i}\right)  \left(
\prod_{i\in I\setminus J}\mathbf{a}_{i}\right)  \right)  =\left[
x^{m}\right]  \left(  \prod_{i\in M}\mathbf{a}_{i}\right)
\]
(by Definition \ref{def.fps.xneq}). Applying this to $m=n$, we obtain
\[
\left[  x^{n}\right]  \left(  \left(  \prod_{i\in J}\mathbf{a}_{i}\right)
\left(  \prod_{i\in I\setminus J}\mathbf{a}_{i}\right)  \right)  =\left[
x^{n}\right]  \left(  \prod_{i\in M}\mathbf{a}_{i}\right)  .
\]
Comparing this with (\ref{pf.prop.fps.union-mulable.b.2}), we obtain%
\begin{equation}
\left[  x^{n}\right]  \left(  \prod_{i\in I}\mathbf{a}_{i}\right)  =\left[
x^{n}\right]  \left(  \left(  \prod_{i\in J}\mathbf{a}_{i}\right)  \left(
\prod_{i\in I\setminus J}\mathbf{a}_{i}\right)  \right)  .
\label{pf.prop.fps.union-mulable.b.at}%
\end{equation}

Now, forget that we fixed $n$. We thus have proved that each $n\in\mathbb{N}$
satisfies (\ref{pf.prop.fps.union-mulable.b.at}). In other words, each
coefficient of the FPS $\prod_{i\in I}\mathbf{a}_{i}$ equals the corresponding
coefficient of $\left(  \prod_{i\in J}\mathbf{a}_{i}\right)  \left(
\prod_{i\in I\setminus J}\mathbf{a}_{i}\right)  $. Hence, we have%
\[
\prod_{i\in I}\mathbf{a}_{i}=\left(  \prod_{i\in J}\mathbf{a}_{i}\right)
\cdot\left(  \prod_{i\in I\setminus J}\mathbf{a}_{i}\right)  .
\]
This proves Proposition \ref{prop.fps.union-mulable} \textbf{(b)}.
\end{proof}
\end{fineprint}

\begin{fineprint}
\begin{proof}
[Detailed proof of Proposition \ref{prop.fps.prod-mulable}.]\textbf{(a)} Fix
$n\in\mathbb{N}$. We know that the family $\left(  \mathbf{a}_{i}\right)
_{i\in I}$ is multipliable. Hence, there exists an $x^{n}$-approximator $U$
for $\left(  \mathbf{a}_{i}\right)  _{i\in I}$ (by Lemma
\ref{lem.fps.mulable.approx}). Consider this $U$.

We also know that the family $\left(  \mathbf{b}_{i}\right)  _{i\in I}$ is
multipliable. Hence, there exists an $x^{n}$-approximator $V$ for $\left(
\mathbf{b}_{i}\right)  _{i\in I}$ (by Lemma \ref{lem.fps.mulable.approx},
applied to $\mathbf{b}_{i}$ instead of $\mathbf{a}_{i}$). Consider this $V$.

We know that $U$ is an $x^{n}$-approximator for $\left(  \mathbf{a}%
_{i}\right)  _{i\in I}$. In other words, $U$ is a finite subset of $I$ that
determines the first $n+1$ coefficients in the product of $\left(
\mathbf{a}_{i}\right)  _{i\in I}$ (by the definition of an $x^{n}%
$-approximator). Hence, in particular, $U$ is finite. Similarly, $V$ is
finite. Moreover, $U\subseteq I$ (since $U$ is a subset of $I$); similarly,
$V\subseteq I$.

Let $M=U\cup V$. This set $M$ is finite (since $U$ and $V$ are finite).
Moreover, using $U\subseteq I$ and $V\subseteq I$, we obtain $M=\underbrace{U}%
_{\subseteq I}\cup\underbrace{V}_{\subseteq I}\subseteq I\cup I=I$. Hence, $M$
is a finite subset of $I$.

Now, let $N$ be a finite subset of $I$ satisfying $M\subseteq N\subseteq I$.
We shall show that
\[
\left[  x^{n}\right]  \left(  \prod_{i\in N}\left(  \mathbf{a}_{i}%
\mathbf{b}_{i}\right)  \right)  =\left[  x^{n}\right]  \left(  \prod_{i\in
M}\left(  \mathbf{a}_{i}\mathbf{b}_{i}\right)  \right)  .
\]

Indeed, let $m\in\left\{  0,1,\ldots,n\right\}  $. We have $U\subseteq U\cup
V=M\subseteq N$. The set $N$ is a finite subset of $I$ and satisfies
$U\subseteq N$. Now, recall that $U$ determines the first $n+1$ coefficients
in the product of $\left(  \mathbf{a}_{i}\right)  _{i\in I}$. Hence, $U$
determines the $x^{m}$-coefficient in the product of $\left(  \mathbf{a}%
_{i}\right)  _{i\in I}$ (since $m\in\left\{  0,1,\ldots,n\right\}  $). In
other words, every finite subset $T$ of $I$ satisfying $U\subseteq T\subseteq
I$ satisfies%
\begin{equation}
\left[  x^{m}\right]  \left(  \prod_{i\in T}\mathbf{a}_{i}\right)  =\left[
x^{m}\right]  \left(  \prod_{i\in U}\mathbf{a}_{i}\right)
\label{pf.prop.fps.prod-mulable.4}%
\end{equation}
(by the definition of what it means to \textquotedblleft determine the $x^{m}%
$-coefficient in the product of $\left(  \mathbf{a}_{i}\right)  _{i\in I}%
$\textquotedblright). We can apply this to $T=N$ (since $N$ is a finite subset
of $I$ satisfying $U\subseteq N\subseteq I$), and thus obtain%
\[
\left[  x^{m}\right]  \left(  \prod_{i\in N}\mathbf{a}_{i}\right)  =\left[
x^{m}\right]  \left(  \prod_{i\in U}\mathbf{a}_{i}\right)  .
\]
However, we can also apply (\ref{pf.prop.fps.prod-mulable.4}) to $T=M$ (since
$M$ is a finite subset of $I$ satisfying $U\subseteq M\subseteq I$), and thus
obtain%
\[
\left[  x^{m}\right]  \left(  \prod_{i\in M}\mathbf{a}_{i}\right)  =\left[
x^{m}\right]  \left(  \prod_{i\in U}\mathbf{a}_{i}\right)  .
\]
Comparing these two equalities, we obtain%
\begin{equation}
\left[  x^{m}\right]  \left(  \prod_{i\in N}\mathbf{a}_{i}\right)  =\left[
x^{m}\right]  \left(  \prod_{i\in M}\mathbf{a}_{i}\right)  .
\label{pf.prop.fps.prod-mulable.5a}%
\end{equation}
The same argument (applied to $\mathbf{b}_{i}$ and $V$ instead of
$\mathbf{a}_{i}$ and $U$) yields%
\begin{equation}
\left[  x^{m}\right]  \left(  \prod_{i\in N}\mathbf{b}_{i}\right)  =\left[
x^{m}\right]  \left(  \prod_{i\in M}\mathbf{b}_{i}\right)  .
\label{pf.prop.fps.prod-mulable.5b}%
\end{equation}

Forget that we fixed $m$. We thus have proved the equalities
(\ref{pf.prop.fps.prod-mulable.5a}) and (\ref{pf.prop.fps.prod-mulable.5b})
for each $m\in\left\{  0,1,\ldots,n\right\}  $. Hence, Lemma
\ref{lem.fps.prod.irlv.cong-mul} (applied to $a=\prod_{i\in N}\mathbf{a}_{i}$
and $b=\prod_{i\in M}\mathbf{a}_{i}$ and $c=\prod_{i\in N}\mathbf{b}_{i}$ and
$d=\prod_{i\in M}\mathbf{b}_{i}$) yields that%
\[
\left[  x^{m}\right]  \left(  \left(  \prod_{i\in N}\mathbf{a}_{i}\right)
\left(  \prod_{i\in N}\mathbf{b}_{i}\right)  \right)  =\left[  x^{m}\right]
\left(  \left(  \prod_{i\in M}\mathbf{a}_{i}\right)  \left(  \prod_{i\in
M}\mathbf{b}_{i}\right)  \right)
\]
for each $m\in\left\{  0,1,\ldots,n\right\}  $. Applying this to $m=n$, we
obtain%
\[
\left[  x^{n}\right]  \left(  \left(  \prod_{i\in N}\mathbf{a}_{i}\right)
\left(  \prod_{i\in N}\mathbf{b}_{i}\right)  \right)  =\left[  x^{n}\right]
\left(  \left(  \prod_{i\in M}\mathbf{a}_{i}\right)  \left(  \prod_{i\in
M}\mathbf{b}_{i}\right)  \right)  .
\]
In view of
\[
\prod_{i\in N}\left(  \mathbf{a}_{i}\mathbf{b}_{i}\right)  =\left(
\prod_{i\in N}\mathbf{a}_{i}\right)  \left(  \prod_{i\in N}\mathbf{b}%
_{i}\right)  \ \ \ \ \ \ \ \ \ \ \left(
\begin{array}
[c]{c}%
\text{by the standard rules for finite}\\
\text{products, since }N\text{ is a finite set}%
\end{array}
\right)
\]
and%
\[
\prod_{i\in M}\left(  \mathbf{a}_{i}\mathbf{b}_{i}\right)  =\left(
\prod_{i\in M}\mathbf{a}_{i}\right)  \left(  \prod_{i\in M}\mathbf{b}%
_{i}\right)  \ \ \ \ \ \ \ \ \ \ \left(
\begin{array}
[c]{c}%
\text{by the standard rules for finite}\\
\text{products, since }M\text{ is a finite set}%
\end{array}
\right)  ,
\]
this rewrites as%
\[
\left[  x^{n}\right]  \left(  \prod_{i\in N}\left(  \mathbf{a}_{i}%
\mathbf{b}_{i}\right)  \right)  =\left[  x^{n}\right]  \left(  \prod_{i\in
M}\left(  \mathbf{a}_{i}\mathbf{b}_{i}\right)  \right)  .
\]

Forget that we fixed $N$. We thus have shown that every finite subset $N$ of
$I$ satisfying $M\subseteq N\subseteq I$ satisfies $\left[  x^{n}\right]
\left(  \prod_{i\in N}\left(  \mathbf{a}_{i}\mathbf{b}_{i}\right)  \right)
=\left[  x^{n}\right]  \left(  \prod_{i\in M}\left(  \mathbf{a}_{i}%
\mathbf{b}_{i}\right)  \right)  $. In other words, $M$ determines the $x^{n}%
$-coefficient in the product of $\left(  \mathbf{a}_{i}\mathbf{b}_{i}\right)
_{i\in I}$ (by the definition of what it means to \textquotedblleft determine
the $x^{n}$-coefficient in the product of $\left(  \mathbf{a}_{i}%
\mathbf{b}_{i}\right)  _{i\in I}$\textquotedblright). Hence, the $x^{n}%
$-coefficient in the product of $\left(  \mathbf{a}_{i}\mathbf{b}_{i}\right)
_{i\in I}$ is finitely determined (by the definition of \textquotedblleft
finitely determined\textquotedblright, since $M$ is a finite subset of $I$).

Forget that we fixed $n$. We thus have proved that for each $n\in\mathbb{N}$,
the $x^{n}$-coefficient in the product of $\left(  \mathbf{a}_{i}%
\mathbf{b}_{i}\right)  _{i\in I}$ is finitely determined. In other words, each
coefficient in the product of $\left(  \mathbf{a}_{i}\mathbf{b}_{i}\right)
_{i\in I}$ is finitely determined. In other words, the family $\left(
\mathbf{a}_{i}\mathbf{b}_{i}\right)  _{i\in I}$ is multipliable (by the
definition of \textquotedblleft multipliable\textquotedblright). This proves
Proposition \ref{prop.fps.prod-mulable} \textbf{(a)}. \medskip

\textbf{(b)} Proposition \ref{prop.fps.prod-mulable} \textbf{(a)} shows that
the family $\left(  \mathbf{a}_{i}\mathbf{b}_{i}\right)  _{i\in I}$ is multipliable.

Let $n\in\mathbb{N}$. In our above proof of Proposition
\ref{prop.fps.prod-mulable} \textbf{(a)}, we have seen the following:

\begin{itemize}
\item There exists an $x^{n}$-approximator $U$ for $\left(  \mathbf{a}%
_{i}\right)  _{i\in I}$.

\item There exists an $x^{n}$-approximator $V$ for $\left(  \mathbf{b}%
_{i}\right)  _{i\in I}$.
\end{itemize}

Consider these $U$ and $V$. Let $M=U\cup V$. Thus, $M=U\cup V\supseteq U$ and
$M=U\cup V\supseteq V$. In our above proof of Proposition
\ref{prop.fps.prod-mulable} \textbf{(a)}, we have seen the following:

\begin{itemize}
\item The set $M$ is a finite subset of $I$.

\item The set $M$ determines the $x^{n}$-coefficient in the product of
$\left(  \mathbf{a}_{i}\mathbf{b}_{i}\right)  _{i\in I}$.
\end{itemize}

We recall that the relation $\overset{x^{n}}{\equiv}$ on $K\left[  \left[
x\right]  \right]  $ is an equivalence relation (by Theorem
\ref{thm.fps.xneq.props} \textbf{(a)}). Thus, this relation $\overset{x^{n}%
}{\equiv}$ is transitive and symmetric.

Now, the definition of the infinite product $\prod_{i\in I}\left(
\mathbf{a}_{i}\mathbf{b}_{i}\right)  $ (namely, Definition
\ref{def.fps.multipliable} \textbf{(b)}) yields that%
\begin{equation}
\left[  x^{n}\right]  \left(  \prod_{i\in I}\left(  \mathbf{a}_{i}%
\mathbf{b}_{i}\right)  \right)  =\left[  x^{n}\right]  \left(  \prod_{i\in
M}\left(  \mathbf{a}_{i}\mathbf{b}_{i}\right)  \right)
\label{pf.prop.fps.prod-mulable.b.2}%
\end{equation}
(since $M$ is a finite subset of $I$ that determines the $x^{n}$-coefficient
in the product of $\left(  \mathbf{a}_{i}\mathbf{b}_{i}\right)  _{i\in I}$).
On the other hand, the set $U$ is an $x^{n}$-approximator for $\left(
\mathbf{a}_{i}\right)  _{i\in I}$. Thus, Proposition
\ref{prop.fps.infprod-approx-xneq} \textbf{(b)} (applied to $U$ instead of
$M$) yields
\[
\prod_{i\in I}\mathbf{a}_{i}\overset{x^{n}}{\equiv}\prod_{i\in U}%
\mathbf{a}_{i}%
\]
(since the family $\left(  \mathbf{a}_{i}\right)  _{i\in I}$ is multipliable).
Since the relation $\overset{x^{n}}{\equiv}$ is symmetric, we thus obtain%
\begin{equation}
\prod_{i\in U}\mathbf{a}_{i}\overset{x^{n}}{\equiv}\prod_{i\in I}%
\mathbf{a}_{i}. \label{pf.prop.fps.prod-mulable.b.3}%
\end{equation}
Moreover, $M$ is a finite subset of $I$ satisfying $U\subseteq M$ (since
$M\supseteq U$) and therefore $U\subseteq M\subseteq I$; hence, Proposition
\ref{prop.fps.infprod-approx-xneq} \textbf{(a)} (applied to $U$ and $M$
instead of $M$ and $J$) yields%
\[
\prod_{i\in M}\mathbf{a}_{i}\overset{x^{n}}{\equiv}\prod_{i\in U}%
\mathbf{a}_{i}\overset{x^{n}}{\equiv}\prod_{i\in I}\mathbf{a}_{i}%
\ \ \ \ \ \ \ \ \ \ \left(  \text{by (\ref{pf.prop.fps.prod-mulable.b.3}%
)}\right)  .
\]
Hence,%
\begin{equation}
\prod_{i\in M}\mathbf{a}_{i}\overset{x^{n}}{\equiv}\prod_{i\in I}%
\mathbf{a}_{i} \label{pf.prop.fps.prod-mulable.b.4}%
\end{equation}
(since the relation $\overset{x^{n}}{\equiv}$ is transitive). The same
argument (applied to $\left(  \mathbf{b}_{i}\right)  _{i\in I}$ and $V$
instead of $\left(  \mathbf{a}_{i}\right)  _{i\in I}$ and $U$) yields%
\begin{equation}
\prod_{i\in M}\mathbf{b}_{i}\overset{x^{n}}{\equiv}\prod_{i\in I}%
\mathbf{b}_{i}. \label{pf.prop.fps.prod-mulable.b.5}%
\end{equation}
From (\ref{pf.prop.fps.prod-mulable.b.4}) and
(\ref{pf.prop.fps.prod-mulable.b.5}), we obtain%
\[
\left(  \prod_{i\in M}\mathbf{a}_{i}\right)  \left(  \prod_{i\in M}%
\mathbf{b}_{i}\right)  \overset{x^{n}}{\equiv}\left(  \prod_{i\in I}%
\mathbf{a}_{i}\right)  \left(  \prod_{i\in I}\mathbf{b}_{i}\right)
\]
(by (\ref{eq.thm.fps.xneq.props.b.*}), applied to $a=\prod_{i\in M}%
\mathbf{a}_{i}$ and $b=\prod_{i\in I}\mathbf{a}_{i}$ and $c=\prod_{i\in
M}\mathbf{b}_{i}$ and $d=\prod_{i\in I}\mathbf{b}_{i}$). In view of
\[
\left(  \prod_{i\in M}\mathbf{a}_{i}\right)  \left(  \prod_{i\in M}%
\mathbf{b}_{i}\right)  =\prod_{i\in M}\left(  \mathbf{a}_{i}\mathbf{b}%
_{i}\right)  \ \ \ \ \ \ \ \ \ \ \left(
\begin{array}
[c]{c}%
\text{by the properties of finite products,}\\
\text{since the set }M\text{ is finite}%
\end{array}
\right)  ,
\]
this rewrites as
\[
\prod_{i\in M}\left(  \mathbf{a}_{i}\mathbf{b}_{i}\right)  \overset{x^{n}%
}{\equiv}\left(  \prod_{i\in I}\mathbf{a}_{i}\right)  \left(  \prod_{i\in
I}\mathbf{b}_{i}\right)  .
\]
In other words,
\[
\text{each }m\in\left\{  0,1,\ldots,n\right\}  \text{ satisfies }\left[
x^{m}\right]  \left(  \prod_{i\in M}\left(  \mathbf{a}_{i}\mathbf{b}%
_{i}\right)  \right)  =\left[  x^{m}\right]  \left(  \left(  \prod_{i\in
I}\mathbf{a}_{i}\right)  \left(  \prod_{i\in I}\mathbf{b}_{i}\right)  \right)
\]
(by Definition \ref{def.fps.xneq}). Applying this to $m=n$, we obtain
\[
\left[  x^{n}\right]  \left(  \prod_{i\in M}\left(  \mathbf{a}_{i}%
\mathbf{b}_{i}\right)  \right)  =\left[  x^{n}\right]  \left(  \left(
\prod_{i\in I}\mathbf{a}_{i}\right)  \left(  \prod_{i\in I}\mathbf{b}%
_{i}\right)  \right)  .
\]
In view of (\ref{pf.prop.fps.prod-mulable.b.2}), this rewrites as%
\begin{equation}
\left[  x^{n}\right]  \left(  \prod_{i\in I}\left(  \mathbf{a}_{i}%
\mathbf{b}_{i}\right)  \right)  =\left[  x^{n}\right]  \left(  \left(
\prod_{i\in I}\mathbf{a}_{i}\right)  \left(  \prod_{i\in I}\mathbf{b}%
_{i}\right)  \right)  . \label{pf.prop.fps.prod-mulable.b.at}%
\end{equation}

Now, forget that we fixed $n$. We thus have proved that each $n\in\mathbb{N}$
satisfies (\ref{pf.prop.fps.prod-mulable.b.at}). In other words, each
coefficient of the FPS $\prod_{i\in I}\left(  \mathbf{a}_{i}\mathbf{b}%
_{i}\right)  $ equals the corresponding coefficient of $\left(  \prod_{i\in
I}\mathbf{a}_{i}\right)  \left(  \prod_{i\in I}\mathbf{b}_{i}\right)  $.
Hence, we have%
\[
\prod_{i\in I}\left(  \mathbf{a}_{i}\mathbf{b}_{i}\right)  =\left(
\prod_{i\in I}\mathbf{a}_{i}\right)  \left(  \prod_{i\in I}\mathbf{b}%
_{i}\right)  .
\]
This proves Proposition \ref{prop.fps.prod-mulable} \textbf{(b)}.
\end{proof}
\end{fineprint}

\begin{fineprint}
In order to prove Proposition \ref{prop.fps.div-mulable}, we need a lemma (an
analogue of Lemma \ref{lem.fps.prod.irlv.cong-mul} for division instead of multiplication):

\begin{lemma}
\label{lem.fps.prod.irlv.cong-div}Let $a,b,c,d\in K\left[  \left[  x\right]
\right]  $ be four FPSs such that $c$ and $d$ are invertible. Let
$n\in\mathbb{N}$. Assume that
\[
\left[  x^{m}\right]  a=\left[  x^{m}\right]  b\ \ \ \ \ \ \ \ \ \ \text{for
each }m\in\left\{  0,1,\ldots,n\right\}  .
\]
Assume further that%
\[
\left[  x^{m}\right]  c=\left[  x^{m}\right]  d\ \ \ \ \ \ \ \ \ \ \text{for
each }m\in\left\{  0,1,\ldots,n\right\}  .
\]
Then,
\[
\left[  x^{m}\right]  \dfrac{a}{c}=\left[  x^{m}\right]  \dfrac{b}%
{d}\ \ \ \ \ \ \ \ \ \ \text{for each }m\in\left\{  0,1,\ldots,n\right\}  .
\]

\end{lemma}

\begin{proof}
[Proof of Lemma \ref{lem.fps.prod.irlv.cong-div}.]We have assumed that
\[
\left[  x^{m}\right]  a=\left[  x^{m}\right]  b\ \ \ \ \ \ \ \ \ \ \text{for
each }m\in\left\{  0,1,\ldots,n\right\}  .
\]
In other words, $a\overset{x^{n}}{\equiv}b$ (by the definition of the relation
$\overset{x^{n}}{\equiv}$). Similarly, $c\overset{x^{n}}{\equiv}d$. Hence,
(\ref{eq.thm.fps.xneq.props.e./}) yields $\dfrac{a}{c}\overset{x^{n}}{\equiv
}\dfrac{b}{d}$. In other words,%
\[
\left[  x^{m}\right]  \dfrac{a}{c}=\left[  x^{m}\right]  \dfrac{b}%
{d}\ \ \ \ \ \ \ \ \ \ \text{for each }m\in\left\{  0,1,\ldots,n\right\}
\]
(by the definition of the relation $\overset{x^{n}}{\equiv}$). This proves
Lemma \ref{lem.fps.prod.irlv.cong-div}.
\end{proof}

Now we can easily prove Proposition \ref{prop.fps.div-mulable}:

\begin{proof}
[Detailed proof of Proposition \ref{prop.fps.div-mulable}.]This can be proved
similarly to Proposition \ref{prop.fps.prod-mulable}, with the obvious changes
(replacing multiplication by division, and requiring each $\mathbf{b}_{i}$ to
be invertible). Of course, instead of Lemma \ref{lem.fps.prod.irlv.cong-mul},
we need to use Lemma \ref{lem.fps.prod.irlv.cong-div}.
\end{proof}
\end{fineprint}

\begin{fineprint}
In order to eventually prove Proposition \ref{prop.fps.prods-mulable-subfams},
we shall first prove a slightly stronger auxiliary statement:

\begin{lemma}
\label{lem.fps.prods-mulable-subfams-appr}Let $\left(  \mathbf{a}_{i}\right)
_{i\in I}\in K\left[  \left[  x\right]  \right]  ^{I}$ be a family of
invertible FPSs. Let $J$ be a subset of $I$. Let $n\in\mathbb{N}$. Let $U$ be
an $x^{n}$-approximator for $\left(  \mathbf{a}_{i}\right)  _{i\in I}$. Then,
$U\cap J$ is an $x^{n}$-approximator for $\left(  \mathbf{a}_{i}\right)
_{i\in J}$.
\end{lemma}

\begin{proof}
[Proof of Lemma \ref{lem.fps.prods-mulable-subfams-appr}.]The set $U$ is an
$x^{n}$-approximator for $\left(  \mathbf{a}_{i}\right)  _{i\in I}$. In other
words, $U$ is a finite subset of $I$ that determines the first $n+1$
coefficients in the product of $\left(  \mathbf{a}_{i}\right)  _{i\in I}$ (by
the definition of an $x^{n}$-approximator). Hence, in particular, $U$ is
finite. Moreover, $U\subseteq I$ (since $U$ is a subset of $I$).

Let $M=U\cap J$. Thus, $M=U\cap J\subseteq U$, so that the set $M$ is finite
(since $U$ is finite). Moreover, $M$ is a subset of $J$ (since $M=U\cap
J\subseteq J$).

Now, let $N$ be a finite subset of $J$ satisfying $M\subseteq N\subseteq J$.
We shall show that
\[
\prod_{i\in N}\mathbf{a}_{i}\overset{x^{n}}{\equiv}\prod_{i\in M}%
\mathbf{a}_{i}.
\]

Indeed, the set $N\cup U$ is finite (since $N$ and $U$ are finite) and is a
subset of $I$ (since $\underbrace{N}_{\subseteq J\subseteq I}\cup
\underbrace{U}_{\subseteq I}\subseteq I\cup I=I$); it also satisfies
$U\subseteq N\cup U\subseteq I$. Now, we recall that $U$ is an $x^{n}%
$-approximator for $\left(  \mathbf{a}_{i}\right)  _{i\in I}$. Hence,
Proposition \ref{prop.fps.infprod-approx-xneq} \textbf{(a)} (applied to $U$
and $N\cup U$ instead of $M$ and $J$) yields
\begin{equation}
\prod_{i\in N\cup U}\mathbf{a}_{i}\overset{x^{n}}{\equiv}\prod_{i\in
U}\mathbf{a}_{i} \label{pf.lem.fps.prods-mulable-subfams-appr.ab}%
\end{equation}
(since $N\cup U$ is a finite subset of $I$ satisfying $U\subseteq N\cup
U\subseteq I$).

On the other hand, we have assumed that $\left(  \mathbf{a}_{i}\right)  _{i\in
I}$ is a family of invertible FPSs. Thus, for each $i\in I$, the FPS
$\mathbf{a}_{i}$ is invertible, so that its inverse $\mathbf{a}_{i}^{-1}$ is
well-defined. We obviously have%
\begin{equation}
\prod_{i\in U\setminus J}\mathbf{a}_{i}^{-1}\overset{x^{n}}{\equiv}\prod_{i\in
U\setminus J}\mathbf{a}_{i}^{-1}
\label{pf.lem.fps.prods-mulable-subfams-appr.cd}%
\end{equation}
(since Theorem \ref{thm.fps.xneq.props} \textbf{(a)} yields that the relation
$\overset{x^{n}}{\equiv}$ is reflexive).

We have now proved (\ref{pf.lem.fps.prods-mulable-subfams-appr.ab}) and
(\ref{pf.lem.fps.prods-mulable-subfams-appr.cd}). Hence,
(\ref{eq.thm.fps.xneq.props.b.*}) (applied to $a=\prod_{i\in N\cup
U}\mathbf{a}_{i}$ and $b=\prod_{i\in U}\mathbf{a}_{i}$ and $c=\prod_{i\in
U\setminus J}\mathbf{a}_{i}^{-1}$ and $d=\prod_{i\in U\setminus J}%
\mathbf{a}_{i}^{-1}$) yields%
\begin{equation}
\left(  \prod_{i\in N\cup U}\mathbf{a}_{i}\right)  \left(  \prod_{i\in
U\setminus J}\mathbf{a}_{i}^{-1}\right)  \overset{x^{n}}{\equiv}\left(
\prod_{i\in U}\mathbf{a}_{i}\right)  \left(  \prod_{i\in U\setminus
J}\mathbf{a}_{i}^{-1}\right)  .
\label{pf.lem.fps.prods-mulable-subfams-appr.ac=bd}%
\end{equation}

On the other hand, the sets $N$ and $U\setminus J$ are disjoint\footnote{since
$\underbrace{N}_{\subseteq J}\cap\left(  U\setminus J\right)  \subseteq
J\cap\left(  U\setminus J\right)  =\varnothing$ and thus $N\cap\left(
U\setminus J\right)  =\varnothing$}, and their union is $N\cup\left(
U\setminus J\right)  =N\cup U$\ \ \ \ \footnote{This follows from%
\[
N\cup\underbrace{U}_{=\left(  U\cap J\right)  \cup\left(  U\setminus J\right)
}=N\cup\underbrace{\left(  U\cap J\right)  }_{=M}\cup\left(  U\setminus
J\right)  \subseteq\underbrace{N\cup M}_{\substack{=N\\\text{(since
}M\subseteq N\text{)}}}\cup\left(  U\setminus J\right)  =N\cup\left(
U\setminus J\right)  .
\]
}. Hence, the set $N\cup U$ is the union of its two disjoint subsets $N$ and
$U\setminus J$. Thus, we can split the product $\prod_{i\in N\cup U}%
\mathbf{a}_{i}$ as follows:%
\[
\prod_{i\in N\cup U}\mathbf{a}_{i}=\left(  \prod_{i\in N}\mathbf{a}%
_{i}\right)  \left(  \prod_{i\in U\setminus J}\mathbf{a}_{i}\right)  .
\]
Multiplying both sides of this equality by $\prod_{i\in U\setminus
J}\mathbf{a}_{i}^{-1}$, we obtain%
\begin{align}
\left(  \prod_{i\in N\cup U}\mathbf{a}_{i}\right)  \left(  \prod_{i\in
U\setminus J}\mathbf{a}_{i}^{-1}\right)   &  =\left(  \prod_{i\in N}%
\mathbf{a}_{i}\right)  \underbrace{\left(  \prod_{i\in U\setminus J}%
\mathbf{a}_{i}\right)  \left(  \prod_{i\in U\setminus J}\mathbf{a}_{i}%
^{-1}\right)  }_{=\prod_{i\in U\setminus J}\left(  \mathbf{a}_{i}%
\mathbf{a}_{i}^{-1}\right)  }\nonumber\\
&  =\left(  \prod_{i\in N}\mathbf{a}_{i}\right)  \left(  \prod_{i\in
U\setminus J}\underbrace{\left(  \mathbf{a}_{i}\mathbf{a}_{i}^{-1}\right)
}_{=1}\right)  =\left(  \prod_{i\in N}\mathbf{a}_{i}\right)
\underbrace{\left(  \prod_{i\in U\setminus J}1\right)  }_{=1}\nonumber\\
&  =\prod_{i\in N}\mathbf{a}_{i}. \label{pf.prop.fps.prods-mulable-subfams.7a}%
\end{align}

However, the set $U$ is the union of its two disjoint subsets $U\cap J$ and
$U\setminus J$. Thus, we can split the product $\prod_{i\in U}\mathbf{a}_{i}$
as follows:%
\[
\prod_{i\in U}\mathbf{a}_{i}=\left(  \prod_{i\in U\cap J}\mathbf{a}%
_{i}\right)  \left(  \prod_{i\in U\setminus J}\mathbf{a}_{i}\right)  .
\]
Multiplying both sides of this equality by $\prod_{i\in U\setminus
J}\mathbf{a}_{i}^{-1}$, we obtain%
\begin{align}
\left(  \prod_{i\in U}\mathbf{a}_{i}\right)  \left(  \prod_{i\in U\setminus
J}\mathbf{a}_{i}^{-1}\right)   &  =\underbrace{\left(  \prod_{i\in U\cap
J}\mathbf{a}_{i}\right)  }_{\substack{=\prod_{i\in M}\mathbf{a}_{i}%
\\\text{(since }U\cap J=M\text{)}}}\underbrace{\left(  \prod_{i\in U\setminus
J}\mathbf{a}_{i}\right)  \left(  \prod_{i\in U\setminus J}\mathbf{a}_{i}%
^{-1}\right)  }_{=\prod_{i\in U\setminus J}\left(  \mathbf{a}_{i}%
\mathbf{a}_{i}^{-1}\right)  }\nonumber\\
&  =\left(  \prod_{i\in M}\mathbf{a}_{i}\right)  \left(  \prod_{i\in
U\setminus J}\underbrace{\left(  \mathbf{a}_{i}\mathbf{a}_{i}^{-1}\right)
}_{=1}\right)  =\left(  \prod_{i\in M}\mathbf{a}_{i}\right)
\underbrace{\left(  \prod_{i\in U\setminus J}1\right)  }_{=1}\nonumber\\
&  =\prod_{i\in M}\mathbf{a}_{i}. \label{pf.prop.fps.prods-mulable-subfams.7b}%
\end{align}

In view of (\ref{pf.prop.fps.prods-mulable-subfams.7a}) and
(\ref{pf.prop.fps.prods-mulable-subfams.7b}), we can rewrite the relation
(\ref{pf.lem.fps.prods-mulable-subfams-appr.ac=bd}) as follows:%
\[
\prod_{i\in N}\mathbf{a}_{i}\overset{x^{n}}{\equiv}\prod_{i\in M}%
\mathbf{a}_{i}.
\]
In other words, each $m\in\left\{  0,1,\ldots,n\right\}  $ satisfies
\begin{equation}
\left[  x^{m}\right]  \left(  \prod_{i\in N}\mathbf{a}_{i}\right)  =\left[
x^{m}\right]  \left(  \prod_{i\in M}\mathbf{a}_{i}\right)
\label{pf.lem.fps.prods-mulable-subfams-appr.at}%
\end{equation}
(by Definition \ref{def.fps.xneq}).

Forget that we fixed $N$. We have thus shown that every finite subset $N$ of
$J$ satisfying $M\subseteq N\subseteq J$ satisfies
(\ref{pf.lem.fps.prods-mulable-subfams-appr.at}) for each $m\in\left\{
0,1,\ldots,n\right\}  $.

Now, let $m\in\left\{  0,1,\ldots,n\right\}  $. Then, every finite subset $N$
of $J$ satisfying $M\subseteq N\subseteq J$ satisfies%
\[
\left[  x^{m}\right]  \left(  \prod_{i\in N}\mathbf{a}_{i}\right)  =\left[
x^{m}\right]  \left(  \prod_{i\in M}\mathbf{a}_{i}\right)
\]
(by (\ref{pf.lem.fps.prods-mulable-subfams-appr.at})). In other words, the set
$M$ determines the $x^{m}$-coefficient in the product of $\left(
\mathbf{a}_{i}\right)  _{i\in J}$ (by the definition of \textquotedblleft
determining the $x^{m}$-coefficient in a product\textquotedblright).

Forget that we fixed $m$. We thus have shown that $M$ determines the $x^{m}%
$-coefficient in the product of $\left(  \mathbf{a}_{i}\right)  _{i\in J}$ for
each $m\in\left\{  0,1,\ldots,n\right\}  $. In other words, $M$ determines the
first $n+1$ coefficients in the product of $\left(  \mathbf{a}_{i}\right)
_{i\in J}$. In other words, $M$ is an $x^{n}$-approximator for $\left(
\mathbf{a}_{i}\right)  _{i\in J}$ (by the definition of an \textquotedblleft%
$x^{n}$-approximator\textquotedblright, since $M$ is a finite subset of $J$).
In other words, $U\cap J$ is an $x^{n}$-approximator for $\left(
\mathbf{a}_{i}\right)  _{i\in J}$ (since $M=U\cap J$). This proves Lemma
\ref{lem.fps.prods-mulable-subfams-appr}.
\end{proof}
\end{fineprint}

\begin{fineprint}
\begin{proof}
[Detailed proof of Proposition \ref{prop.fps.prods-mulable-subfams}.]Let $J$
be a subset of $I$. We shall show that the family $\left(  \mathbf{a}%
_{i}\right)  _{i\in J}$ is multipliable.

Fix $n\in\mathbb{N}$. We know that the family $\left(  \mathbf{a}_{i}\right)
_{i\in I}$ is multipliable. Hence, there exists an $x^{n}$-approximator $U$
for $\left(  \mathbf{a}_{i}\right)  _{i\in I}$ (by Lemma
\ref{lem.fps.mulable.approx}). Consider this $U$.

Set $M=U\cap J$. Lemma \ref{lem.fps.prods-mulable-subfams-appr} yields that
$U\cap J$ is an $x^{n}$-approximator for $\left(  \mathbf{a}_{i}\right)
_{i\in J}$. In other words, $M$ is an $x^{n}$-approximator for $\left(
\mathbf{a}_{i}\right)  _{i\in J}$ (since $M=U\cap J$). In other words, $M$ is
a finite subset of $J$ that determines the first $n+1$ coefficients in the
product of $\left(  \mathbf{a}_{i}\right)  _{i\in J}$ (by the definition of an
$x^{n}$-approximator). Thus, the set $M$ determines the first $n+1$
coefficients in the product of $\left(  \mathbf{a}_{i}\right)  _{i\in J}$.
Hence, in particular, this set $M$ determines the $x^{n}$-coefficient in the
product of $\left(  \mathbf{a}_{i}\right)  _{i\in J}$. Therefore, the $x^{n}%
$-coefficient in the product of $\left(  \mathbf{a}_{i}\right)  _{i\in J}$ is
finitely determined (by the definition of \textquotedblleft finitely
determined\textquotedblright, since $M$ is a finite subset of $J$).

Forget that we fixed $n$. We thus have proved that for each $n\in\mathbb{N}$,
the $x^{n}$-coefficient in the product of $\left(  \mathbf{a}_{i}\right)
_{i\in J}$ is finitely determined. In other words, each coefficient in the
product of $\left(  \mathbf{a}_{i}\right)  _{i\in J}$ is finitely determined.
In other words, the family $\left(  \mathbf{a}_{i}\right)  _{i\in J}$ is
multipliable (by the definition of \textquotedblleft
multipliable\textquotedblright).

Forget that we fixed $J$. We thus have shown that the family $\left(
\mathbf{a}_{i}\right)  _{i\in J}$ is multipliable whenever $J$ is a subset of
$I$. In other words, any subfamily of $\left(  \mathbf{a}_{i}\right)  _{i\in
I}$ is multipliable. This proves Proposition
\ref{prop.fps.prods-mulable-subfams}.
\end{proof}
\end{fineprint}

\begin{fineprint}
Our next goal is to prove Proposition \ref{prop.fps.prods-mulable-rules.SW1}.
First, however, let us restate a piece of Theorem \ref{thm.fps.xneq.props}
\textbf{(f)} in more convenient language:

\begin{lemma}
\label{lem.fps.prods-mulable-rules.SW1.lem1}Let $n\in\mathbb{N}$. Let $V$ be a
finite set. Let $\left(  c_{w}\right)  _{w\in V}\in K\left[  \left[  x\right]
\right]  ^{V}$ and $\left(  d_{w}\right)  _{w\in V}\in K\left[  \left[
x\right]  \right]  ^{V}$ be two families of FPSs such that
\[
\text{each }w\in V\text{ satisfies }c_{w}\overset{x^{n}}{\equiv}d_{w}.
\]
Then, we have%
\[
\prod_{w\in V}c_{w}\overset{x^{n}}{\equiv}\prod_{w\in V}d_{w}.
\]

\end{lemma}

\begin{proof}
[Proof of Lemma \ref{lem.fps.prods-mulable-rules.SW1.lem1}.]This is just
(\ref{eq.thm.fps.xneq.props.e.*}), with the letters $S$, $s$, $a_{s}$ and
$b_{s}$ renamed as $V$, $w$, $c_{w}$ and $d_{w}$.
\end{proof}

\begin{proof}
[Detailed proof of Proposition \ref{prop.fps.prods-mulable-rules.SW1}.]We
shall subdivide our proof into several claims:

\begin{statement}
\textit{Claim 1:} Let $w\in W$. Then, the family $\left(  \mathbf{a}%
_{s}\right)  _{s\in S;\ f\left(  s\right)  =w}$ is multipliable.
\end{statement}

[\textit{Proof of Claim 1:} This was an assumption of Proposition
\ref{prop.fps.prods-mulable-rules.SW1}.] \medskip

Let us set%
\begin{equation}
\mathbf{b}_{w}:=\prod_{\substack{s\in S;\\f\left(  s\right)  =w}%
}\mathbf{a}_{s}\ \ \ \ \ \ \ \ \ \ \text{for each }w\in W.
\label{pf.prop.fps.prods-mulable-rules.SW1.bw=}%
\end{equation}
This is well-defined, because for each $w\in W$, the product $\prod
_{\substack{s\in S;\\f\left(  s\right)  =w}}\mathbf{a}_{s}$ is well-defined
(since Claim 1 shows that the family $\left(  \mathbf{a}_{s}\right)  _{s\in
S;\ f\left(  s\right)  =w}$ is multipliable).

For each $w\in W$, define the subset
\[
S_{w}:=\left\{  s\in S\ \mid\ f\left(  s\right)  =w\right\}
\]
of $S$. Then, for each $w\in W$, the family $\left(  \mathbf{a}_{s}\right)
_{s\in S_{w}}=\left(  \mathbf{a}_{s}\right)  _{s\in S;\ f\left(  s\right)
=w}$ is multipliable (by Claim 1), and we have%
\begin{equation}
\mathbf{b}_{w}=\prod_{\substack{s\in S;\\f\left(  s\right)  =w}}\mathbf{a}%
_{s}=\prod_{s\in S_{w}}\mathbf{a}_{s}
\label{pf.prop.fps.prods-mulable-rules.SW1.bw=2}%
\end{equation}
(since $\left\{  s\in S\ \mid\ f\left(  s\right)  =w\right\}  =S_{w}$).
Moreover, we have%
\[
S=\bigcup_{w\in W}S_{w}%
\]
(since each $s\in S$ satisfies $f\left(  s\right)  =w$ for some $w\in W$, and
thus belongs to $S_{w}$ for this $w$).

Now, let $n\in\mathbb{N}$. Lemma \ref{lem.fps.mulable.approx} (applied to $S$
and $\left(  \mathbf{a}_{s}\right)  _{s\in S}$ instead of $I$ and $\left(
\mathbf{a}_{i}\right)  _{i\in I}$) shows that there exists an $x^{n}%
$-approximator for $\left(  \mathbf{a}_{s}\right)  _{s\in S}$. Pick such an
$x^{n}$-approximator, and call it $U$. Then, $U$ is an $x^{n}$-approximator
for $\left(  \mathbf{a}_{s}\right)  _{s\in S}$; in other words, $U$ is a
finite subset of $S$ that determines the first $n+1$ coefficients in the
product of $\left(  \mathbf{a}_{s}\right)  _{s\in S}$ (by the definition of an
$x^{n}$-approximator).

The set $U$ is finite. Thus, its image $f\left(  U\right)  =\left\{  f\left(
u\right)  \ \mid\ u\in U\right\}  $ is finite as well. Now, we claim the following:

\begin{statement}
\textit{Claim 2:} Let $V$ be a finite subset of $W$ such that $f\left(
U\right)  \subseteq V$. Then,%
\[
\prod_{w\in V}\mathbf{b}_{w}\overset{x^{n}}{\equiv}\prod_{s\in U}%
\mathbf{a}_{s}.
\]

\end{statement}

[\textit{Proof of Claim 2:} Let $w\in V$. Then, recall that $S_{w}=\left\{
s\in S\ \mid\ f\left(  s\right)  =w\right\}  $. Hence, the family $\left(
\mathbf{a}_{s}\right)  _{s\in S_{w}}=\left(  \mathbf{a}_{s}\right)  _{s\in
S;\ f\left(  s\right)  =w}$ is multipliable (by Claim 1), so Lemma
\ref{lem.fps.mulable.approx} (applied to this family $\left(  \mathbf{a}%
_{s}\right)  _{s\in S_{w}}$ instead of $\left(  \mathbf{a}_{i}\right)  _{i\in
I}$) yields that there exists an $x^{n}$-approximator for this family $\left(
\mathbf{a}_{s}\right)  _{s\in S_{w}}$. Let $D_{w}$ be this $x^{n}%
$-approximator. Thus, $D_{w}$ is a finite subset of $S_{w}$ (by the definition
of an $x^{n}$-approximator). Define the further set%
\[
E_{w}:=D_{w}\cup\left(  U\cap S_{w}\right)  .
\]
Then, the set $E_{w}$ is a subset of $S_{w}$ (since both $D_{w}$ and $U\cap
S_{w}$ are subsets of $S_{w}$) and is finite (since $D_{w}$ and $U$ are
finite) and satisfies $D_{w}\subseteq E_{w}$ (by the definition of $E_{w}$).
Hence, Proposition \ref{prop.fps.infprod-approx-xneq} \textbf{(a)} (applied to
$\left(  \mathbf{a}_{s}\right)  _{s\in S_{w}}$, $D_{w}$ and $E_{w}$ instead of
$\left(  \mathbf{a}_{i}\right)  _{i\in I}$, $M$ and $J$) yields%
\begin{equation}
\prod_{s\in E_{w}}\mathbf{a}_{s}\overset{x^{n}}{\equiv}\prod_{s\in D_{w}%
}\mathbf{a}_{s}. \label{pf.prop.fps.prods-mulable-rules.SW1.c24a}%
\end{equation}
Furthermore, Proposition \ref{prop.fps.infprod-approx-xneq} \textbf{(b)}
(applied to $\left(  \mathbf{a}_{s}\right)  _{s\in S_{w}}$ and $D_{w}$ instead
of $\left(  \mathbf{a}_{i}\right)  _{i\in I}$ and $M$) yields%
\begin{equation}
\prod_{s\in S_{w}}\mathbf{a}_{s}\overset{x^{n}}{\equiv}\prod_{s\in D_{w}%
}\mathbf{a}_{s} \label{pf.prop.fps.prods-mulable-rules.SW1.c24b}%
\end{equation}
(since the family $\left(  \mathbf{a}_{s}\right)  _{s\in S_{w}}$ is
multipliable). Since the relation $\overset{x^{n}}{\equiv}$ is an equivalence
relation (by Theorem \ref{thm.fps.xneq.props} \textbf{(a)}), we can combine
(\ref{pf.prop.fps.prods-mulable-rules.SW1.c24a}) with
(\ref{pf.prop.fps.prods-mulable-rules.SW1.c24b}) to obtain%
\begin{equation}
\prod_{s\in E_{w}}\mathbf{a}_{s}\overset{x^{n}}{\equiv}\prod_{s\in S_{w}%
}\mathbf{a}_{s}. \label{pf.prop.fps.prods-mulable-rules.SW1.c24c}%
\end{equation}

Forget that we fixed $w$. We thus have constructed a finite subset $E_{w}$ of
$S_{w}$ satisfying $\prod_{s\in E_{w}}\mathbf{a}_{s}\overset{x^{n}}{\equiv
}\prod_{s\in S_{w}}\mathbf{a}_{s}$ for each $w\in V$. Hence, Lemma
\ref{lem.fps.prods-mulable-rules.SW1.lem1} (applied to $c_{w}=\prod_{s\in
E_{w}}\mathbf{a}_{s}$ and $d_{w}=\prod_{s\in S_{w}}\mathbf{a}_{s}$) yields
\begin{equation}
\prod_{w\in V}\ \ \prod_{s\in E_{w}}\mathbf{a}_{s}\overset{x^{n}}{\equiv}%
\prod_{w\in V}\ \ \prod_{s\in S_{w}}\mathbf{a}_{s}.
\label{pf.prop.fps.prods-mulable-rules.SW1.c22}%
\end{equation}

Let $E$ be the union $\bigcup\limits_{w\in V}E_{w}$. This is a finite union of
finite sets (since $V$ and each $E_{w}$ are finite), hence itself finite.
Furthermore,%
\[
E=\bigcup\limits_{w\in V}\underbrace{E_{w}}_{\subseteq S_{w}\subseteq
S}\subseteq\bigcup\limits_{w\in V}S\subseteq S,
\]
and it is easy to see that
\[
U\subseteq E
\]
\footnote{\textit{Proof.} Let $u\in U$. Then, $f\left(  u\right)  \in f\left(
U\right)  \subseteq V$. Moreover, $u$ is an element of $S$ (since $u\in
U\subseteq S$) and satisfies $f\left(  u\right)  =f\left(  u\right)  $. Hence,
$u\in\left\{  s\in S\ \mid\ f\left(  s\right)  =f\left(  u\right)  \right\}
=S_{f\left(  u\right)  }$ (by the definition of $S_{f\left(  u\right)  }$).
Thus, $u\in U\cap S_{f\left(  u\right)  }$ (since $u\in U$ and $u\in
S_{f\left(  u\right)  }$). But the definition of $E_{f\left(  u\right)  }$
yields $E_{f\left(  u\right)  }=D_{f\left(  u\right)  }\cup\left(  U\cap
S_{f\left(  u\right)  }\right)  \supseteq U\cap S_{f\left(  u\right)  }$, so
$U\cap S_{f\left(  u\right)  }\subseteq E_{f\left(  u\right)  }$. But
$f\left(  u\right)  \in V$ and thus $E_{f\left(  u\right)  }\subseteq
\bigcup\limits_{w\in V}E_{w}=E$. Hence, $u\in U\cap S_{f\left(  u\right)
}\subseteq E_{f\left(  u\right)  }\subseteq E$.
\par
Forget that we fixed $u$. We thus have shown that $u\in E$ for each $u\in U$.
In other words, $U\subseteq E$.}.

The sets $E_{w}$ for distinct $w\in V$ are disjoint (indeed, each $E_{w}$ is a
subset of $S_{w}=\left\{  s\in S\ \mid\ f\left(  s\right)  =w\right\}  $;
thus, all elements $s\in E_{w}$ satisfy $f\left(  s\right)  =w$, which
uniquely determines $w$; this shows that a single element $s$ cannot belong to
$E_{w}$ for two different values of $w$ simultaneously). Thus, by the
properties of finite products, we have
\[
\prod_{s\in\bigcup\limits_{w\in V}E_{w}}\mathbf{a}_{s}=\prod_{w\in V}%
\ \ \prod_{s\in E_{w}}\mathbf{a}_{s}.
\]
In view of $E=\bigcup\limits_{w\in V}E_{w}$, we can rewrite this as%
\begin{align*}
\prod_{s\in E}\mathbf{a}_{s}  &  =\prod_{w\in V}\ \ \prod_{s\in E_{w}%
}\mathbf{a}_{s}\overset{x^{n}}{\equiv}\prod_{w\in V}\ \ \underbrace{\prod
_{s\in S_{w}}\mathbf{a}_{s}}_{\substack{=\mathbf{b}_{w}\\\text{(by
(\ref{pf.prop.fps.prods-mulable-rules.SW1.bw=2}))}}%
}\ \ \ \ \ \ \ \ \ \ \left(  \text{by
(\ref{pf.prop.fps.prods-mulable-rules.SW1.c22})}\right) \\
&  =\prod_{w\in V}\mathbf{b}_{w}.
\end{align*}
Since the relation $\overset{x^{n}}{\equiv}$ on $K\left[  \left[  x\right]
\right]  $ is symmetric (by Theorem \ref{thm.fps.xneq.props} \textbf{(a)}),
this yields%
\begin{equation}
\prod_{w\in V}\mathbf{b}_{w}\overset{x^{n}}{\equiv}\prod_{s\in E}%
\mathbf{a}_{s}. \label{pf.prop.fps.prods-mulable-rules.SW1.c27}%
\end{equation}

But $U$ is an $x^{n}$-approximator for $\left(  \mathbf{a}_{s}\right)  _{s\in
S}$, whereas $E$ is a finite subset of $S$ (since $E\subseteq S$) satisfying
$U\subseteq E\subseteq S$. Hence, Proposition
\ref{prop.fps.infprod-approx-xneq} \textbf{(a)} (applied to the family
$\left(  \mathbf{a}_{i}\right)  _{i\in I}=\left(  \mathbf{a}_{s}\right)
_{s\in S}$ with $M=U$ and $J=E$) gives%
\begin{equation}
\prod_{s\in E}\mathbf{a}_{s}\overset{x^{n}}{\equiv}\prod_{s\in U}%
\mathbf{a}_{s}. \label{pf.prop.fps.prods-mulable-rules.SW1.c28}%
\end{equation}
Now, recall that the relation $\overset{x^{n}}{\equiv}$ on $K\left[  \left[
x\right]  \right]  $ is transitive (by Theorem \ref{thm.fps.xneq.props}
\textbf{(a)}). Hence, from (\ref{pf.prop.fps.prods-mulable-rules.SW1.c27}) and
(\ref{pf.prop.fps.prods-mulable-rules.SW1.c28}), we obtain $\prod_{w\in
V}\mathbf{b}_{w}\overset{x^{n}}{\equiv}\prod_{s\in U}\mathbf{a}_{s}$. This
proves Claim 2.]

\begin{statement}
\textit{Claim 3:} The set $f\left(  U\right)  $ is an $x^{n}$-approximator for
the family $\left(  \mathbf{b}_{w}\right)  _{w\in W}$.
\end{statement}

[\textit{Proof of Claim 3:} Let $V$ be a finite subset of $W$ satisfying
$f\left(  U\right)  \subseteq V\subseteq W$. We shall show that%
\[
\prod_{w\in f\left(  U\right)  }\mathbf{b}_{w}\overset{x^{n}}{\equiv}%
\prod_{w\in V}\mathbf{b}_{w}.
\]

Indeed, Claim 2 yields%
\begin{equation}
\prod_{w\in V}\mathbf{b}_{w}\overset{x^{n}}{\equiv}\prod_{s\in U}%
\mathbf{a}_{s}. \label{pf.prop.fps.prods-mulable-rules.SW1.c3.pf.4}%
\end{equation}

The same argument (applied to $f\left(  U\right)  $ instead of $V$) shows that%
\begin{equation}
\prod_{w\in f\left(  U\right)  }\mathbf{b}_{w}\overset{x^{n}}{\equiv}%
\prod_{s\in U}\mathbf{a}_{s}
\label{pf.prop.fps.prods-mulable-rules.SW1.c3.pf.5}%
\end{equation}
(since $f\left(  U\right)  $ is a finite subset of $W$ satisfying $f\left(
U\right)  \subseteq f\left(  U\right)  \subseteq W$).

However, the relation $\overset{x^{n}}{\equiv}$ on $K\left[  \left[  x\right]
\right]  $ is symmetric (by Theorem \ref{thm.fps.xneq.props} \textbf{(a)});
thus, (\ref{pf.prop.fps.prods-mulable-rules.SW1.c3.pf.4}) entails%
\[
\prod_{s\in U}\mathbf{a}_{s}\overset{x^{n}}{\equiv}\prod_{w\in V}%
\mathbf{b}_{w}.
\]
Therefore, (\ref{pf.prop.fps.prods-mulable-rules.SW1.c3.pf.5}) becomes
\[
\prod_{w\in f\left(  U\right)  }\mathbf{b}_{w}\overset{x^{n}}{\equiv}%
\prod_{s\in U}\mathbf{a}_{s}\overset{x^{n}}{\equiv}\prod_{w\in V}%
\mathbf{b}_{w}.
\]
Since the relation $\overset{x^{n}}{\equiv}$ on $K\left[  \left[  x\right]
\right]  $ is transitive (by Theorem \ref{thm.fps.xneq.props} \textbf{(a)}),
we thus obtain%
\[
\prod_{w\in f\left(  U\right)  }\mathbf{b}_{w}\overset{x^{n}}{\equiv}%
\prod_{w\in V}\mathbf{b}_{w}.
\]
In other words, each $m\in\left\{  0,1,\ldots,n\right\}  $ satisfies%
\begin{equation}
\left[  x^{m}\right]  \left(  \prod_{w\in f\left(  U\right)  }\mathbf{b}%
_{w}\right)  =\left[  x^{m}\right]  \left(  \prod_{w\in V}\mathbf{b}%
_{w}\right)  \label{pf.prop.fps.prods-mulable-rules.SW1.c3.pf.at}%
\end{equation}
(by Definition \ref{def.fps.xneq}).

Forget that we fixed $V$. We thus have shown that if $V$ is a finite subset of
$W$ satisfying $f\left(  U\right)  \subseteq V\subseteq W$, then each
$m\in\left\{  0,1,\ldots,n\right\}  $ satisfies
(\ref{pf.prop.fps.prods-mulable-rules.SW1.c3.pf.at}).

Now, let $m\in\left\{  0,1,\ldots,n\right\}  $ be arbitrary. Then, every
finite subset $V$ of $W$ satisfying $f\left(  U\right)  \subseteq V\subseteq
W$ satisfies%
\[
\left[  x^{m}\right]  \left(  \prod_{w\in f\left(  U\right)  }\mathbf{b}%
_{w}\right)  =\left[  x^{m}\right]  \left(  \prod_{w\in V}\mathbf{b}%
_{w}\right)  \ \ \ \ \ \ \ \ \ \ \left(  \text{by
(\ref{pf.prop.fps.prods-mulable-rules.SW1.c3.pf.at})}\right)  .
\]
In other words, the set $f\left(  U\right)  $ determines the $x^{m}%
$-coefficient in the product of $\left(  \mathbf{b}_{w}\right)  _{w\in W}$ (by
the definition of \textquotedblleft determining the $x^{m}$-coefficient in a
product\textquotedblright, since $f\left(  U\right)  $ is a finite subset of
$W$).

Forget that we fixed $m$. We thus have shown that the set $f\left(  U\right)
$ determines the $x^{m}$-coefficient in the product of $\left(  \mathbf{b}%
_{w}\right)  _{w\in W}$ for each $m\in\left\{  0,1,\ldots,n\right\}  $. In
other words, the set $f\left(  U\right)  $ determines the first $n+1$
coefficients in the product of $\left(  \mathbf{b}_{w}\right)  _{w\in W}$. In
other words, $f\left(  U\right)  $ is an $x^{n}$-approximator for $\left(
\mathbf{b}_{w}\right)  _{w\in W}$ (by the definition of an \textquotedblleft%
$x^{n}$-approximator\textquotedblright, since $f\left(  U\right)  $ is a
finite subset of $W$). This proves Claim 3.]

\begin{statement}
\textit{Claim 4:} The $x^{n}$-coefficient in the product of $\left(
\mathbf{b}_{w}\right)  _{w\in W}$ is finitely determined.
\end{statement}

[\textit{Proof of Claim 4:} Claim 3 shows that $f\left(  U\right)  $ is an
$x^{n}$-approximator for $\left(  \mathbf{b}_{w}\right)  _{w\in W}$. Thus, the
set $f\left(  U\right)  $ determines the first $n+1$ coefficients in the
product of $\left(  \mathbf{b}_{w}\right)  _{w\in W}$ (by the definition of an
\textquotedblleft$x^{n}$-approximator\textquotedblright). Hence, in
particular, this set $f\left(  U\right)  $ determines the $x^{n}$-coefficient
in the product of $\left(  \mathbf{b}_{w}\right)  _{w\in W}$. Thus, the
$x^{n}$-coefficient in the product of $\left(  \mathbf{b}_{w}\right)  _{w\in
W}$ is finitely determined (by the definition of \textquotedblleft finitely
determined\textquotedblright, since $f\left(  U\right)  $ is a finite subset
of $W$). This proves Claim 4.] \medskip

Now, forget that we fixed $n$. We thus have shown that the $x^{n}$-coefficient
in the product of $\left(  \mathbf{b}_{w}\right)  _{w\in W}$ is finitely
determined for each $n\in\mathbb{N}$. In other words, each coefficient in the
product of $\left(  \mathbf{b}_{w}\right)  _{w\in W}$ is finitely determined.
In other words, the family $\left(  \mathbf{b}_{w}\right)  _{w\in W}$ is
multipliable (by the definition of \textquotedblleft
multipliable\textquotedblright). In view of
(\ref{pf.prop.fps.prods-mulable-rules.SW1.bw=}), we can restate this as
follows: The family $\left(  \prod_{\substack{s\in S;\\f\left(  s\right)
=w}}\mathbf{a}_{s}\right)  _{w\in W}$ is multipliable. \medskip

It remains to prove the equality (\ref{eq.prop.fps.prods-mulable-rules.SW1.eq}).

Let $n\in\mathbb{N}$. Lemma \ref{lem.fps.mulable.approx} (applied to $S$ and
$\left(  \mathbf{a}_{s}\right)  _{s\in S}$ instead of $I$ and $\left(
\mathbf{a}_{i}\right)  _{i\in I}$) shows that there exists an $x^{n}%
$-approximator for $\left(  \mathbf{a}_{s}\right)  _{s\in S}$. Pick such an
$x^{n}$-approximator, and call it $U$. Then, $U$ is a finite subset of $S$ (by
the definition of an $x^{n}$-approximator). Furthermore, Proposition
\ref{prop.fps.infprod-approx-xneq} \textbf{(b)} (applied to $S$, $\left(
\mathbf{a}_{s}\right)  _{s\in S}$ and $U$ instead of $I$, $\left(
\mathbf{a}_{i}\right)  _{i\in I}$ and $M$) yields%
\[
\prod_{i\in S}\mathbf{a}_{i}\overset{x^{n}}{\equiv}\prod_{i\in U}%
\mathbf{a}_{i}.
\]
Renaming the index $i$ as $s$ on both sides of this relation, we obtain%
\begin{equation}
\prod_{s\in S}\mathbf{a}_{s}\overset{x^{n}}{\equiv}\prod_{s\in U}%
\mathbf{a}_{s}. \label{pf.prop.fps.prods-mulable-rules.end.1}%
\end{equation}
However, the relation $\overset{x^{n}}{\equiv}$ on $K\left[  \left[  x\right]
\right]  $ is symmetric (by Theorem \ref{thm.fps.xneq.props} \textbf{(a)});
thus, (\ref{pf.prop.fps.prods-mulable-rules.end.1}) entails%
\begin{equation}
\prod_{s\in U}\mathbf{a}_{s}\overset{x^{n}}{\equiv}\prod_{s\in S}%
\mathbf{a}_{s}. \label{pf.prop.fps.prods-mulable-rules.end.1b}%
\end{equation}

Let $V$ be the set $f\left(  U\right)  =\left\{  f\left(  u\right)
\ \mid\ u\in U\right\}  $. Thus, $V=f\left(  U\right)  \subseteq W$ and
$f\left(  U\right)  =V\subseteq V$, so that $f\left(  U\right)  \subseteq
V\subseteq W$. Moreover, the set $f\left(  U\right)  $ is finite (since $U$ is
finite); in other words, the set $V$ is finite (since $V=f\left(  U\right)  $).

We have already proved (in Claim 3) that the set $f\left(  U\right)  $ is an
$x^{n}$-approximator for the family $\left(  \mathbf{b}_{w}\right)  _{w\in W}%
$. In other words, the set $V$ is an $x^{n}$-approximator for the family
$\left(  \mathbf{b}_{w}\right)  _{w\in W}$ (since $V=f\left(  U\right)  $).
Hence, Proposition \ref{prop.fps.infprod-approx-xneq} \textbf{(b)} (applied to
$W$, $\left(  \mathbf{b}_{w}\right)  _{w\in W}$ and $V$ instead of $I$,
$\left(  \mathbf{a}_{i}\right)  _{i\in I}$ and $M$) yields%
\[
\prod_{i\in W}\mathbf{b}_{i}\overset{x^{n}}{\equiv}\prod_{i\in V}%
\mathbf{b}_{i}%
\]
(since the family $\left(  \mathbf{b}_{w}\right)  _{w\in W}$ is multipliable).
Renaming the index $i$ as $w$ on both sides of this relation, we obtain%
\begin{equation}
\prod_{w\in W}\mathbf{b}_{w}\overset{x^{n}}{\equiv}\prod_{w\in V}%
\mathbf{b}_{w}. \label{pf.prop.fps.prods-mulable-rules.end.2}%
\end{equation}

On the other hand, Claim 2 yields%
\begin{equation}
\prod_{w\in V}\mathbf{b}_{w}\overset{x^{n}}{\equiv}\prod_{s\in U}%
\mathbf{a}_{s}. \label{pf.prop.fps.prods-mulable-rules.end.3}%
\end{equation}

Now, (\ref{pf.prop.fps.prods-mulable-rules.end.2}) becomes%
\begin{align*}
\prod_{w\in W}\mathbf{b}_{w}  &  \overset{x^{n}}{\equiv}\prod_{w\in
V}\mathbf{b}_{w}\overset{x^{n}}{\equiv}\prod_{s\in U}\mathbf{a}_{s}%
\ \ \ \ \ \ \ \ \ \ \left(  \text{by
(\ref{pf.prop.fps.prods-mulable-rules.end.3})}\right) \\
&  \overset{x^{n}}{\equiv}\prod_{s\in S}\mathbf{a}_{s}%
\ \ \ \ \ \ \ \ \ \ \left(  \text{by
(\ref{pf.prop.fps.prods-mulable-rules.end.1b})}\right)  .
\end{align*}
Since the relation $\overset{x^{n}}{\equiv}$ on $K\left[  \left[  x\right]
\right]  $ is transitive (by Theorem \ref{thm.fps.xneq.props} \textbf{(a)}),
we thus obtain%
\[
\prod_{w\in W}\mathbf{b}_{w}\overset{x^{n}}{\equiv}\prod_{s\in S}%
\mathbf{a}_{s}.
\]
In other words,
\[
\text{each }m\in\left\{  0,1,\ldots,n\right\}  \text{ satisfies }\left[
x^{m}\right]  \left(  \prod_{w\in W}\mathbf{b}_{w}\right)  =\left[
x^{m}\right]  \left(  \prod_{s\in S}\mathbf{a}_{s}\right)
\]
(by Definition \ref{def.fps.xneq}). Applying this to $m=n$, we obtain%
\[
\left[  x^{n}\right]  \left(  \prod_{w\in W}\mathbf{b}_{w}\right)  =\left[
x^{n}\right]  \left(  \prod_{s\in S}\mathbf{a}_{s}\right)  .
\]

Forget that we fixed $n$. We thus have shown that%
\[
\left[  x^{n}\right]  \left(  \prod_{w\in W}\mathbf{b}_{w}\right)  =\left[
x^{n}\right]  \left(  \prod_{s\in S}\mathbf{a}_{s}\right)
\ \ \ \ \ \ \ \ \ \ \text{for each }n\in\mathbb{N}.
\]
In other words, each coefficient of the FPS $\prod_{w\in W}\mathbf{b}_{w}$
equals the corresponding coefficient of $\prod_{s\in S}\mathbf{a}_{s}$.
Therefore,%
\[
\prod_{w\in W}\mathbf{b}_{w}=\prod_{s\in S}\mathbf{a}_{s}.
\]
In view of (\ref{pf.prop.fps.prods-mulable-rules.SW1.bw=}), we can rewrite
this as%
\[
\prod_{w\in W}\ \ \prod_{\substack{s\in S;\\f\left(  s\right)  =w}%
}\mathbf{a}_{s}=\prod_{s\in S}\mathbf{a}_{s}.
\]
Thus, (\ref{eq.prop.fps.prods-mulable-rules.SW1.eq}) is proven, and the proof
of Proposition \ref{prop.fps.prods-mulable-rules.SW1} is complete.
\end{proof}
\end{fineprint}

\begin{fineprint}
\begin{proof}
[Detailed proof of Proposition \ref{prop.fps.prods-mulable-rules.fubini1}.]Let
$f:I\times J\rightarrow I$ be the map that sends each pair $\left(
i,j\right)  $ to $i$. We first prove two easy claims:

\begin{statement}
\textit{Claim 1:} Let $w\in I$. Let $J^{\prime}$ be the subset $\left\{  s\in
I\times J\ \mid\ f\left(  s\right)  =w\right\}  $ of $I\times J$. Then, there
is a bijection%
\begin{align*}
J  &  \rightarrow J^{\prime},\\
j  &  \mapsto\left(  w,j\right)  .
\end{align*}

\end{statement}

\begin{proof}
[Proof of Claim 1.]We have
\begin{align*}
J^{\prime}  &  =\left\{  s\in I\times J\ \mid\ f\left(  s\right)  =w\right\}
\\
&  =\left\{  \left(  i,j\right)  \in I\times J\ \mid\ \underbrace{f\left(
i,j\right)  }_{\substack{=i\\\text{(by the definition of }f\text{)}%
}}=w\right\} \\
&  \ \ \ \ \ \ \ \ \ \ \ \ \ \ \ \ \ \ \ \ \left(  \text{here, we have renamed
the index }s\text{ as }\left(  i,j\right)  \right) \\
&  =\left\{  \left(  i,j\right)  \in I\times J\ \mid\ i=w\right\} \\
&  =\left\{  \left(  w,j\right)  \ \mid\ j\in J\right\}
\ \ \ \ \ \ \ \ \ \ \left(  \text{since }i\in I\right)  .
\end{align*}
In other words, the set $J^{\prime}$ consists of all pairs $\left(
w,j\right)  $ with $j\in J$. Hence, there is a bijection%
\begin{align*}
J  &  \rightarrow J^{\prime},\\
j  &  \mapsto\left(  w,j\right)  .
\end{align*}
This proves Claim 1.
\end{proof}

\begin{statement}
\textit{Claim 2:} For each $w\in I$, the family $\left(  \mathbf{a}%
_{s}\right)  _{s\in I\times J;\ f\left(  s\right)  =w}$ is multipliable.
\end{statement}

\begin{proof}
[Proof of Claim 2.]Let $w\in I$. We have assumed that for each $i\in I$, the
family $\left(  \mathbf{a}_{\left(  i,j\right)  }\right)  _{j\in J}$ is
multipliable. Applying this to $i=w$, we see that the family $\left(
\mathbf{a}_{\left(  w,j\right)  }\right)  _{j\in J}$ is multipliable.

Let $J^{\prime}$ be the subset $\left\{  s\in I\times J\ \mid\ f\left(
s\right)  =w\right\}  $ of $I\times J$. Thus, the family $\left(
\mathbf{a}_{s}\right)  _{s\in J^{\prime}}$ is the family $\left(
\mathbf{a}_{s}\right)  _{s\in I\times J;\ f\left(  s\right)  =w}$.

Claim 1 yields that there is a bijection%
\begin{align*}
J  &  \rightarrow J^{\prime},\\
j  &  \mapsto\left(  w,j\right)  .
\end{align*}
Hence, the family $\left(  \mathbf{a}_{s}\right)  _{s\in J^{\prime}}$ is a
reindexing of the family $\left(  \mathbf{a}_{\left(  w,j\right)  }\right)
_{j\in J}$. Since the latter family $\left(  \mathbf{a}_{\left(  w,j\right)
}\right)  _{j\in J}$ is multipliable, we thus conclude that the former family
$\left(  \mathbf{a}_{s}\right)  _{s\in J^{\prime}}$ is also multipliable
(since a reindexing of a multipliable family is still
multipliable\footnote{This is part of Proposition
\ref{prop.fps.prods-mulable-rules.reindex}.}). In other words, the family
$\left(  \mathbf{a}_{s}\right)  _{s\in I\times J;\ f\left(  s\right)  =w}$ is
multipliable (since the family $\left(  \mathbf{a}_{s}\right)  _{s\in
J^{\prime}}$ is the family $\left(  \mathbf{a}_{s}\right)  _{s\in I\times
J;\ f\left(  s\right)  =w}$). This proves Claim 2.
\end{proof}

Thanks to Claim 2, we can apply Proposition
\ref{prop.fps.prods-mulable-rules.SW1} to $S=I\times J$ and $W=I$. This
application yields that%
\begin{equation}
\prod_{s\in I\times J}\mathbf{a}_{s}=\prod_{w\in I}\ \ \prod_{\substack{s\in
I\times J;\\f\left(  s\right)  =w}}\mathbf{a}_{s};
\label{pf.prop.fps.prods-mulable-rules.fubini.1}%
\end{equation}
in particular, it yields that the right hand side of
(\ref{pf.prop.fps.prods-mulable-rules.fubini.1}) is well-defined -- i.e., the
family $\left(  \prod_{\substack{s\in I\times J;\\f\left(  s\right)
=w}}\mathbf{a}_{s}\right)  _{w\in I}$ is also multipliable.

Now, fix $w\in I$. Let $J^{\prime}$ be the subset $\left\{  s\in I\times
J\ \mid\ f\left(  s\right)  =w\right\}  $ of $I\times J$. Thus, the family
$\left(  \mathbf{a}_{s}\right)  _{s\in J^{\prime}}$ is the family $\left(
\mathbf{a}_{s}\right)  _{s\in I\times J;\ f\left(  s\right)  =w}$, and
therefore is multipliable (by Claim 2). Hence, the product $\prod_{s\in
J^{\prime}}\mathbf{a}_{s}$ is well-defined. Furthermore, Claim 1 yields that
there is a bijection%
\begin{align*}
J  &  \rightarrow J^{\prime},\\
j  &  \mapsto\left(  w,j\right)  .
\end{align*}
Thus, we can substitute $\left(  w,j\right)  $ for $s$ in the product
$\prod_{s\in J^{\prime}}\mathbf{a}_{s}$ (because any bijection allows us to
substitute the index in a product\footnote{This is just the claim of
Proposition \ref{prop.fps.prods-mulable-rules.reindex}.}). We thus obtain%
\begin{equation}
\prod_{s\in J^{\prime}}\mathbf{a}_{s}=\prod_{j\in J}\mathbf{a}_{\left(
w,j\right)  } \label{pf.prop.fps.prods-mulable-rules.fubini.2}%
\end{equation}
(and, in particular, the product on the right hand side of this equality is
well-defined, i.e., the family $\left(  \mathbf{a}_{\left(  w,j\right)
}\right)  _{j\in J}$ is multipliable). However, we can replace the product
sign \textquotedblleft$\prod_{s\in J^{\prime}}$\textquotedblright\ by
\textquotedblleft$\prod_{\substack{s\in I\times J;\\f\left(  s\right)  =w}%
}$\textquotedblright\ (since $J^{\prime}=\left\{  s\in I\times J\ \mid
\ f\left(  s\right)  =w\right\}  $). Hence, we can rewrite
(\ref{pf.prop.fps.prods-mulable-rules.fubini.2}) as%
\begin{equation}
\prod_{\substack{s\in I\times J;\\f\left(  s\right)  =w}}\mathbf{a}_{s}%
=\prod_{j\in J}\mathbf{a}_{\left(  w,j\right)  }.
\label{pf.prop.fps.prods-mulable-rules.fubini.3}%
\end{equation}

Forget that we fixed $w$. We thus have proved
(\ref{pf.prop.fps.prods-mulable-rules.fubini.3}) for each $w\in I$.

Now, (\ref{pf.prop.fps.prods-mulable-rules.fubini.1}) becomes%
\[
\prod_{s\in I\times J}\mathbf{a}_{s}=\prod_{w\in I}\ \ \underbrace{\prod
_{\substack{s\in I\times J;\\f\left(  s\right)  =w}}\mathbf{a}_{s}%
}_{\substack{=\prod_{j\in J}\mathbf{a}_{\left(  w,j\right)  }\\\text{(by
(\ref{pf.prop.fps.prods-mulable-rules.fubini.3}))}}}=\prod_{w\in I}%
\ \ \prod_{j\in J}\mathbf{a}_{\left(  w,j\right)  }=\prod_{i\in I}%
\ \ \prod_{j\in J}\mathbf{a}_{\left(  i,j\right)  }%
\]
(here, we have renamed the index $w$ as $i$ in the outer product). Renaming
the index $s$ as $\left(  i,j\right)  $ on the left hand side of this
equality, we can rewrite it as
\[
\prod_{\left(  i,j\right)  \in I\times J}\mathbf{a}_{\left(  i,j\right)
}=\prod_{i\in I}\ \ \prod_{j\in J}\mathbf{a}_{\left(  i,j\right)  }.
\]
A similar argument (but using the map $I\times J\rightarrow J,\ \left(
i,j\right)  \mapsto j$ instead of our map $f:I\times J\rightarrow I,\ \left(
i,j\right)  \mapsto i$) shows that%
\[
\prod_{\left(  i,j\right)  \in I\times J}\mathbf{a}_{\left(  i,j\right)
}=\prod_{j\in J}\ \ \prod_{i\in I}\mathbf{a}_{\left(  i,j\right)  }.
\]
Combining these two equalities, we obtain%
\[
\prod_{i\in I}\ \ \prod_{j\in J}\mathbf{a}_{\left(  i,j\right)  }%
=\prod_{\left(  i,j\right)  \in I\times J}\mathbf{a}_{\left(  i,j\right)
}=\prod_{j\in J}\ \ \prod_{i\in I}\mathbf{a}_{\left(  i,j\right)  }.
\]
(Tracing back our above argument, we see that all products appearing in this
equality are well-defined; indeed, their well-definedness has been shown the
moment they first appeared in our proof.) Proposition
\ref{prop.fps.prods-mulable-rules.fubini1} is thus proved.
\end{proof}

\begin{proof}
[Detailed proof of Proposition \ref{prop.fps.prods-mulable-rules.fubini}.]We
have assumed that $\left(  \mathbf{a}_{\left(  i,j\right)  }\right)  _{\left(
i,j\right)  \in I\times J}\in K\left[  \left[  x\right]  \right]  ^{I\times
J}$ is a multipliable family of invertible FPSs. In other words, $\left(
\mathbf{a}_{s}\right)  _{s\in I\times J}\in K\left[  \left[  x\right]
\right]  ^{I\times J}$ is a multipliable family of invertible FPSs (here, we
have renamed the index $\left(  i,j\right)  $ as $s$). Hence, any subfamily of
$\left(  \mathbf{a}_{s}\right)  _{s\in I\times J}$ is multipliable (by
Proposition \ref{prop.fps.prods-mulable-subfams}, applied to $I\times J$ and
$\left(  \mathbf{a}_{s}\right)  _{s\in I\times J}$ instead of $I$ and $\left(
\mathbf{a}_{i}\right)  _{i\in I}$).

Let $f:I\times J\rightarrow I$ be the map that sends each pair $\left(
i,j\right)  $ to $i$. Let us show three easy claims:

\begin{statement}
\textit{Claim 1:} Let $w\in I$. Let $J^{\prime}$ be the subset $\left\{  s\in
I\times J\ \mid\ f\left(  s\right)  =w\right\}  $ of $I\times J$. Then, there
is a bijection%
\begin{align*}
J  &  \rightarrow J^{\prime},\\
j  &  \mapsto\left(  w,j\right)  .
\end{align*}

\end{statement}

\begin{proof}
[Proof of Claim 1.]This is proved in the exact same way as Claim 1 in our
above proof of Proposition \ref{prop.fps.prods-mulable-rules.fubini1}.
\end{proof}

\begin{statement}
\textit{Claim 2:} For each $i\in I$, the family $\left(  \mathbf{a}_{\left(
i,j\right)  }\right)  _{j\in J}$ is multipliable.
\end{statement}

\begin{proof}
[Proof of Claim 2.]Let $i\in I$. Let $J^{\prime}$ be the subset $\left\{  s\in
I\times J\ \mid\ f\left(  s\right)  =i\right\}  $ of $I\times J$. Then, the
family $\left(  \mathbf{a}_{s}\right)  _{s\in J^{\prime}}$ is a subfamily of
$\left(  \mathbf{a}_{s}\right)  _{s\in I\times J}$. Hence, this family
$\left(  \mathbf{a}_{s}\right)  _{s\in J^{\prime}}$ is multipliable (since any
subfamily of $\left(  \mathbf{a}_{s}\right)  _{s\in I\times J}$ is multipliable).

However, Claim 1 (applied to $w=i$) shows that there is a bijection%
\begin{align*}
J  &  \rightarrow J^{\prime},\\
j  &  \mapsto\left(  i,j\right)  .
\end{align*}
Hence, the family $\left(  \mathbf{a}_{s}\right)  _{s\in J^{\prime}}$ is a
reindexing of the family $\left(  \mathbf{a}_{\left(  i,j\right)  }\right)
_{j\in J}$. Since the former family $\left(  \mathbf{a}_{s}\right)  _{s\in
J^{\prime}}$ is multipliable, we thus conclude that the latter family $\left(
\mathbf{a}_{\left(  i,j\right)  }\right)  _{j\in J}$ is also multipliable
(since a reindexing of a multipliable family is still
multipliable\footnote{This is part of Proposition
\ref{prop.fps.prods-mulable-rules.reindex}.}). This proves Claim 2.
\end{proof}

\begin{statement}
\textit{Claim 3:} For each $j\in J$, the family $\left(  \mathbf{a}_{\left(
i,j\right)  }\right)  _{i\in I}$ is multipliable.
\end{statement}

\begin{proof}
[Proof of Claim 3.]This is analogous to Claim 2.
\end{proof}

Thanks to Claim 2 and Claim 3, we can apply Proposition
\ref{prop.fps.prods-mulable-rules.fubini1}, and conclude that%
\[
\prod_{i\in I}\ \ \prod_{j\in J}\mathbf{a}_{\left(  i,j\right)  }%
=\prod_{\left(  i,j\right)  \in I\times J}\mathbf{a}_{\left(  i,j\right)
}=\prod_{j\in J}\ \ \prod_{i\in I}\mathbf{a}_{\left(  i,j\right)  }%
\]
(and, in particular, all the products appearing in this equality are
well-defined). This proves Proposition
\ref{prop.fps.prods-mulable-rules.fubini}.
\end{proof}
\end{fineprint}

\begin{fineprint}
\begin{proof}
[Detailed proof of Lemma \ref{lem.fps.prod.irlv.inf}.]Theorem
\ref{thm.fps.1+f-mulable} (applied to $I=J$) shows that the family $\left(
1+f_{i}\right)  _{i\in J}$ is multipliable. In other words, each coefficient
in the product of this family $\left(  1+f_{i}\right)  _{i\in J}$ is finitely
determined. In other words, for each $m\in\mathbb{N}$, the $x^{m}$-coefficient
in the product of $\left(  1+f_{i}\right)  _{i\in J}$ is finitely determined.
In other words, for each $m\in\mathbb{N}$, there is a finite subset $M_{m}$ of
$J$ that determines the $x^{m}$-coefficient in the product of $\left(
1+f_{i}\right)  _{i\in J}$. Consider this $M_{m}$.

Let $M=M_{0}\cup M_{1}\cup\cdots\cup M_{n}$. Then, $M$ is a finite subset of
$J$ (since $M_{0},M_{1},\ldots,M_{n}$ are finite subsets of $J$). Moreover, we
claim that
\begin{equation}
\left[  x^{m}\right]  \left(  \prod_{i\in J}\left(  1+f_{i}\right)  \right)
=\left[  x^{m}\right]  \left(  \prod_{i\in M}\left(  1+f_{i}\right)  \right)
\label{pf.lem.fps.prod.irlv.inf.1}%
\end{equation}
for each $m\in\left\{  0,1,\ldots,n\right\}  $.

[\textit{Proof of (\ref{pf.lem.fps.prod.irlv.inf.1}):} Let $m\in\left\{
0,1,\ldots,n\right\}  $. Then, $M_{m}$ is one of the $n+1$ sets in the union
$M_{0}\cup M_{1}\cup\cdots\cup M_{n}$. Hence, $M_{m}\subseteq M_{0}\cup
M_{1}\cup\cdots\cup M_{n}=M$.

However, the subset $M_{m}$ of $J$ determines the $x^{m}$-coefficient in the
product of $\left(  1+f_{i}\right)  _{i\in J}$ (by the definition of $M_{m}$).
In other words, every finite subset $J^{\prime}$ of $J$ satisfying
$M_{m}\subseteq J^{\prime}\subseteq J$ satisfies%
\[
\left[  x^{m}\right]  \left(  \prod_{i\in J^{\prime}}\left(  1+f_{i}\right)
\right)  =\left[  x^{m}\right]  \left(  \prod_{i\in M_{m}}\left(
1+f_{i}\right)  \right)
\]
(by the definition of \textquotedblleft determining the $x^{m}$-coefficient in
a product\textquotedblright). Applying this to $J^{\prime}=M$, we obtain%
\[
\left[  x^{m}\right]  \left(  \prod_{i\in M}\left(  1+f_{i}\right)  \right)
=\left[  x^{m}\right]  \left(  \prod_{i\in M_{m}}\left(  1+f_{i}\right)
\right)
\]
(since $M$ is a finite subset of $J$ satisfying $M_{m}\subseteq M\subseteq
J$). On the other hand, the definition of the product $\prod_{i\in J}\left(
1+f_{i}\right)  $ yields that%
\[
\left[  x^{m}\right]  \left(  \prod_{i\in J}\left(  1+f_{i}\right)  \right)
=\left[  x^{m}\right]  \left(  \prod_{i\in M_{m}}\left(  1+f_{i}\right)
\right)
\]
(since $M_{m}$ is a finite subset of $J$ that determines the $x^{m}%
$-coefficient in the product of $\left(  1+f_{i}\right)  _{i\in J}$).
Comparing these two equalities, we obtain%
\[
\left[  x^{m}\right]  \left(  \prod_{i\in J}\left(  1+f_{i}\right)  \right)
=\left[  x^{m}\right]  \left(  \prod_{i\in M}\left(  1+f_{i}\right)  \right)
.
\]
This proves (\ref{pf.lem.fps.prod.irlv.inf.1}).]

Now, we know that (\ref{pf.lem.fps.prod.irlv.inf.1}) holds for each
$m\in\left\{  0,1,\ldots,n\right\}  $. Thus, we can apply Lemma
\ref{lem.fps.prod.irlv.fg} to $f=\prod_{i\in J}\left(  1+f_{i}\right)  $ and
$g=\prod_{i\in M}\left(  1+f_{i}\right)  $. We thus obtain that%
\begin{equation}
\left[  x^{m}\right]  \left(  a\prod_{i\in J}\left(  1+f_{i}\right)  \right)
=\left[  x^{m}\right]  \left(  a\prod_{i\in M}\left(  1+f_{i}\right)  \right)
\label{pf.lem.fps.prod.irlv.inf.3}%
\end{equation}
for each $m\in\left\{  0,1,\ldots,n\right\}  $.

On the other hand, $M$ is a finite set, so that $\left(  f_{i}\right)  _{i\in
M}$ is a finite family. Furthermore, each $i\in M$ satisfies $\left[
x^{m}\right]  \left(  f_{i}\right)  =0$ for each $m\in\left\{  0,1,\ldots
,n\right\}  $ (by (\ref{eq.lem.fps.prod.irlv.inf.ass}), since $i\in M\subseteq
J$). Thus, we can apply Lemma \ref{lem.fps.prod.irlv.fin} to $M$ instead of
$J$. We therefore obtain that%
\begin{equation}
\left[  x^{m}\right]  \left(  a\prod_{i\in M}\left(  1+f_{i}\right)  \right)
=\left[  x^{m}\right]  a \label{pf.lem.fps.prod.irlv.inf.4}%
\end{equation}
for each $m\in\left\{  0,1,\ldots,n\right\}  $.

Hence, for each $m\in\left\{  0,1,\ldots,n\right\}  $, we have%
\begin{align*}
\left[  x^{m}\right]  \left(  a\prod_{i\in J}\left(  1+f_{i}\right)  \right)
&  =\left[  x^{m}\right]  \left(  a\prod_{i\in M}\left(  1+f_{i}\right)
\right)  \ \ \ \ \ \ \ \ \ \ \left(  \text{by
(\ref{pf.lem.fps.prod.irlv.inf.3})}\right) \\
&  =\left[  x^{m}\right]  a\ \ \ \ \ \ \ \ \ \ \left(  \text{by
(\ref{pf.lem.fps.prod.irlv.inf.4})}\right)  .
\end{align*}
This proves Lemma \ref{lem.fps.prod.irlv.inf}.
\end{proof}
\end{fineprint}

\begin{fineprint}
\begin{proof}
[Detailed proof of Proposition \ref{prop.fps.prodrule-inf-inf}.]The following
proof is an expanded version of the argument given by Mindlack at
\url{https://math.stackexchange.com/a/4123658/} .

This will be a long grind; we thus break it up into several claims. First,
however, let us introduce a few notations:

\begin{itemize}
\item If $J$ is a subset of $I$, then $S^{J}$ shall denote the Cartesian
product $\prod_{i\in J}S_{i}$. This Cartesian product $\prod_{i\in J}S_{i}$
consists of families $\left(  s_{i}\right)  _{i\in J}$, where each $s_{i}$
belongs to the respective set $S_{i}$.

The notation $S^{J}$ should not be misconstrued as being an actual power.
(However, in the particular case when all $S_{i}$ equal one and the same set
$S$, the Cartesian product $S^{J}$ we just defined is indeed the Cartesian
power commonly known as $S^{J}$.)

\item If $J$ is a subset of $I$, then $S_{J}^{I}$ shall denote the set of all
families $\left(  s_{i}\right)  _{i\in I}\in S^{I}$ that satisfy%
\[
\left(  s_{i}=0\text{ for all }i\in I\setminus J\right)  .
\]
This set $S_{J}^{I}$ is in a canonical bijection with $S^{J}$, as elements of
both sets consist of \textquotedblleft essentially the same
data\textquotedblright. To wit, an element of $S^{J}$ is a family that only
has $i$-th entries for $i\in J$, whereas an element of $S_{J}^{I}$ is a family
that has $i$-th entries for all $i\in I$, but subject to the requirement that
the $i$-th entries for all $i\in I\setminus J$ are $0$ (so that only the
$i$-th entries for $i\in J$ carry any information). More rigorously: The map%
\begin{align*}
S_{J}^{I}  &  \rightarrow S^{J},\\
\left(  s_{i}\right)  _{i\in I}  &  \mapsto\left(  s_{i}\right)  _{i\in J}%
\end{align*}
is a bijection (since it merely shrinks the family by removing entries that
are required to be $0$ anyway). We denote this bijection by
$\operatorname*{reduce}\nolimits_{J}$.

\item We define $S_{\operatorname*{fin}}^{I}$ to be the set of all essentially
finite families $\left(  s_{i}\right)  _{i\in I}\in S^{I}$. It is easy to see
that $S_{\operatorname*{fin}}^{I}$ is the union of the sets $S_{J}^{I}$ over
all finite subsets $J$ of $I$.
\end{itemize}

Now, we can begin with our claims:

\begin{statement}
\textit{Claim 1:} The family $\left(  p_{i,k}\right)  _{k\in S_{i}}$ is
summable for each $i\in I$.
\end{statement}

[\textit{Proof of Claim 1:} Let $j\in I$. Then, the pairs $\left(  i,k\right)
\in\overline{S}$ with $i=j$ are precisely the pairs of the form $\left(
j,k\right)  $ with $k\in S_{j}$ and $k\neq0$. In other words, the pairs
$\left(  i,k\right)  \in\overline{S}$ with $i=j$ are precisely the pairs of
the form $\left(  j,k\right)  $ with $k\in S_{j}\setminus\left\{  0\right\}  $.

We assumed that the family $\left(  p_{i,k}\right)  _{\left(  i,k\right)
\in\overline{S}}$ is summable. Hence, its subfamily $\left(  p_{i,k}\right)
_{\left(  i,k\right)  \in\overline{S}\text{ with }i=j}$ is summable as well
(since a subfamily of a summable family is always summable). In other words,
the family $\left(  p_{j,k}\right)  _{k\in S_{j}\setminus\left\{  0\right\}
}$ is summable (since this family is just a reindexing of the family $\left(
p_{i,k}\right)  _{\left(  i,k\right)  \in\overline{S}\text{ with }i=j}$
(because the pairs $\left(  i,k\right)  \in\overline{S}$ with $i=j$ are
precisely the pairs of the form $\left(  j,k\right)  $ with $k\in
S_{j}\setminus\left\{  0\right\}  $)). Thus, the family $\left(
p_{j,k}\right)  _{k\in S_{j}}$ is summable as well (since the summability of a
family does not change if we insert a single entry into it\footnote{Indeed,
the summability of a family is an \textquotedblleft all but finitely many $k$
satisfy something\textquotedblright\ type of statement. If we insert a single
entry into the family, such a statement does not change its validity.}).

Forget that we fixed $j$. We thus have shown that the family $\left(
p_{j,k}\right)  _{k\in S_{j}}$ is summable for each $j\in I$. Renaming $j$ as
$i$ in this statement, we obtain the following: The family $\left(
p_{i,k}\right)  _{k\in S_{i}}$ is summable for each $i\in I$. This proves
Claim 1.]\medskip

Claim 1 shows that the sum $\sum_{k\in S_{i}}p_{i,k}$ is well-defined for each
$i\in I$. Moreover, for each $i\in I$, we have%
\begin{align}
\sum_{k\in S_{i}}p_{i,k}  &  =\underbrace{p_{i,0}}_{\substack{=1\\\text{(by
(\ref{eq.prop.fps.prodrule-inf-inf.pi0=1}))}}}+\sum_{k\in S_{i}\setminus
\left\{  0\right\}  }p_{i,k}\ \ \ \ \ \ \ \ \ \ \left(
\begin{array}
[c]{c}%
\text{here, we have split off}\\
\text{the addend for }k=0\\
\text{from the sum, since }0\in S_{i}%
\end{array}
\right) \nonumber\\
&  =1+\sum_{k\in S_{i}\setminus\left\{  0\right\}  }p_{i,k}.
\label{pf.prop.fps.prodrule-fin-inf.sum-no-1}%
\end{align}

Next, we claim the following:

\begin{statement}
\textit{Claim 2:} The family $\left(  \sum_{k\in S_{i}}p_{i,k}\right)  _{i\in
I}$ is multipliable.
\end{statement}

[\textit{Proof of Claim 2:} The family $\left(  p_{i,k}\right)  _{\left(
i,k\right)  \in\overline{S}}$ is summable (by assumption). We can split its
sum into subsums as follows:%
\[
\underbrace{\sum_{\left(  i,k\right)  \in\overline{S}}}_{\substack{=\sum_{i\in
I}\ \ \sum_{k\in S_{i}\setminus\left\{  0\right\}  }\\\text{(since a pair
}\left(  i,k\right)  \text{ belongs to }\overline{S}\\\text{if and only if it
satisfies }i\in I\\\text{and }k\in S_{i}\setminus\left\{  0\right\}  \text{)}%
}}p_{i,k}=\sum_{i\in I}\ \ \sum_{k\in S_{i}\setminus\left\{  0\right\}
}p_{i,k}.
\]
This shows that the family $\left(  \sum_{k\in S_{i}\setminus\left\{
0\right\}  }p_{i,k}\right)  _{i\in I}$ is summable. Hence, Theorem
\ref{thm.fps.1+f-mulable} (applied to $f_{i}=\sum_{k\in S_{i}\setminus\left\{
0\right\}  }p_{i,k}$) yields that the family $\left(  1+\sum_{k\in
S_{i}\setminus\left\{  0\right\}  }p_{i,k}\right)  _{i\in I}$ is multipliable.
In other words, the family $\left(  \sum_{k\in S_{i}}p_{i,k}\right)  _{i\in
I}$ is multipliable (since (\ref{pf.prop.fps.prodrule-fin-inf.sum-no-1}) shows
that this family is precisely the family $\left(  1+\sum_{k\in S_{i}%
\setminus\left\{  0\right\}  }p_{i,k}\right)  _{i\in I}$). This proves Claim
2.]\medskip

Claim 2 shows that the product $\prod_{i\in I}\ \ \sum_{k\in S_{i}}p_{i,k}$ is
well-defined.\medskip

The family $\left(  p_{i,k}\right)  _{\left(  i,k\right)  \in\overline{S}}$ is
summable (by assumption). In other words, for any $m\in\mathbb{N}$, all but
finitely many $\left(  i,k\right)  \in\overline{S}$ satisfy $\left[
x^{m}\right]  p_{i,k}=0$. In other words, for any $m\in\mathbb{N}$, there
exists a finite subset $T_{m}$ of $\overline{S}$ such that
\begin{equation}
\text{all }\left(  i,k\right)  \in\overline{S}\setminus T_{m}\text{ satisfy
}\left[  x^{m}\right]  p_{i,k}=0. \label{pf.prop.fps.prodrule-fin-inf.Tm-def}%
\end{equation}
Consider these finite subsets $T_{m}$.

For each $n\in\mathbb{N}$, we let $T_{n}^{\prime}$ be the subset $T_{0}\cup
T_{1}\cup\cdots\cup T_{n}$ of $\overline{S}$. This subset $T_{n}^{\prime}$ is
finite (since $T_{0},T_{1},\ldots,T_{n}$ are finite).

For each $n\in\mathbb{N}$, we let $I_{n}$ be the subset
\[
\left\{  i\ \mid\ \left(  i,k\right)  \in T_{n}^{\prime}\right\}
\]
of $I$, and we let $K_{n}$ be the set%
\[
\left\{  k\ \mid\ \left(  i,k\right)  \in T_{n}^{\prime}\right\}  .
\]
These two sets $I_{n}$ and $K_{n}$ are finite (since $T_{n}^{\prime}$ is finite).

The definition of $T_{n}^{\prime}$ shows the following:

\begin{statement}
\textit{Claim 3:} Let $n\in\mathbb{N}$ and $\left(  i,k\right)  \in
\overline{S}\setminus T_{n}^{\prime}$. Then, we have%
\[
\left[  x^{m}\right]  p_{i,k}=0\ \ \ \ \ \ \ \ \ \ \text{for each }%
m\in\left\{  0,1,\ldots,n\right\}  .
\]

\end{statement}

[\textit{Proof of Claim 3:} Let $m\in\left\{  0,1,\ldots,n\right\}  $. Then,
$T_{m}\subseteq T_{n}^{\prime}$ (since $T_{n}^{\prime}$ is defined as the
union $T_{0}\cup T_{1}\cup\cdots\cup T_{n}$, whereas $T_{m}$ is one of the
$n+1$ sets appearing in this union). Hence, $T_{n}^{\prime}\supseteq T_{m}$,
so that $\overline{S}\setminus\underbrace{T_{n}^{\prime}}_{\supseteq T_{m}%
}\subseteq\overline{S}\setminus T_{m}$. Now, $\left(  i,k\right)  \in
\overline{S}\setminus T_{n}^{\prime}\subseteq\overline{S}\setminus T_{m}$.
Therefore, (\ref{pf.prop.fps.prodrule-fin-inf.Tm-def}) shows that $\left[
x^{m}\right]  p_{i,k}=0$. This proves Claim 3.]\medskip

The following is easy to see:

\begin{statement}
\textit{Claim 4:} If $\left(  k_{i}\right)  _{i\in I}\in
S_{\operatorname*{fin}}^{I}$, then the family $\left(  p_{i,k_{i}}\right)
_{i\in I}$ is multipliable.
\end{statement}

[\textit{Proof of Claim 4:} Let $\left(  k_{i}\right)  _{i\in I}\in
S_{\operatorname*{fin}}^{I}$. Thus, $\left(  k_{i}\right)  _{i\in I}$ is an
essentially finite family in $S^{I}$. Now, all but finitely many $i\in I$
satisfy $k_{i}=0$ (since $\left(  k_{i}\right)  _{i\in I}$ is essentially
finite) and thus $p_{i,k_{i}}=p_{i,0}=1$ (by
(\ref{eq.prop.fps.prodrule-inf-inf.pi0=1})). Hence, all but finitely many
entries of the family $\left(  p_{i,k_{i}}\right)  _{i\in I}$ equal $1$. Thus,
this family is multipliable (by Proposition \ref{prop.fps.1-mulable}). This
proves Claim 4.]

Claim 4 shows that the product $\prod_{i\in I}p_{i,k_{i}}$ is well-defined
whenever $\left(  k_{i}\right)  _{i\in I}\in S_{\operatorname*{fin}}^{I}$.
Next, we claim the following:

\begin{statement}
\textit{Claim 5:} Let $n\in\mathbb{N}$. Let $\left(  k_{i}\right)  _{i\in
I}\in S_{\operatorname*{fin}}^{I}$. Assume that some $j\in I$ satisfies
$\left(  j,k_{j}\right)  \in\overline{S}\setminus T_{n}^{\prime}$. Then,
$\left[  x^{n}\right]  \left(  \prod_{i\in I}p_{i,k_{i}}\right)  =0$.
\end{statement}

[\textit{Proof of Claim 5:} We have assumed that some $j\in I$ satisfies
$\left(  j,k_{j}\right)  \in\overline{S}\setminus T_{n}^{\prime}$. Consider
this $j$. Then, the product $\prod_{i\in I}p_{i,k_{i}}$ is a multiple of
$p_{j,k_{j}}$ (since $p_{j,k_{j}}$ is one of the factors of this product).
Moreover, we have $\left[  x^{m}\right]  p_{j,k_{j}}=0$ for each $m\in\left\{
0,1,\ldots,n\right\}  $ (by Claim 3, applied to $\left(  i,k\right)  =\left(
j,k_{j}\right)  $). Hence, Lemma \ref{lem.fps.prod.irlv.mul} (applied to
$u=p_{j,k_{j}}$ and $v=\prod_{i\in I}p_{i,k_{i}}$) shows that we have $\left[
x^{m}\right]  \left(  \prod_{i\in I}p_{i,k_{i}}\right)  =0$ for each
$m\in\left\{  0,1,\ldots,n\right\}  $ (since $\prod_{i\in I}p_{i,k_{i}}$ is a
multiple of $p_{j,k_{j}}$). Applying this to $m=n$, we obtain $\left[
x^{n}\right]  \left(  \prod_{i\in I}p_{i,k_{i}}\right)  =0$. This proves Claim 5.]

\begin{statement}
\textit{Claim 6:} Let $n\in\mathbb{N}$. Let $\left(  k_{i}\right)  _{i\in
I}\in S_{\operatorname*{fin}}^{I}\setminus S_{I_{n}}^{I}$. Then, $\left[
x^{n}\right]  \left(  \prod_{i\in I}p_{i,k_{i}}\right)  =0$.
\end{statement}

[\textit{Proof of Claim 6:} We have $\left(  k_{i}\right)  _{i\in I}\in
S_{\operatorname*{fin}}^{I}\setminus S_{I_{n}}^{I}$. In other words, we have
$\left(  k_{i}\right)  _{i\in I}\in S_{\operatorname*{fin}}^{I}$ but $\left(
k_{i}\right)  _{i\in I}\notin S_{I_{n}}^{I}$. From $\left(  k_{i}\right)
_{i\in I}\notin S_{I_{n}}^{I}$, we see that there exists some $j\in I\setminus
I_{n}$ satisfying $k_{j}\neq0$. Consider this $j$. We have $k_{j}\in S_{j}$
(since $\left(  k_{i}\right)  _{i\in I}\in S_{\operatorname*{fin}}%
^{I}\subseteq S^{I}$) and $j\notin I_{n}$ (since $j\in I\setminus I_{n}$). We
have $\left(  j,k_{j}\right)  \in\overline{S}$ (since $k_{j}\in S_{j}$ and
$k_{j}\neq0$) and $\left(  j,k_{j}\right)  \notin T_{n}^{\prime}$ (because if
we had $\left(  j,k_{j}\right)  \in T_{n}^{\prime}$, then we would have $j\in
I_{n}$ by the definition of $I_{n}$; but this would contradict $j\notin I_{n}%
$). Hence, we have $\left(  j,k_{j}\right)  \in\overline{S}\setminus
T_{n}^{\prime}$. Therefore, Claim 5 yields $\left[  x^{n}\right]  \left(
\prod_{i\in I}p_{i,k_{i}}\right)  =0$. This proves Claim 6.]

\begin{statement}
\textit{Claim 7:} The family $\left(  \prod_{i\in I}p_{i,k_{i}}\right)
_{\left(  k_{i}\right)  _{i\in I}\in S_{\operatorname*{fin}}^{I}}$ is summable.
\end{statement}

[\textit{Proof of Claim 7:} Let $n\in\mathbb{N}$. We shall show that all but
finitely many families $\left(  k_{i}\right)  _{i\in I}\in
S_{\operatorname*{fin}}^{I}$ satisfy
\[
\left[  x^{n}\right]  \left(  \prod_{i\in I}p_{i,k_{i}}\right)  =0.
\]

Indeed, let $\left(  k_{i}\right)  _{i\in I}\in S_{\operatorname*{fin}}^{I}$
be a family such that
\begin{equation}
\left[  x^{n}\right]  \left(  \prod_{i\in I}p_{i,k_{i}}\right)  \neq0.
\label{pf.prop.fps.prodrule-fin-inf.c7.pf.ass}%
\end{equation}
We are going to show that $\left(  k_{i}\right)  _{i\in I}$ satisfies the
following two properties:

\begin{itemize}
\item \textit{Property 1:} All $i\in I\setminus I_{n}$ satisfy $k_{i}=0$.

\item \textit{Property 2:} All $i\in I_{n}$ satisfy $k_{i}\in K_{n}%
\cup\left\{  0\right\}  $.
\end{itemize}

\noindent These two properties together will restrict the family $\left(
k_{i}\right)  _{i\in I}$ to finitely many possibilities (since the sets
$I_{n}$ and $K_{n}$ are finite).

If we had $\left(  k_{i}\right)  _{i\in I}\notin S_{I_{n}}^{I}$, then we would
have $\left(  k_{i}\right)  _{i\in I}\in S_{\operatorname*{fin}}^{I}\setminus
S_{I_{n}}^{I}$ (since $\left(  k_{i}\right)  _{i\in I}\in
S_{\operatorname*{fin}}^{I}$) and thus $\left[  x^{n}\right]  \left(
\prod_{i\in I}p_{i,k_{i}}\right)  =0$ (by Claim 6), which would contradict
(\ref{pf.prop.fps.prodrule-fin-inf.c7.pf.ass}). Hence, we cannot have $\left(
k_{i}\right)  _{i\in I}\notin S_{I_{n}}^{I}$. Thus, we have $\left(
k_{i}\right)  _{i\in I}\in S_{I_{n}}^{I}$. In other words, all $i\in
I\setminus I_{n}$ satisfy $k_{i}=0$. This proves Property 1.

If there was some $j\in I_{n}$ that satisfies $k_{j}\notin K_{n}\cup\left\{
0\right\}  $, then this $k_{j}$ would satisfy $\left(  j,k_{j}\right)
\in\overline{S}\setminus T_{n}^{\prime}$\ \ \ \ \footnote{\textit{Proof.} Let
$j\in I_{n}$ be such that $k_{j}\notin K_{n}\cup\left\{  0\right\}  $. We must
show that $\left(  j,k_{j}\right)  \in\overline{S}\setminus T_{n}^{\prime}$.
\par
Indeed, we have $k_{j}\notin K_{n}\cup\left\{  0\right\}  $; thus,
$k_{j}\notin K_{n}$ and $k_{j}\neq0$. However, $j\in I_{n}\subseteq I$ and
$k_{j}\in S_{j}$ (since $\left(  k_{i}\right)  _{i\in I}\in
S_{\operatorname*{fin}}^{I}\subseteq S^{I}$) and therefore $\left(
j,k_{j}\right)  \in\overline{S}$ (by the definition of $\overline{S}$, since
$k_{j}\neq0$). If we had $\left(  j,k_{j}\right)  \in T_{n}^{\prime}$, then we
would have $k_{j}\in K_{n}$ (by the definition of $K_{n}$), which would
contradict $k_{j}\notin K_{n}$. Thus, we have $\left(  j,k_{j}\right)  \notin
T_{n}^{\prime}$. Combining this with $\left(  j,k_{j}\right)  \in\overline{S}%
$, we obtain $\left(  j,k_{j}\right)  \in\overline{S}\setminus T_{n}^{\prime}%
$. This completes our proof.}, and therefore we would have $\left[
x^{n}\right]  \left(  \prod_{i\in I}p_{i,k_{i}}\right)  =0$ (by Claim 5, since
$j\in I_{n}\subseteq I$); but this would contradict
(\ref{pf.prop.fps.prodrule-fin-inf.c7.pf.ass}). Hence, there exists no $j\in
I_{n}$ that satisfies $k_{j}\notin K_{n}\cup\left\{  0\right\}  $. In other
words, all $i\in I_{n}$ satisfy $k_{i}\in K_{n}\cup\left\{  0\right\}  $. This
proves Property 2.

Now, we have shown that our family $\left(  k_{i}\right)  _{i\in I}$ satisfies
Property 1 and Property 2.

Forget that we fixed $\left(  k_{i}\right)  _{i\in I}$. We thus have shown
that any family $\left(  k_{i}\right)  _{i\in I}\in S_{\operatorname*{fin}%
}^{I}$ that satisfies $\left[  x^{n}\right]  \left(  \prod_{i\in I}p_{i,k_{i}%
}\right)  \neq0$ must satisfy Property 1 and Property 2. In other words, any
such family must belong to the set of all families $\left(  k_{i}\right)
_{i\in I}\in S_{\operatorname*{fin}}^{I}$ that satisfy Property 1 and Property
2. However, the latter set is finite\footnote{\textit{Proof.} We must show
that Property 1 and Property 2 leave only finitely many options for the family
$\left(  k_{i}\right)  _{i\in I}$. Indeed, Property 1 shows that all entries
$k_{i}$ with $i\in I\setminus I_{n}$ are uniquely determined; meanwhile,
Property 2 ensures that the remaining entries (of which there are only
finitely many, since the set $I_{n}$ is finite) must belong to the finite set
$K_{n}\cup\left\{  0\right\}  $ (this set is finite, since $K_{n}$ is finite).
Therefore, a family $\left(  k_{i}\right)  _{i\in I}$ that satisfies Property
1 and Property 2 is uniquely determined by finitely many of its entries
(namely, by its entries $k_{i}$ with $i\in I_{n}$), and there are finitely
many choices for each of them (since they must belong to the finite set
$K_{n}\cup\left\{  0\right\}  $). Hence, there are only finitely many such
families (namely, at most $\left\vert K_{n}\cup\left\{  0\right\}  \right\vert
^{\left\vert I_{n}\right\vert }$ many options). In other words, the set of all
families $\left(  k_{i}\right)  _{i\in I}\in S_{\operatorname*{fin}}^{I}$ that
satisfy Property 1 and Property 2 is finite.}. Hence, there are only finitely
many families $\left(  k_{i}\right)  _{i\in I}\in S_{\operatorname*{fin}}^{I}$
that satisfy $\left[  x^{n}\right]  \left(  \prod_{i\in I}p_{i,k_{i}}\right)
\neq0$ (since we have shown that any such family must belong to the finite set
of all families $\left(  k_{i}\right)  _{i\in I}\in S_{\operatorname*{fin}%
}^{I}$ that satisfy Property 1 and Property 2). In other words, all but
finitely many families $\left(  k_{i}\right)  _{i\in I}\in
S_{\operatorname*{fin}}^{I}$ satisfy $\left[  x^{n}\right]  \left(
\prod_{i\in I}p_{i,k_{i}}\right)  =0$.

Forget that we fixed $n$. We thus have shown that for each $n\in\mathbb{N}$,
all but finitely many families $\left(  k_{i}\right)  _{i\in I}\in
S_{\operatorname*{fin}}^{I}$ satisfy $\left[  x^{n}\right]  \left(
\prod_{i\in I}p_{i,k_{i}}\right)  =0$. In other words, the family $\left(
\prod_{i\in I}p_{i,k_{i}}\right)  _{\left(  k_{i}\right)  _{i\in I}\in
S_{\operatorname*{fin}}^{I}}$ is summable. This proves Claim 7.]

Claim 7 shows that the sum $\sum_{\left(  k_{i}\right)  _{i\in I}\in
S_{\operatorname*{fin}}^{I}}\ \ \prod_{i\in I}p_{i,k_{i}}$ is well-defined. We
shall next focus on proving (\ref{eq.prop.fps.prodrule-inf-inf.eq}).

\begin{statement}
\textit{Claim 8:} Let $n\in\mathbb{N}$. Then,
\[
\left[  x^{n}\right]  \left(  \sum_{\left(  k_{i}\right)  _{i\in I}\in
S_{\operatorname*{fin}}^{I}}\ \ \prod_{i\in I}p_{i,k_{i}}\right)  =\left[
x^{n}\right]  \left(  \sum_{\left(  k_{i}\right)  _{i\in I_{n}}\in S^{I_{n}}%
}\ \ \prod_{i\in I_{n}}p_{i,k_{i}}\right)  .
\]

\end{statement}

[\textit{Proof of Claim 8:} The set $I_{n}$ is finite (as we know). Thus,
$S_{I_{n}}^{I}\subseteq S_{\operatorname*{fin}}^{I}$%
\ \ \ \ \footnote{\textit{Proof.} Let $\left(  k_{i}\right)  _{i\in I}\in
S_{I_{n}}^{I}$. Thus, $\left(  k_{i}\right)  _{i\in I}$ is a family in $S^{I}$
that satisfies $k_{i}=0$ for all $i\in I\setminus I_{n}$ (by the definition of
$S_{I_{n}}^{I}$). However, $I_{n}$ is finite. Thus, $k_{i}=0$ for all but
finitely many $i\in I$ (since $k_{i}=0$ for all $i\in I\setminus I_{n}$). In
other words, the family $\left(  k_{i}\right)  _{i\in I}$ is essentially
finite. In other words, $\left(  k_{i}\right)  _{i\in I}\in
S_{\operatorname*{fin}}^{I}$ (by the definition of $S_{\operatorname*{fin}%
}^{I}$).
\par
Forget that we fixed $\left(  k_{i}\right)  _{i\in I}$. We thus have shown
that $\left(  k_{i}\right)  _{i\in I}\in S_{\operatorname*{fin}}^{I}$ for each
$\left(  k_{i}\right)  _{i\in I}\in S_{I_{n}}^{I}$. In other words, $S_{I_{n}%
}^{I}\subseteq S_{\operatorname*{fin}}^{I}$.}. Hence, the set
$S_{\operatorname*{fin}}^{I}$ is the union of its two disjoint subsets
$S_{I_{n}}^{I}$ and $S_{\operatorname*{fin}}^{I}\setminus S_{I_{n}}^{I}$.

Furthermore, for each $\left(  k_{i}\right)  _{i\in I}\in S_{I_{n}}^{I}$, we
have%
\[
k_{i}=0\ \ \ \ \ \ \ \ \ \ \text{for all }i\in I\setminus I_{n}%
\]
(by the definition of $S_{I_{n}}^{I}$) and therefore%
\[
p_{i,k_{i}}=p_{i,0}=1\ \ \ \ \ \ \ \ \ \ \text{for all }i\in I\setminus I_{n}%
\]
(by (\ref{eq.prop.fps.prodrule-inf-inf.pi0=1})) and thus%
\[
\prod_{i\in I\setminus I_{n}}\underbrace{p_{i,k_{i}}}_{=1}=\prod_{i\in
I\setminus I_{n}}1=1
\]
and therefore%
\begin{align}
\prod_{i\in I}p_{i,k_{i}}  &  =\left(  \prod_{i\in I_{n}}p_{i,k_{i}}\right)
\underbrace{\left(  \prod_{i\in I\setminus I_{n}}p_{i,k_{i}}\right)  }%
_{=1}\nonumber\\
&  \ \ \ \ \ \ \ \ \ \ \ \ \ \ \ \ \ \ \ \ \left(
\begin{array}
[c]{c}%
\text{here, we have split the product into two}\\
\text{parts, since the set }I_{n}\text{ is a subset of }I
\end{array}
\right) \nonumber\\
&  =\prod_{i\in I_{n}}p_{i,k_{i}}.
\label{pf.prop.fps.prodrule-fin-inf.c8.pf.1}%
\end{align}

Now,%
\begin{align*}
&  \left[  x^{n}\right]  \left(  \sum_{\left(  k_{i}\right)  _{i\in I}\in
S_{\operatorname*{fin}}^{I}}\ \ \prod_{i\in I}p_{i,k_{i}}\right) \\
&  =\sum_{\left(  k_{i}\right)  _{i\in I}\in S_{\operatorname*{fin}}^{I}%
}\left[  x^{n}\right]  \left(  \prod_{i\in I}p_{i,k_{i}}\right) \\
&  =\sum_{\left(  k_{i}\right)  _{i\in I}\in S_{I_{n}}^{I}}\left[
x^{n}\right]  \underbrace{\left(  \prod_{i\in I}p_{i,k_{i}}\right)
}_{\substack{=\prod_{i\in I_{n}}p_{i,k_{i}}\\\text{(by
(\ref{pf.prop.fps.prodrule-fin-inf.c8.pf.1}))}}}+\sum_{\left(  k_{i}\right)
_{i\in I}\in S_{\operatorname*{fin}}^{I}\setminus S_{I_{n}}^{I}}%
\underbrace{\left[  x^{n}\right]  \left(  \prod_{i\in I}p_{i,k_{i}}\right)
}_{\substack{=0\\\text{(by Claim 6)}}}\\
&  \ \ \ \ \ \ \ \ \ \ \ \ \ \ \ \ \ \ \ \ \left(
\begin{array}
[c]{c}%
\text{here, we have split the sum, since the set }S_{\operatorname*{fin}}%
^{I}\text{ is}\\
\text{the union of its two disjoint subsets }S_{I_{n}}^{I}\text{ and
}S_{\operatorname*{fin}}^{I}\setminus S_{I_{n}}^{I}%
\end{array}
\right) \\
&  =\sum_{\left(  k_{i}\right)  _{i\in I}\in S_{I_{n}}^{I}}\left[
x^{n}\right]  \left(  \prod_{i\in I_{n}}p_{i,k_{i}}\right)  +\underbrace{\sum
_{\left(  k_{i}\right)  _{i\in I}\in S_{\operatorname*{fin}}^{I}\setminus
S_{I_{n}}^{I}}0}_{=0}\\
&  =\sum_{\left(  k_{i}\right)  _{i\in I}\in S_{I_{n}}^{I}}\left[
x^{n}\right]  \left(  \prod_{i\in I_{n}}p_{i,k_{i}}\right)  =\sum_{\left(
k_{i}\right)  _{i\in I_{n}}\in S^{I_{n}}}\left[  x^{n}\right]  \left(
\prod_{i\in I_{n}}p_{i,k_{i}}\right) \\
&  \ \ \ \ \ \ \ \ \ \ \ \ \ \ \ \ \ \ \ \ \left(
\begin{array}
[c]{c}%
\text{here, we have substituted }\left(  k_{i}\right)  _{i\in I_{n}}\text{ for
}\left(  k_{i}\right)  _{i\in I_{n}}\text{ in the}\\
\text{sum, since the map }S_{I_{n}}^{I}\rightarrow S^{I_{n}},\ \left(
k_{i}\right)  _{i\in I}\mapsto\left(  k_{i}\right)  _{i\in I_{n}}\\
\text{is a bijection (indeed, this map is the map}\\
\text{we have called }\operatorname*{reduce}\nolimits_{I_{n}}\text{)}%
\end{array}
\right) \\
&  =\left[  x^{n}\right]  \left(  \sum_{\left(  k_{i}\right)  _{i\in I_{n}}\in
S^{I_{n}}}\ \ \prod_{i\in I_{n}}p_{i,k_{i}}\right)  .
\end{align*}
This proves Claim 8.]

\begin{statement}
\textit{Claim 9:} Let $n\in\mathbb{N}$. Let $i\in I\setminus I_{n}$. Then,
\[
\left[  x^{m}\right]  \left(  \sum_{k\in S_{i}\setminus\left\{  0\right\}
}p_{i,k}\right)  =0\ \ \ \ \ \ \ \ \ \ \text{for each }m\in\left\{
0,1,\ldots,n\right\}  .
\]

\end{statement}

[\textit{Proof of Claim 9:} We have $i\in I\setminus I_{n}$. In other words,
$i\in I$ and $i\notin I_{n}$.

Let $m\in\left\{  0,1,\ldots,n\right\}  $.

Let $k\in S_{i}\setminus\left\{  0\right\}  $. Thus, by the definition of
$\overline{S}$, we have $\left(  i,k\right)  \in\overline{S}$ (since $k\in
S_{i}\setminus\left\{  0\right\}  $ entails $k\in S_{i}$ and $k\neq0$). On the
other hand, $\left(  i,k\right)  \notin T_{n}^{\prime}$ (since otherwise, we
would have $\left(  i,k\right)  \in T_{n}^{\prime}$ and thus $i\in I_{n}$ (by
the definition of $I_{n}$), which would contradict $i\notin I_{n}$). Combining
$\left(  i,k\right)  \in\overline{S}$ with $\left(  i,k\right)  \notin
T_{n}^{\prime}$, we obtain $\left(  i,k\right)  \in\overline{S}\setminus
T_{n}^{\prime}$. Therefore, Claim 3 yields $\left[  x^{m}\right]  p_{i,k}=0$.

Now, forget that we fixed $k$. We thus have shown that
\begin{equation}
\left[  x^{m}\right]  p_{i,k}=0\ \ \ \ \ \ \ \ \ \ \text{for each }k\in
S_{i}\setminus\left\{  0\right\}  .
\label{pf.prop.fps.prodrule-fin-inf.c9.pf.1}%
\end{equation}
Now,%
\[
\left[  x^{m}\right]  \left(  \sum_{k\in S_{i}\setminus\left\{  0\right\}
}p_{i,k}\right)  =\sum_{k\in S_{i}\setminus\left\{  0\right\}  }%
\underbrace{\left[  x^{m}\right]  p_{i,k}}_{\substack{=0\\\text{(by
(\ref{pf.prop.fps.prodrule-fin-inf.c9.pf.1}))}}}=\sum_{k\in S_{i}%
\setminus\left\{  0\right\}  }0=0.
\]
This proves Claim 9.]

\begin{statement}
\textit{Claim 10:} Let $n\in\mathbb{N}$. Then,%
\[
\left[  x^{n}\right]  \left(  \prod_{i\in I}\ \ \sum_{k\in S_{i}}%
p_{i,k}\right)  =\left[  x^{n}\right]  \left(  \prod_{i\in I_{n}}%
\ \ \sum_{k\in S_{i}}p_{i,k}\right)  .
\]

\end{statement}

[\textit{Proof of Claim 10:} The set $I_{n}$ is a subset of $I$. Hence, the
set $I$ is the union of its two disjoint subsets $I_{n}$ and $I\setminus
I_{n}$. Thus, we can split the product $\prod_{i\in I}\ \ \sum_{k\in S_{i}%
}p_{i,k}$ as follows:\footnote{Here, we are tacitly using the fact that any
subfamily of the family $\left(  \sum_{k\in S_{i}}p_{i,k}\right)  _{i\in I}$
is multipliable. But this can be proved in the exact same way as we proved
Claim 2.}%
\begin{align}
\prod_{i\in I}\ \ \sum_{k\in S_{i}}p_{i,k}  &  =\left(  \prod_{i\in I_{n}%
}\ \ \sum_{k\in S_{i}}p_{i,k}\right)  \left(  \prod_{i\in I\setminus I_{n}%
}\ \ \underbrace{\sum_{k\in S_{i}}p_{i,k}}_{\substack{=1+\sum_{k\in
S_{i}\setminus\left\{  0\right\}  }p_{i,k}\\\text{(by
(\ref{pf.prop.fps.prodrule-fin-inf.sum-no-1}))}}}\right) \nonumber\\
&  =\left(  \prod_{i\in I_{n}}\ \ \sum_{k\in S_{i}}p_{i,k}\right)  \left(
\prod_{i\in I\setminus I_{n}}\left(  1+\sum_{k\in S_{i}\setminus\left\{
0\right\}  }p_{i,k}\right)  \right)  .
\label{pf.prop.fps.prodrule-fin-inf.c10.pf.1}%
\end{align}

However, the family $\left(  \sum_{k\in S_{i}\setminus\left\{  0\right\}
}p_{i,k}\right)  _{i\in I}$ is summable (as we have seen in the proof of Claim
2). Hence, its subfamily $\left(  \sum_{k\in S_{i}\setminus\left\{  0\right\}
}p_{i,k}\right)  _{i\in I\setminus I_{n}}$ is summable as well (since a
subfamily of a summable family is always summable). Moreover, Claim 9 shows
that each $i\in I\setminus I_{n}$ satisfies $\left[  x^{m}\right]  \left(
\sum_{k\in S_{i}\setminus\left\{  0\right\}  }p_{i,k}\right)  =0$ for each
$m\in\left\{  0,1,\ldots,n\right\}  $. Hence, Lemma
\ref{lem.fps.prod.irlv.inf} (applied to $a=\prod_{i\in I_{n}}\ \ \sum_{k\in
S_{i}}p_{i,k}$ and $J=I\setminus I_{n}$ and $f_{i}=\sum_{k\in S_{i}%
\setminus\left\{  0\right\}  }p_{i,k}$) yields that%
\begin{align*}
&  \left[  x^{m}\right]  \left(  \left(  \prod_{i\in I_{n}}\ \ \sum_{k\in
S_{i}}p_{i,k}\right)  \left(  \prod_{i\in I\setminus I_{n}}\left(
1+\sum_{k\in S_{i}\setminus\left\{  0\right\}  }p_{i,k}\right)  \right)
\right) \\
&  =\left[  x^{m}\right]  \left(  \prod_{i\in I_{n}}\ \ \sum_{k\in S_{i}%
}p_{i,k}\right)  \ \ \ \ \ \ \ \ \ \ \text{for each }m\in\left\{
0,1,\ldots,n\right\}  .
\end{align*}
Applying this to $m=n$, we find%
\[
\left[  x^{n}\right]  \left(  \left(  \prod_{i\in I_{n}}\ \ \sum_{k\in S_{i}%
}p_{i,k}\right)  \left(  \prod_{i\in I\setminus I_{n}}\left(  1+\sum_{k\in
S_{i}\setminus\left\{  0\right\}  }p_{i,k}\right)  \right)  \right)  =\left[
x^{n}\right]  \left(  \prod_{i\in I_{n}}\ \ \sum_{k\in S_{i}}p_{i,k}\right)
.
\]
In view of (\ref{pf.prop.fps.prodrule-fin-inf.c10.pf.1}), this rewrites as%
\[
\left[  x^{n}\right]  \left(  \prod_{i\in I}\ \ \sum_{k\in S_{i}}%
p_{i,k}\right)  =\left[  x^{n}\right]  \left(  \prod_{i\in I_{n}}%
\ \ \sum_{k\in S_{i}}p_{i,k}\right)  .
\]
This proves Claim 10.]

\begin{statement}
\textit{Claim 11:} Let $n\in\mathbb{N}$. Then,%
\[
\prod_{i\in I_{n}}\ \ \sum_{k\in S_{i}}p_{i,k}=\sum_{\left(  k_{i}\right)
_{i\in I_{n}}\in S^{I_{n}}}\ \ \prod_{i\in I_{n}}p_{i,k_{i}}.
\]

\end{statement}

[\textit{Proof of Claim 11:} The set $I_{n}$ is finite. For any $i\in I_{n}$,
the family $\left(  p_{i,k}\right)  _{k\in S_{i}}$ is summable (by Claim 1).
Hence, Proposition \ref{prop.fps.prodrule-fin-infJ} (applied to $N=I_{n}$)
yields
\[
\prod_{i\in I_{n}}\ \ \sum_{k\in S_{i}}p_{i,k}=\sum_{\left(  k_{i}\right)
_{i\in I_{n}}\in\prod_{i\in I_{n}}S_{i}}\ \ \prod_{i\in I_{n}}p_{i,k_{i}}%
=\sum_{\left(  k_{i}\right)  _{i\in I_{n}}\in S^{I_{n}}}\ \ \prod_{i\in I_{n}%
}p_{i,k_{i}}%
\]
(since $\prod_{i\in I_{n}}S_{i}=S^{I_{n}}$). This proves Claim 11.]\medskip

Now, for each $n\in\mathbb{N}$, we have%
\begin{align*}
&  \left[  x^{n}\right]  \left(  \prod_{i\in I}\ \ \sum_{k\in S_{i}}%
p_{i,k}\right) \\
&  =\left[  x^{n}\right]  \left(  \prod_{i\in I_{n}}\ \ \sum_{k\in S_{i}%
}p_{i,k}\right)  \ \ \ \ \ \ \ \ \ \ \left(  \text{by Claim 10}\right) \\
&  =\left[  x^{n}\right]  \left(  \sum_{\left(  k_{i}\right)  _{i\in I_{n}}\in
S^{I_{n}}}\ \ \prod_{i\in I_{n}}p_{i,k_{i}}\right) \\
&  \ \ \ \ \ \ \ \ \ \ \ \ \ \ \ \ \ \ \ \ \left(  \text{since Claim 11 yields
}\prod_{i\in I_{n}}\ \ \sum_{k\in S_{i}}p_{i,k}=\sum_{\left(  k_{i}\right)
_{i\in I_{n}}\in S^{I_{n}}}\ \ \prod_{i\in I_{n}}p_{i,k_{i}}\right) \\
&  =\left[  x^{n}\right]  \left(  \sum_{\left(  k_{i}\right)  _{i\in I}\in
S_{\operatorname*{fin}}^{I}}\ \ \prod_{i\in I}p_{i,k_{i}}\right)
\ \ \ \ \ \ \ \ \ \ \left(  \text{by Claim 8}\right)  .
\end{align*}
That is, any coefficient of the FPS $\prod_{i\in I}\ \ \sum_{k\in S_{i}%
}p_{i,k}$ equals the corresponding coefficient of $\sum_{\left(  k_{i}\right)
_{i\in I}\in S_{\operatorname*{fin}}^{I}}\ \ \prod_{i\in I}p_{i,k_{i}}$.
Hence,%
\[
\prod_{i\in I}\ \ \sum_{k\in S_{i}}p_{i,k}=\sum_{\left(  k_{i}\right)  _{i\in
I}\in S_{\operatorname*{fin}}^{I}}\ \ \prod_{i\in I}p_{i,k_{i}}=\sum
_{\substack{\left(  k_{i}\right)  _{i\in I}\in\prod_{i\in I}S_{i}\\\text{is
essentially finite}}}\ \ \prod_{i\in I}p_{i,k_{i}}%
\]
(since $S_{\operatorname*{fin}}^{I}$ is the set of all essentially finite
families $\left(  k_{i}\right)  _{i\in I}\in\prod_{i\in I}S_{i}$). In
particular, the family $\left(  \prod_{i\in I}p_{i,k_{i}}\right)  _{\left(
k_{i}\right)  _{i\in I}\in\prod_{i\in I}S_{i}\text{ is essentially finite}}$
is summable. This proves Proposition \ref{prop.fps.prodrule-inf-inf}.
\end{proof}
\end{fineprint}

\begin{fineprint}
Our proof of Proposition \ref{prop.fps.subs.rule-infprod} will use the
\textbf{finite} analogue of Proposition \ref{prop.fps.subs.rule-infprod},
which is easy:

\begin{lemma}
\label{lem.fps.subs.rule-infprod-fin}Let $I$ be a \textbf{finite} set. If
$\left(  f_{i}\right)  _{i\in I}\in K\left[  \left[  x\right]  \right]  ^{I}$
is a family of FPSs, and if $g\in K\left[  \left[  x\right]  \right]  $ is an
FPS satisfying $\left[  x^{0}\right]  g=0$, then $\left(  \prod_{i\in I}%
f_{i}\right)  \circ g=\prod_{i\in I}\left(  f_{i}\circ g\right)  $.
\end{lemma}

\begin{proof}
[Proof of Lemma \ref{lem.fps.subs.rule-infprod-fin}.]This follows by a
straightforward induction on $\left\vert I\right\vert $. (The base case is the
case when $\left\vert I\right\vert =0$, and relies on the fact that
$\underline{1}\circ g=\underline{1}$ for any $g\in K\left[  \left[  x\right]
\right]  $. The induction step relies on Proposition \ref{prop.fps.subs.rules}
\textbf{(b)}. The details are left to the reader, who must have seen dozens of
such proofs by now.)
\end{proof}

\begin{proof}
[Detailed proof of Proposition \ref{prop.fps.subs.rule-infprod}.]Let $\left(
f_{i}\right)  _{i\in I}\in K\left[  \left[  x\right]  \right]  ^{I}$ be a
multipliable family of FPSs. Let $g\in K\left[  \left[  x\right]  \right]  $
be an FPS satisfying $\left[  x^{0}\right]  g=0$.

We shall first show the following auxiliary claim:

\begin{statement}
\textit{Claim 1:} Let $M$ be an $x^{n}$-approximator for $\left(
f_{i}\right)  _{i\in I}$. Then, the set $M$ determines the $x^{n}$-coefficient
in the product of $\left(  f_{i}\circ g\right)  _{i\in I}$.
\end{statement}

[\textit{Proof of Claim 1:} The set $M$ is an $x^{n}$-approximator for
$\left(  f_{i}\right)  _{i\in I}$. In other words, $M$ is a finite subset of
$I$ that determines the first $n+1$ coefficients in the product of $\left(
f_{i}\right)  _{i\in I}$ (by the definition of \textquotedblleft$x^{n}%
$-approximator\textquotedblright).

Let $J$ be a finite subset of $I$ satisfying $M\subseteq J\subseteq I$. Then,
Lemma \ref{lem.fps.subs.rule-infprod-fin} (applied to $J$ instead of $I$)
yields%
\begin{equation}
\left(  \prod_{i\in J}f_{i}\right)  \circ g=\prod_{i\in J}\left(  f_{i}\circ
g\right)  . \label{pf.prop.fps.subs.rule-infprod.3}%
\end{equation}
Also, Lemma \ref{lem.fps.subs.rule-infprod-fin} (applied to $M$ instead of
$I$) yields%
\begin{equation}
\left(  \prod_{i\in M}f_{i}\right)  \circ g=\prod_{i\in M}\left(  f_{i}\circ
g\right)  . \label{pf.prop.fps.subs.rule-infprod.4}%
\end{equation}
However, Proposition \ref{prop.fps.infprod-approx-xneq} \textbf{(a)} (applied
to $\mathbf{a}_{i}=f_{i}$) yields%
\[
\prod_{i\in J}f_{i}\overset{x^{n}}{\equiv}\prod_{i\in M}f_{i}.
\]
We also have $g\overset{x^{n}}{\equiv}g$ (since the relation $\overset{x^{n}%
}{\equiv}$ is an equivalence relation). Hence, Proposition
\ref{prop.fps.xneq.comp} (applied to $a=\prod_{i\in J}f_{i}$ and
$b=\prod_{i\in M}f_{i}$ and $c=g$ and $d=g$) yields%
\[
\left(  \prod_{i\in J}f_{i}\right)  \circ g\overset{x^{n}}{\equiv}\left(
\prod_{i\in M}f_{i}\right)  \circ g.
\]
In view of (\ref{pf.prop.fps.subs.rule-infprod.3}) and
(\ref{pf.prop.fps.subs.rule-infprod.4}), this rewrites as%
\[
\prod_{i\in J}\left(  f_{i}\circ g\right)  \overset{x^{n}}{\equiv}\prod_{i\in
M}\left(  f_{i}\circ g\right)  .
\]
In other words,
\[
\text{each }m\in\left\{  0,1,\ldots,n\right\}  \text{ satisfies }\left[
x^{m}\right]  \left(  \prod_{i\in J}\left(  f_{i}\circ g\right)  \right)
=\left[  x^{m}\right]  \left(  \prod_{i\in M}\left(  f_{i}\circ g\right)
\right)
\]
(by the definition of the relation $\overset{x^{n}}{\equiv}$). Applying this
to $m=n$, we obtain%
\[
\left[  x^{n}\right]  \left(  \prod_{i\in J}\left(  f_{i}\circ g\right)
\right)  =\left[  x^{n}\right]  \left(  \prod_{i\in M}\left(  f_{i}\circ
g\right)  \right)  .
\]

Forget that we fixed $J$. We thus have shown that every finite subset $J$ of
$I$ satisfying $M\subseteq J\subseteq I$ satisfies%
\[
\left[  x^{n}\right]  \left(  \prod_{i\in J}\left(  f_{i}\circ g\right)
\right)  =\left[  x^{n}\right]  \left(  \prod_{i\in M}\left(  f_{i}\circ
g\right)  \right)  .
\]
In other words, the set $M$ determines the $x^{n}$-coefficient in the product
of $\left(  f_{i}\circ g\right)  _{i\in I}$ (by the definition of what it
means to \textquotedblleft determine the $x^{n}$-coefficient in the product of
$\left(  f_{i}\circ g\right)  _{i\in I}$\textquotedblright). This proves Claim
1.] \medskip

Now, let $n\in\mathbb{N}$. Lemma \ref{lem.fps.mulable.approx} (applied to
$\mathbf{a}_{i}=f_{i}$) shows that there exists an $x^{n}$-approximator for
$\left(  f_{i}\right)  _{i\in I}$. Consider this $x^{n}$-approximator for
$\left(  f_{i}\right)  _{i\in I}$, and denote it by $M$. Thus, $M$ is an
$x^{n}$-approximator for $\left(  f_{i}\right)  _{i\in I}$; in other words,
$M$ is a finite subset of $I$ that determines the first $n+1$ coefficients in
the product of $\left(  f_{i}\right)  _{i\in I}$ (by the definition of
\textquotedblleft$x^{n}$-approximator\textquotedblright). Claim 1 shows that
the set $M$ determines the $x^{n}$-coefficient in the product of $\left(
f_{i}\circ g\right)  _{i\in I}$. Hence, there is a finite subset of $I$ that
determines the $x^{n}$-coefficient in the product of $\left(  f_{i}\circ
g\right)  _{i\in I}$ (namely, $M$). In other words, the $x^{n}$-coefficient in
the product of $\left(  f_{i}\circ g\right)  _{i\in I}$ is finitely determined.

Forget that we fixed $n$. We thus have shown that for each $n\in\mathbb{N}$,
the $x^{n}$-coefficient in the product of $\left(  f_{i}\circ g\right)  _{i\in
I}$ is finitely determined. In other words, each coefficient in the product of
$\left(  f_{i}\circ g\right)  _{i\in I}$ is finitely determined. In other
words, the family $\left(  f_{i}\circ g\right)  _{i\in I}$ is multipliable (by
the definition of \textquotedblleft multipliable\textquotedblright). \medskip

It remains to prove that $\left(  \prod_{i\in I}f_{i}\right)  \circ
g=\prod_{i\in I}\left(  f_{i}\circ g\right)  $.

In order to do so, we again fix $n\in\mathbb{N}$. Lemma
\ref{lem.fps.mulable.approx} (applied to $\mathbf{a}_{i}=f_{i}$) shows that
there exists an $x^{n}$-approximator for $\left(  f_{i}\right)  _{i\in I}$.
Consider this $x^{n}$-approximator for $\left(  f_{i}\right)  _{i\in I}$, and
denote it by $M$. Thus, $M$ is an $x^{n}$-approximator for $\left(
f_{i}\right)  _{i\in I}$; in other words, $M$ is a finite subset of $I$ that
determines the first $n+1$ coefficients in the product of $\left(
f_{i}\right)  _{i\in I}$ (by the definition of \textquotedblleft$x^{n}%
$-approximator\textquotedblright). Moreover, Proposition
\ref{prop.fps.infprod-approx-xneq} \textbf{(b)} (applied to $\mathbf{a}%
_{i}=f_{i}$) yields%
\begin{equation}
\prod_{i\in I}f_{i}\overset{x^{n}}{\equiv}\prod_{i\in M}f_{i}.
\label{pf.prop.fps.subs.rule-infprod.1}%
\end{equation}
We also have $g\overset{x^{n}}{\equiv}g$ (since the relation $\overset{x^{n}%
}{\equiv}$ is an equivalence relation). Hence, Proposition
\ref{prop.fps.xneq.comp} (applied to $a=\prod_{i\in I}f_{i}$ and
$b=\prod_{i\in M}f_{i}$ and $c=g$ and $d=g$) yields%
\begin{equation}
\left(  \prod_{i\in I}f_{i}\right)  \circ g\overset{x^{n}}{\equiv}\left(
\prod_{i\in M}f_{i}\right)  \circ g. \label{pf.prop.fps.subs.rule-infprod.2}%
\end{equation}
However, Lemma \ref{lem.fps.subs.rule-infprod-fin} (applied to $M$ instead of
$I$) yields%
\[
\left(  \prod_{i\in M}f_{i}\right)  \circ g=\prod_{i\in M}\left(  f_{i}\circ
g\right)  .
\]
In view of this, we can rewrite (\ref{pf.prop.fps.subs.rule-infprod.2}) as
\[
\left(  \prod_{i\in I}f_{i}\right)  \circ g\overset{x^{n}}{\equiv}\prod_{i\in
M}\left(  f_{i}\circ g\right)  .
\]
In other words,
\[
\text{each }m\in\left\{  0,1,\ldots,n\right\}  \text{ satisfies }\left[
x^{m}\right]  \left(  \left(  \prod_{i\in I}f_{i}\right)  \circ g\right)
=\left[  x^{m}\right]  \left(  \prod_{i\in M}\left(  f_{i}\circ g\right)
\right)
\]
(by the definition of the relation $\overset{x^{n}}{\equiv}$). Applying this
to $m=n$, we obtain%
\begin{equation}
\left[  x^{n}\right]  \left(  \left(  \prod_{i\in I}f_{i}\right)  \circ
g\right)  =\left[  x^{n}\right]  \left(  \prod_{i\in M}\left(  f_{i}\circ
g\right)  \right)  . \label{pf.prop.fps.subs.rule-infprod.6}%
\end{equation}

However, Claim 1 shows that the set $M$ determines the $x^{n}$-coefficient in
the product of $\left(  f_{i}\circ g\right)  _{i\in I}$. Hence, the definition
of the infinite product $\prod_{i\in I}\left(  f_{i}\circ g\right)  $
(specifically, Definition \ref{def.fps.multipliable} \textbf{(b)}) yields%
\[
\left[  x^{n}\right]  \left(  \prod_{i\in I}\left(  f_{i}\circ g\right)
\right)  =\left[  x^{n}\right]  \left(  \prod_{i\in M}\left(  f_{i}\circ
g\right)  \right)  .
\]
Comparing this with (\ref{pf.prop.fps.subs.rule-infprod.6}), we obtain%
\[
\left[  x^{n}\right]  \left(  \left(  \prod_{i\in I}f_{i}\right)  \circ
g\right)  =\left[  x^{n}\right]  \left(  \prod_{i\in I}\left(  f_{i}\circ
g\right)  \right)  .
\]

Forget that we fixed $n$. We thus have shown that each $n\in\mathbb{N}$
satisfies%
\[
\left[  x^{n}\right]  \left(  \left(  \prod_{i\in I}f_{i}\right)  \circ
g\right)  =\left[  x^{n}\right]  \left(  \prod_{i\in I}\left(  f_{i}\circ
g\right)  \right)  .
\]
In other words, the FPSs $\left(  \prod_{i\in I}f_{i}\right)  \circ g$ and
$\prod_{i\in I}\left(  f_{i}\circ g\right)  $ agree in all their coefficients.
Hence, $\left(  \prod_{i\in I}f_{i}\right)  \circ g=\prod_{i\in I}\left(
f_{i}\circ g\right)  $. This completes our proof of Proposition
\ref{prop.fps.subs.rule-infprod}.
\end{proof}
\end{fineprint}

\subsection{\label{sec.details.domino}Domino tilings}

\begin{fineprint}
In Subsection \ref{subsec.gf.weighted-set.domino}, we have given a
classification of faultfree domino tilings of the rectangle $R_{n,3}$. Let us
now outline how this classification can be proved. First, we introduce names
for the faultfree domino tilings that appear in our classification:

\begin{definition}
\label{def.gf.weighted-set.domino.Rn3.ABC}\textbf{(a)} For each even positive
integer $n$, we let $A_{n}$ be the domino tiling%
\[%
\qquad\begin{tikzpicture}[scale=0.8]
\draw(0, 0) rectangle (2, 1);
\draw(2, 1) rectangle (4, 0);
\draw(6, 1) rectangle (4, 0);
\draw(0, 1) rectangle (1, 3);
\draw(1, 1) rectangle (3, 2);
\draw(1, 3) rectangle (3, 2);
\draw(3, 1) rectangle (5, 2);
\draw(3, 3) rectangle (5, 2);
\node(X) at (8, 2) {\Large$\cdots\cdots\cdots$};
\node(X) at (8, 0.5) {\Large$\cdots$};
\draw(10, 0) rectangle (12, 1);
\draw(12, 0) rectangle (14, 1);
\draw(11, 3) rectangle (13, 2);
\draw(11, 2) rectangle (13, 1);
\draw(13, 1) rectangle (14, 3);
\end{tikzpicture}%
\ \ \ \ \ \ \ \ \ \ \text{of }R_{n,3}.
\]
Formally speaking, this domino tiling is the set partition of $R_{n,3}$
consisting of the following dominos:

\begin{itemize}
\item the horizontal dominos $\left\{  \left(  2i-1,\ 1\right)  ,\ \left(
2i,\ 1\right)  \right\}  $ for all $i\in\left[  n/2\right]  $, which fill the
bottom row of $R_{n,3}$, and which we call the \textbf{basement dominos};

\item the vertical domino $\left\{  \left(  1,2\right)  ,\ \left(  1,3\right)
\right\}  $ in the first column, which we call the \textbf{left wall};

\item the vertical domino $\left\{  \left(  n,2\right)  ,\ \left(  n,3\right)
\right\}  $ in the last column, which we call the \textbf{right wall};

\item the horizontal dominos $\left\{  \left(  2i,\ 2\right)  ,\ \left(
2i+1,\ 2\right)  \right\}  $ for all $i\in\left[  n/2-1\right]  $, which fill
the middle row of $R_{n,3}$ (except for the first and last columns), and which
we call the \textbf{middle dominos};

\item the horizontal dominos $\left\{  \left(  2i,\ 3\right)  ,\ \left(
2i+1,\ 3\right)  \right\}  $ for all $i\in\left[  n/2-1\right]  $, which fill
the top row of $R_{n,3}$ (except for the first and last columns), and which we
call the \textbf{top dominos}.
\end{itemize}

\textbf{(b)} For each even positive integer $n$, we let $B_{n}$ be the domino
tiling%
\[%
\qquad\begin{tikzpicture}[scale=0.8]
\draw(0, 3) rectangle (2, 2);
\draw(2, 2) rectangle (4, 3);
\draw(6, 2) rectangle (4, 3);
\draw(0, 2) rectangle (1, 0);
\draw(1, 2) rectangle (3, 1);
\draw(1, 0) rectangle (3, 1);
\draw(3, 2) rectangle (5, 1);
\draw(3, 0) rectangle (5, 1);
\node(X) at (8, 1) {\Large$\cdots\cdots\cdots$};
\node(X) at (8, 2.5) {\Large$\cdots$};
\draw(10, 3) rectangle (12, 2);
\draw(12, 3) rectangle (14, 2);
\draw(11, 0) rectangle (13, 1);
\draw(11, 1) rectangle (13, 2);
\draw(13, 2) rectangle (14, 0);
\end{tikzpicture}%
\ \ \ \ \ \ \ \ \ \ \text{of }R_{n,3}.
\]
This is the reflection of the tiling $A_{n}$ across the horizontal axis of
symmetry of $R_{n,3}$. \medskip

\textbf{(c)} We let $C$ denote the domino tiling
\[%
\begin{tikzpicture}[scale=0.8]
\draw(0, 0) rectangle (2, 1);
\draw(0, 1) rectangle (2, 2);
\draw(0, 2) rectangle (2, 3);
\end{tikzpicture}%
\ \ \ \ \ \ \ \ \ \ \text{of }R_{2,3}.
\]

\end{definition}

Our classification now can be stated as follows:

\begin{proposition}
\label{prop.gf.weighted-set.domino.Rn3.ABC}The faultfree domino tilings of
height-$3$ rectangles are precisely the tilings%
\[
A_{2},A_{4},A_{6},A_{8},\ldots,\ \ \ \ \ \ \ \ \ \ B_{2},B_{4},B_{6}%
,B_{8},\ldots,\ \ \ \ \ \ \ \ \ \ C
\]
we have just defined. More concretely: \medskip

\textbf{(a)} The faultfree domino tilings of a height-$3$ rectangle that
contain a vertical domino in the \textbf{top} two squares of the first column
are $A_{2},A_{4},A_{6},A_{8},\ldots$. \medskip

\textbf{(b)} The faultfree domino tilings of a height-$3$ rectangle that
contain a vertical domino in the \textbf{bottom} two squares of the first
column are $B_{2},B_{4},B_{6},B_{8},\ldots$. \medskip

\textbf{(c)} The only faultfree domino tiling of a height-$3$ rectangle that
contains \textbf{no} vertical domino in the first column is $C$.
\end{proposition}

\begin{proof}
[Proof of Proposition \ref{prop.gf.weighted-set.domino.Rn3.ABC} (sketched).]It
clearly suffices to prove parts \textbf{(a)}, \textbf{(b)} and \textbf{(c)}.
We begin with the easiest part, which is \textbf{(c)}: \medskip

\textbf{(c)} Clearly, $C$ is a faultfree domino tiling of $R_{2,3}$ that
contains no vertical domino in the first column. It remains to show that $C$
is the only such tiling.

Indeed, let $T$ be any faultfree domino tiling of a height-$3$ rectangle that
contains no vertical domino in the first column. Thus, its first column must
be filled with three horizontal dominos, which all must protrude into the
second column and thus cover that second column as well. If $T$ had any
further column, then $T$ would have a fault between its $2$-nd and its $3$-rd
column, which is impossible for a faultfree tiling. Thus, $T$ must consist
only of the three horizontal dominos already mentioned. In other words, $T=C$.
This completes the proof of Proposition
\ref{prop.gf.weighted-set.domino.Rn3.ABC} \textbf{(c)}. \medskip

\textbf{(a)} It is straightforward to check that the tilings $A_{2}%
,A_{4},A_{6},A_{8},\ldots$ are faultfree (indeed, the basement dominos prevent
faults between the $\left(  2i-1\right)  $-st and $\left(  2i\right)  $-th
columns, whereas the top dominos prevent faults between the $\left(
2i\right)  $-th and $\left(  2i+1\right)  $-st columns). Thus, they are
faultfree domino tilings of a height-$3$ rectangle that contain a vertical
domino in the \textbf{top} two squares of the first column. It remains to show
that they are the only such tilings.

Indeed, let $T$ be any faultfree domino tiling of a height-$3$ rectangle that
contains a vertical domino in the \textbf{top} two squares of the first
column. Let $n$ be the width of this rectangle (so that the rectangle is
$R_{n,3}$). We shall show that $n$ is even, and that $T=A_{n}$.

We know that $T$ contains a vertical domino in the \textbf{top} two squares of
the first column. In other words, $T$ contains the left wall (where we are
using the terminology from Definition \ref{def.gf.weighted-set.domino.Rn3.ABC}
\textbf{(a)}). The remaining square in the first column is the square $\left(
1,1\right)  $, and it must thus be covered by the basement domino $\left\{
\left(  1,1\right)  ,\ \left(  2,1\right)  \right\}  $ (since no other domino
would fit). Hence, $T$ contains the basement domino $\left\{  \left(
1,1\right)  ,\ \left(  2,1\right)  \right\}  $.

We shall now prove the following claims:

\begin{statement}
\textit{Claim 1:} For each positive integer $i<n/2$, the tiling $T$ contains
the basement domino $\left\{  \left(  2i-1,\ 1\right)  ,\ \left(
2i,\ 1\right)  \right\}  $, the middle domino $\left\{  \left(  2i,\ 2\right)
,\ \left(  2i+1,\ 2\right)  \right\}  $ and the top domino $\left\{  \left(
2i,\ 3\right)  ,\ \left(  2i+1,\ 3\right)  \right\}  $.
\end{statement}

\begin{proof}
[Proof of Claim 1.]We induct on $i$:

\textit{Base case:} We must prove that Claim 1 holds for $i=1$, provided that
$1<n/2$. So let us assume that $1<n/2$. Thus, $2<n$, so that $n\geq3$. We must
prove that Claim 1 holds for $i=1$; in other words, we must prove that $T$
contains the basement domino $\left\{  \left(  1,1\right)  ,\ \left(
2,1\right)  \right\}  $, the middle domino $\left\{  \left(  2,2\right)
,\ \left(  3,2\right)  \right\}  $ and the top domino $\left\{  \left(
2,3\right)  ,\ \left(  3,3\right)  \right\}  $. For the basement domino, we
have already proved this. For the middle and the top domino, we argue as
follows: If $T$ had a vertical domino in the $2$-nd column, then this domino
would cover the top two squares of that column (since the bottom square is
already covered by the basement domino $\left\{  \left(  1,1\right)
,\ \left(  2,1\right)  \right\}  $), and thus $T$ would have a fault between
the $2$-nd and $3$-rd columns (since the basement domino $\left\{  \left(
1,1\right)  ,\ \left(  2,1\right)  \right\}  $ also ends at the $2$-nd
column), which would contradict the faultfreeness of $T$. Hence, $T$ has no
vertical domino in the $2$-nd column. Thus, the top two squares of the $2$-nd
column of $T$ must be covered by horizontal dominos. These horizontal dominos
must both protrude into the $3$-rd column (since the corresponding squares in
the $1$-st column are already covered by the left wall), and thus must be the
middle domino $\left\{  \left(  2,2\right)  ,\ \left(  3,2\right)  \right\}  $
and the top domino $\left\{  \left(  2,3\right)  ,\ \left(  3,3\right)
\right\}  $. Hence, we have shown that $T$ contains the middle domino
$\left\{  \left(  2,2\right)  ,\ \left(  3,2\right)  \right\}  $ and the top
domino $\left\{  \left(  2,3\right)  ,\ \left(  3,3\right)  \right\}  $. This
completes the base case.

\textit{Induction step:} Let $j$ be a positive integer such that $j+1<n/2$.
Assume (as the induction hypothesis) that Claim 1 holds for $i=j$. In other
words, $T$ contains the basement domino $\left\{  \left(  2j-1,\ 1\right)
,\ \left(  2j,\ 1\right)  \right\}  $, the middle domino $\left\{  \left(
2j,\ 2\right)  ,\ \left(  2j+1,\ 2\right)  \right\}  $ and the top domino
$\left\{  \left(  2j,\ 3\right)  ,\ \left(  2j+1,\ 3\right)  \right\}  $. We
must now show that Claim 1 holds for $i=j+1$ as well, i.e., that $T$ also
contains the basement domino $\left\{  \left(  2j+1,\ 1\right)  ,\ \left(
2j+2,\ 1\right)  \right\}  $, the middle domino $\left\{  \left(
2j+2,\ 2\right)  ,\ \left(  2j+3,\ 2\right)  \right\}  $ and the top domino
$\left\{  \left(  2j+2,\ 3\right)  ,\ \left(  2j+3,\ 3\right)  \right\}  $.

Indeed, we first recall that $j+1<n/2$, so that $2\left(  j+1\right)  <n$ and
therefore $n>2\left(  j+1\right)  =2j+2$, so that $n\geq2j+3$. This shows that
the rectangle $R_{n,3}$ has a $\left(  2j+3\right)  $-th column (along with
all columns to its left).

Now, the square $\left(  2j+1,\ 1\right)  $ cannot be covered by a vertical
domino in $T$, since this vertical domino would collide with the middle domino
$\left\{  \left(  2j,\ 2\right)  ,\ \left(  2j+1,\ 2\right)  \right\}  $
(which we already know to belong to $T$). Thus, this square must be covered by
a horizontal domino. This horizontal domino cannot protrude into the $\left(
2j\right)  $-th column (since it would then collide with the basement domino
$\left\{  \left(  2j-1,\ 1\right)  ,\ \left(  2j,\ 1\right)  \right\}  $,
which we already know to belong to $T$), and thus must be the basement domino
$\left\{  \left(  2j+1,\ 1\right)  ,\ \left(  2j+2,\ 1\right)  \right\}  $. So
we have shown that $T$ contains the basement domino $\left\{  \left(
2j+1,\ 1\right)  ,\ \left(  2j+2,\ 1\right)  \right\}  $. It remains to show
that $T$ also contains the middle domino $\left\{  \left(  2j+2,\ 2\right)
,\ \left(  2j+3,\ 2\right)  \right\}  $ and the top domino $\left\{  \left(
2j+2,\ 3\right)  ,\ \left(  2j+3,\ 3\right)  \right\}  $.

If $T$ had a vertical domino in the $\left(  2j+2\right)  $-nd column, then
this domino would cover the top two squares of that column (since the bottom
square is already covered by the basement domino $\left\{  \left(
2j+1,\ 1\right)  ,\ \left(  2j+2,\ 1\right)  \right\}  $), and thus $T$ would
have a fault between the $\left(  2j+2\right)  $-nd and $\left(  2j+3\right)
$-rd columns (since the basement domino $\left\{  \left(  2j+1,\ 1\right)
,\ \left(  2j+2,\ 1\right)  \right\}  $ also ends at the $\left(  2j+2\right)
$-nd column), which would contradict the faultfreeness of $T$. Hence, $T$ has
no vertical domino in the $\left(  2j+2\right)  $-nd column. Thus, the top two
squares of the $\left(  2j+2\right)  $-nd column of $T$ must be covered by
horizontal dominos. These horizontal dominos must both protrude into the
$\left(  2j+3\right)  $-rd column (since the corresponding squares in the
$\left(  2j+1\right)  $-st column are already covered by the middle domino
$\left\{  \left(  2j,\ 2\right)  ,\ \left(  2j+1,\ 2\right)  \right\}  $ and
the top domino $\left\{  \left(  2j,\ 3\right)  ,\ \left(  2j+1,\ 3\right)
\right\}  $), and thus must be the middle domino $\left\{  \left(
2j+2,\ 2\right)  ,\ \left(  2j+3,\ 2\right)  \right\}  $ and the top domino
$\left\{  \left(  2j+2,\ 3\right)  ,\ \left(  2j+3,\ 3\right)  \right\}  $.
Hence, we have shown that $T$ contains the middle domino $\left\{  \left(
2j+2,\ 2\right)  ,\ \left(  2j+3,\ 2\right)  \right\}  $ and the top domino
$\left\{  \left(  2j+2,\ 3\right)  ,\ \left(  2j+3,\ 3\right)  \right\}  $.
This completes the induction step.

Thus, Claim 1 is proved by induction.
\end{proof}

\begin{statement}
\textit{Claim 2:} The number $n$ is even.
\end{statement}

\begin{proof}
[Proof of Claim 2.]Assume the contrary. Thus, $n$ is odd. But $n>1$ (since
$R_{1,3}$ cannot be tiled by dominos), so that $n-1>0$. Hence, $\left(
n-1\right)  /2$ is a positive integer (since $n$ is odd). Therefore, Claim 1
(applied to $i=\left(  n-1\right)  /2$) shows that the tiling $T$ contains the
basement domino $\left\{  \left(  n-2,\ 1\right)  ,\ \left(  n-1,\ 1\right)
\right\}  $, the middle domino $\left\{  \left(  n-1,\ 2\right)  ,\ \left(
n,\ 2\right)  \right\}  $ and the top domino $\left\{  \left(  n-1,\ 3\right)
,\ \left(  n,\ 3\right)  \right\}  $. But $T$ must also include a domino that
contains the square $\left(  n,\ 1\right)  $ (since $\left(  n,\ 1\right)  \in
R_{n,3}$). This domino cannot be horizontal (since it would then collide with
the basement domino $\left\{  \left(  n-2,\ 1\right)  ,\ \left(
n-1,\ 1\right)  \right\}  $, which we know to belong to $T$), and cannot be
vertical either (since it would then collide with the middle domino $\left\{
\left(  n-1,\ 2\right)  ,\ \left(  n,\ 2\right)  \right\}  $). This is clearly
absurd. This contradiction shows that our assumption was wrong, so that Claim
2 is proved.

(Alternatively, Claim 2 can be obtained from a parity argument: Since $T$ is a
tiling of $R_{n,3}$ by dominos, the total \# of squares of $R_{n,3}$ must be
even (since each domino covers exactly $2$ squares). But this total \# is
$3n$. Thus, $3n$ must be even, so that $n$ must be even.)
\end{proof}

Now, Claim 2 shows that $n$ is even, so that $n/2$ is a positive integer.
Furthermore, we know that the tiling $T$ contains the left wall, the first
$n/2-1$ basement dominos
\[
\left\{  \left(  1,1\right)  ,\ \left(  2,1\right)  \right\}
,\ \ \ \ \ \left\{  \left(  3,1\right)  ,\ \left(  4,1\right)  \right\}
,\ \ \ \ \ \ldots,\ \ \ \ \ \left\{  \left(  n-3,\ 1\right)  ,\ \left(
n-2,\ 1\right)  \right\}
\]
(by Claim 1), all the middle dominos
\[
\left\{  \left(  2,2\right)  ,\ \left(  3,2\right)  \right\}
,\ \ \ \ \ \left\{  \left(  4,2\right)  ,\ \left(  5,2\right)  \right\}
,\ \ \ \ \ \ldots,\ \ \ \ \ \left\{  \left(  n-2,\ 2\right)  ,\ \left(
n-1,\ 2\right)  \right\}
\]
(also by Claim 1) and all the top dominos%
\[
\left\{  \left(  2,3\right)  ,\ \left(  3,3\right)  \right\}
,\ \ \ \ \ \left\{  \left(  4,3\right)  ,\ \left(  5,3\right)  \right\}
,\ \ \ \ \ \ldots,\ \ \ \ \ \left\{  \left(  n-2,\ 3\right)  ,\ \left(
n-1,\ 3\right)  \right\}
\]
(also by Claim 1). This leaves only the four squares $\left(  n-1,\ 1\right)
$, $\left(  n,\ 1\right)  $, $\left(  n,\ 2\right)  $ and $\left(
n,\ 3\right)  $ unaccounted for, but there is only one way to tile them:
namely, with the last basement domino $\left\{  \left(  n-1,\ 1\right)
,\ \left(  n,\ 1\right)  \right\}  $ and the right wall $\left\{  \left(
n,2\right)  ,\ \left(  n,3\right)  \right\}  $. Thus, $T$ must contain all the
basement dominos, all the middle dominos, all the top dominos and both walls
(left and right). Since these dominos cover all the squares of $R_{n,3}$, this
entails that $T$ consists precisely of these dominos. In other words,
$T=A_{n}$. The proof of Proposition \ref{prop.gf.weighted-set.domino.Rn3.ABC}
\textbf{(a)} is thus complete. \medskip

\textbf{(b)} Proposition \ref{prop.gf.weighted-set.domino.Rn3.ABC}
\textbf{(b)} is just Proposition \ref{prop.gf.weighted-set.domino.Rn3.ABC}
\textbf{(a)} reflected across the horizontal axis of symmetry of $R_{n,3}$.
\end{proof}
\end{fineprint}

\subsection{\label{sec.details.gf.lim}Limits of FPSs}

\begin{fineprint}
Let us now prove a few facts stated in Section \ref{sec.gf.lim}.

\begin{proof}
[Detailed proof of Lemma \ref{lem.fps.lim.xn-equiv}.]We have $\lim
\limits_{i\rightarrow\infty}f_{i}=f$. In other words, the sequence $\left(
f_{i}\right)  _{i\in\mathbb{N}}$ coefficientwise stabilizes to $f$. In other
words, for each $n\in\mathbb{N}$,%
\begin{equation}
\text{the sequence }\left(  \left[  x^{n}\right]  f_{i}\right)  _{i\in
\mathbb{N}}\text{ stabilizes to }\left[  x^{n}\right]  f
\label{pf.lem.fps.lim.xn-equiv.stab}%
\end{equation}
(by the definition of \textquotedblleft coefficientwise
stabilizing\textquotedblright).

Now, let $n\in\mathbb{N}$. Furthermore, let $k\in\left\{  0,1,\ldots
,n\right\}  $. Then, the sequence $\left(  \left[  x^{k}\right]  f_{i}\right)
_{i\in\mathbb{N}}$ stabilizes to $\left[  x^{k}\right]  f$ (by
(\ref{pf.lem.fps.lim.xn-equiv.stab}), applied to $k$ instead of $n$). In other
words, there exists some $N\in\mathbb{N}$ such that%
\[
\text{all integers }i\geq N\text{ satisfy }\left[  x^{k}\right]  f_{i}=\left[
x^{k}\right]  f
\]
(by the definition of \textquotedblleft stabilizing\textquotedblright). Let us
denote this $N$ by $N_{k}$. Thus,
\begin{equation}
\text{all integers }i\geq N_{k}\text{ satisfy }\left[  x^{k}\right]
f_{i}=\left[  x^{k}\right]  f. \label{pf.lem.fps.lim.xn-equiv.Nk}%
\end{equation}

Forget that we fixed $k$. Thus, for each $k\in\left\{  0,1,\ldots,n\right\}
$, we have defined an integer $N_{k}\in\mathbb{N}$ for which
(\ref{pf.lem.fps.lim.xn-equiv.Nk}) holds. Altogether, we have thus defined
$n+1$ integers $N_{0},N_{1},\ldots,N_{n}\in\mathbb{N}$.

Let us set $P:=\max\left\{  N_{0},N_{1},\ldots,N_{n}\right\}  $. Then, of
course, $P\in\mathbb{N}$.

Now, let $i\geq P$ be an integer. Then, for each $k\in\left\{  0,1,\ldots
,n\right\}  $, we have $i\geq P=\max\left\{  N_{0},N_{1},\ldots,N_{n}\right\}
\geq N_{k}$ (since $k\in\left\{  0,1,\ldots,n\right\}  $) and therefore
$\left[  x^{k}\right]  f_{i}=\left[  x^{k}\right]  f$ (by
(\ref{pf.lem.fps.lim.xn-equiv.Nk})). In other words, each $k\in\left\{
0,1,\ldots,n\right\}  $ satisfies $\left[  x^{k}\right]  f_{i}=\left[
x^{k}\right]  f$. Renaming the variable $k$ as $m$ in this result, we obtain
the following:%
\[
\text{Each }m\in\left\{  0,1,\ldots,n\right\}  \text{ satisfies }\left[
x^{m}\right]  f_{i}=\left[  x^{m}\right]  f\text{.}%
\]
In other words, $f_{i}\overset{x^{n}}{\equiv}f$ (by the definition of $x^{n}%
$-equivalence\footnote{i.e., Definition \ref{def.fps.xneq}}).

Forget that we fixed $i$. We thus have shown that all integers $i\geq P$
satisfy $f_{i}\overset{x^{n}}{\equiv}f$. Hence, there exists some integer
$N\in\mathbb{N}$ such that
\[
\text{all integers }i\geq N\text{ satisfy }f_{i}\overset{x^{n}}{\equiv}f
\]
(namely, $N=P$). This proves Lemma \ref{lem.fps.lim.xn-equiv}.
\end{proof}

\begin{proof}
[Detailed proof of Proposition \ref{prop.fps.lim.sum-prod}.]Recall that
$\lim\limits_{i\rightarrow\infty}f_{i}=f$. Hence, Lemma
\ref{lem.fps.lim.xn-equiv} shows that there exists some integer $N\in
\mathbb{N}$ such that
\[
\text{all integers }i\geq N\text{ satisfy }f_{i}\overset{x^{n}}{\equiv}f.
\]
Let us denote this $N$ by $K$. Hence, all integers $i\geq K$ satisfy
$f_{i}\overset{x^{n}}{\equiv}f$.

Thus, we have found an integer $K\in\mathbb{N}$ such that%
\begin{equation}
\text{all integers }i\geq K\text{ satisfy }f_{i}\overset{x^{n}}{\equiv}f.
\label{prop.fps.lim.sum-prod.K}%
\end{equation}
Similarly, using $\lim\limits_{i\rightarrow\infty}g_{i}=g$, we can find an
integer $L\in\mathbb{N}$ such that
\begin{equation}
\text{all integers }i\geq L\text{ satisfy }g_{i}\overset{x^{n}}{\equiv}g.
\label{prop.fps.lim.sum-prod.L}%
\end{equation}
Consider these $K$ and $L$. Let us furthermore set $P:=\max\left\{
K,L\right\}  $. Thus, $P\in\mathbb{N}$.

We shall now show that each integer $i\geq P$ satisfies $\left[  x^{n}\right]
\left(  f_{i}g_{i}\right)  =\left[  x^{n}\right]  \left(  fg\right)  $.

Indeed, let $i\geq P$ be an integer. Then, (\ref{prop.fps.lim.sum-prod.K})
yields $f_{i}\overset{x^{n}}{\equiv}f$ (since $i\geq P=\max\left\{
K,L\right\}  \geq K$), whereas (\ref{prop.fps.lim.sum-prod.L}) yields
$g_{i}\overset{x^{n}}{\equiv}g$ (since $i\geq P=\max\left\{  K,L\right\}  \geq
L$). Hence, we obtain%
\[
f_{i}g_{i}\overset{x^{n}}{\equiv}fg
\]
(by (\ref{eq.thm.fps.xneq.props.b.*}), applied to $a=f_{i}$, $b=f$, $c=g_{i}$
and $d=g$). In other words,
\[
\text{each }m\in\left\{  0,1,\ldots,n\right\}  \text{ satisfies }\left[
x^{m}\right]  \left(  f_{i}g_{i}\right)  =\left[  x^{m}\right]  \left(
fg\right)
\]
(by the definition of $x^{n}$-equivalence). Applying this to $m=n$, we find%
\[
\left[  x^{n}\right]  \left(  f_{i}g_{i}\right)  =\left[  x^{n}\right]
\left(  fg\right)  .
\]

Forget that we fixed $i$. We thus have shown that all integers $i\geq P$
satisfy $\left[  x^{n}\right]  \left(  f_{i}g_{i}\right)  =\left[
x^{n}\right]  \left(  fg\right)  $. Hence, there exists some $N\in\mathbb{N}$
such that%
\[
\text{all integers }i\geq N\text{ satisfy }\left[  x^{n}\right]  \left(
f_{i}g_{i}\right)  =\left[  x^{n}\right]  \left(  fg\right)
\]
(namely, $N=P$). In other words,%
\[
\text{the sequence }\left(  \left[  x^{n}\right]  \left(  f_{i}g_{i}\right)
\right)  _{i\in\mathbb{N}}\text{ stabilizes to }\left[  x^{n}\right]  \left(
fg\right)
\]
(by the definition of \textquotedblleft stabilizes\textquotedblright).

Forget that we fixed $n$. We thus have shown that for each $n\in\mathbb{N}$,
the sequence $\left(  \left[  x^{n}\right]  \left(  f_{i}g_{i}\right)
\right)  _{i\in\mathbb{N}}$ stabilizes to $\left[  x^{n}\right]  \left(
fg\right)  $. In other words, the sequence $\left(  f_{i}g_{i}\right)
_{i\in\mathbb{N}}$ coefficientwise stabilizes to $fg$ (by the definition of
\textquotedblleft coefficientwise stabilizing\textquotedblright). In other
words, $\lim\limits_{i\rightarrow\infty}\left(  f_{i}g_{i}\right)  =fg$.

The same argument -- but using (\ref{eq.thm.fps.xneq.props.b.+}) instead of
(\ref{eq.thm.fps.xneq.props.b.*}) -- shows that $\lim\limits_{i\rightarrow
\infty}\left(  f_{i}+g_{i}\right)  =f+g$. Thus, the proof of Proposition
\ref{prop.fps.lim.sum-prod} is complete.
\end{proof}

\begin{proof}
[Detailed proof of Proposition \ref{prop.fps.lim.sum-quot}.]First, let us show
that $g$ is invertible:

\begin{statement}
\textit{Claim 1:} The FPS $g\in K\left[  \left[  x\right]  \right]  $ is invertible.
\end{statement}

\begin{proof}
[Proof of Claim 1.]We have assumed that $\lim\limits_{i\rightarrow\infty}%
g_{i}=g$. In other words, the sequence $\left(  g_{i}\right)  _{i\in
\mathbb{N}}$ coefficientwise stabilizes to $g$. In other words, for each
$n\in\mathbb{N}$, the sequence $\left(  \left[  x^{n}\right]  g_{i}\right)
_{i\in\mathbb{N}}$ stabilizes to $\left[  x^{n}\right]  g$ (by the definition
of \textquotedblleft coefficientwise stabilizing\textquotedblright). Applying
this to $n=0$, we see that the sequence $\left(  \left[  x^{0}\right]
g_{i}\right)  _{i\in\mathbb{N}}$ stabilizes to $\left[  x^{0}\right]  g$. In
other words, there exists some $N\in\mathbb{N}$ such that%
\[
\text{all integers }i\geq N\text{ satisfy }\left[  x^{0}\right]  g_{i}=\left[
x^{0}\right]  g
\]
(by the definition of \textquotedblleft stabilizes\textquotedblright).
Consider this $N$. Thus, all integers $i\geq N$ satisfy $\left[  x^{0}\right]
g_{i}=\left[  x^{0}\right]  g$. Applying this to $i=N$, we obtain $\left[
x^{0}\right]  g_{N}=\left[  x^{0}\right]  g$.

However, each FPS $g_{i}$ is invertible (by assumption). Hence, in particular,
the FPS $g_{N}$ is invertible. By Proposition \ref{prop.fps.invertible}
(applied to $a=g_{N}$), this entails that its constant term $\left[
x^{0}\right]  g_{N}$ is invertible in $K$. In other words, $\left[
x^{0}\right]  g$ is invertible in $K$ (since $\left[  x^{0}\right]
g_{N}=\left[  x^{0}\right]  g$).

But the FPS $g$ is invertible if and only if its constant term $\left[
x^{0}\right]  g$ is invertible in $K$ (by Proposition
\ref{prop.fps.invertible}, applied to $a=g$). Hence, $g$ is invertible (since
its constant term $\left[  x^{0}\right]  g$ is invertible in $K$). This proves
Claim 1.
\end{proof}

It remains to prove that $\lim\limits_{i\rightarrow\infty}\dfrac{f_{i}}{g_{i}%
}=\dfrac{f}{g}$.

Recall that $\lim\limits_{i\rightarrow\infty}f_{i}=f$. Hence, Lemma
\ref{lem.fps.lim.xn-equiv} shows that there exists some integer $N\in
\mathbb{N}$ such that
\[
\text{all integers }i\geq N\text{ satisfy }f_{i}\overset{x^{n}}{\equiv}f.
\]
Let us denote this $N$ by $K$. Hence, all integers $i\geq K$ satisfy
$f_{i}\overset{x^{n}}{\equiv}f$.

Thus, we have found an integer $K\in\mathbb{N}$ such that%
\begin{equation}
\text{all integers }i\geq K\text{ satisfy }f_{i}\overset{x^{n}}{\equiv}f.
\label{pf.prop.fps.lim.sum-quot.K}%
\end{equation}
Similarly, using $\lim\limits_{i\rightarrow\infty}g_{i}=g$, we can find an
integer $L\in\mathbb{N}$ such that
\begin{equation}
\text{all integers }i\geq L\text{ satisfy }g_{i}\overset{x^{n}}{\equiv}g.
\label{pf.prop.fps.lim.sum-quot.L}%
\end{equation}
Consider these $K$ and $L$. Let us furthermore set $P:=\max\left\{
K,L\right\}  $. Thus, $P\in\mathbb{N}$.

We shall now show that each integer $i\geq P$ satisfies $\left[  x^{n}\right]
\dfrac{f_{i}}{g_{i}}=\left[  x^{n}\right]  \dfrac{f}{g}$.

Indeed, let $i\geq P$ be an integer. Then, (\ref{pf.prop.fps.lim.sum-quot.K})
yields $f_{i}\overset{x^{n}}{\equiv}f$ (since $i\geq P=\max\left\{
K,L\right\}  \geq K$), whereas (\ref{pf.prop.fps.lim.sum-quot.L}) yields
$g_{i}\overset{x^{n}}{\equiv}g$ (since $i\geq P=\max\left\{  K,L\right\}  \geq
L$). Since both FPSs $g_{i}$ and $g$ are invertible (indeed, $g_{i}$ is
invertible by assumption, while $g$ is invertible by Claim 1), we thus obtain
\[
\dfrac{f_{i}}{g_{i}}\overset{x^{n}}{\equiv}\dfrac{f}{g}%
\]
(by Theorem \ref{thm.fps.xneq.props} \textbf{(e)}, applied to $a=f_{i}$,
$b=f$, $c=g_{i}$ and $d=g$). In other words,
\[
\text{each }m\in\left\{  0,1,\ldots,n\right\}  \text{ satisfies }\left[
x^{m}\right]  \dfrac{f_{i}}{g_{i}}=\left[  x^{m}\right]  \dfrac{f}{g}%
\]
(by the definition of $x^{n}$-equivalence). Applying this to $m=n$, we find%
\[
\left[  x^{n}\right]  \dfrac{f_{i}}{g_{i}}=\left[  x^{n}\right]  \dfrac{f}%
{g}.
\]

Forget that we fixed $i$. We thus have shown that all integers $i\geq P$
satisfy $\left[  x^{n}\right]  \dfrac{f_{i}}{g_{i}}=\left[  x^{n}\right]
\dfrac{f}{g}$. Hence, there exists some $N\in\mathbb{N}$ such that%
\[
\text{all integers }i\geq N\text{ satisfy }\left[  x^{n}\right]  \dfrac{f_{i}%
}{g_{i}}=\left[  x^{n}\right]  \dfrac{f}{g}%
\]
(namely, $N=P$). In other words,%
\[
\text{the sequence }\left(  \left[  x^{n}\right]  \dfrac{f_{i}}{g_{i}}\right)
_{i\in\mathbb{N}}\text{ stabilizes to }\left[  x^{n}\right]  \dfrac{f}{g}%
\]
(by the definition of \textquotedblleft stabilizes\textquotedblright).

Forget that we fixed $n$. We thus have shown that for each $n\in\mathbb{N}$,
the sequence $\left(  \left[  x^{n}\right]  \dfrac{f_{i}}{g_{i}}\right)
_{i\in\mathbb{N}}$ stabilizes to $\left[  x^{n}\right]  \dfrac{f}{g}$. In
other words, the sequence $\left(  \dfrac{f_{i}}{g_{i}}\right)  _{i\in
\mathbb{N}}$ coefficientwise stabilizes to $\dfrac{f}{g}$ (by the definition
of \textquotedblleft coefficientwise stabilizing\textquotedblright). In other
words, $\lim\limits_{i\rightarrow\infty}\dfrac{f_{i}}{g_{i}}=\dfrac{f}{g}$.
This proves Proposition \ref{prop.fps.lim.sum-quot} (since we already have
shown that $g$ is invertible).
\end{proof}
\end{fineprint}

\begin{fineprint}
\begin{proof}
[Proof of Theorem \ref{thm.fps.lim.sum-lim}.]The family $\left(  f_{n}\right)
_{n\in\mathbb{N}}$ is summable. In other words, the family $\left(
f_{k}\right)  _{k\in\mathbb{N}}$ is summable (since this is the same family as
$\left(  f_{n}\right)  _{n\in\mathbb{N}}$). In other words, for each
$n\in\mathbb{N}$,
\begin{equation}
\text{all but finitely many }k\in\mathbb{N}\text{ satisfy }\left[
x^{n}\right]  f_{k}=0 \label{pf.thm.fps.lim.sum-lim.1}%
\end{equation}
(by the definition of \textquotedblleft summable\textquotedblright).

Let us define%
\begin{equation}
g_{i}:=\sum_{k=0}^{i}f_{k}\ \ \ \ \ \ \ \ \ \ \text{for each }i\in\mathbb{N}.
\label{pf.thm.fps.lim.sum-lim.gi=}%
\end{equation}
Let us furthermore set%
\begin{equation}
g:=\sum_{k\in\mathbb{N}}f_{k}. \label{pf.thm.fps.lim.sum-lim.g=}%
\end{equation}

Now, fix $n\in\mathbb{N}$. Then, all but finitely many $k\in\mathbb{N}$
satisfy $\left[  x^{n}\right]  f_{k}=0$ (by (\ref{pf.thm.fps.lim.sum-lim.1})).
In other words, there exists a finite subset $J$ of $\mathbb{N}$ such that
\begin{equation}
\text{all }k\in\mathbb{N}\setminus J\text{ satisfy }\left[  x^{n}\right]
f_{k}=0. \label{pf.thm.fps.lim.sum-lim.2}%
\end{equation}
Consider this subset $J$. The set $J$ is a finite set of nonnegative integers,
and thus has an upper bound (since any finite set of nonnegative integers has
an upper bound). In other words, there exists some $m\in\mathbb{N}$ such that%
\begin{equation}
\text{all }k\in J\text{ satisfy }k\leq m. \label{pf.thm.fps.lim.sum-lim.3}%
\end{equation}
Consider this $m$.

Let $k\in\mathbb{N}$ be such that $k\geq m+1$. If we had $k\in J$, then we
would have $k\leq m$ (by (\ref{pf.thm.fps.lim.sum-lim.3})), which would
contradict $k\geq m+1>m$. Thus, we cannot have $k\in J$. Hence, $k\in
\mathbb{N}\setminus J$ (since $k\in\mathbb{N}$ but not $k\in J$). Therefore,
(\ref{pf.thm.fps.lim.sum-lim.2}) yields $\left[  x^{n}\right]  f_{k}=0$.

Forget that we fixed $k$. We thus have shown that%
\begin{equation}
\left[  x^{n}\right]  f_{k}=0\ \ \ \ \ \ \ \ \ \ \text{for each }%
k\in\mathbb{N}\text{ satisfying }k\geq m+1. \label{pf.thm.fps.lim.sum-lim.4}%
\end{equation}

Now, let $i$ be an integer such that $i\geq m$. Then, $i\in\mathbb{N}$ (since
$i\geq m\geq0$) and $i+1\geq m+1$ (since $i\geq m$). From $g=\sum
_{k\in\mathbb{N}}f_{k}$, we obtain%
\begin{align}
\left[  x^{n}\right]  g  &  =\left[  x^{n}\right]  \left(  \sum_{k\in
\mathbb{N}}f_{k}\right)  =\sum_{k\in\mathbb{N}}\left[  x^{n}\right]
f_{k}\ \ \ \ \ \ \ \ \ \ \left(  \text{by (\ref{eq.def.fps.summable.sum}%
)}\right) \nonumber\\
&  =\sum_{k=0}^{i}\left[  x^{n}\right]  f_{k}+\sum_{k=i+1}^{\infty
}\underbrace{\left[  x^{n}\right]  f_{k}}_{\substack{=0\\\text{(by
(\ref{pf.thm.fps.lim.sum-lim.4})}\\\text{(since }k\geq i+1\geq m+1\text{))}%
}}=\sum_{k=0}^{i}\left[  x^{n}\right]  f_{k}+\underbrace{\sum_{k=i+1}^{\infty
}0}_{=0}\nonumber\\
&  =\sum_{k=0}^{i}\left[  x^{n}\right]  f_{k}.
\label{pf.thm.fps.lim.sum-lim.5}%
\end{align}
On the other hand, from $g_{i}=\sum_{k=0}^{i}f_{k}$, we obtain%
\[
\left[  x^{n}\right]  g_{i}=\left[  x^{n}\right]  \left(  \sum_{k=0}^{i}%
f_{k}\right)  =\sum_{k=0}^{i}\left[  x^{n}\right]  f_{k}.
\]
Comparing this with (\ref{pf.thm.fps.lim.sum-lim.5}), we obtain $\left[
x^{n}\right]  g_{i}=\left[  x^{n}\right]  g$.

Forget that we fixed $i$. We thus have shown that
\[
\text{all integers }i\geq m\text{ satisfy }\left[  x^{n}\right]  g_{i}=\left[
x^{n}\right]  g.
\]
Hence, there exists some $N\in\mathbb{N}$ such that
\[
\text{all integers }i\geq N\text{ satisfy }\left[  x^{n}\right]  g_{i}=\left[
x^{n}\right]  g
\]
(namely, $N=m$). In other words, the sequence $\left(  \left[  x^{n}\right]
g_{i}\right)  _{i\in\mathbb{N}}$ stabilizes to $\left[  x^{n}\right]  g$ (by
the definition of \textquotedblleft stabilizes\textquotedblright).

Forget that we fixed $n$. We thus have shown that for each $n\in\mathbb{N}$,%
\[
\text{the sequence }\left(  \left[  x^{n}\right]  g_{i}\right)  _{i\in
\mathbb{N}}\text{ stabilizes to }\left[  x^{n}\right]  g.
\]
In other words, the sequence $\left(  g_{i}\right)  _{i\in\mathbb{N}}$
coefficientwise stabilizes to $g$. In other words, $\lim\limits_{i\rightarrow
\infty}g_{i}=g$. In view of (\ref{pf.thm.fps.lim.sum-lim.gi=}) and
(\ref{pf.thm.fps.lim.sum-lim.g=}), we can rewrite this as%
\[
\lim\limits_{i\rightarrow\infty}\sum_{k=0}^{i}f_{k}=\sum_{k\in\mathbb{N}}%
f_{k}.
\]
Renaming the summation index $k$ as $n$ everywhere in this equality, we can
rewrite it as%
\[
\lim\limits_{i\rightarrow\infty}\sum_{n=0}^{i}f_{n}=\sum_{n\in\mathbb{N}}%
f_{n}.
\]
This proves Theorem \ref{thm.fps.lim.sum-lim}.
\end{proof}

\begin{proof}
[Proof of Theorem \ref{thm.fps.lim.prod-lim}.]The family $\left(
f_{n}\right)  _{n\in\mathbb{N}}$ is multipliable. In other words, the family
$\left(  f_{k}\right)  _{k\in\mathbb{N}}$ is multipliable (since this is the
same family as $\left(  f_{n}\right)  _{n\in\mathbb{N}}$). In other words,
each coefficient in the product of this family is finitely determined (by the
definition of \textquotedblleft multipliable\textquotedblright).

Let us define%
\begin{equation}
g_{i}:=\prod_{k=0}^{i}f_{k}\ \ \ \ \ \ \ \ \ \ \text{for each }i\in\mathbb{N}.
\label{pf.thm.fps.lim.prod-lim.gi=}%
\end{equation}
Let us furthermore set%
\begin{equation}
g:=\prod_{k\in\mathbb{N}}f_{k}. \label{pf.thm.fps.lim.prod-lim.g=}%
\end{equation}

Now, fix $n\in\mathbb{N}$. Then, the $x^{n}$-coefficient in the product of the
family $\left(  f_{k}\right)  _{k\in\mathbb{N}}$ is finitely determined (since
each coefficient in the product of this family is finitely determined). In
other words, there is a finite subset $M$ of $\mathbb{N}$ that determines the
$x^{n}$-coefficient in the product of $\left(  f_{k}\right)  _{k\in\mathbb{N}%
}$ (by the definition of \textquotedblleft finitely
determined\textquotedblright). Consider this subset $M$.

The set $M$ is a finite set of nonnegative integers, and thus has an upper
bound (since any finite set of nonnegative integers has an upper bound). In
other words, there exists some $m\in\mathbb{N}$ such that%
\begin{equation}
\text{all }k\in M\text{ satisfy }k\leq m. \label{pf.thm.fps.lim.prod-lim.3}%
\end{equation}
Consider this $m$.

Now, let $i$ be an integer such that $i\geq m$. Then, $i\in\mathbb{N}$ (since
$i\geq m\geq0$) and $m\leq i$. From (\ref{pf.thm.fps.lim.prod-lim.3}), we see
that all $k\in M$ satisfy $k\leq m\leq i$ and therefore $k\in\left\{
0,1,\ldots,i\right\}  $. In other words, $M\subseteq\left\{  0,1,\ldots
,i\right\}  $.

However, the set $M$ determines the $x^{n}$-coefficient in the product of
$\left(  f_{k}\right)  _{k\in\mathbb{N}}$. In other words, every finite subset
$J$ of $\mathbb{N}$ satisfying $M\subseteq J\subseteq\mathbb{N}$ satisfies%
\[
\left[  x^{n}\right]  \left(  \prod_{k\in J}f_{k}\right)  =\left[
x^{n}\right]  \left(  \prod_{k\in M}f_{k}\right)
\]
(by the definition of \textquotedblleft determines the $x^{n}$-coefficient in
the product of $\left(  f_{k}\right)  _{k\in\mathbb{N}}$\textquotedblright).
Applying this to $J=\left\{  0,1,\ldots,i\right\}  $, we obtain%
\begin{equation}
\left[  x^{n}\right]  \left(  \prod_{k\in\left\{  0,1,\ldots,i\right\}  }%
f_{k}\right)  =\left[  x^{n}\right]  \left(  \prod_{k\in M}f_{k}\right)
\label{pf.thm.fps.lim.prod-lim.4}%
\end{equation}
(since $\left\{  0,1,\ldots,i\right\}  $ is a finite subset of $\mathbb{N}$
satisfying $M\subseteq\left\{  0,1,\ldots,i\right\}  \subseteq\mathbb{N}$).

Moreover, the definition of the product of the multipliable family $\left(
f_{k}\right)  _{k\in\mathbb{N}}$ (Definition \ref{def.fps.multipliable}
\textbf{(b)}) shows that%
\[
\left[  x^{n}\right]  \left(  \prod_{k\in\mathbb{N}}f_{k}\right)  =\left[
x^{n}\right]  \left(  \prod_{k\in M}f_{k}\right)
\]
(since $M$ is a finite subset of $\mathbb{N}$ that determines the $x^{n}%
$-coefficient in the product of $\left(  f_{k}\right)  _{k\in\mathbb{N}}$).
Comparing this with (\ref{pf.thm.fps.lim.prod-lim.4}), we find%
\[
\left[  x^{n}\right]  \left(  \prod_{k\in\left\{  0,1,\ldots,i\right\}  }%
f_{k}\right)  =\left[  x^{n}\right]  \left(  \prod_{k\in\mathbb{N}}%
f_{k}\right)  .
\]

In view of $\prod_{k\in\left\{  0,1,\ldots,i\right\}  }f_{k}=\prod_{k=0}%
^{i}f_{k}=g_{i}$ (by (\ref{pf.thm.fps.lim.prod-lim.gi=})) and $\prod
_{k\in\mathbb{N}}f_{k}=g$ (by (\ref{pf.thm.fps.lim.prod-lim.g=})), we can
rewrite this as%
\[
\left[  x^{n}\right]  g_{i}=\left[  x^{n}\right]  g.
\]

Forget that we fixed $i$. We thus have shown that
\[
\text{all integers }i\geq m\text{ satisfy }\left[  x^{n}\right]  g_{i}=\left[
x^{n}\right]  g.
\]
Hence, there exists some $N\in\mathbb{N}$ such that
\[
\text{all integers }i\geq N\text{ satisfy }\left[  x^{n}\right]  g_{i}=\left[
x^{n}\right]  g
\]
(namely, $N=m$). In other words, the sequence $\left(  \left[  x^{n}\right]
g_{i}\right)  _{i\in\mathbb{N}}$ stabilizes to $\left[  x^{n}\right]  g$ (by
the definition of \textquotedblleft stabilizes\textquotedblright).

Forget that we fixed $n$. We thus have shown that for each $n\in\mathbb{N}$,%
\[
\text{the sequence }\left(  \left[  x^{n}\right]  g_{i}\right)  _{i\in
\mathbb{N}}\text{ stabilizes to }\left[  x^{n}\right]  g.
\]
In other words, the sequence $\left(  g_{i}\right)  _{i\in\mathbb{N}}$
coefficientwise stabilizes to $g$. In other words, $\lim\limits_{i\rightarrow
\infty}g_{i}=g$. In view of (\ref{pf.thm.fps.lim.prod-lim.gi=}) and
(\ref{pf.thm.fps.lim.prod-lim.g=}), we can rewrite this as%
\[
\lim\limits_{i\rightarrow\infty}\prod_{k=0}^{i}f_{k}=\prod_{k\in\mathbb{N}%
}f_{k}.
\]
Renaming the product index $k$ as $n$ everywhere in this equality, we can
rewrite it as%
\[
\lim\limits_{i\rightarrow\infty}\prod_{n=0}^{i}f_{n}=\prod_{n\in\mathbb{N}%
}f_{n}.
\]
This proves Theorem \ref{thm.fps.lim.prod-lim}.
\end{proof}

\begin{proof}
[Proof of Corollary \ref{cor.fps.lim.fps-as-pol}.]Let a FPS $a\in K\left[
\left[  x\right]  \right]  $ be written in the form $a=\sum_{n\in\mathbb{N}%
}a_{n}x^{n}$ with $a_{n}\in K$. Then, Corollary \ref{cor.fps.sumakxk} shows
that the family $\left(  a_{n}x^{n}\right)  _{n\in\mathbb{N}}$ is summable.
Hence, Theorem \ref{thm.fps.lim.sum-lim} (applied to $f_{n}=a_{n}x^{n}$)
yields
\[
\lim\limits_{i\rightarrow\infty}\sum_{n=0}^{i}a_{n}x^{n}=\sum_{n\in\mathbb{N}%
}a_{n}x^{n}=a.
\]
In other words, $a=\lim\limits_{i\rightarrow\infty}\sum_{n=0}^{i}a_{n}x^{n}$.
This proves Corollary \ref{cor.fps.lim.fps-as-pol}.
\end{proof}
\end{fineprint}

\begin{fineprint}
\begin{proof}
[Proof of Theorem \ref{thm.fps.lim.sum-lim-conv}.]Let us define%
\begin{equation}
g_{i}:=\sum_{k=0}^{i}f_{k}\ \ \ \ \ \ \ \ \ \ \text{for each }i\in\mathbb{N}.
\label{pf.thm.fps.lim.sum-lim-conv.gi=}%
\end{equation}
We have assumed that the limit $\lim\limits_{i\rightarrow\infty}\sum_{n=0}%
^{i}f_{n}$ exists. Let us denote this limit by $g$. Thus,%
\begin{align}
g  &  =\lim\limits_{i\rightarrow\infty}\sum_{n=0}^{i}f_{n}=\lim
\limits_{i\rightarrow\infty}\underbrace{\sum_{k=0}^{i}f_{k}}_{\substack{=g_{i}%
\\\text{(by (\ref{pf.thm.fps.lim.sum-lim-conv.gi=}))}}%
}\ \ \ \ \ \ \ \ \ \ \left(
\begin{array}
[c]{c}%
\text{here, we have renamed the}\\
\text{summation index }n\text{ as }k
\end{array}
\right) \nonumber\\
&  =\lim\limits_{i\rightarrow\infty}g_{i}.
\label{pf.thm.fps.lim.sum-lim-conv.g=lim}%
\end{align}
In other words, the sequence $\left(  g_{i}\right)  _{i\in\mathbb{N}}$
coefficientwise stabilizes to $g$. In other words, for each $n\in\mathbb{N}$,%
\begin{equation}
\text{the sequence }\left(  \left[  x^{n}\right]  g_{i}\right)  _{i\in
\mathbb{N}}\text{ stabilizes to }\left[  x^{n}\right]  g.
\label{pf.thm.fps.lim.sum-lim-conv.gstab}%
\end{equation}

Let $n\in\mathbb{N}$. Then, the sequence $\left(  \left[  x^{n}\right]
g_{i}\right)  _{i\in\mathbb{N}}$ stabilizes to $\left[  x^{n}\right]  g$ (by
(\ref{pf.thm.fps.lim.sum-lim-conv.gstab})). In other words, there exists some
$N\in\mathbb{N}$ such that
\begin{equation}
\text{all integers }i\geq N\text{ satisfy }\left[  x^{n}\right]  g_{i}=\left[
x^{n}\right]  g. \label{pf.thm.fps.lim.sum-lim-conv.1}%
\end{equation}
Consider this $N$. Let $M$ be the set $\left\{  0,1,\ldots,N\right\}  $; this
is a finite subset of $\mathbb{N}$.

Let $i\in\mathbb{N}\setminus M$. Thus,
\begin{align*}
i  &  \in\mathbb{N}\setminus M=\mathbb{N}\setminus\left\{  0,1,\ldots
,N\right\}  \ \ \ \ \ \ \ \ \ \ \left(  \text{since }M=\left\{  0,1,\ldots
,N\right\}  \right) \\
&  =\left\{  N+1,\ N+2,\ N+3,\ \ldots\right\}  .
\end{align*}
In other words, $i$ is a nonnegative integer with $i\geq N+1$. Hence, $i-1\geq
N$. Therefore, (\ref{pf.thm.fps.lim.sum-lim-conv.1}) (applied to $i-1$ instead
of $i$) yields $\left[  x^{n}\right]  g_{i-1}=\left[  x^{n}\right]  g$. But
(\ref{pf.thm.fps.lim.sum-lim-conv.1}) also yields $\left[  x^{n}\right]
g_{i}=\left[  x^{n}\right]  g$ (since $i\geq N+1\geq N$). Comparing these two
equalities, we find
\begin{equation}
\left[  x^{n}\right]  g_{i-1}=\left[  x^{n}\right]  g_{i}.
\label{pf.thm.fps.lim.sum-lim-conv.2}%
\end{equation}

However, (\ref{pf.thm.fps.lim.sum-lim-conv.gi=}) yields
\begin{equation}
g_{i}=\sum_{k=0}^{i}f_{k}=f_{i}+\sum_{k=0}^{i-1}f_{k}.
\label{pf.thm.fps.lim.sum-lim-conv.3}%
\end{equation}
Moreover, (\ref{pf.thm.fps.lim.sum-lim-conv.gi=}) (applied to $i-1$ instead of
$i$) yields $g_{i-1}=\sum_{k=0}^{i-1}f_{k}$. In view of this, we can rewrite
(\ref{pf.thm.fps.lim.sum-lim-conv.3}) as%
\[
g_{i}=f_{i}+g_{i-1}.
\]
Hence, $\left[  x^{n}\right]  g_{i}=\left[  x^{n}\right]  \left(
f_{i}+g_{i-1}\right)  =\left[  x^{n}\right]  f_{i}+\underbrace{\left[
x^{n}\right]  g_{i-1}}_{\substack{=\left[  x^{n}\right]  g_{i}\\\text{(by
(\ref{pf.thm.fps.lim.sum-lim-conv.2}))}}}=\left[  x^{n}\right]  f_{i}+\left[
x^{n}\right]  g_{i}$. Subtracting $\left[  x^{n}\right]  g_{i}$ from both
sides of this equality, we find $0=\left[  x^{n}\right]  f_{i}$. In other
words, $\left[  x^{n}\right]  f_{i}=0$.

Forget that we fixed $i$. We thus have shown that all $i\in\mathbb{N}\setminus
M$ satisfy $\left[  x^{n}\right]  f_{i}=0$. Hence, all but finitely many
$i\in\mathbb{N}$ satisfy $\left[  x^{n}\right]  f_{i}=0$ (since all but
finitely many $i\in\mathbb{N}$ belong to $\mathbb{N}\setminus M$ (because the
set $M$ is finite)). Renaming the index $i$ as $k$ in this statement, we
obtain the following: All but finitely many $k\in\mathbb{N}$ satisfy $\left[
x^{n}\right]  f_{k}=0$.

Forget that we fixed $n$. We thus have shown that for each $n\in\mathbb{N}$,
\[
\text{all but finitely many }k\in\mathbb{N}\text{ satisfy }\left[
x^{n}\right]  f_{k}=0.
\]
In other words, the family $\left(  f_{n}\right)  _{n\in\mathbb{N}}$ is
summable (by the definition of \textquotedblleft summable\textquotedblright).

Hence, Theorem \ref{thm.fps.lim.sum-lim} yields $\lim\limits_{i\rightarrow
\infty}\sum_{n=0}^{i}f_{n}=\sum_{n\in\mathbb{N}}f_{n}$. In other words,
$\sum_{n\in\mathbb{N}}f_{n}=\lim\limits_{i\rightarrow\infty}\sum_{n=0}%
^{i}f_{n}$. This completes the proof of Theorem \ref{thm.fps.lim.sum-lim-conv}.
\end{proof}

\begin{proof}
[Proof of Theorem \ref{thm.fps.lim.prod-lim-conv}.]Let us define%
\begin{equation}
g_{i}:=\prod_{k=0}^{i}f_{k}\ \ \ \ \ \ \ \ \ \ \text{for each }i\in\mathbb{N}.
\label{pf.thm.fps.lim.prod-lim-conv.gi=}%
\end{equation}

We have assumed that the limit $\lim\limits_{i\rightarrow\infty}\prod
_{n=0}^{i}f_{n}$ exists. Let us denote this limit by $g$. Thus,%
\begin{align}
g  &  =\lim\limits_{i\rightarrow\infty}\prod_{n=0}^{i}f_{n}=\lim
\limits_{i\rightarrow\infty}\underbrace{\prod_{k=0}^{i}f_{k}}%
_{\substack{=g_{i}\\\text{(by (\ref{pf.thm.fps.lim.prod-lim-conv.gi=}))}%
}}\ \ \ \ \ \ \ \ \ \ \left(
\begin{array}
[c]{c}%
\text{here, we have renamed the}\\
\text{product index }n\text{ as }k
\end{array}
\right) \nonumber\\
&  =\lim\limits_{i\rightarrow\infty}g_{i}.
\label{pf.thm.fps.lim.prod-lim-conv.g=lim}%
\end{align}

Let $n\in\mathbb{N}$ be arbitrary. Lemma \ref{lem.fps.lim.xn-equiv} (applied
to $g_{i}$ and $g$ instead of $f_{i}$ and $f$) shows that there exists some
integer $N\in\mathbb{N}$ such that
\begin{equation}
\text{all integers }i\geq N\text{ satisfy }g_{i}\overset{x^{n}}{\equiv}g
\label{pf.thm.fps.lim.prod-lim-conv.gi-equiv-g}%
\end{equation}
(since $g=\lim\limits_{i\rightarrow\infty}g_{i}$). Consider this $N$.

Let $M=\left\{  0,1,\ldots,N\right\}  $. Thus, $M$ is a finite subset of
$\mathbb{N}$. We shall show that this subset $M$ determines the $x^{n}%
$-coefficient in the product of $\left(  f_{k}\right)  _{k\in\mathbb{N}}$.

For this purpose, we first observe that $M=\left\{  0,1,\ldots,N\right\}  $,
so that%
\begin{equation}
\prod_{k\in M}f_{k}=\prod_{k\in\left\{  0,1,\ldots,N\right\}  }f_{k}%
=\prod_{k=0}^{N}f_{k}=g_{N} \label{pf.thm.fps.lim.prod-lim-conv.=gN}%
\end{equation}
(since $g_{N}$ was defined to be $\prod_{k=0}^{N}f_{k}$). Applying
(\ref{pf.thm.fps.lim.prod-lim-conv.gi-equiv-g}) to $i=N$, we obtain
$g_{N}\overset{x^{n}}{\equiv}g$ (since $N\geq N$). Therefore, $g\overset{x^{n}%
}{\equiv}g_{N}$ (since $\overset{x^{n}}{\equiv}$ is an equivalence relation
and thus symmetric).

Next, we shall show two claims:

\begin{statement}
\textit{Claim 1:} For each integer $j>N$, we have $g\overset{x^{n}}{\equiv
}gf_{j}$.
\end{statement}

\begin{proof}
[Proof of Claim 1.]Let $j>N$ be an integer. Thus, $j\geq N+1$, so that
$j-1\geq N$. Hence, (\ref{pf.thm.fps.lim.prod-lim-conv.gi-equiv-g}) (applied
to $i=j-1$) yields $g_{j-1}\overset{x^{n}}{\equiv}g$. However,
(\ref{pf.thm.fps.lim.prod-lim-conv.gi-equiv-g}) (applied to $i=j$) yields
$g_{j}\overset{x^{n}}{\equiv}g$ (since $j\geq N+1\geq N$). Thus,
$g\overset{x^{n}}{\equiv}g_{j}$ (since $\overset{x^{n}}{\equiv}$ is an
equivalence relation and thus symmetric).

Also, we obviously have $f_{j}\overset{x^{n}}{\equiv}f_{j}$ (since
$\overset{x^{n}}{\equiv}$ is an equivalence relation and thus reflexive).

However, the definition of $g_{j-1}$ yields $g_{j-1}=\prod_{k=0}^{j-1}f_{k}$.
Furthermore, the definition of $g_{j}$ yields%
\[
g_{j}=\prod_{k=0}^{j}f_{k}=\underbrace{\left(  \prod_{k=0}^{j-1}f_{k}\right)
}_{=g_{j-1}}f_{j}=g_{j-1}f_{j}.
\]
But recall that $g_{j-1}\overset{x^{n}}{\equiv}g$ and $f_{j}\overset{x^{n}%
}{\equiv}f_{j}$. Hence, $g_{j-1}f_{j}\overset{x^{n}}{\equiv}gf_{j}$ (by
(\ref{eq.thm.fps.xneq.props.b.*}), applied to $a=g_{j-1}$, $b=g$, $c=f_{j}$
and $d=f_{j}$). In view of $g_{j}=g_{j-1}f_{j}$, we can rewrite this as
$g_{j}\overset{x^{n}}{\equiv}gf_{j}$.

Now, recall that $g\overset{x^{n}}{\equiv}g_{j}$. Thus, $g\overset{x^{n}%
}{\equiv}g_{j}\overset{x^{n}}{\equiv}gf_{j}$, so that $g\overset{x^{n}%
}{\equiv}gf_{j}$ (since $\overset{x^{n}}{\equiv}$ is an equivalence relation
and thus transitive). This proves Claim 1.
\end{proof}

\begin{statement}
\textit{Claim 2:} If $U$ is any finite subset of $\mathbb{N}$ satisfying
$M\subseteq U$, then $g\overset{x^{n}}{\equiv}\prod_{k\in U}f_{k}$.
\end{statement}

\begin{proof}
[Proof of Claim 2.]We induct on the nonnegative integer $\left\vert U\setminus
M\right\vert $:

\textit{Base case:} Let us check that Claim 2 holds when $\left\vert
U\setminus M\right\vert =0$.

Indeed, let $U$ be any finite subset of $\mathbb{N}$ satisfying $M\subseteq U$
and $\left\vert U\setminus M\right\vert =0$. Then, $U\setminus M=\varnothing$
(since $\left\vert U\setminus M\right\vert =0$), so that $U\subseteq M$.
Combined with $M\subseteq U$, this yields $U=M$. Thus, $\prod_{k\in U}%
f_{k}=\prod_{k\in M}f_{k}=g_{N}$ (by (\ref{pf.thm.fps.lim.prod-lim-conv.=gN}%
)). But we have already proved that $g\overset{x^{n}}{\equiv}g_{N}$. In other
words, $g\overset{x^{n}}{\equiv}\prod_{k\in U}f_{k}$ (since $\prod_{k\in
U}f_{k}=g_{N}$). Thus, Claim 2 is proved under the assumption that $\left\vert
U\setminus M\right\vert =0$. This completes the base case.

\textit{Induction step:} Let $s\in\mathbb{N}$. Assume (as the induction
hypothesis) that Claim 2 holds when $\left\vert U\setminus M\right\vert =s$.
We must now prove that Claim 2 holds when $\left\vert U\setminus M\right\vert
=s+1$.

Indeed, let $U$ be any finite subset of $\mathbb{N}$ satisfying $M\subseteq U$
and $\left\vert U\setminus M\right\vert =s+1$. Then, the set $U\setminus M$ is
nonempty (since $\left\vert U\setminus M\right\vert =s+1>s\geq0$). Hence,
there exists some $u\in U\setminus M$. Consider this $u$. From $u\in
U\setminus M$, we obtain $\left\vert \left(  U\setminus M\right)
\setminus\left\{  u\right\}  \right\vert =\left\vert U\setminus M\right\vert
-1=s$ (since $\left\vert U\setminus M\right\vert =s+1$). In view of $\left(
U\setminus M\right)  \setminus\left\{  u\right\}  =\left(  U\setminus\left\{
u\right\}  \right)  \setminus M$, we can rewrite this as $\left\vert \left(
U\setminus\left\{  u\right\}  \right)  \setminus M\right\vert =s$. Moreover,
$U\setminus\left\{  u\right\}  $ is a finite subset of $\mathbb{N}$ (since $U$
is a finite subset of $\mathbb{N}$), and satisfies $M\subseteq U\setminus
\left\{  u\right\}  $ (since $u\in U\setminus M$ and thus $u\notin M$ and
therefore $M\setminus\left\{  u\right\}  =M$, so that $M=\underbrace{M}%
_{\subseteq U}\setminus\left\{  u\right\}  \subseteq U\setminus\left\{
u\right\}  $). Hence, by the induction hypothesis, we can apply Claim 2 to
$U\setminus\left\{  u\right\}  $ instead of $U$ (since $\left\vert \left(
U\setminus\left\{  u\right\}  \right)  \setminus M\right\vert =s$). As a
result, we obtain $g\overset{x^{n}}{\equiv}\prod_{k\in U\setminus\left\{
u\right\}  }f_{k}$.

However, $u\in U\setminus M\subseteq U$, so that $\prod_{k\in U}f_{k}=\left(
\prod_{k\in U\setminus\left\{  u\right\}  }f_{k}\right)  f_{u}$ (here, we have
split off the factor for $k=u$ from the product).

We have $g\overset{x^{n}}{\equiv}\prod_{k\in U\setminus\left\{  u\right\}
}f_{k}$ (as we have proved above) and $f_{u}\overset{x^{n}}{\equiv}f_{u}$
(since $\overset{x^{n}}{\equiv}$ is an equivalence relation and thus
reflexive). Hence, $gf_{u}\overset{x^{n}}{\equiv}\left(  \prod_{k\in
U\setminus\left\{  u\right\}  }f_{k}\right)  f_{u}$ (by
(\ref{eq.thm.fps.xneq.props.b.*}), applied to $a=g$, $b=\prod_{k\in
U\setminus\left\{  u\right\}  }f_{k}$, $c=f_{u}$ and $d=f_{u}$). In view of
$\prod_{k\in U}f_{k}=\left(  \prod_{k\in U\setminus\left\{  u\right\}  }%
f_{k}\right)  f_{u}$, we can rewrite this as $gf_{u}\overset{x^{n}}{\equiv
}\prod_{k\in U}f_{k}$.

But $u\in\underbrace{U}_{\subseteq\mathbb{N}}\setminus\underbrace{M}%
_{=\left\{  0,1,\ldots,N\right\}  }\subseteq\mathbb{N}\setminus\left\{
0,1,\ldots,N\right\}  =\left\{  N+1,\ N+2,\ N+3,\ \ldots\right\}  $, so that
$u\geq N+1>N$. Hence, Claim 1 (applied to $j=u$) yields $g\overset{x^{n}%
}{\equiv}gf_{u}$. Therefore, $g\overset{x^{n}}{\equiv}gf_{u}\overset{x^{n}%
}{\equiv}\prod_{k\in U}f_{k}$. Therefore, $g\overset{x^{n}}{\equiv}\prod_{k\in
U}f_{k}$ (since $\overset{x^{n}}{\equiv}$ is an equivalence relation and thus transitive).

Forget that we fixed $U$. We thus have shown that if $U$ is any finite subset
of $\mathbb{N}$ satisfying $M\subseteq U$ and $\left\vert U\setminus
M\right\vert =s+1$, then $g\overset{x^{n}}{\equiv}\prod_{k\in U}f_{k}$. In
other words, Claim 2 holds when $\left\vert U\setminus M\right\vert =s+1$.
This completes the induction step. Thus, Claim 2 is proved by induction.
\end{proof}

Now, let $J$ be a finite subset of $\mathbb{N}$ satisfying $M\subseteq
J\subseteq\mathbb{N}$. Then, $g\overset{x^{n}}{\equiv}\prod_{k\in J}f_{k}$ (by
Claim 2, applied to $U=J$). In other words,%
\[
\text{each }m\in\left\{  0,1,\ldots,n\right\}  \text{ satisfies }\left[
x^{m}\right]  g=\left[  x^{m}\right]  \left(  \prod_{k\in J}f_{k}\right)
\]
(by the definition of $x^{n}$-equivalence). Applying this to $m=n$, we find%
\[
\left[  x^{n}\right]  g=\left[  x^{n}\right]  \left(  \prod_{k\in J}%
f_{k}\right)  .
\]
The same argument can be applied to $M$ instead of $J$ (since $M\subseteq
M\subseteq\mathbb{N}$), and thus yields%
\[
\left[  x^{n}\right]  g=\left[  x^{n}\right]  \left(  \prod_{k\in M}%
f_{k}\right)  .
\]
Comparing these two equalities, we find
\[
\left[  x^{n}\right]  \left(  \prod_{k\in J}f_{k}\right)  =\left[
x^{n}\right]  \left(  \prod_{k\in M}f_{k}\right)  .
\]

Forget that we fixed $J$. We have now shown that every finite subset $J$ of
$\mathbb{N}$ satisfying $M\subseteq J\subseteq\mathbb{N}$ satisfies%
\[
\left[  x^{n}\right]  \left(  \prod_{k\in J}f_{k}\right)  =\left[
x^{n}\right]  \left(  \prod_{k\in M}f_{k}\right)  .
\]
In other words, the set $M$ determines the $x^{n}$-coefficient in the product
of $\left(  f_{k}\right)  _{k\in\mathbb{N}}$ (by the definition of
\textquotedblleft determines the $x^{n}$-coefficient in the product of
$\left(  f_{k}\right)  _{k\in\mathbb{N}}$\textquotedblright). Hence, the
$x^{n}$-coefficient in the product of $\left(  f_{k}\right)  _{k\in\mathbb{N}%
}$ is finitely determined (by the definition of \textquotedblleft finitely
determined\textquotedblright, since $M$ is a finite subset of $\mathbb{N}$).

Forget that we fixed $n$. We thus have shown that for each $n\in\mathbb{N}$,
the $x^{n}$-coefficient in the product of $\left(  f_{k}\right)
_{k\in\mathbb{N}}$ is finitely determined. In other words, each coefficient in
the product of $\left(  f_{k}\right)  _{k\in\mathbb{N}}$ is finitely
determined. In other words, the family $\left(  f_{k}\right)  _{k\in
\mathbb{N}}$ is multipliable (by the definition of \textquotedblleft
multipliable\textquotedblright). In other words, the family $\left(
f_{n}\right)  _{n\in\mathbb{N}}$ is multipliable (since this is the same
family as $\left(  f_{k}\right)  _{k\in\mathbb{N}}$).

Hence, Theorem \ref{thm.fps.lim.prod-lim} yields $\lim\limits_{i\rightarrow
\infty}\prod_{n=0}^{i}f_{n}=\prod_{n\in\mathbb{N}}f_{n}$. In other words,
$\prod_{n\in\mathbb{N}}f_{n}=\lim\limits_{i\rightarrow\infty}\prod_{n=0}%
^{i}f_{n}$. This completes the proof of Theorem
\ref{thm.fps.lim.prod-lim-conv}.
\end{proof}
\end{fineprint}

\subsection{\label{sec.details.gf.laure}Laurent power series}

\begin{fineprint}
Next, we shall prove some claims made in Section \ref{sec.gf.laure}. We shall
not go into too many details here, since the proofs are mostly analogous to
the corresponding proofs in the theory of \textquotedblleft
usual\textquotedblright\ FPSs (with nonnegative exponents), which we have
already seen. However, every once in a while, some minor change is required;
we will mainly focus on these changes.

\begin{proof}
[Proof of Theorem \ref{thm.fps.laure.laupol-ring} (sketched).]First, we need
to prove that the set $K\left[  x^{\pm}\right]  $ is closed under addition,
under scaling, and under the multiplication introduced in Definition
\ref{def.fps.laure.laupol}. The proof of this is analogous to Theorem
\ref{thm.fps.pol.ring} (with the obvious changes, such as replacing the
$\sum_{i=0}^{n}$ sums by $\sum_{i\in\mathbb{Z}}$ sums).

Next, we need to prove that the multiplication on $K\left[  x^{\pm}\right]  $
is associative, commutative, distributive and $K$-bilinear. This is analogous
to the corresponding parts of Theorem \ref{thm.fps.ring}. Likewise, we can
show that the element $\left(  \delta_{i,0}\right)  _{i\in\mathbb{Z}}$ (which
is our analogue of the FPS $\underline{1}$) is a neutral element for this
multiplication. Hence, $K\left[  x^{\pm}\right]  $ is a commutative
$K$-algebra with unity $\left(  \delta_{i,0}\right)  _{i\in\mathbb{Z}}$.

It remains to show that the element $x$ is invertible in this $K$-algebra. But
this is easy: Set $\overline{x}:=\left(  \delta_{i,-1}\right)  _{i\in
\mathbb{Z}}$, and show (by direct calculation) that $x\overline{x}%
=\overline{x}x=1$. (For example, let us prove that $x\overline{x}=1$. Indeed,
from $x=\left(  \delta_{i,1}\right)  _{i\in\mathbb{Z}}$ and $\overline
{x}=\left(  \delta_{i,-1}\right)  _{i\in\mathbb{Z}}$, we obtain%
\[
x\overline{x}=\left(  \delta_{i,1}\right)  _{i\in\mathbb{Z}}\cdot\left(
\delta_{i,-1}\right)  _{i\in\mathbb{Z}}=\left(  c_{n}\right)  _{n\in
\mathbb{Z}},\ \ \ \ \ \ \ \ \ \ \text{where}\ \ \ \ \ \ \ \ \ \ c_{n}%
=\sum_{i\in\mathbb{Z}}\delta_{i,1}\delta_{n-i,-1}%
\]
(by Definition \ref{def.fps.laure.laupol}). However, for each $n\in\mathbb{Z}%
$, we have%
\begin{align*}
c_{n}  &  =\sum_{i\in\mathbb{Z}}\delta_{i,1}\delta_{n-i,-1}=\underbrace{\delta
_{1,1}}_{=1}\delta_{n-1,-1}+\sum_{\substack{i\in\mathbb{Z};\\i\neq
1}}\underbrace{\delta_{i,1}}_{\substack{=0\\\text{(since }i\neq1\text{)}%
}}\delta_{n-i,-1}\\
&  \ \ \ \ \ \ \ \ \ \ \ \ \ \ \ \ \ \ \ \ \left(
\begin{array}
[c]{c}%
\text{here, we have split off the addend for }i=1\\
\text{from the sum}%
\end{array}
\right) \\
&  =\delta_{n-1,-1}+\underbrace{\sum_{\substack{i\in\mathbb{Z};\\i\neq
1}}0\delta_{n-i,-1}}_{=0}=\delta_{n-1,-1}=\delta_{n,0}%
\ \ \ \ \ \ \ \ \ \ \left(
\begin{array}
[c]{c}%
\text{since }n-1=-1\text{ holds}\\
\text{if and only if }n=0
\end{array}
\right)  .
\end{align*}
Thus, $\left(  c_{n}\right)  _{n\in\mathbb{Z}}=\left(  \delta_{n,0}\right)
_{n\in\mathbb{Z}}=1$, so that $x\overline{x}=\left(  c_{n}\right)
_{n\in\mathbb{Z}}=1$, as we wanted to prove. The proof of $\overline{x}x=1$ is
analogous, or follows from $x\overline{x}=1$ by commutativity.)
\end{proof}
\end{fineprint}

\begin{fineprint}
Proposition \ref{prop.fps.laure.a=sumaixi} is an analogue of Corollary
\ref{cor.fps.sumakxk}. Recall that we proved the latter corollary using Lemma
\ref{lem.fps.xa} and Proposition \ref{prop.fps.xk}. Thus, in order to prove
the former proposition, we need to find analogues of Lemma \ref{lem.fps.xa}
and Proposition \ref{prop.fps.xk} for Laurent polynomials.

The analogue of Lemma \ref{lem.fps.xa} for Laurent polynomials is the following:

\begin{lemma}
\label{lem.fps.laure.xa}Let $\mathbf{a}=\left(  a_{n}\right)  _{n\in
\mathbb{Z}}$ be a Laurent polynomial in $K\left[  x^{\pm}\right]  $. Then,
\[
x\cdot\mathbf{a}=\left(  a_{n-1}\right)  _{n\in\mathbb{Z}}%
\ \ \ \ \ \ \ \ \ \ \text{and}\ \ \ \ \ \ \ \ \ \ x^{-1}\cdot\mathbf{a}%
=\left(  a_{n+1}\right)  _{n\in\mathbb{Z}}.
\]

\end{lemma}

\begin{proof}
[Proof of Lemma \ref{lem.fps.laure.xa} (sketched).]From $x=\left(
\delta_{i,1}\right)  _{i\in\mathbb{Z}}$ and $\mathbf{a}=\left(  a_{n}\right)
_{n\in\mathbb{Z}}=\left(  a_{i}\right)  _{i\in\mathbb{Z}}$, we obtain%
\[
x\cdot\mathbf{a}=\left(  \delta_{i,1}\right)  _{i\in\mathbb{Z}}\cdot\left(
a_{i}\right)  _{i\in\mathbb{Z}}=\left(  c_{n}\right)  _{n\in\mathbb{Z}%
},\ \ \ \ \ \ \ \ \ \ \text{where}\ \ \ \ \ \ \ \ \ \ c_{n}=\sum
_{i\in\mathbb{Z}}\delta_{i,1}a_{n-i}%
\]
(by Definition \ref{def.fps.laure.laupol}). However, for each $n\in\mathbb{Z}%
$, we have%
\begin{align*}
c_{n}  &  =\sum_{i\in\mathbb{Z}}\delta_{i,1}a_{n-i}=\underbrace{\delta_{1,1}%
}_{=1}a_{n-1}+\sum_{\substack{i\in\mathbb{Z};\\i\neq1}}\underbrace{\delta
_{i,1}}_{\substack{=0\\\text{(since }i\neq1\text{)}}}a_{n-i}%
\ \ \ \ \ \ \ \ \ \ \left(
\begin{array}
[c]{c}%
\text{here, we have split off the}\\
\text{addend for }i=1\\
\text{from the sum}%
\end{array}
\right) \\
&  =a_{n-1}+\underbrace{\sum_{\substack{i\in\mathbb{Z};\\i\neq1}}0a_{n-i}%
}_{=0}=a_{n-1}.
\end{align*}
Thus, $\left(  c_{n}\right)  _{n\in\mathbb{Z}}=\left(  a_{n-1}\right)
_{n\in\mathbb{Z}}$, so that $x\cdot\mathbf{a}=\left(  c_{n}\right)
_{n\in\mathbb{Z}}=\left(  a_{n-1}\right)  _{n\in\mathbb{Z}}$. Thus we have
proved%
\begin{equation}
x\cdot\mathbf{a}=\left(  a_{n-1}\right)  _{n\in\mathbb{Z}}.
\label{pf.lem.fps.laure.xa.1}%
\end{equation}

It remains to prove that $x^{-1}\cdot\mathbf{a}=\left(  a_{n+1}\right)
_{n\in\mathbb{Z}}$. Here we use a trick: Let $\mathbf{b}=\left(
a_{n+1}\right)  _{n\in\mathbb{Z}}$; this is again a Laurent polynomial in
$K\left[  x^{\pm}\right]  $. Hence, (\ref{pf.lem.fps.laure.xa.1}) (applied to
$\mathbf{b}$ and $a_{n+1}$ instead of $\mathbf{a}$ and $a_{n}$) yields%
\[
x\cdot\mathbf{b}=\left(  a_{\left(  n-1\right)  +1}\right)  _{n\in\mathbb{Z}%
}=\left(  a_{n}\right)  _{n\in\mathbb{Z}}%
\]
(since $a_{\left(  n-1\right)  +1}=a_{n}$ for each $n\in\mathbb{Z}$). In other
words, $x\cdot\mathbf{b}=\mathbf{a}$ (since $\mathbf{a}=\left(  a_{n}\right)
_{n\in\mathbb{Z}}$). Dividing this equality by the invertible element $x$, we
find $\mathbf{b}=x^{-1}\cdot\mathbf{a}$, so that $x^{-1}\cdot\mathbf{a}%
=\mathbf{b}=\left(  a_{n+1}\right)  _{n\in\mathbb{Z}}$. This completes the
proof of Lemma \ref{lem.fps.laure.xa}.
\end{proof}

The analogue of Proposition \ref{prop.fps.xk} for Laurent polynomials is the following:

\begin{proposition}
\label{prop.fps.laure.xk}We have%
\[
x^{k}=\left(  \delta_{i,k}\right)  _{i\in\mathbb{Z}}%
\ \ \ \ \ \ \ \ \ \ \text{for each }k\in\mathbb{Z}.
\]

\end{proposition}

\begin{proof}
[Proof of Proposition \ref{prop.fps.laure.xk} (sketched).]This is similar to
Proposition \ref{prop.fps.xk}, which we proved by induction on $k$. Here, too,
we can use induction, but (since $k$ can be negative) we have to use
\textquotedblleft two-sided induction\textquotedblright, which contains both
an induction step from $k$ to $k+1$ and an induction step from $k$ to $k-1$.
(See \cite[\S 2.15]{detnotes} for a detailed explanation of two-sided induction.)

Both induction steps rely on Lemma \ref{lem.fps.laure.xa}. (Specifically, the
step from $k$ to $k+1$ uses the $x\cdot\mathbf{a}=\left(  a_{n-1}\right)
_{n\in\mathbb{Z}}$ part of Lemma \ref{lem.fps.laure.xa}, whereas the step from
$k$ to $k-1$ uses the $x^{-1}\cdot\mathbf{a}=\left(  a_{n+1}\right)
_{n\in\mathbb{Z}}$ part.)
\end{proof}

Now, Proposition \ref{prop.fps.laure.a=sumaixi} is easy:

\begin{proof}
[Proof of Proposition \ref{prop.fps.laure.a=sumaixi} (sketched).]Analogous to
Corollary \ref{cor.fps.sumakxk}, but using Proposition \ref{prop.fps.laure.xk}
instead of Proposition \ref{prop.fps.xk}.
\end{proof}

Next, let us prove Theorem \ref{thm.fps.laure.lauser-ring}:

\begin{proof}
[Proof of Theorem \ref{thm.fps.laure.lauser-ring} (sketched).]This is
analogous to the proof of Theorem \ref{thm.fps.laure.laupol-ring}. The only
(slightly) different piece is the proof that the set $K\left(  \left(
x\right)  \right)  $ is closed under multiplication. So let us prove this:

Let $\left(  a_{n}\right)  _{n\in\mathbb{Z}}$ and $\left(  b_{n}\right)
_{n\in\mathbb{Z}}$ be two elements of $K\left(  \left(  x\right)  \right)  $.
We must prove that their product $\left(  a_{n}\right)  _{n\in\mathbb{Z}}%
\cdot\left(  b_{n}\right)  _{n\in\mathbb{Z}}$ belongs to $K\left(  \left(
x\right)  \right)  $ as well. This product is defined by%
\[
\left(  a_{n}\right)  _{n\in\mathbb{Z}}\cdot\left(  b_{n}\right)
_{n\in\mathbb{Z}}=\left(  c_{n}\right)  _{n\in\mathbb{Z}}%
,\ \ \ \ \ \ \ \ \ \ \text{where}\ \ \ \ \ \ \ \ \ \ c_{n}=\sum_{i\in
\mathbb{Z}}a_{i}b_{n-i}.
\]
Thus, we need to prove that $\left(  c_{n}\right)  _{n\in\mathbb{Z}}$ belongs
to $K\left(  \left(  x\right)  \right)  $. In other words, we must prove that
the sequence $\left(  c_{-1},c_{-2},c_{-3},\ldots\right)  $ is essentially
finite (by the definition of $K\left(  \left(  x\right)  \right)  $).

So let us prove this now. We know that the sequence $\left(  a_{-1}%
,a_{-2},a_{-3},\ldots\right)  $ is essentially finite (since $\left(
a_{n}\right)  _{n\in\mathbb{Z}}\in K\left(  \left(  x\right)  \right)  $).
Thus, there exists some negative integer $p$ such that%
\begin{equation}
\text{all }i\leq p\text{ satisfy }a_{i}=0.
\label{pf.thm.fps.laure.lauser-ring.3}%
\end{equation}
Similarly, there exists some negative integer $q$ such that%
\begin{equation}
\text{all }j\leq q\text{ satisfy }b_{j}=0
\label{pf.thm.fps.laure.lauser-ring.4}%
\end{equation}
(since $\left(  b_{n}\right)  _{n\in\mathbb{Z}}\in K\left(  \left(  x\right)
\right)  $). Consider these $p$ and $q$. Now, set $r:=p+q$. We shall show that
all negative integers $n\leq r$ satisfy $c_{n}=0$.

Indeed, let $n\leq r$ be any integer. Then, each integer $i\geq p$ satisfies
$\underbrace{n}_{\leq r=p+q}-\underbrace{i}_{\geq p}\leq p+q-p=q$ and
therefore
\begin{equation}
b_{n-i}=0 \label{pf.thm.fps.laure.lauser-ring.4b}%
\end{equation}
(by (\ref{pf.thm.fps.laure.lauser-ring.4}), applied to $j=n-i$). Now,%
\[
c_{n}=\sum_{i\in\mathbb{Z}}a_{i}b_{n-i}=\sum_{\substack{i\in\mathbb{Z}%
;\\i<p}}\underbrace{a_{i}}_{\substack{=0\\\text{(by
(\ref{pf.thm.fps.laure.lauser-ring.3}))}}}b_{n-i}+\sum_{\substack{i\in
\mathbb{Z};\\i\geq p}}a_{i}\underbrace{b_{n-i}}_{\substack{=0\\\text{(by
(\ref{pf.thm.fps.laure.lauser-ring.4b}))}}}=\underbrace{\sum_{\substack{i\in
\mathbb{Z};\\i<p}}0b_{n-i}}_{=0}+\underbrace{\sum_{\substack{i\in
\mathbb{Z};\\i\geq p}}a_{i}0}_{=0}=0.
\]

Forget that we fixed $n$. We thus have shown that all $n\leq r$ satisfy
$c_{n}=0$. Hence, the sequence $\left(  c_{-1},c_{-2},c_{-3},\ldots\right)  $
is essentially finite. As explained, this completes our proof of the claim
that the set $K\left(  \left(  x\right)  \right)  $ is closed under multiplication.

The rest of Theorem \ref{thm.fps.laure.lauser-ring} is proved just like
Theorem \ref{thm.fps.laure.laupol-ring}.
\end{proof}
\end{fineprint}

\subsection{\label{sec.details.sign.intro}Cancellations in alternating sums}

\begin{fineprint}
We shall now prove Lemma \ref{lem.sign.cancel2} and Lemma
\ref{lem.sign.cancel3}. We start with the latter lemma, since the former will
then follow trivially from it.

\begin{proof}
[Detailed proof of Lemma \ref{lem.sign.cancel3}.]The set $\mathcal{X}$ is
finite (since it is a subset of the finite set $\mathcal{A}$). Thus,
$\left\vert \mathcal{X}\right\vert =n$ for some $n\in\mathbb{N}$. Consider
this $n$.

Let $\left[  n\right]  $ be the set $\left\{  1,2,\ldots,n\right\}  $. Then,
$\left\vert \left[  n\right]  \right\vert =n$. Comparing this with $\left\vert
\mathcal{X}\right\vert =n$, we obtain $\left\vert \mathcal{X}\right\vert
=\left\vert \left[  n\right]  \right\vert $. Hence, there exists a bijection
$\alpha:\mathcal{X}\rightarrow\left[  n\right]  $. Consider this $\alpha$.

Now, define two subsets $\mathcal{U}$ and $\mathcal{W}$ of $\mathcal{X}$ by%
\begin{align*}
\mathcal{U}  &  :=\left\{  I\in\mathcal{X}\ \mid\ \alpha\left(  f\left(
I\right)  \right)  <\alpha\left(  I\right)  \right\}  ;\\
\mathcal{W}  &  :=\left\{  I\in\mathcal{X}\ \mid\ \alpha\left(  f\left(
I\right)  \right)  >\alpha\left(  I\right)  \right\}  .
\end{align*}

Then, $f\left(  \mathcal{U}\right)  \subseteq\mathcal{W}$%
\ \ \ \ \footnote{\textit{Proof.} Let $J\in f\left(  \mathcal{U}\right)  $.
Thus, $J=f\left(  K\right)  $ for some $K\in\mathcal{U}$. Consider this $K$.
We have $K\in\mathcal{U}=\left\{  I\in\mathcal{X}\ \mid\ \alpha\left(
f\left(  I\right)  \right)  <\alpha\left(  I\right)  \right\}  $; in other
words, $K$ is an $I\in\mathcal{X}$ satisfying $\alpha\left(  f\left(
I\right)  \right)  <\alpha\left(  I\right)  $. In other words, $K$ is an
element of $\mathcal{X}$ and satisfies $\alpha\left(  f\left(  K\right)
\right)  <\alpha\left(  K\right)  $. However, we have $f\circ
f=\operatorname*{id}$ (since $f$ is an involution) and thus $\left(  f\circ
f\right)  \left(  K\right)  =\operatorname*{id}\left(  K\right)  =K$, so that
$K=\left(  f\circ f\right)  \left(  K\right)  =f\left(  \underbrace{f\left(
K\right)  }_{=J}\right)  =f\left(  J\right)  $. However, from $J=f\left(
K\right)  $, we obtain $\alpha\left(  J\right)  =\alpha\left(  f\left(
K\right)  \right)  <\alpha\left(  \underbrace{K}_{=f\left(  J\right)
}\right)  =\alpha\left(  f\left(  J\right)  \right)  $, so that $\alpha\left(
f\left(  J\right)  \right)  >\alpha\left(  J\right)  $.
\par
Now, $J$ is an element of $\mathcal{X}$ and satisfies $\alpha\left(  f\left(
J\right)  \right)  >\alpha\left(  J\right)  $. In other words, $J$ is an
$I\in\mathcal{X}$ satisfying $\alpha\left(  f\left(  I\right)  \right)
>\alpha\left(  I\right)  $. In other words, $J\in\left\{  I\in\mathcal{X}%
\ \mid\ \alpha\left(  f\left(  I\right)  \right)  >\alpha\left(  I\right)
\right\}  $. In other words, $J\in\mathcal{W}$ (since $\mathcal{W}=\left\{
I\in\mathcal{X}\ \mid\ \alpha\left(  f\left(  I\right)  \right)
>\alpha\left(  I\right)  \right\}  $).
\par
Forget that we fixed $J$. We thus have shown that $J\in\mathcal{W}$ for each
$J\in f\left(  \mathcal{U}\right)  $. In other words, $f\left(  \mathcal{U}%
\right)  \subseteq\mathcal{W}$.} and $f\left(  \mathcal{W}\right)
\subseteq\mathcal{U}$\ \ \ \ \footnote{The \textit{proof} of $f\left(
\mathcal{W}\right)  \subseteq\mathcal{U}$ is completely analogous to the proof
of $f\left(  \mathcal{U}\right)  \subseteq\mathcal{W}$ we just gave; the only
changes are that all \textquotedblleft$<$\textquotedblright\ signs have to be
replaced by \textquotedblleft$>$\textquotedblright\ signs and vice versa, and
that all \textquotedblleft$\mathcal{U}$\textquotedblright s have to be
replaced by \textquotedblleft$\mathcal{W}$\textquotedblright s and vice
versa.}. Now, the map%
\begin{align*}
g:\mathcal{U}  &  \rightarrow\mathcal{W},\\
I  &  \mapsto f\left(  I\right)
\end{align*}
is well-defined (since each $I\in\mathcal{U}$ satisfies $f\left(  I\right)
\in f\left(  \mathcal{U}\right)  \subseteq\mathcal{W}$), and the map
\begin{align*}
h:\mathcal{W}  &  \rightarrow\mathcal{U},\\
I  &  \mapsto f\left(  I\right)
\end{align*}
is also well-defined (since each $I\in\mathcal{W}$ satisfies $f\left(
I\right)  \in f\left(  \mathcal{W}\right)  \subseteq\mathcal{U}$). Consider
these two maps $g$ and $h$. It is clear that $g\circ h=\operatorname*{id}%
$\ \ \ \ \footnote{\textit{Proof.} Let $I\in\mathcal{W}$. Recall that $f$ is
an involution; thus, $f\circ f=\operatorname*{id}$. Now, $\left(  g\circ
h\right)  \left(  I\right)  =g\left(  h\left(  I\right)  \right)  =f\left(
h\left(  I\right)  \right)  $ (by the definition of $g$). However, $h\left(
I\right)  =f\left(  I\right)  $ (by the definition of $h$). Thus, $\left(
g\circ h\right)  \left(  I\right)  =f\left(  \underbrace{h\left(  I\right)
}_{=f\left(  I\right)  }\right)  =f\left(  f\left(  I\right)  \right)
=\underbrace{\left(  f\circ f\right)  }_{=\operatorname*{id}}\left(  I\right)
=\operatorname*{id}\left(  I\right)  =I=\operatorname*{id}\left(  I\right)  $.
\par
Forget that we fixed $I$. We thus have shown that $\left(  g\circ h\right)
\left(  I\right)  =\operatorname*{id}\left(  I\right)  $ for each
$I\in\mathcal{W}$. In other words, $g\circ h=\operatorname*{id}$.} and $h\circ
g=\operatorname*{id}$\ \ \ \ \footnote{The \textit{proof} of this is analogous
to the proof of $g\circ h=\operatorname*{id}$ we just gave.}. Hence, the maps
$g$ and $h$ are mutually inverse, and thus are bijections.

We have%
\begin{equation}
\operatorname*{sign}I+\operatorname*{sign}\left(  g\left(  I\right)  \right)
=0\ \ \ \ \ \ \ \ \ \ \text{for all }I\in\mathcal{U}
\label{pf.lem.sign.cancel3.signg}%
\end{equation}
\footnote{\textit{Proof of (\ref{pf.lem.sign.cancel3.signg}):} Let
$J\in\mathcal{U}$. Then, $J\in\mathcal{U}\subseteq\mathcal{X}$ and $g\left(
J\right)  =f\left(  J\right)  $ (by the definition of $g$). Now, recall our
assumption saying that $\operatorname*{sign}\left(  f\left(  I\right)
\right)  =-\operatorname*{sign}I$ for all $I\in\mathcal{X}$. Applying this to
$I=J$, we obtain $\operatorname*{sign}\left(  f\left(  J\right)  \right)
=-\operatorname*{sign}J$. In view of $g\left(  J\right)  =f\left(  J\right)
$, this rewrites as $\operatorname*{sign}\left(  g\left(  J\right)  \right)
=-\operatorname*{sign}J$. In other words, $\operatorname*{sign}%
J+\operatorname*{sign}\left(  g\left(  J\right)  \right)  =0$.
\par
Forget that we fixed $J$. We thus have shown that $\operatorname*{sign}%
J+\operatorname*{sign}\left(  g\left(  J\right)  \right)  =0$ for all
$J\in\mathcal{U}$. Renaming $J$ as $I$ in this statement, we obtain that
$\operatorname*{sign}I+\operatorname*{sign}\left(  g\left(  I\right)  \right)
=0$ for all $I\in\mathcal{U}$. This proves (\ref{pf.lem.sign.cancel3.signg}%
).}. Furthermore, we have%
\begin{equation}
\operatorname*{sign}I=0\ \ \ \ \ \ \ \ \ \ \text{for all }I\in\mathcal{X}%
\text{ satisfying }\alpha\left(  f\left(  I\right)  \right)  =\alpha\left(
I\right)  . \label{pf.lem.sign.cancel3.signV}%
\end{equation}
\footnote{\textit{Proof of (\ref{pf.lem.sign.cancel3.signV}):} Recall our
assumption saying that
\begin{equation}
\operatorname*{sign}I=0\ \ \ \ \ \ \ \ \ \ \text{for all }I\in\mathcal{X}%
\text{ satisfying }f\left(  I\right)  =I.
\label{pf.lem.sign.cancel3.signV.pf.1}%
\end{equation}
\par
Now, let $I\in\mathcal{X}$ satisfy $\alpha\left(  f\left(  I\right)  \right)
=\alpha\left(  I\right)  $. The map $\alpha$ is injective (since $\alpha$ is a
bijection). Thus, from $\alpha\left(  f\left(  I\right)  \right)
=\alpha\left(  I\right)  $, we obtain $f\left(  I\right)  =I$. Therefore,
(\ref{pf.lem.sign.cancel3.signV.pf.1}) yields $\operatorname*{sign}I=0$.
\par
Forget that we fixed $I$. We thus have shown that $\operatorname*{sign}I=0$
for all $I\in\mathcal{X}$ satisfying $\alpha\left(  f\left(  I\right)
\right)  =\alpha\left(  I\right)  $. This proves
(\ref{pf.lem.sign.cancel3.signV}).}.

However, each $I\in\mathcal{X}$ satisfies exactly one of the three conditions
\textquotedblleft$\alpha\left(  f\left(  I\right)  \right)  <\alpha\left(
I\right)  $\textquotedblright, \textquotedblleft$\alpha\left(  f\left(
I\right)  \right)  =\alpha\left(  I\right)  $\textquotedblright\ and
\textquotedblleft$\alpha\left(  f\left(  I\right)  \right)  >\alpha\left(
I\right)  $\textquotedblright\ (because $\alpha\left(  f\left(  I\right)
\right)  $ and $\alpha\left(  I\right)  $ are two integers). Hence, we can
split the sum $\sum_{I\in\mathcal{X}}\operatorname*{sign}I$ as follows:%
\begin{align*}
\sum_{I\in\mathcal{X}}\operatorname*{sign}I  &  =\sum_{\substack{I\in
\mathcal{X};\\\alpha\left(  f\left(  I\right)  \right)  <\alpha\left(
I\right)  }}\operatorname*{sign}I+\sum_{\substack{I\in\mathcal{X}%
;\\\alpha\left(  f\left(  I\right)  \right)  =\alpha\left(  I\right)
}}\underbrace{\operatorname*{sign}I}_{\substack{=0\\\text{(by
(\ref{pf.lem.sign.cancel3.signV}))}}}+\sum_{\substack{I\in\mathcal{X}%
;\\\alpha\left(  f\left(  I\right)  \right)  >\alpha\left(  I\right)
}}\operatorname*{sign}I\\
&  =\sum_{\substack{I\in\mathcal{X};\\\alpha\left(  f\left(  I\right)
\right)  <\alpha\left(  I\right)  }}\operatorname*{sign}I+\underbrace{\sum
_{\substack{I\in\mathcal{X};\\\alpha\left(  f\left(  I\right)  \right)
=\alpha\left(  I\right)  }}0}_{=0}+\sum_{\substack{I\in\mathcal{X}%
;\\\alpha\left(  f\left(  I\right)  \right)  >\alpha\left(  I\right)
}}\operatorname*{sign}I\\
&  =\underbrace{\sum_{\substack{I\in\mathcal{X};\\\alpha\left(  f\left(
I\right)  \right)  <\alpha\left(  I\right)  }}}_{\substack{=\sum
_{I\in\mathcal{U}}\\\text{(since }\left\{  I\in\mathcal{X}\ \mid
\ \alpha\left(  f\left(  I\right)  \right)  <\alpha\left(  I\right)  \right\}
=\mathcal{U}\text{)}}}\operatorname*{sign}I+\underbrace{\sum_{\substack{I\in
\mathcal{X};\\\alpha\left(  f\left(  I\right)  \right)  >\alpha\left(
I\right)  }}}_{\substack{=\sum_{I\in\mathcal{W}}\\\text{(since }\left\{
I\in\mathcal{X}\ \mid\ \alpha\left(  f\left(  I\right)  \right)
>\alpha\left(  I\right)  \right\}  =\mathcal{W}\text{)}}}\operatorname*{sign}%
I\\
&  =\sum_{I\in\mathcal{U}}\operatorname*{sign}I+\underbrace{\sum
_{I\in\mathcal{W}}\operatorname*{sign}I}_{\substack{=\sum_{I\in\mathcal{U}%
}\operatorname*{sign}\left(  g\left(  I\right)  \right)  \\\text{(here, we
have substituted }g\left(  I\right)  \\\text{for }I\text{ in the sum, since
the}\\\text{map }g:\mathcal{U}\rightarrow\mathcal{W}\text{ is a bijection)}%
}}=\sum_{I\in\mathcal{U}}\operatorname*{sign}I+\sum_{I\in\mathcal{U}%
}\operatorname*{sign}\left(  g\left(  I\right)  \right) \\
&  =\sum_{I\in\mathcal{U}}\underbrace{\left(  \operatorname*{sign}%
I+\operatorname*{sign}\left(  g\left(  I\right)  \right)  \right)
}_{\substack{=0\\\text{(by (\ref{pf.lem.sign.cancel3.signg}))}}}=\sum
_{I\in\mathcal{U}}0=0.
\end{align*}
However, the set $\mathcal{A}$ is the union of its two disjoint subsets
$\mathcal{X}$ and $\mathcal{A}\setminus\mathcal{X}$ (since $\mathcal{X}%
\subseteq\mathcal{A}$). Thus, we can split the sum $\sum_{I\in\mathcal{A}%
}\operatorname*{sign}I$ as follows:%
\[
\sum_{I\in\mathcal{A}}\operatorname*{sign}I=\underbrace{\sum_{I\in\mathcal{X}%
}\operatorname*{sign}I}_{=0}+\sum_{I\in\mathcal{A}\setminus\mathcal{X}%
}\operatorname*{sign}I=\sum_{I\in\mathcal{A}\setminus\mathcal{X}%
}\operatorname*{sign}I.
\]
This proves Lemma \ref{lem.sign.cancel3}.
\end{proof}

\begin{proof}
[Detailed proof of Lemma \ref{lem.sign.cancel2}.]The map $f$ has no fixed
points (by assumption). In other words, there exist no $I\in\mathcal{X}$
satisfying $f\left(  I\right)  =I$. Hence, we have $\operatorname*{sign}I=0$
for all $I\in\mathcal{X}$ satisfying $f\left(  I\right)  =I$ (because
non-existing objects satisfy any possible claim; this is known as being
\textquotedblleft vacuously true\textquotedblright). Thus, Lemma
\ref{lem.sign.cancel3} yields $\sum_{I\in\mathcal{A}}\operatorname*{sign}%
I=\sum_{I\in\mathcal{A}\setminus\mathcal{X}}\operatorname*{sign}I$. This
proves Lemma \ref{lem.sign.cancel2}.
\end{proof}
\end{fineprint}

\subsection{\label{sec.details.det.comb}Determinants in combinatorics}

\begin{fineprint}
\begin{proof}
[Proof of Proposition \ref{prop.lgv.jordan-2}.]If $q$ is any path, then the
\emph{length} $\ell\left(  q\right)  $ of $q$ is defined to be the number of
arcs of $q$.

We shall now prove Proposition \ref{prop.lgv.jordan-2} by strong induction on
$\ell\left(  p\right)  +\ell\left(  p^{\prime}\right)  $:

\textit{Induction step:} Fix a nonnegative integer $N$. Assume (as the
induction hypothesis) that Proposition \ref{prop.lgv.jordan-2} holds whenever
$\ell\left(  p\right)  +\ell\left(  p^{\prime}\right)  <N$. We must now prove
that Proposition \ref{prop.lgv.jordan-2} holds when $\ell\left(  p\right)
+\ell\left(  p^{\prime}\right)  =N$.

So let $A$, $B$, $A^{\prime}$, $B^{\prime}$, $p$ and $p^{\prime}$ be as in
Proposition \ref{prop.lgv.jordan-2}, and let us assume that $\ell\left(
p\right)  +\ell\left(  p^{\prime}\right)  =N$. We must prove that $p$ and
$p^{\prime}$ have a vertex in common.

Assume the contrary. Thus, $p$ and $p^{\prime}$ have no vertex in common.

Recall that each arc of the lattice is either an east-step or a north-step.
Thus, the x-coordinates of the vertices of a path are always weakly increasing
(i.e., if $\left(  v_{0},v_{1},\ldots,v_{n}\right)  $ is a path, then
$\operatorname*{x}\left(  v_{0}\right)  \leq\operatorname*{x}\left(
v_{1}\right)  \leq\cdots\leq\operatorname*{x}\left(  v_{n}\right)  $), and so
are the y-coordinates. Hence, the existence of a path $p$ from $A$ to
$B^{\prime}$ shows that $\operatorname{x}\left(  A\right)  \leq
\operatorname{x}\left(  B^{\prime}\right)  $ and $\operatorname{y}\left(
A\right)  \leq\operatorname{y}\left(  B^{\prime}\right)  $. Similarly, the
existence of a path $p^{\prime}$ from $A^{\prime}$ to $B$ yields
$\operatorname{x}\left(  A^{\prime}\right)  \leq\operatorname{x}\left(
B\right)  $ and $\operatorname{y}\left(  A^{\prime}\right)  \leq
\operatorname{y}\left(  B\right)  $.

We are in one of the following three cases:

\textit{Case 1:} We have $\operatorname*{y}\left(  A^{\prime}\right)
>\operatorname*{y}\left(  A\right)  $.

\textit{Case 2:} We have $\operatorname*{x}\left(  A^{\prime}\right)
<\operatorname*{x}\left(  A\right)  $.

\textit{Case 3:} We have neither $\operatorname*{y}\left(  A^{\prime}\right)
>\operatorname*{y}\left(  A\right)  $ nor $\operatorname*{x}\left(  A^{\prime
}\right)  <\operatorname*{x}\left(  A\right)  $.

We shall derive a contradiction in each of these cases.

Let us first consider Case 1. In this case, we have $\operatorname*{y}\left(
A^{\prime}\right)  >\operatorname*{y}\left(  A\right)  $. Thus,
$\operatorname*{y}\left(  A\right)  <\operatorname*{y}\left(  A^{\prime
}\right)  \leq\operatorname*{y}\left(  B\right)  \leq\operatorname*{y}\left(
B^{\prime}\right)  $ (since $\operatorname*{y}\left(  B^{\prime}\right)
\geq\operatorname*{y}\left(  B\right)  $), so that $\operatorname*{y}\left(
A\right)  \neq\operatorname*{y}\left(  B^{\prime}\right)  $ and therefore
$A\neq B^{\prime}$. This shows that the path $p$ has at least two vertices
(since $p$ is a path from $A$ to $B^{\prime}$). Let $P$ be its second vertex.
Hence, $P$ lies on a path from $A$ to $B^{\prime}$ (namely, on the path $p$).
Therefore, $\operatorname{x}\left(  A\right)  \leq\operatorname{x}\left(
P\right)  \leq\operatorname{x}\left(  B^{\prime}\right)  $ (since the
x-coordinates of the vertices of a path are always weakly increasing) and
$\operatorname{y}\left(  A\right)  \leq\operatorname{y}\left(  P\right)
\leq\operatorname{y}\left(  B^{\prime}\right)  $ (since the y-coordinates of
the vertices of a path are always weakly increasing). Let $r$ be the path from
$P$ to $B^{\prime}$ obtained by removing the first arc from $p$ (in other
words, let $r$ be the part of $p$ from the point $P$ onwards)\footnote{Here is
an illustration (with $r$ drawn extra-thick):%
\[%
\begin{tikzpicture}
\draw[densely dotted] (0,0) grid (7.2,7.2);
\draw[->] (0,0) -- (0,7.2);
\draw[->] (0,0) -- (7.2,0);
\foreach\x/\xtext in {0, 1, 2, 3, 4, 5, 6, 7}
\draw(\x cm,1pt) -- (\x cm,-1pt) node[anchor=north] {$\xtext$};
\foreach\y/\ytext in {0, 1, 2, 3, 4, 5, 6, 7}
\draw(1pt,\y cm) -- (-1pt,\y cm) node[anchor=east] {$\ytext$};
\node
[circle,fill=white,draw=black,text=black,inner sep=1pt] (A') at (1,3) {$A^{\prime
}$};
\node
[circle,fill=white,draw=black,text=black,inner sep=1pt] (A) at (2,1) {$A$};
\node
[circle,fill=white,draw=black,text=black,inner sep=1pt] (P) at (3,1) {$P$};
\node
[circle,fill=white,draw=black,text=black,inner sep=1pt] (B') at (5,6) {$B^{\prime
}$};
\node
[circle,fill=white,draw=black,text=black,inner sep=1pt] (B) at (6,4) {$B$};
\begin{scope}[thick,>=stealth,darkred]
\draw(A) edge[->] (P);
\draw(2.5,1) node[anchor=north] {$p$};
\draw(P) edge[ultra thick, ->] (3,2);
\draw(3,2) edge[ultra thick, ->] (4,2);
\draw(4,2) edge[ultra thick, ->] (4,3);
\draw(4,3) edge[ultra thick, ->] (4,4);
\draw(4,4) edge[ultra thick, ->] (4,5);
\draw(4,5) edge[ultra thick, ->] (5,5);
\draw(5,5) edge[ultra thick, ->] (B');
\draw(5,5.4) node[anchor=east] {$r$};
\end{scope}
\begin{scope}[thick,>=stealth,dbluecolor]
\draw(A') edge[->] (2,3);
\draw(1.6,3) node[anchor=north] {$p^{\prime}$};
\draw(2,3) edge[->] (2,4);
\draw(2,4) edge[->] (3,4);
\draw(3,4) edge[->] (4,4);
\draw(4,4) edge[->] (5,4);
\draw(5,4) edge[->] (B);
\draw(5.4,4) node[anchor=north] {$p^{\prime}$};
\end{scope}
\end{tikzpicture}%
\]
}. Thus, every vertex of $r$ is a vertex of $p$. Hence, the paths $r$ and
$p^{\prime}$ have no vertex in common (since $p$ and $p^{\prime}$ have no
vertex in common). Also, $\ell\left(  r\right)  =\ell\left(  p\right)
-1<\ell\left(  p\right)  $ and thus $\underbrace{\ell\left(  r\right)
}_{<\ell\left(  p\right)  }+\,\ell\left(  p^{\prime}\right)  <\ell\left(
p\right)  +\ell\left(  p^{\prime}\right)  =N$. Moreover, $\operatorname{x}%
\left(  A^{\prime}\right)  \leq\operatorname{x}\left(  A\right)
\leq\operatorname{x}\left(  P\right)  $. Thus, if we had $\operatorname{y}%
\left(  A^{\prime}\right)  \geq\operatorname{y}\left(  P\right)  $, then we
could apply Proposition \ref{prop.lgv.jordan-2} to $P$ and $r$ instead of $A$
and $p$ (by the induction hypothesis, since $\ell\left(  r\right)
+\ell\left(  p^{\prime}\right)  <N$). We would consequently conclude that the
paths $r$ and $p^{\prime}$ have a vertex in common; this would contradict the
fact that the paths $r$ and $p^{\prime}$ have no vertex in common. Hence, we
cannot have $\operatorname{y}\left(  A^{\prime}\right)  \geq\operatorname{y}%
\left(  P\right)  $. Thus, $\operatorname{y}\left(  A^{\prime}\right)
<\operatorname{y}\left(  P\right)  $. Hence, $\operatorname{y}\left(
A^{\prime}\right)  \leq\operatorname{y}\left(  P\right)  -1$ (since
$\operatorname{y}\left(  A^{\prime}\right)  $ and $\operatorname{y}\left(
P\right)  $ are integers). But $P$ is the next vertex after $A$ on the path
$p$. Hence, there is an arc from $A$ to $P$. If this arc was an east-step,
then we would have $\operatorname{y}\left(  P\right)  =\operatorname{y}\left(
A\right)  $, which would contradict $\operatorname{y}\left(  A\right)
\leq\operatorname{y}\left(  A^{\prime}\right)  <\operatorname{y}\left(
P\right)  $. Hence, this arc cannot be an east-step. Thus, this arc must be a
north-step. Therefore, $\operatorname{y}\left(  P\right)  =\operatorname{y}%
\left(  A\right)  +1$. Hence, $\operatorname{y}\left(  A^{\prime}\right)
\leq\operatorname{y}\left(  P\right)  -1=\operatorname{y}\left(  A\right)  $
(since $\operatorname{y}\left(  P\right)  =\operatorname{y}\left(  A\right)
+1$). But this contradicts $\operatorname{y}\left(  A^{\prime}\right)
>\operatorname{y}\left(  A\right)  $. Thus, we have found a contradiction in
Case 1.

An analogous argument can be used to find a contradiction in Case 2. In fact,
there is a symmetry inherent in Proposition \ref{prop.lgv.jordan-2}, which
interchanges Case 1 with Case 2. Namely, if we reflect all points and paths
across the $x=y$ line (i.e., if we replace each point $\left(  i,j\right)  $
by $\left(  j,i\right)  $), and if we rename $A$, $B$, $A^{\prime}$,
$B^{\prime}$, $p$ and $p^{\prime}$ as $A^{\prime}$, $A$, $B^{\prime}$, $B$,
$p^{\prime}$ and $p$ (respectively), then Case 1 becomes Case 2 and vice
versa\footnote{Importantly, this reflection preserves our digraph (in fact, it
transforms north-steps into east-steps and vice versa).}. Thus, Case 2 spawns
a contradiction just like Case 1 did.

Finally, let us consider Case 3. In this case, we have neither
$\operatorname*{y}\left(  A^{\prime}\right)  >\operatorname*{y}\left(
A\right)  $ nor $\operatorname*{x}\left(  A^{\prime}\right)
<\operatorname*{x}\left(  A\right)  $. In other words, we have
$\operatorname*{y}\left(  A^{\prime}\right)  \leq\operatorname*{y}\left(
A\right)  $ and $\operatorname*{x}\left(  A^{\prime}\right)  \geq
\operatorname*{x}\left(  A\right)  $. Combining $\operatorname*{x}\left(
A^{\prime}\right)  \geq\operatorname*{x}\left(  A\right)  $ with
$\operatorname*{x}\left(  A^{\prime}\right)  \leq\operatorname*{x}\left(
A\right)  $, we obtain $\operatorname*{x}\left(  A^{\prime}\right)
=\operatorname*{x}\left(  A\right)  $. Combining $\operatorname*{y}\left(
A^{\prime}\right)  \leq\operatorname*{y}\left(  A\right)  $ with
$\operatorname*{y}\left(  A^{\prime}\right)  \geq\operatorname*{y}\left(
A\right)  $, we obtain $\operatorname*{y}\left(  A^{\prime}\right)
=\operatorname*{y}\left(  A\right)  $. Now, the vertices $A$ and $A^{\prime}$
have the same x-coordinate (since $\operatorname{x}\left(  A^{\prime}\right)
=\operatorname{x}\left(  A\right)  $) and the same y-coordinate (since
$\operatorname{y}\left(  A^{\prime}\right)  =\operatorname{y}\left(  A\right)
$). Hence, these two vertices are equal. In other words, $A=A^{\prime}$.
Hence, the vertex $A$ belongs to the path $p^{\prime}$ (since the vertex
$A^{\prime}$ belongs to the path $p^{\prime}$). However, the vertex $A$
belongs to the path $p$ as well. Thus, the paths $p$ and $p^{\prime}$ have a
vertex in common (namely, $A$). This contradicts the fact that $p$ and
$p^{\prime}$ have no vertex in common. Thus, we have found a contradiction in
Case 3.

We have now found contradictions in all three Cases 1, 2 and 3. Hence, our
assumption must have been false. We thus conclude that $p$ and $p^{\prime}$
have a vertex in common. Now, forget that we fixed $A$, $B$, $A^{\prime}$,
$B^{\prime}$, $p$ and $p^{\prime}$. We thus have proven that if $A$, $B$,
$A^{\prime}$, $B^{\prime}$, $p$ and $p^{\prime}$ are as in Proposition
\ref{prop.lgv.jordan-2}, and if $\ell\left(  p\right)  +\ell\left(  p^{\prime
}\right)  =N$, then $p$ and $p^{\prime}$ have a vertex in common. In other
words, Proposition \ref{prop.lgv.jordan-2} holds when $\ell\left(  p\right)
+\ell\left(  p^{\prime}\right)  =N$. This completes the induction step. Hence,
Proposition \ref{prop.lgv.jordan-2} is proven.

(This proof is essentially the first proof from
\url{https://math.stackexchange.com/questions/2870640/} .)
\end{proof}

\begin{noncompile}
TODO: There is a simpler variant of the above proof, which replaces the
induction by the extremal principle. Namely, label the vertices on path $p$ as
$P_{0},P_{1},\ldots,P_{u}$, and the vertices on path $p^{\prime}$ as
$P_{0}^{\prime},P_{1}^{\prime},\ldots,P_{v}^{\prime}$. Thus, $P_{0}=A$ and
$P_{u}=B^{\prime}$ and $P_{0}^{\prime}=A^{\prime}$ and $P_{v}^{\prime}=B$.

A pair $\left(  i,j\right)  \in\left\{  0,1,\ldots,u\right\}  \times\left\{
0,1,\ldots,v\right\}  $ will be called \emph{good} if $\operatorname*{x}%
\left(  P_{i}\right)  \geq\operatorname*{x}\left(  P_{j}^{\prime}\right)  $
and $\operatorname*{y}\left(  P_{i}\right)  \leq\operatorname*{y}\left(
P_{j}^{\prime}\right)  $. There is at least one good pair (namely, $\left(
0,0\right)  $). Now pick a good pair $\left(  i,j\right)  $ with maximum
$i+j$. We are in one of the three cases \textquotedblleft$\operatorname*{x}%
\left(  P_{i}\right)  >\operatorname*{x}\left(  P_{j}^{\prime}\right)
$\textquotedblright\ and \textquotedblleft$\operatorname*{y}\left(
P_{i}\right)  <\operatorname*{y}\left(  P_{j}^{\prime}\right)  $%
\textquotedblright\ and \textquotedblleft neither of the two\textquotedblright%
. In Case 3, we get $P_{i}=P_{j}^{\prime}$, so we have found an intersection.
In Case 1, we argue that $j\neq v$ (since $\operatorname*{x}\left(
P_{j}^{\prime}\right)  <\operatorname*{x}\left(  P_{i}\right)  \leq
\operatorname*{x}\left(  P_{u}\right)  =\operatorname*{x}\left(  B^{\prime
}\right)  \leq\operatorname*{x}\left(  B\right)  =\operatorname*{x}\left(
P_{v}^{\prime}\right)  $) and therefore $P_{j+1}^{\prime}$ is well-defined,
and argue that $\left(  i,j+1\right)  $ is still good (since
$\operatorname*{x}\left(  P_{i}\right)  >\operatorname*{x}\left(
P_{j}^{\prime}\right)  $ and $\operatorname*{x}\left(  P_{j+1}^{\prime
}\right)  \leq\operatorname*{x}\left(  P_{j}^{\prime}\right)  +1$ lead to
$\operatorname*{x}\left(  P_{i}\right)  \geq\operatorname*{x}\left(
P_{j+1}^{\prime}\right)  $), which contradicts the maximality of $i+j$ in our
choice of $\left(  i,j\right)  $. Similarly for Case 2.
\end{noncompile}

\begin{proof}
[Some details for the proof of Proposition \ref{prop.lgv.kpaths.count}.]In our
above proof of Proposition \ref{prop.lgv.kpaths.count}, we have omitted some
details. Let us elaborate on them here. Specifically, we shall show that the
map $f:\mathcal{X}\rightarrow\mathcal{X}$ (constructed in the proof of
Proposition \ref{prop.lgv.kpaths.count}) is an involution. In other words, we
shall show that applying $f$ twice in succession to any given ipat $\left(
\sigma,\mathbf{p}\right)  \in\mathcal{X}$ returns $\left(  \sigma
,\mathbf{p}\right)  $.

So let $\left(  \sigma,\mathbf{p}\right)  \in\mathcal{X}$ be any ipat. Let
$\left(  \sigma^{\prime},\mathbf{q}\right)  :=f\left(  \sigma,\mathbf{p}%
\right)  $ and $\left(  \sigma^{\prime\prime},\mathbf{r}\right)  :=f\left(
\sigma^{\prime},\mathbf{q}\right)  $. We must thus show that $\left(
\sigma^{\prime\prime},\mathbf{r}\right)  =\left(  \sigma,\mathbf{p}\right)  $.

In order to see this, we write the $k$-tuples $\mathbf{p}$ and $\mathbf{q}$ as
$\mathbf{p}=\left(  p_{1},p_{2},\ldots,p_{k}\right)  $ and $\mathbf{q}=\left(
q_{1},q_{2},\ldots,q_{k}\right)  $. We recall that $f\left(  \sigma
,\mathbf{p}\right)  =\left(  \sigma^{\prime},\mathbf{q}\right)  $ was defined
by the following procedure:

\begin{itemize}
\item We say that a point $u$ is \emph{crowded} if it is contained in at least
two of our paths $p_{1},p_{2},\ldots,p_{k}$.

\item We pick the smallest $i\in\left[  k\right]  $ such that $p_{i}$ contains
a crowded point.

\item Then, we pick the first crowded point $v$ on $p_{i}$.

\item Then, we pick the largest $j\in\left[  k\right]  $ such that $v$ belongs
to $p_{j}$.

\item Call the part of $p_{i}$ that comes after $v$ the \emph{tail} of $p_{i}%
$. Call the part of $p_{i}$ that comes before $v$ the \emph{head} of $p_{i}$.
Call the part of $p_{j}$ that comes after $v$ the \emph{tail} of $p_{j}$. Call
the part of $p_{j}$ that comes before $v$ the \emph{head} of $p_{j}$.

\item Then, we exchange the tails of the paths $p_{i}$ and $p_{j}$ (while
leaving all other paths unchanged).

\item We let $\mathbf{q}=\left(  q_{1},q_{2},\ldots,q_{k}\right)  $ be the
resulting path tuple, and we set $\sigma^{\prime}:=\sigma t_{i,j}$.

\item The map $f$ transforms $\left(  \sigma,\mathbf{p}\right)  $ into
$\left(  \sigma^{\prime},\mathbf{q}\right)  $.
\end{itemize}

We call this the \textquotedblleft first procedure\textquotedblright.

An analogous procedure computes $f\left(  \sigma^{\prime},\mathbf{q}\right)
=\left(  \sigma^{\prime\prime},\mathbf{r}\right)  $ from $\left(
\sigma^{\prime},\mathbf{q}\right)  $:

\begin{itemize}
\item We say that a point $u$ is \emph{crowded} if it is contained in at least
two of our paths $q_{1},q_{2},\ldots,q_{k}$. (We shall soon see that this
notion of \textquotedblleft crowded\textquotedblright\ agrees with the notion
of \textquotedblleft crowded\textquotedblright\ in the previous procedure.)

\item We pick the smallest $i\in\left[  k\right]  $ such that $q_{i}$ contains
a crowded point. (We shall soon see that this is the same $i$ that was defined
in the previous procedure.)

\item Then, we pick the first crowded point $v$ on $q_{i}$. (We shall soon see
that this is the same $v$ that was defined in the previous procedure.)

\item Then, we pick the largest $j\in\left[  k\right]  $ such that $v$ belongs
to $q_{j}$. (We shall soon see that this is the same $j$ that was defined in
the previous procedure.)

\item Call the part of $q_{i}$ that comes after $v$ the \emph{tail} of $q_{i}%
$. Call the part of $q_{i}$ that comes before $v$ the \emph{head} of $q_{i}$.
Call the part of $q_{j}$ that comes after $v$ the \emph{tail} of $q_{j}$. Call
the part of $q_{j}$ that comes before $v$ the \emph{head} of $q_{j}$.

\item Then, we exchange the tails of the paths $q_{i}$ and $q_{j}$ (while
leaving all other paths unchanged).

\item We let $\mathbf{r}=\left(  r_{1},r_{2},\ldots,r_{k}\right)  $ be the
resulting path tuple, and we set $\sigma^{\prime\prime}:=\sigma^{\prime
}t_{i,j}$.

\item The map $f$ transforms $\left(  \sigma^{\prime},\mathbf{q}\right)  $
into $\left(  \sigma^{\prime\prime},\mathbf{r}\right)  $.
\end{itemize}

We call this the \textquotedblleft second procedure\textquotedblright. We want
to show that this second procedure undoes the first procedure (so that
$\left(  \sigma^{\prime\prime},\mathbf{r}\right)  $ is $\left(  \sigma
,\mathbf{p}\right)  $ again). For this purpose, we make the following observations:

\begin{enumerate}
\item The crowded points as defined in the second procedure are precisely the
crowded points as defined in the first procedure. [\textit{Proof:} When we
exchanged the tails of the paths $p_{i}$ and $p_{j}$, the crowded points
stayed crowded, while the non-crowded points stayed non-crowded, since the
number of our paths through each given point stayed unchanged.]

\item The $i$ defined in the second procedure is precisely the $i$ defined in
the first procedure. [\textit{Proof:} Consider the latter $i$. Then, $i$ is
the smallest $i\in\left[  k\right]  $ such that $p_{i}$ contains a crowded
point\footnote{We can speak of \textquotedblleft crowded
points\textquotedblright\ without having to specify whether we are talking
about the first or the second procedure, because we have shown above that the
crowded points as defined in the second procedure are precisely the crowded
points as defined in the first procedure.}. Thus, the path $p_{i}$ contains a
crowded point -- namely, $v$ -- whereas the paths $p_{1},p_{2},\ldots,p_{i-1}$
don't. Hence, the path $q_{i}$ contains a crowded point (since the definition
of $q_{i}$ ensures that $q_{i}$ still contains $v$), whereas the paths
$q_{1},q_{2},\ldots,q_{i-1}$ don't (since the definition of $\mathbf{q}$ shows
that the paths $q_{1},q_{2},\ldots,q_{i-1}$ are identical with the paths
$p_{1},p_{2},\ldots,p_{i-1}$). Hence, our $i$ is the smallest $i\in\left[
k\right]  $ such that $q_{i}$ contains a crowded point. In other words, our
$i$ is the $i$ defined in the second procedure.]

\item The $v$ defined in the second procedure is precisely the $v$ defined in
the first procedure. [\textit{Proof:} Consider the latter $v$. Then, $v$ is
the first crowded point on $p_{i}$. Therefore, $v$ is crowded, but no point on
$p_{i}$ before $v$ is crowded. Hence, all points on $p_{i}$ before $v$ are
non-crowded. But the path $q_{i}$ agrees with $p_{i}$ up until the point $v$
inclusively (by the definition of $q_{i}$). Hence, the points on $q_{i}$
before $v$ are the same as the points on $p_{i}$ before $v$. Thus, the former
are non-crowded (since the latter are). Thus, $v$ is the first crowded point
on $q_{i}$. In other words, our $v$ is the $v$ defined in the second procedure.]

\item The $j$ defined in the second procedure is precisely the $j$ defined in
the first procedure. [\textit{Proof:} Consider the latter $j$. Then, $j$ is
the largest $j\in\left[  k\right]  $ such that $v$ belongs to $p_{j}$. Hence,
$v$ belongs to $p_{j}$ but does not belong to any of $p_{j+1},p_{j+2}%
,\ldots,p_{k}$. Therefore, $v$ belongs to $q_{j}$ (by the construction of
$q_{j}$) but does not belong to any of $q_{j+1},q_{j+2},\ldots,q_{k}$ (since
the definition of $\mathbf{q}$ shows that the paths $q_{j+1},q_{j+2}%
,\ldots,q_{k}$ are identical to the paths $p_{j+1},p_{j+2},\ldots,p_{k}$). In
other words, our $j$ is the largest $j\in\left[  k\right]  $ such that $v$
belongs to $q_{j}$. In other words, our $j$ is the $j$ defined in the second procedure.]

\item When we exchange the tails of $q_{i}$ and $q_{j}$ in the second
procedure, we are simply undoing the corresponding step of the first
procedure, thus recovering the original paths $p_{i}$ and $p_{j}$.
[\textit{Proof:} We already have shown that the numbers $i$ and $j$ defined in
the second procedure are precisely the numbers $i$ and $j$ defined in the
first procedure. Thus, the tails of $p_{i}$ and $p_{j}$ that we exchanged in
the first procedure are precisely the tails of $q_{j}$ and $q_{i}$ in the
second procedure, and consequently, when we exchange them a second time, they
end up back on their original paths.]

\item We have $\mathbf{r}=\mathbf{p}$. [\textit{Proof:} The tuple $\mathbf{q}$
is obtained from $\mathbf{p}$ by exchanging the tails of $p_{i}$ and $p_{j}$,
whereas the tuple $\mathbf{r}$ is obtained from $\mathbf{q}$ by exchanging the
tails of $q_{i}$ and $q_{j}$. As we just showed, the latter exchange merely
undoes the former exchange. Hence, $\mathbf{r}=\mathbf{p}$.]

\item We have $\sigma^{\prime\prime}=\sigma$. [\textit{Proof:} We already know
that the numbers $i$ and $j$ defined in the second procedure are precisely the
numbers $i$ and $j$ defined in the first procedure. Now, $\sigma^{\prime
\prime}$ was defined by $\sigma^{\prime\prime}=\sigma^{\prime}t_{i,j}$,
whereas $\sigma^{\prime}$ was defined by $\sigma^{\prime}=\sigma t_{i,j}$.
Hence, $\sigma^{\prime\prime}=\underbrace{\sigma^{\prime}}_{=\sigma t_{i,j}%
}t_{i,j}=\sigma\underbrace{t_{i,j}t_{i,j}}_{=t_{i,j}^{2}=\operatorname*{id}%
}=\sigma$.]
\end{enumerate}

Combining $\sigma^{\prime\prime}=\sigma$ with $\mathbf{r}=\mathbf{p}$, we
conclude $\left(  \sigma^{\prime\prime},\mathbf{r}\right)  =\left(
\sigma,\mathbf{p}\right)  $, and this completes our proof of the fact that $f$
is an involution.
\end{proof}
\end{fineprint}

\subsection{\label{sec.details.sf.sp}Definitions and examples of symmetric
polynomials}

\begin{fineprint}
\begin{proof}
[Detailed proof of Lemma \ref{lem.sf.simples-enough}.]Let $\sigma\in S_{N}$.
We shall prove that $\sigma\cdot f=f$.

Let us follow Convention \ref{conv.perm.simple}; thus, we shall refer to the
simple transpositions $s_{1},s_{2},\ldots,s_{N-1}$ in $S_{N}$ as
\textquotedblleft simples\textquotedblright.

Theorem \ref{thm.perm.len.redword1} \textbf{(a)} (applied to $n=N$) shows that
we can write $\sigma$ as a composition (i.e., product) of $\ell\left(
\sigma\right)  $ simples. Thus, in particular, we can write $\sigma$ as a
finite product of simples.\footnote{Alternatively, we can derive this from
Corollary \ref{cor.perm.generated} (which is more well-known than Theorem
\ref{thm.perm.len.redword1} \textbf{(a)}):
\par
Corollary \ref{cor.perm.generated} (applied to $n=N$) shows that the symmetric
group $S_{N}$ is generated by the simples $s_{1},s_{2},\ldots,s_{N-1}$. Hence,
each element of $S_{N}$ is a (finite) product of simples and their inverses.
Since the inverses of the simples are simply these simples themselves (because
each $i\in\left[  N-1\right]  $ satisfies $s_{i}^{-1}=s_{i}$), we can simplify
this statement as follows: Each element of $S_{N}$ is a (finite) product of
simples. Applying this to the element $\sigma$, we conclude that $\sigma$ is a
(finite) product of simples.} In other words, there exist finitely many
elements $k_{1},k_{2},\ldots,k_{p}\in\left[  N-1\right]  $ such that
$\sigma=s_{k_{1}}s_{k_{2}}\cdots s_{k_{p}}$. Consider these $k_{1}%
,k_{2},\ldots,k_{p}$.

Now,%
\begin{align*}
\underbrace{\sigma}_{=s_{k_{1}}s_{k_{2}}\cdots s_{k_{p}}}\cdot\,f  &  =\left(
s_{k_{1}}s_{k_{2}}\cdots s_{k_{p}}\right)  \cdot f=s_{k_{1}}\cdot s_{k_{2}%
}\cdot\cdots\cdot s_{k_{p-1}}\cdot\underbrace{s_{k_{p}}\cdot f}%
_{\substack{=f\\\text{(by (\ref{eq.lem.sf.simples-enough.ass}))}}}\\
&  =s_{k_{1}}\cdot s_{k_{2}}\cdot\cdots\cdot s_{k_{p-2}}\cdot
\underbrace{s_{k_{p-1}}\cdot f}_{\substack{=f\\\text{(by
(\ref{eq.lem.sf.simples-enough.ass}))}}}=s_{k_{1}}\cdot s_{k_{2}}\cdot
\cdots\cdot s_{k_{p-3}}\cdot\underbrace{s_{k_{p-2}}\cdot f}%
_{\substack{=f\\\text{(by (\ref{eq.lem.sf.simples-enough.ass}))}}}\\
&  =\cdots=f.
\end{align*}
(To be fully rigorous, this is really an induction argument: We are showing
(by induction on $i$) that $s_{k_{1}}s_{k_{2}}\cdots s_{k_{i}}f=f$ for each
$i\in\left\{  0,1,\ldots,p\right\}  $; the induction base is obvious (since
$s_{k_{1}}s_{k_{2}}\cdots s_{k_{0}}=\left(  \text{empty product in }%
S_{N}\right)  =\operatorname*{id}$), while the induction step relies on
(\ref{eq.lem.sf.simples-enough.ass}). It is straightforward to fill in the
details of this induction.)

Forget that we fixed $\sigma$. We thus have proved that $\sigma\cdot f=f$ for
all $\sigma\in S_{N}$. In other words, the polynomial $f$ is symmetric. This
proves Lemma \ref{lem.sf.simples-enough}.
\end{proof}
\end{fineprint}

\subsection{\label{sec.details.sf.m}$N$-partitions and monomial symmetric
polynomials}

\begin{fineprint}
\begin{proof}
[Proof of Proposition \ref{prop.sf.sigma-pol-coeff}.]Let us write the
polynomial $f$ as%
\begin{equation}
f=\sum_{\left(  b_{1},b_{2},\ldots,b_{N}\right)  \in\mathbb{N}^{N}}%
f_{b_{1},b_{2},\ldots,b_{N}}x_{1}^{b_{1}}x_{2}^{b_{2}}\cdots x_{N}^{b_{N}},
\label{pf.prop.sf.sigma-pol-coeff.1}%
\end{equation}
where the coefficients $f_{b_{1},b_{2},\ldots,b_{N}}$ belong to $K$. Thus, the
coefficient of any monomial $x_{1}^{b_{1}}x_{2}^{b_{2}}\cdots x_{N}^{b_{N}}$
in $f$ is $f_{b_{1},b_{2},\ldots,b_{N}}$. In other words,%
\begin{equation}
\left[  x_{1}^{b_{1}}x_{2}^{b_{2}}\cdots x_{N}^{b_{N}}\right]  f=f_{b_{1}%
,b_{2},\ldots,b_{N}} \label{pf.prop.sf.sigma-pol-coeff.2}%
\end{equation}
for each $\left(  b_{1},b_{2},\ldots,b_{N}\right)  \in\mathbb{N}^{N}$.

The map $\sigma$ is a permutation of $\left[  N\right]  $ (since $\sigma\in
S_{N}$). Hence, the map%
\begin{align}
\mathbb{N}^{N}  &  \rightarrow\mathbb{N}^{N},\nonumber\\
\left(  b_{1},b_{2},\ldots,b_{N}\right)   &  \mapsto\left(  b_{\sigma\left(
1\right)  },b_{\sigma\left(  2\right)  },\ldots,b_{\sigma\left(  N\right)
}\right)  \label{pf.prop.sf.sigma-pol-coeff.map}%
\end{align}
is a bijection (indeed, this map simply permutes the entries of any given
$N$-tuple using the permutation $\sigma$; thus, it can be undone by permuting
them using $\sigma^{-1}$). Moreover, the map $\sigma$ itself is a bijection
(since it is a permutation).

The definition of the action of $S_{N}$ on $\mathcal{P}$ yields%
\begin{align*}
\sigma\cdot f  &  =f\left[  x_{\sigma\left(  1\right)  },x_{\sigma\left(
2\right)  },\ldots,x_{\sigma\left(  N\right)  }\right] \\
&  =\sum_{\left(  b_{1},b_{2},\ldots,b_{N}\right)  \in\mathbb{N}^{N}}%
f_{b_{1},b_{2},\ldots,b_{N}}x_{\sigma\left(  1\right)  }^{b_{1}}%
x_{\sigma\left(  2\right)  }^{b_{2}}\cdots x_{\sigma\left(  N\right)  }%
^{b_{N}}\\
&  \ \ \ \ \ \ \ \ \ \ \ \ \ \ \ \ \ \ \ \ \left(
\begin{array}
[c]{c}%
\text{here, we have substituted }x_{\sigma\left(  1\right)  },x_{\sigma\left(
2\right)  },\ldots,x_{\sigma\left(  N\right)  }\\
\text{for }x_{1},x_{2},\ldots,x_{N}\text{ on both sides of
(\ref{pf.prop.sf.sigma-pol-coeff.1})}%
\end{array}
\right) \\
&  =\sum_{\left(  b_{1},b_{2},\ldots,b_{N}\right)  \in\mathbb{N}^{N}%
}f_{b_{\sigma\left(  1\right)  },b_{\sigma\left(  2\right)  },\ldots
,b_{\sigma\left(  N\right)  }}x_{\sigma\left(  1\right)  }^{b_{\sigma\left(
1\right)  }}x_{\sigma\left(  2\right)  }^{b_{\sigma\left(  2\right)  }}\cdots
x_{\sigma\left(  N\right)  }^{b_{\sigma\left(  N\right)  }}%
\end{align*}
(here, we have substituted $\left(  b_{\sigma\left(  1\right)  }%
,b_{\sigma\left(  2\right)  },\ldots,b_{\sigma\left(  N\right)  }\right)  $
for the summation index $\left(  b_{1},b_{2},\ldots,b_{N}\right)  $ in the
sum, since the map (\ref{pf.prop.sf.sigma-pol-coeff.map}) is a bijection).
Thus,%
\begin{align*}
\sigma\cdot f  &  =\sum_{\left(  b_{1},b_{2},\ldots,b_{N}\right)
\in\mathbb{N}^{N}}f_{b_{\sigma\left(  1\right)  },b_{\sigma\left(  2\right)
},\ldots,b_{\sigma\left(  N\right)  }}\underbrace{x_{\sigma\left(  1\right)
}^{b_{\sigma\left(  1\right)  }}x_{\sigma\left(  2\right)  }^{b_{\sigma\left(
2\right)  }}\cdots x_{\sigma\left(  N\right)  }^{b_{\sigma\left(  N\right)  }%
}}_{\substack{=\prod_{i\in\left\{  1,2,\ldots,N\right\}  }x_{\sigma\left(
i\right)  }^{b_{\sigma\left(  i\right)  }}=\prod_{i\in\left\{  1,2,\ldots
,N\right\}  }x_{i}^{b_{i}}\\\text{(here, we have substituted }i\text{ for
}\sigma\left(  i\right)  \text{ in the product,}\\\text{since the map }%
\sigma:\left[  N\right]  \rightarrow\left[  N\right]  \text{ is a bijection)}%
}}\\
&  =\sum_{\left(  b_{1},b_{2},\ldots,b_{N}\right)  \in\mathbb{N}^{N}%
}f_{b_{\sigma\left(  1\right)  },b_{\sigma\left(  2\right)  },\ldots
,b_{\sigma\left(  N\right)  }}\underbrace{\prod_{i\in\left\{  1,2,\ldots
,N\right\}  }x_{i}^{b_{i}}}_{=x_{1}^{b_{1}}x_{2}^{b_{2}}\cdots x_{N}^{b_{N}}%
}\\
&  =\sum_{\left(  b_{1},b_{2},\ldots,b_{N}\right)  \in\mathbb{N}^{N}%
}f_{b_{\sigma\left(  1\right)  },b_{\sigma\left(  2\right)  },\ldots
,b_{\sigma\left(  N\right)  }}x_{1}^{b_{1}}x_{2}^{b_{2}}\cdots x_{N}^{b_{N}}.
\end{align*}
Hence, for each $\left(  b_{1},b_{2},\ldots,b_{N}\right)  \in\mathbb{N}^{N}$,
we have%
\begin{align}
f_{b_{\sigma\left(  1\right)  },b_{\sigma\left(  2\right)  },\ldots
,b_{\sigma\left(  N\right)  }}  &  =\left(  \text{the coefficient of }%
x_{1}^{b_{1}}x_{2}^{b_{2}}\cdots x_{N}^{b_{N}}\text{ in }\sigma\cdot f\right)
\nonumber\\
&  =\left[  x_{1}^{b_{1}}x_{2}^{b_{2}}\cdots x_{N}^{b_{N}}\right]  \left(
\sigma\cdot f\right)  . \label{pf.prop.sf.sigma-pol-coeff.5}%
\end{align}

Now, let $\left(  a_{1},a_{2},\ldots,a_{N}\right)  \in\mathbb{N}^{N}$ be
arbitrary. Then, (\ref{pf.prop.sf.sigma-pol-coeff.2}) (applied to $\left(
b_{1},b_{2},\ldots,b_{N}\right)  =\left(  a_{\sigma\left(  1\right)
},a_{\sigma\left(  2\right)  },\ldots,a_{\sigma\left(  N\right)  }\right)  $)
yields%
\[
\left[  x_{1}^{a_{\sigma\left(  1\right)  }}x_{2}^{a_{\sigma\left(  2\right)
}}\cdots x_{N}^{a_{\sigma\left(  N\right)  }}\right]  f=f_{a_{\sigma\left(
1\right)  },a_{\sigma\left(  2\right)  },\ldots,a_{\sigma\left(  N\right)  }%
}=\left[  x_{1}^{a_{1}}x_{2}^{a_{2}}\cdots x_{N}^{a_{N}}\right]  \left(
\sigma\cdot f\right)
\]
(by (\ref{pf.prop.sf.sigma-pol-coeff.5}), applied to $\left(  b_{1}%
,b_{2},\ldots,b_{N}\right)  =\left(  a_{1},a_{2},\ldots,a_{N}\right)  $). This
proves Proposition \ref{prop.sf.sigma-pol-coeff}.
\end{proof}
\end{fineprint}

\subsection{\label{sec.details.sf.schur}Schur polynomials}

\begin{fineprint}
\begin{proof}
[Detailed proof of Lemma \ref{lem.sf.skew-diag.convexity}.]Write the
$M$-partitions $\lambda$ and $\mu$ as $\lambda=\left(  \lambda_{1},\lambda
_{2},\ldots,\lambda_{M}\right)  $ and $\mu=\left(  \mu_{1},\mu_{2},\ldots
,\mu_{M}\right)  $. Then, the definition of $Y\left(  \lambda/\mu\right)  $
yields
\begin{equation}
Y\left(  \lambda/\mu\right)  =\left\{  \left(  i,j\right)  \ \mid\ i\in\left[
M\right]  \text{ and }j\in\mathbb{Z}\text{ and }\mu_{i}<j\leq\lambda
_{i}\right\}  . \label{pf.lem.sf.skew-diag.convexity.Ylm=}%
\end{equation}

We know that $\left(  a,b\right)  $ is an element of $Y\left(  \lambda
/\mu\right)  $. Hence,
\[
\left(  a,b\right)  \in Y\left(  \lambda/\mu\right)  =\left\{  \left(
i,j\right)  \ \mid\ i\in\left[  M\right]  \text{ and }j\in\mathbb{Z}\text{ and
}\mu_{i}<j\leq\lambda_{i}\right\}  .
\]
From this, we obtain $a\in\left[  M\right]  $ and $b\in\mathbb{Z}$ and
$\mu_{a}<b\leq\lambda_{a}$.

We know that $\left(  e,f\right)  $ is an element of $Y\left(  \lambda
/\mu\right)  $. Hence,
\[
\left(  e,f\right)  \in Y\left(  \lambda/\mu\right)  =\left\{  \left(
i,j\right)  \ \mid\ i\in\left[  M\right]  \text{ and }j\in\mathbb{Z}\text{ and
}\mu_{i}<j\leq\lambda_{i}\right\}  .
\]
From this, we obtain $e\in\left[  M\right]  $ and $f\in\mathbb{Z}$ and
$\mu_{e}<f\leq\lambda_{e}$.

Now, from $a\leq c$, we obtain $c\geq a\geq1$ (since $a\in\left[  M\right]
$). Also, we have $c\leq e\leq M$ (since $e\in\left[  M\right]  $). Combining
this with $c\geq1$, we obtain $1\leq c\leq M$, so that $c\in\left[  M\right]
$ (since $c\in\mathbb{Z}$).

Now, $\mu_{1}\geq\mu_{2}\geq\cdots\geq\mu_{M}$ (since $\mu$ is an
$M$-partition). In other words, if $u$ and $v$ are two elements of $\left[
M\right]  $ satisfying $u\leq v$, then $\mu_{u}\geq\mu_{v}$. Applying this to
$u=a$ and $v=c$, we obtain $\mu_{a}\geq\mu_{c}$ (since $a\leq c$). Hence,
$\mu_{c}\leq\mu_{a}<b\leq d$.

Also, $\lambda_{1}\geq\lambda_{2}\geq\cdots\geq\lambda_{M}$ (since $\lambda$
is an $M$-partition). In other words, if $u$ and $v$ are two elements of
$\left[  M\right]  $ satisfying $u\leq v$, then $\lambda_{u}\geq\lambda_{v}$.
Applying this to $u=c$ and $v=e$, we obtain $\lambda_{c}\geq\lambda_{e}$
(since $c\leq e$). Hence, $d\leq f\leq\lambda_{e}\leq\lambda_{c}$ (since
$\lambda_{c}\geq\lambda_{e}$).

Combining this with $\mu_{c}<d$, we obtain $\mu_{c}<d\leq\lambda_{c}$.

Thus, we know that $c\in\left[  M\right]  $ and $d\in\mathbb{Z}$ and $\mu
_{c}<d\leq\lambda_{c}$. In other words, $\left(  c,d\right)  \in\left\{
\left(  i,j\right)  \ \mid\ i\in\left[  M\right]  \text{ and }j\in
\mathbb{Z}\text{ and }\mu_{i}<j\leq\lambda_{i}\right\}  $. In view of
(\ref{pf.lem.sf.skew-diag.convexity.Ylm=}), this rewrites as $\left(
c,d\right)  \in Y\left(  \lambda/\mu\right)  $. This proves Lemma
\ref{lem.sf.skew-diag.convexity}.
\end{proof}

\begin{proof}
[Detailed proof of Lemma \ref{lem.sf.skew-ssyt.increase}.]We shall first prove
parts \textbf{(a)} and \textbf{(c)}, and then quickly derive the rest from
them. \medskip

\textbf{(a)} The skew tableau $T$ is semistandard. Hence, we have%
\begin{equation}
T\left(  i,j\right)  \leq T\left(  i,j+1\right)
\label{pf.lem.sf.skew-ssyt.increase.a.1}%
\end{equation}
for any $\left(  i,j\right)  \in Y\left(  \lambda/\mu\right)  $ satisfying
$\left(  i,j+1\right)  \in Y\left(  \lambda/\mu\right)  $. (Indeed, this is
one of the requirements placed on $T$ in Definition \ref{def.sf.skew-ssyt}.)

Now, let $\left(  i,j_{1}\right)  $ and $\left(  i,j_{2}\right)  $ be two
elements of $Y\left(  \lambda/\mu\right)  $ satisfying $j_{1}\leq j_{2}$. We
must prove that $T\left(  i,j_{1}\right)  \leq T\left(  i,j_{2}\right)  $.

Let $k\in\left\{  j_{1},j_{1}+1,\ldots,j_{2}-1\right\}  $ be arbitrary. Then,
$j_{1}\leq k\leq j_{2}-1$. From $k\leq j_{2}-1$, we obtain $k+1\leq j_{2}$, so
that $k\leq k+1\leq j_{2}$. Also, $j_{1}\leq k\leq k+1$.

Thus, we have $i\leq i\leq i$ and $j_{1}\leq k\leq j_{2}$. Hence, Lemma
\ref{lem.sf.skew-diag.convexity} (applied to $\left(  a,b\right)  =\left(
i,j_{1}\right)  $ and $\left(  e,f\right)  =\left(  i,j_{2}\right)  $ and
$\left(  c,d\right)  =\left(  i,k\right)  $) yields $\left(  i,k\right)  \in
Y\left(  \lambda/\mu\right)  $. Therefore, the entry $T\left(  i,k\right)  $
of $T$ is well-defined.

Also, we have $i\leq i\leq i$ and $j_{1}\leq k+1\leq j_{2}$. Hence, Lemma
\ref{lem.sf.skew-diag.convexity} (applied to $\left(  a,b\right)  =\left(
i,j_{1}\right)  $ and $\left(  e,f\right)  =\left(  i,j_{2}\right)  $ and
$\left(  c,d\right)  =\left(  i,k+1\right)  $) yields $\left(  i,k+1\right)
\in Y\left(  \lambda/\mu\right)  $. Therefore, the entry $T\left(
i,k+1\right)  $ of $T$ is well-defined.

Now, (\ref{pf.lem.sf.skew-ssyt.increase.a.1}) (applied to $j=k$) yields
$T\left(  i,k\right)  \leq T\left(  i,k+1\right)  $ (since $\left(
i,k\right)  \in Y\left(  \lambda/\mu\right)  $ and $\left(  i,k+1\right)  \in
Y\left(  \lambda/\mu\right)  $).

Forget that we fixed $k$. We thus have shown that for each $k\in\left\{
j_{1},j_{1}+1,\ldots,j_{2}-1\right\}  $, the inequality $T\left(  i,k\right)
\leq T\left(  i,k+1\right)  $ holds (and both entries $T\left(  i,k\right)  $
and $T\left(  i,k+1\right)  $ are well-defined). In other words, we have%
\[
T\left(  i,j_{1}\right)  \leq T\left(  i,j_{1}+1\right)  \leq T\left(
i,j_{1}+2\right)  \leq\cdots\leq T\left(  i,j_{2}-1\right)  \leq T\left(
i,j_{2}\right)  .
\]
Hence, $T\left(  i,j_{1}\right)  \leq T\left(  i,j_{2}\right)  $. This proves
Lemma \ref{lem.sf.skew-ssyt.increase} \textbf{(a)}. \medskip

\textbf{(c)} The skew tableau $T$ is semistandard. Hence, we have%
\begin{equation}
T\left(  i,j\right)  <T\left(  i+1,j\right)
\label{pf.lem.sf.skew-ssyt.increase.c.1}%
\end{equation}
for any $\left(  i,j\right)  \in Y\left(  \lambda/\mu\right)  $ satisfying
$\left(  i+1,j\right)  \in Y\left(  \lambda/\mu\right)  $. (Indeed, this is
one of the requirements placed on $T$ in Definition \ref{def.sf.skew-ssyt}.)

Now, let $\left(  i_{1},j\right)  $ and $\left(  i_{2},j\right)  $ be two
elements of $Y\left(  \lambda/\mu\right)  $ satisfying $i_{1}<i_{2}$. We must
prove that $T\left(  i_{1},j\right)  <T\left(  i_{2},j\right)  $.

Let $k\in\left\{  i_{1},i_{1}+1,\ldots,i_{2}-1\right\}  $ be arbitrary. Then,
$i_{1}\leq k\leq i_{2}-1$. From $k\leq i_{2}-1$, we obtain $k+1\leq i_{2}$, so
that $k\leq k+1\leq i_{2}$. Also, $i_{1}\leq k\leq k+1$.

Thus, we have $i_{1}\leq k\leq i_{2}$ and $j\leq j\leq j$. Hence, Lemma
\ref{lem.sf.skew-diag.convexity} (applied to $\left(  a,b\right)  =\left(
i_{1},j\right)  $ and $\left(  e,f\right)  =\left(  i_{2},j\right)  $ and
$\left(  c,d\right)  =\left(  k,j\right)  $) yields $\left(  k,j\right)  \in
Y\left(  \lambda/\mu\right)  $. Therefore, the entry $T\left(  k,j\right)  $
of $T$ is well-defined.

Also, we have $i_{1}\leq k+1\leq i_{2}$ and $j\leq j\leq j$. Hence, Lemma
\ref{lem.sf.skew-diag.convexity} (applied to $\left(  a,b\right)  =\left(
i_{1},j\right)  $ and $\left(  e,f\right)  =\left(  i_{2},j\right)  $ and
$\left(  c,d\right)  =\left(  k+1,j\right)  $) yields $\left(  k+1,j\right)
\in Y\left(  \lambda/\mu\right)  $. Therefore, the entry $T\left(
k+1,j\right)  $ of $T$ is well-defined.

Now, (\ref{pf.lem.sf.skew-ssyt.increase.c.1}) (applied to $i=k$) yields
$T\left(  k,j\right)  <T\left(  k+1,j\right)  $ (since $\left(  k,j\right)
\in Y\left(  \lambda/\mu\right)  $ and $\left(  k+1,j\right)  \in Y\left(
\lambda/\mu\right)  $).

Forget that we fixed $k$. We thus have shown that for each $k\in\left\{
i_{1},i_{1}+1,\ldots,i_{2}-1\right\}  $, the inequality $T\left(  k,j\right)
<T\left(  k+1,j\right)  $ holds (and both entries $T\left(  k,j\right)  $ and
$T\left(  k+1,j\right)  $ are well-defined). In other words, we have%
\[
T\left(  i_{1},j\right)  <T\left(  i_{1}+1,j\right)  <T\left(  i_{1}%
+2,j\right)  <\cdots<T\left(  i_{2}-1,j\right)  <T\left(  i_{2},j\right)  .
\]
Hence, $T\left(  i_{1},j\right)  <T\left(  i_{2},j\right)  $ (since
$i_{1}<i_{2}$). This proves Lemma \ref{lem.sf.skew-ssyt.increase}
\textbf{(c)}. \medskip

\textbf{(b)} Let $\left(  i_{1},j\right)  $ and $\left(  i_{2},j\right)  $ be
two elements of $Y\left(  \lambda/\mu\right)  $ satisfying $i_{1}\leq i_{2}$.
We must prove that $T\left(  i_{1},j\right)  \leq T\left(  i_{2},j\right)  $.
If $i_{1}<i_{2}$, then this follows from Lemma \ref{lem.sf.skew-ssyt.increase}
\textbf{(c)}. Hence, for the rest of this proof, we WLOG assume that we don't
have $i_{1}<i_{2}$. Thus, we have $i_{1}\geq i_{2}$. Combining this with
$i_{1}\leq i_{2}$, we obtain $i_{1}=i_{2}$. Thus, $T\left(  i_{1},j\right)
=T\left(  i_{2},j\right)  $, so that $T\left(  i_{1},j\right)  \leq T\left(
i_{2},j\right)  $. This proves Lemma \ref{lem.sf.skew-ssyt.increase}
\textbf{(b)}. \medskip

\textbf{(d)} Let $\left(  i_{1},j_{1}\right)  $ and $\left(  i_{2}%
,j_{2}\right)  $ be two elements of $Y\left(  \lambda/\mu\right)  $ satisfying
$i_{1}\leq i_{2}$ and $j_{1}\leq j_{2}$. Then, $i_{1}\leq i_{2}\leq i_{2}$ and
$j_{1}\leq j_{1}\leq j_{2}$. Hence, Lemma \ref{lem.sf.skew-diag.convexity}
(applied to $\left(  a,b\right)  =\left(  i_{1},j_{1}\right)  $ and $\left(
e,f\right)  =\left(  i_{2},j_{2}\right)  $ and $\left(  c,d\right)  =\left(
i_{2},j_{1}\right)  $) yields $\left(  i_{2},j_{1}\right)  \in Y\left(
\lambda/\mu\right)  $. Therefore, the entry $T\left(  i_{2},j_{1}\right)  $ of
$T$ is well-defined. Hence, Lemma \ref{lem.sf.skew-ssyt.increase} \textbf{(b)}
(applied to $j=j_{1}$) yields $T\left(  i_{1},j_{1}\right)  \leq T\left(
i_{2},j_{1}\right)  $. Furthermore, Lemma \ref{lem.sf.skew-ssyt.increase}
\textbf{(a)} (applied to $i=i_{2}$) yields $T\left(  i_{2},j_{1}\right)  \leq
T\left(  i_{2},j_{2}\right)  $. Thus,%
\[
T\left(  i_{1},j_{1}\right)  \leq T\left(  i_{2},j_{1}\right)  \leq T\left(
i_{2},j_{2}\right)  .
\]
This proves Lemma \ref{lem.sf.skew-ssyt.increase} \textbf{(d)}. \medskip

\textbf{(e)} Let $\left(  i_{1},j_{1}\right)  $ and $\left(  i_{2}%
,j_{2}\right)  $ be two elements of $Y\left(  \lambda/\mu\right)  $ satisfying
$i_{1}<i_{2}$ and $j_{1}\leq j_{2}$. Then, $i_{1}\leq i_{2}\leq i_{2}$ and
$j_{1}\leq j_{1}\leq j_{2}$. Hence, Lemma \ref{lem.sf.skew-diag.convexity}
(applied to $\left(  a,b\right)  =\left(  i_{1},j_{1}\right)  $ and $\left(
e,f\right)  =\left(  i_{2},j_{2}\right)  $ and $\left(  c,d\right)  =\left(
i_{2},j_{1}\right)  $) yields $\left(  i_{2},j_{1}\right)  \in Y\left(
\lambda/\mu\right)  $. Therefore, the entry $T\left(  i_{2},j_{1}\right)  $ of
$T$ is well-defined. Hence, Lemma \ref{lem.sf.skew-ssyt.increase} \textbf{(c)}
(applied to $j=j_{1}$) yields $T\left(  i_{1},j_{1}\right)  <T\left(
i_{2},j_{1}\right)  $. Furthermore, Lemma \ref{lem.sf.skew-ssyt.increase}
\textbf{(a)} (applied to $i=i_{2}$) yields $T\left(  i_{2},j_{1}\right)  \leq
T\left(  i_{2},j_{2}\right)  $. Thus,%
\[
T\left(  i_{1},j_{1}\right)  <T\left(  i_{2},j_{1}\right)  \leq T\left(
i_{2},j_{2}\right)  .
\]
This proves Lemma \ref{lem.sf.skew-ssyt.increase} \textbf{(e)}.
\end{proof}
\end{fineprint}

\begin{fineprint}
\begin{proof}
[Detailed proof of Lemma \ref{lem.sf.tab-greater-i}.]Let $\left(  i,j\right)
\in Y\left(  \lambda\right)  $. Set $p:=T\left(  i,j\right)  $.

All entries of $T$ are elements of $\left[  N\right]  $ (by the definition of
a tableau), and thus are positive integers. Thus, in particular, the $i$
entries $T\left(  1,j\right)  ,T\left(  2,j\right)  ,\ldots,T\left(
i,j\right)  $ are positive integers\footnote{Here, we are tacitly using the
fact that the boxes $\left(  1,j\right)  ,\left(  2,j\right)  ,\ldots,\left(
i,j\right)  $ all belong to $Y\left(  \lambda\right)  $ (so that the
corresponding entries $T\left(  1,j\right)  ,T\left(  2,j\right)
,\ldots,T\left(  i,j\right)  $ are well-defined). This fact can be checked as
follows: Let $u\in\left[  i\right]  $. Thus, $u\leq i$. Now, write $\lambda$
in the form $\lambda=\left(  \lambda_{1},\lambda_{2},\ldots,\lambda
_{N}\right)  $. Thus, $\lambda_{1}\geq\lambda_{2}\geq\cdots\geq\lambda_{N}$
(since $\lambda$ is an $N$-partition). Hence, $\lambda_{u}\geq\lambda_{i}$
(since $u\leq i$), so that $\lambda_{i}\leq\lambda_{u}$ and therefore $\left[
\lambda_{i}\right]  \subseteq\left[  \lambda_{u}\right]  $. However, $\left(
i,j\right)  \in Y\left(  \lambda\right)  $. In other words, $i\in\left[
N\right]  $ and $j\in\left[  \lambda_{i}\right]  $ (by the definition of the
Young diagram $Y\left(  \lambda\right)  $). Hence, $u\leq i\leq N$ (since
$i\in\left[  N\right]  $), so that $u\in\left[  N\right]  $, and furthermore
$j\in\left[  \lambda_{i}\right]  \subseteq\left[  \lambda_{u}\right]  $. Now,
from $u\in\left[  N\right]  $ and $j\in\left[  \lambda_{u}\right]  $, we
obtain $\left(  u,j\right)  \in Y\left(  \lambda\right)  $. Forget that we
fixed $u$. We thus have shown that $\left(  u,j\right)  \in Y\left(
\lambda\right)  $ for each $u\in\left[  i\right]  $. In other words, the boxes
$\left(  1,j\right)  ,\left(  2,j\right)  ,\ldots,\left(  i,j\right)  $ all
belong to $Y\left(  \lambda\right)  $, qed.}. Moreover, the tableau $T$ is
semistandard; thus, its entries increase strictly down each column. Hence, in
particular, we have%
\[
T\left(  1,j\right)  <T\left(  2,j\right)  <\cdots<T\left(  i,j\right)  .
\]
Thus, the $i$ numbers $T\left(  1,j\right)  ,T\left(  2,j\right)
,\ldots,T\left(  i,j\right)  $ are distinct. Moreover, all these numbers are
positive integers (as we have seen above) and are $\leq p$ (since $T\left(
1,j\right)  <T\left(  2,j\right)  <\cdots<T\left(  i,j\right)  =p$); thus,
they all belong to the set $\left[  p\right]  $. This shows that there are $i$
distinct numbers in the set $\left[  p\right]  $ (namely, the $i$ numbers
$T\left(  1,j\right)  ,T\left(  2,j\right)  ,\ldots,T\left(  i,j\right)  $);
in other words, the set $\left[  p\right]  $ has at least $i$ elements. In
other words, we have $\left\vert \left[  p\right]  \right\vert \geq i$. Since
$\left\vert \left[  p\right]  \right\vert =p=T\left(  i,j\right)  $, this
rewrites as $T\left(  i,j\right)  \geq i$. This proves Lemma
\ref{lem.sf.tab-greater-i}.
\end{proof}

\begin{proof}
[Detailed proof of Lemma \ref{lem.sf.alternant-0}.]Write the $N$-tuple
$\alpha\in\mathbb{N}^{N}$ as $\alpha=\left(  \alpha_{1},\alpha_{2}%
,\ldots,\alpha_{N}\right)  $. Then, Definition \ref{def.sf.alternants}
\textbf{(b)} yields
\begin{equation}
a_{\alpha}=\det\left(  \left(  x_{i}^{\alpha_{j}}\right)  _{1\leq i\leq
N,\ 1\leq j\leq N}\right)  . \label{pf.lem.sf.alternant-0.aalpha=}%
\end{equation}

\textbf{(a)} Assume that the $N$-tuple $\alpha$ has two equal entries. In
other words, the $N$-tuple $\left(  \alpha_{1},\alpha_{2},\ldots,\alpha
_{N}\right)  $ has two equal entries (since $\alpha=\left(  \alpha_{1}%
,\alpha_{2},\ldots,\alpha_{N}\right)  $). In other words, there exist two
elements $u,v\in\left[  N\right]  $ such that $u<v$ and $\alpha_{u}=\alpha
_{v}$. Consider these $u,v$. Now, from $\alpha_{u}=\alpha_{v}$, we conclude
that the $u$-th and the $v$-th columns of the $N\times N$-matrix $\left(
x_{i}^{\alpha_{j}}\right)  _{1\leq i\leq N,\ 1\leq j\leq N}$ are equal. Hence,
this $N\times N$-matrix $\left(  x_{i}^{\alpha_{j}}\right)  _{1\leq i\leq
N,\ 1\leq j\leq N}$ has two equal columns (since $u<v$).

However, if an $N\times N$-matrix $A$ has two equal columns, then $\det A=0$
(by Theorem \ref{thm.det.colop} \textbf{(c)}\footnote{more precisely: by the
analogue of Theorem \ref{thm.det.rowop} \textbf{(c)} for columns instead of
rows}). Applying this to $A=\left(  x_{i}^{\alpha_{j}}\right)  _{1\leq i\leq
N,\ 1\leq j\leq N}$, we obtain%
\[
\det\left(  \left(  x_{i}^{\alpha_{j}}\right)  _{1\leq i\leq N,\ 1\leq j\leq
N}\right)  =0
\]
(since the $N\times N$-matrix $\left(  x_{i}^{\alpha_{j}}\right)  _{1\leq
i\leq N,\ 1\leq j\leq N}$ has two equal columns). In view of
(\ref{pf.lem.sf.alternant-0.aalpha=}), this rewrites as $a_{\alpha}=0$. This
proves Lemma \ref{lem.sf.alternant-0} \textbf{(a)}. \medskip

\textbf{(b)} Write the $N$-tuple $\beta\in\mathbb{N}^{N}$ as $\beta=\left(
\beta_{1},\beta_{2},\ldots,\beta_{N}\right)  $. Then, Definition
\ref{def.sf.alternants} \textbf{(b)} yields
\begin{equation}
a_{\beta}=\det\left(  \left(  x_{i}^{\beta_{j}}\right)  _{1\leq i\leq
N,\ 1\leq j\leq N}\right)  . \label{pf.lem.sf.alternant-0.abeta=}%
\end{equation}

However, the $N$-tuple $\beta$ is obtained from $\alpha$ by swapping two
entries. In other words, the $N$-tuple $\left(  \beta_{1},\beta_{2}%
,\ldots,\beta_{N}\right)  $ is obtained from $\left(  \alpha_{1},\alpha
_{2},\ldots,\alpha_{N}\right)  $ by swapping two entries (since $\beta=\left(
\beta_{1},\beta_{2},\ldots,\beta_{N}\right)  $ and $\alpha=\left(  \alpha
_{1},\alpha_{2},\ldots,\alpha_{N}\right)  $). Thus, the matrix $\left(
x_{i}^{\beta_{j}}\right)  _{1\leq i\leq N,\ 1\leq j\leq N}$ is obtained from
the matrix $\left(  x_{i}^{\alpha_{j}}\right)  _{1\leq i\leq N,\ 1\leq j\leq
N}$ by swapping two columns\footnote{Indeed, if the $N$-tuple $\left(
\beta_{1},\beta_{2},\ldots,\beta_{N}\right)  $ is obtained from $\left(
\alpha_{1},\alpha_{2},\ldots,\alpha_{N}\right)  $ by swapping the $u$-th and
the $v$-th entries, then the matrix $\left(  x_{i}^{\beta_{j}}\right)  _{1\leq
i\leq N,\ 1\leq j\leq N}$ is obtained from the matrix $\left(  x_{i}%
^{\alpha_{j}}\right)  _{1\leq i\leq N,\ 1\leq j\leq N}$ by swapping the $u$-th
and the $v$-th columns.}.

However, if we swap two columns of an $N\times N$-matrix $A$, then $\det A$
gets multiplied by $-1$ (by Theorem \ref{thm.det.colop} \textbf{(b)}%
\footnote{more precisely: by the analogue of Theorem \ref{thm.det.rowop}
\textbf{(b)} for columns instead of rows}). In other words, if $A$ and $B$ are
two $N\times N$-matrices such that $B$ is obtained from $A$ by swapping two
columns, then $\det B=-\det A$. Applying this to $A=\left(  x_{i}^{\alpha_{j}%
}\right)  _{1\leq i\leq N,\ 1\leq j\leq N}$ and $B=\left(  x_{i}^{\beta_{j}%
}\right)  _{1\leq i\leq N,\ 1\leq j\leq N}$, we obtain%
\[
\det\left(  \left(  x_{i}^{\beta_{j}}\right)  _{1\leq i\leq N,\ 1\leq j\leq
N}\right)  =-\det\left(  \left(  x_{i}^{\alpha_{j}}\right)  _{1\leq i\leq
N,\ 1\leq j\leq N}\right)
\]
(since the $N\times N$-matrix $\left(  x_{i}^{\beta_{j}}\right)  _{1\leq i\leq
N,\ 1\leq j\leq N}$ is obtained from the matrix $\left(  x_{i}^{\alpha_{j}%
}\right)  _{1\leq i\leq N,\ 1\leq j\leq N}$ by swapping two columns). In view
of (\ref{pf.lem.sf.alternant-0.aalpha=}) and
(\ref{pf.lem.sf.alternant-0.abeta=}), this rewrites as $a_{\beta}=-a_{\alpha}%
$. This proves Lemma \ref{lem.sf.alternant-0} \textbf{(b)}.
\end{proof}

\begin{proof}
[Some details omitted from the proof of Lemma \ref{lem.sf.stemb-lem}.]In our
above proof of Lemma \ref{lem.sf.stemb-lem}, we have omitted certain arguments
-- namely, the proofs of the equalities
(\ref{pf.lem.stemb-lem.o2.pf.gammak+1=alk}) and
(\ref{pf.lem.stemb-lem.o2.pf.gammai=ali}) in the proof of Observation 2. Let
us now show these proofs:\footnote{In both of these proofs, we will use the
notations that were introduced in the proof of Observation 2 in the proof of
Lemma \ref{lem.sf.stemb-lem}.}

[\textit{Proof of (\ref{pf.lem.stemb-lem.o2.pf.gammak+1=alk}):} The definition
of $\operatorname*{cont}\left(  T^{\ast}\right)  $ yields
\begin{align*}
\left(  \operatorname*{cont}\left(  T^{\ast}\right)  \right)  _{k+1}  &
=\left(  \text{\# of }\left(  k+1\right)  \text{'s in }T^{\ast}\right) \\
&  =\left(  \text{\# of }\left(  k+1\right)  \text{'s in }%
\underbrace{\operatorname{col}_{<j}\left(  T^{\ast}\right)  }%
_{\substack{=\beta_{k}\left(  \operatorname{col}_{<j}T\right)  \\\text{(by
(\ref{pf.lem.stemb-lem.T*1}))}}}\right)  +\left(  \text{\# of }\left(
k+1\right)  \text{'s in }\underbrace{\operatorname{col}_{\geq j}\left(
T^{\ast}\right)  }_{\substack{=\operatorname{col}_{\geq j}T\\\text{(by
(\ref{pf.lem.stemb-lem.T*2}))}}}\right) \\
&  \ \ \ \ \ \ \ \ \ \ \ \ \ \ \ \ \ \ \ \ \left(  \text{by
(\ref{pf.lem.stemb-lem.o2.pf.numiinT*}), applied to }i=k+1\right) \\
&  =\underbrace{\left(  \text{\# of }\left(  k+1\right)  \text{'s in }%
\beta_{k}\left(  \operatorname{col}_{<j}T\right)  \right)  }%
_{\substack{=\left(  \text{\# of }k\text{'s in }\operatorname{col}%
_{<j}T\right)  \\\text{(by (\ref{pf.thm.sf.skew-schur-symm.num-k+1}), applied
to }\operatorname{col}_{<j}T\text{ instead of }T\text{)}}}+\underbrace{\left(
\text{\# of }\left(  k+1\right)  \text{'s in }\operatorname{col}_{\geq
j}T\right)  }_{\substack{=b_{k+1}-\nu_{k+1}\\\text{(by
(\ref{pf.lem.stemb-lem.o2.pf.bk+1=}))}}}\\
&  =\left(  \text{\# of }k\text{'s in }\operatorname{col}_{<j}T\right)
+\underbrace{b_{k+1}}_{\substack{=b_{k}+1\\\text{(since }b_{k}+1=b_{k+1}%
\text{)}}}-\,\nu_{k+1}\\
&  =\left(  \text{\# of }k\text{'s in }\operatorname{col}_{<j}T\right)
+b_{k}+1-\nu_{k+1}.
\end{align*}
However, $\gamma=\nu+\operatorname*{cont}\left(  T^{\ast}\right)  +\rho$, so
that%
\begin{align}
\gamma_{k+1}  &  =\left(  \nu+\operatorname*{cont}\left(  T^{\ast}\right)
+\rho\right)  _{k+1}\nonumber\\
&  =\nu_{k+1}+\underbrace{\left(  \operatorname*{cont}\left(  T^{\ast}\right)
\right)  _{k+1}}_{=\left(  \text{\# of }k\text{'s in }\operatorname{col}%
_{<j}T\right)  +b_{k}+1-\nu_{k+1}}+\underbrace{\rho_{k+1}}%
_{\substack{=N-\left(  k+1\right)  \\\text{(by the definition of }\rho
\text{)}}}\nonumber\\
&  =\nu_{k+1}+\left(  \text{\# of }k\text{'s in }\operatorname{col}%
_{<j}T\right)  +b_{k}+1-\nu_{k+1}+N-\left(  k+1\right) \nonumber\\
&  =\left(  \text{\# of }k\text{'s in }\operatorname{col}_{<j}T\right)
+b_{k}+N-k. \label{pf.lem.stemb-lem.o2.pf.gammak+1=}%
\end{align}

On the other hand, the definition of $\operatorname*{cont}T$ yields
\begin{align*}
\left(  \operatorname*{cont}T\right)  _{k}  &  =\left(  \text{\# of }k\text{'s
in }T\right) \\
&  =\left(  \text{\# of }k\text{'s in }\operatorname{col}_{<j}T\right)
+\underbrace{\left(  \text{\# of }k\text{'s in }\operatorname{col}_{\geq
j}T\right)  }_{\substack{=b_{k}-\nu_{k}\\\text{(by
(\ref{pf.lem.stemb-lem.o2.pf.bk=}))}}}\\
&  \ \ \ \ \ \ \ \ \ \ \ \ \ \ \ \ \ \ \ \ \left(  \text{by
(\ref{pf.lem.stemb-lem.o2.pf.numiinT}), applied to }i=k\right) \\
&  =\left(  \text{\# of }k\text{'s in }\operatorname{col}_{<j}T\right)
+b_{k}-\nu_{k}.
\end{align*}
However, $\alpha=\nu+\operatorname*{cont}T+\rho$, so that%
\begin{align*}
\alpha_{k}  &  =\left(  \nu+\operatorname*{cont}T+\rho\right)  _{k}\\
&  =\nu_{k}+\underbrace{\left(  \operatorname*{cont}T\right)  _{k}}_{=\left(
\text{\# of }k\text{'s in }\operatorname{col}_{<j}T\right)  +b_{k}-\nu_{k}%
}+\underbrace{\rho_{k}}_{\substack{=N-k\\\text{(by the definition of }%
\rho\text{)}}}\\
&  =\nu_{k}+\left(  \text{\# of }k\text{'s in }\operatorname{col}%
_{<j}T\right)  +b_{k}-\nu_{k}+N-k\\
&  =\left(  \text{\# of }k\text{'s in }\operatorname{col}_{<j}T\right)
+b_{k}+N-k.
\end{align*}
Comparing this with (\ref{pf.lem.stemb-lem.o2.pf.gammak+1=}), we obtain
$\gamma_{k+1}=\alpha_{k}$. This proves
(\ref{pf.lem.stemb-lem.o2.pf.gammak+1=alk}).]

[\textit{Proof of (\ref{pf.lem.stemb-lem.o2.pf.gammai=ali}):} Let $i\in\left[
N\right]  $ be such that $i\neq k$ and $i\neq k+1$. We must prove that
$\gamma_{i}=\alpha_{i}$.

The definition of $\operatorname*{cont}\left(  T^{\ast}\right)  $ yields
\begin{align*}
\left(  \operatorname*{cont}\left(  T^{\ast}\right)  \right)  _{i}  &
=\left(  \text{\# of }i\text{'s in }T^{\ast}\right) \\
&  =\left(  \text{\# of }i\text{'s in }\underbrace{\operatorname{col}%
_{<j}\left(  T^{\ast}\right)  }_{\substack{=\beta_{k}\left(
\operatorname{col}_{<j}T\right)  \\\text{(by (\ref{pf.lem.stemb-lem.T*1}))}%
}}\right)  +\left(  \text{\# of }i\text{'s in }\underbrace{\operatorname{col}%
_{\geq j}\left(  T^{\ast}\right)  }_{\substack{=\operatorname{col}_{\geq
j}T\\\text{(by (\ref{pf.lem.stemb-lem.T*2}))}}}\right)
\ \ \ \ \ \ \ \ \ \ \left(  \text{by (\ref{pf.lem.stemb-lem.o2.pf.numiinT*}%
)}\right) \\
&  =\underbrace{\left(  \text{\# of }i\text{'s in }\beta_{k}\left(
\operatorname{col}_{<j}T\right)  \right)  }_{\substack{=\left(  \text{\# of
}i\text{'s in }\operatorname{col}_{<j}T\right)  \\\text{(by
(\ref{pf.thm.sf.skew-schur-symm.num-i}),}\\\text{applied to }%
\operatorname{col}_{<j}T\text{ instead of }T\text{)}}}+\left(  \text{\# of
}i\text{'s in }\operatorname{col}_{\geq j}T\right) \\
&  =\left(  \text{\# of }i\text{'s in }\operatorname{col}_{<j}T\right)
+\left(  \text{\# of }i\text{'s in }\operatorname{col}_{\geq j}T\right)  .
\end{align*}
On the other hand, the definition of $\operatorname*{cont}T$ yields
\[
\left(  \operatorname*{cont}T\right)  _{i}=\left(  \text{\# of }i\text{'s in
}T\right)  =\left(  \text{\# of }i\text{'s in }\operatorname{col}%
_{<j}T\right)  +\left(  \text{\# of }i\text{'s in }\operatorname{col}_{\geq
j}T\right)
\]
(by (\ref{pf.lem.stemb-lem.o2.pf.numiinT})). Comparing these two equalities,
we obtain $\left(  \operatorname*{cont}\left(  T^{\ast}\right)  \right)
_{i}=\left(  \operatorname*{cont}T\right)  _{i}$.

However, $\gamma=\nu+\operatorname*{cont}\left(  T^{\ast}\right)  +\rho$, so
that%
\[
\gamma_{i}=\left(  \nu+\operatorname*{cont}\left(  T^{\ast}\right)
+\rho\right)  _{i}=\nu_{i}+\underbrace{\left(  \operatorname*{cont}\left(
T^{\ast}\right)  \right)  _{i}}_{=\left(  \operatorname*{cont}T\right)  _{i}%
}+\,\rho_{i}=\nu_{i}+\left(  \operatorname*{cont}T\right)  _{i}+\rho_{i}.
\]

However, $\alpha=\nu+\operatorname*{cont}T+\rho$, so that%
\[
\alpha_{i}=\left(  \nu+\operatorname*{cont}T+\rho\right)  _{i}=\nu_{i}+\left(
\operatorname*{cont}T\right)  _{i}+\rho_{i}.
\]
Comparing these two equalities, we obtain $\gamma_{i}=\alpha_{i}$. This proves
(\ref{pf.lem.stemb-lem.o2.pf.gammai=ali}).]
\end{proof}
\end{fineprint}


\bigskip

\bigskip

\bigskip

\begin{thebibliography}{999999999}                                                                                        %


\bibitem[17f-hw7s]{17f-hw7s}Darij Grinberg, \textit{UMN Fall 2017 Math 4990
homework set \#7 with solutions}, \url{https://darijgrinberg.gitlab.io/t/17f/hw7os.pdf}

\bibitem[18f-hw2s]{18f-hw2s}Darij Grinberg, \textit{UMN Fall 2018 Math 5705
homework set \#2 with solutions}, \url{https://darijgrinberg.gitlab.io/t/18f/hw2s.pdf}

\bibitem[18f-hw4s]{18f-hw4s}Darij Grinberg, \textit{UMN Fall 2018 Math 5705
homework set \#4 with solutions}, \url{https://darijgrinberg.gitlab.io/t/18f/hw4s.pdf}

\bibitem[18f-hw4se]{18f-hw4se}Jacob Elafandi, \textit{Math 5705: Enumerative
Combinatorics, Fall 2018: Homework 4: solutions to exercises 1, 2, 3,
7}.\newline\url{https://darijgrinberg.gitlab.io/t/18f/hw4s-elafandi.pdf}

\bibitem[18f-mt3s]{18f-mt3s}Darij Grinberg, \textit{Math 5705: Enumerative
Combinatorics, Fall 2018: Midterm 3 with solutions}.\newline\url{https://darijgrinberg.gitlab.io/t/18f/mt3s.pdf}

\bibitem[19fla]{19fla}Darij Grinberg, \textit{Math 201-003: Linear Algebra,
Fall 2019}.\newline\url{https://darijgrinberg.gitlab.io/t/19fla/}

\bibitem[19s]{19s}Darij Grinberg, \textit{Introduction to Modern Algebra (UMN
Spring 2019 Math 4281 notes)}, 29 June 2019.\newline\url{https://darijgrinberg.gitlab.io/t/19s/notes.pdf}

\bibitem[19s-mt3s]{19s-mt3s}Darij Grinberg, \textit{Math 4281: Introduction to
Modern Algebra, Spring 2019: Midterm 3 with solutions}.\newline\url{https://darijgrinberg.gitlab.io/t/19s/mt3s.pdf}

\bibitem[19fco]{19fco}Darij Grinberg, \textit{Enumerative Combinatorics: class
notes (Drexel Fall 2019 Math 222 notes)}, 11 September 2022.\newline\url{https://darijgrinberg.gitlab.io/t/19fco/n/n.pdf}

\bibitem[20f]{20f}Darij Grinberg, \textit{Math 235: Mathematical Problem
Solving}, 22 March 2021.\newline\url{https://darijgrinberg.gitlab.io/t/20f/mps.pdf}

\bibitem[22fco]{22fco}Darij Grinberg, \textit{Math 220: Enumerative
Combinatorics, Fall 2022}.\newline\url{https://darijgrinberg.gitlab.io/t/22fco/}

\bibitem[23wa]{23wa}Darij Grinberg, \textit{An introduction to the algebra of
rings and fields}, 20 September 2026,
\href{https://arxiv.org/abs/2508.13165v2}{arXiv:2508.13165v2}.\newline\url{https://darijgrinberg.gitlab.io/t/23wa/23wa.pdf}

\bibitem[AdeGri26]{AdeGri26}\href{https://arxiv.org/abs/2611.07365}{Benjamin
Adenbaum, Darij Grinberg, \textit{The q-Poincar\'{e} sum and separability in a
Coxeter group}, arXiv:2611.07365.}\newline\url{https://darijgrinberg.gitlab.io/algebra/separable.pdf}

\bibitem[Aigner07]{Aigner07}%
\href{https://doi.org/10.1007/978-3-540-39035-0}{Martin Aigner, \textit{A
Course in Enumeration}, Graduate Texts in Mathematics \#238, Springer 2007.}

\bibitem[AndEri04]{AndEri04}%
\href{https://doi.org/10.1017/CBO9781139167239}{George E. Andrews, Kimmo
Eriksson, \textit{Integer Partitions}, Cambridge University Press 2004}.

\bibitem[Anders19]{Anders19}Kevin Anderson, \textit{Combinatorics of Symmetric
Functions}, summer research report, 28 August 2019.\newline\url{https://hadi-salmasian.com/report/report-KevinAnderson.pdf}

\bibitem[AndFen04]{AndFen04}%
\href{https://doi.org/10.1007/978-0-8176-8154-8}{Titu Andreescu, Zuming Feng,
\textit{A Path to Combinatorics for Undergraduates: Counting Strategies},
Springer 2004.}

\bibitem[Andrew16]{Andrews317}George E. Andrews, \textit{Euler's Partition
Identity -- Finite Version}, The Mathematics Student \textbf{85} (2016), pp.
99--102.\newline See
\url{https://web.archive.org/web/20230102204605/http://www.personal.psu.edu/gea1/pdf/317.pdf}
for a preprint.

\bibitem[Andrew76]{Andrew76}%
\href{https://doi.org/10.1017/CBO9780511608650}{George E. Andrews, \textit{The
Theory of Partitions}, Cambridge University Press 1976.}

\bibitem[ApaKau13]{ApaKau13}%
\href{https://doi.org/10.1016/j.exmath.2013.01.004}{Ainhoa Aparicio Monforte,
Manuel Kauers, \textit{Formal Laurent series in several variables},
Expositiones Mathematicae, \textbf{31}(4), pp. 350--367}.

\bibitem[Armstr19]{Armstrong-561fa18sp19}Drew Armstrong, \textit{Abstract
Algebra I (Fall 2018) and Abstract Algebra II (Spring 2019) lecture notes},
2019.\newline\url{https://www.math.miami.edu/~armstrong/561fa18.php}\newline\url{https://www.math.miami.edu/~armstrong/562sp19.php}

\bibitem[Artin10]{Artin}Michael Artin, \textit{Algebra}, 2nd edition, Pearson
2010.\newline See \url{https://math.mit.edu/~poonen/701/comments_algebra.pdf}
for unofficial errata.

\bibitem[BaBoSi09]{BaBoSi09}\href{https://arxiv.org/abs/0910.0963v1}{Marilena
Barnabei, Flavio Bonetti, Matteo Silimbani, \textit{The descent statistic over
123-avoiding permutations}, arXiv:0910.0963v1.}

\bibitem[Bell06]{Bell06}\href{https://arxiv.org/abs/math/0510054v2}{Jordan
Bell, \textit{Euler and the pentagonal number theorem}, arXiv:math/0510054v2}.

\bibitem[Benede25]{Benede25}Bruno Benedetti, \textit{Introduction to Abstract
Algebra \textquotedblleft Rings First\textquotedblright}, lecture notes, 6
March 2025.\newline\url{https://web.archive.org/web/20250531154025/https://www.math.miami.edu/~bruno/algebra2.pdf}

\bibitem[BenQui03]{BenQui03}\href{https://bookstore.ams.org/dol-27}{Arthur T.
Benjamin and Jennifer J. Quinn, \textit{Proofs that Really Count: The Art of
Combinatorial Proof}, The Mathematical Association of America, 2003}.

\bibitem[BenQui04]{BenQui-fib}%
\href{https://math.hmc.edu/benjamin/wp-content/uploads/sites/5/2019/06/Proofs-that-Really-Count.pdf}{Arthur
T. Benjamin and Jennifer J. Quinn, \textit{Proofs that Really Count: The Magic
of Fibonacci Numbers and More}, Mathematical Adventures for Students and
Amateurs, (David F. Hayes and Tatiana Shubin, editors), Spectrum Series of
MAA, pp. 83--98, 2004.}

\bibitem[BenQui08]{BenQui08}%
\href{https://web.archive.org/web/20170809070707/https://www.math.hmc.edu/~benjamin/papers/DIE.pdf}{Arthur
T. Benjamin and Jennifer J. Quinn, \textit{An Alternate Approach to
Alternating Sums: A Method to DIE for}, The College Mathematics Journal,
Volume 39, Number 3, May 2008, pp. 191-202(12).}

\bibitem[Berndt06]{Berndt06}\href{https://bookstore.ams.org/stml-34/}{Bruce C.
Berndt, \textit{Number Theory in the Spirit of Ramanujan}, Student
Mathematical Library \#34, AMS 2006}.\newline See
\url{https://berndt.web.illinois.edu/spiritcorrections.pdf} for errata.

\bibitem[Berndt17]{Berndt17}Bruce C. Berndt, \textit{Spring 2017, MATH 595.
Theory of Partitions}, lecture notes, 2017.\newline\url{https://berndt.web.illinois.edu/math595-tp.html}

\bibitem[Bharga00]{Bharga00}%
\href{https://doi.org/10.1080/00029890.2000.12005273}{Manjul Bhargava,
\textit{The Factorial Function and Generalizations}, The American Mathematical
Monthly \textbf{107}, No. 9 (Nov., 2000), pp. 783--799}.

\bibitem[BjoBre05]{BjoBre05}%
\href{https://doi.org/10.1007/3-540-27596-7}{Anders Bjorner, Francesco Brenti,
\textit{Combinatorics of Coxeter Groups}, Springer 2005}.\newline See
\url{https://web.archive.org/web/20250402180202/https://www.mat.uniroma2.it/~brenti/correct.ps}
for errata.

\bibitem[BluCos16]{BluCos16}\href{https://arxiv.org/abs/1301.7116v5}{Ben
Blum-Smith, Samuel Coskey, \textit{The Fundamental Theorem on Symmetric
Polynomials: History's First Whiff of Galois Theory}, arXiv:1301.7116v5.}

\bibitem[Bona22]{Bona22}Miklos B\'{o}na, \textit{Combinatorics of
Permutations}, 3rd edition, Taylor\&Francis 2022.\newline\url{https://doi.org/10.1201/9780429274107}

\bibitem[Bona25]{Bona25}Miklos B\'{o}na, \textit{Introduction to Enumerative
and Analytic Combinatorics}, 3rd edition, CRC Press 2025.

\bibitem[Bourba02]{Bourba02}%
\href{https://annas-archive.pk/md5/6bee19481829ab1bce2fa2ef34eb30b2}{Nicolas
Bourbaki, \textit{Lie Groups and Lie Algebras: Chapters 4--6}, Springer 2002}.

\bibitem[Bourba03]{Bourba03}%
\href{https://doi.org/10.1007/978-3-642-61698-3}{Nicolas Bourbaki,
\textit{Algebra II: Chapters 4--7}, Springer 2003}.

\bibitem[Bourba68]{Bourba68}%
\href{https://doi.org/10.1007/978-3-642-59309-3}{Nicolas Bourbaki,
\textit{Theory of Sets}, Springer 1968}.

\bibitem[Bourba74]{Bourba74}%
\href{https://annas-archive.pk/md5/3270565f6d0052635a1550883588204c?&check=1}{Nicolas
Bourbaki, \textit{Algebra I: Chapters 1--3}, Addison-Wesley 1974}.

\bibitem[Brande14]{Brande14}\href{https://arxiv.org/abs/1410.6601v1}{Petter
Br\"{a}nd\'{e}n, \textit{Unimodality, log-concavity, real-rootedness and
beyond}, arXiv:1410.6601v1}.

\bibitem[Bresso99]{Bresso99}%
\href{https://doi.org/10.1017/CBO9780511613449}{David M. Bressoud,
\textit{Proofs and Confirmations: The Story of the Alternating Sign Matrix
Conjecture}, Cambridge University Press 1999}.\newline See
\url{http://web.archive.org/web/20220531023303/https://www.macalester.edu/~bressoud/books/PnC/PnCcorrect.html}
for errata.

\bibitem[Brewer14]{Brewer14}Thomas S. Brewer, \textit{Algebraic properties of
formal power series composition}, PhD thesis at University of Kentucky,
2014.\newline\url{https://uknowledge.uky.edu/cgi/viewcontent.cgi?article=1021&context=math_etds}

\bibitem[BruRys91]{BruRys91}%
\href{https://doi.org/10.1017/CBO9781107325708}{Richard A. Brualdi and Herbert
J. Ryser, \textit{Combinatorial Matrix Theory}, Cambridge University Press
1991.}

\bibitem[BruSch83]{BruSch83}%
\href{https://doi.org/10.1016/0024-3795(83)80049-4}{Richard A. Brualdi and
Hans Schneider, \textit{Determinantal Identities: Gauss, Schur, Cauchy,
Sylvester, Kronecker, Jacobi, Binet, Laplace, Muir, and Cayley}, Linear
Algebra and its Applications \textbf{52--53}, July 1983, pp. 769--791}.

\bibitem[Camero16]{Camero16}Peter J. Cameron, \textit{Combinatorics 1: The art
of counting (vol. 1 of St Andrews Notes on Advanced Combinatorics)}, 28 March
2016, \url{https://cameroncounts.wordpress.com/lecture-notes/} .\newline See
\url{https://darijgrinberg.gitlab.io/algebra/acnotes1-errata.pdf} for corrections.

\bibitem[Cohen08]{Cohen08}Arjeh M. Cohen, \textit{Coxeter groups: Notes of a
MasterMath course, Fall 2007}, January 24, 2008.\newline\url{http://web.archive.org/web/20210725090759/http://arpeg.nl/wp-content/uploads/2016/01/CoxNotes.pdf}

\bibitem[Cohn04]{Cohn04}\href{https://arxiv.org/abs/math/0407093v1}{Henry
Cohn, \textit{Projective geometry over }$\mathbb{F}_{1}$\textit{ and the
Gaussian binomial coefficients}, American Mathematical Monthly \textbf{111}
(2004), pp. 487--495, arXiv:math/0407093v1}.

\bibitem[Comtet74]{Comtet74}%
\href{https://doi.org/10.1007/978-94-010-2196-8}{Louis Comtet,
\textit{Advanced Combinatorics: The Art of Finite and Infinite Expansions}, D.
Reidel Publishing Company, 1974.}

\bibitem[Conrad-UI]{Conrad-ui}Keith Conrad, \textit{Universal identities}, 13
February 2021.\newline\url{https://kconrad.math.uconn.edu/blurbs/linmultialg/univid.pdf}

\bibitem[Dodgso67]{Dodgso67}%
\href{https://archive.org/details/anelementarytre03carrgoog}{Charles L.
Dodgson, \textit{Elementary Treatise on Determinants with their Applications
to simultaneous linear equations and algebraical geometry}, Macmillan 1867}.

\bibitem[Doyle19]{Doyle19}\href{https://arxiv.org/abs/1904.06573v1}{Peter G.
Doyle, \textit{Frobenius's last proof}, arXiv:1904.06573v1}.\newline See
\url{https://darijgrinberg.gitlab.io/algebra/doyle-frob-errata.pdf} for corrections.

\bibitem[Dumas08]{Dumas08}Fran\c{c}ois Dumas, \textit{An introduction to
noncommutative polynomial invariants}, lecture notes (Cimpa-Unesco-Argentina
\textquotedblleft Homological methods and representations of non-commutative
algebras\textquotedblright, Mar del Plata, Argentina March 6 - 17,
2006).\newline\url{https://lmbp.uca.fr/~fdumas/fichiers/CIMPA.pdf}

\bibitem[DumFoo04]{DumFoo04}David S. Dummit, Richard M. Foote,
\textit{Abstract Algebra}, 3rd edition, Wiley 2004. ISBN:
978-0-471-43334-7.\newline See
\url{https://site.uvm.edu/ddummit/files/2024/01/errata_3rd_edition_6_12_2022-1.pdf}
for errata.

\bibitem[EdeStr04]{EdelStrang}%
\href{https://web.archive.org/web/20200910121631/https://www.maa.org/sites/default/files/pdf/upload_library/22/Ford/Edelman189-197.pdf}{Alan
Edelman and Gilbert Strang, \textit{Pascal Matrices}, American Mathematical
Monthly, Vol. 111, No. 3 (March 2004), pp. 189--197.}

\bibitem[Edward22]{Edward22}%
\href{https://doi.org/10.1007/978-3-030-98558-5}{Harold M. Edwards,
\textit{Essays in Constructive Mathematics}, 2nd edition, Springer 2022}.

\bibitem[Egge19]{Egge19}\href{https://bookstore.ams.org/stml-91/}{Eric S.
Egge, \textit{An Introduction to Symmetric Functions and Their Combinatorics},
AMS 2019}.\newline See \url{https://www.ericegge.net/cofsf/index.html} for
corrections and addenda.

\bibitem[EGHetc11]{EGHetc11}%
\href{https://math.mit.edu/~etingof/reprbook.pdf}{Pavel Etingof, Oleg Golberg,
Sebastian Hensel, Tiankai Liu, Alex Schwendner, Dmitry Vaintrob, Elena
Yudovina, \textit{Introduction to Representation Theory}, with historical
interludes by Slava Gerovitch, Student Mathematical Library \textbf{59}, AMS
2011, updated version 2018}.

\bibitem[Erdos42]{Erdos42}%
\href{https://users.renyi.hu/~p_erdos/1942-02.pdf}{P. Erd\"{o}s, \textit{On an
elementary proof of some asymptotic formulas in the theory of partitions},
Annals of Mathematics \textbf{43} (1942), pp. 437--450}.

\bibitem[Euler48]{Euler48}%
\href{https://scholarlycommons.pacific.edu/euler-works/101/}{Leonhard Euler,
\textit{Introductio in analysin infinitorum, tomus 1}, Lausann\ae \ 1748}.

\bibitem[Fink17]{Fink17}Alex Fink, \textit{Enumerative Combinatorics}, module
taught at the London Taught Course Centre, 2017.\newline\url{https://webspace.maths.qmul.ac.uk/a.fink/enumcombi/}

\bibitem[FlaSed09]{FlaSed09}Philippe Flajolet, Robert Sedgewick,
\textit{Analytic Combinatorics}, Cambridge University Press 2009.\newline%
\url{https://algo.inria.fr/flajolet/Publications/book.pdf}\newline See
\url{https://ac.cs.princeton.edu/errata/} for errata.

\bibitem[FoaHan04]{FoaHan04}Dominique Foata, Guo-Niu Han, \textit{The q-series
in combinatorics; permutation statistics}, preliminary version, 5 May
2011.\newline\url{https://irma.math.unistra.fr/~guoniu/papers/p56lectnotes2.pdf}

\bibitem[Ford21]{Ford21}Timothy J. Ford, \textit{Abstract Algebra}, draft of a
book, 14 August 2024.\newline\url{https://tim4datfau.github.io/Timothy-Ford-at-FAU/preprints/Algebra_Book.pdf}

\bibitem[Fulton97]{Fulton97}%
\href{https://doi.org/10.1017/CBO9780511626241}{William Fulton, \textit{Young
Tableaux, With Applications to Representation Theory and Geometry, }London
Mathematical Society Student Texts \textbf{35}, Cambridge University Press
1997}.\newline See \url{https://mathoverflow.net/questions/456463} for errata.

\bibitem[GaiGup77]{GaiGup77}\href{https://doi.org/10.1137/0132025}{P. Gaiha,
S. K. Gupta, \textit{Adjacent Vertices on a Permutohedron}, SIAM Journal on
Applied Mathematics \textbf{32} (1977), issue 2, pp. 323--327}.

\bibitem[Galvin17]{Galvin}David Galvin, \textit{Basic discrete mathematics},
13 December 2017.\newline%
\url{https://www3.nd.edu/~dgalvin1/60610/60610_S21/index.html} \newline(Follow
the Overleaf link and compile main.tex and Course-notes.tex. See also
\url{https://web.archive.org/web/20180205122609/http://www-users.math.umn.edu/~dgrinber/comb/60610lectures2017-Galvin.pdf}
for an archived old version.)

\bibitem[Gashar98]{Gashar98}%
\href{https://doi.org/10.1006/eujc.1998.0212}{Vesselin Gasharov, \textit{A
Short Proof of the Littlewood--Richardson Rule}, Europ. J. Combinatorics
(1998) \textbf{19}, pp. 451--453}.

\bibitem[Gauss08]{Gauss08}\href{https://eudml.org/doc/203313}{Carl Friedrich
Gau\ss , \textit{Summatio quarumdam serierum singularium}, Comm. soc. reg.
sci. Gottingensis rec. \textbf{1} (1811)}.

\bibitem[Gauss16]{Gauss16}C. F. Gauss, \textit{Demonstratio nova altera
theorematis omnem functionem algebraicam rationalem integram unius variabilis
in factores reales primi vel secundi gradus resolvi posse}, Comm. Recentiores
\textbf{3} (1816), pp. 107--142.\newline English translation at
\url{https://en.wikisource.org/wiki/Translation:Demonstratio_nova_altera_theorematis_omnem_functionem_algebraicam_rationalem_integram_unius_variabilis_in_factores_reales_primi_vel_secundi_gradus_resolvi_posse}
.

\bibitem[GesVie85]{GesVie85}%
\href{https://doi.org/10.1016/0001-8708(85)90121-5}{Ira Gessel, G\'{e}rard
Viennot, \textit{Binomial Determinants, Paths, and Hook Length Formulae},
Advances in Mathematics \textbf{58} (1985), pp. 300-321}.

\bibitem[GesVie89]{GesVie89}Ira M. Gessel, X. G. Viennot,
\textit{Determinants, Paths, and Plane Partitions}, 1989 preprint.\newline\url{https://peeps.unet.brandeis.edu/~gessel/homepage/papers/pp.pdf}

\bibitem[Ghys17]{Ghys17}\href{https://arxiv.org/abs/1612.06373v4}{\'{E}tienne
Ghys, \textit{A singular mathematical promenade}, ENS Editions, 2017,
arXiv:1612.06373v4.}

\bibitem[Godsil06]{Godsil06}Chris Godsil, \textit{Lecture Notes on
Combinatorics}, version 5 December 2006.\newline\url{https://web.archive.org/web/20070824060559/http://www.math.uwaterloo.ca/~dgwagner/MATH249/enum.pdf}

\bibitem[Goodma15]{Goodman}Frederick M. Goodman, \textit{Algebra: Abstract and
Concrete}, edition 2.6, 1 May 2015.\newline\url{http://homepage.math.uiowa.edu/~goodman/algebrabook.dir/book.2.6.pdf}

\bibitem[GouJac83]{GouJac83}I. P. Goulden, D. M. Jackson,
\textit{Combinatorial Enumeration}, John Wiley \& Sons 1983, reprinted by
Dover 2004.

\bibitem[Grinbe09]{Gri-19.9}%
\href{https://darijgrinberg.gitlab.io/19-9ML.pdf}{Darij Grinberg,
\textit{Solution to Problem 19.9 from \textquotedblleft Problems from the
Book\textquotedblright}}.\newline\url{https://darijgrinberg.gitlab.io/solutions.html}

\bibitem[Grinbe10]{GriHyp}Darij Grinberg, \textit{A hyperfactorial
divisibility}, version of 27 July 2015.\newline\url{https://darijgrinberg.gitlab.io/other.html}

\bibitem[Grinbe15]{detnotes}Darij Grinberg, \textit{Notes on the combinatorial
fundamentals of algebra}, 15 September 2022.\newline%
\url{https://darijgrinberg.gitlab.io/primes2015/sols.pdf} \newline The
numbering of theorems and formulas in this link might shift when the project
gets updated; for a \textquotedblleft frozen\textquotedblright\ version whose
numbering is guaranteed to match that in the citations above, see
\url{https://github.com/darijgr/detnotes/releases/tag/2022-09-15c} or
\href{https://arxiv.org/abs/2008.09862v3}{arXiv:2008.09862v3}.

\bibitem[Grinbe17]{logexp}Darij Grinberg, \textit{Why the log and exp series
are mutually inverse}, 11 May 2018.\newline\url{https://darijgrinberg.gitlab.io/t/17f/logexp.pdf}

\bibitem[Grinbe18]{diamond}Darij Grinberg, \textit{The diamond lemma and its
applications (talk)}, 20 May 2018.\newline\url{https://darijgrinberg.gitlab.io/algebra/diamond-talk.pdf}

\bibitem[Grinbe19]{trach}Darij Grinberg, \textit{The trace Cayley-Hamilton
theorem}, 23 October 2025,
\href{https://arxiv.org/abs/2510.20689v1}{arXiv:2510.20689v1}.\newline\url{https://darijgrinberg.gitlab.io/algebra/trach.pdf}

\bibitem[Grinbe20]{Grinbe20}Darij Grinberg, \textit{Alternierende Summen:
Aufgaben und L\"{o}sungen}, 28 June 2022.\newline\url{https://darijgrinberg.gitlab.io/algebra/aimo2020-altsum-lsg.pdf}

\bibitem[Grinbe21]{regpol}Darij Grinberg, \textit{Regular elements of a ring,
monic polynomials and \textquotedblleft lcm-coprimality\textquotedblright}, 22
May 2021.\newline\url{https://darijgrinberg.gitlab.io/algebra/regpol.pdf}

\bibitem[Grinbe23]{23s}\href{https://arxiv.org/abs/2308.04512v2}{Darij
Grinberg, \textit{An introduction to graph theory}, arXiv:2308.04512v2.}

\bibitem[GriRei20]{GriRei}Darij Grinberg, Victor Reiner, \textit{Hopf algebras
in Combinatorics}, version of 27 July 2020,
\href{http://www.arxiv.org/abs/1409.8356v7}{\texttt{arXiv:1409.8356v7}}.
\newline See also
\url{https://darijgrinberg.gitlab.io/algebra/HopfComb-sols.pdf} for a version
that gets updated.

\bibitem[GrKnPa94]{GKP}Ronald L. Graham, Donald E. Knuth, Oren Patashnik,
\textit{Concrete Mathematics, Second Edition}, Addison-Wesley 1994.\newline
See \url{https://www-cs-faculty.stanford.edu/~knuth/gkp.html} for errata.

\bibitem[Gruson18]{Gruson18}%
\href{https://archive.org/details/abhandlungenderp14akad/page/36/}{Gruson,
\textit{Neuer analytischer Lehrsatz}, Abhandlungen der mathematischen Klasse
der K\"{o}niglich-Preussischen Akademie der Wissenschaften 1814--1815 (Berlin
1818), pp. 36--41.}

\bibitem[Guicha20]{Guicha20}David Guichard, \textit{An Introduction to
Combinatorics and Graph Theory}, 23 April 2021.\newline\url{https://www.whitman.edu/mathematics/cgt_online/book/}

\bibitem[Hazewi08]{HW-Witt}\href{https://arxiv.org/abs/0804.3888v1}{Michiel
Hazewinkel, \textit{Witt vectors. Part 1}, arXiv:0804.3888v1.}\newline See
\url{https://darijgrinberg.gitlab.io/algebra/typos1web.pdf} for errata.

\bibitem[Hellel08]{Hellel08}Geir T. Helleloid, \textit{Algebraic
Combinatorics}, 11 November 2008.\newline\url{https://annas-archive.pk/md5/421e2dabbcd43e900bc280af5a122fe6?&check=1}

\bibitem[Henric74]{Henric74}%
\href{https://annas-archive.pk/md5/007B2FC5BFD81BC1F14ACF68FE70FEAD?check=1}{Peter
Henrici, \textit{Applied and Computational Complex Analysis, volume 1}, Wiley
1974}.

\bibitem[Henric77]{Henric77}Peter Henrici, \textit{Applied and Computational
Complex Analysis, volume 2}, Wiley 1977.

\bibitem[Hirsch17]{Hirsch17}%
\href{https://doi.org/10.1007/978-3-319-57762-3}{Michael D. Hirschhorn,
\textit{The Power of q: A Personal Journey}, Springer 2017}.\newline See
\url{https://link.springer.com/chapter/10.1007/978-3-319-57762-3_44} for errata.

\bibitem[Hirsch87]{Hirsch87}%
\href{https://web.maths.unsw.edu.au/~mikeh/webpapers/paper25.pdf}{Michael D.
Hirschhorn, \textit{A simple proof of Jacobi's four-square theorem},
Proceedings of the American Mathematical Society \textbf{101} (1987), pp.
436--438}.

\bibitem[Hopkin17]{Hopkins-RSK}Sam Hopkins, \textit{RSK via local
transformations}, \href{https://arxiv.org/abs/2510.22082}{arXiv:2510.22082}%
.\newline\url{https://www.samuelfhopkins.com/docs/rsk.pdf}

\bibitem[Hulpke25]{Hulpke25}Alexander Hulpke, \textit{Combinatorics}, 10
December 2025.\newline\url{https://www.math.colostate.edu/~hulpke/lectures/m501/notes.pdf}

\bibitem[Johnso20]{Johnso20}\href{https://bookstore.ams.org/mbk-134}{Warren
Pierstorff Johnson, \textit{Introduction to q-analysis}, American Mathematical
Society 2020}.

\bibitem[Joyner08]{Joyner08}W. D. Joyner, \textit{Mathematics of the Rubik's
cube}, 19 August 2008.\newline%
\url{https://web.archive.org/web/20160304122348/http://www.permutationpuzzles.org/rubik/webnotes/}
(link to the PDF at the bottom).

\bibitem[KacChe02]{KacChe02}%
\href{https://doi.org/10.1007/978-1-4613-0071-7}{Victor Kac, Pokman Cheung,
\textit{Quantum Calculus}, Springer 2002}.

\bibitem[Kitaev11]{Kitaev11}%
\href{https://doi.org/10.1007/978-3-642-17333-2}{Sergey Kitaev,
\textit{Patterns in Permutations and Words}, Springer 2011}.

\bibitem[KlaPol79]{KlaPol79}%
\href{https://doi.org/10.1016/0012-365x(80)90098-9}{David Klarner, Jordan
Pollack, \textit{Domino tilings of rectangles with fixed width}, Discrete
Mathematics, \textbf{32}(1), pp. 45--52}.

\bibitem[KliSch97]{KliSch97}%
\href{https://doi.org/10.1007/978-3-642-60896-4}{Anatoli Klimyk, Konrad
Schm\"{u}dgen, \textit{Quantum groups and their representations}, Springer
1997}.\newline See
\url{https://www.math.uni-leipzig.de/~schmuedgen/Artikel-Schm-Klimyk-Korr.pdf}
for errata.

\bibitem[Knapp16]{Knapp16}Anthony W. Knapp, \textit{Basic Algebra}, Digital
Second Editions By Anthony W. Knapp, 2017,
\url{http://www.math.stonybrook.edu/~aknapp/download.html} .

\bibitem[Knuth1]{Knuth-TAoCP1}Donald Ervin Knuth, \textit{The Art of Computer
Programming, volume 1: Fundamental Algorithms}, 3rd edition, Addison--Wesley
1997.\newline See \url{https://www-cs-faculty.stanford.edu/~knuth/taocp.html}
for errata.

\bibitem[Knuth2]{Knuth-TAoCP2}Donald Ervin Knuth, \textit{The Art of Computer
Programming, volume 2: Seminumerical Algorithms}, 3rd edition, Addison--Wesley
1998.\newline See \url{https://www-cs-faculty.stanford.edu/~knuth/taocp.html}
for errata.

\bibitem[Knuth3]{Knuth-TAoCP3}Donald Ervin Knuth, \textit{The Art of Computer
Programming, volume 3: Sorting and Searching}, 2nd edition, Addison--Wesley
1998.\newline See \url{https://www-cs-faculty.stanford.edu/~knuth/taocp.html}
for errata.

\bibitem[Koch16]{Koch16}Dick Koch, \textit{The Pentagonal Number Theorem and
All That}, 26 August 2016.\newline\url{https://darkwing.uoregon.edu/~koch/PentagonalNumbers.pdf}

\bibitem[KraPro10]{KraPro10}Hanspeter Kraft, Claudio Procesi,
\textit{Classical invariant theory: A primer}, July 1996.\newline%
\url{https://kraftadmin.wixsite.com/hpkraft}\newline See
\url{https://darijgrinberg.gitlab.io/algebra/KP-errata-web.pdf} for unofficial errata.

\bibitem[Kratte17]{Kratte17}%
\href{https://arxiv.org/abs/1503.05930v3}{Christian Krattenthaler,
\textit{Lattice Path Enumeration}, arXiv:1503.05930v3, published in: Handbook
of Enumerative Combinatorics, M. B\'{o}na (ed.), Discrete Math. and Its Appl.,
CRC Press, Boca Raton-London-New York, 2015, pp. 589--678.}

\bibitem[Kratte99]{Krattenthaler}%
\href{https://www.mat.univie.ac.at/~kratt/artikel/detsurv.html}{Christian
Krattenthaler, \textit{Advanced Determinant Calculus}, S\'{e}minaire
Lotharingien Combin. 42 (1999) (The Andrews Festschrift), paper B42q, 67 pp.},
\href{http://arxiv.org/abs/math/9902004v3}{arXiv:math/9902004v3}.

\bibitem[Krishn86]{Krishn86}V. Krishnamurthy, \textit{Combinatorics: Theory
and Applications}, Ellis Horwood Ltd. 1986.

\bibitem[Krob95]{Krob95}Daniel Krob, \textit{El\'{e}ments de combinatoire},
version 1.0, 1995.\newline\url{https://web.archive.org/web/20220119141152/http://krob.cesames.net/IMG/ps/combi.ps}

\bibitem[Lando03]{Lando03}\href{https://www.ams.org/books/stml/023/}{Sergei K.
Lando, \textit{Lectures on Generating Functions}, Student Mathematical Library
\textbf{23}, AMS 2003}.

\bibitem[Laue15]{Laue-det}Hartmut Laue, \textit{Determinants}, version 17 May
2015,\newline\url{https://web.archive.org/web/20250401010219/https://www.math.uni-kiel.de/algebra/laue/homepagetexte/det.pdf}

\bibitem[LLPT95]{LLPT95}D. Laksov, A. Lascoux, P. Pragacz, and A. Thorup,
\textit{The LLPT Notes}, edited by A. Thorup, 28 March 2018,\newline\url{https://web.archive.org/web/20260220073406/https://noter.math.ku.dk/sympol.pdf}

\bibitem[Loehr11]{Loehr-BC}\href{https://doi.org/10.1201/b17376}{Nicholas A.
Loehr, \textit{Bijective Combinatorics}, Chapman \& Hall/CRC 2011.}\newline
See \url{https://personal.math.vt.edu/nloehr/bijbook.html} for errata.

\bibitem[Macdon95]{Macdon95}%
\href{https://math.berkeley.edu/~corteel/MATH249/macdonald.pdf}{Ian G.
Macdonald, \textit{Symmetric Functions and Hall Polynomials}, Oxford
Mathematical Monographs, 2nd edition, Oxford Science Publications 1995.}

\bibitem[Manive01]{Manive01}Laurent Manivel, \textit{Symmetric Functions,
Schubert Polynomials and Degeneracy Loci}, SMF/AMS Texts and Monographs
\textbf{6}, AMS 2001.

\bibitem[Martin13]{Martin-Dyck}Jeremy L. Martin, \textit{Counting Dyck Paths},
11 September 2013.\newline\url{https://web.archive.org/web/20250726052421/https://jlmartin.ku.edu/~jlmartin/courses/math724-F13/count-dyck.pdf}

\bibitem[Martin21]{Martin21}Jeremy L. Martin, \textit{Lecture Notes on
Algebraic Combinatorics}, 19 November 2026.\newline%
\url{https://jeremymartinmath.github.io/CombinatoricsNotes.pdf}


\bibitem[Melcze24]{Melcze24}Stephen Melczer, \textit{An Invitation to
Enumeration}, lecture notes 2024.\newline\url{https://enumeration.ca/}

\bibitem[MenRem15]{MenRem15}%
\href{https://doi.org/10.1007/978-3-319-23618-6}{Anthony Mendes, Jeffrey
Remmel, \textit{Counting with Symmetric Functions}, Springer 2015.}

\bibitem[MetRot83]{MetRot83}%
\href{https://doi.org/10.1016/0001-8708(83)90035-X}{N. Metropolis, Gian-Carlo
Rota, \textit{Witt Vectors and the Algebra of Necklaces}, Advances in
Mathematics \textbf{50} (1983), pp. 95--125.}

\bibitem[MiRiRu87]{MiRiRu87}%
\href{https://doi.org/10.1007/978-1-4419-8640-5}{Ray Mines, Fred Richman, Wim
Ruitenburg, \textit{A Course in Constructive Algebra}, Springer 1988}.

\bibitem[Muir30]{Muir}Thomas Muir, \textit{The theory of determinants in the
historical order of development}, 5 volumes (1906--1930), later reprinted by
Dover.\newline\url{https://web.archive.org/web/20180328002253/http://www-igm.univ-mlv.fr/~al/}

\bibitem[MuiMet60]{MuiMet60}%
\href{https://annas-archive.pk/md5/6ec2fc61344d589535172f73e2c5a553?&check=1}{Thomas
Muir, \textit{A Treatise on the Theory of Determinants}, revised and enlarged
by William H. Metzler, Dover 1960}.

\bibitem[Mulhol21]{Mulhol21}Jamie Mulholland, \textit{Permutation Puzzles: A
Mathematical Perspective}, 12 January 2021.\newline%
\url{https://www.sfu.ca/~jtmulhol/permutationpuzzles/notes/permutation-puzzles-book.pdf}


\bibitem[Ness61]{Ness61}%
\href{https://annas-archive.pk/md5/6500a6e9b13e18852fc3c6f5fbd422af}{Wilhelm
Ness, \textit{Proben aus der elementaren additiven Zahlentheorie}, Otto Salle
Verlag, Frankfurt am Main / Hamburg 1961}.

\bibitem[Newste19]{Newste19}Clive Newstead, \textit{An Infinite Descent into
Pure Mathematics}, version 0.4, 1 January 2020.\newline\url{https://infinitedescent.xyz}

\bibitem[Niven69]{Niven69}Ivan Niven, \textit{Formal Power Series}, The
American Mathematical Monthly \textbf{76}, No. 8 (Oct., 1969), pp.
871--889.\newline\url{https://web.archive.org/web/20211114061642/https://www.maa.org/sites/default/files/pdf/upload_library/22/Ford/IvanNiven.pdf}

\bibitem[OlvSha18]{OlvSha}Peter J. Olver, Chehrzad Shakiban, \textit{Applied
Linear Algebra}, 2nd edition, Springer 2018.\newline%
\url{https://doi.org/10.1007/978-3-319-91041-3} \newline See
\url{https://www-users.cse.umn.edu/~olver/ala.html} for errata.

\bibitem[Pak06]{Pak06}\href{https://doi.org/10.1007/s11139-006-9576-1}{Igor
Pak, \textit{Partition bijections, a survey}, Ramanujan J \textbf{12} (2006),
pp. 5--75}.\newline See
\url{https://www.math.ucla.edu/~pak/papers/research.htm} for a preprint and updates.

\bibitem[Potter50]{Potter50}\href{https://doi.org/10.2307/2306202}{H. S. A.
Potter, \textit{On the latent roots of quasi-commutative matrices}, Amer.
Math. Monthly \textbf{57} (1950), pp. 321--322.}

\bibitem[Prasad15]{Prasad-rep}%
\href{https://doi.org/10.1017/CBO9781139976824}{Amritanshu Prasad,
\textit{Representation Theory: A Combinatorial Viewpoint}, Cambridge
University Press 2015}.

\bibitem[Prasol94]{Prasolov}%
\href{http://www2.math.su.se/~mleites/books/prasolov-1994-problems.pdf}{Viktor
V. Prasolov, \textit{Problems and Theorems in Linear Algebra}, Translations of
Mathematical Monographs, vol. \#134, AMS 1994}.

\bibitem[Proces07]{Proces07}%
\href{https://doi.org/10.1007/978-0-387-28929-8}{Claudio Procesi, \textit{Lie
Groups: An Approach through Invariants and Representations}, Springer 2007}.

\bibitem[Quinla21]{Quinlan21}Rachel Quinlan, \textit{MA3343 Groups, Semester
2020-2021},\newline\url{http://www.maths.nuigalway.ie/~rquinlan/groups/}

\bibitem[Riorda79]{Riorda79}John Riordan, \textit{Combinatorial Identities},
corrected reprint, Robert E. Krieger Publishing Company 1979.

\bibitem[Robins05]{Robins05}%
\href{https://doi.org/10.1016/j.laa.2003.10.006}{Donald W. Robinson,
\textit{The classical adjoint}, Linear Algebra and its Applications
\textbf{411} (2005), pp. 254--276}.

\bibitem[Sagan01]{Sagan01}Bruce Sagan, \textit{The Symmetric Group}, Graduate
Texts in Mathematics \textbf{203}, 2nd edition 2001.\newline%
\url{https://doi.org/10.1007/978-1-4757-6804-6}\newline See
\url{https://users.math.msu.edu/users/bsagan/Books/Sym/errata.pdf} for errata.

\bibitem[Sagan19]{Sagan19}Bruce Sagan, \textit{Combinatorics: The Art of
Counting}, Graduate Studies in Mathematics \textbf{210}, 21 September
2020.\newline%
\url{https://users.math.msu.edu/users/bsagan/Books/Aoc/final.pdf}\newline See
\url{https://users.math.msu.edu/users/bsagan/Books/Aoc/errata.pdf} for errata.

\bibitem[Sam19]{Sam19}Steven V. Sam, \textit{Notes for Math 184:
Combinatorics}, 9 December 2019.\newline\url{https://mathweb.ucsd.edu/~ssam/old/19F-184/notes.pdf}

\bibitem[Sam21]{Sam21}Steven V. Sam, \textit{Notes for Math 188: Algebraic
Combinatorics}, 18 November 2024.\newline\url{https://mathweb.ucsd.edu/~ssam/old/24F-188/notes-188.pdf}

\bibitem[Sambal22]{Sambal22}\href{https://arxiv.org/abs/2205.00879v7}{Benjamin
Sambale, \textit{An invitation to formal power series}, arXiv:2205.00879v7}.

\bibitem[Savage22]{Savage22}Alistair Savage, \textit{Symmetric Functions},
lecture notes, 4 October 2026.\newline\url{https://alistairsavage.ca/symfunc/notes/Savage-SymmetricFunctions.pdf}

\bibitem[Schwar16]{Schwar16}Rich Schwartz, \textit{The Cauchy-Binet Theorem},
9 February 2016.\newline\url{https://www.math.brown.edu/reschwar/M123/cauchy.pdf}

\bibitem[Sills18]{Sills18}%
\href{https://www.routledge.com/An-Invitation-to-the-Rogers-Ramanujan-Identities/Sills/p/book/9780367657611}{Andrew
V. Sills, \textit{An Invitation to the Rogers--Ramanujan Identities}, CRC
Press 2018.}\newline See
\url{https://web.archive.org/web/20260629111711/http://home.dimacs.rutgers.edu/~asills/}
for errata.

\bibitem[Smid09]{Smid09}Vita Smid, \textit{Inclusion-Exclusion Principle:
Proof by Mathematical Induction}, 2 December 2009.\newline\url{https://web.archive.org/web/20230523111347/https://faculty.math.illinois.edu/~nirobles/files453/iep_proof.pdf}

\bibitem[Smith95]{Smith95}\href{https://doi.org/10.1201/9781439864470}{Larry
Smith, \textit{Polynomial Invariants of Finite Groups}, A K Peters 1995.}

\bibitem[SmiTut24]{SmiTut24}%
\href{https://doi.org/10.1007/978-3-031-50341-2}{Evgeny Smirnov, Anna
Tutubalina, \textit{Symmetric Functions: A Beginner's Course}, Springer 2024.}

\bibitem[Spivey19]{Spivey19}%
\href{https://doi.org/10.1201/9781351215824}{Michael Z. Spivey, \textit{The
Art of Proving Binomial Identities}, CRC Press 2019}.\newline See
\url{https://web.archive.org/web/20240506043422/https://mathcs.pugetsound.edu/~mspivey/Errata.html}
for errata.

\bibitem[Stanko94]{Stanko94}%
\href{https://doi.org/10.1016/0012-365X(94)90242-9}{Zvezdelina E. Stankova,
\textit{Forbidden subsequences}, Discrete Mathematics \textbf{132} (1994), pp.
291--316}.

\bibitem[Stanle11]{Stanley-EC1}Richard P. Stanley, \textit{Enumerative
Combinatorics, volume 1}, Second edition, Cambridge University Press
2012.\newline See \url{http://math.mit.edu/~rstan/ec/} for a draft (2011) and errata.

\bibitem[Stanle15]{Stanle15}%
\href{https://doi.org/10.1017/CBO9781139871495}{Richard P. Stanley,
\textit{Catalan Numbers}, 1st edition, Cambridge University Press
2015}.\newline See \url{http://math.mit.edu/~rstan/catalan/} for errata.

\bibitem[Stanle18]{Stanle18}%
\href{https://doi.org/10.1007/978-3-319-77173-1}{Richard P. Stanley,
\textit{Algebraic Combinatorics: Walks, Trees, Tableaux, and More}, 2nd
edition, Springer 2018}.\newline See
\url{http://www-math.mit.edu/~rstan/algcomb/index.html} for errata.

\bibitem[Stanle24]{Stanley-EC2}Richard P. Stanley, \textit{Enumerative
Combinatorics, volume 2}, Second edition, Cambridge University Press
2024.\newline See \url{http://math.mit.edu/~rstan/ec/} for errata.

\bibitem[Stanle89]{Stanle89}%
\href{https://doi.org/10.1111/j.1749-6632.1989.tb16434.x}{Richard P. Stanley,
\textit{Log-Concave and Unimodal Sequences in Algebra, Combinatorics, and
Geometry}, Annals of the New York Academy of Sciences, \textbf{576} (1 Graph
Theory), pp. 500--535}.

\bibitem[Stembr02]{Stembr02}\href{https://doi.org/10.37236/1666}{John R.
Stembridge, \textit{A Concise Proof of the Littlewood-Richardson Rule},
Electronic Journal of Combinatorics \textbf{9} (2002), \#N5}.

\bibitem[Strick13]{Strick13}Neil Strickland, \textit{MAS201 Linear Mathematics
for Applications}, lecture notes, 11 February 2020.\newline%
\url{https://github.com/NeilStrickland/linear_maths/blob/master/notes/linear_maths.pdf}
%


\bibitem[Strick20]{Strick20}Neil Strickland, \textit{MAS334 Combinatorics},
lecture notes and solutions, 6 December 2020.\newline\url{https://strickland1.org/courses/MAS334/}

\bibitem[Stucky15]{Stucky15}Eric Stucky, \textit{An Exposition of Kasteleyn's
Solution of the Dimer Model}, senior thesis at Harvey Mudd College,
2015.\newline\url{https://scholarship.claremont.edu/hmc_theses/89/}

\bibitem[Talask12]{Talask12}\href{https://arxiv.org/abs/1202.3128v1}{Kelli
Talaska, \textit{Determinants of weighted path matrices}, arXiv:1202.3128v1.}

\bibitem[Tignol16]{Tignol16}\href{https://doi.org/10.1142/9719}{Jean-Pierre
Tignol, \textit{Galois' Theory of Algebraic Equations}, 2nd edition, World
Scientific 2016}.

\bibitem[Uecker17]{Ueckerdt}Torsten Ueckerdt, \textit{Lecture Notes
Combinatorics (2017)}, 30 May 2017.\newline%
\url{https://web.archive.org/web/20230102204620/http://www.math.kit.edu/iag6/lehre/combinatorics2017s/media/script.pdf}\newline
See
\url{https://darijgrinberg.gitlab.io/algebra/ueckerdt-script2017-errata.pdf}
for an unofficial list of errata.

\bibitem[Vorobi02]{Vorobi02}%
\href{https://doi.org/10.1007/978-3-0348-8107-4}{Nicolai N. Vorobiev,
\textit{Fibonacci Numbers}, Translated from the Russian by Mircea Martin,
Springer 2002 (translation of the 6th Russian edition)}.

\bibitem[Wagner05]{Wagner05}Carl G. Wagner, \textit{Basic Combinatorics}, 14
February 2005.\newline\url{https://web.math.utk.edu/~cwagner/papers/comb.pdf}


\bibitem[Wagner08]{Wagner08}David G. Wagner, \textit{C\&O 330: Introduction to
Combinatorial Enumeration}, version 20 June 2012.\newline\url{https://melczer.ca/330/WagnerNotes.pdf}

\bibitem[Wagner17]{Wagner17}Stephan Wagner, \textit{Combinatorics}, 19 June
2017.\newline\url{https://web.archive.org/web/20221222034640/https://math.sun.ac.za/swagner/NotesComb.pdf}

\bibitem[Wagner20]{Wagner20}\href{https://bookstore.ams.org/amstext-49/}{Carl
G. Wagner, \textit{A First Course in Enumerative Combinatorics}, Pure and
Applied Undergraduate Texts \textbf{49}, AMS 2020}.

\bibitem[Warner90]{Warner90}Seth Warner, \textit{Modern Algebra: two volumes
bound as one}, Dover 1990.

\bibitem[White10]{White10}Dennis White, \textit{Math 4707: Inclusion-Exclusion
and Derangements}, 18 October 2010.\newline\url{https://www-users.cse.umn.edu/~reiner/Classes/Derangements.pdf}

\bibitem[Wildon19]{Wildon19}Mark Wildon, \textit{Introduction to
Combinatorics}, 14 September 2020.\newline\url{http://www.ma.rhul.ac.uk/~uvah099/Maths/CombinatoricsWeb.pdf}

\bibitem[Wildon20]{Wildon20}Mark Wildon, \textit{An involutive introduction to
symmetric functions}, 8 May 2020.\newline%
\url{http://www.ma.rhul.ac.uk/~uvah099/Maths/Sym/SymFuncs2020.pdf}\newline See
\url{https://darijgrinberg.gitlab.io/algebra/symfuncs2017-2020-05-08-errata.pdf}
for an unofficial list of errata.

\bibitem[Wilf04]{Wilf04}Herbert S. Wilf, \textit{generatingfunctionology}, 2nd
edition 2004.\newline\url{https://www.math.upenn.edu/~wilf/DownldGF.html}

\bibitem[Wilf09]{Wilf09}Herbert S. Wilf, \textit{Lectures on Integer
Partitions}, 2009.\newline\url{https://www.math.upenn.edu/~wilf/PIMS/PIMSLectures.pdf}

\bibitem[Xantch17]{Xantch17}\href{https://arxiv.org/abs/1104.1931v4}{Qimh
Richey Xantcha, \textit{Binomial Rings: Axiomatisation, Transfer and
Classification}, arXiv:1104.1931v4.}

\bibitem[Zabroc03]{Zabroc03}Mike Zabrocki, \textit{F. Franklin's proof of
Euler's pentagonal number theorem}, 28 February 2003.\newline\url{https://garsia.math.yorku.ca/~zabrocki/math4160w03/eulerpnt.pdf}

\bibitem[Zeilbe85]{Zeilbe}Doron Zeilberger, \textit{A combinatorial approach
to matrix algebra}, Discrete Mathematics 56 (1985), pp. 61--72, corrected
transcription 2026.\newline\url{https://sites.math.rutgers.edu/~zeilberg/mamarimY/DM85dg.pdf}

\bibitem[Zeilbe98]{zeilberger-twotime}Doron Zeilberger, \textit{Dodgson's
Determinant-Evaluation Rule proved by Two-Timing Men and Women}, The
Electronic Journal of Combinatorics, vol. 4, issue 2 (1997) (The Wilf
Festschrift volume), R22.\newline\url{https://doi.org/10.37236/1337} \newline
Also available as
\href{http://arxiv.org/abs/math/9808079v1}{arXiv:math/9808079v1}.

\bibitem[Zelevi81]{Zelevi81}%
\href{https://doi.org/10.1016/0021-8693(81)90128-9}{A. V. Zelevinsky,
\textit{A generalization of the Littlewood--Richardson rule and the
Robinson--Schensted--Knuth correspondence}, J. Algebra \textbf{69} (1981), pp.
82--94}.

\bibitem[Zeng93]{Zeng93}%
\href{https://doi.org/10.1016/0024-3795(93)90371-T}{Jiang Zeng, \textit{A
bijective proof of Muir's identity and the Cauchy-Binet formula}, Linear
Algebra and its Applications \textbf{184}, 15 April 1993, pp. 79--82}.
\end{thebibliography}
\end{document}